\documentclass[12pt, reqno]{amsbook}

\usepackage{chngcntr}
\counterwithin{table}{section}

 \usepackage{caption}
 
\usepackage{pagenote}
\makepagenote

\renewcommand*{\pagenotesubhead}[2]{}

\usepackage[utf8]{inputenc}

\usepackage{latexsym} 
\usepackage{amsxtra}
\usepackage{amssymb}
\usepackage{bm}
\usepackage{amsthm}
\usepackage{mathtools}
\usepackage{xspace}
\usepackage{graphicx}
\usepackage[margin=1in]{geometry}
\usepackage{xfrac}
\usepackage{imakeidx}
\makeindex
\let\oldindex\index
\renewcommand*{\index}[1]{\oldindex{#1}\ignorespaces}
\makeindex[name=symbols,title={Index of Symbols}]

 \usepackage{caption}
\usepackage[T1]{fontenc}
\usepackage{gfsbodoni}

\usepackage[section]{placeins}

\usepackage{verbatim}

\DeclareFontFamily{U}{mathx}{\hyphenchar\font45}
\DeclareFontShape{U}{mathx}{m}{n}{
      <5> <6> <7> <8> <9> <10>
      <10.95> <12> <14.4> <17.28> <20.74> <24.88>
      mathx10
      }{}
\DeclareSymbolFont{mathx}{U}{mathx}{m}{n}
\DeclareFontSubstitution{U}{mathx}{m}{n}
\DeclareMathAccent{\widecheck}{0}{mathx}{"71}
\DeclareMathAccent{\wideparen}{0}{mathx}{"75}

\input xy
\xyoption{all} \CompileMatrices

\begin{document}

\def\Log{\operatorname{Log}}
\def\ERi{\operatorname{ERi}}
\def\Mert{M_2}
\def\Mertt{M_3}
\def\ee{\mathbf{e}}
\def\dd{\mathbf{d}}
\def\dom{\operatorname{dom}}
\def\HH{{\mathbb{H}}}
\def\orb{{\operatorname{orb}}}
\def\diam{{\operatorname{diam}}}
\def\II{{\mathfrak{I}}}
\def\PO{{\operatorname{PO}}}
\def\Cl{{\operatorname{Cl}}}
\def\Max{{\operatorname{-Max}}}
\def\XX{{\bf{X}}}
\def\YY{{\bf{Y}}}
\def\BBB{{\mathcal B}}
\def\inv{{\operatorname{inv}}}
\def\emph{\it}
\def\Int{{\operatorname{Int}}}
\def\Spec{\operatorname{Spec}}
\def\Bin{{\operatorname{B}}}
\def\n{\operatorname{b}}
\def\N{{\operatorname{GB}}}
\def\BC{{\operatorname{BC}}}
\def\dlog{\frac{d \log}{dT}}
\def\Sym{\operatorname{Sym}}
\def\Nr{\operatorname{Nr}}
\def\lbrack{{\{}}
\def\rbrack{{\}}}
\def\burnside{\operatorname{B}}
\def\Sym{\operatorname{Sym}}
\def\Hom{\operatorname{Hom}}
\def\Inj{\operatorname{Inj}}
\def\Aut{{\operatorname{Aut}}}
\def\Mor{{\operatorname{Mor}}}
\def\Map{{\operatorname{Map}}}
\def\CMap{{\operatorname{CMap}}}
\def\GMaps{G{\operatorname{-Maps}}}
\def\Fix{{\operatorname{Fix}}}
\def\res{{\operatorname{res}}}
\def\ind{{\operatorname{ind}}}
\def\inc{{\operatorname{inc}}}
\def\coind{{\operatorname{cnd}}}
\def\Equiv{{\mathcal{E}}}
\def\W{\operatorname{W}}
\def\F{\operatorname{F}}
\def\witt{\operatorname{gh}}
\def\ngh{\operatorname{ngh}}
\def\Fm{{\operatorname{Fm}}}
\def\bij{{\iota}}
\def\mk{{\operatorname{mk}}}
\def\km{{\operatorname{mk}}}
\def\VV{{\bf{V}}}
\def\ff{{\bf{f}}}
\def\ZZ{{\mathbb Z}}
\def\pp{{\mathbb P}}
\def\ss{{\mathbb S}}
\def\Zhat{{\widehat{\mathbb Z}}}
\def\CC{{\mathbb C}}
\def\PP{{\mathbf p}}
\def\DD{{\mathbb D}}
\def\EE{{\mathbb E}}
\def\MM{{\mathbb M}}
\def\JJ{{\mathbb J}}
\def\NN{{\mathbb N}}
\def\RR{{\mathbb R}}
\def\QQ{{\mathbb Q}}
\def\FF{{\mathbb F}}
\def\mm{{\mathfrak m}}
\def\nn{{\mathfrak n}}
\def\jj{{\mathfrak j}}
\def\aaa{{{{\mathfrak a}}}}
\def\bbb{{{{\mathfrak b}}}}
\def\ppp{{{{\mathfrak p}}}}
\def\qqq{{{{\mathfrak q}}}}
\def\PPP{{{{\mathfrak P}}}}
\def\BB{{\mathfrak B}}
\def\jj{{\mathfrak J}}
\def\LL{{\mathfrak L}}
\def\qq{{\mathfrak Q}}
\def\rr{{\mathfrak R}}
\def\cc{{\mathfrak S}}
\def\TT{{\mathcal{T}}}
\def\SS{{\mathcal S}}
\def\UU{{\operatorname{U}}}
\def\AA{\operatorname{Arith}}
\def\BB{{\mathcal B}}
\def\Primes{{\mathcal P}}
\def\genS{{\langle S \rangle}}
\def\genT{{\langle T \rangle}}
\def\bT{\mathsf{T}}
\def\bD{\mathsf{D}}
\def\bC{\mathsf{C}}
\def\VV{{\bf V}}
\def\ff{{\bf f}}
\def\uu{{\bf u}}
\def\aa{{\bf{a}}}
\def\bb{{\bf{b}}}
\def\zero{{\bf 0}}
\def\rad{\operatorname{rad}}
\def\End{\operatorname{End}}
\def\id{\operatorname{id}}
\def\mod{\operatorname{mod}}
\def\im{\operatorname{im}}
\def\ker{\operatorname{ker}}
\def\coker{\operatorname{coker}}
\def\ord{\operatorname{ord}}
\def\li{\operatorname{li}}
\def\Ei{\operatorname{Ei}}
\def\Ein{\operatorname{Ein}}
\def\Ri{\operatorname{Ri}}
\def\Rie{\operatorname{Rie}}
\def\degl{\operatorname{ldeg}}
\def\dege{\operatorname{ledeg}}
\def\lerank{\operatorname{lerank}}

\newtheorem{theorem}{Theorem}[section]
\newtheorem{proposition}[theorem]{Proposition}
\newtheorem{corollary}[theorem]{Corollary}
\newtheorem{lemma}[theorem]{Lemma}

\theoremstyle{definition}
\newtheorem{definition}[theorem]{Definition}
\newtheorem{excursion}[theorem]{Excursion}
\newtheorem{conjecture}[theorem]{Conjecture}
\newtheorem{problem}[theorem]{Problem}
\newtheorem{outstandingproblem}[theorem]{(Outstanding) Problem}
\newtheorem{example}[theorem]{Example}

\newtheorem{zremark}[theorem]{Remark}
\newenvironment{remark}{\par\footnotesize\zremark}{\endzremark}

\makeatletter
\newcounter{savechapter}
\newcounter{apdxchapter}
\renewcommand\appendix{\par
  \setcounter{savechapter}{\value{chapter}}%
  \setcounter{chapter}{\value{apdxchapter}}%
  \gdef\thechapter{\@Alph\c@chapter}}
\newcommand\unappendix{\par
  \setcounter{apdxchapter}{\value{chapter}}%
  \setcounter{chapter}{\value{savechapter}}%
  \gdef\thechapter{\@arabic\c@chapter}}
\makeatother

 \newenvironment{map}[1]
   {$$#1:\begin{array}{rcl}}
   {\end{array}$$
   \\[-0.5\baselineskip]
 }

 \newenvironment{map*}
   {\[\begin{array}{rcl}}
   {\end{array}\]
   \\[-0.5\baselineskip]
 }

 \newenvironment{nmap*}
   {\begin{eqnarray}\begin{array}{rcl}}
   {\end{array}\end{eqnarray}
   \\[-0.5\baselineskip]
 }

 \newenvironment{nmap}[1]
   {\begin{eqnarray}#1:\begin{array}{rcl}}
   {\end{array}\end{eqnarray}
   \\[-0.5\baselineskip]
 }

\newcommand{\eq}{eq.\@\xspace}
\newcommand{\eqs}{eqs.\@\xspace}
\newcommand{\diagram}{diag.\@\xspace}

\newtheoremstyle{mystyle}
  {}
  {}
  {\itshape}
  {}
  {\bfseries}
  {.}
  { }
  {}

\newtheoremstyle{mydefinitionstyle}
  {}
  {}
  {}
  {}
  {\bfseries}
  {.}
  { }
  {}

\theoremstyle{mystyle}

\renewcommand{\thesection}{\thechapter.\arabic{section}}

\numberwithin{figure}{section}
\numberwithin{equation}{section}

\newcounter{Lcount}
\newcounter{Pcount}

\newenvironment{exercises}{
  \begin{list}{\arabic{Lcount}.}
    {\usecounter{Lcount}\leftmargin=1em}}
  {\end{list}}

\newenvironment{exparts}{
  \begin{list}{\alph{Pcount})}
    {\usecounter{Pcount}\leftmargin=1.5em}}
  {\end{list}}

\newcommand\bigsqcapp{\mathop{\mathchoice
	{\vcenter{\hbox{\huge $ \sqcap $}}}	
	{\vcenter{\hbox{\LARGE $ \sqcap $}}}	
	{\vcenter{\hbox{$ \sqcap $}}}		
	{\vcenter{\hbox{\small $ \sqcap $}}}}}	

\newcommand\bigsqcupp{\mathop{\mathchoice
	{\vcenter{\hbox{\huge$ \sqcup $}}}	
	{\vcenter{\hbox{\LARGE$ \sqcup $}}}	
	{\vcenter{\hbox{$ \sqcup $}}}		
	{\vcenter{\hbox{\small $ \sqcup $}}}}}	

\title{\huge Analytic Number Theory and \\ Algebraic Asymptotic Analysis}

\author{ \LARGE Jesse Elliott}

\frontmatter

\maketitle

\renewcommand{\contentsname}{Table of Contents}

\tableofcontents

\chapter*{Dedication}

\begin{center} This book is dedicated to my mother, Nancy (Elliott) Cappello.
\end{center}

\chapter*{Preface}

\section{Overview}

This monograph introduces a new framework for measuring and comparing the asymptotic behavior of complex functions of a  real variable, one that proves especially effective within number theory. At its core is a formalism built around two invariants, \textbf{degree} and \textbf{logexponential degree},  that quantify asymptotic growth.   These invariants extend classical notations such as $O$, $o$, and $\sim$ that were developed by du Bois-Reymond, Bachmann, Landau, Hardy, and others, into a more expressive and numerically structured analytic-algebraic language of asymptotic comparison.

This  framework is situated within a broader field that we refer to as  \textit{algebraic asymptotic analysis}: the study of asymptotic behavior using tools from both algebra and analysis.    Algebraic asymptotic analysis builds on several  other fields, including asymptotic analysis \cite{debruijn} \cite{est} \cite{har3} \cite{har4} \cite{olver} \cite{wong}; Karamata theory \cite{bgt} \cite{kar1} \cite{kar2} \cite{seneta}; the theory of continued fractions and moments in the tradition of Stieltjes \cite{akh} \cite{cuyt} \cite{stie}; and asymptotic differential algebra \cite{asch} \cite{bos} \cite{rosenl} \cite{dmm} \cite{dmm2} \cite{edgar} \cite{rosenl1} \cite{vdh}.

The book is structured in three parts. Part~1 is a survey of classical analytic number theory in both mathematical and historical contexts.  Anticipating Parts~2 and 3, it includes foundational material on asymptotic analysis, arithmetic functions, Dirichlet series, special functions,  and the analytic theory of primes,  prioritizing conceptual breadth, clarity, and visual insight over formal proof.   Part~2 presents the central theoretical framework of the book, developing the theory of logexponential degree and its interactions with asymptotic analysis and real asymptotic differential algebra---notably, through connections to ordered exponential fields,  Hardy fields,  Hardy's ordered differential field of logarithmico-exponential functions \cite{har3} \cite{har4}, and the ordered differential field  of well-based logarithmic-exponential  transseries \cite[Appendix A]{asch} \cite{dmm} \cite{dmm2} \cite{edgar}.  Finally, Part~3 applies Parts~1 and 2 to a broad range of functions in number theory.   Specifically, it uses logexponential degree as a unifying organizing principle for analyzing asymptotic behavior and motivating conjectures related to the following: the Riemann hypothesis,  the Lindel\"of hypothesis,  and the density hypothesis; the ordinates of the Riemann zeta  function zeros and their gaps; the prime counting function, weighted prime counting functions,  the Mertens function,  and summatory functions more broadly; the prime listing function and prime gaps; the Dirichlet divisor problem; least primitive roots, quadratic non-residues, and primes in arithmetic progressions; Roth's theorem on irrationality measure; and the abc conjecture.

At the foundation of the new formalism is the notion of degree, which appears implicitly throughout analytic number theory.  Precisely, we define the {\bf degree} $\deg f \in [-\infty,\infty]$, for any complex-valued function $f$ defined on a subset of $\RR$ that is not bounded above, to be the infimum of all $t \in \RR$ such that $|f(x)| \ll x^t$ as $x \to \infty$. Equivalently, one has $\deg f = \limsup_{x \to \infty} L_f(x)$, where $L_f(x) = \frac{\log |f(x)|}{\log x}$ is the unique function from $\RR_{>0} \cap \dom f$ to $[-\infty,\infty)$ satisfying $|f(x)| = x^{L_f(x)}$, where $\log 0 = -\infty$ and $0 = x^{-\infty}$. The corresponding {\bf lower degree} $\underline{\deg} f$ is the supremum of all $t \in \RR$ such that $|f(x)| \gg x^t$ and is given equivalently by $\liminf_{x \to \infty} L_f(x)$, and also by $-\deg(1/f)$ if $f$ is not eventually $0$. Intuitively, degree can be seen as a degree of freedom, representing the maximal power a function can asymptotically achieve. Likewise, lower degree represents a degree of variation—the minimal power it can asymptotically maintain.

A few notable examples illustrate the significance of the degree invariant. The Riemann hypothesis, for instance, is equivalent to the statement that $\deg(\li - \pi) = \tfrac{1}{2}$, where $\pi(x)$ is the prime counting function and  $\li(x) = \int_0^x \frac{dt}{\log t}$ is the logarithmic integral function.  Likewise, it is equivalent to the statement that $\deg M = \tfrac{1}{2}$, where $M(x) = \sum_{n \leq x} \mu(n)$ is the Mertens function. More generally, one has $\deg(\li - \pi) = \Theta = \deg M(x)$, where $\Theta \in [\tfrac{1}{2},1]$ is the supremum of the real parts of the nontrivial zeros of the Riemann zeta function.  
The equality $\deg(\li - \pi) = \Theta$ means that the bound $\li(x) - \pi(x) = O(x^t)$ holds for all $t > \Theta$ and fails for all $t < \Theta$.     Intuitively,  this says that, to first order,   the prime counting function fluctuates around its smooth approximation $\li(x)$ by a $\Theta$-power law, and thus it is like a ``law of gravity'' for the primes.   The  constant $\Theta$, which we call the {\bf Riemann constant},  thus encodes a fundamental relationship between the prime counting function and the zeros of  Riemann zeta function.   The best case scenario is the Riemann hypothesis, $\Theta = \tfrac{1}{2}$; the worst case scenario is $\Theta = 1$, which we call the {\bf anti-Riemann hypothesis}.   Although a disproof of the anti-Riemann hypothesis would represent major progress toward a proof of the Riemann hypothesis,  currently there is no guarantee that the anti-Riemann hypothesis is false, let alone that the Riemann hypothesis is true. This reveals an asymmetry in the problem of settling the Riemann hypothesis,  namely,  that  a mere negative answer would still leave the exact value of $\Theta = \deg(\li - \pi)$ a mystery.  It also hints at what is widely known but rarely stated: many statements equivalent to the Riemann hypothesis can   be recast as unconditional statements about the Riemann constant $\Theta$.   This gives weight to the claim that the closer that $\Theta$ is to $\tfrac{1}{2}$, the ``more true'' the Riemann hypothesis is.   In this respect,  the Riemann hypothesis acts less like a binary proposition, and more like a fuzzy one whose truth value is a decreasing function $\text{RH}(\Theta)$ of $\Theta$---for example, $\text{RH}(\Theta) = 2(1-\Theta) \in [0,1]$, where $1-\Theta$ is the infimum of the real parts of the nontrivial zeros of $\zeta(s)$.   

Since the functions $\li-\pi$ and $M$ are likely not correlated to more than first order,  the equality $\deg(\li -\pi) = \deg M$  suggests the need for a finer measure of asymptotic behavior,  one capable of quantifying such behavior not just relative to power functions, but more broadly to Hardy's logarithmico-exponential functions \cite{har4}.   A {\bf logarithmico-exponential function} is a real function that is defined on a neighborhood of $\infty$ and can be built from all real constants and the functions $\id$, $\exp$, and $\log$ using the operations of addition, multiplication, division,  and composition.   On the ubiquity of such functions, Hardy remarked that ``the only scales of infinity that are of any practical importance in analysis are those which may be constructed by means of the logarithmic and exponential functions'' \cite[p.\ 22]{har3}. While this assessment might be slightly overstated, the vast corpus of number theory literature reveals that logarithmico-exponential functions indeed offer natural and precise benchmarks against which to compare a vast array of asymptotic behavior.

Logexponential degree captures such behavior by ordering functions naturally and linearly.  It is defined precisely as follows. Let $f$ be a real function whose domain is a subset of $[-\infty,\infty)$ that is not bounded above, and let $T(f)$ denote the real function
\[
T(f)(x) = \left.
\begin{cases}
    f(e^x) e^{-(\deg f) x} & \text{if } \deg f \neq \pm \infty \\
    \max( \log |f(x)|, 0)  & \text{if } \deg f = \infty \\
    -\dfrac{1}{\log |f(x)|} & \text{if } \deg f = -\infty,
\end{cases}
\right.
\]
where again we interpret $\log 0 = -\infty$. The {\bf logexponential degree of $f$}, denoted $\dege f$, is the infinite sequence
\[
\dege f = (\deg f, \deg T(f), \deg T(T(f)), \deg T(T(T(f))), \ldots) \in \prod_{n = 0}^\infty [-\infty,\infty],
\]
where the product is equipped with the lexicographic (total) ordering.   We write $\dege_k f$ for the $k$th coordinate of $\dege f$, so in particular $\dege_0 f = \deg f$.  
Note, for example, that
$$\dege(x^{1/2} (\log x)^{-1}\log \log \log x)  = (1/2,-1,0,1, 0,0,0,\ldots)$$  and $$\dege (x{e^{- 0.2098(\log x)^{3/5}(\log \log x)^{-1/5}}}) = (1,-\infty,-3/5,1/5,0,0,0,\ldots).$$ 
Both of these logarithmico-exponential functions have appeared in the study of the asymptotic behavior of the prime counting function $\pi(x)$ \cite{litt} \cite{ford}.

 Though seemingly {\it ad hoc}, the logexponential degree formalism lends itself to a natural axiomatic characterization: see Theorem \ref{thm:ledege_axioms}.     Importantly,  one has the implications
\[
\dege f < \operatorname{ledeg} g \quad \Rightarrow \quad f = o(g) \quad \Rightarrow \quad f = O(g) \quad \Rightarrow \quad \operatorname{ledeg} f \leq \dege g,
\]
where the first implication holds if $g$  is logarithmico-exponential.  What distinguishes this hierarchy is that \( o \) and \( O \) are qualitative and partially ordered, while $\dege$ is quantitative and linearly ordered---not as a single number, but a lexicographically ordered sequence encoding an infinite tower of power-law growth rates.    
Moreover, the transitions in the hierarchy are not arbitrary: each level in the \( \operatorname{ledeg} \) sequence reflects a ``slower'' scale than the prior levels,  all arbitrated by the map $T$.

The recursive definition of $\dege$ reflects extensive trial and refinement; earlier versions the author attempted failed to preserve the analytic-algebraic properties necessary for global coherence across the degree hierarchy.    Once established, the formalism opens a path to recasting a vast number of classical problems in more general terms.   For example, it uncovers two problems that transcend the Riemann hypothesis: compute the logexponential degree $\dege(\li - \pi)$ and of $\dege M$. Several well-known results in the literature imply the following  constraints on the constants $\Theta_1 = \dege_1(\li - \pi)$ and $\text{M}_1 = \dege_1 M$: unconditionally, one has $\Theta_1 \in [-\infty, 1]$; if the Riemann hypothesis is true, then $\Theta_1 \in [-1, 1]$; if the Riemann hypothesis is false, then $\Theta_1 \in [-\infty, -1]$; if the anti-Riemann hypothesis is true, then $\Theta_1 = \text{M}_1 = -\infty$ and $\dege_2(\li - \pi)= \dege_2 M \in [-1, -\tfrac{3}{5}]$; and if the Riemann hypothesis is true and $\text{M}_1 = 0$, then all of the zeros of $\zeta(s)$ are simple.   Moreover,  assuming a 1979 conjecture of Montgomery, one has $\dege (\li-\pi) = (\tfrac{1}{2},-1,0,2,0,0,0,\ldots)$ \cite[Conjecture, p.\ 16]{mont1}.  Likewise, assuming an unpublished conjecture of Gonek, later supported by Ng  \cite[(20)]{ng},  one has  $\dege M = (\tfrac{1}{2}, 0, 0, \frac{5}{4}, 0,0,0,\ldots)$.   These inferences all motivate the conjecture---backed  also by modest numerical data---that $\Theta_1 = -1$ and $\text{M}_1 = 0$.

One of the strengths of the logexponential degree formalism is its ability to express structural relationships among number-theoretic functions---relationships that are often difficult to formulate or prove using traditional methods. For example, one can relate the logexponential degree of various functions to that of \(\operatorname{li}(x) - \pi(x)\), and vice versa. In particular, we accomplish this for each of the functions
\[
p_n - \operatorname{li}^{-1}(n),\quad x - \vartheta(x),\quad x - \psi(x),\quad e^{\gamma} \prod_{p \leq x} \left(1 - \frac{1}{p} \right) - \frac{1}{\log x},
\]
\[
\sum_{p \leq x} \frac{1}{p} - \log \log x - M,\quad e^{-M} \prod_{p \leq x} e^{1/p} - \log x ,\quad \text{and}  \quad  \log x - \sum_{p \leq x} \frac{\log p}{p} - B,
\]
where \(\vartheta\) and \(\psi\) denote the first and second Chebyshev functions, respectively, and \(\gamma\), \(M\), and \(B\) are well-known constants.  

The formalism also applies to Diophantine approximation. In Chapter~13, we show that the irrationality measure \( \mu(\alpha) \) of any irrational number \( \alpha \) can be interpreted as the degree $\mu(\alpha) = \deg \frac{1}{\mu_\alpha(x)} =  - \underline{\deg} \, \mu_\alpha(x)$, where
$$\mu_\alpha(x) =  \min\left\{ \left|\alpha - \frac{a}{b} \right| : a, b \in \mathbb{Z},\ 1 \leq b \leq x \right\}.$$  It follows that Roth's theorem is equivalent to $\deg \frac{1}{\mu_\alpha(x)} = 2$ for all algebraic irrationals $\alpha$.  Naturally, we define the {\bf logexponential irrationality degree} of $\alpha$ to be the logexponential degree of $\frac{1}{\mu_\alpha(x)}$. 
 We show that, for almost all real numbers $\alpha$,  the logexponential irrationality degree of $\alpha$ is equal to $(2,1,1,1,\ldots)$.  This leads us to speculate that this law also holds for the numbers $\pi$,  $\gamma$, $\log 2$,  and any other ``natural'' real number for which the terms of its continued fraction is expected to obey the asymptotic running geometric mean of Khinchin's constant---including all algebraic numbers of degree greater than two.   

Yet another application is to the {\it abc conjecture}: one can show that it is equivalent to the statement that the  lower degree $\underline{\deg} \, \underline{\operatorname{ABC}}(n)$ of the function
\begin{align*}
\underline{\operatorname{ABC}}(n)  & =  \min\{\operatorname{rad}(ab(a + b)): a, b \in \mathbb{Z}_{>0},\ \gcd(a, b) = 1,  \ a + b = n\} \\ & \leq \rad(n-1)\rad(n)
\end{align*}
is equal to $1$, where $\operatorname{rad}(n) = \prod_{p | n} p$.   Known bounds on $\underline{\operatorname{ABC}}(n)$ put constraints and suggest conjectures on what its logexponential degree could be.

Such applications demonstrate the utility of the  logexponential degree formalism across a vast array of functions in number theory. To manage these effectively,  we organize them around a core set of primitives,  such as $\li(x) - \pi(x)$, $M(x)$, and $\frac{1}{\mu_\alpha(x)}$,  in terms of which we express the logexponential degree of many other number-theoretic functions. This strategy both clarifies existing dependencies and reveals new ones, offering a more coherent view of how various number-theoretic  functions and error terms in their natural approximations interrelate.

The theoretical contributions of this monograph reveal  connections  between the  logexponential degree formalism and real asymptotic differential algebra.   For example, we apply the formalism  to the study of Hardy fields,  and vice versa, where a {\bf Hardy field} is a subfield of the ring of germs of all real functions at $\infty$ that is closed under differentiation.   We also define logexponential degree formally as a canonical map on the ordered differential field $\mathbb{T}$ of well-based logarithmic-exponential transseries.   We then compute  not only which  sequences are equal to $\dege f$ some function $f$,  but also the image of $\dege$ on both $\mathbb{T}$ and the field $\mathbb{L}$ of logarithmico-exponential functions (and we find that the two images coincide).   As one might expect,  for example, the logexponential degree of any transseries, hence of any logarithmico-exponential function,  ends in a trail of $0$s.  

At its core, algebraic asymptotic analysis retains those aspects of real asymptotic differential algebra governing asymptotic growth and comparison that persist even in the absence of a derivation. Since differentiation still enters into some of its applications---particularly through Karamata theory---connections between logexponential degree and Hardy fields are not unexpected. However, although both frameworks provide canonical models of asymptotic behavior, logexponential degree isolates that behavior without relying on derivational structure.

Among the many branches of mathematics, number theory is perhaps unique in that many of its central problems are simple to state yet notoriously difficult to solve. It is likely that newer and more powerful methods are needed to tackle its longstanding problems.  Even the slightest improvement in an error bound often demands extraordinary ingenuity. Rather than attempting to sharpen such bounds directly, this text offers a framework for relating them—so that an improvement in one may yield corresponding improvements in others. While our results are unconditional, they also reflect and organize the best known error bounds to date, drawing on the accumulated knowledge of analytic number theory,  both classical and contemporary.

As a result, the degree framework serves not merely as a technical tool, but as a language for measuring both mathematical and epistemological uncertainty. More precisely, it provides a structure for quantifying, relating, and tracking asymptotic errors across a network of approximations, treating them as functionally and relationally dependent, often in ways previously expressed only heuristically.   Ultimately, this offers a refinement in how we describe asymptotic behavior: not merely as a set of binary outcomes, but as a structured space of epistemic and metaphysical possibilities that expresses much of what we understand and  seek to understand in precise numerical form.

\section{Prerequisites and target audience}

 The target audience for this book is anyone  interested in analytic number theory---particularly those with the equivalent of a Bachelor's degree in mathematics---as well as researchers and graduate students in analytic number theory,  asymptotic analysis, and asymptotic differential algebra.   This book connects those three fields but does not assume prior knowledge of any of them. 

The mathematical prerequisites are  elementary number theory, abstract algebra, real and complex analysis,  and topology, all at the beginning graduate or advanced undergraduate level.    Many standard results from analytic number theory are stated without proof,  with citations provided.   These results are used as  black boxes: techniques like contour integration are rarely unpacked in detail.  The goal is not to reconstruct known results, but to develop a complementary framework focused on asymptotic structure.  

Part 1 of the text is an extenstive  introductory  survey of analytic number theory written for the non-specialist.  It is included not only to make the text more accessible, but also to lay out the results our theory relies on.  Classical proofs appear in standard graduate or advanced undergraduate textbooks (e.g., \cite{bate} \cite{borg} \cite{dav} \cite{gros2} \cite{kou} \cite{mont} \cite{nark} \cite{over} \cite{plyman} \cite{tit} \cite{zud}), and proofs of more recent results are always referenced in the bibliography.   Table \ref{theoremproofs} in Section 1.3 lists several important classical theorems concerning the distribution of the primes,  with extensive references to proofs showcasing diverse methods and levels of generality.   This helps guide the reader into the deeper theory with less risk of privileging one methodological tradition over others.

\section{Outline of contents}

This text is divided into three parts.   Part 1 ({\it A survey of analytic number theory}) is a survey of analytic number theory at the  advanced undergraduate  or beginning graduate level.  The expert in analytic number theory may skip this part of the book.

\begin{enumerate}
\item[Ch.\ 1.]  {\it A brief history of primes.}  Topics include: prime and composite numbers, the fundamental theorem of arithmetic, Mersenne primes, twin primes, prime $k$-tuples conjecture,  Riemann zeta function, arithmetic functions, summatory functions, Dirichlet series, prime listing function $p_n$, prime counting function $\pi(x)$, prime number theorem,  Euler--Mascheroni constant, harmonic numbers, Meissel--Mertens constant, Riemann hypothesis, logarithmic integral function $\li(x)$, Riemann's function $\Ri(x)$, Riemann--von Mangoldt explicit formula, prime number theorem with error bound, Riemann constant $\Theta$.   
\item[Ch.\ 2.]   {\it Asymptotic analysis.}    Topics include: asymptotic relations ($O$, $o$, $\asymp$,  $\sim$, $\Omega_+$, $\Omega_-$, $\Omega_{\pm}$), asymptotic expansions, asymptotic expansions of $\pi(x)$, Euler--Maclaurin formula, Stirling's approximation, the degree $\deg f$ of a real function $f$,  slowly varying and regularly varying functions, Karamata's integral representation theorem and integral theorem.   Section  2.3 introduces the new notion of the {\it degree} of a real function,  and Section 2.4 relates it to the study of regularly varying functions.  
\item[Ch.\ 3.] {\it Arithmetic functions.}    Topics include: elementary complex functions,  formal Dirichlet series, Dirichlet convolution,  ring of arithmetic functions, the rings of multiplicative and additive arithmetic functions, M\"obius inversion theorem, Bell series,  summatory functions,  inversion theorems, Abel's summation formula,  partial summation,    Mertens function $M(x)$,  first Chebyshev function $\vartheta(x)$, second Chebyshev function $\psi(x)$,  Riemann's prime counting function $\Pi(x)$,  von Mangoldt function $\Lambda(n)$,  average value and average order, Dirichlet hyperbola method,  Dirichlet series,  abscissa of convergence, Euler products.  
\item[Ch.\ 4.] {\it Special functions in analytic number theory.}  Topics include: the gamma function $\Gamma(s)$,  log-gamma function $\log\Gamma(s)$,  digamma function $\Psi(s)$,  Riemann zeta function $\zeta(s)$, functional equation, zeros of $\zeta(s)$, Riemann hypothesis, Riemann xi function $\xi(s)$, prime zeta function $P(s)$, and special functions $\Ein(s)$, $\Ei(s)$, $E_1(s)$,  $\ERi(s)$, and $\Ri(x)$.  
\item[Ch.\ 5.]  {\it The analytic theory of primes.}   Topics include:  Riemann von--Mangoldt explicit formulas for $\pi_0(x)$, $\Pi_0(x)$ and $\psi_0(x)$,   prime number theorem with error bound,  Riemann--Siegel theta function and $Z$ function,  nontrivial zeros of $\zeta(s)$, the functions $N(T)$ and $S(T)$,  Riemann--von Mangoldt formula,  Lambert $W$ function, Montgomery's pair correlation conjecture, abstract analytic number theory, arithmetic semigroup,  Dedekind zeta function.  
\end{enumerate}
 Part 2 ({\it Algebraic asymptotic analysis})  extends results in algebraic  asymptotic analysis,  introduces the notion of logexponential degree, and studies asymptotic continued fraction expansions and their applications to analytic number theory.  Most of the results in Part 2 are new.
\begin{enumerate}
\item[Ch.\ 6.] {\it Logexponential degree}. The results in Chapter 6 comprise the main new analytic-algebraic tools introduced in this book and form the basis for the investigation in Parts 2 and 3.    The chapter introduces the iterated logarithmic degree and logexponential degree formalisms and states and proves many of their properties, relating them to the various asymptotic relations $O$, $o$, $\asymp$, $\sim$, etc.,  to the operations $+$, $-$, $\cdot$, $/$, and $\circ$, on functions, and to Hardy's ordered field $\mathbb{L}$ of all {\it logarithmico-exponential functions} \cite{har3} \cite{har4}.   Section 6.1  furthers our study of the notions of the  degree $\deg f$ and lower degree  $\underline{\deg} \, f$ of a real function $f$ and provides some further uses of degree in real and complex analysis.  Section 6.3 establishes various fundamental properties of logexponential degree that are used in the remainder of the text. 
\item[Ch.\ 7.] {\it Asymptotic algebra}.  Chapter 7 explores applications of  the field of asymptotic algebra, e.g., asymptotic differential algebra, to the study of logexponential degree, and vice versa.    We use logexponential degree to provide universal properties for some important Hardy fields,  including  Hardy's ordered differential field $\mathbb{L}$, and we introduce the notion of the {\it logexponential degree} of a logarithmic-exponential transseries.    We also provide several axiomatic characterizations of the degree map.    Chapter 7 is not used in any later chapters,  the only exceptions being the definitions  in Section 7.3 of a {\it Hardian} function and of the ordered differential field $\mathbb{H} \supsetneq \mathbb{L}$ of all  {\it universally Hardian} functions,  along with Theorem \ref{hardintth} and Propositions \ref{hardianexactlog} and \ref{oexppropstrong}.  
\item[Ch.\ 8.]  {\it Asymptotic continued fraction expansions}.      Chapter 8, which is based on \cite{ell0} and \cite{ell}, is a discussion of asymptotic continued fraction expansions and their applications to the prime counting function and related functions.  We show, for example, that,   for each positive integer $n$,  two well-known continued fraction expansions of the exponential integral function $E_n(z)$ correspondingly yield two (divergent) asymptotic continued fraction expansions of the prime counting function.    Chapter 8 is not used in any later chapters, with the sole exception of Proposition \ref{mertensprop}. 
\end{enumerate}
 Part 3 ({\it Applications of algebraic asymptotic analysis to number theory}) applies Parts 1  and 2 to the study of  various important functions arising in number theory.    The main goal of Part 3 is to expess the logexponential degree of the various functions  of arising in  number theory in terms of  the logexponential degree of $f$ for as  few ``logexponential primitives'' $f$ (e.g., $f = \li-\pi$) as possible.
\begin{enumerate}
\item[Ch.\ 9.]  {\it The prime counting function $\pi(x)$ and related functions}. Chapter 9 uses the degree formalisms to study the prime counting function $\pi(x)$ and various functions closely related to $\pi(x)$, including  the first and second Chebyshev functions $\vartheta(x)$ and $\psi(x)$ and Riemann's prime counting function $\Pi(x)$.   
\item[Ch.\ 10.] {\it Summatory functions}. Chapter 10 uses the degree formalisms to study the summatory function  $\sum_{n \leq x} f(n)$ of various arithmetic functions $f(n)$, including the M\"obius function $\mu(n)$, the Liouville lambda function $\lambda(n)$,  the divisor function $d(n)$,  and Euler's totient $\phi(n)$.
\item[Ch.\ 11.] {\it The Riemann zeta function $\zeta(s)$}.  Chapter 11 uses the degree formalisms to study the Riemann zeta function $\zeta(s)$, the Riemann zeta zero counting functions $N(T)$ and $S(T)$,  the ordinates $\gamma_n$ of the zeros of $\zeta(s)$, and the gaps $\gamma_{n+1}-\gamma_n$ between them. 
\item[Ch.\ 12.] {\it Primes in intervals, the $n$th prime, and the $n$th prime gap}.   Chapter 12 uses the degree formalisms to study primes in intervals, the prime listing function $p_n$, the prime gap function $g_n =  p_{n+1}-p_n$, and the maximal prime gap function $G(x) = \max_{p_k \leq x} g_k$.    Section 12.1 explores the problem of determining which real  functions $h$ satisfy $\pi(x+h(x)) - \pi(x) \sim \frac{h(x)}{\log x} \ (x \to \infty)$.
\item[Ch.\ 13.] {\it  Diophantine approximation and continued fractions}.    Chapter  13 provides applications of degree and logexponential degree to  Diophantine approximation and regular continued fractions, e.g., to the study of irrationality measure,  Markov constants, L\'evy constants,  $Q$-order of convergence, rates of convergence,  quadratic irrationals, badly approximable numbers, well approximable numbers, and very well approximable numbers.   We show, for example, that the irrationality measure of any irrational number $\alpha$ is equal to $\deg \frac{n}{\Vert n \alpha \Vert}$, where $\Vert x\Vert$ denotes the distance from $x \in \RR$ to the nearest integer, and we provide several equivalent characterizations of the logexponential degree of $\frac{n}{\Vert n \alpha \Vert}$.   We also pose some conjectures generalizing Roth's theorem concerning the rational approximation of algebraic numbers.   No prior knowledge of Diophantine approximation is assumed in this chapter.
\item[Ch.\ 14.] {\it Conjectures}.  Chapter 14, which is largely data-driven, uses the degree formalisms to express evidence, both numerically and graphically, and in a novel way, for some of the conjectures discussed in this book, including the Riemann hypothesis and various extensions of the Riemann hypothesis.
\end{enumerate}

\section{Motivation and detailed summary}

In this somewhat lengthy technical section, we provide a detailed summary of the text, along with our primary motivations for introducing the notions of degree and logexponential degree.   

Let $\pi: \RR_{\geq 0} \longrightarrow \RR$ denote the function that for any $x \geq 0$ counts the number of primes less than or equal to $x$: $$\pi(x) =   \# \{p \leq x: p \mbox{ is prime}\}, \quad \forall x \geq 0.$$  The function $\pi(x)$ is known as the {\bf prime counting function}. The celebrated {\bf prime number theorem}, proved by de la Vall\'ee Poussin \cite{val1} and Hadamard \cite{had}  in 1896,  states that
\begin{align*}
\pi(x) \sim \frac{x}{\log x} \ (x \to \infty),
\end{align*}
where $\log x$ is the natural logarithm.    It is known, however, that the {\bf logarithmic integral function} $$\li(x) = \int_0^x \frac{dt}{\log t}, \quad \forall x \geq 0,$$
where the Cauchy principal value of the integral is assumed, provides a better approximation to $\pi(x)$ than any algebraic function of $\log x$.  The {\bf prime number theorem with error bound},\index{prime number theorem with error bound}  proved by de la Vall\'ee Poussin in 1899 \cite{val2}, states that  the error $\li(x)-\pi(x)$ in the approximation $\pi(x) \approx \li(x)$ is bounded above by
\begin{align}\label{PNTET10}
\li(x)-\pi(x) = O \left(\frac{x}{e^{C\sqrt{\log x}}} \right) \ (x \to \infty)
\end{align}
for some constant $C  > 0$.  This has since been improved to
\begin{align}\label{bestPNT0}
\li(x)  - \pi(x) = O\left(\frac{x}{e^{ A(\log x)^{3/5}(\log \log x)^{-1/5}}}\right) \ (x \to \infty),
\end{align}
where  $A  = 0.2098$ \cite{ford}, which is the strongest known $O$ bound on the error $\li(x)-\pi(x)$ to date.  Proofs of such bounds on the error are based on the {\it Riemann--von Mangoldt explicit formula for $\pi(x)$} in terms of the zeros of the {\it Riemann zeta function} $\zeta(s)$ \cite{rie} and  rather sophisticated methods for verifying zero-free regions for $\zeta(s)$ in the {\bf critical strip} $\{s \in \CC: 0 \leq \operatorname{Re} s \leq 1\}$.    As is well known, Riemann proved, in his landmark paper  \cite{rie} of 1859, that  the zeros of $\zeta(s)$, besides the negative even integers, are all non-real and lie in the critical strip.  The non-real zeros of $\zeta(s)$,  or, equivalently, the zeros of $\zeta(s)$ that lie in the critical strip, are known as the {\bf nontrivial zeros of $\zeta(s)$}. The celebrated {\bf Riemann hypothesis}, conjectured by Riemann in his paper, states that all nontrivial zeros of $\zeta(s)$ lie on the {\bf critical line} $\{s \in \CC: \operatorname{Re} s = \frac{1}{2}\}$.    

The problem of settling the Riemann hypothesis is widely regarded as one of the most important, if not the most important, unsolved problems in mathematics today.   One reason for this is that there are hundreds of statements known to be equivalent to the Riemann hypothesis (many of which are collected in \cite{broughan} \cite{broughan2}).  Thus far, none of them stand out as a ``correct'' approach, i.e., as an approach that is most likely to lead to an eventual proof or disproof of the conjecture.   Many avenues that were once thought promising ultimately were found to be just another reformuation of the same longstanding problem.   (In more recent times, mathematicians and physicists have been seeking ways to think about the zeros of $\zeta(s)$ from the perspective of quantum physics \cite{bend} \cite{berry} \cite[Part 5]{miller} \cite{wolch} \cite{wolf}.)   Probably the most important known equivalent of the Riemann hypothesis was found in 1901, when  von Koch proved \cite{koch} that the  Riemann hypothesis is equivalent to the error bound
\begin{align}\label{vk0}
\li(x)- \pi(x) = O(\sqrt{x} \log x) \ (x \to \infty),
\end{align}
which to date is  the strongest bound on the error $\li(x)- \pi(x)$ that is widely conjectured to hold.      

Now, let
\begin{align}\label{delta}
\Theta = \sup\{\operatorname{Re} \rho: \rho \in \CC\backslash \RR, \, \zeta(\rho) = 0\}\index[symbols]{.f td@$\Theta$}
\end{align}
denote the supremum of the real parts of the nontrivial zeros of $\zeta(s)$.   Well-known results from Riemann's paper \cite{rie} concerning the zeros of $\zeta(s)$ imply that
$$\tfrac{1}{2} \leq \Theta \leq 1,$$
and that the Riemann hypothesis is equivalent to $\Theta = \frac{1}{2}$.
It is also well known \cite[Theorem 15.2 and Exercise 13.1.1.1]{mont} that von Koch's  equivalent (\ref{vk0}) of the Riemann hypothesis generalizes to the fact that $\Theta$ is  given by 
\begin{align}\label{ess}
\Theta = \min\left\{t \in \RR: \li(x)- \pi(x) = O(x^t \log x) \ (x \to \infty)\right\}
\end{align}
and also by
\begin{align}\label{ess2}
\Theta = \inf\left\{t \in \RR: \li(x)- \pi(x) = O(x^t) \ (x \to \infty)\right\}.
\end{align}
Thus, the constant $\Theta$  carries vital information about both  the Riemann zeta function  and the distribution of the prime numbers.   Indeed,  many known equivalents of the Riemann hypothesis can be generalized to unconditional results concerning  the constant $\Theta$, e.g.,  quintessentially,   the Riemann hypothesis equivalent (\ref{vk0}) generalizes to the expressions (\ref{ess}) and (\ref{ess2}) for $\Theta$.   Throughout this book, and following several other authors,  like Ingham in \cite{ing2}, we use $\Theta$ to denote the constant defined by  (\ref{delta}) above.   Because $\Theta$ is such an important constant, we also  give it a name:  the {\bf Riemann constant}.\index{Riemann constant $\Theta$}    The  results noted above form the basis for the following research program, which is one of our main concerns in Part 3.

\begin{problem}
Given a known equivalent of the Riemann hypothesis, generalize the equivalence to an unconditional statement regarding the Riemann constant $\Theta$.
\end{problem}

The results above also serve to motivate a very natural notion of the {\it degree} of a real function.   Let $f: X \longrightarrow \RR$ be any real function whose domain $X$ is a subset of $\RR$ that is not bounded above.   We define the {\bf degree of $f$} to be the extended real number
$$\deg f = \inf\{t \in \RR: f(x) = O(x^t) \ (x \to \infty)\} \in  \overline{\RR} = [-\infty,\infty],$$ 
which, equivalently, can be defined by
$$\deg f = \limsup_{x \to \infty} \frac{\log |f(x)|}{\log x}.$$
 This notion of degree  extends the usual definition of the degree of a polynomial.   
Given the definition of degree above, statement (\ref{ess2}) concerning the Riemann constant $\Theta$ is equivalent to $$\Theta = \deg (\li-\pi).$$
Statements of the form $$f(x) = O(x^{d + \varepsilon}) \ (x \to \infty), \quad \forall \varepsilon > 0,$$ 
and of the form
$$f(x) = o(x^{d + \varepsilon}) \ (x \to \infty), \quad \forall \varepsilon > 0,$$ 
appear throughout analytic number theory,  and it is common but unstated knowledge to analytic number theorists that both of the statements above are equivalent to
$$\limsup_{x \to \infty}  \frac{\log|f(x)|}{\log x} \leq d.$$
Thus, according to our definition of degree,  all three of the above statements are equivalent to  $\deg f \leq d$.   Loosely speaking,  $\deg f$ is a measure of the degree of freedom of $f$,  and thus the Riemann constant $\Theta$ is a measure of the degree of freedom of $\li(x)-\pi(x)$, i.e., of the degree of freedom that $\pi(x)$ has relative to $\li(x)$.   The smaller the Riemann constant, the less the degree of freedom,  and the more constrained $\pi(x)$ is.

Another important example of a generalization of a Riemann hypothesis equivalent to a statement about the Riemann constant $\Theta$ concerns the {\bf Mertens function} $$M(x) = \sum_{n \leq x} \mu(n), \quad \forall x \geq 0.$$
It is widely known \cite[Theorem 4.16]{broughan} that the Riemann hypothesis is equivalent to 
$$M(x) = O(x^{1/2+\varepsilon}) \ (x \to \infty)$$
for all $\varepsilon > 0$,  which, using our degree terminology, is equivalent to  $\deg M \leq \frac{1}{2}$.   More generally, it is known, but seldom disclosed in the literature, that $$\deg M = \Theta.$$  See Theorem \ref{MLProp} for a sketch of the proof.   Yet another illustrative example concerns the {\bf Riesz function}\index{Riesz function}
$$F(x) = \sum_{k = 1}^\infty \frac{(-1)^{k-1}x^k}{(k-1)! \zeta(2k)} = x \sum_{k = 1}^\infty \frac{\mu(k)}{k^2}e^{-x/k^2},$$
which is an analytic function on all of $\RR$, as the given Taylor series has radius of convergence $\infty$.   In \cite{riesz}, Riesz proved that the Riemann hypothesis is equivalent to $\deg F \leq \frac{1}{4}$, which is known as the {\bf Riesz criterion}\index{Riesz criterion}.   A careful reading of the proof reveals that, unconditionally, one has $$\deg F = \frac{\Theta}{2},$$ which, of course, implies Riesz's equivalence.

 A theme of this book is that the degree formalism and its generalizations to {\it iterated logarithmic degree} and {\it logexponential degree} are useful for motivating and investigating important questions about various well-studied number-theoretic functions, including the prime counting function $\pi(x)$, the prime listing function $p_n$, the prime gap function $g_n = p_{n+1}-p_n$,  the M\"obius function $\mu(n)$, the Liouville lambda function $\lambda(n)$, the divisor function $d(n)$, the prime Omega function $\Omega(n)$, the prime  omega function $\omega(n)$, Euler's totient $\phi(n)$,  the Riemann zeta function $\zeta(s)$,   the Riemann zeta function zero ordinate listing function $\gamma_n$, and the Riemann zeta  zero ordinate gap function $\gamma_{n+1}-\gamma_n$.

The degree formalism  alone is   useful for formulating and gathering numerical and graphical evidence for conjectures.  When one is uncertain what the degree of a given function $f$ is (e.g., see Example \ref{degexamples}),  it can be helpful to compute or to graph the function
$$L_f(x) = \frac{\log|f(x)|}{\log x}\index[symbols]{.g k@$L_f(x)$}$$
(which is the unique function $g(x)$ satisfying $|f(x)| = x^{g(x)}$) for as many and as large values of $x \in \dom f$ as is feasible.  Based on such information, one can then try  to conjecture what  the value  of
$$\deg f = \limsup_{x \to \infty} L_f(x)$$
might be.  

Let us apply this method to the function $\li(x)-\pi(x)$.   Let $\Ri(x)$ denote Riemann's approximation $$\Ri(x)=\sum_{n=1}^\infty \frac{ \mu(n)}{n} \li(x^{1/n})$$ to $\pi(x)$ \cite{rie},  which is studied in Section 4.5, where $\mu(n)$ denotes the {\it M\"{o}bius function}.   One has $\Ri(x) \sim \li(x) \ (x \to \infty)$ and
$\li(x)-\Ri(x) \sim \frac{\sqrt{x}}{\log x} \ (x \to \infty)$.
In Figure \ref{eureka100} below, we provide graphs of the functions $$L_{\li-\pi}(e^x) = \frac{\log|\li(e^x)-\pi(e^x)|}{x} \quad \text{and} \quad L_{\li-\Ri}(e^x) = \frac{\log|\li(e^x)-\Ri(e^x)|}{x},$$
which are just the graphs of the functions $L_{\li-\pi}(x)$ and $L_{\li-\Ri}(x)$  but on a lin-log scale.    The constant $\Theta$ is exactly the lim sup of the blue curve as $x \to \infty$, that is, one has
$$\Theta  = \limsup_{x \to \infty} L_{\li-\pi}(e^x).$$
Since $\li(x)-\Ri(x) \sim \frac{\sqrt{x}}{\log x} \ (x \to \infty)$, one has
$\lim_{x \to \infty} L_{\li-\Ri}(e^x) = \tfrac{1}{2}$.  Thus, the black curve in Figure \ref{eureka100} tends to $\tfrac{1}{2}$ as $x \to \infty$.   Since the Riemann hypothesis is equivalent to  $\deg(\li-\pi) = \tfrac{1}{2}$, the Riemann hypothesis holds if and only if the blue curve has a lim sup of $\tfrac{1}{2}$, if and only if the blue curve minus the black curve has a lim sup of $0$, as $x \to \infty$.  This provides a new way to visualize the Riemann constant $\Theta$, along with the Riemann hypothesis.

\begin{figure}[h!]
\includegraphics[width=70mm]{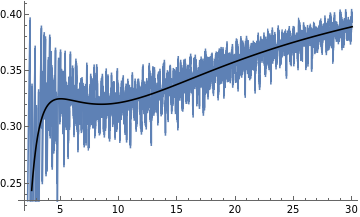}
    \caption{\centering Graphs of $L_{\li-\pi}(e^x)$ and $L_{\li-\Ri}(e^x)$}
\label{eureka100}
\end{figure}

 It should be stressed that, at the writing of this book, there is no guarantee that the Riemann hypothesis is true, i.e., that $\Theta = \frac{1}{2}$.  If the Riemann hypothesis is false, then the  most obvious next best guess for $\Theta$ is $\Theta = 1$, which we dub the {\bf anti-Riemann hypothesis}.  This hypothesis holds that the prime number theorem with error bound is nearly the best of its kind, and thus the hypothesis seems to be the most compelling alternative to the Riemann hypothesis.  A third alternative is that the truth is somewhere in between, i.e., $\frac{1}{2} < \Theta < 1$---perhaps, say, $\Theta = \frac{2}{3}$ or $\Theta = \frac{\pi}{4}$---and the closer $\Theta$ is to $\frac{1}{2}$, the closer the Riemann hypothesis is to being true.   In particular, the problem of computing the Riemann constant $\Theta$ is a natural generalization of the problem of settling the Riemann hypothesis (which asks only whether or not $\Theta = \frac{1}{2}$).

Let us assume for the moment that the Riemann hypothesis is true.   Based on values of $\pi(x)$ that have been either computed or estimated,  it might appear that the  error bound (\ref{vk0})  from von Koch's Riemann hypothesis equivalent should be improved (conjecturally) to $\li(x)- \pi(x) = O\left( \frac{\sqrt{x}}{\log x} \right) \ (x \to \infty)$.   In 1910, for example, Hardy wrote that ``there is reason to anticipate that'' this error bound holds \cite[p.\ 48]{har3}.   More recently,  in 1994, Riesel wrote: ``Judging only from the values [given in a table] we might even try to estimate the order of magnitude of $\li(x)-\pi(x)$ and find it to be about $\sqrt{x}/\log x$.  However, {\it for large values of $x$, this is completely wrong!}'' \cite[p.\ 52]{ries}.   Indeed, the $O$ bound above fails to hold because of Littlewood's 1914 result \cite{litt} that 
\begin{align}\label{litt0}
\li(x)-\pi(x) = \Omega_{\pm} \left(\frac{\sqrt{x}\, \log \log \log x}{\log x} \right) \ (x \to \infty),
\end{align}
where  one writes $f(x) = \Omega_{\pm}(g(x)) \ (x \to a)$
if $\limsup_{x \to a} \frac{f(x)}{|g(x)|} $ is positive and $\liminf_{x \to a} \frac{f(x)}{|g(x)|}$ is negative (both possibly infinite).      Today, Littlewood's result stands as  one of   a great number of occurrences of iterated logarithms  in analytic number theory.   Another example, proved by Rankin in 1938 \cite{rankin}, is the result that  $$p_{n+1}-p_n \neq o\left(\frac{\log n \log \log n\log \log \log \log n}{(\log \log \log n)^{2}}\right) \ (n \to \infty),$$
which  was  strengthened to 
\begin{align}\label{fgkmtgap}
p_{n+1}-p_n \neq o\left(\frac{\log n \log \log n\log \log \log \log n}{\log \log \log n}\right) \ (n \to \infty)
\end{align}
in 2016 by Ford, Green, Konyagin, Maynard, and Tao \cite{fgktm}.

One of the surprising consequences of Littlewood's result is that  $\li(x)-\pi(x)$ changes sign an infinite number of times, even though no one knows any specific value of $x \geq 2$ for which $\li(x)-\pi(x) < 0$.   Even Gauss had been misled into believing that $\li(x) -\pi(x)$ is positive for all $x \geq 2$.   We now know that the infimum of all $x \geq 2$ such that $\li(x)-\pi(x) < 0$, known as {\bf Skewes' number}, is at least $10^{20}$ and at most $s = 1.3971671494 \cdot 10^{316}$ \cite{kul}.   A lesson one learns from this  is that numerical considerations sometimes mean very little in analytic number theory, and this is especially true when iterated logarithms  are lurking in the background.  Nevertheless, despite these warnings, we make the following conjecture (albeit contingently upon the Riemann hypothesis).

\begin{conjecture}\label{eurekaconjecture} 
For any real number $t$, one has
$$\li(x)-\pi(x) = o\left(\sqrt{x}\, (\log x)^{t} \right) \ (x \to \infty)$$
if (and only if) $t > -1$.
\end{conjecture}

The conjecture above expresses the feeling that, if the Riemann hypothesis is true, then the magnitude of $\li(x)-\pi(x)$ should be sufficiently close to $\frac{\sqrt{x}}{\log x}$ without contradicting Littlewood's result.  We also assert, with somewhat less confidence, the following generalization of Conjecture \ref{eurekaconjecture}.

\begin{conjecture}\label{eurekaconjecture2} 
For any real number $t$, one has 
$$\li(x)-\pi(x) = o\left(\frac{\sqrt{x}\, (\log \log x)^{t}}{\log x} \right) \ (x \to \infty)$$
if (and only if) $t > 0$.   Consequently, there exists a $\delta_3 \in [1,\infty]$ such that
$$\li(x)-\pi(x) = o\left(\frac{\sqrt{x}\, (\log \log \log x)^{t}}{\log x} \right) \ (x \to \infty)$$
for all $t > \delta_3$ but for no $t < \delta_3$.
\end{conjecture}

It is easy to see that Conjecture \ref{eurekaconjecture2}  implies Conjecture \ref{eurekaconjecture} and, by von Koch's 1901 result, either conjecture implies the Riemann hypothesis.    We should thus qualify these conjectures as contingent on the Riemann hypothesis.  The idea motivating Conjectures \ref{eurekaconjecture} and  \ref{eurekaconjecture2}   is that there is an infinite sequence  $$\Theta = \delta_0, \delta_1, \delta_2, \delta_3, \ldots \in \overline{\RR}$$ of  invariants of $\li(x)-\pi(x)$, where $\delta_k$, roughly, is the ``degree of $\log^{\circ k}x$ occurring in $\li(x)-\pi(x)$.''  Conjecture \ref{eurekaconjecture} (resp., Conjecture \ref{eurekaconjecture2}) is equivalent to the statement that the sequence $\delta_0, \delta_1, \delta_2, \delta_3, \ldots$ begins $\frac{1}{2}, -1, \delta_2, \delta_3, \ldots$ (resp., $\frac{1}{2}, -1, 0, \delta_3, \ldots$).  Here, one has $$\delta_k = \degl_k (\li-\pi),\index[symbols]{.f tg@$\delta_k$}$$ where  we define $\degl_k f$ recursively for any real function $f$ whose domain is a subset of $\RR$ that is not bounded above, as follows.  First, we let  $f_{[0]} = f$.    Suppose that $f_{[k]}$ is defined, and set $d_k = \dege_k f = \deg f_{[k]}.$  We then let
$$f_{[k+1]}(x) =   \left.
 \begin{cases}
    f_{[k]}(e^x) e^{-d_k x}& \text{if } d_k \neq \pm \infty \\
    f_{[k]}(x) & \text{if } d_k =\pm \infty.
 \end{cases}
\right.$$
This defines $ f_{[k]}$ and $\degl_k f \in \overline{\RR}$  for all  nonnnegative integers $k$.
We call $\degl_k f$ the {\bf (iterated) logarithmic degree of $f$ of order $k$}.

By definition, the constant  $\delta_1 =  \degl_1 (\li-\pi)$ is the infimum of all $t \in \RR$ such that  
 $$\li(x)-\pi(x) = \displaystyle O\left(x^\Theta (\log x)^t \right) \ (x \to \infty),$$ and  (\ref{ess}) implies that $\delta_1 \leq 1$.     If  the Riemann hypothesis is true, then Littlewood's result (\ref{litt0}) requires that $\delta_1 \geq -1$.    Thus,  Conjecture \ref{eurekaconjecture}  is equivalent to the conjecture that $\Theta = \frac{1}{2}$ (the Riemann hypothesis) and $\delta_1 = -1$.  Likewise,  Conjecture \ref{eurekaconjecture2}  is equivalent to the conjecture that $\Theta = \frac{1}{2}$, $\delta_1 = -1$, and $\delta_2 = 0$.  On the other hand, it is known \cite{gros} that, if the Riemann hypothesis is false, then $\li(x)-\pi(x) =  O\left(\frac{x^\Theta}{\log x}\right) \ (x \to \infty)$ and thus $\delta_1 \leq -1$.  It follows that, if $\delta_1 > -1$, then the Riemann hypothesis is true, while if $\delta_1 < -1$, then the Riemann hypothesis is false.  Any result that were to imply $\delta_1 \neq -1$, then, would have to settle the Riemann hypothesis.   However,  it is conceivable that improvements on the unconditional inequalities $\frac{1}{2} \leq \Theta \leq 1$ and $\delta_1 \leq 1$ could be proved absent a proof or disproof of the Riemann hypothesis.

To address the possibility that the Riemann hypothesis (and Conjectures \ref{eurekaconjecture}  and \ref{eurekaconjecture2}) may be false, we also assert the following conjecture.  

\begin{conjecture}\label{eitheror}
Either the Riemann hypothesis or the anti-Riemann hypothesis is true.  Equivalently, either $\deg(\li-\pi) =  \frac{1}{2}$ or  $\deg(\li-\pi) = 1$.
\end{conjecture}

In the  event that the conjecture above is false, the constant $\Theta = \deg(\li-\pi)$ would be quite a noteworthy constant, indeed!

 In Chapter 13, we  show that  Conjectures \ref{eurekaconjecture}  and \ref{eurekaconjecture2} are modestly supported by numerical evidence.     In Chapter 9, we show  that a 1979 conjecture of Montgomery \cite[Conjecture, p.\ 16]{mont1}, namely, that
\begin{align*}
\limsup_{x \to \infty} \frac{x-\psi(x)}{\sqrt{x}\, (\log \log \log x)^2} = \frac{1}{2\pi} = -\liminf_{x \to \infty} \frac{x-\psi(x)}{\sqrt{x}\, (\log \log \log x)^2},
\end{align*}
where $$\psi(x) =  \sum_{k = 1}^\infty \sum_{p^k \leq x} \log p \sim x \ (x \to \infty)$$ is the {\bf second Chebyshev function},  implies that Conjecture \ref{eurekaconjecture2}  holds with $\delta_3 = 2$ and that the entire sequence of iterated logarithmic degrees is given by $\frac{1}{2}, -1, 0, 2, 0, 0, 0, \ldots$.  Although  Montgomery's conjecture implies Conjecture \ref{eurekaconjecture2} and the latter conjecture provides some support for the former, the motivations for these two conjectures are very different from one another.  Indeed, we came to formulate Conjectures \ref{eurekaconjecture} and \ref{eurekaconjecture2} several months before we learned of Montgomery's conjecture, which is a far more ambitious conjecture that was motivated on very different grounds.    There are uncountably many possible scenarios in which Conjecture \ref{eurekaconjecture2} might hold while Montogomery's conjecture fails, including the current best case scenario, implied by the conjecture \cite[(11)]{stoll} of Stoll and Demichel,   according to which the iterated logarithmic degree sequence for $\li(x)-\pi(x)$ is given by $\frac{1}{2},-1,0,1,0,0,0,\ldots$.  Furthermore, as discussed throughout this book, whether or not Conjectures \ref{eurekaconjecture} and \ref{eurekaconjecture2} are true, the degree concerns that motivate the conjecture illuminate many other fundamental questions regarding various number-theoretic functions besides the prime counting function.

A fundamental example of such a problem concerns the Mertens function $M(x)$.   Unlike the situation with $\degl_1(\li-\pi)$, no upper bound is known for $\degl_1 M$, not even on condition of the Riemann hypothesis.   Nevertheless, the conjecture \cite[(20)]{ng} of Gonek and Ng implies that the entire sequence $\degl_k M$  is given by $\frac{1}{2}, 0, 0, \frac{5}{4}, 0,0,0,\ldots$, while  conjectures of Good and Churchhouse \cite{good} and L\'evy \cite{lev} imply that the sequence is given by $\frac{1}{2}, 0, \frac{1}{2}, 0,0,0,\ldots$.  This provides modest support for the following conjecture.

\begin{conjecture}\label{Mconjecture}
One has $\deg M = \frac{1}{2}$ and $\degl_1 M = 0$.  Equivalently,  one has $M(x) = o(\sqrt{x} \,  (\log x)^t)\ (x \to \infty)$ if (and only if) $t > 0$.
\end{conjecture}

Note that, by Theorem \ref{multzero} of Section 10.2, if Conjecture \ref{Mconjecture} holds, then the Riemann hypothesis holds and all of the zeros of the Riemann zeta function are simple.

Now, let  $f$ be any real function whose domain is a subset of $\RR$ that is not bounded above.   In order to better handle situations where $\degl_k f = \pm \infty$ for some $k$, we define a natural refinement $\dege f$ of $\degl f$, as follows.   Let $f_{(0)} = f$.  Suppose that $f_{(k)}$ is defined, and set $d_k = \dege_k f = \deg f_{(k)}.$  We then let
$$f_{(k+1)}(x) =   \left.
 \begin{cases}
    f_{(k)}(e^x) e^{-d_k x}& \text{if } d_k \neq \pm \infty \\
    \max( \log |f_{(k)}(x)|, 0)  & \text{if } d_k = \infty \\
 \displaystyle   -\frac{1}{\log |f_{(k)}(x)|} & \text{if } d_k =- \infty.
 \end{cases}
\right.$$
This defines $\dege_k f \in \overline{\RR}$  for all  nonnnegative integers $k$.
We call $\dege_k f$ the {\bf logexponential degree of $f$ of order $k$}.
Also, we write
$$\dege f = (\dege_0 f, \dege_1 f, \dege_2 f, \ldots) \in \prod_{n = 0}^ \infty \overline{\RR}$$
and
$$\degl f = (\degl_0 f, \degl_1 f, \degl_2 f, \ldots) \in \prod_{n = 0}^ \infty \overline{\RR},$$  
which we call the {\bf logexponential degree of $f$} and the {\bf (iterated) logarithmic degree of $f$}, respectively.
We also endow the set   $\prod_{n = 0}^ \infty \overline{\RR}$ with the lexicographic (total) ordering.  Note that, if $\degl_k f \neq \pm \infty$ for all $k < n$, then $\dege_k f = \degl_k f$ for all $k \leq n$.  However, if $\degl_{n} f = \pm \infty$, then $\dege_{n+1} f$ might not equal $\degl_{n+1} f$, and the definition of $f_{(n+1)}(x)$ as above is designed to ``tame'' the function $f_{(n)}$ by applying a log to $|f_{(n)}|$ appropriately.   Thus, $\degl f$ is tantamount to the truncation of $\dege f$ at the $n$th coordinate for the smallest $n$, if any, such that $\dege_n f= \pm \infty$.   Note  that $$\dege(f+g) \leq  \max(\dege f, \dege g)$$ provided that $\dom f \cap \dom g  $ is not bounded above.  Moreover, $f(x) = O(g(x)) \ (x \to \infty)$ implies  $\dege f \leq \dege g$, which in turn implies $\degl f \leq \degl g$.  

 The seemingly {\it ad hoc} definition  of $\dege$ is motivated by the following examples.

\begin{example}\label{exex0} \
\begin{enumerate}
\item For every positive integer $n$, let $T(n) = \#\{ab: a,b \in \{1,2,3,\ldots,n\}\}$ denote the number of distinct integers in the $n \times n$ multiplication table.  In 2008, Ford proved \cite{ford3} that
$$T(n) \asymp \frac{n^2}{(\log n)^c(\log \log n)^{3/2}} \ (n \to \infty),$$
where 
$$c =1-{\frac {1+\log \log 2}{\log 2}}=0.086071332055\ldots.$$
From this it follows that $\dege T = (2,-c, -\tfrac{3}{2}, 0,0,0,\ldots)$.
\item By (\ref{bestPNT0}), one has $\dege (\li-\pi) \leq (1,-\infty,-\tfrac{3}{5},\tfrac{1}{5},0,0,0,\ldots)$.  Moreover, Littlewood's result (\ref{litt0}) implies that $\dege (\li-\pi) \geq (\tfrac{1}{2},-1,0,1,0,0,0,\ldots)$.
\item Walfisz showed in  \cite{wal} that
$$M(x) = O\left(\frac{x}{e^{ A(\log x)^{3/5}(\log \log x)^{-1/5}}}\right) \ (x \to \infty),$$
for some $A > 0$.  From this it follows that $\dege M \leq (1,-\infty,-\tfrac{3}{5},\tfrac{1}{5},0,0,0,\ldots)$.
\item The arithmetic function $d(n) = \sum_{d \mid n} 1$ is called the {\bf divisor function}, since $d(n)$ equals the number of positive divisors of $n$.  \cite[Theorem 317]{har} states that
\begin{align*}
\limsup_{n \to \infty} \frac{\log d(n)\log \log n}{\log n} = \log 2,
\end{align*}
which implies that $\dege  d(n) = (0,\infty,1,-1,0,0,0,\ldots)$.
\item Let $f(x)$ be the function that on the interval $[N,N+1)$  assumes the value $e^x$ for even integers $N$ and $e^{-x^2}$ for odd integers $N$.  One has $\dege f = (\infty, 1,0,0,0,\ldots)$, since the function $ \max( \log |f(x)|, 0) $ has logexponential degree $(1,0,0,0,\ldots)$.  The function $\log |f(x)|$, on the other hand, has logexponential degree $(2,0,0,0,\ldots)$. 
\item Let $f(x)$ be the function that on the interval $[N,N+1)$  assumes the value $e^{-x}$ for even integers $N$ and $e^{-x^2}$ for odd integers $N$.  One has $\dege f = (-\infty, -1,0,0,0,\ldots)$, since the function $-\frac{1}{ \log |f(x)|} $ has logexponential degree $(-1,0,0,0,\ldots)$.  The function $\log |f(x)|$, on the other hand, has logexponential degree $(2,0,0,0,\ldots)$. 
\end{enumerate}
\end{example}

Examples (5) and (6) above provide some explanation as to why we defined $f_{(k+1)}(x)$ as we did in the definition of $\dege f$, rather than, say, as $\log |f(x)|$, when $\deg f_{(k)} = \pm \infty$.  Our results in Chapter 6,  which comprise the main new analytic-algebraic tools introduced in this book,   demonstrate that the logexponential degree formalism is quite natural, and, in some sense, inevitable.   Theorem \ref{thm:ledege_axioms}, for example, provides a natural axiomatization of logexponential degree.    In Chapter 7, we characterize the degree map ``universally''  in various ways,  and we apply logexponential degree to the study of Hardy fields, logarithmic-exponential transseries, and real asymptotic differential algebra more broadly.   In the remainder of this section, we provide  some evidence for the main thesis of this text, namely, that the degree and logexponential degree formalisms are useful in analytic number theory.

Our quintessential example concerns the prime counting function, to which we  associate the invariants $$\Theta_k = \dege_k(\li-\pi) \in \overline{\RR}$$ for all $k$.  The {\bf Riemann constants} $\Theta_k$\index{Riemann constant $\Theta_k$}\index[symbols]{.f th@$\Theta_k$}  provide fine-tuned information about the difference $\li-\pi$, more so than do the constants $\delta_k = \degl_k (\li-\pi)$.    For several important number-theoretic functions $f$, we are able to express all of the constants $\dege_k f$ in terms of the constants $\Theta_k$ and vice versa.   Examples include the following.

\begin{example}\label{EE1}  Let $k$ be a nonnegative integer.
\begin{enumerate}
\item One has
$$\dege_k(x - \psi(x)) = \begin{cases} \Theta_k  & \text{ if } k \neq 1 \\
  \Theta_k +1 &  \text{ if } k = 1,
\end{cases} $$
where $\psi(x)$ is the second Chebyshev function.
\item Let $p_n$ for every positive integer $n$ denote the $n$th prime.  It is well known that the prime number theorem is equivalent to $p_n \sim n \log n \ (n \to \infty)$.  However, the function $\li^{-1}n \sim n \log n$ is a much better approximation to $p_n$ than is $n \log n$, where $\li^{-1}: \RR \longrightarrow (1, \infty)$ is the inverse of the restriction of the logarithmic integral function $\li$ to the interval $(1,\infty)$.   We show in Section 12.2 that
$$\dege_k(p_n - \li^{-1}n) = \begin{cases} \Theta_k  & \text{ if } k \neq 1 \\
  \Theta_k + \Theta+1 &  \text{ if } k = 1,
\end{cases} $$
and thus, in particular, $\Theta = \deg(p_n - \li^{-1}n)$.
\item In Section 10.1,  we prove that
$$\dege_k\left(\sum_{p\leq x}\frac{1}{p} -\log \log x - M\right) = \begin{cases} \Theta_k-1  & \text{ if } k = 0\\
  \Theta_k  &  \text{ if } k \geq 1,
\end{cases} $$
where 
\begin{align}\label{MMCC}
M = \lim_{x \to \infty} \left(\sum_{p\leq x}\frac{1}{p} -\log \log x\right) = 0.261497212847\ldots
\end{align}
is the {\bf Meissel--Mertens constant}.  Similarly, we prove that
$$\dege_k \left(e^{\gamma} \prod_{p \leq x}\left(1-\frac{1}{p}\right) - \frac{1}{ \log x}\right) = \begin{cases} \Theta_k-1  & \text{ if } k \leq 1\\
  \Theta_k  &  \text{ if } k \geq 2.
\end{cases} $$   
\item  Since one has $\Ri(x) \sim \li(x) \ (x \to \infty)$ and
$\li(x)-\Ri(x) \sim \frac{\sqrt{x}}{\log x} \ (x \to \infty)$, Littlewood's theorem (\ref{litt0}) yields 
$$\dege_k (\li-\pi) = \dege_k(\Ri-\pi)$$
for all $k$.  Thus, although numerical evidence suggests that  Riemann's function $\Ri(x)$ is a better approximation of $\pi(x)$ than is $\li(x)$,  it is no better in the long run,  at least with respect to the logexponential degree formalism.  
\end{enumerate}
\end{example}

The operations $\degl$ and $\dege$ not only facilitate the proofs of various degree relationships, as in the examples above, but they also yield some  new relations between functions: one has the irreversible implications
\begin{align*}
\degl f < \degl g  \quad & \Longrightarrow \quad \dege f  < \dege g  \\
 \quad & \Longrightarrow \quad  f(x) = o(g(x)) \ (x \to \infty) & \\
 \quad & \Longrightarrow \quad f(x) = O(g(x)) \ (x \to \infty)  \\ 
 \quad & \Longrightarrow \quad \dege f  \leq \dege g \\
\quad & \Longrightarrow \quad \degl f \leq \degl g,
\end{align*} 
where the second implication holds provided that the function $g$ is sufficiently nice, e.g., if $g$  can be can be built from all real constants and the functions $\id$, $\exp$,  and $\log$  using the operations $+$, $\cdot$, $/$, and $\circ$.  Examples throughout this text, including Example \ref{EE1} above, show that the six relations above are all distinct.  For example,  one has
$$\dege(\li-\pi) = \dege(\Ri-\pi) = \dege \left(\frac{p_n-\li^{-1}n}{(\log n)^{\Theta+1}} \right) = \dege \left(x\sum_{p\leq x}\frac{1}{p} -x\log \log x - Mx\right),$$ but none of these four functions is $O$ of the others.  
Thus, the notions of degree, iterated logarithmic degree, and logexponential degree serve not only to quantify the asymptotic growth of real functions,  including those arising prominently in number theory,  but they also provide natural benchmarks for comparing the asymptotic behavior of such functions.

The following examples  provide further indication that computing the degree, no less the logexponential degree, of a given number-theoretic function is quite often a difficult problem.

\begin{example}\label{degexamples} \
\begin{enumerate}
\item Let $g_n = p_{n+1}-p_n$ denote the {\bf $n$th prime gap}.
It is known \cite{baker} that
$$g_n = O ((n \log n)^{21/40}) \ (n \to \infty).$$
In \cite{cramer1}, Cram\'er proved on condition of the Riemann hypothesis that
$$g_n = O (n^{1/2}(\log n)^{3/2}) \ (n \to \infty).$$
Thus, one has
$$\dege g_n \leq (\tfrac{21}{40},\tfrac{21}{40},0,0,0,\ldots),$$
while also
$$\dege g_n \leq (\tfrac{1}{2},\tfrac{3}{2},0,0,0,\ldots)$$
on condition of the Riemann hypothesis.   It is conjectured, however, that $g_n = O((\log n)^t)  \ (n \to \infty)$ for all $t > 2$, or, equivalently, that $\deg g_n = 0$ and $\dege_1 g_n \leq 2$.  Computing $\deg g_n$, no less $\dege_k g_n$ for all $k$, is a difficult open problem.
\item  Let $L(r)$ equal the number of lattice points in $\ZZ^2$ lying on or  inside the circle in $\RR^2$ of radius $r$ centered at the origin.   It is straightfoward to show that $L(r) \sim \pi r^2 \ (r \to \infty)$.  Let $H(r) = L(r)-\pi r^2$.  Gauss proved that $|H(r)
|\leq 2\sqrt{2}\pi r$ for all $r$ and therefore $H(r) = O(r) \ (r \to \infty)$. The {\bf Gauss circle problem} is (equivalent to) the problem  of computing  $\deg H$.  Hardy and Landau proved that $$H(r) \neq o(r^{1/2}(\log r)^{1/4}) \ (r \to \infty),$$ from which it follows that
$$\dege H \geq (\tfrac{1}{2},\tfrac{1}{4},0,0,0,\ldots).$$  In 2003, Huxley proved \cite{hux2} that  $H(r) = O(r^{131/208} (\log r)^{18627/8320} ) \ (r \to \infty)$ and therefore $\deg H \in [\tfrac{1}{2},\tfrac{131}{208}]$ and 
$$\dege H \leq (\tfrac{131}{208},\tfrac{18627}{8320},0,0,0,\ldots).$$    It is widely conjectured that $\deg H = \tfrac{1}{2}$.  
\item Using what we now call the  {\it Dirichlet hyperbola method} (Proposition \ref{hyperbola}), Dirichlet proved that
$$\sum_{n\leq x}d(n) = x \log x + (2 \gamma-1)x + O(x^{t}) \ (x \to \infty)$$
for $t = \frac{1}{2}$.  The {\bf Dirichlet divisor problem} is the problem of determining the infimum of all such $t$ for which the $O$ bound  above  holds, which is equal to the degree $\deg D$ of the function $D(x) = \sum_{n\leq x}d(n) - x \log x - (2 \gamma-1)x$.  Hardy proved in 1914 that $$D(x) = \Omega_{\pm}(x^{1/4}) \ (x \to \infty).$$  It follows that $\deg D \in [\frac{1}{4},\frac{1}{2}]$ and 
$$\dege D \geq (\tfrac{1}{4},0,0,0,\ldots).$$  In 2003, Huxley  used the same methods in his approach to the Gauss circle problem to show \cite{hux2} that
 $$D(x) = O(x^{131/416}(\log x)^{26947/8320}) \ (x \to \infty)$$ and therefore  $\deg D \in [\frac{1}{4},\frac{131}{416}]$ and
$$\dege D \leq (\tfrac{131}{416},\tfrac{26947}{8320},0,0,0,\ldots).$$  
 It is widely conjectured that $\deg D = \tfrac{1}{4}$.   
\item An integer $a$ is said to be a {\bf primitive root modulo $p$}, where $p$ is prime, if $a$ generates the (cyclic) group $(\ZZ/p\ZZ)^*$ of units in the field $\ZZ/p\ZZ$.   Let $g(p)$ denote the smallest positive primitive root mod $p$.  In 1962, Burgess proved  \cite{burg} that $$\deg g(p) \leq \tfrac{1}{4}.$$  
 Moreover, Fridlander (1949) and Sali\'e (1950) proved that $g(p) \neq o(\log p) \ (p \to \infty)$, and therefore $$\dege g(p) \geq (0,1,0,0,0,\ldots).$$   Shoup, in 1990--92, proved that $g(p) = O((\log p)^6) \ (p \to \infty)$, and therefore
$$\dege g(p) \leq (0,6,0,0,0,\ldots),$$ provided that the {\it extended Riemann hypothesis} is true \cite[Theorem 1.3]{shoup}.  It is thus widely conjectured that $\deg g(p) = 0$.
\item An integer $a$ not divisible by $p$, where $p \neq 2$ is prime, is said to be a {\bf quadratic residue modulo $p$} if  there exists an integer $x$ such that $x^2 \equiv a \ (\text{mod } p)$, and a  {\bf quadratic non-residue modulo $p$} if $a$ is not a quadratic residue modulo $p$.   The smallest positive quadratic non-residue $n_p$ modulo $p$ is always a prime number less than $p$.  In 1957, Burgess proved  \cite{burg1} that $$\deg n_p \leq \tfrac{1}{4\sqrt{e}}.$$  {\bf Vinogradov's conjecture} is (equivalent to) the statement that $\deg n_p = 0$.    In 2001,  Wedeniwski proved in his doctoral thesis that $n_p < \frac{3}{2} (\log p)^2$ for all primes $p >3$,  provided that the {\it generalized Riemann hypothesis} is true.   A lower bound for the growth of $ n_p$, due to Chowla, is $$n_p \neq o(\log p) \ (p \to \infty),$$ which Montgomery, in 1971, proved can be strengthened to $$n_p \neq o(\log p \log \log p) \ (p \to \infty)$$ on condition of the generalized  Riemann hypothesis.  Thus, under the generalized Riemann hypothesis, one has
$$(0,1,1,0,0,0,\ldots) \leq \dege n_p \leq (0,2,0,0,0,\ldots).$$
See \cite{mcgown} for an account of these results, in historical context, along with some explicit inequalities for $ n_p$.

\item  Given relatively prime positive integers $a$ and $d$ with $1 \leq a \leq d$, let $p(a,d)$ denote the smallest prime in the arithmetic progression $a,a+d, a+2d,a+3d,a+4d,\ldots$.   {\bf Linnik's theorem}\index{Linnik's theorem}  states that there exist positive constants $C$ and $L$ such that $p(a,d) \leq Cd^L$ for all relatively prime positive integers $a$ and $d$ with $1 \leq a \leq d$ \cite{linnik} \cite{linnik2}.  Equivalently, it says that the function $$p(d) = \max\{p(a,d): 1 \leq a \leq d, \, \gcd(a,d) = 1\}$$ is $O(d^L)$ for some $L$, that is, $\deg p(d)$ is finite.  It is known that $L = 5$ is admissible, and thus
$$\dege p(d)  \leq (5,0,0,0,\ldots).$$  
 According to \cite[p.\ 282]{rib}, in 1961--62, Prachar and Schinzel proved that 
$$p(d) \neq o \left( \frac{d \log d \log \log d \log \log \log \log d}{(\log \log \log d)^2}\right) \ (d \to \infty),$$
from which it follows that
$$\dege p(d)  \geq (1,1,1,-2,1,0,0,0,\ldots).$$   The constant $\deg p(d) \in [1,5]$ is known as {\bf Linnik's constant} \cite[p.\ 279]{rib}. 
In 1934, Chowla conjectured that $\deg   p(d) = 1$ \cite{chowla}.  
In 1963, Kanold conjectured that $p(d) \leq d^2$
for all $d$, or equivalently that, for all positive integers $a$ and $d$ with $\gcd(a,d) = 1$, there is at least one prime among the numbers $a,a+d,a+2d,\ldots,a+(d-1)d$ \cite{kanold}.  According to \cite[pp.\ 280--283]{rib}, Schinzel and Sierpi\'nski (in 1958) and Kanold (in 1963)  conjectured that $\deg p(d)= 2$, while Heath-Brown (in 1978) conjectured that $p(d) = O(d (\log d)^2) \ (d \to \infty)$ and therefore $\dege p(d)\leq (1,2,0,0,0,\ldots)$.    
In 1992, Heath-Brown proved \cite{heath}  that the generalized Riemann hypothesis implies that $p(d) = O (\phi(d)^2 (\log d)^2) \ (d \to \infty)$,
where $\phi$ is {\it Euler's totient}, and therefore that $\dege p(d)  \leq (2,2,0,0,0,\ldots)$.
Thus, the conjecture that $\deg  p(d) \leq 2$ is implied by any of the conjectures noted above.  In 1990, Granville and Pomerance conjectured the lower bound  $p(d) \gg \phi(d) (\log d)^2 \ (d \to \infty)$ \cite{granpom}, which, if true,  would yield $\dege p(d)\geq (1,2,0,0,0,\ldots)$.
\item The {\bf irrationality measure $\mu(\alpha)$} of an irrational number $\alpha$ is the infimum of all $t >0$ such that  there are only finitely many pairs $(a,b)$ of integers $a$ and $b$  with $b> 0$ and
$$\left|\alpha -\frac{a}{b} \right| < \frac{1}{b^t}.$$   We show in Section 13.2 that, for any irrational number $\alpha$,  one has
$$\mu(\alpha) = \deg \frac{n}{\Vert n\alpha \Vert} = \deg \frac{1}{ \min\left\{ \left| \alpha-\frac{a}{b}\right|: \text{$a,b \in\ZZ$ and $1 \leq b \leq x$} \right\}} < \infty,$$
where $\Vert x \Vert$ for any $x \in \RR$ is the distance from $x$ to the nearest integer.   Although it is known that $\mu(\alpha) = 2$ for all irrationals $\alpha$ outside a set of Lebesgue measure $0$,  and also that $\mu(\alpha) = 2$ for all algebraic numbers $\alpha$ and for $\alpha = e$, the irrationality measure of many important transcendental numbers, like $\pi$ and $\log 2$,  is unknown.    The tightest known bounds for $\mu(\pi)$, for example, are $\mu(\pi) \in [2,7.103205334137\ldots]$ \cite{zei}.  In Section 13.4, we  prove several equivalent characterizations of the logexponential degree $\dege \frac{n}{\Vert n\alpha \Vert} = \dege \frac{1}{ \min\left\{ \left| \alpha-\frac{a}{b}\right|: \text{$a,b \in\ZZ$ and $1 \leq b \leq x$} \right\}}$, which we prove is equal to $(2,1,1,1,\ldots)$ for all irrationals $\alpha$ outside a set of Lebesgue measure $0$.
\item The {\bf radical} of a positive integer $n$, denoted $\operatorname{rad}(n)$,  is the product $\prod_{p |n} p $ of the distinct prime factors of $n$.    The {\bf abc conjecture}\index{abc conjecture} states that, for every $t > 1$, there exist only finitely many triples $(a, b, c)$ of mutually prime positive integers such that  $a + b =  c>\operatorname {rad} (abc)^t$.   It is an extremely important conjecture in Diophantine analysis and has many deep consequences, both known and conjectured. Let $\operatorname{ABC}(n)$ denote the maximal value of $c = a + b$ over all relatively prime positive integers $a, b$ such that $\operatorname{rad}(abc) = n$; that is, let
\[
\operatorname{ABC}(n) = \max\{a + b : a, b \in \mathbb{Z}_{>0},\ \gcd(a, b) = 1,\ \operatorname{rad}(ab(a + b)) = n\}.
\]
It is known that $\operatorname{ABC}(n)  < \infty$ for all $n$.   Moreover, the abc conjecture is equivalent to $\deg \operatorname{ABC} = 1$ and thus also to 
\[
 \dege \operatorname{ABC}  \leq (1, \infty, 1, 0, 0, 0, \ldots).
\]
To date, the tightest known bounds for $\dege  \operatorname{ABC}$ are
\[
(1, 0, \infty, \tfrac{1}{2}, -1, 0, 0, 0, \ldots) \leq \dege \operatorname{ABC}   \leq (\infty, \tfrac{1}{3}, 3, 0, 0, 0, \ldots),
\]
A conjecture of van Frankenhuysen (1995),  supported by Robert,  Stewart, and Tenenbaum (2014),  along with heuristics and numerical evidence,  implies that
\[
 \dege \operatorname{ABC} = (1, 0, \infty, \tfrac{1}{2}, -\tfrac{1}{2}, 0, 0, 0, \ldots).
\]
This leads us to conjecture,  more modestly, that $\deg \text{ABC}(n) = 1$ and $\dege_1  \text{ABC}(n) = 0$,  that is, that
$$\text{ABC}(n) = O(n (\log n)^t) \ (n \to \infty)$$
if (and only if) $t > 0$.  Equivalently,  this conjecture states that
$$\dege  \text{ABC}(n) \leq (1, 0, \infty, 1, 0, 0, 0, \ldots),$$
or, equivalently still, that for every $\varepsilon > 0$,  there exist only finitely many triples $(a,b,c)$ of mutually (or pairwise) relatively prime positive integers with $a+b = c > \operatorname{rad}(abc) \cdot( \log \operatorname{rad}(abc))^\varepsilon$.  
Note that the abc conjecture is also equivalent to $ \deg \frac{1}{ \underline{\operatorname{ABC}}(n)} = -1$, where
\[
\underline{\operatorname{ABC}}(n) =  \min\{ \operatorname{rad}(abc): a, b \in \mathbb{Z}_{>0},\ \gcd(a, b) = 1,\   a+b  = c  = n\}.
\]
\end{enumerate} 
\end{example}

One moral to be drawn from these examples is that, in number theory, an unknown logexponential degree is usually expected to be much closer to its best known lower  bound than its best known upper bound, to several terms.    In other words,  we expect that our current machinery is better at saying how bad our best approximations are than it is in saying how good they truly are; intuitively, this means that we expect that lots of  ``cancellations'' are happening  in our equations and inequalities that we do not yet know how to explain.  

 A research program initiated in Part 3 of this text is to express the logexponential degree of the various real functions arising in number theory in terms of $\dege f$ for as  few ``logexponential primitives'' $f$ as possible.   Additionally, given a real number-theoretic function $f$, one seeks a ``nice'' function $g$,  i.e., a function $g$ within a specified class $\mathcal{F}$ of ``well-behaved'' functions defined in a neighborhood of $\infty$, so that $\dege (f-g) < \dege g$,  which is a stronger condition  than $f(x) \sim g (x) \ (x \to \infty)$.  Ideally,  one would like to find a $g \in \mathcal{F}$ that minimizes $\dege (f-g)$, but, currently,  most such problems of interest are intractable.

One might rightfully demand an explanation as to why logexponential degree is relevant to number theory at all.  
One explanation for this is that the notion of logexponential degree is intimately connected to Hardy's field of all {\bf logarithmico-exponential functions} \cite{har3} \cite{har4},  that is, the field of all germs of real functions that are defined on a neighborhood of $\infty$ and can be can be built from all real constants and the functions $\id$, $\exp$,  and $\log$  using the operations $+$, $\cdot$, $/$, and $\circ$.   Indeed,  several of our results in Chapters 6 and 7 (e.g., Theorem \ref{infpropexp} and Proposition \ref{oexppropstrong}) show that the logexponential degree of a real function captures information about the asymptotic growth rate of the function in comparison to the logarithmico-exponential functions.    Moreover,  as suggested by Hardy \cite[p.\ 22]{har3}, the logarithmico-exponential functions offer precise benchmarks against which to compare the order of growth of nearly any function  in number theory that one might be interested in.   Nevertheless,  how well such functions can approximate the more complicated functions of number theory,  such as $\li-\pi$ and $p_{n+1}-p_n$, remains to be seen.

Broadly speaking,  occurrences of iterated logarithms in analytic number theory asymptotics  can be classified into one of two types:  those that are essential, and those that are  artifacts of our existing machinery for bounding various number-theoretic functions.  Essential occurrences of the two-fold iterated logarithm $\log \log$ abound.  Examples include the {\it law of the iterated logarithm} \cite{feller}, Example \ref{sumex}(3), Remark \ref{prnt}, the asymptotic expansion of $\frac{p_n}{\log n}$  in Example \ref{pnas}, examples (1) and (4) of Example  \ref{exex0}, examples (4), (7), (10), and (11) of Example \ref{exex}, and the definition (\ref{MMCC}) of the Meissel--Mertens constant $M$.  The occurrence of $\log \log$ in the last of these examples is motivated by {\it Cram\'er's model of the primes},  a heuristic model of the primes in which the ``probability'' that an integer $n >1$ is prime is $\frac{1}{\log n}$. Indeed,  Cram\'er's model  suggests that the sum $\sum_{p \leq x} \frac{1}{p}$ can be approximated by $\sum_{1<n \leq x} \frac{1}{n \log n}$, hence  also by $\int_e^x \frac{dt}{t \log t} = \log \log x$.   Cram\'er's model also suggests two  readily verified essential occurrences of $\log \log \log$, namely, that the limit 
\begin{align*}
\lim_{x \to \infty}\left( \sum_{p \leq x} \frac{1}{p \log \log p} - \log \log \log x \right), 
\end{align*}
and thus also the limit
\begin{align*}
\lim_{x \to \infty}\left( \sum_{p \leq x} \frac{1}{p\sum_{ p' \leq p} \frac{1}{p'}} - \log \log \log x \right),\end{align*}
exists.  Likewise, the limits
\begin{align*}
\lim_{x \to \infty}\left( \sum_{p \leq x} \frac{1}{p \log \log p \log \log \log p} - \log \log \log \log x \right)
\end{align*}
and
\begin{align*}
\lim_{x \to \infty}\left( \sum_{p \leq x} \frac{1}{p \left( \sum_{{ p' \leq p}} \frac{1}{p'} \right) \left( \sum_{ p' \leq p} \frac{1}{p'\sum_{ p'' \leq p'} \frac{1}{p''}} \right)} - \log  \log \log \log x \right)
\end{align*}
exist, and so on.  Admittedly, these  last examples are somewhat artificial.   Some natural occurrences of $\log \log \log$, on the other hand, are {\it conjectured} to be essential, e.g., by Montgomery's  conjecture \cite[Conjecture, p.\ 16]{mont1} and by the conjecture \cite[(20)]{ng} of Gonek and Ng.  However, the majority of currently known occurrences of $\log \log \log$, $\log \log \log \log$, and beyond are  unlikely to be essential occurrences (although Example \ref{exex}(11) is another exception).   P.\ Nielsen notes in \cite{nielsen} that the iterated logs in (\ref{fgkmtgap}), for example, come from ``the current state of the art sieve methods, together with bounding techniques.  When you solve for the best fit functions to undo some of the exponentiation that occurs in calculations, the logs just fall out.  In these types of problems, it is not inconceivable (and actually occurs quite regularly) that one new idea is applied to the problem, and the asymptotic changes (sometimes involving more multi-log factors, to account for the small additional room for improvement that was gained).''  

 Although the logexponential degree formalism does not explain {\it why} iterated logarithms arise in analytic number theory,  or why logarithmico-exponential functions are so prevalent, it at least provides a new way of delineating such occurrences and establishing their various interrelationships in order to explain  {\it   how they proliferate}.   Perhaps neither degree nor logexponential degree have been defined or studied before because our ignorance even of $\deg  f$ is so great for so many  number-theoretic functions $f$ that we care about.  Such problems draw attention to the rather humbling fact that our ignorance is much greater than that revealed by problems like the Riemann hypothesis, the Dirichlet divisor problem, and the Gauss circle problem.  One at least can  be hopeful that future solutions to these difficult problems may at the same time shed light on the accompanying logexponential  degress of higher order---as would happen, for example, if something close to Montgomery's conjecture \cite[Conjecture, p.\ 16]{mont1}  or the conjecture \cite[(20)]{ng} of Gonek and Ng were someday  to be resolved in the positive.

 \begin{remark}[The ANTEDB]
 In late January 2025, the author learned of a recent project of Terence Tao, Timothy Trudgian,  and Andrew Yang,  namely,  the Analytic Number Theory Exponent Database, or ANTEDB (currently housed at https://github.com/teorth/expdb),  which they describe as an ``ongoing project to systematically record known theorems for various exponents appearing in analytic number theory, as well as the relationships between them.''  Their notion of exponent is more general, though informal, than our notion of degree. While this book relates asymptotic invariants regardless of tractability, the ANTEDB project 
catalogs and improves the best known explicit bounds on invariants that are amenable to study with our current machinery.
 \end{remark}

\subsection*{Key properties of logexponential degree}

Below we state the key properties of $\dege$ that are proved in Sections 6.3, 7.5, and 7.6 and that  are used throughout Part 3.  Since the proofs in Part 2 are somewhat technical,  the reader may choose to take these properties as given and read Part 3  after just skimming, or even skipping, Part 2.

First, we require some notation.  The set $\prod_{n=0}^\infty \overline{\RR}$  is equipped with the lexicographic order and endowed with the product topology where each factor $ \overline{\RR}$ has the discrete topology.   Let $$\prod_{n = 0}^{ \infty *}\overline{\RR} \subsetneq \prod_{n = 0}^{\infty}\overline{\RR}$$ denote the set of all sequences $\dd = (\dd_i)$ in $\prod_{n = 0}^{\infty}\overline{\RR}$ satisfying the following four conditions for all nonnegative integers $n$.
\begin{enumerate}
\item If $\dd_n = \infty$, then $\dd \geq (\dd_0, \dd_1, \ldots, \dd_n, 0,1,0,0,0,\ldots).$
\item If $\dd_n = -\infty$, then $\dd \leq (\dd_0, \dd_1, \ldots, \dd_n, 0,-1,0,0,0,\ldots).$
\item If $\dd_n$ is finite and $\dd_{n+1} = \infty$, then $\dd \leq (\dd_0, \dd_1, \ldots, \dd_{n+1},1,0,0,0,\ldots).$
\item If $\dd_n$ is finite and $\dd_{n+1} = -\infty$, then $\dd \geq (\dd_0, \dd_1, \ldots, \dd_{n+1}, -1,0,0,0,\ldots).$
\end{enumerate}

For all $\dd, \ee \in \prod_{n = 0}^\infty\overline{\RR}$, we define:
$$\dd\oplus \ee  =  \dd+ \ee\index[symbols]{.g m@$\dd \oplus \ee$}  $$
if  $\dd_k$ and $\ee_k$ are finite for all $k$, and
$$\dd\oplus \ee =  (\dd_0+\ee_0, \ldots, \dd_{n-1}+\ee_{n-1},f_{0},f_{1},f_{2},\ldots)  
$$
if $n$ is the smallest nonnegative integer such that $\dd_n$ and $\ee_n$ are not both finite, where 
$$
(f_0,f_1,f_2,\ldots) = \begin{cases} \max( (\dd_{n},\dd_{n+1},\ldots),(\ee_{n},\ee_{n+1},\ldots)) & \text{if }  \dd_{n} = \infty \text{ or }  \ee_{n} = \infty \\
 \min( (\dd_{n},\dd_{n+1},\ldots),(\ee_{n},\ee_{n+1},\ldots)) & \text{otherwise}
\end{cases}$$
as computed in  $\prod_{n = 0}^\infty\overline{\RR}$.  The operation $\oplus$ is a binary operation on $\prod_{n = 0}^\infty\overline{\RR}$ and restricts to a binary operation on $\prod_{n = 0}^{\infty*}\overline{\RR}$.

We say that $f$ has {\bf exact degree} if $\lim_{x \to \infty} \frac{\log |f(x)|}{\log x}$ exists or is $\pm \infty$.    We say that $f$ has {\bf exact logexponential degree} if  $T^{\circ n}(f)$ has exact degree for all $n$.   If $f$ is not eventually zero, we let $\underline{\dege}\, f = -\dege (1/f)$; otherwise, we let $\underline{\dege}\, f  = (-\infty, -\infty, -\infty,\ldots)$.
 
 Let $f$ and $g$ be real functions with $\dom f \cap \dom g$ containing a subset $X$ of $\RR$ that is not  bounded above.   One has the following.
\begin{enumerate}
\item \textbf{Degree compatibility:}  If $\dege f = \dd$, then $\dd_0 = \deg f$.
\item  \textbf{Exp shifting:} Let $\dd = \dege f$.  If $\deg f = 0$, then $\dege (f \circ \exp) = (\dd_1, \dd_2, \ldots)$.     If $\deg f > 0$ and $\lim_{x \to \infty} f(x) = \infty$, then $\dege (\exp \circ f) = (\infty, \dd_0, \dd_1, \ldots)$.  
\item  \textbf{Log shifting:}  Let $\dd = \dege f$, and let $\log^+ x = \max(\log x,0)$.      If $\deg f  =  \infty$, then $\dege ( \log^+ \circ |f|) = (\dd_1, \dd_2,\dd_3, \ldots)$.   If $\deg f  =  -\infty$, and $\dege ( 1/\log \circ |f|) = (\dd_1, \dd_2,\dd_3, \ldots)$.   If $\deg f \in \RR$, then $\dege (f \circ \log) = (0, \dd_0, \dd_1, \dd_2, \ldots)$.
        \item  \textbf{Exactness equivalence:} $f$ has exact logexponential degree if and only if $\dege f = \underline{\dege} \,  f$.
    \item \textbf{Exactness of $\mathbb{L}$:}  One has $T(f) \in \mathbb{L}$ for all $f \in \mathbb{L}$, and every $f \in \mathbb{L}$ has exact degree;  thus, every $f \in \mathbb{L}$ has exact logexponential degree.
    \item \textbf{$O$ compatibility:} If $f(x) = O(g(x)) \ (x \to \infty)$, then   $\dege f \leq \dege g.$
    \item \textbf{$o$ compatibility:} Suppose that $g$ is eventually defined on $\dom f$.  If $\dege f < \dege g$ and $g$ has exact logexponential degree, then  $f(x) = o(g(x)) \ (x \to \infty)$.  More generally, if $\dege f < \underline{\dege}\, g$, then $f(x) = o(g(x)) \ (x \to \infty)$. 
    \item \textbf{Nonarchimedean property:} One has $\dege(f + g) \leq \max(\dege f, \dege g)$.
      \item \textbf{Submultiplicativity:}  One has $\dege(f g) \leq \dege f \oplus \dege g$.
            \item \textbf{Exact multiplicativity:} If $g$ has exact logexponential degree, and if exactly one of $\dege_n f$ and $\dege_n g$ is finite for the least $n$, if any, for which at least of them is finite, then $\dege(f g) = \dege f \oplus \dege g$.  
    \item \textbf{Composition:} Suppose that $f$ and $g$ are both defined on a neighborhood of $\infty$ and satisfy the following conditions.
    \begin{enumerate}
        \item $f$ has finite degree $d$.
        \item $g$ has positive degree and is eventually positive.
                \item $g(x) \asymp r(x) \ (x \to \infty)$ for some $r \in \mathfrak{L}$.
        \item Either $f$ has exact logexponential degree or $g$ is eventually continuous and increasing.
    \end{enumerate}
    Then one has
    \[
    \dege(f \circ g) = \dege f \oplus d \cdot \big(\dege g + (-1, 0, 0, 0, \ldots)\big).
    \]
    \item   \textbf{Compositional inversion:}  If $f$ is increasing and unbounded  of finite positive  degree with $\dege f = \dd$ (resp.,  $\underline{\dege}\, f = \dd$), then the inverse function $f^{-1}$ exists, is increasing and unbounded, and has lower logexponential degree $\underline{\dege}(f^{-1}) = \dd'$  (resp.,  logexponential degree $\dege  f^{-1} = \dd'$) given by
\[
\dd' = \left( \tfrac{1}{\dd_0},-\tfrac{\dd_1}{\dd_0},\ -\tfrac{\dd_2}{\dd_0}, -\tfrac{\dd_3}{\dd_0},\ldots \right),
\]
where each coordinate is given as above until the first $k$,  if any, such that $\dd_k = \pm\infty$, after which tail of $\dd'$ is exactly the negated tail of $\dd$, that is,  $\dd_j'= -\dd_j$ for all $j \geq k$.
            \item \textbf{Restriction compatibility:} $\dege f|_X \leq \dege f$.
        \item  \textbf{Arithmetic function extendibility:} If $f$ is an arithmetic function, then $\dege f = \dege f(\lfloor x \rfloor)  = \dege f(\lceil x \rceil)$.
    \item \textbf{Majoration:}  Suppose that $f$ is unbounded, but bounded  on $[N,x] \cap X$ for all $x \geq N$, and let
$$\widetilde{f}(x) = \sup_{t \in [N,x] \cap X} |f(t)|, \quad \forall x \in [N,\infty) \cap X.$$  Then 
 $ \dege \widetilde{f} = \dege f $.
\item  \textbf{Admissibility:} One has $\dege  \RR^{\RR_\infty} = \prod_{n = 0}^{ \infty *}\overline{\RR}$.   In fact, for any sequence $\dd \in \prod_{n= 0}^{\infty}\overline{\RR}$, one has $\dd = \dege f$ for some $f \in \RR^{\RR_\infty}$ if and only if $\dd \in \prod_{n = 0}^{ \infty *}\overline{\RR}$, if and only if $\dd = \dege f$ for some positive,  monotonic,  infinitely differentiable function $f$ on $\RR_{>0}$ of exact logexponential degree.
 \item \textbf{Completeness and denseness:}  The poset $\dege  \RR^{\RR_\infty}$ is order-wise dense and complete, and  $\dege \mathbb{L}_{>0}$ is  both order-wise dense and topologically dense in $\dege  \RR^{\RR_\infty}$.
    \item \textbf{Reducibility to  $\mathbb{L}$:}  One has
\begin{align*}
\dege f & = \inf\{\dege r:r \in \mathbb{L}_{> 0}, \, f(x) = O(r(x)) \ (x \to \infty) \} \\ 
 & = \inf\{\dege r:r \in \mathbb{L}_{> 0}, \, \forall x \gg 0\, |f(x)|\leq r(x) \} \\
& = \inf\{\dege r:r \in \mathbb{L}_{> 0}, \, f(x) = o(r(x)) \ (x \to \infty) \} \\
& = \inf\{\dege r:r \in \mathbb{L}_{> 0}, \, \dege f < \dege r \},
\end{align*}
where the infima (exist and) are computed in  $\prod_{n = 0}^ {\infty *}\overline{\RR}$.
    \item \textbf{Stabilization on $\mathbb{L}$:} If $f \in \mathbb{L}$, then $\dege f$ eventually stabilizes to $(0, 0, \ldots)$, that is, there exists an $N$ such that $\dege_k f = 0$ for all $k \geq N$.
\end{enumerate}

 See Tables \ref{tableofsymbols} and \ref{tableofsymbols2} for a list of some important symbols and terms defined and used in this text.  Tables  \ref{tab1}, \ref{tab1a}, and  \ref{tab1b} provide a list of many of the degree theorems  and conjectures in analytic number theory that are discussed in this book.    Note that, if the third column in a given row is empty, then the given result is unconditional.  Functions in the tables whose degrees are unknown include $\li(x)-\pi(x)$, the Mertens function $M(x)$, the functions $\sum_{n \leq x} d(n) - x \log x - (2 \gamma-1)x$, $\sum_{n \leq x} \mu^2(n) -\frac{6}{\pi^2} x$,  and $\sum_{n \leq x} \phi(n) -\frac{3}{\pi^2}x^2$,  the function $\zeta(\sigma + ix)$ for a fixed $\sigma \in \RR$,  and the prime gap function $g_n = p_{n+1}-p_n$. The logexponential degrees of all other functions appearing in the tables are expressed in terms of the logexponential degrees of those particular functions, along with the Riemann zeta zero ordinate gap function $\gamma_{n+1}-\gamma_n$ and the function $S(T)$.

\newpage

\begin{table}[!htbp]
\tiny  \caption{\centering Definitions of some important symbols and terms}
\begin{tabular}{|l|l|l|} \hline
Symbol & Definition & Name \\ \hline\hline
$\overline{\RR}$ & $\RR \cup \{\infty,-\infty\}$ & set of all extended real numbers \\ \hline
$\RR^{\RR_\infty}$  & the set of all real functions $f$ with $\sup \dom f = \infty$ &   \\ \hline
$\mathbb{L}$  & the field of (germs at $\infty$ of) all real functions that can  & the field of (germs at $\infty$ of) all  \\ 
& be built from all constants and the functions $\id$, $\exp$,    &   logarithmico-exponential functions  \\ 
& and $\log$  using the operationas $+$, $\cdot$, $/$, and $\circ$  & \\ \hline
$F^{\circ k}$ & $k$th iterate of $F$  & \\ \hline
$f(x) = O( g(x)) $ & $\exists M > 0: |f(x)| \leq M|g(x)|$ for all $x$ in the intersection & $f(x)$ is big $O$ of $g(x)$ as $x \to a$ \\  
  $\quad (x \to a)$ & of $\dom f$ with some punctured neighborhood of $a$ &  \\ \hline
  $f(x) \ll g(x) $ & $f(x) = O( g(x))  \ (x \to a)$  &  \\  
  $\quad (x \to a)$  & & \\ \hline
$f(x) \gg g(x) $ & $\exists M > 0: |f(x)| \geq M|g(x)|$ for all $x$ in the intersection &  \\  
  $\quad (x \to a)$ & of $\dom f$ with some punctured neighborhood of $a$ &  \\ \hline
      $f(x) \asymp g(x) $ & $f(x) \ll g(x) \ (x \to a)$ and  $f(x) \gg g(x) \ (x \to a)$  &  \\  
  $\quad (x \to a)$  & & \\ \hline
 $f(x) = o( g(x))$ & $\forall M > 0: |f(x)| \leq M|g(x)|$ for all $x$ in the intersection & $f(x)$ is little $o$ of $g(x)$ as $x \to a$ \  \\  
   $\quad (x \to a)$ & of $\dom f$ with some punctured neighborhood of $a$ &  \\ \hline
   $f(x) \sim g(x) $  &  $f(x)-g(x) = o(g(x)) \ (x \to a)$  & $f(x)$ is asymptotic to $g(x)$ as $x \to a$ \\ 
   $\quad (x \to a)$  &   & \\ \hline
   $f(x) = \Omega_{+}(g(x))$  & $\exists M > 0:$ for all $x$ in a punctured neighborhood of $a$  &  $f(x)$ is Omega plus of $g(x)$ as $x \to a$ \\
 $\quad (x \to a)$ &   $\exists y \neq a$ closer to $a$ than $x$ such that $f(y) > M| g(y)|$ & \\ \hline
    $f(x) = \Omega_{-}(g(x))$  & $\exists M > 0:$ for all $x$ in a punctured neighborhood of $a$  &  $f(x)$ is Omega minus of $g(x)$ as $x \to a$ \\
     $\quad (x \to a)$ &   $\exists y \neq a$ closer to $a$ than $x$ such that $f(y) < -M| g(y)|$ & \\ \hline
    $f(x) = \Omega_{\pm}(g(x)) $ &  $f(x) = \Omega_{+}(g(x)) \ (x \to a)$ and   & 
    $f(x)$ is Omega plus minus of $g(x)$  \\
 $\quad (x \to a)$ & $f(x) = \Omega_{-}(g(x)) \ (x \to a)$  & as $x \to a$ \\ \hline
$\deg f$  & $\inf\{t \in \RR: f(x) = O(x^t) \ (x \to \infty)\}$ & degree of $f$    \\ 
  & $= \limsup_{x \to \infty} \frac{\log |f(x)|}{\log x}$  &  \\ \hline
$\underline{\deg} \, f$ & $\sup\{t \in \RR: f(x) \gg x^t \ (x \to \infty)\}$ &  lower degree of $f$  \\ 
  & $= \liminf_{x \to \infty} \frac{\log |f(x)|}{\log x}$  &  \\ \hline 
$\overline{\underline{\deg}} \, f$ & $\lim_{x \to \infty} \frac{\log |f(x)|}{\log x}$ &  exact degree of $f$  \\ \hline
$T(f)$ & $\left.
 \begin{cases}
    f(e^x) e^{-(\deg f) x}& \text{if } \deg f  \neq \pm \infty \\
    \max( \log |f(x)|, 0)  & \text{if } \deg f = \infty \\
 \displaystyle   -\frac{1}{\log |f(x)|} & \text{if } \deg f =- \infty
 \end{cases}
\right.$ & \\ \hline
$\dege f$ & $(\deg f, \deg T(f), \deg T(T(f)), \deg T(T(T(f))), \ldots)$ & logexponential degree of $f$ \\ \hline
$\dege_k f$ & $\deg T^{\circ k}(f)$ = $k$th coordinate of $\dege f$ & logexponential degree of $f$   \\ 
&  & of order $k$ \\ \hline
$f_{(k)}$  & $T^{\circ k}(f)$  & \\ \hline
$f(x) \simeq \sum_{n=1}^{\infty} a_n \varphi_{n}(x)$  &  $f(x) = \sum_{k=1}^{n} a_k \varphi_{k}(x) + o(\varphi_{n}(x)) \ (x \to a)$  for all  & $f(x)$ has the asymptotic expansion  \\ 
$\quad (x \to a)$ & $\quad n \in \ZZ_{>0}$ & $ \sum_{n=1}^{\infty} a_n \varphi_{n}(x)$ at $a$ with  respect  \\ 
& & to $\{\varphi_n\}$ \\ \hline
$S_f(x)$ & $\sum_{n \leq x} f(n)$ & summatory function of $f$  \\ \hline
$D_f(X)$ & $\sum_{n = 1}^\infty f(n)n^{-X}$ & formal Dirichlet series of $f$  \\ \hline
$D_f(s)$ & $\sum_{n = 1}^\infty \frac{f(n)}{n^s}$ & Dirichlet series of $f$  \\ \hline
$(f* g)(n)$ &  $\sum_{ab = n} f(a)g(b)$  & Dirichlet convolution of $f$ and $g$ \\ \hline 
$\zeta(s)$ & meromorphic continuation of $\sum_{n = 1}^\infty \frac{1}{n^s}$ to $\CC$ & Riemann zeta function \\ \hline
$P(s)$ & $\sum_{p} \frac{1}{p^s}$  & prime zeta function \\ \hline
$ \Gamma (s)$  & meromorphic continuation of $\int_{0}^{\infty }x^{s}e^{-x}\,\frac{dx}{x}$  to $\CC$ & gamma function  \\ \hline
$\Psi(s)$  &   $\frac{\Gamma'(s)}{\Gamma(s)}$ & digamma function \\ \hline 
$\Ei(s)$ &  $ \int_{-\infty}^s \frac{e^{z}}{z} dz$ & exponential integral function \\  \hline
$\li(x)$ & $ \int_0^x \frac{dt}{\log t} = \Ei(\log x)$ & logarithmic integral function \\ \hline
$\ERi(s)$ &  $ \sum_{n=1}^\infty \frac{ \mu(n)}{n} \Ei\left(\frac{s}{n} \right)$ &  \\  \hline
$\Ri(x)$ &   $\sum_{n=1}^\infty \frac{ \mu(n)}{n} \li(x^{1/n}) = \ERi(\log x)$ & Riemann's function \\ \hline
$W(x)$  & inverse of the restriction of $xe^x$ to $[-1,\infty)$  & Lambert $W$ function \\ \hline
$\pi(x)$ & $ \sum_{p\leq x} 1  = \#\{p \leq x: p \text{ is prime}\}$ &  prime counting function   \\ \hline
$\PP(x)$ & $\frac{\pi(x)}{x}$ &   prime density function \\ \hline
$\vartheta(x)$ & $\sum_{p \leq x} \log p$ &  first Chebyshev function  \\ \hline
 $\psi(x)$ &  $\sum_{k = 1}^\infty \sum_{p^k \leq x} \log p =  \sum_{k = 1}^\infty \vartheta(x^{1/k})$  &   second Chebyshev function  \\ \hline
 $p_n$ & $n$th prime number ($= \inf\{x \in \RR: \pi(x) \geq n\}$) & prime listing function \\ \hline
$g_n$ & $p_{n+1}-p_n$ & $n$th prime gap \\ \hline
$G(x)$ &  $\max_{p_{k} \leq x}  g_k = \max_{k \leq \pi(x)}  g_k$ & maximal prime gap function \\ \hline
\end{tabular}\label{tableofsymbols}
\end{table}

\newpage

\begin{table}[!htbp]
\tiny \caption{\centering Definitions of some important symbols and terms}
\begin{tabular}{|l|l|l|} \hline
Symbol & Definition & Name \\ \hline\hline
$\gamma_n$ & ordinate of  the $n$th nontrivial zero of  $\zeta(s)$ &  \\ \hline
$\tau_n$ & $\frac{\gamma_n}{2\pi}$ &  \\ \hline
$\widehat{\gamma}_n$  & $\tau_n \log \frac{\tau_n}{e}+ \frac{11}{8}$ & \\ \hline
$N(T)$ &  $\#\{\rho \in \CC: \zeta(\rho) = 0,\, \operatorname{Im} \rho \in (0,T]\}$ & \\ \hline
$S(T)$ &  $\frac{1}{\pi}\operatorname{arg}\zeta\left(\frac{1}{2}+i T\right)$, where the argument is chosen to be  &  \\ 
 & $0$  at $\infty+iT$   and to vary continuously on the line    &  \\ 
  & from $\infty+iT$ to $1/2+iT$  &  \\ \hline
  $v_p(n)$ & exponent of $p$ in the prime factorization of $n$  &   $p$-adic valuation \\ \hline
  $\Omega(n)$   & $\sum_{p|n} v_p(n)$ &  prime Omega function   \\ \hline
 $\omega(n)$  & $ \sum_{p| n} 1$ &  prime omega function  \\ \hline
  $\mu(n)$ & $\begin{cases} (-1)^{\Omega(n)} = (-1)^{\omega(n)} & \quad  \text{if }  n \text{ is squarefree} \\
 0  & \quad \text{if } n \text{ is not squarefree}  
\end{cases}$  &  M\"obius function  \\ \hline
$\lambda(n)$ & $(-1)^{\Omega(n)}$ &   Liouville function  \\ \hline
$\phi(n)$ & $\# (\ZZ/n\ZZ)^* = \#\{k \in \ZZ_{>0} : k \leq n, \, \gcd(k,n) = 1\}$  & Euler's totient \\ \hline
     $\sigma_a(n)$ & $\sum_{d \mid n} d^a$ &   generalized sum of divisors function \\ \hline
 $d(n)$ & $\sigma_0(n) = \#\{d \in \ZZ_{>0}: d \mid n\}$ &  divisor function  \\ \hline
  $\sigma(n)$ & $\sigma_1(n) = \sum_{d \mid n} d$ &   sum of divisors function \\ \hline
$H_n$ & $\sum_{k = 1}^n \frac{1}{n}$ & $n$th harmonic number  \\ \hline
$M(x)$ & $ \sum_{n\leq x}\mu(n)$ & Mertens function \\ \hline
$L(x)$ & $ \sum_{n\leq x}\lambda(n)$ & summatory Liouville function \\ \hline
$\Theta$   & $\sup\{\operatorname{Re} s: s \in \CC,\, \zeta(s)= 0\} = \deg (\li-\pi)$ &  Riemann constant   \\ \hline
$\Theta_k$ & $\dege_k(\li-\pi)$ & $k$th Riemann constant \\ \hline
$\Theta_{-1}$ &  $ \sup\left\{t \in \RR: \li(x)-\pi(x)  \ll xe^{-(\log x)^{t}} \ (x \to \infty) \right\}$ & anti-Riemann constant \\  \hline
$\gamma$ & $\lim_{x \to \infty} \left( \sum_{n \leq x} \frac{1}{n}-\log x \right ) = \lim_{s \to 1} \left( \zeta(s)-\frac{1}{s-1}\right)$ &  Euler--Mascheroni constant  \\ \hline
$M$ &   $\lim_{x \to \infty} \left( \sum_{p \leq x} \frac{1}{p}-\log \log x \right)$ &  Meissel--Mertens constant  \\ \hline
$H$ &    $- \sum_p \left(\frac{1}{p} +\log\left(1- \frac{1}{p} \right) \right)  = \int_1^\infty  \frac{\li(t)-\pi(t)}{t^{2}} \, dt = \gamma-M$ &   Mertens constant \\ \hline
$B$  &  $\lim_{x \to \infty} \left( \log x - \sum_{p \leq x} \frac{\log p}{p} \right)$   &    \\ \hline
$\mu$ & unique positive zero of $\li(x)$  & Ramanujan--Soldner constant  \\ \hline
$\lfloor \alpha \rfloor$ & $\max\{n \in \ZZ: n \leq \alpha\}$ & floor of $\alpha$  \\ \hline
$\lceil \alpha \rceil$ & $\min\{n \in \ZZ: n \geq \alpha\}$ & ceiling of $\alpha$  \\ \hline
$\{\alpha\}$ & $\alpha-\lfloor \alpha \rfloor $  & fractional part of $\alpha$ \\   \hline
$\Vert \alpha \Vert$ & $\min\{|\alpha-n|: n \in \ZZ\} = \min(\{\alpha\},1-\{\alpha\})$ &  distance from $\alpha$ to the nearest integer \\ \hline
   $S(\alpha)$ &   $ \left.
 \begin{cases}
   \frac{1}{\alpha-\lfloor \alpha \rfloor} & \text{if } \alpha \notin \ZZ\cup \{ \infty\}   \\
   \infty  & \text{if } \alpha \in \ZZ\cup \{ \infty\}.\index[symbols]{.u  d@$S(\alpha)$}
 \end{cases}
\right.$  &  regular continued fraction operator \\ \hline
    $a_n(\alpha)$  &  $S^{\circ n}(\alpha)$ & $n$th term of regular continued fraction of $\alpha$ \\ \hline
    $p_n(\alpha)$, $q_n(\alpha)$  &  $\frac{p_n(\alpha)}{q_n(\alpha)} =  [a_0(\alpha), a_1(\alpha),\ldots,a_n(\alpha)]$,  $\gcd(p_n(\alpha),q_n(\alpha)) = 1$, &  \\
       &  $\quad q_n(\alpha) > 0$ &   \\  \hline
$\ord_1 \alpha$ &  $ \#((1,\alpha)\ZZ/\ZZ)$  & order of $\alpha$ modulo 1 \\ \hline
$| \alpha|_1$ &  $\frac{1}{\ord_1 \alpha} = \inf \, (1,\alpha)\ZZ \cap \RR_{>0}$  &   \\ \hline
$\alpha \gg_1 f$ &  $\left|\alpha -r \right| < f(\ord_1 r)$ for only finitely many $r \in \QQ$ & \\ \hline
$\mu(\alpha)$ &  $\inf\{t \in \RR: \alpha \gg_1 n^{-t}\}$ & irrationality measure of $\alpha$ \\ \hline
$M(\alpha)$ & $ \inf \left\{ t \in \RR_{>0}: \alpha  \gg_1 \frac{1}{t} n^{-2}\right\}$  & Markov constant of $\alpha$ \\ \hline
$m(\alpha)$ & $ \inf \left\{ t \in \RR_{>0}: \alpha  \gg_1 \frac{1}{t} n^{-\mu(\alpha)}\right\}$  & relativized Markov constant of $\alpha$ \\ \hline
$\lambda(\alpha)$ &  $\lim_{n \to \infty} \frac{1}{n} \log q_n(\alpha)$ & L\'evy constant of $\alpha$ \\  \hline
$\overline{\lambda}(\alpha;\beta)$ &  $\limsup_{n \to \infty} \frac{\log q_n(\alpha)}{\log q_n(\beta)}$ &  upper relative L\'evy constant of $\alpha$  \\
 & & with respect to $\beta$ \\  \hline
 $\underline{\lambda}(\alpha;\beta)$ &  $\liminf_{n \to \infty} \frac{\log q_n(\alpha)}{\log q_n(\beta)}$ &  lower relative L\'evy constant of $\alpha$  \\
  & &with respect to $\beta$ \\  \hline
  $Q\{\alpha_n\}_n$ &  $\sup\{ c \in [1,\infty): |\alpha-\alpha_{n+1}| = O\left( | \alpha-\alpha_n|^c \right) \ (n \to \infty) \}$ & $Q$-order of  convergence of $\{\alpha_n\}_n$ \\ \hline
$Q(\alpha)$  & $Q\left \{ \frac{p_n(\alpha)}{q_n(\alpha)}\right \}_n$  & \\ \hline 
$R(\alpha)$ &   $\limsup_{n \to \infty } \frac{ \left|\alpha-\frac{p_{n}(\alpha)} {q_{n}(\alpha)}\right|}{ \left|\alpha-\frac{p_{n+1}(\alpha)} {q_{n+1}(\alpha)}\right|} $ & \\ \hline
 $B_k(\alpha)$  &  $ \limsup_{n \to \infty} \frac{q_{n+k}(\alpha)}{q_{n}(\alpha)} $ & \\ \hline
$\pmb{\mu}( \alpha)$ &  $\inf\{\dege f: f \in \mathbb{L}_{>0} \text{ and } \alpha \gg_1  1/f|_{\ZZ_{>0}}\}$  
&  logexponential irrationality degree of $\alpha$ \\ \hline
\end{tabular}\label{tableofsymbols2}
\end{table}

\newpage

\begin{table}[htb!]
\caption{\centering Degree theorems and conjectures in analytic number theory}\tiny 
\centering 
\begin{tabular}{|l|l|l|} \hline 
Degree  & Equals   & Assuming  \\  \hline\hline
$\deg(\li-\pi)$ & $\Theta = \sup_{\rho: \, \zeta(\rho) = 0} \operatorname{Re} \rho \in [\frac{1}{2},1]$ & \\ \hline
$\deg(\li-\pi)$ & $\frac{1}{2}$  & RH (Riemann hypothesis) \\ \hline  
$\deg(\li-\pi)$ & $1$  & ARH (Anti-Riemann hypothesis) \\ \hline  
$\dege_1(\li-\pi)$ & $\leq 1$  & \\ \hline 
$\dege_1(\li-\pi)$ & $ \in [-1,1]$  & RH\\ \hline 
$\dege_1(\li-\pi)$ & $\in [-\infty,-1]$  & $\lnot$RH\\ \hline 
$\dege_1(\li-\pi)$ & $-\infty$  &  ARH \\ \hline    
$\dege_2(\li-\pi)$ & $ \in [-1,-\frac{3}{5}]$  &  ARH \\ \hline    
$\dege(\li-\pi)$ & $\geq (\frac{1}{2},-1,0,1,0,0,0,\ldots)$ & \\ \hline 
$\dege(\li-\pi)$ & $\leq (\Theta,1,0,0,0,0,\ldots)$ & \\ \hline 
$\dege(\li-\pi)$ & $\leq (1,-\infty,-\tfrac{3}{5},\tfrac{1}{5},0,0,0,\ldots)$ & \\ \hline 
$\dege(\li-\pi)$  &  $\leq (\Theta,-1,0,0,0,\ldots)$ & $\lnot$RH\\ \hline
$\dege(\li-\pi)$  & $(\Theta, -1, 0, 0, 0,\ldots)$  &   $\lnot$RH, and $\exists \rho: \zeta(\rho) = 0, \, \operatorname{Re} \rho = \Theta$\\ \hline 
$\dege(\li-\pi)$ & $\geq (1,-\infty,-1,0,0,0,\ldots)$  &  ARH  \\ \hline
$\dege(\li-\pi)$ & $(\frac{1}{2}, -1, \Theta_2, \Theta_3, \ldots)$, $\Theta_2 \geq 0$ &  Conjecture \ref{eurekaconjecture} \\ \hline 
$\dege(\li-\pi)$  & $(\frac{1}{2}, -1, 0, \Theta_3, \Theta_4, \ldots)$,  $\Theta_3 \geq 1$ &  Conjecture \ref{eurekaconjecture2}  \\ \hline 
$\dege(\id -\vartheta)$   &  $\dege(\li-\pi)+(0,1,0,0,\ldots)$  &   \\ \hline  
$\dege(\id-\psi)$  & $\dege(\id-\vartheta)$  &  \\ \hline 
$\dege(\sum_\rho \frac{x^\rho}{\rho})$ & $\dege(\id-\psi)$  & \\ \hline 
$\dege(\psi -\vartheta)$   &  $(\frac{1}{2},0,0,0,\ldots)$  &   \\ \hline  
$\dege(\li-\Ri)$  &  $(\frac{1}{2},-1,0,0,0,\ldots)$  &  \\ \hline  
$\dege(\Ri-\pi)$  &  $\dege(\li-\pi)$  &  \\ \hline  
$\dege(\li -\Pi)$   &  $\dege(\li-\pi)$  &   \\ \hline  
$\dege(\Ri-\Pi)$  &  $\dege(\li-\pi)$  &   \\ \hline  
$\dege(\id -\psi)$   &  $(\frac{1}{2},0,0,2,0,0,0,\ldots)$  &  (\ref{MMC3}) (Montgomery)  \\ \hline  
$\dege(\li - \pi) $  & $(\frac{1}{2},-1,0,2,0,0,0,\ldots)$  &    (\ref{MMC3}) (Montgomery)  \\ \hline  
$\dege(\li - \pi) $  & $(\frac{1}{2},-1,0,1,0,0,0,\ldots)$  &    (\ref{dem}) (Stoll and Demichel) \\ \hline  
$\deg M$ & $\deg(\li-\pi)$  & \\ \hline 
$\dege_1 M$ & $\in [m-1,\infty]$  &  RH, and $\zeta(s)$ has a zero of order $m$ \\ \hline    
$\dege M$ & $\geq (\frac{1}{2},0,0,0,\ldots)$  &  \\ \hline    
$\dege M$ & $\leq (1,-\infty,-\tfrac{3}{5},\tfrac{1}{5},0,0,0,\ldots)$  &  \\ \hline      
$\dege M$ & $\geq (1,-\infty,-1,0,0,0,\ldots)$  &  ARH \\ \hline    
$\dege M$ & $(\frac{1}{2},0,0,\frac{5}{4},0,0,0,\ldots)$  & (\ref{gonek}) (Gonek and Ng) \\ \hline 
$\dege M$ & $(\frac{1}{2},0,\frac{1}{2},0,0,0,\ldots)$  & (\ref{good}) (Good, Churchhouse, L\'evy)  \\  \hline
$\dege M$ & $\geq (\frac{1}{2},0,0,\frac{1}{2},0,0,\ldots)$  & (\ref{kotc}) (Kotnik and van de Lune)  \\   \hline
$\dege M$ & $(\frac{1}{2},0,d_2,d_3,\ldots)$  &  Conjecture \ref{Mconjecture} \\ \hline    
$\displaystyle \dege\left(\sum_\rho \frac{x^\rho}{\rho \zeta'(\rho)}\right)$  & $\displaystyle \dege M$  & RH, and the zeros of $\zeta(s)$ are simple \\ \hline 
$\displaystyle \dege\left(\pi(e^{t} n) -  \sum_{\mu e^{-t}  \leq k < n} \frac{e^{t} }{H_k -\gamma+t} \right)$,  $t \in \RR$  & $\dege(\li-\pi)$ &   \\  \hline  
$\displaystyle \dege  \left(\sum_{p \leq x} \frac{1}{p}  -  \log \log x-M\right)$  & $\displaystyle \dege(\li-\pi) + (-1,0,0,0,\ldots)$ & \\ \hline
$\displaystyle \dege \left( \sum_{p \leq x} \log\left(1- \frac{1}{p}\right) +\log\log x +\gamma\right)$ & 
$\displaystyle \dege(\li-\pi) + (-1,0,0,0,\ldots)$ & \\ \hline
$\displaystyle \dege\left(\frac{1}{s} \sum_{p \leq x} \log\left(1- \frac{s}{p}\right) +\log\log x +G(s) \right)$ &  $\dege(\li-\pi)+(-1,0,0,0,\ldots)$  &   \\ \hline
$\displaystyle \dege \left(e^{\gamma} \prod_{p \leq x}\left(1-\frac{1}{p}\right) - \frac{1 }{\log x}\right)$ & 
$\displaystyle \dege(\li-\pi) + (-1,-1,0,0,0,\ldots)$ & \\ \hline
$\displaystyle \dege \left(e^{-\gamma} \prod_{p \leq x}\left(1-\frac{1}{p}\right)^{-1} -  \log x\right)$ & 
$\displaystyle \dege(\li-\pi) + (-1,1,0,0,0,\ldots)$ & \\ \hline
\end{tabular}\label{tab1}
\end{table}

\newpage

\begin{table}[htb!]
\caption{\centering Degree theorems and conjectures in analytic number theory}\tiny 
\centering 
\begin{tabular}{|l|l|l|} \hline 
Degree  & Equals   & Assuming  \\  \hline\hline
$\displaystyle \dege\left(\log x - \sum_{p \leq x} \frac{\log p}{p}-B \right)$ &  $\dege(\li-\pi)+(-1,1,0,0,0,\ldots) $  &   \\ \hline
$\displaystyle \dege\left(e^{-M}\prod_{p\leq x} e^{1/p}-\log x  \right)$ &  $\dege(\li-\pi)+(-1,1,0,0,0,\ldots) $  &   \\ \hline
$\displaystyle\dege \left(\sum_{n \leq x} d(n) - x \log x - (2 \gamma-1)x\right)$ & $\geq (\tfrac{1}{4},0,0,0,\ldots)$ & \\ \hline
$\displaystyle\dege \left(\sum_{n \leq x} d(n) - x \log x - (2 \gamma-1)x\right)$ & $\leq (\tfrac{131}{416},\tfrac{26947}{8320},0,0,0,\ldots)$ & \\ \hline
$\displaystyle \deg \left(\sum_{n \leq x} d(n) - x \log x - (2 \gamma-1)x\right)$ & $\frac{1}{4}$ & widely conjectured \\ \hline
$\displaystyle \deg\left(\sum_{n \leq x} \mu^2(n) -\frac{6}{\pi^2} x\right)$ & $\in [\frac{1}{4}, \frac{1}{2}]$ & \\ \hline
$\displaystyle \deg\left(\sum_{n \leq x} \mu^2(n) -\frac{6}{\pi^2} x\right)$  & $\in [\frac{1}{4}, \frac{11}{35}]$ & RH \\ \hline
$\displaystyle \dege\left(\sum_{n \leq x} \phi(n) -\frac{3}{\pi^2} x^2\right)$  & $\geq (1,0,\tfrac{1}{2},0,0,0,\ldots)$  &  \\ \hline
$\displaystyle \dege\left(\sum_{n \leq x} \phi(n) -\frac{3}{\pi^2} x^2\right)$  & $\leq (1,\tfrac{2}{3}, \tfrac{1}{3}, 0,0,0,\ldots)$ &  \\ \hline
$\displaystyle \dege\left(\sum_{n \leq x} \phi(n) -\frac{3}{\pi^2} x^2\right)$  & $(1,0, 1, 0,0,0,\ldots)$ &  (\ref{montR1})  and (\ref{montR2}) (Montgomery) \\ \hline
$\displaystyle \dege\left(\sum_{n \leq x}\frac{\phi(n)}{n} -\frac{6}{\pi^2} x\right)$  & $\displaystyle\dege\left(\sum_{n \leq x} \phi(n) -\frac{3}{\pi^2} x^2\right)+(-1,0,0,0,\ldots)$ &  \\ \hline
$\displaystyle \deg \zeta(\sigma + ix), \ \sigma \in \RR$ & $0$ if $\sigma  \geq \frac{1}{2}$, and $\frac{1}{2}-\sigma$ if $\sigma \leq \frac{1}{2}$  & Lindel\"of hypothesis \\  \hline
  $\deg \zeta(\frac{1}{2} + it)$   &  $\leq \frac{13}{84}$   & \\ \hline
    $\deg \zeta(\frac{1}{2} + it)$   &  $0$   & Lindel\"of hypothesis  \\ \hline
$\dege \zeta(\frac{1}{2} + it)$   &  $\geq (0,\infty, \tfrac{1}{2}, -\tfrac{1}{2},\tfrac{1}{2},0,0,0,\ldots)$   & \\ \hline
 $\dege \zeta(\tfrac{1}{2}+it)$   & $\leq (0,\infty, 1, -1,0,0,0,\ldots)$    & RH  \\ \hline
  $\dege \zeta(\tfrac{1}{2}+it)$   & $(0,\infty, \tfrac{1}{2}, \tfrac{1}{2},0,0,0,\ldots)$    &   (\ref{fghcb}) (Farmer, Gonek, Hughes)  \\ \hline
  $\dege \zeta(1 + it)$ & $\geq (0,0,1,0,0,0,\ldots)$   & \\  \hline
$\dege \zeta(1 + it)$ & $(0,0,1,0,0,0,\ldots)$   & RH \\ \hline
$\dege S$  &  $\leq (0,1,0,0,0,\ldots) $ & \\  \hline
$\dege S$  &  $\geq (0,\tfrac{1}{3},-\tfrac{1}{3},0,0,0,\ldots) $ & \\  \hline
$\dege S$  &  $\leq (0,1,-1,0,0,\ldots) $ & RH \\  \hline
$\dege S$  &  $\geq (0,\tfrac{1}{2},-\tfrac{1}{2}, \tfrac{1}{2},0,0,0,\ldots) $ & RH \\  \hline
$\dege S$  &  $(0,\tfrac{1}{2},\tfrac{1}{2}, 0,0,0,\ldots) $ &   (\ref{fghc}) (Farmer, Gonek, Hughes) \\  \hline
$\dege S(\gamma_n)$ & $\dege S$  &  \\ \hline
$\dege (n-\widehat{\gamma}_n)$ & $\dege S$  &  \\ \hline
$\dege (2\pi r_n-\gamma_n)$ &  $\dege S + (0,-1,0,0,0,\ldots)$ &  \\ \hline
$\dege(\gamma_{n+1}-\gamma_n)$ & $\geq (0,-1,0,0,0,\ldots)$ & \\ \hline 
$\dege(\gamma_{n+1}-\gamma_n)$ & $\leq (0,0,0-1,0,0,0,\ldots)$ & \\ \hline 
$\dege(\gamma_{n+1}-\gamma_n)$ & $\leq (0,0,-1,0,0,0,\ldots)$ & RH \\ \hline 
$\dege(\gamma_{n+1}-\gamma_n)$ & $(0,-\frac{1}{2},0,0,0,\ldots)$ & (\ref{arousbour}) (Arous and Bourgade) \\ \hline
$\dege_1(\gamma_{n+1}-\gamma_n)$ & $>-1$ & Conjecture \ref{gapc} \\ \hline
$\dege \delta_n$  &    $\dege (\gamma_{n+1}-\gamma_n) + (0,1,0,0,0,\ldots)$  & \\ \hline
$\dege (\widehat{\gamma}_{n+1}-\widehat{\gamma}_n)$  &    $\dege (\gamma_{n+1}-\gamma_n) + (0,1,0,0,0,\ldots)$  &  \\ \hline
$\dege(p_n-\li^{-1}n)$ & $\dege(\li-\pi) + (0,\Theta+1,0,0,0,\ldots)$ &   \\  \hline
$\dege(p_n-\Ri^{-1}n)$ & $\dege(\li-\pi) + (0,\Theta+1,0,0,0,\ldots)$ &   \\  \hline
\end{tabular}\label{tab1a}
\end{table}

\newpage

\begin{table}[htb!]
\caption{\centering Degree theorems and conjectures in analytic number theory}\tiny 
\centering 
\begin{tabular}{|l|l|l|} \hline 
Degree  & Equals   & Assuming  \\  \hline\hline
$\dege(\li(p_n)-n)$  & $\dege(\li-\pi)+(0,\Theta,0,0,0,\ldots)$ &  \\ \hline
$\dege(\Ri(p_n)-n)$  & $\dege(\li-\pi)+(0,\Theta,0,0,0,\ldots)$ &  \\ \hline
$\dege(x-p_{\li(x)})$  & $\dege(\li-\pi) + (0,1,0,0,0,\ldots)$ &  \\ \hline
$\dege(x-p_{\Ri(x)})$  & $\dege(\li-\pi) + (0,1,0,0,0,\ldots)$ &  \\ \hline
$\deg g_n$ & $\leq \frac{1}{2}$   & density hypothesis \\  \hline
$\deg g_n$ & $0$ &   Conjecture \ref{gapconj} (Piltz) \\  \hline 
$\dege g_n$ & $\leq (\tfrac{21}{40}, \tfrac{21}{40},0,0,0,\ldots)$ &   \\  \hline
$\dege g_n$ & $\geq (0, 1,1,-1,1,0,0,0,\ldots)$ &   \\  \hline
$\dege g_n$ & $\leq (\frac{1}{2}, \frac{3}{2},0,0,0,\ldots)$ &  RH  \\  \hline
$\dege g_n$ & $\leq (0,2,\infty,1,0,0,0,\ldots)$ &  (\ref{pintzcon}) (Pintz) \\  \hline
 $\dege G$ &  $\dege g_n+(0,-\deg g_n,0,0,0,\ldots)$ & \\ \hline
  $\deg G(e^x)$  & 2 & Conjecture \ref{mayncram} (Maynard)   \\  \hline 
  $\dege (x-p_{\pi(x)})$ &  $\dege G$ & \\ \hline
  $\deg \frac{n}{\Vert n \alpha \Vert}$  & $\mu(\alpha) $ &  \\ \hline
$\dege \frac{n}{\Vert n \alpha \Vert}$ & $\pmb{\mu}(\alpha) $ &  \\ \hline
 $\displaystyle \deg  \frac{1}{\min \left\{\left| \alpha-\frac{a}{b}\right|: \text{$a,b \in\ZZ$ and $1 \leq b \leq x$}\right\}}$ & $\mu(\alpha) $ &  \\ \hline
$\displaystyle \dege  \frac{1}{\min \left\{\left| \alpha-\frac{a}{b}\right|: \text{$a,b \in\ZZ$ and $1 \leq b \leq x$}\right\}}$ & $\pmb{\mu}(\alpha) $ &  \\ \hline
    $-\sup\{\deg f: f: \ZZ_{>0}  \longrightarrow \RR_{>0} \text{ and } \alpha \gg_1  f\}$   & $\mu(\alpha) $ &  \\ \hline
    $-\sup\{\dege f: f: \ZZ_{>0}  \longrightarrow \RR_{>0} \text{ and } \alpha \gg_1  f\}$   & $\pmb{\mu}(\alpha) $ &  \\ \hline
$\pmb{\mu}(\alpha)$ & $ = (2,1,1,1,\ldots)$ for almost all $\alpha \in \RR$ &  \\ \hline
$\pmb{\mu}(\alpha)$ & $ \geq (2,0,0,0,\ldots)$ for all irrational $\alpha$ &   \\  \hline
$\pmb{\mu}(\alpha)$ & $ = (2,0,0,0,\ldots)$ for all real algebraic &   \\ 
 & $\quad$  numbers $\alpha$ of degree $2$ &  \\ \hline
$\pmb{\mu}(e)$ & $(2,1,-1,0,0,0,\ldots)$ & \\ \hline
$\pmb{\mu}(\pi)$, $\pmb{\mu}(\log 2)$, $\pmb{\mu}(\gamma)$ & $(2,1,1,1,\ldots)$ &  Conjecture \ref{diophconj0} \\ \hline
$\pmb{\mu}(\alpha)$ & $ \leq (2,1,1,1,\ldots)$ for all real algebraic & Conjecture \ref{diophconja}  \\  & $\quad$  numbers $\alpha$  &  \\ \hline
$\pmb{\mu}(\alpha)$ & $ = (2,1,1,1,\ldots)$ for all real algebraic & Conjecture \ref{diophconj3}  \\ 
 & $\quad$  numbers $\alpha$ of degree $> 2$ & \\ \hline
 $\underline{\deg} \, \min\{\operatorname{rad}(ab(a + b)) : a, b \in \ZZ_{> 0},  $ & $1$ &   abc conjecture  \\ 
  $\quad$ $\gcd(a, b) = 1,  \ a+b = n\}$ & & \\ \hline
         $\underline{\dege} \, \min\{\operatorname{rad}(ab(a + b)) : a, b \in \ZZ_{> 0},  $ &  $  \leq   (1, 0, -\infty, -\tfrac{1}{2}, 1, 0, 0, 0, \ldots)$  &     \\ 
  $\quad$ $ \gcd(a, b) = 1, \ a+b = n\}$ &  \\ \hline
     $\underline{\dege} \,\min\{\operatorname{rad}(ab(a + b)) : a, b \in \ZZ_{> 0},  $ & $  \geq (-\infty, -\tfrac{1}{3}, -3, 0, 0, 0, \ldots)$ &     \\ 
  $\quad$ $ \gcd(a, b) = 1, \ a+b = n\}$ && \\ \hline
         $\underline{\dege} \, \min\{\operatorname{rad}(ab(a + b)) : a, b \in \ZZ_{> 0}, $ &  $ (1, 0, -\infty, -\tfrac{1}{2}, \tfrac{1}{2}, 0, 0, 0, \ldots)$  &    1995 conjecture of \\ 
  $\quad$ $\gcd(a, b) = 1,  \ a+b = n\}$ & &van Frankenhuysen    \\ \hline
\end{tabular}\label{tab1b}
\end{table}

\section{Notation and conventions}

In this section, we discuss some notation and conventions that are assumed throughout the text.  An index of symbols and an index of terms appear at the end of the book.

\subsection*{Numbers}

As is standard, we let $\ZZ$, $\QQ$, $\RR$, and $\CC$\index[symbols]{.a  A@$\ZZ$}\index[symbols]{.a  B@$\QQ$}\index[symbols]{.a  C@$\RR$}\index[symbols]{.a  D@$\CC$} denote the ring of all integers, the field of all rational numbers, the field of all real numbers, and the field of all complex numbers, respectively.  We let $\ZZ_{\geq 0}$ denote the set of all nonnegative integers, and we let $\ZZ_{> 0}$ denote the set of all positive integers.   We also use self-explanatory notations like $\RR_{\geq 1}$, $\QQ_{< 0}$, and so on.   We also let $\pp$ denote the set of all prime numbers.\index[symbols]{.a  P@$\pp$}

For any $x \in \RR$, a {\bf neighborhood of $x$}\index{neighborhood} is a subset of $\RR$ containing an interval of the form $(a,b)$ for some $a,b \in \RR$ with $a<x<b$, while a  {\bf punctured neighborhood of $x$}\index{punctured neighborhood} is a subset of $\RR$ containing an interval of the form $(a,b)\backslash \{x\} = (a,x)\cup (x,b)$ for some $a,b \in \RR$ with $a<x<b$.  A {\bf (punctured) neighborhood of $\infty$}\index{neighborhood of $\infty$} is a subset of $\RR$ containing  the interval $(a, \infty)$  for some $a \in \RR$, while a {\bf (punctured) neighborhood of $-\infty$} is a subset of $\RR$ containing the interval $(-\infty,a)$  for some $a \in \RR$.  The structure $$\overline{\RR} = \RR \cup \{\infty,-\infty\} = [-\infty,\infty]$$ of all {\bf extended real numbers}\index{extended real number}\index[symbols]{.a  ZZ@$\overline{\RR}$}\index[symbols]{.a  ZZ@$\overline{\RR}$} ordered by $\leq$ is complete as a totally ordered set.  Endowed with the order topology, it is a compact  topological space, isomorphic to $[0,1] \subsetneq \RR$.  It is ``almost'' a field, with exceptions made for  the usual indeterminate forms, e.g., $\infty+(-\infty)$, $0\cdot \infty$, and $\frac{\infty}{\infty}$, which are undefined.   
If $X$ is a subset of $\RR$ or  $\overline{\RR}$, then we write $\overline{X}$\index[symbols]{.d NP@$\overline{X}$} for the closure of $X$ in $\overline{\RR}$.  Thus, for example, one has $\infty \in \overline{X}$ if and only $X$ contains an unbounded set of positive real numbers (i.e., $\sup X = \infty$) or $\infty \in X$.  If $X \subseteq \RR$ is nonempty, then $\overline{X}$ is the union of the closure of $X$ in $\RR$ with $\{\sup X, \inf X\}$.

We  use the expression ``for all $x \gg 0$'' to mean ``for all sufficiently large $x$,'' i.e., ``there exists an $N > 0$ such that for all $x \geq N$.''\index[symbols]{.f a@$\gg 0$}  Also, we say that a given property expressed in terms of $x$ holds ``eventually'' if it holds for all $x \gg 0$.  However, in some contexts these notions are also used relativized to the domain of any of the functions involved.  

\subsection*{Sets and functions}

The cardinality of a set $X$ is denoted $\# X$.\index[symbols]{.b a@$\#X$}  For any sets $X$ and $Y$,  we let $Y^X$\index[symbols]{.c edd@$Y^X$}  denote the set of all functions from $X$ to $Y$.   The identity function on a set $X$ is denoted $\id = \id_X$.\index[symbols]{.b a@$\id$} We write $\dom f$, $\operatorname{codom} f$ and $\operatorname{im} f$\index[symbols]{.b ca@$\dom f$}\index[symbols]{.c  da@$\operatorname{codom} f$}\index[symbols]{.c ea@$\operatorname{im} f$} for the domain, codomain, and range of a function $f$.   The expression ``$f$ is defined on $X$'' means that $X$ is a subset of the domain of $f$.    If $f$ is defined on $X$, then we write $f\vert_X$  for the restriction $f\vert_X: X \longrightarrow \operatorname{codom}  f$\index[symbols]{.c ed@$f\vert_X$}   of $f$ to the set $X$.   We sometimes abuse notation by writing $f = g$ when $f$ and $g$ are functions with $\dom g \supseteq \dom f$ and $f(x) = g(x)$ for all $x \in \dom f$.   If $f$ and $g$ are functions,  then we assume that the composition $f \circ g$ has domain $\{x \in \dom g: g(x) \in \dom f\} = g^{-1}(\dom f \cap \operatorname{codom} g)  = g^{-1} (\dom f \cap \operatorname{im} g)$ and codomain $\operatorname{codom} f$.
For all nonnegative integers $k$, we denote the $k$th iterate $f \circ f \circ \cdots \circ f$ of a function $f$ by $f^{\circ k}$.\index[symbols]{.d  fa@$f^{\circ k}$} Thus, for example, the $k$th iterate $\log^{\circ k}$ of $\log$ has domain $(\exp^{\circ (k-1)}(0),\infty)$ for all positive integers $k$.  If the compositional inverse $f^{-1}$ of $f$ exists, then we write $f^{\circ (-k)} = (f^{-1})^{\circ k}$.   
 If $*$ is an associative binary operation, written multiplicatively, on a set $S$, then  for all $x \in S$ we denote by $x^{* n}$ the $n$-fold $*$-product $x * x* \cdots * x$.

 A {\bf real function}\index{real function} is a function with domain a subset of $\RR$ and codomain a subset of $\RR$;  a {\bf complex function}\index{complex function} is a function with domain a subset of $\CC$ and codomain a subset of $\CC$; and a {\bf complex-valued function of a real variable}\index{complex-valued function of a real variable}  is a function with domain a subset of $\RR$ and codomain a subset of $\CC$.  For any such functions $f$ and $g$, unless otherwise stated, we assume that the functions $f+g$, $f-g$, and $f \cdot g$ have domain $\dom f  \cap \dom g$,   and the function $f/g$  has domain $\dom f \cap \{x \in  \dom g: g(x) \neq 0\} = (\dom f \cap  \dom g)\backslash g^{-1}(\{0\})$.

For functions $f$ of a real variable, we often require the assumption $\infty \in \overline{\dom f}$, or equivalently, $\sup \dom f = \infty$, i.e., the domain of $f$ is a subset of $\RR$ that is not bounded above.   We also let $\RR^{\RR_\infty}$  denote the set of all real functions $f$ with $\infty \in \overline{\dom f}$.  More generally, for any $a \in \overline{\RR}$, we let $\RR^{\RR_a}$\index[symbols]{.a rrr@$\RR^{\RR_a}$} denote the set of all real functions $f$ such that $a$ is a limit point of $\dom f$.

 We consider all limits, limits superior, and limits inferior to be relativized to the domain of the function involved.  Thus, for example, for all $f \in \RR^{\RR_\infty}$ we write $\lim_{x \to \infty} f(x) = L$ if for every $\varepsilon > 0$ there exists an $N > 0$ such that $|f(x)-L|< \varepsilon$ for all  $x \in \dom f$ with $x> N$.   
Specifically for the purpose of studying degree, we extend the definitions of $\lim_{x \to a} f(x)$, $\limsup_{x \to a} f(x)$, and $\liminf_{x \to a} f(x)$ to functions $f$ of a real variable  having  values in $\overline{\RR}$.  The extended definitions for ``extended real valued functions'' still require that $a \in \overline{\dom f}$, but where $\dom f$ may include values of $x$ such that $f(x) = \pm \infty$.  Thus, for example, if $f(x) = \infty$ for some unbounded set of values of $x > 0$, then $\limsup_{x \to \infty} f(x) = \infty$.  

We also set the conventions $\log |0|  = \log(0^+) = -\infty$,  $\log \infty = \infty$, and $e^{-\infty} = 0$.  Thus, if $f$ is a real function, then $\log |f(x)|$ is a function of a real variable, with the same domain as $f$, that assumes values in the extended reals, where $\log |f(x)|$ is finite if and only if $f(x) \neq 0$.

For all $x \in \RR$, the {\bf floor of $x$},\index{floor function $\lfloor x \rfloor$}\index[symbols]{.d M@$\lfloor x \rfloor$}  denoted $\lfloor x \rfloor$, is the largest integer less than or equal to $x$,   the {\bf ceiling of $x$},\index{ceiling $\lceil x \rceil$}\index[symbols]{.d N@$\lceil x \rceil$} denoted $\lceil x \rceil$, is the  smallest integer greater than or equal to $x$ (and is equal to $-\lfloor -x \rfloor$),  and the {\bf fractional part of $x$},\index{fractional part $\lceil x \rceil$}\index[symbols]{.d N@$\{ x\}$} denoted $\{ x \}$, is equal to $x - \lfloor x \rfloor$.

\subsection*{Other conventions}

All rings are assumed commutative with identity.  If $R$ is a ring, then $R^+$ denotes the group $R$ under addition, while $R^\times$ denotes the monoid $R$ under multiplication and $R^*$ denotes the group of units of $R$ under multiplication.

 We use the bold letters, like $\mathbf{d}$ and  $\mathbf{e}$, to denote tuples with coordinates $\mathbf{d}_n$ and $\mathbf{e}_n$, respectively,  indexed by the nonnegative integers.
 
In all finite or infinite sums and products with dummy variable ``$p$,'' the values of $p$ are restricted to the primes, taken in succession.  Thus, for example,   $\sum_{p} \frac{1}{p^2}$ represents the sum $\frac{1}{2^2} + \frac{1}{3^2} + \frac{1}{5^2}+ \frac{1}{7^2} +\frac{1}{11^2}+\cdots$.

A term (or phrase) that is being defined  precisely is written in boldface.   Technical terms for which we do not provide the definition are written in italics, and in such cases their definition can be found in Wikipedia.    A term is also written in italics if it is defined precisely later in the text (in which case it will appear in the index), or if it is being singled out for emphasis. 

Ultimately,  theorems, propositions, lemmas, and corollaries are all theorems, and their labelling as such is a matter of taste.    A  ``Remark'' is a relevant but parenthetical comment that can be skipped over without loss of continuity.    A ``Problem'' is a newly posed problem that the author was unable to solve, while an ``(Outstanding) Problem'' is also a Problem, but its solution necessitates the solution of some longstanding open problem, such  as the Riemann hypothesis.

\section{Acknowledgments} 

I would like to thank my family, especially my mother, Nancy (Elliott) Cappello, for all of their support through some difficult times.   This book would not have been possible without them.   I would also like to thank my former teachers who helped inspire my passion for mathematics,  writing,  and teaching, most especially, John Foley and Jean Jonker of Holyoke High School,  and my former MIT professors James Munkres,  Michael Artin,  Haynes Miller,  Gian-Carlo Rota,  and George Boolos.

I would like to thank the late Doug Stoll for assisting me with some of the data analysis and graphs concerning the function $V(x)$ provided in Sections 1.3 and 13.2, for correcting some of my miscalculations in Mathematica regarding the function $G(s)$  in Section 8.5, and for carrying out several time-consuming computations.   I would  also like to thank Grant Molnar for reading and providing corrections and extensive comments on an early draft of the book.   Finally,  I would like to thank the following former students of mine for collecting and analyzing data and for creating some of the graphs in Section 5.3 and Chapter 14, as a part of semester-long undergraduate research classes in Spring 2022: Noah Browne and Paul Kime at CSU Channel Islands, and Alexa Alcala, Aditya Baireddy, Sudhanva Kulkarni,  Lucas Salim, Christopher Silbermann,  and Haolin Zhang at UC Berkeley.   

All graphs in this book were created using Mathematica.

\section{About the author}

Jesse Elliott is a professor of mathematics and philosophy at California State University, Channel Islands.  He was born and raised in Massachusetts and received a BS in Mathematics in 1995 from the Massachusetts Institute of Technology.   In 2003, he completed a PhD in Mathematics from the University of California, Berkeley,  under the advisorship of Hendrik W.\ Lenstra, Jr.    
His areas of research are algebra, number theory, and the foundations and philosophy of mathematics.   His first book,  {\it Rings, Modules, and Closure Operations},  published by Springer in 2019,  is a research monograph on the applications of closure operations to the study of rings and modules.   

\section{Publication}

This book is scheduled to be published as Volume 13 of World Scientific's Monographs in Number Theory Series.  Due to a severe illness, I am not certain that I will be able to complete the lengthy publication process.  World Scientific has  graciously agreed to allow me to post this unofficial draft to the arXiv.   Once the process  (hopefully) is complete, the book will be available at  the link below.

\bigskip

https://www.worldscientific.com/worldscibooks/10.1142/13521\#/t=aboutBook.

\mainmatter

\part{A survey of analytic number theory}

\chapter{A brief history of primes}

In this chapter, we provide a brief history of the  mathematical study of the prime numbers.   Nearly all of the results described in this chapter are well known  by analytic number theorists and are stated without proof.  In subsequent chapters, some of these results are proved.  

\section{The prime numbers,  algebraically}

A positive integer is said to be {\bf prime}\index{prime} if it is greater than $1$ and has two and only two positive divisors, namely, $1$ and itself, and a positive integer is said to be {\bf composite}\index{composite} if it is greater than $1$ but not prime and thus has a positive divisor other than $1$ and itself.  Prime numbers have been the subject of intense study for at least two millenia.  In Book 7 of Euclid's {\it Elements}, for example, one finds an ingenious proof of the fundamental theorem of arithmetic, and another of the infinitude of the  set of all primes.   The {\bf fundamental theorem of arithmetic}\index{fundamental theorem of arithmetic} is the statement that every positive integer $N$ can be written as a product of primes (an empty product in the case $N = 1$, and a product of one prime in the case where $N$ is prime), and, moreover, the prime factorization of $N$ is unique if the prime factors are listed in nondecreasing order.  The prime numbers, therefore, are like ``multiplicative atoms,'' and positive integers are  like ``molecules'' that can be built up uniquely from such atoms.  In other words, the prime numbers are the ``multiplicative building blocks'' of the positive integers.   The list of positive integers in terms of their prime factorizations proceeds as follows:
$$1, \  2, \  3, \  2\cdot 2, \  5, \  2\cdot 3, \  7, \  2\cdot 2 \cdot 2, \  3 \cdot 3, \  2 \cdot 5, \  11, \  2 \cdot 2 \cdot 3, \  13, \  2 \cdot 7, \  3 \cdot 5, \ 2\cdot 2 \cdot 2 \cdot 2, \ 17,$$
$$2 \cdot 3 \cdot 3, \ 19, \ 2 \cdot 2 \cdot 5, \ 3 \cdot 7, \ 2 \cdot 11, \ 23, \ 2 \cdot 2 \cdot 2 \cdot 3, \ 5 \cdot 5, \ 2 \cdot 13, \ 3 \cdot 3 \cdot 3, \ 2 \cdot 2 \cdot 7,  \ 29, \  2 \cdot 3 \cdot 5,   \ldots.$$
Analogously, the number $1$ is the one and only ``additive atom'' of the nonnegative integers, and the list of the nonnegative integers in terms of their ``additive building blocks'' proceeds as follows:
$$0, \ 1, \ 1+1, \ 1+1+1, \ 1+1+1+1, \ 1+1+1+1+1, \ldots.$$

The 25 prime numbers less than 100 are
$$2, \  3, \  5, \  7, \  11, \  13, \  17, \  19, \  23, \  29, \  31, \  37, \  41, \  43, \  47, \  53, \  59, \  61, \  67, \  71, \  73, \  79, \  83, \  89, \  97.$$
The number $1$ is neither prime nor composite: it is not considered to be prime, because if it were prime then the number $2$ could be factored as a product $2 \cdot 1 \cdot 1$ of three primes, or four primes, etc., which would invalidate the fundamental theorem of arithmetic.   Of course, $2$ is the only even prime, just as $3$ is the only prime divisible by $3$, and $5$ is the only prime divisible by $5$, and so on.  It follows that $2$ and $3$ are the only consecutive integers both of which are prime.  Since among any three consecutive odd integers, exactly one is divisible by $3$,  the numbers $3$, $5$, and $7$ are the only three consecutive odd integers all three of which are prime.  

See {\it The PrimePages: prime number research \& records} (https://t5k.org) for the latest records concerning the prime numbers.  A  {\bf Mersenne}\index{Mersenne number}  number is a number of the form $2^k-1$ for some positive integer $k$, or, equivalently, a number whose binary expansion is a sequence of all $1$s.  
Currently, the largest known and verified prime number, found in 2018, is the Mersenne prime $2^{82589933}-1$, which has $24862048$ digits, or, in binary, $82589933$ bits.    Note that, if the Mersenne number $2^k-1$ is prime, then $k$ must be prime, since $2^a-1$ divides $2^{ab}-1$ for all $a,b \in \ZZ_{>0}$.  Note also that $82589933$ is the 
$4811740$th prime.  Fast algorithms are known for testing Mersenne numbers for primality, so much so that the size of known Mersenne primes vastly outstrips the size of known non-Mersenne primes.  Currently there are 51 known Mersenne primes, the smallest $18$  of which are $2^p-1$ for $p$ equal to
$$2, \  3, \  5, \  7, \  13, \  17, \  19, \  31, \  61, \  89, \  107, \  127, \  521, \  607, \  1279, \  2203, \  2281, \  3217.$$
The 17 largest known Mersenne primes were all discovered  by the Great Internet Mersenne Prime Search (GIMPS), a distributed computing projec founded in 1996.  The first and smallest of these 17 Mersenne primes, namely, $2^{1398269}-1$, was discovered on November 13, 1996.   To this day,  no one has yet proved or disproved   {\bf Mersenne's conjecture}\index{Mersenne's conjecture}  that there are infinitely many Mersenne primes.   Nevertheless,  the consensus among number theorists is that Mersenne's conjecture is very likely to be true.

One reason that the prime numbers are so fascinating is the fact that, not only are there infinitely many of them, but  various statements as to how they are distributed are extremely beautiful and profound, and yet, more often than not, are extremely difficult to prove.  There are many easily formulated conjectures about the prime numbers, such as Mersenne's conjecture, that have been open for centuries.   The centuries-old {\bf twin prime conjecture}\index{twin prime conjecture} is the statement that there are infinitely many twin primes,  where {\bf twin primes}\index{twin primes} are pairs of consecutive odd integers both of which are prime.  Examples of small twin primes are $11$ and $13$, and $101$ and $103$.  Currently, the largest known twin primes, discovered in 2016, are $$2996863034895 \cdot 2^{1290000} \pm 1,$$ which have $388342$ digits.  Very little progress was made on the twin prime conjecture until J.\ R.\ Chen proved in 1966 that there are infinitely many prime numbers $p$ such that $p+2$ is either prime or semiprime, where a {\bf semiprime}\index{semiprime} is a positive integer that is a product of two primes.  More recently, in 2013, the groundbreaking work of  Y.\ Zhang, along with subsequent work by J.\ Maynard, T.\ Tao and others, resulted in a proof that there are infinitely many primes that differ at most by 246, or, in analytic terms,
$$2 \leq \liminf_{n \to \infty} (p_{n+1}-p_n) \leq 246,$$
where $p_n$ denotes the $n$th prime.  The twin prime conjecture, on the other hand, is equivalent to the  statement that
$$\liminf_{n \to \infty} (p_{n+1}-p_n) = 2.$$

Not only is the twin prime conjecture widely believed true, but it is believed also that there are  infinitely many {\bf prime triplets},\index{prime triplet} that is, triples of primes of the form $(p,p+2,p+6)$ or $(p,p+4,p+6)$.  For example, the  six prime triplets of the first form are
$$(5,7,11), \, (11,13,17), \, (17,19,23),\, (41,43,47), \, (101, 103, 107), \, (107, 109, 113),$$
and the  first six prime triplets of the second form are
$$(7, 11, 13),\, (13, 17, 19),\,  (37, 41, 43),\, (67, 71, 73),\, (97, 101, 103),\, (103, 107, 109).$$
As of October 2020, the largest known proven prime triplet contains primes with 20008 digits, namely the primes $(p,p+2,p+6)$ with $$p  = 411128692197\cdot 2^{66420}-1.$$

One can generalize the notions of twin primes and prime triplets as follows. If a given $k$-tuple has infinitely many translations whose coordinates are all prime,  then there cannot exist a prime $p$ such that the set of all of the $k$-tuple's coordinates contains every possible residue modulo $p$: if such a prime $p$ were to exist, then, for any positive integer $n$, one of the coordinates formed by adding $n$ to the $k$-tuple would be divisible by $p$, so there could only be finitely many translations whose coordinates are all prime, namely, only those that include $p$.  A $k$-tuple that satisfies this condition is said to be {\bf admissible}.\index{admissible $k$-tuple}  The {\bf prime $k$-tuples conjecture}\index{prime $k$-tuples conjecture}   states that, for any admissible $k$-tuple $(m_1,m_2,\ldots,m_k)$ (e.g., $(0,2)$, $(0,2,6)$, $(0,4,6)$, and $(0,2,6,8)$), there are infinitely many positive integers  $n$ such that the $k$-tuple $(n+m_1, n+m_2,\ldots,n+m_k)$ is a $k$-tuple of primes, or {\bf prime $k$-tuple}.\index{prime $k$-tuple}

Besides the celebrated {\it G\"odel's incompleteness theorems}, there are  many reasons why some questions about the integers are so easy to formulate and yet so difficult to settle.  One such reason can be explicated as follows. Note first that the concept of primality, for example, is a purely  {\it multiplicative} concept in that it is defined solely in terms of multiplication, i.e., in terms of the monoid $\ZZ_{> 0}$ under multiplication.  By contrast, an {\it additive} concept is one that is defined solely in terms of addition, i.e., in terms of the monoid $\ZZ_{\geq 0}$ under addition.  Of course, some concepts, like being even, are both additive and multiplicative, and multiplication itself can be defined recursively in terms of addition via the distributive law.  Nevertheless, this distinction persists, since the monoids $\ZZ_{>0}$ and $\ZZ_{\geq 0}$  are not even close to being isomorphic to one another,  despite the fact that their respective operations are related via the distributive law.  Notice, then, that the notion of a twin prime combines an additive concept (differing by $\pm 2$) with a multiplicative concept (primality), and the twin prime conjecture thus asks how a particular set of multiplicative concepts relate to a particular set of additive concepts. Other famous examples of this type of problem include Mersenne's conjecture, {\it Goldbach's conjecture}, and the {\it abc conjecture} (all still open), as well as the {\it Taylor--Wiles theorem} (formerly {\it Fermat's last theorem}), {\it Mih\u{a}ilescu's theorem} (formerly {\it Catalan's conjecture}), and the {\it Green--Tao theorem} and its generalization, the {\it Tao--Ziegler theorem}.  Many important questions about the distribution of primes are of this type,  e.g.,  questions about how the prime numbers are distributed additively, or linearly, across the number line.  Thus, our understanding---and lack thereof---of how multiplicative properties of the integers relate to additive properties of the integers goes far beyond the distributive law.   Somewhat ironically,  one has group isomorphisms  $\exp : \RR \longrightarrow \RR_{>0}$ and $\log : \RR_{> 0} \longrightarrow \RR$, where  $\RR$ is a group under addition and $\RR_{>0}$ is a group under multiplication,  so that profound questions of this particular sort  do not arise in the realm of the real numbers.

\begin{remark}[Euclid's proof of the fundamental theorem of arithmetic]
The existence part of the fundamental theorem of arithmetic is easy to establish using the fact that the ordered set $\ZZ_{>0}$ of all positive integers is {\bf well-ordered},\index{well-ordered} that is, every nonempty set of positive integers has a smallest element.  The well-orderedness of the ordered set $\ZZ_{>0}$ is known as the {\bf well-ordering principle},\index{well-ordering principle} and it is equivalent to the principle of induction, given the other axioms of arithmetic.   Suppose, to obtain a contradiction, that not every positive integer is a product of prime numbers.  Then, since  $\ZZ_{>0}$ is well-ordered, there must exist a smallest positive integer $N$ that is not a product of prime numbers.   Since we allow products of zero or one primes,  the integer $N$ must be greater than $1$ and not prime.  It follows that $N$ is composite, and so it must have a positive integer divisor $d$ other than $1$ and $N$, so that $1 < d < N$.  The {\bf complementary divisor}\index{complementary divisor} of a positive divisor $a$ of an integer $n$ is  the integer $n/a$, which is also a positive divisor of $n$, since $n = (n/a)a$.  Clearly the complementary divisor $N/d$ of the divisor $d$ of $N$   also satisfies $1< N/d < N$.  Thus, since $d$ and $N/d$ are both positive integers less than $N$, they must both be equal to a product of prime numbers.  But then $N = (N/d)d$ is a product of two numbers that are expressible as a product of prime numbers.  It follows that $N$ itself is expressible as a product of prime numbers, and this is our desired contradiction.  

Euclid's proof of the uniqueness of prime factorizations rests crucially on  {\bf Euclid's lemma},\index{Euclid's lemma} which states that, if a positive integer $d$ divides the product $ab$ of two integers $a$ and $b$, and if $\gcd(d,a) = 1$, then the integer $d$ must divide $b$.   From Euclid's lemma one can deduce that an integer $p >1$ is prime if and only if, for all integers $a$ and $b$, if $p$ divides the product $ab$, then $p$ must divide $a$ or $p$ must divide $b$.  This property of primes is exactly what is needed to show that if one has two prime factorizations $q_1 q_2\cdots q_m =  N = q_1' q_2' \cdots q_n'$ of the same positive integer $N$, then $q_1$ must equal $q_i'$ for some $i = 1,2,\ldots,n$, whence it can be cancelled from both sides of the equation, and the argument repeated, to deduce that the prime factorizations of a given integer are essentially unique.  Finally, the proof of Euclid's lemma rests on the fact that, for any integers $a$ and $b$, there exist integers $s$ and $t$ such that $\gcd(a,b) = sa+tb$,  which has a simple constructive proof called the {\it extended Euclidean algorithm}. 
\end{remark}

\begin{remark}[The infinitude of the primes]
If $q_1, q_2, \ldots q_n$ is any finite list of primes, not necessarily distinct, then the integer $N = 1+q_1 q_2 \cdots q_n > 1$ leaves a remainder of $1$ when divided by $q_i$ for any $i$, and thus $N$ is not divisble by any of the primes $q_1, q_2, \ldots q_n$.  However, every number greater than $1$ is divisible by some prime, which follows from the existence of prime factorizations.  Thus, $N$ has a prime factor $p$, and any  such prime factor $p$ of $N$ (e.g., the smallest prime factor of $N$) cannot be among the list of primes $q_1, q_2, \ldots q_n$, because none of those primes divide $N$.  This is Euclid's ingenious ``constructive'' proof that, for any finite list of primes, there is a prime that is not on that list.  Such a proof is preferable to a proof by contradiction, which is not constructive, but which is often given in textbooks.   It is also worth noting that,  around Euclid's time, an {\it actual infinity}  had been regarded by many influential thinkers, including Zeno of Elea and Aristotle, with much suspicion, but a {\it potential infinity} had been deemed acceptable, which is probably why Euclid phrased his proposition (Book IX Proposition 20) as ``Prime numbers are more than any assigned multitude of prime numbers''\ (English translation) rather than as the modern paraphrase ``There are infinitely many prime numbers.''  
\end{remark}

In the terminology of abstract algebra, the fundamental theorem of arithmetic states that  the commutative monoid $\ZZ_{> 0}$ of all positive integers under multiplication is the {\it free commutative monoid} generated by the set of all prime numbers.    Similarly, the group $\QQ_{>0}$ of all positive rational numbers under multiplication is the {\it  free abelian group} generated by the set of all prime numbers.  Analogously, the commutative monoid $\ZZ_{\geq 0}$ of all nonnegative integers under addition is the free commutative monoid generated by $1$.   

 In category theory, a {\it universal property} of an object is a property that uniquely characterizes  the object up to isomorphism as an object of some category.   For example,  the ring $\ZZ$ of all integers is  characterized uniquely up to isomorphism as an {\it initial object}  in the category of rings, that is, for any ring $R$ (with identity), there is a unique ring homomorphism from $\ZZ$ to $R$, and any ring $\ZZ'$ with the same property is isomorphic to $\ZZ$ via a unique isomorphism $\ZZ \longrightarrow \ZZ'$.    Similarly, the field $\QQ$ is characterized uniquely up to isomorphism as an initial object  in the category of fields of characteristic $0$,  or as an initial object  in the category of {\it ordered fields}, and the ordered field $\RR$  of all real numbers is characterized uniquely up to isomorphism as a {\it terminal object} in the category of {\it archimedean ordered fields}.  Generally speaking, a universal property of a mathematical object offers  evidence that the given object is worthy of study.   
 
\begin{remark}[Rings, prime ideals, and unique factorization domains]
 Let $R$ be a commutative ring (with identity).  The ring $R$ is an {\bf integral domain}\index{integral domain} if $R$ is a subring of some field, or, equivalently, if $0$ is the one and only zerodivisor of $R$.  An ideal $\ppp$  of $R$ is  {\bf maximal}\index{maximal ideal}  if the quotient ring $R/\ppp$ is a field, and an ideal $\ppp$ is {\bf prime}\index{prime ideal}  if $R/\ppp$ is an integral domain.   Equivalently, a maximal ideal  of $R$ is an ideal that is contained by exactly two ideals, namely, itself and $(1) = R$, while a prime ideal of $R$ is equivalently an ideal $\ppp$ of $R$ properly contained in $(1)$ such that $\ppp  \supseteq  \aaa \bbb$ implies $\ppp \supseteq  \aaa$ or $\ppp \supseteq \bbb$, for all ideals $\aaa$ and $\bbb$ of $R$.   For example, the ideal $(0)$ is prime if and only if $R$ is an integral domain, and the ideal $(0)$ is maximal if and only if $R$ is a field, if and only if $R$ has exactly two ideals, namely, $(1) \supsetneq (0)$.  The maximal ideals of the integral domain $\ZZ$ are precisely the principal ideals $(p)$ generated by $p$ for some prime number $p$, and $\ZZ$ has exactly one other prime ideal, namely, $(0)$.  An element $p$ of $R$ is said to be {\bf prime}\index{prime element of a ring} if the principal ideal $(p)$ of $R$ generated by $p$ is a prime ideal of the ring $R$, or, equivalently,  if the quotient ring $R/(p)$ is an integral domain, or equivalently still, if $p$ is a nonunit of $R$ such that $p$ divides a product in $R$ if and only if it divides at least one of the factors.   For example, the prime elements of the integral domain $\ZZ$ are precisely $0$ and $\pm p$ for $p$ prime.  Prime factorizations of a given non-zerodivisor of $R$, when such a factorization exists, are always unique up to reordering and associates, where two elements of $R$ are said to be {\bf associates}\index{associates} if they generate the same principal ideal.   An element $r$ of an integral domain $R$ is said to be {\bf irreducible} if it is a nonzero non-unit and if $r = ab$ implies either $a$ or $b$ is a unit of $R$, for all $a,b \in R$.   In an arbitrary integral domain,  a nonzero element may not have a factorization as a product of irreducibles, and such factorizations  need not be unique up to  reordering and  associates when they do exist.  For example, the two factorizations $2 \cdot 2 = (2i) \cdot (-2i)$ of $4$ in the integral domain $\ZZ[2i]$ are  distinct irreducible factorizations of $4$, since none of the four numbers can be factored any further (except by using a factor of $\pm 1$, which are the only units in $\ZZ[2i]$), and since $2$ is associate to nether $2i$ nor $-2i$ in the ring $\ZZ[2i]$, since $i \notin \ZZ[2i]$.   A {\bf unique factorization domain}\index{unique factorization domain} is an integral domain $R$ in which every nonzero nonunit element factors (necessarily uniquely up to reordering and associates) as a product of primes.   The fundamental theorem of arithmetic is equivalent to the statement that the ring $\ZZ$ is a unique factorization domain.  
\end{remark}

\section{The prime numbers, asymptotically}

 {\it Analytic number theory} is the branch of number theory that uses methods from real and complex analysis to deepen our understanding of the integers.   Its mere existence is a startling testament to the interconnectedness of mathematics.   Historically, one of the first applications of analysis to number theory was unveiled in 1737, when Euler related the study of prime numbers to what is now known as the {\bf Riemann zeta function}  $$\zeta(x) = \sum_{n = 1}^\infty \frac{1}{n^x},$$ which  Euler studied as a function defined for all real $x > 1$.  Essentially, Euler showed that the fundamental theorem of arithmetic is equivalent to the {\bf Euler product representation}
$$\zeta(x) = \prod_{p \text{ prime }} \frac{1}{1-\frac{1}{p^x}}, \quad \forall x > 1,$$
of $\zeta(x)$.  Thus,  Euler translated a fundamental algebraic property of the monoid $\ZZ_{>0}$ into a particular analytic  property of the function $\zeta: \RR_{> 1} \longrightarrow \RR_{> 1}$.  As a consequence of Euler's result, if there were only finitely many primes, then the product $$\prod_{p \text{ prime }} \frac{1}{1-\frac{1}{p}}  = \zeta(1)$$ would be finite, but that would contradict $$\zeta(x) \to \zeta(1) = \sum_{n = 1}^\infty \frac{1}{n} = \infty \mbox{ \  as \ } x \to 1^+.$$  This provided the first {\it analytic proof} of the infinitude of the primes.

  In the same year, Euler proved that the sum $\sum_{p \text{ prime }} \frac{1}{p}$ of the reciprocals of all prime numbers diverges, which yielded yet another analytic proof of the infinitude of the primes.    A modernization of Euler's argument shows that
\begin{align}\label{ZL}
\lim_{x \longrightarrow 1^+} \frac{P(x)}{\log  \zeta(x)} = 1,
\end{align}
where $$P(x) = \sum_{p \text{ prime}} \frac{1}{p^{x}}$$
for $x> 1$ is the {\bf prime zeta function},\index{prime zeta function $P(s)$}  and where $\log$ denotes the logarithm to the base $e$.  This in particular implies that  $P(1) = \sum_{p \text{ prime }} \frac{1}{p}$ diverges.

\begin{figure}[ht!]
\includegraphics[width=60mm]{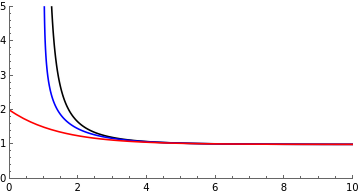}
\caption{\centering Graphs of $\zeta(x)$ (in black), $1+P(x)$ (in blue), and $1+2^{-x}$ (in red) on $[0,10]$}
 \label{zetap1b}
\end{figure}

Figure \ref{zetap1b} provides a graph of the functions $\zeta(x)$, $1+P(x)$, and $1+2^{-x}$ on the interval $[0,10]$.  Notice from the graphs that it appears that one has
$$\lim_{x \to \infty} \zeta(x) = 1$$
and $$\lim_{x \to \infty} P(x) = 0,$$
which are not difficult statements to prove.  One can prove, more generally, that
$$\zeta(x)-1 \sim P(x) \sim 2^{-x}  \ (x \to \infty),$$
where one writes
$$f(x) \sim g(x) \ (x \to a)$$ whenever $f(x)$ and $g(x)$ are functions such that $$\lim_{x \to a} \frac{f(x)}{g(x)} = 1.$$   Equivalently,  the condition $f(x) \sim g(x) \ (x \to a)$ means that $f(x)$ is a ``good approximation'' for $g(x)$ for $x$ near $a$ (or for large $x$ if $a = \infty$) in the sense that the percentage error of this approximation, given as a fraction by the ratio $\frac{f(x)-g(x)}{g(x)}$, tends to $0$ as $x \to a$:
$$\lim_{x \to a} \left(\frac{f(x)-g(x)}{g(x)}\right) = \lim_{x \to a}\left(\frac{f(x)}{g(x)} - 1\right) = 0.$$
Note that the symbols ``$f(x) \sim g(x) \ (x \to a)$'' are read ``$f(x)$ is asymptotic to $g(x)$ as $x$ approaches $a$.''  

The asymptotic notation $\sim$ is used  routinely in analytic number theory.   One also writes
$$f(x) = o(g(x)) \ (x \to a)$$
whenever $$\lim_{x \to a} \frac{f(x)}{g(x)} = 0.$$      This means that $f(x)$ is ``infinitely smaller'' than $g(x)$ for $x$ near $a$ (or for large $x$ if $a = \infty$).
Note, in particular, that $$f(x) \sim g(x) \ (x \to a) \iff f(x)-g(x) = o(g(x)) \ (x \to a).$$  It is also convenient to write
 $$f(x) = O(g(x)) \ (x \to a)$$
whenever $$\limsup_{x \to a} \left| \frac{f(x)}{g(x)}\right| < \infty,$$   
or, equivalently, whenever the function $\left| \frac{f(x)}{g(x)}\right|$ is bounded on some punctured neighborhood of $a$.  
More general definitions of these relations are provided in Section 2.1, but for this chapter the definitions above suffice.  The symbols ``$f(x) = o(g(x)) \ (x \to a)$'' are read ``$f(x)$ is little $o$ of $g(x)$ as $x$ approaches $a$,''  and the symbols ``$f(x) = O(g(x)) \ (x \to a)$'' are read ``$f(x)$ is big $O$ of $g(x)$ as $x$ approaches $a$.''

\begin{figure}[ht!]
\includegraphics[width=80mm]{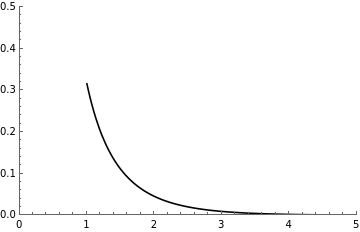}
\caption{\centering Graph of $\log \zeta(x)-P(x)$ on $(1,5]$}
 \label{zetazeta}
\end{figure}

In the notations defined above, (\ref{ZL}) is expressed as
$$P(x) \sim \log \zeta(x)  \ (x \to 1^+).$$
In fact, the approximation $$P(x) \approx \log \zeta(x)$$ is excellent for all $x > 1$, as can be witnessed from the graph in Figure \ref{zetazeta} of the difference $\log \zeta(x)-P(x)$, which decreases rapidly to $0$ on  the interval $(1,\infty)$, and where the limit $$H = \lim_{x \to 1^+}\left( \log \zeta(x)-P(x)\right)  = 0.315718452053\ldots\index{Mertens constant $H$}\index[symbols]{.f tb@$H$}
$$ is an important constant, which we call the {\bf Mertens constant}, discussed further later in this section.

Euler was interested in the function $\zeta(x)$ in part because of the renowned {\bf Basel problem},\index{Basel problem} posed by P.\ Mengoli in 1650, which was the problem of computing the exact value of the convergent series $\sum_{n = 1}^\infty \frac{1}{n^2}$.  In 1734, Euler solved the problem by proving, non-rigorously by today's standards, that
$$\zeta(2) = \sum_{n = 1}^\infty \frac{1}{n^2} = \frac{\pi^2}{6} = 1.644934066848\ldots.$$
From the Euler product representation of $\zeta(x)$, which Euler proved just three years later, it follows that
$$\frac{\pi^2}{6} = \prod_{p \text{ prime }} \frac{p^2}{(p-1)(p+1)} = \frac{2}{1} \cdot \frac{2}{3} \cdot \frac{3}{2} \cdot \frac{3}{4} \cdot\frac{5}{4} \cdot \frac{5}{6} \cdot\frac{7}{6} \cdot \frac{7}{8} \cdot\frac{11}{10} \cdot \frac{11}{12}\cdot \frac{13}{12} \cdot \frac{13}{14}\cdot \frac{17}{16} \cdot \frac{17}{18} \cdots.$$

After Euler's work, another giant leap occurred  in  1837,  with Dirichlet's   proof of {\bf Dirichlet's theorem on primes in arithmetic progression} \cite{diri},\index{Dirichlet's theorem on primes in arithmetic progression} which many say kickstarted {\it rigorous} analytic number theory.  

\begin{theorem}[{Dirichlet's theorem on primes in arithmetic progression \cite{diri}}]
For any positive integers $a$ and $b$, there are infinitely many primes of the form $a+bn$ for $n \in \ZZ_{>0}$ if (and only if) $a$ and $b$ are relatively prime. 
\end{theorem} 

Actually, Dirichlet proved a much stronger theorem: if $a$ and $b$ are relatively prime positive integers, then the proportion of all prime numbers less than or equal to $x$ that are congruent to $a$ modulo $b$ tends to $\frac{1}{\phi(b)}$ as $x$ tends to $\infty$.  Here, $\phi(n)$ denotes {\bf Euler's totient},\index{Euler's totient $\phi(n)$}\index[symbols]{.rt  A@$\phi(n)$} which equals the number of positive integers from $1$ to $n$ that are relatively prime to $n$.  Thus, Dirichlet's theorem says that the primes for any integer $b$  are equally distributed modulo $b$ among the $\phi(b)$ distict congruence classes of integers relatively prime to $b$.  A well-known generalization of Dirichlet's theorem is that the sum 
$$\sum_{p \text{ prime } \atop p \equiv a \, (\text{mod b})} \frac{1}{p}$$ diverges, and in fact the limit
$$\lim_{x \to \infty} \left(\sum_{p \leq x \atop p \equiv a \, (\text{mod b})} \frac{1}{p} - \frac{1}{\phi(b)} \sum_{p \leq x} \frac{1}{p}\right)$$
converges (and both terms in the limit differ from $\frac{1}{\phi(b)}\log \log x$ by a limiting constant), for all relatively prime positive integers $a$ and $b$.

Dirichlet proved his theorem by studying a certain class of {\it Dirichlet series}, where the {\bf Dirichlet series}\index{Dirichlet series $D_f(s)$}\index[symbols]{.l E@$D_f(s)$} of a function $f: \ZZ_{> 0} \longrightarrow \CC$ is the function 
$$D_f(s) = \sum_{n = 1}^\infty \frac{f(n)}{n^s},$$
of a complex variable $s$.   Specifically,  Dirichlet's proof employed the Dirichlet series of functions $f$ called {\it Dirichlet characters}, whose corresponding Dirichlet series are called {\it Dirichlet $L$-functions}.    Although Dirichlet considered Dirichlet series and Dirichlet $L$-functions as functions of a real variable, Riemann was the first to study a particular Dirichlet series as a function of a complex variable,   namely, the series
$$\zeta(s)= \sum_{n = 1}^\infty \frac{1}{n^s},$$
which converges (absolutely) to an analytic function, called the {\bf Riemann zeta function}, on the right half plane $\{s \in \CC: \operatorname{Re} s > 1\}$ \cite{rie}.   Dirichlet's and Riemann's work later inspired Jensen (in 1884 and 1888) and Cahen (in 1894) to undertake the general  study of Dirichlet series of a complex variable.  Today, Dirichlet series are a highly-developed tool for studying functions $f: \ZZ_{> 0} \longrightarrow \CC$ and are  one of the main topics of Chapter 3.

On a more simplistic level, using functions to study the primes is possible for the basic reason that there are numerous functions in number theory that carry all of the information about the primes.    Perhaps the most obvious example is the {\bf prime listing function}\index{prime listing function $p_n$}\index[symbols]{.rt P@$p_n$}  $p_-: \ZZ_{> 0} \longrightarrow \ZZ_{> 0}$, where $p_n$ for all $n$ denotes the $n$th prime, so that, for example, $p_1 = 2$ and $p_{25} = 97$.  Another obvious example is the characteristic function $\chi_{\pp}:  \ZZ_{> 0} \longrightarrow \{0,1\}$ of the set $\pp$ of all prime numbers, where for any subset $X$ of $\ZZ_{> 0}$ the {\bf characteristic function $\chi_X: \ZZ_{> 0} \longrightarrow \{0,1\}$ of $X$\index{characteristic function $\chi_X(n)$}\index[symbols]{.rt Eb@$\chi_X(n)$}} is the function defined by
$$\chi_X(n) = \begin{cases}   1 & \quad \text{if } n \in X \\
  0 & \quad \text{if } n \in \ZZ_{>0}\backslash X. 
\end{cases}$$
Both of the examples above are discrete functions, defined on $\ZZ_{>0}$.    Functions $f: \ZZ_{>0} \longrightarrow \CC$,  as those above, are called {\bf arithmetic functions}.\index{arithmetic function}  Essentially, they are just sequences of complex numbers indexed by the positive integers.  

Besides forming its Dirichlet series $D_f(s)$, another way to ``analysis-ize'' an arithmetic function $f$  is to form its {\bf summatory function}\index{summatory function $S_f(x)$}\index[symbols]{.k Fr@$S_f(x)$}  $S_f$, which is the function $S_f: \RR_{\geq 0} \longrightarrow \CC$ defined by
$$S_f(x) = \sum_{n \leq x} f(n),\quad \forall x \geq 0.$$
The summatory function $S_f$ of $f$  is a step function that is constant on the interval  $[n-1,n)$, and changes by the value $+f(n)$ at $n$, for every positive integer $n$.   In many ways, to be made apparent in Chapter 3, the association $f \longmapsto S_f$ is a discrete analogue of integration.    The summatory function $S_{\chi_{\pp}}$ of the characteristic function $\chi_{\pp}$ of the set $\pp$ of all primes is called the {\bf prime counting function}\index{prime counting function $\pi(x)$} $\pi(x)$.\index[symbols]{.s B@$\pi(x)$}  Specifically,
$$\pi(x) = \sum_{p \leq x \atop p \text{ is prime}} 1 = \#\{p \leq x: p \text{ is prime}\}$$
for any $x \geq 0$ is equal to the number of primes less than or equal to $x$.   For example, one has $\pi(100) = 25$ because there are $25$ primes less than or equal to $100$, the largest being $p_{25} = 97$. The function $\pi(x)$ is a nondecreasing step function, continuous from the right, that jumps up by $1$ at every prime  and is constant, assuming the value $n$, over the interval $[p_n, p_{n+1})$, for all $n$.  Thus, the $n$th prime $p_n$ can be recovered from $\pi(x)$ as the $n$th smallest discontinuity of $\pi(x)$, or as the smallest real number $x$ such that $\pi(x) \geq n$, and one has
$$p_n =  \inf\{x \in \RR: \pi(x) \geq n\} = \min\{x \in \RR: \pi(x) \geq n\}$$
for all positive integers $n$.  Since $\pi(p_n) = n$ for all $n$, the function $\pi(x): \RR_{\geq 0} \longrightarrow \ZZ_{\geq 0}$ is a left inverse to the function $p_n : \ZZ_{\geq 0} \longrightarrow \RR_{\geq 0}$, where one defines $p_0 = 0$.   In the reverse direction, the composition $p_{\pi(x)}$ for any $x$ equals the largest prime number less than or equal to $x$, e.g., $p_{\pi(100)} = 97$.   Generally speaking, information about one of the two functions $p_n$ and $\pi(x)$ yields corresponding information about the other. 

Other ``primes-equivalent'' functions worth mentioning  are the {\bf prime density function}\index{prime density function $\PP(x)$}$$\PP(x) = \frac{\pi(x)}{x}, \quad \forall x > 0,\index[symbols]{.s BA@$\PP(x)$}$$ and the related {\bf prime probability function}\index{prime probability function $\PP(\lfloor x \rfloor)$} $$\PP(\lfloor x \rfloor) = \frac{\pi(x)}{\lfloor x \rfloor}, \quad \forall x \geq 1.$$   For any $x\geq 1$, the number $\PP(\lfloor x \rfloor)$ equals the probability that a randomly selected positive integer less than or equal to $x$ is prime.  Thus, for example, one has
$$\PP(100) =  \frac{25}{100} = \frac{1}{4}.$$

Given that larger primes are rarer than smaller primes, one might expect that
\begin{align}\label{Plim}
\lim_{x \to \infty} \PP(x)  = 0,
\end{align}
or, equivalently, 
$$\pi(x) = o(x) \ (x \to \infty),$$
and indeed this is true.   If $X$ is any set of positive integers, then one defines the {\bf natural density $\delta(X)$ of $X$}\index{natural density $\delta(X)$} to be the limit
$$\delta(X) = \lim_{x \to \infty}  \frac{\#\{n \in X: n \leq x\}} {\#\{n \in \ZZ_{> 0}: n \leq x\}} = \lim_{x \to \infty}  \frac{\#\{n \in X: n \leq x\}} {\lfloor x \rfloor}  = \lim_{x \to \infty}  \frac{\#\{n \in X: n \leq x\}} {x},$$
provided that any of the limits above exist.  Thus, (\ref{Plim}) states that the natural density of the set of all primes is $0$.  (By contrast, if $a$ and $b$ are positive integers, then the natural density of the set of all positive integers congruent to $a$ modulo $b$ is equal to $\frac{1}{b} > 0$.)    Note also that (\ref{Plim})  implies that
$$\lim_{n \to \infty} \frac{n}{p_n} = \lim_{n \to \infty} \frac{\pi(p_n)}{p_n} = 0.$$  Below is a relatively simple proof of (\ref{Plim}). 

\begin{proof}[Proof of (\ref{Plim})]  Let $n$ be a fixed positive integer, and let $x \geq 1$.  Since all prime numbers are either less than or equal to $n$ or relatively prime to $n$, and since there are  at most $\left \lceil \frac{x}{n} \right \rceil \phi(n)$ integers less than or equal to $x$ that are relatively prime to $n$, one has
$$\pi(x) \leq n + \left\lceil  \frac{x}{n} \right \rceil \phi(n),$$
and therefore
$$\frac{\pi(x)}{x} \leq \frac{n}{x} + \frac{1}{x} \left\lceil  \frac{x}{n} \right \rceil \phi(n),$$
whence
$$0 \leq \limsup_{x \to \infty} \frac{\pi(x)}{x} \leq \lim_{x \to \infty}\left( \frac{n}{x} + \frac{1}{x} \left\lceil  \frac{x}{n} \right \rceil \phi(n)\right) = \frac{\phi(n)}{n} = \rho(n),$$
where $$\rho(n) = \frac{\phi(n)}{n} = \prod_{p |n} \left(1-\frac{1}{p} \right)\index[symbols]{.rt  A@$\rho(n)$}$$
 is equal to the probability that a randomly chosen integer from $1$ to $n$ is relatively prime to $n$.   But one has
$$\rho(n!) =  \prod_{p |n!} \left(1-\frac{1}{p} \right) = \prod_{p \leq n} \left(1-\frac{1}{p} \right)$$
and therefore
$$\lim_{n \to \infty} \rho(n!) = \lim_{n \to \infty}\prod_{p \leq n} \left(1-\frac{1}{p} \right) =  \prod_{p } \left(1-\frac{1}{p} \right) = \frac{1}{ \prod_{p } \left(1-\frac{1}{p} \right) ^{-1}} = \frac{1}{\zeta(1)} = \frac{1}{\infty} = 0,$$
where the limits above can be made rigorous.  It follows that
$$0 \leq \limsup_{x \to \infty} \frac{\pi(x)}{x} \leq \lim_{n \to \infty}  \rho(n!) = 0,$$
which immediately implies (\ref{Plim}). 
\end{proof}

As of the writing of this book, the record for computing $\pi(x)$ is for $x = 10^{29}$.  For $x = 10^{29}$, it is known
 (On-Line Encyclopedia of Integer Sequences, February 2022) that the exact value of $\PP(x)$ is  $$\PP(x) = \frac{1520698109714272166094258063}{10^{29}} = \frac{1}{65.759271587961\ldots},$$
which is much larger than one might naively expect: the odds that a randomly chosen positive integer having at most 29 digits is prime are better than $1$ in $66$.   Thus, although the function $\PP(x)$ tends to $0$,  apparently  it does so very slowly.   

In analysis, when a given limit converges, it is natural to ask how fast it converges.  Thus, (\ref{Plim}) raises a more refined question: How quickly does $\PP(x)$ tend to $0$?  Like $1/x$?  Like $1/\sqrt{x}$?   An answer to this  challenging question is provided by the celebrated {\bf prime number theorem}.\index{prime number theorem}

\begin{theorem}[{Prime number theorem \cite{val1} \cite{had}}]
One has
\begin{align*}
\pi(x) \sim \frac{x}{\log x} \ (x \to \infty),
\end{align*}
or, equivalently,
\begin{align*}
\PP(x) \sim \frac{1}{\log x} \ (x \to \infty).
\end{align*}
\end{theorem}

Since $\frac{1}{\log x}$ tends to $0$ as $x \to \infty$, the prime number theorem is a precise mathematical statement to the effect that {\it $\PP(x)$ tends to $0$ like $\frac{1}{\log x}$, as $x$ tends to $\infty$}.     
The prime number theorem is easily seen to be equivalent to each of the following statements.
\begin{enumerate}
\item  $\displaystyle \lim_{x \to \infty} \frac{\pi(x)}{\frac{x}{\log x}}  = 1$.
\item  $\displaystyle  \lim_{x \to \infty} \frac{\PP(x)}{\frac{1}{\log x}}  = 1$.
\item $\displaystyle \pi(x) \sim \frac{x}{\log x -1} \ (x \to \infty)$.
\item $\displaystyle \PP(x) \sim \frac{1}{\log x-1} \ (x \to \infty)$.
\item $\displaystyle  \lim_{x \to \infty} x^{\PP(x)} = e$.
\end{enumerate}
The last of these statements reveals that famous constant $e$, which was discovered by Bernoulli in 1683 and  is given by
$$e = \lim_{x \to \infty} \left(1+\frac{1}{x}\right)^x = \sum_{n = 0}^\infty \frac{1}{n!}  =  2.718281828459\ldots,$$
``encodes'' information about the asymptotic distribution of the primes.

The (informal) approximation $\pi(x) \approx \frac{x}{\log x}$ was first conjectured by Gauss in 1792 or 1793 at the age 15 or 16, according to his own recollection  in his famous letter to the astronomer Encke in 1849  \cite[pp.\ 444--447]{gau}.   The first actual published statement of such an approximation was made by Legendre in 1798, which he refined further in 1808 \cite{gold}.  Legendre's and Gauss' speculations are discussed further in Remark \ref{legendre}.   Chebyshev (if not Gauss) seems to have been the first to have come close to a precise statement of the theorem.   In 1848, he proved a result, namely,  \cite[II-\'eme Th\'eor\`eme]{cheb1}, that immediately implies that
$$\liminf_{x \to \infty}  \frac{\pi( x)}{\frac{x}{\log x}} \leq 1 \leq \limsup_{x \to \infty}  \frac{\pi( x)}{\frac{x}{\log x}},$$
and, therefore,  if the limit $\lim_{x \to \infty} \frac{\pi( x)}{\frac{x}{\log x}}$ exists, then it must equal $1$.   Then,  in 1850, he  proved  \cite{cheb1}
 that
\begin{align}\label{chebin}
c_1\frac{x}{\log x}< \pi(x) < c_2 \frac{x}{\log x}, \quad \forall x \gg 0,
\end{align}
for all $$c_1 < A = \log (2^{1/2}3^{1/3}5^{1/5}30^{-1/30}) = 0.9212920229\ldots$$
and for  $$c_2 = \tfrac{6}{5} A = 1.1055504275\ldots,$$
and thus
$$A \leq \liminf_{x \to \infty}  \frac{\pi( x)}{\frac{x}{\log x}} \leq 1 \leq \limsup_{x \to \infty}  \frac{\pi( x)}{\frac{x}{\log x}} \leq \frac{6}{5}A.$$
Of course, the prime number theorem itself is equivalent to
$$\liminf_{x \to \infty}  \frac{\pi( x)}{\frac{x}{\log x}} = 1 = \limsup_{x \to \infty}  \frac{\pi( x)}{\frac{x}{\log x}}$$
and to
\begin{align*}
C_1\frac{x}{\log x}< \pi(x) < C_2 \frac{x}{\log x}, \quad \forall x \gg 0,
\end{align*}
for all $C_1$ and $C_2$ with $C_1<1<C_2$.    Thus, Chebyshev came very close to proving the prime number theorem, and there is little doubt that he had conjectured the theorem in its precise form.

In 1896, almost 100 years after Legendre's aforementioned published work, the prime number theorem was finally proved,  by de la Vall\'ee Poussin \cite{val1} and Hadamard \cite{had}.   Although it is has been claimed many times that de la Vall\'ee Poussin and Hadamard proved  the theorem ``independently,''  ``separately'' is a more accurate description, as the two proofs were very far from being independent from one another: just one year prior,  in 1895,  von Mangoldt  had rigorously proved \cite{mang0} several unproved assertions that Riemann made in his landmark paper \cite{rie} of 1859, and it was those results that finally made the proofs of the prime number theorem attainable.  Specific details of these developments are discussed  in the next section.   

Today, many ``effective'' versions of the prime number theorem are known.  For example, in  1961, Rosser and Schoenfeld proved \cite{ross} that
\begin{align*}
 \frac{x}{\log x - \frac{1}{2}} < \pi(x) < \frac{x}{ \log x - \frac{3}{2}}, \quad \forall x \geq 67,
\end{align*}
which is easily seen to imply the prime number theorem.    See \cite[p.\ 36]{ing2} and \cite[Theorem 4.5]{apos} for simple proofs that the prime number theorem is equivalent to
$$p_n \sim n \log n \ (n \to \infty),$$
where $p_n$ denotes the $n$th prime.

\begin{figure}[ht!]
\includegraphics[width=80mm]{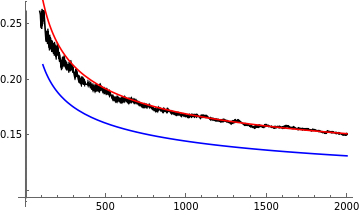}
\caption{\centering Graphs of $\PP(x)$ (in black), $\frac{1}{\log x}$ (in blue), and $\frac{1}{\log x-1}$ (in red) on $[1,2000]$}
 \label{overflow1b}
\end{figure}

\begin{figure}[ht!]
\includegraphics[width=80mm]{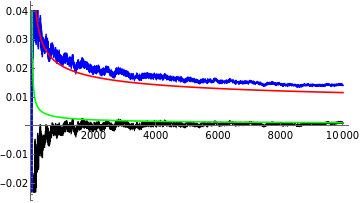}
\caption{\centering Graphs of $\PP(x)-\frac{1}{\log x}$ (in blue), $\frac{1}{(\log x)^2}$ (in red), $\PP(x)-\frac{1}{\log x-1}$ (in black), and $\frac{1}{(\log x)^3}$ (in green)  on $[1,10000]$}
 \label{overflow1bc}
\end{figure}

A graphical illustration of the prime number theorem is provided in Figure \ref{overflow1b}, which provides the graph of the functions $\PP(x)$, $\frac{1}{\log x}$, and $\frac{1}{\log x-1}$.
As Figure \ref{overflow1b} suggests, the function $\frac{1}{\log x-1}$ is in fact a much better approximation to $\PP(x)$ than is $\frac{1}{\log x}$.  In simple terms, one can show that
$$\left| \PP(x) - \frac{1}{\log x - 1}\right| < \left|\PP(x) - \frac{1}{\log x}\right|, \quad \forall x \geq 97.$$
More dramatically, one has
$$\PP(x)-\frac{1}{\log x} \sim \frac{1}{(\log x)^2} \ (x \to \infty),$$
while
$$\PP(x)-\frac{1}{\log x-1} \sim \frac{1}{(\log x)^3} \ (x \to \infty),$$
so that the error in the latter approximation is substantially less than the error in the former for large $x$.  
For example, one has $\frac{1}{\log (10^{29})} = \frac{1}{66.774967\ldots}$, so that $\frac{1}{\log (10^{29})-1} =  \frac{1}{65.774967\ldots}$,
while, as we noted earlier, $\PP(10^{29}) = \frac{1}{65.759271\ldots}$.  Figure  \ref{overflow1bc} provides  graphs of the functions $\PP(x)-\frac{1}{\log x}$, $\frac{1}{(\log x)^2}$, $\PP(x)-\frac{1}{\log x-1}$, and $\frac{1}{(\log x)^3}$  on $[1,10000]$.

One of the morals that one draws from the numerical examples above   is that  numbers even as large as $10^{29}$ are {\it small} when it comes to the study of the primes and that the true ``size'' of a prime $p$ for the purpose of studying prime asymptotics is better measured by $\log p$ than by $p$ itself.  This is in large part due to the prime number theorem $\PP(x) \sim \frac{1}{\log x} \ (x \to \infty)$ and the remarkable fact that it is only the zeroth term in the infinite ``asymptotic expansion'' 
\begin{align}\label{pias}
\PP(x) - \sum_{k= 0}^{n-1} \frac{k!}{(\log x)^{k+1}} \sim \frac{n!}{(\log x)^{n+1}} \ (x \to \infty), \quad \forall n \geq 0,
\end{align}
of $\PP(x)$.  Also, more practically, with respect to computers, the ``size'' of a prime $p$ is measured by how many bits it has in its binary expansion, which is given exactly by $\lfloor \log_2 p \rfloor+1$ and which is larger than $\log p$ asymptotically by a factor of $\frac{1}{\log 2} = 1.4426950408\ldots$.   A random prime number in the interval $[1,10^{29}]$ can be expressed using at most $97$ bits, and 97 is not a very large number, so, by this measure, neither is $10^{29}$.  As mentioned earlier, the largest known prime, $2^{82589933}-1$, has $24862048$ digits, and, expressed in binary notation, is the sequence $(11 1\ldots 1)_2$ of  $82589933$ bits.    In fact, there is a deterministic algorithm to test if an integer $n$ is prime that runs on $O((\log n)^{t})$ operations for any  $t > 6$, which is  ``polynomial time'' in  the number of bits of $n$, not in $n$ itself.   The best known algorithms to compute $p_n$, $\pi(n)$, and $\PP(n)$, on the other hand, are polynomial time in $n$, which is ``exponential time'' in the number of bits of $n$.  Unfortunately, this makes analyzing these functions from a purely computational perspective rather difficult.

As we have noted, the number $e$ encodes information about the distribution of the primes, and, as we have also seen, so does the number $\zeta(2) = \frac{\pi^2}{6}$.  There are many other such constants, including the  Mertens constant $H$ mentioned earlier.  After $e$ and $\pi$, arguably the most important in analytic number theory is the {\bf Euler--Mascheroni constant}\index{Euler--Mascheroni constant $\gamma$}\index[symbols]{.f t@$\gamma$}
$$\gamma = \lim_{x \to \infty}\left( \sum_{n \leq x} \frac{1}{n}-\log x\right) = 0.577215664901\ldots,$$
discovered by Euler in 1735, and the related constant
 $$e^\gamma = \lim_{n \to \infty}  \left(\frac{1}{n}\exp H_n \right)  = \lim_{x \to \infty} \left(\frac{1}{x}\prod_{n \leq x} e^{1/n}\right) =  1.781072417990\ldots,$$ where $H_n$ is the {\bf $n$th harmonic  number} $$H_n =\sum_{k = 1}^n \frac{1}{k}.\index{harmonic number $H_n$}\index[symbols]{.rt H@$H_n$}$$
 To prove that the limits above exist and are so related is a straightforward calculus exercise.  The harmonic numbers 
satisfy
$$H_n -\gamma - \log n  \sim \frac{1}{2n}  \ (n \to \infty)$$
and in fact
$$\frac{1}{2n+1} <  H_n-\gamma-\log n < \frac{1}{2n}, \quad \forall n \geq 1,$$
and thus $ \gamma+\log n $ is an excellent approximation of $H_n$ for large $n$.  In particular, since $$H_n \sim \log n \ (n \to \infty),$$ the prime number theorem is equivalent to
$$\pi(n) \sim \frac{n}{H_n} \ (n \to \infty),$$
where $$\frac{n}{H_n} = \left( \frac{ 1^{-1}+ 2^{-1}+ 3^{-1}+\cdots+ n^{-1} }{n} \right)^{-1}$$
is the harmonic mean of the numbers $1, 2, 3,\ldots,n$.  See \cite{ell3}, Proposition \ref{harmonicprime2}, and Remark \ref{harmonicprime} for further relationships between the prime counting function and the harmonic numbers.

The Euler--Mascheroni constant $\gamma$ also features prominently in {\it Mertens' theorems}, which are three well-known theorems that were proved by Mertens in 1874 \cite{mertens}, twenty-two years before the first proofs of the prime number theorem.   {\bf Mertens' third theorem}\index{Mertens' third theorem} is the identity
$$\lim_{x \to \infty}\left(\sum_{p\leq x} \log\left(1- \frac{1}{p}\right) +\log\log x\right) = -\gamma,$$
or, equivalently,
$$e^{\gamma} \prod_{p \leq x} \left(1-\frac{1}{p} \right)  \sim \frac{1}{\log x} \ (x \to \infty).$$  From Mertens' third theorem,  it follows immediately that the prime number theorem (which  in 1874 was a mere conjecture) is equivalent to the asymptotic relation
$$\PP(x)  \sim e^{\gamma} \prod_{p \leq x} \left(1-\frac{1}{p} \right) \ (x \to \infty),$$
or, equivalently, 
$$e^\gamma = \lim_{x\to \infty }  \frac{\PP(x)}{\prod_{{p\leq x}}\left(1-{\frac{1}{p}}\right)}.$$
These reformulations of the prime number theorem are profoundly beautiful.  The function 
$$\mathbf{P}(x) = \prod_{{p\leq x}}\left(1-{\frac{1}{p}}\right)$$ approximates the probability that a randomly selected number is not divisible by any prime less than or equal to $x$,  but assuming, erroneously of course, that these are independent events.   The constant $e^\gamma$ is therefore an asymptotic measure of how independent these events {\it aren't}.   Moreover, $\mathbf{P}(x) $ for any $x> 0$ is the limit
$$\mathbf{P}(x)  = \lim_{t \to \infty} \frac{\#\{n \in (0,t]: \text {all prime factors of $n$ are $> x$} \}}{t},$$
which is also the limit as $N \to \infty$ of the proportion of the integers in the interval $[1,N]$ having no prime factor less than or equal to $x$ \cite{lag2}.

{\bf Mertens' second theorem}\index{Mertens' second theorem} states that the limit
$$M = \lim_{x \to \infty} \left( \sum_{p \leq x} \frac{1}{p} - \log \log x \right) = 0.261497212847\ldots,$$  
known as the {\bf Meissel--Mertens constant},\index{Meissel--Mertens constant $M$}\index[symbols]{.f ta@$M$} exists.   This theorem is much stronger than Euler's result that the sum $\sum_p \frac{1}{p}$ diverges.  Moreover, it suggests the following analogy: the relationship between the divergent series $\sum_{n = 1}^\infty \frac{1}{n}$, the function $\log x$, and the constant $\gamma$ is analogous to the relationship between the divergent series $\sum_{p} \frac{1}{p}$, the function $\log \log x$, and the constant $M$.   Mertens' second theorem  is one of the first appearances, if not the first appearance, of an iterated logarithm function in number theory.  It has, as a consequence, the estimates
$$\sum_{p \leq 10^{100}} \frac{1}{p} \approx  \log \log (10^{100}) + M  = \log 100 + \log \log 10 + M \approx 5.700699\ldots$$
and
$$\sum_{p \leq 10^{10^{100}}} \frac{1}{p} \approx  \log \log (10^{10^{100}}) + M  = 100 \log 10 + \log \log 10 + M \approx 231.354038\ldots,$$
the values of which are miniscule when compared to a {\bf googol} $10^{100}$ and a {\bf googolplex} $10^{10^{100}}$, respectively.  Thus, although the series $\sum_{p \leq x} \frac{1}{p}$ diverges, it does so extraordinarily slowly.  This is one of the first instances in analytic number theory of a result showing that some patterns that emerge among the primes do so only at astronomically large numbers.  Moreover, the ``culprit'' in this particular instance is the iterated logarithm $\log \log x$.

The occurrence of $\log \log x$ in Mertens' second theorem can be motivated by {\it Cram\'er's model of the primes},  which, in very loose terms,  uses the prime number theorem to model the ``probability'' that an integer $n >1$ is prime as $\frac{1}{\log n}$ \cite{cramer1}.   Under Cram\'er's model, the sum $\sum_{p \leq x} \frac{1}{p}$ is approximated by the sum $\sum_{1 < n \leq x} \frac{1}{n \log n}$, and a straightforward calculus exercise, using the  fact that $$\int_e^x \frac{dt}{t \log t} = \log \log x, \quad \forall x >1,$$ shows that the limit
$$R = \lim_{x \to \infty} \left( \sum_{1 < n \leq x} \frac{1}{n \log n} -\log \log x \right)$$
exists, where the constant
$$R = 0.794678645452\ldots$$
has its sequence of digits given by Sequence A361972 of the On-Line Encyclopedia of Integer Sequences (OEIS).  It follows that Mertens' second theorem is equivalent to the existence of the limit
$$\lim_{x \to \infty} \left(\sum_{1 < n \leq x} \frac{1}{n \log n}-\sum_{p \leq x} \frac{1}{p} \right) = R-M = 0.533181432605\ldots.$$

 It is clear that Mertens' second theorem is also equivalent to
$$e^{-M} \prod_{p \leq x} e^{1/p} \sim \log x \ (x \to \infty).$$
Thus, given Mertens' second theorem, the prime number theorem is equivalent to $$\PP(x) \sim e^M \prod_{p \leq x} e^{-1/p} \ (x \to \infty),$$  
where $e^M = 1.371244130303\ldots$.  
{\bf Mertens' first theorem}\index{Mertens' first theorem} states that
\begin{align}\label{mert1}
\left| \log x - \sum_{p \leq x} \frac{\log p}{p}\right | < 2, \quad \forall x \geq 2.
\end{align}
In fact,  it is known that the limit
$$B = \lim_{x \to \infty} \left( \log x - \sum_{p \leq x} \frac{\log p}{p} \right) = 1.332582275733\ldots$$
exists, which in turn implies the prime number theorem.  It is also known that
\begin{align*}
\gamma = \lim_{x \to \infty} \left( \log x - \sum_{p \leq x} \frac{\log p}{p-1} \right),
\end{align*}
and therefore
$$B  = \gamma+ \sum_{p}  \frac{\log p}{p(p-1)}  = \gamma+  \sum_{n = 2}^\infty \sum_p \frac{\log p}{p^n} =  \gamma-\sum_{n = 2}^\infty P'(n).$$

As noted previously, another important constant relating to the prime numbers that was introduced by Mertens is the Mertens constant $H$, which is also given by
\begin{align*}
H  & =- \sum_p \left(\frac{1}{p} +\log\left(1- \frac{1}{p} \right) \right)  \\
 & =  \sum_p \left(\frac{1}{2p^2} + \frac{1}{3p^3} + \frac{1}{4p^4} + \cdots \right) \\
 & =  \sum_{n = 2}^\infty \frac{P(n)}{n}  \\
 & = 0.315718452053\cdots,
\end{align*}
where also 
$$e^{-H} = \prod_p \left(1-\frac{1}{p}\right)e^{1/p} = 0.729264744257\ldots.$$
An immediate consequence of Mertens' third theorem and the first expression above for the constant $H$ is the relationship
$$\gamma = M+H.$$
These three constants encode important information about the primes, and they crop up in many seemingly unrelated contexts in analytic number theory---much how the constants $e$ and $\pi$ make appearances in many seemingly unrelated contexts, as, for instance, in the normal distribution.    Notably,  for example (see Section 3.6), one has
$$\gamma = \lim_{x \rightarrow 1^+} \left(\zeta(x) - \frac{1}{x-1}\right) = \sum_{n = 2}^\infty \frac{(-1)^n \zeta(n)}{n} =   \int_1^\infty \frac{1-\{x\}}{x^{2}} \, dx$$
and
$$H =   \lim_{x \rightarrow 1^+} \left( \log \frac{1}{x-1} - P(x) \right)    =   \int_1^\infty  \frac{\li(x)-\pi(x)}{x^{2}} \, dx.$$

\begin{remark}[Algebraic numbers and transcendental numbers]\label{eulerconst}
The most remarkable example, in the author's opinion, of an equation in mathematics  that combines an assortment of fundamental constants is {\bf Euler's equation}\index{Euler's equation}
$$1+e^{i \pi} = 0,$$
which follows immediately from the well-known {\bf Euler's formula}\index{Euler's formula}
$$e^{ix} = \cos x+ i \sin  x, \quad \forall x \in \RR,$$
by plugging in $x = \pi$.  (See Remark \ref{EFpf} for a ``Calculus 1'' proof of Euler's formula.)  Arguably, the five mathematical constants appearing in Euler's formula are the five most fundamental constants in all of mathematics: $0$ and $1$ are the most important constants in arithmetic, $e$ the most important in analysis, $\pi$ the most important in geometry, and $i$ is one of the most important constants in both algebra and complex analysis.     Another startling aspect of Euler's equation $1+e^{i \pi} = 0$ 
 is that it uses each of the three fundamental arithmetical operations---addition, multiplication, and exponentiation---exactly once.    The renowned physicist Richard Feynman wrote in his {\it Lectures on Physics} that Euler's formula is ``our jewel'' and is ``one of the most remarkable, almost astounding, formulas in all of mathematics'' \cite[pp.\ 22-1, 22-10]{feyn}, and a poll of readers conducted in 1990 by The Mathematical Intelligencer named Euler's equation as the ``most beautiful theorem in mathematics'' \cite{wells}.  More importantly,  Euler's equation is useful for deriving an important property of the number $\pi$: the  equation,  along with the proof of the {\it transcendence} of the number $e$ in 1873 by Hermite, formed the bases for the first proof of the transcendence of $\pi$ by  Lindemann in 1882.   

Here, a complex number is said to be {\bf transcendental}\index{transcendental number} if is not {\it algebraic},\index{algebraic number} where a complex number is said to be {\bf algebraic} if it is a root of a nonzero polynomial with rational coefficients.  For example, trivially every rational number is algebraic, and $\sqrt{2}$, though irrational, is algebraic since it is a root of the polynomial $x^2-2$.  One can show that the set $\overline{\QQ}$ of all algebraic numbers is a countable subfield of the field $\CC$, and the transcencence of both $e$ and $\pi$ can then be expressed simply as $e, \pi \notin \overline{\QQ}$.   Note that a complex number $\alpha$ is transcendental if and only if the  unique surjective ring homomorphism $\QQ[X]\longrightarrow \QQ[\alpha]$ sending $X$ to $\alpha$ is an isomorphism, if and only if the field extension $\QQ(\alpha)$ of $\QQ$ is infinite dimensional as a vector space over $\QQ$.   More generally, complex numbers $\alpha_1, \alpha_2, \ldots, \alpha_n$ are said to be {\bf algebraically independendent}\index{algebraically independent} if the unique ring homomorphism $\QQ[X_1,X_2,\ldots,X_n] \longrightarrow \CC$ sending $X_i$ to $\alpha_i$ for all $i$ is injective, that is, if $f(\alpha_1, \alpha_2, \ldots, \alpha_n) \neq 0$ for all nonzero polynomials $f \in \QQ[X_1,X_2,\ldots,X_n]$ with rational coefficients in $n$ variables.  It is an important open problem in the field of {\it transcendental number theory}---whose primary goal is develop techniques for proving that various  collections of constants in mathematics are algebraically dependent or independent---to prove or disprove the conjecture that the numbers $e$ and $\pi$ are algebraically independent.   It is also believed that the three constants $\gamma$, $M$, and $H$ are transcendental, but, in fact, no one has yet proved that any of the three constants is irrational, much less transcendental.   
\end{remark}

\begin{remark}[Mathematical constants in physics]\label{physics1}
The mathematical constants $e$, $\pi$, and $i$ are important not only in mathematics, but also in other sciences,  physics especially.  The constants $\pi$ and $i$, for example, appear in {\it Schr\"odinger's equation} ${ {\hat {H}}\vert \Psi (t)\rangle = i\frac{h}{2\pi} {\frac {d}{dt}}\vert \Psi (t)\rangle}$, which is the quantum counterpart of {\it Newton's second law} $\vec{F} =  \frac{d}{dt}\vec{p}$, where $h  = 6.62607015\cdot 10^{-34} \operatorname{J}\cdot \operatorname{Hz}^{-1}$ is a universal physical constant known as {\bf Planck's constant}.  Even the constant $\gamma$ arises in physics, as, for example, in {\it dimensional regularization of Feynman diagrams}  in quantum field theory.    For applications of special values of the Riemann zeta function to physics, see  Remark \ref{physics}.
\end{remark}

\begin{remark}[Mersenne primes]
The {\bf Lenstra--Pomerance--Wagstaff conjecture},\index{Lenstra--Pomerance--Wagstaff conjecture} conjectured independently by  Lenstra and Pomerance \cite{pome},  modifying a prior conjecture of Gillies \cite{gill},  states that the number of Mersenne primes less than or equal to $x$ is asymptotic to $e^\gamma \log_2 \log_2x$, or, equivalently,  that the number of Mersenne primes $2^p-1$ with $p \leq x$ is asymptotic to $e^\gamma \log_2 x$.  Equivalently still,  the conjecture states that the proportion of primes $p \leq x$ for which   $2^p-1$ is also prime is asymptotic to  $\frac{e^\gamma}{ \log 2} \frac{(\log x)^2}{x}$.
\end{remark}

\begin{remark}[Twin primes and Sophie Germain primes]
The {\it Hardy--Littlewood conjecture for prime $k$-tuples}, applied to the twin primes,  states  that
$${\displaystyle \pi_{2}(x)\sim 2C_{2}{\frac {x}{(\log x)^{2}}}\sim 2C_{2}\int_{0}^{x}{dt \over (\log t)^{2}}} \ (x \to \infty),$$
where $\pi_2(x)$ denotes the number of primes $p \leq x$ such that $p+2$ is also prime, where 
$$C_{2}=\prod_{\textstyle {p\;{\rm {prime}} \atop p \neq 2}}\left(1-{\frac {1}{(p-1)^{2}}}\right) =  0.660161815846\ldots \index{twin prime constant $C_2$}$$
is the {\bf twin prime constant}, and where $2C_2 = 1.320323631693\ldots$.  The conjecture is equivalent to
$$q_n \sim \frac{1}{2C_2} n (\log n)^2 \ (n \to \infty),$$
where $q_n$ is the sequence that enumerates  in succession the first of each twin prime pair, so that the sequence is $3,5,11,17,29,41,\ldots$.  By the prime number theorem, the conjecture is also equivalent to
$$\frac{\pi_2(x)}{\pi(x)} \sim  \frac{2C_2}{\log x}  \sim 2C_2 \PP(x) \ (x \to \infty),$$
which means that the probability $\frac{\pi_2(x)}{\pi(x)}$ that a randomly chosen prime number less than or equal $x$ is the smaller of two twin primes is asymptotic to $2C_2 \PP(x)$, i.e., $2C_2$ times the probability that a randomly chosen number less than or equal to $x$ is prime.   In other words, the conjecture says that the likelihood that $n+2$ is prime is increased asymptotically by a factor of $2C_2$ if it is known also that $n$ is prime, and thus, if the conjecture is true, then $2C_2$ can be thought of as a ``prime coupling constant.''  The  Hardy--Littlewood conjecture for twin primes is the likely analogue of the prime number theorem for twin primes, but since it implies that there are infinitely many twin primes, the proof is likely to be far more difficult, perhaps involving new methods that have not been discovered yet.  However, in 2007, J.\ Wu made some  progress on the conjecture by proving \cite{wu} that
\begin{align*}
\pi_{2}(x) \leq 3.3996 \cdot 2 C_2\frac{x}{(\log x)^{2}},  \quad \forall x \gg 0.
\end{align*}
Wu's theorem is a partial analogue for $\pi_2(x)$ of Chebyshev's  1850 result (\ref{chebin}). Wu's result  improved upon a 1914 result of V.\ Brun, which had implied that the sum
$$B_2 = \sum \limits_{p\,:\,p+2\in \mathbb {P} }{\left({{\frac {1}{p}}+{\frac {1}{p+2}}}\right)}={{\frac {1}{3}}+{\frac {1}{5}}}+{{\frac {1}{5}}+{\frac {1}{7}}}+{{\frac {1}{11}}+{\frac {1}{13}}}+\cdots,$$
known now as {\bf  Brun's constant},\index{Brun's constant} is finite, which stands in sharp contrast to the divergence of the sum of the reciprocals of the primes. It is known that $1.8304< B_2 < 2.347$, and various heuristic estimates of $B_2$ have been made, e.g., in  2002, P.\ Sebah and P.\ Demichel used all twin primes up to $10^{16}$ to extrapolate the  estimate $B_2\approx  1.902160583$.

A prime $p$ is a {\bf Sophie Germain}\index{Sophie Germain prime} prime if $2p+1$  is also prime.  It is conjectured that the number of Sophie Germain primes less than or equal to $x$ is asymptotic to 
$$2C_{2}{\frac {x}{(\log x)^{2}}}\sim 2C_{2}\int_{0}^{x}{dt \over (\log t)^{2}}\sim  \pi_{2}(x)\ (x \to \infty).$$
Both of these asymptotics, as well as the Hardy--Littlewood conjecture for prime $k$-tuples,  follow from the far more general {\it Bateman--Horn conjecture} \cite{batem}.
\end{remark}

\section{The prime numbers, analytically}

The prime number theorem is one of the major mathematical achievements of the 19th century.  The fact that the theorem was  proved by two mathematicians both in the same year (1896) was no coincidence.   As mentioned earlier,  much progress had been made by Chebyshev  in 1848 and 1850.
Nine years later, major groundwork for the eventual 1896 proofs was laid down by Riemann in his landmark eight-page paper of 1859 \cite{rie}, which was the only work he ever wrote on number theory.   Riemann first was able to show that the function $$\zeta(s) = \sum_{n = 1}^\infty \frac{1}{n^s}$$ on the right half plane $\{s \in \CC: \operatorname{Re} s > 1\}$ extends (uniquely) to a meromorphic function on all of $\CC$ with a single (simple) pole at $s = 1$ with residue $1$.  In fact, the extended function $\zeta(s)-\frac{1}{s-1}$ is entire, with limiting value
 $$ \lim_{s \to 1} \left( \zeta(s)-\frac{1}{s-1} \right) = \gamma$$ at $s = 1$.   Note that $\frac{1}{s-1}$ is the meromorphic continuation to $\CC$ of  the function
$$\int_{1}^\infty \frac{dx}{x^s} = \frac{1}{s-1}$$
on $\{s \in \CC: \operatorname{Re} s> 1\}$,  where the given integral is employed in the integral test from calculus to determine the region of absolute convergence of the sum $\sum_{n = 1}^\infty \frac{1}{n^s}$.  In particular, one has
$$\gamma = \lim_{s \to 1 \atop {\operatorname{Re} s> 1}} \left( \sum_{n = 1}^\infty \frac{1}{n^s}-\int_{1}^\infty \frac{dx}{x^s} \right).$$

Today, it is the meromorphic function $\zeta(s)$ on all of $\CC$ that we refer to as the {\bf Riemann zeta function}.\index{Riemann zeta function $\zeta(s)$}  One may think of the Riemann zeta function as a simultaneously arithmetic, algebraic, and analytic generating function for the sequence of positive integers.  Both proofs of the prime number theorem involved very deep results and conjectures from  Riemann's paper, specifically concerning a remarkable explicit formula that Riemann conjectured for the prime counting function $\pi(x)$, expressed in terms of zeros of the Riemann zeta function $\zeta(s)$.    One of the most important outstanding open problems in all of mathematics today is to prove or disprove the famous {\bf Riemann hypothesis},\index{Riemann hypothesis} conjectured by Riemann in his 1859 paper, which is the statement that the zeros of $\zeta(s)$ besides the negative even integers $-2,-4,-6,-8,\ldots$ all have real part $\frac{1}{2}$, that is, they all lie on the {\bf critical line}\index{critical line} $$\{s \in \CC: \operatorname{Re} s = \tfrac{1}{2}\} = \{\tfrac{1}{2}+it: t\in \RR\}.$$  Riemann was able to show in his paper that  the real zeros of $\zeta(s)$ are precisely the negative even integers and that all of the non-real zeros of  $\zeta(s)$ lie in the {\bf criticial strip}\index{critical strip} $\{s \in \CC: 0 \leq \operatorname{Re} s \leq 1\}$.   He also showed that they are discretely ordered in the critical strip, occur in conjugate pairs $\rho$ and $\overline{\rho}$ (which are reflections of each other over the real axis), as well as in pairs $\rho$ and $1-\overline{\rho}$ (which are reflections of each other over the critical line).  Thus, the non-real zeros of $\zeta(s)$, also called the {\bf nontrivial zeros}\index{nontrivial zeros of $\zeta(s)$} of $\zeta(s)$, can be listed in order of nondecreasing absolute imaginary part.  In 1914, Hardy proved that there are infinitely many  zeros of $\zeta(s)$ on the critical line \cite{har2}.  In more recent times, mathematicians, using sophisticated computers and algorithms,  have calculated at least the first ten trillion zeros of $\zeta(s)$  to very high precision, and they all lie on the critical line \cite{gour}.  This provides  some evidence, though circumstantial, that the Riemann hypothesis is true.

\begin{table}[!htbp]
  \caption{\centering Some noteworthy real values of $\zeta(s)$}
  \footnotesize
\begin{tabular}{|c|c|} \hline
$s$ & $\zeta(s)$ \\ \hline\hline
$2n$, where $n \in \ZZ_{> 0}$   & ${\frac {(-1)^{n+1}B_{2n}(2\pi )^{2n}}{2(2n)!}}$ \\ \hline
$-n$, where $n \in \ZZ_{\geq 0}$ & ${\frac {(-1)^{n}B_{n+1}}{n+1}}$  \\ \hline
$-2n$, where $n \in \ZZ_{> 0}$ & $0$  \\ \hline
$-2n+1$, where $n \in \ZZ_{> 0}$ & $-{\frac {B_{2n}}{2n}}$  \\ \hline
$0$ & $-\frac{1}{2}$ \\ \hline
$1$ & $\infty$ \\  \hline
$2$ & $\frac{\pi^2}{6} = 1.644934\ldots$ \\ \hline
$3$ & $1.202056\ldots$ \\ \hline
$4$ & $\frac{\pi^4}{90} = 1.082323\ldots$ \\ \hline
$5$ & $1.036927\ldots $ \\ \hline
$6$ & $\frac{\pi^6}{945}= 1.017343\ldots$ \\ \hline
$7$ & $1.008349\ldots$ \\ \hline
$8$ & $\frac{\pi^8}{9450}= 1.004077\ldots$ \\ \hline
$-1$ & $-\frac{1}{12}$ \\ \hline
$-2$ & $0$ \\ \hline
$-3$ & $\frac{1}{120}$ \\ \hline
$-4$ & $0$ \\ \hline
$-5$ & $-\frac{1}{252}$ \\ \hline
$-6$ & $0$ \\ \hline
$-7$ & $\frac{1}{240}$ \\ \hline
$-8$ & $0$ \\ \hline
$\tfrac{1}{2} \pm 14.134725\ldots i$ & $0$  \\ \hline
$\tfrac{1}{2}\pm 21.022040 \ldots i$ & $0$  \\ \hline
$\tfrac{1}{2}\pm 25.010858  \ldots i$ & $0$  \\ \hline
$\tfrac{1}{2}\pm 30.424876 \ldots i$ & $0$  \\ \hline
$\tfrac{1}{2}\pm  32.935062 \ldots i$ & $0$  \\ \hline
$\tfrac{1}{2}\pm 37.586178 \ldots i$ & $0$  \\ \hline
$\tfrac{1}{2}\pm 40.918719 \ldots i$ & $0$  \\ \hline
$\tfrac{1}{2}\pm 43.327073 \ldots i$ & $0$  \\ \hline
$\tfrac{1}{2}\pm 48.005150 \ldots i$ & $0$  \\ \hline
$\tfrac{1}{2}\pm 49.773832 \ldots i$ & $0$  \\ \hline
\end{tabular}\label{zetavalues}
\end{table}

Some noteworthy real values of $\zeta(s)$, discovered by Euler, Riemann, and others, are provided in Table \ref{zetavalues}, including the (exactly) 20 nontrivial zeros of $\zeta(s)$ with absolute imaginary part  less than or equal to $50$.  In the table, $B_n$ for all nonnegative integers $n$ denotes the {\bf $n$th Bernoulli number},\index{Bernoulli number $B_n$}\index[symbols]{.rt K@$B_n$} which can be defined via their generating function
$$ \sum_{n=0}^{\infty }\frac {B_{n}}{n!}X^{n} = \frac{X}{e^X-1} = \left(\sum_{n = 0}^\infty \frac{1}{(n+1)!}X^n \right)^{-1},$$
or explicitly by
$$B_n = \sum_{k = 0}^n \frac{1}{k+1} \sum_{j = 0}^k (-1)^j {k \choose j} j^n,$$
whence the $B_n$ are rational numbers for all $n$.   Moreover, one has $B_n = 0$ if and only if $n$ is an odd integer greater than $1$, and $B_n >0$ if and only if $n = 0$ or $n$ is congruent to $2$ modulo $4$.
The Bernoulli numbers have importance in mathematics well beyond just the Riemann zeta function, e.g., they are employed in the {\it Euler--Maclaurin formula} (Theorem \ref{EMthm}).  

\begin{figure}[ht!]
\includegraphics[width=80mm]{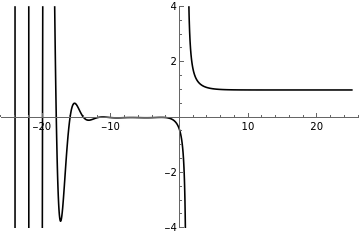}
\caption{\centering Graph of $\zeta(t)$ on $[-25,25]$}
   \label{zeta1}
\end{figure}

\begin{figure}[ht!]
\includegraphics[width=80mm]{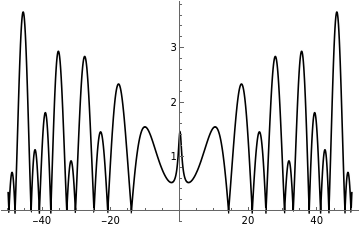}
\caption{\centering  Graph of $|\zeta(\tfrac{1}{2}+it)|$ on $[-50,50]$}
   \label{zeta4}
\end{figure}

\begin{figure}[ht!]
\includegraphics[width=80mm]{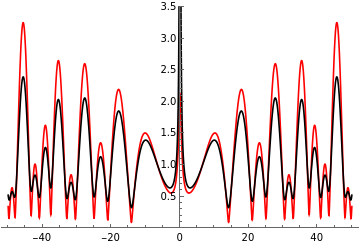}
\caption{\centering Graph of $|\zeta(1+it)|$ (in black) and $|\zeta(\tfrac{2}{3}+it)|$ (in red)  on $[-50,50]$}
\label{zeta23}
\end{figure}

\begin{figure}[ht!]
\includegraphics[width=90mm]{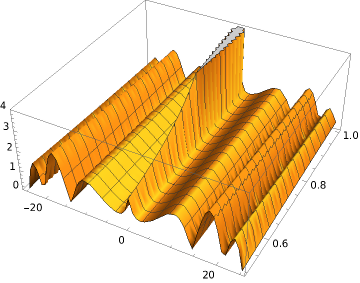}
\caption{\centering Graph of $|\zeta(a+it)|$ on $[-25.010858  \ldots, 25.010858  \ldots]\times[\tfrac{1}{2},1]$}
\label{zeta6}
\end{figure}

\begin{figure}[ht!]
\includegraphics[width=120mm]{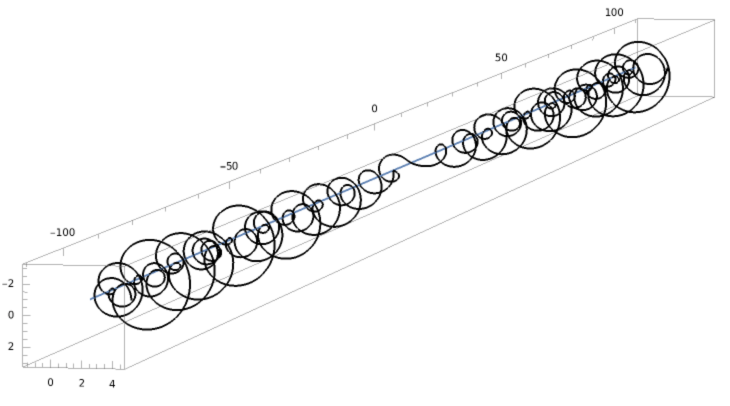}
\caption{\centering Graph of $(t,\zeta(\tfrac{1}{2}+it))$ and $(t,0+i0)$ on $[-100,100]$}
  \label{zeta3d}
\end{figure}

\begin{figure}[ht!]
\includegraphics[width=50mm]{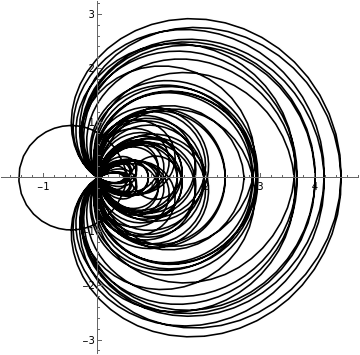}
\caption{Parametric plot of $\zeta(\tfrac{1}{2}+it)$ on $[-100,100]$}
   \label{zeta1_12a}
\end{figure}

 Figure \ref{zeta1} is a graph of the function $\zeta(t)$ on $[-25,25]$.    Figure \ref{zeta4} is a graph of the function $|\zeta(\tfrac{1}{2}+it)|$, to be contrasted with the graphs of $|\zeta(\tfrac{2}{3}+it)|$ and $|\zeta(1+it)|$ in Figure \ref{zeta23},  on $[-50,50]$.  Note that $\zeta(s) = 0$ if and only if $|\zeta(s)| = 0$, so that these graphs provide a ``snapshot''  visualization of the Riemann hypothesis.  A better visualization would be a moving graph of $|\zeta(a+it)|$ with a ``draggable'' parameter $a \in [0,1]$, but since we cannot provide such a graph here, instead we provide the 3D graph of
 $|\zeta(a+it)|$ for $(a, t) \in[-25.010858  \ldots, 25.010858  \ldots]\times[\tfrac{1}{2},1]$ as in  Figure \ref{zeta6},   whose domain contains the first six nontrivial zeros of $\zeta(s)$, i.e., the first three with positive imaginary part and the first three with negative imaginary part.   Figure \ref{zeta3d} is a 3D plot of the graph of $(t,\zeta(\tfrac{1}{2}+it))$ on $[-100,100]$, and Figure  \ref{zeta1_12a} is its projection onto the complex plane.  These graphs provide a way of visualizing the function $\zeta(s)$ and its zeros.  Further graphs for this purpose are provided in Section 4.2.

\begin{figure}[ht!]
\includegraphics[width=65mm]{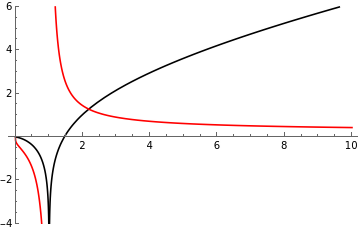}
\caption{\centering Graph of $\li(x)$ (in black) and its derivative $\frac{1}{\log x}$ (in red) on $[0,10]$}
\label{lilog}
\end{figure}

\begin{figure}[ht!]
\includegraphics[width=65mm]{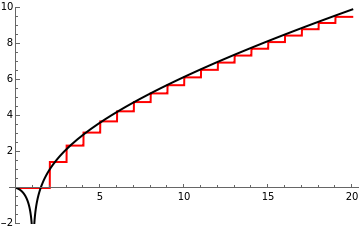}
\caption{\centering Graph of $\li(x)$ (in black) and the ``logarithmic sum'' function $\sum_{1< n \leq x } \frac{1}{\log n}$ (in red) on $[0,20]$}
 \label{ls}
\end{figure}

To describe the explicit formula for $\pi(x)$ in terms of the zeros of $\zeta(s)$, we first travel back a little to 1838, when Dirichlet observed that $\pi(x)$ can be well approximated by the {\bf logarithmic integral function}\index{logarithmic integral function $\li(x)$}  
$$\li(x) = \int_0^x \frac{dt}{\log t},$$
where the Cauchy principal value of the integral is assumed.  The function $\li(x)$ is a particular antiderivative of $\frac{1}{\log x}$.  Figure \ref{lilog} provides graphs on the interval $[0,10]$  of both $\li(x)$ and its derivative $\frac{1}{\log x}$.
Since $$\frac{d}{dx}  \frac{x}{\log x}  =   \frac{\log x -1}{(\log x)^2} \sim  \frac{1}{\log x}  = \frac{d}{dx} \li(x) \ (x \to \infty),$$ L'H\^opital's rule implies that  $$\li(x) \sim \frac{x}{\log x} \ (x \to \infty).$$   It follows that the prime number theorem is equivalent to
$$\pi(x) \sim \li(x) \ (x \to \infty).$$
However, $\li(x)$ is a significantly better approximation to $\pi(x)$ than is $\frac{x}{\log x}$ or any other known rational function of $x$ and $\log x$.    Note that the logarithmic integral approximation of $\pi(x)$ is motivated by Cram\'er's model  of the primes, since, under that model, $\pi(x)$ is approximated by $\sum_{1< n \leq x } \frac{1}{\log n}$, and an easy calculus exercise shows that the limit
$$\lim_{n \to \infty} \left(\li(x) -\sum_{1< n \leq x } \frac{1}{\log n} \right) > 0$$
exists: see Figure \ref{ls} for a graph of the ``logarithmic sum'' function $\sum_{1< n \leq x } \frac{1}{\log n}$ in comparison to $\li(x)$.   However, it was Riemann's  brilliant insight,  described below, that explained more precisely how the functions $\pi(x)$ and $\li(x)$ are related to one another.

Let $$\pi_0(x) = \lim_{\varepsilon \to 0} \frac{\pi(x+\varepsilon)+ \pi(x-\varepsilon)}{2} = \frac{\pi(x^+)+\pi(x^-)}{2}, \quad \forall x \geq 0.\index[symbols]{.s BAA@$\pi_0(x)$}$$
The function $\pi_0(x)$ is equal to the prime counting  function $\pi(x)$ except at its discontinuities (namely, at the primes), where $\pi_0(x)$ assumes the average of the limit from the right and the limit from the left of $\pi(x)$.
The {\bf Riemann--von Mangoldt explicit formula for $\pi_0(x)$}\index{Riemann--von Mangoldt explicit formula for $\pi_0(x)$} states that
\begin{align}\label{pi0s}
\pi_0(x) = \Ri(x) - \sum_{n=1}^\infty \sum_\rho \frac{\mu(n)}{n}\Ei\left(\frac{\rho \log x}{n}\right), \quad \forall x > 1,
\end{align}
where the inner sum is over all of the zeros $\rho$ of the Riemann zeta function, with the nontrivial zeros taken in order of increasing absolute value of the imaginary part and repeated to multiplicity, and with the real zeros  summed in the natural order $-2,-4,-6,-8,\ldots$.  Here, $\mu(n)$ is the {\it M\"obius function} and $\Ei(s)$ is the {\it complex exponential integral function}, defined in Sections 3.3 and 4.5, respectively, and {\bf Riemann's function} $\Ri(x)$ is defined by
$$\Ri(x) = \sum_{n=1}^\infty \frac{ \mu(n)}{n} \li(x^{1/n}) = \sum_{n=1}^\infty \frac{ \mu(n)}{n} \Ei\left(\frac{\log x}{n}\right),$$  
where the latter equality follows from the fact that $\li(x) = \Ei(\log x)$ for all $x > 0$.   In particular, the Riemann--von Mangoldt explicit formula implies that all of the information about the prime numbers is encoded by the zeros of $\zeta(s)$, along with  the functions $\Ei(s)$ and $\mu(n)$.  See \cite{hut} for an implementation of the Riemann von--Mangoldt explicit formula in Sage.  A more in-depth discussion of the formula and its relation to the prime number theorem is provided in Sections 5.1 and 5.2.

In his landmark paper, Riemann wrote down an equivalent form of the explicit formula (\ref{pi0s}) (specifically, (\ref{realexplicit0})).  However, Riemann's formula was not proved rigorously until 1895, by von Mangoldt.  Then, just one year later, the prime number theorem was finally proved.      Both  1896 proofs of the prime number theorem relied on showing that the prime number theorem is equivalent to $\zeta(s)$ having no zeros on the line $\{s \in \CC: \operatorname{Re} s  = 1\}$, and then proving that in fact $\zeta(s)$ has no such zeros.  The Riemann--von Mangoldt explicit formula is precisely what made the first part of those proofs possible.   Thus, both Riemann and von Mangoldt played major roles in the eventual proofs of the prime number theorem.

\begin{figure}[ht!]
\includegraphics[width=80mm]{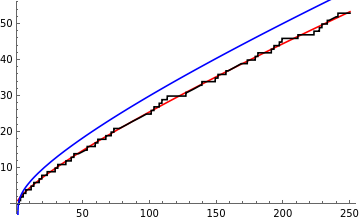} 
\caption{\centering Graph of $\pi(x)$ (in black), $\Ri(x)$ (in red),  and $\li(x)$ (in blue)  on $[1, 250]$}
 \label{RP}
\end{figure}

\begin{figure}[ht!]
\includegraphics[width=80mm]{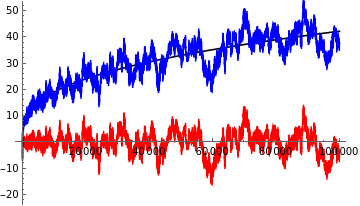}
\caption{\centering Graph of  $\li(x)-\Ri(x)$ (in black), $\Ri(x)-\pi(x)$ (in red) and $\li(x)-\pi(x)$ (in blue) on $[1, 100000]$}
\label{PLR}
\end{figure}

The ``main term'' of the  explicit formula (\ref{pi0s}) for $\pi_0(x)$ is the smooth function $\Ri(x)$, which captures the ``size'' of $\pi(x)$ in the sense that
$$\pi(x) \sim \pi_0(x) \sim \Ri(x) \sim \li(x) \ (x \to \infty)$$
The terms $\frac{\mu(n)}{n}\left( \Ei\left(\frac{\rho \log x}{n}\right) + \Ei\left(\frac{\overline{\rho} \log x}{n}\right)\right)$,  grouped together in conjugate pairs, are  smooth oscillatory terms.  Along with the somewhat neglible terms $\frac{\mu(n)}{n}\Ei\left(\frac{-2k \log x}{n}\right)$ over the real zeros $\rho = -2k$,  together they precisely capture all of the deviations of $\pi_0(x)$ from $\Ri(x)$, in that their sum over all $n$ and $\rho$ is equal to $\Ri(x)-\pi_0(x)$.    Figure \ref{RP} provides a graph of the functions $\pi(x)$, $\li(x)$, and $\Ri(x)$ on the interval $[1,250]$.   Figure \ref{PLR} compares the differences $\li(x)-\Ri(x)$, $\Ri(x)-\pi(x)$ and $\li(x)-\pi(x)$ on the interval $[1,100000]$.  One has the relationship
$$\li(x)-\pi(x) =  (\Ri(x)-\pi(x)) +( \li(x)-\Ri(x)),$$
where the  smooth function
\begin{align}\label{liri}
\li(x)-\Ri(x) \sim \frac{\li(x^{1/2})}{2} \sim  \frac{x^{1/2}}{2\log (x^{1/2})} = \frac{\sqrt{x}}{\log x} \ (x \to \infty)
\end{align}
tracks the function $\li(x)-\pi(x)$ much more closely than it does the function $\Ri(x)-\pi(x)$, at least for small $x$.

In 1899, just a few years after his proof of the prime number theorem, de la Vall\'ee Poussin proved a much stronger version of the prime number theorem \cite{val2}, which we call the {\bf prime number theorem with error bound}.\index{prime number theorem with error bound} 

\begin{theorem}[{Prime number theorem with error bound \cite{val2}}]\label{PNTET1} There exists a constant $C > 0$ such that
\begin{align*}
\li(x)-\pi(x) = O \left(\frac{x}{e^{C\sqrt{\log x}}} \right) \ (x \to \infty)
\end{align*}
\end{theorem}

Since $x^{t} = o( e^{c\sqrt{x}})  \ (x \to  \infty)$ for all $t \in \RR$ and all $c > 0$, the prime number theorem with error bound implies that
\begin{align}\label{PNTET}
 \li(x) -\pi(x) =  o\left(\frac{x}{(\log x)^{t}}\right) \ (x \to \infty)
\end{align}
for all $t \in \RR$.   By contrast,  the prime number theorem is equivalent to
$$  \li(x) -\pi(x) =  o(\li(x)) \ (x \to \infty),$$
and thus to
$$\ \li(x)  - \pi(x) =  o\left(\frac{x}{\log x} \right) \ (x \to \infty).$$
Thus, the prime number theorem with error bound is a substantial improvement over the prime number theorem.

\begin{align}\label{vk}
\li(x)-\pi(x) = O(\sqrt{x} \log x) \ (x \to \infty).
\end{align}
In 1976 \cite{schoen}, Schoenfeld proved an effective version of von  Koch's result, namely, that the Riemann hypothesis is equivalent to
$$|\li(x)-\pi(x)| \leq \frac{1}{8\pi} \sqrt{x} \log x, \quad \forall x \geq 2657.$$

Of course, there is no guarantee that the Riemann hypothesis is true.   
Thus, we follow an  unspoken tradition (following Ingham in \cite{ing2}, for example) and let
 $$\Theta = \sup\{\operatorname{Re} \rho : \rho \in \CC\backslash \RR,\, \zeta(\rho) = 0\}$$
denote the supremum of the real parts of the nontrivial zeros of the Riemann zeta function, which we call the {\bf Riemann constant}.\index{Riemann constant $\Theta$} Equivalently, one has
$$\Theta =  \inf\{t \in \RR:  \forall s \in \CC\backslash \RR \ ( \operatorname{Re} s> t  \Rightarrow \zeta(s) \neq  0 )\}.$$
Since $\zeta(s)$ has zeros (in fact, infinitely many) on the critical line,  and all of its nontrivial zeros lie in the critical strip and are symmetric about the critical line, one has
$$\tfrac{1}{2} \leq \Theta \leq 1$$
and
$$0 \leq 1-\Theta =  \inf\{\operatorname{Re} \rho : \rho \in \CC\backslash \RR,\, \zeta(\rho) = 0\} \leq \tfrac{1}{2}.$$
It follows that the vertical strip $\{s \in \CC: \operatorname{Re} s \in [1-\Theta,\Theta]\}$ is the smallest closed vertical strip containing all of the nontrivial zeros of $\zeta(s)$.  Note also that $\Theta-\tfrac{1}{2}$ is the supremum of the distances of the nontrivial zeros of $\zeta(s)$ to the critical line.
It follows from these observations that Riemann hypothesis is equivalent to $\Theta = \tfrac{1}{2}$, to $1-\Theta = \tfrac{1}{2}$, and to $\Theta = 1-\Theta$.  In particular, the problem of settling the Riemann hypothesis generalizes as follows.

\begin{outstandingproblem}\label{thetaprob}
Compute the Riemann constant $\Theta$.
\end{outstandingproblem}

 One of the main reasons that the Riemann constant is of fundamental importance  in number theory is that von Koch's Riemann hypothesis equivalent (\ref{vk}) generalizes to the following unconditional result.

\begin{theorem}[{\cite[Theorems 30 and 31]{ing2} \cite[Theorem 12.3]{ivic} \cite[Theorem 15.2 and Exercise 13.1.1.1]{mont}}]\label{RC}
One has
\begin{align*}
\Theta = \inf\left\{t \in \RR: \li(x)- \pi(x) = O(x^t) \ (x \to \infty)\right\}
\end{align*}
and
\begin{align*}
\Theta = \min\left\{t \in \RR: \li(x)- \pi(x) = O(x^t \log x) \ (x \to \infty)\right\}.
\end{align*}
\end{theorem}

It follows from Theorem \ref{RC} that the Riemann constant is the unique real number  $\Theta$ such that the $O$ bound $\li(x)- \pi(x) = O(x^t) \ (x \to \infty)$ holds for all $t > \Theta$ and fails for all $t < \Theta$.   Therefore, since $\Theta \geq \tfrac{1}{2}$,  one has
$$\li(x)-\pi(x) \neq O(x^{t}  ) \ (x \to \infty), \quad \forall t < \tfrac{1}{2}.$$
 It also follows that,  if the Riemann hypothesis is false, that is, if some zero $\rho$ of $\zeta(s)$ lies to the right of the critical line,  then $\Theta \geq \operatorname{Re} \rho> \tfrac{1}{2}$ and therefore
$$\li(x)-\pi(x) \neq O (x^{t}) \ (x \to \infty)$$
for all $t < \operatorname{Re}\rho$, including $t = \tfrac{1}{2}$ and the uncountably many real numbers strictly between $\tfrac{1}{2}$ and $\operatorname{Re}\rho$.  Thus, the Riemann hypothesis says that the approximation $\li(x)$ of $\pi(x)$ is about as close to $\pi(x)$ as is {\it prima facie} possible.    More generally, the Riemann constant $\Theta$ sets precise limits on how well $\li(x)$ approximates $\pi(x)$ for large $x$, and,  moreover, the smaller $\Theta$ is, the better the approximation.  The worst case scenario for the error function $\li(x)-\pi(x)$, then, is $\Theta = 1$, which we dub the {\bf anti-Riemann hypothesis}.\index{anti-Riemann hypothesis}   In  spirit opposite to the Riemann hypothesis, the anti-Riemann hypothesis says that de la Vall\'ee Poussin's prime number theorem with error bound is close to the best error bound  possible.

To this day, proofs of the strongest known bounds on the error $\li(x) - \pi(x)$  are based on the Riemann--von Mangoldt explicit formula for certain functions related to $\pi(x)$ and  rather sophisticated methods for verifying zero-free regions for $\zeta(s)$ in the critical strip.   Even the most current of methods have not proved $\Theta < 1$ and allow us only to bound the zeros of $\zeta(s)$ asymptotically away from the line $\{s \in \CC: \operatorname{Re} s = 1\}$   The largest known zero-free region of the critical strip yields the following best known error bound to date.

\begin{theorem}[{Prime number theorem with error bound \cite{ford}}]\label{bestPNT}
One has
$$\li(x)  - \pi(x) = O\left(xe^{ - A(\log x)^{3/5}(\log \log x)^{-1/5}}\right) \ (x \to \infty),$$
where  $A  = 0.2098$.
\end{theorem}  

 Any such result  superceding de la Vall\'ee Poussin's  prime number theorem with error bound (Theorem \ref{PNTET1}) of 1899 is also called a {\bf prime number theorem with error bound}. \index{prime number theorem with error bound}  Note that the exponent  $1$ of $x$ appearing in the $O$ bound  is precisely the best known upper bound  of $\Theta$, and the improvement  over de la Vall\'ee Poussin's $O$ bound $O\left(xe^{-C\sqrt{\log x}}\right)$ is due to an enlargement of the zero-free region of $\zeta(s)$.   Such improvements are extraordinarily difficult to carry out,  and yet,  by some measures,  they have not brought us much closer to proving the Riemann hypothesis since 1899.   Thus,  it would seem that either newer and substantially stronger techniques are needed, or the Riemann hypothesis is false (or both, even, if $\frac{1}{2} < \Theta < 1$).

Nevertheless,  the situation is not hopeless, as,  subsequent to von Koch's 1901 result,  hundreds of other statements have been shown to be equivalent to the Riemann hypothesis, e.g., those collected in \cite{broughan} \cite{broughan2}.   Thus, even if the Riemann hypothesis is eventually proved false, the negation of hundreds of statements will immediately have been proved true.   For these and other reasons, the problem of settling the Riemann hypothesis is widely regarded as one of the most important, if not the most important, unsolved problems in mathematics today.   Absent a solution to Problem \ref{thetaprob},  the following research program is therefore warranted.

\begin{problem}
Given a known equivalent of the Riemann hypothesis, generalize the equivalence to an unconditional statement regarding the Riemann constant $\Theta$.
\end{problem}

Several widely known  examples of such unconditional statements, along with several new ones,  are discussed in Part 3.

 Let us assume,  for the moment, that the Riemann hypothesis is true.   It might appear, then,  based on values of $\pi(x)$ that have been computed or estimated,  that the conjectural error bound (\ref{vk})   can be improved to $\li(x)- \pi(x) = O\left( \frac{\sqrt{x}}{\log x} \right) \ (x \to \infty)$.   For example,  Hardy wrote in 1910 that ``there is reason to anticipate that'' this error bound holds \cite[p.\ 48]{har3}.   However, as Riesel noted in 1994:
\begin{quote}
Judging only from the values [given in a table] we might even try to estimate the order of magnitude of $\li(x)-\pi(x)$ and find it to be about $\sqrt{x}/\log x$.  However, {\it for large values of $x$, this is completely wrong!}'' \cite[p.\ 52]{ries}. 
\end{quote}
 Indeed, although the big $O$ estimate above is suggested by numerical data, it is known to be false because of  Littlewood's celebrated 1914 result,  Theorem \ref{litt} below,  which we call {\bf Littlewood's theorem}\index{Littlewood's theorem}  \cite{litt}  \cite[Theorem 35]{ing2} \cite[Theorem 15.11]{mont}  \cite[Theorem 6.20]{nark} \cite[Theorem 6.3]{plyman}.    One writes $f(x) = \Omega_{\pm}(g(x)) \ (x \to a)$
if $\limsup_{x \to a} \frac{f(x)}{|g(x)|} $ is positive and $\liminf_{x \to a} \frac{f(x)}{|g(x)|}$ is negative (both possibly infinite).

\begin{theorem}[{Littlewood's theorem \cite{litt}}]\label{litt}
One has
\begin{align*}
\li(x)-\pi(x) = \Omega_{\pm} \left(\frac{\sqrt{x}\, \log \log \log x}{\log x} \right) \ (x \to \infty).
\end{align*}
\end{theorem}

 As a consequence of Littlewood's  theorem, one has
\begin{align*}
\li(x)-\pi(x) \neq o \left(\frac{\sqrt{x}\, \log \log \log x}{\log x} \right) \ (x \to \infty),
\end{align*}
and therefore also
\begin{align*}
\li(x)-\pi(x) \neq O \left(\frac{\sqrt{x}}{\log x} \right) \ (x \to \infty).
\end{align*}
Littlewood's theorem is unconditional and provides the strongest known lower  bound on $\li(x)-\pi(x)$ to date.

\begin{figure}[h!]
\includegraphics[width=80mm]{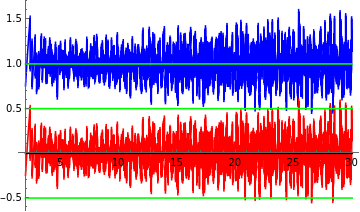}
    \caption{\centering Graphs of $\frac{\Ri(e^x)-\pi(e^x)}{\li(e^x)-\Ri(e^x)}$ (in red) and $\frac{\li(e^x)-\pi(e^x)}{\li(e^x)-\Ri(e^x)} = 1+\frac{\Ri(e^x)-\pi(e^x)}{\li(e^x)-\Ri(e^x)}$ (in blue) and $1$, $\frac{1}{2}$, and $-\frac{1}{2}$ (in green) on $[2,30]$}
\label{LiRi}
\end{figure}

\begin{figure}[h!]
\includegraphics[width=80mm]{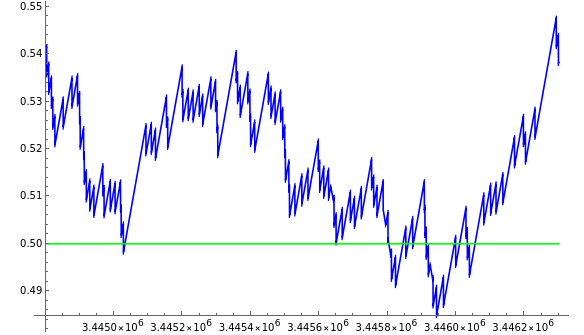}
    \caption{\centering Graph of $\frac{\li(x)-\pi(x)}{\li(x)-\Ri(x)}$ (in blue) and $\frac{1}{2}$ (in green) on $[3444800,3446300]$}
\label{LiRi2}
\end{figure}

Back in 1914, Littlewood's result was astonishing, not only because it was quite possibly the first occurrence of the function $\log \log \log x$ in relationship to the prime numbers, but also because it had the immediate consequence that
\begin{align*}
\li(x)-\pi(x) = \Omega_{\pm} (1) \ (x \to \infty)
\end{align*}
and therefore the function $\li(x)-\pi(x)$ changes sign an infinite number of times.  To this day, no one knows any specific value of $x \geq \mu$ for which $\li(x)-\pi(x) < 0$, where $\mu = 1.451369234883\ldots$ is the unique positive zero of $\li(x)$, known as the {\bf Ramanujan--Soldner constant}.\index{Ramanujan--Soldner constant $\mu$}\index[symbols]{.f tc@$\mu$}  (It is easy to see that $\li(x)-\pi(x) < 0$ for $0< x < \mu$.)     It is currently known that the infimum $\mathfrak{s}$ of all $x \geq \mu$ such that $\li(x)-\pi(x) < 0$, known as {\bf Skewes' number},\index{Skewes' number} is at least $10^{20}$ and at most $e^{727.9513}$.  For  $x \approx \mathfrak{s}$, or for any other number $x$ where $\li(x)-\pi(x)$ is approximately $0$, the function $\Ri(x)-\pi(x)$ is approximately equal to $\Ri(x)-\li(x) \approx -\frac{\sqrt{x}}{\log x} \ll 0$.  Thus, the  trend observed earlier in Figures \ref{RP} and \ref{PLR}  that $\Ri(x)$ is  a better approximation of $\pi(x)$ than $\li(x)$ is violated whenever $\li(x) - \pi(x) \approx 0$, or more generally, precisely when
$$\frac{\li(x)-\pi(x)}{\li(x)-\Ri(x)}  < \frac{1}{2},$$ that is,
$$\frac{\Ri(x)-\pi(x)}{\li(x)-\Ri(x)}  < -\frac{1}{2},$$
where the difference between the two functions above is equal to $1$, and where we assume also that $x > \nu$, where $\nu = 2.68945880489\ldots$ is the largest positive root of $\li(x)-\Ri(x)$. Figure \ref{LiRi} provides a graph of the two functions on a lin-log scale, that is, with $e^x$ substituted for $x$, on the interval $[2,30]$, where the first function is graphed in blue and the second in red.   By Littlewood's theorem, the functions in blue and red attain arbitrarily large positive and large negative values, even when further divided by $\log \log \log x$.   Skewes' number is the precise point at which the blue graph first dips down below the $x$-axis, i.e., when the red graph first dips down below the  line $y = -1$.  Moreover, $\li(e^x)$ is a better approximation of $\pi(e^x)$ than is $\Ri(e^x)$ precisely for those $x$ for which the blue graph lies below the line $y = \tfrac{1}{2}$, which is also where the red graph lies below the line $y = -\tfrac{1}{2}$.  One can see from the graph that, at least on the interval $[2,30]$, the function $\Ri(e^x)$ is the overall ``winner'' in this competition to  approximate $\pi(x)$.   A  search on Mathematica carried out by D.\ Stoll  \cite{stolll} reveals that the first interval for $x > \nu$ on which $\li(x)$ is a better approximation to $\pi(x)$ is $[3445027,3445031.758\ldots)$, where  $3445027 = p_{246588} = e^{15.052442\ldots}$ is a prime number with subsequent prime $p_{246589} = 3445093$, and therefore the function $\frac{\li(x)-\pi(x)}{\li(x)-\Ri(x)}$ is continuous and increasing on the interval $[3445027,3445093)$, attaining the value $0.498060\ldots$ at the left endpoint and approaching the value $0.524956\ldots$ at the right endpoint (from the left).  The next interval on which $\li(x)$ wins over $\Ri(x)$ is the rather small interval $[3445649,3445649.000498\ldots)$, where $3445649 = p_{246629} = e^{15.052622\ldots}$ is prime.   Figure \ref{LiRi2} provides a graph of  the function $\frac{\li(x)-\pi(x)}{\li(x)-\Ri(x)}$ on the interval $[3444800,3446300]$, which includes the first eight intervals, including the first two mentioned above, on which $\li(x)$ wins over $\Ri(x)$.

Unlike the number $3445027$, Skewes' number is extremely difficult to calculate, not only because it is very large with respect to our ability to compute $\pi(x)$, but also because of the following surprising result.  Let $X$ denote the set $$X = \left \{x\in [\log 2,\infty): \operatorname{Li}(e^x)\geq \pi(e^x)\right\},$$ where $$\operatorname{Li}(x) = \int_2^x \frac{1}{\log t} \, dt = \li(x)-\li(2) = \li(x)-1.045163780117\ldots.$$
It is known \cite{rubin}  that
$$\lim_{x \to \infty} \frac{1}{x-\log 2}\int_{t \in [\log 2, x] \cap X} \, dt = 0.99999973\ldots.$$
In other words, the {\it  asymptotic density} of the set $X$ of all $x\geq \log 2$ such that  $\li(e^x)-\li(2)\geq \pi(e^x)$ is equal to $0.99999973\ldots$.   Thus, counterexamples to $\li(x)-\li(2)\geq \pi(x)$, hence also counterexamples to  $\li(x)\geq \pi(x)$, are extremely rare.       A general lesson one learns from this is that numerical considerations sometimes do not carry much weight in analytic number theory,  especially when iterated logarithms might be lurking in the background.   

 See Table \ref{theoremproofs} for a list of some major  classical theorems in analytic number theory concerning the distribution of the primes, with their original sources, along with references to more recent proofs in the  literature.   The uninitiated reader seeking expertise in analytic number theory should consider undertaking a careful study of some of those proofs.

\begin{table}[!htbp]
  \caption{\centering Proofs of classical theorems concerning the distribution of the primes}
    \footnotesize
\begin{tabular}{|l|l|} \hline
Theorem  & Proofs \\ \hline\hline
Prime number theorem  & \cite[Chapters 4 and 13]{apos} \cite[Chapter 5]{bate}  \cite[Chapter 4]{edw} \\ 
 \cite{val1} \cite{had}  &  \cite[Chapters 4 and 5]{kon} \cite[Chapter 9]{gros2} \cite[Chapter II]{ing2}    \\ 
&  \cite[Chapter 6]{kno}  \cite{korevaar} \cite[Chapter 8]{mont} \cite[Chapter 5]{nark}  \\
 &  \cite[Chapter VII]{newman}  \cite[Chapter 5]{plyman}  \cite[Chapter III]{tit} \\
 &   \cite{zag}  \cite[Chapter 2]{zud} \\ \hline
Prime number theorems  &  \cite[Chapter 8]{bate} \cite[Section 3.5]{borg}  \cite{deb} \cite[Chapter 18]{dav} \\ 
 with error bound   \cite{val2}  &   \cite[Chapter 5]{edw}   \cite{ford}   \cite[Chapter III]{ing2} \cite[Chapter 12]{ivic}    \\
  & \cite[Chapter 8]{kou} \cite{liu3} \cite[Chapter 6]{mont}    \\
   &   \cite[Chapters 5 and 6]{nark} \cite[Chapter 6]{over}  \cite[Chapter 4]{patt}  \\ \hline
Riemann--von Mangoldt  &  \cite[Chapter 8]{bate} \cite[Section 3.5]{borg}  \cite[Chapter 17]{dav}   \\
  explicit formulas  &  \cite[Chapter 3]{edw}  \cite[Chapter IV]{ing2}  \cite[Chapter 12]{ivic}    \\ 
 \cite{rie} \cite{mang0} &  \cite[Chapters 5 and 8]{kou} \cite[Chapter 6]{nark} \cite[Chapter 10]{over} \\
 &  \cite[Chapter 3]{patt}  \\ \hline
  Littlewood's theorem   &  \cite[Chapter 11]{bate}  \cite[Theorem 4.13]{broughan} \cite[Theorem 35]{ing2}  \\
 \cite{litt} &  \cite[Theorem 15.11]{mont}  \cite[Theorem 6.20]{nark}  \cite[Chapter 6]{plyman}   \\ \hline
 Dirichlet's theorem on  & \cite[Chapter 7]{apos} \cite[Chapter 9]{bate} \cite[Section 3.4]{borg}  \\
 primes in arithmetic  & \cite[Chapter 22]{dav} \cite[Chapters 13 and 14]{kon}    \\
 progression \cite{diri} & \cite[Chapter 16]{ireland} \cite[Chapters 4 and 11]{mont} \\
 &  \cite[Chapter 2]{nark}  \cite[Chapter 7]{over} \cite[Chapter VI]{serre} \\
  &  \cite[Chapter 5]{zud} \\  \hline
\end{tabular}\label{theoremproofs}
\end{table}

\begin{remark}[The Riemann zeta function in physics]\label{physics}
The Riemann zeta function has several applications to physics, most notably, to {\it   zeta function regularization} in quantum field theory.  For example,  the value $\zeta(-3) = \frac{1}{120}$ is used in the derivation of the {\it Casimir effect}, and $\zeta(s)$ more broadly is used in the regularization of the energy–momentum tensor in curved spacetime, e.g., in the calculation of the vacuum expectation value of the energy of a particle field.  Additionally, the constant $\zeta(4) = \frac{\pi^4}{90}$ arises in the theory of black-body radiation, in the derivation of the {\it Stefan--Boltzmann law} from {\it Planck's law}, where $${\sigma ={\frac {2\pi ^{5}k^{4}}{15h^{3}c^{2}}}=5.670374\ldots \cdot 10^{-8}\,\mathrm {W\,m^{-2}\,K^{-4}}}$$ is the {\bf Stefan--Boltzmann constant}, and where $k = 1.380649\ldots \cdot 10^{-23} \operatorname{J}\cdot \operatorname{K}^{-1}$ is  the {\bf Boltzmann constant} and $c = 299792458 \, \operatorname{m}\cdot \operatorname{s}^{-1}$ is the speed of light in a vacuum.  An extensive list of appearances of the Riemann zeta function in physics can be found at http://empslocal.ex.ac.uk/people/staff/mrwatkin/zeta/physics1.htm (accessed by the author on 1 April 2024).
\end{remark}

\begin{figure}[ht!]
\includegraphics[width=70mm]{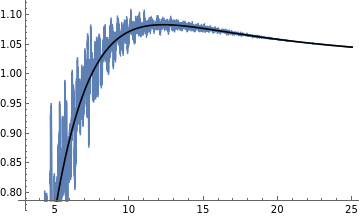}
\caption{\centering Graph of $x-\frac{e^x}{\pi(e^x)}$ and $x-\frac{e^x}{\Ri(e^x)}$ on $[3,25]$} \label{eureka0caaa}
\end{figure}

\begin{figure}[ht!]
\includegraphics[width=70mm]{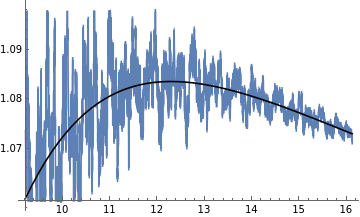}
\caption{\centering Graph of $x-\frac{e^x}{\Ri(e^x)}$ and $x-\frac{e^x}{\pi(e^x)}$  on $[\log(10^4), \log(10^7)]$}
 \label{primes15}
\end{figure}

\begin{figure}[ht!]
\includegraphics[width=70mm]{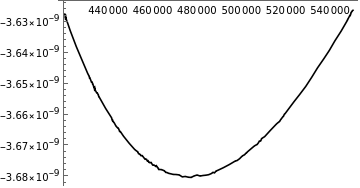}
\caption{\centering  Graph of $\frac{d}{dx}\left(\log x-\frac{x}{\Ri(x)}\right)$ on $[420000, 550000]$}
 \label{primes12}
\end{figure}

\begin{remark}[Legendre's constant]\label{legendre}
The first actual published statement of something close to the prime number theorem was made by Legendre in 1798 \cite[p. 19]{leg1}, which he refined further in 1808 \cite[pp.\ 394--398]{leg2}.   Legendre expressed his
1798  conjecture  as follows (English translation):
\begin{quote}
It is probable that the strict formula which gives the value
of $b$ [i.e., $\pi(a)$] when $a$ is very large is of the form $\frac{a}{A \log a +B}$,  $A$ and $B$ denoting constant coefficients and $\log a$ denoting a hyperbolic logarithm.  The exact determination of these coefficients would be a curious problem and worthy of exercising the astuteness of analysts.
\end{quote}
    In 1808, he wrote (English translation): 
\begin{quote}
Although the sequence of prime numbers is extremely irregular, one can however find, with a very satisfying precision, how many of these numbers there are from $1$ to a given limit $x$.   The formula that resolves this question is
$$y = \frac{x}{\log x - 1.08366}.$$
\ldots [Table of values of $\pi(x)$ given and each compared with $y$] \ldots It is impossible that a formula could represent a series of numbers of such a great extent,
and one subject to frequent anomalies, [completely] accurately.  There is therefore no doubt, not only
that the general law is represented by a function of the form $\frac{x}{A \log x +B}$, but that the coefficients $A$ and $B$ indeed have values very close to $A = 1\ldots$, $B = -1.08366$.
\end{quote}
Due to a lack of rigor by today's standards, it is difficult to give a faithful and precise interpretation of Legendre's conjectures.  

 Following Chebyshev and Gauss, we let $A$ denote the unique function $[2,\infty) \longrightarrow \RR$ such that
$$\pi(x) = \frac{x}{\log x - A(x)}\index[symbols]{.s BAb@$A(x)$} $$ 
for all $x \geq 2$, so that $$A(x) = \log x - \frac{1}{\PP(x)}.$$  A generous interpretation of Legendre's 1808 conjecture is that the limit $$L= \lim_{x \to \infty} A(x)$$ exists and is approximately equal to $1.08366$.    The limit $L$ is now referred to as {\bf Legendre's constant}.\index{Legendre's constant}  Clearly,  the existence of the limit $L$ implies the prime number theorem.   In 1848, Chebyshev proved \cite{cheb} that if Legendre's constant exists then it must equal $1$; more specifically, he proved that
$$\liminf_{x \to \infty} A(x) \leq 1 \leq \limsup_{x \to \infty} A(x).$$
In fact, the asymptotic relation (\ref{pias}) for $n = 1$, that is, the relation
$$\PP(x) -\frac{1}{\log x} \sim \frac{1}{(\log x)^2} \ (x \to \infty),$$
is equivalent to $$ A(x)  =\left(\PP(x) -\frac{1}{\log x}\right)\frac{\log x }{\PP(x)} \sim \frac{1}{(\log x)^2} (\log x)^2\sim 1 \ (x \to \infty)$$
and thus to the existence of $L$, and each of these equivalent statements  implies the prime number theorem.  

According to F.\ L.\ Bauer, ``Legendre's mistake can be explained easily: the largest tables of primes were those of J.\ H.\ Lambert (1770), going up to $10^5$, and of G.\ Vega (1796), going up to $4 \cdot 10^5$'' \cite{bau}.  Using Riemann's approximation $\Ri(x)$ to $\pi(x)$, we can provide a  more detailed explanation for Legendre's approximation $ 1.08366$  of the constant $L = 1$.  Figure \ref{eureka0caaa} compares  the functions $A(e^x) = x-\frac{e^x}{\pi(e^x)}$ and  $x-\frac{e^x}{\Ri(e^x)}$  on the interval $[3,25]$, and Figure \ref{primes15} compares the two functions on a smaller interval.  It is likely not a mere coincidence that the function  $\log x-\frac{x}{\Ri(x)}$  appears to attain a global maximum of approximately $1.08356$ at $x \approx 216811 \approx e^{12.2871}$, with a very small derivative nearby that appears to attain a local (and perhaps even global) minimum of only about $-3.68 \cdot 10^{-9}$ somewhat near the point $(4.75 \cdot 10^5, 1.0828)$.  See Figure \ref{primes12} for a graph of the derivative of $\log x-\frac{x}{\Ri(x)}$ near its apparent local minimum.  Since  $\log x-\frac{x}{\Ri(x)}$ approximates the ``center''  of the graph of  $A(x)$, the misleading features of the function $\log x-\frac{x}{\Ri(x)}$ described above  would explain why Legendre arrived at something close to $1.08356$, namely, his approximation $1.08366$, as an approximation of the constant $L$.    Even at $x = 10^{29}$ one has $A(x) =  1.015696108866115\ldots$ and $\log x-\frac{x}{\Ri(x)} = 1.015696108866137\ldots$,  thus indicating that convergence  of $A(x)$ to the limit $1$ is  slow.
Regarding this very matter, in his 1849 letter to Encke \cite[pp.\ 444--447]{gau}, Gauss made some prescient remarks, based on his extensive calculations  (English translation):
\begin{quote}
It appears that, with increasing $n$, the (average) value of $A$ decreases; however, I dare not conjecture whether the limit as $n$ approaches infinity is $1$ or a number different from 1.  I cannot say that there is any justification for expecting a very simple limiting value; on the other hand, the excess of $A$ over $1$ might well be a quantity of the order of $\frac{1}{\log n}$.
\end{quote}
Gauss' speculations  turned out to be correct, in that
\begin{align}\label{aass}
A(x) -1 \sim \frac{1}{\log x} \ (x \to \infty),
\end{align} 
and this explains why the convergence of $A(x)$ to Legendre's constant is so slow.   Indeed, by the prime number theorem and the asymptotic relation (\ref{pias}) for $n = 3$, one has
\begin{align*}
A(x) -1  & =\left(\left(\log x-1 \right)\PP(x)-1\right)\frac{1}{\PP(x)} \\
& =\left(\left(\log x-1  \right)\left(\frac{1}{\log x} + \frac{1}{(\log x)^2}  +\frac{2}{(\log x)^3}+ O\left( \frac{1}{(\log x)^4} \right)\right)-1\right)\frac{1}{\PP(x)} \\
& =\left( \frac{1}{(\log x)^2} + O \left(\frac{1}{(\log x)^3} \right) \right)\frac{1}{\PP(x)}  \\
 & =  \frac{1}{\log x} + O\left( \frac{1}{(\log x)^2} \right) \ (x \to \infty),
\end{align*} 
which implies (\ref{aass}).  
\end{remark}

\chapter{Asymptotic analysis}

In this chapter, we discuss {\it asymptotic relations} and {\it asymptotic expansions}.  We  introduce the notion of  the {\it degree} of a real function (which we study in much greater detail in Section 6.1),  and we apply it to the study of {\it slowly varying} functions and {\it regularly varying} functions, two natural notions that have applications to analysis, probability theory, and analytic number theory.   We also discuss the {\it Euler--Maclaurin formula}, {\it Karamata's integral representation theorem}, and {\it Karamata's insttegral theorem}.

\section{Asymptotic relations}

Let $a \in \overline{\RR}$, and let $f, g \in \RR^{\RR_a}$ (where $\RR^{\RR_a}$ denotes the set of all real functions $f$ such that $a$ is a limit point of $\dom f$).   Assume also that $\dom g$ contains the intersection of $\dom f$ with some punctured neighborhood of $a$.
One defines the following.
\begin{enumerate}
\item $f(x) = O( g(x)) \ (x \to a)$, also written $f(x) \ll g(x) \ (x \to a)$, if for some $M > 0$ one has $|f(x)| \leq M|g(x)|$ for all $x$ in the intersection of $\dom f$ with some punctured neighborhood of $a$.\index[symbols]{.d aa@$O$}\index[symbols]{.e bb@$\ll$}
\item $f(x) \gg g(x) \ (x \to a)$ if for some $M > 0$ one has $|f(x)| \geq M|g(x)|$ for all $x$ in the intersection of $\dom f$ with some punctured neighborhood of $a$.\index[symbols]{.e  cc@$\gg$}
\item $f(x) \asymp g(x) \ (x \to a)$ if $f(x) \ll g(x) \ (x \to a)$ and $f(x) \gg g(x) \ (x \to a)$.\index[symbols]{.e ee@$\asymp$}
\item $f(x) = o( g(x)) \ (x \to a)$ if  for all  $M > 0$ one has $|f(x)| \leq M|g(x)|$ for all $x$ in the intersection of $\dom f$ with some punctured neighborhood of $a$.\index[symbols]{.d cc@$o$}
\item $f(x) \sim g(x) \ (x \to a)$ if  $f(x)-g(x) = o(g(x)) \ (x \to a)$.\index[symbols]{.e hh@$\sim$}
\item $f(x) = \Omega_{+}(g(x)) \ (x \to a)$ if  there exists an $M > 0$ such that for every $x$ in a punctured neighborhood of $a$ there exists a $y \neq a$ closer to $a$ than $x$ such that $f(y)> M| g(y)|$.\index[symbols]{.e ii@$\Omega_+$}
\item $f(x) = \Omega_{-}(g(x)) \ (x \to a)$ if  there exists an $M > 0$ such that for every $x$ in a punctured neighborhood of $a$ there exists a $y \neq a$ closer to $a$ than $x$ such that $f(y) > -M |g(y)|$.\index[symbols]{.e jj@$\Omega_{-}$}
\item $f(x) = \Omega_{\pm}(g(x)) \ (x \to a)$ if $f(x) = \Omega_{+}(g(x)) \ (x \to a)$ and $f(x) = \Omega_{-}(g(x)) \ (x \to a)$.\index[symbols]{.e kk@$\Omega_{\pm}$}
\end{enumerate}

Note that both of the conditions $f(x) = o( g(x)) \ (x \to a)$ and $f(x)  \sim g(x) \ (x \to a)$ are stronger than the condition $f(x) = O(g(x)) \ (x \to a)$, and all three conditions require that $\dom g$ contains the intersection of $\dom f$ with some punctured neighborhood of $a$.   The relations $\sim$, $O$, and $o$ are transitive, in the obvious sense.  Moreover, if $f(x) = O( g(x)) \ (x \to a)$ and $g(x) = o(h(x)) \ (x \to a)$, or if  $f(x) = o( g(x)) \ (x \to a)$ and $g(x) = O(h(x)) \ (x \to a)$,  then $f(x) = o( h(x)) \ (x \to a)$.  The relations $\asymp$ and $\sim$ are  symmetric on functions $f, g \in \RR^{\RR_a}$ such that $U \cap \dom f = U \cap \dom g$ for some punctured neighborhood $U$ of $a$.  Moreover, a function $f \in  \RR^{\RR_a}$ satisfies  $f(x) \sim 0 \ (x \to a)$ if and only if $f$ is zero on the intersection of $\dom f$ with some  punctured neighborhood of $a$.   Note, also, that,  if $f(x) = \Omega_{+}(g(x)) \ (x \to a)$ or if $f(x) = \Omega_{-}(g(x)) \ (x \to a)$, then  $f(x) \neq  o(g(x)) \ (x \to a)$.

As above, let $f, g \in \RR^{\RR_a}$, where $\dom g$ contains the intersection of $\dom f$ with some punctured neighborhood of $a$.  If also $g$ is nonzero on the intersection of $\dom f$ with some  punctured neighborhood of $a$, then the relations (1)--(8) defined above simplify as in the following proposition.

\begin{proposition}
Let $f, g \in \RR^{\RR_a}$, where $g$ is defined and nonzero on the intersection of $\dom f$ with some punctured neighborhood of $a$.   One has the following.
\begin{enumerate}
\item $f(x) = O( g(x)) \ (x \to a)$  if and only if  
$$\limsup_{x \to a} \left|\frac{f(x)}{g(x)}\right| < \infty.$$
\item $f(x) \gg g(x) \ (x \to a)$  if and only if  
$$\liminf_{x \to a} \left|\frac{f(x)}{g(x)}\right| > 0.$$
\item $f(x) \asymp g(x) \ (x \to a)$ if and only if 
$$\limsup_{x \to a} \left|\frac{f(x)}{g(x)}\right| < \infty \text{ and } \liminf_{x \to a} \left|\frac{f(x)}{g(x)}\right| > 0.$$
\item $f(x) = o( g(x)) \ (x \to a)$  if and only if
$$\lim_{x \to a} \frac{f(x)}{g(x)} = 0,$$
if and only if 
$$\limsup_{x \to a} \left| \frac{f(x)}{g(x)} \right| = 0.$$
\item $f(x) \sim g(x) \ (x \to a)$  if and only if
$$\lim_{x \to a} \frac{f(x)}{g(x)} = 1.$$
\item $f(x) = \Omega_{+}(g(x)) \ (x \to a)$ if  and only if 
$$\limsup_{x \to a} \frac{f(x)}{|g(x)|}> 0.$$
\item $f(x) = \Omega_{-}(g(x)) \ (x \to a)$ if  and only if 
$$\liminf_{x \to a} \frac{f(x)}{|g(x)|}< 0.$$
\item $f(x) = \Omega_{\pm}(g(x)) \ (x \to a)$  if and only if
$$\limsup_{x \to a} \frac{f(x)}{|g(x)|}> 0 \text{ and } \liminf_{x \to a} \frac{f(x)}{|g(x)|}< 0 .$$
\end{enumerate}
\end{proposition}

When we write $f(x) \neq O( g(x)) \ (x \to a)$, or $f(x)  \not \asymp g(x) \ (x \to a)$,  etc., for any of the relations above, we mean not only that the negation of the given relation holds, but also that $\dom g$ contains the intersection of $\dom f$ with some punctured neighborhood of $a$.    Thus, for example, if $f$ is defined on all of $\RR$, then $f|_\ZZ(x) = O(f(x)) \ (x \to \infty)$, but it is neither the case that $f(x) = O(f|_\ZZ(x)) \ (x \to \infty)$ nor that $f(x) \neq O(f|_\ZZ(x)) \ (x \to \infty)$.   We follow the convention that $O(g(x))$ (resp., $o(g(x))$) denotes some function $f$ such that $f(x) = O(g(x) ) \ (x \to a)$ (resp., $f(x) = o(g(x) ) \ (x \to a)$), assuming that $a$ is understood by context.   Our domain conventions can be summarized by stipulating that a given asymptotic relation must make sense for all $x \neq a$ near $a$ in the domain of the function appearing on the {\it left} of the given relation.  The rationale for this convention is that it is natural to interpret any asymptotic relation  as expressing a property of the function appearing on the left.    Thus, for example, when we write $f(x) = h(x)+ O( g(x)) \ (x \to a)$,  we mean  that $f(x) -h(x) =  O( g(x)) \ (x \to a)$ and that both $\dom g$ and $\dom h$ contain the intersection of $\dom f$ with some punctured neighborhood of $a$.     In a similar vein,  the relations $f(x) \ll g(x) \ (x \to a)$ and $g(x) \gg f(x) \ (x \to a)$ are equivalent if and only if $U \cap \dom f = U \cap \dom g$ for some punctured neighborhood $U$ of $a$, and likewise for the relations $f(x) \asymp g(x) \ (x \to a)$ and $g(x) \asymp f(x) \ (x \to a)$ and for the relations $f(x) \sim g(x) \ (x \to a)$ and $g(x) \sim f(x) \ (x \to a)$.  Of course, if all functions involved are defined on some punctured neighborhood of $a$,  or if they all have the same domain, then these domain restrictions are automatic.

\begin{example} Let $f,g \in \RR^{\RR_\infty}$ and $t, u, A\in \RR$ with $A \neq 0$. One has the following.
\begin{enumerate}
\item $f(x) = O(1) \ (x \to \infty)$ if and only if $f(x)$ is eventually bounded.
\item $f(x) = o(1) \ (x \to \infty)$ if and only if $\lim_{x \to \infty} f(x)= 0$.
\item  $\sin x = O(1) \ (x \to \infty)$.
\item $\sin x \neq o(1) \ (x \to \infty)$.
\item $\sin x =  \Omega_\pm(1) \ (x \to \infty)$.
\item $x^t  = O(x^u) \ (x \to \infty)$ if and only if $t \leq u$.
\item $x^t  = o(x^u) \ (x \to \infty)$ if and only if $t < u$.
\item $x^t  = o(e^{x}) \ (x \to \infty)$ for all $t \in \RR$.
\item $\log x = o(x^t) \ (x \to \infty)$ if and only if $t > 0$, if and only if $\log x = O(x^t) \ (x \to \infty)$.
\item  If $f$ and $g$ are polynomials, then $f(x) = O(g(x)) \ (x \to \infty)$ if and only if $\deg f \leq \deg g$.
\item If $f$ and $g$ are polynomials, then $f(x) = o(g(x)) \ (x \to \infty)$ if and only if $\deg f < \deg g$.
\item  If $f$ and $g$ are polynomials, then $f(x) \sim g(x) \ (x \to \infty)$ if and only if $f$ and $g$ have the same leading term.
\item $f(x) \sim A x^n \ (x \to \infty)$ if $f$ is a polynomial of degree $n$ with leading term $A x^n$.
\item $f(x) \sim A \ (x \to \infty)$ if and only if $\lim_{x \to \infty} f(x) = A$.
\end{enumerate}
\end{example}

We note the following properties of the $O$ and $o$ relations.

\begin{proposition} Let $a \in \overline{\RR}$, let $f_1, f_2, g_1, g_2, g \in \RR^{\RR_a}$ with $a$ a limit point of $\dom f_1 \cap \dom f_2$, and let $r_1, r_2 \in \RR$.  One has the following.
\begin{enumerate}
\item If $$ f_{1}(x)=O(g_{1}(x)) \ (x \to a) {\text{ and }} f_{2}(x)=O(g_{2}(x)) \ (x \to a),$$ then 
$$ f_{1}(x)+f_{2}(x)=O(\max(|g_{1}(x)|,|g_{2}(x)|)) \ (x \to a)$$ and
$$f_{1}(x)f_{2}(x)=O(g_{1}(x)g_{2}(x)) \ (x \to a).$$
\item If  
$$ f_{1}(x)=O(g(x)) \ (x \to a) \text{ and }  f_{2}(x)=O(g(x)) \ (x \to a),$$
then  $$r_1 f_1(x) + r_2 f_2(x)  = O(g(x)) \ (x \to a).$$
\item If $$ f_{1}(x)=o(g_{1}(x)) \ (x \to a) {\text{ and }} f_{2}(x)=o(g_{2}(x)) \ (x \to a),$$ then 
$$ f_{1}(x)+f_{2}(x)=o(\max(|g_{1}(x)|,|g_{2}(x)|)) \ (x \to a)$$ and
$$f_{1}(x)f_{2}(x)=o(g_{1}(x)g_{2}(x)) \ (x \to a).$$
\item If  
$$ f_{1}(x)=o(g(x)) \ (x \to a) \text{ and }  f_{2}(x)=o(g(x)) \ (x \to a),$$
then  $$r_1 f_1(x) + r_2 f_2(x)  = o(g(x)) \ (x \to a).$$
\end{enumerate}
\end{proposition}

In this book, we are interested primarily in the case where $a = \infty$.  However, we also require generalizations of the $O$, $o$, $\asymp$, and $\sim$ relations  to complex functions, where we then assume that $a \in \CC \cup \{\infty\}$.  To obtain these definitions, one  replaces the real absolute value in the  definitions  above with the complex absolute value, and one replaces punctured real neighborhoods with punctured complex neighborhoods, and then appropriate analogues of the results above hold for complex functions defined on a subset of $\CC$ containing $a$ as a limit point.

\section{Asymptotic expansions}

Let $a \in \overline{\RR}$ (resp, $a \in \CC \cup \{\infty\}$).  An  {\bf asymptotic sequence  at $a$}\index{asymptotic sequence} is a sequence $\{\varphi_n\}_{n = 1}^\infty$  of real (resp., complex) functions $\varphi_n$ such that $$\varphi_{n+1}(x) = o(\varphi_n(x)) \ (x \to a)$$ for all $n \geq 1$ (so $a$ is a limit point of $\dom \varphi_n$ for each $n$).  Let $\{\varphi_n\}$ be an asymptotic sequence at $a$, let $f$ be a real (resp., complex) function, and let $\{a_n\}$ be a sequence of real (resp., complex)  numbers.  The function $f$ is said to have the {\bf asymptotic expansion}\index{asymptotic expansion}\index[symbols]{.f s@$\simeq$}
\begin{align}\label{asympa}
f(x) \simeq \sum_{n=1}^{\infty} a_n \varphi_{n}(x) \ (x \to a)
\end{align}
{\bf (at $a$ with respect to $\{\varphi_n\}$)} if 
\begin{align}\label{asympb}
f(x) = \sum_{k=1}^{n} a_k \varphi_{k}(x) + o(\varphi_{n}(x)) \  (x \to a)
\end{align}
for all positive integers $n$ (which requires that $\varphi_n$ for each positive integer $n$ be defined on the intersection of $\dom f$ with some punctured neighborhood of $a$).

\begin{remark}\label{asrem} \
\begin{enumerate}
\item It is clear that, for a given positive integer $n$, condition (\ref{asympb}) implies that
\begin{align}\label{asymp3}
f(x) = \sum_{k=1}^{n-1} a_k \varphi_{k}(x) + O(\varphi_{n}(x)) \  (x \to a),
\end{align}
which in turn implies that $f(x) = \sum_{k=1}^{n-1} a_k \varphi_{k}(x) + o(\varphi_{n-1}(x)) \  (x \to a)$
if $n \geq 2$.   Consequently, the asymptotic expansion (\ref{asympa}) holds if and only if either (\ref{asympb}) or (\ref{asymp3}) holds for infinitely many integers $n \geq 1$, if and only if (\ref{asymp3}) holds for all $n \geq 1$.
\item  For any positive integer $n$,  if $a_n \neq 0$  and $\varphi_k$ for each $k < n$ is defined on the intersection of $\dom f$ with some punctured neighborhood of $a$, then (\ref{asympb}) is equivalent to
$$f(x) - \sum_{k=1}^{n-1}a_{k}\varphi_{k}(x) \sim a_n \varphi_n(x) \ (x \to a).$$
\end{enumerate}
\end{remark}

\begin{example} \
\begin{enumerate}
\item  Important examples of asymptotic expansions over $\RR$ and $\CC$ follow from Taylor's theorem: if $f$ is real or complex function defined in a neighborhood of some number $a$, then one has an asymptotic expansion  $f(x) \simeq \sum_{n = 0}^\infty a_n (x-a)^n \ (x \to a)$
of $f$  at $a$ with respect to the asymptotic sequence $\{(x-a)^n\}$ if $f$ is infinitely differentiable at $a$, in which case  $a_n = \frac{f^{(n)}(a)}{n!}$ for all $n$.   This also applies  to $f$ over $\CC$ at  $\infty$ (resp., over $\RR$ at $\infty$, over $\RR$ at $-\infty$) by considering the function $f(\tfrac{1}{x})$ with respect to the asymptotic sequence $\left\{ \frac{1}{x^k}\right\}$ at $0$ (resp., at $0^+$, at $0^-$).   
\item  One can generalize the notion of an asymptotic expansion ({\it of infinite order}) to the notion of an {\it asymptotic expansion  of finite order $N$}, and  then $f$  has an asymptotic expansion $f(x) \simeq \sum_{n = 0}^N a_n (x-a)^n \ (x \to a)$ of order $N+1$ at $a$  with respect to the asymptotic sequence $\{(x-a)^n\}$ if $f$ is $N$-times differentiable at $a$,  in which case  $a_n = \frac{f^{(n)}(a)}{n!}$ for all $n \leq N$.
\item Two asymptotic expansions at $a$ with respect to  $\{(x-a)^n\}$ or at  $\infty$ with respect to $\left\{ \frac{1}{x^n}\right\}$ can be added, subtracted, multiplied, divided, and composed just like formal power series.  
\end{enumerate}
\end{example}

Excellent references on the theory of asymptotic expansions include \cite{erd} \cite{est} \cite{olver} \cite{wong}.
 
The following result provides some necessary and sufficient conditions for two functions to have the same asymptotic expansion with respect to a given asymptotic sequence.

\begin{proposition}[{\cite[Lemma 2.4]{ell0}}]\label{asympprop}
Let $a \in \CC \cup \{\infty\}$, let $\{\varphi_n\}$ be an asymptotic sequence at $a$, and let $f$ and $g$ be complex functions such that $\dom f$ contains the intersection of $\dom g$ with some punctured neighborhood of $a$.  A given asymptotic expansion of $f$  at $a$ with respect to $\{\varphi_n\}$ is also an asymptotic expansion of  $g$ at $a$ with respect to $\{\varphi_n\}$  if and only if $$g(x) = f(x) + o(\varphi_n(x)) \ (x \to a)$$
for all (or, equivalently, for infinitely many) positive integers $n$,
if only if $$g(x) = f(x) + O(\varphi_n(x)) \ (x \to a)$$
for all (or, equivalently, for infinitely many) positive integers $n$.
\end{proposition}

\begin{proof}
Suppose that $f(x) \simeq \sum_{n = 1}^\infty a_n \varphi_n(x) \ (x \to a)$ is an asymptotic expansion of $f$, or, equivalently, that $f(x)=  \sum_{k = 1}^n a_k \varphi_k(x) + o(\varphi_n(x)) \ (x \to a)$
for all positive integers $n$.   Let $n$ be any positive integer. If one has
$g(x) =   \sum_{k = 1}^n a_k \varphi_k(x) + o(\varphi_n(x)) \ (x \to a),$
then subtracting we see that
$g(x)- f(x) = o(\varphi_n(x)) \ (x \to a)$ on $\dom g$, and therefore $g(x)= f(x) + o(\varphi_n(x)) \ (x \to a)$.  The converse is also clear.  Similar statements hold for the $O$ relation. The proposition follows.
\end{proof}

The following result provides a natural condition under which two asymptotic sequences can be viewed as  {\it asymptotically equivalent}.

\begin{proposition}[{\cite[Lemma 2.5]{ell0}}]\label{asympprop2a}
Let $a \in \CC \cup \{\infty\}$, let $\{\varphi_n\}$ be an asymptotic sequence at $a$, and let $\{\psi_n\}$ be  a sequence of complex functions such that, for all positive integers $n$, one has $\psi_n(x) -\varphi_n(x) = o(\varphi_N(x)) \ (x \to a)$ for all $N \geq n$ and $U_n \cap \dom \varphi_n   = U_n \cap \dom \psi_n$ for some punctured neighborhood $U_n$ of $a$.   Then $\{\psi_n\}$ is an asymptotic sequence at $a$ with  $\varphi_n(x) -\psi_n(x) = o(\psi_N(x)) \ (x \to a)$ for all  positive integers $n$ and $N \geq n$.  Moreover, any asymptotic expansion 
$$f(x) \simeq \sum_{n = 1}^\infty a_n \varphi_n(x) \ (x \to a)$$
of a complex-valued function $f$ at $a$ with respect to $\{\varphi_n\}$ is equivalent to the asymptotic expansion $$f(x) \simeq \sum_{n = 1}^\infty a_n \psi_n(x) \ (x \to a)$$
of $f$ at $a$ with respect to $\{\psi_n\}$.
\end{proposition}

\begin{proof}
Let $N$ be a positive integer.
For all $n$ one has $\psi_n(x) -\varphi_n(x) = o(\varphi_n(x)) \ (x \to a)$ and therefore $\psi_n(x) \sim \varphi_n(x) \ (x \to a).$   It follows that $\varphi_n(x) -\psi_n(x) = o(\varphi_N(x)) = o(\psi_N(x)) \ (x \to a)$ for all $n \leq N$ and that $\{\psi_n\}$ is an asymptotic sequence  at $a$.  If the first asymptotic expansion of $f$ holds, then one has
\begin{align*} f(x) - \sum_{n = 1}^N a_n \psi_n(x)  & = \left( f(x) - \sum_{n = 1}^N a_n \varphi_n(x) \right) +  \sum_{n = 1}^N a_n (\psi_n(x)-\varphi_n(x)) \\
  &=  o( \varphi_N(x))  \ (x \to a) \\ 
  &=  o( \psi_N(x))  \ (x \to a),
\end{align*}
for all $N$, and therefore second asymptotic expansion of $f$ also holds.  By symmetry, the reverse implication holds as well.  The proposition follows.
\end{proof}

In 1848  \cite[p.\ 153]{cheb}, Chebyshev noted the asymptotic expansion
\begin{align}\label{asex}
\frac{\li(x)}{x} \simeq \sum_{n = 0}^\infty {\frac {n!}{(\log x)^{n+1}}} \ (x \to \infty)
\end{align}
of the function $\frac{\li(x)}{x}$ with respect to the asymptotic sequence $\left\{\frac{1}{(\log x)^{n+1}}\right\}$.  This asymptotic expansion is now well known (see \cite[Section 10.3]{stop}, for example) and follows from  the fact that  
\begin{align}\label{lias}
\li(x) - \sum_{k = 0}^{n-1} \frac{k!x}{(\log x)^{k+1}} = \int_e^x \frac{n! \, dt} {(\log t)^{n+1}} +C_n \sim \frac{n!x}{(\log x)^{n+1}} \ (x \to \infty)
\end{align}
for all nonnegative integers $n$, where $C_n$ is a constant  and the given equality is proved by repeated integration by parts.  (Since $\frac{d}{dx} \frac{x}{(\log x)^{n+1}} = \frac{\log x -n-1}{(\log x)^{n+2}} \sim \frac{1} {(\log x)^{n+1}} \ (x \to \infty)$, the  asymptotic relation in  (\ref{lias}) above follows from L'H\^opital's rule.)   From (\ref{PNTET}), (\ref{asex}), and Proposition \ref{asympprop}, we see that $\PP(x)$ has the same asymptotic expansion as $\frac{\li(x)}{x}$, namely,
\begin{align}\label{asex2}
\PP(x) \simeq \sum_{n = 0}^\infty {\frac {n!}{(\log x)^{n+1}}} \ (x \to \infty).
\end{align}
In fact, assuming (\ref{asex}), it follows easily   that (\ref{PNTET})  and (\ref{asex2}) are equivalent.  Thus, the asymptotic expansion (\ref{asex2}) carries essentially the same information as the weak version (\ref{PNTET}) of the prime number theorem with error bound.  Note that the series $\sum_{k = 0}^\infty {\frac {k!}{(\log x)^{k+1}}}$ is divergent for all $x$, and the definition of asymptotic expansions equates (\ref{asex2}) with the statement
\begin{align*}
\PP(x) - \sum_{k = 0}^{n-1} {\frac {k!}{(\log x)^{k+1}}} \sim \frac{n!}{(\log x)^{n+1}} \ (x \to \infty), \quad \forall  n \geq 1.
\end{align*}

It is noteworthy that one can ``rationalize'' the asymptotic expansion (\ref{asex2}) by substituting $e^x$ for $x$: one has
\begin{align*}
\PP(e^x) - \sum_{k = 0}^{n-1} {\frac {k!}{x^{k+1}}} \sim \frac{n!}{x^{n+1}} \ (x \to \infty), \quad \forall  n \geq 1,
\end{align*}
and therefore the function $\PP(e^x)$ has the asymptotic expansion
\begin{align*}
\PP(e^x) \simeq \sum_{n= 0}^\infty {\frac {n!}{x^{n+1}}} \ (x \to \infty)
\end{align*}
with respect to the asymptotic sequence $\left\{ \frac{1}{x^{n+1}} \right\}$.  
For this and many other reasons,  we consider the function $\PP(e^x) = \frac{\pi(e^x)}{e^x} \sim \frac{1}{x} \ (x \to \infty)$ to be a ``rationalized''  version of the function $\pi(x)$.

\begin{example}[{\cite{ell0}}]
Let $n$ be a nonnegative integer.   The number $D_n = n!\sum_{k = 0}^n \frac{(-1)^k}{k!}$ is equal to the number of {\it derangements} of any $n$-element set, and the number $A_n = n! \sum_{k = 0}^n \frac{1}{k!}$ is equal to the number of {\it arrangements} of any $n$-element set.  The sequence $D_n$ is $1,0,1,2,9,44,265,1854,\ldots$, and the sequence  $A_n$ is $1,2,5,16,65,326,1957,13700,\ldots$, and they are Sequences A000166 and A000522, respectively,  of the On-Line Encyclopedia of Integer Sequences (OEIS).  Using the binomial theorem, one can show that the following asymptotic expansions are equivalent to (\ref{asex2}):
\begin{enumerate} 
\item $\displaystyle\PP(x) \simeq \sum_{n=0}^\infty \frac{D_{n}}{(\log x-1 )^{n+1}} \ (x \to \infty)$.
\item $\displaystyle\PP(ex) \simeq \sum_{n=0}^\infty \frac{D_{n}}{(\log x )^{n+1}} \ (x \to \infty)$.
\item $\displaystyle\PP(x) \simeq \sum_{n=0}^\infty \frac{A_{n}}{(\log x+1)^{n+1}} \ (x \to \infty)$.
\item $\displaystyle\PP(x/e) \simeq \sum_{n=0}^\infty \frac{A_{n}}{(\log x )^{n+1}} \ (x \to \infty)$.
\end{enumerate}
In particular, one has the asymptotic expansion
$$\PP(x/e) \simeq \frac{1}{\log x} +  \frac{2}{(\log x)^2} + \frac{5}{(\log x)^3} + \frac{16}{(\log x)^4} + \frac{65}{(\log x)^5}  + \cdots \ (x \to \infty).$$
At the same time, squaring the asymptotic expansion (\ref{asex2}) of $\PP(x)$ yields
$$\PP(x)^2 \simeq \frac{1}{(\log x)^2} +  \frac{2}{(\log x)^3} + \frac{5}{(\log x)^4} + \frac{16}{(\log x)^5} + \frac{64}{(\log x)^6}  + \cdots \ (x \to \infty),$$
It follows that $$ \PP(x/e) - \PP(x)^2 \log x  \sim \frac{1}{(\log x)^5} \ (x \to \infty).$$
Consequently, one has $\PP(x)^2 \log x  < \PP(x/e)$ for all sufficiently large $x$, which is equivalent to Ramanujan's famous inequality
$$\pi(x)^2 < \frac{ex}{\log x} \pi(x/e), \quad  \forall x \gg 0.$$
It is known \cite{axler} that the Riemann hypothesis implies that the smallest integer $N$ such that Ramanujan's  inequality holds for all $x \geq N$ is equal to $38358837683$, and, moreover, the inequality holds unconditionally for all $x \geq e^{9032}$.

For any sequence $a_n$ of complex numbers and any $z \in \CC$, the sequence $b_n = \sum_{k = 0}^n {n \choose k} a_k z^{n-k}$ is called the {\bf $z$-binomial transform of $a_n$}. \index{binomial transform} The $0$-binomial transform of the sequence $n!$ is the sequence $n!$, which is equal to the number of permutations of any $n$-element set.  The sequence $D_n$ is the $(-1)$-binomial transform of the sequence $n!$, and  sequence $A_n$ is the $1$-binomial transform of the sequence $n!$. For every nonnegative integer $n$,  let $r_n(X)$ denote the monic integer polynomial
$$r_n(X) = n! \sum_{k = 0}^n \frac{X^k}{k!}  =  \sum_{k = 0}^n \frac{n!}{k!} {X^k}\in \ZZ[X].$$ 
For any $z \in \CC$, the sequence $r_n(z)$ is the $z$-binomial transform of the sequence $n!$, and one has $r_n(0) = n!$, $r_n(-1) = D_{n}$, and $r_n(1) = A_{n}$ for all $n$, so the family of sequences $r_n(z)$ interpolates those three sequences.   
For all $t \in \RR$, one has the asymptotic expansion
$$\PP(e^{-t }x) \simeq \sum_{n=0}^\infty \frac{r_n(t)}{(\log x)^{n+1}} \ (x \to \infty).$$
Note that, since $r_n(X) \in \ZZ[X]$, one has $r_n(k) \in \ZZ$ for all $n$ and all $k \in \ZZ$.  The integer sequence $r_n(2)$ is OEIS Sequence A010842, and thus, for example, $r_n(2)$ is the number of ways to split the set $\{1,2,\ldots,n\}$ into two disjoint subsets $S$ and $T$ and linearly order $S$ and then choose a subset of $T$.   Also, the integer sequence $r_n(-2)$ is OEIS Sequence A000023.   Note also that $r_n(z) \sim   n! e^z =  r_n(0)e^z  \ (n \to \infty)$ for all $z \in \CC$.  Thus, for example, one has $D_n  \sim \,  n! e^{-1} \ (n \to \infty)$ and $A_n   \sim \,  n!e \ (n \to \infty),$ which are well-known asymptotics for the sequences $D_n$ and $A_n$.  
\end{example}

\begin{example}[{\cite{pan}}]
Another interesting reformulation of the asymptotic expansion (\ref{asex2}) is  a result  proved by Panaitopol in 2000: reciprocating the asymptotic expansion of $\PP(x)$ yields the asymptotic expansion
\begin{align*}
A(x) = \log x- \frac{1}{\PP(x)}\simeq \sum_{n = 0}^\infty {\frac {k_n}{(\log x)^{n}}} \ (x \to \infty),
\end{align*}
where $\{k_n\}$ is the sequence with generating function $$\sum_{n = 0}^\infty k_n X^n = 
\frac{1}{X}-\frac{1}{X\sum_{n = 0}^\infty n!X^n},$$
and where $A(x)$ is the unique function defined on $[2,\infty)$ such that
$$\pi(x)  = \frac{x}{\log x - A(x)}, \quad \forall x \geq 2.$$
 It follows that  $\{k_n\}$ is OEIS Sequence A233824, and therefore $k_n$ for any $n \geq 0$ is the number of subgroups of index $n$ of the free group on two generators, and the sequence  $\{k_n\}$ has its first several terms given by  $0, 1, 3, 13, 71, 461, 3447, 29093, \ldots$.  It  also follows that  $k_n = I_{n+1}$ for all $n$,  where $\{I_n\}$ is OEIS Sequence A003319 and $I_n$ for any nonnegative integer $n$ is equal to the number of indecomposable permutations of $\{1,2,3\ldots, n\}$, where a permutation of $\{1,2,3,\ldots,n\}$ is said to be {\bf indecomposable}\index{indecomposable permutation} if it does not fix $\{1,2,3,\ldots,j\}$ for any  $1 \leq j < n$.
\end{example}

\begin{example}[{\cite{rey}}]\label{pnas}
 The function $\frac{p_n}{n \log n}$, where $p_n$ denotes the $n$th prime, has a (divergent) asymptotic expansion of the form 
$$\frac{p_n}{n \log n} \simeq 1+ \sum_{k= 1}^n  \frac{P_k(\log \log n)}{(\log n)^k} \ (n \to \infty),$$
starting
$$ 1+\frac{\log \log n-1}{\log n}+{\frac {\log \log n-2}{(\log n)^2}}-{\frac {(\log \log n)^{2}-6\log \log n+11}{2(\log n)^{3}}}+O\left({\frac {(\log \log n)^{3}}{(\log n)^{4}}}\right),$$
 where  $P_0 = 1$, $P_1(x) = x-1$, and $P_k(x)$ for all $k \geq 2$ is a polynomial of degree at most $k-1$ that can be computed recursively as in \cite{rey}.
\end{example}

Recall that $B_n$ denotes the $n$th Bernoulli number (which equals $0$ if $n > 1$ is odd, since the defining generating function $\frac{X}{e^X-1}$, plus $\frac{X}{2}-1$, of the Bernoulli numbers is an even function).  For every nonnegative integer $n$, let  $B_n(T)$ denote   the {\bf $n$th  Bernoulli polynomial},\index{Bernoulli polynomial $B_n(T)$} which are defined collectively by their generating function
$$\sum_{n = 0}^\infty B_n(T)\frac{X^n}{n!} = \frac{Xe^{TX}}{e^X-1}.$$
Explicitly,  one has
$$B_n(T) = \sum_{k = 0}^n {n \choose k} B_{n-k}T^k$$
for all $n$, and thus $B_n(T)$ has constant term $B_n(0) = B_n$ for all $n$.  Moreover, one has $B_n(1) = B_n$ for all $n \neq 1$ and $B_1(1) = -B_1 = \frac{1}{2}$.    
The first six Bernoulli polynomials, for example, are given by
\begin{align*}
B_0(T) & =1 \\
B_1(T) & =T-\tfrac{1}{2} \\
B_2(T) & =T^2-T+\tfrac{1}{6} \\
B_3(T) & =T^3-\tfrac{3}{2}T^2+\tfrac{1}{2}T \\
B_4(T) & =T^4-2T^3+T^2-\tfrac{1}{30} \\
B_5(T) & =T^5-\tfrac{5}{2}T^4+\tfrac{5}{3}T^3-\tfrac{1}{6}T. 
\end{align*}
 The following theorem is known as the {\bf Euler--Maclaurin formula}.\index{Euler--Maclaurin formula}

\begin{theorem}[Euler--Maclaurin formula {\cite[Proposition 1.3]{kon} \cite[(A.24)]{ivic}}]\label{EMthm}
Let $f$ be an $N$-times continuously  differentiable complex-valued function on $[a,b]$, where $a,b \in \ZZ$ with $a< b$ and $N$ is a positive integer.  One has 
$$ \sum_{k=a}^{b}f(k) =  \int_{a}^{b}f(x)\,dx+{\frac {f(b)+f(a)}{2}} + \sum_{k=1}^{\lfloor N/2\rfloor }\,{\frac {B_{2k}}{(2k)!}}\left(f^{(2k-1)}(b)-f^{(2k-1)}(a)\right) +R,$$
where
$$R =(-1)^{N+1} \frac{1}{N!}\int_a^b B_{N}(\{x\})f^{(N)}(x) \, dx$$
and
$$|R| \leq \frac{2\zeta(N)}{(2\pi)^{N}} \int_a^b |f^{(N)}(x)| \, dx.$$
\end{theorem}

\begin{corollary}
Let $f$ be an infinitely differentiable complex-valued function on $[a,\infty)$, where $a \in \ZZ$.   Suppose that $f^{(n)}$ has a constant sign on $(a,\infty)$ for infinitely many nonnegative integers $n$.  Then one has the asymptotic expansion
$$ \sum_{k=a}^{b}f(k) \simeq \int_{a}^{b}f(x)\,dx+{\frac {f(b)+f(a)}{2}} + \sum_{k=1}^{\infty }\,{\frac {B_{2k}}{(2k)!}}\left(f^{(2k-1)}(b)-f^{(2k-1)}(a)\right) \ (b \to \infty).$$
\end{corollary}

Loosely speaking, the  Euler--Maclaurin formula can be viewed as an extension of the trapezoid rule by the inclusion of correction terms. The combinatorist G.-C.\ Rota  described the Euler--Maclaurin formula as ``one of the most remarkable formulas of mathematics''  that ``has proved very useful for over 200 years'' \cite[p.\ 11]{rota}.  It is particularly useful for deriving asymptotic expansions of functions that are important in analytic number theory and combinatorics.

\begin{example}  For any $s \in \CC$ with $\operatorname{Re} s > 1$, the Euler--Maclaurin formula can be shown to yield the asymptotic expansion
$$ \sum_{k = 1}^n \frac{1}{k^{s}} \simeq \zeta(s) + \frac{n^{1-s}}{1-s}+ \frac{1}{2}n^{-s}  -\sum_{k = 1}^\infty  \frac{(s)_{2k-1} }{n^s} \frac{B_{2k}}{(2k)!} \frac{1}{n^{2k-1}}   \ (n\to \infty),$$
where
$$(s)_n = s(s+1)(s+2)\cdots(s+n-1)$$ denotes the {\bf Pochhammer symbol}.\index[symbols]{.rt Ec@$(s)_n$}\index{Pochhammer symbol $(s)_n$}   
For $s = 2$, this simplifies to 
$$ \sum_{k=1}^{n}{\frac {1}{k^{2}}}\simeq \zeta (2)-{\frac {1}{n}}+{\frac {1}{2n^{2}}}-\sum_{k=1}^\infty{\frac {B_{2k}}{n^{2k+1}}} \ (n \to \infty).$$
For $s = 1$, the Euler--Maclaurin formula yields the asymptotic expansion
$$H_n  \simeq \log n +  \gamma  +\frac{1}{2n}-  \sum_{k = 1}^\infty \frac{B_{2k}}{2kn^{2k}} \ (n \to \infty),$$
where $H_n = \sum_{k = 1}^n \frac{1}{k}$ denotes the $n$th harmonic  number.
\end{example}

\begin{example}[{\cite{nemes}}]
{\bf Stirling's approximation}\index{Stirling's approximation} is the asymptotic relation
$$n! \sim  {\sqrt {2\pi n}}\left({\frac {n}{e}}\right)^{n} \ (n \to \infty)$$
for the factorial function $n!$.   The Euler--Maclaurin formula yields the asymptotic expansion
$${\displaystyle \log(n!) \simeq n\log n-n+{\tfrac {1}{2}}\log 2\pi n+\sum_{k=1}^{\infty}{\frac {B_{2k}}{2k(2k-1)n^{2k-1}}} \ (n \to \infty),}$$
whose first few terms are
 $$\log(n!)\simeq n\log n-n+{\tfrac {1}{2}}\log 2\pi n+{\frac {1}{12n}}-{\frac {1}{360n^{3}}}+{\frac {1}{1260n^{5}}}-{\frac {1}{1680n^{7}}}+\cdots \ (n \to \infty).$$
 (The constant ${\tfrac {1}{2}}\log 2\pi$ is derived from the {\it Wallis product} for $\pi$.)
From this one obtains the asymptotic expansion
$$\frac{n!}{ {\sqrt {2\pi n}}\left({\frac {n}{e}}\right)^{n}} \simeq 1+{\frac {1}{12n}}+{\frac {1}{288n^{2}}}-{\frac {139}{51840n^{3}}}-{\frac {571}{2488320n^{4}}}+\cdots \ (n \to \infty).$$
An explicit but rather complicated formula for the coefficients in the expansion above is  provided in \cite{nemes}.
\end{example}

\begin{figure}[ht!]
\includegraphics[width=70mm]{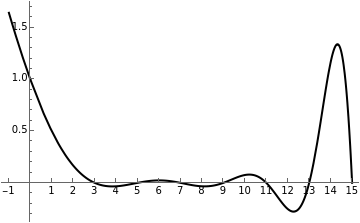}
\caption{\centering Graph of the Bernoulli function $B(x) = -x\zeta(1-x)$ on $[-1,15]$} 
\label{Bernoulli}
\end{figure}

\begin{remark}[Conventions regarding the Bernoulli number $B_1$ \cite{lusch}]
Numerous results,  including the Euler--Maclaurin formula, suggest that the alternative convention $B_1 = \frac{1}{2}$ is more natural than the standard convention  $B_1 =- \frac{1}{2}$ that we have employed \cite{lusch}.    The Bernoulli numbers so revised are denoted $B_n^+$ (where $B_n^+ = B_n(1)$ for all $n$,  so that $B_1^+ = -B_1$ and $B_n^+ = B_n$ for all $n \neq 1$) and have generating function
$$\sum_{n  = 0}^\infty \frac{B_n^+}{n!} X^n = \frac{X}{1-e^{-X}}.$$
The most pursuasive of the various arguments presented in \cite{lusch} in favor of the alternative convention  is that  the {\bf Bernoulli function} $$B(s) = -s\zeta(1-s)$$  is entire,  with $B(0) = 1$,  and one has $B(n) = B_n^+$ for all nonnegative integers $n$, and thus the Bernoulli function interpolates the revised Bernoulli numbers $B_n^+$.  Moreover, one has $$\zeta(s) =- \frac{B(1-s)}{1-s}$$
 for all $s \in \CC\backslash \{1\}$,  and therefore
 $$\zeta(-n) =- \frac{B_{n+1}^+}{n+1}, \quad \forall n \geq 0,$$
 which is slightly cleaner than the formula
  $$\zeta(-n) =(-1)^n \frac{B_{n+1}}{n+1}, \quad \forall n \geq 0.$$
  See  Figure \ref{Bernoulli} for a graph of the Bernoulli function on $[-1,15]$.  Note that $B(-1) = \zeta(2)$.
\end{remark}

\section{The degree $\deg f$ of a real function $f$}

Let $f \in \RR^{\RR_\infty}$, that is, let $f$ be a real function defined on a subset of $\RR$ that is not bounded above.   We define
$$\deg f = \inf\{t \in \RR: f(x) = o(x^t) \ (x \to \infty)\} \in \overline{\RR}.$$
We call $\deg f$ the {\bf (upper) degree of $f$}.\index{degree $\deg$}\index[symbols]{.g d@$\deg f$}

The following proposition provides an equivalent definition of degree.

\begin{proposition}\label{limsupdef}
For all $f \in \RR^{\RR_\infty}$, one has
$$\deg f = \limsup_{x \to \infty}  \frac{\log|f(x)|}{\log x}.$$
\end{proposition}

\begin{proof}
 Let $t \in \RR$.  One has
$f(x) = o\left(x^{s} \right) \ (x \to \infty)$
for all $s > t$ if and only if
for all $s > t$ there exists an $N > 0$ such that
$|f(x)| \leq x^{s},$
or equivalently
$\frac{\log|f(x)|}{\log x}  \leq s,$
for all $x \geq N$ in $\operatorname{dom} f$, where our convention is that $\log|0| = -\infty$.   The proposition follows readily from this (and from our convention that a lim sup of an extended real valued function is defined even if the values are $\pm \infty$ for any set of values of $x$).
\end{proof}

Using the proposition, one can easily verify the following examples of degree for some functions arising in calculus.

\begin{example}\label{degreeexamples} Let $a,b,c \in \RR$.
\begin{enumerate}
\item If $f \in \RR[x]$ is any real polynomial, then $\deg f$ coincides with the ususal degree of $f$ as a polynomial.
\item If $f/g \in \RR(x)$ is a rational function, where $f, g \in \RR[x]$ and $g$ is nonzero, then $\deg(f/g) = \deg f -\deg g$.
\item $\deg x^a = a$.
\item $\deg  \, (x^a+c)^b = b\max(a,0)$.
\item $\deg \exp = \infty$.
\item $\deg \, (\log)^a = 0$.
\item $\deg \, (\sin)^a  = \deg \, (\cos)^a = 0$.
\item $\deg \, \tan = \infty$.
\item $\deg e^{a\sqrt{x}}$ is equal to $\infty$ if $a > 0$, and to $-\infty$ if $a < 0$.
\item $\deg e^{a\sqrt{\log x}}= 0$.
\end{enumerate}
\end{example}

Statements of the form  $$f(x) = o(x^{d + \varepsilon}) \ (x \to \infty), \quad \forall \varepsilon > 0,$$ 
and  of the form
$$f(x) = O(x^{d + \varepsilon}) \ (x \to \infty), \quad \forall \varepsilon > 0,$$ appear throughout analytic number theory, and it is common but unstated knowledge to analytic number theorists that both of the statements above are equivalent to
$$\limsup_{x \to \infty}  \frac{\log|f(x)|}{\log x} \leq d.$$  By our definition of degree,  then, all three of the above statements are equivalent to  $\deg f \leq d$.    Expressions of the form $\limsup_{x \to \infty}  \frac{\log|f(x)|}{\log x}$  appear in the study of Dirichlet series, as shown in Theorem \ref{dirichlet}, and a great many other examples are provided throughout this book.  

 The following corollary of Proposition \ref{limsupdef} implies that the notions of limit superior, limit inferior, and degree are all interdefinable.  

\begin{corollary}\label{lims} For all $f \in \RR^{\RR_\infty}$, one has
$$\limsup_{ x \to \infty} f(x) =  \deg x^{f(x)}$$
and
$$\liminf_{ x \to \infty} f(x)  = - \limsup_{x \to \infty}\, (-f(x)) = - \deg x^{-f(x)}.$$
\end{corollary}

We  note the following proposition, whose proof is straightforward.

\begin{proposition}
Let $f, g \in \RR^{\RR_\infty}$.  One has the following.
\begin{enumerate}
\item $\deg f$ is the unique $d \in \overline{\RR}$ such that $f(x) = o(x^{t}) \ (x \to \infty)$  for all $t > d$ but for no $t < d$.
\item $\deg f$ is the unique $d \in \overline{\RR}$ such that $f(x) = O(x^{t}) \ (x \to \infty)$  for all $t > d$ but for no $t < d$.
\item $\deg f = \infty$ if and only if  $f(x) \neq o(x^{t}) \ (x \to \infty)$  for all $t \in \RR$, if and only if  $f(x) \neq O(x^{t}) \ (x \to \infty)$  for all $t \in \RR$.
\item $\deg f = -\infty$ if and only if  $f(x) = o(x^{t}) \ (x \to \infty)$  for all $t \in \RR$, if and only if  $f(x)= O(x^{t}) \ (x \to \infty)$  for all $t \in \RR$.
\item If $f$ is eventually bounded (on its domain), then $\deg f \leq 0$.
\item If $\deg f < 0$, then $\displaystyle \lim_{x \to \infty} f(x) = 0$ (on its domain).
\item If $f(x) = O(g(x)) \ (x \to \infty)$, then $\deg f  \leq \deg g$.
\item If  $f(x) = O(g(x)) \ (x \to \infty)$ and $g(x) = O(f(x)) \ (x \to \infty)$, then $\deg f = \deg g$. 
\end{enumerate}
\end{proposition}

Let  $f  \in \RR^{\RR_\infty}$.   The  {\bf exact degree of $f$},\index{exact degree $\overline{\underline{\deg}}\, f$}\index[symbols]{.g db@$\overline{\underline{\deg}}  \, f$} written $\overline{\underline{\deg}}\, f$, is the limit
$$\overline{\underline{\deg}} \, f= \lim_{x \to \infty} \frac{\log |f(x)|}{\log x},$$
provided that the limit exists or is $\pm \infty$,  in which case we say that $f$ has {\bf exact degree}.   Note that if $\deg f = -\infty$, then $f$ has exact degree.   We also define the  {\bf lower degree}\index{lower degree $\underline{\deg}\, f$}\index[symbols]{.g da@$\underline{\deg} \, f$} $\underline{\deg} \, f$ of $f$ to be
$$\underline{\deg} \, f = \liminf_{x \to \infty} \frac{\log |f(x)|}{\log x}.$$
Clearly one has $\underline{\deg} \, f \leq \deg f$,
with equality if and only if $f$ has exact degree.   If $f$ is not eventually nonzero (on its domain), then  $\underline{\deg} \, f = -\infty$; otherwise,  $f$ is eventually nonzero (on its domain),  and then one has
$$\underline{\deg} \, f = - \limsup_{x \to \infty} \frac{-\log|f(x)|}{\log x} = -\deg (1/f).$$

   The following are some number-theoretic examples of degree, lower degree, and exact degree.

\begin{example}\label{divisorf}  \
\begin{enumerate}
\item As noted in Theorem \ref{RC} (for which some details of the proof are provided in Section 5.2), it known that the Riemann constant $\Theta$ is equal to the infimum of all $t \in \RR$ such that 
$$\li(x)-\pi(x) = O(x^t) \ (x \to \infty).$$
 In other words, one has
$$\Theta = \deg(\li-\pi).$$
It follows that  $$\tfrac{1}{2} \leq \deg(\li-\pi) \leq 1,$$ and the Riemann hypothesis is equivalent to $\deg(\li-\pi) = \frac{1}{2}$.  Note, however, that $\li-\pi$ does  have not exact degree $\Theta$.  In fact, one has
$\li(x)-\pi(x) = 0$ for an unbounded set of positive real numbers $x$, and therefore
$$\underline{\deg} (\li-\pi) = -\infty.$$
\item Since
$$\li(x)-\Ri(x) \sim \frac{\sqrt{x}}{\log x} \ (x \to \infty),$$
one has
$$\overline{\underline{\deg}}(\li-\Ri) = \tfrac{1}{2},$$
that is, the function $\li-\Ri$ has exact degree $\tfrac{1}{2}$.
\item  \cite[Theorem 317]{har} states that
\begin{align*}
\limsup_{n \to \infty} \frac{\log d(n)\log \log n}{\log n} = \log 2,
\end{align*}
where $$d(n) = \sum_{d | n} 1$$ is the {\bf divisor function}.\index{divisor function $d(n)$}\index[symbols]{.rt  A@$d(n)$}  (The lim sup is attained by the sequence $p_1 p_2 \cdots p_n$.)  By Corollary \ref{lims}, this statement can be re-interpreted as $$\deg d(n)^{\log \log n} = \deg\, (\log n)^{\log d(n)} =  \log 2.$$
On the other hand, since there are infinitely many prime numbers, one has
\begin{align*}
\liminf_{n \to \infty} d(n) = 2
\end{align*}
and therefore
\begin{align*}
\liminf_{n \to \infty} \frac{\log d(n)\log \log n}{\log n} = 0,
\end{align*}
whence
\begin{align*}
\underline{\deg} \, d(n)^{\log \log n} = \underline{\deg}\, (\log n)^{\log d(n)} = 0.
\end{align*}
It also follows that
\begin{align*}
\overline{\underline{\deg}}  \, d(n) = \lim_{n \to \infty} \frac{\log d(n)}{\log n} = 0.
\end{align*}
\item It is known that $$\limsup_{n \to \infty} \frac{\sigma(n)}{n \log \log n} = e^\gamma,$$ 
where $$\sigma(n) = \sum_{d|n } d$$ is the {\bf sum of divisors function}\index{sum of divisors function $\sigma(n)$}\index[symbols]{.rt  B@$\sigma(n)$} and  $\gamma$ is the Euler--Mascheroni constant \cite[(25)]{gron}.   (The lim sup is attained by the sequence $(p_1 p_2 \cdots p_n)^{\lfloor \log p_n\rfloor}$.)
By Corollary \ref{lims}, the lim sup above can be re-interpreted as the statement $$\deg n^{\sigma(n)/(n \log \log n)} = e^\gamma.$$  Moreover, since there are infinitely many prime numbers, one has
$$\liminf_{n \to \infty} \frac{\sigma(n)}{n} = \liminf_{n \to \infty} \frac{\sigma(n)}{1+n} = 1.$$
It follows that
$$\overline{\underline{\deg}}  \, \sigma(n) = 1.$$
\item Let $\mathbb{L}$ denote the field of all {\bf logarithmico-exponential functions}\index{logarithmico-exponential functions} \cite{har3} \cite{har4},  that is, the field of all (germs of) real functions defined on a neighborhood of $\infty$ that can be can be built from all real constants and the functions $\id$, $\exp$,  and $\log$  using the operations $+$, $\cdot$, $/$, and $\circ$.  By Proposition  \ref{Kfield} of Section 6.3, every function in $\mathbb{L}$ has exact degree.
\end{enumerate}
\end{example}

Finally, we note the following.

\begin{proposition}
Let $f$ be a real function that is nonzero and differentiable on a neighborhood of $\infty$, and suppose that the limit
$$d = \lim_{x \to \infty} \frac{x f'(x)}{f(x)} \in \overline{\RR}$$ exists or is $\pm \infty$.  Then $f$ has exact degree $d$.     
\end{proposition}

\begin{proof}
By L'H\^opital's rule,  one has
$$\overline{\underline{\deg}} \, f   = \lim_{x \to \infty} \frac{\log |f(x)|}{\log x} = \lim_{x \to \infty} \frac{ f'(x)/f(x)}{1/x} = d.$$
\end{proof}

See Section 6.1 for a more thorough study of degree, lower degree, and exact degree.

\section{Slowly varying and regularly varying functions}

In this section we study {\it slowly varying} functions and {\it regularly varying} functions,  which were introduced by Karamata in \cite{kar1}  \cite{kar2}.   The study of such functions  is called {\it Karamata theory}.  Excellent references for Karamata theory and for the proofs of the theorems stated in this section are \cite{bgt} \cite[Chapter IV]{korevaar}  \cite{seneta}.  When necessary, to avoid Lebesgue integration, we restrict our discussion to continuous functions.

A real function $f$  defined  and either positive or negative on some neighborhood of $\infty$ is said to be {\bf slowly varying}\index{slowly varying} if 
$$\lim_{x \to \infty} \frac{f(cx)}{f(x)}  = 1$$
for all $c > 0$. 

\begin{example} Let $f$ and $g$ be real functions  defined on a neighborhood of $\infty$, let $c,a \in \RR$ be nonzero, and let $k$ and $l$ be nonnegative integers.
\begin{enumerate}
\item If $\lim_{x \to \infty} f(x) \neq 0$ exists, then $f$ is slowly varying.
\item $\sin x$ is of degree $0$ but is not slowly varying.
\item If $f$ is slowly varying and $f(x) \sim c g(x) \ (x \to \infty)$, then $g$ is slowly varying.
\item If $f$ and $g$ are slowly varying, then $fg$ is slowly varying.
\item If $f$ is eventually positive, then $f$ is slowly varying if and only if $cf^a$ is slowly varying.
\item $c(\log^{\circ k} x)^a$ is slowly varying (resp., of degree $0$) if and only if $k \geq 1$.
\item $e^{c(\log x)^a}$ is slowly varying  (resp., of degree $0$) if and only if $a < 1$.
\item $\exp^{\circ l}({c(\log^{\circ k}x)^a})$ is slowly varying  (resp., of degree $0$) if and only if $k > l$, or $k = l$ and $a< 1$.
\item $f(x) = 2+ \sin \log \log x$ is slowly varying with $\limsup_{x \to \infty} f(x) =3$ and $\liminf_{x \to \infty} f(x) = 1$.
\item $f(x) = \exp((\log x)^{1/2} \sin( (\log x)^{1/2}))$ is slowly varying with $\limsup_{x \to \infty} f(x) = \infty$ and $\liminf_{x \to \infty} f(x) = 0$.
\end{enumerate}
\end{example}

Note that, if $f$ is a real function that is eventually positive, and if $g$ is defined by $g(x) = \log f(e^x)$, so that $f(x) = e^{g(\log x)}$, then $f$ is slowly varying if and only if $$\lim_{ x\to \infty} (g(x+c)-g(x) ) = 0$$
for all $c \in \RR$.   This equivalence is used in proofs of the following result, which is known as the {\bf Karamata's integral representation theorem}.\index{Karamata's integral representation theorem}

\begin{theorem}[Karamata's integral representation theorem]\label{karintrep}
Let  $f$ be a real function defined on a neighborhood of $\infty$.  Then $f$ is slowly varying and continuous on a neighborhood of $\infty$ if and only if there exists an $N > 0$ and continuous functions $C$ and $\eta$  on $[N,\infty)$ such that  $\lim_{x \to \infty} C(x) \neq 0$ exists, $\lim_{x \to \infty} x \eta(x) = 0$, and $$f(x) = C(x) \exp \int_{N}^x \eta(t) \, dt$$
for all $x \in [N,\infty)$.
\end{theorem}

\begin{corollary}\label{svd0}
Let  $f$ be a real function defined on a neighborhood of $\infty$.  If $f$ is slowly varying and continuous on a neighborhood of $\infty$, then $\overline{\underline{\deg}} \, f = 0$.
\end{corollary}

\begin{proof}
By the theorem, L'H\^opital's rule,  and the fundamental theorem of calculus,  one has
$$\lim_{x \to \infty}  \frac{\log|f(x)|}{\log x} =  \lim_{x \to \infty} \frac{\log|C(x)|+\int_{N}^x \eta(t) \, dt}{\log x}  =  \lim_{x \to \infty} \frac{\int_{N}^x \eta(t) \, dt}{\log x}  = \lim_{x \to \infty} x\eta(x) = 0,$$
whence $\overline{\underline{\deg}}\,  f = 0$.
\end{proof}

Note that sufficiency of the equivalent condition in Karamata's integral representation theorem is easy to prove.   Indeed, assuming  $f(x) = C(x) \exp \int_{N}^x \eta(t) \, dt$, where $C$ and $\eta$ are as in the theorem,  one has
$$\left| \int^{cx}_x \eta(t) \, dt \right| \leq  \int^{cx}_x |\eta(t)| \, dt \leq (cx-x) \max_{t \in [x,cx]} |\eta(t)|= o(1) \ (x\to \infty)$$
for all $x \geq N$ and all $c > 1$,  and therefore
\begin{align*}
\lim_{x \to \infty} \frac{f(cx)}{f(x)}  =   \lim_{x \to \infty} \frac{C(cx)}{C(x)}\, \exp \lim_{x \to \infty}  \int^{cx}_x \eta(t) \, dt  = 1 \exp 0 = 1.
\end{align*}

Let $f$ be a real function defined and either positive or negative on some neighborhood of $\infty$.  The function $f$ is said to be {\bf regularly varying}\index{regularly varying} if the limit $\lim_{x \to \infty} \frac{f(cx)}{f(x)}$ exists and is finite and positive for all $c > 0$ (or, equivalently, for all $c > 1$).    If  there exists a constant $r \in \RR$ such that $$\lim_{x \to \infty} \frac{f(cx)}{f(x)} = c^r,$$ or, equivalently,
$$f(cx) \sim c^r f(x) \ (x \to \infty),$$  for all $c > 0$, then $f$ is regularly varying, $r$ is called the {\bf index of regular variation of $f$},\index{index of regular variation} and $f$ is said to be {\bf regularly varying of index $r$}.\index{regularly varying} 
The following result is clear.

\begin{proposition}\label{regvarp}  Let $f$ and $g$ be real functions defined on a neighborhood of $\infty$, and let $r,s,A \in \RR$ with $A \neq 0$.
\begin{enumerate}
\item  The power function $Ax^r$ is regularly varying of index $r$, and in fact it is the unique function $f(x)$ on $(0,\infty)$ with $f(cx) = c^r f(x)$ for all $c,x > 0$ and $f(1) = A$.
\item  If $f$ is regularly varying  (resp.,  regularly varying of index $r$), and if $f(x) \sim g(x) \ (x \to \infty)$, then $g$ is regularly varying  (resp.,  regularly varying of index $r$).  
\item  If $f$ is regularly varying  (resp.,  regularly varying of index $r$) and $g$ is regularly varying  (resp.,  regularly varying of index $s$), then $fg$ is is regularly varying  (resp.,  regularly varying of index $r+s$).
\end{enumerate}
\end{proposition}

Thus, for example,  the prime number theorem implies that the prime counting function is regularly varying of index $1$.

Karamata's integral representation theorem and Proposition \ref{regvarp} together yield the following.

\begin{corollary}\label{regvarcor}
Let $f$ be a real function that is nonzero and continuously differentiable on a neighborhood of $\infty$, and suppose that the limit
$$r = \lim_{x \to \infty} \frac{x f'(x)}{f(x)} \in \RR$$ exists.  Then $f$ is regularly varying of index $r$.
\end{corollary}

\begin{proof}
Let $g(x) = x^{-r} f(x)$ and $\eta(x) =  \frac{g'(x)}{g(x)}$, so that
$$ \lim_{x \to \infty} x\eta(x) = \lim_{x \to \infty} \frac{x g'(x)}{g(x)}=  \lim_{x \to \infty} \left(-r+\frac{x f'(x)}{f(x)}\right)= 0,$$
where $\eta(x)$ is continuous on a neighborhood $[N,\infty)$ of $\infty$, and where
$$g(x) = g(N)  \exp \int_{N}^x \eta(t)\, dt$$
for all $x \in [N,\infty)$.  From Karamata's integral representation theorem, then, it follows that $g$ is slowly varying,  whence $f$ is regularly varying of index $r$  by Proposition \ref{regvarp}.
\end{proof}

By the following theorem, which is known as {\bf Karamata's characterization theorem}\index{Karamata's characterization theorem}, any eventually continuous regularly varying function is regularly varying of index $r$ for a unique $r \in \RR$.  

\begin{theorem}[Karamata's characterization theorem]\label{regdeg}
Let  $f$ be a real function defined on a neighborhood of $\infty$. Then $f$ is regularly varying and continuous on a neighborhood of $\infty$  if and only if there exists a (unique)  $r \in \RR$ and a slowly varying function $F$ continuous on a neighborhood of $\infty$ such that $f(x) = x^r F(x)$ for all $x \gg 0$.  If those conditions hold, then $f$ is regularly varying of index $r$.
\end{theorem}

\begin{corollary}\label{regdegcor}
 Let  $f$ be a real function defined on a neighborhood of $\infty$.    If $f$ is regularly varying and continuous on a neighborhood of $\infty$, then  $\overline{\underline{\deg}} \, f \in \RR$ exists and is equal to the index of regular variation of $f$.
\end{corollary}

\begin{proof}
By the theorem,  the function $f$ is regularly varying of index $r$, and $x^{-r} f(x)$ is slowly varying, for some $r \in \RR$.  By Corollary \ref{svd0}, then, one has  $-r+ \overline{\underline{\deg}} \, f = \overline{\underline{\deg}}(x^{-r} f(x)) = 0$, whence $r = \overline{\underline{\deg}}\, f$. 
\end{proof}

By the characterization theorem,  Karamata's integral representation theorem generalizes as follows.

\begin{theorem}\label{karintrep2}
Let  $f$ be a real function defined on a neighborhood of $\infty$.  Then $f$ is regularly varying and continuous on a neighborhood of $\infty$ if and only if there exists an $N > 0$ and continuous functions $C$ and $\eta$ on $[N,\infty)$ such that $\lim_{x \to \infty} C(x) \neq 0$ and $r = \lim_{x \to \infty}x\eta(x)$ exist and $$f(x) = C(x) \exp \int_{N}^x \eta(t) \, dt$$
for all $x \in [N,\infty)$.  Moreover,  if those equivalent conditions hold, then $f$ is regularly varying of index $r$.
\end{theorem}

\begin{proof}
Necessity follows by applying Karamata's integral representation theorem to $f_1(x) = x^{-r}f(x)$, where $r$ is as provided by Theorem \ref{regdeg}, and sufficiency follows by applying Karamata's integral representation theorem to  $\eta_1(x) = \eta(x)-\frac{r}{x}$, where $r = \lim_{x \to \infty}x\eta(x)$.
\end{proof}

The following result characterizes the eventually monotonic  functions that are regularly varying.

\begin{proposition}\label{regvarfun}
 Let  $f$ be an eventually nonzero and monotonic real function defined on a neighborhood of $\infty$.  Then $f$ is regularly varying of index $r \in \RR$ if and only if $f(g(x)) \sim c^r f(x) \ (x \to \infty)$ for all $c >0$ and all functions $g$ defined on a neighborhood of $\infty$ with $g(x) \sim cx \ (x \to \infty)$.
\end{proposition}

\begin{proof}
Necessity is clear.  Let $h(x) =g(x)-cx$, so that $h(x) = o(x) \ (x \to \infty)$.   Suppose that $f$ is eventually nondecreasing.  It follows that, for every $\varepsilon > 0$,  one has
$$c^r(1-\varepsilon)^r  \sim \frac{ f(cx-c\varepsilon x)}{f(x)} \leq \frac{f(x+h(x))}{f(x)} \leq \frac{f(cx+c\varepsilon x)}{f(x)} \sim c^r (1+\varepsilon)^r \ (x \to \infty)$$
for all $x \gg 0$, and therefore $$\frac{f(g(x))}{f(x)} = \frac{f(x+h(x))} {f(x) } \sim c^r \  (x \to \infty).$$ 
The proof when $f$ is eventually nonincreasing is similar, with the inequalities reversed.
\end{proof}

The following result, known as {\bf Karamata's integral theorem}\index{Karamata's integral theorem},  which we use on several occasions in Part 3, is one of the main reasons that the notions of slowly varying and regularly varying functions are so useful.  It can be described as an asymptotic mean value theorem for continuous regularly varying functions.

\begin{theorem}[Karamata's integral theorem]\label{karam}
Let $f$ be a real function that is nonzero and continuous  on $[N,\infty)$, where $N > 0$.  
\begin{enumerate}
\item $f$ is regularly varying of index $r > -1$ if and only if there is a constant $s > 0$ such that
$$\frac{1}{x}\int_N^x f(t) \, dt \sim \frac{ f(x)}{s} \ (x \to \infty).$$  Moreover, if those equivalent conditions hold,  then one has
$s = r+1 = |r+1|$.
\item $f$ is regularly varying of index $r < -1$ if and only if the integral $\int_N^\infty f(t) \, dt$ exists and there is a constant $s > 0$ such that
$$\frac{1}{x}\int_x^\infty f(t) \, dt \sim \frac{f(x)}{s} \ (x \to \infty).$$  Moreover, if those  equivalent  conditions hold,  then one has
$s = - (r+1)  = |r+1|$.
\item If $f$ is regularly varying of index $-1$ and $\int_N^\infty f(t) \, dt = \pm \infty$, then 
$$f(x) = o\left(\frac{1}{x} \int_N^x f(t) \, dt \right) \ (x \to \infty).$$
\item If $f$ is regularly varying of index $-1$ and $\int_N^\infty f(t) \, dt$ exists, then 
$$f(x) = o\left( \frac{1}{x}\int_x^\infty f(t) \, dt \right) \ (x \to \infty).$$
\end{enumerate}
\end{theorem}

\begin{corollary}
Let $f$ be a real function that is nonzero and continuous  on $[N,\infty)$, where $N > 0$.   Then $f$ is slowly
 varying if and only if 
$$\frac{1}{x}\int_N^x f(t) \, dt \sim f(x) \ (x \to \infty).$$
\end{corollary}

\begin{example}\label{exkar}
For all nonnegative integers $n$, the asymptotic relation
\begin{align*}
\int_e^x \frac{n!\, dt} {(\log t)^{n+1}} +C_n \sim \frac{n!x}{(\log x)^{n+1}} \ (x \to \infty)
\end{align*}
stated previously in (\ref{lias}) follows from Karamata's integral theorem and the fact that  the function $\frac{1}{(\log x)^{n+1}}$ is slowly varying.
\end{example}

\begin{example}
Since the function $\frac{1}{\log x}$ is slowly varying, one has
$$\li(x) = \int_\mu^x \frac{1}{\log t} \, dt \sim \frac{x}{\log x} \ (x \to \infty),$$
where $\mu$ denotes the Ramanujan--Soldner constant.  More generally, for any $r > -1$, the function  $\frac{x^r}{\log x}$ is regularly varying of index $r > -1$, and therefore
$$\li(x^{r+1}) -  \li(\mu^{r+1})= \int_\mu^x \frac{t^r}{\log t} \, dt \sim \frac{x^{r+1}}{(r+1)\log x} \ (x \to \infty)$$
for all $x > 1$.  For $r = -1$, one has
$$\int_e^x \frac{t^{-1}}{\log t}  \, dt = \log \log x.$$
Finally, for $r < -1$, the function  $\frac{x^r}{\log x}$ is regularly varying of index $r < -1$, and therefore
$$\li(x^{r+1}) =- \int_x^\infty \frac{t^r}{\log t} \, dt \sim \frac{x^{r+1}}{(r+1)\log x} \ (x \to \infty)$$
for all $x > 1$. 
\end{example}

Karamata's integral theorem  can be viewed  in the following broader context.   Let $f$ be a function that is  Riemann integrable on $[N,x]$  for all $x >N$, where $N \in \RR$.  The {\bf average value $\operatorname{Avg}_f[y,x]$ of $f$ on the interval $[y,x]$}, for any $x > y \geq N$, is defined by
$$\operatorname{Avg}_f[y,x] = \frac{1}{x-y}\int_y^x f(t)\, dt.$$
Note that the power functions $f(x) = Ax^r$ for  $A, r \in \RR$ with $r > -1$ and $A \neq 0$ are characterized as those  real functions $f(x)$,  continuous and nonzero on $(0,\infty)$, such that $f(0^+) \in \{0,\pm \infty\}$ and
$$\operatorname{Avg}_f[0,x] = \frac{f(x)}{s}, \quad \forall x >  0,$$
for some $ s >0$ (where, necessarily,  $s = r+1$).  Indeed, by the fundamental theorem of calculus, the integral equation above  implies $sf = (xf)' = f+xf'$ and therefore $\frac{f'}{f} = \frac{s-1}{x}$.  Statement (1) of Karamata's integral theorem states that the eventually continuous regularly varying functions (of some index $r > -1$)  are characterized as those real functions $f(x)$,  continuous and nonzero on a neighborhood $[N,\infty)$  on $\infty$, such that
$$\operatorname{Avg}_f[N,x] \sim \frac{ f(x)}{s} \ (x \to \infty)$$  for some $ s >0$ (where, necessarily,  $s = r+1$).   Loosely speaking, then,  the continuous regularly varying functions are those continuous functions  whose running average asymptotically  behaves like the running average of a power  function.

Finally,  we note that all of the continuity  hypotheses utilized in this section can be weakened considerably: all of the results stated in this section can be generalized to {\it  (Lebesgue) measurable} functions $[N,\infty) \longrightarrow \RR$, as is typically done in the literature on Karamata theory.  For example, if  $f: [N,\infty) \longrightarrow \RR$ is measurable, then $f$ is regulary varying of index $\overline{\underline{\deg}} \, f \neq \pm \infty$ if $f$ is regularly varying.  Likewise, Karamata's integral representation theorem and  integral theorem extend to this more general setting, where the functions $f$ and $C$ are assumed measurable instead of continuous, and where the Riemann integrals are replaced with Lebesgue integrals.  In particular, since the derivative of a differentiable function is measurable,  one need only assume in Corollary \ref{regvarcor} that $f$ is differentiable (rather than continuously differentiable) on a neighborhood of $\infty$.  Note also that, if $f$ and $g$ are eventually measurable and regularly varying of index $r$ and $s$, respectively,  then $|f|+|g|$ is regularly varying of index $\max(r,s)$, and $f \circ g$ is regularly varying of index $rs$ provided that $\lim_{x \to \infty} g(x) = \infty$.

\chapter{Arithmetic functions}

In this chapter, we study arithmetic functions from the point of view of abstract algebra, e.g., from the point of view groups, rings, and fields, and then from the point of view of real and complex analysis.  The main algebraic tools introduced to study arithmetic functions are {\it Dirichlet convolution}, {\it formal Dirichlet series}, and {\it  formal Euler products}.  Many examples of both {\it multiplicative} and {\it additive} arithmetic functions and their corresponding formal Dirichlet series are provided.    One important analytic tool is a well-known theorem relating the summatory function of an arithmetic function to its Dirichlet series.  Also discussed in this chapter are  {\it Abel's summation formula}, various inversion theorems analogous to the M\"obius inversion theorem, the {\it Dirichlet hyperbola method},  and the notions of the {\it average value} and {\it average order} of an arithmetic function.  Applications are given to the prime counting function and its variants known as the {\it first and second Chebyshev functions} and the {\it Riemann prime counting function}, as well as to the {\it M\"obius function},  the {\it prime Omega function},  the {\it prime omega function}, the {\it divisor function},  and Euler's totient.

\section{Elementary complex functions}

There are several  equivalent ways to define the {\bf complex exponential function}\index{complex exponential function $\exp s = e^s$} $\exp s = e^s$. \index[symbols]{.sz  A@$\exp s$}\index[symbols]{.sz  Aa@$e^s$} The most popular is via its Taylor series $$e^s = \sum_{n = 0}^\infty \frac{s^n}{n!}.$$  Equivalently, $e^s$ is the unique entire function $\exp$ such that $\exp(0) = 0$ and $$\frac{d}{ds} \exp s = \exp s,$$ or, equivalently, the unique entire function  $\exp$ such that $\exp^{(n)}(0) =  1$ for all nonnegative integers $n$.    The author is partial to the definition
\begin{align}\label{EXP}
e^s = \lim_{n \rightarrow \infty} \left(1 + \frac{s}{n}\right)^n
\end{align}
for any complex number $s$ for which the limit exists, for, as we show in Remark \ref{EFpf}, it leads to a simple  ``Calculus 1'' derivation of {\bf Euler's formula} $$e^{a+bi} = e^a(\cos b + i\sin b), 
\quad \forall a,b \in \RR.$$
Note that, together with Euler's formula,  the identity $\cos^2x + \sin^2x = 1$ implies that
$$|e^s| = e^{\operatorname{Re} s}, \quad \forall s \in \CC,$$
and the sum laws for $e^x$, $\cos x$, and $\sin x$ imply that
$$e^{s+t} = e^s e^t, \quad \forall s, t \in \CC.$$

The argument of a complex number is unique only up to an additive  integer multiple of $2\pi$.  The {\bf principal branch of the argument of $s$},\index{principal branch $\operatorname{Arg} s$ of the argument of $s$}\index[symbols]{.sz  C@$\operatorname{Arg} s$} denoted $\operatorname{Arg} s$, is taken to be in the interval $(-\pi,\pi]$, where one defines $\operatorname{Arg} 0 = 0$.   Note that
$$e^{i \operatorname{Arg} s}  = \cos\operatorname{Arg} s+ i \sin \operatorname{Arg} s = \frac{s}{|s|}$$ 
for all $s \in \CC\backslash\{0\}$
and therefore
$$s = |s| e^{i \operatorname{Arg} s}$$
for all $s \in \CC$.   The {\bf principal branch  $\Log s$ of the  complex natural logarithm}\index{principal branch $
\Log s$ of the complex natural logarithm}\index[symbols]{.sz  B@$\Log s$} is the unique complex number $z$ such that $e^z = s$ and  $\operatorname{Im} z \in (-\pi,\pi]$.  For all $s \in \CC\backslash \{0\}$, one has
$$\Log s = \log |s|+i \operatorname{Arg} s$$
and
 $$\exp \Log s = s.$$
However, one has
 $$\Log \exp s =   \operatorname{Re} s+ i \operatorname{Arg} \exp s = s+2 \pi i\left\lfloor \frac{1}{2}-\frac{\operatorname{Im} s }{2\pi} \right \rfloor $$
for all $s \in \CC$, and one has $\Log \exp  s = s $ if and only if $\operatorname{Im} s \in (-\pi,\pi]$.
One also has
$$\frac{d}{ds} \Log s = \frac{1}{s}$$
for all $s \in \CC \backslash (-\infty, 0]$.   For all $s_1, s_2, \ldots, s_n \in \CC$, one has
\begin{align*}
\Log (s_1 s_2 \ldots s_n) & = \log |s_1|+\log |s_2| +\cdots+ \log |s_n| + i \operatorname{Arg} (s_1 s_2 \ldots s_n) \\
 & = \Log s_1+\Log s_2 +\cdots+ \Log s_n + i \left(\operatorname{Arg} (s_1 s_2 \ldots s_n) -\sum_{k = 1}^n \operatorname{Arg} s_k \right).
\end{align*}
Thus, for example, one has
$$\Log( st) = \Log s + \Log t + i (\operatorname{Arg}( st) - \operatorname{Arg}s-\operatorname{Arg} t),$$
which equals $\Log s+\Log t$ if  $\operatorname{Re}s$ and $\operatorname{Re} t$ are positive.  Note also that
$$\frac{\Log(s)+\Log(1/s)}{2} =  \begin{cases} \displaystyle
0  & \text{if } s\in \CC\backslash (-\infty,0] \\
 \pi i & \text{if } s \in (-\infty,0)
\end{cases}$$
for all $s \in \CC$.

 Note that one defines
$$x^s = e^{s \log x}$$
for all $x > 0$ and all $s \in \CC$.   More generally,  we let
$$a^s = e^{s \Log a}$$
for all $a,s \in \CC$ with $a \neq 0$.
For all $s \in\CC$, one defines
$$\cosh s = \frac{e^{s}+e^{-s}}{2},$$
$$\sinh s = \frac{e^{s}-e^{-s}}{2},$$
$$\tanh s = \frac{\sinh s }{\cosh s},$$
$$\cos s = \cosh i s,$$
$$\sin s = \frac{1}{i}\sinh i s$$
and
$$\tan s = \frac{\sin s}{\cos s} = \frac{1}{i}\tanh i s.$$
It follows that
$$e^s = \cos s + i \sin s = \cosh s + \sinh s$$
and 
$$1 = e^s e^{-s} = \cos^2 s + \sin^2 s = \cosh^2 s-\sinh^2 s$$
for all $s \in \CC$.

\begin{remark}[Derivation of Euler's formula]\label{EFpf}  In this remark, we provide a ``Calculus 1''  derivation Euler's formula from the definition (\ref{EXP}) of $e^s$. First, note that the complex numbers $$w_n =  \left(1 + \frac{s}{n}\right)^n$$ are well defined, and they converge to $w \in \CC$ if and only if the real part $a_n$ and imaginary part $b_n$ of $w_n = a_n + ib_n$ converge respectively to the real and imaginary parts of $w$.  Alternatively, the $w_n$ converge to $w$ if and only if the  modulus $r_n$ and  argument $\theta_n$ of $w_n = r_n (\cos \theta_n + i \sin \theta_n)$ converge respectively to the modulus and argument of $w$, provided that one counts angles as equal that differ by an integer multiple of $2\pi$.  We exploit this equivalence to compute the given limit: to compute the limit of $(1+s/n)^n = r_n(\cos \theta_n + i \sin \theta_n)$ as $n \rightarrow 
\infty$, we compute the limit of the modulus $r_n$ and the limit of the argument $\theta_n$.    Note first that the complex number $1 + s/n$, where $s = a+bi$, has modulus $((1+a/n)^2 + (b/n)^2)^{1/2}$ and argument $\arctan((b/n)(1+a/n)^{-1})$, where the given arctangent is well-defined if $1+s/n$ is close to $1$, hence if $n$ is sufficiently large.  Raising any complex number to the $n$th power raises its modulus to the $n$th power and multiplies its argument by $n$.  Therefore, the complex number $w_n  = (1+s/n)^n$ has modulus  $$r_n = ((1+a/n)^2 + (b/n)^2)^{n/2}$$ and argument $$\theta_n = n\arctan((b/n)(1+a/n)^{-1}).$$
Using these expressions for $r_n$ and $\theta_n$, we compute their respective limits as $n \to \infty$:
\begin{eqnarray*}
\lim_{n \rightarrow \infty}  r_n & = & 
\lim_{n \rightarrow \infty} \left(\left(1+\frac{a}{n}\right)^2+\frac{b^2}{n^2} \right)^{n/2} \\ & = &  \exp \left(\lim_{n \rightarrow \infty} \frac{n}{2} \log \left( \left(1+\frac{a}{n}\right)^2 + \frac{b^2}{n^2} \right)\right) \\
  & = & \exp \left(\frac{1}{2}\lim_{x \rightarrow 0} \frac{\log ( (1 +ax)^2+ b^2x^2 )}{x} \right) \\
  & = & \exp \left(\left.{\frac{1}{2} \frac{d}{dx}  \log ( (1 +ax)^2+ b^2x^2)}\right|_{x=0} \right)  \\
  & = & \exp \left(\left.{\frac{a(1+ax) + b^2x}{(1 +ax)^2+ b^2x^2}} \right|_{x=0} \right) \\
  & = & \exp \, a,
\end{eqnarray*}
and likewise
\begin{eqnarray*}
\lim_{n \rightarrow \infty}  \theta_n & = & 
\lim_{n \rightarrow \infty} n \arctan \left(\frac{b/n}{1+a/n}\right)  \\
  & = & \lim_{x \rightarrow 0} \frac{\arctan \left(\frac{bx}{1+ax}\right)}{x} \\
  & = & \left.{\frac{d}{dx} \arctan \left(\frac{bx}{1+ax}\right)}\right|_{x=0}  \\
  & = & \left.{\frac{b}{(1+ax)^2 + b^2x^2} }\right|_{x = 0}  \\
  & = & b.
\end{eqnarray*}
Thus, the limit of $(1+ s/n)^n$ as $n$ approaches infinity exists and is equal to that complex number with modulus $e^a$ and argument $b$, which is nothing other than $e^a(\cos b + i \sin b)$.  In other words, one has $$e^{a+bi} = \lim_{n \rightarrow \infty} \left(1+ \frac{a+bi}{n} \right)^n = e^a(\cos b + i \sin b),$$   
that is, 
\begin{align*}
e^s = e^{\operatorname{Re}s} (\cos \operatorname{Im} s + i \sin \operatorname{Im} s)
\end{align*}
for all $s \in \CC$, which is a generalization of Euler's formula to all complex numbers.   The derivation above of Euler's formula is more elementary than the usual ``Calculus 2'' proof using Taylor series in that it relies only upon definitions of the complex exponential function and derivatives  in terms of limits, the derivative formulas for the natural logarithm and arctangent functions, the chain rule, and some basic properties of complex numbers and the trigonometric functions cosine and sine.  
\end{remark}

\section{Arithmetic functions and formal Dirichlet series}

Recall that an {\bf arithmetic function} is a function $f: \ZZ_{> 0} \longrightarrow \CC$.  
Important arithmetic functions abound in number theory.  We introduced several in Section 1.2, namely, $p_n$, $p_{n+1}-p_n$,  $\chi_{\pp}(n)$, $\phi(n)$, and $\rho(n) = \frac{\phi(n)}{n}$.  We also mentioned the all-important {\it M\"obius function} $\mu$,   which we define later in this section.   The prime counting function $\pi(x)$ restricted to the positive integers is also an arithmetic function.  More generally, any complex-valued function $F(x)$ defined on $[1,\infty)$ restricts to an arithmetic function  $F(n)$ and can be studied as such.

Any arithmetic function $f$ can be alternatively conceptualized as the sequence $$\{f(n)\}_{n = 1}^\infty = (f(1),f(2),f(3),\ldots)$$ of complex numbers indexed by the positve integers, in that this correspondence defines a natural bijection
$$\CC^{\ZZ_{>0}} \longrightarrow \prod_{n = 1}^ \infty \CC.$$
Moreover, the bijection is an isomorphism of rings if one endows the set $\CC^{\ZZ_{>0}}$ of all arithmetic functions with pointwise addition $+$ and pointwise multiplication $\cdot$ and one endows the set $ \prod_{n = 1}^ \infty \CC$ with coordinatewise addition  and coordinatewise multiplication.  Both of these rings have many zerodivisors (those functions assuming the value $0$ at some $n$ and those sequences with some coordinate equal to $0$, respectively) and many invertible elements (precisely those that are not zerodivisors).    However, the operation $\cdot$ of pointwise multiplication on  the set $\CC^{\ZZ_{>0}}$  of all arithmetic functions is by no means the only important operation that endows that set with a ring structure.  Three other important operations are described below.

  There are many ways to form a  generating function of a sequence of complex numbers, hence of an arithmetic function, where in combinatorics  the generating function usually lies in  the field $\CC((T)) = \CC[[T]][1/T]$ of {\bf formal Laurent series}\index{formal Laurent series} over $\CC$, which is the quotient field of the integral domain $\CC[[T]]$ of formal power series.  The simplest way to create such a generating function from a function $f: \ZZ_{\geq 0} \longrightarrow \CC$ on the nonnegative integers is via its {\bf formal Taylor series}\index{formal Taylor series}
$$F_f(T) = \sum_{n = 0}^\infty f(n)T^n \in \CC[[T]].$$
 For example, for the function $a^{\id}$ with $a^{\id}(n) = a^n$ for all $n$, where $a \in \CC$ is fixed, one has
$$F_{a^{\id}}(T) = 1+aT+a^2T^2+a^3T^3\cdots  = \frac{1}{1-aT}.$$
Note that, for all functions $f,g \in \CC^{\ZZ_{\geq 0}}$, one has
$$F_{f+g}(T) = F_f (T)+ F_g(T)$$
and
$$F_{f \star g}(T) = F_f(T)F_g(T),$$
where $f \star g$ is the {\bf additive convolution}\index{additive convolution}  of  $f$ and $g$, which  is defined by
\begin{align*} 
(f \star g)(n) & = \sum_{a+b = n} f(a)g(b) \\
& = \sum_{k = 0}^{n} f(k)g(n-k) \\
& = f(0)g(n)+f(1)g(n-1)+\cdots+ f(n-1)g(1)+f(n)g(0)
\end{align*}
for all nonnegative integers $n$, where the first sum is over all pairs $(a,b)$ of nonnegative integers $a$ and $b$ such that $a+b = n$.  Thus, the bijection 
$$F_{-}(T): \CC^{\ZZ_{\geq 0}} \longrightarrow  \CC[[T]]$$ acting by
$$f \longmapsto F_f(T)$$
is an isomorphism of integral domains if one endows the set $\CC^{\ZZ_{\geq 0}}$ with the operations of (pointwise) addition and additive convolution.    This ring structure is more interesting than that of $\prod_{n = 0}^\infty \CC$, as the ring $\CC[[T]]$  is a {\bf discrete valuation ring},\index{discrete valuation ring} that is, a principal ideal domain with a unique maximal ideal, namely,  in this case, the ideal $(T)$ generated by $T$.  A discrete valuation ring is equivalently a unique factorization domain with at most one prime (or irreducible) element up to associates, and the prime elements of $\CC[[T]]$  are precisely the elements of the set $(T)\backslash (T^2)$.

Another generating function of a function $f\in \CC^{\ZZ_{\geq 0}}$ that is often used is the {\bf exponential generating function} $\sum_{n = 0}^\infty \frac{f(n)}{n!} T^n$.  This is just the formal Taylor series of the arithmetic function  $g(n) = \frac{f(n)}{n!}$, but then multiplication of generating functions corresponds to a different convolution, namely, the {\bf biomial convolution}\index{binomial convolution}
$$(f \square g)(n) =  \sum_{a+b = n} {n \choose a, b} f(a)g(b) = \sum_{k = 0}^{n} {n \choose k} f(k)g(n-k)$$
{\bf of $f$ and $g$}.   Addition and binomial convolution  endow a ring structure on $\CC^{\ZZ_{\geq 0}}$ that is also isomorphic to $\CC[[T]]$, where the map to $\CC[[T]]$ sends $f$ to $\sum_{n = 0}^\infty \frac{f(n)}{n!}T^n$.

The two types of generating functions defined above, and their corresponding convolution products, are the most prevalent in combinatorics.  In number theory, one of the most important generating functions associated to an arithmetic function $f$ is its {\bf formal Dirichlet series}\index{formal Dirichlet series $D_f(X)$}\index[symbols]{.k R@$D_f(X)$} 
$$D_f(X) = \sum_{n = 1}^\infty f(n)n^{-X}.$$ 
Like the expression ``$X^n$,'' the expression ``$n^{-X}$'' is just a formal symbol, and, as with formal power series, we do not consider issues of convergence.   The {\it Dirichlet series of $f$} is obtained from the formal Dirichlet series $D_f(X)$ of $f$ by ``substituting'' a complex variable $s$ for the formal ``variable'' $X$, i.e., by substituting $n^{-s}$ for the monomial $n^{-X}$.  Dirichlet series, for which matters of convergence are relevant, are studied in Section 3.8.

In analogy with the equations $X^m X^n = X^{m+n}$ and $(X^m)^n = X^{mn}$, we write $$m^{-X} n^{-X} = (mn)^{-X}$$
and
$$(m^{-X})^n = (m^n)^{-X} = m^{-nX}.$$
  Formal Dirichlet series, by definition, are added pointwise, so that
$$D_{f}(X) + D_{g}(X)= \sum_{n = 1}^\infty(f(n)+g(n))n^{-X} = D_{f+g}(X),$$
and they are multiplied as follows: 
\begin{align*}
D_{f}(X) D_{g}(X) & = \sum_{a = 1}^\infty \sum_{b = 1}^\infty f(a)g(b) a^{-X} b^{-X}\\
& = \sum_{n = 1}^\infty \sum_{ab = n} f(a)g(b) (ab)^{-X}\\
& =  \sum_{n = 1}^\infty \left(\sum_{ab = n} f(a) g(b)\right) n^{-X} \\
& =  \sum_{n = 1}^\infty ((f*g)(n) )n^{-X} \\
& = D_{f * g}(X),
 \end{align*}
where $f * g$ is the {\bf Dirichlet convolution},  or {\bf Dirichlet product},\index{Dirichlet convolution $f*g$}\index{Dirichlet product $f*g$} of  $f$ and $g$, defined by
\begin{align*} 
(f * g)(n) & = \sum_{ab = n} f(a)g(b) \\
& = \sum_{d|n} f(d)g\left(\frac{n}{d}\right),\index[symbols]{.ia  t@$f*g$} 
\end{align*}
where the respective sums are over all ordered pairs $(a,b)$ of positive integers $a$ and $b$ such that $ab = n$, and  over all positive divisors $d$ of $n$.    Thus, for example, one has
$$(f*g)(6) = f(1)g(6)+f(2)g(3)+f(3)g(2)+f(6)g(1).$$
The operations $+$ and $*$ make the set $\CC^{\ZZ_{> 0}}$ into an integral domain called the {\bf ring $\operatorname{Arith}_\CC$ of arithmetic functions (over $\CC$)}.  Let us denote by $\operatorname{Dir}_\CC$ the set of all formal Dirichlet series $D_f(X)$ for $f \in \operatorname{Arith}_\CC$.  Then  $\operatorname{Dir}_\CC$ is also an integral domain under the operations defined above, and the map
$$D_{-}(X): \operatorname{Arith}_\CC \longrightarrow \operatorname{Dir}_\CC$$
acting by 
$$f \longmapsto D_f(X)$$
that associates to each arithmetic function its formal Dirichlet series is an isomorphism of integral domains.
Note that the multiplicative identity of the ring $\operatorname{Arith}_\CC$ is the function $\iota = \delta_{1,-}$,\index[symbols]{.q  A@$\iota(n)$}   where $\delta_{i,j}$ denotes the Kronecker delta, so that 
$$\iota(n) = \begin{cases} 1 & \quad  \text{if }  n = 1 \\
 0  & \quad \text{if } n > 1. 
\end{cases}$$
The corresponding element of the ring $\operatorname{Dir}_\CC$ is denoted $1 = 1^{-X}$.

 By the fundamental theorem of arithmetic, one can write
$$n^{-X} = \prod_{p} (p^{-X})^{v_p(n)},$$
where $v_p(n)$ for any prime $p$ denotes the exponent of $p$ in the unique prime factorization of $n$.
Let $$X_k = p_k^{-X}$$ for all positive integers $k$, where $p_k$ denotes the $k$th prime.  Then the map
$$\operatorname{Arith}_\CC \longrightarrow \CC[[X_1,X_2,X_3,\ldots]]$$
acting by $$f  \longmapsto \sum_{n = 1}^\infty f(n) \prod_{k= 1}^\infty X_k^{v_{p_k}(n)}$$
is an isomorphism of integral domains, where in the formal power series ring $\CC[[X_1,X_2,X_3,\ldots]]$ one assumes that each monomial $\prod_{k = 1}^\infty X_k^{n_k}$ involves only finitely many variables $X_k$, but where the formal sums  of such monomials in the ring can have arbitrary coefficients.   Thus, one can conceptualize arithmetic functions and formal Dirichlet series as formal power series in a countably infinite number of variables (one for each prime).  The formal Taylor series generating function $F_f(X)$ of $f \in \CC^{\ZZ_{\geq 0}}$ reflects the structure of the positive integers under addition, which is the free commutative semigroup generated by $1$, while the formal Dirichlet series generating function $D_f(X)$ of $f \in \CC^{\ZZ_{>0}}$ reflects the  structure of the positive integers under multiplication, which by the fundamental theorem of arithmetic is the free  commutative monoid generated by the prime numbers.

It is well known that the integral domain $\CC[[X_1,X_2,X_3,\ldots]]$ is a unique factorization domain.  Thus, so is the integral domain ${\operatorname{Arith}}_\CC$.    Moreover, the integral domain $\CC[[X_1,X_2,X_3,\ldots]]$ has a unique maximal ideal, namely, the set of all series with constant coefficient $0$, and thus ${\operatorname{Arith}}_\CC$ has the unique maximal ideal $$\mm_\CC = \{f \in {\operatorname{Arith}}_\CC: f(1) = 0\}.$$    It follows that the invertible elements of the integral domain ${\operatorname{Arith}}_\CC$ are precisely those lying in the group
$${\operatorname{Arith}}_\CC^* = {\operatorname{Arith}}_\CC\backslash \mm_\CC = \{f \in {\operatorname{Arith}}_\CC: f(1) \neq 0\}.$$
under (Dirichlet) multiplication.   For any arithmetic function $f \in \operatorname{Arith}_\CC^*$, we denote its inverse with respect to $*$, called the {\bf Dirichlet inverse of $f$},\index{Dirichlet inverse $f^{*-1}$}\index[symbols]{.ia  t@$f^{*-1}$}  by $f^{*-1}$.  For any $f \in  \operatorname{Arith}_\CC^*$, one has
$$D_{f^{*-1}}(X) = \frac{1}{D_f(X)}.$$

Now,  if $R$ is any ring (commutative with identity), then any  ideal $\aaa$ of $R$ induces a topology on $R$ called the {\it $\aaa$-adic topology}.    In particular, the ring $\CC[[X_1,X_2,X_3,\ldots]]$ has a topology induced by its unique maximal ideal, and it is complete with respect to that topology.  The same thus holds for the rings ${\operatorname{Arith}}_\CC$ and ${\operatorname{Dir}}_\CC$, where the topology on ${\operatorname{Arith}}_\CC$ is the $\mm_\CC$-adic topology.   Equivalently, the topology on  ${\operatorname{Arith}}_\CC$  is just the product topology on $\prod_{n = 1}^\infty \CC$, where $\CC$ has the discrete topology (which of course is not the usual topology on $\CC$).
In fact, there is a {\it  norm} $$\Vert \cdot\Vert : {\operatorname{Arith}}_\CC \longrightarrow \RR_{\geq 0}$$ on ${\operatorname{Arith}}_\CC$, where
$$\Vert f\Vert  = \frac{1}{\inf\{n \in \ZZ_{> 0} : f(n) \neq 0\}} \in \{0\} \cup \{1/n : n \in \ZZ_{> 0}\}  \subseteq \RR_{\geq 0}\index[symbols]{.ia  w@${\Vert f\Vert }$}$$ 
for all $f \in {\operatorname{Arith}}_\CC$, that induces the $\mm_\CC$-adic topology, where a {\bf norm}\index{norm}  on an integral domain $A$ is a map $| \cdot|  : A \longrightarrow \RR_{\geq 0}$ satisfying the following conditions for all $f, g \in A$.
\begin{enumerate}
\item $| f|  = 0$ if and only if $f = 0$.
\item $| f+g|  \leq |f|+|g|.$
\item $| f g|  = | f|  | g| $.
\end{enumerate}
In this case, $d(f,g) = \vert f-g\vert $ defines a metric, hence a topology, on the ring $A$, called the {\bf norm metric}\index{norm metric}  and {\bf norm topology},\index{norm topology} respectively, induced by the norm $| \cdot |$.    It is clear that addition, subtraction, and multiplication are then continuous as functions from $A\times A \longrightarrow A$, where $A \times A$ has the product topology.    A norm $|\cdot|$ on $A$ is {\bf nonarchimedean} if $| f+g|  \leq \max(| f| ,| g| )$  for all $f, g \in A$, which is much stronger than condition (2) above.     It is straightforward to check that the map $\Vert \cdot \Vert$ is a nonarchimedean norm on the integral domain $\operatorname{Arith}_\CC$.    Two  arithmetic functions  $f, g \in \operatorname{Arith}_\CC$ are close in the induced metric $d(f,g) = \Vert f-g\Vert $ if $\Vert f-g\Vert  \ll 1$, that is,  if $f(n) = g(n)$ for all $n \leq N$ for some $N \gg 0$.   Moreover, the ring ${\operatorname{Arith}}_\CC$ is {\it complete} with respect to the norm metric in the sense that every Cauchy sequence converges.  Note also that $$\operatorname{Arith}_\CC^*  = \{f \in \operatorname{Arith}_\CC: \Vert f\Vert  = 1\}$$ and $$\mm_\CC  = \{f \in \operatorname{Arith}_\CC: \Vert f\Vert  < 1\},$$ and thus the units of $\operatorname{Arith}_\CC$ are the elements of largest norm.

Let $R$ be a subring of $\CC$.  Then the set ${\operatorname{Arith}}_R$  of all $R$-valued arithmetic functions is a subring of  the integral domain ${\operatorname{Arith}}_\CC$, called the {\bf ring of arithmetic functions over $R$},\index{ring of arithmetic functions $\operatorname{Arith}_R$}\index[symbols]{.ia  r@$\operatorname{Arith}_R$} and is isomorphic to the ring $R[[X_1,X_2,X_3,\ldots]]$ of formal power series over $R$ in a countably infinite number of variables.  Moreover, the ideal $$\mm_R = \mm_\CC \cap {\operatorname{Arith}}_R\index[symbols]{.ia  rc@$\mm_R$}$$   is the intersection of all of the maximal ideals  of the integral domain ${\operatorname{Arith}}_R$, and the set
$$1+\mm_R = \{f \in {\operatorname{Arith}}_R: f(1) = 1\}\index[symbols]{.ia  rd@$1+\mm_R$}$$  is a subgroup of the group
$${\operatorname{Arith}}_R^* = \{f \in \operatorname{Arith}_R: f(1) \in R^*\}\index[symbols]{.ia  ra@$\operatorname{Arith}_R^*$}$$
of all units of the integral domain ${\operatorname{Arith}}_R$, where $R^*$ denotes the group of units of the integral domain $R$.

\section{Multiplicative and additive arithmetic functions}

Let $R$ be a subring of $\CC$.  The subset
$${\operatorname{Mult}}_R =  \{f \in {\operatorname{Arith}}_R: f \text{ is multiplicative}\}$$  
of $1+\mm_R$, where an arithmetic function $f$ is said to be {\bf multiplicative}\index{multiplicative} if
$$f(1) = 1, \text{ and } f(mn) = f(m)f(n), \quad \forall m,n \in \ZZ_{>0}: \, \gcd(m,n) = 1,$$
is a subgroup of the group $1+\mm_R$.  To gain familiarity with Dirichlet convolution, it is a good exercise to prove that ${\operatorname{Mult}}_R$ is indeed a subgroup of the group $1+\mm_R$.  The  multiplicative arithmetic functions $f$ are freely and uniquely determined by their values $f(p^k)$ at the prime powers $p^k$ for $p$ prime and $k \geq 1$, and for such $f$ one has
$$f(n) = \prod_{p} f(p^{v_p(n)})$$
for all $n$.\index[symbols]{.ia  s@$\operatorname{Mult}_R$}  

Many important arithmetic functions that arise in number theory are multiplicative.  In fact, some are even {\it completely multiplicative}, where an arithmetic function $f$ is said to be {\bf completely multiplicative}\index{completely multiplicative} if
$$f(1) = 1, \text{ and } f(mm) = f(m)f(n), \quad \forall m,n \in \ZZ_{>0}.$$  The completely multiplicative functions $f$ are freely and uniquely determined by their values $f(p)$ at the prime numbers $p$, and for such $f$ one has
$$f(n) = \prod_{p} f(p)^{v_p(n)}$$
for all $n$.   Unlike the set of all multiplicative functions, the set of all completely multiplicative functions is not closed under Dirichlet multiplication or Dirichlet inversion.    However, both sets are closed under pointwise multiplication (and also under pointwise reciprocation, for those whose values are all nonzero).

The most obvious completely multiplicative functions are those of the form $\id^a$ for some $a \in\CC$, e.g., the {\bf zeta (arithmetic) function}\index{zeta arithmetic function $\zeta(n)$} $\zeta = \id^0$,\index[symbols]{.q  B@$\zeta(n)$}  which is identically $1$ and has formal Dirichlet series 
$$\zeta(X) = D_{\zeta}(X) = \sum_{n = 1}^\infty n^{-X},$$
which is called the {\bf formal Riemann zeta function}.\index{formal Riemann zeta function $\zeta(X)$}  (We abuse notation and use $\zeta$ to denote an arithmetic function, its formal Dirichlet series, and the usual Riemann zeta function, but in any given context it is clear which of the three is intended.)   Since the arithmetic function $\zeta$ is completely multiplicative, its Dirichlet inverse $\zeta^{*-1}$, denoted $\mu$, exists, and is multiplicative.  The function $\mu$ is called the {\bf M\"obius function}\index{Mobius function @ M\"obius function $\mu(n)$} \index[symbols]{.q M@$\mu(n)$}  and is extremely important in number theory and combinatorics.  For all positive integers $n$, one has
\begin{align}\label{mob}
\mu(n) =  \begin{cases} (-1)^{\Omega(n)} = (-1)^{\omega(n)} & \quad  \text{if }  n \text{ is squarefree} \\
 0  & \quad \text{if } n \text{ is not squarefree}, 
\end{cases}
\end{align}
where 
$$\Omega(n) =  \sum_{p|n} v_p(n) =  \sum_p v_p(n)$$
denotes the number of prime factors of $n$ counting multiplicity,  where
$$\omega(n) = \sum_{p| n} 1 = \sum_{p:\, v_p(n) > 0} 1 \leq \Omega(n)$$
denotes the number of distinct prime factors of $n$, and where $n$ is said to be {\bf squarefree}\index{squarefree} if any of the following equivalent conditions hold.
\begin{enumerate}
\item $1$ is the only perfect square that divides $n$.
\item No square of a prime divides $n$.
\item $n$ is a product of distinct primes.
\item $\Omega(n) = \omega(n)$.
\end{enumerate}
The formal Dirichlet series of $\mu$ is given by 
$$D_\mu(X) = \frac{1}{\zeta(X)}.$$

  The functions $\Omega$ and $\omega$ are also  important arithmetic functions, called the {\bf prime Omega function} and the  {\bf prime omega function}, respectively.\index{prime Omega function $\Omega(n)$}\index{prime omega function $\omega(n)$}\index[symbols]{.rr  N@$\Omega(n)$}\index[symbols]{.rr  O@$\omega(n)$}  Those two arithmetic functions, along with $0$, $\iota$, $\zeta$, and $\mu$, are among the most important arithmetic functions in the algebraic theory of arithmetic functions.   All six of these functions are integer-valued and thus lie in the ring ${\operatorname{Arith}}_\ZZ$.

The functions $\omega$ and $\Omega$ are {\it additive} and {\it completely additive}, respectively, where an arithmetic function $f$ is said to be {\bf additive}\index{additive} if
$$f(mn) = f(m)+f(n), \quad \forall m,n \in \ZZ_{>0}: \, \gcd(m,n) = 1,$$
and where $f$ is said to be {\bf completely additive}\index{completely additive} if
$$f(mn) = f(m)+f(n), \quad \forall m,n \in \ZZ_{>0}.$$  Note that, if $f$ is additive, then necessarily $f(1) = 0$.  The additive functions $f$ are freely and uniquely determined by their values $f(p^k)$ at the prime powers $p^k$ for $p$ prime and $k \geq 1$, and for such $f$ one has
$$f(n) = \sum_{p} f(p^{v_p(n)})$$
for all $n$.  Likewise, the completely additive functions $f$ are freely and uniquely determined by their values $f(p)$ at the prime numbers $p$, and for such $f$ one has
$$f(n) = \sum_{p} {v_p(n)} f(p)$$
for all $n$.  Moreover, for any prime $p$,  the function $v_p$ itself is  a completely additive arithmetic function,  known as the {\bf $p$-adic valuation},\index[symbols]{.rr Na@$v_p(n)$}\index{adic valuation $v_p(n)$}  and the sum $\sum_{p} v_p$ over all primes is equal to $\Omega$.  

\begin{remark}[Combinatorial interpretation of $\omega(n)$ and $\Omega(n)$]
For any positive integer $n$, the {\bf divisor lattice of $n$}\index{divisor lattice} is the finite partially ordered set $D(n) = \{d \in \ZZ_{>0}: d | n\}$ of all positive divisors of $n$,  partially ordered by the divides relation $\mid$.  Note that $D(n)$ is a lattice, where least upper bounds coincide with lcms, and greatest lower bounds coincide  with gcds.  Let $P$ be any finite poset.  The {\bf  order dimension of $P$} is the minimal number of totally ordered sets such that $P$ embeds into their product with the componentwise ordering.    The  {\bf height of $P$} is the largest possible cardinality of any totally ordered subset of $P$.   For any positive integer $n$, the positive integer $\omega(n)$ coincides with the order dimension of the poset $D(n)$, while the positive integer $\Omega(n)$ coincides with its height.   Note in fact that all maximal totally ordered subsets of $D(n)$ have cardinality $\Omega(n)$.  The well-studied function $d(n) = \# D(n)$  is called the {\bf divisor function}: see Example \ref{divex}.
\end{remark}

\begin{example}
Let $a_p \in \CC$ for all primes $p$.  Then the arithmetic function
$$f(n) = \prod_{p|n} a_p, \quad \forall n \in \ZZ_{> 0},$$
is multiplicative, and the arithmetic function
$$g(n) = \sum_{p|n} a_p, \quad \forall n\in \ZZ_{> 0},$$ is  additive.  For example, if $a_p = 1-\frac{1}{p}$ for all $p$, then $f(n) = \frac{\phi(n)}{n}$ for all $n$, where $\phi$ denotes  Euler's totient.  On the other hand, if $a_p = 1$ for all $p$, then $g = \omega$.  Arithmetic functions $f$ of the first type are said to be {\bf strongly multiplicative},\index{strongly multiplicative} and arithmetic functions $g$ of the second type are said to be {\bf strongly additive}.\index{strongly additive}   An arithmetic function $h$ is strongly multiplicative (resp., strongly additive) if and only if it is multiplicative (resp., additive) and $h(p^n) = h(p)$ for all primes $p$ and all positive integers $n$. Like the completely multiplicative functions and the completely additive functions, the strongly multiplicative functions and the strongly additive functions are freely and uniquely determined by their values at the prime numbers.
\end{example}

 Let $f \in \operatorname{Arith}_\CC$ be an arithmetic function.  The arithmetic function $$\widehat{f} = \zeta* f,$$ given  by
$$\widehat{f}(n)= \sum_{d| n} f(n)$$
for all $n$, is called the {\bf Dirichlet transform of $f$},\index{Dirichlet transform $\widehat{f}$} and the arithmetic function $$\widecheck{f} = \mu* f\index[symbols]{.ia  ta@$\widehat{f}$},$$ given by
$$\widecheck{f}(n) = \sum_{d| n} \mu\left(\frac{n}{d}\right)f(d) = \sum_{d| n} \mu(d)f\left(\frac{n}{d}\right)\index[symbols]{.ia  tb@$\widecheck{f}$}$$
for all $n$, is called the {\bf M\"obius transform of $f$}.\index{Mobius transform @ M\"obius transform $\widecheck{f}$}  Note that, for any arithmetic function $f$, one has 
$$D_{\widehat{f}}(X) = \zeta(X) D_f(X)$$
and
$$D_{\widecheck{f}}(X) =\frac{ D_f(X)}{\zeta(X)}.$$
Since the operation $*$ is associative and $\zeta$ and $\mu$ are inverses with respect to $*$, one has the following result, which is known as the {\bf M\"obius inversion theorem (for Dirichlet convolution)}.\index{Mobius inversion @ M\"obius inversion theorem}

\begin{proposition}[M\"obius inversion theorem for Dirichlet convolution]\label{mobius1}
One has
$$f = \widecheck{\widehat{f}}= \widehat{\widecheck{f}}$$
for all arithmetic functions $f$.   Moreover, if any one of the functions $f$, $\widehat{f}$, and $\widecheck{f}$ is multiplicative, then all three of them are. 
\end{proposition}

 This is one of many ``inversion'' theorems that are important in number theory and combinatorics.  See Section 3.5 for several others.

Below are examples of  some widely-studied multiplicative arithmetic functions and their corresponding formal Dirichlet series.

\begin{example}\label{divex}  \
\begin{enumerate} 
\item For any $a \in \CC$, let  $\sigma_a = \widehat{\id^a}$ denote the Dirichlet transform of the multiplicative function $\id^a$, which is also multiplicative, and which is called the {\bf generalized sum of divisors function}.\index{generalized sum of divisors function $\sigma_a(n)$}\index[symbols]{.rt  C@$\sigma_a(n)$} Thus, one has
$$\sigma_a(n)= \sum_{d|n} d^a$$
for all  $n$.   The function $d(n) = \sigma_0(n)$ is called the {\bf divisor function},\index{divisor function $d(n)$}\index[symbols]{.rt  A@$d(n)$}  since $d(n)= \sum_{d | n} 1$ equals the number of positive divisors of $n$ for all $n$.  The function $\sigma(n) = \sigma_1(n)$ is called the {\bf sum of divisors function},\index{sum of divisors function $\sigma(n)$}\index[symbols]{.rt  B@$\sigma(n)$} since $\sigma(n)= \sum_{d | n} d$ equals the sum of the positive divisors of $n$ for all $n$.  Note that
$$\sigma_{-a}(n) =  \sum_{d|n} d^{-a} = n^{-a} \sum_{d|n} \left(\frac{n}{d}\right)^{a} =n^{-a} \sum_{d|n} d^a = \frac{\sigma_a(n)}{n^a}$$
for all $n$ and all $a \in  \CC$, so, for example,
$$\sigma_{-1}(n) = \frac{\sigma(n)}{n}$$
for all $n$. The formal Dirichlet series of $\id^a$ is written
$$D_{\id^a}(X) = \sum_{n = 1}^\infty n^a n^{-X}  = \sum_{n = 1}^\infty n^{-(X-a)} = \zeta(X-a).$$
More generally,  for any arithmetic function $f$ and any $a \in \CC$, we write
$$D_{\id^a \cdot f}(X) = \sum_{n = 1}^\infty n^a f(n) n^{-X}  = \sum_{n = 1}^\infty f(n) n^{-(X-a)} = D_f(X-a).$$
Using this notation,  one has 
$$D_{\sigma_a}(X) = \zeta(X) \zeta(X-a).$$
Thus, for example, one has
$$D_{d}(X) = \zeta(X)^2$$
and
$$D_{\sigma}(X) = \zeta(X)\zeta(X-1).$$
\item It is not difficult to show that $\widehat{\phi} = \id$, that is,  that
$$\sum_{d|n} \phi(d) = n$$
for all $n$, where $\phi$ denotes Euler's totient.  From M\"obius inversion, it then follows that $\phi = \widecheck{\id}$, that is,
$$\phi(n) = \sum_{d|n} d \mu\left(\frac{n}{d}\right)=  n \sum_{d|n} \frac{\mu(d)}{d}$$ 
and thus
$$\frac{\phi(n)}{n} = \prod_{p|n} \left( 1-\frac{1}{p} \right)= \sum_{d|n} \frac{\mu(d)}{d}$$ 
for all $n$, both of which are multiplicative functions.   The formal Dirichlet series of $\phi$ and $\frac{\phi}{\id}$ are given by
$$D_\phi(X) = \frac{\zeta(X-1)}{\zeta(X)}$$
and
$$D_{\frac{\phi}{\id}}(X) = D_\phi(X+1) =  \frac{\zeta(X)}{\zeta(X+1)},$$
respectively.
\end{enumerate}
\end{example}

For any arithmetic function $f$ and any prime $p$, the {\bf Bell series of $f$ at $p$}\index{Bell series $B_{f,p}(X)$}\index[symbols]{.k ZA@$B_{f,p}(X)$} is the formal Dirichlet  series
$$B_{f,p}(X) = \sum_{k = 0}^\infty f(p^k)p^{-kX} \in  \operatorname{Dir}_\CC.$$  
For all primes $p$, the map
$$B_{-,p}(X) : \operatorname{Dir}_\CC \longrightarrow \operatorname{Dir}_\CC$$
is a ring homomorphism.  In particular, one has
$$B_{f*g,p}(X)  = B_{f,p}(X) B_{g,p}(X)$$
For all arithmetic functions $f$ and $g$.
For any arithmetic function $f$,  there is a unique multiplicative arithmetic function $\operatorname{Mult}(f)$ that agrees with $f$ at every prime power greater than $1$.  If $f \in 1+\mm_\CC$, that is, if $f(1) = 1$, then, on the level of Dirichlet series, one has
$$D_{\operatorname{Mult}(f)}(X) = \prod_{p} B_{f,p}(X),$$
where the product converges in the norm topology.  

\begin{proposition}[Euler products of multiplicative and completely multiplicative functions]\label{eulerprod} Let $f \in 1+\mm_\CC$.
\begin{enumerate}
\item $f$ is multiplicative if and only if $f = \operatorname{Mult}(f)$, if and only if 
$$D_f(X) = \prod_{p} B_{f,p}(X).$$  
\item $f$ is completely multiplicative if and only if $f$ is multiplicative and 
$$B_{f,p}(X) = \sum_{k = 0}^\infty f(p)^k(p^{-X})^k = \frac{1}{1-f(p)p^{-X}}$$
for all primes $p$.
\item $f$ is completely multiplicative if and only if $$D_f(X) = \prod_{p} \frac{1}{1-f(p)p^{-X}},$$   if and only if $$ D_{f^{*-1}}(X) =\prod_p(1-f(p)p^{-X}),$$   if and only if $$D_{f^{*-1}}(X) =  D_{\mu \cdot f}(X),$$ if and only if $$f^{*-1} = \mu \cdot f.$$
\end{enumerate}
\end{proposition}

For any multiplicative arithmetic function $f$,  the expression $D_f(X) = \prod_{p} B_{f,p}(X)$ is called the {\bf Euler product representation of $D_f(X)$}.\index{Euler product representation}   For example, one has
$$\zeta(X) = \prod_p  \frac{1}{1-p^{-X}},$$
and, more generally,
$$D_{\id^a}(X) = \prod_p \frac{1}{1-p^ap^{-X}}$$
for all $a \in \CC$.  One also has
$$D_{\mu}(X) = \frac{1}{\zeta(X)} = \prod_p (1-p^{-X}),$$
from which one can immediately derive the expression (\ref{mob}) for $\mu(n)$ stated earlier.

The following proposition is easily verified.

\begin{proposition}\label{additivelem}
An arithmetic function $f$ is additive if and only if $\widecheck{f}(n) \neq 0$ implies that $n$ is a prime power greater than $1$, for any positive integer $n$.  
\end{proposition}

Using the proposition above, one can verify the following analogue of Proposition \ref{eulerprod} for  the additive and completely additive arithmetic functions.

\begin{proposition} Let $f$ be an arithmetic function.
\begin{enumerate}
\item For all primes $p$, one has
$$ B_{\widecheck{f},p}(X)  = B_{\mu,p}(X)B_{f,p}(X) = (1-p^{-X}) B_{f,p}(X).$$
\item $f  \in \mm_\CC$ is additive if and only if
$$\frac{D_f(X)} {\zeta(X)} = D_{\widecheck{f}}(X) = \sum_{p} B_{\widecheck{f},p}(X),$$
if and only if 
$$D_{f}(X)  = \zeta(X)  \sum_{p} B_{\widecheck{f},p}(X),$$
 if and only if
$$D_{f}(X)  = \zeta(X)  \sum_{p} (1-p^{-X}) B_{f,p}(X)$$
(where the given sums converges in the norm topology).
\item  $f$ is completely additive if and only if $f$ is additive and
$$ B_{\widecheck{f},p}(X) =   f(p) \sum_{k = 1}^\infty p^{-kX}$$
for all primes $p$, if and only if $f$ is additive and
$$ B_{f,p}(X) =   f(p) \sum_{k = 1}^\infty k p^{-kX},$$
for all primes $p$, 
where also
$$  \sum_{k = 1}^\infty p^{-kX} = \frac{p^{-X}}{1-p^{-X}}$$
and
$$  \sum_{k = 1}^\infty kp^{-kX} = \frac{p^{-X}}{(1-p^{-X})^2}.$$
\item  $f$ is completely additive if and only if
$$D_{f}(X)  = \zeta(X)  \sum_{p} f(p) \frac{p^{-X}}{1-p^{-X}}.$$
\end{enumerate}
\end{proposition}

\begin{example}
$$D_{\omega}(X) =   \zeta(X) \sum_p p^{-X} =  \zeta(X)P(X)$$
and 
$$D_{\Omega}(X) = \zeta(X) \sum_p   \sum_{k = 1}^\infty p^{-kX} = \zeta(X) \sum_{k = 1}^\infty P(kX),$$
where 
$$P(X) = D_{\chi_{\pp}}(X) = \sum_p p^{-X}$$
 is the {\bf formal prime zeta function}, \index{formal prime zeta function $P(X)$} where  $\chi_{\pp}$ is the characteristic function of the set $\pp$  of all primes,
and where one defines  $$D_f(kX)    = \sum_{n = 1}^\infty f(n)(n^k)^{-X}$$
for any arithmetic function $f$ and any positive integer $k$. 
\end{example}

If $f \in \operatorname{Arith}_\CC$ is additive (resp., completely additive), then for any $a \in \CC\backslash \{0\}$ the function $a^f$ given by
$$(a^f)(n) = a^{f(n)} = e^{f(n)\Log a}$$
is multiplicative (resp., completely multiplicative).  On the other hand, if $f \in \operatorname{Arith}_\CC$ is multiplicative (resp., completely multiplicative), then for any $a \in \CC $ the function $a \Log f$ is additive (resp., completely additive), provided that $\operatorname{Re}f(n) >0 $ for all $n$ (which guarantees that $\Log (f(m)f(n)) = \Log f(m) + \Log f(n)$ for all $m$ and $n$).  The following two examples are examples of the first construction.  

\begin{example}
The completely multiplicative arithmetic function $\lambda = (-1)^\Omega$, that is, the arithmetic function $\lambda$ defined by
$$\lambda(n) = (-1)^{\Omega(n)}$$ for all $n$, is called the {\bf Liouville lambda function}.\index{Liouville lambda function $\lambda(n)$}\index[symbols]{.r  N@$\lambda(n)$}  Note that $\lambda(n)$ is  equal to $1$ if $n$ has an even number of prime factors, and equal to $-1$ if $n$ has an odd number of prime factors, counting multiplicities.    It is not difficult to show that $\widehat{\lambda} = \chi_{\ss}$, where $\ss$ is the set of all perfect squares, that is, 
$$\sum_{d|n} \lambda(d) =  \begin{cases} 1 & \quad  \text{if }  n \text{ is a perfect square} \\
 0  & \quad  \text{otherwise}.
\end{cases}$$
Therefore, by M\"obius inversion, one has $\lambda = \widecheck{\chi_{\ss}}$ and therefore
$$\lambda(n) = \sum_{d|n} \mu\left(\frac{n}{d}\right)\chi_{\ss}(d) = \sum_{d^2|n} \mu\left(\frac{n}{d^2}\right).$$
Moreover, the Dirichlet inverse of $\lambda$ is given by
$$\lambda^{*-1} = \widehat{(\chi_{\ss})^{*-1}} =  |\mu| = \mu^2 = \mu \cdot \lambda = \chi_{\mathbb{SF}},$$
where $\mathbb{SF}$ is the set of all squarefree integers.  It follows that 
$$D_{\chi_{\ss}}(X) = \zeta(2X),$$
$$D_{\lambda}(X) =\frac{1}{D_{|\mu|}(X)} = \frac{\zeta(2X)}{\zeta(X)},$$
and thus
 $$D_{|\mu|}(X) = \frac{\zeta(X)}{\zeta(2X)}.$$
\end{example}

\begin{example}
The multiplicative arithmetic function $2^\omega$ at $n$ equals the number $\sum_{d|n} |\mu(d)|$ of squarefree divisors of $n$ for all $n$, and thus $2^\omega = \widehat{|\mu|} = \widehat{ \widehat{(\chi_{\ss})^{*-1}} }$ and $|\mu| =\widecheck{2^\omega}$.    It follows that
$$D_{2^\omega}(X) = \frac{\zeta(X)^2}{\zeta(2X)}.$$
\end{example}

An arithmetic function $f \in 1+\mm_\CC$ is completely multiplicative if and only if $$f^{*-1} = \mu \cdot f,$$
if and only if 
$$f \cdot(g*h) = (f \cdot g)*(f \cdot h)$$
for all $g, h \in \operatorname{Arith}_\CC$.   The extensive remark below shows that there is a natural structure of a commutative ring on the group $\operatorname{Mult}_\CC$, where addition in the ring is Dirichlet convolution $*$,  where multiplication is denoted $\odot$, where $\zeta$ is the $1$ of the ring and thus $\mu$ is the $-1$ of the ring, and where $f \odot g = f \cdot g$ if either $f$ or $g$ is completely multiplicative.     This provides a ring-theoretic interpretation of the arithmetic functions $\zeta$ and $\mu$, and  also of the two identities noted above for completely multiplicative arithmetic functions $f$.

\begin{remark}[The rings of multiplicative and additive arithmetic functions {\cite{ell2}}]\label{muladd}   The definitions of $\operatorname{Arith}_R = R^{\ZZ_{> 0}}$,  $\operatorname{Mult}_R$, and  $\operatorname{Add}_R$ generalize in the obvious way to any commutative ring $R$.     In 2008, the author proved that there exists a uniqiue functorial ring structure, that is functorial in the obvious way (so that a ring homomorphism $\varphi: R \longrightarrow S$ sends $f \in \operatorname{Mult}_R$ to $\varphi \circ f \in \operatorname{Mult}_S$), on the group $\operatorname{Mult}_R$ satisfying properties (1)--(4) below for all $f, g,h \in \operatorname{Mult}_R$, where multiplication in the ring $\operatorname{Mult}_R$ is denoted $\odot$.
\begin{enumerate}
\item Addition in the ring $\operatorname{Mult}_R$ is Dirichlet convolution $*$.
\item The binary operation   $\odot$ on $\operatorname{Mult}_R$ is of {\bf Dirichlet type}: for any positive integer $n$, the value of $f \odot g$ at $n$ depends on $f(d)$ and $g(d)$ only  for the positive integers $d$ dividing $n$.
\item If $f$ and $g$ are completely multiplicative (or, in fact, if either $f$ or $g$ is completely multiplicative), then $f \odot g = f \cdot g$.
\item $\zeta$ is the $1$ of the ring  $\operatorname{Mult}_R$, that is, it is the identity element  with respect to $\odot$.  Equivalently, one has $\zeta \odot f = f$ for all $f$.
\end{enumerate}
One also has the following.
\begin{enumerate}
\item[(5)]  $\mu$ is  the $-1$ of the ring $\operatorname{Mult}_R$, that is, it is the additive inverse of the $1$ of the ring.  Equivalently, one has $\mu \odot f = f^{*-1}$ for all $f$.  Thus, if $f$ is completely multiplicative, then $\mu \cdot f = \mu\odot f = f^{*-1}$.
\item[(6)] Since $\mu \odot \mu = \zeta \neq \mu \cdot \mu$ if $R$ is not the trivial ring, one need not have $f \odot g = f \cdot g$.
\item[(7)] $f \odot(g*h) = (f \odot g)*(f \odot h)$.
\item[(8)] If $f$ is completely multiplicative, then
$$f \cdot(g*h)  =  f \odot(g*h) = (f \odot g)*(f \odot h) =  (f \cdot g)*(f \cdot h).$$
\end{enumerate}
It is not easy  to give an explicit  equation for $(f \odot g)(n)$ for all $n$ in terms of $f(d)$ and $g(d)$ for the positive integers $d$ dividing $n$.  Instead, we can at least give an implicit definition of the operation $\odot$, as follows.  First, we describe a natural ring structure on the set $\operatorname{Add}_R$.  By Proposition \ref{additivelem}, the set $$\mu* \operatorname{Add}_R =  \{\widecheck{f} = \mu*f: f \in \operatorname{Add}_R\}$$ naturally forms a ring under pointwise addition and pointwise multiplication isomorphic to the ring $ \prod_{n = 1}^\infty \prod_p R$ under coordinatewise addition and coordinatewise multiplication.  Since the map
$$\widecheck{-} = \mu*-: \operatorname{Add}_R \longrightarrow \mu* \operatorname{Add}_R$$ is a bijection with inverse $\widehat{-} = \zeta*-$,  we can endow $\operatorname{Add}_R$ with the unique ring structure so that the given bijection is a ring isomorphism.  Equivalently, we endow $\operatorname{Add}_R$  with the unique ring structure so that the bijection
$$ \begin{aligned}   \operatorname{Add}_R &  \longrightarrow \prod_p \prod_{n = 1}^\infty R \\   f & \longmapsto  ((f(p^n)-f(p^{n-1}))_{n = 1}^\infty)_p
\end{aligned}$$
is an isomorphism of rings, where the direct product has coordinatewise addition and coordinatewise multiplication.   Explicitly, addition in the ring $\operatorname{Add}_R$  is pointwise addition $+$, and  if we denote multiplication in the ring $\operatorname{Add}_R$ by $\otimes$, then one has
\begin{align*}
{f \otimes g}  & = \zeta *( (\mu * f) \cdot( \mu * g))  \\ & ={ \widehat{{\widecheck{f}} \cdot{ \widecheck{g}}}}\index[symbols]{.ia tw@$f \otimes g$}
\end{align*}
for all $f, g \in \operatorname{Add}_R$.  Clearly the multiplicative identity of the ring $\operatorname{Add}_R$ is the completely additive function $\Omega$.    One has $$ \mu* \operatorname{Add}_R = \{f \in \operatorname{Arith}_R: f = \chi_{\mathcal{P}}\cdot f\},$$ where $\mathcal{P}$ is the set of all prime powers greater than $1$.  Thus, the function $\chi_{\mathcal{P}}$, with formal Dirichlet series
$$D_{\chi_{\mathcal{P}}}(X) = \sum_p \sum_{k = 1}^\infty (p^k)^{-X} =  \sum_{k = 1}^\infty P(kX) = \frac{D_\Omega(X)}{\zeta(X)},$$
  is the multiplicative identity element of the ring $ \mu* \operatorname{Add}_R$ (with respect to $\cdot$), and, by Proposition \ref{additivelem}, an arithmetic function $f$ is additive if and only  if $\widecheck{f} = \chi_{\mathcal{P}} \cdot \widecheck{f}$, if and only if $f =\widehat{ \chi_{\mathcal{P}} \cdot \widecheck{f}}$.

It turns out that the map
\begin{align*}
\widehat{\operatorname{DLog}}: \begin{aligned}   \operatorname{Mult}_R & \longrightarrow \operatorname{Add}_R \\  f & \longmapsto  \zeta*(\Omega\cdot f) * f^{*-1}\index[symbols]{.ia  u@$\widehat{\operatorname{DLog}}(f)$}
\end{aligned}
\end{align*}
is a homomorphism of abelian groups, it is injective if $R$ is $\ZZ$-torsion-free, and it is an isomorphism if $R$ is a field.  The operation $\odot$ is defined to be the unique operation of muliplication on $\operatorname{Mult}_R$ so that the map $\widehat{\operatorname{DLog}}$ is a homomorphism of rings for every ring $R$.  Equivalently, $*$ and $\odot$ are the unique operations  on $\operatorname{Mult}_R$ so that the map
\begin{align*}
\operatorname{DLog}: \begin{aligned}   \operatorname{Mult}_R &  \longrightarrow \mu* \operatorname{Add}_R \\   f & \longmapsto  (\Omega\cdot f) * f^{*-1}\index[symbols]{.ia  ua@$\operatorname{DLog}(f)$}
\end{aligned}
\end{align*}
is a ring homomorphism for every ring $R$.   Thus, one has 
$$f * g =  \operatorname{DLog}^{-1}(\operatorname{DLog}(f) + \operatorname{DLog}(g))$$
and 
$$f \odot g =  \operatorname{DLog}^{-1}(\operatorname{DLog}(f) \cdot  \operatorname{DLog}(g))\index[symbols]{.ia tv@$f \odot g$}$$
for all $f, g \in \operatorname{Mult}_R$, if $R$ is $\ZZ$-torsion-free.   The rings  $\operatorname{Mult}_R$ and $\operatorname{Add}_R$ are called the {\bf ring of multiplicative arithmetic functions}, and the {\bf ring of additive arithmetic functions}, respectively, {\bf over $R$}.\index{ring of multiplicative arithmetic functions $\operatorname{Mult}_R$}\index{ring of additive arithmetic functions $\operatorname{Add}_R$} 

The operation $\operatorname{DLog}$ on $\operatorname{Mult}_R$, or more generally its extension $f \longmapsto  (\Omega\cdot f) * f^{*-1}$ to the group  $\operatorname{Arith}_R^*$,  is like a logarithmic derivative, because the operation $\partial$ on the ring $\operatorname{Arith}_R$ acting by   $$\partial: f \longmapsto \Omega \cdot f,$$  which is called the {\bf standard derivation of $\operatorname{Arith}_R$},  is analogous to differentiation.\index{standard derivation $\partial$  of $\operatorname{Arith}_R$}\index[symbols]{.ia  u@$\partial(f)$}    As a {\it derivation of $\operatorname{Arith}_R$ over $R$}, the map $\partial$ satisfies the following conditions for all $f, g \in \operatorname{Arith}_R$.
\begin{enumerate}
\item $\partial(a \iota) = 0$ for all $a \in R$.
\item  $\partial(f+g) = \partial(f)+\partial(g)$.
\item $\partial(f*g) = f *\partial(g) + \partial(f) *g$.
\end{enumerate}
More generally, if $h$ is any arithmetic function in $\operatorname{Arith}_R$, then $h$ is completely additive if and only if the operation $\partial_h: f \longmapsto h \cdot f$ on $\operatorname{Arith}_R$ is a derivation of $\operatorname{Arith}_R$ over $R$, that is, if and only if $\partial_h$ satisfies the three conditions above.  Thus, one has $\partial = \partial_\Omega$, and the ring homomorphism $$\widehat{\operatorname{DLog}}: \operatorname{Mult}_R \longrightarrow \operatorname{Add}_R$$ from the ring $\operatorname{Mult}_R$ of all multiplicative arithmetic functions  (under the operations $*$ and $\odot$) to the ring $\operatorname{Add}_R$ of all additive arithmetic functions (under the operations $+$ and $\otimes$) acts by $$\widehat{\operatorname{DLog}}: f \longmapsto  \widehat{\operatorname{DLog}(f)} =  \widehat{\partial(f)*f^{*-1}}.$$

It is clear that the map
\begin{align*}
 \begin{aligned}  \operatorname{Mult}_R & \longrightarrow \prod_p (1+TR[[T]]) \\   f &  \longmapsto \left (\sum_{k= 0}^\infty f(p^k)T^k\right)_p
\end{aligned}
\end{align*}
is an isomorphism of abelian groups.  The abelian group $\Lambda(R)  = 1+TR[[T]]$ has the structure of a  {\it lambda ring} and is known as the  {\it universal lambda ring over $R$} \cite{yau}.  The group isomophism above is a ring isomorphism when $1+TR[[T]]$ is endowed with this ring  structure and $\prod_p (1+TR[[T]])$ the direct product ring structure.  By transport of structure, then, this isomorphism provides an alternative way to construct the ring structure on the group $\operatorname{Mult}_R$.   It follows that the ring  $ \operatorname{Mult}_R $, too, has the structure of a lambda ring, being naturally isomorphic to the ring $\prod_p\Lambda(R) \cong \Lambda( \prod_p R)$.  Moreover, the ring homomorphism $\widehat{\operatorname{DLog}}$  is analogous to the {\it ghost homomorphism} from the ring $\Lambda(R)$ to the ring $\prod_{n= 1}^\infty R$.   Note that the ring $\Lambda(R)$ can be recovered, up to isomorphism,  from the ring $\operatorname{Mult}_R$ as its subring
$$\{f \in \operatorname{Mult}_R: f(p^n) = f(q^n) \text{ for all primes } p,q \text{ and all } n \in \ZZ_{> 0} \} \cong \Lambda(R).$$

See \cite[Section 11]{ell2}  for an application of the results above to the study of schemes of finite type over $\ZZ$.  Another application is to the study of powers of arithmetic functions, described below.

For all $f \in 1+\mm_\CC$ and all $a \in \CC$, define
$$f^{*a} = \sum_{n = 0}^\infty {a \choose n} (f-1)^n,$$
where $$ {a \choose n} = \frac{a(a-1)(a-2)\cdots (a-n+1)}{n!}$$
for all $n$.  The arithmetic function $f^{*a} \in 1+\mm_\CC$  is well-defined, since the series above converges in the norm topology.  Moreover, if $a \in \QQ$, then $f^{*a}$  is interpreted algebraically as the unique $a$th Dirichlet power of $f$, so, for example,  $f^{*2} = f*f$, while $f^{*1/2}$ is the unique arithmetic function $g  \in 1+\mm_\CC$ such that $g * g = f$, and $f^{*-1}$ is the Dirichlet inverse of $f$.  The ``scalar multiplication'' operation $(a,f) \longmapsto f^a$ gives the abelian group $1+\mm_\CC$ the structure of a complex vector space.  Moreover, the group $\operatorname{Mult}_\CC$ is a $\CC$-subspace of $1+\mm_\CC$ and in fact has the structure of a $\CC$-algebra.  For $f \in \operatorname{Mult}_\CC$, an equivalent definition of   $f^{*a}$ is 
$$f^{*a} = \zeta^{*a} \odot f,$$ 
where, by the generalized binomial theorem, one has
$$D_{\zeta^{*a}}(X) = \prod_p \left(1-p^{-X}\right)^{-a} = \prod_p \left(\sum_{n = 0}^\infty (-1)^n {-a \choose n}p^{-nX}\right),$$
and therefore $\zeta^{*a}$ is the unique multiplicative function such that
\begin{align*}
\zeta^{*a}(p^n) & = (-1)^n {-a \choose n}  = \frac{(a)_n}{n!}
\end{align*}
for all prime powers $p^n$, where $(a)_n = a(a+1)(a+2)\cdots(a+n-1)$ denotes the Pochhammer symbol.   The map $a \longrightarrow \zeta^{*a}$ is an embedding of the field $\CC$ in the ring $\operatorname{Mult}_\CC$.
It follows that $d = \zeta * \zeta$ is the $2$, or $1+1$, of the ring $\operatorname{Mult}_\CC$.  More generally, for any positive integer $k$, the function $d_k = \zeta^{*k}$ is the $k = 1+1+\cdots+1$ of the ring $\operatorname{Mult}_\CC$, where $$d_k(n) = \sum_{n = a_1 a_2 \cdots a_k} 1$$ for any $n$ is the number of ways of writing $n$ as a product of $k$ factors (where order matters), and, for any $a \in \CC$, the series $\zeta^{*a}$ is the ``$a$'' of the $\CC$-algebra $\operatorname{Mult}_\CC$.  Moreover, $a\Omega$, which is the completely additive function that weights every prime $p$ by $a$,  is the ``$a$'' of the $\CC$-algebra $\operatorname{Add}_\CC$.

An immediate consequence of the isomorphisms $\widehat{\operatorname{DLog}}$ and $\operatorname{DLog}$ is that, for any multiplicative arithmetic function $f$, the function
$$f^{*a} =  \widehat{\operatorname{DLog}}^{-1}(a \, \widehat{\operatorname{DLog} }f) =  \operatorname{DLog}^{-1}(a \operatorname{DLog} f)$$
is multiplicative for all $a \in \CC$.  
Moreover, if $f$ is completely multiplicative, then $f^{*a} = \zeta^{*a} \cdot f$, and therefore $f^{*a}$ is the unique multiplicative function such that
$$f^{*a}(p^n) = \frac{(a)_n}{n!} f(p)^n$$
for all prime powers $p^n$.  Thus, for example, one has
$$\mu(p^n) = \zeta^{*-1}(p^n) =   \frac{(-1)_n}{n!}   = \begin{cases}  1 & \quad \text{if }  n = 0 \\
   -1  & \quad \text{if } n = 1 \\
   0 & \quad \text{if } n > 1 
\end{cases}$$
and, if $f$ is completely multiplicative and $a = -1$, then
\begin{align*}
f^{*-1}(p^n)  & = \begin{cases}  1 & \quad \text{if }  n = 0 \\
   -f(p)  & \quad \text{if } n = 1 \\
   0 & \quad \text{if } n > 1 
\end{cases} \\
& = \mu(p^n) f(p^n),
\end{align*}
from which we recover the fact that $f^{*-1} = \mu \cdot f$.
\end{remark}

In Section 3.8, we study Dirichlet series as complex functions rather than as formal Dirichlet series.

\section{Summatory functions}

As noted in Section 1.2, there are two natural ways to ``analysis-ize'' an arithmetic function $f$.   One way,  discussed in Section 3.8,  is to form its {\bf Dirichlet series} $D_f(s) = \sum_{n = 1}^\infty \frac{f(n)}{n^s}$ as a function of a complex variable $s$.   Another way  is to form the {\bf summatory function $S_f$ of $f$},\index{summatory function $S_f(x)$}\index[symbols]{.k Fr@$S_f(x)$}  which is the function $S_f: \RR_{\geq 0} \longrightarrow \CC$ given by
$$S_f(x) = \sum_{n \leq x} f(n)$$ for all $x \geq 0$, which is defined, more generally, for any function $f$ whose domain contains $\ZZ_{>0}$.   The summatory function $S_f$ of $f$ is a step function that is constant on the interval  $[n-1,n)$, and jumps by the value $+f(n)$ at $n$, for every positive integer $n$.  
Note, in particular, that $S_f$ is completely determined by the relations \begin{align*}S_f(0) = 0,\end{align*}
\begin{align*}S_f(x) = S_f(\lfloor x \rfloor), \quad \forall x \geq 0,\end{align*}
and \begin{align*}f(n) = S_f(n)-S_f(n-1), \quad \forall n \in \ZZ_{> 0}.\end{align*}

\begin{example} Let $x \geq 0$.  One has the following.
\begin{enumerate}
\item  $S_\zeta(x) = \lfloor x \rfloor$, where $\zeta(n) = 1$  for all $n$.
\item $S_{\frac{1}{\id}}(x)  = \sum_{n \leq x} \frac{1}{n} = H_{\lfloor x \rfloor}$, where  $H_n$ for any nonnegative integer $n$ denotes  the $n$th harmonic number $H_n =  \sum_{k = 1}^n \frac{1}{k}.$
\item  $S_{\chi_{\pp}}(x) = \pi(x)$, where $\chi_{\pp}$ is the characteristic function of the set $\pp$ of all prime numbers. 
\item  $S_\mu(x)  = M(x)$, where $\mu$ is the M\"obius function and $M$ is the famous {\bf Mertens function}\index{Mertens function $M(x)$}\index[symbols]{.st M@$M(x)$}
$$M(x) = \sum_{n \leq x} \mu(n).$$
\item  $S_\lambda(x)  = L(x)$, where $\lambda$ is the Liouville lambda function and $L$ is the  {\bf summatory Liouville  function}\index{summatory Liouville function $L(x)$} \index[symbols]{.st N@$L(x)$}
$$L(x) = \sum_{n \leq x} \lambda(n).$$
\end{enumerate}
\end{example}

If $F: \RR_{\geq 0} \longrightarrow \CC$ is any function, then  $\Delta(F)$ denotes the arithmetic function defined by
\begin{align*}\Delta(F)(n) = F(n)-F(n-1), \quad \forall n \in \ZZ_{> 0}.\end{align*}
One then has
$$\Delta(S_f) = f$$
for any arithmetic function $f$.  Moreover, for any function  $F: \RR_{\geq 0} \longrightarrow \CC$, one has
\begin{align*}S_{\Delta(F)}(x) =   F(\lfloor x \rfloor) -F(0)
\end{align*}
and therefore
\begin{align*}S_{\Delta(F)} = F \quad \Longleftrightarrow  \quad F(0) = 0 \text{ and } F(x) =F(\lfloor x \rfloor), \quad \forall x \geq 0.\end{align*}
The operations $S_{-}$ and $\Delta$ are discrete analogues of integration and differentiation, respectively, and the facts noted above are the discrete analogues of the fundamental theorem of calculus.   Generally, for  a given ``naturally'' defined  arithmetic function $f$, one expects the function $S_f(x)$ to be ``less erratic'' than $f$,  and more amenable to asymptotic  analysis, analogous to how an antiderivative of a function, when it exists, is  ``smoother'' than the original function.  One reason for this is that $\frac{1}{x}S_f(x) = \frac{1}{x}\sum_{n \leq x}f(n)$ represents a running average of $f$, and any erratic behavior in a function tends to be somewhat abated in its running average.

There are many well-known techniques in analytic number theory for studying summatory functions, several of which  are discussed in subsequent sections of this chapter.  One of the most elementary is the following result, which generalizes  \cite[Theorem 1]{apostol} and the method outlined in \cite[Section 1.5 pp.\ 19--20]{iwan}.  It implies that the summatory function $S_g(x)$ of any nonnegative monotonic arithmetic function $g$ is well-approximated by the function $\int_1^x f(t) \, dt$ for any nonnegative monotonic extension $f: [1,\infty) \longrightarrow \RR$ of  $g$ (e.g., for $f(x) = g(\lfloor x\rfloor)$).

\begin{proposition}\label{iwanncor}
Let $f: [N,\infty) \longrightarrow \RR$ be any nonnegative monotonic function on $[N,\infty)$, where $N$ is a positive integer,  and let
$$S_{f,N}(x) = \sum_{N \leq k \leq x} f(k)$$
and  $$F(x)= \int_N^x f(t) \, dt$$
for all $x \geq N$.
\begin{enumerate}
\item If $f$ is nondecreasing on $[N,\infty)$,  then 
$$f(N) \leq S_{f,N}(n) -F(n) \leq f(n)$$
for all  integers $n \geq N$  and
$$\left|S_{f,N}(x) -F(x) \right|  \leq f(x)$$
for all $x \geq N$.
\item If $f$ is nonincreasing on $[N,\infty)$,  then the limits $L = \lim_{x \to \infty} f(x) \geq 0$ and
$$0\leq \gamma_{f,N} =\lim_{n \to \infty} \left(\sum_{k = N}^{n-1}f(k)- \int_N^n f(t) \, dt \right) \leq f(N)$$
exist, and one has
$$0 \leq S_{f,N}(n) -F(n) - \gamma_{f,N} \leq f(n)$$
for all  integers $n \geq N$ and
$$\left|S_{f,N}(x) - F(x) - \gamma_{f,N} \right |  \leq f(\lfloor x\rfloor)$$
for all $x \geq N$,  and therefore
$$0 \leq \limsup_{x \to \infty} \left|S_{f,N}(x) -F(x) - \gamma_{f,N} \right | \leq L.$$
\item If $f$ is nondecreasing on $[N,\infty)$ and $f(x) = o(F(x)) \ (x \to \infty)$, then  $$S_{f,N}(x) \sim F(x) \ (x \to \infty).$$
\item If $f$ is nonincreasing  on $[N,\infty)$ and  $\lim_{x \to \infty} f(x) = 0$, then 
$$\gamma_{f,N} =\lim_{x \to \infty} \left(S_{f,N}(x)- \int_N^x f(t) \, dt \right).$$
\item If $f$ is nonincreasing on $[N,\infty)$, $\lim_{x \to \infty} f(x) = 0$,  and $\int_N^\infty f(t)\, dt$ exists,
then `
$$\sum_{k = N}^\infty f(k) = \gamma_{f,N} +\int_N^\infty f(t) \, dt,$$
and, if also $f(\lfloor x \rfloor) = o\left(\int_x^\infty f(t)\, dt\right) \ (x \to \infty)$,  then
$$\sum_{k > x} f(k) \sim \int_x^\infty f(t)\, dt \ (x \to \infty).$$
\item If $f$ is increasing (resp., decreasing) on $[N,\infty)$, then the inequalities in (1) (resp., (2)) are strict for all $n >N$ and all $x > N$.
\end{enumerate}
\end{proposition}

\begin{proof}
Note that any monotonic function on a closed interval has only countably many discontinuities and is therefore Riemann integrable.   By a simple linear change of variables, we may assume without loss of generality that $N = 1$.   Suppose that $f$ is nondecreasing.  Then one has
$$0 \leq \sum_{k = 2}^{n}f(k)- \int_1^n f(t) \, dt \leq \sum_{k = 2}^n (f(k)-f(k-1)) = f(n)-f(1)$$
for all positive integers $n$, and thus
$$f(1)-\int_{\lfloor x \rfloor}^x f(t) \, dt \leq S_f(x)- \int_1^x f(t) \, dt \leq f(\lfloor x \rfloor)-\int_{\lfloor x \rfloor}^x f(t) \, dt$$
for all $x \geq 1$.
Therefore, since
$$\{x\}f(\lfloor x \rfloor) \leq \int_{\lfloor x \rfloor}^x f(t) \, dt \leq \{x\} f(x),$$
one has
$$-f(x) \leq f(1) -\{x\} f(x) \leq  S_f(x)- \int_1^x f(t) \, dt \leq f(\lfloor x \rfloor)-\{x\}f(\lfloor x \rfloor) \leq f(x),$$
for all $x \geq 1$.  Statement (1) follows.

Suppose now, that $f$ is nonincreasing.  Since $f$ is nonnegative,  the limit  $\lim_{x \to \infty} f(x)$ exists, and  \cite[pp.\ 409--410]{apostol} (which is a simple geometric argument) implies that the limit $\gamma_{f,N}$ exists.  Let
$$E_f(x) = S_f(x)- \int_1^x f(t) \, dt- \gamma_{f,N}$$ for all $x \geq 1$.  Again by \cite[pp.\ 409--410]{apostol}, one has
$$0 \leq E_f(n) \leq f(n)$$
for all positive integers $n$, and therefore
$$-f( \lfloor x \rfloor ) \leq -\{x\}f( \lfloor x \rfloor )\leq E_f(x) = E_f(\lfloor x \rfloor)-\int_{\lfloor x \rfloor}^x f(t) \, dt \leq f(\lfloor x \rfloor)-\{x\}f(x) \leq f(\lfloor x \rfloor)$$ 
for all $x \geq 1$.   Statement (2) follows.  Finally,   statements (3)--(5) follow  from (1) and  (2),  and statement (6) follows by modifying the proof of (1) and (2).
\end{proof}

Note that, if $f: [1,\infty) \longrightarrow \RR$ is defined on $[1,\infty)$, then $S_{f,N}(x) = S_f(x)-S_f(N-1)$ for all $x \geq 1$,  and $S_f = S_{f,1}$.   If $\gamma_{f,1}$ exists, then we write $\gamma_f = \gamma_{f,1}$.

\begin{figure}[ht!]
\includegraphics[width=80mm]{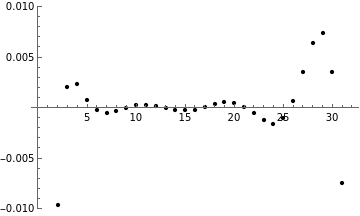}
\caption{\centering Plot of $\gamma_k$ for $k  = 2,3,4,\ldots, 31$}
  \label{stielt}
\end{figure}

\begin{example}\label{stiec} \
\begin{enumerate}
\item One has
$$0 < H_n - \log n - \gamma < \frac{1}{n}$$
for all positive integers $n$ and
$$\left|H_{\lfloor x \rfloor} -\log x - \gamma \right |  < \frac{1}{\lfloor x\rfloor}$$
for all $x \geq 1$, where $\gamma = \gamma_{\frac{1}{\id}}$.  Consequently, one has
$$H_{\lfloor x \rfloor} = \log x + \gamma +  O\left( \frac{1}{x}\right) \ (x \to \infty).$$
\item  \cite{farr}.  Example (1) generalizes as follows: for any real number $a \geq 0$, one has
$$\sum_{n \leq x} \frac{(\log n)^a}{n} =  \frac{(\log x)^{a+1}}{a+1}  + \gamma_a + O\left(\frac{(\log x)^a}{x}\right) \ (x \to \infty),$$
where 
\begin{align*}
\gamma_a & = \lim_{x \to \infty} \left(\sum_{n \leq x} \frac{(\log n)^a}{n} -  \int_1^x \frac{(\log t)^a}{t}\, dt  \right) \\& =  \lim_{x \to \infty} \left(\sum_{n \leq x} \frac{(\log n)^a}{n} -  \frac{(\log x)^{a+1}}{a+1} \right),
\end{align*}
where ($0^0 = 1$ and) $\gamma_0 = \gamma$ and  
$$\lim_{a \to 0^+} \gamma_a  = \gamma -1.$$  For any nonnegative integer $k$, the constant $\gamma_k$ is known as the {\bf $k$th  Stieltjes constant} \cite[(1.12)]{ivic}.\index{Stieltjes constant $\gamma_k$}\index[symbols]{.f tf@$\gamma_k$} See Figure \ref{stielt} for a plot of $\gamma_k$ for $k = 2,3,4,\ldots,31$.   The constant $\gamma_a$ for $a \geq 0$ is known as the {\bf $a$th fractional Stieltjes constant},\index{fractional Stieltjes constant $\gamma_a$}\index[symbols]{.f tfa@$\gamma_a$}  and the function $\gamma_a$ is discontinuous at $a = 0$ but infinitely differentiable on $(0,\infty)$ and thus smoothly interpolates the constant $\gamma_0-1 = \gamma-1$ and the Stieltjes constants $\gamma_k$ for all positive integers $k$.  
 Moreover,  for all $a \geq 0$ and all $s \in \CC$ with $\operatorname{Re} s > 1$, one has $$ \sum_{n =1}^\infty \frac{(\log n)^a}{n^s} =\frac{\Gamma(a+1)}{(s-1)^{a+1}} +\sum_{n = 0}^\infty \frac{(-1)^n \gamma_{a+n}}{n!}(s-1)^n,$$
where $\Gamma(s)$ denotes the {\it gamma function} (defined in Section 4.1).  See Example \ref{sumex}(3)  and (\ref{stcon}) for further applications of the Stieltjes constants $\gamma_k$ to analytic number theory.
\item  \cite[p.\ 20]{iwan}.  For any nonnegative integer $k$, one has
$$\sum_{n \leq x} (\log n)^k  = x P_k(\log x) + O((\log x)^k) \ (x \to \infty),$$
where 
$$P_k(X) =  \sum_{l = 0}^k (-1)^{k-l} \frac{k!}{l!} X^l \in \ZZ[X].$$
Thus, for example, one has
$$\sum_{n \leq x} \log n  = x \log x -x + O(\log x) \ (x \to \infty).$$
\end{enumerate}
\end{example}

From Proposition \ref{iwanncor} and Karamata's integral theorem, we obtain the following.

\begin{corollary}\label{karcor}
Let $f: [N,\infty)\longrightarrow \RR$ be continuous, nonnegative, monotonic, and  regularly varying of index $r \in \RR$, where $N$ is a positive integer.  If $r > -1$, then one has
$$\sum_{N \leq k\leq x} f(k) \sim \int_N^x f(t)\, dt \sim \frac{x f(x)}{r+1} \gg 1  \ (x \to \infty).$$  On the other hand, if $r < -1$, then $$ \sum_{k = N}^\infty f(k) = \gamma_{f ,N}+ \int_N^\infty f(t) \, dt$$ exists, and one has
$$\sum_{k > x} f(k)   \sim  \int_x^\infty f(t)\, dt  \sim -\frac{x f(x)}{r+1} \ll 1 \ (x \to \infty).$$ 
\end{corollary}

\begin{remark}[Measure-theoretic interpretation of summatory functions]\label{measarith}
 The summatory function $S_f(x)$ of an arithmetic function $f$ can be interpreted as a {\bf discrete integral of $f$}\index{discrete integral} in the sense that $$S_f(x) = \int_{[1,x]} f(t) \, d \nu(t)$$ for all $x \geq 1$, where $f(t)$ is any extension  of $f$ to $[1,\infty)$,  and where $$\nu = \sum_{n = 1}^\infty \delta_n$$ is the unique {\it discrete measure on the power set of $\RR$ that is supported on $\ZZ_{>0}$ with all weights equal to $1$}, called the {\bf integer measure},\index{integer measure  $\nu$}  where
$\delta_a$ for all $a \in \RR$ is the {\bf Dirac meaure}\index{Dirac measure $\delta_a$} defined by
$$\delta_a(X) = \begin{cases} 1 & \quad \text{if } a \in X \\
 0 & \quad \text{if } a \notin X 
\end{cases}
$$
for any subset $X$ of $\RR$.      In general, one has
$$\nu(X) = |X\cap \ZZ_{>0}|$$
and, more generally, 
$$\int_X f \, d\nu = \sum_{n =1}^\infty f(n) \chi_X(n) = \sum_{n \in X \cap \ZZ_{>0}} f(n),$$
for any subset $X$ of $\RR$ and any function $f$ for which the sum $\sum_{n \in X\cap \ZZ_{>0}} |f(n)|$ is finite, where $\chi_X$ is the characteristic function of $X$.   Thus, for example, the summatory function $$S_{\frac{1}{\id}}(x) = H_{\lfloor x \rfloor} =  \int_{[1,x]} \frac{1}{t} \, d\nu(t)$$ of $\frac{1}{\id}$ is a discrete integral of $\frac{1}{\id}$ and is in this sense the {\bf discrete natural logarithm}. 
\end{remark}

\section{Inversion theorems}

As mentioned in Section 3.3, Proposition \ref{mobius1} (the M\"obius inversion theorem for Dirichlet convolution) is one of many inversion theorems.    In this section, we discuss a few others.

We begin with a useful application of Proposition \ref{mobius1}.

\begin{example}
For any prime power $q > 1$, there is a field of $q$ elements, which is unique up to unique isomorphism, and it is denoted $\FF_q$, or  sometimes $\operatorname{GF}(q)$, where ``GF'' stands for ``Galois field.''   Any finite field must have characteristic $p$ for some prime number $p$ and therefore must be a finite dimensional vector space over the field $\FF_p = \ZZ/p\ZZ$.  Therefore any finite field has $q$ elements for some prime power $q>1$, and so the finite fields, up to isomorphism, are
$$\FF_2, \, \FF_3, \, \FF_4, \, \FF_5, \, \FF_7, \, \FF_8, \, \FF_9, \, \FF_{11}, \, \FF_{13}, \, \FF_{16}, \, \FF_{17}, \, \FF_{19}, \, \FF_{23}, \, \FF_{25}, \, \FF_{27}, \,\FF_{29}, \, \FF_{31}, \, \FF_{32}, \, \FF_{36}, \, \FF_{37}, \, \ldots.$$ 
Moreover, $\FF_{q'}$ is isomorphic to a (unique) subfield of $\FF_{q}$ if and only if $q'$ divides $q$.  Fix any prime power $q> 1$.   Then the product of all monic irreducible polymomials in the polynomial ring $\FF_p[X]$ of degree dividing $n$ is equal to $X^{q^n}-X$, for any integer $n$.  Therefore, if we denote by
$N_n(q)$ the number of monic irreducible polymomials in the polynomial ring $\FF_q[X]$ of degree $n$, then one has
$$q^n = \sum_{d|n} dN_d(q)$$
for all $n$.
From M\"obius inversion, then, it follows that
$$N_n(q) = \frac{1}{n} \sum_{d|n} \mu\left(\frac{n}{d}\right) q^d$$
for all $n$.  The polynomial $$N_n(X) = \frac{1}{n} \sum_{d|n} \mu\left(\frac{n}{d}\right) X^d \in \QQ[X]$$ is called the  {\bf  $n$th necklace polynomial}\index{necklace polynomials}  and also arises naturally in  other contexts, both in algebra and in combinatorics.   The name ``necklace polynomial'' was proposed by Metropolis and Rota in \cite{rota2} because the $N_n(k)$ for any positive integers $n$ and $k$ is equal to the number of necklaces formed  by $n$ colored beads, from $k$ available colors, that have the property that all of their cyclic rotations are distinct.  From this combinatorial interpretation it follows that polynomial $N_n(X)$ assumes integer values at all of the integers.  For example, if $n = p$ is prime, then $N_p(X) = \frac{1}{p}(X^p-X)$, and thus  the fact that $N_p(X)$ is integer-valued at the integers  is equivalent to Fermat's little theorem.  More generally, if $n = p^k > 1$ is a power of a prime $p$, then $$N_{p^k}(X) = \frac{1}{p^k}\left(X^{p^k}-X^{p^{k-1}}\right),$$
is integer-valued at the integers, which is equivalent to
$$a^{p^k}\equiv a^{p^{k-1}}  \ (\text{mod } p^k), \quad \forall a \in \ZZ.$$
\end{example}

The following result, which (to distinguish it from several other inversion theorems) we call the {\bf  first Dirichlet inversion theorem},\index{Dirichlet inversion theorem} is an obvious generalization of Proposition \ref{mobius1}.

\begin{proposition}[First Dirichlet inversion theorem]
Let $h$ be any arithmetic function possessing a Dirichlet inverse $h^{*-1}$, that is, satisfying $h(1) \neq 0$.  
For any arithmetic functions $f$ and $g$, one has
$$g(n)=\sum_{d|n} h(d)f\left({\frac {n}{d}}\right), \quad  \forall n\geq 1,$$
if and only if
$$f(n)=\sum_{d|n}h ^{*-1}(d)g\left({\frac {n}{d}}\right), \quad  \forall n\geq 1.$$
\end{proposition}

Let $f$ and $g$ be arithmetic functions, and let $F, G \in \CC^{\RR_{\geq 1}}$ be arbitrary functions from $\RR_{\geq 1}$ to $\CC$.  We define
$$(f*F)(x) = \sum_{n\leq x}f (n)F\left({\frac {x}{n}}\right), \quad  \forall x\geq 1,$$
which is also a function in $\CC^{\RR_{\geq 1}}$.  Note that 
$$f*(F+G) = f*(F+G),$$
$$(f+g)*F = f*F+g*F,$$
$$\iota * F = F,$$
and $$f*1 = S_f.$$
Moreover,  for all $x \geq 1$ one has
\begin{align*}
(g*(f*F))(x) & = \sum_{n\leq x}g (n)(f*F)\left({\frac {x}{n}}\right)  \\ 
& = \sum_{n\leq x}\sum_{m\leq \frac{x}{n}}g (n)f (m)F\left({\frac {x}{nm}}\right) \\
& =  \sum_{r\leq x}\sum_{nm = r} g (n)f (m)F\left({\frac {x}{r}}\right) \\
& = ((g*f)*F)(x),
\end{align*}
and therefore
$$g*(f*F) = (g*f)*F.$$
It follows that the abelian group $\CC^{\RR_{\geq 1}}$   under addition is a left module over the integral domain $\operatorname{Arith}_\CC$  under the scalar multiplication $(f,F) \longmapsto f *F$.   

As  another consequence of the associative property proved above, one has
$$S_{f*g} = (f* g)* 1 = f*(g*1) = f*S_g,$$ and therefore
$$S_{f*g}(x) = \sum_{n \leq x} f(n) S_g\left(\frac{x}{n}\right).$$

\begin{example}\label{piiex}
One has $\pi(x) = S_{\chi_{\pp}}(x)$, where $\chi_{\pp}$ is the characteristic function of the set $\pp$ of all prime numbers. Since also $\omega = \chi_{\pp} * \zeta$ and $\chi_{\pp} = \omega* \mu$,  it follows that
$$S_\omega(x) = \sum_{p \leq x} \left\lfloor \frac{x}{p}\right\rfloor,$$
$$\pi(x) =\sum_{n \leq x} \omega(n) M\left(\frac{x}{n} \right),$$
and
$$\pi(x) =\sum_{n \leq x} \mu(n) S_\omega\left(\frac{x}{n} \right) = \sum_{n \leq x} \sum_{p\leq \frac{x}{n}} \mu(n)\left\lfloor \frac{x}{pn}\right\rfloor = \sum_{n \leq x}\left( \sum_{p|n} \mu\left(\frac{n}{p}\right) \right)\left\lfloor \frac{x}{n}\right\rfloor.$$
The last equation implies that $\pi = f * \lfloor-\rfloor$, where $f(n) = \sum_{p|n} \mu\left(\frac{n}{p}\right)$ for all $n$.
\end{example}

As yet another consequence of the associative property proved above, we obtain the following result, which we call the {\bf second Dirichlet inversion theorem}.\index{Dirichlet inversion theorem}

\begin{proposition}[Second Dirichlet inversion theorem]\label{mobius2}
Let $h$ be any arithmetic function possessing a Dirichlet inverse $h^{*-1}$, that is, satisfying $h(1) \neq 0$, and let $F$ and $G$ be functions from $\RR_{\geq  1}$ to $\CC$.  Then one has
$$G(x)=\sum_{n\leq x}h(n)F\left({\frac {x}{n}}\right), \quad \forall x\geq 1,$$
if and only if
$$F(x)=\sum_{n\leq x}h^{*-1}(n)G\left({\frac {x}{n}}\right), \quad \forall x\geq 1.$$
\end{proposition}

\begin{example}
By Example \ref{piiex} and Proposition \ref{mobius2}, one has
$$\sum_{p \leq x} \left\lfloor \frac{x}{p}\right\rfloor = S_\omega (x) =\sum_{n \leq x} \pi\left(\frac{x}{n} \right)$$
for all $x \geq 0$.
\end{example}

\begin{corollary}
Let $h$ be any completely multiplicative arithmetic function, and let $F$ and $G$ be functions from $\RR_{\geq  1}$ to $\CC$.  Then one has
$$G(x)=\sum_{n\leq x}h (n)F\left({\frac {x}{n}}\right), \quad  \forall x\geq 1,$$
if and only if
$$F(x)=\sum_{n\leq x}\mu(n)h(n)G\left({\frac {x}{n}}\right), \quad  \forall x\geq 1.$$
\end{corollary}

Applying the corollary to the function $h = \zeta$ that is identically $1$  on $\ZZ_{> 0}$, we obtain the following.

\begin{corollary}
Let $F$ and $G$ be functions from $\RR_{\geq  1}$ to $\CC$.  Then one has
$$G(x)=\sum_{n\leq x}F\left({\frac {x}{n}}\right), \quad  \forall x\geq 1,$$
if and only if
$$F(x)=\sum_{n\leq x}\mu(n) G\left({\frac {x}{n}}\right), \quad  \forall x\geq 1.$$
\end{corollary}

Applying Proposition \ref{mobius2} to the function $F$ that is identically $1$  on $\RR_{\geq 1}$, we obtain the following.

\begin{corollary}
Let $f$ be any arithmetic function with $f(1) \neq 0$.  
Then one has
$${1=\sum_{n\leq x}f ^{*-1}(n)S_f\left({\frac {x}{n}}\right)}$$
for all $x \geq 1$.  Consequently, if $f$ is completely multiplicative, then
$${1=\sum_{n\leq x}\mu(n)f(n)S_f\left({\frac {x}{n}}\right)}$$
for all $x \geq 1$.
\end{corollary}

\begin{example}
One has
$$ \sum_{n \leq x} \mu(n)\left \lfloor \frac{x}{n}\right \rfloor = 1 =  \sum_{n \leq x} M\left(\frac{x}{n} \right)$$
and
$$1 =  \sum_{n \leq x}   \frac{\mu(n)}{n} H_{\lfloor x/n\rfloor}$$
for all $x \geq 1$.
\end{example}

Let $f$ and $g$ be arithmetic functions.  Applying  Proposition \ref{mobius2} to the functions $F(x) = f(\lfloor x \rfloor)$ and $G(x) = g(\lfloor x \rfloor)$, we obtain the following.

\begin{corollary}
Let $f$, $g$, and $h$ be arithmetic functions with $h(1) \neq 0$.  
Then one has
$$g(n) = \sum_{m = 1}^n h(m)f\left(\left\lfloor \frac{n}{m}\right\rfloor\right), \quad   \forall n\geq 1,$$
if and only if 
$$f(n) = \sum_{m = 1}^n h^{*-1}(m) g\left(\left\lfloor \frac{n}{m}\right\rfloor\right), \quad \forall n\geq 1.$$
\end{corollary}

\begin{remark}[Generalization of second Dirichlet inversion theorem]\label{secondgen}
The second Dirichlet inversion theorem has the following generalization to complex functions.  Let $\Delta$ be any subset of  $\{s \in \CC: |s| \geq 1\}$ such that $\frac{s}{n} \in \Delta$ for all $s \in \Delta$ and all positive integers $n \leq |s|$.
Let $h$ be any arithmetic function satisfying $h(1) \neq 0$, and let $F$ and $G$ be functions from $\Delta$ to $\CC$.  Then one has
$$G(s)=\sum_{n\leq |s|}h(n)F\left({\frac {s}{n}}\right), \quad \forall s \in \Delta,$$
if and only if
$$F(s)=\sum_{n\leq |s|}h^{*-1}(n)G\left({\frac {s}{n}}\right), \quad \forall s\in \Delta.$$
The proof is the same, in that it follows from associativity of the convolution
$$(f*F)(s) = \sum_{n\leq |s|}f (n)F\left({\frac {s}{n}}\right), \quad  \forall s \in \Delta.$$
\end{remark}

Finally, we note the following.

\begin{proposition}\label{thirdD}
Let $f$ and $h$ be arithmetic functions with $h(1) \neq 0$,  and suppose that
$$\sum_{n = 1}^\infty \sum_{m = 1}^\infty |h^{*-1}(n) h(m) f(mn)| < \infty.$$
Let $g$ be the (well-defined) arithmetic function defined by
$$g(n) = \sum_{m = 1}^\infty h(m) f(mn), \quad \forall n \in \ZZ_{> 0}.$$
Then one has
$$f(1) = \sum_{n = 1}^\infty h^{*-1}(n) g(n).$$   More generally, if  $n$ is any positive integer such that
$$\sum_{m = 1}^\infty \sum_{l = 1}^\infty |h^{*-1}(m) h(l) f(lmn)| < \infty,$$
then one has
\begin{align*}
f(n ) = \sum_{m = 1}^\infty h^{*-1}(m) g(mn).
\end{align*}
\end{proposition}

\begin{proof}
It suffices to prove the last statement of the proposition.   To do so, simply note that
\begin{align*}
\sum_{m = 1}^\infty h^{*-1}(m) g(mn) & = \sum_{m = 1}^\infty \sum_{l = 1}^\infty h^{*-1}(m)h(l) f(lmn) \\
& = \sum_{r = 1}^\infty \left(\sum_{lm = r} h^{*-1}(m)h(l) \right) f(rn) \\
& = f(n),
\end{align*}
where the rearrangement of the series follows from the assumption of the absolute convergence of the double sum.
\end{proof}

\section{Asymptotics of summatory functions and prime counting functions}

{\bf Abel's summation formula},\index{Abel's summation formula} stated below, is an extremely useful tool for studying summatory functions.

\begin{theorem}[Abel's summation formula, e.g., {\cite[Proposition 1.4]{kon} \cite[Theorem 1.4]{borg}}]\label{abels}
Let $f$ be an arithmetic function, let $y \leq x$ be nonnegative real numbers, and let $\alpha$ be a continuously differentiable complex-valued function on $[y,x]$.  Then one has
\begin{align*}
S_{\alpha \cdot f} (x)-S_{\alpha \cdot f} (y) = \sum_{y< n \leq x}  \alpha(n) f(n) = \alpha(x)S_f(x)-\alpha(y)S_f(y) - \int_y^x  \alpha'(t)S_f(t) \, dt.
\end{align*}  In particular, if $y \in [0,1)$, then one has
\begin{align*}
S_{ \alpha \cdot f}(x) = \alpha(x) S_f(x)- \int_1^x \alpha'(t)S_f(t)  \, dt.
\end{align*}
\end{theorem}

Abel's summation formula can be proved using the fundamental theorem of calculus \cite[Proposition 1.4]{kon}  or by applying integration by parts for {\it Riemann--Stieltjes integration} to the functions $S_f$ and $\alpha$ \cite[Section 1.1.3]{borg}.  It  can be generalized to the case where $\alpha$ is only assumed to be continuous on $[y,x]$, provided that  the integral is interpreted as a Riemann–Stieltjes integral: in this more general situation, one has
$${\displaystyle S_{\alpha \cdot f} (x)-S_{\alpha \cdot f} (y) = \alpha (x)S_f(x)-\alpha (y)S_f(y)-\int_{y}^{x}S_f(t)\,d\alpha (t).}$$

\begin{remark}[Summation by parts]\label{sumbyparts}
{\bf Summation by parts}, or {\bf partial summation},\index{summation by parts}\index{partial summation}  is a discrete analogue of Abel's summation formula.  It states that, if $f$ and $g$ are arithmetic functions, then
$$S_{ g \cdot f} (n) = g(n+1) S_f(n)  - \sum_{k = 1}^n (g(k+1)-g(k))S_f(k)$$
for all positive integers $n$.
\end{remark}

Applying Abel's summation formula and integration by parts to the arithmetic function $f =\zeta$, we obtain the following.

\begin{corollary}\label{abelcor}
Let $y \leq x$ be nonnegative real numbers, and let $\alpha$ be a continuously differentiable complex-valued function on $[y,x]$. Then one has
\begin{align*}
S_{ \alpha} (x)-S_{\alpha} (y)  = \sum_{y< n \leq x} \alpha(n) &= \lfloor x \rfloor \alpha(x) -\lfloor y \rfloor \alpha(y) - \int_y^x \lfloor t \rfloor  \alpha'(t)  \, dt \\ 
&=\int_y^x \alpha(t)\, dt - \{x\} \alpha(x) +\{y\} \alpha(y) + \int_y^x \{ t\}  \alpha'(t)  \, dt.
\end{align*}
  In particular, if $y \in [0,1)$, then one has
\begin{align*}
S_{\alpha}(x) = \lfloor x \rfloor \alpha(x) - \int_1^x \lfloor t \rfloor \alpha'(t) \, dt = \int_1^x \alpha(t)\, dt - \{x\} \alpha(x) +\alpha(1)+ \int_1^x \{ t\}  \alpha'(t)  \, dt.
\end{align*}
\end{corollary}

A simple application of Corollary \ref{abelcor} yields the weak form of Stirling's approximation provided in the following example.

\begin{example}
For the function $\alpha(x) = \log x$,  Corollary \ref{abelcor} yields
\begin{align*}
\log (\lfloor x \rfloor!) = \sum_{n \leq x} \log n &  = \int_1^x \log t \, dt  -\{x\} \log x+ \int_1^x\frac{ \{t\}}{t}\, dt \\
  & = x \log x -x +1-\{x\} \log x+ \int_1^x\frac{ \{t\}}{t}\, dt  \\
 & = x\log x -x + O(\log x) \ (x\to \infty)
\end{align*}
for all $x \geq 1$.
\end{example}

Another application of  Corollary \ref{abelcor}  is to the Euler--Mascheroni constant $\gamma$ and the  harmonic numbers $H_n$.

\begin{example}\label{harmex}
For the function $\alpha(x) = \frac{1}{x}$,  Corollary \ref{abelcor}  yields
\begin{align*}
 H_{\lfloor x \rfloor} & =  \frac{\lfloor x \rfloor}{x}  + \int_1^x \frac{\lfloor t \rfloor }{t^2} \, dt \\
 & = \log x -\frac{\{x\}}{x}+ 1- \int_1^x \frac{\{t\}}{t^2} \, dt
\end{align*}
for all $x \geq 1$.  It follows that
$$H_{\lfloor x \rfloor} - \log x  = -\frac{\{x\}}{x}+ 1- \int_1^x \frac{\{t\}}{t^2} \, dt \to 1- \int_1^\infty \frac{\{t\}}{t^2} \, dt$$
as $x \to \infty$, since the given integral converges by the comparison test.  This proves that the Euler--Mascheroni constant $\gamma$ exists and is given by $$\gamma =  1- \int_1^\infty \frac{\{t\}}{t^2} \, dt.$$  
It follows, then, that
$$-\frac{1}{x} \leq -\frac{\{x\}}{x}\leq  H_{\lfloor x \rfloor} - \log x-\gamma  = -\frac{\{x\}}{x}+\int_x^\infty \frac{\{t\}}{t^2} \, dt \leq  \int_x^\infty \frac{\{t\}}{t^2} \, dt \leq \frac{1}{x},$$
and therefore
$$ H_{\lfloor x \rfloor} - \log x-\gamma  =  O\left(\frac{1}{x}\right) \ (x \to \infty)$$
for all $x \geq 1$.  Consequently,  one has
$$H_{\lfloor x \rfloor} \sim \log x \ (x \to \infty).$$  
\end{example}

An application analogous to Example \ref{harmex} of  Abel's summation formula yields the following result relating Mertens' second theorem to the prime counting function.   Recall that  Mertens' second theorem states that the limit
$$M = \lim_{x \to \infty} \left( \sum_{p \leq x} \frac{1}{p} - \log \log x \right) = 0.261497212847\ldots,$$  
known as the Meissel--Mertens constant, exists.  Let
$$\Mert(x) =\sum_{p \leq x} \frac{1}{p}  -  \log \log x-M$$
for all $x> 1$.  

\begin{proposition}\label{psthm}
For all $x > 1$ one has
\begin{align*}
 \sum_{p \leq x} \frac{1}{p} & = \frac{\pi(x)}{x} + \int_1^x \frac{ \pi(t)}{t^2}\, dt,
\end{align*}
\begin{align*}
\pi(x)  = x \sum_{p \leq x}\frac{1}{p} - \int_1^x \left(\sum_{p \leq t}\frac{1}{p} \right)  dt,
\end{align*}
and
\begin{align*}
 \Mert(x)+ \frac{\li(x)-\pi(x)}{x} &=  \int_x^\infty \frac{ \li(t)-\pi(t)}{t^2}\, dt  \\ &= H- \int_1^x \frac{ \li(t)-\pi(t)}{t^2}\, dt  \\
 & =  \frac{H}{x} + \frac{1}{x} \int_1^x M_2(t) \, dt, 
\end{align*}
where
$$H=   \int_1^\infty  \frac{\li(t)-\pi(t)}{t^{2}} \, dt = \gamma -M = 0.315718452053\ldots.$$  
\end{proposition}

\begin{proof}
Abel's summation formula and integration by parts yield
\begin{align*}
 \sum_{p \leq x} \frac{1}{p} & = \frac{\pi(x)}{x} + \int_1^x \frac{ \pi(t)}{t^2}\, dt  \\ 
& = \log \log x+\gamma+ \frac{\pi(x)-\li(x)}{x} + \int_1^x \frac{ \pi(t)-\li(t)}{t^2}\, dt
\end{align*}
and therefore
\begin{align*}
 \Mert(x)& = - \frac{\li(x)-\pi(x)}{x}+ \gamma-M- \int_1^x \frac{ \li(t)-\pi(t)}{t^2}\, dt 
\\
 & = - \frac{\li(x)-\pi(x)}{x} + \int_x^\infty \frac{ \li(t)-\pi(t)}{t^2}\, dt,
\end{align*}
where  the integral
$$H :=   \int_1^\infty  \frac{\li(t)-\pi(t)}{t^{2}} \, dt = \gamma-M$$
exists by Mertens' second theorem and the fact that 
$$\lim_{x \to \infty} \frac{\li(x)-\pi(x)}{x}  = 0$$ 
which  in turn follows from  $\pi(x) = o(x)\ (x \to \infty)$ and $\li(x) = o(x)\ (x \to \infty)$.

Now, let $\pi_{-1}(x) =  \sum_{p \leq x} \frac{1}{p}$.  Abel's summation formula and integration by parts yield
\begin{align*}
\pi(x) &  = x \pi_{-1}(x)- \int_1^x \pi_{-1}(t)\, dt \\ &  =  x( \pi_{-1}(x)-\log\log x) -\gamma+ \li(x)- \int_1^x\left( \pi_{-1}(t)-\log \log t\right)  dt 
\end{align*}
and therefore
\begin{align*}
 \frac{\li(x)-\pi(x)}{x} &  = -(\pi_{-1}(x)- \log \log x) + \frac{1}{x} \left(\gamma+\int_1^x\left( \pi_{-1}(t)-\log \log t\right)  dt \right)  \\
& = -M_2(t)-M + \frac{\gamma}{x}+  \frac{1}{x} \int_1^x M_2(t) \, dt + \frac{x-1}{x}M \\
& = -M_2(t) + \frac{H}{x}+  \frac{1}{x} \int_1^x M_2(t) \, dt.
\end{align*}
The proposition follows.  
\end{proof}

Proposition \ref{psthm} implies that Mertens' second theorem is essentially equivalent to the convergence of the integral  $\int_1^\infty  \frac{\li(t)-\pi(t)}{t^{2}} \, dt$.   Mertens actually proved the stronger result that
$$\left|\sum_{p \leq x} \frac{1}{p}  -  \log \log x- M\right| < \frac{4}{\log x}, \quad\forall x \geq 2,$$
which one can show follows from Abel's summation formula and Mertens' first theorem (\ref{mert1}).
It follows that
\begin{align*}
\sum_{p \leq x} \frac{1}{p}  =  \log \log x+ M + O\left( \frac{1}{\log x} \right) \ (x \to \infty),
\end{align*}
which, given Proposition \ref{psthm} and Chebyshev's 1850 result below,  is equivalent to
$$\int_x^\infty \frac{ \li(t)-\pi(t)}{t^2}\, dt = O\left( \frac{1}{\log x} \right) \ (x \to \infty).$$

\begin{theorem}[{\cite{cheb1}}]\label{chebasymp}
One has $$\pi(x) \asymp \frac{x}{\log x} \ (x \to \infty).$$
\end{theorem}

For the remainder of this section, we assume Chebyshev's theorem above.  We will not assume the full strength of  the prime number theorem until Chapter 5.

Let $\vartheta(x)$ denote the {\bf first Chebyshev function},\index{first Chebyshev function $\vartheta(x)$} \index[symbols]{.s C@$\vartheta(x)$}   which is the function defined by $$\vartheta(x) = \sum_{p \leq x} \log p$$ for all $x \geq 0$.  The function $\vartheta(x)$ is the summatory function $S_{\chi_{\pp} \cdot \log}(x)$ of the arithmetic function $\chi_{\pp} \cdot \log$ and is a weighted prime counting function that weights each prime $p$ by $\log p$.  
The following proposition, which, like Proposition \ref{psthm}, follows from Abel's summation formula,  reveals  important and well-known connections  between the functions $\vartheta$ and $\pi$. 

\begin{proposition}\label{pitheta}
Let $x> 0$.  One has
$$\vartheta(x) = \pi(x) \log x -\int_0^x \frac{\pi(t)}{t} \, dt $$
and
\begin{align*}
\li(x)-\pi(x) -\frac{x-\vartheta(x)}{\log x}= \frac{1}{\log x}\int_0^x \frac{\li(t)-\pi(t)}{t} \, dt,
\end{align*}
where the Cauchy principal value of the latter integral is assumed.  One also has
$$\pi(x) = \frac{\vartheta(x)}{\log x} + \int_\mu^x \frac{\vartheta(t)}{t \log^2 t} \, dt$$
and
\begin{align*}
\li(x)-\pi(x) -\frac{x-\vartheta(x)}{\log x}= \li(N)-\pi(N) -\frac{N-\vartheta(N)}{\log N}+   \int_N^x \frac{t-\vartheta(t)}{t \log^2 t} \, dt
\end{align*}
for all $N > 1$.
\end{proposition}

\begin{corollary}\label{pitheta2}
One has
$$\vartheta(x) = \pi(x) \log x + O \left( \frac{x}{\log x} \right) \ (x \to \infty)$$
and therefore
$$\vartheta(x) \sim \pi(x) \log x \ (x \to \infty).$$ 
Consequently, the prime number theorem is equivalent to
$$\vartheta(x) \sim x \ (x \to \infty).$$ 
\end{corollary}

\begin{proof}
By  the proposition and Theorem \ref{chebasymp}, one has
\begin{align*}
 \pi(x) \log x -\vartheta(x)  = \int_0^x \frac{\pi(t)}{t} \, dt   = O\left(  \int_0^x \frac{dt}{\log t}  \right)  =   O\left( \frac{x}{\log x}\right) = O(\pi(x)) \ (x \to \infty).
\end{align*}
The corollary follows.
\end{proof}

Let  $\Lambda(n)$ be the {\bf von Mangoldt function},\index{von Mangoldt function $\Lambda(n)$}\index[symbols]{.rt  E@$\Lambda(n)$}   which by definition  is given by
$$\Lambda(n) =  \begin{cases} \log p & \quad \text{if } n > 1 \text{ is a power of some prime } p \\
 0 &  \quad \text{otherwise}.
\end{cases}$$
Equivalently, one defines $\Lambda = \widecheck{\log} = \mu*\log$.
The summatory function of the arithmetic function $\Lambda(n)$ is the {\bf second Chebyshev function}\index{second Chebyshev function $\psi(x)$} \index[symbols]{.s D@$\psi(x)$}  $\psi(x)$, given by $$\psi(x) =S_{\Lambda}(x) =  \sum_{k = 1}^\infty \sum_{p^k \leq x} \log p =  \sum_{k = 1}^\infty \vartheta(x^{1/k})$$ for all $x \geq 0$.  Thus, $\psi(x)$ is a weighted prime power counting function that weights each  prime power $p^n$ by $\log p$.    Since $\vartheta(x) = \psi(x) = 0$ if $x < 2$ and $x^{1/n} < 2$ if $n > \log_2 x$, one has
$$\psi(x) = \sum_{1 \leq n \leq  \log_2 x} \vartheta(x^{1/n}), \quad \forall x > 0.$$
By M\"obius inversion, then, one has
$$\vartheta(x) =  \sum_{1 \leq n \leq  \log_2 x} \mu(n) \psi(x^{1/n}) = \sum_{n = 1}^\infty \mu(n) \psi(x^{1/n}), \quad \forall x > 0.$$
The following is an analogue of Corollary \ref{pitheta2} for $\psi(x)$.

\begin{proposition}\label{pitheta2b}
Let $x \geq 2$ and $L(x) = \lfloor \log_2 x \rfloor -1$.  
\begin{enumerate}
\item One has
$$\vartheta(x^{1/2}) \leq  \psi(x)-\vartheta(x) \leq \vartheta(x^{1/2}) + L(x) \vartheta(x^{1/3})$$
and
$$\psi(x^{1/2}) \leq  \psi(x)-\vartheta(x) \leq \psi(x^{1/2}) + L(x) \psi(x^{1/3}).$$
 \item One has $$ \psi(x)-\vartheta(x) = O( \sqrt{x} ) \ (x \to \infty)$$
and therefore
$$\psi(x) \sim \vartheta(x) \ (x \to \infty).$$ 
\item The prime number theorem is equivalent to
$$\psi(x) \sim x \ (x \to \infty).$$ 
\end{enumerate}
\end{proposition}

 Let $\Pi(x)$ denote the {\bf Riemann prime counting function},\index{Riemann prime counting function $\Pi(x)$}  \index[symbols]{.s G@$\Pi(x)$}    defined by
$$\Pi(x) = \sum_{n = 1}^\infty \sum_{p^n \leq x} \frac{1}{n} = \sum_{n = 1}^\infty \frac{1}{n}\pi(x^{1/n}).$$
The function $\Pi(x)$ is a weighted prime power counting function that weights each prime power $p^n$ by $\frac{1}{n}$.   Since $\pi(x) = \Pi(x) = 0$ if $x < 2$, and $\sqrt[n]{x} < 2$ if $n > \log_2 x$, one has
\begin{align}\label{Pi2}
\Pi(x) = \sum_{n \leq \log_2 x} \frac{1}{n}\pi(\sqrt[n]{x}), \quad \forall x > 0.
\end{align}
By M\"obius inversion, one has
$$\pi(x) = \sum_{n \leq \log_2 x} \frac{ \mu(n)}{n} \Pi(\sqrt[n]{x}) = \sum_{n=1}^\infty \frac{ \mu(n)}{n} \Pi(\sqrt[n]{x}), \quad \forall x > 0.$$

\begin{proposition}\label{pitheta2c}
One has
 $$ \Pi(x)-\pi(x) = O\left (\frac{\sqrt{x}}{\log x} \right) \ (x \to \infty)$$
and therefore
$$\Pi(x) \sim \pi(x) \ (x \to \infty).$$
Consequently, the prime number theorem is equivalent to 
$$\Pi(x) \sim \frac{x}{\log x} \ (x \to \infty).$$
\end{proposition}  

\begin{lemma}\label{RAlem}
One has 
$$\pi(x) = o(\pi(x^t) ) \ (x \to \infty)$$
for all $t > 1$. 
\end{lemma}

\begin{proof}
This follows immediately from Theorem \ref{chebasymp}.
\end{proof}

The following result provides  various asymptotic relations between the weighted prime and prime power counting functions discussed in this section.

\begin{proposition}\label{RAprop2}
One has the following asymptotic expansions.
\begin{enumerate}
\item  $\displaystyle \Pi(x) \simeq \sum_{n = 1}^\infty \frac{1}{n}\pi(\sqrt[n]{x})  \ (x \to \infty)$.
\item $\displaystyle \pi(x) \simeq \sum_{n = 1}^\infty \frac{\mu(n)}{n}\Pi(\sqrt[n]{x})  \ (x \to \infty)$.
\item  $\displaystyle \psi(x) \simeq \sum_{n = 1}^\infty \vartheta(x^{1/n})  \ (x \to \infty)$.
\item $\displaystyle \vartheta(x) \simeq  \sum_{n = 1}^\infty \mu(n) \psi(x^{1/n})\ (x \to \infty)$.
\end{enumerate}
\end{proposition}

\begin{proof}  By (\ref{Pi2}), for any positive integer $N$, one has $$\frac{1}{N}\pi(x^{1/N}) \leq  \Pi(x)-\sum_{k = 1}^{N-1} \frac{1}{n}\pi(x^{1/n}) \leq \frac{1}{N}\pi(x^{1/N}) + \frac{1}{N+1}(\log_2 x) \pi(x^{1/(N+1)})$$
for all $x > 2^N$, and therefore
$$1 \leq \frac{ \Pi(x)-\sum_{n = 1}^{N-1} \frac{1}{n}\pi(x^{1/n})}{\frac{1}{N}\pi(x^{1/N})} \leq 1  + \frac{\frac{1}{N+1}(\log_2 x )\pi(x^{1/(N+1)})}{\frac{1}{N}\pi(x^{1/N})} \to 1$$
as $x \to \infty$, by Lemma \ref{RAlem}.  It follows that
$$\lim_{x \to \infty} \frac{ \Pi(x)-\sum_{n = 1}^{N-1} \frac{1}{n}\pi(x^{1/n})}{\frac{1}{N}\pi(x^{1/N})} = 1.$$
This proves statement (1), and the other three statements are proved in like fashion.
\end{proof}

The prime number theorem in a stronger form  yields, via Abel's summation formula, the following.

\begin{proposition}\label{mert1th}
Assume the prime number theorem in the stronger form 
$$\vartheta(x) = x + O \left(\frac{x}{(\log x)^t} \right) \ (x  \to \infty),$$
for some $t > 1$.  One has
$$\sum_{p \leq x} \frac{\log p}{p} -\log x = \frac{\vartheta(x)}{x}-\int_1^x \frac{t-\vartheta(t)}{t^2}\, dt = -B + o(1) \ (x \to \infty)$$
for some constant $B$, where
$$B = \lim_{x\to \infty}\left(\log x -\sum_{p \leq x} \frac{\log p}{p}\right)  = -1+\int_1^\infty \frac{t-\vartheta(t)}{t^2}\, dt = \int_1^\infty \frac{t-1-\vartheta(t)}{t^2}\, dt.$$
\end{proposition}

Note that $B= 1.332582275733\ldots$.

\begin{remark}[Equivalents of the prime number theorem]
\cite[pp.\ 12--13]{ing2} provides a simple proof that the three functions $\frac{\pi(x)}{\frac{x}{\log x}}$, $\frac{\psi(x)}{x}$, and $\frac{\vartheta(x)}{x}$ all have the same lim sup as $x \to \infty$, all at most $4 \log 2$, and the same lim inf as $x \to \infty$, all at least $\log 2$.  \cite[Theorem 6]{ing2} proves more, namely, that the lim sup is at least $1$ and the lim inf is at most $1$.  It follows that, if the limit of any of the three functions exists as $x \to \infty$, then all three limits  must exist and be equal to $1$.  These results are substantially easier to prove than the prime number theorem.
\end{remark}

\section{Average value and average order}

An arithmetic function $f$ is said to have {\bf average value}, or {\bf mean value}, $c \in \CC$\index{average value}\index{mean value} if
$$\lim_{x \to \infty} \frac{1}{x} \sum_{n \leq x} f(n) = c,$$
or, equivalently, if
$$\lim_{N \to \infty} \frac{1}{N} \sum_{n =1}^N f(n) = c.$$

\begin{example}\label{sumex0} \
\begin{enumerate}
\item  An arithmetic function $f$  has mean value $0$ if and only if $S_f(x) = o(x) \ (x \to \infty)$.  For example, since $M(x) = o(x) \ (x \to \infty)$, where $M = S_\mu$ is the Mertens function, the M\"obius function $\mu$ has average value $0$.  
\item Consider the number $$S_{|\mu|}(x) = \sum_{n \leq x} |\mu(n)|$$ of squarefree positive integers less than or equal to $x$.   It is well known \cite[Theorem 2.2]{mont} that 
\begin{align}\label{Qbound}
Q(x) = O(x^{1/2}) \ (x \to \infty),
\end{align}
 where $$Q(x) = S_{|\mu|}(x) -  \frac{x}{\zeta(2)},$$ 
and where $\zeta(2) = \frac{\pi^2}{6}$.  It follows that
$$\lim_{x \to \infty} \frac{1}{x}\sum_{n \leq x} |\mu(n)| = \frac{6}{\pi^2},$$
that is, the function $|\mu(n)|$ has average value $\frac{6}{\pi^2}$. Thus,  the asymptotic density of the set of all squarefree numbers in the positive integers is equal to  $\frac{6}{\pi^2}$.  
\item The multiplicative arithmetic function
$$\rho(n) = \frac{\phi(n)}{n} = \prod_{p |n} \left(1-\frac{1}{p} \right),$$
where $\rho(n)$ for any positive integer $n$ is equal to the probability that a randomly chosen integer from $1$ to $n$ is relatively prime to $n$, satisfies
\begin{align}\label{rhobound}
S_\rho(x) = \frac{6}{\pi^2}x + O(\log x) \ (x \to \infty)
\end{align}
 \cite[Theorem 2.1]{mont}.  Therefore, like the function $|\mu|$, the function $\rho = \frac{\phi}{\id}$ has average value $\frac{6}{\pi^2}$. 
\end{enumerate}
\end{example}

Let $f$ and $g$ be arithmetic functions.  The function $f$ is said to have {\bf average order $g$}\index{average order} if
$$S_f(x) \sim S_g(x) \ (x \to \infty),$$ or, equivalently, if
$$\frac{1}{x} \sum_{n \leq x} f(n) \sim \frac{1}{x} \sum_{n \leq x} g(n) \ (x \to \infty),$$ or, equivalently still, if
$$\frac{1}{N} \sum_{n = 1}^N f(n) \sim \frac{1}{N} \sum_{n = 1}^N g(n) \ (N \to \infty).$$  For example, $f$ has average value $c \neq 0$ if and only if $f$ has average order $c \neq 0$.   Note that the relation ``$f$ has average order $g$'' is an equivalence relation on the set of all arithmetic functions.   In the case where $f$ is real-valued, one often seeks an approximating function $g: \RR_{\geq 1} \to \RR$ that is continuous and monotone. 

\begin{example}[{\cite{tenen} \cite{finch}}]\label{sumex} \
\begin{enumerate}
\item Using what we now call the  {\it Dirichlet hyperbola method}, described later in this section, Dirichlet proved that
\begin{align}\label{Sdbound}
S_d(x) = x \log x + (2 \gamma-1)x + O(\sqrt{x}) \ (x \to \infty),
\end{align}
where $d$ denotes the divisor function.
Since
$$\sum_{n\leq x}\log n=x\log x-x+O(\log x) \sim \sum_{n \leq x} d(n) \ (x \to \infty),$$
it follows that the arithmetic function $d(n)$ has average order $\log n$.
\item The arithmetic function $2^\omega$ at $n$ equals the number $\sum_{d|n} |\mu(d)|$ of squarefree divisors of $n$ for all $n$.  Mertens proved that 
$$S_{2^\omega} (x) = \frac{1}{\zeta(2)}x \log x + \left(\frac{2\gamma-1}{\zeta(2)}- \frac{2\zeta'(2)}{\zeta^2(2)}\right) x+O(x^{1/2}\log x) \ (x \to \infty).$$
It follows that the arithmetic function $2^\omega$ has average order $\frac{1}{\zeta(2)} \log n$.
\item 
The summatory function of the additive arithmetic function $\omega$ satisfies
\begin{align}\label{Somegabound}
S_\omega(x) = x \log \log x + M x + O\left(\frac{x}{\log x} \right) \ (x \to \infty),
\end{align}
where $M$ is the Meissel--Mertens constant, 
and, since $\log \log$ is slowly varying, by Corollary \ref{karcor} one has $S_{\log \log}(x) \sim x \log \log x \ (x \to \infty)$, and therefore $\omega$ has average order $\log \log n$.    More generally, one has the asymptotic expansions
$${\frac {1}{x}}S_\omega(x) \simeq \log \log x+M+\sum_{n = 0}^\infty \left(\sum_{k=0}^{n}{\frac {\gamma_{k}}{k!}}-1\right){\frac {n!}{(\log x)^{n+1}}} \ (x \to \infty)$$
and
$${\frac {1}{x}}S_\Omega(x) \simeq \log \log x+K+\sum_{n = 0}^\infty \left(\sum_{k=0}^{n}{\frac {\gamma_{k}}{k!}}-1\right){\frac {n!}{(\log x)^{n+1}}} \ (x \to \infty),$$
where the $\gamma_k$ are the Stieltjes constants, and where
$$K = M + \sum_p \frac{1}{p(p-1)}  = \gamma+\sum_{n = 2}^\infty  \left(1-\frac{1}{n} \right)P(n) =1.034653881897\ldots$$
\cite[Section 1.4.3]{finch}.   It follows that the arithmetic function $\Omega$ has average order $\log \log n$ and the arithmetic function $\Omega-\omega \geq 0$ has average value $$K-M = \sum_p \frac{1}{p(p-1)} = 0.773156669049\ldots$$
\end{enumerate}
\end{example}

Now we present a few known techniques for verifying some of  the examples above.  

\begin{proposition}[{\cite[pp.\ 35--36]{mont} \cite[Theorem 6.13]{kon}}]\label{montprop}
Let $f$ be an arithmetic function.    For all $ x \geq 0$, one has
\begin{align*}
S_{\widehat{f}}(x) = \sum_{n \leq x} f(n)\left\lfloor \frac{x}{n} \right \rfloor =  x\sum_{n \leq x} \frac{f(n)}{n} -  \sum_{n \leq x} f(n)\left\{ \frac{x}{n} \right \} = x\sum_{n \leq x} \frac{f(n)}{n} + O\left( S_{|f|}(x) \right)  \ (x \to \infty).
\end{align*}
Consequently, one has  $S_{\widehat{f}} = f * \lfloor - \rfloor$, and, if 
$$S_{|f|}(x) = o\left( x\sum_{n \leq x} \frac{f(n)}{n}\right) \ (x \to \infty),$$
then
\begin{align*}
S_{\widehat{f}}(x) \sim x\sum_{n \leq x} \frac{f(n)}{n} \ (x \to \infty).
\end{align*}
Moreover, if the sum $s = \sum_{n = 1}^\infty \frac{f(n)}{n}$ converges and $S_{|f|}(x) = o(x) \ (x \to \infty)$, which both hold if the sum  $s = \sum_{n = 1}^\infty \frac{f(n)}{n}$ converges absolutely, then  the arithmetic function
${\widehat{f}}$ has average value $s$.
\end{proposition}

\begin{proof}
Since $\lfloor x \rfloor = x + O(1) \ (x \to \infty)$, one has
\begin{align*}
S_{\widehat{f}}(x)  & = \sum_{k \leq x} \sum_{n|k} f(n)  \\
& = \sum_{n \leq x} f(n) \sum_{k\leq x, \, n|k} 1 \\
& = \sum_{n \leq x} f(n)\left\lfloor \frac{x}{n} \right \rfloor \\
 & =  x\sum_{n \leq x} \frac{f(n)}{n} -  \sum_{n \leq x} f(n)\left\{ \frac{x}{n} \right \}  \\
& = x\sum_{n \leq x} \frac{f(n)}{n} + O\left( \sum_{n \leq x} |f(n)| \right)  \ (x \to \infty).
\end{align*}
Moreover,  if the sum $s = \sum_{n = 1}^\infty \frac{f(n)}{n}$ converges absolutely, then one has
\begin{align*}
\sum_{n \leq x} |f(n)| & = \sum_{n \leq \sqrt{x}}n\frac{|f(n)|}{n} + \sum_{\sqrt{x}< n \leq x }n\frac{|f(n)|}{n} \\
& \leq \sqrt{x}\sum_{n \leq \sqrt{x}}\frac{|f(n)|}{n}  + x  \sum_{\sqrt{x}< n \leq x }\frac{|f(n)|}{n} \\
& = O(\sqrt{x}) + o(x) \\
& = o(x) \ (x \to \infty).
\end{align*}
This completes the proof.
\end{proof}

Proposition \ref{montprop} is easily seen to yield (\ref{Qbound}), (\ref{rhobound}), and even (\ref{Somegabound}).    (See, for example, \cite[Theorems 2.1 and 2.2]{mont}).)  Applied to the arithmetic function $f = \zeta$, the proposition yields
\begin{align*}
S_{d}(x) & = \sum_{n \leq x} \left\lfloor \frac{x}{n} \right \rfloor = xH_{\lfloor x \rfloor} + O\left( {\lfloor x \rfloor} \right)   = x\log x+O(x) \ (x \to \infty) \ (x \to \infty),
\end{align*}
which implies  (\ref{Sdbound})  for $t = 1$  but  is too crude an estimate to yield (\ref{Sdbound}) for any $t < 1$.  For this purpose, one can use the following result, known as the {\bf Dirichlet hyperbola method},\index{Dirichlet hyperbola method} which has a simple verification \cite[p.\ 37]{mont}.

\begin{proposition}[{Dirichlet hyperbola method  \cite[p.\ 37]{mont}}]\label{hyperbola}
Let $f$ and $g$ be arithmetic functions.  For all $x, y \in \RR$ with $0 < y \leq x$, one has
$$S_{f*g}(x) = \sum_{n \leq y} f(n) S_g\left( \frac{x}{n} \right) + \sum_{n \leq \frac{x}{y}} g(n) S_f\left( \frac{x}{n} \right)-S_f(y)S_g\left( \frac{x}{y} \right).$$
\end{proposition}

Applying the Dirichlet hyperbola method to $f = g = \zeta$, we obtain
$$S_d(x) = \sum_{n \leq y}\left\lfloor  \frac{x}{n} \right \rfloor + \sum_{n \leq \frac{x}{y}}  \left\lfloor \frac{x}{n} \right\rfloor-\lfloor y \rfloor \left\lfloor \frac{x}{y} \right \rfloor.$$
The first term is
$$\sum_{n \leq y}\left\lfloor  \frac{x}{n} \right \rfloor  = x H_{\lfloor y \rfloor} + O(y)  = x \log y + \gamma x + O\left(\frac{x}{y}+y \right),$$
and the error bound is minimized by taking  $y =  \sqrt{x}$.  Since $\lfloor \sqrt{x} \rfloor^2 = x+ O( \sqrt{x} ) \ (x \to \infty)$, for this choice of $y$ one has
\begin{align*}
S_d(x) & = 2\sum_{n \leq \sqrt{x}}\left\lfloor  \frac{x}{n} \right \rfloor -\lfloor \sqrt{x} \rfloor^2 \\
& =  2(x \log \sqrt{x} + \gamma x ) + O( \sqrt{x} ) -x + O( \sqrt{x} ) \\
& = x \log x  + (2 \gamma - 1) x + O( \sqrt{x} ) \ (x \to \infty),
\end{align*}
thus yielding Dirichlet's result.

The asymptotic relation (\ref{Somegabound}) for $S_\omega$ follows  from Proposition \ref{montprop}, applied to $f = \chi_{\pp}$ and $g = \zeta$:
\begin{align*}
S_\omega(x)  & = \sum_{p \leq x} \left \lfloor \frac{x}{p}\right \rfloor  \\
& = x\sum_{p \leq x} \frac{1}{p} -\sum_{p \leq x} \left\{\frac{x}{p}\right \} \\
&  = x\sum_{p \leq x} \frac{1}{p} + O(\pi(x)) \\
& = x \left (\log \log x + M + O\left( \frac{1}{\log x}\right) \right) + O(\pi(x)) \\
& = x \log \log x + M x +   O\left( \frac{x}{\log x}\right) \ (x \to \infty),
\end{align*}
where the penultimate equation follows from Mertens' second theorem and the last equation follows from Chebyshev's result that $\pi(x)  = O \left (\frac{x}{\log x} \right) \ (x \to \infty)$ (and thus the proof does not require the full strength of the prime number theorem).  It is interesting to observe that $S_\omega(x) = \sum_{p \leq x}  \lfloor \frac{x}{p}  \rfloor$, while
$S_d(x) =  \sum_{n \leq x}\lfloor \frac{x}{n} \rfloor$, and thus the $O$ bound (\ref{Somegabound}) for $S_\omega$ is a ``primes'' analogue of Dirichlet's $O$ bound (\ref{Sdbound}) for $S_d$.  The result for $S_\omega$ generalizes as follows.

\begin{proposition}[{\cite[Theorem 6.19]{kon}}]\label{erdd}
Let $f$ be a real-valued additive arithmetic function such that $f(p) = 1$ for all primes $p$ and $f(p^{n}) -f(p^{n-1})$ is uniformly bounded for all primes $p$ and all  positive integers $n$.  Then one has
$$S_f(x) = x \log \log x +(M+ C_f)x + O\left (\frac{x}{\log x} \right) \ (x \to \infty),$$
where $C_f$ is the constant
$$C_f = \sum_{n = 2}^\infty \sum_p \frac{f(p^n)-f(p^{n-1})}{p^n}.$$
\end{proposition}

For example,  the proposition applies to both $f = \omega$ and $f= \Omega$, where $C_\omega = 0$ and $C_\Omega$ is given by
$$C_\Omega =   \sum_{n = 2}^\infty \sum_p \frac{1}{p^n} =     \sum_p \sum_{n = 2}^\infty \frac{1}{p^n} =  \sum_p \frac{1}{p(p-1)},$$
thus yielding (\ref{Somegabound}) and
\begin{align*}
S_\Omega(x) = x \log \log x + (M +C_\Omega)x + O\left(\frac{x}{\log x} \right) \ (x \to \infty).
\end{align*}

\begin{remark}[Probabilistic number theory]\label{prnt}
{\it Probabilistic number theory} was founded  in the 1930s by Erd\H{o}s, Wintner, and Kac.  One of its key ideas is that the prime numbers can be treated as independent random variables.  The  {\bf Erd\H{o}s--Kac theorem},\index{Erd\H{o}s--Kac theorem} also known as the {\bf fundamental theorem of probabilistic number theory},\index{fundamental theorem of probabilistic number theory}   states that, for any strongly additive real-valued arithmetic function $f$ with $|f(p)| \leq 1$ for all primes $p$, and for all real numbers $a \leq b$, one has
$$\lim_{{x\rightarrow \infty }}\left({\frac  {1}{x}}\cdot \#\left\{n\leq x:a\leq {\frac  {f(n)-A(n)}{B(n)}}\leq b\right\}\right)=\Phi (a,b),$$
where 
$$A(n)=\sum_{{p\leq n}}{\frac  {f(p)}{p}},\qquad B(n)={\sqrt  {\sum_{{p\leq n}}{\frac  {f(p)^{2}}{p}}}},$$
and where $\Phi(a,b)$ is the normal distribution
$$\Phi(a,b)= \frac{1}{\sqrt{2\pi}}\int_a^b e^{-t^2/2} \, dt$$
with mean $0$ and standard deviation $1$ \cite{erdos}.    From this result it follows that
$$\lim_{x \rightarrow \infty}  \left ( \frac {1}{x} \cdot \#\left\{ n \leq x : a \le \frac{\omega(n) - \log \log n}{\sqrt{\log \log n}} \le b \right\} \right ) = \Phi(a,b)$$
 for all real numbers $a \leq b$.
\end{remark}

\section{Dirichlet series}

Along with forming its summatory function, a natural way to ``analysis-ize'' an arithmetic function $f$ is to form the {\bf Dirichlet series $D_f(s)$ of $f$},\index{Dirichlet series $D_f(s)$}\index[symbols]{.l E@$D_f(s)$}    which is the complex function defined by
$$D_f(s) = \sum_{n = 1}^\infty \frac{f(n)}{n^s}$$
for all $s \in \CC$ (which may or may not coverge for a given $s \in \CC$).  The Dirichlet series of $f$ is obtained from the formal Dirichlet series $D_f(X)$ of $f$ by ``substituting'' a complex variable $s$ for the formal ``variable'' $X$, i.e., by substituting $n^{-s}$ for the monomial $n^{-X} \in \operatorname{Dir}_\CC$.  

\begin{example}\label{direx} \
\begin{enumerate}
\item The Dirichlet series of the arithmetic function $\zeta(n) = 1$ converges precisely on the set $\{s \in \CC: \operatorname{Re} s > 1\}$, and on that domain one has
$$D_{\zeta} (s) = \zeta(s)= \sum_{n = 1}^\infty \frac{1}{n^s}.$$   
\item The Dirichlet series of the arithmetic function $\frac{1}{\id}$ converges precisely on the set $\{s \in \CC: \operatorname{Re} s > 0\}$, and on that domain one has
$$D_{\frac{1}{\id}} (s) = \zeta(s+1)= \sum_{n = 1}^\infty \frac{1}{n^{s+1}}.$$    More generally, for any $a \in \CC$, 
the Dirichlet series of the arithmetic function $\id^a$ converges precisely on the set $\{s \in \CC: \operatorname{Re} s > 1+\operatorname{Re} a\}$, and on that domain one has
$$D_{\id^a} (s) = \zeta(s-a)= \sum_{n = 1}^\infty \frac{1}{n^{s-a}}.$$   
\item The Dirichlet series  $$D_{\chi_{\pp}} (s) = P(s)= \sum_{p}^\infty \frac{1}{p^s}$$  of the characteristic function $\chi_{\pp}$ of the set $\pp$ of all prime numbers is the {\bf prime zeta function} $P(s)$, and\index{prime zeta function $P(s)$}\index[symbols]{.t  B@$P(s)$}   it converges precisely on the set $\{s \in \CC: \operatorname{Re} s > 1\}$.
\end{enumerate}
\end{example}

As  was first proved by Riemann in \cite{rie} and is now well known, the Dirichlet series $\zeta(s) = D_{\zeta} (s)$ on $\{s \in \CC: \operatorname{Re} s> 1\}$ extends to a meromorphic function on all of $\CC$ possessing a single pole, namely, a simple pole at $s = 1$ with residue $1 = \lim_{s \to 1} (s-1) \zeta(s)$.   The meromorphic continuation of $D_{\zeta} (s)$ is also denoted $\zeta(s)$ and is called the {\bf Riemann zeta function}.\index{Riemann zeta function $\zeta(s)$}
 It is well known \cite[Theorem 1.3]{ivic} that the Stieltjes constants $\gamma_k$, as defined in Example \ref{stiec}(2), are  the real constants occurring in the Laurent expansion
\begin{align}\label{stcon}
 \zeta (s)={\frac {1}{s-1}}+\sum_{k=0}^{\infty }{\frac {(-1)^{k}\gamma_k }{k!}}(s-1)^{k}
 \end{align}
of $\zeta(s)$ at $1$.   In particular, one has
$$\gamma = \lim_{s \to 1} \left( \zeta(s)-\frac{1}{s-1}\right).$$  See Section 4.2 for a survey of results on the Riemann zeta function.   Many treatises are devoted to the topic, e.g.,  \cite{edw} \cite{ivic} \cite{patt} \cite{tit}.

\begin{proposition}[{\cite[pp.\ 9--10]{flaj}}]\label{flajprop}
Let $F(x)$ be a complex-valued function that is  bounded and Riemann integrable on $[1,x]$ for all $x > 1$, with $\deg F = \sigma < \infty$, and let
$$G(s) = \int_1^\infty \frac{ F(x)}{x^{s+1}} \,dx.$$
Then $G(s)$ is analytic in the right half plane $\{s \in \CC: \operatorname{Re} s > \sigma\}$, and one has
$$\frac{d}{ds} G(s) = -\int_1^\infty \frac{ F(x)\log x}{x^{s+1}} \,dx$$
for all $s$ in the given right half plane.
\end{proposition}

The equation for $\frac{d}{ds} G(s)$ in the proposition states that  $\frac{d}{ds} G(s)$ can be obtained from $G(s) = \int_1^\infty \frac{ F(x)}{x^{s+1}} \,dx$ by differentiating the integrand with respect to $s$.  Note that \cite[pp.\ 9--10]{flaj} generalizes the proposition to    complex-valued functions $F(s)$ that are Lebesgue integrable on all compact subsets of $[1,\infty)$.

Let $f$ be an arithmetic function,  let $s \in \CC$,  and let $N$ be a positive integer.   Applying Abel's summation formula to $\alpha(x) = x^{-s}$, we see that
\begin{align}\label{abelrem}
\sum_{n = 1}^N \frac{f(n)}{n^s}  = \frac{S_f(N)}{N^{s}}+s\int_1^N \frac{S_f(x)}{ x^{s+1}}\, dx.
\end{align}  
This forms the basis for proofs of Theorem \ref{dirichlet} below.  The theorem,  which is of fundamental importance to analytic number theory, shows that the summatory function $S_f(x)$ and the Dirichlet series $D_f(s)$ of  an arithmetic function $f$ are  intimately related to one another.  We have replaced the limits superior in the traditional formulations of the result with degrees, as appropriate.    For any arithmetic function $f$,  the extended real number $\sigma_{\mathrm{c}}(f)$ provided by the theorem---namely,  the unique extended real number for which the Dirichlet series $D_f(s)$ converges for $\operatorname{Re}s > \sigma_{\mathrm{c}}(f)$ and diverges for $\operatorname{Re} s < \sigma_{\mathrm{c}}(f)$---is called
the {\bf abscissa of convergence of $D_f(s)$}.\index{abscissa of convergence $\sigma_{\mathrm{c}}(f)$}\index[symbols]{.l F@$\sigma_{\mathrm{c}}(f)$}   Likewise, the extended real number $\sigma_{\mathrm{a}}(f)$ provided by the theorem is called the {\bf abscissa of absolute convergence of $D_f(s)$}.\index[symbols]{.l Fa@$\sigma_{\mathrm{a}}(f)$}\index{abscissa of absolute convergence $\sigma_{\mathrm{a}}(f)$}

\begin{theorem}[{\cite[Theorem 1.3]{mont} \cite[Chapter II Sections 3 and 4]{widd}}]\label{dirichlet}
Let $f: \ZZ_{> 0} \longrightarrow \CC$ be an arithmetic function. 
\begin{enumerate}
\item There exists a unique $\sigma_{\mathrm{c}}  = \sigma_{\mathrm{c}}(f) \in \overline{\RR}$ such that  the series $\sum_{n =1}^\infty f(n)n^{-s}$ converges for $\operatorname{Re}s > \sigma_{\mathrm{c}}$ and diverges for $\operatorname{Re} s < \sigma_{\mathrm{c}}$.
\item If $D_f(0) = \sum_{n = 1}^\infty f(n)$ is divergent, then 
$$\sigma_{\mathrm{c}} = \deg S_f \geq 0.$$
\item If $D_f(0) = \sum_{n = 1}^\infty f(n)$ is convergent, then 
$$\sigma_{\mathrm{c}} = \deg (D_f(0)-S_f) \leq 0,$$
where $D_f(0) -S_f(x) = \sum_{n > x} f(n)$ for all $x \geq 0$.
\item One has
\begin{align*}
{D_f(s)} = s\int_1^\infty  \frac{S_f(x)}{x^{s+1}} \,dx
\end{align*} 
for all $s \in \CC$ with $\operatorname{Re} s > \max(\sigma_{\mathrm{c}},0)$.
\item If $S_f(x)$ is bounded, then $\sigma_{\mathrm{c}} \leq 0$.
\item If $\sigma_{\mathrm{c}} < 0$, then $S_f(x)$ is bounded.
 \item There exists a unique $\sigma_{\mathrm{a}} = \sigma_{\mathrm{a}}(f)  \in \overline{\RR}$ such that  the series $\sum_{n =1}^\infty f(n)n^{-s}$ converges absolutely if $\operatorname{Re}s > \sigma_{\mathrm{a}}$ but does not converge absolutely if $\operatorname{Re} s < \sigma_{\mathrm{a}}$.
\item One has $$\sigma_{\mathrm{c}} \leq \sigma_{\mathrm{a}} \leq \sigma_{\mathrm{c}} + 1.$$
\item If $\sigma_{\mathrm{c}} = \sigma_{\mathrm{a}}$, then both are equal to $1+\deg f$.
\item If $f(n) \in \RR_{\geq 0}$ for all sufficiently large $n$, then $\sigma_{\mathrm{c}} = \sigma_{\mathrm{a}} = 1+\deg f$.
\item $D_f(s)$ is analytic on $\{s \in \CC: \operatorname{Re} s > \sigma_{\mathrm{c}}\}$, the Dirichlet series $\sum_{n = 1}^\infty \frac{f(n)\log n}{n^s}$ has the same abscissa of convergence and abscissa of absolute convergence as $D_f(s)$, and one has
$$\frac{d}{ds} D_f(s) = -\sum_{n = 1}^\infty \frac{f(n)\log n}{n^s}$$
for all $s$ for which the given series on the right converges.  
\end{enumerate}
\end{theorem}

 See \cite[Chapter II]{widd} for a generalization of the theorem to {\it Laplace--Stieltjes transforms}, which generalize both  Dirichlet series and Laplace transforms.

The following useful result relates Dirichlet series with the operation of Dirichlet convolution.

\begin{proposition}[{\cite[Theorem 11.5]{apos} \cite[Theorem 6.9]{kon}}]\label{Dconjprop}
Let $f$ and $g$ be arithmetic functions.  For all $s \in \CC$, if $D_f(s)$  and $D_g(s)$ are absolutely convergent, then $D_{f*g}(s)$ is absolutely converent and
$$D_{f*g}(s) = D_f(s) D_g(s).$$
In particular, one has 
$$\sigma_{\mathrm{a}}(f*g) \leq \max(\sigma_{\mathrm{a}}(f),\sigma_{\mathrm{a}}(g)).$$
\end{proposition}

\begin{example}
For all $s, a \in \CC$, one has
$$D_\phi(s) = \frac{\zeta(s-1)}{\zeta(s)}$$
if $\operatorname{Re} s >2$
and
$$D_{\sigma_a}(s) = \zeta(s) \zeta(s-a)$$
if  $\operatorname{Re} s > 1+\operatorname{Re} a$.
\end{example}

\begin{example}
Let $f$ be the arithmetic function with $f(1) =1$ and $f(2) = -1$, and $f(n) = 0$ for all $n \geq 3$.   Then 
$D_f(s) = 1-2^{-s}$ is absolutely convergent for all $s \in \CC$, so that $\sigma_{\mathrm{c}}(f) = \sigma_{\mathrm{a}}(f)  = -\infty$.  On the other hand,
$D_{f^{*-1}}(s) = \frac{1}{1-2^{-s}}$ converges, and converges absolutely, precisely for all $s \in \CC$ with $\operatorname{Re} s > 0$, so that $\sigma_{\mathrm{c}}(f^{*-1}) = \sigma_{\mathrm{a}}(f^{*-1}) = 0$.
\end{example}

 Below are applications of Theorem \ref{dirichlet} and Proposition \ref{Dconjprop} to  seven important arithmetic functions and their corresponding Dirichlet series.

\begin{example}\label{direxample}  \
\begin{enumerate}
\item Since $S_\zeta(x) = {\lfloor x \rfloor}$ and $D_\zeta(s) = \zeta(s)$, one has
$${\zeta(s)} = s \int_1^\infty  \frac{\lfloor x \rfloor}{x^{s+1}} \, dx$$
and
$$\frac{1}{s-1} =  \int_1^\infty  \frac{x}{x^{s+1}} \, dx$$
for all $s \in \CC$ with $\operatorname{Re} s > 1 = \deg \, \lfloor x \rfloor$.   It follows that
\begin{align}\label{zetaexp}
{\zeta(s)} -\frac{1}{s-1} = 1-s \int_1^\infty  \frac{\{x \}}{x^{s+1}} \, dx = s \int_1^\infty  \frac{1-\{x \}}{x^{s+1}} \, dx 
\end{align}
for all $s \in \CC$ with $\operatorname{Re} s > 1$.  Moreover, the function ${\zeta(s)} -\frac{1}{s-1}$ is entire since it has limit $\gamma$ as $s \to 1$, and the two integrals above converge for all $s \in \CC$ with $\operatorname{Re} s > 0$.  By uniqueness of analytic continuation,  (\ref{zetaexp}) is valid for all such $s$.  In particular, one has
$$\gamma = 1- \int_1^\infty  \frac{\{x \}}{x^{2}} \, dx = \int_1^\infty  \frac{1-\{x \}}{x^{2}} \, dx.$$
\item Since $S_{\frac{1}{\id}}(x) = H_{\lfloor x \rfloor}$ and $D_{\frac{1}{\id}}(s) = \zeta(s+1)$, one has
$${\zeta(s+1)} = s \int_1^\infty  \frac{H_{\lfloor x \rfloor}}{x^{s+1}} \, dx,$$
for all $s \in \CC$ with $\operatorname{Re} s > 0 = \deg \, H_{\lfloor x \rfloor}$.   It follows that
$${\zeta(s)}  = \gamma+\frac{1}{s-1}+(s-1)\int_1^\infty \frac{H_{\lfloor x \rfloor}-\gamma - \log x}{x^{s}} \, dx$$
for all $s \in \CC$ with $\operatorname{Re} s > 1$.  Moreover, since
$$H_{\lfloor x \rfloor}-\gamma - \log x  = O\left( \frac{1}{x}\right) \ (x \to \infty),$$
the integral on the right converges to an analytic function on  $\{s \in \CC: \operatorname{Re} s > 0\}$, and thus the expression above for  $\zeta(s)$  is valid for all such $s$.    Taking the limit as $s \to 0^+$, we find that
 $$\lim_{s \to 0^+} \int_1^\infty \frac{H_{\lfloor x \rfloor}-\gamma - \log x}{x^s} \, dx = \gamma-1-\zeta(0) = \gamma-\frac{1}{2} = 0.077215664901\ldots.$$
Moreover, one has
$$\lim_{s \to 1} \frac{\zeta(s)-\gamma-\frac{1}{s-1}}{s-1} = -\gamma_1,$$
where $\gamma = \gamma_0$  and $\gamma_1  = -0.072815845483\ldots$ are the zeroth and first Stieltjes constants.  
It follows that
$$\int_1^\infty \frac{H_{\lfloor x \rfloor}-\gamma - \log x}{x} \, dx = -\gamma_1 = 0.072815845483\ldots.$$
Likewise, one has
$$\int_1^\infty \frac{H_{\lfloor x \rfloor}-\gamma - \log x}{x^2} \, dx = \zeta(2)-\gamma -1  = 0.067718401946\ldots.$$
\item Since $S_\mu(x) = M(x) = O(x) \ (x \to \infty)$ and $D_\mu(s) = \frac{1}{\zeta(s)}$, one has
$$\frac{1}{\zeta(s)} = s \int_1^\infty  \frac{M(x)}{x^{s+1}} \, dx,$$
for all $s \in \CC$ with $\operatorname{Re} s > 1$.  Taking a limit as $s \to 1^+$, we find that
$$\int_1^\infty  \frac{M(x)}{x^2} \, dx = 0.$$
It is known that  the special value
$$D_\mu(1) = \sum_{n = 1}^\infty \frac{\mu(n)}{n} = 0$$
is ``elementarily'' equivalent to the prime number theorem \cite[p.\ 285]{nark}.
\item Recall that $\chi_{\pp}$ denotes the characteristic function of the set $\pp$ of all primes.  The Dirichlet series $$D_{\chi_{\pp}}(s) = \sum_p \frac{1}{p^s} = P(s)$$  is the prime zeta function $P(s)$, and  the summatory function $S_{\chi_{\pp}}(x)$  is the prime counting function $\pi(x)$.  Therefore, one has
$${P(s)} = s \int_1^\infty  \frac{\pi(x)}{x^{s+1}} \, dx$$
for all $s \in \CC$ with $\operatorname{Re} s > 1 = \deg \pi$.   
It is also known that
$$\Log \frac{1}{s-1} = s \int_1^\infty  \frac{\li(x)}{x^{s+1}} \, dx,$$
and therefore
\begin{align}\label{zetaP}
\Log \frac{1}{s-1}- P(s)=  s \int_1^\infty  \frac{\li(x)-\pi(x)}{x^{s+1}} \, dx,
\end{align}
for all $s \in \CC$ with $\operatorname{Re} s > 1$.  
By Example \ref{primez} later in this section, the function $\Log \frac{1}{s-1}- P(s)$ can be analytically continued  to $\{s\in \CC:\operatorname{Re} s > \Theta\}\cup\{1\}$, with  $$\lim_{s \to 1} \left(\Log \frac{1}{s-1}- P(s)\right) = \sum_{n = 2}^\infty \frac{P(n)}{n} = \gamma-M = H = 0.315718452053\ldots.$$
 Moreover,  the function $\int_1^\infty \frac{\li(x)-\pi(x)}{x^{s+1}} \, dx$ is analytic,  and  thus (\ref{zetaP}) holds,  on $\{s \in \CC: \operatorname{Re} s > \Theta\}$.
\item The first Chebyshev function $\vartheta(x)$ is the summatory function $S_{\chi_{\pp} \cdot \log}(x)$ of the completely multiplicative function $\chi_{\pp} \cdot \log$, which has Dirichlet series
$$D_{\chi_{\pp} \cdot \log}(s) = \sum_p \frac{\log p}{p^s} = -P'(s).$$  Recall from Corollary \ref{pitheta2},  that
$$\vartheta(x) \sim \pi(x) \log x \ (x \to \infty),$$
and therefore the prime number theorem is equivalent to
$$\vartheta(x) \sim x \ (x \to \infty).$$  
By Theorem \ref{dirichlet}, one has
$$P'(s) = - s \int_1^\infty\frac{\vartheta(x)}{x^{s+1}} \, dx,$$
and therefore
$$P'(s)+ \frac{1}{s-1} =  P'(s)+\frac{s}{s-1} -1=  -1+s \int_1^\infty\frac{x-\vartheta(x)}{x^{s+1}} \, dx,$$
 for all $s \in \CC$ with $\operatorname{Re} s > 1 = \deg \vartheta$.   Moreover, by example (4), one has
$$P'(s)+ \frac{1}{s-1} =- \frac{d}{ds}\left( \Log \frac{1}{s-1}-P(s) \right)= -\frac{d}{ds} \left( s \int_1^\infty  \frac{\li(x)-\pi(x)}{x^{s+1}} \, dx\right),$$
which can be analytically continued  to $\{s\in \CC:\operatorname{Re} s > \Theta\} \cup \{1\}$.
Furthermore, by  Proposition \ref{mert1th} (or \cite[(4.21)]{ross}), one has
$$-1+\int_1^\infty \frac{x-\vartheta(x)}{x^2}\, dt = \lim_{x\to \infty}\left(\log x -\sum_{p \leq x} \frac{\log p}{p}\right) =  B.$$
It follows \cite[(5.2)]{ross} that
$$\lim_{s\to 1}\left(P'(s)+ \frac{1}{s-1} \right)  = B.$$
Consequently, one has
\begin{align*}
B & = \lim_{s\to 1^+}\left(P'(s)+ \frac{1}{s-1} \right) \\
&  = - \lim_{s\to 1^+} \frac{d}{ds} \left( s \int_1^\infty  \frac{\li(x)-\pi(x)}{x^{s+1}} \, dx\right) \\
& = - \lim_{s\to 1^+} \left(  \int_1^\infty  \frac{\li(x)-\pi(x)}{x^{s+1}} \, dx +  \frac{d}{ds}  \int_1^\infty  \frac{\li(x)-\pi(x)}{x^{s+1}} \, dx \right) \\
& = -H +\lim_{s\to 1^+} \int_1^\infty  \frac{(\li(x)-\pi(x))\log x}{x^{s+1}} \, dx  \\
& = -H + \int_1^\infty  \frac{(\li(x)-\pi(x))\log x}{x^{2}} \, dx, 
\end{align*}
whence 
$$B+H = \int_1^\infty  \frac{(\li(x)-\pi(x))\log x}{x^{2}} \, dx$$
and therefore
$$B=  \int_1^\infty  \frac{(\li(x)-\pi(x))(\log x-1)}{x^{2}} \, dx.$$
\item The second Chebyshev function $\psi(x)$ is the summatory function of  the von Mangoldt function
$$\Lambda(n) =  \begin{cases} \log p & \quad \text{if } n > 1 \text{ is a power of some prime } p \\
 0 &  \quad \text{otherwise},
\end{cases}$$
where $\Lambda = \widecheck{\log} = \mu*\log$, so that $\log = \widehat{\Lambda} = \zeta* \Lambda$ and 
$$D_\Lambda(s) = \frac{D_{\log}(s)}{\zeta(s)} = - \frac{\zeta'(s)}{\zeta(s)}.$$  By Theorem \ref{dirichlet}, then, one has
\begin{align}\label{ingz}
\frac{\zeta'(s)}{\zeta(s)} =- s \int_1^\infty\frac{\psi(x)}{x^{s+1}} \, dx,
\end{align}
and therefore
$$\frac{\zeta'(s)}{\zeta(s)}+ \frac{1}{s-1} =  \frac{\zeta'(s)}{\zeta(s)}+\frac{s}{s-1} -1=  -1+s \int_1^\infty\frac{x-\psi(x)}{x^{s+1}} \, dx,$$
 for all $s \in \CC$ with $\operatorname{Re} s > 1 = \deg \psi$.   Moreover,  the function $\frac{\zeta'(s)}{\zeta(s)}+ \frac{1}{s-1}$ can be meromorphically continued to $\CC$, and its meromorphic extension is analytic on $\{s\in \CC:\operatorname{Re} s > \Theta\}$ and at $s = 1$, with
$$\gamma = \lim_{s\to 1}\left(\frac{\zeta'(s)}{\zeta(s)}+ \frac{1}{s-1} \right)  = -1+ \int_1^\infty\frac{x-\psi(x)}{x^{2}} \, dx.$$
\item   Recall that $\Pi(x)$ denotes the Riemann prime counting function.   Since
$$\Pi(x) = \sum_{n \leq x} \frac{\Lambda(n)}{\log n}$$ for all $x \geq 0$, the function $\Pi(x)$ is the summatory function $S_{\frac{\Lambda}{\log}}$ of the completely multiplicative arithmetic function $\frac{\Lambda}{\log}$, which has Dirichlet series
$$D_{\frac{\Lambda}{\log}} (s) = \log \zeta(s),$$
where $\log \zeta(s)$ is the unique analytic continuation of $\log \zeta(x)$ on $\RR_{> 1}$ to $\{s \in \CC:\operatorname{Re} s > 1 \}$.   Therefore, by Theorem \ref{dirichlet}, one has
\begin{align}\label{ingz2}
\log \zeta(s) = s\int_1^\infty \frac{\Pi(x)}{x^{s+1}} \, dx,
\end{align}
whence
\begin{align}\label{Psz}
\log \zeta(s) -P(s) = s\int_1^\infty \frac{\Pi(x)-\pi(x)}{x^{s+1}} \, dx,
\end{align}
for all $s \in \CC$ with $\operatorname{Re} s > 1 = \deg \Pi$.   In fact, by Example \ref{primez} later in this section,  the function $\log \zeta(s)-P(s)$ can be analytically continued to $\{s \in \CC: \operatorname{Re} s > \frac{1}{2}\}$, and, since $\deg(\Pi(x)-\pi(x)) = \frac{1}{2}$,  it follows that (\ref{Psz}) holds for all $s$ with $\operatorname{Re} s > \frac{1}{2}$.  Thus,  since 
$$\lim_{s\to 1} \left(\log \zeta(s)- \Log \frac{1}{s-1}\right) = 0,$$
one has
$$\int_1^\infty \frac{\Pi(x)-\pi(x)}{x^{2}} \, dx = \lim_{s\to 1} \left(\log \zeta(s)-P(s)\right) = \lim_{s\to 1} \left(\Log \frac{1}{s-1}-P(s)\right) =  H.$$
It follows that
$$\int_1^\infty \frac{\Pi(x)-\li(x)}{x^{2}} \, dx =  \int_1^\infty \frac{\Pi(x)-\pi(x)}{x^{2}} \, dx -  \int_1^\infty  \frac{\li(x)-\pi(x)}{x^{2}} \, dx = H-H = 0.$$
\end{enumerate}
\end{example}

\begin{remark}[Mellin transforms, Laplace transforms, and Fourier transforms]
The {\bf Mellin transform}\index{Mellin transform} ${\mathcal M}_f(s)$ \cite{flaj} of a locally Lebesgue integrable complex-valued function $f$ on $(0,\infty)$ is the complex function $${\mathcal M}_f(s) = \int_0^\infty x^s f(x) \, \frac{dx}{x}.$$    The identity in statement (4) of Theorem \ref{dirichlet} can be expressed as
$$D_f(s) = s {\mathcal M}_{S_f}(-s).$$  The {\bf two-sided Laplace transform}\index{two-sided Laplace transform} $\mathcal{B}_g(s)$ of  a locally Lebesgue integrable complex-valued function $g$ on $\RR$ is the complex function $$\mathcal{B}_g(s) = \int_{-\infty}^\infty e^{-st}g(t)\, dt.$$   One has
$$\mathcal{B}_g(s) = \mathcal{M}_{g (\log x)}(-s)$$
and
$$\mathcal{M}_f(s) = \mathcal{B}_{f(e^{x})}(-s).$$
Thus, the identity in statement (4)  of Theorem \ref{dirichlet} can also be expressed as
$$D_f(s) = s {\mathcal B}_{S_f(e^{x})}(s).$$ 
One may also define the {\it Fourier transform} in terms of the Mellin transform and vice versa.  
\end{remark}

The following  theorem is an analytic version both of its analogue for formal Dirichlet series and of the fundamental theorem of arithmetic.

\begin{theorem}[{\cite[Theorem 11.7]{apos}}]
Let $f$ be a multiplicative arithmetic function.   One has
$$D_f(s) = \prod_p \sum_{k = 0}^\infty \frac{f(p^k)}{p^{ks}}$$
for all $s \in \CC$ with $\operatorname{Re} s  >\sigma_{\mathrm{a}}$.  Consequently, if $f$ is completely multipicative, then one has 
$$D_f(s) = \prod_p \frac{1}{1-\frac{f(p)}{p^s}}$$
for all such $s$.
\end{theorem}

Note that the series  $B_{f,p}(s)= \sum_{k = 0}^\infty \frac{f(p^k)}{p^{ks}}$ is the Bell series $B_{f,p}(X)$  of $f$ ``evaluated at'' $X = s$.
The product representation of $D_f(s)$ in the theorem is known as the  {\bf Euler product representation of  $D_f(s)$}. \index{Euler product representation}

\begin{figure}[ht!]
\includegraphics[width=70mm]{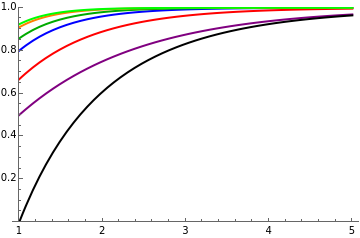} 
\caption{\centering Graph of $\frac{1}{\zeta(s)}$ (in black) and its Euler factors $1-{\frac {1}{p^{s}}}$ for all primes $p \leq 13$, on $[1,5]$}
\label{eulerprodu}
\end{figure}

\begin{example}
One has the Euler product representations
$$\zeta(s) = \prod_{p} \frac{1}{1-{\frac {1}{p^{s}}}}$$
and 
$$D_{\mu}(s) = \frac{1}{\zeta(s)} = \prod_{p} \left(1-{\frac {1}{p^{s}}} \right)$$
of $\zeta(s)$ and $\frac{1}{\zeta(s)}$,  respectively,  both
valid for all $s \in \CC$ with $\operatorname{Re} s > 1$.    See Figure \ref{eulerprodu} for a graph of  $\frac{1}{\zeta(s)}$ and its Euler factors $1-{\frac {1}{p^{s}}}$ for all primes $p \leq 13$, on $[1,5]$.
\end{example}

\begin{example}\label{primez}
Let $s \in \CC$ with $\operatorname{Re} s > 1$.  From the Euler product representation
$$\zeta(s) = \prod_{p}\left(1-{\frac {1}{p^{s}}}\right)^{-1}$$
of $\zeta(s)$,   it follows \cite{frob} that
\begin{align*}
\log \zeta (s) & =-\sum_{p }\Log \left(1-{\frac {1}{p^{s}}}\right) \\
  & =  \sum_p \sum_{n = 1}^\infty \frac{1}{np^{ns}} \\
  & =  \sum_{n = 1}^\infty \frac{1}{n} \sum_p  \frac{1}{p^{ns}} \\
  & =\sum_{n = 1}^\infty{\frac {P(ns)}{n}},
\end{align*}
where $\log \zeta(s)$ is the unique analytic continuation of $\log \zeta(x)$ on $\RR_{> 1}$ to $\{s \in \CC:\operatorname{Re} s > 1 \}$, and where $P(s)$ is the prime zeta function.
Thus, by Proposition \ref{thirdD} applied to $f(n) = P(ns)$ and $h(n) = \frac{1}{n}$,  one has
\begin{align*} P(s)=\sum_{n = 1}^\infty {\frac {\mu (n)}{n}}\log \zeta (ns).
\end{align*}
It follows that the functions $P(s)-\log \zeta(s)$ and $P(s)-\Log \frac{1}{s-1}$ can be analytically continued to $\{s \in \CC: \operatorname{Re} s> \frac{1}{2}\}$ and to $\{s \in \CC: \operatorname{Re} s> \Theta\}\cup \{1\}$,  respectively, with
$$\lim_{s \to 1} \left(\Log \frac{1}{s-1}-P(s)\right)=  \lim_{s \to 1} \left(\log \zeta(s)-P(s)\right) = \sum_{n = 2}^\infty \frac{P(n)}{n},$$
where
\begin{align*}
  \sum_{n = 2}^\infty \frac{P(n)}{n} & =  \sum_p \left(\frac{1}{2p^2} + \frac{1}{3p^3} + \frac{1}{4p^4} + \cdots \right) \\
  & =- \sum_p \left(\frac{1}{p} +\log\left(1- \frac{1}{p} \right) \right) \\
& = \gamma-M \\
& = H.
\end{align*}
\end{example}

The nonvanishing of $\zeta(s)$ on the line $\{s \in \CC: \operatorname{Re} s = 1\}$, in conjunction with the
 following theorem,  which is known as the {\bf Wiener--Ikehara theorem},\index{Wiener--Ikehara theorem} yields the well-known proof of the prime number theorem given in Example \ref{liw} below.

\begin{theorem}[{Wiener--Ikehara theorem \cite{ike} \cite{wie}}]\label{wiener}
Let $f$ be a nonnegative real-valued  arithmetic function.  Suppose that $\sigma_{\mathrm{c}} = \sigma_{\mathrm{c}}(f)$ is finite.  Then the Dirichlet series $D_f(s)$ has a singularity at $s = \sigma_{\mathrm{c}}$.  Moreover, if $\sigma_{\mathrm{c}} > 0$ and the function $D_f(s) - \frac{a}{s-\sigma_{\mathrm{c}}}$ for some $a >0$ extends to a continuous function on $\{s \in \CC: \operatorname{Re} s \geq \sigma_{\mathrm{c}}\}$, then one has
$$S_f(x) \sim \frac{a}{\sigma_{\mathrm{c}}}x^{\sigma_{\mathrm{c}}} \ (x \to \infty).$$
\end{theorem}

\begin{example}\label{liw}
Applying the theorem to $f = \Lambda$, since the Dirichlet series $D_\Lambda = -\frac{\zeta'}{\zeta}$ has abscissa of convergence $\sigma_{\mathrm{c}}(\Lambda) = 1$, if the function $\zeta(s)$ does not vanish on $\{s \in \CC: \operatorname{Re} s  = 1\}$, then the function $-\frac{\zeta'(s)}{\zeta(s)}-\frac{1}{s-1}$ extends to a continuous function on $\{s \in \CC: \operatorname{Re} s \geq 1\}$, and therefore
$$\psi(x) = S_\Lambda(s) \sim \tfrac{1}{1}x^1 = x \ (x \to \infty),$$
which implies the prime number theorem.
Thus, Theorem \ref{wiener} and the nonvanishing of $\zeta(s)$ on the line $\{s \in \CC: \operatorname{Re} s = 1\}$ implies the prime number theorem.  In Section 5.2, we discuss some more extensive  (well-known) zero-free regions for $\zeta(s)$ and their consequences for the error $\li(x)-\pi(x)$ in the approximation $\li(x)$ of $\pi(x)$.
\end{example}

The Wiener--Ikehara theorem is one of many {\it Tauberian theorems}  that can be used to supply elegant proofs of the prime number theorem (without error bounds).   See \cite{korevaar} for an excellent treatment of  various Tauberian theorems, e.g., those of Wiener--Ikehara,  Hardy--Littlewood,  Karamata, and Lambert, and the corresponding proofs of the prime number theorem.

Examples of Dirichlet series and summatory functions are provided in Table \ref{table1a} at the end of this section.  For each arithmetic function $f$ in the table, we provide the simplest  function $F(x)$ known to be asymptotic to the function
$$S_{f,0}(x) = \begin{cases} D_f(0)-S_f(x)  = \sum_{n > x} f(n) & \quad \text{if } D_f(0) = \sum_{n = 1}^\infty f(n) \text{ exists} \\
    S_f(x)  & \quad \text{otherwise},
\end{cases}$$
so that, by Theorem \ref{dirichlet}, one has $\sigma_{\mathrm{c}}(f) = \deg S_{f,0}  = \deg F$.

\begin{table}[!htbp]
  \caption{\centering Dirichlet series and summatory function of arithmetic functions $f(n)$}
    \footnotesize
\begin{tabular}{|l||l|l|l|l|} \hline
$f(n)$ &  $D_f(s)$ &  $S_{f,0}(x)$ & $S_{f,0}(x) \sim$ & $\sigma_{\mathrm{c}}(f)$  \\  \hline\hline
$1$ & $\zeta(s)$  &  $\lfloor x \rfloor$ & $x$   &  $1$  \\  \hline
$(-1)^{n-1}$ & $\left(1-2^{1-s}\right) \zeta(s)$  &  $\sum_{n\leq x} (-1)^n$  & $\sum_{n\leq x} (-1)^n$     &   $0$ \\  \hline
$\chi_{\pp}(n)$ & $P(s)$  &  $\pi(x)$ &  $\frac{x}{\log x}$     &  $1$   \\  \hline
$\chi_{\pp}(n)\log n$ & $-P'(s)$  &  $\vartheta(x)$ &  $x$    &  $1$  \\  \hline
$\Lambda(n)$ & $-\frac{\zeta'(s)}{\zeta(s)}$ &  $\psi(x)$ & $x$     &  $1$  \\  \hline
$\frac{\Lambda(n)}{\log n}$ & $\log \zeta(s)$ &  $\Pi(x)$ & $\frac{x}{\log x}$     &   $1$ \\  \hline
$n^a$, $a > -1$ &  $\zeta(s-a)$ & $\sum_{n \leq x}n^a$ & $\frac{x^{a+1}}{a+1}$    & $a+1$   \\  \hline
$\frac{1}{n}$ &  $\zeta(s+1)$ & $\sum_{n \leq x} \frac{1}{n}$ & $\log x$     &  $0$   \\  \hline
$n^a$, $a < -1$ &  $\zeta(s-a)$ & $\sum_{n > x}n^a$ &  $- \frac{x^{a+1}}{a+1} $    &  $a+1$   \\  \hline
$\chi_{\pp}(n)n^a$, & $P(s-a)$  &  $\sum_{p \leq x}p^a$ & $\frac{x^{a+1}}{(a+1)\log x}$   &  $a+1$  \\ 
   $\quad a > -1$&   &   &   &   \\  \hline
$\chi_{\pp}(n)\frac{1}{n}$ & $P(s+1)$  &  $\sum_{p \leq x}\frac{1}{p}$ &  $\log \log x$    &  $0$  \\  \hline
$\chi_{\pp}(n)n^a$, & $P(s-a)$  &  $\sum_{p > x}p^a$ & $  -\frac{x^{a+1}}{(a+1)\log x}$  & $a+1$  \\ 
   $\quad a < -1$ &   &   &    &   \\  \hline
$(\log n)^k$, $k \in \ZZ_{\geq 0}$ & $(-1)^k\zeta^{(k)}(s)$  &  $\sum_{n \leq x} (\log n)^k$ & $x (\log x)^k$     &     $1$ \\  \hline
$n^a(\log n)^k$, & $(-1)^k\zeta^{(k)}(s-a)$  &  $\sum_{n \leq x} n^a(\log n)^k$ & $\frac{ x^{a+1}(\log x)^k}{a+1}$     &   $a+1$ \\ 
   $\quad a > -1,  \ k \in \ZZ_{\geq 0}$ &   &    &  &    \\  \hline
$\frac{1}{n}(\log n)^k$, $k \in \ZZ_{\geq 0}$ & $(-1)^k\zeta^{(k)}(s+1)$  &  $\sum_{n \leq x} \frac{1}{n}(\log n)^k$ & $\frac{ (\log x)^{k+1}}{k+1}$    &  $0$  \\  \hline
$n^a(\log n)^k$, & $(-1)^k\zeta^{(k)}(s-a)$  &  $\sum_{n > x} n^a(\log n)^k$ & $ -\frac{ x^{a+1}(\log x)^k}{a+1}$  & $a+1$  \\ 
   $\quad a < -1,  \ k \in \ZZ_{\geq 0}$ &   &  \quad   &  &     \\  \hline
$\sigma_a(n)$, $a > 1$ & $\zeta(s)\zeta(s-a)$ &  $\sum_{n \leq x} \sigma_a(n)$ & $\frac{\zeta(a+1)x^{a+1}}{a+1}$    &   $a+1$  \\  \hline
$\sigma(n) = \sigma_1(n)$ & $\zeta(s)\zeta(s-1)$ &  $\sum_{n \leq x} \sigma(n)$ & $\frac{\zeta(2)x^{2}}{2}$    &     $2$  \\  \hline
$\sigma_a(n)$, $0<a<1$ & $\zeta(s)\zeta(s-a)$ &  $\sum_{n \leq x} \sigma_a(n)$ & $\frac{\zeta(a+1)x^{a+1}}{a+1}$     &   $a+1$   \\  \hline
$d(n) = \sigma_0(n)$ & $\zeta(s)^2$ &  $\sum_{n \leq x} d(n)$ & $x \log x$      & $1$   \\  \hline
$\sigma_{a}(n)$, $a< 0$  & $\zeta(s)\zeta(s-a)$ &  $\sum_{n \leq x} \sigma_{a}(n)$ & $\zeta(1-a)x$     &  $1$  \\  \hline
$\omega(n)$ & $\zeta(s)P(s)$ &  $\sum_{n \leq x} \omega(n)$ & $x \log \log x$     &  $1$  \\  \hline
$\Omega(n)$ & $\zeta(s)\sum_{n = 1}^\infty P(ns)$ &  $\sum_{n \leq x} \Omega(n)$ & $x \log \log x$     & $1$   \\  \hline
$\mu(n)$ & $\frac{1}{\zeta(s)}$ &  $M(x)$ & $M(x)$    &  $\Theta$  \\  \hline
$\lambda(n)$ & $\frac{\zeta(2s)}{\zeta(s)}$ &  $L(x)$ &  $L(x)$   & $\Theta$    \\  \hline
$|\mu(n)|$ & $\frac{\zeta(s)}{\zeta(2s)}$ &  $\sum_{n \leq x}|\mu(n)| $ & $\frac{6}{\pi^2} x$    &  $1$  \\  \hline
$\phi(n)$ & $\frac{\zeta(s-1)}{\zeta(s)}$ &  $\sum_{n \leq x} \phi(n)$ & $\frac{3}{\pi^2} x^2$    &   $2$  \\  \hline
$2^{\omega(n)}$ & $\frac{\zeta(s)^2}{\zeta(2s)}$ &  $\sum_{n \leq x} 2^{\omega(n)}$ & $\frac{6}{\pi^2} x\log x$     &  $1$  \\  \hline
\end{tabular}\label{table1a}
\end{table}

\chapter{Special functions  in analytic number theory}

{\it Special functions} are particular mathematical functions that have established names and notations due to their importance in mathematics, physics, and other fields.  In this chapter, we state some important and well-known facts about some important special functions in analytic number theory,  including the gamma function $\Gamma(s)$, the Riemann zeta function $\zeta(s)$,  the Riemann xi function $\xi(s)$ an Riemann Xi function $\Xi(s)$, the prime zeta function $P(s)$,  the complementary exponential integral function $\Ein(s)$, the exponential integral functions $\Ei(s)$ and $E_1(s)$,  and Riemann's function $\Ri(x)$.  We also introduce an entire function $\ERi(s)$ such that $\Ri(x) = \ERi(\log x)$ for all $x > 0$.

\section{The gamma function $\Gamma(s)$}

The {\bf gamma function}\index{gamma function $\Gamma(s)$}\index[symbols]{.szw  A@$\Gamma(s)$}  $\Gamma(s)$ for $s$ in the right open half plane $\{s \in \CC: \operatorname{Re} s > 0\}$ is defined by
$$ \Gamma (s)=\int_{0}^{\infty }x^{s}e^{-x}\,\frac{dx}{x}.$$  
The {\bf Mellin transform}\index{Mellin transform} of a function $\rho(x)$ on $[0,\infty)$ is defined as $\int_{0}^{\infty }x^{s}\rho(x)\,\frac{dx}{x}$, and thus the gamma function is the Mellin transform of the {\bf exponential distribution} $\rho(x) = e^{-x}$, on $[0,\infty)$.  The gamma function  is of fundamental importance in analytic number theory, real and complex analysis, probability theory, and combinatorics.  References on the gamma function, and for the facts noted in this section, abound: see \cite[Chapter 1]{andrews} \cite[Chapter 1]{lebed} \cite[Chapter 2 Sections 1 and 2]{olver}, for example.

\begin{figure}[ht!]
\includegraphics[width=80mm]{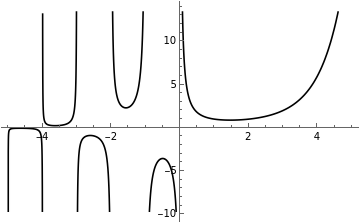}
\caption{\centering Graph of $\Gamma(x)$ on $[-5,5]$}
 \label{Gamma}
\end{figure}

\begin{figure}[ht!]
\includegraphics[width=80mm]{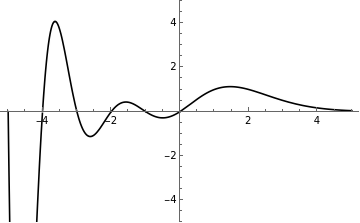}
\caption{\centering Graph of $\frac{1}{\Gamma(x)}$ on $[-5,5]$}
\label{GammaR}
\end{figure}

Any of the formulas
$$\Gamma(s) = \lim_{n \to \infty} \frac{n!\, n^s}{s(s+1)(s+2)\cdots (s+n)}, \quad \forall s \in \CC\backslash \ZZ_{\leq 0},$$
$$\Gamma(s) =\frac{1}{s}\prod_{n=1}^\infty \frac{\left(1+{\frac {1}{n}}\right)^{s}}{1+\frac{s}{n}}, \quad \forall s \in \CC\backslash \ZZ_{\leq 0},$$
and
$$\Gamma(s) =  \frac{e^{-\gamma s}}{s} \prod_{n  = 1}^\infty\left( 1+\frac{s}{n}\right)^{-1}e^{s/n}, \quad \forall s \in \CC\backslash \ZZ_{\leq 0},$$
can be used to extend the gamma function to a unique meromorphic function on $\CC$ with (simple) poles precisely at the nonpositive integers, with residues $$\operatorname{Res}(\Gamma, -n) = \lim_{s \to -n} (s+n)\Gamma(s) = \frac{(-1)^n}{n!}$$ for all nonnegative integers $n$.    
Like the complex exponential function $\exp s$, the gamma function $\Gamma(s)$ has no zeros on $\CC$.  It follows that the {\bf reciprocal gamma function}\index{reciprocal gamma function $\frac{1}{\Gamma(s)}$} $\frac{1}{\Gamma(s)}$ is entire, with (simple) zeros precisely at the nonnegative integers.   Moreover, the expression
$$\frac{1}{\Gamma(s)} =  se^{\gamma s} \prod_{n  = 1}^\infty\left( 1+\frac{s}{n}\right)e^{-s/n},$$
is valid for all $s \in \CC$ and is called the {\bf Weierstrass factorization  of $\frac{1}{\Gamma(s)}$}.  
Figure \ref{Gamma} provides a graph of $\Gamma(x)$, and Figure \ref{GammaR} a graph of $\frac{1}{\Gamma(x)}$, on the interval $[-5,5]$.

  The gamma function satisfies the {\bf functional equation}\index{functional equation for $\Gamma(s)$} $$\Gamma(s+1) = s\Gamma(s), \quad \forall s \in \CC\backslash \ZZ_{\leq 0}.$$
Since $\Gamma(1) = 1$, it follows by induction that $$\Gamma(n+1) = n!$$ for all nonnegative integers $n$, so that the gamma function interpolates the sequence of factorials.     The  {\bf Bohr--Mollerup theorem}\index{Bohr--Mollerup theorem} states that $\Gamma(s)$ is the unique meromorphic function $f(s)$ on $\CC$ satisfying the following properties.
\begin{enumerate}
\item $f(x) > 0$  for all $x \in\RR_{>0}$.
\item $\log f(x)$ is a convex function on $\RR_{>0}$.
\item $f(x+1) = xf(x)$ for all $x\in\RR_{>0}$.
\item $f(1) = 1$.
\end{enumerate}

Two other important functional equations satisfied by the gamma function are  {\bf Euler's reflection formula}\index{Euler's reflection formula}
$$s \Gamma(s)\Gamma(-s) = \Gamma (s)\Gamma (1-s)={\frac {\pi }{\sin \pi s}}, \quad \forall s \in \CC \backslash \ZZ,$$  
and {\bf Legendre's duplication formula}\index{Legendre's duplication formula}
$$\Gamma (s)\Gamma \left(s+{\tfrac {1}{2}}\right)=\Gamma\left(\tfrac{1}{2}\right)2^{1-2s}\Gamma (2s), \quad  \forall s \in \CC \backslash \tfrac{1}{2}\ZZ_{\leq 0},$$
where $\Gamma\left(\frac{1}{2}\right)= \sqrt{\pi}$, and where
$$\sin \pi s = \pi s \prod_{n \in \ZZ\backslash \{0\}} \left(1-\frac{s}{n} \right)e^{s/n} = \pi s \prod_{n = 1}^\infty \left(1-\frac{s^2}{n^2} \right), \quad  \forall s \in \CC,$$
is the {\bf Weierstrass factorization of $\sin \pi s$}.\index{Weierstrass factorization of $\sin \pi s$}

 {\bf Stirling's approximation for $\Gamma(s)$}\index{Stirling's approximation} is the asymptotic relation
$${\displaystyle \Gamma (s)={\sqrt {\frac {2\pi }{s}}}\,{\left({\frac {s}{e}}\right)}^{s}\left(1+O\left({\frac {1}{s}}\right)\right)} \ (s \to \infty)$$
on $\{s \in \CC: |\operatorname{Arg}(s)| < \pi-\varepsilon\}$ for any fixed $\varepsilon > 0$.  This generalizes the well-known asymptotic relation
$${n!={\sqrt {{2\pi n }}}\,{\left({\frac {n}{e}}\right)}^{n}\left(1+O\left({\frac {1}{n}}\right)\right)} \ (n \to \infty).$$

One has 
$$\Gamma'(1) = -\gamma$$
and
$$\Gamma''(1) = \gamma^2 + \zeta(2)$$
and therefore $$\lim_{s \to 0}\left(  \Gamma(s)-\frac{1}{s}\right) = \lim_{s \to 0}\left( \frac{ \Gamma(1+s)-\Gamma(1)}{s}\right) = \Gamma'(1) =  - \gamma.$$
Thus, the function $ \Gamma(s)-\frac{1}{s}$ is analytic at $0$ (and its derivative at $0$ is $\frac{1}{2} \Gamma''(1)$).    More generally,  for any nonnegative integer $n$, one has
$$\Gamma'(n+1) = n!(H_n -\gamma)$$
and 
$$\lim_{s \to -n}\left(  \Gamma(s)-\frac{(-1)^n}{n!(s+n)}\right) =\frac{(-1)^n}{n!}(H_n-\gamma),$$
where $H_n$ denotes the $n$th harmonic number.

The  {\bf log-gamma function}\index{log-gamma function}  $\log \Gamma(s)$ is the unique 
analytic continuation of the function $\log (\Gamma(x))$ on $\RR_{>0}$  to $\CC \backslash(-\infty,0]$.  It can be extended to $\CC \backslash \ZZ_{\leq 0}$ by setting
$$\log \Gamma(s) =  -\gamma s -\Log s+ \sum_{n  = 1}^\infty\left(\frac{s}{n} -\Log \left( 1+\frac{s}{n}\right)\right), \quad \forall s \in \CC \backslash \ZZ_{\leq 0},$$ and then it satisfies
$$\operatorname{Re} \log\Gamma(s) =  \operatorname{Re} \Log(\Gamma(s)), \quad \forall s \in \CC \backslash \ZZ_{\leq 0}.$$
The log-gamma function has the asymptotic expansion
$${\displaystyle \log \Gamma(s) \simeq (s-\tfrac{1}{2})\Log s-s+{\tfrac {1}{2}}\log 2\pi+\sum_{k=1}^{\infty}{\frac {B_{2k}}{2k(2k-1)s^{2k-1}}} \ (s \to \infty)} $$
on $\{s \in \CC:  |\operatorname{Arg}(s)| < \pi-\varepsilon\}$ for any fixed $\varepsilon > 0$.
One also has
$$\left.{\frac{d^n}{ds^n} \log \Gamma(s) }\right|_{s = 1} = (-1)^{n} n!\, \zeta(n), \quad \forall n > 1.$$

The {\bf digamma function}\index{digamma function $\Psi(x)$} is defined by $$\Psi(s) = \frac{\Gamma'(s)}{\Gamma(s)},  \quad \forall s \in \CC\backslash \ZZ_{\leq 0},\index[symbols]{.szz   A@$\Psi(s)$}$$
and, like the gamma function, is meromorphic on $\CC$ with (simple) poles precisely at the nonpositive integers, but with residues
$$\operatorname{Res}(\Psi, -n) = \lim_{s \to -n} (s+n)\Psi(s) =-1$$
for all nonnegative integers $n$.      One has
$$\Psi(s) = \frac{d}{ds} \log \Gamma(s), \quad \forall s \in \CC \backslash(-\infty,0],$$
and $$\log\Gamma(s) = \int_{1}^s \Psi(z) \, dz,  \quad \forall s \in \CC \backslash(-\infty,0],$$
where the integral is along any path from $1$ to $s$ that does not cross $(-\infty,0]$.

From the functional equation for the gamma function follows, by logarithmic differentiation, the {\bf functional equation}
$$\Psi(s+1) = \frac{1}{s}+\Psi(s), \quad \forall s \in \CC\backslash \ZZ_{\leq 0},$$
 for $\Psi(s)$. \index{functional equation for $\Psi(s)$}  Since $\Psi(1) = \Gamma'(1) = -\gamma$ and $H_0 = 0$, it follows from the functional equation for the digamma function that the harmonic numbers $H_n$ are interpolated by the complex function
\begin{align}\label{Hpsi}
H_s = \Psi(s+1)+ \gamma = \sum_{k = 1}^\infty \left(\frac{1}{k}-\frac{1}{s+k}\right) = \lim_{n \to \infty} \left( H_n - \sum_{k = 1}^n \frac{1}{s+k}\right),  \quad \forall s \in \CC\backslash \ZZ_{< 0},\index[symbols]{.szz  A@$H_s$}
\end{align}
which satisfies the   functional equation
$$H_s = \frac{1}{s}+ H_{s-1}  \quad \forall s \in \CC\backslash \ZZ_{< 0}.$$
 One has the Taylor series expansion
$${\displaystyle  H_s =\sum_{n=1}^{\infty }(-1)^{n+1}\zeta (n+1)s^{n}}$$
valid for all $s \in \CC$ with $|s| < 1$, and, in particular,
$$\left.{\frac{d^n}{ds^n} H_s }\right|_{s = 0} = (-1)^{n+1} n!\, \zeta(n+1), \quad \forall n \geq 1.$$
One also has
$$H_s -\gamma -\Log s \sim \frac{1}{2s} \ (s \to \infty)$$
on  $\{s \in \CC: |\operatorname{Arg}(s)| < \pi-\varepsilon\}$ for  any  $\varepsilon > 0$.
More generally,  one has the (divergent) asymptotic expansion
$$H_s -\gamma -\Log s - \frac{1}{2s} \simeq  \sum_{k = 1}^\infty \frac{-B_{2k}}{2ks^{2k}} \ (s \to \infty)$$
on $\{s \in \CC: |\operatorname{Arg}(s)| < \pi-\varepsilon\}$.

\section{The Riemann zeta function $\zeta(s)$}

Recall from Section 3.8 that the Dirichlet series $$\zeta(s) = \sum_{n = 1}^\infty \frac{1}{n^s}$$
converges, and converges absolutely, precisely for all complex numbers $s$ with $\operatorname{Re} s > 1$.   In  Example \ref{direxample}(1), we showed how to extend the definition of $\zeta(s)$ to all $s \in \CC \backslash \{1\}$ with $\operatorname{Re} s > 0$.   In his seminal paper \cite{rie} of 1859,  Riemann proved that the  function $\zeta(s)$ extends uniquely to an analytic function on $\CC\backslash \{1\}$, with a simple pole at $s = 1$ with residue $1$.  The resulting meromorphic function $\zeta(s)$ on $\CC$ is called the  {\bf Riemann zeta function}.\index{Riemann zeta function $\zeta(s)$}\index[symbols]{.t  A@$\zeta(s)$}    References on the Riemann zeta function, and for the facts noted in this section,  include \cite[Chapters 6 and 8]{bate} \cite{edw} \cite{ivic} \cite{patt} \cite{tit} \cite[Chapters 2 and 3]{zud}.

Also in  \cite{rie}, Riemann proved the {\bf functional equation}\index{functional equation for $\zeta(s)$}
$$ \zeta (s)=2^{s}\pi ^{s-1} \sin \left({\frac {\pi s}{2}}\right) \Gamma (1-s) \zeta (1-s), \quad \forall s \in \CC,$$
or, equivalently,
$$2\cos\left(\frac{\pi s}{2}\right) \zeta(s) = \frac{(2\pi)^{s}}{\Gamma(s)}\zeta(1-s),  \quad \forall s \in \CC,$$
where both sides of the equations assume their limiting values as $s \to 0$ and $s \to 1$.  For more contemporary proofs of the functional equation, see any of \cite[Theorem 12.7]{apos} \cite[Sections 1.6 and 1.7]{edw} \cite[Theorem 15]{ing2} \cite[Section 1.2]{ivic2}.

 Since every negative  integer $s = n$ is a simple zero of the reciprocal gamma function $\frac{1}{\Gamma(s)}$, while $\cos\frac{\pi s}{2} \neq 0$ if $s = n$ is even, the functional equation implies that every negative even integer is a simple zero of $\zeta(s)$.  In fact, these {\bf trivial zeros of $\zeta(s)$}  are all of the real zeros of $\zeta(s)$.   (This can be verified by showing that $\zeta(s)  > 1$ for all $s \in (1,\infty)$ and $\zeta(s) < 0$ for all $s \in [0,1)$ and then applying the functional equation.)  Riemann showed that the remaining zeros of $\zeta(s)$  lie in the  critical strip $\{s \in \CC: 0 \leq \operatorname{Re} s \leq 1\}$.   The {\bf Riemann hypothesis} is the conjecture,  posed by Riemann in his paper of 1859, that all of the  nontrivial (i.e., nonreal) zeros of $\zeta(s)$ have real part $\frac{1}{2}$.    The functional equation implies that, if $\rho$ is a zero of $\zeta(s)$ in the critical strip, then so is $1-\rho$ (as are $\overline{\rho}$ and $1-\overline{\rho}$).   We let $\Theta$ denote the  {\bf Riemann constant}
$$\Theta = \sup\{\operatorname{Re} s : s \in \CC, \, \zeta(s) = 0\}.$$ Since the nontrivial zeros lie in the critical strip are symmetric about the critical line,  and they are also infinite in number, one has $\tfrac{1}{2} \leq \Theta \leq 1$,
and the Riemann hypothesis is equivalent to $\Theta = \tfrac{1}{2}$.

Note that, by Theorem \ref{dirichlet}, one has
$$\zeta'(s) = -\sum_{n = 1}^\infty \frac{\log n}{n^s}$$
and
$$\frac{\zeta'(s)}{\zeta(s)} = -\sum_{n = 1}^\infty \frac{\Lambda(n)}{n^s}$$
for all $s$ with $\operatorname{Re} s > 1$.   
Since $\zeta(s)$ is meromorphic with a single simple pole at $s = 1$, the function $(s-1)\zeta(s)$ is entire.  More generally,  the function $\zeta(s)-\frac{1}{s-1}$ is entire, since one has $$\gamma = \lim_{s \to 1} \left( \zeta(s)-\frac{1}{s-1} \right).$$  It follows that there are real constants $\gamma_k$, where $\gamma_0 = \gamma$,  such that
$$\displaystyle \zeta (s)={\frac {1}{s-1}}+\sum_{k=0}^{\infty }{\frac {(-1)^{k}\gamma_k }{k!}}(s-1)^{k}, \quad \forall s \in \CC\backslash\{1\}.$$
In fact,  $\gamma_k$ for every nonnegative integer $k$ is precisely the $k$th Stieltjes constant, defined in Example \ref{stiec}(2) \cite[Theorem 1.3]{ivic}.    This provides one of the simplest expressions for the extension of the function $\zeta(s)$ on $\{s \in \CC: \operatorname{Re} s > 1\}$ to an analytic function on $\CC\backslash \{1\}$ (but defining $\zeta(s)$  this way requires verifying that the radius of convergence of the given Taylor series is infinite).   

 Notably, one has
$${\displaystyle \zeta (2n)=(-1)^{n+1}{\frac {(2\pi )^{2n}B_{2n}}{2(2n)!}}\!}$$
 and
$${\displaystyle \zeta^{\prime }(-2n)=(-1)^{n}{\frac {(2n)!}{2(2\pi )^{2n}}}\zeta (2n+1)}$$
(the latter of which can be deduced from the functional equation) for all positive integers $n$,  and also
$$\zeta(-n)  =(-1)^{n}{\frac {B_{n+1}}{n+1}} = \begin{cases} \displaystyle   -\frac{1}{2} & \quad \text{if } n = 0 \\
\displaystyle   0  & \quad \text{if } n > 0 \text{ is even} \\
  \displaystyle -\frac{B_{n+1}}{n+1} & \quad \text{if } n > 0 \text{ is odd}
\end{cases}$$
 \cite[Sections 1.4 and 1.5]{edw}. 
It is also known \cite[Sections 3.8 and 6.8]{edw} that
$$\zeta'(0) = -\frac{1}{2} \log 2 \pi = -\int_0^1 \log\Gamma(x) \, dx$$
and therefore
$$\frac{\zeta'(0)}{\zeta(0)} = \log 2 \pi.$$  

For all $s = x+ia$ with $a \in \RR$ fixed,  one has 
$$\zeta(x+ia) \sim 1 \ (x \to \infty)$$
and in fact
$$\zeta(x+ia) - 1  = 2^{-s}+3^{-s} + 4^{-s} + \cdots \sim 2^{-s} = 2^{-ia}2^{-x} \ (x \to \infty).$$
Moreover,  the functional equation and Stirling's approximation imply that
\begin{align*}
\zeta(1-s) & = 2\cos\left(\frac{\pi s}{2}\right) (2\pi)^{-s}\Gamma(s)\zeta(s) \\ & \sim 2\cos\left(\frac{\pi s}{2}\right){(2\pi)^{-s}} \Gamma(s) \\
& \sim 2\cos\left(\frac{\pi s}{2}\right){\sqrt {\frac {2\pi }{s}}}\,{\left({\frac {s}{2 \pi e}}\right)}^{s} \ (x \to \infty).
\end{align*}
This explains the extreme oscillating behavior of $\zeta(x)$ as $x \to -\infty$ as witnessed in Figure \ref{zeta1}: as $x \to \infty$, the factor $(2\pi)^{-x}\Gamma(x) \sim {\sqrt {\frac {2\pi }{x}}}\,{\left({\frac {x}{2 \pi e}}\right)}^{x}$  of  $\zeta(1-x)$ grows super-exponentially, while the cosine factor oscillates between $1$ and $-1$.

 Another important fact is that $\zeta(s)$ can be expressed in terms of its zeros.  One has the factorization
\begin{align}\label{ingz3}
(s-1)\zeta (s)=\pi ^{\frac {s}{2}}\frac{1}{\Gamma \left(1+{\frac {s}{2}} \right)}\xi(s)
\end{align}
for a unique entire function $\xi(s)$ called the {\bf Riemann xi function}.\index{Riemann xi function $\xi(s)$}\index[symbols]{.t  C@$\xi(s)$}  In 1893 \cite{had1}, Hadamard made rigorous Riemann's argument in \cite{rie} that
\begin{align}\label{ingz4}
\xi(s) = \frac{1}{2}\prod_{\rho }\left(1-{\frac {s}{\rho }}\right),
\end{align}
where the product is over all nontrivial zeros $\rho$ of $\zeta(s)$, repeated to multiplicity, and the factors for a pair of zeros of the form $\rho$ and $1-\rho$ are combined.  This result (also proved, for example, in \cite[Sections 1.8--1.10]{edw})  was an important step in the 1896 proofs of the prime number theorem three years later.  Note that the functions $(s-1)\zeta(s)$, $\pi^{\frac{s}{2}}$, $\frac{1}{\Gamma\left(1 +\frac{s}{2}\right)}$, and $\xi(s)$ appearing in equation (\ref{ingz3}) are all entire.   Moreover, the functional equation for $\zeta(s)$ is equivalent to
$$\xi(s) = \xi(1-s)$$ for all $s \in \CC$.  Equivalently, one has
$$\Xi(s) = \Xi(-s)$$ for all $s \in \CC$,
where $$\Xi(s) = \xi(\tfrac{1}{2}+is)$$ is the {\bf Riemann Xi function}.\index{Riemann Xi function $\Xi(s)$}\index[symbols]{.t  B@$\Xi(s)$}    The functions $\xi(s)$ and $\Xi(s)$ are real for real values  of $s$, and the function $\Xi(s)$ is also real for purely imaginary values of $s$.   Moreover, it is clear from the facts mentioned above that the Riemann hypothesis is equivalent to the statement that all of the zeros of the function $\Xi(s)$ are real (and are equal to $\operatorname{Im} \rho$, where $\rho$ ranges over the nontrival zeros $\rho$ of $\zeta(s)$).

\begin{figure}[ht!]
\includegraphics[width=75mm]{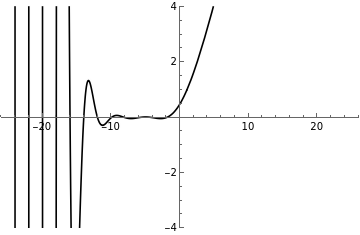}
\caption{\centering Graph of $(x-1)\zeta(x)$ on $[-25,25]$}
\label{zeta3}
\end{figure}

\begin{figure}[ht!]
\includegraphics[width=75mm]{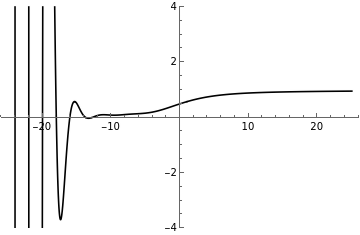}
\caption{\centering Graph of $\zeta(x)-\frac{1}{x-1}$ on $[-25,25]$}
 \label{zeta2}
\end{figure}

\begin{figure}[ht!]
\includegraphics[width=75mm]{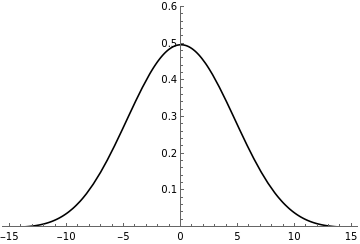}
\caption{\centering Graph of $\Xi(x) = \xi(\tfrac{1}{2}+xi)$  on $[-15,15]$}
 \label{Xi2}
\end{figure}

\begin{figure}[ht!]
\includegraphics[width=75mm]{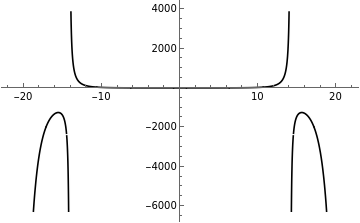}
\caption{\centering Graph of and $\frac{1}{\Xi(x)}$ on $[-26,26]$}
 \label{Xi3}
\end{figure}

\begin{figure}[ht!]
\includegraphics[width=35mm]{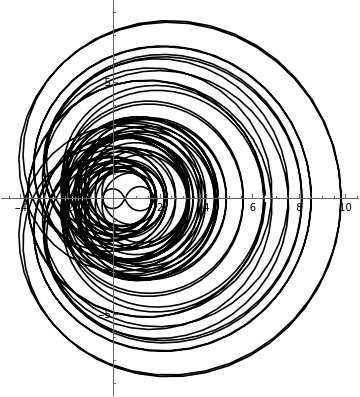}
\caption{\centering Parametric plot of $\zeta(0+it)$ on $[-100,100]$}   
\label{zeta1_0}
\end{figure}

\begin{figure}[ht!]
\includegraphics[width=35mm]{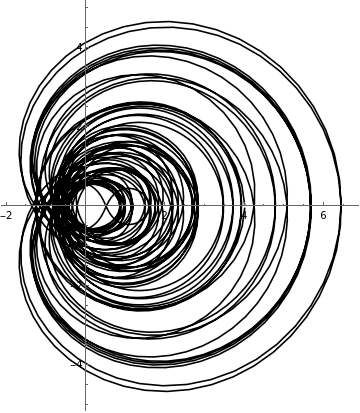}
\caption{\centering Parametric plot of $\zeta(\tfrac{1}{4}+it)$ on $[-100,100]$}
\label{zeta1_14}
\end{figure}

\begin{figure}[ht!]
\includegraphics[width=40mm]{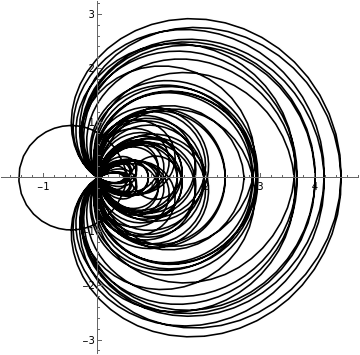}
\caption{\centering Parametric plot of $\zeta(\tfrac{1}{2}+it)$ on $[-100,100]$}
  \label{zeta1_12}
\end{figure}

\begin{figure}[ht!]
\includegraphics[width=60mm]{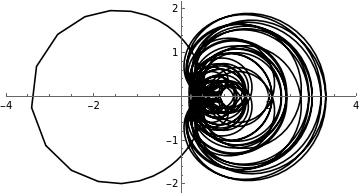} 
\caption{\centering Parametric plot of $\zeta(\tfrac{3}{4}+it)$ on $[-100,100]$}
  \label{zeta1_34}
\end{figure}

\begin{figure}[ht!]
\includegraphics[width=45mm]{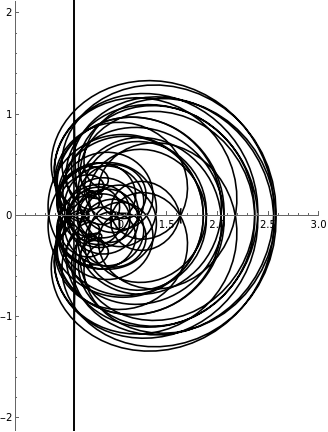}
\caption{\centering Parametric plot of $\zeta(1+it)$ on $[-100,100]$}
  \label{zeta1_1}
\end{figure}

\begin{figure}[ht!]
\includegraphics[width=55mm]{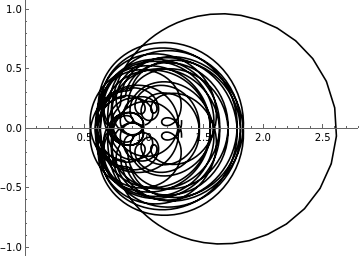}
\caption{\centering Parametric plot of $\zeta(\tfrac{3}{2}+it)$ on $[-100,100]$}
   \label{zeta1_32}
\end{figure}

\begin{figure}[ht!]
\includegraphics[width=55mm]{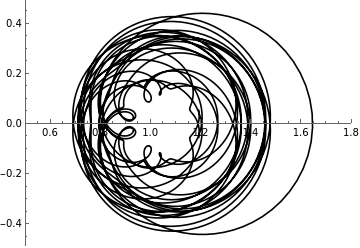}
\caption{\centering Parametric plot of $\zeta(2+it)$ on $[-100,100]$}
   \label{zeta1_2}
\end{figure}

\begin{figure}[ht!]
\includegraphics[width=40mm]{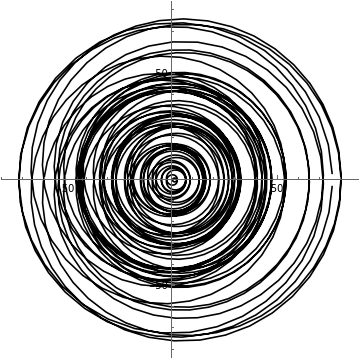}
\caption{\centering Parametric plot of $\zeta(-1+it)$ on $[-100,100]$}
   \label{zeta1_-1}
\end{figure}

\begin{figure}[ht!]
\includegraphics[width=45mm]{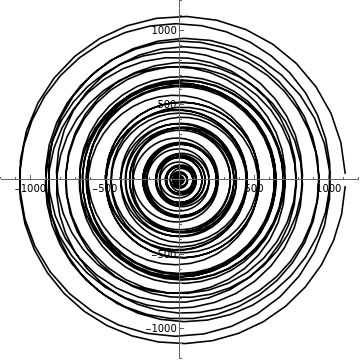}
\caption{\centering Parametric plot of $\zeta(-2+it)$ on $[-100,100]$}
   \label{zeta1_-2}
\end{figure}

 See Figures  \ref{zeta3} and \ref{zeta2} for graphs of the entire functions and $(x-1)\zeta(x)$ and  $\zeta(x)-\frac{1}{x-1}$, respectively, on $[-25,25]$,  Figure \ref{Xi2} for a graph of the entire function $\Xi(x) = \xi(\tfrac{1}{2}+xi)$ on $[-15,15]$, and  Figure \ref{Xi3} for a graph of the function $\frac{1}{\Xi(x)}$ on $[-26,26]$.   It is a widely held conjecture that all of the (nontrivial) zeros of $\zeta(s)$ are simple, or, equivalently, that all of the zeros of $\Xi(s)$ are simple.  This conjecture implies that the function $\Xi(x)$ changes signs between each successive imaginary part of the nontrivial zeros of $\zeta(s)$.  Finally, Figures \ref{zeta1_0} through \ref{zeta1_-2} are parametric plots of $|\zeta(a+ti)|$ for various fixed $a$ and $t \in [-100,100]$.  These graphs are meant to exhibit various zeros and zero-free regions of $\zeta(s)$.

\nopagebreak

\section{The prime zeta function $P(s)$}

Recall from Section 3.8 that the prime zeta function $P(s)$ is defined by $$P(s) = \sum_p \frac{1}{p^s}$$
for all $s \in \CC$ with $\operatorname{Re} s > 1$.   There are many known  relationships between the functions $P(s)$ and $\zeta(s)$.   For example, we have seen that
$$\zeta(s) P(s) = \sum_{n = 1}^\infty \frac{\omega(n)}{n^s},$$
$$\zeta(s)\sum_{n = 1}^\infty P(ns) = \sum_{n = 1}^\infty \frac{\Omega(n)}{n^s},$$
$${\displaystyle \log \zeta (s)=\sum_{n = 1}^\infty{\frac {P(ns)}{n}}},$$
and
$$P(s)=\sum_{n = 1}^\infty \frac {\mu(n)}{n}\log \zeta (ns),$$
for all $s \in \CC$ with $\operatorname{Re} s > 1$.   The last of these equations allows one to extend $P(s)$  to a  function on $\{s \in \CC: \operatorname{Re} s > 0\}$,\index[symbols]{.t  B@$P(s)$}\index{prime zeta function $P(s)$}    but with logarithmic singularities at $s$ for all squarefree integers $n$ such that $ns$ is a pole or zero of the Riemann zeta function, i.e., at $\frac{1}{n}$ and at $\frac{\rho}{n}$ for every squarefree integer $n$ and for every nontrivial zero $\rho$ of $\zeta(s)$.  Thus,  for every squarefree integer $n$, the function $P(s)$ has at least one singularity on the line $\{s \in \CC: \operatorname{Re} s = \frac{1}{n}\}$ if $n$ is odd (namely, $\frac{1}{n}$)  and infinitely many singularities on the lines $\{s \in \CC: \operatorname{Re} s = \frac{1}{n}\}$  and $\{s \in \CC: \operatorname{Re} s = \frac{1}{2n}\}$ if $n$ is even (namely, $\frac{2\rho}{n}$ and $\frac{\rho}{n}$, respectively, for each zero  $\rho$ of $\zeta(s)$ with real part $\frac{1}{2}$). Moreover, if the Riemann hypothesis holds, then these are all of the singularities of $P(s)$.

By Example \ref{primez}, the function $P(s)-\log \zeta(s)$ can be analytically continued to $\{s \in \CC: \operatorname{Re} s> \frac{1}{2}\}$.  If one defines
$$a_n = \prod_{p^k |n} \frac{1}{k} $$ for all $n$,
then one has
$$\exp P(s) = \sum_{n= 1}^\infty \frac{a_n}{n^s} $$
for all $s$ with $\operatorname{Re} s > 1$  \cite[Lemma 2.7]{li1}.   Thus,  like the function $\log \zeta(s)$, the prime zeta function $P(s)$ is a logarithm of a Dirichlet series.

\begin{figure}[ht!]
\includegraphics[width=70mm]{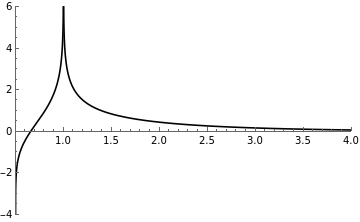}
\caption{\centering Graph of $\operatorname{Re}P(x)$ (which equals $P(x)$ for $x > 1$) on $(\tfrac{1}{2},4]$}
\label{primezeta0}
\end{figure}

\begin{figure}[ht!]
\includegraphics[width=70mm]{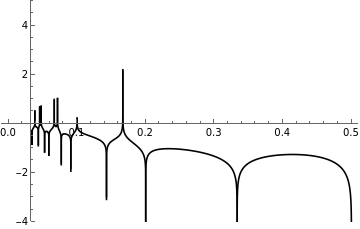}
\caption{\centering Graph of $\operatorname{Re}P(x)$ on $[\tfrac{1}{30},\tfrac{1}{2})$}
   \label{primezeta1}
\end{figure}

\begin{figure}[ht!]
\includegraphics[width=70mm]{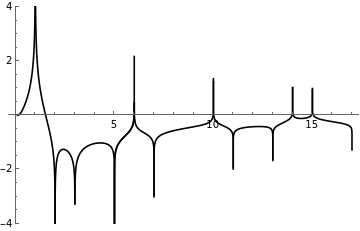}
\caption{\centering Graph of $\operatorname{Re}P(1/x)$ on $(0,17)$}
   \label{primezeta2}
\end{figure}

It is straighforward to show, using the fact that $\operatorname{Im}(\log x) = \pi$ if $x < 0$, that
$$\operatorname{Im} P(x) = \pi \sum_{n \leq \lfloor 1/x \rfloor} \frac{\mu(n)}{n}$$
for all $x \in (0,1)$ such that $P(x)$ is defined,  that is, such that $x \neq \frac{1}{n}$ for all squarefree positive integers $n$.  Since one has
$$\sum_{n = 1}^\infty \frac{\mu(n)}{n} = 0,$$
 it follows that
$$\lim_{x \to 0^+} \operatorname{Im} P(x) = 0$$
 on the domain of $P(x)$.       Moreover, for all $s = x+ia$ with $a$ fixed,  one has 
$$P(x+ia)  = 2^{-s}+3^{-s} + 5^{-s} + 7^{-s} + \cdots \sim 2^{-s}  = 2^{-ia}2^{-x} \sim \zeta(x+ia) - 1 \ (x \to \infty)$$
and therefore 
$$\lim_{x \to \infty} P(x+ia) = 0,$$
and thus, in particular,
$$\lim_{x \to \infty} P(x) = 0.$$
Figures \ref{primezeta0} and \ref{primezeta1} provide graphs of the function $\operatorname{Re} P(x)$ (which equals $P(x)$ for $x > 1$) on  $(\tfrac{1}{2},4]$ and $[\tfrac{1}{30},\tfrac{1}{2})$, respectively.  On this scale,  the graph was impossible to see for $x < \tfrac{1}{30}$.  Thus, it makes more sense to graph the function $\operatorname{Re} P(1/x)$, which in Figure \ref{primezeta2} is graphed on the interval $(0,17)$.  If one sets $\operatorname{Re} P(1/0) = 0$, then this function is infinitely differentiable on $[0,\infty)$ except at the squarefree positive integers, and all of its derivatives at $0^+$ are $0$.

\section{The functions $\Ein(s)$, $\Ei(s)$, and $E_1(s)$}

It is clear that the function
$$\frac{1-e^{-s}}{s} = \sum_{n=1}^{\infty }{\frac {(-1)^{n+1}s^{n-1}}{n!}}$$
is entire, where $\frac{1-e^{-s}}{s}$ assumes the limiting value  $1$ at $0$.  It follows that its antiderivative
$$ \Ein(s) =\int_{0}^{s}\frac{1-e^{-z}}{z} \, dz=\sum_{n=1}^{\infty }{\frac {(-1)^{n+1}s^{n}}{n! \, n}}, \quad \forall s \in \CC,$$
which is called the {\bf  complementary exponential integral function} \cite[Chapter 2 Section 3]{olver},\index{complementary exponential integral function $\Ein(s)$} is also entire.\index[symbols]{.t  F@$\Ein(s)$}  It is readily checked that
$$\Ein(s) = e^{-s}\sum_{n=1}^{\infty } H_{n}{\frac {s^{n}}{n!}}. \quad \forall s \in \CC.$$
The function $\Ein(s)$ is the unique entire function with $\Ein(0)= 0$ and with derivative
 $$\frac{d}{ds} \Ein(s) = \frac{1-e^{-s}}{s}, \quad \forall s \in \CC.$$
Equivalently, $\Ein(s)$ is the unique entire function with $\Ein(0) = 0$ and $$\Ein^{(n)} (0)= \frac{(-1)^{n+1}}{n}$$ for all positive integers $n$.  The function $\Ein(s)$ has the special values
$$\Ein(1) = -\int_0^1 \frac{\log x}{\exp x} \, dx =  0.796599599297\ldots$$
and
$$\Ein(-1)  = \int_0^1 \exp x  \log x  \, dx = \int_1^e \log \log x \, dx =  -1.317902151454\ldots,$$
and in fact one has
$$\Ein(s) = -s\int_0^1 e^{-sx} \log x  \, dx, \quad \forall s \in \CC.$$

The {\bf exponential  integral function} $E_1(s)$\index{exponential integral function $E_1(s)$} is defined by
$$E_1(s) = - \gamma-\Log s+\Ein(s),  \quad \forall s \in \CC\backslash\{0\}.\index[symbols]{.t  Ga@$E_1(s)$}$$
  An equivalent definition is
$$E_1(s) = \int_{s}^\infty \frac{e^{-z}}{z} dz, \quad \forall s \in \CC \backslash (-\infty, 0],$$
where the integral is along any path of integration not crossing $(-\infty, 0]$, with
$$E_1(x) := \lim_{\varepsilon \to 0^+} E_1(x+\varepsilon i), \quad \forall x  \in  (-\infty, 0),$$ and $$\lim_{\varepsilon \to 0^-} E_1(x+\varepsilon i) = \overline{E_1(x)}, \quad \forall x  \in  (-\infty, 0)$$  \cite[Chapter 2 Section 3]{olver} \cite[Chapter 14]{cuyt}.  The function $E_1(x)$ for nonzero real $x$ is given by
\begin{align*} E_1(x) =   \left.
  \begin{cases}
   -\li(e^{-x}) & \text{if } x > 0 \\
    -\li(e^{-x})-\pi i   & \text{if } x < 0,
 \end{cases}
\right.
\end{align*}
where $\li(x)$ is the logarithmic integral function.

\begin{figure}[ht!]
\includegraphics[width=70mm]{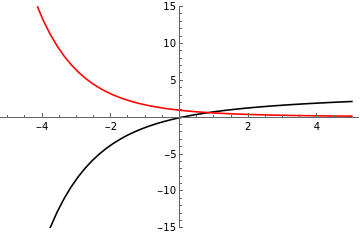}
\caption{\centering Graph of $\Ein(x)$ (in black) and $\frac{d}{dx}\Ein(x) = \frac{1-e^{-x}}{x}$  (in red) on $[-5,5]$}
 \label{Ein}
\end{figure}

\begin{figure}[ht!]
\includegraphics[width=70mm]{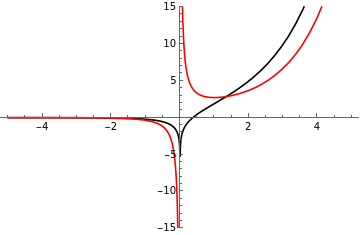}
\caption{\centering Graph of $\Ei(x)$ (in black) and $\frac{d}{dx}\Ei(x) = \frac{e^{x}}{x}$  (in red) on $[-5,5]$} \label{Ei}
\end{figure}

The related function
\begin{align*} 
\Ei(s) & =  -E_1(-s)-\frac{\log(-s)+\log(-1/s)}{2}
\\   & =   \left.
 \begin{cases}
     -E_1(-s)& \text{if } s \in \CC \backslash [0,\infty) \\
      -E_1(-s)-\pi i  & \text{if }  s \in (0,\infty)
 \end{cases} 
 \right.  
\\   & =   \left.
 \begin{cases}
    \gamma+\Log(-s)-\Ein(-s) & \text{if } s \in \CC \backslash [0,\infty) \\
       \gamma+\Log s   -\Ein(-s)  & \text{if }  s \in (0,\infty)
 \end{cases} 
 \right.  
\end{align*}
is called the {\bf exponential integral function $\Ei(s)$}.\index{exponential integral function $\Ei(s)$}\index[symbols]{.t  G@$\Ei(s)$} 
An equivalent definition is
$$\Ei(s) = \int_{-\infty}^s \frac{e^{z}}{z} dz, \quad \forall s \in \CC \backslash [0,\infty),$$
where the integral is along any path of integration not crossing $[0,\infty)$, with
$$\Ei(x) := \lim_{\varepsilon \to 0} \frac{ \Ei(x+\varepsilon i)+\Ei(x-\varepsilon i)}{2}, \quad \forall x \in (0, \infty),$$ 
so that
$$\Ei(x) = \li(e^x), \quad \forall x \in \RR\backslash \{0\},$$
and thus
$$\frac{d}{dx} \Ei(x) = \frac{e^x}{x}, \quad \forall x \in \RR\backslash \{0\},$$
Note also that
$$\frac{d}{ds} E_1(s) = \frac{1-e^{-s}}{s}- \frac{1}{s}  = -\frac{e^{-s}}{s}, \quad \forall s \in \CC \backslash (-\infty, 0].$$
and
$$\frac{d}{ds} \Ei(s) = \frac{1}{s} -\frac{1-e^{s}}{s} = \frac{e^s}{s},  \quad \forall s \in \CC \backslash [0,\infty).$$
See Figure \ref{Ein} (resp., Figure \ref{Ei}) for a graph of the function  $\Ein(x)$ (resp., $\Ei(x)$) and its derivative on the interval $[-5,5]$.

We also note the following.

\begin{proposition}[{\cite[(14.1.13)]{cuyt} \cite[Chapter 1 Example 1]{wong}}]\label{EIlem}
One has
$$E_1 (s) \sim \frac{e^{-s}}{s} \ (s \to \infty),$$
and, more generally, the asymptotic expansion
$$e^{s}E_1 (s) \simeq  \sum_{n = 0}^\infty \frac{(-1)^n n!}{s^{n+1}}  \ (s \to \infty)$$
on $\{s \in \CC: |\operatorname{Arg} s| \leq \pi-\varepsilon\}$ for all $\varepsilon > 0$.  Moreover, one has
$$\Ei (s) \sim \frac{e^{s}}{s}  \ (s \to \infty),$$
and, more generally, the asymptotic expansion
$$e^{-s} \Ei (s) \simeq  \sum_{n = 0}^\infty \frac{n!}{s^{n+1}}  \ (s \to \infty),$$
on  $\{s \in \CC: |\operatorname{Arg} s| \geq \varepsilon\}$ for all $\varepsilon > 0$ and on $(0,\infty)$.
\end{proposition}

\begin{remark}[Caveats regarding Mathematica's functions $\Ei(s)$ and $\li(s)$]\label{lieiremark}
Note that Mathematica's exponential integral function $\Ei(s)$,  defined by
\begin{align*} 
\Ei(s)  = \gamma+\Log s-\frac{\Log(s)+\Log(1/s)}{2}-\Ein(-s), \quad \forall s \in \CC\backslash \{0\},
\end{align*}
is different than ours.   In particular,  Mathematica's function $\Ei(s)$, minus our function $\Ei(s)$, is precisely 
\begin{align*} 
\Log s-\frac{\Log(s)+\Log(1/s)}{2} - \Log (-s)+\frac{\Log(-s)+\Log(-1/s)}{2} = \pi i \operatorname{sgn} (\operatorname{Im} s), \quad \forall s \in \CC\backslash \{0\}.
\end{align*}
Under our definition of $\Ei(s)$, for example, one has
 $$\lim_{x \to \infty}\Ei(a+ x i ) = 0, \quad \forall a \in \RR,$$ 
while, under Mathematica's definition,  one has
 $\lim_{x \to \infty}\Ei(a\pm x i ) = \pm \pi i$ for all $a \in \RR$, and  the analogue of Proposition \ref{EIlem} does not hold.  Moreover, Mathematica's version of $\Ei(s)$ is less suitable for use in the Riemann--von Mangoldt explicit formula (as described in Sections 1.3 and 5.1).  
For documentation on Mathematica's  function $\Ei(s)$, see https://functions.wolfram.com/GammaBetaErf/ExpIntegralEi/ (accessed by the author on 1 April 2024).

In traditional expressions of  the Riemann--von Mangoldt explicit formula,  what should be written as
$\Ei(\rho \log x)$ for $\rho \in \CC$ is written, somewhat misleadingly, as $\li(x^\rho)$.   Note, in particular,  that
$\Ei(\rho \log x)$ is not a function of $x^\rho$, while $\li(x^\rho)$ appears to be a function of $x^\rho$, at least under a naive intepretation of the expression $\li(x^\rho)$.   For example, Mathematica extends the logarithmic integral function $\li(x)$ to a complex function by setting $\li(s) =  \Ei(\Log s)$, thus interpreting $\li(x^\rho)$ as  $\Ei(\Log (x^\rho))$, rather than as $\Ei(\rho \log x)$. To avoid confusion, we avoid using the expression $\li(s)$ for nonreal $s$ altogether and instead utilize the more general function $\Ei(s)$.
\end{remark}

\section{The functions $\ERi(s)$ and $\Ri(x)$}

Let  us define $\ERi(s)$ to be the function
$$\ERi(s)= \begin{cases} \displaystyle \sum_{n=1}^\infty \frac{ \mu(n)}{n} \Ei\left(\frac{s}{n} \right) & \text{if } s \neq 0 \\
1 & \text{if } s = 0.
\end{cases}$$\index[symbols]{.t  Gc@$\ERi(s)$}

\begin{proposition}\label{ERiprop}
The function $\ERi(s)$ is entire and  has the Taylor series representation
$$\ERi(s) = 1+ \sum_{n = 1}^\infty \frac{s^n}{n! \, n \, \zeta(n+1)}$$
on all of $\CC$.
 Thus, $\ERi(s)$ is the unique entire function such that $\ERi(0) = 1$ and $\ERi^{(n)}(0) = \frac{1}{n \, \zeta(n+1)}$ for all positive integers $n$.  Moreover, one has
$$\ERi(s) = 1  -  \sum_{n = 1}^\infty \frac{\mu(n)}{n}\Ein(-s/n)$$
for all $s \in \CC$.
\end{proposition}

\begin{proof}
Let $s \in \CC\backslash \{0\}$.  Using the facts that $\sum_{n = 1}^\infty \frac{\mu(n)}{n}  = 0$ and $\sum_{n = 1}^\infty \frac{\mu(n)\log n}{n} = -1$ and $\Log(xs) = \log x+ \Log s$ for all $x > 0$, one has
\begin{align*}
\ERi(s) & = \sum_{n = 1}^\infty \frac{\mu(n)}{n}  \Ei\left(\frac{s}{n} \right) \\
& =  \sum_{n = 1}^\infty \frac{\mu(n)}{n}  \left( \gamma+ \Log( -s/n)-\frac{\Log(-s/n)+\Log(-n/s)}{2} - \Ein(-s/n) \right) \\
& =  \sum_{n = 1}^\infty \frac{\mu(n)}{n}  \left( \gamma+ \Log(- s)-\log n-\frac{\Log(-s)+\Log(-1/s)}{2} - \Ein(-s/n) \right) \\
& =   \left( \gamma+\frac{\Log(-s)-\Log(-1/s)}{2}\right)\sum_{n = 1}^\infty \frac{\mu(n)}{n}   - \sum_{n = 1}^\infty \frac{\mu(n)\log n}{n} -  \sum_{n = 1}^\infty \frac{\mu(n)}{n}\Ein(-s/n) \\
& =  1  -  \sum_{n = 1}^\infty \frac{\mu(n)}{n}\Ein(-s/n) \\
  & = 1 + \sum_{n = 1}^\infty \sum_{k = 1}^\infty  \frac{\mu(n)}{n} \frac{(s/n)^k}{k! \, k}  \\
   &= 1+ \sum_{k = 1}^\infty \sum_{n = 1}^\infty \frac{\mu(n)}{n^{k+1}} \frac{s^k}{k! \, k} \\
   &= 1+ \sum_{k = 1}^\infty \frac{s^k}{k! \, k \, \zeta(k+1)},
\end{align*}
where the last three sums are absolutely convergent, and are therefore equal, since
$$ \sum_{k = 1}^\infty \sum_{n = 1}^\infty \frac{1}{n^{k+1}} \frac{s^k}{k! \, k} =  \sum_{k = 1}^\infty \frac{\zeta(k+1)}{k! \, k} s^k $$
converges for all  $s \in \CC$. 
\end{proof}

\begin{figure}[htb!]
\includegraphics[width=70mm]{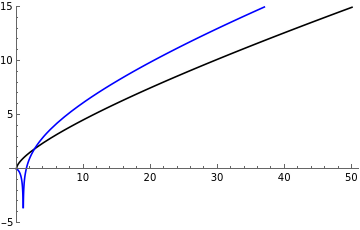}
\caption{\centering Graph of $\Ri(x)$ (in black) and $\li(x)$ (in blue) on $[0,50]$}
 \label{liri2b}
\end{figure}

\begin{figure}[htb!]
\includegraphics[width=70mm]{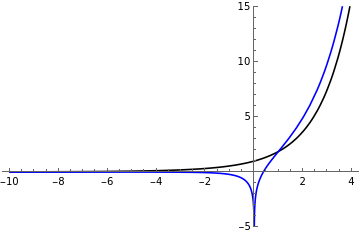}
\caption{\centering Graph of  $\ERi(x) = \Ri(e^x)$ (in black) and $\Ei(x) = \li(e^x)$ (in blue)  on $[-10,4]$}  \label{liri1}
\end{figure}

Let $\Ri(x)$\index[symbols]{.t  Gd@$\Ri(x)$}  denote the function
$$\Ri(x) =  \begin{cases} \displaystyle   \ERi(\log x)  = \sum_{n=1}^\infty \frac{ \mu(n)}{n} \li(x^{1/n}) & \quad \text{if } x >0  \\
 0 & \quad \text{if } x = 0.
\end{cases}$$
By Proposition \ref{ERiprop}, the function $\Ri(x)$ has the series representation
$$\Ri(x) = 1+ \sum_{n = 1}^\infty \frac{(\log x)^n}{n \cdot n!\, \zeta(n+1)}.$$
The series above is known as  {\bf Gram series representation of $\Ri(x)$}\index{Gram series representation of $\Ri(x)$}.   See Figure \ref{liri2b} for a comparison of the graphs of $\Ri(x)$ and $\li(x)$, and see Figure \ref{liri1} for a comparison of the graphs of $\ERi(x) = \Ri(e^x)$ and $\Ei(x) = \li(e^x)$.

From the graph of $\Ri(x)$ one might guess that $\Ri(x)$ has only one zero, namely, at $x = 0$.  However, it is a rather shocking fact  \cite{born}  that it has many zeros, and very likely an infinite number of zeros, just to the right of $x = 0$, where the {\it largest zero} of the function $\Ri(x)$ is 
$$x = 1.828643269752522610 0973252\ldots \cdot 10^{-14828}.$$
Equivalently, the  largest zero of the function $\ERi(x) = \Ri(e^x)$ is $$x = -34142.12818 46064 64992 63004 60679 8\ldots.$$   This amazing fact provides yet another stark warning  against basing conjectures  in analytic number theory on numerical and graphical evidence alone.

Now, for all $c > 0$, let
\begin{align*}
\operatorname{ER}(s,c) = \sum_{n \leq c|s|} \frac{\mu(n)}{n} \Ei\left(\frac{s}{n} \right),  \quad \forall s \in \CC\backslash\{0\}.
\end{align*}
so that
\begin{align*}
\ERi(s) = \lim_{c \to \infty} \operatorname{ER}(s,c), \quad \forall s \in \CC\backslash\{0\}.
\end{align*}
 By  M\"obius inversion (specifically, by Remark \ref{secondgen}), one has
\begin{align*}
\Ei(s) = \sum_{n \leq c|s|} \frac{1}{n} \operatorname{ER}\left(\frac{s}{n},c \right)
\end{align*}
for all $c > 0$ and all $s \in \CC$ with  $|s| \geq \frac{1}{c}$.

\begin{proposition}\label{ERIlem}
Let $c>0$ and $\varepsilon > 0$, and let
\begin{align*} h(x) =   \left.
  \begin{cases}
    \frac{1}{\sqrt{x}} e^{(\log x)^{1/2+\varepsilon}} & \text{if the Riemann hypothesis holds} \\
     e^{-(\log x)^{3/5-\varepsilon}}   & \text{if the anti-Riemann hypothesis holds} \\
    x^{\Theta-1+\varepsilon} & \text{otherwise.} 
 \end{cases}
 \right.
\end{align*}
One has
$$\ERi(s) = \operatorname{ER}(s,c) + O(h(c|s|)) \ (s \to \infty).$$
In particular, one has
$$\ERi(s) = \operatorname{ER}(s,c) + O\left(\frac{1}{(\log |s|)^t} \right) \ (s \to \infty)$$
for all $t >0$.
\end{proposition}

\begin{proof}
By the well-known results (\ref{sbd}) and Corollary \ref{walf} of Section 10.2 concerning the Mertens function $M(x)$, one has
$$M(x) = O(xh(x)) \ (x \to \infty).$$
By Abel's summation formula (Theorem \ref{abels}) and Karamata's integral theorem (Theorem \ref{karam}), it follows that
$$\sum_{n > x} \frac{\mu(n)}{n}  = O\left(h(x)\right) \ (x \to \infty)$$
and in fact
$$\sum_{n > x} \frac{\mu(n)}{n^{k+1}}  = O\left(\frac{h(x)}{x^k}\right) \ (x \to \infty)$$
for all nonnegative integers $k$, where the $O$ constant does not depend on $k$.
Similarly, by  \cite[Theorem 1.4]{rama} and  Karamata's integral theorem, one has
$$\sum_{n > x} \frac{\mu(n)}{n} \log(x/n) = O(h(x))  \ (x \to \infty).$$
Therefore, by the series representation
$$\Ei(s) =   \gamma+  \Log(- s)-\frac{\Log (-s)+\Log(-1/s)}{2} + \sum_{k =1}^\infty \frac{s^k}{k! \, k}, \quad \forall s \in \CC\backslash\{0\},$$
one has
\begin{align*}
\sum_{n > c|s|} \frac{\mu(n)}{n} \Ei(s/n) &  = \sum_{n > c|s|} \frac{\mu(n)}{n} \left(\gamma + \Log(-s/n)-\frac{\Log(-s)+\Log(-1/s)}{2} \right) +   \sum_{n > c|s|} \sum_{k =1}^\infty \frac{\mu(n)}{n^{k+1}} \frac{ s^k}{k! \, k}  \\
 &  = \sum_{n > c|s|} \frac{\mu(n)}{n} \Log( -s/n) + O(h(c|s|))+  \sum_{k =1}^\infty \left(\sum_{n > c|s|} \frac{\mu(n)}{n^{k+1}}\right) \frac{ s^k}{k! \, k}  \\
  &  = \sum_{n > c|s|} \frac{\mu(n)}{n} \Log(- s/n) + O(h(c|s|))+  O\left(\sum_{k =1}^\infty \frac{h(c|s|)}{c^k|s|^{k}} \frac{ s^k }{k! \, k} \right) \\
    &  = \sum_{n > c|s|} \frac{\mu(n)}{n} \Log( -s/n) + O(h(c|s|)) \\
  &  = \sum_{n > c|s|} \frac{\mu(n)}{n} \log(c |s|/n) + \sum_{n > c|s|} \frac{\mu(n)}{n}(i \operatorname{Arg}(-s)-\log c) + O(h(c|s|)) \\
  &  = O(h(c|s|)) \ (s \to \infty).
\end{align*}
The proposition follows.
\end{proof}

\begin{proposition}\label{RAprop0}
One has the asymptotic expansion
\begin{align*}
\ERi(s) \simeq \sum_{n = 1}^\infty \frac{\mu(n)}{n} \Ei \left(\frac{s}{n}\right) \ (s \to \infty)
\end{align*}
on  $\{s \in \CC: |\operatorname{Arg} s| \geq \varepsilon\}$ for all $\varepsilon > 0$ and on $(0,\infty)$.
\end{proposition}

\begin{proof}
Let $\varepsilon > 0$, and let $N > 1$ be a fixed positive integer.  By  Propositions \ref{ERIlem} and \ref{EIlem}, for all  $s \in \{s \in \CC: |\operatorname{Arg} s| \geq \varepsilon\}$ (resp., $s \in (0,\infty)$) with $|s|> N$, one has
\begin{align*}\ERi(s) & = \sum_{n \leq |s|} \frac{\mu(n)}{n} \Ei\left(\frac{s}{n} \right) + o(1)  \\
  & = \sum_{n = 1}^{N} \frac{\mu(n)}{n} \Ei\left(\frac{s}{n} \right) + \sum_{N< n \leq |s| } \frac{\mu(n)}{n} \Ei\left(\frac{s}{n} \right) + o(1) \\
& = \sum_{n = 1}^{N} \frac{\mu(n)}{n} \Ei\left(\frac{s}{n} \right) + O\left(  \frac{e^{s/(N+1)}}{s/(N+1)} |s|\right) + o(1) \\
& = \sum_{n = 1}^{N} \frac{\mu(n)}{n} \Ei\left(\frac{s}{n} \right) + O\left(  e^{s/(N+1)}\right) + o(1)  \\
& = \sum_{n = 1}^{N-1} \frac{\mu(n)}{n} \Ei\left(\frac{s}{n} \right) + O\left( \frac{\mu(N)}{N} \Ei\left(\frac{s}{N} \right) \right) \ (s \to \infty).
\end{align*}
The proposition follows.
\end{proof}

\begin{corollary}\label{elll00}
One has
\begin{align*} \Ei(s) \sim \ERi(s)  \sim \frac{e^s}{s} \ (s \to \infty),
\end{align*}
\begin{align*} \Ei(s)-\ERi(s) \sim \frac{e^{s/2}}{s} \ (s \to \infty),
\end{align*}
and
\begin{align*}
\Ei(s)-\ERi(s) - \frac{1}{2}\Ei\left(\frac{s}{2} \right) \sim \frac{1}{3}\Ei\left(\frac{s}{3} \right) \sim \frac{e^{s/3}}{s} \ (s \to \infty)
\end{align*}
on  $\{s \in \CC: |\operatorname{Arg} s| \geq \varepsilon\}$ for all $\varepsilon > 0$ and on $(0,\infty)$.
\end{corollary}

Since $\lim_{n \to \infty} \operatorname{ERi}\left(\frac{s}{n} \right)  =  \ERi(0) =  1$, the sum $\sum_{n = 1}^\infty \frac{1}{n} \operatorname{ERi}\left(\frac{s}{n}\right)$ diverges, for all $s \in \CC$.   Nevertheless, one has the following.

\begin{proposition}\label{RApropaaa}
One has the (divergent) asymptotic expansion
\begin{align*}
\Ei(s) \simeq \sum_{n = 1}^\infty \frac{1}{n} \ERi \left(\frac{s}{n}\right) \ (s \to \infty)
\end{align*}
on  $\{s \in \CC: |\operatorname{Arg} s| \geq \varepsilon\}$ for all $\varepsilon > 0$ and on $(0,\infty)$.
\end{proposition}

\begin{proof}
Let $\varepsilon > 0$, and let $N > 1$ be a fixed positive integer.  For all $s \in \{s \in \CC: |\operatorname{Arg} s| \geq \varepsilon\}$ (resp., $s \in (0,\infty)$) with $|s|> N$, one has
\begin{align*}
\Ei(s) - \sum_{n = 1}^N  \frac{1}{n} \ERi \left( \frac{s}{n} \right) &  =  \sum_{N < n \leq |s|} \frac{1}{n} \ERi \left( \frac{s}{n} \right)+ \sum_{n \leq |s|}  \frac{1}{n} \left(\operatorname{ER}\left( \frac{s}{n} ,1\right) - \ERi \left( \frac{s}{n} \right)\right)\\
 & = \sum_{N < n \leq |s|} \frac{1}{n} \ERi \left( \frac{s}{n} \right) + O\left( \sum_{n \leq |s|} \frac{1}{n} \right) \\
 & = \sum_{N < n \leq |s|} \frac{1}{n} \ERi \left( \frac{s}{n} \right) + O(\log |s|) \\
  & =  \frac{1}{N+1} \ERi\left( \frac{s}{N+1} \right) +  \sum_{N+1 < n \leq |s|} \frac{1}{n} \ERi \left( \frac{s}{n} \right) + O(\log |s|)  \\
 & =  \frac{1}{N+1} \ERi\left( \frac{s}{N+1} \right) + O\left(  \frac{e^{s/(N+2)}}{s/(N+2)} |s|\right) + O(\log |s|) \\
 & \sim  \frac{1}{N+1} \ERi\left( \frac{s}{N+1} \right) \ (s \to \infty).
\end{align*}
The proposition follows.
\end{proof}

\begin{figure}[htb!]
\includegraphics[width=70mm]{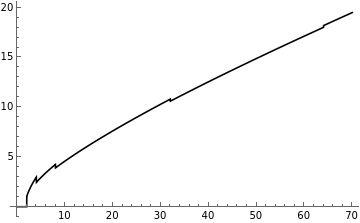}
\caption{\centering Graph of  $R(x)$ on $[0,70]$}
 \label{RR}
\end{figure}

\begin{figure}[htb!]
\includegraphics[width=70mm]{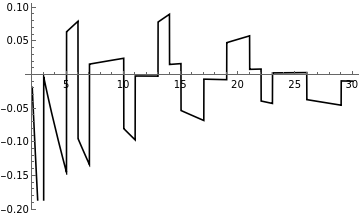}
\caption{\centering Graph of  $\Ri(2^x)- R(2^x)$ on $[2,30]$}
\label{RR1}
\end{figure}

Now, let $$R(x) = \operatorname{ER}(\log x,{\tfrac{1}{\log 2}}), \quad \forall x > 1,$$ that is,
\begin{align*}
R(x) = \sum_{n \leq \log_2 x} \frac{\mu(n)}{n}\li(x^{1/n}), \quad \forall x > 1.
\end{align*}
By Proposition \ref{mobius2}, it follows that
\begin{align}\label{RRRR}
\li(x) = \sum_{n \leq \log_2 x} \frac{1}{n} R(x^{1/n}), \quad \forall x \geq 2.
\end{align}
Note that the function $R(x)$ is piecewise continuous, with discontinuities precisely at each power $2^n$ of $2$ for which $n$ is squarefree, where the function leaps from the left to the right by $+\frac{\mu(n)}{n}\li(2)$.   Figure \ref{RR} provides a graph of the function $R(x)$ on $[0,70]$, and  Figure \ref{RR1} provides a graph of the function $\Ri(2^x)-R(2^x)$ on $[2,30]$.  As shown in Section 5.1, the functions $\Ri(x)$ and $R(x)$ are the ``main terms''  in two explicit formulas for the prime counting function.

\begin{corollary}\label{RRR}
Let  $\varepsilon > 0$, and let
\begin{align*} h(x) =   \left.
  \begin{cases}
    \frac{1}{\sqrt{x}} e^{(\log x)^{1/2+\varepsilon}} & \text{if the Riemann hypothesis holds} \\
     e^{-(\log x)^{3/5-\varepsilon}}   & \text{if the anti-Riemann hypothesis holds} \\
    x^{\Theta-1+\varepsilon} & \text{otherwise.} 
 \end{cases}
 \right.
\end{align*}
One has
$$\Ri(x) = R(x) + O(h(\log_2 x)) \ (x \to \infty).$$
In particular, one has
$$\Ri(x) = R(x) + O\left(\frac{1}{(\log \log x )^t} \right) \ (x \to \infty)$$
for all $t >0$. 
\end{corollary}

\begin{corollary}\label{RAprop}
One has the asymptotic expansion
\begin{align*}
\Ri(x) \simeq \sum_{n = 1}^\infty \frac{\mu(n)}{n} \li (x^{1/n}) \ (x \to \infty)
\end{align*}
and the (divergent) asymptotic expansion 
\begin{align*}
\li(x) \simeq \sum_{n = 1}^\infty \frac{1}{n} \Ri (x^{1/n}) \ (x \to \infty).
\end{align*}
\end{corollary}

\begin{corollary}\label{elll}
One has
\begin{align*} \li(x) \sim \Ri(x)  \ (x \to \infty),
\end{align*}
\begin{align*} \li(x)-\Ri(x) \sim \frac{\sqrt{x} }{\log x} \ (x \to \infty),
\end{align*}
and
\begin{align*}
\li(x)-\Ri(x) - \frac{1}{2}\li(\sqrt{x}) \sim \frac{1}{3}\li(\sqrt[3]{x}) \sim \frac{\sqrt[3]{x}}{\log x} \ (x \to \infty).
\end{align*}
\end{corollary}

\chapter{The analytic theory of primes}

In this chapter, we discuss some well-known deep  results, including the {\it Riemann--von Mangoldt explicit formulas}, relating the prime counting function $\pi(x)$ and the related Chebyshev functions $\psi(x)$ and $\vartheta(x)$ to the Riemann zeta function $\zeta(s)$, the exponential integral function $\Ei(s)$, and the  function $\ERi(s)$.   We use these results to deduce a strong version of  the prime number theorem with error bound.    We also discuss some well-known results, such as the {\it Riemann--von Mangoldt  formula},  concerning the zeros of the Riemann zeta function.    The final section of the chapter generalizes the prime number theorem to the setting of {\it abstract analytic number theory} \cite{kno}.

\section{The explicit formulas for $\pi_0(x)$, $\Pi_0(x)$, and $\psi_0(x)$}

Recall from Section 3.6 that the {\bf Riemann prime counting function} is defined by
$$\Pi(x) = \sum_{n = 1}^\infty \frac{1}{n}\pi(\sqrt[n]{x}) = \sum_{n = 1}^\infty \sum_{p^n \leq x} \frac{1}{n}, \quad \forall x \geq 0.$$
It is a weighted prime power counting function, where each power $p^n > 1$ of a prime $p$ is weighted by $\frac{1}{n}$.  
Let
$$\pi_0(x) = \lim_{\varepsilon \to 0} \frac{\pi(x+\varepsilon)+ \pi(x-\varepsilon)}{2}$$
and
$$\Pi_0(x) = \lim_{\varepsilon \to 0} \frac{\Pi(x+\varepsilon)+\Pi(x-\varepsilon)}{2} = \sum_{n = 1}^\infty \frac{1}{n}\pi_0(x^{1/n}).$$
\index[symbols]{.s GA@$\Pi_0(x)$}   By M\"obius inversion, one has
\begin{align}\label{piPi}
\pi_0(x) =\sum_{ n \leq \log_2 x} \frac{\mu(n)}{n}\Pi_0(x^{1/n}) =\sum_{n=1}^\infty \frac{\mu(n)}{n}\Pi_0(x^{1/n}).
\end{align}
The {\bf  Riemann--von Mangoldt explicit formula for $\Pi_0(x)$} \cite{mang0} \cite[Chapter 3]{edw} \cite[Chapter 4]{plyman}
\index{Riemann--von Mangoldt explicit formula for $\Pi_0(x)$}  states that
\begin{align}\label{RVM}
\Pi_0(x) = \li(x) - \sum_\rho \Ei(\rho \log x) +\log \xi(0), \quad \forall x > 1,
\end{align}
where the sum runs over all of the zeros $\rho$ of the Riemann zeta function $\zeta(s)$, with the nontrivial zeros taken in order of increasing absolute value of the imaginary part and repeated to multiplicity, and where $\xi(0) = \frac{1}{2}$ and thus $\log \xi(0) = -\log 2$.   Note  that $\li(x) = \Ei(\log x)$ for all $x > 0$, but that $ \Ei(\rho \log x)$ in general is not equal to $\li(x^\rho)$ as defined in Mathematica, nor is it a function of $x^\rho$, as noted in Remark \ref{lieiremark}.

The essence of Riemann's derivation of (\ref{RVM}) in \cite{rie}, which von Mangoldt made rigorous in 1895 \cite{mang0}, is the inversion of the expression (\ref{ingz2}) for $\log \zeta(s)$ in terms of $\Pi(x)$ and  the use of (\ref{ingz3})  and (\ref{ingz4}) to express $\log \zeta(s)$ in terms of the log and log-gamma functions and the zeros of $\zeta(s)$.

By (\ref{piPi}) and the  Riemann--von Mangoldt explicit formula for $\Pi_0(x)$, for all $x >1$ one has
\begin{align*}
\pi_0(x)  & =\sum_{n=1}^\infty \frac{\mu(n)}{n}\Pi_0(x^{1/n}) \\  
& = \sum_{n=1}^\infty  \frac{\mu(n)}{n} \left( \li(x^{1/n}) - \sum_\rho \Ei\left(\frac{\rho \log x}{n}\right) - \log 2  \right) \\
 & = \Ri(x) - \sum_{n=1}^\infty \sum_\rho \frac{\mu(n)}{n}\Ei\left(\frac{\rho \log x}{n}\right),
\end{align*}
 where the inner sum runs over all of the zeros $\rho$ of the Riemann zeta function $\zeta(s)$, with the nontrivial zeros taken in order of increasing absolute value of the imaginary part and repeated to multiplicity, and
where we have used the fact that $\sum_{n = 1}^\infty \frac{\mu(n)}{n} = 0$.  We call the resulting equation
\begin{align}\label{piexplicit}
\pi_0(x) = \Ri(x) - \sum_{n=1}^\infty \sum_\rho \frac{\mu(n)}{n}\Ei\left(\frac{\rho \log x}{n}\right), \quad \forall x>1,
\end{align}
the {\bf  Riemann--von Mangoldt explicit formula for $\pi_0(x)$}\index{Riemann--von Mangoldt explicit formula for $\pi_0(x)$}.   It follows from the explicit formulas noted above that the function $\li(x)$ is properly considered an approximation for $\Pi(x)$, while Riemann's function $\Ri(x)$ is the analogous approximation for $\pi(x)$.   

The sum $\sum_\rho \Ei(\rho \log x)$ in the explicit formula for $\Pi_0(x)$ over the trivial zeros $\rho = -2k$ of $\zeta(s)$ is given by
$$\sum_{k = 1}^\infty \li(x^{-2k}) =- \sum_{k = 1}^\infty \int_x^\infty {\frac {t^{-2k}}{t\log t}}\,dt =
 -\int_x^\infty {\frac {1}{t\log t}}\left(\sum_{k = 1}^\infty t^{-2k}\right)\,dt = - \int_x^\infty {\frac {dt }{t(t^{2}-1)\log t}}.$$
From the identity above, one can show  \cite{kul}  \cite[(32)]{ries1} that 
 \begin{align*}
 \sum_{n=1}^\infty  \frac{\mu(n)}{n}\left(  \sum_{k = 1}^\infty\li(x^{-2k/n})+\frac{n}{2\log x} \right)   =- \frac{1}{\pi} \arctan \frac{\pi}{\log x}.
 \end{align*} Thus, the explicit formulas can be re-expressed as
\begin{align}\label{realexplicit0}
\Pi_0(x) = \li(x) - \sum_\rho \Ei(\rho \log x) - \log 2+\int_x^\infty \frac{dt}{t(t^2-1) \log t}, \quad \forall x > 1,
\end{align}
and
\begin{align}\label{realexplicit}
\pi_0(x)   = \Ri(x) - \sum_{n=1}^\infty \frac{\mu(n)}{n}\left( \sum_\rho \Ei\left(\frac{\rho \log x}{n}\right) -\frac{n}{2\log x} \right) + \frac{1}{\pi} \arctan \frac{\pi}{\log x}, \quad \forall x > 1,
\end{align}
where now the sums run over all {\it nontrivial} zeros $\rho$ of $\zeta(s)$.   In \cite{kul} (which is an online source that is not peer-reviewed),  Kulsha claims that the identity
\begin{align}\label{realexplicit3}
\sum_{n = 1}^\infty  \frac{\mu(n)}{n}\left(\sum_\rho  \Ei\left(\frac{\rho \log x}{n}\right) -\frac{n}{2\log x} \right) =\sum_\rho \ERi(\rho \log x)+\frac{1}{\log x} 
\end{align}
``is obvious if we allow generalized summation $\sum_{n = 1}^\infty \mu(n) = \frac{1}{\zeta(0)} = -2$,'' but in private communication with the author he did not fill in any further details.   Given (\ref{realexplicit}), Kulsha's conjectural identity (\ref{realexplicit3}) is equivalent to the following.

\begin{conjecture}[{\cite{kul}}]\label{piexpconj}
One has
\begin{align*}
\pi_0(x) = \Ri(x) - \sum_\rho \ERi(\rho \log x)-\frac{1}{\log x} + \frac{1}{\pi} \arctan \frac{\pi}{\log x}, \quad \forall x > 1,
\end{align*}
where the sum runs over all of the nontrivial zeros $\rho$ of the Riemann zeta function $\zeta(s)$ taken in order of increasing absolute value of the imaginary part and repeated to multiplicity.  
\end{conjecture}

\begin{remark}[{On Conjecture \ref{piexpconj}}]\label{piexpconjrem}
Conjecture \ref{piexpconj}, or, rather, something close to it,  seems to  originate with Zagier's  claim in \cite[p.\ 14]{zag2} that the identity $\pi(x) = \Ri(x) - \sum_\rho \Ri(x^\rho)$ holds, with the sum apparently running over all zeros of $\zeta(s)$, and the nontrivial zeros taken in conjugate pairs and in nondecreasing order of imaginary part.   We interpret this claim more precisely as
\begin{align}\label{realexplicit4}
\pi_0(x) = \Ri(x) - \sum_\rho \ERi(\rho \log x), \quad \forall x > 1,
\end{align}
writing $\pi_0(x)$ instead of $\pi(x)$, and writing $\ERi(\rho \log x)$ instead of $\Ri(x^\rho)$ (which Mathematica interprets as $\ERi(\Log (x^\rho))$.   Zagier's claim is repeated in  \cite[p.\ 6]{hut} and at  https://en.wikipedia.org/wiki/Prime-counting\_function.
 In \cite[p.\ 224]{rib}, Ribenboim  claims that that the identity $\pi(x) = \Ri(x) - \sum_\rho \Ri(x^\rho)$ is ``exact,''  but with the sum running over  the {\it nontrivial} zeros of $\zeta(s)$.  MathWorld repeats this claim at https://mathworld.wolfram.com/RiemannPrimeCountingFunction.html,  but notes  that ``no proof of the equality $\ldots$ appears to exist in the literature.''  Regarding Ribenboim's claim and the more likely Conjecture \ref{piexpconj}, the author agrees with the comment in \cite[p.\ 249]{borw} of Borwein, Bradley, and Crandall that ``such a proof should be nontrivial, as the conditionally convergent series involved are problematic.''  Note, for example,  that, in \cite[Theorem, p.\ 248]{stop}, Stopple claims that
$$\pi(x) = \Ri(x)+ \sum_\rho \Ri(x^\rho) + \sum_{n = 1}^\infty \frac{\mu(n)}{n} \int_{x^{1/n}}^\infty \frac{dt}{t(t^2-1) \log t},$$
with the first sum running over the nontrivaial zeros of $\zeta(s)$.
However, the second sum is divergent (since $\sum_{n=1}^{\infty} \frac{\mu(n)}{n} ( - \frac{n}{2 \log x} + \int_{x^{1/n}}^{\infty} \frac{dt}{t (t^2-1) \log t}) = \frac{1}{\pi}\arctan\frac{\pi}{\log x}$),  the term $\sum_\rho \Ri(x^\rho)$  has a sign error,  and its convergence not been established.

Note that, since $$\Ri(x) - \sum_\rho \ERi(\rho \log x) =  \Ri(x) -  \sum_\rho \sum_{n=1}^\infty \frac{\mu(n)}{n}\Ei\left(\frac{\rho \log x}{n}\right),$$ whether the given series are convergent or not,  the conjecture (\ref{realexplicit4}) above is equivalent to the two sums in the explicit formula (\ref{piexplicit}) being interchangeable, which is tacitly assumed in the various claims discussed above.
Moreover, conjecture (\ref{realexplicit4}) requires convergence of the sum $\sum_\rho \ERi(\rho \log x)$ over both the nontrivial zeros and the trivial zeros of $\zeta(s)$ separately.  We have already remarked that convergence  over the nontrivial zeros is not obvious.  It has also been claimed, in  \cite[p.\ 6]{hut} and at  https://en.wikipedia.org/wiki/Prime-counting\_function, that the sum $\sum_\rho \ERi(\rho \log x) =  \sum_{k = 1}^\infty  \Ri(x^{-2k})$  over the trivial zeros $\rho  = -2k$ converges and is given by
\begin{align}\label{realexplicit5}
\sum_{k = 1}^\infty  \Ri(x^{-2k}) = \frac{1}{\log x} - \frac{1}{\pi} \arctan \frac{\pi}{\log x}, \quad \forall x>1.
\end{align}
However, the behavior of $\Ri(x)$ immediately to the right of $x = 0$ is alarmingly complex \cite{born}, so much so that it is also not obvious even that the sum $\sum_{k = 1}^\infty  \Ri(x^{-2k})$ converges for all $x > 1$.  Of course, if  Conjecture \ref{piexpconj} holds, then  the conjectures  (\ref{realexplicit4}) and  (\ref{realexplicit5}) are equivalent.  Unfortunately, the author  does not see how to establish any of these claims.  Given the various missteps and errors of omission in the literature,  any purported proofs of these conjectures, including Kulsha's (\ref{realexplicit3}), should be held to high standards of rigor.
\end{remark}

\begin{remark}[Alternative explicit formulas for $\pi_0(x)$]
An alternative explicit formula for $\pi_0(x)$ is
\begin{align*}
\pi_0(x)  & =\sum_{n \leq \log_2 x} \frac{\mu(n)}{n}\Pi_0(x^{1/n}) \\  
& = \sum_{n\leq \log_2 x}  \frac{\mu(n)}{n} \left( \li(x^{1/n}) - \sum_\rho \Ei\left(\frac{\rho \log x}{n}\right) + \int_{x^{1/n}}^\infty \frac{dt}{t(t^2-1) \log t} - \log 2  \right) \\
& =R(x) - \sum_{n\leq \log_2 x} \frac{\mu(n)}{n}  \left(\sum_\rho   \Ei\left(\frac{\rho \log x}{n}\right) - \int_{x^{1/n}}^\infty \frac{dt}{t(t^2-1) \log t}+ \log 2  \right)
\end{align*}
for all $x > 1$, where $\rho$ ranges over all nontrivial zeros of $\zeta(s)$, and where the piecewise smooth function
\begin{align*}
R(x) = \sum_{n \leq \log_2 x} \frac{\mu(n)}{n}\li(x^{1/n}), \quad \forall x > 1,
\end{align*}
is related to $\li(x)$ and $\Ri(x)$ as in (\ref{RRRR}) and Corollary \ref{RRR}, respectively.  By similar reasoning, for any positive integer $N$ one has
\begin{align*}
\pi_0(x)  = \sum_{n = 1}^N  \frac{\mu(n)}{n} \left( \li(x^{1/n}) - \sum_\rho \Ei\left(\frac{\rho \log x}{n}\right) + \int_{x^{1/n}}^\infty \frac{dt}{t(t^2-1) \log t} - \log 2  \right)  
\end{align*}
for all $x \in  (1,2^{N+1})$ (which is equivalent to \cite[(17)]{ries1}), where  $\rho$ ranges over all nontrivial zeros of $\zeta(s)$.
\end{remark}

\begin{figure}[ht!]
\includegraphics[width=55mm]{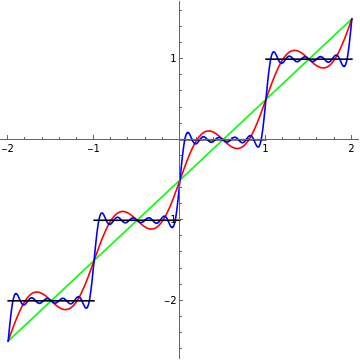}
\caption{\centering Graph of  $\lfloor x\rfloor$ (in black), $x-\frac{1}{2}$ (in green), and $x - \frac{1}{2}+\sum_{n = 1}^N \frac{\sin 2 \pi n x}{\pi n}$  for $N = 1$ (in red) and $N = 6$ (in blue) on $[-2,2]$}
\label{floor}
\end{figure}

\begin{figure}[ht!]
\includegraphics[width=70mm]{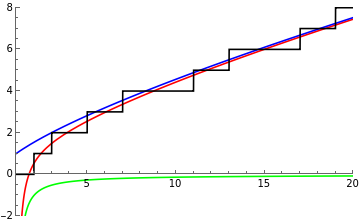}
\caption{\centering Graph of  $\pi(x)$ (in black), $\Ri(x)$ (in blue), $\Ri(x) - \frac{1}{\log x} +\frac{1}{\pi} \arctan \frac{\pi}{\log x}$ (in red), and $- \frac{1}{\log x} +\frac{1}{\pi} \arctan \frac{\pi}{\log x} $ (in green) on $[1,20]$}
\label{correction2}
\end{figure}

The explicit formula for $\pi_0(x)$ is analogous to a {\it Fourier series expansion} of a piecewise continuous and bounded periodic function.   For the sake of analogy, we discuss the floor function
$$\lfloor x \rfloor = x+ \{x \},$$
where the fractional part
$$\{x\} = x-\lfloor x \rfloor \in [0,1)$$ of $x$ is a periodic function of period $1$. The main term of the expression above for $\lfloor x \rfloor$ is the smooth function $x$, which approximates the ``size'' of the function $\lfloor x \rfloor$ in the sense that
$$\lfloor x \rfloor \sim x \ (x \to \infty),$$ 
thus analogous to how $\Ri(x)$ approximates $\pi(x)$.
Analogous to the function $\Ri(x)-\pi_0(x)$, the subordinate term $ \{x \}$, or, rather, the function
$$\{x\}_0 = \lim_{\varepsilon \to 0} \frac{\{x+\varepsilon\}+ \{x-\varepsilon\}}{2} = x-\lfloor x \rfloor_0$$
that assumes the value $\{x\}$ at the non-integers and $\frac{1}{2}$ at the integers, where 
$$\lfloor x \rfloor_0  = \lim_{\varepsilon \to 0} \frac{\lfloor x+\varepsilon \rfloor+ \lfloor x-\varepsilon \rfloor}{2},$$
has the Fourier series expansion
$$\{x\}_0 = \frac{1}{2}-\sum_{n = 1}^\infty \frac{\sin 2 \pi n x }{\pi n} $$
valid for all $x \in \RR$. (For example, one has $\frac{1}{4}  = \frac{1}{2}-\sum_{n = 1}^\infty \frac{\sin\frac{\pi n}{2}}{\pi n} $ and therefore $\frac{\pi}{4} = \sum_{n = 1}^\infty \frac{\sin \frac{\pi n}{2}}{n}  = 1-\frac{1}{3}+\frac{1}{5}-\frac{1}{7}+\cdots.$) 
Thus, one has the expansion
$$\lfloor x \rfloor_0 = x-\{x\}_0 = x  + \lfloor 0 \rfloor_0 +\sum_{n = 1}^\infty \frac{ \sin 2 \pi n x}{\pi n}$$
of the function $\lfloor x \rfloor_0$, with smooth main term $x + \lfloor 0 \rfloor_0  = x - \frac{1}{2}$ and smooth oscillatory terms $\frac{\sin 2\pi n x}{ \pi n}$, one for every positive root $n$ of the function $\{x\}$, or  one for every positive root $n$ of the function $\sin \pi x$.
 For each positive integer $N$, let $${\lfloor x\rfloor}_N = x - \frac{1}{2}+\sum_{n = 1}^N \frac{\sin 2 \pi n x}{\pi n},$$ 
so that
$$\lfloor x \rfloor_0 = \lim_{N \to \infty} {\lfloor x\rfloor}_N.$$   On each interval $[n,n+1]$ for $n \in \ZZ$, the oscillating function ${\lfloor x\rfloor_N}$ has exactly $N$ local maxima and $N$ local minima.   See Figure \ref{floor} for the graph of  the functions $\lfloor x\rfloor$, $x-\frac{1}{2}$,  and ${\lfloor x\rfloor}_N$ for $N = 1$ and $N = 5$, on the interval $[-2,2]$. 
Notice that the function $x-\frac{1}{2}$, which is the function $x$ but with the slight ``correction term'' $\lfloor 0 \rfloor_0  = -\frac{1}{2}$, traces the ``center'' of the function $\lfloor x\rfloor_0$.  This is analogous to how the function $\Ri(x)$, or, slightly more accurately for small $x$ in particular, the function
$$\Ri(x)- \frac{1}{\log x} +\frac{1}{\pi} \arctan \frac{\pi}{\log x},$$
 traces the ``center'' of the function $\pi_0(x)$.  The approximating function above is the function $\Ri(x)$ but with  the negative ``correction term''
\begin{align*}
 - \frac{1}{\log x} +\frac{1}{\pi} \arctan \frac{\pi}{\log x} 
& = \sum_{n = 1}^\infty \frac{(-1)^{n}\pi^{2n}}{(2n+1)(\log x)^{2n+1}} \\
& = - \frac{\pi^2}{3(\log x)^3} +  \frac{\pi^4}{5(\log x)^5}(1+o(1)) \ (x \to \infty),
\end{align*}
which increases to $0$ as $x \to \infty$, and where the given Taylor series in $\frac{1}{\log x}$ converges for $x > e$.  As we have noted, the correction term  comes, conjecturally, from the contribution of the trivial zeros of $\zeta(s)$.  The  term is negligible if $x$ isn't too small, e.g., its value at $x = 20$ is just $-0.0762\ldots$.  This is exhibited in Figure \ref{correction2}, which provides a graph of the functions $\pi(x)$, $\Ri(x)$, $\Ri(x) - \frac{1}{\log x} +\frac{1}{\pi} \arctan \frac{\pi}{\log x}$, and $- \frac{1}{\log x} +\frac{1}{\pi} \arctan \frac{\pi}{\log x} $ on the interval $[1,20]$. 

\begin{figure}[ht!]
\includegraphics[width=70mm]{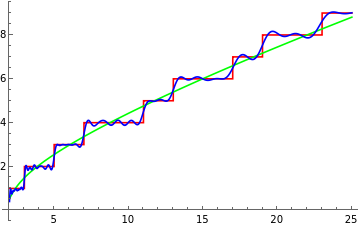} 
\caption{\centering Graph of $\pi(x)$ (in red), $\Ri(x)- \frac{1}{\log x} +\frac{1}{\pi} \arctan \frac{\pi}{\log x}$ (in green),  and $\Ri(x)- \frac{1}{\log x} +\frac{1}{\pi} \arctan \frac{\pi}{\log x}-\sum_{n = 1}^5 \sum_{k = 1}^{20} \frac{\mu(n)}{n} 2\operatorname{Re} \Ei\left( \frac{\rho_k\log x}{n} \right)$ (in blue)  on $[2,25]$ }
 \label{20zeros}
\end{figure}

\begin{figure}[ht!]
\includegraphics[width=70mm]{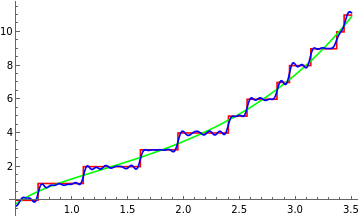} 
\caption{\centering Graph of $\pi(e^x)$ (in red), $\Ri(e^x)- \frac{1}{x} +\frac{1}{\pi} \arctan \frac{\pi}{x}$ (in green),  and $\Ri(e^x)- \frac{1}{x} +\frac{1}{\pi} \arctan \frac{\pi}{x}-\sum_{n = 1}^5 \sum_{k = 1}^{20}  \frac{\mu(n)}{n}2\operatorname{Re} \Ei\left( \frac{\rho_k x}{n} \right)$ (in blue)  on $[0.5,3.5]$ }
 \label{20zeros2}
\end{figure}

\begin{figure}[ht!]
\includegraphics[width=70mm]{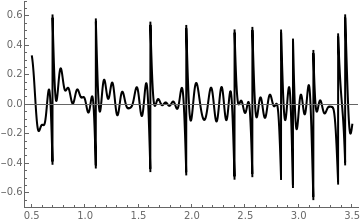} 
\caption{\centering Graph of $\pi(e^x)-\left(\Ri(e^x)- \frac{1}{x} +\frac{1}{\pi} \arctan \frac{\pi}{x}- \sum_{n = 1}^5 \sum_{k = 1}^{20} \frac{\mu(n)}{n}2\operatorname{Re} \Ei\left( \frac{\rho_k x}{n} \right)\right)$  on $[0.5,3.5]$ }
 \label{20zerosa}
\end{figure}

\begin{figure}[ht!]
\includegraphics[width=70mm]{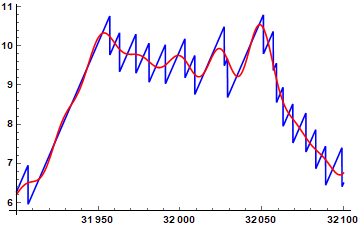} 
\caption{\centering Graph of $-\pi(x)+\Ri(x)- \frac{1}{\log x} +\frac{1}{\pi} \arctan \frac{\pi}{\log x}$ (in blue) and $ \sum_{n = 1}^{10} \sum_{k = 1}^{10000}\frac{\mu(n)}{n} 2\operatorname{Re}\Ei\left( \frac{\rho_k\log x}{n} \right)$ (in red) on $[31900,32100]$ }
 \label{10000zeros}
\end{figure}

Note that $$\Ei\left( \frac{\rho\log x}{n}\right)+ \Ei\left( \frac{\overline{\rho}\log x}{n}\right) = 2\operatorname{Re}\Ei\left( \frac{\rho\log x}{n} \right)$$
for all $\rho \in \CC$ and all $x > 0$.  Let $\rho_k$ denote the $k$th nontrivial zero of $\zeta(s)$ with positive imaginary part, repeated to multiplicity, and listed in order of absolute value.   The explicit formula for $\pi_0(x)$ then yields the approximations
\begin{align*}
\pi_0(x) \approx \Ri(x) - \frac{1}{\log x} +\frac{1}{\pi} \arctan \frac{\pi}{\log x}-\sum_{n = 1}^M  \sum_{k = 1}^{N} \frac{\mu(n)}{n}2\operatorname{Re} \Ei\left( \frac{\rho_k\log x}{n} \right),
\end{align*}
for positive integers $M$ and $N$, which are to be chosen  suitably large  for $x$ in a given range.  Note that Conjecture \ref{piexpconj} implies equality for $M = N = \infty$, but with the two sums interchanged.  Figure \ref{20zeros} provides a graph of the function $\pi(x)$ and the approximations $$\Ri(x)- \frac{1}{\log x} +\frac{1}{\pi} \arctan \frac{\pi}{\log x}$$ and $$\Ri(x)- \frac{1}{\log x} +\frac{1}{\pi} \arctan \frac{\pi}{\log x}- \sum_{n = 1}^5 \sum_{k = 1}^{20} \frac{\mu(n)}{n} 2\operatorname{Re}\Ei\left( \frac{\rho_k\log x}{n} \right)$$ of $\pi(x)$.  Figure \ref{20zeros2} provides the same graph as the three functions in Figure  \ref{20zeros} but on a lin-log scale, and Figure \ref{20zerosa} provides a graph of the error in the latter approximation of $\pi_0(x)$.  Figure \ref{10000zeros}, which is courtesy of D.\ Stoll, and which took 180 hours to compute, compares the graph $-\pi(x)+\Ri(x)- \frac{1}{\log x} +\frac{1}{\pi} \arctan \frac{\pi}{\log x}$ and $\sum_{n = 1}^{10} \sum_{k = 1}^{10000}  \frac{\mu(n)}{n} \Ei\left( \frac{\rho_k\log x}{n} \right)$ on $[31900,32100]$.

We mention the analogy with the explicit formula for the floor function not only because it aids in understanding the Riemann--von Mangoldt explicit formula, but also because there are important connections between the Riemann zeta function and the floor function, both of which can be seen as functions encoding the sequence of positive integers.  One connection between the Riemann zeta function and the floor function is the identity 
\begin{align}\label{zetafloor}
\zeta(s) = s \int_1^\infty  \frac{\lfloor x \rfloor}{x^{s+1}} \, dx,
\end{align}
which holds for all  $s$ with $\operatorname{Re} s > 1$.  Since
$$\int_1^\infty \frac{x}{x^{s+1}} \, dx = \int_1^\infty \frac{1}{x^{s}}\, dx  = \frac{1}{s-1}$$
for all such $s$, one has
$$\zeta(s) - \frac{1}{s-1} = 1- s \int_1^\infty \frac{\{x\}}{x^{s+1}} \, dx.$$
The identity above is more useful than (\ref{zetafloor})  because it is valid for all $s \in \CC$ with $\operatorname{Re} s > 0$, even for $s = 1$ if one takes the limit as $s \to 1$ on both sides of the equation, so that
$$\gamma = 1-  \int_1^\infty \frac{\{x\}}{x^{2}} \, dx = \int_1^\infty \frac{1-\{x\}}{x^{2}} \, dx.$$
Thus, the identity allows one to {\it define} $\zeta(s)$ for all $s \in \CC$ with $\operatorname{Re} s > 0$  by the equation
$$\zeta(s) =  \frac{s}{s-1}- s \int_1^\infty \frac{\{x\}}{x^{s+1}} \, dx.$$
It follows from this definition that the nontrival zeros of $\zeta(s)$ (none of which have real part $0$ or $1$) are precisely those complex numbers $\rho$ with $\operatorname{Re} \rho > 0$ (and $\operatorname{Re} \rho < 1$) such that
$$(\rho-1)\int_1^\infty \frac{\{x\}}{x^{\rho+1}} \, dx = 1.$$

Now,  let
$$\psi_0(x) = \lim_{\varepsilon \to 0} \frac{\psi(x+\varepsilon)+ \psi(x-\varepsilon)}{2}.$$  The
{\bf  Riemann--von Mangoldt explicit formula for $\psi_0(x)$}\index[symbols]{.s DA@$\psi_0(x)$}\index{Riemann--von Mangoldt explicit formula for $\psi_0(x)$} \cite{mang0} \cite[Chapter 3]{edw} \cite[Chapter IV]{ing2} \cite[Chapter 10]{over}  \cite[Chapter 4]{plyman} states that
\begin{align}\label{RVM2}
\psi_{0}(x)=x-\sum_{\rho }{\frac {x^{\rho }}{\rho }}-{\frac {\zeta '(0)}{\zeta (0)}}, \quad \forall x > 1,
\end{align}
where the sum runs over all of the zeros $\rho$ of $\zeta(s)$, with the nontrivial zeros taken in order of increasing absolute value of the imaginary part and repeated to multiplicity, and where
$${\frac {\zeta '(0)}{\zeta (0)}} = \log 2 \pi.$$  Since
$$\sum_{k=1}^{\infty }{\frac {x^{-2k}}{-2k}}={\frac {1}{2}}\log \left(1-x^{-2}\right)$$
is the corresponding sum over the trivial zeros of $\zeta(s)$, the explicit formulas may be re-expressed as
$$\psi_{0}(x)=x-\sum_{\rho }{\frac {x^{\rho }}{\rho }}-{\frac {\zeta '(0)}{\zeta (0)}}-{\frac {1}{2}}\log(1-x^{-2}),$$
where now the sum is over all of the nontrivial zeros $\rho$ of $\zeta(s)$.

Von Mangoldt's derivation of (\ref{RVM}) in \cite{mang0} relies on the  inversion of the expression (\ref{ingz}) for $\frac{\zeta'(s)}{\zeta(s)}$ in terms of $\psi(x)$. Since $\frac{\zeta'(s)}{\zeta(s)}$ is meromorphic, while $\log \zeta(s)$ has logarithmic singularities, the explicit formula for $\psi_0(x)$ is more readily verified than that for $\Pi_0(x)$.  Moreover, the explicit formula for $\psi_0(x)$ is much simpler than that for the functions $\pi_0(x)$ and $\Pi_0(x)$, making it is easier to handle, e.g., in bounding the sum of oscillatory terms involving the zeros of $\zeta(s)$.  For these reasons, it is typical for analytic number theorists to first prove a given result for $\psi(x)$ (such as the prime number theorem in the form $\psi(x) \sim x \ (x \to \infty)$)  and then deduce a corresponding result for $\pi(x)$ by known relationships between the two functions, as described in  Proposition \ref{pitheta} and Example \ref{direxample}, for example.

\begin{remark}[Explict formula for $\vartheta_0(x)$]
Since by M\"obius inversion one has
$$\vartheta_0(x) =  \sum_{1 \leq n \leq  \log_2 x} \mu(n) \psi_0(x^{1/n}) = \sum_{n = 1}^\infty \mu(n) \psi_0(x^{1/n}),$$
where
$$\vartheta_0(x) = \lim_{\varepsilon \to 0} \frac{\vartheta(x+\varepsilon)+ \vartheta(x-\varepsilon)}{2},$$
one has
$$\vartheta_{0}(x)=  \sum_{1 \leq n \leq  \log_2 x} \mu(n) x^{1/n} -  \sum_{1 \leq n \leq  \log_2 x} \mu(n) \sum_{\rho }  {\frac {x^{\rho/n}}{\rho }}-{\frac {\zeta '(0)}{\zeta (0)}} M(\log_2 x), \quad \forall x > 1,$$
where $M(x) = \sum_{n \leq x} \mu(n)$ is the Mertens function.
\end{remark}

\begin{figure}[htb!]
\includegraphics[width=70mm]{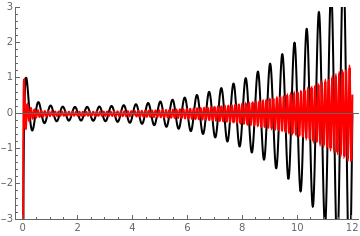}
\caption{\centering Graph of $\Ei(\rho_k x) + \Ei(\overline{\rho_k} x) = 2 \operatorname{Re} \Ei(\rho_k x)$ on $[0,12]$ for $k = 1$ (in black) and $k = 10$ (in red)}
 \label{lirho}
\end{figure}

\begin{figure}[htb!]
\includegraphics[width=70mm]{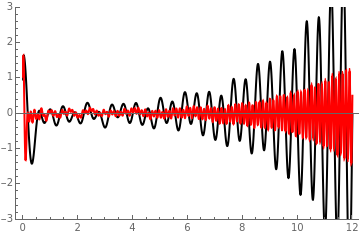}
\caption{\centering Graph of $\ERi(\rho_k x) + \ERi(\overline{\rho_k} x) = 2 \operatorname{Re} \ERi(\rho_k x)$ on $[0,12]$ for $k = 1$ (in black) and $k = 10$ (in red)}
 \label{rirho1}
\end{figure}

\begin{figure}[htb!]
\includegraphics[width=70mm]{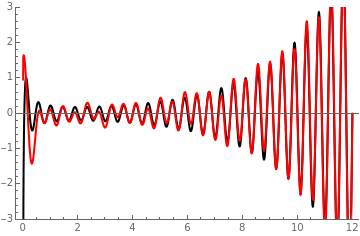}
\caption{\centering Graph of  $\Ei(\rho_1 x) + \Ei(\overline{\rho_1} x) = 2 \operatorname{Re} \Ei(\rho_1 x)$ (in black) and $\ERi(\rho_1 x) + \ERi(\overline{\rho_1} x) = 2 \operatorname{Re} \ERi(\rho_1 x)$ (in red) on $[0,12]$ 
\label{liri4}}
\end{figure}

\begin{figure}[htb!]
\includegraphics[width=70mm]{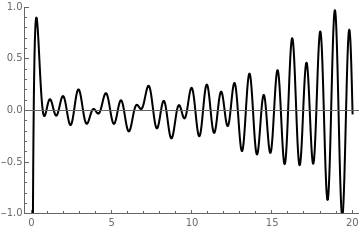}
\caption{\centering Graph of $(\Ei(\rho_1 x) + \Ei(\overline{\rho_1} x)) -(\ERi(\rho_1 x) + \ERi(\overline{\rho_1} x))$ on $[0,20]$}
 \label{liminusri}
\end{figure}

In exponential form, the explicit formulas  are
$$\Pi_0(e^x) = \li(e^x) - \sum_\rho \Ei(\rho x) +\log \xi(0), \quad \forall x > 0,$$
$$\pi_0(e^x)  = \Ri(e^x) - \sum_{n=1}^\infty \sum_\rho \frac{\mu(n)}{n}\Ei\left(\frac{\rho x}{n}\right), \quad \forall x>0,$$
$$\psi_{0}(e^x)=e^x-\sum_{\rho }{\frac {e^{\rho x }}{\rho }}-{\frac {\zeta '(0)}{\zeta (0)}}, \quad \forall x> 0.$$
These forms are ``better'' than the linear forms because their terms exhibit more easily discernible behavior, e.g., the oscillations are linear in $x$ rather than logarithmic in $x$.   Figures \ref{lirho} and \ref{rirho1} provide graphs of the function
$$\Ei(\rho x) + \Ei(\overline{\rho} x) = 2 \operatorname{Re} \Ei(\rho x)$$
and
$$\ERi(\rho x) + \ERi(\overline{\rho} x) = 2 \operatorname{Re} \ERi(\rho x),$$ 
respectively, on the interval $[0,12]$ for the first nontrivial zero $\rho = \rho_1 = \tfrac{1}{2} + 14.134725\ldots i$ and the tenth nontrivial zero $\rho = \rho_{10} = \frac{1}{2} + 49.773832\ldots i$ of $\zeta(s)$, where $\rho_k$ for any positive integer $k$ denotes the $k$th nontrivial zero of $\zeta(s)$ with positive imaginary part, repeated to multiplicity.   Figure \ref{liri4} compares the two graphs for $\rho = \rho_1$, while Figure \ref{liminusri} provides a graph of their difference on $[0,20]$.

\begin{figure}[htb!]
\includegraphics[width=70mm]{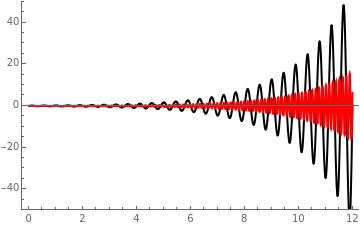}
\caption{\centering Graph of $\frac{e^{\rho_k x}}{\rho_k} + \frac{e^{\overline{\rho_k }x}}{\overline{\rho_k}}  = 2\operatorname{Re}\frac{e^{\rho_k x}}{\rho_k}$ on $[0,12]$ for $k = 1$ (in black) and $k = 10$ (in red) }
\label{Psi0}
\end{figure}

For all $x  \in \RR$  and all  $\rho \in \CC$ with nonzero real and imaginary part, one has
\begin{align*}
\frac{e^{\rho x}}{\rho} + \frac{e^{\overline{\rho }x}}{\overline{\rho}}  & = 2\operatorname{Re}\frac{e^{\rho x}}{\rho} \\
& = \frac{2e^{x\operatorname{Re} \rho}}{|\rho|^2}\left((\operatorname{Re} \rho) \cos (x \operatorname{Im} \rho) + (\operatorname{Im} \rho) \sin (x \operatorname{Im} \rho) \right) \\
& = \frac{2e^{x\operatorname{Re} \rho}}{|\rho|} \cos (x \operatorname{Im} \rho - \arctan(\operatorname{Im} \rho /\operatorname{Re} \rho )).
\end{align*}
Therefore, the function $\frac{e^{\rho x}}{\rho} + \frac{e^{\overline{\rho }x}}{\overline{\rho}}$ is just the complex exponential function $e^{x\operatorname{Re} \rho}$ multiplied by the cosine wave with amplitude $\frac{2}{|\rho|}$,  period $\frac{2\pi}{ \operatorname{Im} \rho}$, and phase shift $- \arctan(\operatorname{Im} \rho /\operatorname{Re} \rho )$.   Figure \ref{Psi0} provides a graph of the function on the interval $[0,12]$ for both $\rho = \rho_1$ and $\rho = \rho_{10}$. 

By the remarks above, one has
\begin{align}
\sum_{\rho }{\frac {e^{\rho x}}{\rho }} & = \sum_{\rho: \, \operatorname{Im}  \rho> 0} \frac{2e^{x\operatorname{Re} \rho}}{|\rho|} \cos (x \operatorname{Im} \rho - \arctan(\operatorname{Im} \rho /\operatorname{Re} \rho ))  \nonumber \\
& =2e^{\Theta x} \sum_{\rho: \, \operatorname{Im}  \rho> 0} \frac{e^{-(\Theta-\operatorname{Re} \rho)x}}{|\rho|} \cos (x \operatorname{Im} \rho - \arctan(\operatorname{Im} \rho /\operatorname{Re} \rho )).\label{psieq}
\end{align}
Moreover, the Riemann hypothesis is of course equivalent to $\Theta = \operatorname{Re} \rho$, hence to $e^{-(\Theta-\operatorname{Re} \rho)x} = 1$,  for all $\rho$.

\begin{remark}[Infinite sums involving Riemann zeta function zeros]
Riemann proved that the value $\sum_\rho \frac{1}{\rho}$ of the function $\sum_\rho \frac{x^\rho}{\rho}$ at $x = 1$, where the sums are over the nontrivial zeros of $\zeta(s)$, is conditionally convergent to the value
$$\sum_\rho \frac{1}{\rho}  = \sum_{\rho:\, \operatorname{Re} \rho >0}  \frac{2 \operatorname{Re \rho}}{|\rho|^2}  = \tfrac{1}{2}\gamma+1-\tfrac{1}{2} \log 4 \pi = 0.023095708966\ldots.$$  Thus, from the explicit formula for $\psi_0(x)$, one has
$$\lim_{x \to 1^+} \sum_\rho \frac{x^\rho}{\rho}   = \infty \neq \sum_\rho \frac{1}{\rho}.$$
Moroeover,  by the functional equation for $\zeta(s)$, one has
$$\sum_\rho \frac{1}{\rho(1-\rho)} = \sum_\rho \left (\frac{1}{\rho} +\frac{1}{1-\rho}\right)  =2 \sum_\rho \frac{1}{\rho} = \gamma+2- \log 4 \pi =  0.046191417932\ldots,$$
and, moreover, 
$$\sum_\rho \frac{ 2\operatorname{Re} \rho}{|\rho|^2}  = \sum_\rho \left( \frac{1}{\rho}+\frac{1}{\overline{\rho}} \right) = 2\sum_\rho \frac{1}{\rho} = \gamma+2- \log 4 \pi,$$
or, equivalently, 
$$\sum_{\rho}  \frac{\operatorname{Re} \rho}{|\rho|^2} = \sum_\rho \frac{1}{\rho} = \tfrac{1}{2}\gamma+1-\tfrac{1}{2} \log 4 \pi.$$
It follows that the Riemann hypothesis is equivalent to
$$\sum_\rho \frac{1}{|\rho|^2}   = \gamma+2- \log 4 \pi.$$
Note, however, that the sums $\sum_\rho \frac{1}{|\rho|}$  and $\sum_\rho \frac{1}{\operatorname{Im} \rho}$ both diverge.
\end{remark}

\begin{remark}[Approximations of the explicit formulas]
By Proposition \ref{EIlem}, for any fixed $\rho \in \CC$, one has
\begin{align*} \Ei(\rho \log x)  \sim \frac{x^\rho}{\rho \log x}  \ (x \to \infty),
\end{align*} 
and, more generally, the asymptotic expansion
$$ \Ei (\rho \log x) \simeq  \sum_{n = 0}^\infty \frac{n!\, x^{\rho} }{(\rho \log x)^{n+1}}  \ (x \to \infty).$$
Likewise, for any fixed $x >1$, one has
$$\Ei (\rho_k \log x) \sim \frac{x^{\rho_k}}{\rho_k\log x}  \ (k \to \infty),$$
and, more generally, the asymptotic expansion
$$ \Ei (\rho_k \log x) \simeq  \sum_{n = 0}^\infty \frac{n!\, x^{\rho_k} }{(\rho_k \log x)^{n+1}}  \ (k \to \infty).$$
For example, for $x = e$ and $k = 1$, one has
$$\Ei(\rho_1) = 0.115940896079\ldots+ 0.004272165141\ldots i$$
and
$$\frac{e^{\rho_1}}{\rho_1} = 0.116507260809\ldots + 0.003836493969\ldots i.$$
Using three terms of the asymptotic expansion of $\Ei (\rho_k \log x)$ above, one obtains the better approximation
$$\frac{e^{\rho_1}}{\rho_1} + \frac{e^{\rho_1}}{\rho_1^2} + \frac{2e^{\rho_1}}{\rho_1^3}  = 0.115910338967\ldots +
0.004506810004\ldots i.$$
The approximations above are useful for approximating the terms in the explicit formulas for  $\Pi_0(x)$ and $\pi_0(x)$ \cite{stoll} and ascertaining their qualitative behavior \cite[pp.\ 970--972]{ries1}.
\end{remark}

\section{The prime number theorem with error bound}

It is known \cite[Theorem 12.5]{mont} that the Riemann--von Mangoldt explicit formula for $\psi_0(x)$ generalizes to the truncated form
\begin{align}\label{psiest}
\psi_{0}(x)=x-\sum_{\rho: \, |\operatorname{Im} \rho| \leq  T} \frac{x^\rho}{\rho} -{\frac {\zeta '(0)}{\zeta (0)}}-{\frac {1}{2}}\log(1-x^{-2}) + O \left( \min \left(1, \frac{x}{\langle x \rangle T} \right) \log x + \frac{x (\log xT)^2}{T}\right)
\end{align}
uniformly for $x \geq c$ and $T \geq 2$ (so the $O$ constant does not depend on $x$ or $T$), where $c > 1$,  and where $\langle x \rangle$ denotes the distance from $x$ to the nearest prime power not equal to $x$.  Taking a limit as $T \to \infty$ yields the explicit formula for $\psi_0(x)$.  Moreover, although the series $\sum_\rho \frac{1}{|\rho|}$ diverges, it is known \cite[Theorem 25b]{ing2} \cite[(13.1)]{mont} that
\begin{align}\label{rhoest}
\sum_{\rho: \, |\operatorname{Im}\rho| \leq T} \frac{1}{|\rho|} = O((\log T)^2) \ (T \to \infty).
\end{align}
(This is an easy consequence of the {\it Riemann--von Mangoldt formula} (\ref{RVMF}) of the next section.) 
From (\ref{psiest}) and (\ref{rhoest}), one deduces the following.

\begin{theorem}[{\cite[Theorems 30 and 31]{ing2} \cite[Theorem 12.3]{ivic}}]\label{psiestprop}
The Riemann constant $\Theta$ is the smallest real number $t$ such that 
$$\psi(x) = x + O(x^{t} (\log x)^2) \ (x \to \infty).$$
\end{theorem}

\begin{proof}
Taking $T = x$ in (\ref{psiest}) and applying (\ref{rhoest}), we obtain
\begin{align*}
x-\psi(x) & =\sum_{\rho:\, |\operatorname{Im} \rho| \leq  x} \frac{x^\rho}{\rho}+  O (\log x)^2) \ (x \to \infty) \\
  & =   x^{\Theta} \sum_{\rho:\, |\operatorname{Im} \rho| \leq  x} \frac{x^{-(\Theta-\rho)}}{\rho}+  O (\log x)^2) \ (x \to \infty) \\
 & =O\left( x^{\Theta} \sum_{\rho: \, |\operatorname{Im}\rho| \leq x} \frac{1}{|\rho|}  \right) + O (\log x)^2) \ (x \to \infty) \\
 &  = O(x^{\Theta} (\log x)^2) \ (x \to \infty).
\end{align*}
On the other hand,  suppose that $x-\psi(x) = O(x^{t} (\log x)^2) \ (x \to \infty)$, where $t \in \RR$.  By Example \ref{direxample}(6), one has
$$\frac{\zeta'(s)}{\zeta(s)}+ \frac{s}{s-1} =  s \int_1^\infty\frac{x-\psi(x)}{x^{s+1}} \, dx$$
for all $s$ with $\operatorname{Re} s > 1$.  By the hypothesized $O$ bound on $x-\psi(x)$, the integral above  is analytic  on $\{s \in \CC: \operatorname{Re} s > t\}$, and thus so is the function $\frac{\zeta'(s)}{\zeta(s)}+ \frac{s}{s-1}$.   This implies that $\zeta(s)$ has no zeros in $\{s \in \CC: \operatorname{Re} s > t\}$, hence that $t \geq \Theta$.     This completes the proof.
\end{proof}

\begin{corollary}\label{montlem}
The Riemann constant $\Theta$ is the smallest real number $t$ such that 
$$\vartheta(x) = x + O(x^t(\log x)^2) \ (x \to \infty).$$
Likewise,  $\Theta$ is the smallest real number $t$ such that 
$$\pi(x) = \li(x) + O(x^t\log x) \ (x \to \infty).$$
\end{corollary}

\begin{proof}
From the theorem,  the statement for $\vartheta(x)$ follows from  Proposition  \ref{pitheta2b}(2),  and then the statement for $\pi(x)$ follows from Proposition \ref{pitheta} and Karamata's integral theorem.
\end{proof}

From the theorem and corollary above,  one also deduces the following.

\begin{corollary}\label{maindegree}
One has
$$\Theta = \deg(\li-\pi) = \deg(\id -\psi) = \deg(\id-\vartheta).$$
\end{corollary}

\begin{proof}
By Theorem \ref{psiestprop}, one has $\psi(x) = x + O(x^{t}) \ (x \to \infty)$ for all $t > \Theta$ but for no $t < \Theta$.   It follows that $\deg(\id-\psi) = \Theta$.  The proofs for $\id-\vartheta$ and $\li-\pi$ are similar.
\end{proof}

\begin{corollary}\label{maindegree2}
One has $\frac{1}{2} \leq \deg(\li-\pi) \leq 1$, and the Riemann hypothesis holds if and only if $\deg(\li-\pi) = \frac{1}{2}$.
\end{corollary}

Note that the prime number theorem follows from Corollary \ref{montlem} provided that $\Theta <1$, i.e., provided that the anti-Riemann hypothesis is false.  In fact,   the conjecture $\Theta < 1$ is much stronger than any known $O$ bound on $\li-\pi$.   We now show how (\ref{psiest}), together with a result concerning zero-free regions of $\zeta(s)$,  is used to derive a strong  version of the prime number theorem with error bound.   

\begin{figure}[ht!]
\includegraphics[width=55mm]{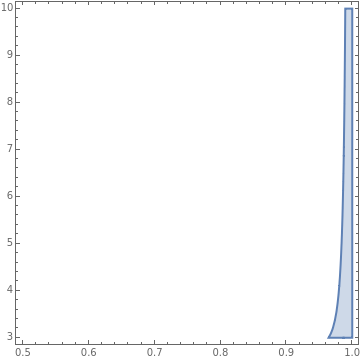}
\caption{\centering  Zero-free region $\{s = \sigma+it: \sigma \geq 1-\frac{1}{57.54}(\log t)^{-2/3} (\log \log t)^{-2/3} \text{ and } t \geq 3\}$  \newline  of $\zeta(s)$  \cite[Theorem 5]{ford2}, in $\{s: \operatorname{Im} s \leq 10\}$}
\label{ZeroFreeRegion}
\end{figure}

To this end,  one seeks zero-free regions for $\zeta(s)$ of the form
\begin{align}\label{ZFR}
\{s = \sigma+it: \sigma \geq 1-c(\log t)^{-a/(a+1)} (\log \log t)^{-1/(a+1)} \text{ and } t \geq t_0\},
\end{align}
where $a, c, t_0 > 0$.   The best admissible value of $a$ known is $a = 2$.  For example,  K.\ Ford proved in 2002  \cite[Theorem 5]{ford2}  that $\zeta(s)$ has an explicit zero-free region of the form above for $a = 2$, with accompanying explicit values $c =  \frac{1}{57.54}$ and $t_0 = 3$.    Note that the function
$ 1-c(\log t)^{-2/3} (\log \log t)^{-1/3}$ is increasing on  $(e,\infty)$,  with limit $1$ as $t \to \infty$, and, for $c = \frac{1}{57.54}$ and $t =  3$,  it assumes the value $0.9640188077\ldots$.   See Firgure \ref{ZeroFreeRegion} for a  plot of Ford's zero-free region of $\zeta(s)$.

\begin{theorem}
 Given any zero-free region (\ref{ZFR})  for $\zeta(s)$, where $a, c, t_0 > 0$, 
there exists a constant $C> 0$ such that
\begin{align*}
\psi(x) =x + O\left( x\exp \left(-C(\log x)^{(a+1)/(2a+1)} (\log \log x)^{-1/(2a+1)} \right) \right) \ (x \to \infty).
\end{align*}
\end{theorem}

\begin{proof}
For any zero $\rho = \sigma+it$ of $\zeta(s)$ with $t_0 \leq |t| \leq T$, we have
$$\sigma \leq 1-c(\log T)^{-a/(a+1)} (\log \log T)^{-1/(a+1)}.$$
Let $\varepsilon  \in (0,1)$.  Applying (\ref{psiest}) to $x^\varepsilon \leq T \leq x$, along with (\ref{rhoest}), we obtain
\begin{align*}
x-\psi(x) & = O\left(x^{1-c(\log T)^{-a/(a+1)} (\log \log T)^{-1/(a+1)}} \sum_{\rho: \, |\operatorname{Im}\rho| \leq T} \frac{1}{|\rho|} \right) + O\left( \frac{x(\log x)^2}{T}\right)  \\
& = O\left( x\exp\left(-c \log x (\log T)^{-a/(a+1)} (\log \log T)^{-1/(a+1)} \right)( \log T)^2 \right) + O\left( \frac{x(\log x)^2}{T}\right)  \\
& = O\left( x\exp\left(-d \log x (\log T)^{-a/(a+1)} (\log \log T)^{-1/(a+1)} \right) \right) + O\left( \frac{x(\log x)^2}{T}\right)
\end{align*}
for some $d > 0$.  To make the two $O$ terms above approximately equal to one another,  let
$$T =\exp \left( (\log x)^{(a+1)/(2a+1)} (\log \log x)^{-1/(2a+1)}  \right).$$ It then follows that
\begin{align*}
x-\psi(x)  = O\left( x\exp \left(-C(\log x)^{(a+1)/(2a+1)} (\log \log x)^{-1/(2a+1)} \right) \right) \ (x \to \infty)
\end{align*}
for some $C > 0$. 
\end{proof}

Applying the theorem to $a = 2$ (using, for example,  \cite[Theorem 5]{ford2}), we obtain the {\bf prime number theorem with error bound}\index{prime number theorem with error bound} in the  following form.

\begin{theorem}[{Prime number theorem with error bound \cite{ford}}]\label{PNTEB}
There exists a constant $C> 0$ such that
\begin{align*}
\psi(x) =x + O\left(x\exp \left( -C(\log x)^{3/5} (\log \log x)^{-1/5}  \right) \right)  \ (x \to \infty).
\end{align*}
\end{theorem}

Note that the constant $C$ in Theorem \ref{PNTEB} can be made explicit, e.g., $C = 0.2098$ \cite{ford}.

Let $\alpha \in (0,1]$.   As in the proof of Theorem \ref{PNTEB}, one can show that, if there exist constants $C_1> 0$ and $t_0 > 0$ such that 
$$\zeta(s) \neq 0 \text{ for all } s = \sigma + it \text{ with } \sigma \geq 1-C_1 (\log t)^{1-1/\alpha} \text{ and } t \geq t_0,$$
then there exists a constant $C_2 > 0$ such that
$$\psi(x) = x + O \left(x e^{-C_2(\log x)^\alpha} \right) \ (x \to \infty).$$   In fact, this follows from a more general theorem of Ingham from 1932,  namely, \cite[Theorem 22]{ing2}.  The converse, which is  more difficult to establish,  was proved by Tur\'an in 1950 \cite[Theorem]{turan}.    In 1961, Sta\'s generalized Tur\'an's result by proving  a partial converse \cite[Satz I]{stas} of Ingham's general result \cite[Theorem 22]{ing2}.  For a survey of these results, along with a slight sharpening of Ingham's theorem, see  \cite{pintz0}.  Together, they yield, for example, the following.

\begin{theorem}[{\cite[Theorem 22]{ing2}  \cite[Satz I$^\prime$]{stas}}]\label{antiR2}
Let $\alpha \in (0,1)$ and $\beta \in \RR$.    The following conditions are equivalent.
\begin{enumerate}
\item There exist constants $C_1> 0$ and $t_0 > 0$ such that $\zeta(s) \neq 0$ for all $s = \sigma + it$ with $\sigma \geq 1-C_1 (\log t)^{1-1/\alpha}(\log \log t)^{\beta/\alpha}$  and $t \geq t_0$.
\item There exists a constant $C_2 > 0$ such that any of the following equivalent conditions  holds.
\begin{enumerate}
\item $\psi(x) = x + O \left(x e^{-C_2(\log x)^\alpha (\log\log x)^\beta} \right) \ (x \to \infty)$.
\item $\vartheta(x) = x + O \left(x e^{-C_2(\log x)^\alpha (\log\log x)^\beta} \right) \ (x \to \infty)$.
\item $\pi(x) = \li(x) + O \left(x (\log x)^{-1} e^{-C_2(\log x)^\alpha(\log\log x)^\beta} \right) \ (x \to \infty)$.
\end{enumerate}
\item There exists a constant $C_3 > 0$ such that
$$\pi(x) = \li(x) + O \left(x e^{-C_3(\log x)^\alpha(\log\log x)^\beta} \right) \ (x \to \infty).$$
\end{enumerate}
\end{theorem}

\begin{proof}
The implication $(1) \Rightarrow (2)(a)$ follows from Ingham's \cite[Theorem 22]{ing2}, and the reverse implication follows from  Sta\'s's  \cite[Satz I$^\prime$]{stas}.   Moreover, the equivalence of (2)(a) and (2)(b) follows from Proposition \ref{pitheta2b}(2),  and the equivalence of statements (2)(b) and (2)(c) follows from Proposition \ref{pitheta} and Karamata's integral theorem (or, see Theorem \ref{lithetapsi}(3)).  Finally, the equivalence of statements (2)(c) and (3) is clear.
\end{proof}

\section{The zeros of $\zeta(s)$ and the Riemann--von Mangoldt formula}

For any $T > 0$, let $$N(T) = \#\{\rho \in \CC: \zeta(\rho) = 0,\, \operatorname{Im} \rho \in (0,T]\}\index[symbols]{.t  Gi@$N(T)$} $$ denote the number of zeros of $\zeta(s)$, counting multiplicities, with imaginary part lying in the interval $(0,T]$,  that is, with imaginary part greater than $0$ and less than or equal to $T$.  For example, one has $N(50) = 10$, since there are exactly 10 zeros of $\zeta(s)$ with  imaginary part lying in the interval $(0,50]$.  The {\bf Riemann--von Mangoldt formula},\index{Riemann--von Mangoldt formula}\index{Riemann--von Mangoldt formula}  conjectured by Riemann in 1859 \cite{rie} and proved by von Mangoldt  in 1905 \cite{mang}, states that
\begin{align}\label{RVMF}
N(T)={\frac  {T}{2\pi }}\log {{\frac  {T}{2\pi }}}-{\frac  {T}{2\pi }}+O(\log T) \ (T \to \infty),
\end{align}
 or, equivalently,
$$N(2 \pi T)=T \log T-T+O(\log T)  = \int_0^T \log t \, dt + O(\log T) \ (T \to \infty).$$ 
A much stronger  version of the formula is stated in Theorem \ref{RMF} below.  

\begin{figure}[ht!]
\includegraphics[width=80mm]{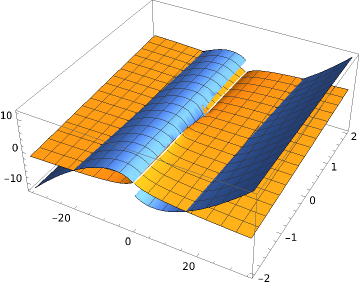}
\caption{\centering  Graph of $\operatorname{Re} \theta(s)$ (in blue) and $\operatorname{Im} \theta(s)$ (in gold) for $\operatorname{Re} s \in [-35,35]$ and $\operatorname{Im} s \in [-2,2]$}
    \label{Siegel1}
\end{figure}

\begin{figure}[ht!]
\includegraphics[width=80mm]{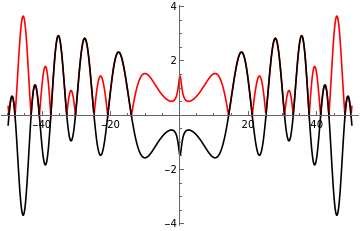}
\caption{\centering Graph  $Z(t)$ (in black) and $|\zeta(\frac{1}{2}+it)|$ for $t \in [-50,50]$}
 \label{Siegel2}
\end{figure}

\begin{figure}[ht!]
\includegraphics[width=80mm]{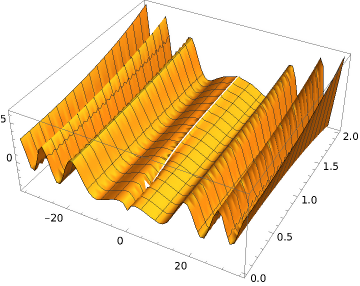}
\caption{\centering Graph of $\operatorname{Re} Z(s)$ for $\operatorname{Re} s \in [-35,35]$ and $\operatorname{Im} s \in [0,2]$}
   \label{Siegel3}
\end{figure}

The {\bf Riemann--Siegel theta function}\index{Riemann--Siegel theta function $\theta(s)$}\index[symbols]{.t  Ge@$\theta(s)$} \cite[Chapter 1]{ivic2} is the function
$$\theta(t) = \arg \Gamma \left( \frac{1}{4}+\frac{it}{2}\right)-\frac{\log \pi}{2}t = \operatorname{Im} \log \Gamma \left( \frac{1}{4}+\frac{it}{2} \right)-\frac {\log \pi}{2}t,$$
where the argument is chosen in such a way that $\theta(t)$ is continuous and $\theta(0) = 0$.  The second expression allows  one to extend the function $\theta$ to a complex function, namely, by defining
$${\theta (s) = -{\frac {i}{2}}\left(\log \Gamma \left({\frac {1}{4}}+{\frac {is}{2}}\right)-\log \Gamma \left({\frac {1}{4}}-{\frac {is}{2}}\right)\right)-{\frac {\log \pi }{2}s}}.$$ 
Since the log-gamma function $\log \Gamma(s)$ is analytic on $\CC\backslash (-\infty,0]$,  the function $\theta(s)$ so defined  is analytic on $\CC\backslash(\frac{ i}{2}\RR_{\geq 1} \cup \frac{ -i}{2}\RR_{\geq 1})$.
The function $\theta(s)$ is an odd function and has the asymptotic expansion
$$\theta(t) \simeq \frac{t}{2}\log \frac{t}{2\pi} - \frac{t}{2} - \frac{\pi}{8}+\frac{1}{48t}+ \frac{7}{5760t^3}+ \frac{31}{80640t^5}+\cdots \ (t \to \infty)$$
over $\RR$, where the numerators and denominators of the coefficients of $\frac{1}{t^{2n+1}}$ are OEIS Sequences A036282 and A114721, respectively.   The {\bf Riemann--Siegel $Z$ function},\index{Riemann--Siegel $Z$ function $Z(s)$} or {\bf Hardy's $Z$ function},\index{Hardy's $Z$ function $Z(s)$}\index[symbols]{.t  Gf@$Z(s)$} is defined by
$$Z(s) = e^{i \theta(s)} \zeta\left(\frac{1}{2}+ is \right)$$
\cite[Chapter 1]{ivic2}.  The function $Z(s)$ an even function, analytic in the region  $\CC\backslash(\frac{ i}{2}\RR_{\geq 1} \cup \frac{ -i}{2}\RR_{\geq 1})$, and it satisfies
$$Z(s) \in \RR, \quad \forall s \in \RR,$$
and
$$|Z(s)| = \left | \zeta\left(\frac{1}{2}+ is\right) \right|, \quad \forall s \in \CC.$$
The real zeros of $Z(s)$ are precisely the imaginary parts of the nontrivial zeros of $\zeta(s)$, and the zeros of $Z(s)$ are  $-i(\rho-\frac{1}{2})$, where $\rho$ ranges over the zeros of $\zeta(s)$.  Thus, the Riemann hypothesis is true if and only if all of the zeros of $Z(s)$ in the region $\{s \in \CC: \operatorname{Im} s \in (-\frac{1}{2},\frac{1}{2})\}$ are real.    See Figures \ref{Siegel1}, \ref{Siegel2}, and \ref{Siegel3} for graphs of the functions $\theta(s)$ and $Z(s)$.

As is customary, we  let
$$S(T) =  \frac{1}{\pi}\operatorname{arg}\zeta\left(\frac{1}{2}+ i T\right),\index[symbols]{.t  Gj@$S(T)$}$$
where the argument of $\zeta\left(\frac{1}{2}+i T\right)$ is chosen to be $0$ at $\infty+iT$ and to vary continuously on the line from $\infty+iT$ to $1/2+iT$ (or, equivalently,  to be $0$ at $2$ and to vary  continuously on the line from $2$ to $2+iT$ followed by the line from $2+iT$ to $1/2+iT$).  

\begin{theorem}[{Riemann--von Mangoldt formula  \cite[Chapter 1]{ivic2} \cite[Chapter IX]{tit}}]\label{RMF}
Let $T > 0$ with $T \neq \frac{\operatorname{Im}\rho}{2\pi}$ for all nontrivial zeros $\rho$ of $\zeta(s)$. One has
$$N(2 \pi T) =1+ \frac{1}{\pi} \theta(2\pi T) +S(2\pi T).$$
Moreover, one has 
$$S(T) = O(\log T) \ (T \to \infty),$$ as well as the asymptotic expansion
\begin{align*}
1+ \frac{1}{\pi} \theta(2\pi T)  \simeq T \log T - T+ \frac{7}{8}+\frac{1}{96 \pi^2 T}+ \frac{7}{11340\pi^4 T^3}+ \frac{31}{161280 \pi^6 T^5}+\cdots \ (T \to \infty),
\end{align*}
and, in particular,
\begin{align*}
1+ \frac{1}{\pi} \theta(2\pi T) =T \log T-T+\frac{7}{8}+\frac{1+o(1)}{96\pi^2 T} \ (T \to \infty).
\end{align*}
\end{theorem}

\begin{figure}[ht!]
\includegraphics[width=80mm]{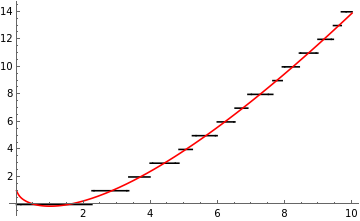}
\caption{\centering Graph of $N(2\pi t)$ (in black) and $1+ \frac{1}{\pi} \theta(2\pi t)$ (in red) for $t \in [0,10]$}
\label{RiemannSiegel}
\end{figure}

\begin{figure}[ht!]
\centering
\caption{\centering Graph of $1+\frac{1}{\pi}\theta(2\pi t)-(t \log t-t+\frac{7}{8})$ (in black) and $\frac{1}{96\pi^2 t}$ (in red) on $[0,0.5]$    \label{Siegel0a}}
\includegraphics[width=80mm]{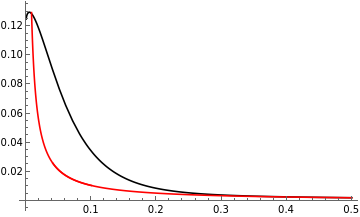}
\end{figure}

See Figure \ref{RiemannSiegel} for a graph of the function $N(2\pi t)$ and its smooth approximation $1+ \frac{1}{\pi} \theta(2\pi t)$, and see Figure \ref{Siegel0a} for a graph of the functions $1+\frac{1}{\pi}\theta(2\pi t)-(t \log t-t+\frac{7}{8})$ and $\frac{1}{96\pi^2 t}$.

At the values of $T$ that are the imaginary part of some zero of $\zeta(s)$, one defines $S(T) = \frac{S(T^+)+S(T^-)}{2}$ and re-defines $N(T)$ so that $N(T) = \frac{N(T^+)+N(T^-)}{2}$, so that the  Riemann--von Mangoldt formula $$N(2 \pi T) =1+ \frac{1}{\pi} \theta(2\pi T) +S(2\pi T)$$ holds for all $T > 0$.  

Now, let $$0<\gamma_1 \leq \gamma_2 \leq \gamma_3 \leq \cdots\index[symbols]{.ru N2@$\gamma_n$}$$ denote the ordinates of the  zeros of $\zeta(s)$ in the upper half plane, repeated to multiplicity, and listed so that
$$0<\gamma_1 \leq \gamma_2 \leq \gamma_3 \leq \cdots$$ 
and $\gamma_n \leq \gamma_{n+1}$ if $\sigma_n = \sigma_{n+1}$.
The notation $\gamma_n$ here for the {\bf $n$th zeta zero ordinate}\index{zeta zero ordinate $\gamma_n$} is standard, despite the fact $\gamma_n$ is also used to denote the $n$th Stieltjes constant.  For the remainder of this text, the notation $\gamma_n$ is used as defined above.  Two well-known corollaries of the Riemann--von Mangoldt formula are the $O$ bound (\ref{rhoest}) used in the proofs of Proposition \ref{psiestprop} and Theorem \ref{PNTEB}, and the asymptotic relation
$$\frac{\gamma_n}{2 \pi} \sim \frac{ n}{\log n} \ (n \to \infty).$$
 See  \cite[pp.\ 69--71]{ing2}, \cite[Section 1.4]{ivic}, or \cite[Section 1.4]{ivic2} for the proofs.

The asymptotic above suggests that the constants $$\tau_n  = \frac{\gamma_n}{2\pi}\index[symbols]{.ru N3@$\tau_n$}$$ are a natural renormalization of  the $\gamma_n$.  Indeed, one also has $$N(2 \pi T) = \#\{n: \tau_n \in (0,T]\} = T\log T-T + O(\log T) \ (T \to \infty)$$ for all $T > 0$ with $T \neq \tau_n$ for all $n$.  See Table \ref{Ntable} for the values of $N(2 \pi n)$ and $\tau_n$ for $n = 1,2,3,\ldots, 30$.    It is  noteworthy that $N(2 \pi n)$ is asymptotic to $p_n \sim n \log n\sim \int_0^n \log t \, dt$ and $\tau_{n}$ is asymptotic to $ \pi(n) \sim \frac{n}{\log n} \sim \int_0^n \frac{dt}{\log t}$,  especially given that the nontrivial zeros of $\zeta(s)$ are like ``eigen'' or ``characteristic'' values  of the distribution of the primes, where the function $\tau_n$ is the ``eigen'' analogue of the function $p_n$, and where $N(2\pi x)$ is the corresponding analogue of the function $\pi(x)$.   Further aspects of this duality are discussed at the end of Section 12.2, and, particularly, in Table \ref{duality}.

\begin{table}[!htbp]
 \caption{\centering Values and approximations of $N(2\pi T) \approx \int_1^{n} \log t \, dt$ and $\tau_n = \frac{\gamma_n}{2 \pi}  \approx \frac{n -11/8}{W((n-11/8)/e)}$ for   $n = 1,2,3,\ldots,30$}
  \footnotesize
\begin{tabular}{|c|c|c|c|c|c|} \hline
$n$ & $\displaystyle N(2\pi n)$   & $\displaystyle  \int_1^{n} \log t \, dt$    &  $\displaystyle  \tau_n  = \frac{\gamma_n}{2 \pi}$ & $\displaystyle   \frac{n -11/8}{W((n-11/8)/e)}$    \\ \hline 
\hline
$1$ &  $0$    & $0$  &       $2.249611375552\ldots$   &  $2.3111\ldots$ \\ \hline
$2$ &  $0$    &  $0.3862\ldots$    &      $3.345761522384\ldots$     &  $3.2875\ldots$    \\   \hline
$3$ &  $1$    &   $1.2958\ldots$   &      $3.980601614847\ldots$    &   $4.0573\ldots$   \\   \hline
$4$ &  $3$     &  $2.5451\ldots$   &      $4.842269428389\ldots$    &  $4.7332\ldots$     \\   \hline
$5$ &   $4$   &   $4.0471\ldots$   &      $5.241777852724\ldots$    &    $5.3515\ldots$  \\   \hline
$6$ &   $6$   &  $5.7505\ldots$   &      $5.982026045909\ldots$     &   $5.9297\ldots$   \\   \hline
$7$ &  $8$     &   $7.6213\ldots$  &      $6.512416395771\ldots$    &   $6.4777\ldots$    \\   \hline
$8$ &  $10$    &  $9.6355\ldots$   &      $6.895717882362\ldots$       &  $7.0019\ldots$    \\   \hline
$9$ &   $12$   &   $11.7750\ldots$  &      $7.640257056610\ldots $       &  $7.5066\ldots$    \\   \hline
$10$ & $14$      &  $14.0258\ldots$   &      $7.921751475449 \ldots$      &  $7.9949\ldots$    \\   \hline
$11$ &  $16$    &   $16.3768\ldots$  &      $8.430488500345 \ldots$      & $8.4693\ldots$     \\   \hline
$12$ &  $18$     & $18.8188\ldots$    &      $8.983699339977 \ldots$      &  $8.9316\ldots$    \\   \hline
$13$ &  $21$    &   $21.3443\ldots$  &      $9.445375410906 \ldots$      &   $9.3832\ldots$   \\   \hline
$14$ & $24$     &  $23.9468\ldots$    &      $9.681678249263 \ldots$      &   $9.8252\ldots$   \\   \hline
$15$ & $26$     & $26.6207\ldots$     &      $10.362983242540 \ldots$      &   $10.2588\ldots$   \\   \hline
$16$ &  $29$    &   $29.3614\ldots$    &      $10.676083427436 \ldots$      &  $10.6846\ldots$     \\   \hline
$17$ &  $32$    &  $32.1646\ldots$  &      $11.068653606587 \ldots$      &   $11.1033\ldots$   \\   \hline
$18$ &  $35$    &   $35.0266\ldots$  &      $11.469844378476 \ldots$      &   $11.5155\ldots$   \\   \hline
$19$ &  $38$    &   $37.9443\ldots$  &      $12.048775740002 \ldots$      &   $11.9219\ldots$   \\   \hline
$20$ &  $41$    &  $40.9146\ldots$    &      $12.277982630995  \ldots$      &  $12.3227\ldots$    \\   \hline
$21$ &  $44$    &   $43.9349\ldots$   &      $12.626935406408 \ldots$      &   $12.7184\ldots$   \\   \hline
$22$ &  $47$    &  $47.0029\ldots$    &      $13.195596946559 \ldots$      &    $13.1092\ldots$  \\   \hline
$23$ &  $50$    &  $ 50.1163\ldots$  &      $13.486072563177 \ldots$      &    $13.4957\ldots$  \\   \hline
$24$ &  $53$    & $ 53.2732\ldots$    &      $13.914164605845  \ldots$      &  $13.8779\ldots$    \\   \hline
$25$ &  $56$    &   $56.4718\ldots$   &      $14.134409040292  \ldots$      &  $14.2561\ldots$    \\   \hline
$26$ &   $60$   &   $59.7105\ldots$  &      $14.720542964867 \ldots$      &   $14.6306\ldots$   \\   \hline
$27$ &  $63$    &  $62.9875\ldots$    &      $15.064229274340  \ldots$      &  $15.0016\ldots$    \\   \hline
$28$ &  $66$    &   $66.3017\ldots$   &      $15.258285334780  \ldots$      &  $15.3692\ldots$    \\   \hline
$29$ &  $70$    &   $69.6515\ldots$   &      $15.729473091500 \ldots$      &   $15.7336\ldots$   \\   \hline
$30$ &  $73$    &  $73.0359\ldots$   &      $16.125236811010  \ldots$      &   $16.0949\ldots$   \\   \hline
\end{tabular}\label{Ntable}
\end{table}

\begin{figure}[ht!]
\includegraphics[width=80mm]{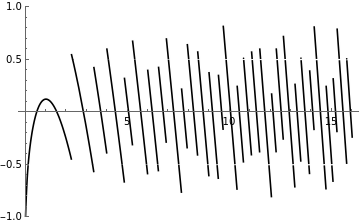}
\caption{\centering Graph of $\frac{1}{\pi}\operatorname{Arg}\zeta(\frac{1}{2}+2\pi i t)$ on $[0,18]$}
\label{Steps}
\end{figure}

\begin{figure}[ht!]
\includegraphics[width=80mm]{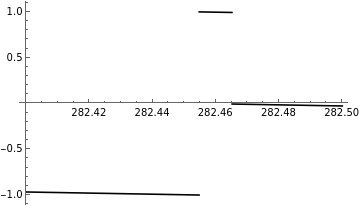}
\caption{\centering Graph of $\frac{1}{\pi}\operatorname{Arg}\zeta(\frac{1}{2}+i t)$ on $[282.4,282.5]$}
\label{zet}
\end{figure}

We now discuss the function $S(T)$.  By the Riemann--von Mangoldt formula, one has
$$S(T) = 2R(T)+ \frac{1}{\pi}\operatorname{Arg}\zeta\left(\frac{1}{2}+ i T\right), \quad \forall T > 0,$$
where $\operatorname{Arg}$ is the principal branch of $\arg$, and where $R(T) \in \ZZ$ and
$$R(T) = O(\log T) \ (T \to \infty).$$
 The function $R(T)$ is piecewise integer-valued, assuming half-integer values at the $\gamma_n$.   
See Figure \ref{Steps} for a graph of the function $$ \frac{1}{\pi}\operatorname{Arg}\zeta\left(\frac{1}{2}+2\pi i T\right) = S(2\pi T)-2R(2\pi T), \quad \forall T \in \RR_{>0} \backslash\{\tfrac{\gamma_n}{2\pi}: n \in \ZZ_{>0}\}.$$   Since  $\frac{1}{\pi}\operatorname{Arg}\zeta(\frac{1}{2}+i T) \in (-1,1],$ the size of $S(T)$ is captured by the function $R(T)$.   
The trend over intervals where the function $R(T)$ vanishes is that the function is decreasing on the interval $(\gamma_n, \gamma_{n+1})$, jumping up by $1$ unit when passing over $\gamma_n$,  from the left to the right.   
Using Mathematica, one can verify that $R(T) = -1$ for $T \in (282.4547207\ldots,\gamma_{127})$, where
$$\gamma_{127} = 282.46511476505209623302\ldots.$$  One can also verify that $282.4547207\ldots$ is the infimum of all $T > 0$ such that $R(T) \neq 0$.   See Figure \ref{zet} for a graph of $\frac{1}{\pi}\operatorname{Arg}\zeta(\frac{1}{2}+i T)$ on the interval $[282.4,282.5]$.  The fact that $R(T) = -1$ for $T \in (282.4547207\ldots,\gamma_{127})$ is reflected in the portion of the graph immediately to the right of $282.4547207\ldots$, where the graph  jumps up by $+2$.   Although $S(T)$ and $R(T)$ are unbounded above and below,  the average order of $|S(T)|$ is only on the order of $(\log \log T)^{1/2}$ (see (\ref{STav}) below),  and finding explicit values of $T$, no less the smallest value of $T$, for which $|R(T)| > 1$, appears to be a difficult problem \cite[Section 4.3.4]{gour}.

Like the function $\li(x)-\pi(x)$, much is known, but more is unknown, about the order of growth of $S(T)$.      It is known, for example, that
$$S(T) = \Omega_{\pm} \left( \frac{(\log T)^{1/3}}{(\log \log T)^{1/3}} \right) \ (T \to \infty),$$
and 
$$\int_0^T S(t) \, dt = O(\log T) \ (T\to \infty).$$
The latter of these bounds implies that $S(T)$ has mean value $0$, while  the variance from the mean satisfies
$$\frac{1}{T} \int_0^T S(t)^2 \, dt =  \frac{1}{2\pi^2}\log \log T + O((\log \log T)^{-1/2}) \ (T \to  \infty),$$
with $O$ bound $O((\log \log T)^{-1})$ on condition of the Riemann hypothesis.  See \cite[pp.\ 18--19]{ivic2} and \cite{tsang} for a discussion of these results.  It is also known that $$\frac{1}{T} \int_0^T |S(t)|^\lambda \, dt \sim \pi^{-\lambda-1/2} \Gamma\left( \frac{\lambda+1}{2} \right) (\log \log T)^{\lambda/2} \ (T \to \infty
)$$ for all $\lambda >-1$ \cite{selb}, and, in particular,
\begin{align}\label{STav}
\frac{1}{T} \int_0^T |S(t)| \, dt \sim \pi^{-3/2}  (\log \log T)^{1/2}\ (T \to \infty
).
\end{align}
Moreover,  on condition of the Riemann hypothesis, one has 
$$S(T) \leq \left( \frac{1}{4}+o(1)  \right) \frac{\log T}{\log \log T},$$
$$\int_0^T S(t) \, dt = O \left (\frac{\log T}{(\log \log T)^2} \right) \ (T \to \infty),$$
$$|S(T)| = \Omega_{+ }\left( \frac{(\log T)^{1/2}(\log \log \log T)^{1/2}}{(\log \log T)^{1/2}} \right) \ (T \to \infty),$$
and
$$\int_0^T S(t) \, dt = \Omega_+ \left( \frac{(\log T)^{1/2}(\log \log \log T)^{1/2}}{(\log \log T)^{3/2}} \right) \ (T \to \infty)$$ \cite[Theorems 1 and 2]{carn} \cite[Theorem 2]{bond}.

\begin{figure}[ht!]
\includegraphics[width=80mm]{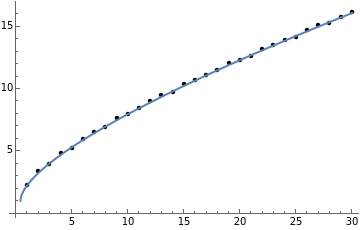}
\caption{\centering Plot of $\frac{\gamma_n}{2\pi}$ for $n = 1,2,3,\ldots,30$ and graph of $\frac{x-11/8}{W((x-11/8)/e)}$ for $x \in [3/8,30]$}
   \label{TauW}
\end{figure}

\begin{figure}[ht!]
\includegraphics[width=80mm]{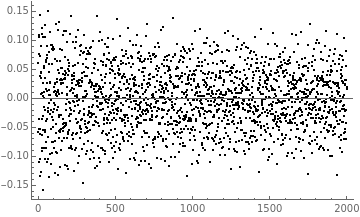}
\caption{\centering Plot of $\frac{n-11/8}{W((n-11/8)/e)}-\frac{\gamma_n}{2\pi}$ for $n = 1,2,3,\ldots,2000$}
\label{ZetaZeros1}
\end{figure}

\begin{figure}[ht!]
\includegraphics[width=80mm]{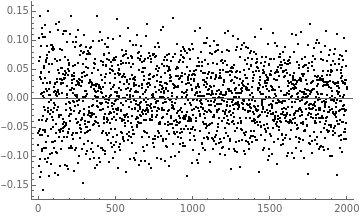}
\caption{Plot of  $\frac{1}{2}\left(\frac{n-7/8}{W((n-7/8)/e)}+\frac{n-7/8-1}{W((n-7/8-1)/e)}\right) - \frac{\gamma_{n}}{2\pi}$  for $n = 1,2,3,\ldots,2000$}
   \label{ZetaZeros}
\end{figure}

\begin{figure}[ht!]
\includegraphics[width=80mm]{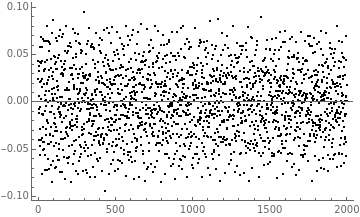}
\caption{\centering Plot of $\frac{n-7/8}{W((n-7/8)/e)}-\frac{1}{2}\left(\frac{\gamma_{n+1}}{2\pi}+\frac{\gamma_{n}}{2\pi}\right)$ for $n = 1,2,3,\ldots,2000$}
   \label{ZetaZeroa}
\end{figure}

The function $S(T)$ captures the fluctuations in $N(T)$ from its approximate mean $1+ \frac{1}{\pi} \theta(T)$,  just as $\pi(x)-\Ri(x)$ captures the fluctuations in $\pi(x)$ from its approximate mean $\Ri(x)$.   In particular, $S(T)$ captures all of the information about the $\gamma_n$.  Notice in  Table \ref{Ntable} that $N(2\pi n)$ is equal to the nearest integer to $\int_1^n \log t \, dt = t\log t - t+1$ for all  $n = 2,3,4\ldots,30$ except for $n = 15$.  Let $a \in \RR$.   The inverse of the function $$\int_0^x \log t \, dt +a= x\log x - x+a$$ is the function
$$\frac{x-a}{W((x-a)/e)} =eL((x-a)/e),$$ where $$W: [-1/e,\infty) \longrightarrow [-1,\infty)$$ is the {\bf (real and principal branch of the) Lambert $W$ function},\index{Lambert $W$ function}\index{Lambert $W$ function}\index[symbols]{.t  Ge@$W(x)$} which is the inverse of the restriction of the function $xe^x$ to the domain $[-1,\infty)$ \cite{corless}, and where $$L:  \begin{aligned}  [-1/e,\infty) &   \longrightarrow [1/e,\infty)  
 \\   x &  \longmapsto \frac{x}{W(x)} =  e^{W(x)} \end{aligned}$$ is the  inverse of the restriction of the function $x \log x$ to the domain $[1/e,\infty)$.  One can see from  Table \ref{Ntable} and Figure \ref{TauW} that  $\frac{n-11/8}{W((n-11/8)/e)}$ is a good approximation of $\tau_n$ for $n = 1,2,3,\ldots,30$.   For instance, the approximation $\frac{-3/8}{W(-(3/8)/e)}$ of $\tau_1$ is correct to within $2.74\%$.  Numerical evidence, Figure \ref{RiemannSiegel}, and Theorem \ref{tauprop} below suggest that $a = \frac{11}{8} = \frac{\frac{7}{8}+\frac{7}{8}+1}{2}$ is optimal.  Figure \ref{ZetaZeros1} shows a plot of the difference $\frac{n-11/8}{W((n-11/8)/e)}-\tau_n$ for $n = 1,2,3,\ldots,2000$, and  Figure \ref{ZetaZeros} shows a plot of the difference $\frac{1}{2}\left(\frac{n-7/8}{W((n-7/8)/e)}+\frac{n-7/8-1}{W((n-7/8-1)/e)}\right) - \tau_n$ for $n = 1,2,3,\ldots,2000$.  These approximations are based on the approximation
$$\frac{1}{2}(\tau_{n+1}-\tau_n) \approx \frac{n-7/8}{W((n-7/8)/e)},$$
which is suggested by Figure \ref{RiemannSiegel} (assuming that the red curve crosses each black horizontal line near its center, on average).  The error $\frac{n-7/8}{W((n-7/8)/e)}-\frac{1}{2}(\tau_{n+1}-\tau_n)$ in the latter approximation is plotted in Figure \ref{ZetaZeroa}. 
Note that the difference
$$\frac{x-11/8}{W((x-11/8)/e)} -\frac{1}{2}\left (\frac{x-7/8}{W((x-7/8)/e)} +\frac{x-7/8-1}{W((x-7/8-1)/e)} \right) \sim \frac{1}{8x (\log x)^2} \ (x \to \infty)$$
is very small.   Consequently, by the following theorem,  the approximation  $\frac{n-11/8}{W((n-11/8)/e)}$ of $\tau_n$ is  just as good as the approximation $\frac{1}{2}\left (\frac{n-7/8}{W((n-7/8)/e)} +\frac{n-7/8-1}{W((n-7/8-1)/e)} \right)$.

\begin{theorem}\label{tauprop}
Let $a \in \RR$, and let 
$$r_n = r_n(a) =   \frac{n-a}{W((n-a)/e)}.$$
Then one has
\begin{align*}
r_n- \frac{\gamma_n}{2\pi} & =  \frac{S(\gamma_n)+\frac{11}{8}-a}{\log r_n} + O \left( \frac{1}{r_n \log r_n}\right)   \\
& = O\left(\frac{S(\gamma_n)}{\log \gamma_n} \right)  \\
& = O(1) \ (n \to \infty).
\end{align*}
\end{theorem}

\begin{proof}
One has
\begin{align*}
n-1-\frac{1}{\pi}\theta(2\pi r_n) & = n- r_n \log \frac{r_n}{e}-\frac{7}{8} - \frac{1}{96\pi^2 r_n}+ O\left( \frac{1}{r_n^3} \right)\\
 & = a-\frac{7}{8}- \frac{1}{96\pi^2 r_n}+ O\left( \frac{1}{r_n^3} \right) \ (n \to \infty),
\end{align*}
while also
$$S(\gamma_n) = -\frac{1}{2}+n-1-\frac{1}{\pi}\theta(\gamma_n)$$ 
and therefore
$$\frac{1}{\pi}\theta(\gamma_n)-\frac{1}{\pi}\theta(2\pi r_n) =a-\frac{11}{8} -S(\gamma_n) -\frac{1}{96 \pi^2 r_n}  + O\left( \frac{1}{r_n^3} \right)\ (n \to \infty).$$
Since $\frac{1}{\pi}\theta(2\pi t)$ is approximated by $t\log \frac{t}{e} - \frac{1}{8}+\frac{1}{96 \pi^2 t}$ to within $O(t^{-3})$, where the latter function has derivative $\log t -\frac{1}{96\pi^2 t^2}$,
by the mean value theorem it follows that
\begin{align*}
(\tau_n-r_n)\left(\log u_n  -\frac{1}{96 \pi^2 u_n^2}\right) &  = \left(\tau_n \log  \frac{\tau_n}{e}   + \frac{1}{96 \pi^2 \tau_n}\right)-\left(r_n \log  \frac{r_n}{e} +  \frac{1}{96 \pi^2 r_n} \right) \\ 
&  = \frac{1}{\pi}(\theta(\gamma_n)-\theta(2\pi r_n)) + O \left( \frac{1}{r_n^3} \right) \\
& =a-\frac{11}{8} -S(\gamma_n) -\frac{1}{96 \pi^2 r_n}+ O \left( \frac{1}{r_n^3}\right)  \ (n \to \infty)
\end{align*}
for some $u_n$ lying between $r_n$ and $\tau_n$.  Using the facts that
$$S(T) = O (\log T) \ (T \to \infty)$$
and
$$\log r_n \sim \log \tau_n \sim \log u_n \sim \log n \ (n \to \infty),$$ 
we obtain
$$\tau_n - r_n = O(1) \ (n \to \infty)$$
and thus
$$u_n - r_n = O(1) \ (n \to \infty).$$
It follows that
\begin{align*}
\log u_n -\frac{1}{96 \pi^2 u_n^2}= \log r_n + O\left(\frac{1}{r_n} \right)  \ (n \to \infty).
\end{align*}
Therefore, one has
\begin{align*}
(\tau_n-r_n)\left(\log r_n +O\left(\frac{1}{r_n}\right)\right)  =a-\frac{11}{8} -S(\gamma_n) -\frac{1}{96 \pi^2 r_n}+ O \left( \frac{1}{r_n^3}\right)  \ (n \to \infty),
\end{align*}
whence
\begin{align}\label{lindelo}
\tau_n-r_n = \left(\frac{a-\frac{11}{8} -S(\gamma_n) -\frac{1}{96 \pi^2 r_n}+ O \left( \frac{1}{r_n^3}\right)}{\log r_n}\right)\left(1+ O \left( \frac{1}{r_n \log r_n}\right) \right) \ (n \to \infty).
\end{align}
The theorem follows.
\end{proof}

Since $S(T)$ has mean value $0$ in the sense that
$$\frac{1}{T} \int_0^T S(t) \, dt = O\left(\frac{\log T}{T} \right)  = o(1) \ (T\to \infty),$$
it seems reasonable to conjecture that the $S(\gamma_n)$ have mean value $0$ in the sense that
$$\frac{1}{n}\sum_{k = 1}^n S(\gamma_k) = o(1) \ (n \to \infty).$$
If this is true, then $S(\gamma_n)+\frac{11}{8}-a$ has mean value $\frac{11}{8}-a$, which, given Theorem \ref{tauprop},  implies that $a =\frac{11}{8}$ is optimal in the approximation $\gamma_n \approx r_n$, for which one has
\begin{align*}
2\pi r_n-\gamma_n = 2\pi \frac{S(\gamma_n)}{\log r_n} + O \left( \frac{1}{r_n \log r_n}\right)  \ (n \to \infty).
\end{align*}

\begin{problem}
Do the $S(\gamma_k)$ have mean value $0$?
\end{problem}

For all positive integers $n$,  we define the {\bf normalized (zeta zero) ordinates} $\widehat{\gamma}_n$\index{normalized (zeta zero) ordinate $\widehat{\gamma}_n$} of the ordinates $\gamma_n$ by 
$$\widehat{\gamma}_n = \tau_n \log \frac{\tau_n}{e}+ \frac{11}{8}.\index[symbols]{.ru N4@$\widehat{\gamma}_n$}$$ 
We also define the related normalizations
$$\widehat{N}(x) =N\left(2\pi  \frac{x-7/8}{W((x-7/8)/e)}\right),\index[symbols]{.t  Gj@$\widehat{N}(x)$}$$
$$\widehat{S}(x) =S\left(2\pi  \frac{x-7/8}{W((x-7/8)/e)}\right),\index[symbols]{.t  Gk@$\widehat{S}(x)$}$$
and
$$\widehat{Z}(x) = Z\left(2\pi  \frac{x-7/8}{W((x-7/8)/e)}\right)\index[symbols]{.t  Gl@$\widehat{Z}(x)$}$$
of the functions $N(T)$, $S(T)$, and $Z(T)$,  respectively,
where $W$ is the Lambert $W$ function.  Note that the function $\frac{x-7/8}{W((x-7/8)/e)}$ is the inverse of the function $x \log x-x+\tfrac{7}{8}$, and both functions are increasing on their (restricted) domains.   It follows that
$$\widehat{N}(x) = \# \left\{n: 0 < \widehat{\gamma}_{n}-\tfrac{1}{2} \leq x \right\}$$
for all $x \in  \RR_{> 0}\backslash \{\widehat{\gamma}_n-\tfrac{1}{2}: n \in \ZZ_{>0}\}$.  Moreover, if the Riemann hypothesis holds, then $\widehat{\gamma}_n-\tfrac{1}{2}$ is  the $n$th largest zero of $\widehat{Z}(x)$ on $[0,\infty)$.    Note the asymptotic expansion
\begin{align*}
\widehat{N}(x) \simeq x + \widehat{S}(x)+\frac{1}{96\pi^2  \frac{x-7/8}{W((x-7/8)/e)}}+ \frac{7}{11340\pi^4 \left( \frac{x-7/8}{W((x-7/8)/e)}\right)^3}+ \frac{31}{161280 \pi^6\left ( \frac{x-7/8}{W((x-7/8)/e)}\right)^5}+\cdots
\end{align*}
and, in particular,
$$\widehat{N}(x)  = x + \widehat{S}(x)+ (1+o(1)) \frac{\log x}{96 \pi^2 x}   \ (x \to \infty),$$
where also $\widehat{S}(x) = O(\log x)  \ (x \to \infty)$.

\begin{figure}[ht!]
\includegraphics[width=80mm]{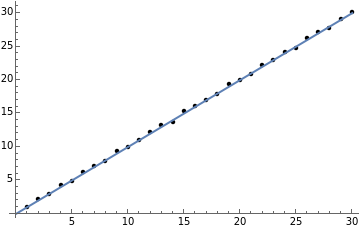}
\caption{\centering Plot of $\widehat{\gamma}_n$ for $n= 1,2,3,\ldots, 30$ and graph of  $\id(x) = x$ on $[0,30]$}
\label{GammaHatX}
\end{figure}

\begin{figure}[ht!]
\includegraphics[width=80mm]{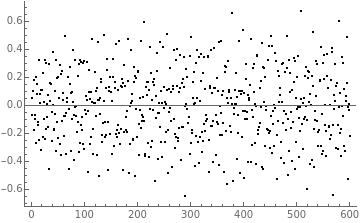}
\caption{\centering Plot of $n-\widehat{\gamma}_n$ for $n= 1,2,3,\ldots, 600$}
 \label{nln}
\end{figure}

See Figure  \ref{GammaHatX} for a plot of $\widehat{\gamma}_n$ for $n= 1,2,3,\ldots, 30$ and a graph of  $\id(x) = x$ on $[0,30]$, and see Figure \ref{nln} for a plot of the difference $n-\widehat{\gamma}_n$ for $n= 1,2,3,\ldots, 600$.  The following result shows that $n$ is a good approximation of $\widehat{\gamma}_n$.   It is an analogue of Theorem \ref{tauprop} for the $\widehat{\gamma}_n$, but its proof is simpler, which suggests that  the $\widehat{\gamma}_n$ are more directly amenable to study than the $\gamma_n$.

\begin{proposition}\label{lnnew}
One has
\begin{align*}
n - \widehat{\gamma}_n  & \simeq S(\gamma_n) + \frac{1}{96 \pi^2 \tau_n}+ \frac{7}{11340\pi^4 \tau_n^3}+ \frac{31}{161280 \pi^6 \tau_n^5}+\cdots  \ (n \to \infty)\\
& = S(\gamma_n) +\frac{1+o(1)}{96 \pi^2\tau_n}  \ (n \to \infty) \\
 & = S(\gamma_n)   + (1+o(1)) \frac{\log n}{96 \pi^2 n} \ (n \to \infty)  \\
&  = O(\log n) \ (n \to \infty).
\end{align*}
\end{proposition}

\begin{proof}
By Theorem \ref{RMF}, one has
\begin{align*}
n & =  N(\gamma_n) +\tfrac{1}{2}\\
& = N(2 \pi \tau_n) +\tfrac{1}{2}\\
& \simeq S(\gamma_n) + \widehat{\gamma}_n + \frac{1}{96 \pi^2 \tau_n}+ \frac{7}{11340\pi^4 \tau_n^3}+ \frac{31}{161280 \pi^6 \tau_n^5}+\cdots  \ (n \to \infty)\\
& = S(\gamma_n) + \widehat{\gamma}_n+\frac{1+o(1)}{96 \pi^2\tau_n}   \ (n \to \infty) \\
 & = S(\gamma_n) + \widehat{\gamma}_n   + (1+o(1)) \frac{\log n}{96 \pi^2 n}   \ (n \to \infty).
\end{align*}
 The proposition follows.
\end{proof}

\begin{proposition}\label{sntn}
For all positive integers $n$, let
$$r_n = \frac{n-11/8}{W((n-11/8)/e)}.$$ 
One has
$$n - \widehat{\gamma}_n \sim (r_n-\tau_n)\log \tau_n \ (n \to \infty).$$ More precisely, one has
$$0 \leq n - \widehat{\gamma}_n-(r_n -\tau_n)\log \tau_n  \leq (1+o(1)) \frac{(r_n-\tau_n)^2}{\tau_n} \ (n \to \infty).$$
\end{proposition}

\begin{proof}
Let $f(x) = x\log \frac{x}{e}+\frac{11}{8}$, so that  $f(\tau_n) = \widehat{\gamma}_n$, $f(r_n) = n$, and $f'(x) = \log x$.
By the mean value theorem, there exists a $u_n$ lying between $\tau_n$ and $r_n$ such that 
$$n - \widehat{\gamma}_n = f(r_n) - f(\tau_n) \sim (r_n -\tau_n)\log u_n \ (n \to \infty).$$
The first  claim follows, since $\log u_n \sim \log \tau_n \sim \log r_n \ (n \to \infty)$.  
Next, since $\frac{r_n - \tau_n}{\tau_n} \to 0$ as  $n \to \infty$, one  has
$$|\log \tau_n-\log u_n| \leq | \log r_n-\log \tau_n| = \left|\log\left(1 +\frac{r_n-\tau_n}{\tau_n} \right) \right| \sim \left|\frac{r_n-\tau_n}{\tau_n}\right| \ (n \to \infty),$$
and therefore, since also
$$n - \widehat{\gamma}_n-(r_n -\tau_n)\log \tau_n = (r_n -\tau_n)(\log u_n-\log \tau_n) \geq 0 \ (n \to \infty),$$
it follows that
$$0 \leq n- \widehat{\gamma}_n-(r_n -\tau_n)\log \tau_n  \leq (1+o(1)) \frac{(r_n-\tau_n)^2}{\tau_n} \ (n \to \infty).$$
The proposition follows.
\end{proof}

See Table \ref{lnvalues} for a list of values of $\frac{1}{n-\widehat{\gamma}_n}$ and $\frac{1}{(r_n-\tau_n)\log \tau_n}$ for various $n$.  The data in the table exhibits the conjecture that
$$\lim_{n \to \infty} \left(  \frac{1}{(r_n-\tau_n)\log \tau_n}-\frac{1}{n-\widehat{\gamma}_n}\right) = 0,$$
which holds provided that $\tau_n \neq r_n$ for all $n \gg 0$.  In fact, Proposition \ref{sntn} implies that
$$0 \leq \frac{1}{(r_n-\tau_n)\log \tau_n}-\frac{1}{n-\widehat{\gamma}_n} \leq (1+o(1))  \frac{1}{\tau_n (\log \tau_n)^2}  \ (n \to \infty)$$
provided that $\tau_n \neq r_n$ for all $n \gg 0$.

\begin{table}[!htbp]
 \footnotesize  \caption{\centering Table of $\frac{1}{n-\widehat{\gamma}_n}$  and $\frac{1}{(r_n-\tau_n)\log \tau_n}$ for various $n$}
\begin{tabular}{|r|r|r|} \hline
$n$ & $\displaystyle \frac{1}{n-\widehat{\gamma}_n} \quad$     &  $\displaystyle \frac{1}{(r_n-\tau_n)\log \tau_n} \quad $  \\ \hline \hline
$1$ & $+19.71526144673401056154\ldots$     &  $+20.04483812568006946446\ldots$  \\ \hline
$2$ & $-14.30703388801728527409\ldots$     &    $-14.20321866424249532694 \ldots$   \\  \hline
$3$ & $+9.37494309862990575103\ldots$     &    $+9.43989502957890698116 \ldots$   \\  \hline
$4$ & $-5.85320825639670374797\ldots$     &    $-5.81109012207778107267 \ldots$   \\  \hline
$5$ & $+5.46656140931184211744\ldots$     &    $+5.50086172566936059062 \ldots$   \\  \hline
$6$ & $-10.70907325228702936779\ldots$     &    $-10.68280964256289874907 \ldots$   \\  \hline
$7$ & $-15.39242101342964707573\ldots$     &    $-15.37048208068774567801 \ldots$   \\  \hline
$8$ & $+4.85974926820033887027\ldots$     &    $+4.87902162159842157520 \ldots$   \\  \hline
$9$ & $-3.69439395541925938817\ldots$     &    $-3.67840445621719254779 \ldots$   \\  \hline
$10$ & $+6.58774685109342091006\ldots$     &    $+6.60240476862318634039 \ldots$   \\  \hline
$11$ & $+12.06100732809558158889\ldots$     &    $+12.07402302848027647876 \ldots$   \\  \hline
$12$ & $-8.75505285753604597378\ldots$     &    $-8.74346777883453395125 \ldots$   \\  \hline
$13$ & $-7.17084681659511459533\ldots$     &    $-7.16031001417998468879 \ldots$   \\  \hline
$14$ & $+3.05818889085855308905\ldots$     &    $+3.06812765439980107028 \ldots$   \\  \hline
$15$ & $-4.11251374323512307297\ldots$     &    $-4.10364002054709520885 \ldots$   \\  \hline
$16$ & $+49.86704521721930167706\ldots$     &    $+49.87539364992112101443 \ldots$   \\  \hline
$17$ & $+12.00277843216467699984\ldots$     &    $+12.01058085117508425587 \ldots$   \\  \hline
$18$ & $+8.95896680160964698747\ldots$     &    $+8.96627485472449696526 \ldots$   \\  \hline
$19$ & $-3.17277050042100725006\ldots$     &    $-3.16603385988764717873 \ldots$   \\  \hline
$20$ & $+8.91550638155334653584\ldots$     &    $+8.92196907457112139586 \ldots$   \\  \hline
$10^2$ & $-3.22614633250428601805\ldots$     &    $-3.22513623565617638448 \ldots$   \\  \hline
$10^3$ & $+12.16604268013329528233\ldots$     &    $+12.16611801684375851694 \ldots$   \\  \hline
$10^4$ & $-5.57820925609988932271\ldots$     &    $-5.57820338502181397472 \ldots$   \\  \hline
$10^5$ & $+10.53613007653254099926\ldots$     &    $+10.53613055247819481992 \ldots$   \\  \hline
$10^6$ & $-14.02630395955963681475\ldots$     &    $-14.02630391975951795638 \ldots$   \\  \hline
$10^7$ & $+4.79124001591285692567\ldots$     &    $+4.79124001932233517057 \ldots$   \\  \hline
$10^8$ & $+31.85036629856406953094\ldots$     &    $+31.85036629886171286854 \ldots$   \\  \hline
$10^9$ & $+1.64237795987680678667\ldots$     &    $+1.64237795990318454539 \ldots$   \\  \hline
$10^{10}$ & $-2.26405998909815729826\ldots$     &    $-2.26405998909579096837 \ldots$   \\  \hline
$10^{11}$ & $+1.40942189440702258896\ldots$     &    $+1.40942189440723701467 \ldots$   \\  \hline
$10^{12}$ & $+2.14781207913575065516\ldots$     &    $+2.14781207913577024942 \ldots$   \\  \hline
$10^{13}$ & $-10.58795756132271662229\ldots$     &    $-10.58795756132271481896 \ldots$   \\  \hline
$10^{14}$ & $+1.77283144533498110966\ldots$     &    $+1.77283144533498127664 \ldots$   \\  \hline
$10^{15}$ & $-2.37473921460691433900\ldots$     &    $-2.37473921460691432346 \ldots$   \\  \hline
$10^{16}$ & $+2.22795385571023829615\ldots$     &    $+2.22795385571023829760 \ldots$   \\  \hline
$10^{17}$ & $+7.06289642244433025703\ldots$     &    $+7.06289642244433025717 \ldots$   \\  \hline
$10^{18}$ & $-13.67446709093282300622\ldots$     &    $-13.67446709093282300621 \ldots$   \\  \hline
$10^{19}$ & $+6.70242133201323746240\ldots$     &    $+6.70242133201323746240 \ldots$   \\  \hline
\end{tabular}\label{lnvalues}
\end{table}

Now we address the Riemann zeta zero gaps.   One defines the {\bf normalized (zeta zero ordinate) spacings}\index{normalized (zeta zero ordinate) spacing $\delta_n$}
$$\delta_n = (\tau_{n+1}-\tau_n) \log \tau_n.$$\index[symbols]{.ru N6@$\delta_n$}
In 1973, Montgomery conjectured the following.

\begin{conjecture}[{\cite[p.\ 185]{mont3}}]\label{normsp}
One has
$$\limsup_{n \to \infty} \delta_n =  \infty$$
and
$$\liminf_{n \to \infty} \delta_n = 0.$$
\end{conjecture}

Odlyzko  observed  in \cite[p.\ 277]{odl2} that the normalized spacings $\delta_n$ have mean value $1$ in the sense that
\begin{align}\label{odlm}
\sum_{n = N+1}^{N+M} \delta_n = M + O(\log(NM)) \ (N,M \to \infty).
\end{align}
In 2012, S.\ Feng and X.\ Wu proved that
$$\limsup_{n \to \infty} \delta_n \geq 2.7327$$
and
$$\liminf_{n \to \infty} \delta_n \leq 0.5154.$$
From the lim sup statement above it follows that
\begin{align}\label{tau1}
\tau_{n+1}-\tau_n \neq o\left(\frac{1}{\log n}\right) \ (n \to \infty).
\end{align}
On the other hand, Littlewood proved unconditionally (see \cite[Theorem 9.12]{tit}) that
\begin{align}\label{tau2}
\tau_{n+1}-\tau_n = O\left(\frac{1}{\log \log \log \tau_n} \right) \ (n \to \infty),
\end{align}
while Goldston and Gonek, sharpening another result of Littlewood, proved that
$$\limsup_{n \to \infty} \, (\tau_{n+1}-\tau_n) \log \log \tau_n \leq \tfrac{1}{2},$$
that is,
\begin{align}\label{tau3}
\tau_{n+1}-\tau_n \leq \frac{\frac{1}{2}+o(1)}{\log  \log \tau_n}  \ (n \to \infty),
\end{align}
on condition of the Riemann hypothesis \cite[Corollary 1]{gogo}.

The following proposition relates the normalized spacings $\delta_n$ to the spacings $\widehat{\gamma}_{n+1}-\widehat{\gamma}_n$ of the  normalized zeros $\widehat{\gamma}_n$.

\begin{proposition}\label{deltal}
One has
$$\widehat{\gamma}_{n+1}-\widehat{\gamma}_n \sim \delta_n \sim (\tau_{n+1}-\tau_n)\log n \ (n \to \infty).$$
More precisely, one has
\begin{align*}
\widehat{\gamma}_{n+1}-\widehat{\gamma}_n & = 1-(S(\gamma_{n+1})-S(\gamma_n)) +\frac{\tau_{n+1}-\tau_n}{96 \pi^2 \tau_n \tau_{n+1}}+ O \left(\frac{1}{\tau_n^3}\right)  \ (n \to \infty) \\
& = \delta_n + O\left(\frac{(\tau_{n+1}-\tau_n)^2}{\tau_n}\right) \ (n \to \infty) \\
& = \delta_n + O\left(\frac{\log n}{n (\log \log \log n)^2}\right) \ (n \to \infty).
\end{align*}
In fact, one has
\begin{align*}
\delta_n  \leq \widehat{\gamma}_{n+1}-\widehat{\gamma}_n  \leq \delta_n + \frac{(\tau_{n+1}-\tau_n)^2}{\tau_n}
\end{align*}
for all positive integers $n$.   Moreover, if the Riemann hypothesis holds, then one has
\begin{align*}
\delta_n  \leq \widehat{\gamma}_{n+1}-\widehat{\gamma}_n  \leq \delta_n + (\tfrac{1}{4}+o(1)) \frac{\log n}{n (\log \log n)^2} \ (n \to \infty).
\end{align*}
\end{proposition}

\begin{proof}
By Proposition \ref{lnnew}, one has
\begin{align*}
\widehat{\gamma}_{n+1}-\widehat{\gamma}_n -1  = -(S(\gamma_{n+1})-S(\gamma_n))+\frac{\tau_{n+1}-\tau_n}{96 \pi^2 \tau_n \tau_{n+1}}+ O \left(\frac{1}{\tau_n^3}\right)  \ (n \to \infty).
\end{align*}
Let $f(t) = t\log (t/e)+\frac{11}{8}$, so that $f'(t) = \log t$.  By the mean value theorem, for every positive integer $n$ there exists a $u_n \in [\tau_n,\tau_{n+1}]$ such that
\begin{align*}
\widehat{\gamma}_{n+1}-\widehat{\gamma}_n   = (\tau_{n+1}-\tau_n) \log u_n.
\end{align*}
Since $\log u_n \in [\log \tau_n, \log \tau_{n+1}]$ and
$$ \log \tau_{n+1}-\log \tau_n = \log\left(1 +\frac{\tau_{n+1}-\tau_n}{\tau_n} \right) \leq \frac{\tau_{n+1}-\tau_n}{\tau_n},$$
it follows that
\begin{align*}
\delta_n & \leq \widehat{\gamma}_{n+1}-\widehat{\gamma}_n   \\
 & \leq (\tau_{n+1}-\tau_n)\left(\log \tau_n +\frac{\tau_{n+1}-\tau_n}{\tau_n}\right) \\
& = \delta_n + \frac{(\tau_{n+1}-\tau_n)^2}{\tau_n}.
\end{align*}
Noting also (\ref{tau3}), the proposition follows.
\end{proof}

See Table \ref{lnvaluesa}  for a list of values of $\frac{1}{\delta_n-1}$, $\frac{1}{\widehat{\gamma}_{n+1}-\widehat{\gamma}_n-1}$, and $\frac{1}{\widehat{\gamma}_{n+1}-\widehat{\gamma}_n-\delta_n}$ for various $n$.  
Note that $\widehat{\gamma}_{n+1}-\widehat{\gamma}_n-\delta_n < \frac{1}{n}$ for all $n$ in the table.  If the conjecture
$$\limsup_{n \to \infty} \, (\gamma_{n+1}-\gamma_n) \sqrt{\frac{\log \gamma_n}{32}} = 1$$ of
G.\ B.\ Arous and P.\ Bourgade holds \cite{arous}, then one has
$$(\tau_{n+1}-\tau_n)^2 = O\left( \frac{1}{\log \tau_n}\right)  \ (n \to \infty)$$
and therefore
\begin{align*}
0 \leq \widehat{\gamma}_{n+1}-\widehat{\gamma}_n - \delta_n = O\left(\frac{1}{\tau_n \log \tau_n}\right) = O\left(\frac{1}{n}\right) \ (n \to \infty),
\end{align*}
by Proposition \ref{deltal}, where all three of the $O$ bounds above can assume the implicit constant $\frac{32}{4 \pi^2}+o(1) = 0.8105694691\ldots+o(1)$.  The data in Table \ref{lnvaluesa} suggests the  conjecture that
$$\lim_{n \to \infty} \left(\frac{1}{\delta_n-1}- \frac{1}{\widehat{\gamma}_{n+1}-\widehat{\gamma}_n-1}\right) = 0.$$

\begin{table}[!htbp]
 \caption{\centering Table of $\frac{1}{\delta_n-1}$, $\frac{1}{\widehat{\gamma}_{n+1}-\widehat{\gamma}_n-1}$, and $\frac{1}{\widehat{\gamma}_{n+1}-\widehat{\gamma}_n-\delta_n}$ for various $n$}  \footnotesize 
\begin{tabular}{|r|r|r|r|} \hline
$n$  & $\displaystyle \frac{1}{\delta_n-1} \quad$  & $\displaystyle \frac{1}{\widehat{\gamma}_{n+1}-\widehat{\gamma}_n-1} \quad$ & $\displaystyle \frac{1}{\widehat{\gamma}_{n+1}-\widehat{\gamma}_n-\delta_n} \quad$ \\ \hline \hline
$1$  &  $-8.9856891221963539\ldots$   &  $+8.2906491421063345\ldots$  & $4.3120940718150668\ldots$  \\ \hline
$2$  &  $-4.28619394606081674\ldots$    &  $-5.6637006566621591\ldots$   & $17.6229409846494304\ldots$ \\  \hline
$3$  &   $+5.2538578871861330 \ldots$   &  $+3.6034245437176360\ldots$ & $11.4708543273273740\ldots$  \\  \hline
$4$  &   $-2.7040035261286828 \ldots$   &  $-2.8266407639029349\ldots$  & $62.3240276070317937\ldots$ \\  \hline
$5$  &   $+4.4181298882565574 \ldots$   &  $+3.6191350630232028\ldots$  & $20.0124059463214975\ldots$ \\  \hline
$6$  &   $-19.5086468326425003 \ldots$   &  $-35.1967380108975351\ldots$  & $43.7682777154803746\ldots$ \\  \hline
$7$  &   $-3.5485662581776234 \ldots$   &  $-3.6935945982884427\ldots$  & $90.3751994462474824\ldots$  \\  \hline
$8$  &   $+2.2850292741403123 \ldots$   &  $+2.0988458869520216\ldots$   & $25.7591419192741935\ldots$ \\  \hline
$9$  &   $ -2.3386317634534486\ldots$   &  $-2.3669907468196058\ldots$  & $195.1945763654789423\ldots$ \\  \hline
$10$  &   $+18.9077313457072797 \ldots$   &  $+14.5169160832988024\ldots$  & $62.5127528414006595\ldots$  \\  \hline
$11$  &   $+5.5752193971525683\ldots$   &  $+5.0727541970451081\ldots$   & $56.2859230655267172\ldots$  \\  \hline
$12$  &   $+73.6970120987157428 \ldots$   &  $+39.6294050711331499\ldots$ & $85.7286143587246766\ldots$   \\  \hline
$13$  &   $ -2.1304880898672305\ldots$   &  $-2.1438779470267999\ldots$  & $341.1168900338339691\ldots$ \\  \hline
$14$  &   $+1.8290810863901267 \ldots$   &  $+1.7539207083092768\ldots$ & $42.6829038984551597\ldots$ \\  \hline
$15$  &   $+1.8290810863901267 \ldots$   &  $-3.7991957092606679\ldots$  & $213.5396805957398935\ldots$ \\  \hline
$16$  &   $-14.2062723455891276 \ldots$   &  $-15.8075976541891212\ldots$ & $140.2382365399825668\ldots$  \\  \hline
$17$  &   $-28.1766341874751666 \ldots$   &  $-35.3282353025346048\ldots$ & $139.1899165786024738\ldots$ \\  \hline
$18$  &   $+2.4246467634726636 \ldots$   &  $+2.3430070133187266\ldots$ & $69.5857638089923711\ldots$ \\  \hline
$19$  &   $ -2.3282208612304208\ldots$   &  $-2.3400238032178214\ldots$  & $461.5876482527182659\ldots$ \\  \hline
$20$  &   $-8.0068174603311085 \ldots$   &  $-8.3346465407566333\ldots$ & $203.5633732114223383\ldots$ \\  \hline
$10^2$  &   $ -3.5619658723459903\ldots$   &  $-3.5685894190197872\ldots$ & $1919.0917417777986960\ldots$ \\  \hline
$10^3$  &   $-7.0177169246659111 \ldots$   &  $-7.0204455885530378\ldots$  & $18055.5399504932805630\ldots$ \\  \hline
$10^4$  &   $+46.2846407132271848 \ldots$   &  $+46.2715181276575658\ldots$  & $163204.1628101496369075\ldots$ \\  \hline
$10^5$  &   $+1.5463268306131056 \ldots$   &  $+1.5463237447063188\ldots$  & $774852.2753397137795964\ldots$  \\  \hline
$10^6$  &   $+7.1949205422820425 \ldots$   &  $+7.1949178694384820\ldots$   & $19367711.2803112925965866\ldots$ \\  \hline
$10^7$  &   $-2.1633699390951802 \ldots$   &  $-2.1633699437096690\ldots$  & $1014233582.7536971885448601\ldots$ \\  \hline
$10^8$  &   $-3.5156002438513547 \ldots$   &  $-3.5156002457349145\ldots$ & $6561748046.8380965510577709\ldots$ \\  \hline
$10^9$  &   $+1.9702723351780313 \ldots$   &  $+1.9702723349453130\ldots$  & $16681000179.0303197276858354\ldots$ \\  \hline
$10^{10}$  &   $-13.5300800629068466 \ldots$   &  $-13.5300800632783673\ldots$ & $492739924570.5018453462276246\ldots$  \\  \hline
$10^{11}$  &   $+1.5452148222874096 \ldots$   &  $+1.5452148222860205\ldots$  & $1718905452507.2077476987420989\ldots$ \\  \hline
\end{tabular}\label{lnvaluesa}
\end{table}

The following renowned conjecture concerning the nontrivial zeros of $\zeta(s)$, known as {\bf Montgomery's pair correlation conjecture},\index{Montgomery's pair correlation conjecture} was first conjectured by H.\ Montgomery in 1973 \cite{mont3}.  (We re-express the conjecture in terms of the constants $\tau_n$.)   Let
$$\operatorname{sinc} x =  \begin{cases} \displaystyle \frac{\sin x}{x}  & \text{if } x \neq 0 \\
1 & \text{if } x = 0.
\end{cases}$$

\begin{conjecture}[{Montgomery's pair correlation conjecture \cite{mont3}}]
For all $0 \leq a \leq b$, one has
$$\frac{\# \left\{(\tau_m,\tau_n): 0 < \tau_m,\tau_n \leq T, \ \frac{a}{\log T} < \tau_m-\tau_n \leq \frac{b}{\log T} \right \} }{\# \left\{\tau_n: 0 < \tau_n \leq T \right\} } \sim  \int_a^b \left( 1- \operatorname{sinc}^2(\pi t) \right)  dt$$
as $T \to \infty$.
\end{conjecture}

\begin{figure}[ht!]
\includegraphics[width=70mm]{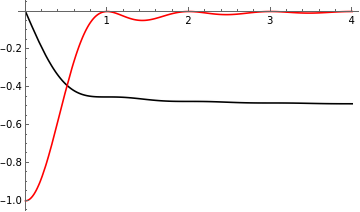}
\caption{\centering Graph of $T(x) = -\int_0^x \operatorname{sinc}^2(\pi t)  \, dt$ (in black) and $\frac{d}{dx} T(x) = - \operatorname{sinc}^2(\pi x)$ (in red) on $[0,4]$}
 \label{SX}
\end{figure}

Note that the denominator in the conjecture can be replaced with $T \log T$.  Note also that
$$ \int_0^x \left( 1- \operatorname{sinc}^2(\pi t) \right) dt= x-\frac{1}{\pi}\operatorname{Si} 2 \pi x +\frac{1-\cos 2 \pi x}{2 \pi^2 x},$$
where
$$\operatorname{Si} x = \int_0^x \operatorname{sinc} t \, dt = \operatorname{Im} \Ein(ix) \sim \frac{\pi}{2} \ (x \to \infty),$$
and therefore
$$ \lim_{x \to \infty} T(x) = -\frac{1}{2},$$
where 
$$T(x) = - \int_0^x \operatorname{sinc}^2(\pi t) \, dt = -x+\int_0^x\left( 1-\operatorname{sinc}^2(\pi t) \right) dt$$
is the unique antiderivative of the function $-\operatorname{sinc}^2(\pi x)$ with $T(0) = 0$.     Figure \ref{SX} is a graph of the functions $T(x)$ and $\frac{d}{dx}T(x)$ on the interval $[0,4]$.   Note, for example, that 
$$T(1/2) = -0.386847504951\ldots,$$
$$T(1) =-0.451411666790\ldots,$$
and
$$T(2) = -0.474969669883\ldots,$$
while $\frac{d}{dx}T(x) = 0$ precisely for $x \in \ZZ\backslash \{0\}$.

In \cite[p.\ 277]{odl2}, Odlyzko re-expressed (without proof) Montgomery's pair correlation conjecture in terms of the  normalized spacings $\delta_n$, as follows.

\begin{conjecture}[{Odlyzko's reformulation of Montgomery's pair correlation conjecture \cite[p.\ 277]{odl2}}]\label{odlc}
For all $0 \leq a \leq b$, one has
\begin{align*}
\frac{\# \left\{(n,k): n,k \in \ZZ, \, 1\leq n \leq N, \, k \geq 0,  \, \sum_{j = n}^{n+k} \delta_j \in (a,b] \right \} }{N} \sim   \int_a^b \left( 1- \operatorname{sinc}^2(\pi t) \right)  dt.
\end{align*}
as $N \to \infty$. 
\end{conjecture}

Note that, since the $\delta_n$ have mean value $1$ in the sense of (\ref{odlm}), the denominator $N$ in the conjecture can be replaced with  $\#\{n  \in \ZZ_{> 0}: \delta_1+\delta_2+\cdots+\delta_n \leq N\}$.  More importantly,  the conjecture implies Montgomery's conjecture that $$\liminf_{n \to \infty} \delta_n = 0,$$ since if $L = \liminf_{n \to \infty} \delta_n > 0$, then the limit of the ratio on the left is $0$ as $N \to \infty$ for $0 = a < b < L$.    
By Proposition \ref{deltal}, the latter conjecture is equivalent to
 $$\liminf_{n \to \infty} (\widehat{\gamma}_{n+1}-\widehat{\gamma}_n) = 0.$$
 
\begin{problem}
Deduce Odlyzko's Conjecture \ref{odlc} from  Montgomery's pair correlation conjecture and vice versa.
\end{problem}

The following conjecture is an analogue of Odlyzko's Conjecture \ref{odlc}  for the normalized zeros $\widehat{\gamma}_n$.   Our   Conjecture \ref{montpairnorm} is similar to Odlyzko's conjecture but attempts to re-express Montgomery's pair correlation conjecture in terms of the {\it spacings of the  normalized zeros}  $\widehat{\gamma}_n \approx n$ rather than the normalized spacings $\delta_n$.

\begin{conjecture}\label{montpairnorm}
For all $0\leq a\leq b$, one has
\begin{align*}
\frac{\# \left\{(\widehat{\gamma}_{m},\widehat{\gamma}_{n}): 0 < \widehat{\gamma}_m, \widehat{\gamma}_n \leq T, \ a < \widehat{\gamma}_{m}-\widehat{\gamma}_{n} \leq b \right \} }{\# \left\{\widehat{\gamma}_n: 0 < \widehat{\gamma}_{n} \leq T \right\} } \sim   \int_a^b \left( 1- \operatorname{sinc}^2(\pi t) \right)  dt.
\end{align*}
as $T \to \infty$. 
\end{conjecture}

\begin{figure}[ht!]
\includegraphics[width=80mm]{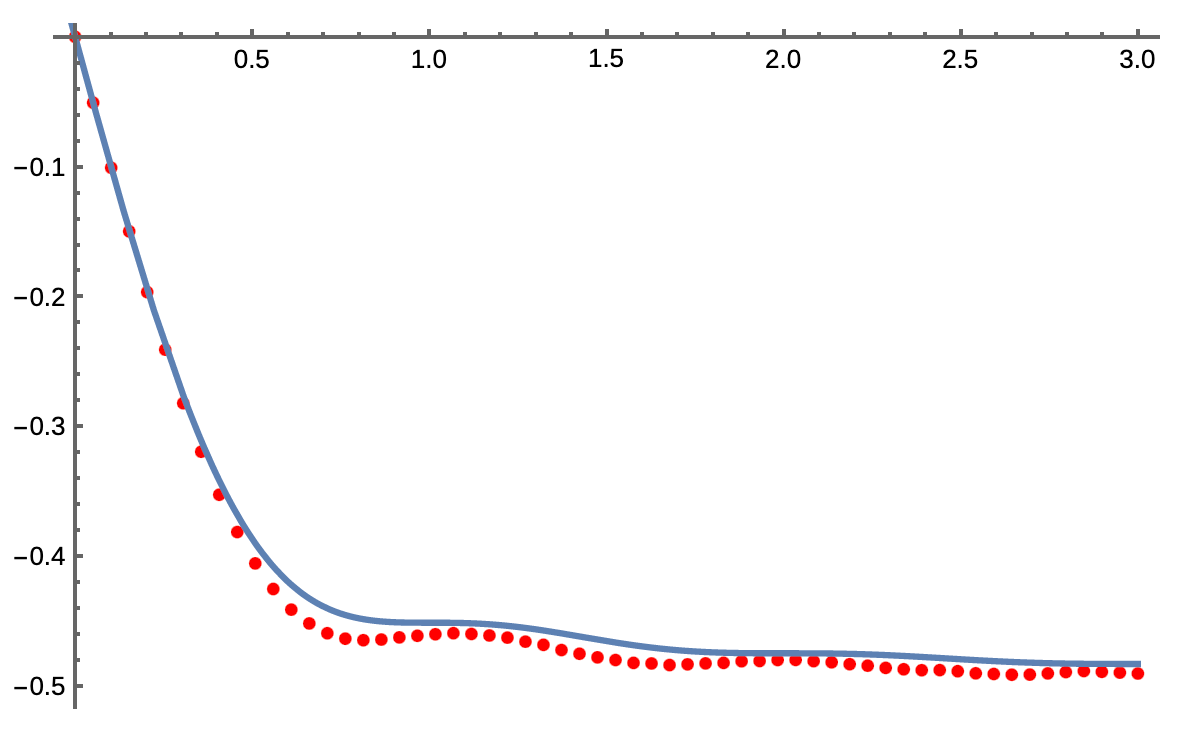}
\caption{\centering Graph of $-b+\int_0^b \left( 1- \operatorname{sinc}^2(\pi t) \right)  dt$, and plot of $-b+\frac{\# \left\{(\widehat{\gamma}_{m},\widehat{\gamma}_{n}): 0 < \widehat{\gamma}_m, \widehat{\gamma}_n \leq T, \ 0 < \widehat{\gamma}_{m}-\widehat{\gamma}_{n} \leq b \right \} }{\# \left\{\widehat{\gamma}_n: 0 < \widehat{\gamma}_{n} \leq T \right\} } $ for  $b = k/20$, where $k = 0,1,2,\ldots, 60$ and $T = 10^5$}
\label{mont1}
\end{figure}

\begin{figure}[ht!]
\includegraphics[width=80mm]{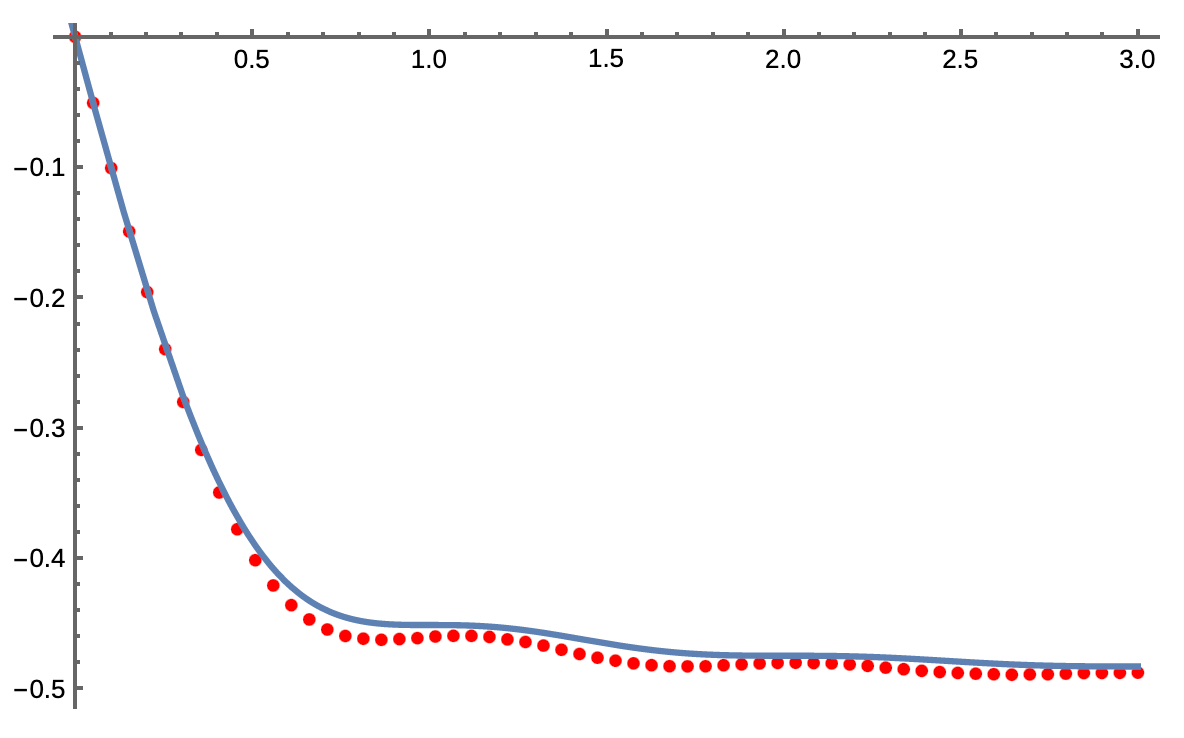}
\caption{\centering Graph of $-b+\int_0^b \left( 1- \operatorname{sinc}^2(\pi t) \right)  dt$, and plot of $-b+\frac{\# \left\{(\widehat{\gamma}_{m},\widehat{\gamma}_{n}): 0 < \widehat{\gamma}_m, \widehat{\gamma}_n \leq T, \ 0 < \widehat{\gamma}_{m}-\widehat{\gamma}_{n} \leq b \right \} }{\# \left\{\widehat{\gamma}_n: 0 < \widehat{\gamma}_{n} \leq T \right\} } $ for  $b = k/20$, where $k = 0,1,2,\ldots, 60$ and $T = 2\cdot 10^6$}
\label{mont2}
\end{figure}

\begin{figure}[ht!]
\includegraphics[width=80mm]{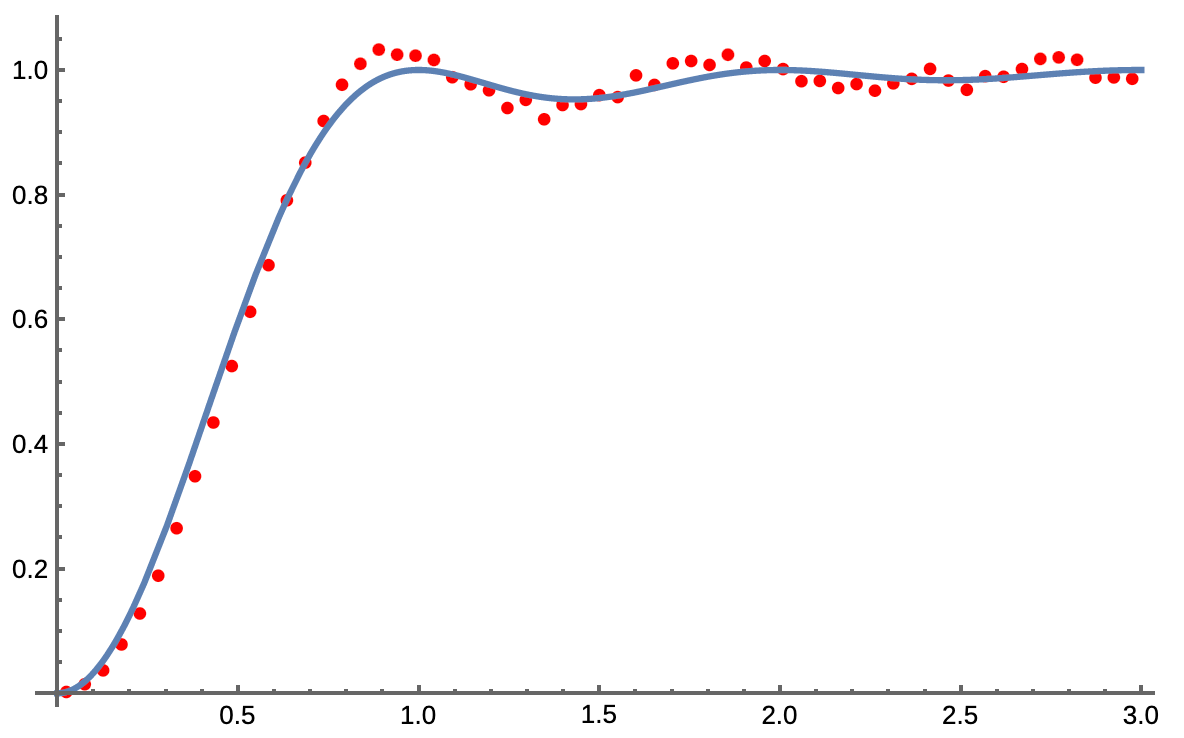}
\caption{\centering Graph of $1-\operatorname{sinc}^2(\pi x)$, and ``slope plot'' of $\frac{\# \left\{(\widehat{\gamma}_{m},\widehat{\gamma}_{n}): 0 < \widehat{\gamma}_m, \widehat{\gamma}_n \leq T, \ 0 < \widehat{\gamma}_{m}-\widehat{\gamma}_{n} \leq b \right \} }{\# \left\{\widehat{\gamma}_n: 0 < \widehat{\gamma}_{n} \leq T \right\} }$ for  $b = k/20$, where $k = 0,1,2,\ldots, 60$ and $T = 10^5$  (cf.\   \cite[Figure 1]{odl2}) }
\label{mont3}
\end{figure}

For $a = 0$ and $b = k/20$, where $k = 0,1,2,\ldots, 60$, we computed the ratio on left-hand side of the asymptotic of Conjecture \ref{montpairnorm}, with $b$ subtracted, for $T = 10^5$ and for $T = 2 \cdot 10^6$.  The results are shown in Figures \ref{mont1} and \ref{mont2}, respectively, which provide  evidence for Conjecture \ref{montpairnorm}.   Figure \ref{mont3} is  precisely Odlyzko's \cite[Figure 1]{odl2} (which he provides  as  numerical support for his Conjecture \ref{odlc}),  but using the normalized zeros instead of the normalized spacings.  (Figures \ref{mont1}, \ref{mont2}, and \ref{mont3} are courtesy of Lucas Salim.)  Along with Proposition \ref{deltal}, this gives us hope that the following conjecture is true.

\begin{conjecture}
Conjecture \ref{montpairnorm} and Odlyzko's Conjecture \ref{odlc} are equivalent.
\end{conjecture}

Interestingly, Proposition \ref{deltal} provides some intuitive support for Conjecture \ref{montpairnorm}.  If one replaces $\widehat{\gamma}_n$ with $n$ in the conjecture, then the quantity on the left is asymptotic to $b-a = \int_a^b  dt$.  Thus the integral $-\int_a^b  \operatorname{sinc}^2(\pi t) \, dt $ reflects the deviations of $\widehat{\gamma}_n$ from $n$.  Since $\widehat{\gamma}_n$ is close to $n$, one expects that the difference $\widehat{\gamma}_{n+1}-\widehat{\gamma}_n$ is close to $1$, on average.  Moreover, one expects that $\widehat{\gamma}_{m}-\widehat{\gamma}_n$ is close to $m-n$, but that the degree of its closeness to $m-n$ decays on average as $m-n$ increases.  Both of these expectations can be visualized in the graph of the density function  $-\operatorname{sinc}^2(\pi t)$ in  Figure \ref{SX}, where it is seen that the function has a global minimum of $-1$ at $0$ and a global maximum of $0$ at the positive integers, while the local minima in between the integers get smaller negative and approach $0$, occurring  at numbers that approach the half odd integers from below, as $x$ approaches $\infty$.  

\begin{figure}[ht!]
\includegraphics[width=80mm]{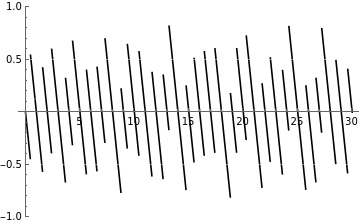}
\caption{\centering Graph of $\widehat{A}(x) = \frac{1}{\pi}\operatorname{Arg}\zeta\left(\frac{1}{2}+2\pi i \frac{x-7/8}{W((x-7/8)/e)}\right)$ on $[0,30]$ }
  \label{ZetaArg}
\end{figure}

\begin{figure}[ht!]
\includegraphics[width=80mm]{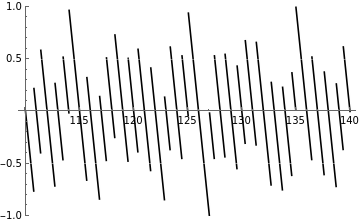}  
\caption{\centering Graph of $\widehat{A}(x) = \frac{1}{\pi}\operatorname{Arg}\zeta\left(\frac{1}{2}+2\pi i \frac{x-7/8}{W((x-7/8)/e)}\right)$ on $[110,140]$ }
\label{ZetaArg2}
\end{figure}

Let us now define
$$A(T) = \frac{1}{\pi}\operatorname{Arg}\zeta\left(\frac{1}{2}+i T\right)\index[symbols]{.t  Gm@$A(T)$}$$
and
$$\widehat{A}(x) = A\left(2\pi \frac{x-7/8}{W((x-7/8)/e)}\right).\index[symbols]{.t  Gn@$\widehat{A}(x)$}$$
Figures  \ref{ZetaArg} and  \ref{ZetaArg2} provide a graph of the function $\widehat{A}(x)$  on $[0,30]$ and $[110,140]$, respectively, to be contrasted with the graph of $A(2\pi x)$ in Figure \ref{Steps}. 
 Notice that  each continuous piece of the graph of $\widehat{A}(x)$ is approximately a line of slope $-1$.
The reason for this is that, since
$$e^{i\pi x} \zeta(\tfrac{1}{2}+i\theta^{-1}(\pi x)) = Z(\theta^{-1}(\pi x)) \in \RR,$$  the function
\begin{align}\label{atheta}
A(\theta^{-1}(\pi x)) = \frac{1}{\pi} \operatorname{Arg}\left(e^{-i\pi x}  Z(\theta^{-1}(\pi x)) \right) =  1-2\left\{\frac{1+x+ I(\theta^{-1}(\pi x))}{2} \right\}
\end{align}
is exactly piecewise linear of slope $-1$, where
$$I(x) = \left.
\begin{cases} 0 & \text{if } Z(x) \geq 0 \\
   -1   & \text{if } Z(x) < 0.
 \end{cases}\right.  $$
In fact, we have
$$A(2\pi x) = 1-2\left\{\frac{1+\frac{1}{\pi}\theta(2 \pi x)+ I(2\pi x)}{2} \right\} $$
and therefore
$$\widehat{A}(x) = 1-2\left\{\frac{1+\frac{1}{\pi}\theta\left(2 \pi \frac{x-7/8}{W((x-7/8)/e)}\right) +I\left(2 \pi \frac{x-7/8}{W((x-7/8)/e)}\right)}{2} \right\},$$
where also, since
$$1+\frac{1}{\pi}\theta(2 \pi t) = t \log (t/e) +(1+o(1)) \frac{1}{96 \pi^2 t} \ (t \to \infty),$$
one has
$$1+\frac{1}{\pi}\theta\left(2 \pi \frac{x-7/8}{W((x-7/8)/e)}\right) = x +(1+o(1)) \frac{1}{96 \pi^2 \frac{x-7/8}{W((x-7/8)/e)}}  \ (x \to \infty)$$
and where the function $1+\frac{1}{\pi}\theta\left(2 \pi \frac{x-7/8}{W((x-7/8)/e)}\right) -x$ is positive and decreasing on $[0,\infty)$ and its value at $0$ is a mere $0.000685715\ldots$ (and thus $2\pi \frac{x-7/8}{W((x-7/8)/e)} \approx \theta^{-1}((x-1)\pi)$ is a near inverse to $\frac{1}{\pi} \theta(x)+1$).  
Since also
$$ 1-2\left\{\frac{x+y}{2} \right\} + \{x\}+\{y\} = \left \lceil 1-2\left\{\frac{x+y}{2} \right\} + \{x\}\right\rceil   \in \{0,1,2\}, \quad \forall x,y \in \RR,$$
 one has the following.

\begin{proposition}\label{ahatprop}
Let $$E(x) = \left\lceil \widehat{A}(x) + \{x\} \right \rceil$$
for all $x > 0$.  The function $E: \RR_{> 0} \longrightarrow \{0,1,2\}$ is piecewise continuous, and one has
$$\{-\widehat{A}(x) - \{x\}\}  = -(\widehat{A}(x) + \{x\} -E(x)) \sim   \frac{1}{96\pi^2 \frac{x-7/8}{W((x-7/8)/e)}} \ (x \to \infty).$$ 
Consequently, one  has
$$\{-\widehat{A}(n)\} = -(\widehat{A}(n)  -E(n)) \sim  \frac{1}{96\pi^2 \tau_n} =  \frac{1}{48\pi \gamma_n}\ (n \to \infty).$$
\end{proposition}

The sequence $E(n)$ of the proposition is apparently related to the well-known sequence of {\it Gram points}.  For any positive integer $n$, the {\bf $n$th Gram point}\index{Gram point $G_n$} \cite[Section 2]{trudg}
 $$G_n = \theta^{-1}((n-1)\pi) \approx 2\pi \frac{n-7/8}{W((n-7/8)/e)} = 2\pi r_n(7/8),$$
is the unique solution to
$$\theta(G_n) = (n-1)\pi.$$
(Note that $g_{n} = G_{n+1}$, instead of $G_n$, is often called the $n$th Gram point, starting at $n = 0$.)  The $n$th Gram point $G_n$ satisfies
$$ (-1)^{n-1} \zeta(\tfrac{1}{2}+iG_n) = e^{i(n-1)\pi} \zeta(\tfrac{1}{2}+iG_n) = Z(G_n) \in \RR$$ 
and thus
$$ \zeta(\tfrac{1}{2}+iG_n) =  (-1)^{n-1}Z(G_n) \in \RR.$$   Thus, all Gram points are points $t = G_n >0$ where $\zeta(\tfrac{1}{2}+it)$ crosses the real axis.  It is known \cite{tit2} that $(-1)^{n-1}Z(G_n) $ has mean value $2$ in the sense that
$$\lim_{N \to \infty} \frac{1}{N}\sum_{n = 1}^N (-1)^{n-1}Z(G_n)   = 2.$$
One says that the gram point $G_n$ is  {\bf good}\index{good Gram point} if $(-1)^{n-1}Z(G_n) >0$.   Gram points are useful for finding zeros of $\zeta(s)$ on the critical line because, by continuity, the function $Z(x)$ necessarily has a zero between any two consecutive Gram points that are both good or both not good.  The sequence of values of $n$ for which the gram point $G_n$ is not good is $$127, \ 135, \ 196, \ 212, \ 233, \ 255, \ 289, \ 368, \ 378, \ 380, \ 398, \ 401,\ldots,$$ which is OEIS Sequence A114856 offset by $+1$.  It is known that the upper density (defined as an appropriate lim sup) of the set of all good Gram points within the set of all Gram points is less than $1$ \cite{trudg}.  However, one expects that ``most'' Gram points are good, e.g., the lower density (defined as an appropriate lim inf) is likely  greater than $1/2$.  By (\ref{atheta}), one has
\begin{align*}
A(G_n) & = 1-2\left\{ \frac{n+I(G_n)}{2}  \right\} \\
& =  \left.
  \begin{cases} 0 & \text{if $n$ is odd and $Z(G_n) \geq 0$, or if $n$ is even and $Z(G_n) < 0$}  \\
   1   & \text{if $n$ is even and $Z(G_n) \geq 0$, or if $n$ is odd and $Z(G_n) < 0$} 
 \end{cases}
 \right. \\
 & =  \left.
  \begin{cases} 0 & \text{if $G_n$ is good}  \\
   1   & \text{otherwise},
 \end{cases}
 \right.
\end{align*}
where the last equality holds provided that no $G_n$ is a zero of $Z(x)$.  This suggests that the rule
$$E(n) =  \left.
  \begin{cases}
   0 & \text{if $G_n$ is good}  \\
   1   & \text{otherwise}.
 \end{cases}
\right.$$
for which the author could find no exceptions, holds for `nearly all' $n$.    Note, however, that the function $E(x)$ does assume the value $2$, for example, at  $x = 7299022.0139$, and the apparent growing density of values of $x$ for which $E(x) = 2$ as  $x \to \infty$ suggests that the rule for  $E(n)$ above may have some exceptions.

\begin{remark}[Normalizations using the Riemann--Siegel theta function]\label{altnorm}
An alternative normalization $\widecheck{\gamma}_n$ of the $\gamma_n$, employed in \cite{reyna}, is $$\widecheck{\gamma}_n = 1+\frac{1}{\pi}\theta(\gamma_n)  = n-\frac{1}{2}-S(\gamma_n) = \widehat{\gamma}_n  -\frac{1}{2}+ \frac{1+o(1)}{48 \pi \gamma_n}\ (n \to \infty).$$  The corresponding normalizations $\widecheck{S}(x)$, $\widecheck{N}(x)$, $\widecheck{Z}(x)$, and $\widecheck{A}(x)$  of the functions $S(T)$, $N(T)$, $Z(T)$, and $A(T)$ compose the latter functions, respectively, with the  inverse $\theta^{-1}((x-1)\pi)$ of the function $1+\frac{1}{\pi}\theta(x)$,  rather than the inverse $2\pi\frac{x-7/8}{W((x-7/8)/e)}$ of the function $\frac{x}{2\pi}\log\frac{x}{2\pi e}+ \frac{7}{8}$.   (Since $\theta(x)$ is increasing for $x \geq 6.289835\ldots$, with $\theta(6.289835\ldots)  = -3.530972\ldots <-\pi$, the function $\theta^{-1}((x-1)\pi)$ is well defined for $x \geq 0$.)  Thus, for example, one has
$$ \widecheck{N}(x) = x + \widecheck{S}(x),$$
$$\widecheck{Z}(x) = -e^{i\pi x}\zeta(\tfrac{1}{2}+ i\theta^{-1}((x-1)\pi)),$$
$$\widecheck{A}(x) = 1-2\left\{\frac{x+ I(\theta^{-1}((x-1)\pi))}{2} \right\},$$
and
$$\widecheck{A}(n) = A(G_n) = 1-2\left\{ \frac{n+I(G_n)}{2}  \right\}.$$
One can prove analogues of our results with respect to these alternative normalizations, which, admittedly,  are more natural, but not as readily computable.
\end{remark}

Section 11.2 provides further analysis of the functions $N(T)$, $S(T)$, $R(T)$, $\tau_n$, $\tau_{n+1}-\tau_n$, $\widehat{\gamma}_n$, and $\widehat{\gamma}_{n+1}-\widehat{\gamma}_n$.

\section{Primes in abstract analytic number theory}

A complex number is said to be {\bf algebraic}\index{algebraic number} if it is a root of a nonzero polynomial with rational coefficients.  It is well known that the set $\overline{\QQ}$\index[symbols]{.a rrrr@$\overline{\QQ}$} of all algebraic numbers is a countable subfield of the field $\CC$.  A complex number is said to be an {\bf algebraic integer}\index{algebraic integer}\index{algebraic integer} if it is a root of a monic polynomial with integer coefficients.  It is also well known that the set $\mathcal{O}_\CC$\index[symbols]{.a rrrrr@$\mathcal{O}_\CC$} of all algebraic integers is a subring of the field $\overline{\QQ}$ with $\mathcal{O}_\CC \cap \QQ  = \ZZ$.  Moreover, every algebraic number is of the form $\alpha/n$ for some $\alpha \in \mathcal{O}_\CC$ and some positive integer $n$.  In particular, $\overline{\QQ}$ is the field of fractions of the integral domain $\mathcal{O}_\CC$.  It is known that the integral domain  $\mathcal{O}_\CC$ is a {\bf B\'ezout domain},\index{Bezout @B\'ezout domain} that is, all of its finitely generated ideals are principal.  {\it Algebraic number theory} is largely the study of the field $\overline{\QQ}$ of all algebraic numbers and its various subrings and subfields.

For any subfield $K$ of $\CC$, one lets $\mathcal{O}_K = \mathcal{O}_\CC \cap K$\index[symbols]{.a rrrrs@$\mathcal{O}_K$} denote the ring of all algebraic integers in $K$.\index{ring $\mathcal{O}_K$ of all algebraic integers in $K$}   For example, one has $\mathcal{O}_\QQ = \ZZ$ and $\mathcal{O}_{\QQ(i)} = \ZZ[i]$.  If $K$ is finite dimensional as a vector space of $\QQ$, then $K$ is said to be a {\bf number field}, and $\mathcal{O}_K$ is said to be a {\bf number ring}.\index{number field}\index{number ring}   Every number field $K$ is of the form $\QQ(\alpha)$ for some algebraic number $\alpha$, and conversely every such field is a number field.  Moreover,  every number ring $\mathcal{O}_K$ is of the form $\ZZ[\alpha_1, \alpha_2,\ldots,\alpha_n]$ for some algebraic integers $\alpha_1, \alpha_2, \ldots, \alpha_n$,  but not every such ring is a number ring.  For example, for $K = \QQ(\sqrt{5})$, one has $\mathcal{O}_K = \ZZ[\frac{1+\sqrt{5}}{2}]$, and $\ZZ[\sqrt{5}]$ is not a number ring.  It is known that, for any number field $K$, the  number ring $\mathcal{O}_K$ is a {\it Dedekind domain} with finite {\it ideal class group}.  The {\bf norm  $N_{K/\QQ}(\aaa)$ of $\aaa$}\index{norm $N_{K/\QQ}(\aaa)$ of a nonzero ideal $\aaa$ of $\mathcal{O}_K$} for any nonzero ideal $\aaa$ of $\mathcal{O}_K$ is simply the cardinality $\# (\mathcal{O}_K/\aaa)$ of the quotient ring   $\mathcal{O}_K/\aaa$, which is known to be finite for all nonzero ideals $\aaa$.   Moreover, one has
$$N_{K/\QQ}(\aaa \bbb)  = N_{K/\QQ}(\aaa )N_{K/\QQ}(\bbb)$$
for all nonzero ideals $\aaa$ and $\bbb$ of $\mathcal{O}_K$.   For example, one has $N_{\QQ/\QQ} ((n))= |n|$ for every nonzero (principal) ideal $(n)$ of $\ZZ$.

For all $x > 0$, let $\pi_K(x)$ denote the number of prime ideals of ${\mathcal O}_K$ of norm less than or equal to $x$.  The {\bf Landau prime ideal theorem},\index{Landau prime ideal theorem} proved by Landau in 1903 \cite{land1}, states that 
$$\pi_K(x) = \li(x) + O\left(xe^{-c_K {\sqrt  {\log x }}}\right)\ (x \to \infty)$$
for some constant $c_K > 0$ depending on $K$ (where the $O$ constant also depends on $K$).  Landau's  theorem implies that
\begin{align}\label{landauPNT}
\pi_K(x) = \li(x) + O(x(\log x)^{-t}) \ (x \to \infty)
\end{align}
for all $t > 0$.  

 An {\bf arithmetic semigroup}\index{arithmetic semigroup} is a commutative monoid $G$ (written multiplicatively) that is freely generated by a countable set $P$, called the set of {\bf primes} of $G$, and that is equipped with a map $|\!-\!|: G \longrightarrow \RR_{>0}$, called the {\bf norm} map of $G$, such that the following conditions hold.
\begin{enumerate}
\item  $|ab| = |a||b|$ for all $a,b \in G$. 
\item $|p| > 1$ for all $p \in P$.
\item  For all $x > 0$, the number $$N_G(x) = \#\{a \in G: |a| \leq x\}$$ of elements of $G$ of norm at most $x$ is finite.
\end{enumerate}
In fact, condition (3) can be replaced with
\begin{enumerate}
\item[$(3')$]  For all $x > 0$, the number $$\pi_G(x) = \#\{p \in P: |p| \leq x\}$$ of primes of $G$ of norm at most $x$ is finite.
\end{enumerate} 
For all $x \geq 1$, we let $$\PP_G(x) = \frac{\pi_G(x)}{N_G(x)}$$ denote the probability that a randomly selected element of $G$ of norm less than or equal to $x$ is prime.   The   {\bf zeta function of $G$} is the complex function $$\zeta_G(s) = \sum_{a \in G} \frac{1}{|a|^{s}}.\index{zeta function $\zeta_G(s)$ of an arithmetic semigroup $G$}$$

\begin{example}\
\begin{enumerate}
\item  The commutative monoid $\ZZ_{>0}$ of all positive integers under multiplication is freely generated by the prime numbers, and in fact $\ZZ_{>0}$  is an arithmetic semigroup, with the norm map acting by $n \longmapsto n$,, with $\pi_{\ZZ_{>0}}(x) = \pi(x)$ and $N_{\ZZ_{>0}}(x) = \lfloor x \rfloor$ for all $ x > 0$, and with $\zeta_{\ZZ_{>0}}(s) = \zeta(s)$.
\item   Let $K$ be a number field.    Let $G_K$ be the commutative monoid, under ideal multiplication, of all nonzero ideals of $\mathcal{O}_K$.  Then $G_K$ is a freely generated by the prime ideals of $\mathcal{O}_K$, and in fact $G_K$ is arithmetic semigroup, with the norm map acting by $\mathfrak{a} \longmapsto N_{K/\QQ}(\mathfrak{a}) = \#( \mathcal{O}_K/\mathfrak{a})$,  and with $\pi_{G_K}(x) = \pi_K(x)$ for all $x >0$.  Moreover, $\PP_{G_K}(x)$ for all $x \geq 1$ is the probability that a randomly selected ideal of ${\mathcal O}_K$ of norm less than or equal to $x$ is prime, and $\zeta_{G_K}(s)$ is 
the {\bf Dedekind zeta function}\index{Dedekind zeta function $\zeta_K(s)$ of a number field $K$}   $$\zeta_K(s) = \sum_{\mathfrak{a} \in {G_K}} \frac{1}{(N_{K/\QQ}(\mathfrak{a}))^s}$$ {\bf of $K$}.
  \end{enumerate}
\end{example}

 The following is one of several {\it abstract prime number theorems}.  

\begin{theorem}[{\cite[Theorem 1]{deb} \cite{ny}}]\label{abstractPNT}\index{abstract prime number theorem}
 Let $G$ be any  arithmetic semigroup.  
For any  $\delta > 0$, one has
\begin{align*}
\pi_G(x) = \li(x^\delta) + O(x^\delta (\log x)^{-t}) \ (x \to \infty)
\end{align*}
for all $t > 0$ if and only if there exists an $A >0$ such that
\begin{align}\label{Nhyp}
N_G(x) = Ax^\delta+ O(x^\delta (\log x)^{-t}) \ (x \to \infty)
\end{align}
for all $t>0$.  Moreover, if the conditions above hold, then 
one has  $$A \delta = \operatorname{Res}_{s = \delta} \zeta_G(s) = \lim_{\operatorname{Re}{s} \to \delta^+} (s-\delta) \zeta_G(s),$$ 
where the series defining $\zeta_G(s) = \sum_{a \in G} \frac{1}{|a|^{s}}$ is absolutely convergent for $\operatorname{Re} s > \delta$ and has meromorphic continuation to $\operatorname{Re} s  \geq \delta$ with a unique (simple) pole at $s = \delta$.  
\end{theorem}

 Let $G$ be any  arithmetic semigroup.   For any $\delta > 0$, we let $G_\delta$ denote the commutative monoid $G$ together with norm map $|\!-\!|^{\delta}: G \longrightarrow \RR_{>0}$, where $|\!-\!|$ is the norm map of the arithmetic semigroup $G$.  It is clear that $G_\delta$ is also an arithmetic semigroup, and one has
$$N_{G_\delta}(x) =  \#\{a \in G: |a|^\delta \leq x\} =  \#\{a \in G: |a| \leq x^{1/\delta}\} = N_G(x^{1/\delta}),$$
and likewise $\pi_{G_\delta}(x)  =  \pi_G(x^{1/\delta}),$
for all $x > 0$.    This simple trick allows one to reduce the proof of Theorem \ref{abstractPNT} to the case where $\delta = 1$, which is assumed throughout \cite{deb} and \cite{ny}.

An arithmetic semigroup $G$ satisfies {\bf Axiom A}  \cite[Chapter 4]{kno}\index{Axiom A}  if there exist positive constants $A$ and $\delta$ and a nonnegative constant $\eta < \delta$ such that  $$N_G(x) = Ax^\delta+O(x^\eta) \ (x \to \infty).$$   If $G$ is an arithmetic semigroup satisfying Axiom A, then clearly the hypothesis (\ref{Nhyp}) of Theorem \ref{abstractPNT} holds, and, moreover,  by \cite[Chapter 4 Proposition 2.8 and Corollary 2.14]{kno}, the zeta function $\zeta_G(s)$ has analytic continuation to   $\{s \in \CC: \operatorname{Re} s > \eta, \, s \neq \delta\}$, where one has
$$\sum_{a \in G:\, |a| \leq x} \frac{1}{|a|^\delta} = A\delta \log x + \gamma_G + O(x^{\eta-\delta}) \ (x \to \infty),$$
where $$\gamma_G = \lim_{s \to \delta} \left(\zeta_G(s)-\frac{A\delta}{s-\delta}\right) = A+ \delta\int_1^\infty \frac{N_G(t)-At^\delta}{t^\delta} \, \frac{dt}{t}.$$
  \cite[Chapters 4 and 5]{kno} provides a wide variety of natural examples of arithmetic semigroups satisfying Axiom A, including the following defined for any number field.

\begin{example}[{\cite{kno}}]\label{axa}  Let $K$ be a number field.
\begin{enumerate}
\item Let ${G_K}$ be the monoid, under ideal multiplication, of all nonzero ideals of $\mathcal{O}_K$.  Then ${G_K}$ is an arithmetic semigroup satisfying Axiom A, with the norm map acting by $\mathfrak{a} \longmapsto N_{K/\QQ}(\mathfrak{a}) = \#( \mathcal{O}_K/\mathfrak{a})$, with $\delta = 1$, and with $$A = \operatorname{Res}_{s = 1} \zeta_K(s).$$  For example, if $K  = \QQ(i)$, then $A = \pi/4$, and if $K = \QQ(\sqrt{2})$, then $A =\frac{1}{\sqrt{2}}\log(1+\sqrt{2})$.  Also, if $K = \QQ$, then $\zeta_{G_K}(s) = \zeta(s)$ is the Riemann zeta function, $A = 1$, and $\PP_{G_K}(x) =  \PP(\lfloor x \rfloor)$ for all $x \geq 1$, where $\PP(x) = \frac{\pi(x)}{x}$.
\item Let $G$ be the monoid,  under direct product, of all isomorphism classes of $\mathcal{O}_K$-modules of finite cardinality.  Then $G$ is an arithmetic semigroup satisfying Axiom A, with the norm map acting by $M \longmapsto \#M$, and with $\delta = 1$ and $$\zeta_G(s) = \prod_{n =1}^\infty \zeta_K(ns),$$ so that $$A = \operatorname{Res}_{s = 1} \zeta_K(s) \prod_{n = 2}^\infty \zeta_K(n).$$  Moreover, $\PP_G(x)$ for all $x \geq 1$ is the probability that an $\mathcal{O}_K$-module of cardinality less than or equal to $x$, randomly selected among a set of representatives of the isomorphism classes, is {\it indecomposable}.  Note that, if $K = \QQ$, then $\mathcal{O}_K = \ZZ$ and an $\mathcal{O}_K$-module is equivalently an abelian group.
\end{enumerate}
\end{example}

Note that Theorem \ref{abstractPNT} and Example \ref{axa}(1) together  yield the weak version (\ref{landauPNT}) of  Landau's prine number theorem.

{\it Abstract analytic number theory} is the study of arithmetic semigroups and various natural axioms that one can impose on them.    The many  results and examples of abstract analytic number theory cut across all fields of mathematics, even topology and category theory \cite{kno}.    Although much of analytic number theory beyond the prime number theorems with error term, e.g., the Riemann--von Mangoldt explicit formulas and the Riemann hypothesis, have no known generalizations to abstract analytic number theory,  the results that do so generalize apply to a wide variety of examples that lie well beyond the scope of analytic number theory.

\begin{remark}[Zeta functions]\label{zetafunctions}
``Zeta functions'' are defined for many objects in mathematics besides number fields and arithmetic semigroups, e.g., for special types of the following classes of objects: partially ordered sets, graphs, categories, schemes over $\ZZ$,  algebraic varieties over a field, compact Lie groups, compact Riemannian manifolds, dynamical systems, quadratic forms,  and error-correcting codes. The Riemann zeta function is an important special case of many of these constructions.   See \cite{lag3} \cite{li} for excellent surveys of zeta functions arising in number theory. 
\end{remark}

\part{Algebraic asymptotic analysis}

\chapter{Logexponential degree}

In this chapter,  we introduce a generalization of $\deg f$ that we call the {\it (iterated) logarithmic degree} $\degl f$ of $f$, along with a generalization of $\degl f$ that we call the {\it logexponential degree} $\dege f$ of $f$, in Sections 6.2 and 6.3, respectively.  We prove many properties of these resepective formalisms, including how they relate to the relations $O$, $o$, $\gg$, $\asymp$, and $\sim$ on functions and how they relate to the basic operations $+$, $-$, $\cdot$, $/$, and $\circ$ on functions.   The definitions and results in this chapter, especially those in Sections 6.1 and  6.3,  form the foundation for Parts 2 and 3.

\section{The degree of a real function}

Let $f \in \RR^{\RR_\infty}$, that is, let $f$ be a real function whose domain is a subset of $\RR$ that is not bounded above.   Recall from Section 2.3 that we define the {\bf (upper) degree} of $f$ to be the extended real number
$$\deg f = \inf\{t \in \RR: f(x) = o(x^t) \ (x \to \infty)\} \in \overline{\RR}.$$
See Example \ref{degreeexamples} for some simple examples from calculus and Example \ref{divisorf} for some important number-theoretic examples.
Also, see Proposition \ref{limsupdef} for a proof of the alternative expression 
$$\deg f = \limsup_{x \to \infty} \frac{\log |f(x)|}{\log x}$$
for the degree of $f$, which has as a corollary that
$$\limsup_{ x \to \infty} f(x) = \deg x^{f(x)}$$
and
$$\liminf_{ x \to \infty} f(x) = -\limsup_{x \to \infty}\, (-f(x)) =- \deg x^{-f(x)}.$$
In this section, we provide some useful properties of the degree formalism.  All of the proofs are fairly straightforward and are left to the reader.

Let  $f  \in \RR^{\RR_\infty}$.   Recall that the  {\bf lower degree} $\underline{\deg} \, f$  of $f$ is defined by
$$\underline{\deg} \, f = \liminf_{x \to \infty} \frac{\log |f(x)|}{\log x},$$  and the  {\bf exact degree} $\overline{\underline{\deg}} \, f$ of $f$ is the limit
$$\overline{\underline{\deg}}\, f= \lim_{x \to \infty} \frac{\log |f(x)|}{\log x},$$
provided that the limit exists or is $\pm \infty$, in which case we say that $f$ has {\bf exact degree}.   Note that if $\deg f = -\infty$ or if $\underline{\deg}\, f = \infty$, then $f$ has exact degree.   Moreover, if $f$ is not eventually nonzero on its domain, then $\underline{\deg}\, f = -\infty$, and thus  $f$ has exact degree if and only if $\deg f = -\infty$.    

The following proposition collects several  properties and equivalent characterizations of degree.

\begin{proposition}\label{firstprop0}
For all $f \in \RR^{\RR_\infty}$, one has the following.
\begin{enumerate}
\item $\deg f$ is the unique $d \in \overline{\RR}$ such that $f(x) = o(x^{t}) \ (x \to \infty)$  for all $t > d$ but for no $t < d$.
\item $\deg f$ is the unique $d \in \overline{\RR}$ such that $f(x) = O(x^{t}) \ (x \to \infty)$  for all $t > d$ but for no $t < d$.
\item One has
\begin{align*}
\deg f  &   = \inf \{t \in \RR: f(x) = o(x^t) \ (x \to \infty)\} \\
  & = \inf \{t \in \RR: f(x) = O(x^t) \ (x \to \infty)\} \\
  & = \sup \{t \in \RR: f(x) \neq o(x^t) \ (x \to \infty)\} \\
  & = \sup \{t \in \RR: f(x) \neq O(x^t) \ (x \to \infty)\}.
\end{align*}
\item For all $t \in \RR$, one has $\deg f  \leq t$
if and only if  $f(x) = o\left(x^s \right) \ (x \to \infty)$
for all $s > t$, if and only if 
$f(x) = O\left(x^s \right) \ (x \to \infty)$
for all $s > t$.
\item For all $t \in \RR$, one has $\deg f  < t$
if and only if  $f(x) = o\left(x^s \right) \ (x \to \infty)$
for some $s < t$, if and only if 
$f(x) = O\left(x^s \right) \ (x \to \infty)$
for some $s < t$.
\item $\deg f = \infty$ if and only if $f(x) \neq o(x^t) \ (x \to \infty)$ for all $t \in \RR$.
\item $\deg f = -\infty$ if and only if $f(x) = o(x^t) \ (x \to \infty)$ for all $t \in \RR$.
\item $\deg f  \neq \pm \infty$ if and only if $f(x) = x^d F(x)$, for all $x \gg 0$ in $\dom f$,  for some $d \in \RR$ and for some real function $F$ with $\deg F = 0$ (and, moreover, one has $d = \deg f$).
\item $f$ has finite exact degreee $d \in \RR$ if and only if $f (x) = x^{d+o(1)} \ (x \to \infty)$.
\item $\deg f = \deg |f| = \deg Cf$ for all $C \neq 0$.
\item $\overline{\underline{\deg }} \, 0  = -\infty$ and $\overline{\underline{\deg }}\, 1 = 0$.
\item If $f$ is eventually bounded (i.e., if $\limsup_{x \to \infty}|f(x)| < \infty$), then $\deg f \leq 0$.
\item If $\deg f < 0$, then $\displaystyle \lim_{x \to \infty} f(x) = 0$.
\item $\deg af(x) = \deg f(ax) = \deg f(x)$ and $\deg f(x^a) = \deg f(x)^a = a \deg f(x)$  for all $a > 0$.
\item If $f(x)$ is differentiable and nonzero for all $x \gg 0$ in $\dom f$, then $$\displaystyle \liminf_{x \to \infty} \frac{x f'(x)}{f(x)} \leq \underline{\deg}\, f  \leq \deg f \leq \limsup_{x \to \infty} \frac{x f'(x)}{f(x)}.$$
\item Suppose that $\displaystyle L = \lim_{x \to \infty} \frac{x f'(x)}{f(x)}$ on $\dom f$ exists or is $\pm \infty$.   Then $f$ has exact degree $L$.  Moreover, if $L \neq 0,\pm \infty$ then $\deg f' = -1+\deg f$; if $L = 0$ then $\deg f' \leq -1$; and if $L =   \infty$ then $\deg f' =\infty$.
\item If $f  = \sum_{a \in \RR} r_a x^a \in \RR[x^{a}: a \in \RR]$ is a polynomial in the power functions $x^a$ for $a \in \RR$, then $\overline{\underline{\deg }}\, f = \max\{a \in \RR: r_a \neq 0\}$.
\item If $f  = r/s \in \RR(x^{a}: a \in \RR)$, where $r,s \in \RR[x^{a}: a \in \RR]$ and $s \neq 0$, then $\overline{\underline{\deg }}\, f  = \overline{\underline{\deg }} \, r -\overline{\underline{\deg }}\, s$.
\end{enumerate}
\end{proposition}

The following proposition shows how degree relates to various  relations on functions.    (Recall that by convention we assume that each of the asymptotic relations between $f$ and $g$ in the proposition implies that $\dom g$ contains the intersection of $\dom f$ with $[N,\infty)$ for some real number $N$.)

\begin{proposition}\label{firstprop1}
Let $f, g \in \RR^{\RR_\infty}$, and let $\underline{g} = g|_{\dom f \cap \dom g}$.
\begin{enumerate}
\item If $f(x) = O(g(x)) \ (x \to \infty)$, then $\deg f \leq \deg  \underline{g} \leq \deg g$.
\item If $f(x) \gg g(x) \ (x \to \infty)$, then $\deg f \geq \deg \underline{g}$. 
\item If $f(x) \asymp g(x) \ (x \to \infty)$, then $\deg f = \deg \underline{g}$, and $f$ has exact degree if and only if $\underline{g}$ has exact degree. 
\item  If $f(x) \asymp g(x) \ (x \to \infty)$ and $[N,\infty) \cap \dom f = [N , \infty) \cap \dom g$ for some  $N > 0$, then $\deg f = \deg g$. 
\item  $\deg f \leq  \deg g$ if and only if for all $t$ such that $g(x) = o\left(x^t \right) \ (x \to \infty)$ (or $g(x) = O\left(x^t \right) \ (x \to \infty)$) one has $f(x) = o\left(x^s \right) \ (x \to \infty)$  (or $f(x) = O\left(x^t \right) \ (x \to \infty)$) for all $s > t$.
\item   $\deg f <  \deg g$ if and only if there exist $s$ and $t$ with $s >t$ such that $f(x) = o\left(x^t \right) \ (x \to \infty)$ (or $f(x) = O\left(x^t \right) \ (x \to \infty)$) and $g(x) \neq  o\left(x^s \right) \ (x \to \infty)$ (or $g(x) \neq  O\left(x^s \right) \ (x \to \infty)$).
\end{enumerate}
\end{proposition}

The next proposition shows how degree relates to the  operations of addition, multiplication, and division of functions.

\begin{proposition}\label{firstprop2}
Let $f$ and $g$ be real functions such that both $\underline{f} = f|_{\dom f \cap \dom g}$ and $\underline{g} = g|_{\dom f \cap \dom g}$ are in $\RR^{\RR_\infty}$.
\begin{enumerate}
\item One has $\deg(f+ g) \leq \max(\deg \underline{f},\deg \underline{g})$, with equality if $\deg \underline{f} \neq \deg \underline{g}$.
\item One has $\deg(fg) \leq \deg \underline{f}+\deg \underline{g}$ provided that the sum of degrees is defined in $\overline{\RR}$,
with equality if also  either $\underline{f}$ or $\underline{g}$ has exact degree.
\item One has $\deg(f/g) \leq \deg \underline{f}+\deg(1/\underline{g})$ provided that  $\underline{g}$ is eventually nonzero and the sum of degrees is defined in $\overline{\RR}$, with equality if also  either $\underline{f}$ or $\underline{g}$ has exact degree.
\item One has $\deg(f/g) \geq \deg \underline{f}-\deg \underline{g}$  (resp., $\deg(fg) \geq \deg \underline{f} -\deg(1/\underline{g})$) provided that  $\underline{g}$ is eventually nonzero and the difference  of degrees is defined in $\overline{\RR}$, with equality if  also $\underline{g}$ has exact degree.
\item  If $g$ is eventually nonzero, then $\underline{\deg} \, g  = -\deg(1/g)$, while,  if $g$ is not eventually nonzero, then $\underline{\deg} \, g = -\infty$.
\item  $g$ has exact degree if and only if  $\deg g =-\infty$ or else $g$ is eventually nonzero and $\deg(1/g) = -\deg g$.
\item If $\deg \underline{f} < \deg \underline{g}$ and $\underline{g}$ has exact degree,  or more generally if  $\deg \underline{f} < \underline{\deg} \, \underline{g}$, then  $\deg(f/g) < 0$ and therefore $\underline{f}(x) = o(\underline{g}(x)) \ (x \to \infty)$.
\item If $\deg (f-g) <  \max(\deg \underline{f},\deg \underline{g})$ and  $\underline{g}$ has exact degree, or more generally if $\deg(f-g) < \underline{\deg} \, \underline{g}$, then $\deg\left(f/g-1 \right) < 0$ and therefore $\underline{f}(x) \sim \underline{g}(x) \ (x \to \infty)$.
\item If $L = \displaystyle \lim_{x \to \infty} \frac{\log |f(x)|}{\log |g(x)|} \geq 0$ exists and $L > 0$ or $\deg \underline{g}$ is finite, then $\deg \underline{f} = L \deg \underline{g}$.
\end{enumerate}
\end{proposition}

\begin{example}
Since $\deg \sin x= 0$ and $\underline{\deg} \, \sin x= -\infty$, the  function $\sin x$ does not have exact degree.  Moreover, one has $\deg (1/x) = -1 < 0 = \deg \sin x$, and yet $1/x \neq o(\sin x) \ (x \to \infty)$ and in fact $1/x \neq O(\sin x) \ (x \to \infty)$.
\end{example}

Let $a \in \overline{\RR}$ and  $f, g \in \RR^{\RR_a}$,  where $\dom g$ contains the intersection of $\dom f$ with some punctured neighborhood of $a$.   Let us write
$f(x)\ggg g(x) \ (x \to a)$ if for all  $M > 0$ one has $|f(x)| \geq M|g(x)|$ for all $x$ in the intersection of $\dom f$ with some punctured neighborhood of $a$.  Supposing that $\dom g$ contains the intersection of $\dom f$ with some punctured neighborhood of $a$,  the condition $f(x)\ggg g(x) \ (x \to a)$ is equivalent to $g|_{\dom f}(x) = o(f(x)) \ (x \to a)$, just as the condition  $f(x) \gg g(x) \ (x\to a)$ is equivalent to $g|_{\dom f}(x) = O(f(x)) \ (x \to a)$.

The following are analogues, for lower degree, of the results above.

\begin{proposition}\label{firstprop0lower}
For all $f \in \RR^{\RR_\infty}$, one has the following.
\begin{enumerate}
\item $\underline{\deg} \, f$ is the unique $d \in \overline{\RR}$ such that $f(x) \ggg x^t \ (x \to \infty)$  for all $t < d$ but for no $t > d$.
\item $\underline{\deg} \, f$ is the unique $d \in \overline{\RR}$ such that $f(x) \gg x^t \ (x \to \infty)$  for all $t < d$ but for no $t > d$.
\item One has
\begin{align*}
\underline{\deg} \, f  &   = \sup \{t \in \RR: f(x) \ggg x^t  \ (x \to \infty)\} \\
  & = \sup \{t \in \RR: f(x) \gg x^t \ (x \to \infty)\} \\
  & = \inf \{t \in \RR: f(x) \not \ggg  x^t  \ (x \to \infty)\} \\
  & = \inf \{t \in \RR: f(x) \not \gg x^t  \ (x \to \infty)\}.
\end{align*}
\item For all $t \in \RR$, one has $\underline{\deg} \, f  \geq t$
if and only if  $f(x) \ggg x^s \ (x \to \infty)$
for all $s < t$, if and only if 
$f(x) \gg x^s \ (x \to \infty)$
for all $s < t$.
\item For all $t \in \RR$, one has $\underline{\deg} \, f  > t$
if and only if  $f(x) \ggg x^s \ (x \to \infty)$
for some $s > t$, if and only if 
$f(x) \gg x^s (x \to \infty)$
for some $s > t$.
\item $\underline{\deg} \, f = \infty$ if and only if $f(x) \gg x^t  \ (x \to \infty)$ for all $t \in \RR$.
\item $\underline{\deg} \, f = -\infty$ if and only if $f(x) \not \gg x^t \ (x \to \infty)$ for all $t \in \RR$.
\item $\underline{\deg} \, f  \neq \pm \infty$ if and only if $f(x) = x^d F(x)$, for all $x \gg 0$ in $\dom f$,  for some $d \in \RR$ and for some real function $F$ with $\underline{\deg} \, F = 0$ (and, moreover, one has $d = \underline{\deg} \, f$).
\item $\underline{\deg} \, f = \underline{\deg} \, |f| = \underline{\deg} \, Cf$ for all $C \neq 0$.
\item If $|f|$ is eventually bounded away from $0$ (i.e., if $\liminf_{x \to \infty} |f(x) | >0)$, then $\underline{\deg} \, f \geq 0$.
\item If $\underline{\deg} \, f > 0$, then $\displaystyle \lim_{x \to \infty} |f(x) |= \infty$.
\item $\underline{\deg} \, af(x) = \underline{\deg} \, f(ax) = \underline{\deg} \, f(x)$ and $\underline{\deg} \, f(x^a) = \underline{\deg} \, f(x)^a = a \, \underline{\deg} \, f(x)$  for all $a > 0$.
\end{enumerate}
\end{proposition}

\begin{proposition}\label{firstprop1lower}
Let $f, g \in \RR^{\RR_\infty}$, and let $\underline{g} = g|_{\dom f \cap \dom g}$.
\begin{enumerate}
\item If $f(x) = O(g(x)) \ (x \to \infty)$, then $\underline{\deg} \, f \leq \underline{\deg} \,\underline{g}$.
\item If $f(x) \gg g(x) \ (x \to \infty)$, then $\underline{\deg} \, f \geq \underline{\deg} \, \underline{g} \geq \underline{\deg}\, g$. 
\item If $f(x) \asymp g(x) \ (x \to \infty)$, then $\underline{\deg} \, f = \underline{\deg} \, \underline{g}$.
\item  If $f(x) \asymp g(x) \ (x \to \infty)$ and $[N,\infty) \cap \dom f = [N \cap \infty) \cap \dom g$ for some  $N > 0$, then $\underline{\deg} \, f = \underline{\deg} \, g$. 
\item  $\underline{\deg} \, f \leq  \underline{\deg} \, g$ if and only if for all $t$ such that $f(x) \ggg x^t  \ (x \to \infty) \ (x \to \infty)$ (or $f(x) \gg x^t  \ (x \to \infty)$) one has $g(x) \ggg x^s  \ (x \to \infty)$  (or $g(x) \gg x^s  \ (x \to \infty)$) for all $s < t$.
\item   $\underline{\deg} \, f <  \underline{\deg} \, g$ if and only if there exist $s$ and $t$ with $s <t$ such that $g(x) \ggg x^t  \ (x \to \infty)$ (or $g(x) \gg x^t  \ (x \to \infty)$) and $f(x) \not \ggg x^s  \ (x \to \infty)$ (or $f(x) \not \gg x^s  \ (x \to \infty)$).
\item   $\deg \, f <  \underline{\deg} \, g$ if and only if there exist $s$ and $t$ with $s <t$ such that $g(x) \ggg x^t  \ (x \to \infty)$ (or $g(x) \gg x^t  \ (x \to \infty)$) and $f(x)\lll x^s  \ (x \to \infty)$ (or $f(x) \ll x^s  \ (x \to \infty)$).
\end{enumerate}
\end{proposition}

\begin{proposition}\label{firstprop2lower}
Let $f$ and $g$ be real functions such that both $\underline{f} = f|_{\dom f \cap \dom g}$ and $\underline{g} = g|_{\dom f \cap \dom g}$ are in $\RR^{\RR_\infty}$.
\begin{enumerate}
\item One has $\underline{\deg} \,(fg) \geq \underline{\deg} \, \underline{f}+\underline{\deg} \, \underline{g}$ provided that the sum of  lower degrees is defined in $\overline{\RR}$,
with equality if also  either $\underline{f}$ or $\underline{g}$ has exact degree.
\item One has $\underline{\deg} \,(f/g) \geq \underline{\deg} \, \underline{f}+\underline{\deg} \,(1/\underline{g})$ provided that  $\underline{g}$ is eventually nonzero and the sum of degrees is defined in $\overline{\RR}$, with equality if also  either $\underline{f}$ or $\underline{g}$ has exact degree.
\item One has $\underline{\deg} \,(f/g) \leq \underline{\deg} \, \underline{f}-\underline{\deg} \, \underline{g}$  (resp., $\underline{\deg} \,(fg) \leq \underline{\deg} \, \underline{f} -\underline{\deg} \,(1/\underline{g})$) provided that  $\underline{g}$ is eventually nonzero and the difference  of lower degrees is defined in $\overline{\RR}$, with equality if  also $\underline{g}$ has exact degree.
\end{enumerate}
\end{proposition}

Because of the corresponding properties of limits superior, limits inferior, and limits,  degree, lower degree, and exact degree are interrelated as follows.

\begin{proposition}\label{fresupperlower}
Let $f \in \RR^{\RR_\infty}$.  One has
\begin{align*}
\underline{\deg}\, f & = \min\left\{\deg f|_X:  X \subseteq \dom f \text{ and }  \sup X = \infty \right\} \\
  & =  \min\left\{\overline{\underline{\deg}} \, f|_X:  X \subseteq \dom f \text{ and }  \overline{\underline{\deg}} \, f|_X \text{ exists} \right\}   \\
   & =  \min\left\{\overline{\underline{\deg}} \, f|_X:  X = \{x_n\} \subseteq \dom f,  \, x_1 < x_2 < \cdots, \, \text{and }  \overline{\underline{\deg}} \, f|_X \text{ exists} \right\} 
\end{align*}
and
\begin{align*}
\deg f & = \max\left\{\underline{\deg}\, f :   X \subseteq \dom f \text{ and }  \sup X = \infty \right\} \\
  & =  \max\left\{\overline{\underline{\deg}} \, f|_X:  X \subseteq \dom f \text{ and }  \overline{\underline{\deg}} \, f|_X \text{ exists} \right\} \\
     & =  \max\left\{\overline{\underline{\deg}} \, f|_X:  X = \{x_n\} \subseteq \dom f,  \, x_1 < x_2 < \cdots, \, \text{and }  \overline{\underline{\deg}} \, f|_X \text{ exists} \right\}.
\end{align*}
\end{proposition}

Note that, although degree and lower degree are dual notions, there is a slight asymmetry between them owing to the fact the absolute value function on $\RR$ has a minimum value but no maximum value:  one cannot define $\deg$ in terms of $\underline{\deg}$ as one can define $\underline{\deg}$  in terms of $\deg$.   Indeed,  although one has $\underline{\deg}\, f = -\deg(1/f)$ and $\deg f = -\underline{\deg}(1/f)$  if $f$ is eventually nonzero, the condition that $f$ not be eventually nonzero requires that $\underline{\deg} \, f= -\infty$ but does not necessitate anything about $\deg f$.

It is important to note that the degree of a function $\RR^{\RR_\infty}$, or the various  asymptotic relations $O$, $o$, $\asymp$, etc., depend only on the {\it germ}  class of each function  involved, where two functions $f,g \in \RR^{\RR_\infty}$ {\bf define the same germ at $\infty$} if they are eventually equal, that is, if there exists an $N> 0$ such that $[N,\infty) \cap \dom f = [N,\infty) \cap \dom g$ and $f(x) = g(x)$ for all $x \in  [N,\infty) \cap \dom f$.  Defining the same germ at $\infty$ is an equivalence relation on the set $\RR^{\RR_\infty}$, and the {\bf germ of $f$ at $\infty$}\index{germ of a real function at $\infty$} is defined to be the equivalence class $[f]_\infty$ of $f$, for all $f \in \RR^{\RR_\infty}$.  Obviously $\deg \, [f]_\infty$ is well-defined for all germs $[f]_\infty$  at $\infty$ of functions $f$ in $\RR^{\RR_\infty}$.  By statements  (1) and (2)  of Proposition \ref{firstprop2}, the correspondence $ f \longmapsto e^{\deg f}$ defines for any subset $X$ of $\RR$ with $\sup X = \infty$ a {\it nonarchimedean pseudonorm} on the ring of germs at $\infty$ of all real functions on $X$.  For simplicity, we sometimes blur the distinction between a function in $\RR^{\RR_\infty}$ and its germ at  $\infty$.  Effectively, this means that we sometimes treat two functions in $\RR^{\RR_\infty}$ as equal if they are eventually equal.

\begin{example} \
\begin{enumerate}
\item Let $f(x)$ be the function that on the interval $[N,N+1)$  assumes the value $e^x$ for even integers $N$ and $e^{-x}$ for odd integers $N$.   Then $\deg f = \deg \frac{1}{f} = \infty$, while $f \frac{1}{f} = 1$ and therefore $\deg (f  \frac{1}{f}) = 0$.
\item Let $f(x)$ be the function that on the interval $[N,N+1)$  assumes the value $1$ for even integers $N$ and $e^{-x}$ for odd integers $N$.   Then $\deg f = 0$ and  $\deg \frac{1}{f} = \infty$,  while $f \frac{1}{f} = 1$ and therefore $\deg (f  \frac{1}{f}) = 0$.
\item Let $f(x) = \sin(x^2)$.  Then $\deg f^{(n)} = n$ for all nonnegative integers $n$.  In particular, one has $\deg f' > \deg f$.
\end{enumerate}
\end{example}

The following proposition relates degree to the operation of composition of functions.

\begin{proposition}\label{circle1} Let $f$ and $g$ be real functions such that  $ \lim_{x \to \infty} g(x) = \infty$ (on $\dom g$) and $f \circ g \in \RR^{\RR_\infty}$, and let $\underline{f} = f|_{\dom f \cap \im g}$.  One has the following.
\begin{enumerate}
\item Suppose that $f$ has positive exact degree.  Then $\deg(f \circ g) = \deg f \deg g|_{\dom (f\circ g)}\leq \deg f \deg g$, and equality holds provided that $[N,\infty) \cap \im g \subseteq \dom f$ for some $N > 0$.
\item Suppose that $g$ has finite exact degree.  Suppose also that either $g$ has positive degree or $\underline{f}$ has finite degree.  Then $\deg(f \circ g) = \deg \underline{f} \deg g \leq \deg f \deg g$, and equality holds if $[N,\infty) \cap \dom f \subseteq \im g$ for some $N > 0$ (e.g., if $g$ is defined and continuous on some neighborhood of $\infty$) or if $f$ has exact degree.
\item Suppose that $|f(x)| \geq 1$ for all $x \gg 0$ in $\dom f \cap \im g$.  Suppose also that either $ \deg \underline{f}> 0$ or $\deg g < \infty$.  Then $\deg(f \circ g) \leq \deg \underline{f} \deg g \leq \deg f \deg g$, and the first inequality is an equality (resp., both inequalities are equalities)  if  $\underline{f}$ (resp., $f$) has exact degree.
\end{enumerate}
\end{proposition}

One defines $\deg f= \deg |f|$ for any complex-valued function $f = \operatorname{Re} f+ i \operatorname{Im}f$ of a real variable with $\sup \dom f = \infty$.  All of the definitions and results in this section generalize appropriately to this more general setting.   In particular, if $f$ is any arithmetic function, then $\deg f$ is well-defined.

\begin{remark}[{Degree at $a^+$ and $a^-$ for $a \in \RR$}]
One defines $\deg_{-\infty} f = \deg f(-x)$ if $-\infty = \inf \dom f.$
Let $a \in \RR$.  One naturally defines the {\bf degree $\deg_{a^+} f$ of $f$ at $a^+$}\index{degree $\deg$} (resp., {\bf the degree $\deg_{a^-} f$ of $f$ at $a^-$}) of any real function $f$ with $a \in \overline{(a,\infty) \cap \dom f}$ (resp.,  $a \in \overline{(-\infty,a) \cap \dom f}$)  by
 $$\deg_{a^+} f =\deg f (a+1/x)$$
and
$$\deg_{a^-} f = \deg f (a-1/x),$$
respectively.   Thus, one has
$$\deg_{a^+} f = \limsup_{x \to \infty} \frac{\log |f(a+1/x)|}{\log x} = \limsup_{x \to a^+}  \frac{\log |f(x)|}{\log \frac{1}{x-a}} =  -\liminf_{\varepsilon \to 0^+}  \frac{\log |f(a+\varepsilon)|}{\log \varepsilon}$$
and
$$\deg_{a^-} f = \limsup_{x \to \infty} \frac{\log |f(a-1/x)|}{\log x} = \limsup_{x \to a^-}  \frac{\log |f(x)|}{\log \frac{1}{a-x}} =  -\liminf_{\varepsilon  \to 0^+}  \frac{\log |f(a-\varepsilon)|}{\log \varepsilon}.$$
For all $f \in \RR^{\RR_a}$, one defines the {\bf degree $\deg_a f$ of $f$ at $a$}\index{degree $\deg_a$}\index[symbols]{.g dd@$\deg_a  f$}  to be
$$\deg_{a} f =   -\liminf_{\varepsilon  \to 0}  \frac{\log |f(a+\varepsilon)|}{\log| \varepsilon|} = \max(\deg_{a^+} f, \deg_{a^-} f).$$
\end{remark}

\begin{remark}[Upper and lower leading coefficient]\label{lc}
One can generalize the definition of the leading coefficient of a polynomial as follows.  Let $f \in \RR^{\RR_\infty}$ with $\deg f \neq \pm \infty$.   Define   $$\operatorname{lc}^+ f = \limsup_{x \to \infty} f(x)x^{-\deg f}\in \overline{\RR},\index{leading coefficient $\operatorname{lc}^\pm f$}\index[symbols]{.g d@$\operatorname{lc}^+ f$}$$
and $$\operatorname{lc}^- f = \liminf_{x \to \infty} f(x)x^{-\deg f} \in \overline{\RR},\index[symbols]{.g d@$\operatorname{lc}^- f$}$$
and also $\operatorname{lc} f = \operatorname{lc}^+ |f|$.\index[symbols]{.g d2@$\operatorname{lc} f$}  Here, ``lc'' stands for ``leading coefficient.'' For example, one has $\operatorname{lc}^+ \sin x = 1$ and $\operatorname{lc}^- \sin x = -1$.   Note that   $\operatorname{lc} \log = \infty$  and $\operatorname{lc}(1/\log) = 0$ even though $\deg \log  = \deg (1/\log) = 0$.   In general, one has   $f(x) = O(x^{\deg f}) \ (x \to \infty)$ if and only if $\operatorname{lc} f$ is finite, if and only if $\operatorname{lc}^+ f$ and $\operatorname{lc}^- f$ are both finite, and   one has $f(x) = o(x^{\deg f}) \ (x \to \infty)$ if and only if $\operatorname{lc} f = 0$, if and only if $\operatorname{lc}^+ f = \operatorname{lc}^- f = 0$.
\end{remark}

In the two remarks below, we make explicit several examples of implicit uses of degree in real and complex analysis. Specifically, we show that the degree formalism has applications to entire functions and various notions of generalized dimension.

\begin{remark}[The order and type of an entire function]
A reference for the definitions and claims made in this remark is {\cite[Chapter 2]{boas}}.  Let $f$ be a nonconstant entire function of a single complex variable, and let 
$$M_f(x) = \sup\{|f(z)|: z \in \CC, \, |z|\leq x\}, \quad \forall x \geq 0.$$ The {\bf order (at $\infty$)}\index{order of an entire function} of $f$ is defined to be 
$$\rho = \limsup_{x \to \infty} \frac{\log \log M_f(x)}{\log x} \in [0,\infty].$$
Thus, one has
$$\rho = \deg (\log M_f(x) )=  \deg (\sup\{\log |f(z)|:  z \in \CC, \, |z|\leq x\}).$$
If $\rho < \infty$, then one defines the {\bf type (at $\infty$)}\index{type of an entire function of finite order} of $f$ to be
$$\sigma =\limsup_{x\to \infty }{\frac {\log M_f(x)}{x^{\rho }}} \in [0,\infty].$$
(If $\rho = 0$, then $\sigma = \infty$.)  Thus, one has
$$\sigma = \operatorname{lc}^+ ( \log M_f(x)),$$
where $ \operatorname{lc}^+$ is defined as in Remark \ref{lc}.

Let $f(z)=\sum_{n=0}^{\infty }a_{n}z^{n}$ be an arbitrary complex power series.   The {\bf Cauchy--Hadamard theorem}\index{Cauchy--Hadamard theorem} states that the radius of convergence $R$  of the series is given by
$${\displaystyle {R}= \frac{1}{\limsup_{n\to \infty } |a_{n}|^{1/n}}} = \frac{1}{\deg n^{|a_{n}|^{1/n}}}.$$
The function $f(z)$ is entire if and only if  $R = \infty$, if and only if $\deg n^{|a_{n}|^{1/n}}  = 0$.  If $f(z)$ is entire and nonconstant, then one has
$$\rho = \limsup_{n \to \infty} \frac{n \log n}{-\log |a_n|} =  -\frac{1}{\limsup_{n \to \infty} \frac{\log |a_n|^{1/n}}{\log n}} = -\frac{1}{\deg |a_n|^{1/n}}.$$
If also $0< \rho < \infty$, i.e., $-\infty < \deg |a_n|^{1/n} < 0$, then one has
$$(e  \rho \sigma )^{1/\rho} = \limsup_{n \to \infty} n^{1/\rho} |a_n|^{1/n} = \operatorname{lc}^+ |a_n|^{1/n},$$
or, equivalently,
$$\sigma= \frac{1}{e \rho}\left(\operatorname{lc}^+ |a_n|^{1/n}\right)^\rho.$$
  More generally, for any fixed $z_0 \in \CC$ one has
$$\rho = \frac{1}{1-\deg |f^{(n)}(z_0)|^{1/n}}$$
and, if $0< \rho < \infty$, also
$$( \rho \sigma)^{1/\rho} =  e^{1-1/\rho}\operatorname{lc}^+  |f^{(n)}(z_0)|^{1/n} = e^{\deg |f^{(n)}(z_0)|^{1/n}}\operatorname{lc}^+  |f^{(n)}(z_0)|^{1/n}.$$
Thus, one can express the various formulas for the order and type of an entire function more simply if one avails oneself of the more elementary notions of degree and leading coefficient.

In the following examples,  $f(z)$ is a nonconstant entire function of order $\rho \in (0,\infty)$ and type $\sigma \in [0,\infty]$, and $c$ is a nonzero complex number and $k$ a positive integer.  
\begin{enumerate}
\item Any nonconstant polynomial is of order $0$ and type $\infty$.
\item  $f(z)$ has no zeros if and only if $\frac{1}{f(z)}$ is entire, if and only if $f(z) =e^{g(z)}$ for some polynomial $g(z)  \in \CC[z]$, in which case $\rho = \deg g$ and $\sigma$ is the absolute value of  the leading coefficient of $g$.
\item  $cf(z)^k$ is entire of order $\rho$ and type $k\sigma$.
\item $f(cz^k)$ is entire of order $k \rho$ and type $|c|^\rho \sigma$.
\item If $f(z)$ is even, then  $f(c\sqrt{z})$ is entire of order $\frac{1}{2}\rho$ and type $|c|^{\rho/2}\sigma$.
\item $f(z)+f(-z)$ is even and entire of order $\rho'\leq \rho$, and therefore $f(c\sqrt{z})+f(-c\sqrt{z})$ is entire of order $\frac{1}{2}\rho' \leq \frac{1}{2}\rho$.
\item $f(z)$ has the same order and type as all of its derivatives.
\item $e^z$ is entire of order $1$ and type $1$.
\item  $e^{cz^k}$  is entire of order $k$ and type $|c|$.
\item $e^{e^z}$ is entire of order $\infty$.
\item  The functions $(s-1)\zeta(s)$, $\zeta(s)-\frac{1}{s-1}$, and $\frac{1}{\Gamma(s)}$ are entire of order $1$ and type $\infty$.
\end{enumerate}
\end{remark}

\begin{remark}[Generalized dimensions and fractals]
References for the definitions and claims in this remark are \cite{olsen} and \cite{falc}.  Let $\mu$ be a  Borel probability measure on a metric space $M$, and let $X$ be the support of $\mu$.  Let $q \in \RR$.  The {\bf lower generalized dimension},\index{lower generalized dimension} or {\bf lower $q$-R\'enyi dimension},\index{lower $q$-R\'enyi dimension} $\underline{D}_\mu(q)$ of $\mu$  is defined by
\begin{align*}
\underline{D}_\mu(q) & = \liminf_{\varepsilon \to 0^+} \frac{\log \int_X \mu(B(x,\varepsilon))^{q-1}d \mu(x)}{(q-1)\log \varepsilon} 
 = -\frac{1}{q-1} \deg\int_X \mu(B(x,1/r))^{q-1}d \mu(x)
\end{align*}
if $q \neq 1$,  and
\begin{align*}
\underline{D}_\mu(1) =  \liminf_{\varepsilon \to 0^+} \frac{\int_X \log \mu(B(x,\varepsilon))d \mu(x)}{\log \varepsilon} = -\deg \exp  \int_X \log \mu(B(x,1/r))d \mu(x).
\end{align*}
if $q = 1$.   The  {\bf upper generalized dimension},\index{upper generalized dimension}  or {\bf upper $q$-R\'enyi dimension},\index{upper $q$-R\'enyi dimension} $\overline{D}_\mu(q)$ of $\mu$ is defined by replacing  the limits inferior above with limits superior  and the degrees above with lower degrees.    The R\'enyi dimensions were introduced by R\'enyi in 1960 as a tool for analyzing various problems in information theory \cite{ren}.
 One of the most significant aspects of R\'enyi dimensions is their relationship with {\it Hausdorff multifractal spectra}.  For any probability measure $\mu$ on $\RR^d$, the {\bf local dimension of $\mu$ at $x \in \RR^d$}\index{local dimension} is defined by
$$\operatorname{dim}_{\operatorname{loc}}(x;\mu) = \lim_{\varepsilon \to 0^+}  \frac{\log \mu(B(x,\varepsilon))}{\log \varepsilon} =  -\overline{\underline{\deg}} \, \mu(B(x,1/r)).$$  One defines the {\bf Hausdorff multifractal spectrum function}\index{Hausdorff multifractal spectrum function} $f_\mu$ of $\mu$ as the Hausdorff dimension of
the level sets of the local dimension of $\mu$, that is,
$$f_\mu(\alpha)= \dim\{x\in \RR^d: \operatorname{dim}_{\operatorname{loc}}(x;\mu) = \alpha\}, \quad \forall \alpha \geq 0,$$
where $\dim$ denotes the {\it Hausdorff dimension}.   The {\bf Legendre transform}\index{Legendre transform} $\varphi^*$ of a function $\varphi: \RR \to \RR$ is defined by $\varphi^*(x) = \inf_{y \in \RR} (xy + \varphi(y))$. In the 1980s, it was conjectured that for ``nice'' measures $\mu$ the upper and lower R\'enyi dimensions coincide, i.e., $\underline{D}_\mu(q) = \overline{D}_\mu(q)$, for all $q \in \RR$, and that the multifractal spectrum function $f_\mu$ of $\mu$ coincides with the Legendre
transform $\gamma_\mu^*$ of the function $\gamma_\mu: \RR \longrightarrow \RR$ defined by $$\gamma_\mu(q) = (1-q)\underline{D}_\mu(q) = (1-q)\overline{D}_\mu(q), \quad \forall q \in \RR.$$  Such  measures $\mu$ are said to satisfy the {\bf multifractal formalism}.
Since the 1990s several analysts and physicists have computed the multifractal spectrum and verified the multifractal formalism  for various measures in the mathematical literature.

Let $X$ be a subset of $\RR^d$.   For any $\varepsilon > 0$, let $N(X,\varepsilon)$ denote the number of boxes of side length $\varepsilon$ required to cover the set $X$.  The {\bf upper box dimension}, {\bf entropy dimension}, {\bf Kolmogorov dimension}, {\bf Kolmogorov capacity}, {\bf limit capacity}, or {\bf upper Minkowski dimension}, of $X$ is defined to be
$$\overline{\dim}_{\operatorname{box}} (X) = \limsup_{\varepsilon \to 0^+} \frac{\log N(X,\varepsilon)}{\log (1/\varepsilon)} = \deg N(X,1/r).$$  The {\bf lower box dimension} (resp., {\bf box dimension}\index{box dimension}) of $X$ is defined by replacing  the limit superior above with a limit inferior (resp., a limit) and the degree above with a lower degree (resp., exact degree).  For example, the famous {\it (middle-thirds) Cantor set} has box dimension $\frac{\log 2}{\log 3}$, the {\it Koch curve} has box dimension $\frac{\log 4}{\log 3}$, and the {\it Sierpi\'nski triangle} has box dimension $\frac{\log 3}{\log 2}$.

Many other notions of dimension, including  {\it correlation dimension}, {\it topological entropy}, and {\it uncertainty exponents} in dynamical systems and chaos theory, and  {\it information dimension}  in information theory, can also be expressed in terms of degrees.
\end{remark}

\section{The iterated logarithmic degree  $\degl f$ of a real function $f$}

As we demonstrate throughout this book, many functions in number theory are comparable, or are conjectured to be comparable, in some way (via $O$, $o$, $\sim$, etc.),  to a function $f(x)$ of the form $$f(x) = x^{d_0} (\log x)^{d_1} (\log^{\circ 2} x)^{d_2}\cdots (\log^{\circ k} x)^{d_k}$$
with $d_0, d_1, d_2, \ldots, d_k \in \RR$.   The set $\mathfrak L$ of all (germs of) such functions $f(x)$ is an abelian group\index[symbols]{.h  e@$\mathfrak L$}  (and, in fact, a real vector space) canonically isomorphic to $\bigoplus_{n = 0}^{ \infty}\RR$.  We denote the canonical isomorphism $\mathfrak L \longrightarrow \bigoplus_{n = 0}^{ \infty}\RR$ by $\degl$, where $\degl (\log^{\circ n} x)$ for all $n$ is the $n$th elementary unit vector of $\bigoplus_{n = 0}^{ \infty}\RR$.  The {\it group ring} $\RR[{\mathfrak L}]$\index[symbols]{.h  r@$\RR[{\mathfrak L}]$}
is an integral domain with quotient field $\RR({\mathfrak L})$,\index[symbols]{.h  rr@$\RR({\mathfrak L})$}  and the functions in $\mathfrak L$ are the {\bf monomials} of the ring $\RR[{\mathfrak L}]$\index{monomial of the field $\RR({\mathfrak L})$}  and field $\RR({\mathfrak L})$.  The field $\RR({\mathfrak L})$ 
is (isomorphic to) a subfield of the ring $\mathcal{C}^\infty$ of germs at $\infty$ of all infinitely differentiable functions on $\RR$, and like $\mathcal{C}^\infty$ it is closed under differentiation and is therefore a {\it differential $\mathbb{R}$-algebra}.

 We give the abelian group $\bigoplus_{n = 0}^{ \infty}\RR$ the lexicographic order and  $\mathfrak L$ the order $\preceq$ induced by the isomorphism with $\bigoplus_{n = 0}^{ \infty}\RR$.    Both groups are totally ordered abelian groups, and one has  $f \preceq g$ in  $\mathfrak L$ if and only if $f(x) = O(g(x)) \ (x \to \infty)$.  Moreover, one has $f \prec g$ if and only if  $f(x) = o(g(x)) \ (x \to \infty)$, if and only if $g-f$ is eventually positive.  The total ordering $\preceq$ of $\mathfrak{L}$  is a precise expression of the idea that $\log^{\circ l}x$ is ``infinitesimal'' with respect to $\log^{\circ k}x$ for all $l > k \geq 0$.

 Since $\mathfrak L$ is a totally ordered abelian group, the field $\RR({\mathfrak L})$ is a {\it valued field  with value group} $\bigoplus_{n = 0}^{ \infty}\RR$ (or $\mathfrak L$) and {\it valuation}
$$\degl: \RR({\mathfrak L}) \longmapsto \left(\bigoplus_{n = 0}^{ \infty}\RR\right) \cup \{\infty\},$$
where $$\degl \left (\sum_{f \in {\mathfrak L}} a_f f \right)= \max\{\degl f: a_f \neq 0\}$$
for all $\sum_{f \in {\mathfrak L}} a_f f \in  \RR[{\mathfrak L}]$, and where $\degl  (f/g) = \degl f-\degl g$ for all $f,g \in \RR[{\mathfrak L}]$ with $g \neq 0$.    (See Section 7.1 for the definitions of the terms in italics.)  The induced preorder $f \preceq g$  on  $\RR({\mathfrak L})$, defined by $\degl f \leq \degl g$, coincides with $O$, and $f \prec g$ coincides with $o$.  Moreover, any function in $\RR({\mathfrak L})$ is eventually positive, eventually $0$, or eventually negative, so that $\RR({\mathfrak L})$ is a totally ordered field under the total ordering $\leq_\infty$, where $f \leq_\infty g$ if $g-f$ is eventually positive or eventually $0$.  Furthermore, $f \preceq g$  is equivalent to $|f| \preceq |g|$ and is implied by $|f| \leq_\infty |g|$.

The ring $\RR[x^a: a \in \RR] \cong \RR[\RR^+]$ is a subring of $\RR[{\mathfrak L}]$, and the valuation $\degl$ restricted to $\RR[x^a: a \in \RR]$ coincides with the degree map $\deg$ defined in the previous section, making its quotient field $\RR(x^a: a \in \RR)$ a valued field with value group $\RR$.  We now extend the map $\degl$ to all functions in $\RR^{\RR_\infty}$ as we did with $\deg$ in the previous section.   

Let $f \in \RR^{\RR_\infty}$.  We define $f_{[k]}$ and $\degl_k f = \deg f_{[k]} \in \overline{\RR}$ for all nonnegative integers $k$, recursively, as follows:  let $f_{[0]} = f$ and $\degl_0 f = \deg f$,  suppose that $f_{[k]}$ is defined, and set $d_k = \degl_k f = \deg f_{[k]}$  and
$$f_{[k+1]}(x) =   \left.
  \begin{cases}
   f_{[k]}(e^x)e^{-d_k x} & \text{if } d_k \neq \pm \infty \\
   f_{[k]}(x)   & \text{if } d_k = \pm \infty.
 \end{cases}
\right.$$
We write $$\degl f = (\degl_0 f, \degl_1 f, \degl_2 f, \ldots) \in \prod_{n = 0}^ \infty \overline{\RR}.$$   We call $\degl f$ the {\bf (iterated) logarithmic degree of $f$}\index{iterated logarithmic degree $\degl$}\index{logarithmic degree $\degl$}\index[symbols]{.g e@$\degl f$}  and $\degl_k f$ the {\bf (iterated) logarithmic degree of $f$ of order $k$}.   Note that, if $\degl_k f = \pm \infty$ for some $k$, then $\degl_l f = \degl_k f$ for all $l \geq k$.    Clearly $\degl f$ coincides with $\degl f$ as defined earlier for any $f$ in the totally ordered normed field $\RR({\mathfrak L})$.

\begin{example}\
\begin{enumerate}
\item $\degl\left(\frac{\sqrt{x} \, \log \log \log x}{\log x} \right) = (\frac{1}{2}, -1, 0, 1, 0, 0, 0, \ldots)$. 
\item  $e^{(\log \log x)^2} = (\log x)^{\log \log x}$ and $e^{ \sqrt{\log x}}$  are examples of functions  $f$ such that $\deg f = 0$ while $\degl_1 f = \infty$, and thus $\degl f = (0,\infty,\infty,\infty,\ldots)$.
\item By Example \ref{divisorf}(3) one has $\deg d(n) = 0$, and by \cite[Theorem 314]{har} one has $d(n) \neq O((\log n)^t) \ (n \to \infty)$ for all $t$, which can be expressed as the equation $\degl_1 d(n) = \infty$.  Therefore, one has $\degl d(n) = (0,\infty,\infty,\infty,\ldots)$.
\end{enumerate}
\end{example}

In the next section, we redefine $f_{[k+1]}(x)$ in the case where $d_k = \pm \infty$ to yield an invariant, denoted $\dege f$, that is more refined than $\degl f$.

Throughout this book, we endow the set $\prod_{n = 0}^ \infty \overline{\RR}$ with  the lexicographic ordering.  It follows that $\prod_{n = 0}^ \infty \overline{\RR}$ is a totally ordered set that is complete, i.e., it possesses both the supremum and infimum of any of its subsets.  In this section, we also endow the set $\prod_{n = 0}^ \infty \overline{\RR}$ with total binary operations of coordinatewise addition and subtraction, where the sums $\infty+(-\infty)$ and $(-\infty)+ \infty$ are defined to be $\infty$, while the differences $\infty- \infty$ and $-\infty- (-\infty)$ are defined to be $-\infty$.    This  makes addition and subtraction defined on all of $\prod_{n = 0}^ \infty \overline{\RR}$, although the two operations are not inverses of each other since for example $(\infty+\infty)-\infty = -\infty$ and $(\infty-\infty)+\infty = \infty$.  The only reason to set these particular conventions is so that $\degl(fg) \leq \degl f+ \degl g$ and $\degl(f/g) \geq \degl f-\degl g$  can both hold generally as stated in Proposition \ref{diffprop} below.  

Note that $\degl f \in \prod_{n = 0}^{ * \infty} \overline{\RR}$ for all functions $f$ defined on a neighborhood of $\infty$, where the {\bf restricted product} $\prod_{n = 0}^{ * \infty} \overline{\RR}$\index{restricted product $\prod_{n = 0}^{ * \infty} \overline{\RR}$}\index[symbols]{.g ga@$\prod_{n = 0}^ {* \infty} \overline{\RR}$} is defined to be the union of $\prod_{n = 0}^{ \infty } \RR$ with the set of all  sequences in $\prod_{n = 0}^ \infty \overline{\RR}$ that are either eventually $\infty$ or eventually  $-\infty$.   The  induced total ordering on the restricted product is also complete.  Moreover, the restricted product contains $\bigoplus_{n = 0}^{ \infty}\RR$ as a dense subset, i.e., for any two distinct elements of $\prod_{n = 0}^ {*\infty} \overline{\RR}$  there is an element of $\bigoplus_{n = 0}^{ \infty}\RR$ lying strictly between them.

Let $N$ be a positive integer, and suppose that $\degl_k f \neq \pm \infty$ for all $k < N$.  The functions $f_{[k]}$ are then determined by the relations $f_{k}(x) = x^{\deg f_{[k]}} f_{[k+1]}(\log x),$
or by the relations
$$f(x) =  x^{\deg f_{[0]}} (\log x)^{\deg f_{[1]}} (\log^{\circ 2}x)^{\deg f_{[2]}}\cdots (\log^{\circ k}x)^{\deg f_{[k]}}f_{[k+1]}(\log^{\circ(k+1)} x),$$
for all $k < N$.   Explicitly, one has $f_{[k]} =  \exp\left(-\sum_{i +j = k-1} \degl_i(f)  \exp^{\circ j} \right)f(\exp^{\circ k})$ for all $k \leq N$.  The following proposition is easily verified.

\begin{proposition}\label{Nprop}
Let $f \in \RR^{\RR_\infty}$, and let $d_k \in \RR$ for all $k \leq N$, where $N$ is a nonnegative integer.   One has $d_k = \degl_k f$ for all $k \leq N$ if and only if $$  f(x) = o\left(x^{d_0} (\log x)^{d_1} (\log^{\circ 2}x)^{d_2} \cdots (\log^{\circ (N-1)} x)^{d_{N-1}}  (\log^{\circ N} x)^{t} \right) \ (x \to \infty)$$ for all $t > d_N$
but for no $t<d_N$; moreover, if those conditions hold, then one has
\begin{align*}
\degl_{N+1} f  & = \limsup_{x \to \infty}  \frac{\log |f(e^x)|-d_0x- d_1 \log x -d_2 \log^{\circ 2}x \cdots - d_{N} \log^{\circ N}x}{\log^{\circ (N+1)} x} \\
   & = \deg\left( f(\exp^{\circ(N+1)})\exp\left(- d_0\exp^{\circ N} x - d_1 \exp^{\circ (N-1)} x - \cdots - d_{N}x \right) \right).
\end{align*}
\end{proposition}

The following proposition generalizes Proposition \ref{firstprop1}.

\begin{proposition}\label{asprop}
Let $f, g \in \RR^{\RR_\infty}$, and let $\underline{g} = g|_{\dom f \cap \dom g}$.
\begin{enumerate}
\item If $f(x) = O(g(x)) \ (x \to \infty)$, then $\degl f \leq \degl \underline{g} \leq \degl g$.
\item If $f(x) \gg g(x) \ (x \to \infty)$, then $\degl f \geq \degl \underline{g}$.
\item If $f(x) \asymp g(x) \ (x \to \infty)$, then $\degl f = \degl \underline{g}$ and $f_{[k]}(x) \asymp \underline{g}_{[k]}(x) \ (x \to \infty)$ for all $k$.
\item If $f$ is eventually bounded, then  $\degl f \leq (0,0,0,\ldots)$.
\item If $\degl f < (0,0,0,\ldots)$, then $\displaystyle \lim_{x \to \infty} f(x)  = 0$.
\end{enumerate}
\end{proposition}

\begin{proof}
Suppose that $f(x) = O(g(x)) \ (x \to \infty)$.  Then $\deg f \leq \deg \underline{g}$.  If $\deg f < \deg \underline{g}$, or if $\deg f = -\infty$ or $\deg \underline{g} = \infty$, then $\degl f \leq \degl \underline{g}$.  Thus we may suppose that $-\infty < \deg f = \deg \underline{g} < \infty$.  It follows that 
$$f_{[1]}(x) = f(e^x)e^{-(\deg f)x} = O(\underline{g}(e^x)e^{-(\deg \underline{g})x}) = O(\underline{g}_{[1]}(x)) \ (x \to \infty).$$
Now proceed by induction.   This proves (1).  Statements (2)--(5) are also easily proved by induction.
\end{proof}

Because of Proposition \ref{asprop}(1), we think of $\degl f \leq \degl g$ as a total preorder that is substantially coarser than the preorder $f(x) = O(g(x)) \ (x \to \infty)$  (and substantially finer than the total preorder $\deg f \leq \deg g$), although the two preorders agree on $\RR(\mathfrak{L})$.  The following proposition provides a more precise relationship between these two preorders.

\begin{proposition}\label{infprop}
Let $f$ be a real function defined on a neighborhood of $\infty$.    One has
\begin{align*}
\degl f & = \inf\{\degl r:r \in \mathfrak{L}, \, \degl f \leq \degl r \} \\ 
& = \inf\{\degl r:r \in \mathfrak{L}, \, f(x) = O(r(x)) \ (x \to \infty) \} \\ 
 & = \inf\{\degl r:r \in \mathfrak{L}, \, \forall x \gg 0\, |f(x)|\leq r(x) \} \\
& = \inf\{\degl r:r \in \mathfrak{L}, \, f(x) = o(r(x)) \ (x \to \infty) \} \\
& = \inf\{\degl r:r \in \mathfrak{L}, \, \degl f < \degl r \},
\end{align*}
where the infima are computed in the restricted product $\prod_{n = 0}^ {* \infty} \overline{\RR}$.
\end{proposition}

\begin{proof}
If $r \in \mathfrak{L}$, then clearly $\degl (f/r) = \degl f-\degl r$, so that $f(x) = o(r(x)) \ (x \to \infty)$ if $\degl f < \degl r$.  Therefore, Proposition \ref{asprop} implies that
\begin{align*}
\degl f & \leq \inf\{\degl r:r \in \mathfrak{L}, \, \degl f \leq \degl r \} \\ 
& \leq \inf\{\degl r:r \in \mathfrak{L}, \, f(x) = O(r(x)) \ (x \to \infty) \} \\ 
 & \leq \inf\{\degl r:r \in \mathfrak{L}, \, \forall x \gg 0\, |f(x)|\leq r(x) \} \\
& \leq \inf\{\degl r:r \in \mathfrak{L}, \, f(x) = o(r(x)) \ (x \to \infty) \} \\
& \leq \inf\{\degl r:r \in \mathfrak{L}, \, \degl f < \degl r\},
\end{align*}
Moreover, the equality $\degl f =  \inf\{\degl r:r \in \mathfrak{L}, \, \degl f < \degl r\}$ follows from the fact that  the ordered set $\prod_{n = 0}^{*\infty} \overline{\RR}$ is complete and
$\bigoplus_{n = 0}^{ \infty}\RR = \{\degl r: r \in \mathfrak{L}\}$ is dense in the poset $\prod_{n = 0}^{*\infty} \overline{\RR}$.
\end{proof}

 The next proposition relates $\degl$ to the operations of addition, multiplication, and division of functions.

\begin{lemma}\label{difflem}
Let $f,g \in \RR^{\RR_\infty}$ with $\dom f = \dom g$.  Then one has
$$\max(\degl f,\degl g) = \degl \max(|f(x)|,|g(x)|).$$
\end{lemma}

\begin{proof}
Since $0 \leq |f(x)| \leq  \max(|f(x)|,|g(x)|)$ and $0 \leq  |g(x)| \leq  \max(|f(x)|,|g(x)|)$ for all $x$ (in $\dom f = \dom g$), the inequality   $\max(\degl f,\degl g) \leq \degl \max(|f(x)|,|g(x)|)$ is immediate.  To prove the reverse inequality, we may suppose without loss of generality that $\degl f \geq \degl g$.  Let $r \in \mathfrak{L}$ with $\degl f < \degl r$.  Since $\degl g < \degl r$, one has $|f(x)| \leq r(x)$ and $|g(x)| \leq r(x)$,  and therefore $\max(|f(x)|,|g(x)|) \leq r(x)$, for all $x \gg 0$, whence
$\degl\max(|f(x)|,|g(x)|) \leq \degl r$.  Taking the infimum over all $r$ as chosen and applying Proposition \ref{infprop},  we conclude that $\degl \max(|f(x)|,|g(x)|) \leq \degl f =  \max(\degl f,\degl g)$.
\end{proof}

\begin{proposition}\label{diffprop}
Let $f$ and $g$ be real functions such that both $\underline{f} = f|_{\dom f \cap \dom g}$ and $\underline{g} = g|_{\dom f \cap \dom g}$ are in $\RR^{\RR_\infty}$.  One has the following.
\begin{enumerate}
\item $\degl(f + g) \leq \max (\degl \underline{f},\degl \underline{g})$, with equality if $\degl \underline{f} \neq  \degl \underline{g}$.  
\item If $\degl(f+g) \neq \max(\degl \underline{f}, \degl \underline{g})$, then $\degl \underline{f} = \degl \underline{g}$.   
\item $\degl(fg) \leq \degl \underline{f} + \degl \underline{g}$.
\item $\degl(f/g) \geq \degl \underline{f}-  \degl \underline{g}$ (resp., $\degl (fg) \geq \degl \underline{f} - \degl (1/g)$) if $g$ is eventually nonzero on its domain.
\end{enumerate}
\end{proposition}

\begin{proof}
Since $f+g = \underline{f}+\underline{g}$ and $f g = \underline{f}\, \underline{g}$, we may suppose without loss of generality that $f = \underline{f}$ and $g = \underline{g}$.   Since $$|f(x)+g(x)| \leq |f(x)|+|g(x)| \leq 2 \max(|f(x)|,|g(x)|)$$ for all $x \in \dom f = \dom g$, by  Lemma \ref{difflem} one has 
$$\degl (f+g) \leq \degl  \max(|f(x)|,|g(x)|) = \max(\degl f,\degl g).$$  This proves statement (1), and statement (2) is an immediate consequence of (1).  Since $(fg)_{(1)} = f_{[1]} g_{[1]}$ in the case where $\deg (fg) = \deg f + \deg g$, we see that statements (3) and (4) follow by induction from Proposition \ref{firstprop2} and our conventions for adding and subtracting $\infty$ and $-\infty$.
\end{proof}

Next, we wish to determine some conditions on $f$ or $g$ that imply  $\degl(fg) = \degl f + \degl g$.  To do so, we need to generalize the notion of exact degree appropriately.

Let $f \in \RR^{\RR_\infty}$.    We  say that $f$ has {\bf finite (iterated) logarithmic degree to order $n$}\index{finite (iterated) logarithmic degree} if $\degl_k f \neq \pm \infty$ for all $k \leq n$, and  $f$ has {\bf finite (iterated) logarithmic degree (to order $\infty$)} if $f$ has finite logarithmic degree to order $n$ for all $n$.    We say that $f$ has {\bf exact (iterated) logarithmic degree to order $n$}\index{exact (iterated) logarithmic degree} if $f_{[k]}$ has exact degree for all $k \leq n$, where the functions $f_{[k]}$ are defined recursively as in the definition of $\degl f$.  We say that $f$ has {\bf exact (iterated) logarithmic degree (to order $\infty$)} if $f$ has exact degree to order $n$ for all $n$, that is, if $f_{[k]]}$ has exact degree for all $k$.  Thus, $f$ has finite  logarithmic degree to order $0$ if and only if $f$ has finite degree; and $f$ has exact logarithmic degree to order $0$ if and only if $f$ has exact degree.    As a degenerate case, if a function $f$ is eventually $0$, i.e., if $f(x) \sim 0 \ (x \to \infty)$, then one has $\degl f = (-\infty,-\infty,-\infty,\ldots$), but, since our convention is that  $\log |0| = -\infty$ and $\lim_{x \to \infty} (-\infty) = -\infty$, we consider $f$ to have exact logarithmic degree.

\begin{example} \item
\begin{enumerate}
\item All functions in  $\RR(\mathfrak L)\backslash \{0\}$ have finite and exact logarithmic degree.
\item If $f \in \RR^{\RR_\infty}$ and $f(x) \asymp r(x) \ (x \to \infty)$  for some $r \in \RR(\mathfrak L) \backslash \{0\}$,  then $f$ has finite and exact logarithmic degree $\degl r$. 
\item $\sin x$ has finite logarithmic degree $(0,0,0,\ldots)$ but does not have  exact logarithmic degree to any order.
\item $e^x$ has exact logarithmic degree $(\infty,\infty,\infty,\ldots)$ but does not have finite logarithmic degree to any order.
\item $e^x \sin x$ has iterated logarithmic degree $(\infty,\infty,\infty,\ldots)$ and does not have finite or exact logarithmic degree to any order.
\item  $e^{\sqrt{\log^{\circ k}x}}$ has exact logarithmic degree  $(0,0,0,\ldots, 0, \infty, \infty,\infty,\ldots)$, but it has finite logarithmic degree only to order $k$, for all nonnegative integers $k$.
\item $(\log x)^{\sin \log \log x}$ has iterated logarithmic degree $(0,1,0,0,0,\ldots)$ and, by Corollary \ref{regvarcor}, is slowly varying, hence of exact degree $0$, but it does not have exact logarithmic degree to order $n$ for any $n > 0$.  
\end{enumerate}
\end{example}

For any sequences $\dd = (\dd_0,\dd_1,\dd_2,\ldots)$ and  $\ee = (\ee_0,\ee_1,\ee_2,\ldots)$ in $\prod_{n = 0}^\infty\overline{\RR}$, we let
$\mathcal{S}(\dd,\ee)$ denote the smallest nonnegative integer $k$ such that $\dd_k \neq \ee_k$, where we let $\mathcal{S}(\dd,\ee)= \infty$ if no such integer $k$ exists.\index[symbols]{.g n@$\mathcal{S}(\dd,\ee)$}  
Also, for  any functions $f, g\in \RR^{\RR_\infty}$, we let $\mathcal{L}(f,g) = \mathcal{S}(\degl f, \degl g).$\index[symbols]{.g o@$\mathcal{L}(f,g)$}    Thus, $\mathcal{L}(f,g)$ is the smallest nonnegative integer $k$ such that $\degl_k f \neq \degl_k g$, where $\mathcal{L}(f,g) = \infty$ if no such integer $k$ exists.

\begin{proposition}\label{mmm1}
Let $f \in \RR^{\RR_\infty}$, and let $n$ be a nonnegative integer.  One has the following.
\begin{enumerate} 
\item $f$ has finite logarithmic degree to order $n$ if and only if there  exists an $r \in \mathfrak{L}$ such that $\mathcal{L}(f,r) > n$, in which case $r$ is unique subject to the condition $\degl_k r = 0$ for all $k > n$.
\item $f$ has finite and exact logarithmic degree to order $n$ if and only if the limit
$$d =  \lim_{x \to \infty} \frac{ \log \left|\frac{f(x)}{r(x)} \right|}{\log  \log^{\circ n}x}$$
exists for some $r \in \mathfrak{L}$.  In that case, one has 
$\mathcal{L}(f,r) \geq n$ and $\degl_n f = d+\degl_n r.$  Moreover, $r$ is unique subject to the condition $\degl_k r= 0$ for all $k \geq n$, and in that case one has $\degl_{n} f = d$.
\end{enumerate}
\end{proposition}

\begin{proposition}\label{mmm2} Let $f, g\in \RR^{\RR_\infty}$ with  $f(x) \asymp g(x) \ (x \to \infty)$, let $\underline{g} = g|_{\dom f \cap \dom g}$, and let $n$ be a nonnegative integer.   Then $\degl f = \degl \underline{g}$, and $f$ has exact  (resp., finite) logarithmic degree to order $n$ if and only if $\underline{g}$ has exact (resp., finite) logarithmic degree to order $n$.
\end{proposition}

\begin{proof}
This follows by induction from Propositions \ref{firstprop1}(3) and \ref{asprop}(3).
\end{proof}

The following result  provides equivalent characterization of all real functions that have exact logarithmic degree to order $n \geq 0$.

\begin{proposition} \label{exact}
Let $g \in \RR^{\RR_\infty}$ with $g$ not eventually $0$, and let $n$ be a nonnegative integer.  The following conditions are equivalent.
\begin{enumerate}
\item $g$ has  exact logarithmic degree to order $n$.
\item $\degl_k(fg) - \degl_k g= \degl_k f $ for all $k \leq n$ for all real functions $f$ such that $[N,\infty) \cap \dom f = [N , \infty) \cap \dom g$ for some $N > 0$ and such that $\degl_k(fg) - \degl_k g$ is defined in $\overline{\RR}$ for all $k \leq n$.
\item $\degl_k f = -\degl_k g$ for all $k \leq n$ for all real functions $f$ with $\dom f  = \dom g$ and $\degl(fg) = (0,0,0,\ldots)$.
\end{enumerate}
\end{proposition}

\begin{proof}
Statement (1) implies statement (2) (by induction) because $$\limsup_{x \to \infty}\left( \frac{\log |(fg)_{[k]}(x)|}{\log x}-\frac{\log |g_{[k]}(x)|}{\log x} \right) = \limsup_{x \to \infty} \frac{\log | (fg)_{[k]}(x)|}{\log x}-\lim_{x \to \infty} \frac{\log |g_{[k]}(x)|}{\log x}$$
provided that the limit on the right exists and the difference on the right is defined in $\overline{\RR}$.
Clearly statement (2) implies statement (3).  Thus, it remains only to show that (3) implies (1).

 For all $x \in \dom g$, let
$$f(x) = \begin{cases}  \frac{1}{g(x)} & \text{if } g(x) \neq 0 \\ 
  0 & \text{if } g(x) = 0,
\end{cases}$$
so that $\dom f = \dom g$ and
$$f(x)g(x) = \begin{cases}  1 & \text{if } g(x) \neq 0 \\ 
  0 & \text{if } g(x) = 0,
\end{cases}$$
while also
$$\log |f(x)| = \begin{cases} -\log |g(x)| & \text{if } g(x) \neq 0 \\ 
  -\infty & \text{if } g(x) = 0.
\end{cases}$$
It follows that $\degl (f g) = (0,0,0,\ldots)$ and therefore $\degl_k g = -\degl_k f$ for all $k \leq n$, by hypothesis (3).  Now, one has
$$\deg g = -\deg f = -\limsup_{x \to \infty} \frac{\log |f(x)|}{\log x} = - \limsup_{x \to \infty} \left(-\frac{\log |g(x)|}{\log x} \right)= \liminf_{x \to \infty} \frac{\log |g(x)|}{\log x},$$
and therefore $g$ has exact degree, say, $d_0$.  If $d_0$ is not finite, then $g$ has exact logarithmic degree.  Thus we may suppose that  $d_0$ is finite and $n \geq 1$.   Then one has
\begin{align*}
\degl_1 g & = -\degl_1 f \\
 & = -\limsup_{x \to \infty} \frac{\log |f(x)|+d_0 \log x}{\log \log x} \\
& = -\limsup_{x \to \infty} \left(-\frac{\log |g(x)|-d_0 \log x}{\log \log x} \right)\\
& = \liminf_{x \to \infty} \frac{\log |g(x)|-d_0\log x}{\log \log x},
\end{align*}
so that $g$ has exact degree to order $1$.  Clearly this argument can be continued {\it ad infinitum}.
\end{proof}

\begin{proposition}\label{mmm}
Let $f$ and $g$ be real functions such that both $\underline{f} = f|_{\dom f \cap \dom g}$ and $\underline{g} = g|_{\dom f \cap \dom g}$ are in $\RR^{\RR_\infty}$.   Suppose also that $\underline{g}$ has finite and exact  logarithmic degree to order $n \geq 0$.  One has the following.
\begin{enumerate}
\item  $\degl_k(fg) = \degl_k \underline{f} + \degl_k  \underline{g}$ for all $k \leq n$, that is,  $$\mathcal{S}(\degl(fg), \degl  \underline{f} + \degl  \underline{g}) > n.$$
\item $\degl_k(f/g) = \degl_k  \underline{f} -\degl_k  \underline{g}$ for all $k \leq n$, that is,  $$\mathcal{S}(\degl(f/g), \degl  \underline{f} - \degl  \underline{g}) > n.$$
\item For any nonnegative integer $k \leq n$, one has $\mathcal{L}(\underline{f}, \underline{g}) > k$  if and only if $$\limsup_{x \to \infty} \frac{\log|\underline{f}(x)|-\log| \underline{g}(x)|}{\log \log^{\circ k}x} = 0.$$
\item Suppose that $\underline{g}$ has finite and exact logarithmic degree.  Then one has  $\degl \underline{f} = \degl \underline{g}$ if and only if 
$$\limsup_{x \to \infty} \frac{\log|\underline{f}(x)|-\log|\underline{g}(x)|}{\log^{\circ k}x} = 0$$ for all $l$.  More generally,  one has
$$\mathcal{L}(\underline{f},\underline{g}) = \inf \left\{k \in \ZZ_{\geq 0}: \limsup_{x \to \infty} \frac{\log|\underline{f}(x)|-\log|\underline{g}(x)|}{\log \log^{\circ k}x} \neq 0\right \}.$$
\end{enumerate}
\end{proposition}

\begin{proof}
This is more or less a corollary of Proposition \ref{exact}.
\end{proof}

A useful application of Proposition  \ref{mmm} is the following.

\begin{corollary}\label{deglogprop}
Let $f,g \in \RR^{\RR_\infty}$ with $\underline{g} =  g|_{\dom f \cap \dom g} \in \RR^{\RR_\infty}$, and suppose that $\underline{g}$  has finite and exact logarithmic degree and $\dom g$ contains the intersection of $\dom f$ with some neighborhood of $\infty$.  Then one has the following.
\begin{enumerate}
\item If $f(x) \neq o(g(x)) \ (x \to \infty)$, or, equivalently, if $ \limsup_{x \to \infty}  \left| \frac{f(x)}{g(x)} \right| > 0$, then $\degl f \geq \degl  \underline{g}$.
\item  If $f(x) = O(g(x)) \ (x \to \infty)$ but $f(x) \neq o(g(x)) \ (x \to \infty)$, or, equivalently, if $ \limsup_{x \to \infty}  \left| \frac{f(x)}{g(x)} \right|$ is finite and positive, then $\degl f = \degl \underline{g}$.
\end{enumerate}
\end{corollary}

\begin{proof}
To prove (1), first note that, by Proposition \ref{asprop}(6), one has $\degl(f/g) \geq (0,0,0,\ldots)$.  Therefore, by  Proposition \ref{mmm}(2), one has $\degl f-\degl  \underline{g}  =  \degl(f/g) \geq (0,0,0,\ldots)$, and therefore $\degl f\geq \degl  \underline{g}$.  Statement (2), then, follows from statement (1) and Proposition \ref{asprop}(1).
\end{proof}

\begin{example}\label{newex} \
\begin{enumerate}
\item By Corollary \ref{deglogprop} and Example \ref{divisorf}, one has $$\degl \log d(n) = (0,1,-1,0,0,\ldots)$$
and  $$\degl \sigma(n) = (1,0,1,0,0,0\ldots).$$
\item  Recall from Chapter 3 that $\Omega(n)$ for any positive integer $n$ denotes the number of prime factors of $n$, counting multiplicities.  Then one has
$$\Omega(n) \leq \log_2 n$$
for all $n$, with equality if and only if $n$ is a power of $2$.
It follows that
$$\limsup_{n \to \infty} \frac{\Omega(n)}{\log_2 n} = 1.$$
By Corollary  \ref{deglogprop}, then, one has 
$$\degl \Omega(n) = (0,1,0,0,0,\ldots).$$
\item Recall from Chapter 3 that $\omega(n)$ for any positive integer $n$ denote the number of distinct prime factors of $n$.  For any $x \geq 0$, let $$x\# = \prod_{p \leq x} p,\index[symbols]{.s BA@$x\#$}\index{primorial $x\#$}$$  denote the product of all of the primes less than or equal to $x$,  which is called {\bf $x$ primorial}.
Then one has $\omega(x\#) = \pi(x)$ and $\log(x\#) = \vartheta(x)$,
where $$\vartheta(x) = \sum_{p \leq x} \log p, \quad \forall x \geq 0,$$ is the first Chebyshev function.  By Corollary \ref{pitheta2} and the prime number theorem, one has $$\vartheta(x) \sim \pi(x)\log x \sim x \ (x \to \infty).$$  It follows that
$$\omega(x \#) = \pi(x) \sim \frac{\vartheta(x)}{\log x} =  \frac{\log (x\#)}{\log x}\sim \frac{\log (x\#)}{\log \log (x\#)} \ (x \to \infty).$$
From this one can show that
$$\limsup_{n \to \infty}  \frac{\omega(n)}{\frac{\log n}{\log \log n}} = 1.$$
By Corollary \ref{deglogprop}, then, one has
$$\degl \omega(n) = (0,1,-1,0,0,0,\ldots).$$
\end{enumerate}
\end{example}

Finally, we relate logarithmic degree to the operation of composition of functions.

\begin{proposition}\label{fgi}
Let $f, g \in \RR^{\RR_\infty}$, and let $n$ be a nonnegative integer.  Suppose that the following conditions hold.
\begin{enumerate}
\item $g$ has  finite and positive exact degree and is eventually positive.
\item $f \circ g  \in \RR^{\RR_\infty}$. 
\item $\underline{f} = f|_{\dom f \cap \im g}$ has finite degree.
\item  Either $g$ has exact logarithmic degree to order $n$ or $\deg f\geq 0$ and $f$ has exact logarithmic degree to order $n$.
\item  If $n \geq 1$, then $g$ has finite logarithmic degree to order $n-1$ and $\degl_n \underline{f}+ \deg \underline{f} \degl_n g$ is defined in $\overline{\RR}$ (e.g., if $g$ has finite logarithmic degree to order $n$).
\end{enumerate}
Then one has $\deg(f \circ g) = \deg \underline{f} \deg g$ and $$\degl_k(f \circ g) = \degl_k \underline{f} + \deg \underline{f} \degl_k g, \quad \forall k = 1,2,3,\ldots,n.$$
Moreover, if  condition $(5)$ holds for all $n$, then one has
\begin{align*}
\degl(f \circ g) & = \degl \underline{f} + (\deg \underline{f}) (\degl g-(1, 0,0,0,\ldots)) \\
& \leq \degl f+ (\deg f) (\degl g-(1, 0,0,0,\ldots)),
\end{align*}
and equality holds if either   $[N,\infty) \cap \dom f \subseteq \im g$ for some $N > 0$ (e.g., if $g$ is defined and continuous on some neighborhood of $\infty$) or $f$ has exact logarithmic degree, or, more generally, if $\degl \underline{f} = \degl f$.
\end{proposition}

\begin{proof}
Condition (1) implies that $\lim_{x \to \infty}  g(x) = \infty$ and that, for all $k \geq 2$, the limit $\lim_{x \to \infty} (\log^{\circ k} g(x) - \log^{\circ k} x)$ exists and is finite (and equals $0$ if $k \geq 3$), and therefore $\log^{\circ k} g(x) \sim \log^{\circ k} x  \ (x \to \infty)$.
For all $k$, let $d_k = \degl_k \underline{f}$, and $e_k = \degl_k g$, and also let $$c_k = \degl_k \underline{f} + d_0 \degl_k g$$
for all $k \geq 1$.  By Proposition \ref{circle1}(2), the first claim of the proposition is true for $n = 0$ or $k = 0$, that is, one has $\deg (f\circ g) = d_0e_0$.  Thus we may suppose that $n > 0$ and $1 \leq k \leq  n$.   Then one has
\begin{align*}
c_1 & =  \limsup_{x \to \infty}  \frac{\log |\underline{f}(x)|-d_0 \log x}{\log \log x} + d_0 \degl_1 g \\
& =  \limsup_{x \to \infty}  \frac{\log |f(g(x))|-d_0 \log g(x)}{\log \log g(x)} +d_0 \degl_1 g  \\
& =  \limsup_{x \to \infty}  \frac{\log |f(g(x))|-d_0 \log g(x)}{\log \log x} +d_0 \limsup_{x \to \infty}  \frac{\log g(x)-e_0 \log x}{\log \log x} \\
& =  \limsup_{x \to \infty}  \frac{\log |f(g(x))|-d_0 e_0 \log x}{\log \log x} \\
& = \degl_1 (f \circ g),
\end{align*}
where the penultimate equality above follows from condition (4), since by that hypothesis at least one of the two limits superior can be replaced by a limit.  Moreover, if $\deg \underline{f} = \deg f$, then $\degl_1 \underline{f} \leq \degl_1 f$ and therefore
$c_1 \leq \degl_1 f+ d_0 \degl_1 g$.  Now, suppose  that $n \geq 2$.  Then one has
\begin{align*}
c_2 & =  \limsup_{x \to \infty}  \frac{\log |\underline{f}(x)|-d_0 \log x-d_1 \log \log x}{\log\log \log x} + d_0 \degl_2 g \\
& =  \limsup_{x \to \infty}  \frac{\log |f(g(x))|-d_0 \log g(x)-d_1 \log \log g(x)}{\log \log \log g(x)} + d_0 \degl_2 g  \\
& =  \limsup_{x \to \infty}  \frac{\log |f(g(x))|-d_0 \log g(x)-d_1 \log \log g(x)}{\log \log \log x}  \\
  & \qquad \qquad  \qquad +d_0\limsup_{x \to \infty}  \frac{\log g(x)-e_0 \log x- e_1 \log \log x}{\log \log \log x} \\
& =  \limsup_{x \to \infty}  \frac{\log |f(g(x))|-d_0 e_0 \log x-(d_1+d_0e_1) \log\log x-d_1(\log \log g(x)-\log \log x)}{\log \log \log x} \\
& =  \limsup_{x \to \infty}  \frac{\log |f(g(x))|-d_0 e_0 \log x-(d_1+d_0e_1) \log\log x}{\log \log \log x} \\
& = \degl_2 (f \circ g).
\end{align*}
It is clear, then, how to proceed by induction to verify the claim for all $k\leq n$.

Now, if $[N,\infty) \cap \dom f \subseteq \im g$, then $\underline{f}(x) = f(x)$ for all $x \gg 0$ in  $\dom f$, so that $\degl \underline{f} = \degl f$.  On the other hand, if $f$ has exact logarithmic degree,  then
\begin{align*}
\degl_k f  & =  \lim_{x \to \infty}  \frac{\log |f(x)|-d_0 \log x-\cdots -d_k \log^{\circ k} x}{\log \log^{\circ k} x} \\
 & =  \lim_{x \to \infty}  \frac{\log |f(g(x))|-d_0 \log g(x)-\cdots -d_k \log^{\circ k} g(x)}{\log \log^{\circ k} g(x)}  \\
 & = \degl_k \underline{f}
\end{align*}
for all $k$, so that $\degl \underline{f} = \degl f$. This completes the proof.
\end{proof}

\section{The logexponential degree $\dege f$ of a real  function $f$}

In this section, we provide a more refined alternative $\dege$ to $\degl$ that we call {\it logexpondential degree}.  This section forms the basis for Part 3.   We show that  $\dege$ shares many of the  properties possessed by $\deg$ and  $\degl$, which provides further justification for its ostensibly {\it ad hoc} definition.  We also provide a canonical axiomatization of logexponential degree in terms of $O$, $o$, and the logarithmico-exponential functions.  As this section is long and, unavoidably, somewhat technical,   we have divided it into five subsections:  {\it Logexponential degree} (6.3.1),   {\it Relationships with $O$, $o$, and $\asymp$} (6.3.2),  {\it Relationships with logarithmico-exponential functions} (6.3.3),   {\it Relationships with operations on functions} (6.3.4), and {\it Axiomatization of $\dege$} (6.3.5).

\subsection{Logexponential degree}

Let $f \in \RR^{\RR_\infty}$,  and define $f_{(k)}$ and $\dege_k f = \deg f_{(k)} \in \overline{\RR}$ for all nonnegative integers $k$, recursively, as follows.    Let $f_{(0)} = f$.  Suppose that $f_{(k)}$ is defined, and set $d_k = \dege_k f = \deg f_{(k)}.$  We then let
$$f_{(k+1)}(x) =   \left.
 \begin{cases}
    f_{(k)}(e^x) e^{-d_k x}& \text{if } d_k \neq \pm \infty \\
    \max( \log |f_{(k)}(x)|, 0)  & \text{if } d_k = \infty \\
 \displaystyle   -\frac{1}{\log |f_{(k)}(x)|} & \text{if } d_k =- \infty.
 \end{cases}
\right.$$
This defines $\dege_k f$ for all nonnegative integers $k$.  We also set $$\dege f = (\dege_0 f, \dege_1 f, \dege_2 f, \ldots) \in \prod_{n = 0}^\infty \overline{\RR}.$$   We call $\dege f$ the {\bf logexponential degree of $f$}\index{logexponential degree $\dege$}\index{exponential degree $\dege$}\index[symbols]{.g g@$\dege f$}  and $\dege_k f$ the {\bf logexponential degree of $f$ of order $k$}.  

If $f$ has finite logarithmic degree to order $n-1$, then one has $\dege_k f = \degl_k f$ for all $k \leq n$.  However, if $\degl_{n} f = \pm \infty$, then $\dege_{n+1} f$ might not equal $\degl_{n+1} f$.  Indeed, the definition of $f_{(n+1)}(x)$, as $ \max( \log |f_{(n)}(x)|, 0)$ or as $-\frac{1}{\log |f_{(n)}(x)|}$ according to whether  $\degl_{n} f =\infty$ or  $\degl_{n} f = - \infty$, is designed to ``tame'' the function $f_{(n)}$ by applying a log to $|f_{(n)}|$ appropriately.   Trivially, $\dege$ and $\degl$ are related as follows.

\begin{proposition}\label{degeprop0}
Let  $f ,g\in \RR^{\RR_\infty}$, and let the $f_{(k)}$ be defined recursively as in the definition of $\dege$, so that $\dege_k f = \deg f_{(k)}$ for all $ k$.  For all $n \geq 0$, one has the following.
\begin{enumerate}
\item If $f$ has finite logarithmic degree to order $n$, then $\dege_k f = \degl_k f$ for all $k \leq n+1$.
\item If $f$ has finite logarithmic degree, then $\degl f = \dege f$.  Otherwise, there is a smallest $n$ such that  $\dege_n f= \pm \infty$, and then one has $$\degl f = (\dege_0 f, \dege_1 f, \ldots, \dege_{n-1} f, \dege_n f, \dege_n f, \dege_n f, \ldots).$$
\item If $\dege f_{(n)} = \degl g$, then $$\dege f = (\dege_0 f, \dege_1 f, \ldots, \dege_{n-1} f, \degl_0 g, \degl_1 g, \degl_2 g, \ldots).$$
\item If $\degl f < \degl g$, then $\dege f < \dege g$.
\item If $\dege f \leq \dege g$, then $\degl f \leq \degl g$.
\end{enumerate}
\end{proposition}

Note, in particular, that $\degl f$ is tantamount to the truncation of $\dege f$ at the $n$th coordinate for the smallest $n$, if any, such that $\dege_n f= \pm \infty$.   Examples  \ref{exex} and \ref{degeex} show that $\dege$ provides far more refined information than $\degl$.

The seemingly {\it ad hoc} definition of $\dege$  is motivated by the examples below.

\begin{example}\label{exex} \
\begin{enumerate}
\item By Example \ref{divisorf}(3), one has $\dege  d(n) = (0,\infty,1,-1,0,0,0,\ldots)$.
\item By Example \ref{divisorf}(4), one has $\dege  \sigma(n) = (1,0,1,0,0,0,\ldots)$
\item Let $a \in \RR$.  Since $\sigma_{-a}(n) =  \frac{\sigma_a(n)}{n^a}$ for all $n$, one has $$\dege \sigma_{-a} = \dege \sigma_a + (-a,0,0,0,\ldots).$$  In particular, since  $\deg \sigma_a = 0$ if $a \leq 0$, one has $\deg \sigma_a =a$ if $a \geq 0$.
\item Let $a \in (0,1)$.  By \cite[p.\ 122]{gron}, one has
$$\limsup_{n \to \infty} \frac{\log \frac{\sigma_a(n)}{n^a}}{\frac{(\log n)^{1-a}}{\log \log n}}= \frac{1}{1-a}.$$
(The lim sup is attained by the sequence $p_1 p_2 \cdots p_n$.)  It follows that $\dege \sigma_a = (a,\infty,1-a,-1,0,0,0,\ldots)$.  
\item Let $a \in (1,\infty)$.  By \cite[(1)]{gron}, one has
$$\limsup_{n \to \infty} \frac{\sigma_a(n)}{n^a} = \zeta(a).$$
(The lim sup is attained by the sequence $(p_1 p_2 \cdots p_n)^n$.)
It follows that $\dege \sigma_a = (a,0,0,0,\ldots)$.  
\item  By Theorem \ref{bestPNT} and statement (1) of Proposition \ref{aspropexoexp} below, one has $\dege (\li-\pi) \leq (1,-\infty,-\frac{3}{5},\frac{1}{5},0,0,0,\ldots).$
\item Let $a(n)$ denote the number of non-isomorphic abelian groups of order $n$.  In 1970, Kr\"atzel  proved \cite{krat} that
 \begin{align*}
\limsup_{n \to \infty} \frac{\log a(n)\log \log n}{\log n} = \frac{1}{4}\log 5.
\end{align*}
Consequently, one has
$$\dege  a(n) = (0,\infty,1,-1,0,0,0,\ldots).$$
\item Let $G(n)$ denote the number of non-isomorphic  groups of order $n$.  In 1993, Pyber proved \cite{pyber} that
 \begin{align*}
G(n) \leq n^{\frac{2}{27}l(n)^2+O(l(n)^{3/2})},
\end{align*}
where $l(n) \leq \log_2 n$ is the largest exponent of any prime power dividing $n$, and where for $n = p^k$ a power of a prime $p$ one has $l(n) = k = \log_p n$ and the bound above is asymptotically sharp.  It follows that
 \begin{align*}
\limsup_{n \to \infty} \frac{\log G(n)}{(\log n)^3} = \frac{2}{27(\log 2)^2}
\end{align*}
and therefore
$$\dege  G(n) = (\infty, 0, 3,0,0,0,\ldots).$$
\item Let $g(n)$ denote the largest order of an element of the symmetric group $S_n$.  In 1903, Landau proved \cite{land0} that
$$\log g(n) \sim \sqrt{n \log n} \ (n \to \infty).$$
It follows that
$$\dege g(n) = (\infty, \tfrac{1}{2}, \tfrac{1}{2}, 0, 0, 0, \ldots).$$
On the other hand, since $|S_n| = n!$ and $$\log (n!) \sim n \log n = \log (n^n) \ (n \to \infty),$$
one has
$$\dege |S_n|  = \dege (n!) = \dege (n^n) = (\infty, 1, 1, 0, 0, 0,\ldots).$$
See \cite[Section 10.10]{broughan}  and Remark \ref{landaufunction} for connections between the function $g(n)$ and the Riemann hypothesis.
\item   By \cite[Theorem 1]{debruijn}, one has
$$\sum_{n \leq x} \frac{1}{\omega(n)} = e^{(1+o(1))\sqrt{8 \log x/\log \log x}} \  (x \to \infty).$$
It follows that
$$\dege \sum_{n \leq x} \frac{1}{\omega(n)} = (0,\infty,\tfrac{1}{2}, -\tfrac{1}{2},0,0,0,\ldots).$$
\item  The {\bf Carmichael function}\index{Carmichael function $\lambda(n)$} is the arithmetic function $\lambda$ whose value $\lambda(n)$ at $n$ is the exponent of the abelian group $(\ZZ/n\ZZ)^*$, that is, $\lambda(n)$ is the smallest positive integer $k$ such that
$a^k \equiv 1   \ (\text{mod } n)$ for all $a \in \ZZ$ relatively prime to $n$.   By \cite[Theorem 3]{erdpom},  one has
$${\displaystyle \sum_{n\leq x}\lambda(n)={\frac {x^2}{\log x}}e^{C(1+o(1))\log \log x/\log \log \log x}},$$
where
$$C=e^{-\gamma }\prod_{p }\left({1-{\frac {1}{(p-1)^{2}(p+1)}}}\right) =  0.34537\ldots.$$
It follows that
$$\dege \sum_{n\leq x}\lambda(n) = (2,-1,\infty,1,-1,0,0,0,\ldots).$$
\end{enumerate}
\end{example}

\begin{example}\label{degeex} Let $c,d \in \RR$ with $c \neq 0$, and let $n$ be a positive integer.  One has the following.
\begin{enumerate}
\item $\dege  e^{x^2} \sin x = (\infty, 2, 0, 0, 0, \ldots)$.
\item $ \dege  e^{-x^2} = (-\infty, -2,0,0,0,\ldots)$.
\item $\dege   x^x = (\infty, 1,1,0,0,0,\ldots)$.
\item $ \dege  e^{e^x} = (\infty, \infty,1,0,0,0,\ldots)$.
\item $ \dege  e^{-e^x} = (-\infty,- \infty,-1,0,0,0,\ldots)$.
\item $\dege e^{r(x)} = (\infty, \degl_0 r, \degl_1 r, \degl_2 r, \ldots)$ for all $r \in \mathfrak{L}$ with $\degl r > (0,1,0,0,0,\ldots)$.  
\item $\dege e^{r(\log x)} = (0, \infty, \degl_1 r, \degl_2 r, \degl_3 r, \ldots)$ for all $r \in \mathfrak{L}$ with $$(0,1,0,0,0,\ldots)< \degl r < (1,0,0,0,\ldots).$$
\item $\dege e^{-r(x)} = (-\infty, -\degl_0 r, -\degl_1 r, -\degl_2 r, -\degl_3 r, \ldots)$ for all $r \in \mathfrak{L}$ with $\degl r > (0,1,0,0,0,\ldots)$.
\item $\dege e^{-r(x)} = (0,-\infty, -\degl_0 r, -\degl_1 r, -\degl_2 r, -\degl_3 r, \ldots)$ for all $r \in \mathfrak{L}$ with $$(0,1,0,0,0,\ldots)< \degl r < (1,0,0,0,\ldots).$$
\item $\dege   e^{c(\log x)^d} = (\operatorname{sgn} c)(\infty,0, d,0,0,0,\ldots)$ if $d > 1$.
\item $\dege   e^{c(\log x)^d} =  (\operatorname{sgn} c)(0,\infty, d,0,0,0,\ldots)$ if $0<d < 1$.
\item $\dege   e^{c(\log x)^d} = (0,0, 0,\ldots)$ if $d\leq 0$.  
\item $\dege \exp^{\circ \lfloor x \rfloor}x = (\infty,\infty,\infty,\ldots)$.
\item $\dege 1/\exp^{\circ \lfloor x \rfloor}(x) = (-\infty,-\infty,-\infty,\ldots)$.
\item $\dege  \exp^{\circ n }((\log^{\circ n} x)^d) = (\infty, 0,\infty,0,\ldots, \infty,0,d,0,0,0,\ldots)$ if $d > 1$, where ``$d$'' is preceded by $n$ occurrences of ``$\infty,0$.''
\item $\dege  \exp^{\circ \lfloor x \rfloor}((\log^{\circ \lfloor x \rfloor} x)^d)= (\infty,0,\infty,0,\infty,0\ldots)$ if $d >1$.
\item $\dege  x^{\log^{\circ n}x} = (\infty,0,1,0,0,0,\ldots,1,0,0,0\ldots)$ if $n \geq 2$, where the second ``$1$'' is preceded by $n-2$ zeros.
\item Let $f(x)$ be the function that equals ${\log^{\circ N}x}$ on the interval $[\exp^{\circ N}(N),\exp^{\circ (N+1)}(N+1))$ for all positive integers $N$.  Then $\lim_{x \to\infty} f(x) = \infty$ and $f(x) = o(\log^{\circ N} x) \ (x \to \infty)$ for all $N$, and therefore $\degl f = \dege f = (0,0,0,\ldots)$.  It also follows that $f$ has exact logarithmic degree and therefore $\degl \frac{1}{f} = \dege \frac{1}{f} = (0,0,0,\ldots)$.
\item Let $f$ be any function such that $\lim_{x \to \infty} f(x) = \infty$ and $\dege f = \dege \frac{1}{f} = (0,0,0,\ldots)$.  Also, let $r \in \RR^{\RR_\infty}$ be such that $r$ is eventually positive (resp., eventually negative) and $r(x) \neq o(\log x) \ (x \to \infty)$.  One has the following.
\begin{enumerate}
\item $\dege e^{r(x)f(x)} = \pm( \infty,\dege_0 r,\dege_1 r,\dege_2 r,\ldots)$.
\item $\dege e^{r(\log x)f(\log x)} = \pm(0, \infty,\dege_0 r,\dege_1 r,\dege_2 r,\ldots)$ if $\deg r < 1$.
\item $\dege e^{cxf(x)} =  (\operatorname{sgn} c)(\infty,1,0,0,0,\ldots)$.
\item $\dege x^{cf(x)} =  (\operatorname{sgn} c)(\infty,0,1,0,0,0,\ldots)$.
\item $\dege \, (\log x)^{cf(\log x)} =  (\operatorname{sgn} c) (0,\infty,0,1,0,0,0,\ldots)$.
\item $\dege x^{c/f(x)} =  (\operatorname{sgn} c)(0,\infty,1,0,0,0,\ldots)$.
\end{enumerate}
\item Let $f(x)$ be the function that on the interval $[N,N+1)$  assumes the value $e^x$ for even integers $N$ and $e^{-x^2}$ for odd integers $N$.  One has $\dege f = (\infty, 1,0,0,0,\ldots)$, since the function $ \max( \log |f(x)|, 0) $ has logexponential degree $(1,0,0,0,\ldots)$.  The function $\log |f(x)|$, on the other hand, has logexponential degree $(2,0,0,0,\ldots)$. 
\item Let $f(x)$ be the function that on the interval $[N,N+1)$  assumes the value $e^{-x}$ for even integers $N$ and $e^{-x^2}$ for odd integers $N$.  One has $\dege f = (-\infty, -1,0,0,0,\ldots)$, since the function $-\frac{1}{ \log |f(x)|} $ has logexponential degree $(-1,0,0,0,\ldots)$.  The function $\log |f(x)|$, on the other hand, has logexponential degree $(2,0,0,0,\ldots)$. 
\end{enumerate}
\end{example}

It is natural to ask why we do not simply define $f_{(k+1)}(x) = \log|f_{(k)}(x)|$ when $\deg f =  \pm \infty$.   Examples (20) and (21) above indicate why this is not advisable.  Another reason is that few if any of the properties of $\dege$ proved in this section hold under that revised definition.  Note that, if $\dege_k f= \pm \infty$, then for all $n > k$ one has $f_{(n)}(x) \geq 0$ for all $x \gg 0$.  Moreover, if  $\dege_k f= -\infty$, then $\lim_{x \to \infty} f_{(k)}(x) = \lim_{x \to \infty} f_{(k+1)}(x) = 0$, while, if $\dege_k f = \infty$, then
$\limsup_{x \to \infty} f_{(k)}(x) = \limsup_{x \to \infty} f_{(k+1)}(x)= \infty$.

Let $T = (-)_{(1)}:  \RR^{\RR_\infty} \longrightarrow   \RR^{\RR_\infty}$ act by $f \longmapsto f_{(1)}$, so that
$$T(f)(x)=   \left.
 \begin{cases}
    f(e^x) e^{-(\deg f) x}& \text{if } \deg f  \neq \pm \infty \\
    \max( \log |f(x)|, 0)  & \text{if } \deg f = \infty \\
 \displaystyle   -\frac{1}{\log |f(x)|} & \text{if } \deg f =- \infty.
 \end{cases}
\right.$$
for all  $f \in \RR^{\RR_\infty}$.   Then, for all  $f \in \RR^{\RR_\infty}$, one has
$$f_{(k)} = T^{\circ k}(f)$$
for all $k$, and, in particular,  one has
$$(f_{(k)})_{(l)} = f_{(k+l)}  $$
for all $k$ and $l$.   Note, then, that
$$\dege f = (\deg f, \deg T(f), \deg T(T(f)), \ldots),$$
that is, 
$$\dege_k f = \deg T^{\circ k}(f)$$
for all $k$.   This is in analogy with the definition of a regular continued fraction expansion of a real number, where the operator $T$ is the analogue of   the continued fraction operator $S: (-\infty,\infty]  \longrightarrow (1,\infty] $ given by
$$S(x)=   \left.
 \begin{cases}
   \frac{1}{x-\lfloor x \rfloor} & \text{if } x \notin \ZZ\cup \{ \infty\}   \\
   \infty  & \text{if } x \in \ZZ\cup \{ \infty\},
 \end{cases}
\right.$$
and where $\deg$ is the analogue of   the floor function $\lfloor  x \rfloor$.  Indeed,   the regular continued fraction expansion of any $x \in \RR$ is given by 
$$ [\lfloor  x \rfloor, \lfloor S(x) \rfloor,    \lfloor S(S(x)) \rfloor,    \lfloor S(S(S(x)))\rfloor, \ldots],$$
where the notation above  for regular continued fractions  is as defined in Section 13.1.
Just as $\lfloor  x \rfloor$ is an integer measure of the  ``size'' of $x \in\RR$,  so $\deg f$ is an extended real number measure of the  ``size'' of $f \in \RR^{\RR_\infty}$.    The analogy here is not perfect, since any real number can be recovered from its regular continued fraction expansion (as the limit of its sequence of convergents).    Nevertheless, it shows that the definition  of $\dege$ is completely determined by the maps $T$ and $\deg$.

Suppose now that $\dege f$ ends in a trail of $0$s,  i.e., that there exists an $n \geq 1$ such that $\dege f_{(n)} = (0,0,0,\ldots)$.  Let $$\operatorname{cf}^+ f = \limsup_{x \to \infty} f_{(n)} (x) \in [0,\infty]$$ and $$\operatorname{cf}^- f = \liminf_{x \to \infty} f_{(n)} (x) \in [0,\infty],$$ neither of which depends on the choice of $n \gg 0$.     These invariants can be used to encode all of  the explicit constants in Examples \ref{divisorf}, \ref{newex},   and \ref{exex}.  For example, one has $\operatorname{cf}^+ d(n) = \log 2$,  $\operatorname{cf}^+ \sigma(n) = e^\gamma$,  $\operatorname{cf}^+ \sigma_a(n)$ is equal to  $\frac{1}{1-a}$ for $a \in (0,1)$ and to $\zeta(a)$ for $a \in (1,\infty)$,  $\operatorname{cf}^+ \Omega(n) = \frac{1}{ \log 2}$,   $\operatorname{cf}^+ \omega(n) = 1$,  $\operatorname{cf}^+ a(n) = \frac{1}{4}\log 5$,  $\operatorname{cf}^+ G(n)  = \frac{2}{27(\log 2)^2}$,   $\operatorname{cf}^{\pm} g(n) = 1$,   $\operatorname{cf}^{\pm} \sum_{n \leq x} \frac{1}{\omega(n)} = \sqrt{8}$, and $\operatorname{cf}^+ \frac{1}{\phi(n)} = e^\gamma$.

Using {\it  transseries}  \cite[Appendix A]{asch}  \cite{dmm2} and recursion,  it is possible to find,  given  any $\dd \in \prod_{n = 0}^\infty \overline{\RR}$ with $\dd = \dege f$ for some $f \in \RR^{\RR_\infty}$,  a canonical or prototypical germ $F$ at $\infty$ with $\dd = \dege F$: see Corollary \ref{prototype}.   For example, the germ associated to $(1,-\infty,-\tfrac{3}{5},\tfrac{1}{5},0,0,0,\ldots)$ is  $xe^{ -(\log x)^{3/5}(\log \log x)^{-1/5}}$.

By Corollary \ref{maindegree}, one has
$$\Theta = \deg(\li-\pi).$$
In particular,  one has $\frac{1}{2} \leq \deg(\li-\pi) \leq 1$, and the Riemann hypothesis holds if and only if $\deg(\li-\pi) = \frac{1}{2}$.
Thus, the problem of computing $\Theta$ and thereby settling the Riemann hypothesis  generalizes to the following even more difficult problem, which is discussed in Chapter 9 and Section 14.2.

\begin{outstandingproblem}\label{mainproblem00}
Compute $\dege(\li - \pi)$.
\end{outstandingproblem}

One of our main applications of logexponential degree to analytic number theory, in Part 3, is to relate $\dege(\li - \pi)$ to $\dege f$ for several other functions $f$ arising in number theory.

\begin{remark}[Functions of degree $\pm \infty$]
Note that $\dege$ can distinguish between $x$ and $x\log x$,  but it is not refined enough to distinguish between $e^x$ and $xe^x$,  since the latter two functions both have logexponential degree $(\infty,1,0,0,0,\ldots)$.   To better deal with functions of degree $\pm \infty$,  we define the {\bf exponentiallog degree $\operatorname{eldeg} f$} of $f$ as follows.  For any $f \in \RR^{\RR_\infty}$ with $\infty_f := \deg f = \pm \infty$ and $\dege f \neq \pm (\infty, \infty,\infty,\ldots)$,  there is a smallest $n \in \ZZ_{> 0}$ such that $\deg ( f \circ \log^{\circ n}) \neq \infty_f$, and then $\deg ( f \circ \log^{\circ n})   \in (\operatorname{sgn}\infty_f)\RR_{\geq 0}$.    Let $\dd = \dege (f\circ \log^{\circ n})$, so that $\dd_0 \in \RR$.   We let
$$\operatorname{eldeg} f  = (\infty_f, \infty_f, \ldots, \infty_f; \dd_0, \dd_1, \dd_2,\ldots),$$
where the semicolon is preceded by $n$ occurrences of $\infty_f$.  For all other $f \in \RR^{\RR_\infty}$,  that is, for all $f \in \RR^{\RR_\infty}$ with $\deg f \neq \pm \infty$ or $\dege f  = \pm (\infty, \infty,\infty,\ldots)$, we let $\operatorname{eldeg} f = \dege f$.  Thus, for example, we have
$$\operatorname{eldeg} xe^{2x} = (\infty; 2,1,0,0,0,\ldots),$$
$$\operatorname{eldeg} xe^{2x}e^{-e^x} = (-\infty,-\infty;-1,2,1,0,0,0,\ldots),$$
and 
$$\operatorname{eldeg} e^{(\log x)^2} = \operatorname{eldeg} xe^{(\log x)^2} = (\infty;0,\infty,0,2,0,0,0,\ldots).$$
One can show that $\dege f$ can  always be recovered from $\operatorname{eldeg} f$.  As the third example exhibits, any infinities occurring after the semicolon introduce the same level of coarseness in $\operatorname{eldeg} f$ as in $\dege f$.   Moreover,  $\operatorname{eldeg}$ depends on the choice of $e$ as a base, since, for example,  one has
$$\operatorname{eldeg} b^{x} = (\infty; \log b, 0,0,0,\ldots)$$
for all $b > 1$,  whereas $\dege$ is base independent, in the sense of Remark \ref{altbase}.
\end{remark}

\subsection{Relationships with $O$, $o$, and $\asymp$}

\begin{lemma}\label{reslemma}
Let $f \in \RR^{\RR_\infty}$ and $X \subseteq \dom f$ with $\infty \in \overline{X}$.  Then one has
$$\dege f|_X \leq \dege f.$$
\end{lemma}

\begin{proof}
This follows by induction and a property of $\limsup$, namely, that
$$\limsup_{x \to \infty} g|_X(x) \leq \limsup_{x \to \infty} g(x)$$
for all $g \in \RR^{\RR_\infty}$ and $X \subseteq \dom g$ with $\infty \in \overline{X}$.
\end{proof}

The following result is an analogue of Proposition \ref{asprop} for $\dege$.

\begin{proposition}\label{aspropexoexp}
Let $f, g \in \RR^{\RR_\infty}$, and let $\underline{g} = g|_{\dom f \cap \dom g}$.
\begin{enumerate}
\item If $f(x) = O(g(x)) \ (x \to \infty)$, then $\dege f \leq \dege \underline{g} \leq \dege g$.
\item If $f(x) \gg g(x) \ (x \to \infty)$, then $\dege f \geq \dege \underline{g}$.
\item If $f(x) \asymp g(x) \ (x \to \infty)$, then $f_{(k)}(x) \asymp \underline{g}_{(k)}(x) \ (x \to \infty)$ for all $k$ and $\dege f = \dege \underline{g}$.
\item If $f$ is eventually bounded, then  $\dege f \leq (0,0,0,\ldots)$.
\item If $\dege f < (0,0,0,\ldots)$, then $\displaystyle \lim_{x \to \infty} f(x)  = 0$.
\end{enumerate}
\end{proposition}

\begin{proof}
Statements (4 and (5) follow from Propositions  \ref{asprop} and \ref{degeprop0}.

By the lemma, to prove (1), we may suppose without loss of generality that $\dom f = \dom g$.  Suppose that $f(x) = O(g(x)) \ (x \to \infty)$.   If $\dege_k f$ and $\dege_k g$ are finite for all $k$, then statement (1) follows from the proof of Proposition \ref{asprop}.   Thus we may suppose that there is a smallest $n$ such that $\dege_n f$ and $\dege_n g$ are not both finite and that $\dege_k f = \dege_k g$ for all $k < n$.  Clearly then one has $d: = \dege_n f \leq \dege_n g$, and we may suppose without loss of generality that equality holds.  One has $f_{(k)}(x) = O(g_{(k)}(x)) \ (x \to \infty)$ for all $k \leq n$.  Thus, there exists  a $c >0$ such that $|f_{(n)}(x)| \leq e^c|g_{(n)}(x)|$, and therefore  $\log |f_{(n)}(x)| \leq c + \log |g_{(n)}(x)|$, for all $x \gg 0$.

Suppose that $d = \infty$.   Then  $f_{(n+1)}(x) =  \max(\log |f_{(n)}(x)|,0)$ and $\limsup_{x \to \infty} f_{(n+1)}(x) = \infty$, so that $\dege f_{(n+1)} \geq (0,0,0,\ldots)$, and likewise for $g_{(n+1)}$.  Now, one has
\begin{align*}
0 \leq f_{(n+1)}(x)  & = \max(\log |f_{(n)}(x)|,0) \\
 & \leq  \max (c+\log |g_{(n)}(x)|,0) \\
 & \leq \max (c+\log |g_{(n)}(x)|,c) \\
& = c+ g_{(n+1)}(x)
\end{align*}
for all $x \gg 0$.  It follows that $f_{(n+1)}(x)  = O(g_{(n+1)}(x)) \ (x \to \infty)$.  Suppose, on the other hand, that $d = -\infty$. Then  one has
\begin{align*}
0 \leq f_{(n+1)}(x)  & =-\frac{1}{\log |f_{(n)}(x)|} \\
 & \leq-\frac{1}{c+\log |g_{(n)}(x)|}  \\
 & =  -\frac{1}{\log |g_{(n)}(x)|}\left(1 - \frac{c}{c+\log |g_{(n)}(x)|} \right)\\
 & \leq-\frac{2}{\log |g_{(n)}(x)|}  \\
 & = 2 g_{(n+1)}(x)
\end{align*}
for all $x \gg 0$.  Therefore, again one has $f_{(n+1)}(x)  = O(g_{(n+1)}(x)) \ (x \to \infty)$.  Thus, statement (1) follows by induction.  Finally,  statements (2) and (3) follow from statement (1) and its proof.
\end{proof}

We think of $\dege f \leq \dege g$ as a total preorder that is substantially coarser than the preorder $f(x) = O(g(x)) \ (x \to \infty)$  and substantially finer than the total preorder $\degl f \leq \degl g$ (although the three preorders agree on $\RR(\mathfrak{L})$). 

Let $f \in  \RR^{\RR_\infty}$ and $n$ a nonnegative integer.  We say that $f$ has {\bf exact logexponential degree to order $n$}\index{exact logexponential degree} if $f_{(k)}$ has exact degree for all $k \leq n$, where the $f_{(k)}$ are defined recursively as in the definition of $\dege f$.  We say that $f$ has {\bf exact logexponential degree} if $f$ has exact  logexponential degree to order $n$ for all $n$.    For example, if $f$ has finite and exact logarithmic degree to order $n$, then $f$ has finite and exact logexponential degree to order $n$.  If $f$ has exact logexponential degree, then $f$ has exact logarithmic degree, but the converse does not hold.

The following proposition extends Lemma \ref{reslemma}.

\begin{proposition}\label{arithbridge2}
Let $f \in \RR^{\RR_\infty}$ and $X \subseteq \dom f$ with $\infty \in \overline{X}$, and let $n$ be a nonnegative integer.  If $f$ has exact logexponential degree to order $n$, then so does $ f|_X$, and one has
$$\dege_k f|_X = \dege_k f$$
for all $k \leq n$.  Thus, if $f$ has exact logexponential degree, then so does $ f|_X$, and one has
$$\dege f|_X = \dege f.$$
\end{proposition}

\begin{proof}
Let $g =  f|_X$.   Since $f$ has exact degree, one has
$$d := \deg f = \lim_{x \to \infty}  \frac{\log |f(x)|}{\log x} = \lim_{x \to \infty} \frac{\log |g(x)|}{\log x} = \deg g,$$
and thus $g$ has exact degree.  Suppose that $d \neq \pm \infty$.  Then $f_{(1)}(x) = f(e^x)e^{-dx}$ and $g_{(1)}(x) = g(e^x)e^{-dx}$, and the latter function is a restriction of the former.  If $f_{(1)}$ has exact degree, then one has
$$\deg f_{(1)} =  \lim_{x \to \infty}  \frac{\log| f_{(1)}(x)|}{\log x}  =  \lim_{x \to \infty}   \frac{\log |g_{(1)}(x)|}{\log x}  = \deg g_{(1)},$$
and thus $g_{(1)}$ has exact degree.  A similar argument establishes the same equalities when $d = \infty$ or $d = -\infty$.  An obvious inductive argument, then, yields $$\dege_k f = \deg f_{(k)} = \deg g_{(k)} = \dege_k g$$ for all nonnegative integers $k \leq n$.
\end{proof}

Let $f \in \RR^{\RR_\infty}$.  Since $\limsup_{x \to \infty} f(x) = \lim_{x \to \infty} f(x)$ if and only if $\limsup_{x \to \infty} f(x) = \limsup_{x \to \infty} f|_X(x)$ for all  $X \subseteq \dom f$ with $\infty \in \overline{X}$, the converse of Proposition \ref{arithbridge2} also holds, i.e., one has the following.

\begin{corollary}
Let $f \in \RR^{\RR_\infty}$, and let $n$ be a nonnegative integer.   Then $f$  has exact logexponential degree to order $n$ if and only if $\dege_k f|_X = \dege_k f$ for all $X \subseteq \dom f$ with $\infty \in \overline{X}$ and all $k \leq n$.   Thus, $f$  has exact logexponential degree if and only if $\dege f|_X = \dege f$ for all $X \subseteq \dom f$ with $\infty \in \overline{X}$.
\end{corollary}

\begin{corollary}\label{exact0help}
Let $f \in \RR^{\RR_\infty}$, and suppose that $f$ is not eventually nonzero.  Then $f$ has exact  logexponential degree if and only if $\dege f = (-\infty,-\infty,-\infty,\ldots)$.    More generally, $f$ has exact logexponential degree to order $n \geq 0$ if and only if $\dege_k f = -\infty$ for all $k \leq n$.
\end{corollary}

\begin{proposition}\label{exactas}
Let $f, g\in \RR^{\RR_\infty}$ with $f(x) \asymp g(x) \ (x \to \infty)$, and let $n$ be a nonnegative integer.  If $g$ has exact logexponential degree to order $n$, then so does $f$, and one has $\dege_k f = \dege_k g$ for all $k \leq n$.
\end{proposition}

\begin{proof}
Let $\underline{g} =  g|_{\dom f \cap \dom g}$.
By  Proposition \ref{aspropexoexp}(3), one has $f_{(k)}(x) \asymp \underline{g}_{(k)}(x) \ (x \to \infty)$ for all $k$.  Moreover, by Proposition \ref{arithbridge2}, the function $\underline{g}_{(k)} = g_{(k)}|_{\dom f_{(k)} \cap \dom g_{(k)}}$ has exact degree, whence, by Proposition  \ref{firstprop1}(3), the function $f_{(k)}$ also has exact degree, for all $k \leq n$.  Finally, again by Propositions 
 \ref{arithbridge2} and   \ref{firstprop1}(3), one has $\dege_k f = \dege_k \underline{g} = \dege_k g$  for all $k \leq n$.
\end{proof}

 Let $f \in \RR^{\RR_\infty}$,  and define $f_{\{k\}}$ and $\underline{\dege}_k f =\underline{ \deg } \, f_{\{k\}} \in \overline{\RR}$ for all nonnegative integers $k$, recursively, as follows.    Let $f_{\{0\}} = f$.  Suppose that $f_{\{k\}}$ is defined, and set $d_k = \underline{\dege}_k f =\underline{\deg} \, f_{\{k\}}.$  We then let
$$f_{\{k+1\}}(x) =   \left.
 \begin{cases}
    f_{\{k\}}(e^x) e^{-d_k x}& \text{if } d_k \neq \pm \infty \\
     \log |f_{\{k\}}(x)| & \text{if } d_k = \infty \\
 \displaystyle   -\frac{1}{\min(\log |f_{\{k\}}(x)|,0)} & \text{if } d_k =- \infty.
 \end{cases}
\right.$$
We also set $$\underline{\dege}\, f = (\underline{\dege}_0 f, \underline{\dege}_1 f, \underline{\dege}_2 f, \ldots) \in \prod_{n = 0}^\infty \overline{\RR}.$$   We call $\underline{\dege}\, f$ the {\bf lower logexponential degree of $f$}\index{lower logexponential degree $\underline{\dege}$}\index{lower exponential degree $\underline{\dege}$}\index[symbols]{.g g@$\underline{\dege}\, f$}  and $\underline{\dege}_k f$ the {\bf lower logexponential degree of $f$ of order $k$}.

For all $f \in \RR^{\RR_\infty}$, one has
$$\limsup_{x \to \infty}\, (- f(x)) = -\liminf_{x \to \infty} f(x)$$
and
$$\liminf_{x \to \infty}\, (- f(x)) = -\limsup_{x \to \infty} f(x),$$
and  $\lim_{x \to \infty} f(x)$ exists or is $\pm \infty$ if and only if $\limsup_{x \to \infty} f(x) = \liminf_{x \to \infty} f(x)$.   From these basic facts, and from the definitions of degree, lower degree,  $f_{(k)}$, and $f_{\{k\}}$, one readily deduces the following.

\begin{proposition}
 Let $f \in \RR^{\RR_\infty}$.  If $f$ is not eventually nonzero, then $\underline{\dege}\, f = (-\infty, -\infty, -\infty,\ldots)$.   Suppose, on the other hand,  that $f$ is eventually nonzero.   Then one has the following.
 \begin{enumerate}
\item $f_{\{k\}} = 1/(1/f)_{(k)}$ for all $k$.
\item $\underline{\dege}\, f = - \dege(1/f)$.
\item $\underline{\dege}\, f \leq \dege f$, and equality holds if and only if $f$ has exact logexponential degree.
\item $f$ has exact logexponential degree to order $n \geq 0$ if and only if $\dege_k f = \underline{\dege}_k f$ for all $k \leq n$.
\end{enumerate}
\end{proposition}

\begin{corollary}\label{exactloginverse}
Let $r \in \RR^{\RR_{\infty}}$, and let $n$ be a nonnegative integer.   
\begin{enumerate}
\item If $r$ has exact logexponential degree to order $n$ and is eventually nonzero, then, for all nonnegative integers $k \leq n$, one has
$$r_{(k+1)}(x) =   \left.
 \begin{cases}
    r_{(k)}(e^x) e^{-(\deg r_{(k)}) x}& \text{if } \deg r_{(k)} \neq \pm \infty \\
    \log r_{(k)}(x)  & \text{if } \deg r_{(k)} = \infty \\
 \displaystyle   -\frac{1}{\log r_{(k)}(x)} & \text{if } \deg r_{(k)} =- \infty
 \end{cases}
\right.$$ and $(1/r)_{(k)} = 1/r_{(k)}$, and therefore $\dege_k (1/r) = -\dege_k r$.
\item $r$ has exact logexponential degree to order $n$  if and only if  $\dege_k r = -\infty$ for all $k \leq n$ or  $r$ is eventually nonzero and  $\dege_k (1/r) = -\dege_k r$ for all $k \leq n$.
\item $r$ has exact logexponential degree if and only if $\dege r = (-\infty, -\infty,-\infty, \ldots)$ or $r$ is eventually nonzero and  $\dege(1/r) = -\dege r$.
\end{enumerate}
\end{corollary}

\begin{corollary}\label{aspropexoexpunder}
Let $f, g \in \RR^{\RR_\infty}$, and let $\underline{g} = g|_{\dom f \cap \dom g}$.
\begin{enumerate}
\item If $f(x) \gg g(x) \ (x \to \infty)$, then $\underline{\dege} \,  f \geq  \underline{\dege} \,  \underline{g} \geq  \underline{\dege} \,  g$.
\item If $f(x) = O(g(x)) \ (x \to \infty)$, then $\underline{\dege} \,  f \leq \underline{\dege} \,  \underline{g}$.
\item If $f(x) \asymp g(x) \ (x \to \infty)$, then $f_{\{k\}}(x) \asymp \underline{g}_{\{k\}}(x) \ (x \to \infty)$ for all $k$ and $\underline{\dege} \,  f = \underline{\dege} \,  \underline{g}$.
\item If $f$ is eventually bounded away from $0$,  then  $\underline{\dege} \,  f \geq (0,0,0,\ldots)$.
\item If $\underline{\dege} \,  f > (0,0,0,\ldots)$, then $\displaystyle \lim_{x \to \infty} |f(x) | = \infty$.
\end{enumerate}
\end{corollary}

\begin{example}
By \cite[Theorem 328]{har}, one has
$$ \liminf_{n \to\infty} \frac{\phi(n) \log \log n}{n} =e^{-\gamma },$$
where $\phi(n)$ is Euler's totient.  It follows that
$$\underline{\dege} \, \phi(n) = (1,0,-1,0,0,0,\ldots).$$
Moreover, it is easy to check that
$$ \limsup_{n \to\infty} \frac{\phi(n)}{n} =1,$$
and therefore 
$$\dege \phi(n) = (1,0,0,0,\ldots).$$
\end{example}

The following proposition yields, in the subsequent corollary,  a useful property of functions of exact logexponential degree.

\begin{proposition}\label{oexpprop}
Let $f,g \in \RR^{\RR_\infty}$,  where  $\dom g$ contains the intersection of $\dom f$ with some neighborhood of $\infty$,  and let $n$ be a nonnegative integer.
\begin{enumerate}
\item If the smallest nonnegative integer $N$ such that $\dege_N f \neq \underline{\dege}_N g$ exists and is at most $n$,  where also $\dege_N f < \underline{\dege}_N g$,  then $f(x) = o(g(x)) \ (x \to \infty)$.  
\item  $\dege f <\underline{\dege}\, g$, then $f(x) = o(g(x)) \ (x \to \infty)$.  
\end{enumerate}
\end{proposition}

\begin{proof}
It suffices to prove statement (1).   Note that the  function  $g_{(k)}$ is eventually nonzero, for each $k \leq n$.  For all $k < N$, one has $d_k := \dege_k f = \underline{\dege}_k g$.  Moreover, one has $\deg f_{(N)} < \underline{\deg}\, g_{\{N\}}$, so that $f_{(N)}(x) = o(g_{\{N\}}(x)) \ (x \to \infty)$, by Proposition \ref{firstprop2}(7).    Let $m$ be least such that $f_{(m)}(x) = o(g_{\{m\}}(x)) \ (x \to \infty)$, so that $m \leq N$.
We claim that $m = 0$ and therefore  $f(x) = o(g(x)) \ (x \to \infty)$.  Suppose to obtain a contradiction that $m \geq 1$.  Suppose first that $d_{m-1} \neq \pm \infty$.  Then  $g_{\{m-1\}}(x) = x^{d_{m-1}}g_{\{m\}}(\log x)$ and therefore
$$f_{(m-1)}(x)= x^{d_{m-1}}f_{(m)}(\log x) = o( g_{\{m-1\}}(x)) \ (x \to \infty).$$  Suppose, on the other hand, that $d_{m-1} = -\infty$.  Then $\lim_{x \to \infty} \frac{\log |f_{(m-1)}(x)|}{\log x} = -\infty$ and thus
$$\lim_{x \to \infty} \frac{1}{ f_{(m)}(x)} = - \lim_{x \to \infty} \log |f_{(m-1)}(x)| = \infty,$$
so that
\begin{align*}
\frac{|f_{(m-1)}(x)|}{|g_{\{m-1\}}(x)|} \leq |f_{(m-1)}(x)|\max\left(\left| \frac{1}{g_{\{m-1\}}(x)} \right |,1\right) & =  \exp\left(-\frac{1}{f_{(m)}(x)}+\frac{1}{g_{\{m\}}(x)} \right) \\
 & =  \exp\left(\left(\frac{f_{(m)}(x)}{g_{\{m\}}(x)}-1 \right) \frac{1}{f_{(m)}(x)}\right) \\ & \to \exp(-1 \cdot \infty) = 0 \\
\end{align*}
as $x \to \infty$, and therefore   $f_{(m-1)}(x)= o( g_{\{m-1\}}(x)) \ (x \to \infty)$.  
Finally, suppose that $d_{m-1}= \infty$.  Then $\lim_{x \to \infty} \frac{\log |g_{\{m-1\}}(x)|}{ \log x} = \infty$ and $g_{\{m\}}(x)  = \log |g_{\{m-1\}}(x)|$
and thus
$$\lim_{x \to \infty} g_{\{m\}}(x) =  \lim_{x \to \infty} \log |g_{\{m-1\}}(x)|= \infty,$$
so that
\begin{align*}
\frac{|f_{(m-1)}(x)|}{|g_{\{m-1\}}(x)|} \leq \frac{\max(|f_{(m-1)}(x)|,1)}{|g_{\{m-1\}}(x)|} & =  \exp\left(f_{(m)}(x)-g_{\{m\}}(x) \right) \\
 & =  \exp\left(\left(\frac{f_{(m)}(x)}{g_{\{m\}}(x)}-1 \right) g_{\{m\}}(x)\right) \\ & \to \exp(-1 \cdot \infty) = 0, \\
\end{align*}
whence $f_{(m-1)}(x)= o( g_{\{m-1\}}(x)) \ (x \to \infty)$.  Thus, we have shown that $$f_{(m-1)}(x)= o( g_{\{m-1\}}(x)) \ (x \to \infty),$$ in all three cases.  This contradicts the minimality of $m$.  It follows that $m = 0$, as claimed.
\end{proof}

\begin{corollary}\label{oexp}
Let $f,g \in \RR^{\RR_\infty}$,  where  $\dom g$ contains the intersection of $\dom f$ with some neighborhood of $\infty$,  and let $n$ be a nonnegative integer.
\begin{enumerate}
\item If $g$ has exact logexponential degree to order $n$, and  if the smallest nonnegative integer $N$ such that $\dege_N f \neq \dege_N g$ exists and is at most $n$,  where also $\dege_N f < \dege_N g$,  then $f(x) = o(g(x)) \ (x \to \infty)$.  
\item If  $g$ has exact logexponential degree and $\dege f < \dege g$, then $f(x) = o(g(x)) \ (x \to \infty)$.  
\end{enumerate}
\end{corollary}

\begin{corollary}\label{oexpcor0}
Let $f, g \in \RR^{\RR_\infty}$, where $\dom g$ contains the intersection of $\dom f$ with some neighborhood of $\infty$, and let  $\underline{g} = g|_{\dom f \cap \dom g}$.   Suppose that $ f(x) = O(g(x)) \ (x \to \infty)$ but $ f(x) \neq o(g(x)) \ (x \to \infty)$.  Then one has
$$\underline{\dege} \, f \leq  \underline{\dege} \, \underline{g} \leq \dege f \leq \dege \, g.$$
Consequently, if $g$ has exact  logexponential degree, then $\dege f = \dege g$.  Likewise, if $f$ has exact logexponential degree, then $\dege f =  \underline{\dege} \, \underline{g}$.
\end{corollary}

\begin{corollary}\label{oexpcor}
Let $f, g \in \RR^{\RR_\infty}$, where $g$ has exact logexponential degree and $\dom g$ contains the intersection of $\dom f$ with some neighborhood of $\infty$.  One has
\begin{align*}
\degl f < \degl g  \quad & \Longrightarrow \quad \dege f  < \dege g  \\
 \quad & \Longrightarrow \quad  f(x) = o(g(x)) \ (x \to \infty) & \\
 \quad & \Longrightarrow \quad f(x) = O(g(x)) \ (x \to \infty)  \\ 
 \quad & \Longrightarrow \quad \dege f  \leq \dege g \\
\quad & \Longrightarrow \quad \degl f \leq \degl g.
\end{align*}
In particular, if $ f(x) = O(g(x)) \ (x \to \infty)$ but $ f(x) \neq o(g(x)) \ (x \to \infty)$, then $\dege f = \dege g$.
\end{corollary}

Note that none of the implications in the corollary above are reversible.

\begin{remark}[Leading coefficient and logexponential degree]
Assume the notation of Remark \ref{lc}.  Let $f \in \RR^{\RR_\infty}$  with $\deg f = d \neq \pm \infty$.  If $\operatorname{lc} f < \infty$, that is, if $f(x) = O(x^{d}) \ (x \to \infty)$, then $\dege f \leq (d,0,0,0,\ldots)$.  On the other hand, if $\operatorname{lc} f > 0$, that is, if $f(x) \neq  o(x^{d}) \ (x \to \infty)$, then $\dege f \geq (d,0,0,0,\ldots)$.    Consequently, if $0< \operatorname{lc} f < \infty$, then $\dege f  = (d,0,0,0,\ldots)$. 
\end{remark}

\subsection{Relationships with logarithmico-exponential functions}

The following proposition provides constraints on what sequences in $\prod_{n = 0}^\infty\overline{\RR}$ can be of the form $\dege f$ for some $f \in \RR^{\RR_\infty}$.   By Theorem \ref{degeequiv} below,  these constraints are exhaustive.

\begin{proposition}\label{degeprop}
Let  $f \in \RR^{\RR_\infty}$, and let the $f_{(k)}$ be defined recursively as in the definition of $\dege$, so that $\dege_k f = \deg f_{(k)}$ for all $k$.  For all $n \geq 0$, one has the following.
\begin{enumerate}
\item Suppose that $\dege_n f = \infty$.  Then one has $f_{(n+1)}(x) \neq O(\log x) \ (x \to \infty)$ and therefore
$$\dege_{(n+1)} f \geq (0,1,0,0,0,\ldots)$$
and
\begin{align*}
\dege f \geq (\dege_0 f, \ldots, \dege_{n} f, 0,1, 0,0,0,\ldots).
\end{align*}
\item Suppose that $\dege_n f = -\infty$.  Then one has $f_{(n+1)}(x) = o\left(\frac{1}{\log x}\right) \ (x \to \infty)$ and therefore
$$\dege_{n+1} f \leq (0,-1,0,0,0,\ldots)$$
and
\begin{align*}
\dege f \leq (\dege_0 f, \ldots, \dege_{n} f, 0,-1, 0,0,0,\ldots).
\end{align*}
\item  If $\dege_n f$ is finite and $\dege_{n+1} f= \infty$, then one has $f_{(n+2)}(x) = o(x) \ (x \to \infty)$, so that $\dege f_{n+2} \leq (1,0,0,0,\ldots)$ and therefore
\begin{align*}
(\dege_0 f, \ldots, \dege_n f, \infty, 0,1, 0,0,0\ldots)  &  \leq \dege f  \\ & \leq (\dege_0 f, \ldots, \dege_n f, \infty, 1,0,0,\ldots).
\end{align*}
\item  If $\dege_n f$ is finite and $\dege_{n+1} f= -\infty$, then one has $f_{(n+2)}(x) \neq O\left(\frac{1}{x} \right) \ (x \to \infty)$, so that $\dege_{n+2} f \geq (-1,0,0,0,\ldots)$ and therefore
\begin{align*}
(\dege_0 f, \ldots, \dege_n f, -\infty, 0,-1, 0,0,0\ldots) &  \geq \dege f  \\ & \geq (\dege_0 f, \ldots, \dege_n f,- \infty, -1,0,0,\ldots).
\end{align*}
\end{enumerate}
\end{proposition}

\begin{proof}
Suppose that $\dege_n f  =  \infty$.  Then  $$\limsup_{x \to \infty} \frac{\log |f_{(n)}(x)|}{\log x} = \infty$$ and therefore $f_{(n+1)}(x) = \max(\log |f_{(n)}(x)|,0) \neq O(\log x) \ (x \to \infty)$, so that $\dege f_{(n+1)} \geq (0,1,0,0,0,\ldots)$.  Statement (1) follows.

Suppose that $\dege_n f  = -\infty$.   Then
$\lim_{x \to \infty} \frac{-1}{\log |f_{(n)}(x)|} = 0$ and  $$\lim_{x \to \infty} \frac{\log |f_{(n)}(x)|}{\log x} = -\infty$$ and therefore  
 $$\lim_{x \to \infty} \frac{\frac{-1}{\log |f_{(n)}(x)|}}{\frac{1}{\log x}} = 0,$$
so that 
$$f_{(n+1)}(x) =  \frac{-1}{\log |f_{(n)}(x)|} = o \left(\frac{1}{\log x} \right) \ (x \to \infty)$$
and therefore $\dege f_{(n+1)} \leq (0,-1,0,0,0,\ldots)$.  Statement (2) follows.

Suppose now that  $d = \dege_n f = \deg f_{(n)}$ is finite.  For any $t >d$, one has $|f_{(n)}(x)| \leq x^t$ for all $x\gg 0$, and one has $f_{(n+1)}(x) = f_{(n)}(e^x)e^{-dx}$.  If $\dege_{n+1} f = \infty$, then one has
\begin{align*}
0 \leq f_{(n+2)}(x) & =  \log \max (|f_{(n+1)}(x)|,1) \\
& =  \log \max( |f_{(n)}(e^x)e^{-dx}|,1) \\
& \leq \log \max (e^{(t-d)x},1) \\
& = (t-d) x
\end{align*}
for all $t >d$ and all $x \gg 0$, and therefore $0 \leq f_{(n+2)}(x) = o(x) \ (x \to \infty)$.  Statement (3) follows.  Suppose, on the other hand, that $\dege_{n+1} f = -\infty$.  Then for any $t  < d$, one has $|f_{(n)}(x)|> x^t$ on an unbounded set of $x >0$.  It follows that $$0 \leq -\log |f_{(n+1)}(x)| = -\log | f_{(n)}(e^x)| + dx < (d-t)x,$$ so that
$$f_{(n+2)}(x) = -\frac{1}{\log |f_{(n+1)}(x)|} \geq \frac{1}{(d-t)x} > 0,$$
 on an unbounded set of $x> 0$, and therefore   $f_{(n+2)}(x) \neq O\left(\frac{1}{x} \right) \ (x \to \infty)$.  Statement (4) follows.
\end{proof}

 The reader can test these limits for themselves by trying to construct  functions with logexponential degree $(\infty, -\varepsilon, 0,0,0,\ldots)$,  $(\infty, 0, 1-\varepsilon, 0,0,0,\ldots)$, $(0, \infty,1+\varepsilon,0,0,0,\ldots)$, or $(0, \infty,1, \varepsilon,0,0,0,\ldots)$  for some $\varepsilon > 0$, for example.   The impossibility of such  constructions is what led us to Proposition \ref{degeprop}.  
Note, however,  that the various upper and lower bounds in the proposition are attained by examples (19)(e) and (19)(f) of Example \ref{degeex}.   

Proposition \ref{degeprop} motivates the following definition.  Let $$\prod_{n = 0}^{ \infty *}\overline{\RR} \subsetneq \prod_{n = 0}^{\infty}\overline{\RR} \index[symbols]{.g h@$\prod_{n = 0}^ { \infty *} \overline{\RR}$}\index{restricted product $\prod_{n = 0}^{ \infty *}\overline{\RR}$}$$ (to be contrasted with  $\prod_{n = 0}^{* \infty}\overline{\RR}$) denote the {\bf restricted product} consisting of all sequences $\dd$ in $\prod_{n = 0}^{\infty}\overline{\RR}$ satisfying the following four conditions for all nonnegative integers $n$.
\begin{enumerate}
\item If $\dd_n = \infty$, then $\dd \geq (\dd_0, \dd_1, \ldots, \dd_n, 0,1,0,0,0,\ldots).$
\item If $\dd_n = -\infty$, then $\dd \leq (\dd_0, \dd_1, \ldots, \dd_n, 0,-1,0,0,0,\ldots).$
\item If $\dd_n$ is finite and $\dd_{n+1} = \infty$, then $\dd \leq (\dd_0, \dd_1, \ldots, \dd_{n+1},1,0,0,0,\ldots).$
\item If $\dd_n$ is finite and $\dd_{n+1} = -\infty$, then $\dd \geq (\dd_0, \dd_1, \ldots, \dd_{n+1}, -1,0,0,0,\ldots).$
\end{enumerate}
By Proposition \ref{degeprop}, if $\dd = \dege f$ for some $f \in \RR^{\RR_\infty}$, then $\dege f \in \prod_{n = 0}^{ \infty *}\overline{\RR}$.  In fact the converse holds.   Indeed,  we prove the following theorem in Section 7.5.

\begin{theorem}\label{degeequiv}
For any sequence $\dd \in \prod_{n= 0}^{\infty}\overline{\RR}$, one has $\dd = \dege f$ for some $f \in \RR^{\RR_\infty}$ if and only if $\dd \in \prod_{n = 0}^{ \infty *}\overline{\RR}$, if and only if $\dd = \dege f$ for some positive,  monotonic,  infinitely differentiable function $f$ on $\RR_{>0}$ of exact logexponential degree.   Consequently, the image of $\dege$ is equal to $\prod_{n = 0}^{\infty *} \overline{\RR}$.
\end{theorem}

Now, let $$\mathbb{L} = \operatorname{Def}(\RR,\id,\exp, \log,+,\cdot,/, \circ)\index[symbols]{.i  k@$\mathbb{L}$}$$ denote the ring of all {\bf logarithmico-exponential functions}\index{logarithmico-exponential functions} \cite{har3} \cite{har4}, that is, the ring of all (germs of) real functions that are defined on a neighborhood of $\infty$ and can be  built from all real constants and the functions $\id$, $\exp$,  and $\log$  using the operations $+$, $\cdot$, $/$, and $\circ$.   The ring $\mathbb{L}$ contains the field  $\RR(\mathfrak{L})$ as a subfield, and, as a consequence of Theorem \ref{infpropexp} below, it plays the role for $\dege$ that $\RR(\mathfrak{L})$ plays for $\degl$ and $\RR(x^a: a \in \RR)$ plays for $\deg$. 
Note that, by \cite[Theorem, p.\ 24]{har3}, any function in the ring $\mathbb{L}$ is  positive, $0$, or negative, for all $x \gg 0$, and thus $\mathbb{L}$ is a field.

\begin{proposition}\label{Kfield}
Let $r \in \mathbb{L}$.  One has the following.
\begin{enumerate}
\item $r$ is positive, $0$, or negative, for all $x \gg 0$.   
\item $r$ is infinitely differentiable on its domain, and all of its derivatives lie in $\mathbb{L}$.
\item Each of the functions $r_{(k)}$ as in the definition of $\dege r$ has exact degree and is equal to some function in $\mathbb{L}$ on some neighborhood of $\infty$. 
\item $r$ has exact logexponential degree.  
\end{enumerate}
\end{proposition}   

For any real functions $f$ and $g$ defined on a neighborhood of $\infty$, we write $f \leq_\infty g$ if $g-f$ is eventually positive or eventually $0$,  we write $f <_\infty g$ if $g-f$ is eventually positive, and we write $f =_\infty g$ if $g-f$ is eventually $0$.  Thus, for all $f, g\in \mathbb{L}$, one has either $f\leq_\infty g$ or $g \leq_\infty f$, and $f =_\infty g$ holds if and only if both $f\leq_\infty g$ and $g \leq_\infty f$.  The ordering $\leq_\infty $ makes $\mathbb{L}$ into a (totally) ordered field, provided that one identifies any two functions $f$ and $g$ that are eventually equal, i.e., that satisfy $f =_\infty g$.  This totally ordered field is a subfield of the ring of germs of real functions at $\infty$.   Let us also define $$\mathbb{L}_{> 0} = \{f \in \mathbb{L}: 0<_\infty f\}.$$
 Note that, if $r \in \mathbb{L}^*$, then $|r| \in \mathbb{L}_{> 0} $ and $r_{(k)} \in \mathbb{L}^*$ for all $k$,  and, if also $\dege_n r = \pm \infty$ for some $n$, then $r_{(k)} \in \mathbb{L}_{> 0}$ for all $k > n$, where, again, one identifies functions  that are eventually equal.

For any $f,g \in \RR^{\RR_\infty}$, we let $\mathcal{E}(f,g) = \mathcal{S}(\dege f, \dege g)$, that is, we let $\mathcal{E}(f,g)$ denote the infimum of the set of all nonnegative integers $n$ such that $\dege_n f \neq \dege_n g$ (which is $\infty$ if there is no such nonnegative integer $n$).\index[symbols]{.g p@$\mathcal{E}(f,g)$}    We say that $g$  {\bf approximates $f$ to logexponential order $n$}\index{approximates $f$ to logexponential order $n$} if  $\mathcal{E}(f,g) > n$, that is, if $\dege_k f = \dege_k g$ for all $k \leq n$.  The following proposition implies that any function in $\RR^{\RR_\infty}$ can be approximated to any logexponential order by some function in $\mathbb{L}_{> 0}$.

\begin{proposition}\label{degerprop}
Let $\dd \in  \prod_{n = 0}^{ \infty}\overline{\RR}$.  Then $\dd \in  \prod_{n = 0}^{ \infty*}\overline{\RR}$ if and only if, for every nonnegative integer $n$, there exists an $r \in \mathbb{L}_{> 0}$ such that $\dege_k r = \dd_k$ for all $k \leq n$.   In other words,  $ \prod_{n = 0}^{ \infty*}\overline{\RR}$ is equal to the closure of $\dege \mathbb{L}_{>0}$ in  $ \prod_{n = 0}^{ \infty}\overline{\RR}$, where $\overline{\RR}$ is endowed with the discrete topology and the product is endowed with the product topology.   Consequently, for any $f \in \RR^{\RR_\infty}$ and any nonnegative integer $n$, there exists an $r \in \mathbb{L}_{> 0}$ such that $\mathcal{E}(f,r) > n$.
\end{proposition}

\begin{proof}
By Proposition \ref{degeprop} and the definition of $\prod_{n = 0}^{ \infty *}\overline{\RR}$, if for every nonnegative integer $n$ there exists an $r \in \mathbb{L}_{> 0}$ such that $\dege_k r = \dd_k$ for all $k \leq n$, then $\dd \in  \prod_{n = 0}^{ \infty*}\overline{\RR}$.   Conversely, let $\dd \in \prod_{n= 0}^{\infty*}\overline{\RR}$, and let $n$ be a nonnegative integer.  Then the shifted sequence $(\dd_1, \dd_2, \dd_3, \ldots)$ also lies in $\prod_{n= 0}^{\infty*}\overline{\RR}$.  Therefore, by induction we may assume that we have an $s \in  \mathbb{L}_{> 0}$ such that $\dege_k s = \dd_{k+1}$ for all $k$ with $1 \leq k \leq n$, so that, for example, $\deg s = \dd_1$.   We wish to find an $r \in \mathbb{L}_{> 0}$ such that $\dege_k r =  \dd_k$ for all $k \leq n$.  Equivalently, we wish to find an $r \in \mathbb{L}_{> 0}$ such that $\deg r= \dd_0$ and $\dege_k r = \dege_{k-1} s$ for all $k$ with $1 \leq k \leq n$.

Suppose first that $\dd_0 \neq \pm \infty$.    If $\dd_1 \neq \pm \infty$, then $r(x) = x^{\dd_0}s(\log x) \in \mathbb{L}_{> 0}$ satisfies  our requirement for $r$.  Indeed, in this case, one has 
$$\deg s(\log x) = \lim_{x \to \infty} \frac{\log s(\log x)}{\log x} = \lim_{x \to \infty} \frac{\log s(x)}{x} = \lim_{x \to \infty} \frac{\log s(x)}{\log x} \frac{\log x}{x} = 0$$
and therefore $r_{(1)}(x) =(s(\log x))_{(1)} = s(x)$, so that $$\dege s(\log x) = (0,\dege_0 s , \dege_1 s, \dege_2 s, \ldots)$$
and thus
$$\dege r = (\dd_0, \dege_0 s , \dege_1 s, \dege_2 s, \ldots),$$
whence $r$ fulfills our requirement  that $\dege_k r = \dd_k$ for all $k \leq n$ (and in fact $\dege_k r = \dege_{k-1} s$ for all $k \geq 1$).

Suppose, instead, that $\dd_1  = \infty$, so that
$$ \dd \leq (\dd_0, \infty,1,0,0,0,\ldots).$$
Then $s(x) = e^{s_{(1)}(x)}$ for all $x \gg 0$, where $\lim_{x \to \infty} \frac{s_{(1)}(x)}{\log x} = \infty$  and
$$\dege s_{(1)} \leq (1,0,0,0,\ldots).$$   If $ \dege s_{(1)} < (1,0,0,0,\ldots)$, then $s_{(1)}(x) = o(x) \ (x \to \infty)$, so that 
$$\deg s(\log x) = \lim_{x \to \infty} \frac{s_{(1)}(\log x)}{\log x}   =   \lim_{x \to \infty}  \frac{s_{(1)}(x)}{x} = 0$$
and therefore $(s(\log x))_{(1)} = s(x)$ and thus
$$\dege s(\log x) = (0,\dege_0 s, \dege_1 s, \dege_2 s, \ldots).$$ Therefore, in this case, the function
$r(x) = x^{\mathbf{d_0}} s(\log x)$ satisfies  $r_{(1)}(x) = s(x)$ and fulfills our requirement for $r$.  Suppose, on the other hand, that $\dege s_{(1)} = (1,0,0,0,\ldots)$,
so that
$$\dege s = (\infty, 1,0,0,0,\ldots).$$
Then $r(x) = x^{\mathbf{d_0}} x^{1/\log^{\circ(n+1) }x} = x^{\mathbf{d_0}} e^{\log x/\log^{\circ(n+1) }x} \in \mathbb{L}_{> 0}$ has
$$\dege r = (\mathbf{d_0},\infty,1,0,0,0,\ldots, 0,-1,0,0,0,\ldots),$$
where the $-1$ is preceded by $n-1$ zeros, and so $\dege_k  r= \dd_k$ for all $k \leq n$.  

Suppose, on the other hand, that $\dd_1  = -\infty$, so that
$$\dd \geq (\dd_0, -\infty,-1,0,0,0,\ldots).$$
Then $s(x) = e^{-1/s_{(1)}(x)}$ for all $x \gg 0$, where $\lim_{x \to \infty} \frac{1/s_{(1)}(x)}{\log x} = \infty$  and
$$\dege s_{(1)} \geq (-1,0,0,0,\ldots).$$   If $ \dege s_{(1)} > (-1,0,0,0,\ldots)$, then $1/s_{(1)}(x) = o(x) \ (x \to \infty)$, so that 
$$\deg s(\log x) = \lim_{x \to \infty} \frac{-1/s_{(1)}(\log x)}{\log x}   =   \lim_{x \to \infty}  \frac{-1/s_{(1)}(x)}{x} = 0$$
and therefore $(s(\log x))_{(1)} = s(x)$ and thus
$$\dege s(\log x) = (0,\dege_0 s, \dege_1 s, \dege_2 s, \ldots).$$ Therefore, in this case, the function
$r(x) = x^{\mathbf{d_0}} s(\log x)$ satisfies  $r_{(1)}(x) = s(x)$ and fulfills our requirement.  Suppose, on the other hand, that $\dege s_{(1)} = (-1,0,0,0,\ldots)$,
so that
$$\dege s = (-\infty, -1,0,0,0,\ldots).$$
Then $r(x) = x^{\mathbf{d_0}} x^{-1/\log^{\circ(n+1) }x} = x^{\mathbf{d_0}} e^{-\log x/\log^{\circ(n+1) }x} \in \mathbb{L}_{> 0}$ has
$$\dege r = (\mathbf{d_0},-\infty,-1,0,0,0,\ldots, 0,1,0,0,0,\ldots),$$
where the $1$ is preceded by $n-1$ zeros, and so $\dege_k  r= \dd_k$ for all $k \leq n$.  

The proves the proposition in the case where $\dd_0 \neq \pm \infty$.  Suppose now that $\dd_0 =  \infty$, so that
$$\dd \geq (\infty,0,1,0,0,0,\ldots)$$
and therefore
$$\dege s \geq (0,1,0,0,0,\ldots).$$
Suppose that $ \dege s > (0,1,0,0,0,\ldots)$, so that  $s(x) \neq O(\log x) \ (x \to \infty)$ and therefore $\lim_{x \to \infty} \frac{s(x)}{\log x} = \infty$, since $\frac{s(x)}{\log x} \in \mathbb{L}_{> 0}$ is positive and unbounded.  Let $r(x) = e^{s(x)}$.  Then
$$\dege r = \lim_{x \to \infty} \frac{s(x)}{\log x} = \infty$$ and therefore $r_{(1)}(x) = \log e^{s(x)} =  s(x)$, so that
$$\dege r = (\infty, \dege_0 s, \dege_1 s, \dege_2 s, \ldots),$$
whence $\dege_k r = \dd_k$ for all $k\leq n$.  Suppose, on the other hand, that $\dege s = (0,1,0,0,0,\ldots)$. 
Then $r(x) = x^{\log^{\circ(n+1) }x} = e^{(\log x) (\log^{\circ(n+1) }x)} \in \mathbb{L}_{> 0}$ has
$$\dege r = (\infty,0,1,0,0,0,\ldots, 0,1,0,0,0,\ldots),$$
where the second $1$ is preceded by $n-1$ zeros, and so $\dege_k  r= \dd_k$ for all $k \leq n$.  

Finally,  suppose that $\dd_0 =  -\infty$, so that
$$\dd \leq (-\infty,0,-1,0,0,0,\ldots)$$
and therefore
$$\dege s \leq (0,-1,0,0,0,\ldots).$$
Suppose that $ \dege s < (0,-1,0,0,0,\ldots)$.  Then $1/s(x) \neq O(\log x) \ (x \to \infty)$ and therefore $\lim_{x \to \infty} \frac{1/s(x)}{\log x} = \infty$.  Let $r(x) = e^{-1/s(x)}$.  Then
$$\dege r = \lim_{x \to \infty} \frac{-1/s(x)}{\log x} =- \infty$$ and therefore $r_{(1)}(x) =- \frac{1}{\log e^{-1/s(x)}} =  s(x)$ for all $x \gg 0$, so that
$$\dege r = (-\infty, \dege_0 s, \dege_1 s, \dege_2 s, \ldots),$$
whence $\dege_k r = \dd_k$ for all $k\leq n$.  On the other hand, suppose that $\dege s = (0,-1,0,0,0,\ldots)$. 
Then $r(x) = x^{-\log^{\circ(n+1) }x} = e^{-(\log x) (\log^{\circ(n+1) }x)} \in \mathbb{L}_{> 0}$ has
$$\dege r = (-\infty,0,-1,0,0,0,\ldots, 0,-1,0,0,0,\ldots),$$
where the second $-1$ is preceded by $n-1$ zeros, and so $\dege_k  r= \dd_k$ for all $k \leq n$.    This completes the proof.
\end{proof}

In Section 7.6, we use transseries to prove the analogue of Theorem \ref{degeequiv} for the logarithmico-exponential functions, namely, Theorem \ref{chartheorem2}, by computing the  set of all $\dd \in \prod_{n = 0}^{ \infty*}\overline{\RR}$ such that $\dd = \dege f$ for some logarithmico-exponential function $f$.
 
Our next theorem is an analogue of Proposition \ref{infprop} for logexponential degree.   The theorem has many applications,  both to proving various properties of logexponential degree (in Subsection  6.3.4) and to applying  logexponential degree to analytic number theory  (in Part 3).

\begin{lemma}\label{compllem}
The totally ordered set $\prod_{n = 0}^{ \infty*}\overline{\RR}$ is complete.
\end{lemma}

\begin{proof}
Let $\mathcal{S}$ be a subset of $\prod_{n = 0}^{ \infty*}\overline{\RR}$, and let  $\dd$ be the infimum of $ \mathcal{S}$ in $\prod_{n = 0}^{ \infty}\overline{\RR}$.  If $\dd$ lies in $\prod_{n = 0}^{ \infty*}\overline{\RR}$, then clearly $\dd$ is the infimum of $ \mathcal{S}$ in $\prod_{n = 0}^{ \infty*}\overline{\RR}$.  Thus, we may suppose that $\dd$ does not lie in $\prod_{n = 0}^{ \infty*}\overline{\RR}$.   Let $n$  be the smallest nonnegative integer violating conditions (1)--(4) of the definition of   $\prod_{n = 0}^{ \infty*}\overline{\RR}$ preceding Theorem \ref{degeequiv}.  There are four cases to consider.

Suppose first that $\dd_n = \infty$ and $ \dd < (\dd_0, \dd_1, \ldots, \dd_n, 0,1,0,0,0,\ldots)$.  Then there exists some $\dd' \in \mathcal{S}$ such that 
$$\dd \leq (\dd_0, \dd_1 \ldots, \dd_{n-1},\infty, \dd_{n+1}',\dd_{n+1}',\ldots) = \dd' < (\dd_0, \dd_1, \ldots, \dd_{n-1},\infty, 0,1,0,0,0,\ldots).$$
But this contradicts $\dd' \in \prod_{n = 0}^{ \infty*}\overline{\RR}$.  Therefore this case is impossible.

Suppose now that $\dd_n = -\infty$ and $ \dd > (\dd_0, \dd_1, \ldots, \dd_n, 0,-1,0,0,0,\ldots)$. 
 Then one has $$\dd' > (\dd_0, \dd_1, \ldots, \dd_n, 0,-1,0,0,0,\ldots) \in \prod_{n = 0}^{ \infty*}\overline{\RR}$$
for all $\dd'  \in \mathcal{S}$.  Suppose that $\dd''$ is some lower bound of $\mathcal{S}$ in $\prod_{n = 0}^{ \infty*}\overline{\RR}$  such that $\dd'' > (\dd_0, \dd_1, \ldots, \dd_n, 0,-1,0,0,0,\ldots)$.  Then
$$(\dd_0, \dd_1 \ldots, \dd_{n-1},-\infty, \dd_{n+1},\dd_{n+1},\ldots) = \dd \geq \dd'' >(\dd_0, \dd_1, \ldots, \dd_{n-1},-\infty, 0,-1,0,0,0,\ldots),$$
which contradicts $\dd'' \in \prod_{n = 0}^{ \infty*}\overline{\RR}$.  Therefore $(\dd_0, \dd_1, \ldots, \dd_n, 0,-1,0,0,0,\ldots)$ is the greatest lower bound of $\mathcal{S}$ in $\prod_{n = 0}^{ \infty*}\overline{\RR}$.  Thus, $\mathcal{S}$ has an infimum in $\prod_{n = 0}^{ \infty*}\overline{\RR}$ in this case.

Suppose now that $\dd_n \neq \pm \infty$, $\dd_{n+1} = \infty$ and $ \dd > (\dd_0, \dd_1, \ldots, \dd_{n+1},1,0,0,0,\ldots)$.  Then one has $$\dd' > (\dd_0, \dd_1, \ldots, \dd_{n+1}, 1,0,0,0,\ldots) \in \prod_{n = 0}^{ \infty*}\overline{\RR}$$
for all $\dd'  \in \mathcal{S}$.  Suppose that $\dd''$ is some lower bound of $\mathcal{S}$  in $\prod_{n = 0}^{ \infty*}\overline{\RR}$ such that $\dd'' > (\dd_0, \dd_1, \ldots, \dd_{n+1}, 1,0,0,0,\ldots)$.  Then
$$(\dd_0, \dd_1 \ldots, \dd_{n},\infty, \dd_{n+2},\dd_{n+3},\ldots) = \dd \geq \dd'' >(\dd_0, \dd_1, \ldots, \dd_{n},\infty, 1,0,0,0,\ldots),$$
which contradicts $\dd'' \in\prod_{n = 0}^{ \infty*}\overline{\RR}$.  Therefore $(\dd_0, \dd_1, \ldots, \dd_{n+1}, \infty,1,0,0,0,\ldots)$ is the greatest lower bound of $\mathcal{S}$ in $\prod_{n = 0}^{ \infty*}\overline{\RR}$.  Thus, $\mathcal{S}$ has an infimum in $\prod_{n = 0}^{ \infty*}\overline{\RR}$ in this case.

Finally, suppose that $\dd_n \neq \pm \infty$, $\dd_{n+1} =- \infty$ and $ \dd < (\dd_0, \dd_1, \ldots, \dd_{n+1},-1,0,0,0,\ldots)$.  Then there exists some $\dd' \in \mathcal{S}$ such that 
$$\dd \leq (\dd_0, \dd_1 \ldots, \dd_{n},-\infty, \dd_{n+2}',\dd_{n+3}',\ldots) = \dd' < (\dd_0, \dd_1, \ldots, \dd_{n},-\infty, -1,0,0,0,\ldots).$$
But this contradicts $\dd' \in \prod_{n = 0}^{ \infty*}\overline{\RR}$.   Therefore this case is impossible.
\end{proof}

\begin{theorem}\label{infpropexp}
Let $f \in \RR^{\RR_\infty}$.    One has
\begin{align*}
\dege f & = \inf\{\dege r:r \in \mathbb{L}_{> 0}, \, \dege f \leq \dege r \} \\ 
& = \inf\{\dege r:r \in \mathbb{L}_{> 0}, \, f(x) = O(r(x)) \ (x \to \infty) \} \\ 
 & = \inf\{\dege r:r \in \mathbb{L}_{> 0}, \, \forall x \gg 0\, |f(x)|\leq r(x) \} \\
& = \inf\{\dege r:r \in \mathbb{L}_{> 0}, \, f(x) = o(r(x)) \ (x \to \infty) \} \\
& = \inf\{\dege r:r \in \mathbb{L}_{> 0}, \, \dege f < \dege r \},
\end{align*}
where the infima (exist and) are computed in  $\prod_{n = 0}^ {\infty*}\overline{\RR}$.
\end{theorem}

\begin{proof}
The given infima exist by Lemma \ref{compllem}.  By Corollary \ref{oexp}, if $\dege f < \dege r$, where $r \in \mathbb{L}_{> 0}$, then one has $f(x) = o(r(x)) \ (x \to \infty)$.  Moreover, by Proposition \ref{aspropexoexp}, if $f(x) = O(r(x)) \ (x \to \infty)$, then $\dege f \leq \dege r$.  It follows that
\begin{align*}
\dege f & \leq \inf\{\dege r:r \in \mathbb{L}_{> 0}, \, \dege f \leq \dege r \} \\ 
& \leq \inf\{\dege r:r \in \mathbb{L}_{> 0}, \, f(x) = O(r(x)) \ (x \to \infty) \} \\ 
 & \leq \inf\{\dege r:r \in \mathbb{L}_{> 0}, \, \forall x \gg 0\ |f(x)|\leq r(x) \} \\
& \leq \inf\{\dege r:r \in \mathbb{L}_{> 0}, \, f(x) = o(r(x)) \ (x \to \infty) \} \\
& \leq \inf\{\dege r:r \in \mathbb{L}_{> 0}, \, \dege f < \dege r\},
\end{align*}

Suppose to obtain a contradiction that
$$\dege f < \inf\{\dege r:r \in \mathbb{L}_{> 0}, \, \dege f < \dege r\}.$$
Then one has $\dege f < \dd$ for some $\dd  \in \prod_{n = 0}^{ \infty*}\overline{\RR}$ such that $\dd  \leq \dege r$ for all $r \in \mathbb{L}_{> 0}$ such that $\dege f < \dege r$.  Let $N$ be least such that $\dege_{N} f < \dd_{N}$, so that $\dege_{k} f = \dd_k$ for all $k < N$.   Let $t \in \RR$ be such that $\dege_{N} f < t< \dd_{N}$.  We claim that there is an $r \in \mathbb{L}_{> 0}$ such that $\dege_k r = \dege_k f$ for all $k < N$ and $\dege_N r= t$.   Assume that this claim is true.   Then one has $\dege f < \dege r < \dd$, which our desired contradiction.

Thus, we must prove the claim above.  Note first that, by the definition of $\prod_{n = 0}^ {\infty*}\overline{\RR}$, one has
 $$( \dd_0,  \dd_1,  \dd_2, \ldots,  \dd_{N-1}, t, 0,0,0,0,\ldots) \in \prod_{n = 0}^ {\infty*}\overline{\RR}$$
for some $t$ with $\dege_{N} f < t< \dd_{N}$ provided that it is not the case that $ \dd_{N-1} = -\infty$ and $\dege_N f = 0$ and it is not the case that
 $ \dd_{N-2} \neq \pm\infty$,  $ \dd_{N-1} = \infty$, and $\dege_N f = 1$.  Thus, in this case, our claim holds by Proposition \ref{degerprop}.

We may suppose, then, that either $ \dd_{N-1} = -\infty$ and $\dege_N f = 0$, or
 $ \dd_{N-2} \neq \pm\infty$,  $ \dd_{N-1} = \infty$, and $\dege_N f = 1$.   Suppose the first case.
Then
\begin{align*}
\dege f & < \dd = (\dd_0, \dd_1, \ldots, \dd_{N-2},-\infty, \dd_{N}, \dd_{N+1},\dd_{N+2},\ldots) \\ & \leq (\dd_0, \dd_1, \ldots, \dd_{N-2}, -\infty,0,-1,0,0,0,\ldots),
\end{align*}
which implies $0 = \dege_N f  < \dd_{N} \leq 0$, a contradiction. Likewise, the second case is impossible, since
then
\begin{align*}
\dege f & < \dd = (\dd_0, \dd_1, \ldots, \dd_{N-2},\infty, \dd_{N}, \dd_{N+1},\dd_{N+2},\ldots) \\ & \leq (\dd_0, \dd_1, \ldots, \dd_{N-2}, \infty,1,0,0,0,\ldots),\end{align*}
which  implies $1 = \dege_N f   < \dd_{N} \leq 1$, again a contradiction.  This completes the proof.
\end{proof}

\begin{corollary}\label{infpropexplow}
Let $f \in \RR^{\RR_\infty}$.    One has
\begin{align*}
\underline{\dege}\, f & = \sup\{\dege r:r \in \mathbb{L}_{> 0}, \, \underline{\dege}\, f\geq \dege r \} \\ 
& = \sup\{\dege r:r \in \mathbb{L}_{> 0}, \, f(x) \gg r(x) \ (x \to \infty) \} \\ 
 & = \sup\{\dege r:r \in \mathbb{L}_{> 0}, \, \forall x \gg 0\ |f(x)|\geq r(x) \} \\
& = \sup\{\dege r:r \in \mathbb{L}_{> 0}, \, f(x) \ggg r(x) \, (x \to \infty) \} \\
& = \sup\{\dege r:r \in \mathbb{L}_{> 0}, \, \underline{\dege}\, f > \dege r \},
\end{align*}
where the suprema (exist and) are computed in  $\prod_{n = 0}^ {\infty*}\overline{\RR}$.
\end{corollary}

\begin{corollary}
Let $f,g  \in \RR^{\RR_\infty}$.     One has $\dege f \leq \dege g$ if and only if $\dege f \leq \dege r$ for all  $r \in  \mathbb{L}_{> 0}$ such that $g(x) \ll r(x) \ (x \to \infty)$.  Likewise,  one has $\underline{\dege}\,  f \leq \underline{\dege}\,  g$ if and only if $\dege r \leq \underline{\dege}\, g$ for all  $r \in  \mathbb{L}_{> 0}$ such that $f(x) \gg r(x) \ (x \to \infty)$.
\end{corollary}

The following proposition is the analogue of Proposition \ref{fresupperlower} for logexponential degree.

\begin{proposition}
Let $f \in \RR^{\RR_\infty}$.  For all $X \subseteq \dom f$ with   $\sup X = \infty$, one has
$$\underline{\dege}\, f  \leq \underline{\dege}\, f|_X \leq \dege  f|_X \leq \dege f.$$
Moreover, one has
$$\underline{\dege}\, f = \inf\left\{\dege f|_X:  X \subseteq \dom f \text{ and }  \sup X = \infty \right\}$$
and
$$\dege f  = \sup\left\{\underline{\dege}\, f|_X:  X \subseteq \dom f \text{ and }  \sup X = \infty \right\}.$$
\end{proposition}

\begin{proof}
The first statement is clear, as then are the inequalities 
$$\underline{\dege}\, f \leq \inf\left\{\dege f|_X:  X \subseteq \dom f \text{ and }  \sup X = \infty \right\}$$
and
$$\dege f  \geq \sup\left\{\underline{\dege}\, f|_X:  X \subseteq \dom f \text{ and }  \sup X = \infty \right\}.$$
We prove that the second inequality above is an equality; the proof for the other is similar.   Suppose, to obtain a contradiction, that 
$$\dege f > \sup\left\{\underline{\dege}\, f|_X:  X \subseteq \dom f \text{ and }  \sup X = \infty \right\}.$$
By Proposition \ref{degerprop} (or by Proposition  \ref{oexppropstrong} of Section 7.5), there exists an $r \in \mathbb{L}_{>0}$ such that
$$\dege f > \dege r > \sup\left\{\underline{\dege}\, f|_X:  X \subseteq \dom f \text{ and }  \sup X = \infty \right\}.$$
It follows that, for every set $X \subseteq \dom f$ with  $\sup X = \infty$,  one has
$\dege r > \underline{\dege}\, f|_X$, whence it is not the case that $f|_X(x) \gg r(x) \ (x \to \infty)$,  i.e., one has $$\liminf_{x \to \infty} \frac{|f|_X(x)|}{r(x)} = 0.$$  Since this holds for every set $X \subseteq \dom f$ with  $\sup X = \infty$, one has $\limsup_{x \to \infty} \frac{|f(x)|}{r(x)} = 0$  and  therefore $\lim_{x \to \infty} \frac{|f(x)|}{r(x)} = 0$, i.e., $f(x) = o(r(x)) \ (x \to \infty)$.  However, this contradicts $\dege f > \dege r$.
\end{proof}

For any functions $f$ and $g$ on $[0,\infty)$ with $f(x) \geq g(x)$ for all $x$,  there is a function $h$ on $[0,\infty)$ with  $f(x) \geq h(x) \geq  g(x)$ for all $x$ and for which $h(x) = f(x)$ on some unbounded subset of $[0,\infty)$ and $h(x) = g(x)$ on some unbounded subset of $[0,\infty)$. From the proposition, then, it follows that $\dege h = \dege f$ and $\underline{\dege} \, h = \underline{\dege}\, g$.  Moreover, if $f$ and $g$ are continuous and  increasing (resp.,  decreasing), then one can find a corresponding $h$ with the same  properties.

\subsection{Relationships with operations on functions}

 The next proposition relates $\dege$ to the operations of addition, multiplication, and division of functions.  For all $\dd, \ee \in \prod_{n = 0}^\infty\overline{\RR}$, we define
$$\dd\oplus \ee  =  \dd+ \ee\index[symbols]{.g m@$\dd \oplus \ee$}  $$
if  $\dd_k$ and $\ee_k$ are finite for all $k$, and
$$\dd\oplus \ee =  (\dd_0+\ee_0, \ldots, \dd_{n-1}+\ee_{n-1},f_{0},f_{1},f_{2},\ldots)  
$$
if $n$ is the smallest nonnegative integer such that $\dd_n$ and $\ee_n$ are not both finite, where 
$$
(f_0,f_1,f_2,\ldots) = \begin{cases} \max( (\dd_{n},\dd_{n+1},\ldots),(\ee_{n},\ee_{n+1},\ldots)) & \text{if }  \dd_{n} = \infty \text{ or }  \ee_{n} = \infty \\
 \min( (\dd_{n},\dd_{n+1},\ldots),(\ee_{n},\ee_{n+1},\ldots)) & \text{otherwise}
\end{cases}$$
as computed in  $\prod_{n = 0}^\infty\overline{\RR}$.  The operation $\oplus$ is a binary operation on $\prod_{n = 0}^\infty\overline{\RR}$ and restricts to a binary operation on $\prod_{n = 0}^{\infty*}\overline{\RR}$.

\begin{lemma}\label{difflemexp}
Let $f,g \in \RR^{\RR_\infty}$ with $\dom f = \dom g$.  Then one has
$$\max(\dege f,\dege g) = \dege \max(|f(x)|,|g(x)|).$$
\end{lemma}

\begin{proof}
Since $0 \leq |f(x)| \leq  \max(|f(x)|,|g(x)|)$ and $0 \leq  |g(x)| \leq  \max(|f(x)|,|g(x)|)$ for all $x$ (in $\dom f = \dom g$), the inequality   $\max(\dege f,\dege g) \leq \dege \max(|f(x)|,|g(x)|)$ follows from Proposition \ref{aspropexoexp}(1).  To prove the reverse inequality, we may suppose without loss of generality that $\dege f \geq \dege g$.  Let $r \in  \mathbb{L}_{> 0}$ with $\dege r > \dege f$.  By Corollary \ref{oexp}, since $\dege f < \dege r$ and $\dege g < \dege r$, one has $|f(x)| \leq r(x)$ and $|g(x)| \leq r(x)$,  and therefore $\max(|f(x)|,|g(x)|) \leq r(x)$, for all $x \gg 0$, whence
$\dege\max(|f(x)|,|g(x)|) \leq \dege r$.  Taking the infimum over all $r$ as chosen and applying Theorem \ref{infpropexp},  we conclude that $\dege \max(|f(x)|,|g(x)|) \leq \dege f =  \max(\dege f,\dege g)$.
\end{proof}

\begin{theorem}\label{diffpropexp}
Let $f$ and $g$ be real functions such that both $\underline{f} = f|_{\dom f \cap \dom g}$ and $\underline{g} = g|_{\dom f \cap \dom g}$ are in $\RR^{\RR_\infty}$, and let $n$ be a nonnegative integer.  One has the following.
\begin{enumerate}
\item $\dege(f + g) \leq \max (\dege \underline{f},\dege \underline{g})$, with equality if $\dege \underline{f} \neq  \dege \underline{g}$.  
\item If $\dege(f+g) \neq \max(\dege \underline{f}, \dege \underline{g})$, then $\dege \underline{f} = \dege \underline{g}$.   
\item $\dege(fg) \leq \dege \underline{f} \oplus \dege \underline{g}$.
\item If both $\underline{f}$ and $\underline{g}$ have finite logexponential degree to order $n$, then $\dege_k (fg)  \leq \dege_k \underline{f}+ \dege_k \underline{g}$ for the smallest $k \leq n$, if any, for which equality does not hold.
\item Suppose that (a) $\underline{g}$ has finite logexponential degree to order $n-1$ and exact logexponential degree to order $n$; (b) $\underline{f}$ has finite logexponential degree to order $n-1$; (c) exactly one of $\dege_n \underline{g}$ and  $\dege_n \underline{f}$ is finite; and (d) if $\dege_n \underline{g} = \pm \infty$, then $\underline{g}$ has exact logexponential degree.   Then one has $\dege(fg) = \dege \underline{f} \oplus \dege \underline{g}$.
 \item If $\underline{g}$ has exact logexponential degree to order $n$,  and if one of $\dege_k \underline{f}$ and  $\dege_k \underline{g}$ is finite for  the smallest $k \leq n$ such that  $\dege_k \underline{f}$ and  $\dege_k \underline{g}$ are not both finite,  should such a $k$ exist,  then  $\mathcal{S}(\dege(fg), \dege \underline{f} \oplus \dege \underline{g}) > n$.  
 \item If $\underline{g}$ has exact logexponential degree,  and if one of $\dege_k \underline{f}$ and  $\dege_k \underline{g}$ is finite for  the smallest $k$ such that  $\dege_k \underline{f}$ and  $\dege_k \underline{g}$ are not both finite,  should such a $k$ exist, then  $\dege(fg) = \dege \underline{f} \oplus \dege \underline{g}$.
\item If $\underline{g}$ has finite and exact logexponential degree, then 
 $\dege(fg) = \dege \underline{f} \oplus \dege \underline{g}$.
 \item  $\dege (f^k) = \dege f \oplus \dege f \oplus\cdots  \oplus \dege f$ for any positive integer $k$.  More generally,  for any $a > 0$,   if $f$ has finite logexponential degree to order $n$,   then   $\dege_k (|f|^a) = a \dege_k f$ for all $k  \leq n$,  and, if one also has $\dege_{n+1} f = \pm \infty$, then $\dege_k (|f|^a) = \dege_k f$ for  all $k > n$.
\end{enumerate}
\end{theorem}

\begin{proof}
We may suppose without loss of generality that $f = \underline{f}$ and $g = \underline{g}$.   Since $$|f(x)+g(x)| \leq |f(x)|+|g(x)| \leq 2 \max(|f(x)|,|g(x)|)$$ for all $x \in \dom f = \dom g$, by  Lemma \ref{difflemexp} one has 
$$\dege (f+g) \leq \dege  \max(|f(x)|,|g(x)|) = \max(\dege f,\dege g).$$  This proves statement (1), and statement (2) is an immediate consequence of (1).  

We now prove statement (3).  By Proposition \ref{diffprop}, we may assume that $\dege_{k} f$ and $\dege_{k} g$ are not both finite for all $k$.
Let $n$ be least so that $\dege_{n} f$ and $\dege_{n} g$ are not both finite.  We may also suppose that $\degl_k (fg) = \degl_k f + \degl_k g$ for all $k < n$, for otherwise $\dege (fg) < \dege f\oplus \dege g$ by Proposition \ref{diffprop}.  It follows, then, that $(fg)_{(k)} = f_{(k)} g_{(k)}$ for all $k \leq n$.   

First, suppose that $\deg f_{(n)} = \deg g_{(n)} = -\infty$, so that  $\deg \, (fg)_{(n)} =- \infty$.  Then $f_{(n+1)}(x) = -\frac{1}{\log |f_{(n)}(x)|}$ and $g_{(n+1)}(x) = -\frac{1}{\log|g_{(n)}(x)|}$ and 
\begin{align*}
0 \leq (fg)_{(n+1)}(x) & = -\frac{1}{\log|f_{(n)}(x)g_{(n)}(x)|} \\
  &  = -\frac{1}{\log|f_{(n)}(x)|+\log |g_{(n)}(x)|}  \\
    &  = \frac{1}{\frac{1}{f_{(n+1)}(x)}+\frac{1}{g_{(n+1)}(x)}}  \\
        &  \leq \min(f_{(n+1)}(x),g_{(n+1)}(x)).
\end{align*}
for all $x \gg 0$.  It follows, then, that 
\begin{align*}
\dege \,  (fg)_{(n+1)}  \leq \dege \min(f_{(n+1)}(x),g_{(n+1)}(x)) \leq \min (\dege f_{(n+1)}, \dege g_{(n+1)}),
\end{align*}
whence statement (3) holds in this case.

Next,  suppose that $\deg f_{(n)} = -\infty$ and $-\infty < \deg g_{(n)} < \infty$, so that $\deg \, (fg)_{(n)} = \deg (f_{(n)}g_{(n)}) = -\infty$.  One has
\begin{align*}
0 \leq (fg)_{(n+1)}(x)   = -\frac{1}{\log|f_{(n)}(x)|+\log |g_{(n)}(x)|} 
\end{align*}
and therefore
\begin{align*}
 \limsup_{x \to \infty}  \frac{(fg)_{(n+1)}(x) }{f_{(n+1)}(x)} & = \limsup_{x \to \infty} \frac{\log|f_{(n)}(x)|}{\log|f_{(n)}(x)|+\log |g_{(n)}(x)|}   \\ & =  \limsup_{x \to \infty} \frac{1}{1+\frac{\log |g_{(n)}(x)|}{\log|f_{(n)}(x)|}}  \\ 
& = \frac{1}{1+ \liminf_{x \to \infty}  \frac{\log |g_{(n)}(x)|}{\log|f_{(n)}(x)|}}  \\
& = \frac{1}{1- \limsup_{x \to \infty}  \frac{\log |g_{(n)}(x)|/\log x}{-\log|f_{(n)}(x)|/\log x}} \\
& = \frac{1}{1- (\deg g_{(n)}) \left(\frac{1}{\infty} \right) }\\
& = 1,
\end{align*}
whence
$$(fg)_{(n+1)}(x) = O(f_{(n+1)}(x)) \ (x \to \infty).$$
It follows that $$ \dege \, (fg)_{(n+1)} \leq \dege f_{(n+1)}.$$
Thus, statement (3) holds in this case.

Next, suppose that $\deg f_{(n)} = \deg g_{(n)} = \infty$.  If $\deg \, (fg)_{(n)} < \infty$, then $\dege \, (fg)_{(n)} < \max(\dege f_{(n)}, \dege g_{(n)})$.  If, on the other hand, one has $\deg \, (fg)_{(n)} = \infty$, then  $f_{(n+1)}(x) =  \max(\log |f_{(n)}(x)|,0)$ and $g_{(n+1)}(x) =  \max(\log |g_{(n)}(x)|,0)$, and thus $$0\leq (fg)_{(n+1)}(x) =  \max(\log |f_{(n)}(x)g_{(n)}(x)|,0) \leq f_{(n+1)}(x)+g_{(n+1)}(x),$$ so that $\dege \, (fg)_{(n+1)}  \leq \max (\dege f_{(n+1)}, \dege g_{(n+1)})$.  Thus, statement (3) holds in either case.

To complete the proof of (3), suppose that $\deg f_{(n)} = \infty$ and $\deg g_{(n)} < \infty$.    Let $F = f_{(n)}$ and $G = g_{(n)}$.  Note that
$$\dege(FG) = \dege((F+G)^2-F^2-G^2) \leq \max(\dege(F+G)^2, \dege(F^2),\dege( G^2)),$$
where 
$$\dege(F^2) \leq \dege F \oplus \dege F = \dege F,$$
$$\deg(G^2) = 2 \deg G < \infty,$$
and
$$\dege(F+G)^2 \leq \dege (F+G) \oplus \dege (F+G) = \dege (F+G) = \dege F,$$
whence
$$\dege(FG)  \leq \max(\dege(F+G)^2, \dege(F^2),\dege( G^2)) \leq \dege F.$$
Thus, in this case, we have $\dege \, (fg)_{(n)} \leq \dege  f_{(n)}$.  This completes the proof of (3).

Statement (4) follows from Propositions \ref{diffprop} and \ref{degeprop0}.  

Assume the hypotheses of statement (5).    By Proposition \ref{mmm}(1), one has $\dege_k(fg) = \dege_k f + \dege_k g$ and $(fg)_{(k)} = f_{(k)} g_{(k)}$, where  $g_{(k)}$ has exact degree $\dege_k g$, for all $k \leq n$.   Suppose that $d:=\dege_n g$ is finite, so that $\dege_n f = \pm\infty$.   Suppose first that $\dege_n f = -\infty$.  Then one has
$$\lim_{x \to \infty} \frac{\log |g_{(n)}(x)|}{\log |f_{(n)}(x)|} = \lim_{x \to \infty} \frac{d\log x}{\log |f_{(n)}(x)|} = \frac{d}{-\infty}= 0,$$
and therefore
$$\lim_{x \to \infty} \frac{f_{(n+1)}(x)}{(fg)_{(n+1)}(x)}  =  \lim_{x \to \infty} \frac{\log |f_{(n)}(x)g_{(n)}(x)|}{\log |f_{(n)}(x)|}= \lim_{x \to \infty} \left(1+\frac{\log|g_{(n)}(x)|}{\log |f_{(n)}(x)|}\right) = 1.$$  
It follows that $(fg)_{(n+1)}(x) \sim f_{(n+1)}(x) \ (x \to \infty)$, whence $\dege \, (fg)_{(n)} = \dege f_{(n)}$.  
Suppose, on the other hand, that $\dege_n f = \infty$.  Let $N$ be any real number with  $N > -\deg g_{(n)}$ and $N \geq 0$, and let
\begin{align*}
F(x) & = \frac{\max(\log |f_{(n)}(x)|,2N \log x)}{\max(\log |f_{(n)}(x)g_{(n)}(x)|,2N \log x+\log|g_{(n)}(x)|)} \\ 
& = \frac{\max\left(\frac{\log |f_{(n)}(x)|}{\log x},2N\right)}{\max\left(\frac{\log |f_{(n)}(x)|}{\log x}+\frac{\log|g_{(n)}(x)|}{\log x},2N+\frac{\log|g_{(n)}(x)|}{\log x}\right)}.
\end{align*}
Note that $F(x) \geq 0$ for all $x \gg 0$.  If $\deg g_{(n)} > 0$, then we may let $N = 0$, and then $F(x) \leq 1$ for all $x \gg 0$.  Suppose, on the other hand, that
$\deg g_{(n)} \leq 0$, so that $N > -\deg  g_{(n)} \geq 0$.   We  claim  that $F(x) \leq 2$ for all $ x \gg 0$.  One has
$$-N< \frac{\log|g_{(n)}(x)|}{\log x}$$ for all $x \gg 0$.  If $\frac{\log |f_{(n)}(x)|}{\log x} \geq 2N$, then one has
$$F(x)  = \frac{\frac{\log |f_{(n)}(x)|}{\log x}}{\frac{\log |f_{(n)}(x)|}{\log x}+\frac{\log |g_{(n)}(x)|}{\log x}} < \frac{\frac{\log |f_{(n)}(x)|}{\log x}}{\frac{\log |f_{(n)}(x)|}{\log x}-N} \leq 2$$
for all $x \gg 0$.
On the other hand, if  $\frac{\log |f_{(n)}(x)|}{\log x} \leq 2N$, then
$$F(x) = \frac{2N}{2N+ \frac{\log |g_{(n)}(x)|}{\log x}} < \frac{2N}{2N-N} = 2$$
for all $x  \gg 0$.  Thus, we have shown that $0\leq F(x) \leq 2$ for all $x \gg 0$.  It follows that
$$\max(\log |f_{(n)}(x)|,2N \log x) = O\left(\max(\log |f_{(n)}(x)g_{(n)}(x)|,2N \log x+\log|g_{(n)}(x)|) \right) \ (x \to \infty).$$
Therefore, since $$\dege \max(\log |f_{(n)}(x)|,0) \geq (0,1,0,0,0,\ldots) =\dege (2N \log x)$$ and likewise $$\dege \max(\log |f_{(n)}(x)g_{(n)}(x)|,0) \geq  \dege(2N \log x+\log|g_{(n)}(x)|),$$
by Lemma \ref{difflemexp} and Proposition \ref{aspropexoexp}(1) one has
\begin{align*}
\dege f_{(n+1)} & = 
 \dege  \max(\log |f_{(n)}(x)| ,0)  \\ & = \dege \max(\log |f_{(n)}(x)|,2N \log x) 
 \\ & \leq  \dege \max(\log |f_{(n)}(x)g_{(n)}(x)|,2N \log x+\log|g_{(n)}(x)|) \\
  & = \dege \max(\log |f_{(n)}(x)g_{(n)}(x)|,0)  \\
   & = \dege \, (fg)_{(n+1)}.
\end{align*}
Since the reverse inequality was proved in the proof of statement (3), equalities hold.   This proves (5) in the case where $\dege_n g$ is finite.

Suppose,  now, that $\dege_n g = \pm \infty$, so that $\dege_n f$ is finite and $g$ has exact logexponential degree.   It follows that $\dege_n(fg) = \deg (f_{(n)} g_{(n)}) = \dege_n g$.     Suppose that $\dege_n g = -\infty$.
The proof of (3) shows that
$$\limsup_{x \to \infty}  \frac{(fg)_{(n+1)}(x) }{g_{(n+1)}(x)}  = 1,$$
whence $(fg)_{(n+1)}(x) =O( g_{(n+1)}(x)) \ (x \to \infty)$ and $(fg)_{(n+1)}(x)\neq o( g_{(n+1)}(x)) \ (x \to \infty)$.  Therefore, since $g_{(n+1)}$ has exact logexponential degree, it follows from  Corollary \ref{oexpcor} that $\dege \, (fg)_{(n+1)} = \dege g_{(n+1)}$.
Suppose, on the other hand, that $\dege_n g = \infty$, so that
$\lim_{x \to \infty} |g_{(n)}(x)| = \infty$ and $g_{(n+1)}(x) = \log  |g_{(n)}(x)|$ for all $x \gg 0$, while also
$$\limsup_{x \to \infty} \frac{\log |f_{(n)}(x)|}{\log |g_{(n)}(x)|} = \limsup_{x \to \infty} \frac{\log |f_{(n)}(x)|/\log x}{\log |g_{(n)}(x)|/\log x} = \frac{\dege_n f }{\infty}  = 0.$$ It follows that
\begin{align*}
\limsup_{x \to \infty} \frac{(fg)_{(n+1)}(x)}{g_{(n+1)}(x)} &  =  \limsup_{x \to \infty} \frac{\max(\log |f_{(n)}(x)g_{(n)}(x)|,0)}{\log |g_{(n)}(x)|} \\
\\ & = \limsup_{x \to \infty} \max\left(1+\frac{\log |f_{(n)}(x)|}{\log |g_{(n)}(x)|},0\right) \\
& = 1,
\end{align*}
whence $\dege \, (fg)_{(n+1)} = \dege g_{(n+1)}$ also in that case.  Statement (5) follows.

Note that statement (7) is an immediate consequence of statement (5) and Proposition \ref{mmm}(1), and statement (8) is an immediate consequence of statement (7).   Statement (9) is an easy consquence of the fact that $\log  |f|^a = a \log  |f|$.  Finally, to prove (6), note that the proof of (5) implies that $\dege_l (fg) = \dege_l f + \dege_l g$ and $(fg)_{(l)} = f_{(l)}g_{(l)}$ for all $l \leq k$,  while also $(fg)_{(k+1)}(x) =O( g_{(k+1)}(x)) \ (x \to \infty)$ and $(fg)_{(k+1)}(x)\neq o( g_{(k+1)}(x)) \ (x \to \infty)$, where by hypothesis $g_{(k+1)}(x)$ has exact logexponential degree to order $n-(k+1)$,  so that, by Proposition \ref{aspropexoexp}(1) and Corollary \ref{oexp}(1),  one has $$\mathcal{L}((fg)_{(k+1)},  g_{(k+1)}) > n-(k+1),$$
whence
$$\mathcal{S}(\dege(fg), \dege \underline{f} \oplus \dege \underline{g}) > n.$$  This completes the proof.
\end{proof}

\begin{example} \
\begin{enumerate}
\item Let $f(x) = e^{x^2}$ and $g(x) = e^{-x}$, so that $\deg f = \infty$ and $\deg g = -\infty$.  Then $f$ and $g$ have exact logexponential degree,  and one has
$$\dege (fg)  = (\infty,2,0,0,0,\ldots) = \dege f \oplus \dege g.$$
\item Let $f(x) = e^x$ and $g(x) = e^{-x^2}$, so that $\deg f = \infty$ and $\deg g = -\infty$.   Then $f$ and $g$ have exact logexponential degree,  yet one has
$$\dege (fg)  = (-\infty,-2,0,0,0,\ldots) < \dege f \oplus \dege g = (\infty,1,0,0,0,\ldots).$$
\item Let $f(x)$ be defined on $\RR$ so that $f(x) = 1$ for $x \in [2N,2N+1)$ and $f(x) = 0$ for $x \in [2N-1,2N)$ for all integers $N$.   Let $g(x)  =  1-f(x)$.  Then $f(x) + g(x) = 1 = \max(f(x),g(x))$ has exact logexponential degree  $(0,0,0,0\ldots)$,  while also  $\underline{\dege} \, f = \underline{\dege} \, g = (-\infty,-\infty,\infty, \ldots)$.  This example shows that, in general,  one need not have $\underline{\dege}(f+g) \leq \max (\underline{\dege}\, f,\underline{\dege} \, g)$, and the obvious inequality  $\max(\underline{\dege}\, f , \underline{\dege}\, g) \leq \underline{\dege} \max (|f|,|g|)$ need not be an equality.    Nevertheless, since $|f+g| \leq 2\max ( |f|,|g|)$,  in general one has $\underline{\dege}(f+g) \leq \underline{\dege} \max (|f|,|g|)$.   Moreover,  one has $\min(\underline{\dege}\, f, \underline{\dege}\, g) = \underline{\dege} \min (|f|,|g|)$: the proof is similar to that of Lemma \ref{difflemexp}.    
\end{enumerate}
\end{example}

\begin{problem}
Generalize statements (3)--(8) of Theorem \ref{diffpropexp} to the extent possible.
\end{problem}

Next, we relate $\dege$ to the operation of composition of functions.  The most general case possible seems to be rather complicated, so we restrict ourselves to the special cases below, which are sufficient for our purposes.

\begin{theorem}\label{fgie}
Let $f$ and $g$ be real functions defined on a neighborhood of $\infty$.  Suppose that $f$  has finite degree $d$ and exact logexponential degree, $g$ has positive degree and is eventually positive, and $g(x) \asymp  r(x) \ (x \to \infty)$ for some  $r \in \mathfrak{L}$.  Then one has 
$$\dege( f  \circ g )=  \dege f \oplus d (\dege g + (-1,0,0,0,\ldots)),$$
and $f \circ g$ has exact logexponential degree.
\end{theorem}

\begin{proof}
Let $e_k = \dege_k g= \dege_k r$ for all $k$.
Note that 
 $$\log g(x) \sim \log r(x) \sim e_0 \log x \ (x \to \infty)$$
and therefore
 $$\log g(e^x) \sim e_0 x \ (x \to \infty).$$
  It follows by induction, then, that
$$\log^{\circ k} g(x) \sim \log^{\circ k } x \ (x \to \infty)$$ for all $k \geq 2$.  
Thus, for any nonnegative integer $n$, one has
$$h_n(x) := \log^{\circ (n+1)} g(\exp^{\circ (n+1)}) \sim a_n x \ (x \to \infty),$$
where $a_0 = e_0$ and $a_n = 1$ if $n \geq 1$.

Replacing $f$ with $|f|$, we may assume that $f$ is nonnegative and eventually positive.  Since $f$ has finite degree $d$, we may suppose that $f$ has finite logarithmic degree to order $n \geq 0$.  Then, by Proposition \ref{fgi}, one has
$\deg (f \circ g)=  de_0$ and
$$\degl_1 (f \circ g) = \degl_1 f + d \degl_1 g= \degl_1 f + de_1,$$
and, provided that $n \geq  1$, also
$$\degl_k (f \circ g) = \degl_k f + d \degl_k g = \degl_k f, \quad \forall k = 2,3,\ldots,n.$$
Moreover,  since $f$ and $g$ have exact logexponential degree, all of the limits superior in the proof of Proposition \ref{fgi} can be replaced with limits, and thus $f \circ g$ has exact logexponential degree to order $n+1$.  Thus, the conclusion of the theorem holds if $f$ has finite logarithmic degree.  We may suppose, then, without loss of generality, that $\dege_{n+1} f = \deg f_{(n+1)} = \pm \infty$.   

Let $d_k = \dege_k f$  for all $k$. 
One has  $$f_{(n+1)}(\log^{\circ (n+1)} x) = f(x)x^{-d}(\log x)^{-d_1} \cdots (\log^{\circ n}x)^{-d_ {n}},$$
so that 
$$f_{(n+1)}(\log^{\circ (n+1)} g(x)) = f(g(x))g(x)^{-d}(\log g(x))^{-d_1} \cdots (\log^{\circ n}g(x))^{-d_ {n}},$$
and therefore 
\begin{align*}
(f\circ g)_{(n+1)}(\log^{\circ (n+1)} x) & = f(g(x))(\log x)^{-d_1} \cdots (\log^{\circ n}x)^{-d_ {n}} \cdot x^{-de_0} (\log x)^{-de_1}\cdots (\log^{\circ n} x)^{-de_n} \\
 & =  f_{(n+1)}(\log^{\circ (n+1)} g(x))\left( \frac{g(x)}{r_n(x)}\right)^{d} \left( \frac{\log g(x)}{\log x}\right)^{d_1} \cdots\left( \frac{\log^{\circ n} g(x)}{\log^{\circ n}x}\right)^{d_{n}} \\
  & \sim  f_{(n+1)}(\log^{\circ (n+1)} g(x))\left( \frac{g(x)}{r_n(x)}\right)^{d}e_0^{d_1} \\
   & =  f_{(n+1)}(\log^{\circ (n+1)} g(x)) k_n(x),
\end{align*}
where  
$$r_n(x) = x^{e_0}(\log x)^{e_1} \cdots (\log^{\circ n} x)^{e_n}$$
and
$$k_n(x) = \left( \frac{g(x)}{r_n(x)}\right)^{d}e_0^{d_1}.$$
It follows that
$$(f\circ g)_{(n+1)}(x) \sim   f_{(n+1)}(h_n(x))\,  k_n(\exp^{\circ(n+1)}x) \ (x \to \infty),$$
where, as before, 
$$h_n(x) = \log^{\circ (n+1)} g(\exp^{\circ (n+1)}) \sim a_n x \ (x \to \infty).$$
Note that, since $g(x) \asymp  r(x) \ (x \to \infty)$, one has
$$k_n(x) \asymp \prod_{k > n}(\log^{\circ k} x)^{de_k} \ (x \to \infty),$$
and therefore 
$$\log k_n(x) = O( \log^{\circ (n+2)} x) \ (x \to \infty),$$ whence
 $$\log k_n(\exp^{\circ(n+1)}x) = O(\log x) \ (x \to \infty).$$
Since $\overline{\underline{\deg}} \, f_{(n+1)}  = \pm \infty$, one has
$$\lim_{x \to \infty} \frac{\log x}{  \log  f_{(n+1)}(h_n(x))  } = \lim_{x \to \infty} \frac{\log h_n(x)}{  \log  f_{(n+1)}(h_n(x))  }  =  0$$
and therefore
$$\log k_n(\exp^{\circ(n+1)}x)  = O(\log x) = o( \log  f_{(n+1)}(h_n(x)) ) \ (x \to \infty)$$
If $\deg f_{(n+1)} = \infty$, then also $\deg \, (f\circ g)_{(n+1)}  = \infty$ and therefore
\begin{align*}
(f\circ g)_{(n+2)}(x) & = \log \, (f\circ g)_{(n+1)}(x) \\
& \sim   \log  f_{(n+1)}(h_n(x))  + \log  k_n(\exp^{\circ(n+1)}x) \\
& \sim   \log  f_{(n+1)}(h_n(x))  \\
& = f_{(n+2)}(h_n(x))\ (x \to \infty),
\end{align*}
where $f_{(n+2)}$ has finite degree (by Proposition \ref{degeprop}(3)).  On the other hand, if $\deg f_{(n+1)} = -\infty$, then also $\deg \, (f\circ g)_{(n+1)}  = -\infty$ and therefore
\begin{align*}
(f\circ g)_{(n+2)}(x) & = -\frac{1}{\log \, (f\circ g)_{(n+1)}(x)} \\
& \sim  -\frac{1}{  \log  f_{(n+1)}(h_n(x))  + \log  k_n(\exp^{\circ(n+1)}x) }   \\
& \sim  -\frac{1}{  \log  f_{(n+1)}(h_n(x))  }   \\
& =  f_{(n+2)}(h_n(x))\ (x \to \infty),
\end{align*}
where $f_{(n+2)}$ has finite degree (by Proposition \ref{degeprop}(4)).
Thus, we have shown that 
$$(f\circ g)_{(n+2)}(x)  \sim (f_{(n+2)} \circ h_n)(x) \ (x \to \infty),$$
and therefore
$$\dege \, (f\circ g)_{(n+2)}  =  \dege (f_{(n+2)} \circ h_n)$$
for some function $h_n$ with $h_n(x) \sim a_n x \ (x \to \infty)$ for some $a_n > 0$, where $f_{(n+2)}$ has finite degree and exact logexponential degree.   The argument may then be repeated {\it ad infinitum} to show that
$$\dege (f_{(n+2)} \circ h_n) = \dege f_{(n+2)}.$$
Therefore, since
$$\dege \, (f\circ g)_{(n+2)}  =  \dege f_{(n+2)},$$
the identity for $\dege \, (f\circ g)$ stated in the theorem follows from the definition of the operation $\oplus$, and the claim that $f \circ g$ has exact logexponential degree follows by repeated application of Proposition \ref{exactas}.
\end{proof}
 
\begin{lemma}\label{fgie00lemma}
Let $\dd,\dd',\ee \in \prod_{n = 0}^{ \infty*}\overline{\RR}$ with $\ee_n \neq \pm \infty$ for all $n$, and let
$\mathcal{S}$  be a subset of $\prod_{n = 0}^{ \infty*}\overline{\RR}$.  One has the following.
\begin{enumerate}
\item If $\dd \leq \dd'$, then $\dd \oplus \ee \leq \dd' \oplus \ee$.
\item $(\dd\oplus \ee) \oplus (-\ee) = \dd$.
\item $(\inf \mathcal{S}) \oplus \ee = \inf\{\dd\oplus\ee:  \dd \in  \mathcal{S}\}.$
\end{enumerate}
\end{lemma}

\begin{proof}
Statements (1) and (2) are clear, and statement  (3) is an easy consequence of statements  (1) and (2).
\end{proof}

\begin{theorem}\label{fgie00}
Let $f$ and $g$ be real functions defined on a neighborhood of $\infty$.  Suppose that $f$  has finite degree $d$, that $g$ has positive degree and is eventually positive,  continuous, and increasing, and  that $g(x) \asymp  r(x) \ (x \to \infty)$ for some  $r \in \mathfrak{L}$.  Then one has 
$$\dege( f  \circ g )=  \dege f \oplus d (\dege g + (-1,0,0,0,\ldots)).$$
\end{theorem}

\begin{proof}
Let $e_k = \dege_k g= \dege_k r$ for all $k$.  Note that $g$ has exact logexponential degree, by Proposition \ref{exactas}.  Thus, by Proposition \ref{circle1}(2), one has $\deg(f \circ g) = \deg f \deg g$.  Note also that the compositional inverse $g^{-1}$ of $g$ is eventually continuous and increasing, and one has
$$g^{-1}(x) \asymp  x^{1/e_0} (\log x)^{-e_1/e_0}(\log^{ \circ 2} x)^{-e_2/e_0} \cdots (\log^{ \circ n} x)^{-e_n/e_0}  \ (x \to \infty).$$

Suppose first that $\dege f < (d,\infty,1,0,0,0,\ldots)$.
Then there exists an $s \in \mathbb{L}_{> 0}$ with $f(x) = O(s(x)) \ (x \to \infty)$ and $\deg s = d$.   For any such $s$, one has
$f(g(x)) = O(s(g(x))) \ (x \to \infty)$, and therefore
$$\dege (f \circ g) \leq \dege(s \circ g) =  \dege s \oplus d (\dege g + (-1,0,0,0,\ldots)),$$
by Theorem \ref{fgie}.  Taking the infimum over all  $s$ as chosen,  by Theorem \ref{infpropexp} and Lemma \ref{fgie00lemma} we have
$$\dege (f \circ g) \leq \dege f \oplus d ( \dege g + (-1,0,0,0,\ldots)) < (de_0,\infty,1,0,0,0,\ldots).$$
From the inequality above,  it follows that  $(f \circ g)(x) = O(t(x)) \ (x \to \infty)$ for some  $t \in \mathbb{L}_{> 0}$ with $\deg t = \deg(f \circ g) = de_0$. For any such $t$, one has  $f(x) = O(t(g^{-1}(x))) \ (x \to \infty)$, and therefore
\begin{align*}
\dege f & \leq \dege t(g^{-1}(x))  \\ & = \dege t\oplus de_0(1/e_0-1,-e_1/e_0, -e_2/e_0,-e_3/e_0,\ldots) \\
& = \dege t\oplus (d-de_0,-de_1, -de_2,-de_3,\ldots) \\
& = \dege t \oplus -d(\dege g+ (-1,0,0,0,\ldots)),
\end{align*}
again by Theorem \ref{fgie}.  Taking the infimum over all $t$ as chosen, we see that
\begin{align*}
\dege f \leq \dege(f \circ g)\oplus -d(\dege g + (-1,0,0,0,\ldots)).
\end{align*}
Thus, by Lemma \ref{fgie00lemma}, we have
\begin{align*}
\dege (f \circ g)  & \leq  \dege f \oplus d (\dege g + (-1,0,0,0,\ldots)) \\
& \leq (\dege(f \circ g)\oplus -d(\dege g + (-1,0,0,0,\ldots)))\oplus d (\dege g + (-1,0,0,0,\ldots)) \\
& = \dege(f \circ g).
\end{align*}
This proves the theorem in the case where $\dege f < (d,\infty,1,0,0,0,\ldots)$.  

Suppose now that $\dege f = (d,\infty,1,0,0,0,\ldots)$.  Then, since $\deg(f \circ g) = de_0$, one has
$$\dege (f \circ g)   \leq (de_0,\infty,1,0,0,0,\ldots).$$
On the other hand, if $\dege (f \circ g)  < (de_0,\infty,1,0,0,0,\ldots)$, then
\begin{align*}
\dege f & = \dege((f \circ g)\circ g^{-1}) \\ &
= \dege (f \circ g) \oplus de_0(\dege g^{-1} + (-1,0,0,0,\ldots)) \\
& < (d, \infty,1,0,0,0,\ldots),
\end{align*}
which is a contradiction.   It follows that
$$\dege (f \circ g)   = (de_0,\infty,1,0,0,0,\ldots)  = \dege f \oplus d(\dege g + (-1,0,0,0,\ldots)).$$
This completes the proof.
\end{proof}

Theorems  \ref{fgie} and \ref{fgie00} combine as follows.

\begin{theorem}\label{fgieth}
Let $f$ and $g$ be real functions defined on a neighborhood of $\infty$.  Suppose that the following conditions hold.
\begin{enumerate}
\item $f$  has finite degree $d$.
\item $g$ has positive degree and is eventually positive.
\item Either $f$ has exact logexponential degree or $g$ is
eventually continuous and increasing.
\item  $g(x) \asymp  r(x) \ (x \to \infty)$ for some  $r \in \mathfrak{L}$.
\end{enumerate}
Then one has 
$$\dege( f  \circ g )=  \dege f \oplus d (\dege g + (-1,0,0,0,\ldots)).$$
\end{theorem}

\begin{corollary}\label{fgie2}
Let $f$ and $g$ be real functions defined on a neighborhood of $\infty$.  Suppose that $g$ is eventually positive,    $g(x) \asymp x \ (x \to \infty)$, and either $f$ has exact logexponential degree or $g$ is eventually continuous and increasing.
Then one has 
$$\dege( f  \circ g )=  \dege f.$$
\end{corollary}

\begin{proof}
Let $h = f \circ g$. 
If $d \neq \pm \infty$, then by Theorem \ref{fgieth} one has
$$\dege h  = \dege f.$$
On the other hand, if $d = \pm \infty$, then $\deg h = \deg f$ and $h_{(1)} = f_{(1)} \circ g$,  by Proposition \ref{circle1}(2).
An obvious inductive argument, then, yields $$\dege_k h = \deg h_{(k)} = \deg f_{(k)} = \dege_k f$$ for all  nonnegative integers $k$.
\end{proof}

\begin{problem}
Generalize  Theorem \ref{fgieth} to the extent possible.
\end{problem}

\begin{remark}[Base independence of logexponential degree]\label{altbase}
Let $f \in \RR^{\RR_{\infty}}$, and let $b> 1$ and $d \in \RR$.  Then $$f(b^x)b^{-dx} = f(e^x)e^{-dx} \circ  (x \log b)$$
while
  $$\max( \log_{b} |f(x)|, 0) = \frac{1}{\log b} \max( \log |f(x)|, 0)$$
  and
  $$-\frac{1}{\log_{b} |f(x)|}  = (\log b) \left(-\frac{1 }{\log |f(x)|}\right).$$
  Let $b_k > 1$ for all nonnegative integers $k$.  Let $f_0 = f$, and supposing that $f_k$ is defined for some nonnegative integer $k$, let $d_k = \deg f_k$ and
$$f_{k+1}(x) =   \left.
\begin{cases}
    f_{k}(b_k^x) b_k^{-d_k x}  & \text{if } d_k \neq \pm \infty \\
    \max( \log_{b_k} |f_{k}(x)|, 0)  & \text{if } d_k = \infty \\
 \displaystyle   -\frac{1}{\log_{b_k} |f_{k}(x)|}  & \text{if } d_k =- \infty.
 \end{cases}
\right.$$  
By Corollary \ref{fgie2}, Proposition \ref{aspropexoexp}(3), and induction, one has
$$ \dege_k f = d_k = \deg f_k$$
for all $k$.   This shows that $\dege$ is independent of choice of bases.   Our convention to use $b_k = e$ for all $k$  is the most natural choice for functions arising in analytic number theory.   
\end{remark}

\begin{theorem}\label{thm:ledege_inverse}
Let $f$ be an increasing unbounded function of finite positive  degree with $\dege f = \dd$ (resp.,  $\underline{\dege}\, f = \dd$).     
Then the inverse function $f^{-1}$ exists, is  increasing and unbounded, and has lower logexponential degree $\underline{\dege}(f^{-1}) = \dd'$  (resp.,  logexponential degree $\dege  f^{-1} = \dd'$) given by
\[
\dd' = \left( \tfrac{1}{\dd_0},-\tfrac{\dd_1}{\dd_0},\ -\tfrac{\dd_2}{\dd_0}, -\tfrac{\dd_3}{\dd_0},\ldots \right),
\]
where each coordinate is given as above until the first $k$,  if any, such that $\dd_k = \pm\infty$, after which tail of $\dd'$ is exactly the negated tail of $\dd$, that is,  $\dd_j'= -\dd_j$ for all $j \geq k$.
\end{theorem}

\begin{proof}
Let $ \dd = \dege f$ and $\dd' = \underline{\dege}(f^{-1})$,  and let $r \in \mathbb{L}$ with $0< \deg r = \ee_0 < \infty$.  Observe that $f(x) \leq Cr(x)$ for all $x \gg 0$  in $\dom f$ if and only if $ f^{-1}(x) \geq r^{-1}(x/C) $ for all $x \gg 0$ in $\dom f^{-1}$.    Now,  by Karamata's inversion theorem \cite[Theorem 1.5.12]{bgt},  since $r(x) = x^{\ee_0} s(x)$, where $\ee_0 > 0$ and $s(x)$ is positive,  eventually monotonic, and slowly varying,  we have $$r^{-1}(x) \sim x^{1/\ee_0}( s(x^{1/\ee_0}))^{-1/\ee_0} \ (x \to \infty),$$
and therefore, since $s \in  \mathbb{L}$, we have
$$\underline{\dege}(f^{-1})  \geq \dege r^{-1}(x/C) =   \dege r^{-1} = \left( \tfrac{1}{\ee_0},-\tfrac{\ee_1}{\ee_0},\ -\tfrac{\ee_2}{\ee_0}, -\tfrac{\ee_3}{\ee_0},\ldots \right),$$
where $\ee = \dege s$, and where  the tail  is defined appropriately.  Noting that  $1/x$ and $-x$ are decreasing, the inequality \[
\dd' \geq \left( \tfrac{1}{\dd_0},-\tfrac{\dd_1}{\dd_0},\ -\tfrac{\dd_2}{\dd_0}, -\tfrac{\dd_3}{\dd_0},\ldots \right)
\]
 follows since $\underline{\dege}(f^{-1})$ is greater than or equal to the sup of $  \dege r^{-1}$ over all $r \in \mathbb{L}$ with $f = O(r)$.   Reversing the inequalities in this argument, we see that
  \[
\dd \leq \left( \tfrac{1}{\dd'_0},-\tfrac{\dd'_1}{\dd'_0},\ -\tfrac{\dd'_2}{\dd'_0}, -\tfrac{\dd'_3}{\dd'_0},\ldots \right)
\]
and therefore 
 \[
 \left( \tfrac{1}{\dd_0},-\tfrac{\dd_1}{\dd_0},\ -\tfrac{\dd_2}{\dd_0}, -\tfrac{\dd_3}{\dd_0},\ldots \right) \geq \dd',
\]
whence equality holds.   The theorem follows.
\end{proof}

\subsection{Axiomatization of $\dege$}

Recall that $ \RR^{\RR_\infty}$ denotes the set of  real functions whose domain is not bounded above.     The set $\prod_{n=0}^\infty \overline{\RR}$  is equipped with the lexicographic order and with the product topology where each factor $ \overline{\RR}$ has the discrete topology.  The following theorem shows that the map $\dege$ chan be characterized in terms of the relations $O$ and $o$ and the map $\dege$ on $ \mathbb{L}_{>0}$.
  
\begin{theorem}[Axiomatization of $\dege$]\label{thm:ledege_axioms}
Let
\[
\mathbf{D} : \RR^{\RR_\infty} \longrightarrow \prod_{n = 0}^{\infty} \overline{\RR}
\]
be a map.   One has $  \mathbf{D} r = \dege r$ for all $r \in  \mathbb{L}_{>0}$ if and only if $ \mathbf{D}$ satisfies the following four conditions.
\begin{enumerate}
    \item \textbf{Degree compatibility:} For all $r \in \mathbb{L}_{> 0}$, one has   $(\mathbf{D} r)_0 =  \deg r$.
    \item \textbf{Monomial compatibility:} For all $a \in \RR$ and all $r \in \mathbb{L}_{> 0}$ with $\deg r \in \RR$, one has
    \[
    \mathbf{D}(x^a r) = (a, 0, 0, \ldots) + \mathbf{D} r,
    \]
    where addition is coordinatewise.
\item  \textbf{Exp shifting:} For all  $r \in \mathbb{L}_{> 0}$ with $\deg r = 0$ and $\mathbf{D} r = \dd$, one has
\[\mathbf{D} (r \circ \exp) = (\dd_1, \dd_2, \ldots).\]
    \item \textbf{Log shifting:} For all  $r \in \mathbb{L}_{> 0}$ with $\deg r = \infty$ and $\mathbf{D} r = \dd$,  one has
        \[
        \mathbf{D}(\log \circ r) = (\dd_1, \dd_2, \dd_3, \ldots).
        \]
 For all  $r \in \mathbb{L}_{> 0}$ with $\deg r = -\infty$ and $\mathbf{D} r = \dd$,  one has
        \[
        \mathbf{D}(1/\log \circ r) = (\dd_1, \dd_2, \dd_3, \ldots).
        \]
\end{enumerate}
Moreover, one has $\mathbf{D} = \dege$ if and only if  $\mathbf{D}$  satisfies the following  three additional conditions.
\begin{enumerate}
    \item[(5)] \textbf{$O$ compatibility:} For all $f \in  \RR^{\RR_\infty}$ and all $r \in \mathbb{L}_{>0}$ with  $f(x) = O(r(x)) \ (x \to \infty)$,  one has  $\mathbf{D} f \leq \mathbf{D} r$.
    \item[(6)] \textbf{$o$ compatibility:}  For all $f \in  \RR^{\RR_\infty}$ and all $r \in \mathbb{L}_{>0}$ with $\mathbf{D} f < \mathbf{D} r$,  one has  $  f(x) = o(r(x)) \ (x \to \infty)$.
    \item[(7)] \textbf{Admissibility:} The image of  $\mathbf{D}$ is contained in the closure of $\mathbf{D}( \mathbb{L}_{>0})$ in $\prod_{n=0}^\infty \overline{\RR}$.
\end{enumerate}
If these seven conditions hold, then  the image of  $\dege$ is equal to the closure of $\mathbf{D}( \mathbb{L}_{>0})$ in $\prod_{n=0}^\infty \overline{\RR}$, and 
\begin{align*}
\mathbf{D} f& = \inf\{\mathbf{D} r:r \in \mathbb{L}_{> 0}, \, f(x) = O(r(x)) \ (x \to \infty) \} \\ 
& = \inf\{\mathbf{D} r:r \in \mathbb{L}_{> 0}, \, f(x) = o(r(x)) \ (x \to \infty) \} \\
& = \dege f
\end{align*}
for all $f \in  \RR^{\RR_\infty}$,  where the infima are computed in  the closure of $\mathbf{D}( \mathbb{L}_{>0})$ in $\prod_{n=0}^\infty \overline{\RR}$.
\end{theorem}

\begin{proof}
Since $T^{\circ n}(r) \in \mathbb{L}_{>0}$, properties (1)--(4) ensure that $(\mathbf{D} r)_0 = \dege_0 r$, $(\mathbf{D} r)_1 =( \mathbf{D} (T (r)))_0$,  $(\mathbf{D} r)_2 =( \mathbf{D} (T (r)))_1 = ( \mathbf{D} (T(T (r))))_0$,  and so on,  and therefore $\mathbf{D} r = \dege r$,  for all $r \in \mathbb{L}_{>0}$.   By Proposition \ref{degerprop},  it follows that the closure of  $\mathbf{D}( \mathbb{L}_{>0})$ in $\prod_{n=0}^\infty \overline{\RR}$ is equal to $\prod_{n=0}^{\infty *}  \overline{\RR}$, and   therefore, by property (7),   the image of $\mathbf{D}$ is contained in $\prod_{n=0}^{\infty *}  \overline{\RR}$. Then, by properties (5)--(7) and the proof of Theorem \ref{infpropexp},  for all $f \in  \RR^{\RR_\infty}$,  we have
\begin{align*}
\mathbf{D} f & = \inf\{\mathbf{D} r:r \in\mathbb{L}_{> 0} \, \mathbf{D} f \leq \mathbf{D} r \} \\ 
& = \inf\{\mathbf{D} r:r \in \mathbb{L}_{> 0}, \, f(x) = O(r(x)) \ (x \to \infty) \} \\ 
& = \inf\{\mathbf{D} r:r \in \mathbb{L}_{> 0}, \, f(x) = o(r(x)) \ (x \to \infty) \} \\
& = \inf\{\mathbf{D} r:r \in \mathbb{L}_{> 0}, \, \mathbf{D} f < \mathbf{D} r \},
\end{align*}
where the infima (exist and) are computed in  $\prod_{n = 0}^ {\infty *}\overline{\RR}$.   It follows, then, that \begin{align*}
\mathbf{D} f 
& = \inf\{\mathbf{D} r:r \in \mathbb{L}_{> 0}, \, f(x) = O(r(x)) \ (x \to \infty) \} \\ 
& = \inf\{\dege r:r \in \mathbb{L}_{> 0}, \, f(x) = O(r(x)) \ (x \to \infty) \} \\
& = \dege f
\end{align*} 
for all $f \in  \RR^{\RR_\infty}$.   Finally,  by Theorem \ref{degeequiv} (proved  in Section 7.5),  the image of $\dege$ is equal to $\prod_{n = 0}^ {\infty *}\overline{\RR}$.  This completes the proof.
\end{proof}

In  Section 7.6, we compute $\dege \mathbb{L}_{>0}$, and we  show that the map $\dege$ on $\mathbb{L}$ extends to a canonical map on the ordered differential field of well-based logarithmic-exponential transseries.

\section{Further results on degree and logexponential degree}

The following proposition can be useful for studying the logexponential degree of an arithmetic function.  (For example, it is used in the proof of Proposition \ref{lipiprop}.)

\begin{proposition}\label{arithb}
Let $f$ be a real-valued arithmetic function.   One has
$$\dege f = \dege f(\lfloor x \rfloor) = \dege f(\lceil x \rceil).$$
\end{proposition}

\begin{proof}
By Lemma \ref{reslemma}, one has
$$\dege f  \leq  \dege f(\lfloor x \rfloor)$$
and
$$\dege f  \leq  \dege f(\lceil x \rceil).$$
Let $r \in \mathbb{L}_{>0}$ with $f(n) = O(r(n)) \ (n \to \infty)$.  Then one has
$$f(\lfloor x \rfloor) = O(r(\lfloor x \rfloor)) \ (x \to \infty),$$  whence by Corollary \ref{fgie2} one has
$$\dege f(\lfloor x \rfloor) \leq \dege r(\lfloor x \rfloor) = \dege r(x).$$
Taking the infimum over all $r$ as chosen and applying Theorem \ref{infpropexp}, we find that
$$\dege f(\lfloor x \rfloor) \leq \dege f.$$
Similarly, one has
$$\dege f(\lceil x \rceil) \leq \dege f.$$
This completes the proof.
\end{proof}

Next, we show that all functions of finite degree in $\mathbb{L}$ are regularly varying.  This is quite useful in combination with Theorem \ref{infpropexp} and Karamata's integral theorem: for example, see the proofs of Theorem \ref{lithetapsi} in Section 9.3 and Theorem \ref{mertenssecond} in Section 10.1.

\begin{proposition}\label{kregvar}
A function $r \in \mathbb{L}$ is regularly varying if and only if it is of finite degree, and, if those equivalent conditions hold, then $r$ is regularly varying of index $$\overline{\underline{\deg}} \, r = \lim_{x \to \infty}\frac{xr'(x)}{r(x)}.$$
\end{proposition}

\begin{proof}
Let  $r \in \mathbb{L}^*$ with $\deg r \neq \pm \infty$.  Then $r$ is continuously differentiable on $[N,\infty)$ for some real number $N$.  Since $\frac{xr'(x)}{r(x)} \in \mathbb{L}$, the limit $d = \lim_{x \to \infty}\frac{xr'(x)}{r(x)}$
exists or is $\pm \infty$.  But then $d = \deg r$, so that $d \neq \pm \infty$.  Therefore, by Corollary \ref{regvarcor}, the function $r$ is regularly varying of index $d$.
\end{proof}

Next, we note the following.

\begin{proposition}\label{supprop}
Let $f \in \RR^{\RR_\infty}$ with $f$ unbounded on $X = \dom f$ at $\infty$.   Let $N \geq 0$, and suppose that $f$ is bounded  on $[N,x] \cap X$ for all $x \geq N$, and let
$$\widetilde{f}(x) = \sup_{f \in [N,x] \cap X} |f(t)|, \quad \forall x \in [N,\infty) \cap X.$$  Then $\widetilde{f}$ is nonnegative and nondecreasing on $[N,\infty) \cap X$ with $\lim_{x \to \infty} \widetilde{f}(x) = \infty$, and one has
 $ \dege \widetilde{f} = \dege f $.
\end{proposition}

\begin{proof}
Since $0 \leq |f(x)| \leq \widetilde{f}(x)$ for all $x \in [N,\infty) \cap X$, one has $\dege f \leq \dege \widetilde{f}$.  To prove the reverse inequality,  let $r \in \mathbb{L}_{>0}$ with $|f(x)|  \leq r(x)$ for all  $x \in [N',\infty) \cap X$.  Since $f(x)$ is unbounded on $[N',\infty) \cap X$, the function $r$ must be eventually increasing, say,  increasing for $x \geq M \geq \max(N,N')$.  Then, for all sufficiently large $x \in [M,\infty) \cap X$, one has
\begin{align*}
\widetilde{f}(x) & = \max\left(\sup_{t \in [N,M] \cap X} |f(t)|, \sup_{t \in [M,x] \cap X}| f(t) |\right) 
\\ & \leq  \max\left(\sup_{t \in [N,M] \cap X} |f(t)|, \sup_{t \in [M,x] \cap X} r(t) \right) \\ 
& = \max\left(\sup_{t \in [N,M] \cap X} |f(t)|, r(x) \right) \\
& = r(x).
\end{align*}
Therefore, one has $\dege \widetilde{f} \leq \dege r$.
Taking the infimum over all $r$ as chosen, we find that $\dege \widetilde{f} \leq \dege f$.  This completes the proof.
\end{proof}

The corresponding result for lower logexponential degree is as follows.

\begin{proposition}\label{suppropinf}
Let $f \in \RR^{\RR_\infty}$ with $\liminf_{x \to \infty} f(x)= 0$,  and let  $X = \dom f$.   Let $N \geq 0$, and suppose that $f$ is bounded away from $0$  on $[N,x] \cap X$ for all $x \geq N$, and let
$$\widetilde{f}(x) = \inf_{f \in [N,x] \cap X} |f(t)|, \quad \forall x \in [N,\infty) \cap X.$$  Then $\widetilde{f}$ is nonnegative and nonincreasing on $[N,\infty) \cap X$ with $\lim_{x \to \infty} \widetilde{f}(x) = 0$, and one has
 $ \underline{\dege}\, \widetilde{f} = \underline{\dege}\, f $.
\end{proposition}

Regarding nondecreasing nonnegative step functions, we have the following.

\begin{proposition}\label{stepdeg}
Let $f$ be a nondecreasing nonnegative step  function on $[0,\infty)$ that is continuous from the right, and whose set of points of discontinuity  is an unbounded subset of $X = \{x_1,x_2,\ldots\}$, where the sequence $x_n$ is increasing without bound.   Let $U = f|_X$, and let $L$ be defined on $X\backslash\{x_1\}$, with $L(x_n) = f(x_n^-)$, and therefore $L(x_n)= f(x_{n-1})\leq f(x_n) = U(x_n)$, for all $n \geq 1$.    Then one has
$$\underline{\dege} \, f = \underline{\dege}\, L \leq\dege U = \dege f .$$
\end{proposition}

\begin{proof}
Since $U$ is a restriction of $f$, one has $\dege U \leq \dege f$.  Let $r \in \mathbb{L}_{>0}$ with $U(x_n) \leq r(x_n)$ for all $n$.   Since $U$ is nondecreasing,  we may suppose without loss of generality that $r$ is eventually nondecreasing.  Since $r(x)$ is nondecreasing on  $[x_n,x_{n+1})$, one has $f(x) = U(x_n) \leq r(x_n) \leq r(x)$,  for all $n \gg 0$ and all $x \in[x_n,x_{n+1})$.   It follows that $f(x) \leq r(x)$ for all $x \gg 0$, whence $\dege f \leq \dege U$, and equality holds.

Now, let $r \in \mathbb{L}_{>0}$ with $L(x_n) \geq r(x_n)$ for all $n \gg 0$.  Again, we may suppose without loss of generality  that $r$ is eventually nondecreasing.   Since $r(x)$ is nondecreasing on  $[x_{n-1},x_{n})$, one has $f(x) = L(x_n) \geq r(x_n) \geq r(x)$,  for all $n \gg 0$ and all $x \in[x_{n-1},x_{n})$.    It follows that $\underline{\dege} \, f \geq \underline{\dege}\, L$.  Finally, let $r \in \mathbb{L}_{>0}$ with $f(x) \geq r(x)$ for all $x \gg 0$.    Then, by the eventual continuity of $r$, one has $L(x_n) = f(x_n^-) \geq r(x_n)$ for all $n \gg 0$, and therefore $\underline{\dege} \, f \leq \underline{\dege}\, L$.  This completes the proof.
\end{proof}

Thus, for example, if $f$ is any nonnegative arithmetic function, then $\dege S_f = \dege S_f|_{\ZZ_{>0}}$ and $\underline{\dege} \, S_f = \underline{\dege} \, (S_f-f)$.

Finally,  we note the following.

\begin{proposition}\label{strongdegg}
Let $f$ be a real function that is Riemann integrable on $[N,x]$ for all $x>N$,  where $N  >0$ and $\deg f \in \RR\backslash\{-1\}$.   Let  $F(x) = \int_N^x f(t) \, dt$ or  $F(x) = \int_x^\infty f(t) \, dt$ on $(N,\infty)$ according to whether $\deg f > -1$ or $\deg f < -1$.   Let $s \in \RR$.
\begin{enumerate}
\item One has
$$\dege F \leq \dege f + (1,0,0,0,\ldots).$$
\item If $\deg f >-1$ and $s < \deg F$, then
$$\dege \int_N^x \frac{f(t)}{t^s} \, dt \leq \dege F + (-s,0,0,0,\ldots),$$
with equality if $\deg F > 0$.
\item If $\deg f >-1$ and $s > \deg f+1$, then
$$\dege \int_x^\infty \frac{f(t)}{t^s} \, dt \leq \dege F + (-s,0,0,0,\ldots),$$
with equality if $\deg F > 0$.
\item If $\deg f  < -1$ and $s < \deg F$, then
$$\dege \int_N^x \frac{f(t)}{t^s} \, dt = \dege F + (-s,0,0,0,\ldots).$$
\item If $\deg f  < -1$ and $s > \deg f + 1$, then
$$\dege \int_x^\infty \frac{f(t)}{t^s} \, dt = \dege F + (-s,0,0,0,\ldots).$$
\end{enumerate}
\end{proposition}

\begin{proof}
To prove (1), we suppose that $\deg f > -1$.    The proof when $\deg f < -1$ is similar.
Suppose that $f(x) = O(r(x)) \ (x \to \infty)$, where $r$ is regularly varying of index $d$ and both positive and continuous on $[N,\infty)$.  Since $d \geq \deg f >-1$,   Karamata's integral theorem implies that
$$\frac{1}{x}\int_N^x f(t)\, dt \ll \frac{1}{x} \int_N^x r(t)\, dt \asymp  r(x) \ (x \to \infty),$$
and therefore
$$\dege \int_N^x f(t)\, dt + (-1,0,0,0,\ldots) \leq \dege r.$$
Statement (1) follows by taking the infimum of $\dege r$ over all $r$ as chosen.

Suppose, now, that $\deg f > -1$ and $s < \deg F$.  Let
 $$F_s(x) = \int_N^x \frac{f(t)}{t^s} \, dt.$$
By Riemann--Stieltjes integration by parts \cite[Section 1.1.3]{borg} (along with \cite[Chapter I Theorem 6a]{widd}), one has
$$F_s(x) = \frac{F(x)}{x^s} +  s \int_N^x \frac{F(t)}{t^{s+1}} \, dt.$$
Let $r$ be regularly varying of index $d$ and both positive and continuous on $[N,\infty)$.  Suppose that $F(x) \ll r(x) \ (x \to \infty)$.
Then $d-s-1 \geq \deg F-s-1 > -1$, so that, by Karamata's integral theorem, one has
$$F_s(x) \ll \frac{r(x)}{x^s} +  |s| \int_N^x \frac{r(t)}{t^{s+1}} \, dt \ll  \frac{r(x)}{x^s} \ (x \to \infty),$$
and therefore 
$$\dege F_s \leq \dege r+(-s,0,0,0,\ldots).$$
Taking the infimum over all $r$ as chosen,  we see that $$\dege F_s \leq \dege F+(-s,0,0,0,\ldots).$$
Again by Riemann--Stieltjes integration by parts, one has
$$F(x)  = x^s F_s(x) -  s \int_N^x t^{s-1}F_s(x) \, dt.$$
Therefore, if $\deg F_s>-s$,  that is,  if $\deg F_s+s-1 >-1$, then,  by similar reasoning as above,  one has
$$\dege F \leq \dege F_s + (s,0,0,0,\ldots).$$  On the other hand,  
$\deg F_s \leq -s$ implies  $\deg F \leq 0$.    Statement (2) follows.

Suppose, now, that $\deg f > -1$ and $s > \deg f+1$, so that $s > \deg F$.  Let
 $$G_s(x) = \int_x^\infty \frac{f(t)}{t^s} \, dt.$$
By Riemann--Stieltjes integration by parts,  since  $\lim_{x \to \infty} \frac{F(x)}{x^s} = 0$,  one has
$$G_s(x) = -\frac{F(x)}{x^s} +  s \int_x^\infty \frac{F(t)}{t^{s+1}} \, dt.$$
Let $r$ be regularly varying of index $d$ and both positive and continuous on $[N,\infty)$.  Suppose that $F(x) \ll r(x) \ (x \to \infty)$, with $\deg F \leq d < s$.  Then,  by Karamata's integral theorem, one has
$$G_s(x) \ll \frac{r(x)}{x^s} +  |s| \int_N^x \frac{r(t)}{t^{s+1}} \, dt \ll  \frac{r(x)}{x^s} \ (x \to \infty),$$
and therefore 
$$\dege G_s \leq \dege r+(-s,0,0,0,\ldots).$$
Taking the infimum over all $r$ as chosen,  we see that $$\dege G_s \leq \dege F+(-s,0,0,0,\ldots).$$
Again by Riemann--Stieltjes integration by parts, one has
$$F(x)  = N^sG_s(N)-x^s G_s(x) + s \int_N^x t^{s-1}G_s(x) \, dt.$$
Therefore, if $\deg G_s>-s$,  that is,  if $\deg G_s+s-1 >-1$, then,  by similar reasoning as above, one has
$$\dege F \leq \dege G_s + (s,0,0,0,\ldots).$$  On the other hand,  
$\deg G_s \leq -s$ implies  $\deg F \leq 0$.    Statement (3) follows.

Suppose, now, that $\deg f < -1$ and $s< \deg F$, so that $\deg  F <0$.   Let
 $$F_s(x) = \int_N^x \frac{f(t)}{t^s} \, dt.$$
By Riemann--Stieltjes integration by parts, one has
$$F_s(x) = \frac{F(N)}{N^s}-\frac{F(x)}{x^s}  - s \int_N^x \frac{F(t)}{t^{s+1}} \, dt.$$
Let $r$ be regularly varying of index $d$ and both positive and continuous on $[N,\infty)$.  Suppose that $F(x) \ll r(x) \ (x \to \infty)$.
Then $d-s \geq \deg F-s > 0$ and $d-s-1>-1$, so that
$$F_s(x) \ll \left| \frac{F(N)}{N^s} \right|+\frac{r(x)}{x^s} + |s| \int_N^x \frac{r(t)}{t^{s+1}} \, dt \ll  1+\frac{r(x)}{x^s} \ (x \to \infty)$$
and therefore 
$$\dege F_s \leq \max((0,0,0,\ldots), \dege r+(-s,0,0,0,\ldots)) =  \dege r+(-s,0,0,0,\ldots).$$
Taking the infimum over all $r$ as chosen,  we see that  $$\dege F_s \leq \dege F+(-s,0,0,0,\ldots).$$
Thus, one has $\deg F_s \leq \deg F-s < -s$, so that $\deg F_s+s-1 <-1$.
Again by Riemann--Stieltjes integration by parts, one has
$$F(x)  = -x^s F_s(x) + s \int_x^\infty t^{s-1}F_s(x) \, dt.$$
Therefore, by similar reasoning as above, one has
$$\dege F \leq \dege F_s + (s,0,0,0,\ldots).$$  
 Statement (4) follows.  

Finally, suppose that $\deg f < -1$ and $s> \deg f+1$, so that $\deg f-s <-1 $ and $\deg F  \leq \deg f+1 <s$.   Let
 $$G_s(x) = \int_x^\infty \frac{f(t)}{t^s} \, dt.$$
By Riemann--Stieltjes integration by parts, since  $\lim_{x \to \infty} \frac{F(x)}{x^s} = 0$,  one has
$$G_s(x) = \frac{F(x)}{x^s} -  s \int_x^\infty \frac{F(t)}{t^{s+1}} \, dt.$$
Let $r$ be regularly varying of index $d$ and both positive and continuous on $[N,\infty)$.  Suppose that $F(x) \ll r(x) \ (x \to \infty)$,
with $\deg F \leq d < s$.
Then $d-s-1< -1$, so that
$$G_s(x) \ll \frac{r(x)}{x^s} +  |s| \int_x^\infty \frac{r(t)}{t^{s+1}} \, dt \ll  \frac{r(x)}{x^s} \ (x \to \infty).$$
It follows that $$\dege G_s \leq \dege F+(-s,0,0,0,\ldots).$$
Thus, one has $\deg G_s \leq \deg F-s < -s$, so that $\deg G_s+s-1 <-1$.
Again by Riemann--Stieltjes integration by parts,  since  $\lim_{x \to \infty} x^s G_s(x) = 0$, one  has
$$F(x)  = x^s G_s(x) + s \int_x^\infty t^{s-1}G_s(x) \, dt.$$
Therefore,  by similar reasoning as above, one has
$$\dege F \leq \dege G_s + (s,0,0,0,\ldots).$$
Statement (5) follows.   This completes the proof.
\end{proof}

\chapter{Asymptotic algebra}

In this chapter, we provide applications of asymptotic differential algebra (e.g., Hardy fields) to the study of logexponential degree, and vice versa.    Nearly all of this chapter is not used in later chapters,  the only exceptions being the definitions  in Section 7.3 of a {\it Hardian} function and of the ordered field $\mathbb{H} \supsetneq \mathbb{L}$ of all  {\it universally Hardian} functions,  along with Theorem \ref{hardintth} and Propositions \ref{hardianexactlog} and \ref{oexppropstrong}.

\section{Lattice-ordered rings}

In this  section, we generalize the various asymptotic relations studied in Section 2.1 to the setting of {\it lattice-ordered rings},  providing a modest supplement to the well-developed abstract approach  to asymptotic relations taken by asymptotic differential algebra, e.g., in \cite{asch}, and by the theory of ordered exponential fields, e.g., in \cite{kuhl}.  We describe the abstract approach to asymptotic relations, with or without differentials, as {\it asymptotic algebra}.

A {\bf partially ordered abelian group}\index{partially ordered abelian group} is an abelian group $G$ (which we write additively) equipped with a partial ordering $\leq$ on $G$ such that $a \leq a'$ and $b \leq b'$ implies $a+ b \leq a'+b'$ for all $a,a',b,b' \in G$.   A {\bf (totally) ordered  abelian group}\index{ordered abelian group} is a partially ordered abelian group $G$ such that the partial ordering on $G$ is a total ordering.   A {\bf lattice-ordered abelian group}\index{lattice-ordered abelian group} is a partially ordered abelian group $G$ such that the partial ordering on $G$ is a lattice ordering, i.e., such that the least upper bound $a \vee b$ and greatest lower bound $a \wedge b$  of $a$ and $b$ exist in $G$ for all $a, b \in G$.  Note that,  if least upper bounds exists in a partially ordered abelian group $G$, then so do greatest lower bounds, since then $a \wedge b = -((-a)\vee (-b))$ for all $a, b \in G$.   The same holds with  $\wedge$ and $\vee$ interchanged.  Less  obviously,  by \cite[Chapter XIV Theorem 2]{birkh}, a partially ordered abelian group $G$ is lattice-ordered  if and only if $a \vee 0$ exists for all $a \in G$.

Let $G$ be a lattice-ordered abelian group.   For any $a \in G$,  let $|a| = a \vee (-a)$ denote the {\bf absolute value} of $a$, and  let $a^+ = a \vee 0$ and $a^- = (-a) \vee 0 = -(a \wedge 0)$ denote the {\bf positive part} and {\bf  negative part} of $a$, respectively.  One has the following elementary result.

\begin{proposition}[{\cite[Proposition 1.3]{johnson}}]
Let $G$ be a lattice-ordered abelian group, and let $a,b,c \in G$, and let $n$ be a positive integer.  One has the following.
\begin{enumerate}
\item $a+(b  \vee c) = (a+b) \vee (a+c)$ and $a+(b  \wedge c) = (a+b) \wedge (a+c)$.
\item $a+b = (a \vee b)+(a \wedge b)$.
\item $|a+b| \leq |a|+|b|$ and $|a-b| \geq ||a|-|b||$.
\item $n(a \vee b) = (na) \vee (nb)$ and $n(a \wedge b) = (na) \wedge (nb)$.
\item $|na| =n|a|$.
\item If $na \geq 0$, then $a \geq 0$.
\item $a$ has infinite order.
\item $a = a^+-a^-$.
\item $|a| = a^++a^- \geq 0$.
\item $a^+ \wedge a^- = 0$.
\item $|a| = 0$ if and only if $a= 0$.
\end{enumerate}
\end{proposition}

All rings in this chapter are assumed commutative with identity.

A {\bf partially ordered ring}\index{partially ordered ring} is a ring $R$ equipped with a partial ordering $\leq$ on $R$ with $1 \geq 0$ such that 
$R$ is a partially ordered abelian group under addition and $a,b \geq 0$ implies $ab \geq 0$
for all $a,b \in R$ (or, equivalently, $a \geq b$ implies $ac \geq bc$ for all $a,b,c \in R$ with $c \geq 0$.   A {\bf totally ordered ring}\index{totally ordered ring} (resp., {\bf lattice-ordered ring})\index{lattice-ordered ring} is a partially ordered ring whose partial ordering is a total ordering (resp.,  a lattice ordering).  Note that, if $R$ is a lattice-ordered ring, then
$$|ab| \leq |a||b|$$
for all $a,b \in R$.   If $R$ is a totally ordered ring (hence also a lattice-ordered ring), then 
$$|a|= \max(a,-a),$$
$$|ab| = |a||b|,$$
$$a^2 \geq 0,$$
and $a \wedge b = 0$ implies $a = 0$ or $b = 0$, for all $a, b \in R$.   Moreover, none of these four properties need hold in a lattice-ordered ring. 

The theories of  {\it lattice-ordered groups}, or {\it $\ell$-groups},  and {\it lattice-ordered  rings}, or {\it $\ell$-rings}---not necessarily abelian or commutative, respectively---are quite developed.  See \cite{steinberg}, for example.

Let $R$ be a lattice-ordered ring containing a totally ordered subfield $k$.   Let us define the following {\bf asymptotic relations on $R$},\index{asymptotic relations} all defined relative to $k$: for all $a, b \in R$, we let
\begin{enumerate}
\item $a \preceq b$ if  $|a| \leq r|b|$ for some $r \in k_{>0}$;\index[symbols]{.e hh1@$\preceq$}
\item $a \asymp b$ if  $a \preceq b$ and $b \preceq a$;\index[symbols]{.e ee@$\asymp$}
\item $a \prec b$  if  $|a| \leq r|b|$ for all $r \in k_{>0}$;\index[symbols]{.e hh2@$\prec$}
\item $a \sim b$ if $a-b \prec b$.\index[symbols]{.e hh@$\sim$}
\end{enumerate}
Note that, although these four relations depend on the choice of $k$, different choices of $k$ can yield the same relations,  e.g., any subfield of $k = \RR$ yields the same asymptotic relations as does $\RR$ itself.   Also note that, as long as $R$ contains a field, $R$ contains a unique minimal totally ordered subfield,  necessarily isomorphic to $\QQ$.   Thus,  the four relations can be defined as long as $R$ contains a field,  and then its  minimal subfield is a natural choice for $k$.  We call the four asymptotic relations on $R$ defined with respect to  its  minimal subfield the {\bf standard asymptotic relations on $R$}.  Note that,  if $k \subseteq k'$ are totally ordered subfields of $R$, then the relations $\preceq$ and $\asymp$ with respect to $k$ are finer than those with respect to $k'$, while the relations $\prec$ and $\sim$  with respect to $k$ are coarser than those  with respect to $k'$ (because of properties of the existential and universal quantifiers, respectively).

The following proposition notes several properties of these four relations on $R$ and their relationships to one another.

\begin{proposition}\label{lor1}
Let $R$ be a lattice-ordered ring containing a totally ordered subfield $k$, and let $a,b,c \in R$.   One has the following.
\begin{enumerate}
\item $a \preceq b$ if and only if $|a| \preceq |b|$, and likewise for $\prec$ and for $\asymp$.
\item $a \preceq a$, and likewise for $\prec$ and for $\asymp$.
\item If  $a \preceq b$ and  $b \preceq c$, then  $a \preceq c$.
\item If  $a \preceq b$ and  $b \prec c$, then  $a \prec c$.
\item If  $a \prec b$ and  $b \preceq c$, then  $a \prec c$.
\item  $a \asymp b$ if and only if $b \asymp a$.
\item  $a \sim b$ if and only if $b \sim a$.
\item If  $a \prec b$ or $a \asymp b$, then  $a \preceq b$.
\item If $a \sim b$, then $a \asymp b$.
\item If  $a \sim b$ and  $b \sim c$, then  $a \sim c$.
\item $0 \prec a$.
\item $a \preceq 0$ if and only if $a \prec 0$, if and only if $a \prec a$,  if and only if $a \sim 0$,  if and only if $a = 0$.
\end{enumerate}
\end{proposition}

The following proposition relates the asymptotic relations on $R$ to the ring operations on $R$.

\begin{proposition}\label{lor2} Let $R$ be a lattice-ordered ring containing a totally ordered subfield $k$,  let $a_1, a_2, b_1, b_2, a, b \in R$, and let $r_1, r_2 \in k$.  One has the following.
\begin{enumerate}
\item $|a| \vee |b| \leq |a|+|b| \leq  2( |a| \vee |b|)$, and therefore $|a|+|b| \asymp |a| \vee |b|$.
\item If $ a_{1} \preceq b_{1} $ and $a_{2}\preceq b_{2},$ then 
$ a_{1}+a_{2} \preceq |b_{1}| + |b_{2}| $ and
$a_{1}a_{2} \preceq |b_{1}||b_{2}|.$
\item If  
$ a_{1} \preceq b  \text{ and }  a_{2} \preceq b$,
then  $r_1 a_1 + r_2 a_2  \preceq b.$
\item If $ a_{1} \prec b_{1} $ and $a_{2}\prec b_{2}$, then 
$ a_{1}+a_{2} \prec |b_{1}| + |b_{2}| $ and
$a_{1}a_{2} \prec |b_{1}||b_{2}|.$
\item If  
$ a_{1} \prec b  \text{ and }  a_{2} \prec b$,
then  $r_1 a_1 + r_2 a_2  \prec b.$
\item If $ a_{1} \prec b_{1} $ and $a_{2}\preceq b_{2}$, then 
$a_{1}a_{2} \prec |b_{1}||b_{2}|.$
\item If $a_{1} \sim b_{1}$, $a_{2}\sim b_{2}$ and $|b_{1} b_{2}| = |b_{1}||b_{2}|$, then 
$a_{1}a_{2} \sim b_{1}b_{2}.$
\item If $a_1 \sim b$, then $a_1 + a_2 \sim b$ if and only if $a_2 \prec b$.
\end{enumerate}
\end{proposition}

Note that, by statement (1) of the proposition, one can replace $|b_{1}| + |b_{2}| $ in statements (2) and (4) with
$|b_{1}| \vee |b_{2}|$.

Assume now an additional condition on the ring $R$,  namely, that $|ab| = |a||b|$ for all $a,b \in R$.   This holds, for example, if $R$ is totally ordered, or, more generally, if  $R$ is an {\it f-ring} \cite{johnson} \cite{steinberg}.  Then we may define an addition $[a]+[b] = [ab]$ on $\asymp$-equivalence classes $[a],[b] \in  R/{\asymp}$, as well as a partial ordering $[a]\leq [b]$ on $R/{\asymp}$, defined to hold precisely when $b \preceq a$.  Then $R/{\asymp}$ is a lattice-ordered commutative monoid with identity element $0 = [1]$, that is, $R/{\asymp}$ is a commutative monoid with identity element $[1]$ and $\leq$ is a lattice ordering on $R/{\asymp}$ such that $[a_1]\leq [b_1]$ and $[a_2]\leq [b_2]$ implies $[a_1]+[a_2]\leq [b_1]+[b_2]$ for all $a_1,a_2,b_1,b_2 \in R$.   Note that $[|a| \vee |b|] = [a] \vee [b]$ for all $a,b \in R$.   Now, let
$$v: R \longrightarrow R/{\asymp}$$
act by $$v: a \longmapsto [a].$$  Then the map $v$ is surjective and satisfies the following conditions for all $a,b \in R$.
\begin{enumerate}
\item $v(1) = 0$.
\item $v(ab)  = v(a)+v(b)$.
\item $v(a+b) \geq v(a) \wedge v(b)$.
\end{enumerate}
Indeed,  (1) and (2) are obvious, and (3) follows immediately from
$$|a+b| \leq  |a|+|b| \asymp |a| \vee |b|.$$
We denote  $[0]$ by $\infty$, since $[0]+[a]=[0]$ for all $a\in R$, and we let $\Gamma_R = (R/{\asymp}) \backslash \{\infty\}$, so that $v$ restricts to a map
$$v: R\backslash\{0\} \longrightarrow \Gamma_R.$$   Note that $R$ is an integral domain if and only if $R\backslash\{0\}$ is a submonoid of $R^\times$, if and only if $0 \neq \infty$ and $[a]+[b] =  \infty$ implies $[a]= \infty$ or $[b] = \infty$ for all $a,b \in R$, if and only if $\Gamma_R$ is a submonoid of $R/{ \asymp}$.   Moreover, if $R$ is a lattice-ordered field,  or {\bf $\ell$-field} \cite{steinberg}, then $\Gamma_R$ is an abelian group.  Likewise,  if $R$ is totally ordered, then $\Gamma_R$ is totally ordered and $$v(a+b) \geq v(a) \wedge v(b) = \min(v(a),v(b))$$ for all $a,b \in R$.   

An important special case, employed to a great extent in the next section, is where $R$ is a totally ordered field.    In that case, the map $v: R \longrightarrow R/{\asymp} \ =  \Gamma_R \cup \{\infty\}$ is a {\it valuation on $R$} and the relation $\preceq$ is a {\it dominance relation on $R$} \cite[Section 3.1]{asch}.   Let $K$ be any field.  A {\bf valuation on $K$}\index{valuation on a field} is a surjective map $v: K \longrightarrow \Gamma \cup \{\infty\}$, where $\Gamma$ is a (totally) ordered abelian group,  and where one sets $\gamma <\infty = \infty+\infty = \gamma+\infty = \infty+\gamma$ for all $\gamma \in \Gamma$,  such that $v(a) = 0$ if and only if $a = \infty$,  and such that $v(ab)  = v(a)+v(b)$ and $v(a+b) \geq \min(v(a) , v(b))$ for all $a,b \in K$.  Moreover,  the {\bf dominance relation  $\preceq$  on $K$ associated to $v$},  defined by  $a \preceq b$ if $v(a) \geq v(b)$, coincides with our asymptotic relation $\preceq$  defined earlier.  The same holds of the relations $\asymp$, $\prec$, and $\sim$ associated to $v$ as defined in \cite[Section 3.1]{asch}, where, specifically, one writes $a \asymp b$ if $a \preceq b$ and $b \preceq a$,  one writes $a \prec b$ if $a \preceq b$ and $b \not \preceq a$, and one writes $a \sim b$ if $a-b \prec b$.

A {\bf valued field}\index{valued field} is a field $K$ equipped with a valuation  $v: K \longrightarrow  \Gamma \cup \{\infty\}$.   The ordered abelian group $\Gamma_v = \Gamma$ is called the {\bf value group of $v$}. \index{value group of a valuation} Let $K$ be a valued field with valuation $v$.  The subset $$K_v = \{a \in K: v(a) \geq 0\}$$ of $K$ is a subring of $K$ such that $a \in K_v$ or $1/a \in K_v$ for all $a \in K^*$ and is called the {\bf valuation ring of $v$}. \index{valuation ring of a valuation}  It follows that the  ring $K_v$ is  an integral domain with unique  maximal ideal
$$\mm_v = \{a \in K: v(a) > 0\}.$$ 
 It is well known that there are natural one-to-one correspondences between  the valuations on $K$ up to {\it equivalence}, the valuation rings of the valuations on $K$, the {\it valuation rings of $K$}, and the dominance relations on $K$ \cite[Section 3.1]{asch}.

To summarize, we have shown that any totally ordered field  $R$ equipped with a totally ordered subfield $k$ induces  a valuation $v: R \longrightarrow \Gamma_R \cup \{\infty\} = R/{ \asymp}$ on $R$, which is canonical for  $k \cong \QQ$.   Moreover,  the map $v$ and the asymptotic relations $\preceq$, $\asymp$, $\prec$, and $\sim$ associated to $v$ as in  \cite[Section 3.1]{asch} generalize to any lattice-ordered ring $R$ containing a distinguished totally ordered subfield $k$, and they share many of the properties as they do when $R$ is a totally ordered field, e.g., those properties expressed in Propositions \ref{lor1} and \ref{lor2}.

\section{The ring of germs of real functions at $\infty$}

If $R$ is a ring  (commutative with identity) and $X$ is a set, then we denote by $R^X = \prod_{x \in X}R$ the ring of all functions from $X$ to $R$ under   pointwise addition and multiplication.  If $R$ is a partially ordered ring, then $R^X$ is a partially ordered ring  when ordered by the relation $\leq$, where $f \leq g$ if $f(x) \leq g(x)$ for all $x \in X$.      Moreover,  if $R$ is a lattice-ordered ring, then $R^X$ is also a lattice-ordered ring, with $(f \vee g)(x) = f(x) \vee g(x)$ for  all $x \in X$, and likewise for $\wedge$.  In particular, the ring $\RR^X$ is lattice-ordered.
Also when $R = \RR$, one has $f \leq g$ if and only if $g-f = h^2$ for some $h \in \RR^X$. 
We extend any function $f \in \RR^{\RR_\infty}$ to a function on $\RR$ by defining it to be $0$ at all $x \in \RR\backslash \dom f$.    Thus we have a natural surjection  $ \RR^{\RR_\infty} \longrightarrow \RR^{\RR} = \prod_{x \in \RR} \RR$.

The  ring $\mathcal{R}\index[symbols]{.i  ka@$\mathcal{R}$}$ of germs at $\infty$ of all functions in $\RR^\RR$, that is, the ring $\mathcal{R} = \RR^{\RR}/{ =_\infty}$, is a lattice-ordered ring under the relation $\leq \ = \ \leq_{\infty}$.     Another way to construct $\mathcal{R}$ is as the quotient ring $(\RR^{\RR})_\infty  = \RR^{\RR}/(0)_\infty$ by the (radical) ideal $(0)_\infty$ in $\RR^\RR$ of all functions  that vanish on a  neighborhood of $\infty$.    Note that $f \leq g$ in $\mathcal{R}$ if and only if $g-f = h^2$ for some $h \in \mathcal{R}$.    The group of units $\mathcal{R}^*$ of the ring $\mathcal{R}$ is equal to the set of all germs of functions that are eventually nonzero, and it contains $\mathcal{R}_{>0}$ and $\mathcal{R}_{>0} \cup \mathcal{R}_{<0} \cong \mathcal{R}_{>0} \oplus \{1,-1\}$ as proper subgroups.   Moreover, the monoid of nonzerodivisors of $\mathcal{R}$ is equal to $\mathcal{R}^*$, that is, $\mathcal{R}$ is a total quotient ring.

A (commutative) ring $R$ is {\bf von Neumann regular}\index{von Neumann regular ring}  if $I^2 = I$ for all principal ideals (or equivalently, all ideals) $I$ of $R$ \cite{goodearl}.   A ring $R$ is von Neumann regular  if and only if  $R$ is a reduced total quotient ring in which every prime ideal is maximal.   The von Neumann regular rings are a natural generalization of the fields to rings with zerodivisors, since, for example,  a field is equivalently a von Neumann regular ring that is an integral domain.  The class of  von Neumann regular rings is closed under arbitrary direct products and quotient rings.   Thus, the rings $\RR^\RR$ and $\mathcal{R}$ are von Neumann regular.

The absolute value on the lattice-ordered ring $\mathcal{R}$, as described in the previous section, coincides with the usual absolute value on $\mathcal{R}$ and therefore satisfies $|fg| = |f||g|$ for all $f,g \in \mathcal{R}$.   Moreover, since  $\mathcal{R}$ is a lattice-ordered ring containing the totally ordered field $\RR$,  the results of the previous section yield (standard) asymptotic relations $\preceq$,  $\asymp$,  $\prec$, and $\sim$ on $\mathcal{R}$, and it  is clear that these relations coincide with the usual  $O$, $\asymp$, $o$, and $\sim$ relations on $\mathcal{R}$.     In other words,  for all $f, g \in \mathcal{R}$,  one has $f \preceq g$ if and only if  $f(x) = O(g(x)) \ (x \to \infty)$,  and so on.    
(Note that, on $\RR^X$, the former asymptotic relations must hold globally,  e.g.,  one has $f \preceq g$ in $\RR^X$ if and only if there exists an $M \in \RR_{>0}$ such that $|f(x)| \leq M |g(x)|$ for all $x \in X$.  Unfortunately,  then, one has $f \prec g$ in $\RR^X$ if and only if $f = 0$.)  Consequently,  from Propositions \ref{lor1} and \ref{lor2}, we deduce many of the well-known properties of the asymptotic relations $O$, $\asymp$, $o$, and $\sim$ on $\mathcal{R}$  stated in Section 2.1.   Moreover, if $g \in \mathcal{R}^*$, that is, if $g$ is eventually nonzero, then $f \preceq g$ is equivalent to  $\limsup_{x \to \infty} \frac{|f(x)|}{|g(x)|} < \infty$,  while $f \prec g$ is equivalent to  $\lim_{x \to \infty} \frac{|f(x)|}{|g(x)|} = 0$ and $f \sim g$ is equivalent to  $\lim_{x \to \infty} \frac{|f(x)|}{|g(x)|} = 1$.  Likewise, if $f \in \mathcal{R}^*$, then $f \preceq g$ is equivalent to  $\liminf_{x \to \infty} \frac{|g(x)|}{|f(x)|} > 0$, and $f \prec g$  is equivalent to $\lim_{x \to \infty} \frac{|g(x)|}{|f(x)|} = \infty$.    If $R$ is any subring of $\mathcal{R}$,  then we let $\preceq$,  $\asymp$,  $\prec$ and $\sim$ denote the relations on $\mathcal{R}$ restricted to $R$.

Now, one has $e^f \in \mathcal{R}_{>0}$ if $f \in \mathcal{R}$, and $\log f \in \mathcal{R}$ if $f \in \mathcal{R}_{>0}$, and $e^f \leq e^g$ if and only if $f \leq g$, while also $e^{f+g} = e^f e^g$, for all $f, g \in \mathcal{R}$.   Thus,  the map  $$\exp \circ -: \mathcal{R}^+ \longrightarrow \mathcal{R}_{>0}$$ acting by  $f \longmapsto e^f$ is an order isomorphism from the additive partially ordered group $\mathcal{R}^+$ of $\mathcal{R}$ onto the multiplicative partially ordered group $\mathcal{R}_{>0}$, with inverse acting by $f \longmapsto \log f$.

 If $R$ is a partially ordered  ring,  then an {\bf exponential on $R$}\index{exponential on a partially ordered ring} is an  isomorphism $E$  from the partially ordered  additive group $R^+$ onto the partially ordered group $R^* \cap R_{>0}$.    This condition means that $E$ is a group isomorphism and both $E$ and its inverse $E^{-1}$ are increasing.   Let us say that a {\bf (partially) ordered exponential ring}\index{ordered exponential ring} is a partially ordered ring $R$ equipped with an exponential on $R$.   Thus, the pair $(\mathcal{R}, \exp \circ -)$ is an exponential ring.    Note that, if $R$ is a partially ordered field, then $R^* \cap R_{>0} = R_{>0}$.   A {\bf (totally) ordered exponential field}\index{ordered exponential field}  is an ordered exponential ring that is a totally ordered field, that is, it is a totally ordered field $K$ equipped with an exponential $K \longrightarrow K_{>0}$ on $K$ \cite[p.\ 22]{kuhl}.
 
Let $(R,E)$ be an ordered exponential ring.  Then $f^\wedge a := E(aE^{-1}(f))$ behaves much like ``$f^a$'' for all $a \in R$ and all $f \in R^* \cap R_{>0}$  in that it satisfies the obvious laws of exponents:  for all $f,g \in R^* \cap R_{>0}$ and all $a,b \in R$, one has the following.
 \begin{enumerate}
\item  $f^\wedge a \in R^* \cap R_{>0}$.
 \item $f^\wedge 0 =1$.
 \item $1^\wedge a = 1$.
 \item $f^\wedge(a+b) = (f^\wedge a)(f^\wedge b)$.
  \item $f^\wedge (ab) = (f^\wedge a)^\wedge b$.
 \item $(fg)^\wedge a = (f^\wedge a)( g^\wedge a)$.
  \item $f^\wedge a  = f^a$ if $a \in \ZZ$.
  \item $(f^a)^\wedge b = f^{\wedge}  (ab)= (f^{\wedge}b )^a$ if $a \in \ZZ$.
    \item If $a < b$,  then $f^\wedge a < f^\wedge b$ if $f >1$ and $f^\wedge a >  f^\wedge b$ if $f < 1$.
  \item If $f <g$, then $f^\wedge a < g^\wedge a$ if $a > 0$, and $f^\wedge a > g^\wedge a$ if $a < 0$.
 \end{enumerate}
 It follows that, for any fixed $g \in R^* \cap R_{>0}$, the map $g^\wedge- = E( E^{-1}(g) \cdot -)$ acting by $a \longmapsto g^\wedge a$ is an exponential on $R$ with $1 \longmapsto g$.  If $F: R \longrightarrow R^* \cap R_{>0}$ is another exponential on $R$, then $$(E^{-1}\circ F)(r+s) = E^{-1}(F(r)F(s)) =  (E^{-1}\circ F)(r)+(E^{-1}\circ F)(s)$$ for all $r,s \in R$,  and therefore $E^{-1}\circ F$ is an automorphism of the ordered additive group $R^+$ of $R$.  Conversely,  if $\phi$ is any automorphism of  the ordered group $R^+$, then  $E \circ \phi$ is an exponential on $R$.  Thus,  if a partially ordered ring $R$ admits an exponential, then the set of all exponentials on $R$ is a group,  isomorphic to the group of all automorphisms of the ordered group $R^+$, and it contains (a group isomorphic to) the group $R^* \cap R_{>0} \cong R^+$ as a  subgroup, since dilation $x \longmapsto rx$ is an automorphism of the ordered group $R^+$ for some $r \in R$ if and only if $r \in R^* \cap R_{>0}$.  In particular, if $E$ is an exponential on $R$, then $E \circ (r\cdot -)$ is an exponential on $R$ for all $r \in R^* \cap R_{>0}$.  

For the  ordered exponential ring $(\mathcal{R}, \exp \circ -)$,  one has $f^\wedge a = e^{a \log f} = f^a$ for all $f \in \mathcal{R}$ and all $a \in \mathcal{R}_{>0}$.   Thus,  applying the remarks in the previous paragraph to $g = \id$, and to $r = \log$, we see in two distinct ways that the map
$$E_- = \id^- = \exp \circ (\log \cdot -): \mathcal{R}^+ \longrightarrow \mathcal{R}_{>0}$$ acting by  $$f(x) \longmapsto E_f(x) = x^{f(x)} = e^{(\log x)f(x)}$$ is a ``base independent'' exponential on $\mathcal{R}$,  with  inverse $$L_- = (E_-)^{-1} = \frac{\log \circ -}{\log}: \mathcal{R}_{>0} \longrightarrow \mathcal{R}^+,$$ which is also base independent,  acting by 
$$L_-: f(x) \longmapsto L_f(x) =  \frac{\log f(x)}{\log x}.$$

We extend $\mathcal{R}$ to the set $\mathcal{R}_{\pm \infty}$ of germs of  all functions from $\RR$ to $\overline{\RR}$,  and we extend the map $L_-$ to a map $\mathcal{R}_{\geq 0} \longrightarrow \mathcal{R}_{\pm \infty}$ by setting $L_f(x) = -\infty$ if $f(x) = 0$.  Note that  $\limsup$ and $\liminf$ define maps
$$\limsup_{x \to \infty}: \mathcal{R}_{\pm \infty} \longrightarrow \overline{\RR}$$
and
$$\liminf_{x \to \infty}: \mathcal{R}_{\pm \infty} \longrightarrow \overline{\RR},$$   
and the (upper) degree $\deg$ and lower degree $\underline{\deg}$ maps on $\mathcal{R}$ are given respectively by
 the composition
$$\mathcal{R} \overset{|-|}{\longrightarrow} \mathcal{R}_{\geq 0} \overset{L_-}{\longrightarrow} \mathcal{R}_{\pm \infty} \overset{\underset{x \to \infty}{\limsup}}{\longrightarrow} \overline{\RR}$$
and
$$\mathcal{R} \overset{|-|}{\longrightarrow} \mathcal{R}_{\geq 0} \overset{L_-}{\longrightarrow} \mathcal{R}_{\pm \infty} \overset{\underset{x \to \infty}{\liminf}}{\longrightarrow} \overline{\RR}.$$
This is our entry point into understanding the degree map from a ``universal'' perspective,  a topic that is taken up further in Section 7.4.

\section{Real asymptotic differential algebra: Hardy rings and Hardy fields}

Let $\mathcal{C}$\index[symbols]{.i  kaa@$\mathcal{C}$} be the subring of $\mathcal{R}$ consisting of the germs of continuous functions.
The group of units $\mathcal{C}^*$ of  the ring $\mathcal{C}$ is the group of all functions in $\mathcal{C}$ that are eventually positive or eventually negative, that is, that are comparable to $0$, and thus $\mathcal{C}^* = \mathcal{C}_{>0} \cup \mathcal{C}_{<0}$.

\begin{proposition}[{\cite[Section 2]{bos}}]
Let $R$ be a subring of $\mathcal{C}$ containing $\RR$.   The following conditions are equivalent.
\begin{enumerate}
\item $R$ is contained in some subfield of $\mathcal{C}$.
\item $R$ is an integral domain.
\item $R$ is contained in $\mathcal{C}^* \cup \{0\}$.
\item $R$ is contained in $\mathcal{C}_{>0} \cup \mathcal{C}_{<0} \cup \{0\}$.
\item $R$ is totally ordered with respect to the relation $\leq$ on $\mathcal{C}$.
\end{enumerate}
Suppose that these conditions hold.  Then the quotient field of $R$ is the smallest subfield of $\mathcal{C}$ containing $R$.   Suppose also that $R$ contains $\RR$.   Then, for all $f \in R$, the ring  $\RR(f)$ is the quotient field of $\RR[f]$ and is the smallest subfield of $R$ containing $\RR$ and $f$.    Moreover,  an $f \in \mathcal{C}$ belongs to some subfield of $\mathcal{C}$ if and only if $f$ is $\leq$-comparable to all $c \in \RR$, if and only if $\RR[f]$ is totally ordered by $\leq$.
\end{proposition}

\begin{corollary}
Let $K$ be a subfield of $\mathcal{C}$.   Then $$K\backslash\{0\} = K^* \subseteq  \mathcal{C}^* = \mathcal{C}_{>0} \cup \mathcal{C}_{<0}.$$
Thus, every function in $K$ is (eventually) positive,  (eventually) negative, or (eventually) $0$, whence $K$ is a totally ordered field under the ordering $\leq \ = \  \leq_\infty$.   
\end{corollary}

For any nonnegative integer $n$, let $\mathcal{C}^{(n)}$ denote the ring of germs of functions from $\RR$ to $\RR$ that are $n$-times  continuously differentiable on a neighborhood of $\infty$, and let $\mathcal{D} = \mathcal{C}^{(\infty)} =  \bigcap_{n = 1}^\infty \mathcal{C}^{(n)}\index[symbols]{.i  kb@$\mathcal{D}$}$.  The group of units of $\mathcal{D}$  is given by $\mathcal{D}^* = \mathcal{D}_{>0} \cup \mathcal{D}_{<0}$.   Note that a function in $\mathcal{D}$ need not be infinitely differentiable on a neighborhood of $\infty$.   Equivalently,  one has $\mathcal{D} = \mathcal{C}^{(\infty)} \supsetneq \mathcal{C}^\infty$, where $\mathcal{C}^\infty$ is the ring of germs of all infinitely differentiable functions on $\RR$.  Differentiation $D = \frac{d}{dx} = (-)'$ is a {\it derivation} $D: \mathcal{D} \longrightarrow \mathcal{D}$ on $\mathcal{D}$ and makes the ring $\mathcal{D}$ into a {\it differential $\RR$-algebra}.     ({\it Differential algebra},  introduced by J.\ Ritt in 1950, is the study of differential rings and their  applications to the algebraic study of differential equations.)   Let us say that a {\bf Hardy ring}\index{Hardy ring} is a subring of $\mathcal{D}$ that is closed under differentiation.  In other words, a Hardy ring is a subring $R$ of $\mathcal{D}$ such that $D(R) \subseteq R$,  i.e., such that $f' \in R$ for  all $f \in R$.  Thus, $\mathcal{D}$ is the largest Hardy ring, and $\mathcal{C}^\infty$ is also a Hardy ring.   Two obvious classes of examples of Hardy rings are $C$ and $C[x] \cong C[X]$ (generated as a subring of $\mathcal{D}$),  where $C$ is any subring of $\RR$ and  where $x = \id$ is the identity function.   Every Hardy ring $R$ is a {\it differential $C_R$-algebra}, where $C_R = \ker D|_R$ is the subring $C_R = R \cap \RR$ of $\RR$ of all constants in $R$.    Note that the intersection of a collection of Hardy rings is a Hardy ring, and therefore every subring of $\mathcal{D}$ is contained in a smallest Hardy ring.  For convenience,  hereafter we assume that all  subrings $R$ of $\mathcal{R}$ contain $\RR$,  or, equivalently,  have ring of constants $R \cap \RR = \RR$.  

A {\bf Hardy  field}\index{Hardy field} is a Hardy ring that is a field.   Equivalently, a Hardy field is a subfield $K$ of $\mathcal{D}$ containing $\RR$ such that $f' \in K$ for all $f \in K$.    Since every subfield of $\mathcal{C}$ is a totally ordered field ordered by $\leq$,  so is every Hardy field.   Examples of Hardy fields include the ordered fields $$\RR \subsetneq \RR(x) \subsetneq \RR(x^a: a \in \RR) \subsetneq  \RR(\mathfrak{L})\subsetneq \mathbb{L}$$ employed in Chapter 6.    Note that,  if $G$ is any subgroup of $\RR$ containing $\ZZ$,  then $\RR(x^a: r \in G) \cong \RR(G)$ is a Hardy field containing $\RR(x)$, and, if $G'$ is another such group, distinct from $G$, then $\RR(x^a: r \in G) \neq \RR(x^a: r \in G')$.

Let $K$ be a Hardy field.    If $f \in K$, then $f' \in K$, so that $f'$ is (eventually) positive, negative, or $0$,  so that $f$ is (eventually) (strictly) increasing,  (strictly) decreasing, or constant.   It follows that $f \preceq g$ is equivalent to $g \not \prec f$, as long as $f \neq 0$.   Indeed, if $g \neq 0$,  so that $f/g \in K$, then $f  \preceq g$ is equivalent to 
$\lim_{x \to \infty} \frac{f(x)}{g(x)} \in \RR$,  while $f \asymp g$ is equivalent to $\lim_{x \to \infty} \frac{f(x)}{g(x)} \in \RR^*$, and $f \prec g$ is equivalent to $\lim_{x \to \infty} \frac{f(x)}{g(x)} = 0$.     
Since any subfield $K$ of $\mathcal{C}$ is totally ordered with totally ordered subfield $\QQ$,  any such field $K$ is equipped with a natural valuation on $K$, as described in the previous section: the set $K^*/{\asymp}$ of all $\asymp$-equivalence classes $[f]$ of germs $f \in K^*$  is a totally ordered abelian group $\Gamma_K\index[symbols]{.i  kg@$\Gamma_K$}$ under the ordering $\leq$ on $\asymp$-equivalence classes $[f]$, where $[f] \leq [g]$  if $g \preceq f$, and the map $v: K \longrightarrow  \Gamma_K \cup \{\infty\}$ acting by $f \longmapsto v(f) = [f]$, where $0 \longmapsto \infty$, is a valuation on the field $K$.  The valuation ring $$\mathcal{O}_K = K_v =  \{f \in K: v(f) \geq 0\} = \{f \in K: f \preceq 1\}\index[symbols]{.i  kh@$\mathcal{O}_K$}$$ of the valued field $K$ is the ring of all bounded functions $f$ in $K$,  or, equivalently, all functions $f$ in $K$ such that $\lim_{x \to \infty} f(x) \in \RR$.    The unique maximal ideal $\mm_K$ of the valuation ring $\mathcal{O}_K$ is the ideal 
$$\mm_K = \mm_v =  \{f \in K: v(f) > 0\} = \{f \in K: f \prec 1\}\index[symbols]{.i  ki@$\mm_K$}$$ 
of all functions $f$ in $K$ with  $\lim_{x \to \infty} f(x) = 0$.  Thus, one has a canonical isomorphism 
$$\mathcal{O}_K \cong \RR\oplus \mm_K$$
of abelian groups,   and the map to the residue field $\mathcal{O}_K/\mm_K$ is the projection
$$\lim_{x \to \infty} : \mathcal{O}_K \longrightarrow \RR \cong \mathcal{O}_K/\mm_K$$
onto $\RR$, which has a section given by the inclusion $ \RR\longrightarrow \mathcal{O}_K$, and which extends to the map
$$\lim_{x \to \infty} : K \longrightarrow \overline{\RR}$$
with $\mathcal{O}_K$ equal to the preimage of $\RR \subsetneq \overline{\RR}$.
Moreover, the group of units of the valuation domain $\mathcal{O}_K$ is the group
$$\mathcal{O}_K^* = \mathcal{O}_K\backslash \mm_K =  \{f \in K: v(f) =  0\} = \{f \in K: f \asymp 1\}$$
of all functions $f$ in $K$ with $\lim_{x \to \infty} f(x) \in \RR^*$.  Note, then, that $v(f)= v(g)$ if and only if $f/g \in \mathcal{O}_K^*$, if and only if $(f) = (g)$ as principal fractional ideals in $\mathcal{O}_K$.  Thus, $\Gamma_K$ is order-isomorphic to the group of all nonzero principal fractional ideals of $\mathcal{O}_K$ ordered by $\supseteq$.   It follows that, as groups, one has $\Gamma_K \cong K^*/\mathcal{O}_K^*$.  Moreover, the set $K_{>0}$ of all positive elements of $K$ is a subgroup of $K^*$, and one has $K^* \cong K_{>0} \oplus \{1,-1\}$.

We say that a  subring $R$ of $\mathcal{R}$ is {\bf logexponentially closed}\index{logexponentially closed} if $e^f \in R$ for all $f \in R$ and and $\log |f| \in R$ for all $f \in R^*$, that is, if the exponential $\exp \circ -$ on $\mathcal{R}$ restricts/corestricts to an exponential on $R$.   If  a  subring $R$ of $\mathcal{R}$ is  logexponentially closed, then $R$ inherits from $\mathcal{R}$ the structure of an ordered exponential ring $(R,  \exp \circ ( a \cdot -))$ for any $a \in R$.   In this case,  $R$ also has the structure of a ``base independent''  ordered exponential ring with exponential $E_- = \id^- =  \exp\circ (\log \cdot -)$ if and only if  $R \supseteq \RR(x)$,  i.e., if and only if $\id \in R$.     Note that the  Hardy field $\mathbb{L}$ of all logarithmico-exponential functions is the smallest logexponentially closed subfield of $\mathcal{R}$ containing $\RR(x)$.     Similarly, the Hardy field $\RR(\mathfrak{L})$ is the smallest subfield $K$ of $\mathcal{R}$ containing $\RR(x^a : a \in \RR)$ such that $f \circ \log \in K$ for all $f \in K$, and the Hardy field $\RR(x^a : a \in \RR)$ is the smallest subfield $K$ of $\mathcal{R}$ containing $\RR(x)$ such that $f \circ x^a \in K$ for all $f \in K$ and all $a>0$.

The largest Hardy ring $\mathcal{D}$ has several nice closure properties: it is  closed under differentiation,   closed under antidifferentiation,  and logexponentially closed. The rings $\mathcal{R}$ and $\mathcal{C}$ are also logexponentially closed, and $\mathcal{C}$ is also closed under antidifferentiation.   The operation $$\operatorname{Dlog}: \mathcal{D}^* \longrightarrow \mathcal{D}^+\index[symbols]{.i  ke@$\operatorname{Dlog}$}$$
of {\bf logarithmic differentiation},\index{logarithmic differentiation} acting by $$\operatorname{Dlog}: f \longmapsto \frac{f'}{f} = (\log |f|)',$$ is a group homomorphism.   If $R$ is a Hardy ring,  then $\operatorname{Dlog}$ restricts to a group homomorphism  $R^* \longrightarrow R^+$.   

Since the intersection of a collection of Hardy subfields of a Hardy field is a Hardy field, every subset of a Hardy field $K$ is contained in a smallest Hardy subfield of $K$.   However, if $K$ and $L$ are Hardy fields, then there need not be a Hardy field containing both $K$ and $L$.    A {\bf maximal Hardy field}\index{maximal Hardy field} is a Hardy field that is not a proper subfield of any other Hardy field.  By Zorn's lemma, every Hardy field is contained in a maximal Hardy field.  The intersection $\mathbb{H}\index[symbols]{.i  kc@$\mathbb{H}$}$ of all maximal Hardy fields is a Hardy field (denoted by $E$ in \cite{bos}).   As we show later,  one has $\mathbb{L} \subsetneq \mathbb{H}$.   If $K$ is a Hardy field,  then,  following \cite{bos}, we let $\operatorname{E}(K)\index[symbols]{.i  kca@$\operatorname{E}(K)$}$ denote the intersection of all maximal Hardy fields containing $K$.  A Hardy field $K$ is {\bf perfect}\index{perfect Hardy field} if it is the intersection of some collection of maximal Hardy fields, or, equivalently, if $\operatorname{E}(K) = K$ \cite{bos2}.  Thus, $\mathbb{H} = \operatorname{E}(\RR)$ is the smallest perfect Hardy field.  Moreover, if $K$ is a Hardy field, then $\operatorname{E}(K)$ is the smallest perfect Hardy field containing $K$. 

A function $f \in \mathcal{D}$ is {\bf Hardian}\index{Hardian function} if there exists a Hardy field containing $f$.   
If $R$ is a Hardy ring and $f \in \mathcal{D}$, then the $R$-subalgebra $R[f,f',f'',\ldots]$ of $\mathcal{D}$ generated by $f$ and all of its derivatives is the smallest Hardy ring  containing $f$, and there is a unique (surjective) ring homomorphism $$D_{f,R}: R[X_0,X_1,X_2,\ldots] \longrightarrow R[f,f',f'',\ldots]\index[symbols]{.i  ki1@$D_{f,R}$}$$ sending $X_n$ to $f^{(n)}$ for all $n$.   If $K$ is a Hardy field, then a (Hardian) function $f \in \mathcal{D}$ is {\bf Hardy adjoinable to $K$}\index{Hardy adjoinable} if there exists a Hardy field containing both $K$ and $f$.    Thus, $f \in \mathcal{D}$ is Hardian if and only if $f$ is Hardy adjoinable to $\RR$, if and only if $f$ is Hardy adjoinable to the Hardy field $\mathbb{H}$.   If $f$ is Hardy adjoinable to a Hardy field $K$, then there is a smallest Hardy field containing $K$ and $f$, denoted $K\{f\}$.
The following result is clear.

\begin{proposition}[{\cite{bos}}]
Let $f \in \mathcal{D}$, and let $K$ be a Hardy field.  The following conditions are equivalent.
\begin{enumerate}
\item $f$ is Hardy adjoinable to $K$.
\item The image of the map $D_{f,K}$ is contained in $\mathcal{D}^* \cup \{0\}$.  
\item For all $p(X_0, X_1,X_2,\ldots) \in K[X_0,X_1,X_2,\ldots]$, the function $p(f,f',f'',\ldots)$ is either eventually positive, eventually negative, or eventually $0$, i.e.,  $p(f,f',f'',\ldots)$ lies in $\mathcal{D}^* \cup \{0\}$. 
\item The ring $K[f,f',f'',\ldots]$ is an integral domain.
\item The kernel $D_{f,K}^{-1}(0)$ of $D_{f,K}$ is prime, as an ideal of $K[X_0,X_1,X_2,\ldots]$.   
\end{enumerate}
When these conditions hold,  one has
 \begin{align*}
K\{f\}  & = K(f,f',f'',\ldots) \\
& = \{p(f,f',f'',\ldots)/q(f,f',f'',\ldots): p,q \in K[X_0,X_1,X_2,\ldots], \, q(f,f',f'',\ldots) \neq 0\}
\end{align*}
where $K(f,f',f'',\ldots) $ (generated as a subfield of $ \mathcal{D}$) is the quotient field of the integral domain  $K[f,f',f'',\ldots]$, and it is isomorphic to the quotient field of $K[X_0,X_1,X_2,\ldots]/D_{f,K}^{-1}(0)$, i.e., the residue field of the domain $K[X_0,X_1,X_2,\ldots]$ at the prime ideal $D_{f,K}^{-1}(0)$.    
\end{proposition}

We say that a subset $S$ of $ \mathcal{D}$  is {\bf Hardy adjoinable to $K$}\index{Hardy adjoinable} if   there exists a Hardy field containing both $K$ and $S$, and we let $K\{S\}\index[symbols]{.i  kl@$K\{S\}$}$ denote the smallest Hardy field containing $K$ and $S$.   If $K$ and $L$ are Hardy fields, then $L$ is Hardy adjoinable to $K$ if and only if $K$ is Hardy adjoinable to $L$, and when both conditions hold one has $K\{L\} = L\{K\}$. 

\begin{example}  Let $K$ be a Hardy field.   Since $x$, $\exp x$, and  $\log x$ are in $\mathbb{H}$,  the set  $\{x, \exp x, \log x\}$ is Hardy adjoinable to $K$, and one has $K\{x\} = K(x)$, $K\{\exp x\} = K(\exp x)$,  $K\{\log x \} = K(x,\log x)$, and $$K\{\exp x, \log x\} = (K\{\exp x\})\{\log x\} = (K(\exp x))\{\log x\} = K(\exp x,x,\log x).$$  
\end{example}

As we have noted, the Hardy field $\mathbb{H}$ is equal to the set of all functions in $ \mathcal{D}$ that are Hardy adjoinable to every Hardy field.   Thus, we call elements of $\mathbb{H}$ {\bf universally Hardian}, \index{universally Hardian function} and  $\mathbb{H}$ is the Hardy field of all universally Hardian functions (or, rather, germs at $\infty$).       If $K$ is any Hardy field, then an obvious transfinite inductive argument shows that $\mathbb{H}$ is Hardy-ajoinable to $K$, and therefore $K$ is Hardy adjoinable to $\mathbb{H}$.   Thus, every Hardy field is Hardy adjoinable to $\mathbb{H}$.  (Indeed, if $K$ is a Hardy field, then $K$ is contained in some maximal Hardy field $M$, and then $M$ contains both $K$ and $\mathbb{H}$.) Moreover,  $\mathbb{H}$ is the largest Hardy field with this property:  if $K$ is a Hardy field such that every Hardy field is Hardy adjoinable to $K$,  then $K$ is Hardy adjoinable to every Hardy field, so that all elements of $K$ are Hardy adjoinable to every Hardy field,  and thus $K \subseteq \mathbb{H}$.  Thus,  the Hardy field $ \mathbb{H}$ has the universal property that it is the largest Hardy field $K$ such that every Hardy field is Hardy adjoinable to $K$.  Equivalently,  the Hardy field $ \mathbb{H}$ is the largest Hardy field that is Hardy adjoinable to every Hardy field.  Thus, we call $ \mathbb{H}$ the {\bf universally Hardy adjoinable Hardy field}.\index{universally Hardy adjoinable Hardy field $\mathbb{H}$}

Maximal Hardy fields are precisely the Hardy fields $M$ such that $f \in M$ for all Hardian $f \in \mathcal{D}$ that are Hardy adjoinable to $M$.   Moreover, for any Hardy field $K$, the intersection $\operatorname{E}(K)$ of all maximal Hardy fields containing $K$ is equal to the set of all $f \in \mathcal{D}$ that are Hardy adjoinable to every Hardy field containing $K$.

Let $K$ be a subfield of $\mathcal{D}$.  The algebraic closure $K^o\index[symbols]{.i  kn@$K^o$}$ of $K$ in $\mathcal{D}$ is a {\it real closed} subfield of $\mathcal{D}$,  and in fact it is isomorphic to the {\it real closure of $K$},
 and one has $(K^o)^o = K^o$.   Thus, $K$ is real closed if and only if $K^o = K$.  Moreover, if $K$ is a Hardy field,  then $K^o$ is a  (real closed) Hardy field.   It follows that $K$ is real closed if $K$ is a perfect Hardy field.  For proofs of these facts, see \cite[Section 3]{bos}.

\begin{theorem}[{\cite[Theorem 2]{rosenl}}]\label{rosenl}
Let $K$ be a Hardy field,  let $p(X),q(X) \in K[X]$,  and let $y$ be the (germ of) an eventually differentiable function such that $q(y)$ is eventually nonzero and such that $y' = p(y)/q(y)$ on some neighborhood of $\infty$.   Then $K(y)$ is a Hardy field.  Thus, $y$ is Hardy adjoinable to $K$ and $K\{y\} = K(y)$.  Consequently,  for all $f \in K$, one has the following.
\begin{enumerate}
\item $K(e^f)$ is a Hardy field, since $y = e^f$ satisfies $y' = fy$.
\item $K(\log |f|)$ is a Hardy field if $f \neq 0$, since $y = \log |f|$ satisfies $y' = f'/f$ and $f'/f \in K$.
\item $K(\int f(x) \, dx)$ is a Hardy field for any eventual antiderivative $y = \int f(x) \, dx$ of $f$, since $y$ satisfies $y' = f$.
\end{enumerate}
\end{theorem}

\begin{corollary}\label{Hsim}
Let $K$ be a perfect Hardy field.  Let $p(X),q(X) \in K[X]$,  and let $y$ be the (germ of) an eventually differentiable function such that $q(y)$ is eventually nonzero and such that $y' = p(y)/q(y)$ on some neighborhood of $\infty$.   Then $y \in K$.  Thus,  if $f \in \mathcal{D}$ is Hardian (resp.,  lies in $K$),  then $e^f$, any eventual antiderivative of $f$, and $\log |f|$ if $f$ is eventually nonzero,  are also Hardian (resp., lie in $K$).
\end{corollary}

\begin{corollary}
Any Hardy field  is contained in some logexponentially closed Hardy field.
\end{corollary}

Note that, if $f \in \mathcal{D}$ is Hardian and unbounded, then the compositional inverse  $f^{-1} \in \mathcal{D}$ of $f$ exists and is also Hardian \cite[Corollary 6.5]{bos}.   
 
\begin{example}
Let $W$ denote the Lambert $W$ function.   One has
$$W'(x) =  \frac{W(x)}{x(1+W(x))}, \quad \forall x \in (-1/e,0)\cup(0,\infty).$$
It follows that $W \in \RR\{W(x)\} = \RR(x,W(x)) \subseteq \mathbb{H}$, and thus $e^{x/W(x)} \in \mathbb{H}$.   The function $e^{x/W(x)}$ is the compositional inverse of $\log x \log \log x \in \mathfrak{L}$.   It is known (and originally conjectured by Hardy in \cite{har4}) that there is no $f \in \mathbb{L}$ such that $f (x) \sim e^{x/W(x)}$  \cite{dmm} \cite[Theorem 6.49]{kuhl}.  However, one has
$$e^{x/\log x} \prec e^{x/W(x)} \prec (e^{x/\log x})^a$$
for all $a > 1$, and thus $\dege e^{x/W(x)} = \dege e^{x/\log x} = (\infty,1,-1,0,0,0,\ldots)$.  Note also that $W^{-1}(x) = xe^x$ is  in $\RR(x,\exp x)$.  Moreover,  one has $\int_0^x W(t) \, dt = x\left(W(x)-1 + \frac{1}{W(x)}\right)$, so any Hardy field containing $W(x)$ contains $x$ and $\int_0^x W(t) \, dt$.  In fact, any subfield of $\mathcal{D}$ containing $x$ and $W(x)$ contains the $n$-fold derivative and every $n$-fold antiderivative of $W(x)$ for all $n$, since these lie in $\RR(x,W(x))$.
\end{example}

\begin{remark}[The asymptotic equivalence class of $\log$]
Let $K$ be a Hardy field,  let $f \in K_{>0}$ with $f \succ 1$.   \cite[Proposition 6]{rosenl1} states that there exists a $g \in K_{>0}$ with $g \sim \log f$  if and only if there exists a $g \in K_{>0}$ with $g \succ 1$ and $f \succ g^n$ for every positive integer $n$.   This provides a ``universal property'' of the  asymptotic equivalence class of $\log$; more precisely, it provides a characterization of the $\sim$-equivalence class of $\log f$ for any $f$ in a Hardy field $K$ with $f > 0$ and $f \succ 1$.
\end{remark}

Let $R$ be a subring of $\mathcal{D}$.  A function $f \in \mathcal{D}$ is {\bf differentially algebraic over $R$}\index{differentially algebraic} if there exists a nonzero polynomial $p(X_0,X_1,X_2,\ldots) \in R[X_0,X_1,X_2,\ldots]$ such that $p(f,f',f'',\ldots) = 0$, or, equivalently,  if the kernel $D_{f,R}^{-1}(0)$ of the ring homomorphism $D_{f,R}$ is nontrivial.    A function $f \in \mathcal{D}$ is {\bf differentially algebraic} if it is differentially algebraic over $\RR$.    Suppose that $f \in \mathcal{D}$ is differentially algebraic over $R$.  Then for some nonnegative integer $d$ there exists a nonzero polynomial $p(X_0,X_1,\ldots,X_d) \in R[X_0,X_1,\ldots,X_d]$ in the kernel $D_{f,R}^{-1}(0)$ of minimal degree $m > 0$ in $X_d$, and thus the kernel $D_{f,R}^{-1}(0) \cap R[X_0,X_1,X_2,\ldots,X_d]$ of the restricted map $\RR[X_0,X_1,X_2,\ldots,X_d] \longrightarrow S = R[f,f',f'',\ldots,f^{(d)}]$ is a nonzero ideal  of $R[X_0,X_1,X_2,\ldots,X_d]$.   The smallest $d = d_{f,R}\index[symbols]{.i  kp@$d_{f,R}$}$ possible is called the {\bf differential degree of $f$ over $R$}\index{differential degree $d_{f,R}$ of $f$ over $R$}.  We set $d_{g,R} = \infty$ for  any $g \in \mathcal{D}$ that is not differentiably algebraic over $R$.  Differentiating the equation $p = D_{f,R}(p) = p(f,f',f'',\ldots, f^{(d)}) = 0$, we obtain
  $$D(p) = \sum_{k = 0}^d p_k f^{(k+1)} = 0,$$
  where  $p_k = \frac{\partial p}{\partial f^{(k)}} = D_{f, R} \frac{\partial p}{\partial X_k}  \in S$ for all $k \leq d$.
Note then that  $p_d \neq 0$ since the polynomial $p_d$ in $f,f',f'',\ldots,f^{(d)}$ has degree $m-1 < m$ in $f^{(d)}$.   Suppose that $p_d \in \mathcal{D}^*$ (which holds, for example, if $R$ is a Hardy field and $f$ is Hardy adjoinable to $R$).   Then
  $$f^{(d+1)}  = -\frac{1}{p_d}\sum_{k = 0}^{d-1} p_k f^{(k+1)} \in R[f,f',f'',\ldots,f^{(d)}][1/p_d]$$
  and thus
  $$p_k' \in R[f,f',f'',\ldots,f^{(d+1)}] \subseteq R[f,f',f'',\ldots,f^{(d)}][1/p_d]$$
  for all $k \leq d$.  By induction,  then, one has
    $$f^{(k)}  \in R[f,f',f'',\ldots,f^{(d)}][1/p_d]$$
    for all $k$, i.e., 
    $$R[f,f',f'',\ldots] \subseteq R[f,f',f'',\ldots,f^{(d)}][1/p_d].$$  Conversely, if  there exists a $p \in R[f,f',f'',\ldots, f^{(d)}] \cap \mathcal{D}^*$ such that  $$R[f,f',f'',\ldots] \subseteq R[f,f',f'',\ldots,f^{(d)}][1/p]$$ for some $d$,  then   $p^n f^{(d+1)} -q = 0$ for some $q \in R[f,f',f'',\ldots]$ and some nonnegative integer $n$,  so  there exist polynomials
 $P$ and $Q$ in $R[X_0,X_1,\ldots, X_d]$ with $P \neq 0$ such that $0 \neq P X_{d+1}-Q \in R[X_0,X_1,X_2,\ldots]$ lies in the kernel of $D_{f,R}$, and therefore $f$ is differentially algebraic over $R$.  Thus, we have the following.
 
\begin{proposition}
Let $R$ be a subring of $\mathcal{D}$,  let $f \in \mathcal{D}$, and let $d$ be a nonnegative integer.   The following conditions are equivalent.
\begin{enumerate}
\item There  exists  a $p \in R[f,f',f'',\ldots,f^{(d)}] \cap \mathcal{D}^*$ such that  $$R[f,f',f'',\ldots] \subseteq R[f,f',f'',\ldots,f^{(d)}][1/p].$$
\item There  exist $P,Q \in R[X_0,X_1,\ldots, X_d]$ with $P \neq 0$ such that $P X_{d+1}-Q$ lies in  the kernel of $D_{f,R}$ and $D_{f,R}(P) \in \mathcal{D}^*$.
\end{enumerate}
Moreover, if the conditions above hold, then  $f$ is differentially algebraic over $R$.
\end{proposition}

\begin{proposition}\label{diffeq}
Let $K$ be a Hardy field,  let $f \in \mathcal{D}$ be Hardy adjoinable to $K$, and let $d$ be a nonnegative integer.   The following conditions are equivalent.
\begin{enumerate}
\item $f$ is differentially algebraic over $K$  with $d_{f,K} \leq d$.
\item There  exists a nonzero $p \in K[f,f',f'',\ldots,f^{(d)}]$ such that  $$K[f,f',f'',\ldots] \subseteq K[f,f',f'',\ldots,f^{(d)}][1/p].$$
\item There  exist  $P,Q  \in K[X_0,X_1,\ldots, X_d]$ with $P \neq 0$ such that $P X_{d+1}-Q$ lies in  the kernel of $D_{f,R}$ but $P$ does not.
\item There exists an $F \in K(X_0,X_1,\ldots, X_d)$ such that $f^{(d+1)} \in F(f,f',f'',\ldots,f^{(d)})$.
\item $K\{f\} \subseteq K(f,f',f'',\ldots,f^{(d)})$.
\item $K\{f\} = K(f,f',f'',\ldots,f^{(d)})$.
\end{enumerate}
\end{proposition}
 
\begin{corollary}[{\cite[Section 14]{bos2}}]\label{diffeqcor}
Let $f \in \mathcal{D}$ be Hardian, and let $d$ be a nonnegative integer.   The following conditions are equivalent.
\begin{enumerate}
\item $f$ is differentially algebraic of differential degree at most $d$ over $\RR$.
\item There  exists a nonzero $p \in \RR[f,f',f'',\ldots,f^{(d)}]$ such that  $$\RR[f,f',f'',\ldots] \subseteq \RR[f,f',f'',\ldots,f^{(d)}][1/p].$$
\item There  exist  $P,Q  \in \RR[X_0,X_1,\ldots, X_d]$ with $P \neq 0$ such that $P X_{d+1}-Q$ lies in  the kernel of $D_{f,R}$ but $P$ does not.
\item There exists  an $F \in \RR(X_0,X_1,\ldots, X_d)$ such that $f^{(d+1)} \in F(f,f',f'',\ldots,f^{(d)})$.
\item $\RR\{f\} \subseteq \RR(f,f',f'',\ldots,f^{(d)})$.
\item $\RR\{f\} = \RR(f,f',f'',\ldots,f^{(d)})$.
\end{enumerate}
\end{corollary}

 We say that a subring $S$ of $\mathcal{D}$ containing a ring $R$ is {\bf differentially algebraic over $R$}\index{differentially algebraic} if every element of $S$ is differentially algebraic over $R$.   We say that a subring $R$ of $\mathcal{D}$ is {\bf differentially algebraic} if every element of $R$ is differentially algebraic, i.e., if $R$ is differentially algebraic over $\RR$.   A subring $R$ of $\mathcal{D}$ differentially algebraic over a Hardy field $K$ (resp., over $\RR$) if and only if  its elements satisfy the conditions of Prospotion \ref{diffeq} (resp., Corollary \ref{diffeqcor}).   By  \cite[Theorem 14.3]{bos2}, the Hardy field $\mathbb{H}$ is differentially algebraic.     Moreover,  by Theorem \ref{rosenl}, if $f$ is Hardian and $d_{f,\RR} \leq 1$, then $f \in \mathbb{H}$.  \cite[Conjecture 2]{bos} states that, if $f\in \mathcal{D}$ is Hardian and differentially algebraic, then $f \in \mathbb{H}$.  If this conjecture is true, then $\mathbb{H}$ is precisely the set of all differentially algebraic Hardian functions, and thus $\mathbb{H}$ is the largest differentially algebraic Hardy field.    Note that the Hardian elements of $\mathcal{D}$ are those elements $f$ of $\mathcal{D}$ such that $p(f,f',f'',\ldots) \in \mathcal{D}^* \cup\{0\}$ for every  $p \in \RR[X_0,X_1,X_2,\ldots]$, and the differentially algebraic elements of $\mathcal{D}$ are those elements $f $ of $\mathcal{D}$ such that  $p(f,f',f'',\ldots) \in \{0\}$ for some  nonzero  $p \in \RR[X_0,X_1,X_2,\ldots]$.   Thus, \cite[Conjecture 2]{bos} states that $\mathbb{H}$ is precisely the set of all $f \in \mathcal{D}$ for which both conditions are true.    For any  subring $R$  of $\mathcal{R}$,  we say that $R$ is {\bf bounded by some continuous function} if there exists a $u \in \mathcal{C}$ such that $f \leq u$ for all $f \in R$.   By  \cite[Theorem 14.4]{bos2}, if $K$ is a Hardy field bounded by some continuous function, then $\operatorname{E}(K)$ is differentially algebraic over $K$.  It is also conjectured that $\operatorname{E}(K)$ is is differentially algebraic over $K$ for any Hardy field $K$.  
 
Any Hardian function that is differentially algebraic is analytic on some punctured neighborhood of $\infty$.   Thus $\mathbb{H}$ is a subfield of the ring $\mathcal{A}$ of germs of functions that are analytic on some punctured neighborhood of $\infty$.   Furthermore, any function that is analytic at $\infty$ is Hardian.  However,  if a function $r(x) = \sum_{n = 0}^\infty a_n/x^n$ analytic at $\infty$ is differentially algebraic (or in $\mathbb{H}$),  then the (real) coefficients $a_0,a_1,a_2,\ldots$ must be algebraically dependent over $\QQ$.  See \cite{bos} and \cite{bos2} for proofs of these facts.

The following is a  corollary of Proposition \ref{diffeq}.

\begin{corollary}[{cf.\  \cite[Theorem 14.8]{bos2}}]\label{diffalgprop}
Let $K$ be a Hardy field and $f \in \mathcal{D}$, and let $\aaa$ be the kernel of the  ring homomorphism $D_{f,R}: R[X_0,X_1,X_2,\ldots] \longrightarrow R[f,f',f'',\ldots]$.   Then $f$ is Hardy adjoinable to $K$ if and only if $\aaa$ is a prime ideal, and $f$ is differentially algebraic over $K$ if and only if $\aaa$ is nonzero.   Suppose that $f$ is Hardy adjoinable to $K$, so that  $\aaa = \ppp$ is a prime ideal.   
\begin{enumerate} 
\item The ring $K\{f\}$ is isomorphic to the quotient field of $K[X_0,X_1,X_2,\ldots]/\ppp$, that is,  $K\{r\}$ is isomorphic to the residue field of the domain $K[X_0,X_1,X_2,\ldots]$ at $\ppp$.  
\item  If $f$ is not differentially algebraic over $K$, then  $\ppp = 0$ and thus $K\{f\}$ is isomorphic to the field $K(X_0,X_1,X_2,\ldots)$ and  has infinite transcencence degree over $K$, with transcendence basis $f, f',f'',\ldots$.  
\item Suppose that $f$ is differentially algebraic over $K$ with $d = d_{f,K}$.   Then 
 $d$ is least such that $\ppp_d = \ppp \cap K[X_0,X_1,X_2,\ldots,X_d]$ is nonzero, and   $K\{f\} = K(f,f',f'',\ldots, f^{(d)})$ is isomorphic to the residue field  of $K[X_0,X_1,X_2,\ldots,X_d]$ at $\ppp_d$.
Moreover,  $f,f',f'',\ldots,f^{(d-1)}$ are algebraically independent over $K$, and one has $f^{(d)} \in K(f,f',f'',\ldots,f^{(d-1)})^o$, so that the extension $K (f,f',f'',\ldots,f^{(d-1)}) \subseteq K(f,f',f'',\ldots,f^{(d)}) = K\{f\}$ is algebraic.  Therefore $f,f',f'',\ldots,f^{(d-1)}$ is a transcendence basis of $K\{f\}$.  
 \end{enumerate}
 In particular, $K\{f\}$ has transcendence degree $d = d_{f,K}$ over $K$.   
 \end{corollary}
 
Corollary \ref{diffalgprop} stands in perfect analogy with algebraic extensions, where algebraic degree is replaced by transcendence degree and/or differential degree, where $f^n$ is replaced by $f^{(n)}$,  where $K(f)$ is replaced by $K\{f\}$, and where vector space bases of $K(f)$ over $K$ are replaced by transcendence bases of $K\{f\}$ over $K$.
 
 Next, we prove some useful properties of Hardian functions.

\begin{proposition}\label{hardianexactlog}
Let $r \in \mathcal{D}$ be Hardian (resp., in $\mathbb{H}$, in $\mathbb{L}$) with $r$ eventually positive.    For all nonnegative integers $k$, one has
$$r_{(k+1)}(x) =   \left.
 \begin{cases}
    r_{(k)}(e^x) e^{-(\deg r_{(k)}) x}& \text{if } \deg r_{(k)} \neq \pm \infty \\
    \log r_{(k)}(x)  & \text{if } \deg r_{(k)} = \infty \\
 \displaystyle   -\frac{1}{\log r_{(k)}(x)} & \text{if } \deg r_{(k)} =- \infty,
 \end{cases}
\right.$$ and  $r_{(k)}$ is Hardian  (resp., in $\mathbb{H}$, in $\mathbb{L}$)  and eventually positive.  Moreover,  $r$ has exact logexponential degree,  with
$$\dege_k r = \lim_{x \to \infty} \frac{x r_{(k)}'(x)}{r_{(k)}(x)} \in \overline{\RR}$$
for all $k$.
\end{proposition}

\begin{proof}
It is clear that, if $K$ is a Hardy field, then so is $K(x) \circ \exp = \{f  \circ \exp: f \in K(x)\}$.   Moreover,  by \cite[Theorem 5.3]{bos},  if $f$ is Hardian and nonzero, then  $\log |f|$ is Hardian.  It follows that  $r_{(k)}$ is Hardian for all $k$.   Moreover, if $r \in \mathbb{H}$ (resp., $r \in \mathbb{L}$),  then the closure properties of $\mathbb{H}$ and $\mathbb{L}$ imply that  $r_{(k)} \in \mathbb{H}$ (resp., $r_{(k)} \in \mathbb{L}$) for all $k$.  If $f$ is Hardian and nonzero, then $ \frac{x f'(x)}{f(x)}$ is Hardian and therefore $d = \lim_{x \to \infty} \frac{x f'(x)}{f(x)} \in \overline{\RR}$ exists in $\RR$ or is $\pm \infty$, and thus $f$ has exact degree $d$, by Proposition \ref{firstprop0}(16).  Thus, by induction, $r$ has exact logexponential degree.  
\end{proof}

By the proposition and Corollary \ref{regvarcor}, we have the following.

\begin{corollary}
A Hardian function $r$ is regularly varying if and only if it is of finite degree, and, if those equivalent conditions hold, then $r$ is regularly varying of index $$\overline{\underline{\deg}} \, r = \lim_{x \to \infty}\frac{xr'(x)}{r(x)}.$$
\end{corollary}

The following proposition summarizes the relationships between  several types of functions of finite degree that we have considered.

\begin{proposition} For any $d \in \RR$, and for any function $r$ that is defined in a neighborhood of $\infty$, each of the following conditions implies the next.
\begin{enumerate}
\item $r$ is logarithmico-exponential of  degree $d$.
\item $r$ is Hardian  of degree $d$.
\item $r$ is eventually nonzero and differentiable, and  $\lim_{x \to \infty} \frac{xr'(x)}{r(x)} = d$.
\item $r$ is regularly varying of index $d$.
\item $r$ has exact degree $d$.
\end{enumerate}
\end{proposition}
 
We also note the following.
 
 \begin{proposition}[{\cite[Lemma 12.5]{bos2}}]
For any Hardian function $f$ of degree $\pm \infty$, one has
 $$\log |f| \sim \log |f'|$$
 and $$ff'' \sim (f')^2.$$
 \end{proposition}
 
 \begin{corollary}
For any Hardian function $f$ of degree $\pm \infty$, one has $$\dege f' = \dege f.$$
 \end{corollary}

 The following theorem generalizes Hardy's theorem \cite[Theorems 8a--8c]{har4} on the asymptotics of the antiderivative of a logarithmico-exponential function  to the Hardian functions.
 
\begin{theorem}[{cf.\  \cite[Theorems 8a--8c]{har4}}]\label{hardintth}
Let $f$ be a nonzero Hardian function that is defined and continuous on $[N,\infty)$, where $N \in \RR$.    Let  $F_1(x) = \int_N^x f(t) \, dt$ and $F_2(x) = \int_x^\infty f(t) \, dt$ on $[N,\infty)$ according to whether $\int_N^\infty f(t) \, dt$ diverges or converges.  Then one has the following.
\begin{enumerate}
\item If $\deg f = \infty$, then
$$F_1 \sim \frac{f^2}{f'} = \frac{f}{\operatorname{Dlog} f}.$$
\item  If $d = \deg f \in (-1,\infty)$,  then
$$F_1 \sim \frac{x f}{d+1}.$$
\item If $d = \deg f \in (-\infty,-1)$, then
$$F_2 \sim - \frac{x f}{d+1}.$$ 
\item If $\deg f = -\infty$, then
$$F_2 \sim - \frac{f^2}{f'} = -\frac{f}{\operatorname{Dlog} f}.$$
\end{enumerate}
More generally,  suppose that $\dege f \neq (-1,-1,-1,\ldots)$, and let $N$ be the smallest nonnegative integer such that $\dege_N f \neq -1$.  Let $$g =   f_{(N)}(\log^{\circ N}x) = (x \log x \log^{\circ 2}x \cdots \log^{\circ (N-1)}x) f = \frac{f}{\frac{d}{dx} \log^{\circ N}x}.$$
 Then one has the following.
\begin{enumerate}
\item[(5)] Suppose that $\dege_N f = \infty$.  Then
$$F_1 \sim  \frac{f}{\operatorname{Dlog} g}.$$
Moreover, one has 
$$g \prec F_1 \prec g \log^{\circ N} x$$
if $N > 0$, while
$$\dege F_1 = \dege f$$ if $N = 0$.
\item[(6)]  If $d = \dege_N f \in (-1,\infty)$,  then
$$F_1 \sim \frac{g}{d+1}\log^{\circ N}x.$$
\item[(7)] If $d = \dege_N f \in (-\infty,-1)$, then
$$F_2 \sim - \frac{g}{d+1}\log^{\circ N}x .$$ 
\item[(8)] Suppose that $\dege_N f = -\infty$.  Then
$$F_2 \sim - \frac{f}{\operatorname{Dlog} g}.$$
Moreover,  one has
$$g \prec F_2 \prec g \log^{\circ N}x$$
if $N > 0$, while 
$$\dege F_2 = \dege f$$ if $N = 0$.
\end{enumerate}
Consequently,  one has $$\dege F_i = \dege f \oplus (1,1,1,\ldots,1,0,0,0,\ldots) \neq (0,0,0,\ldots)$$
for $i = 1$ in cases (5) and (6) and for $i = 2$ in cases (7) and (8),
where $$(1,1,1,\ldots,1,0,0,0,\ldots) = \dege\frac{1}{ \frac{d}{dx} \log^{\circ (N+1)}x}$$
denotes the vector with $N+1$ $1$s followed by a  tail of $0$s.
\end{theorem}

\begin{proof}
The proof is a straightforward generalization of the proof of \cite[Theorems 8a--8c]{har4}, which establishes the theorem for all $f \in \mathbb{L}$.
\end{proof}

\begin{corollary}\label{karamcor}
Let $f$ be a Hardian function with $\deg f  = 0$.
\begin{enumerate}
\item If $d = \dege_1 f  \in (-\infty,-1)$, then
$$\int_{x}^\infty \frac{f(t)}{t} \, dt  \sim -\frac{f(x) \log x}{d+1} \ (x \to \infty),$$
and therefore
$$ \dege \int_{x}^\infty \frac{f(t)}{t} \, dt  = \dege f + (0,1,0,0,0,\ldots).$$
\item If $\dege_1 f  = -\infty$, then
$$f(x) \prec \int_{x}^\infty \frac{f(t)}{t} \, dt  \prec f(x)\log x,$$
and therefore
$$\dege \int_{x}^\infty \frac{f(t)}{t} \, dt = \dege f = \dege f + (0,1,0,0,0,\ldots).$$
\end{enumerate}
\end{corollary}

\begin{corollary}
Let $f$ be a nonconstant Hardian function.  One has $\dege f' \neq (-1,-1,-1,\ldots)$ if and only if $\dege f \neq (0,0,0,\ldots)$, and, if those coindtions hold,  then the smallest nonnegative integer $N$ such that $\dege_N f' \neq -1$ is equal to the smallest nonnegative integer $N$ such that such that $\dege_N f \neq 0$.  For such a function $f$ and nonnegative integer $n$, we have the following.
\begin{enumerate}
\item  If $d = \dege_N f \neq \pm \infty$,  then
$$f' \sim df \frac{d}{dx}\log^{\circ (N+1)}x.$$
\item If $N = 0$ and $\dege_n f = \pm \infty$,  then
 $$\dege f' = \dege f.$$
\item If $N > 0$ and $\dege_N f = \pm \infty$, then
$$ f \frac{d}{dx} \log^{\circ (N+1)}x  \prec f'  \prec f {\frac{d}{dx} \log^{\circ N}x}.$$
\end{enumerate}
Consequently,   one has $$\dege f' = \dege f \oplus (-1,-1,-1,\ldots,-1,0,0,0,\ldots)$$
and
$$\dege f = \dege f' \oplus (1,1,1,\ldots,1,0,0,0,\ldots),$$
where
$$(-1,-1,-1,\ldots,-1,0,0,0,\ldots) = \dege  \frac{d}{dx} \log^{\circ (N+1)}x$$
denotes the vector with $N+1$ $-1$s followed by a  tail of $0$s,
and where
$$(1,1,1,\ldots,1,0,0,0,\ldots) = \dege\frac{1}{ \frac{d}{dx} \log^{\circ (N+1)}x}$$
denotes the vector with $N+1$ $1$s followed by a  tail of $0$s.
\end{corollary}

\section{Characterizations and generalizations of  the degree map}

In this section, we generalize degree in two  ways, and we use these generalizations to provide some  ``universal properties'' of the degree  map $\deg$.  This section is not used in later sections.

For any two functions  $f$ and $g$ in $\RR^{\RR_\infty}$, we write $f =_\infty g$ if $U \cap \dom f = U \cap \dom g$ and $f|_{U \cap \dom f}  = g|_{U \cap \dom g}$ for some neighborhood $U$ of $\infty$.   We write $f \leq_\infty g$ (resp., $f<_\infty g$) if $U \cap \dom f \subseteq U \cap \dom g$ and $f|_{U \cap \dom f}  \leq g|_{U \cap \dom f}$ (resp.,  $f|_{U \cap \dom f}  < g|_{U \cap \dom f}$, i.e., $f(x) < g(x)$ for all $x \in U \cap \dom f$).  Thus, one has $f =_\infty g$ if and only if $f \leq_\infty g$ and $g \leq_\infty f$.   The relation $=_\infty$ is an equivalence relation on the set $\RR^{\RR_\infty}$, and the relation $\leq_\infty$ induces a partial ordering on the set $\RR^{\RR_\infty}/{=_\infty}$ of $=_\infty$-equivalence classes.     For any subset $\mathcal{F}$ of $\RR^{\RR_\infty}$, we let $\mathcal{F}_\infty$ denote the image of $\mathcal{F}$ in $\RR^{\RR_\infty}/{=_\infty}$, which is thus a poset,  whose  elements are called {\bf germs at $\infty$ of functions in $\mathcal{F}$}.   We write $\leq$ (resp., $<$) for the relation on germs induced by the relation $\leq_\infty$ (resp., $<_\infty$).   

If $Y$ is any poset, then  $Y^{\RR_\infty}$ is defined,  and the definitions above extend  to $Y^{\RR_\infty}$, in the obvious way.   We are interested only in the cases $Y = \RR$ and $Y = \overline{\RR}$.  

Now,  let $\mathcal{F} $ be a subset of $\RR^{\RR_\infty}$, and let $\Phi: \mathcal{F} \longrightarrow \overline{\RR}^{\RR_\infty}$ be any function such that the induced function $\mathcal{F} \longrightarrow (\overline{\RR}^{\RR_\infty})_\infty$ factors through the map $\mathcal{F} \longrightarrow \mathcal{F}_\infty$, that is,  such that $f =_\infty g$ implies  $\Phi(f) =_\infty \Phi(g)$, for all $f,g \in \mathcal{F}$.   Equivalently, we may let $\Phi$ be a function $\Phi: \mathcal{F}_\infty \longrightarrow (\overline{\RR}^{\RR_\infty})_\infty$, and we can blur the distinction between the two.  Given such a map $\Phi$,  we let  $\deg(-,\Phi)\index[symbols]{.i  kq@$\deg(f,\Phi)$}$ denote the composition
$$\mathcal{F} \overset{\Phi}{\longrightarrow} \overline{\RR}^{\RR_\infty} \overset{\underset{x \to \infty}{\limsup}}{\longrightarrow} \overline{\RR},$$ that is, we set
$$\deg(f,\Phi) = \limsup_{x \to \infty} \Phi(f) \in \overline{\RR}$$
for all $f \in  \mathcal{F}$.    Likewise, we let
  $\underline{\deg} (-,\Phi)\index[symbols]{.i  kq@$\underline{\deg} (f,\Phi)$}$ denote the composition
$$\mathcal{F}  \overset{\Phi}{\longrightarrow} \overline{\RR}^{\RR_\infty} \overset{\underset{x \to \infty}{\liminf}}{\longrightarrow} \overline{\RR},$$ that is, we set
$$\underline{\deg} (f,\Phi) = \liminf_{x \to \infty} \Phi(f) \in \overline{\RR}$$
for all $f \in  \mathcal{F}$. 
We call $\deg(f,\Phi)$ (resp., $\underline{\deg}(f,\Phi)$) the {\bf $\Phi$-degree of $f$}\index{$\Phi$-degree of $f$}  (resp., {\bf lower $\Phi$-degree of $f$)}. \index{lower $\Phi$-degree of $f$}   

Of course, if $\Phi(f) = f$, then $\deg(f,\Phi) = \limsup_{x \to \infty} f(x)$ and $\underline{\deg}(f,\Phi)  = \liminf_{x \to \infty}f(x)$, for any $f \in \mathcal{F}$.   The ordinary degree and lower degree maps are precisely $\deg(-,L_-)$ and $\underline{\deg}(-,L_-)$, respectively, where $L_-: \RR^{\RR_\infty} \longrightarrow \overline{\RR}^{\RR_\infty}$ acts by $f \longmapsto L_f(x) = \frac{\log |f(x)|}{\log x}$, where as usual, we define $\log |0| = -\infty$, so that $\dom L_f = \RR_{>0} \cap \dom f$ for all $f \in \RR^{\RR_\infty}$.   

The following lemma is easily verified.

\begin{lemma}
Let $f \in \RR^{\RR_\infty}$.  One has
\begin{align*}
\limsup_{x \to \infty} f(x) & = \inf\{a \in \RR: f \leq_\infty a \} \\
 & = \inf\{a \in \RR: f <_\infty a \}.
\end{align*}
Let $t \in \overline{\RR}$.  One has $\limsup_{x \to \infty} f(x)  \leq t$ if and only if $f\leq_\infty s$ for all $s> t$, if and only if $f<_\infty s$ for all $s> t$.   Likewise,  one has $\limsup_{x \to \infty} f(x) < t$ if and only if $f\leq_\infty s$ for some $s< t$, if and only if $f<_\infty s$ for some $s< t$.
\end{lemma}

As an immediate consequence, we have the following.

\begin{proposition}
Let $\mathcal{F} \subseteq \RR^{\RR_\infty}$ and $\Phi: \mathcal{F}_\infty \longrightarrow (\overline{\RR}^{\RR_\infty})_\infty$ a function, and  let $f \in \mathcal{F}_\infty$.  One has
\begin{align*}
\deg(f,\Phi) & = \inf\{a \in \RR: \Phi(f) \leq a \} \\
 & = \inf\{a \in \RR: \Phi(f) < a \}.
\end{align*}
Let $t \in \overline{\RR}$.  One has $\deg(f,\Phi)  \leq t$ if and only if $\Phi(f) \leq s$ for all $s> t$, if and only if $\Phi(f) < s$ for all $s> t$.   Likewise,  one has $\deg(f,\Phi) < t$ if and only if $\Phi(f) \leq s$ for some $s< t$, if and only if $\Phi(f)< s$ for some $s< t$.
\end{proposition}

The lemma and proposition above have obvious counterparts for limits inferior and lower $\Phi$-degree, respectively.

Let $R$ and $S$ be a rings.  A {\bf logarithm on $R$}\index{logarithm on $R$} is a group homomorphism $L:  R^* \rightarrow R^+$.    More generally, we say that a {\bf logarithm on $R$ in $S$}\index{logarithm on $R$ in $S$} is a group homomorphism $L:  R^* \rightarrow S^+$.    If $L$  is a logarithm on $R$ in $S$, then
$$0 = L(1) = L((-1)^2) = 2L(-1),$$ and therefore, if $S$ has characteristic other than $2$, then  $L(-1) = 0$, whence $L(-f) = L(f)$ for all $f \in R$.  

Let $R$ be a subring of $\mathcal{R}$ and $L: R^* \longrightarrow \mathcal{R}$ a logarithm on $R$ in $\mathcal{R}$.  One has $$\deg(rs,L) \leq \deg(r,L)+\deg(s,L)$$ 
and
$$\underline{\deg}(rs,L) \geq \underline{\deg}(r,L)+\underline{\deg}(s,L)$$ 
for all $r,s\in R$ whenever the respective sums are defined in $\overline{\RR}$.
If $L: R^* \longrightarrow  K^+$ is logarithm on $R$ in some Hardy field $K$,    then  $\deg(-,L) = \underline{\deg}(-,L)$ is the composition
$$R^* \overset{L}{\longrightarrow} K \overset{\underset{x \to \infty}{\lim}}{\longrightarrow} \overline{\RR},$$ that is, 
$$\deg(r,L) =\underline{\deg}(r,L) = \lim_{x \to \infty} L(r)$$
for all $r \in R^*$.  

The ring $\mathcal{C}^{(1)}$  is equipped with the logarithm $$\operatorname{Dlog}: (\mathcal{C}^{(1)} )^* \longrightarrow \mathcal{C}^+$$ on $\mathcal{C}^{(1)}$ in $\mathcal{C}$, and thus
$$\operatorname{Dlog}_h = h \cdot \operatorname{Dlog}:  (\mathcal{C}^{(1)})^*\longrightarrow \mathcal{C}^+$$ 
is a logarithm on $\mathcal{C}^{(1)}$ in  $\mathcal{C}$ for all $ h \in\mathcal{C}$.  We let $$\deg_h = \deg(-,\operatorname{Dlog}_h): (\mathcal{C}^{(1)})^* \longrightarrow \overline{\RR}\index[symbols]{.i  kr@$\deg_h$}$$ denote the corresponding degree map and $$\underline{\deg}_h =\underline{ \deg}(-,\operatorname{Dlog}_h): (\mathcal{C}^{(1)})^* \longrightarrow \overline{\RR}\index[symbols]{.i  kr@$\underline{\deg}_h$}$$the corresponding lower degree map.  Of particular interest is $h = x = \id$, for one has
$$\underline{\deg}_x r = \liminf_{x \to \infty} \frac{x r'(x)}{r(x)} \leq \underline{\deg}\, r \leq  \deg r  \leq   \limsup_{x \to \infty} \frac{x r'(x)}{r(x)} = \deg_x r$$
for all $r \in (\mathcal{C}^{(1)})^*$.  Thus, all are equal if and only if  $\lim_{x \to \infty} \frac{x r'(x)}{r(x)}$ exists or is $\pm \infty$,  both of which hold, for example, if $r$ is Hardian, that is, if $r$ is an element of some Hardy field.  We denote by $LD$ the logarithm
$$LD = \operatorname{Dlog}_x = x \cdot \operatorname{Dlog}:  (\mathcal{C}^{(1)})^*\longrightarrow \mathcal{C}^+$$ 
on $\mathcal{C}^{(1)}$ in $\mathcal{C}$.

We  wish to determine the relationships between the degree maps $\deg_h$ for various $h\in \mathcal{C}_{>0}$.    For all $f \in (\mathcal{C}^{(1)})^*$, one has
$$LD(f(\log x)) = \frac{xf'(\log x)\frac{1}{x}}{f(\log x)} = (\operatorname{Dlog} f(x))\circ \log x$$
and
$$\operatorname{Dlog} (f(e^x)) =  \frac{f'(e^x)e^x}{f(e^x)} = (LD(f(x)) \circ e^x,$$
whence
$$ \deg_1 f(x) = \deg_x f(\log x) $$
and
$$\deg_x f(x) = \deg_1 f(e^x),$$
 Thus, $\deg_x$ and $\deg_1$ both determine the other.  More generally, one has the following.

\begin{proposition}
Let $f \in (\mathcal{C}^{(1)})^*$.  Let $g \in  \mathcal{C}^{(1)}$ with $\lim_{x \to \infty} g(x) = \infty$ and $h = 1/g'  \in \mathcal{C}_{>0}$.   (Equivalently, let $h \in \mathcal{C}_{>0}$ with $\lim_{x \to \infty} g(x) = \infty$, where $g$ is any antiderivative of $1/h$.)  Then $g^{-1} \in  \mathcal{C}^{(1)}$ exists and $\lim_{x \to \infty} g^{-1}(x) = \infty$,
and $1/(g^{-1})' = g' \circ g^{-1} \in  \mathcal{C}_{>0}$, and one has the following.
\begin{enumerate}
\item $\operatorname{Dlog} (f \circ g) \circ g^{-1} = \operatorname{Dlog}_h f = \operatorname{Dlog} (f \circ g^{-1}) \circ g.$
\item $\operatorname{Dlog} (f \circ g)  = (\operatorname{Dlog}_{h} f) \circ g$.
\item $ \operatorname{Dlog} (f \circ g^{-1}) = (\operatorname{Dlog}_h f ) \circ g^{-1}.$
\item $\deg_1 (f \circ g) = \deg_h f = \deg_1 (f \circ g^{-1}).$
\item  $\operatorname{Dlog}_h (f \circ  g) \circ g^{-1}  = \operatorname{Dlog} f  = \operatorname{Dlog}_h (f \circ g^{-1}) \circ g.$
\item $\operatorname{Dlog}_h (f \circ g) = (\operatorname{Dlog} f) \circ g$.
\item  $\operatorname{Dlog}_h (f \circ g^{-1}) = (\operatorname{Dlog} f) \circ g^{-1}$.
\item $\deg_h (f \circ g) = \deg_1 f = \deg_h (f \circ g^{-1}).$
\end{enumerate}
\end{proposition}

Now, one has an injective group homomorphism
$$
\xymatrix{
  \mathcal{G} \ \ar@{>->}[r] &  \Aut(\mathcal{C}^{(1)})}$$
acting by $g \longmapsto - \circ g^{-1}$, where  
$$\mathcal{G} = \{g\in  \mathcal{C}^{(1)}: \lim_{x \to \infty} g(x) = \infty, \, g' \in \mathcal{C}_{>0}\}$$
 is a group  under composition of functions,  and where $\Aut(\mathcal{C}^{(1)})$ is the group of automorphisms of the ring 
$ \mathcal{C}^{(1)}$.   As each logarithm $\operatorname{Dlog}_h  = h \cdot \operatorname{Dlog}: (\mathcal{C}^{(1)})^* \longrightarrow \mathcal{C}$ for $g \in \mathcal{G}$ and $h = 1/g'$ determines the others, so each degree map $\deg_h$ determines the others.  Equivalently, the $\deg_h$ are a one-parameter family of maps from  $\mathcal{C}^{(1)}$ to $\overline{\RR}$ generated by $\deg = \deg_x$, $\deg_1$, or $\deg_h$ for any element  $h$ as above.  For any $g \in \mathcal{G}$ and $h = 1/g'$,  one has the commutative diagram
\begin{align}\label{commdiag}
\begin{gathered}
\xymatrix{ (\mathcal{C}^{(1)})^* \ar[r]^{-\circ g^{-1}} \ar[dr]^{\operatorname{Dlog}_h}   \ar@/^1pc/[dd]^(.66){\deg_h} & (\mathcal{C}^{(1)})^* \ar[r]^{-\circ g} \ar[dr]^{\operatorname{Dlog}}   \ar@/^1pc/[dd]^(.66){\deg_1}|\hole & (\mathcal{C}^{(1)})^*   \ar[dr]^{\operatorname{Dlog}_h}  \ar@/^1pc/[dd]^(.66){\deg_h}|\hole & 
\\  & \mathcal{C}  \ar[dl]^(.33){\underset{x \to \infty}{\limsup}}  \ar[r]^{-\circ g^{-1}} &   \mathcal{C} 
 \ar[dl]^(.33){\underset{x \to \infty}{\limsup}}  \ar[r]^{-\circ g}& { \mathcal{C} \ }    \ar[dl]^(.33){\underset{x \to \infty}{\limsup}}   \\
\overline{\RR}   \ar[r]^{\id} & \overline{\RR}   \ar[r]^{\id} & \overline{\RR} & 
}
\end{gathered}
\end{align}
where the top  horizontal row are group automorphisms,  the middle row are ring automorphisms,  the top three descending arrows  are surjective group homomorphisms with kernel $\RR^*$, all descending arrows are surjective, and the composition across each row is the identity function.    Note, for example,  that one has $g = \log x$ if and only if  $h = 1/g' =x$,  if and only if $g^{-1} = \exp x$,  if and only if $1/(g^{-1})' = \exp(-x)$,  if and only if $\operatorname{Dlog}_h = LD$.  

 For all $g$ in the subgroup $$\mathcal{G} \cap  \mathcal{D} =\{g \in \mathcal{D}: \lim_{x \to \infty} g(x) = \infty, \, g' \in \mathcal{D}_{>0}\}$$
of $\mathcal{G}$,  the commutative diagram (\ref{commdiag}) descends to
\begin{align}\label{commdiag2}
\begin{gathered}
\xymatrix{ \mathcal{D}^* \ar[r]^{-\circ g^{-1}} \ar[dr]^{\operatorname{Dlog}_h}   \ar@/^1pc/[dd]^(.66){\deg_h} &  \mathcal{D}^* \ar[r]^{-\circ g} \ar[dr]^{\operatorname{Dlog}}   \ar@/^1pc/[dd]^(.66){\deg_1}|\hole &  \mathcal{D}^*   \ar[dr]^{\operatorname{Dlog}_h}  \ar@/^1pc/[dd]^(.66){\deg_h}|\hole & 
\\  & \mathcal{D}  \ar[dl]^(.33){\underset{x \to \infty}{\limsup}}  \ar[r]^{-\circ g^{-1}} &   \mathcal{D} 
 \ar[dl]^(.33){\underset{x \to \infty}{\limsup}}  \ar[r]^{-\circ g}& { \mathcal{D} \ }    \ar[dl]^(.33){\underset{x \to \infty}{\limsup}}   \\
\overline{\RR}   \ar[r]^{\id} & \overline{\RR}   \ar[r]^{\id} & \overline{\RR} & 
}
\end{gathered}
\end{align}
where now the given logarithms are logarithms on $\mathcal{D}$.
  We note the following.

  \begin{proposition}
   Let $R$ be a Hardy ring, and let $g \in \mathcal{G} \cap  \mathcal{D} =\{g \in \mathcal{D}: \lim_{x \to \infty} g(x) = \infty, \, g' \in \mathcal{D}_{>0}\}$.    If $1/g' \in  R$, then $R \circ g^{-1}$ is a Hardy ring.     Moreover, the converse holds if  $f' \in R^*$ for some $f \in R$.
 \end{proposition}
 
 \begin{corollary}[{\cite[Lemma 6.4]{bos}}]\label{hardyg}
 Let $K$ be a Hardy field not contained in $\RR$, and let $g  \in \mathcal{G} \cap  \mathcal{D} =\{g \in \mathcal{D}: \lim_{x \to \infty} g(x) = \infty, \, g' \in \mathcal{D}_{>0}\}$.   Then 
 $K \circ g^{-1}$ is a Hardy field if and only if  $g' \in K$.  Thus, if $g \in K$, then  $K \circ g^{-1}$ is a Hardy field containing $\id$.
 \end{corollary}

Thus,  for any Hardy ring $R$ and any $g \in \mathcal{G} \cap  \mathcal{D}$, if $h = 1/g' \in R$  (e.g.,  if $g' \in R^*$), then $R \circ g^{-1}$ is a Hardy ring and  the commutative diagram (\ref{commdiag2}) descends to
$$
\xymatrix{ R^* \ar[r]^{-\circ g^{-1}} \ar[dr]^{\operatorname{Dlog}_h}   \ar@/^1pc/[dd]^(.66){\deg_h} & (R \circ g^{-1})^*  \ar[r]^{-\circ g} \ar[dr]^{\operatorname{Dlog}}   \ar@/^1pc/[dd]^(.66){\deg_1}|\hole & R^*   \ar[dr]^{\operatorname{Dlog}_h}  \ar@/^1pc/[dd]^(.66){\deg_h}|\hole & 
\\  & R  \ar[dl]^(.33){\underset{x \to \infty}{\limsup}}  \ar[r]^{-\circ g^{-1}} &   R \circ g^{-1}
 \ar[dl]^(.33){\underset{x \to \infty}{\limsup}}  \ar[r]^{-\circ g}& {R\ }    \ar[dl]^(.33){\underset{x \to \infty}{\limsup}}   \\
\overline{\RR}   \ar[r]^{\id} & \overline{\RR}   \ar[r]^{\id} & \overline{\RR} & 
}$$
while,  if $1/(g^{-1})' = g' \circ g^{-1} \in R$ (e.g., if $(g^{-1})' \in R^*$),  then $R \circ g$ is a Hardy ring  and the commutative diagram (\ref{commdiag2}) descends to
$$
\xymatrix{ (R \circ g)^* \ar[r]^{-\circ g^{-1}} \ar[dr]^{\operatorname{Dlog}_h}   \ar@/^1pc/[dd]^(.66){\deg_h} & R^*  \ar[r]^{-\circ g} \ar[dr]^{\operatorname{Dlog}}   \ar@/^1pc/[dd]^(.66){\deg_1}|\hole & (R \circ g)^*   \ar[dr]^{\operatorname{Dlog}_h}  \ar@/^1pc/[dd]^(.66){\deg_h}|\hole & 
\\  & R \circ g  \ar[dl]^(.33){\underset{x \to \infty}{\limsup}}  \ar[r]^{-\circ g^{-1}} &   R
 \ar[dl]^(.33){\underset{x \to \infty}{\limsup}}  \ar[r]^{-\circ g}& {R \circ g \ }    \ar[dl]^(.33){\underset{x \to \infty}{\limsup}}   \\
\overline{\RR}   \ar[r]^{\id} & \overline{\RR}   \ar[r]^{\id} & \overline{\RR} & 
}$$
For example,  $R \circ \exp$ is a Hardy ring if $x\in R$, and $R \circ \log$ is a Hardy ring if $\exp(-x) \in R$.  If $R$ is a Hardy field, then the limits superior in the commutative diagrams above can be replaced with limits, and $\deg_x = \deg$ on $R^*$ if $x \in R$.

Now, let $R$ be a subring of $\mathcal{R}$ and $L$ a logarithm on $R$ in $\mathcal{R}$.
Naturally, we wish to impose further conditions on $L$ so that the $L$-degree map satisfies further desirable conditions.  For example,  the condition that $|r| \leq |s|$ implies $\deg(r,L) \leq \deg(s,L)$ for all $r,s \in R$ is desirable, and it holds if and only if $r \leq s$ implies $\limsup_{x \to \infty} L(r) \leq \limsup_{x \to \infty} L(s)$ for all $r,s \in R^* \cap R_{>0}$. This condition holds if $L$ itself is nondecreasing (on $R^* \cap R_{>0}$), but it may hold even if $L$ is not nondecreasing.  Let us say that a  logarithm on $R$ in $\mathcal{R}$ is {\bf semi-increasing} \index{semi-increasing logaritihm} if $r \leq s$ implies $\limsup_{x \to \infty} L(r) \leq \limsup_{x \to \infty} L(s)$ for all $r,s \in R^*\cap R_{>0}$.    Of course,  if $E$ is an exponential on $R$, then $E^{-1}$ is an increasing logarithm on $R$, and any nondecreasing logarithm is semi-increasing.   Moreover, the map $\log \circ -: R^*\cap R_{>0} \longrightarrow \mathcal{R}$ is an increasing logarithm on $R$ in $\mathcal{R}$.    Thus, if $\log f \in R$ for all $f \in R^*\cap R_{>0}$, then the map $\log \circ -: R^*\cap R_{>0} \longrightarrow R$ is an increasing logarithm on $R$.

If $R = K$ is a Hardy field,  then the logarithm $\operatorname{Dlog}$ on $K$ is semi-increasing.  To see this, suppose that $0<r\leq s$,  so that $\frac{r}{s} \to A \in [0,1]$ and therefore  $(\frac{r}{s})' \to 0$ as $x \to \infty$.  If $A > 0$, then  $\operatorname{Dlog}(r) - \operatorname{Dlog}(s) = \operatorname{Dlog}\frac{r}{s} \to \frac{0}{A} =  0$ as $x \to \infty$, so that $$\lim_{x \to \infty} \operatorname{Dlog}(r)(x) = \lim_{x \to \infty} \operatorname{Dlog}(s)(x).$$ 
Suppose, on the other hand, that $A = 0$.    Then $\frac{s}{r}>0$ is unbounded,  so  that $(\frac{s}{r})' > 0$, whence
$\operatorname{Dlog}(s) - \operatorname{Dlog}(r) = \operatorname{Dlog}\frac{s}{r} >0$ and therefore
$\operatorname{Dlog}(s) > \operatorname{Dlog}(r)$, whence 
$$\lim_{x \to \infty} \operatorname{Dlog}(r)(x) \leq \lim_{x \to \infty} \operatorname{Dlog}(s)(x).$$ 
Thus, $\operatorname{Dlog}$ is semi-increasing.
 Note that $\operatorname{Dlog}$ is not nondecreasing if $K$ properly  contains  $\RR$, since then $\mm_K \neq 0$ and thus there is some  $f > 0$ with $f  \in \mm_K$,  whereby $0<1-f<1+f$ and yet $\operatorname{Dlog}(1-f) > 0 > \operatorname{Dlog}(1+f)$.   Nevertheless, the argument above shows that,  if $r,s \in K_{>0}$ with $r \prec s$, then $\operatorname{Dlog}(r) < \operatorname{Dlog}(s)$.

Now, let $K$ be a  Hardy field and $L$ a semi-increasing logarithm on $K$.   We let $\deg(0,L) = -\infty$, so that $\deg(-,L)$ is defined on all of $K$.
For all $r,s \in K$, one has
$$|r+s| \leq |r|+|s| \leq \max(|2r|,|2s|),$$
whence, if $r$, $s$, and $r+s$ are nonzero, then
\begin{align*}
\deg(r+s,L) & \leq \deg(\max(|2r|,|2s|),L) \\
& = \max(\deg(2r,L) ,\deg(2s,L)) \\
& = \max(\deg(r,L) ,\deg(s,L))+\deg(2,L),
\end{align*}
where the last equation holds provided that $\deg(2,L)$ is finite.   This also holds if any of $r$, $s$, and $r+s$ are $0$.
Now, for any $c \in \RR^*$, one has
$$\deg(c,L) = \lim_{x \to \infty} L(|c|).$$
Let $c >0$ with $c \neq 1$.  If $c$ has $L$-degree $0$, then $c^n$ has $L$-degree $0$ for all integers $n$, whence, for any $d > 0$, one has $c^n \leq d < c^{n+1}$ or $c^{n+1} < d \leq c^{n}$ for some integer $n$, and therefore $d$ also has $L$-degree $0$.
Thus all nonzero constants have $L$-degree $0$ if and only if some nonzero constant $c$ with $c \neq 1$ has $L$-degree $0$, if and only if $2$ has $L$-degree $0$, if and only if  $L(c) \prec 1$ for all $c > 0$, if and only if $L(c) \prec 1$ for some $c > 0$ with $c \neq 1$.   If this condition holds, then $\deg(cr,L) = \deg(r,L)$ for all $r \in K$ and all $c \neq 0$, so that, if  $r \preceq s$,  that is, if $|r| \leq c |s|$ for some $c > 0$, then $\deg(r,L)  \leq \deg(s,L)$.    Conversely, if $r \preceq s$ implies $\deg(r,L)  \leq \deg(s,L)$ for all $r,s \in K$, then, since $2  \asymp 1$, one has $\deg(2,L) = \deg(1,L) = 0$, whence  all nonzero constants have $L$-degree $0$.
 Moreover, if this condition holds, then the inequality
\begin{align*}
\deg(r+s,L)\leq  \max(\deg(r,L) ,\deg(s,L))
\end{align*}
holds for all $r,s \in K$.  Conversely, if the  inequality above holds for all $r,s \in K$, then
$$0 = \deg(1,L) \leq \deg(2,L) \leq\max(\deg(1,L),\deg(1,L)) = \deg(1,L) = 0,$$
and thus all nonzero constants have $L$-degree $0$.  We thus have the following.

\begin{proposition}\label{Hproper}
Let $K$ be a  Hardy field and $L$ a semi-increasing logarithm on $K$. The following conditions are equivalent.
\begin{enumerate}
\item All nonzero constants have $L$-degree $0$.
\item Some nonzero constant $c \neq 1$ has $L$-degree $0$.
\item $L(c) \prec 1$ for all $c \in \RR_{>0}$ (or for some $c \in \RR_{>0}$ with $c \neq 1$).
\item $\deg(cr,L) = \deg(r,L)$ for all $r \in K$ and all $c \in \RR^*$.
\item $r \asymp s$ implies $\deg(r,L)= \deg(s,L)$ for all $r,s \in  K$.
\item $r \preceq s$ implies $\deg(r,L) \leq \deg(s,L)$ for all $r,s \in  K$.
\item $\deg(r,L) <\deg(s,L)$ implies $r \prec s$ for all $r,s \in  K$.
\item $\deg(r+s,L)\leq  \max(\deg(r,L) ,\deg(s,L))$ for all $r,s \in K$.
\item The map $\deg(-,L): K \longrightarrow \overline{\RR}$ factors (uniquely) through the valuation $v: K \longrightarrow  \Gamma_K \cup \{\infty\}$ on $K$.
\end{enumerate}
\end{proposition}

Next, we note the following.

\begin{proposition}
Let $R$ be a subring of $\mathcal{R}$ containing $\RR$ and $E: R^+ \longrightarrow R^*\cap R_{>0}$ an exponential on $R$,  let $L= E^{-1}$, and let $r \in R^*$.  One has
\begin{align*}
\deg(r,L) & = \inf\{a \in \RR: L(|r|) \leq a \} \\
& = \inf\{a \in \RR: |r| \leq E(a)\} \\
 & = \inf\{a \in \RR: L(|r|) < a \} \\
& = \inf\{a \in \RR: |r| < E(a)\}.
\end{align*}
In particular, the map $\deg(-,L)$ depends only on the map $E|_\RR: \RR \longrightarrow E(\RR)$, which is an isomorphism of ordered abelian groups with inverse  $\deg(-,L)|_{E(\RR)}: E(\RR) \longrightarrow \RR$.   In particular, one has $\deg(E(a),L) = a$ for all $a \in \RR$.
\end{proposition}

Let $R$ be a subring of $\mathcal{R}$ containing $\RR$.  Since $\deg(-,L)$, for any exponential $E$ on $R$ with $L = E^{-1}$, depends only on the map $E|_\RR$, we may generalize the definition of $\deg(-,L)$  in this case as follows.   Let $E: \RR \longrightarrow \mathcal{R}_{>0}$ be any nondecreasing map of posets.   For all $f \in \mathcal{R}$,  we define
\begin{align*}
\deg_E f   = \inf\{a \in \RR: |f| \leq E(a)\} 
\end{align*}
and
\begin{align*}
\underline{\deg}_E f   = \sup\{a \in \RR: |f| \geq E(a)\} \leq \deg_E f.
\end{align*}
(These definitions extend naturally to all $f \in (\RR^{\RR_\infty})_\infty$, but for the sake of simplicity we restirct our attention to $f \in \mathcal{R} \subsetneq (\RR^{\RR_\infty})_\infty$.)
The following proposition is an easy consequence of the definition of $\deg_E$ and has an obvious dual for $\underline{\deg}_E$.

\begin{proposition}
Let $E: \RR \longrightarrow  \mathcal{R}_{>0}$ be an nondecreasing  map of posets, and let $f,g \in \mathcal{R}$.   One has the following.
\begin{enumerate}
\item $\deg_E E(a) \leq a$ for all $a \in \RR$. 
\item For all $t \in \overline{\RR}$, one has $\deg_E f  \leq t$ if and only if $ |f| \leq E(s)$ for all $s> t$, and  $\deg_E f < t$ if and only if $|f| \leq  E(s)$ for some $s< t$.  
\item  Suppose that $E$ is increasing. Then $\deg_E E(a) = a$ for all $a \in \RR$, and one has
\begin{align*}
\deg_E f = \inf\{a \in \RR: |f| < E(a)\}.
\end{align*}
Moreover,  for all $t \in \overline{\RR}$, one has $\deg_E f  \leq t$ if and only if $ |f| < E(s)$ for all $s> t$, and one has $\deg_E f < t$ if and only if $|f| < E(s)$ for some $s< t$. 
\end{enumerate}
\end{proposition}

In  the case where $E$ is the restriction of an exponential on some subring of $\mathcal{R}$ containing $\RR$, the map $E$ is an embedding of ordered abelian groups.    The next several results establish properties of the maps  $\deg_E$ and $\underline{\deg}_E$ for embeddings $E: \RR \longrightarrow  \mathcal{R}_{>0}$ of ordered abelian groups.

\begin{proposition}
Let $E: \RR \longrightarrow  \mathcal{R}_{>0}$ be a nondecreasing map of posets.  Then $E$ is an embedding of ordered abelian groups if and only if  $E(1) >1$ and $E(a) = E(1)^a$ for all $a \in \RR$.   Therefore, such an embedding $E$ is freely and uniquely determined by $E(1) >1$.   Suppose that $E$ is  an embedding  of ordered abelian groups, and let $f,g \in \mathcal{R}$.   One has the following.
\begin{enumerate}
\item $\deg_E f^a =a \deg_E f$ for all $a \in \RR_{\geq 0}$.
\item $\underline{\deg}_E f  = -\deg_E(1/f)$  and $\deg_E f^a =a \, \underline{\deg}_E f$ and  for all $a  \in \RR_{\leq 0}$  if $f \in  \mathcal{R}^*$.
\item $\deg_E (fg)  \leq \deg_E f+\deg_E g$, provided that the sum is defined in $\overline{\RR}$.
\item $\underline{\deg}_E (fg)  \geq \underline{\deg}_E f+\underline{\deg}_E g$, provided that the sum is defined in $\overline{\RR}$.
\end{enumerate}
\end{proposition}

\begin{proof}
Suppose that  $E: \RR^+ \longrightarrow \mathcal{R}_{>0}$ is any nondecreasing group homomorphism.
For all $a \in \ZZ$,  one has $E(a) = E(1)^a$.  In fact, this extends easily to all $a \in \QQ$.  Now, let $c \in \RR$ and $a,b \in \QQ$ with $a<c<b$.  Then
$$0<E(1)^a =  E(a)\leq E(c) \leq E(b) = E(1)^c.$$
Taking limits through the rationals as $a \to c^-$ and $b\to c^+$,  we see that
$$E(c) = E(1)^c$$ and $E(1) \geq 1$.
Moreover, $E$ is increasing if and only if $E(1)> 1$.  Thus,  an increasing group homomorphism $E: \RR^+ \longrightarrow \mathcal{R}_{>0}$, that is, an embedding $E: \RR^+ \longrightarrow \mathcal{R}_{>0}$ of ordered abelian groups, is freely and uniquely determined by $E(1) >1$.   Statements (1)--(4) are then readily  verified.
\end{proof}

\begin{proposition}
Let $E_1,E_2:  \RR^+ \longrightarrow  \mathcal{R}_{>0}$ be embeddings of ordered abelian groups.  The following conditions are equivalent.
\begin{enumerate}
\item There is a  (unique) $A> 0$  such that $\deg_{E_1} f = A \deg_{E_2} f$ for all $f \in \mathcal{R}$ with $\deg_{E_1} f\geq 0$  or $\deg_{E_2} f\geq 0$.
\item $\deg_{E_1} f = \frac{ \deg_{E_2} f}{\deg_{E_2} E_1(1)}$ for all $f \in \mathcal{R}$ with $\deg_{E_1} f\geq 0$  or $\deg_{E_2} f\geq 0$,   where $\deg_{E_2} E_1(1) \in (0,\infty)$. 
\item  $\deg_{E_1} f = \deg_{E_1} E_2(1) \, \deg_{E_2} f$ for all $f \in \mathcal{R}$ with $\deg_{E_1} f\geq 0$  or $\deg_{E_2} f\geq 0$,  where $\deg_{E_1} E_2(1) \in (0,\infty)$.
\item $ \deg_{E_2} E_1(1) \in (0,\infty)$.
\item $ \deg_{E_1} E_2(1) \in (0,\infty)$.
\end{enumerate}
Likewise, the following conditions are equivalent.
\begin{enumerate}
\item[(6)] There is a  (unique) $B> 0$  such that $\deg_{E_1} f = B \deg_{E_2} f$ for all $f \in \mathcal{R}$ with $\deg_{E_1} f\leq 0$  or $\deg_{E_2} f\leq 0$.
\item[(7)]  $\deg_{E_1} f = \frac{ \deg_{E_2} f}{\underline{\deg}_{E_2} E_1(1)}$ for all $f \in \mathcal{R}$ with $\deg_{E_1} f\leq 0$  or $\deg_{E_2} f\leq 0$,   where $\underline{\deg}_{E_2} E_1(1) \in (0,\infty)$. 
\item[(8)]   $\deg_{E_1} f = \underline{\deg}_{E_1} E_2(1) \deg_{E_2} f$ for all $f \in \mathcal{R}$ with $\deg_{E_1} f\leq 0$  or $\deg_{E_2} f\leq 0$,  where $\underline{\deg}_{E_1} E_2(1) \in (0,\infty)$.
\item[(9)]  $ \underline{\deg}_{E_2} E_1(1) \in (0,\infty)$.
\item[(10)]  $\underline{\deg}_{E_1} E_2(1) \in (0,\infty)$.
\end{enumerate}
\end{proposition}

\begin{proof}
We prove the first half of the theorem; the proof of the latter half is similar.
It is clear that  $(1) \Leftrightarrow (2) \Leftrightarrow (3)$,  $(2) \Rightarrow (4)$, and $(3) \Rightarrow (5)$.   
Suppose that (4) holds,  let $A = \frac{1}{\deg_{E_2} E_1(1)} > 0$,  let $f \in \mathcal{R}$, and let $t \in [0,\infty]$.  Note that 
$ |f| < E_1(s) = E_1(1)^s$ for all $s> t$ if and only if $ |f| < E_2(1)^{s\deg_{E_2} E_1(1) } = E_2(s\deg_{E_2} E_1(1) )$ for all $s >t$, if and only if $|f^{A}| < E_2(s)$ for all $s >t$.   Thus, one has $\deg_{E_1} f \leq t$ if and only if $A\deg_{E_2} f = \deg_{E_2} f^A \leq t$.  Thus  (4) implies (2), and, likewise, (5) implies (3).    This completes the proof.
\end{proof}

\begin{corollary}\label{degEcor}
Let $E:  \RR^+ \longrightarrow  \mathcal{R}_{>0}$ an embedding of ordered abelian groups.  The following conditions are equivalent.
\begin{enumerate}
\item There is a  (unique) $A> 0$  such that $\deg_E f = A \deg  f$ for all $f \in \mathcal{R}$ with $\deg_E f \geq 0$ or $\deg  f\geq 0$.
\item $\deg_E f = \frac{ \deg  f}{\deg E(1)}$ for all $f \in \mathcal{R}$ with $\deg_E f \geq 0$ or $\deg  f\geq 0$, where $\deg E(1) \in (0,\infty)$. 
\item  $\deg_E f = \deg_E x \, \deg  f$ for all $f \in \mathcal{R}$ with $\deg_E f \geq 0$ or $\deg  f\geq 0$,  where $\deg_E x \in (0,\infty)$.
\item $ \deg E(1) \in (0,\infty)$.
\item $\deg_E x \in (0,\infty)$.
\end{enumerate}
Likewise, the following conditions are equivalent.
\begin{enumerate}
\item[(6)]  There is a  (unique) $B> 0$  such that $\deg_E f = B \deg  f$ for all $f \in \mathcal{R}$ with $\deg_E f \geq 0$ or $\deg  f\geq 0$.
\item[(7)]  $\deg_E f = \frac{ \deg  f}{\underline{\deg} \, E(1)}$ for all $f \in \mathcal{R}$ with $\deg_E f \geq 0$ or $\deg  f\geq 0$, where $\underline{\deg} \, E(1) \in (0,\infty)$. 
\item[(8)]   $\deg_E f = \underline{\deg}_E x \deg  f$ for all $f \in \mathcal{R}$ with $\deg_E f \geq 0$ or $\deg  f\geq 0$,  where $\underline{\deg}_E x \in (0,\infty)$.
\item[(9)]  $ \underline{\deg}\, E(1) \in (0,\infty)$.
\item[(10)]  $\underline{\deg}_E x \in (0,\infty)$.
\end{enumerate}
\end{corollary}

\begin{corollary}
Let $E_1,E_2:  \RR^+ \longrightarrow  \mathcal{R}_{>0}$ be embeddings of ordered abelian groups.  The following conditions are equivalent.
\begin{enumerate}
\item There is a  (unique) $A> 0$ such that $\deg_{E_1} f = A \deg_{E_2} f$ for all $f \in \mathcal{R}$.
\item $\deg_{E_1} f = \frac{ \deg_{E_2} f}{\deg_{E_2} E_1(1)}$ for all $f \in \mathcal{R}$, where $\deg_{E_2} E_1(1)\in (0,\infty)$.
\item  $\deg_{E_1} f = \deg_{E_1} E_2(1) \deg_{E_2} f$ for all $f \in \mathcal{R}$, where $\deg_{E_1} E_2(1)  \in (0,\infty)$.
\item $\overline{\underline{\deg}}_{E_2} E_1(1) \in (0,\infty)$.
\item $\overline{\underline{\deg}}_{E_2} E_1(a) \in \RR^*$ for all $a \in \RR$.
\item $\overline{\underline{\deg}}_{E_1} E_2(1) \in (0,\infty)$.
\item $\overline{\underline{\deg}}_{E_1} E_2(a)  \in \RR^*$ for all $a \in \RR$.
\item $\deg_{E_2} E_1(a) = a\deg_{E_2} E_1(1) \in \RR$  for all $a \in \RR$ (or, equivalently, for $a = -1)$.
\item $\deg_{E_2} E_1(a) = a\deg_{E_1} E_2(1) \in \RR$  for all $a \in \RR$ (or, equivalently, for $a = -1)$.
\item There is a  (unique) $A> 0$ such that, for all real numbers $a< b< c$, one has $E_2(a) < E_1(Ab) < E_2(c)$.
\item  There is a  (unique) $A> 0$ such that, for all real numbers $a< b< c$, one has $E_1(a) < E_2(b/A) < E_1(c)$.
\end{enumerate}
Moreover, $A$ is the same constant throughout.
\end{corollary}

\begin{corollary}\label{degEcor2}
Let $E:  \RR^+ \longrightarrow  \mathcal{R}_{>0}$ an embedding of ordered abelian groups.  The following conditions are equivalent.
\begin{enumerate}
\item There is a  (unique) $A> 0$ such that $\deg_E f = A \deg f$ for all $f \in \mathcal{R}$.
\item $\deg_E f = \frac{ \deg f}{\deg E(1)}$ for all $f \in \mathcal{R}$, where $\deg E(1) \in (0,\infty)$.
\item  $\deg_E f = \deg_E x\, \deg f$ for all $f \in \mathcal{R}$, where $\deg_E x  \in (0,\infty)$.
\item $\overline{\underline{\deg}}  \, E(1) \in (0,\infty)$.
\item $\overline{\underline{\deg}}  \, E(a) \in \RR^*$ for all $a \in \RR$.
\item $\overline{\underline{\deg}}_E x \in (0,\infty)$.
\item $\overline{\underline{\deg}}_E x^a  \in \RR^*$ for all $a \in \RR$.
\item $\deg  E(a) = a\deg  E(1) \in \RR$  for all $a \in \RR$ (or, equivalently, for $a = -1)$.
\item $\deg_E x^a = a\deg_E x \in \RR$  for all $a \in \RR$ (or, equivalently, for $a = -1)$.
\item There is a  (unique) $A> 0$ such that, for all real numbers $a< b< c$, one has $x^{a} < E(Ab) < x^{c}$.
\item  There is a  (unique) $A> 0$ such that, for all real numbers $a< b< c$, one has $E(a) < x^{b/A} < E(c)$.
\end{enumerate}	
Moreover, $A$ is the same constant throughout.
\end{corollary}

\begin{corollary}
Let  $E_1,E_2: \RR^+ \longrightarrow \mathcal{R}_{>0}$ be embeddings of ordered abelian groups.  The following conditions are equivalent.
\begin{enumerate}
\item $\deg_{E_1} f = \deg_{E_2} f$ for all $f \in \mathcal{R}$.
\item  $\overline{\underline{\deg}}_{E_1} E_2(1) = 1$.
\item  $\overline{\underline{\deg}}_{E_1} E_2(a) = a$  for all $a \in \RR$.
\item $\deg_{E_1} E_2(a) = a$  for all $a \in \RR$ (or, equivalently, for $a = \pm 1)$.
\item For all real numbers $a<b<c$, one has $E_2(a) < E_1(b)<E_2(c)$.
\end{enumerate}
\end{corollary}

Consequently, we have the following characterization of the ordinary degree map $\deg$.

\begin{corollary}
Let $E:  \RR^+ \longrightarrow  \mathcal{R}_{>0}$ an embedding of ordered abelian groups.  The following conditions are equivalent.
\begin{enumerate}
\item $\deg_E f = \deg f$ for all $f \in \mathcal{R}$.
\item  $\overline{\underline{\deg}}_E x = 1$.
\item  $\overline{\underline{\deg}}_E x^a = a$  for all $a \in \RR$.
\item $\deg_E x^a = a$  for all $a \in \RR$ (or, equivalently, for $a = \pm 1)$.
\item $\overline{\underline{\deg}}\, E(1) = 1$.
\item  $\overline{\underline{\deg}} \, E(a) = a$  for all $a \in \RR$.
\item $\deg E(a) = a$  for all $a \in \RR$ (or, equivalently, for $a = \pm 1)$.
\item For all real numbers $a<b<c$, one has $x^a < E(b) <x^c$.
\item For all real numbers $a<b<c$, one has $E(a) < x^b <E(c)$.
\end{enumerate}
\end{corollary}

 In particular, if $E$ is any embedding $\RR^+ \longrightarrow  \mathcal{R}_{>0}$  of ordered abelian groups, then $\deg = \deg_E: \mathcal{R} \longrightarrow \overline{\RR}$ if and only if $\deg_E x = 1$ and $\deg_E(1/x) = -1$.  Of course, the canonical such embedding acts by $a \longmapsto x^a$.

\begin{example} \
\begin{enumerate}
\item The map $E:   \RR^+ \longrightarrow \mathcal{R}_{>0}$ acting by $a \longmapsto (x \log x)^a$ (so that $E(1) = x \log x$) is an embedding of ordered abelian groups with $\deg_E f = \deg f$ for all $f \in \mathcal{R}$.   
\item  Let $A > 0$.   The map $E:   \RR^+ \longrightarrow \mathcal{R}_{>0}$ acting by $a \longmapsto (x^{1/A})^a = x^{a/A}$ (so that $E(1) = x^{1/A}$)  is an embedding of ordered abelian groups with $\deg_E f = A\deg f$ for all $f \in \mathcal{R}$.   
\end{enumerate}
\end{example}

All of the degree maps satisfying the equivalent conditions (8)--(14) of Corollary \ref{degEcor} or the equivalent conditions of Corollary \ref{degEcor2} also satisfy the equivalent conditions of the following proposition, whose proof is routine.

\begin{proposition}
Let  $E: \RR^+ \longrightarrow \mathcal{R}_{>0}$ be an embedding of ordered abelian groups. The following conditions are equivalent.
\begin{enumerate}
\item $1 \prec E(1)$.
\item $1 \prec E(c)$ for some $c \in \RR_{>0}$.
\item $1 \prec E(c)$ for all $c \in \RR_{>0}$.
\item $E(-1) \prec 1$.
\item  $E(c) \prec 1$ for some $c \in \RR_{<0}$.
\item $E(c) \prec 1$ for all $c \in \RR_{<0}$.
\item $E(c) \prec E(d)$ for all real numbers $c<d$.
\item The map $E|_\RR$ is an order isomorphism from $(\RR,<)$ to  $(E(\RR),\prec)$.
\item  $(E(\RR),\prec)$ is order isomorphic to $(\RR,<)$.
\item  $\deg_E f <0$ implies  $f \prec 1$, for all $f\in  \mathcal{R}$.
\item $\deg_E f  = \inf\{a \in \RR: f \prec E(a) \}$ for all $f \in  \mathcal{R}$.
\item $\deg_E 1 = 0$.
\item $\deg_E c = 0$ for all $c \in \RR^*$.
\item There is an $M > 0$ such that $\deg_E c \leq M$ for all $c \in \RR_{>0}$.
\item $\deg_E(cf)  = \deg_E f$ for all $f \in  \mathcal{R}$ and all  $c \in \RR^*$.
\item $\deg_E (fg)  = \deg_E f$ for all $f,g \in  \mathcal{R}$ with $g \asymp 1$.
\item $\deg_E f  = 0$ for all $f \in  \mathcal{R}$ with $f \asymp 1$.
\item $\deg_E f = \deg_E g$ for all $f,g \in  \mathcal{R}$ with $f \asymp g$.
\item $\deg_E f  \leq 0$ for all $f \in  \mathcal{R}$ with $f \preceq 1$.
\item $\deg_E f  \leq \deg_E g$ for all $f,g \in  \mathcal{R}$ with $f \preceq g$.
\item $\deg_E(f+g)  \leq  \max( \deg_E f,\deg_E g) \  (= \deg_E \max(|f|,|g|))$ for all $f,g \in  \mathcal{R}$.
\item The map $\deg_E: \mathcal{R} \longrightarrow \overline{\RR}$ factors (uniquely) through the valuation $v: \mathcal{R}  \longrightarrow \mathcal{R}/{ \asymp}$ on $\mathcal{R}$.
\end{enumerate}
\end{proposition}

\begin{example} Let $f \in  \mathcal{R}$,  and let  $E: \RR^+ \longrightarrow \mathcal{R}_{>0}$ be an embedding of ordered abelian groups.
\begin{enumerate}
\item Let $E(1) = \log x$.  Then
\begin{align*}
\deg_E f &  = \inf\{a \in \RR: f  \prec (\log x)^a \} \\
& = \deg f(e^x) \\
& = \begin{cases}  \infty &  \text{if } \deg f > 0 \\
  \dege_1 f &  \text{if } \deg f = 0 \\
  -\infty  &  \text{if } \deg f < 0.
\end{cases}
\end{align*}
 \item Let $E(1) = \log^{\circ N} x$, where $N$ is a nonnegative integer.
 Then
 \begin{align*}
\deg_E f &  = \inf\{a \in \RR: f  \prec (\log ^{\circ N})^a \} 
\\ & = \deg (f \circ \exp^{\circ N} x) \\
& = \begin{cases}  \infty &  \text{if } \dege f > (0,0,\ldots,  0,\infty,1,0,0,\ldots) \ \\
  \dege_N f &  \text{if } \deg_k f = 0 \text{ for all } k < N \\
  -\infty  &  \text{if } \deg f < (0,0,\ldots,  0,-\infty,-1,0,0,\ldots),
\end{cases}
\end{align*}
 where the $\pm \infty$ in the given vectors are preceded by exactly  $N$ $0$s.
 \item Let $E(1) = r \in \RR_{>1}$.
 Then
 \begin{align*}
\deg_E f &  = \inf\{a \in \RR: |f|  \leq r^a \}  \\
& =    \limsup_{x \to \infty} \log_r | f(x) | \\
& =   \log_r \limsup_{x \to \infty}| f(x) | 
\end{align*}
and therefore 
$$\deg_E f =  \deg x^{\log_r |f|} = \deg |f|^{\log_r x}$$
and
$$r^{\deg_E f} =  \limsup_{x \to \infty} | f(x) | = \deg x^{|f|}.$$
 \end{enumerate}
\end{example}

The examples above motivate the following.

\begin{proposition}
Let $E:  \RR^+ \longrightarrow  \mathcal{R}_{>0}$ an embedding of ordered abelian groups.    Suppose that
$E(1) \succ 1$ and $E(1)$ is eventually continuous and increasing.   Let $F$ be the (eventual) compositional inverse of $E(1)$.   The one has
$$\deg_E f = \deg (f \circ F)$$
for all $f \in \mathcal{R}$.
\end{proposition}

\section{Logexponentially bounded Hardy fields}

A function $f \in \RR^{\RR_\infty}$ is said to be {\bf transexponential}\index{transexponential}  if $f(x) \neq o(\exp^{\circ n}(x)) \ (x \to \infty)$ for every nonnegative integer $n$.    Of course, $f$ is transexponential if and only if $\dege f = (\infty,\infty,\infty,\ldots)$.  It is known that no functions in $\mathbb{H}$ are transexponential.  In fact, one has the following.

\begin{theorem}[{\cite[Theorem 13.2]{bos2}}]\label{Ebounded}
Let $r \in \mathbb{H}$.   There is an integer $k$ such that $r(x) = o(\exp^{\circ k}(x)) \ (x \to \infty)$.   Moreover,  there is a smallest such integer $k$ if and only if $r$ is unbounded, if and only if $\lim_{x \to \infty} r(x) = \pm \infty$.    Equivalently,  one has $\dege r < (\infty,\infty,\infty,\ldots)$, and $r$ is unbounded if and only if $\dege r > (0,0,0,\ldots)$.
\end{theorem}

Let us say that a Hardy field $K$ is {\bf logexponentially bounded}\index{logexponentially bounded Hardy field} if for every unbounded function $f \in K$ one has
$$(0,0,0,\ldots) < \dege f < (\infty,\infty,\infty,\ldots).$$   Theorem \ref{Ebounded} states, equivalently, that the Hardy field  $\mathbb{H}$ is logexponentially bounded.

\begin{proposition}\label{kcor}
Let $K$ be a logexponentially bounded Hardy field, and let $r,s \in K$ with $s \neq 0$.  
\begin{enumerate}
\item $r(x) = O(1) \ (x \to \infty)$ if and only if $\lim_{x \to \infty} r(x) \in \RR$ exists, if and only if $\dege r \leq (0,0,0,\ldots)$, if and only if $r(x) = o(\log^{\circ k}(x)) \ (x \to \infty)$ for all positive integers $k$.
\item $r(x) = o(1) \ (x \to \infty)$ if and only if $\dege r <(0,0,0,\ldots)$.
\item $r(x) \asymp 1 \ (x \to \infty)$ if and only if $\dege r = (0,0,0,\ldots)$,  if and only if $\lim_{x \to \infty} r(x) \in \RR^*$.
\item $r \preceq s$ if and only if $\lim_{x \to \infty} \frac{r(x)}{s(x)} \in \RR$ exists, if and only if $\dege \frac{r}{s} \leq (0,0,0,\ldots)$.
\item $r \prec s$ if and only if $\dege \frac{r}{s} <(0,0,0,\ldots)$.
\item $r \asymp s$ if and only if $\dege \frac{r}{s} = (0,0,0,\ldots)$.
\item $\dege s >  (-\infty,-\infty,-\infty,\ldots)$.
\end{enumerate}
\end{proposition}

Let $K$ be a logexponentially bounded Hardy field with canonical valuation $v: K \longrightarrow \Gamma_K \cup \{\infty\}$.  By the proposition, if $r,s \in K$ with $s \neq 0$, then one has $v(r) \geq v(s)$ if and only if $v(r/s) \geq 0$, if and only if $\dege \frac{r}{s} \leq (0,0,0,\ldots)$.   Likewise, one has $v(r) = v(s)$ if and only if $v(r/s) = 0$, if and only if $\dege \frac{r}{s} =(0,0,0,\ldots)$.  Thus,  $\Gamma_K$ is the group $K^*$ modulo the congruence relation $\asymp$,  and that relation can be defined in terms of $\dege$,  by $r \asymp s \,  \Longleftrightarrow \,  \dege \frac{r}{s} = (0,0,0,\ldots)$.

Let $S$ be a subset of $\mathcal{R}$.  We let $\prod_{n = 0}^{\infty S}\overline{\RR}$\index[symbols]{.g ha@$\prod_{n = 0}^{\infty S} \overline{\RR}$} denote the image in $\prod_{n = 0}^{\infty *}\overline{\RR}$ of $S$ under $\dege$.   We show later in this section that $\prod_{n = 0}^{\infty S}\overline{\RR} = \prod_{n = 0}^{\infty *}\overline{\RR}$ if $S$ contains $\mathcal{C}^{\infty}$,  and we complete the proof of Theorem \ref{degeequiv}.   Let $\mathcal{A}\index[symbols]{.i  kba@$\mathcal{A}$}$ denote the ring of germs of all functions that are analytic on a punctured neighbohood of $\infty$.   One has the commutative diagram below, which can be extended to contain all subrings of $\mathcal{R}$ and the various inclusions between them.   Note that $\prod_{n = 0}^{\infty S}\overline{\RR}$ as a totally ordered set is order-isomorphic to the set $S/{\equiv_{\dege}}$ of $\equiv_{\dege}$-equivalence classes $[f]$ of all elements $f$ of $S$, where $f \equiv_{\dege} g$ if $\dege f = \dege g$, and where $[f] \leq [g]$ if $\dege f \leq \dege g$.   If $S = K$ is a Hardy field, then $v(f) = v(g)$ implies $\dege f = \dege g$, and so the map onto $\prod_{n = 0}^{\infty K}\overline{\RR}$ factors through the map on $\Gamma_K \cup \{\infty\}$, and the canonical surjective map $\Gamma_K \cup \{\infty\} \longrightarrow \prod_{n = 0}^{\infty K}\overline{\RR}$ is order-reversing.  By Proposition \ref{kcor}, if $K$ is  a logexponentially bounded Hardy field, then $\Gamma_K$ is order-anti-isomorphic to $K^*/{\equiv}$, where $f \equiv g$ if $\dege(f/g) = (0,0,0,\ldots)$.

$$
\xymatrix{ 
\mathcal{R} \ar@{->>}[r]^(.4){\dege}     &{ \prod_{n = 0}^{\infty \mathcal{R}} \overline{\RR} } \ar@{=}[r]   & { \prod_{n = 0}^{\infty *} \overline{\RR}    }  \\
\mathcal{C} \ar@{->>}[rr] \ar@{>->}[u]    & & {\prod_{n = 0}^{\infty \mathcal{C}}\overline{\RR} } \ar@{=}[u]   \\
\mathcal{D} \ar@{->>}[rr] \ar@{>->}[u]    & & {\prod_{n = 0}^{\infty \mathcal{D}} \overline{\RR}}  \ar@{=}[u]   \\
\mathcal{C}^\infty \ar@{->>}[rr] \ar@{>->}[u]    & &{ \prod_{n = 0}^{\infty \mathcal{C}^\infty} \overline{\RR} } \ar@{=}[u]   \\
\mathcal{A} \ar@{->>}[rr] \ar@{>->}[u]    & &  {\prod_{n = 0}^{\infty \mathcal{A}} \overline{\RR}  } \ar@{>->}[u]   \\
\mathbb{H} \ar@{->>}[r]   \ar@{>->}[u] &\Gamma_{\mathbb{H}} \cup\{ \infty\} \ar@{->>}[r]  & { \prod_{n = 0}^{\infty \mathbb{H}} \overline{\RR}   }   \ar@{>->}[u] \\
\mathbb{L} \ar@{->>}[r]   \ar@{>->}[u] &\Gamma_{\mathbb{L}} \cup\{ \infty\} \ar@{->>}[r]  \ar@{>->}[u]  & {\prod_{n = 0}^{\infty \mathbb{L}} \overline{\RR} }     \ar@{>->}[u] \\
\RR(\mathfrak{L}) \ar@{->>}[r] \ar@{>->}[u] & \mathfrak{L} \cup \{ \infty\} \ar[r]^(.36){\cong} \ar@{>->}[u] & {\left(\bigoplus_{n = 0}^\infty\RR \right) \cup \{{\pmb{- \infty}}\}}  \ar@{>->}[u]  \\
\RR(x^a: a \in \RR) \ar@{->>}[r] \ar@{>->}[u] & {\id^\RR \cup \{ \infty\}} \ar[r]^{\cong} \ar@{>->}[u] & {\RR\cup \{-\infty\}}  \ar@{>->}[u]  
}$$

\begin{lemma}\label{abcd}
Let $K$ be a logexponentially bounded Hardy field  containing $\RR(x)$, and let $f \in K^*$.  Then the Hardy field $K(\log |f|)$ is logexponentially bounded.
\end{lemma}

\begin{proof}
Let $g \in K(\log |f|)$ be unbounded, so that there exist nonzero  polynomials $p$ and $q$ in $K[x]$ with $q(g) \neq 0$ such that $g = p(\log |f|)/q(\log |f|)$ and $q(\log |f||) \prec p(\log |f|)$.     Suppose first that $f \not \not \asymp 1$.  Replacing $f(x)$ with $|f(x)|$ or $\frac{1}{|f(x)|}$, we may assume without loss of generality that $\lim_{x \to \infty} f(x) = \infty$.    Then $h(x) = g(f^{-1}(e^x)) = p(x)/q(x) \in K^*$ is unbounded, so that $\log^{\circ n} x \prec h(x) \prec \exp^{\circ n} x$ for some $n$, whnce
$$\log^{\circ (n+1)} f(x) = \log^{\circ n} \log f(x) \prec  g(x) \prec \exp^{\circ n} \log f(x) =  \exp^{\circ (n-1)} f(x).$$
But one also has
$$\log^{\circ m} x \leq  f(x) \leq \exp^{\circ m} x$$ for some $m$, and therefore
$$\log^{\circ (n+m+1)}  x \prec  g(x) \prec   \exp^{\circ (n+m-1)}x.$$
Thus,  if $f \not \asymp 1$, then $K(\log |f|)$ is logexponentially bounded.
Since $\id \in K$,  it follows that $K(\log x)$ is logexponentially bounded.  Thus, we may assume without loss of generality that $\log x \in K$.   Finaly, if $f \asymp 1$, then we may replace $f$ with $x f \not \asymp 1$ and conclude that $K(\log |xf|) = K(\log |f|)$ is logexponentially bounded.  This completes the proof.
\end{proof}

\begin{lemma}\label{klem1}
Let $\dd \in \prod_{n = 0}^{ \infty*}\overline{\RR}$ and $n$ a nonnegative integer.   Let $K$ be a logexponentially bounded Hardy field containing $\RR(x^a: a \in \RR)$.  Then there does not  exist an $r \in K$ satistfying any of the following conditions.
\begin{enumerate}
\item $\dd_n = \infty$ and $\dege r = (\dd_0, \dd_1, \ldots, \dd_n,0,1,0,0,\ldots)$.
\item $\dd_n = -\infty$ and $\dege r = (\dd_0, \dd_1, \ldots, \dd_n,0,-1,0,0,\ldots)$.
\item $\dd_n$ is finite,  $\dd_n = \infty$, and $\dege r = (\dd_0, \dd_1, \ldots, \dd_{n+1},1,0,0,0,\ldots)$.
\item $\dd_n$ is finite,  $\dd_n = -\infty$, and $\dege r = (\dd_0, \dd_1, \ldots, \dd_{n+1},-1,0,0,0,\ldots)$.
\end{enumerate}
\end{lemma}

\begin{proof}
Note that $L \circ \exp$ and $L(\log|f|)$ are logexponentially bounded Hardy fields for any logexponentially bounded Hardy field $L$ and any $f \in L^*$.   It follows that there exists a sequence $K = K_0, K_1, K_2, \ldots$ of logexponentially bounded Hardy fields, where  $K_{n+1}$ is either $K_n \circ \exp$ or $K_n(\log|r_{(n-1)}|)$ for all $n$, such that $r_{(n)} \in K_n$ for all $n$.   Thus,   to prove (1) (or, more precisely, to prove condition (1) impossible),  it suffices to show that $\dege r \neq (\infty,0,1,0,0,\ldots)$ for all $r \in K$.   Let $r \in K$,  and suppose to obtain a contradiction that $\dege r = (\infty,0,1,0,0,0,\ldots)$.  Then  one has $r_{(2)}(x) = \log r(e^x)$ and $\dege  \frac{r_{(2)}(x)}{x} = (0,0,0,\ldots)$,  whence by Proposition \ref{kcor}(3) one has
$$\deg r = \lim_{x \to \infty} \frac {\log r(x)}{\log x} = \lim_{x \to \infty} \frac {\log r(e^x)}{x} \in \RR^*,$$
which contradicts $\deg r = \infty$.   The proof of (2) is similar.  To prove (3), it suffices to show that  $\dege r \neq (0,\infty,1,0,0,0,\ldots)$.   Suppose to obtain a contradiction that $\dege r = (0,\infty,1,0,0,0,\ldots)$.  Then,  again, one has $r_{(2)}(x) = \log r(e^x)$ and $\dege  \frac{r_{(2)}(x)}{x} = (0,0,0,\ldots)$,  whence, again by Proposition \ref{kcor}(3), one has $\deg \in \RR^*$,
which contradicts $\deg r = 0$.   The proof of (4) is similar.
\end{proof}

Let $$\prod_{n = 0}^{ \infty **}\overline{\RR} \subsetneq \prod_{n = 0}^{\infty *}\overline{\RR} \index[symbols]{.g hb@$\prod_{n = 0}^ { \infty **} \overline{\RR}$}$$  denote the {\bf restricted product}\index{restricted product $\prod_{n = 0}^{ \infty**}\overline{\RR}$} consisting of all sequences $\dd$ in $\prod_{n = 0}^{\infty}\overline{\RR}$ satisfying the following conditions for all nonnegative integers $n$.
\begin{enumerate}
\item If $\dd_n = \infty$, then $\dd > (\dd_0, \dd_1, \ldots, \dd_n, 0,1,0,0,0,\ldots).$
\item If $\dd_n = -\infty$, then $\dd < (\dd_0, \dd_1, \ldots, \dd_n, 0,-1,0,0,0,\ldots).$
\item If $\dd_n$ is finite and $\dd_{n+1} = \infty$, then $\dd <(\dd_0, \dd_1, \ldots, \dd_{n+1},1,0,0,0,\ldots).$
\item If $\dd_n$ is finite and $\dd_{n+1} = -\infty$, then $\dd > (\dd_0, \dd_1, \ldots, \dd_{n+1}, -1,0,0,0,\ldots).$
\item $\dd < (\infty,\infty,\infty, \ldots)$.
\item $\dd > (-\infty,-\infty,-\infty, \ldots)$.
\end{enumerate}
By Lemma \ref{klem1} and Propositions \ref{degeprop} and \ref{kcor}(7),  one has the following.

\begin{proposition}\label{kprop0}
For any logexponentially bounded Hardy field  $K$ containing $\RR(x^a: a \in \RR)$,  one has $\dege r \in \prod_{n = 0}^{ \infty **}\overline{\RR}$ for all $r\in  K^*$,  and therefore
$$\prod_{n = 0}^{ \infty K^*}\overline{\RR} \subseteq \prod_{n = 0}^{ \infty **}\overline{\RR}.$$
\end{proposition}

\begin{problem}
Is $K(x)$, and, more generally,  $K(x^a: a \in \RR)$,  logexponentially bounded for any logexponentially bounded Hardy field $K$?  If so, then the hypotheses that $K$ contain $\id$ or $\RR(x^a: a \in \RR)$ in Lemma \ref{abcd} through Proposition \ref{kprop0} can be removed.
\end{problem}

\begin{proposition}\label{kprop1}
Let  $K$ be any logexponentially bounded Hardy field containing $\id$ such that $f \circ \log \in K$ and $\exp \circ f \in K$ for all $f \in K$ (e.g., $K = \mathbb{L}$ and $K = \mathbb{H}$).   Let  $s \in K^*$ and $d \in \overline{\RR}$, and let $\dd = \dege s$.  Then there exists an $r \in K^*$ with  $\deg r = d$ and $r_{(1)} = s$
if and only if $(d, \dd_0,  \dd_1, \dd_2,\dd_3, \ldots) \in  \prod_{n = 0}^{ \infty **}\overline{\RR}$.    Moreover, if those equivalent conditions hold, then they are satisfied by the function
$$r(x) =  \begin{cases}  x^d |s(\log x)|  &  \text{if } d \neq \pm \infty \\
 e^{|s(x)|} &  \text{if } d = \infty \\
  e^{-1/|s(x)|}  &  \text{if } d = -\infty
\end{cases}$$
in $K_{> 0}$. 
\end{proposition}

\begin{proof}
Necessity of the condition follows from Proposition \ref{kprop0}, and sufficiency follows follows from the proof of  Proposition \ref{degerprop}, since Lemma \ref{klem1} and the condition  $(d, \dd_0,  \dd_1, \dd_2,\dd_3, \ldots) \in  \prod_{n = 0}^{ \infty **}\overline{\RR}$ rule out all of the cases in the proof where $r_{(1)} \neq s$.
\end{proof}

\begin{corollary}\label{kprop1cor}
For all $\dd \in \prod_{n= 0}^{\infty **}\overline{\RR}$ such that  $\dd_k = 0$ for all $k \gg 0$, one has $\dd = \dege r$ for some $r \in \mathbb{L}^*$.  Equivalently, one has
$$\left\{\dd \in \prod_{n= 0}^{\infty **}\overline{\RR}:  \dd_k = 0 \text{ for all } k \gg 0\right\}  \subseteq \prod_{n= 0}^{\infty  { \mathbb{L}}^*}\overline{\RR}.$$
\end{corollary}

\begin{proof}
This follows from Proposition \ref{kprop1} and an obvious induction.
\end{proof}

In the next section, we prove the converse of the corollary above,  i.e., we prove that the containment in the corollary is an equality.

\begin{proposition}
Let $r \in \mathbb{L}^*$.  Then $r$ has finite logexponential degree if and only if there exists a (unique) $s \in \mathfrak{L}$ such that $r \asymp s$ (or, equivalently, such that $\dege r = \dege s$), if and only if there is a (unique) $s \in \mathfrak{L}$ and a (unique) $c \in \RR^*$ such that $r \sim cs$.  In particular,  if $r$ has finite logexponential degree,  then $\dege_k r = 0$ for all $k \gg 0$.
\end{proposition}

\begin{proof}
If there exists an $s \in \mathfrak{L}$ such that $r \asymp s$, then clearly $s$ has finite logexponential degree, and $s$ is unique.   Suppose, conversely, that $r$ has finite logexponential degree.  The proof of the existence of $s$ is by induction on the order of $r$, where the {\it order of} $f \in \mathbb{L}$ is defined as in \cite{har3} \cite{har4}.   Note that $(\log^{\circ k} x)^t$ has order at most $k+2$ for all nonnegative integers $k$ and all $t \in \RR$.  By  \cite[Theorem 7]{har4}, if $r$ has order $0$, then $r(x) \asymp x^t$ for some $t \in \QQ$.  If $r$ has order at most $2$, then,  again by  \cite[Theorem 7]{har4}, one has $$r(x) x^{-\deg f} \asymp (\log x)^t$$ for some $t \in \QQ$.   If $r$ has order at most $3$, then, again by \cite[Theorem 7]{har4}, one has  $$r(x) x^{-\deg f} (\log x)^{-\dege_1 f} \asymp (\log \log x)^t$$
for some $t \in \QQ$.  This can be continued indefinitely.  
\end{proof}

Next,  we prove Theorem \ref{degeequiv}.

\begin{lemma}\label{cklem}
There exists a positive,  increasing,  infinitely differentiable function $f$ on $\RR_{>0}$ with $\lim_{x \to \infty} f(x) = \infty$ and
$$\dege f = (0,0,0,\ldots).$$
and any such function $f$ has exact logexponential degree.  Consequently,  for all $c \in \RR^*$, one has
$$\dege e^{cxf(x)} =  (\operatorname{sgn} c)(\infty,1,0,0,0,\ldots)$$
and
$$\dege x^{cf(x)} =  (\operatorname{sgn} c)(\infty,0,1,0,0,0,\ldots),$$
where $e^{cxf(x)}$ and $x^{cf(x)}$ are positive,  monotonic, and  infinitely differentiable  functions on $\RR_{>0}$ of exact logexponential degree.
\end{lemma}

\begin{proof}
The first statement follows from the proof of the italicized statement at the top of page 16 in  \cite{har3}, applied to the functions $\phi_n = \log^{\circ n}$,
and the last statement readily follows.
\end{proof}

\begin{lemma}\label{cklem2}
Let $\dd \in \prod_{n = 0}^{ \infty*}\overline{\RR}$ and $n$ a nonnegative integer.    There exists a positive,  monotonic,  infinitely differentiable function $f$ on $\RR_{>0}$ of exact logexponential degree satisfying any and exactly one of the following conditions.
\begin{enumerate}
\item $\dd_n = \infty$ and $\dege f = (\dd_0, \dd_1, \ldots, \dd_n,0,1,0,0,\ldots)$.
\item $\dd_n = -\infty$ and  $\dege f = (\dd_0, \dd_1, \ldots, \dd_n,0,-1,0,0,\ldots)$.
\item $\dd_n$ is finite,  $\dd_{n+1} = \infty$, and $\dege f = (\dd_0, \dd_1, \ldots, \dd_{n+1},1,0,0,0,\ldots)$.
\item $\dd_n$ is finite, $\dd_{n+1} = -\infty$, and $\dege f = (\dd_0, \dd_1, \ldots, \dd_{n+1},-1,0,0,0,\ldots)$.
\item $\dd_n$ and $\dd_{n+1}$ are finite and $\dege f = (\dd_0, \dd_1, \ldots, \dd_{n+1},0,0,0,0,\ldots)$.
\end{enumerate}
\end{lemma}

\begin{proof}
The proof by induction is similar to that of Proposition \ref{kprop1} and uses Lemma \ref{cklem}.
\end{proof}

\begin{proof}[Proof of Theorem \ref{degeequiv}]
Necessity follows from Proposition \ref{degeprop}.  To prove sufficiency,  let $\dd \in \prod_{n = 0}^{ \infty*}\overline{\RR}$.  We wish to find a positive, monotonic,  infinitely differentiable function $f$ on $\RR_{>0}$ such that $\dd = \dege f$.  We may suppose without loss of generality that $\dd \neq \dege r$ for all $r \in \mathbb{L}_{> 0}$, and that $\dd \neq \pm (\infty, \infty, \infty,\ldots)$.   Moreover, by Lemma \ref{cklem2}, we may suppose without loss of generality that $\dd$ is not eventually $0$.

For each nonnegative integer $n$, let $f_n$ be any positive,  monotonic,  infinitely differentiable function $f$ on $\RR_{>0}$ of exact logexponential degree  as in Lemma \ref{cklem2}.  In cases (1) and (4) of the lemma, one has $\dege f_n < \dege f$, and in cases (2) and (3) of the lemma, one has $\dege f_n > \dege f$.  Suppose that there are infinitely many $n$ for which (1) or (4) hold,  and infinitely many $n$ for which (2) or (3) hold.  We then have positive,  monotonic,  infinitely differentiable functions $\varphi_n$ and $\psi_n$ on $\RR_{>0}$ of exact logexponential degree  such that
$$\dege \varphi_1 < \dege \varphi_2 < \dege  \varphi_3 < \cdots <  \dd < \cdots < \dege  \psi_3 < \dege  \psi_2 < \dege  \psi_1$$
and such that both $\mathcal{S}(\dege \varphi_{n}, \dd)$ and $\mathcal{S}(\dege \psi_{n}, \dd)$ are monotonically increasing to $\infty$ as $n \to \infty$.   By Corollary \ref{oexp}, one has
$$ \varphi_1 \prec  \varphi_2 \prec   \varphi_3 \prec  \cdots \prec   \psi_3 \prec   \psi_2 \prec   \psi_1.$$
By  \cite[pp.\ 16--18]{har3},  then, there exists a positive,  monotonic,  infinitely differentiable function $f$ on $\RR_{>0}$  such that $$\varphi_n \prec f \prec \psi_n$$
for all $n$.  Thus one has
$$\dege \varphi_n \leq \dege f \leq \dege \psi_n$$
for all $n$.  Taking a sup and inf  in $\prod_{n = 0}^{ \infty*}\overline{\RR}$ over all $n$ and $m$, respectively, we see that $\dd \leq \dege f \leq \dd$, and therefore $\dege f = \dd$.  Moreover,  by the squeeze theorem, since the functions $ \varphi_n$ and $ \psi_n$ have exact logexponential degree,  the function $f$ does as well.

 Suppose, alternatively, that there are only finitely many $n$ for which (1) or (4) holds, or that there are only finitely many $n$ for which (2) or (3) holds.  Then $\dd_n$ is finite for all $n \gg 0$, and a similar argument shows that there are  sequences of positive,  monotonic,  infinitely differentiable functions   $\varphi_n$ and $\psi_n$ on $\RR_{>0}$ of exact logexponential degree such that  $$\dege \varphi_1 < \dege \varphi_2 < \dege  \varphi_3 < \cdots <  \dd < \cdots < \dege  \psi_3 < \dege  \psi_2 < \dege  \psi_1$$
and  such that  both $\mathcal{S}(\dege \varphi_{n}, \dd)$ and $\mathcal{S}(\dege \psi_{n}, \dd)$ are monotonically increasing to $\infty$,
  whence the conclusion again follows from \cite[pp.\ 16--18]{har3}.  
\end{proof}

\begin{corollary}
If $S$ is any subset of the ring $\mathcal{R}$ containing $\mathcal{C}^\infty$, then one has $$\prod_{n = 0}^{ \infty S}\overline{\RR} = \prod_{n = 0}^{ \infty *}\overline{\RR}.$$
\end{corollary}

The proof of the following  lemma is straightfoward (though slightly tedious, since there are several cases to check).

\begin{lemma}\label{kprop1lem}
The sets
$$\left\{\dd \in \prod_{n= 0}^{\infty **}\overline{\RR}:  \dd_k = 0 \text{ for all } k \gg 0\right\} \subseteq \prod_{n= 0}^{\infty   {\mathbb{L}}^*}\overline{\RR} $$
are dense in the poset $\prod_{n= 0}^{\infty *}\overline{\RR}$.    
\end{lemma}

From Corollary \ref{kprop1cor}, Lemma \ref{kprop1lem},  and Theorem \ref{degeequiv}, we immediately deduce the following.

\begin{proposition}\label{oexppropstrong}
For all $\dd, \ee \in \prod_{n= 0}^{\infty *}\overline{\RR}$ with $\dd < \ee$, there exists an $r \in \mathbb{L}_{>0}$ such that $\dd < \dege r < \ee$.  Equivalently, for all $f,  g\in \RR^{\RR_\infty}$, one has $\dege f < \dege g$ if and only if there exists an $r \in \mathbb{L}_{>0}$ with $\dege f < \dege r < \dege g$.
\end{proposition}

Thus, we may refine Proposition \ref{oexpprop}(2) as follows.

\begin{corollary}
Let $f,g \in \RR^{\RR_\infty}$,  where  $\dom g$ contains the intersection of $\dom f$ with some neighborhood of $\infty$. 
\begin{enumerate}
\item  If  $\dege f <\underline{\dege}\, g$, then  there exists an $r \in \mathbb{L}_{>0}$ such that $\dege f < \dege r < \underline{\dege}\, g$, and thus $f(x) = o(r(x)) \ (x \to \infty)$ and $r|_{\dom g}(x) = o(g(x)) \ (x \to \infty)$.
\item More generally,  one has $\dege f <\underline{\dege}\, g$ if and only if  there exist  $r,s \in \mathbb{L}_{>0}$ such that $\dege r <\dege s$,  $f(x) = o(r(x)) \ (x \to \infty)$,  and $s|_{\dom g}(x) = o(g(x)) \ (x \to \infty)$.
\end{enumerate}
\end{corollary}

The corollary above describes the  relation  $\dege f <\underline{\dege}\, g$ on functions $f,g \in \RR^{\RR_\infty}$,  where  $\dom g$ contains the intersection of $\dom f$ with some neighborhood of $\infty$,   in terms of the  $o$ relation and the relation $\dege r <\dege s$ on functions $r,s \in \mathbb{L}_{>0}$.  In the next section, we undertake an in depth study of the function $\dege$  on $ \mathbb{L}_{>0}$.

\section{Logexponential degree of transseries}

In this section, we assume familiarity with the ordered differential field $\mathbb{T} = \RR((x^{-1}))^{\operatorname{LE}}$\index[symbols]{.i  kc1@$\mathbb{T}$} of {\it  (well-based) (logarithmic-exponential) transseries},\index{transseries} as described in \cite[Appendix A]{asch} and \cite{dmm2}, for example.  The field $\mathbb{T}$ admits a composition $$\circ: \mathbb{T} \times \mathbb{T}_{>0,\succ 1} \longrightarrow  \mathbb{T}$$ and is
 an  ordered exponential field with exponential $ r \longmapsto \exp \circ r$ and logarithm $r \longmapsto \log \circ r$.  For any nonzero $r \in \mathbb{T}$, we can write 
 $$r = (1+\varepsilon)m,$$
 where $m$ is the leading transmonomial of $r$ and where $\varepsilon \prec 1$.   We may then define
 $$\lim r =    \left. \begin{cases}
    c & \text{if } m \sim c \in \RR^* \\
    0 & \text{if } m \prec 1 \\
    \infty & \text{if } m \succ 1 \text{ and } m >0 \\
       - \infty & \text{if } m \succ 1 \text{ and } m <0
\end{cases}
\right.$$  
and $\lim 0 = 0$.
Thus, we may also define the {\bf degree $\deg r \in \overline{\RR}$ of $r$} by
$$\deg r = \lim \frac{\log \circ |r|}{\log}$$
for all $r \neq 0$,  where we  also set $\deg 0 = -\infty$.
Let $r_{(0)} = r$, and define $r_{(k)} \in \mathbb{T}$ and $\dege_k r = \deg r_{(k)} \in \overline{\RR}$ for all nonnegative integers $k$, recursively, as follows.   Suppose that $r_{(k)}$ is defined, and set $d_k = \dege_k r = \deg r_{(k)}.$  We then let
$$r_{(k+1)}=   \left.
 \begin{cases}
    (x^{-d_k} r_{(k)} ) \circ \exp& \text{if } d_k \neq \pm \infty \\
     \log \circ  |r_{(k)}|  & \text{if } d_k = \infty \\
 \displaystyle   -\frac{1}{\log\circ |r_{(k)}|} & \text{if } d_k =- \infty.
 \end{cases}
\right.$$
We also set $$\dege r = (\dege_0 r, \dege_1 r, \dege_2 r, \ldots) \in \prod_{n = 0}^\infty \overline{\RR}.$$   We call $\dege r$ the {\bf logexponential degree of $r$}\index{logexponential degree of a transseries}\index[symbols]{.g g@$\dege f$} and $\dege_k r$ the {\bf logexponential degree of $r$ of order $k$}.  

Note that the Hardy field $\mathbb{L}$ of all (germs of) logarithmico-exponential functions is canonically isomorphic to a subfield of $\mathbb{T}$, and the $\dege$ map on $\mathbb{T}$ restricts to the $\dege$ map on $\mathbb{L}$ studied previously.

\begin{proposition}\label{Tseriesprop}
Let $r,s \in \mathbb{T}$, with leading transmonomials $m$ and $n$, respectively (where $m = 0$ if and only if $r= 0$, and likewise for $s$).   One has the following.
\begin{enumerate}
\item If $r \preceq s$, then $\dege r \leq \dege s$.
\item If  $\dege r < \dege s$, then $r \prec s$.
\item If $r \asymp s$, then $r_{(k)} \asymp s_{(k)}$ for all $k$ for all $k$ and $\dege r = \dege s$.
\item If $\dege r < (0,0,0,\ldots)$, then $\displaystyle \lim r = 0$.
\item $\dege (rs) = \dege r \oplus \dege s$ if $\dege_k s$ is finite for all $k$.
\item If $r \neq 0$, then $(1/r)_{(k)} = 1/r_{(k)}$ for all $k$ and $\dege (1/r) = -\dege r$.
\item $\dege r = \dege m$.
\item $\dege(r+s) \leq \max(\dege r,\dege s)$, with equality if $r \not \sim -s$ (which holds if and only if $m \not \sim -n$, if and only if $m \neq -n$).
\end{enumerate}
\end{proposition}

\begin{proof}
The  proofs of statements (1)--(6) are similar to the proofs of the corresponding statements concerning the $\dege$ maps on $\RR^{\RR_\infty}$ and $\mathbb{L}$.   Statement (7) follows from statements (3) and (5),  since $r = m(r/m)$ and $r/m \sim 1$.   Finally, statement (8) follows from statements (1), (2), and (7).
\end{proof}

For any $r \in  \mathbb{T}$ with $r > 0$ and $r \succ 1$,  there exist positive integers $m$ and $n$ such that $\log^{\circ m}  \prec r \prec \exp^{\circ n}$.  Equivalently, for any such $r$, one has
$$(0,0,0,\ldots) < \dege r < (\infty, \infty,  \infty, \ldots).$$
Consequently, one has the following analogue of Proposition \ref{kcor}.

\begin{proposition}\label{kcorLE}
Let $r,s \in \mathbb{T}$ with $s \neq 0$.  
\begin{enumerate}
\item $r\preceq 1$ if and only if $\lim r \in \RR$, if and only if $\dege r \leq (0,0,0,\ldots)$, if and only if $r \prec \log^{\circ k}$ for all positive integers $k$.
\item $r \prec 1$ if and only if $\dege r <(0,0,0,\ldots)$.
\item $r \asymp 1$ if and only if $\dege r = (0,0,0,\ldots)$,  if and only if $\lim \in \RR^*$.
\item $r \preceq s$ if and only if $\lim \frac{r}{s} \in \RR$, if and only if $\dege \frac{r}{s} \leq (0,0,0,\ldots)$.
\item $r \prec s$ if and only if $\dege \frac{r}{s} <(0,0,0,\ldots)$.
\item $r \asymp s$ if and only if $\dege \frac{r}{s} = (0,0,0,\ldots)$.
\end{enumerate}
\end{proposition}

For all $r \in \mathbb{T}$, let $\operatorname{lm} (r)$ denote the leading monic transmonomial of $r$.  Note that
$$\operatorname{lm} (rs) = \operatorname{lm} (r) \operatorname{lm} (s) $$
and
$$\operatorname{lm} (1/t) = 1/\operatorname{lm} (t)$$
for all $r,s,t \in \mathbb{T}$ with $t \neq 0$.   Let $m$ be a monic transmonomial,  so that $m =  (x^a  e^r) \circ \log^{(k)}$ for some purely large log-free transseries $r \in \mathbb{T}$, some $a \in \RR$, and some nonnegative integer $k$.   The minimal such $k$ is the {\bf  (logarithmic) depth of $m$}, \index{depth of a transseries} which is defined  more generally for an aribtrary $s \in \mathbb{T}$, and denoted $\operatorname{depth} (s)$, as the sup of the depths of all monomials in the support of $s$ \cite[Section 3]{edgar}.      Moreover,  the {\bf (exponential) height of $s$}, which we denote by $\operatorname{height} (s)$, is defined recursively so that $\operatorname{height} (m) = \operatorname{height} (r)+1$ when $k$ is minimal, and so that $ \operatorname{height} (t)$ is the sup of the heights of all monomials in the support of $t$ for any $t \in \mathbb{T}^*$ \cite[Section 3]{edgar}.   Note that,  if $m =  (x^a  e^r) \circ \log^{(k)}$ as above but $k$ is not minimal, then $\operatorname{depth} (m) < k$ and $\operatorname{height} (m) \leq  \operatorname{height} (r)+1$.
Note also that transseries are defined in such a way that every transseries has finite depth and height, which makes their study amenable to proofs by induction or finite descent.  

Let us define the {\bf (logexponential) breadth of} $r \in \mathbb{T}$, \index{breadth of a transseries} denoted $\operatorname{breadth} (r)$, by
$$\operatorname{breadth} (r)  = (\operatorname{height}(r), \operatorname{depth} (r)) \in \ZZ_{\geq 0} \times \ZZ_{\geq 0},$$
where $\ZZ_{\geq 0} \times \ZZ_{\geq 0}$ is endowed with the product order inherited from $\ZZ_{\geq 0}$.

Now, let $r \in \mathbb{T}$.   Let $\operatorname{lm}_k r = \operatorname{lm}( r_{(k)})$ for all $k \in \ZZ_{\geq 0}$.
By Proposition \ref{kcorLE}(7), one has $\dege_k r = \deg \operatorname{lm}_k r$  for all $k$.  Note that $\operatorname{lm}_0(r) = \operatorname{lm}(r)$ and 
$$\operatorname{lm}_{k+1} r=   \left.
 \begin{cases}
    (x^{-d_k} \operatorname{lm}_k r ) \circ \exp& \text{if } d_k \neq \pm \infty \\
    \operatorname{lm}( \log \circ  \operatorname{lm}_k r) & \text{if } d_k = \infty \\
 \displaystyle  \operatorname{lm}\left(-\frac{1}{\log\circ \operatorname{lm}_k r}\right) =   -\frac{1}{\operatorname{lm}(\log\circ \operatorname{lm}_k r)} & \text{if } d_k =- \infty,
 \end{cases}
\right.$$
 where $d_k = \dege_k r =  \deg \operatorname{lm}_k r$, for all $k$.
 
The following theorem characterizes all sequences of the form $\dege r$ for some $r \in \mathbb{T}^*$ (resp.,  $r \in \mathbb{L}^*$).

\begin{theorem}\label{chartheorem2}
For all $\dd \in \prod_{n = 0}^{ \infty}\overline{\RR}$,  the following conditions are equivalent.
\begin{enumerate}
\item $\dd = \dege r$ for some $r \in  \mathbb{L}^*$.
\item $\dd = \dege r$ for some $r \in  \mathbb{T}^*$.
\item $\dd \in \prod_{n = 0}^{ \infty **}\overline{\RR}$ and $\dd_k = 0$ for all $k \gg 0$.
\end{enumerate}
Moreover,  for all $r \in  \mathbb{T}^*$,  one has $\dege_k r = 0$ for all $k >  \operatorname{height} (\operatorname{lm}(r))+ \operatorname{depth} (\operatorname{lm}(r))$.  More generally,  for all $r \in  \mathbb{T}^*$ and all nonnegative integers $k$, one has $$\operatorname{breadth} (\operatorname{lm} (r_{(k+1)}))  < \operatorname{breadth} (\operatorname{lm}( r_{(k)}))$$ if $r_{(k)} \not \asymp 1$, whence there exists a $k \leq \operatorname{height} (\operatorname{lm}(r))+ \operatorname{depth} (\operatorname{lm}(r))+1$ such that $r_{(n)} \asymp 1$,  and therefore $\dege_n k = 0$, for all $n \geq k$.
\end{theorem}

\begin{proof}
Since $\mathbb{L}$ is isomorphic to a canonical subfield of $\mathbb{T}$ and the $\dege$ map is preserved under this identification,  condition (1) implies condition  (2).  Moreover,  condition (3) implies condition (1) by Corollary \ref{kprop1cor}.  The proof that $r \in \mathbb{T}^*$ implies $\dege r \in \prod_{n = 0}^{ \infty **}\overline{\RR}$ is a straightforward generalization of the proof Proposition \ref{kprop0}.  It remains only to prove the claim in the theorem that $$\operatorname{breadth} (\operatorname{lm} (r_{(k+1)}))  < \operatorname{breadth} (\operatorname{lm}( r_{(k)}))$$ for all $k$ with $r_{(k)} \not \asymp 1$.

 Let $m$ be any monic transmonomial,  so that $m =  (x^a  e^r) \circ \log^{(k)}$ for some purely large log-free transseries $r \in \mathbb{T}$ and some $a \in \RR$,  where $k$ is the depth of $m$.  
Suppose first that $\deg m = \pm \infty$,  or, equivalently,  that $\log \circ m \succ \log$.  Then $\operatorname{height} (m) \geq 1$.    Moreover, one has $$a\log x \circ \log^{(k)}+ r \circ \log^{(k)} = \log \circ m \succ \log \circ \log^{(k)},$$ and therefore $$\operatorname{lm}(\log \circ m) = \operatorname{lm}(r) \circ \log^{(k)}.$$  It follows that
$$\operatorname{height}(-1/\operatorname{lm} (\log \circ m)) = \operatorname{height}(\operatorname{lm} (\log \circ m)) \leq \operatorname{height} (r)  =  \operatorname{height} (m)-1$$
and
$$\operatorname{depth}(-1/\operatorname{lm} (\log \circ m)) = \operatorname{depth} (\operatorname{lm} (\log \circ m)) \leq k  =  \operatorname{depth} (m).$$

Suppose, instead, that $\deg m = d \in \RR$.   If $k = 0$, then $r = 0$ and $m = x^d$.   Suppose, on the other hand, that $k > 0$.  Then $\deg (x^a  \circ \log^{(k)}) = 0$ and thus $\deg (e^r \circ \log^{(k)}) = d$.   If $d \neq 0$, then $r \circ \log^{(k)} \sim d \log x$, whence $r \sim d\exp^{(k-1)}$ and thus $r$ has leading transmonomial $d\exp^{(k-1)}$, say, $r= d\exp^{(k-1)}+s$ for some purely large $s \prec \exp^{(k-1)}$,  so that
$$e^r \circ \log^{(k)} = e^{d\exp^{(k-1)}+s}\circ \log^{(k)} = x^d(e^s \circ \log^{(k)}).$$
Thus, if we let, $m' = x^{-d} m$,  then in each case (namely, $k = 0$, $k >0$ and $d \neq 0$, and $k>0$ and $d = 0$), 
the transseries $m'$ is a monic transmonomial of degree $0$ of the form  $e^t \circ \log^{(k)}$ for some purely large log-free transseries $t$, where $t = 0$, $t = r$, or $t$ is $r$ minus its leading transmonomial.  Thus,  $m'$ has depth and height at most the depth and height, respectively, of $m$.   If $k > 0$, then $m' \circ \exp =  e^t \circ \log^{(k-1)}$ is a monic transmonomial, and one has
$$\operatorname{height}(m' \circ \exp) \leq  \operatorname{height} (t) +1 \leq \operatorname{height} (r) +1 = \operatorname{height} (m)$$
and
$$\operatorname{depth}(m' \circ \exp) \leq k-1 =  \operatorname{depth} (m)-1.$$
On the other hand, if $k = 0$, then $m' = 1$.
This proves our claim and completes the proof.
\end{proof}

\begin{corollary}\label{Ldege}
One has
$$\dege \mathbb{L}^* = \prod_{n= 0}^{\infty   \mathbb{L}^*}\overline{\RR} = \left\{\dd \in \prod_{n= 0}^{\infty **}\overline{\RR}:  \dd_k = 0 \text{ for all } k \gg 0\right\}.$$
\end{corollary}

\begin{example}\
\begin{enumerate}
\item For the transseries $r = e^{e^{(\log \log x)^2}}$, the sequence  $\operatorname{lm}( r_{(k)})$ is $$r, e^{(\log \log x)^2}, e^{(\log x)^2}, (\log x)^2, x^2, 1,1,1\ldots,$$ and the sequence $\operatorname{breadth} (\operatorname{lm}( r_{(k)}))$ is $$(2,2),(1,2),(1,1),(0,1),(0,0),(0,0),(0,0),\ldots.$$
\item  For the transseries $r = e^{\log x \log \log \log x}$, the sequence  $\operatorname{lm}( r_{(k)})$ is $$r, \log x \log \log \log x, x \log \log x,  \log x, x,1,1,1\ldots,$$ and the sequence $\operatorname{breadth} (\operatorname{lm}( r_{(k)}))$ is $$(3,3),(2,3),(2,2),(0,1),(0,0),(0,0),(0,0),\ldots.$$
\item  For the transseries $r = e^{(\log x)^2+( \log \log \log x)^2}$, the sequence  $\operatorname{lm}( r_{(k)})$ is $$r, (\log x)^2,  x^2,1,1,1\ldots,$$ and the sequence $\operatorname{breadth} (\operatorname{lm}( r_{(k)}))$ is $$(3,3),(0,1),(0,0),(0,0),(0,0),\ldots.$$
\end{enumerate}
\end{example}

Theorem \ref{chartheorem2} and Proposition \ref{kprop1} yield the following way to associate  to any $\dd \in \prod_{n = 0}^{ \infty \mathbb{L}}\overline{\RR}$ a canonical  $\dd \, { \uparrow}  \, {\mathbb{L}} \in \mathbb{L}_{>0}$ such that $\dege ( \dd \, {\uparrow} \,  {\mathbb{L}} )= \dd$, or, equivalently,   to any $r \in \mathbb{T}$ a canonical $r\, {\downarrow} \, {\mathbb{L}}  = (\dege r) \, {\uparrow} \, {\mathbb{L}}\in \mathbb{L}_{>0}$  such that $\dege (r \, {\downarrow}\, {\mathbb{L}} )  = \dege r$, where both $\dd \, {\uparrow}\, {\mathbb{L}}$ and $r \, {\downarrow} \, {\mathbb{L}}$ are, in a precise sense, the simplest possible.

\begin{corollary}\label{prototype}
Let $\dd , \ee \in \prod_{n = 0}^{ \infty \mathbb{L}}\overline{\RR}$.    By Theorem \ref{degeprop},  there exists a smallest $N$ such that $\dd_k= 0$ for all $k > N$, and one must then have $\dd_N \in \RR$.   Let $f_N = x^{\dd_N}$.   Assuming that $f_k \in \mathbb{L}_{>0}$ is defined, where $0 < k \leq N$,  let
$$f_{k-1}(x) =  \begin{cases}  x^{\dd_{k-1}} f_k(\log x)  &  \text{if } \dd_{k-1} \neq \pm \infty \\
 e^{f_k(x)} &  \text{if } \dd_{k-1} = \infty \\
  e^{-1/f_k(x)}  &  \text{if } \dd_{k-1} = -\infty
\end{cases}$$
in $\mathbb{L}_{> 0}$.   This defines $f_k \in \mathbb{L}_{> 0}$ for all integers $k$ with $0 \leq k \leq N$.    Let $  \dd \, {\uparrow} \,  {\mathbb{L}} = f_0 \in \mathbb{L}_{> 0}$.     One has the following.
\begin{enumerate}
\item  $(\dd \, {\uparrow} \,  {\mathbb{L}})_{(k)} = f_k$ for all integers $k$ with $0 \leq k \leq N$.
\item $\dege (\dd \, {\uparrow} \,  {\mathbb{L}}) = \dd$.   
\item $\dd \, {\uparrow} \,  {\mathbb{L}} \preceq \ee \, {\uparrow} \,  {\mathbb{L}}$ if and only if $\dd \leq \ee$.
\item $\dd \, {\uparrow} \,  {\mathbb{L}} \prec \ee \, {\uparrow} \,  {\mathbb{L}}$ if and only if $\dd < \ee$.
\item $\dd_k$ is finite for all $k$ if and only if $\dd \in \bigoplus_{n = 0}^{ \infty}\RR$, if and only if $\dd \, {\uparrow} \,  {\mathbb{L}} \in \mathfrak{L}$.
\item The image of $\dd \, {\uparrow} \,  {\mathbb{L}}$ in $\mathbb{T}$ is the unique monic transmonomial $m$ such that $\dege m = \dd$ and $m_{(k)}$ is a monic transmonomial for all $k$.
\end{enumerate}
\end{corollary}

Note that statements (2)--(5) say, equivalently, that the map  $\dd \longmapsto \dd \, {\uparrow} \,  {\mathbb{L}}$ from the poset $\prod_{n = 0}^{ \infty \mathbb{L}}\overline{\RR}$ to the set $\mathbb{L}_{>0}$ preordered by the $O$ relation $\preceq$  is a preorder embedding that is a section of the preorder morphism $\dege:\mathbb{L}_{>0} \longrightarrow \prod_{n = 0}^{ \infty \mathbb{L}}\overline{\RR}$,  and it maps $\bigoplus_{n = 0}^{ \infty}\RR$ onto $\mathfrak{L}$.

\begin{example}\label{prototypeexample}
One has
$$(1,-\infty,-\tfrac{3}{5},\tfrac{1}{5},0,0,0,\ldots) \, {\uparrow} \,  {\mathbb{L}} = xe^{ -(\log x)^{3/5}(\log \log x)^{-1/5}}.$$
Indeed, the sequence $f_3, f_2, f_1, f_0$ as defined in the corollary is
$$x^{1/5}, x^{-3/5}(\log x)^{1/5},e^{-x^{3/5}(\log x)^{-1/5}},xe^{-(\log x)^{3/5}(\log \log x)^{-1/5}}.$$
\end{example}

Next, let us say that a Hardy field $K$ is {\bf transerially logexponential}\index{transerially logexponential Hardy field} if there exists a field embedding $\varphi: K \longrightarrow \mathbb{T}$ such that $\dege f = \dege \varphi(f)$ for all $f \in K$.   For example,  the Hardy field $\mathbb{L}$ of all logarithmico-exponential functions is transerially logexponential, and any subfield of a transerially logexponential Hardy field is also a transerially logexponential Hardy field.  Theorem \ref{chartheorem2} also has the following immediate corollary.

\begin{corollary}
Let $K$ be a transerially logexponential Hardy field.   For all $r \in  K^*$,  one has  $\dege r \in \prod_{n = 0}^{ \infty **}\overline{\RR}$ and $\dege_k r = 0$  for all $k \gg 0$, and therefore
$$\prod_{n= 0}^{\infty  K^*}\overline{\RR}  \subseteq \prod_{n= 0}^{\infty   \mathbb{L}^*}\overline{\RR}.$$ In particular, $K$ is logexponentially bounded.
\end{corollary}

The transerially logexponential condition on Hardy fields is a much stricter finiteness condition than is logexponential boundedness.   Nevertheless, the following proposition provides a large class of  transerially logexponential Hardy fields.

\begin{proposition}\label{lehard}
Let $K$ be a Hardy field containing $\RR(\mathfrak{L})$ such that $\log f \in K$ for all $f \in K_{>0}$ and such that $K$ admits an ordered differential field embedding $\varphi: K \longrightarrow \mathbb{T}$.  Then $\dege f = \dege \varphi(f)$ for all $f \in K$, and therefore $K$ is a transerially logexponential.
\end{proposition}

\begin{proof}
Note that $K' = K \circ \exp$ is a Hardy field  containing $\RR(\mathfrak{L})$, one has $\log f \in K'$ for all $f \in K'_{>0}$, and $K'$ admits an ordered differential field embedding $\varphi': K' \longrightarrow \mathbb{T}$ acting by $f \circ \exp \longmapsto \varphi(f) \circ \exp$.   Let $f \in K_{>0}$.  Let $K_0 = K$ and $\varphi_0 = \varphi$, and, assuming $K_n$ and $\varphi_n: K_n \longrightarrow \mathbb{T}$ are defined,  let $K_{n+1} = K_n$ and $\varphi_{n+1} = \varphi_n$ if $\dege_n f = \pm \infty$, and let $K_{n+1} = K_n' = K_n \circ \exp$ and $\varphi_{n+1} = \varphi_{n}': K_{n+1} \longrightarrow \mathbb{T}$, acting by $f \circ \exp \longmapsto \varphi(f) \circ \exp$, if $\dege_n f \in \RR$.   It follows that each $K_n$ and $\varphi_n$ satisfies the same conditions in the proposition as $K$ and $\varphi$, and, for any $f \in K$, one has $f_{(n)} \in K_n$ for all $n$.  

Since $\varphi$ is the identity on $\QQ$ and respects the relation $\leq$, it must be the identity on $\RR$.
Since $\varphi$ is the identity on $\RR$ and respects the relation $\leq$, it also respects limits.
Let $f \in K_{>0}$.   Then one has $\log f \in K$ and $$(\varphi(\log f))' = \varphi ((\log f)') = \varphi(f'/f) = \varphi(f)'/\varphi(f) = (\log \varphi(f))',$$
whence $$ \log \varphi(f) = \varphi(\log f) +C$$ 
for some $C \in \RR$.    By a similar argument,  for any $a \in \RR$ one has 
$$\varphi(x^a) =  Ax^a \asymp x^a$$
for some $A \in \RR$ with $A >0$.    
Moreover, one has
$$\deg f = \lim_{x \to \infty} \frac{\log f(x)}{\log x} =  \lim \frac{\varphi(\log f)}{\varphi(\log x)} = \lim \frac{\log \varphi(f)-C}{\log x-D} = \lim \frac{\log \varphi(f)}{\log x} = \deg \varphi(f)$$ for some $D \in \RR$,  so that  $\varphi$ is degree preserving.   
Furthermore, if $\lim_{x \to \infty} f(x) \in \{0, \infty\}$ (e.g.,  if $\deg f = \pm \infty$), then
$$ \log \varphi(f) = \varphi(\log f) +C \asymp \varphi(\log f)$$
and
$$-\frac{1}{\log \varphi(f)}  = -\frac{1}{\varphi(\log f)+C}  \asymp -\frac{1}{\varphi(\log f)} = \varphi\left(-\frac{1}{\log f}\right).$$ 
It follows from the facts noted above that
$$\varphi_1(f_{(1)} ) = \varphi'(f_{(1)} ) \asymp \varphi(f)_{(1)},$$
whether  $\deg f \neq \pm \infty$,  $\deg f = \infty$,  or $\deg f = -\infty$.  Applying the same conclusion  to $\varphi_1$ and $f_{(1)}$, we see that
$$\varphi_2(f_{(2)} ) = \varphi_{1}'((f_{(1)} )_{(1)}) \asymp \varphi_1(f_{(1)})_{(1)} \asymp (\varphi(f)_{(1)})_{(1)} = \varphi(f)_{(2)}.$$
 It follows, then, by induction  on $n$, that 
$$\varphi_n(f_{(n)} ) \asymp \varphi(f)_{(n)},$$
and therefore
$$ \dege_n f  =  \deg f_{(n)}  = \deg \varphi_n(f_{(n)} ) = \deg \varphi(f)_{(n)} = \dege_n \varphi(f),$$
for all $n$.
\end{proof}

Thus, for any Hardy field $K$ satisfying the hypotheses of Proposition \ref{lehard}, the map $\dege$ on $K$ is completely determined by the map $\dege$ on $\mathbb{T}$, and one has
$$\prod_{n= 0}^{\infty  K^*}\overline{\RR}  \subseteq \prod_{n= 0}^{\infty   \mathbb{L}^*}\overline{\RR}.$$

\begin{example} \
\begin{enumerate}
\item The Hardy field $K = \mathbb{L}$ satsfies the hypotheses of Proposition \ref{lehard} and is therefore transerially logexponential.  
\item  More generally,  the well-known Hardy field $H(\RR_{\operatorname{an},\operatorname{exp}}) \supsetneq \mathbb{L}$  of \cite[Section 5]{dmm0} \cite{dmm2} satisfies the conditions  of Proposition \ref{lehard}, and therefore any subfield of $H(\RR_{\operatorname{an},\operatorname{exp}})$  is a transerially logexponential Hardy field.
Indeed, by \cite[(3.11), (3.12), (6.29), (6.30)]{dmm2},  the Hardy field  $H(\RR_{\operatorname{an},\operatorname{exp}})$ is real closed,  it is closed under composition  and compositional inverses (i.e., $f \circ g, g^{-1} \in H(\RR_{\operatorname{an},\operatorname{exp}})$ for all $f, g \in H(\RR_{\operatorname{an},\operatorname{exp}})$ with $g>0$ and $g \succ 1$),  and it admits a natural ordered differential field embedding in $\mathbb{T}$ that sends $x$ to $x$, $\log x$ to $\log  x$, and $\exp x$ to $\exp x$, and that respects composition and compositional inverses.   
\item Any {\it transserial Hardy field} $K$ \cite{vdh} containing $\RR(\mathfrak{L})$ such that $\log f \in K$ for all $f \in K_{>0}$  (e.g., $K = \mathbb{L}$)   also satisfies the conditions of Proposition \ref{lehard} and is therefore transerially logexponential.  
\item The smallest Hardy field $K$  satisfying the hypotheses of Proposition \ref{lehard} is the smallest Hardy field  $\mathbb{L}_0$ containing $\RR(\mathfrak{L})$ that is closed under precomposition of positive elements by $\log$, and one has
$$\prod_{n= 0}^{\infty  \mathbb{L}_0^*}\overline{\RR}= \prod_{n= 0}^{\infty  \RR(\mathfrak{L})^*}\overline{\RR}  =  \bigoplus_{n= 0}^{\infty}\overline{\RR}.$$
\end{enumerate}
\end{example}

\begin{problem} \ 
\begin{enumerate}
\item Is the Hardy field $\mathbb{H}$  transerially logexponential?  
\item  Does $\mathbb{H}$  admit an ordered differential field embedding in $\mathbb{T}$? 
\item  Compute $\prod_{n = 0}^{\infty \mathbb{H}} \overline{\RR}$.  Does it properly  contain $\prod_{n = 0}^{\infty \mathbb{L}} \overline{\RR}$?
\end{enumerate}
\end{problem}

\begin{problem}
Axiomatize the logexponential degree map $\operatorname{ledeg}$ as a canonical map from Hardy’s field $\mathbb{L}$ of logarithmico-exponential functions and from the logarithmic-exponential transseries field $\mathbb{T}$ into $\prod_{n=0}^\infty \overline{\mathbb{R}}$, satisfying a minimal set of natural properties, such as compatibility with ordering, asymptotic equivalence, addition, multiplication, and composition.  In particular, determine whether $\operatorname{ledeg}$ can be characterized intrinsically, without reference to the recursive normalization procedure.
\end{problem}

\chapter{Asymptotic continued fraction expansions}

In this chapter, we explore the role of continued fractions as tools for algebraic asymptotic analysis. We define an \textit{asymptotic continued fraction expansion} of a complex-valued function to be a continued fraction whose approximants provide an asymptotic expansion of the function.   Building on \cite{ell0} and \cite{ell}, we establish general properties of such expansions.  We show, for example, that the ``best'' rational function approximations of a function with an asymptotic Jacobi continued fraction expansion are exactly the approximants of the continued fraction.   An elegant application yields asymptotic continued fraction expansions of the  function $\PP(e^x) = \pi(e^x)/e^x$,  induced from  two classical continued fraction expansions of the exponential integral function $-E_n(-z)$, precisely from the region $[0,\infty)$ on which the exansions diverge.   From this, we determine the best rational approximants of $\PP(e^x)$, whose denominators, notably, are the monic Laguerre polynomials rather than powers of $x$.  This chapter is largely self-contained and is not referenced in later chapters, except for Proposition~\ref{mertensprop}.

\section{General asymptotic continued fraction expansions}

Let $a \in\overline{\RR}$ (resp.,  $a \in \CC \cup \{\infty\}$).  Let $f, b_0, b_1, b_2, \ldots$ and $a_1, a_2, \ldots$ be real (resp., complex) functions, each of which has its domain containing $a$ as a limit point.   Consider the formal continued fraction
\begin{align*}
b_0(z) + \frac{a_1(z)}{b_1(z) \, +} \ \frac{a_2(z)}{b_2(z)\, +} \ \frac{a_3(z)}{b_3(z)\, +}  \ \cdots \ = \ b_0(z) +\cfrac{a_1(z)}{b_1(z)+\cfrac{a_2(z)}{b_2(z)+\cfrac{a_3(z)}{b_3(z)+\cdots}}}.
\end{align*}
Let 
\begin{align*}
w_n(z) =  b_0(z) + \frac{a_1(z)}{b_1(z)\, +} \ \frac{a_2(z)}{b_2(z)\, +} \ \cdots \ \frac{a_n(z)}{b_n(z)} \  = \ b_0(z) +\cfrac{a_1(z)}{b_1(z)+\cfrac{a_2(z)}{b_2(z)+\cdots +\cfrac{a_{n-1}(z)}{b_{n-1}(z)+\cfrac{a_{n}(z)}{b_{n}(z)}}}}
\end{align*} for all $n \geq 0$ denote the {\bf $n$th approximant}\index{approximants of a continued fraction} of the given formal continued fraction, where $w_0(z) = b_0(z)$, and where the domain of $w_n(z)$ is the largest possible so that the expression above is defined.  We write 
\begin{align}\label{contlega} 
f(z) \, \simeq \, b_0(z) + \frac{a_1(z)}{b_1(z) \,+} \ \frac{a_2(z)}{b_2(z)\,+} \ \frac{a_3(z)}{b_3(z)\,+} \ \cdots  \ (x \to a)
\end{align}
if  $\{w_n(z)-w_{n-1}(z)\}_{n = 1}^\infty$ is an asymptotic sequence at $a$  and 
$$f(z) -w_0(z) \simeq \sum_{n = 1}^\infty (w_n(z)-w_{n-1}(z))  \ (x \to a)$$
is an  asymptotic expansion of $f(z)-w_0(z)$ at $x = a$ with respect to $\{w_n(z)-w_{n-1}(z)\}$,  where $\dom w_0$ (and consequently $\dom w_n$ for all $n$)  contains the intersection of $\dom f$ with some punctured neighborhood of $a$.   In that case we say that (\ref{contlega}) is an {\bf asymptotic continued fraction expansion of $f$  at $a$.}\index{asymptotic continued fraction expansion}   Thus, (\ref{contlega}) is an asymptotic continued fraction expansion of $f$  at $a$ if and only if 
\begin{align}\label{contas1}
w_{n+1}(z)-w_{n}(z) = o(w_n(z)-w_{n-1}(z)) \ (x \to a)
\end{align}
and
\begin{align}\label{contas2}
f(z) - w_{n}(z) = o( w_{n}(z)-w_{n-1}(z)) \ (x \to a)
\end{align}
for all $n \geq 1$ and $\dom w_n$ for each $n$ contains the intersection of $\dom f$ with some punctured neighborhood of $a$.

\begin{remark}[Remarks on asymptotic continued fraction expansions]\label{crem} With the notation as above,  suppose that $\dom w_n$ for each $n$ contains the intersection of $\dom f$ with some punctured neighborhood of $a$.
\begin{enumerate}
\item By Remark \ref{asrem}, condition (\ref{contas2})  can be replaced with
\begin{align}\label{contas3}
f(z) - w_{n}(z)  = O( w_{n+1}(z)-w_{n}(z)) \ (x \to a)
\end{align}
and with
\begin{align}\label{contas4}
f(z) - w_{n}(z) \sim w_{n+1}(z)-w_{n}(z) \ (x \to a).
\end{align}
More precisely, the asymptotic continued fraction  expansion (\ref{contlega}) holds if and only if (\ref{contas1}) holds for all $n \geq 1$ and any one of the three conditions (\ref{contas2}), (\ref{contas3}), (\ref{contas4}) holds for infinitely many $n \geq 1$, in which case all three of the conditions hold for all  $n \geq 1$.
\item $\{w_{n}(z)-w_{n-1}(z)\}$ is an asymptotic sequence  at $a$ if and only if $w_{l}(z)-w_n(z) \sim w_m(z)-w_n(z) \ (x \to a)$ for all $n,m,l  \geq 0$ with $m > n$ and $l > n$, or equivalently with $m = n+1$ and $l = n+2$.
\item The asymptotic continued fraction expansion (\ref{contlega}) holds if and only if 
$f(z) - w_n(z) \sim w_{m}(z)-w_n(z)  \ (x \to a)$
 for all $n,m  \geq 0$ with $m > n$, or equivalently with $n < m \leq n+2$.
\item  If the asymptotic continued fraction expansion (\ref{contlega}) holds, then
$$f(z) \sim b_0(z) \ (x \to a) \ \Longleftrightarrow\  w_1(z) \sim b_0(z) \ (x \to a) \ \Longleftrightarrow\ f(z) \sim w_n(z) \ (x \to a) \text{ for all } n.$$
\item If $f(z)-b_0(z)$ is nonzero in the intersection of $\dom f$ with some  punctured neighborhood of $a$, then the asymptotic continued fraction  expansion (\ref{contlega}) holds if and only if $\frac{a_1(z)}{f(z) -b_0(z)}  \sim  b_1(z) \ (x \to \infty)$ and one has the asymptotic  continued fraction expansion
\begin{align*}
\frac{a_1(z)}{f(z) -b_0(z)} \, \simeq \,  b_1(z)+ \ \frac{a_2(z)}{b_2(z)\,+} \ \frac{a_3(z)}{b_3(z)\,+} \ \cdots  \ (x \to a).
\end{align*}
\end{enumerate}
\end{remark}

For ease of notation one makes the identification
\begin{align*}
b_0(z) \pm \frac{a_1(z)}{b_1(z) \, -} \ \frac{a_2(z)}{b_2(z)\, -} \ \frac{a_3(z)}{b_3(z)\, -}  \ \cdots \ = \
b_0(z)+ \frac{\pm a_1(z)}{b_1(z) \, +} \ \frac{-a_2(z)}{b_2(z)\, +} \ \frac{-a_3(z)}{b_3(z)\, +}  \ \cdots,
\end{align*}
and one identifies two formal continued fractions if they possess identical approximants.

A {\bf normed ring}\index{normed ring} (resp.,  {\bf normed field})\index{normed field} is an integral domain $R$ (resp., field $R$) equipped with a norm (as defined in Section 3.2) $|\cdot|: R \longrightarrow \RR_{\geq 0}$ on $R$.   The field $\CC((z)) = \CC[[z]][1/z]$ of formal Laurent series in $z$ is a normed field in the {\bf $(z)$-adic norm} $|\cdot|$, defined so that $|0| = 0$ and
$$\left |\sum_{k = n}^\infty a_k z^k \right| = 2^{-k}$$
for all $n \in \ZZ$ and $a_k \in \CC$ with $a_n \neq 0$.   Moreover, the normed field $\CC((z))$ is {\bf complete} in the sense that every Cauchy sequence in $\CC((z))$ converges.  We say that a formal continued fraction
$$G(z) = b_0(z)+  \frac{a_1(z)}{b_1(z) \, +} \  \frac{a_2(z)}{b_2(z) \,+} \  \frac{a_3(z)}{b_3(z) \,+} \ \cdots$$
with $a_n(z),b_n(z) \in \CC((z))$ for all $n$ {\bf converges formally}\index{converges formally in $\CC((z))$} to  $F(z) \in \CC((z))$ if $F(z) = \lim_{n \to \infty} w_n(z)$ in the normed field $\CC((z))$, or equivalently if $F(z) = w_0(z)+ \sum_{n = 1}^\infty( w_n(z)-w_{n-1}(z))$ in the normed field $\CC((z))$, where $w_n(z)$ is the $n$th approximant of $G(z)$.  For further details on formal convergence, see \cite[Chapter V]{lor}.  Note that formal convergence is a special case of the more general notion of the convergence of a continued fraction in a complete normed field $K$ (e.g., $K = \RR$ with the usual absolute value as the norm, or $K = \CC((z))$ with the $(z)$-adic norm).

\section{Asymptotic Jacobi and Steiltjes continued fraction expansions}

In this section  we provide some general results about asymptotic Jacobi and Stieltjes continued fraction expansions.

The {\bf rational degree} $\operatorname{rdeg} w$\index{rational degree of a rational function} of a rational function $w \in K(X)$ over a field $K$ is equal to the maximum of the degree of the numerator and  the degree of denominator of $w$ when $w$ is written as a quotient of two relatively prime polynomials in $K[X]$.  Let $f$ be a complex function defined on an unbounded subset of $\CC$.  We say that a rational  function $w \in \CC(z)$ is a {\bf best rational approximation of $f(z)$  (at $\infty$)}\index{best rational approximation of a function}  if $w$ is the unique rational function $v \in \CC(z)$  of rational degree at most $\operatorname{rdeg} w$  such that $f(z)-v(z) = O(f(z)-w(z)) \ (z \to \infty)$.

 A {\bf Jacobi continued fraction}, or {\bf J-fraction},\index{Jacobi continued fraction}\index{J-fraction} is a continued fraction of the form
\begin{align*}
 \frac{a_1}{z+b_1  \,-} \  \frac{a_2}{z+b_2 \,-} \  \frac{a_3}{z+b_3 \,-}\ \cdots,
\end{align*}
where $a_n, b_n \in \CC$ and $a_n \neq 0$ for all $n$. The following theorem characterizes asymptotic Jacobi continued fraction expansions and shows that, if a given complex function has an asymptotic Jacobi continued fraction expansion, then the best rational approximations of the function are precisely the approximants of the continued fraction.  

\begin{theorem}[{\cite[Theorem 2.7]{ell0}}]\label{wsimJ}
Let $\{a_n\}$ and $\{b_n\}$ be sequences of complex numbers with $a_n \neq 0$ for all $n$, and for all nonnegative integers $n$ let $w_n(z)$ denote the $n$th approximant of the Jacobi continued fraction
\begin{align*}
\frac{a_1}{z+b_1 \,-} \  \frac{a_2}{z+b_2 \,-} \  \frac{a_3}{z+b_3 \,-} \ \cdots.
\end{align*}
One has the following.
\begin{enumerate}
\item $w_n(z)$ is a rational function of rational degree $n$  with $w_n(z) \sim \frac{a_1}{z}  \ (z \to \infty)$ and $$w_n(z)-w_{n-1}(z) \sim \frac{a_1 a_2 \cdots a_n}{z^{2n-1}}  \ (z \to \infty)$$ for all $n \geq 1$.  In particular, $\{w_n(z)-w_{n-1}(z)\}_{n = 1}^\infty$ is an asymptotic sequence at $\infty$.  Moreover, the given continued fraction converges formally to a series  $\sum_{k = 1}^\infty \frac{c_k}{z^k}$ in $(1/z)\CC[[1/z]]$.
\item  Let $f(z)$ be a complex function defined on an unbounded subset of $\CC$.  The following conditions are equivalent.
\begin{enumerate}
\item $f(z)$ has the asymptotic expansion  $f(z) \simeq \sum_{k= 1}^\infty \frac{c_k}{z^k}\ (z \to \infty)$.
\item $f(z)$  has the asymptotic continued fraction expansion
\begin{align*}
f(z) \, \simeq \, \frac{a_1}{z+b_1 \,-} \  \frac{a_2}{z+b_2 \,-} \  \frac{a_3}{z+b_3 \,-} \ \cdots  \ (z \to \infty).
\end{align*}
\item $f(z) - w_n(z) \sim  \frac{a_1 a_2 \cdots a_{n+1}}{z^{2n+1}} \ (z \to \infty)$ for all  nonnegative integers $n$.
\item $f(z) - w_n(z) = O\left(\frac{1}{z^{2n+1}} \right) \ (z \to \infty)$ for infinitely many (or all) nonnegative integers $n$.
\item $w_n(z)$ for every nonnegative integer $n$ is the unique rational function  $w(z) \in \CC(z)$ of rational degree at most $n$ such that $f(z) - w(z) = O\left(\frac{1}{z^{2n+1}} \right) \ (z \to \infty)$.  
\item For every positive integer $n$ one has $f_n(z) \sim \frac{a_{n+1}}{z} \ (z \to \infty)$, where $f_n(z)$ is the function defined by the equation
\begin{align*}
f(z) =  \frac{a_1}{z+b_1 \,-} \  \frac{a_2}{z+b_2 \,-} \  \frac{a_3}{z+b_3 \,-} \ \cdots  \  \frac{a_{n-1}}{z+b_{n-1} \,-} \  \frac{a_n}{z+b_n - f_n(z)},
\end{align*}
or equivalently by the recurrence relation $f_{n+1}(z) = z+b_{n+1}-\frac{a_{n+1}}{f_n(z)}$, where $f_0(z) = f(z)$.
\item For every positive integer $n$, and for some (or all) $d_n \in \CC$, the function 
\begin{align*}
f(z) -  \frac{a_1}{z+b_1 \,-} \  \frac{a_2}{z+b_2 \,-} \  \frac{a_3}{z+b_3 \,-} \ \cdots \   \frac{a_{n-1}}{z+b_{n-1}\, -} \  \frac{a_n}{z +b_n+d_n}
\end{align*}
is asymptotic  at $\infty$ to $\frac{a_1 a_2 \cdots a_n}{z^{2n}}d_n$ if $d_n \neq 0$, or to $\frac{a_1a_2\cdots a_{n+1}}{z^{2n+1}}$ if $d_n = 0$.
\end{enumerate}
\item If the equivalent conditions of statement (2) hold, then the best rational approximations of $f(z)$  are precisely the approximants $w_n(z)$ for $n \geq 0$.
\end{enumerate}
\end{theorem}

 A {\bf Stieltjes continued fraction}, or {\bf S-fraction},\index{Stieltjes continued fraction}\index{S-fraction} is a continued fraction of the form
\begin{align*}
 \cfrac{\frac{a_1}{z}} {1 \, -} \ \cfrac{\frac{a_2}{z}} {1 \, -}  \ \cfrac{\frac{a_3}{z}} {1 \, -} \ \cfrac{\frac{a_4}{z}} {1 \, -} \  \cdots \ = \  \frac{a_1} {z \, -} \ \frac{a_2} {1 \, -}  \ \frac{a_3} {z \, -} \ \frac{a_4} {1 \, -} \  \cdots.
\end{align*}
where $a_n \in \CC$ and $a_n \neq 0$ for all $n$.   Such continued fractions were introduced by Stieltjes in \cite[J.3]{stie}.  The following result is more or less a corollary of Theorem \ref{wsimJ}.

\begin{corollary}[{\cite[Corollary 2.8]{ell0}}]\label{wsim}
Let $\{a_n\}$ be a sequence of nonzero complex numbers, and for all nonnegative integers $n$ let  $w_n(z)$ denote the $n$th approximant of the Stieltjes continued fraction
\begin{align*}
 \cfrac{\frac{a_1}{z}} {1 \, -} \ \cfrac{\frac{a_2}{z}} {1 \, -}  \ \cfrac{\frac{a_3}{z}} {1 \, -} \ \cfrac{\frac{a_4}{z}} {1 \, -} \  \cdots \ = \  \frac{a_1} {z \, -} \ \frac{a_2} {1 \, -}  \ \frac{a_3} {z \, -} \ \frac{a_4} {1 \, -} \  \cdots.
\end{align*}
One has the following.
\begin{enumerate}
\item $w_n(z)$ is a rational function of rational degree ${\left\lfloor \frac{n+1}{2}\right\rfloor}$  with $w_n(z) \sim \frac{a_1}{z}  \ (z \to \infty)$ and $$w_n(z)-w_{n-1}(z) \sim \frac{a_1 a_2 \cdots a_n}{z^{n}}  \ (z \to \infty)$$ for all $n \geq 1$.   In particular, $\{w_n(z)-w_{n-1}(z)\}_{n = 1}^\infty$ is an asymptotic sequence at $\infty$.  Moreover, the $2n$th approximant $w_{2n}(z)$ of the Stieltjes continued fraction above coincides with the $n$th approximant of the Jacobi continued fraction
\begin{align*}
\frac{a_1}{z-a_2 \,-} \  \frac{a_2 a_3}{z-a_3-a_4 \,-} \  \frac{a_4a_5}{z-a_5-a_6 \,-}  \ \cdots,
\end{align*}
and both continued fractions converge formally to a series  $\sum_{k = 1}^\infty \frac{c_k}{z^k}$ in $(1/z)\CC[[1/z]]$.
 \item Let $f(z)$ be a complex function defined on an unbounded subset of $\CC$.    The following conditions are equivalent.
\begin{enumerate}
\item $f(z)$ has the asymptotic expansion  $f(z) \simeq \sum_{k= 1}^\infty \frac{c_k}{z^k}\ (z \to \infty)$.
\item $f(z)$  has the asymptotic continued fraction expansion
\begin{align*}
f(z) \, \simeq \,  \cfrac{\frac{a_1}{z}} {1 \, -} \ \cfrac{\frac{a_2}{z}} {1 \, -}  \ \cfrac{\frac{a_3}{z}} {1 \, -} \ \cfrac{\frac{a_4}{z}} {1 \, -} \  \cdots  \ (z \to \infty).
\end{align*}
\item $f(z)$  has the asymptotic continued fraction expansion
\begin{align*}
f(z) \, \simeq \, \frac{a_1}{z-a_2 \,-} \  \frac{a_2 a_3}{z-a_3-a_4 \,-} \  \frac{a_4a_5}{z-a_5-a_6 \,-} \  \cdots  \ (z \to \infty).
\end{align*}
\item $f(z) - w_n(z) \sim  \frac{a_1 a_2 \cdots a_{n+1}}{z^{n+1}} \ (z \to \infty)$ for all  nonnegative integers $n$. 
\item $f(z) - w_n(z) = O\left(\frac{1}{z^{n+1}} \right) \ (z \to \infty)$ for  infinitely many (or all)   nonnegative integers $n$.
\item $w_{2n}(z)$ for every nonnegative integer $n$ is the unique rational function  $w(z) \in \CC(z)$ of rational degree at most $n$ such that $f(z) - w(z) = O\left(\frac{1}{z^{2n+1}} \right) \ (z \to \infty)$.  
\end{enumerate}
\item If the equivalent conditions of statement (2) hold, then the best rational approximations of $f(z)$  are precisely the even-indexed approximants $w_{2n}(z)$. 
\end{enumerate}
\end{corollary}

\begin{proof}
Statement (1) follows from the well-known formulas \cite[(1.4)--(1.5)]{wall} for continued fraction approximants and from the transformation  \cite[(28.2)--(28.4)]{wall}  from Stieltjes  to Jacobi continued fractions, and statements (2) and (3) follow immediately from statement (1) and Theorem \ref{wsimJ}.
\end{proof}

\begin{remark}[Good rational approximations]
Let us say that a rational  function $w \in \CC(z)$ is a {\bf good rational approximation of $f(z)$  (at $\infty$)}\index{good rational approximation of a function} if $\operatorname{rdeg} v \geq \operatorname{rdeg} w$ for any any $v \in \CC(z)$ such that $f(z)-v(z) = O(f(z)-w(z)) \ (z \to \infty)$.  Clearly any best rational approximation is a good rational approximation.
If $f(z)$ satisfies conditions (2)(a)--(f) of Corollary \ref{wsim}, then all of the approximants $w_n(z)$, not just the even-indexed ones, are good rational approximations of $f(z)$.  Indeed, if $f(z)-v(z) = O(f(z)-w_n(z)) \ (z \to \infty)$, then $w_n(z)-v(z) = O\left(\frac{1}{z^{n+1}}\right) \ (z \to \infty)$; but then $\operatorname{rdeg} v < \operatorname{rdeg} w_n = {\left\lfloor \frac{n+1}{2}\right\rfloor}$ implies $\operatorname{rdeg} w_n+ \operatorname{rdeg} v < n+1$, which leads to a contradiction, whence $\operatorname{rdeg} v \geq \operatorname{rdeg} w_n$.
\end{remark}

\section{Continued fractions with polynomial terms}

The following  theorem shows that any asymptotic continued fraction expansion at $0$ with respect to a continued fraction that has nonzero terms $a_n(z)$ in $z\CC[z]$ and $b_n(z)$ in $\CC[z]\backslash z\CC[z]$ is equivalent to an asymptotic expansion with respect to the asymptotic sequence $\{z^n\}$  at $0$.  

\begin{theorem}[{\cite[Theorem 2.10]{ell0}}]\label{wsimG}
Let $b_0(z) \in \CC[z]$, let $a_n(z), b_n(z) \in \CC[z]$ with $a_n(0) = 0$ to multiplicity $m_n \in \ZZ_{> 0}$ and $b_n(0) =  \beta_n \neq 0$ for all $n \geq 1$, and let  $w_n(z)$ denote the $n$th approximant of the continued fraction
\begin{align*}
 b_0(z)+  \frac{a_1(z)}{b_1(z) \, -} \  \frac{a_2(z)}{b_2(z) \,-} \  \frac{a_3(z)}{b_3(z) \,-} \ \cdots
\end{align*}
One has the following.
\begin{enumerate}
\item $w_n(z)$ is a rational function in $\CC(z)$ with $w_n(z)-b_0(z) \sim \frac{a_1(z)}{\beta_1}  \ (z \to 0)$ and $$w_n(z)-w_{n-1}(z) \sim \frac{a_1(z) a_2(z) \cdots a_n(z)}{(\beta_1\beta_2\cdots \beta_{n-1})^2\beta_n}  \ (z \to 0)$$ for all $n \geq 1$.  Moreover, $\{w_n(z)-w_{n-1}(z)\}_{n = 1}^\infty$ is an asymptotic sequence at $0$, and the given continued fraction converges $(z)$-adically to a series $\sum_{n = 0}^\infty c_n z^n$ in $\CC[[z]]$.  
 \item Let $f(z)$ be a complex function defined on a subset of $\CC$ having $0$ has a limit point.  The following conditions are equivalent.
\begin{enumerate}
\item $f(z)$ has the asymptotic expansion $f(z) \, \simeq \, \sum_{n = 0}^\infty c_n z^n \ (x \to 0)$.
\item $f(z)$  has the asymptotic continued fraction expansion
\begin{align*}
f(z) \, \simeq \,   b_0(z)+  \frac{a_1(z)}{b_1(z) \, -} \  \frac{a_2(z)}{b_2(z) \,-} \  \frac{a_3(z)}{b_3(z) \,-} \ \cdots \   \cdots \ (z \to 0).
\end{align*}
\item $f(z) - w_n(z) \sim  \frac{a_1(z) a_2(z) \cdots a_{n+1}(z)}{(\beta_1\beta_2\cdots \beta_n)^2\beta_{n+1}} \ (z \to 0)$ for all  nonnegative integers $n$. 
\item $f(z) - w_n(z) = O\left( z^{m_1+m_2 + \cdots+ m_{n+1}} \right) \ (z \to 0)$ for  infinitely many (or all)   nonnegative integers $n$.
\item  For every positive integer $n$ one has $f_n(z) \sim \frac{a_{n+1}(z)}{\beta_{n+1}} \ (z \to 0)$, where $f_n(z)$ is the function defined by the equation
\begin{align*}
f(z) =   b_0(z)+  \frac{a_1(z)}{b_1(z) \, -} \  \frac{a_2(z)}{b_2(z) \,-} \  \frac{a_3(z)}{b_3(z) \,-} \ \cdots \   \frac{a_{n}(z)}{b_{n}(z) \,-} \   \frac{f_n(z)}{1},
\end{align*}
or equivalently by the recurrence relation $f_{n+1}(z) = b_{n+1}(z)-\frac{a_{n+1}(z)}{f_n(z)}$, where $f_0(z) = f(z)-b_0(z)$.
\item  For every positive integer $n$, and for some function (or all functions) $g_n(z)$ with $g_n(z) = o(1) \ (x \to 0)$ and $a_{n+1}(z) = o(g_n(z)) \ (x \to 0)$,  the function  
\begin{align*}
f(z) -   b_0(z)-  \frac{a_1(z)}{b_1(z) \, -} \  \frac{a_2(z)}{b_2(z) \,-} \  \frac{a_3(z)}{b_3(z) \,-} \ \cdots \   \frac{a_{n}(z)}{b_{n}(z) \, -} \ \frac{g_n(z)}{1}
\end{align*}
is asymptotic  at $0$ to $-\frac{a_1(z) a_2(z) \cdots a_n(z)}{(\beta_1\beta_2\cdots \beta_n)^2} g_n(z)$.
\end{enumerate}
\end{enumerate}
\end{theorem}

\begin{proof}
We may suppose without loss  of generality that $b_0(z) = 0$.  The numerator $A_n = A_n(z)$ and denominator $B_n = B_n(z)$ of $w_n(z) = \frac{ A_n(z)}{B_n(z)}$ are uniquely determined by the well-known recurrence relations \cite[(1.4)]{wall}  and therefore by induction are polynomials with
$$B_n(z) = b_1(z)b_2(z)\cdots b_n(z) + \, \text{terms each involving some } a_i(z)$$
and
$$A_n(z) = a_1(z)(b_2(z)b_3(z)\cdots b_n(z)+ \, \text{terms each involving some } a_i(z))$$
for all $n$.  Therefore
$$B_n(z) = \beta_1\beta_2\cdots \beta_n+z D_n(z) \quad \text{ and } \quad A_n(z)= a_1(z)(\beta_2\beta_3\cdots \beta_n+ z C_n(z))$$
for some $C_n(z), D_n(z) \in \CC[z]$, and thus $w_n(z) \sim  \frac{a_1(z)\beta_2\beta_3\cdots \beta_n}{\beta_1\beta_2\cdots \beta_n}  = \frac{a_1(z)}{\beta_1} \  (z \to 0).$
Moreover, by \cite[(42.9)]{wall}, one has 
\begin{align*}
w_n(z)-w_{n-1}(z) = \frac{a_1(z) a_2(z) \cdots a_n(z)}{B_n(z)B_{n-1}(z)} = \frac{a_1(z) a_2(z) \cdots a_n(z)}{(\beta_1\beta_2\cdots \beta_n+z D_n(z))(\beta_1\beta_2\cdots \beta_{n-1}+z D_{n-1}(z))},
\end{align*}
and therefore
\begin{align*}
w_n(z)-w_{n-1}(z) 
= a_1(z) a_2(z) \cdots a_n(z)\left(\frac{1}{(\beta_1\beta_2\cdots \beta_{n-1})^2\beta_n} + zG_n(z)\right)
\end{align*}
for some $G_n(z) \in \CC[[z]]$.  It follows that 
\begin{align*}
w_n(z)-w_{n-1}(z) \equiv \frac{a_1(z) a_2(z) \cdots a_n(z)}{(\beta_1\beta_2\cdots \beta_{n-1})^2\beta_n}  \ \left(\operatorname{mod} \, \left(z^{m_1+m_2+\cdots+m_n+1}\right)\right)
\end{align*}
so the continued fraction in the theorem converges $(z)$-adically to a series $F(z) = \sum_{k = 1}^\infty (w_k(z)-w_{k-1}(z))$ in  $\CC[[z]]$.  Moreover, since also $G_n(z)$ is a rational function, it is analytic at $0$, and thus
\begin{align*}
w_n(z)-w_{n-1}(z) - \frac{a_1(z) a_2(z) \cdots a_n(z)}{(\beta_1\beta_2\cdots \beta_{n-1})^2\beta_n}  = O \left(z^{m_1+m_2+\cdots+m_n+1}\right) \ (z \to 0).
\end{align*}
Statement (1) follows.  

By statement (1) and Remark \ref{crem}(1), statements (2)(b)--(d) are equivalent.  Now, one has 
\begin{align*}
F(z)-w_n(z) =  \sum_{k = n+1}^\infty (w_k(z)-w_{k-1}(z)) \equiv w_{n+1}(z)-w_{n}(z) \  \left( \operatorname{mod} \, \left(z^{m_1+m_2+\cdots+m_{n+1}+1}\right)\right),
\end{align*} and therefore $F(z) - w_n(z) = z^{m_1+m_2+\cdots+ m_{n+1}}F_n(z)$
for some invertible $F_n(z) \in \CC[[z]]$.  Thus  one has $$\sum_{k = 1}^{m_1+m_2+\cdots+ m_{n+1}-1} c_kz^k - w_n(z)  = z^{m_1+m_2+\cdots+ m_{n+1}}H_n(z)$$
for some $H_n(z) \in  \CC[[z]]$.  Since then $H_n(z)$ is also a rational function of $z$,  it follows that  $H_n(z)$ is analytic at $0$ and thus
$$\sum_{k = 1}^{m_1+m_2+\cdots+ m_{n+1}-1}  c_kz^k - w_n(z)  = O\left( z^{m_1+m_2+\cdots+ m_{n+1}} \right) \ (z \to 0).$$  Therefore, one has
$$f(z) - \sum_{k = 1}^{m_1+m_2+\cdots+ m_{n+1}-1}c_k z^k  = O\left( z^{m_1+m_2+\cdots+ m_{n+1}} \right) \ (z \to 0), \quad \forall n \geq 0,$$
if and only if 
$$f(z) - w_n(z)  = O\left( z^{m_1+m_2+\cdots+ m_{n+1}} \right) \ (z \to 0), \quad  \forall n \geq 0.$$
This proves that statements (2)(a) and (2)(d) are equivalent.

Next, we verify the equivalence of (2)(c) and (2)(e).  By \cite[(1.3)]{wall} (or \cite[(1.3.2)]{cuyt}) and the definition of the function $f_n(z)$ one has
\begin{align}\label{fn}
f(z) = \frac{A_n(z)-A_{n-1}(z)f_n(z)}{B_n(z)-B_{n-1}(z)f_n(z)}
\end{align}
and therefore by the determinant formula   \cite[(1.5)]{wall}  \cite[(1.3.4)]{cuyt} one has
\begin{align}\label{AB}
f(z) - w_n(z) = \frac{A_n(z)-A_{n-1}(z)f_n(z)}{B_n(z)-B_{n-1}(z)f_n(z)}- \frac{ A_n(z)}{B_n(z)} =  \frac{a_1(z) a_2(z) \cdots a_n(z) f_n(z)}{B_n(z)(B_n(z)-B_{n-1}(z)f_n(z))}.
\end{align}
Thus, if (2)(e) holds, then 
$$f(z) - w_n(z) \sim  \frac{a_1(z) a_2(z) \cdots a_{n}(z)\frac{a_{n+1}(z)}{\beta_{n+1}}}{(\beta_1\beta_2\cdots \beta_{n})(\beta_1\beta_2\cdots \beta_n)}  =  \frac{a_1(z) a_2(z) \cdots a_{n+1}(z)}{(\beta_1\beta_2\cdots \beta_n)^2\beta_{n+1}},$$
whence (2)(c) holds.   Conversely, inverting (\ref{fn}), we have
\begin{align*}
f_n(z) = \frac{A_n(z)-B_n(z)f(z)}{A_{n-1}(z)-B_{n-1}(z)f(z)} = \frac{f(z)-w_n(z)}{f(z)-w_{n-1}(z)}\cdot \frac{B_n(z)}{B_{n-1}(z)},
\end{align*}
so, if (2)(c) holds, then
\begin{align*}
f_n(z) = \frac{f(z)-w_n(z)}{f(z)-w_{n-1}(z)}\cdot \frac{B_n(z)}{B_{n-1}(z)} \sim \frac{\frac{a_1(z) a_2(z) \cdots a_{n+1}(z)}{ ( \beta_1\beta_2\cdots \beta_n)^2 \beta_{n+1} }}{\frac{a_1(z) a_2(z) \cdots a_{n}(z)}{ (\beta_1\beta_2\cdots \beta_{n-1})^2 \beta_n }} \cdot \frac{\beta_1\beta_2\cdots \beta_n}{\beta_1\beta_2\cdots \beta_{n-1}} = \frac{a_{n+1}(z)}{\beta_{n+1}} \ (z \to \infty),
\end{align*} 
 whence (2)(e) holds.  Therefore, conditions (2)(c) and (2)(e) are equivalent.

Finally, we show that (2)(c) and (2)(f) are equivalent.   Let $$G_n(z) =  b_0(z)+ \frac{a_1(z)}{b_1(z) \, -} \  \frac{a_2(z)}{b_2(z) \,-} \  \frac{a_3(z)}{b_3(z) \,-} \ \cdots \   \frac{a_{n}(z)}{b_{n}(z) \,-} \  \frac{g_n(z)}{1}.$$  As in the proof of (\ref{AB}), one has
$G_n(z) - w_n(z) =  \frac{a_1(z)a_2(z) \cdots a_n(z) g_n(z)}{B_n(z)(B_n(z)-B_{n-1}(z)g_n(z))}$,
and therefore
$$G_n(z) - w_n(z)  \sim \frac{a_1(z) a_2(z) \cdots a_n(z) }{(\beta_1 \beta_2 \cdots \beta_n)^2 }g_n(z) =   o(w_n(z)-w_{n-1}(z)) \ (z \to 0)$$
Thus, if (2)(c) holds, then $$f(z) -w_n(x) \sim \frac{a_1(z) a_2(z) \cdots a_n(z) }{(\beta_1 \beta_2 \cdots \beta_n)^2 \beta_{n+1}}a_{n+1}(z) = o(- (G_n(z) - w_n(z)) \ (z \to 0)$$
and therefore
$$f(z)-G_n(z) \sim   - (G_n(z) - w_n(z))  \sim -  \frac{a_1(z) a_2(z) \cdots a_n(z) }{(\beta_1 \beta_2 \cdots \beta_n)^2 }g_n(z)\ (z \to 0),$$   
so that (2)(f) holds.  Conversely, if (2)(f) holds, then
$$f(z)-G_n(z) \sim  -\frac{a_1(z) a_2(z) \cdots a_n(z) }{(\beta_1 \beta_2 \cdots \beta_n)^2 }g_n(z) \sim  -(G_n(z)-w_n(z))  \ (z \to 0)$$
and therefore
$$f(z) - w_n(z) = o ( -(G_n(z)-w_n(z)) ) = o\left( w_n(z)-w_{n-1}(z)\right) \ (z \to 0),$$
so that (2)(b) holds.  Therefore, conditions (2)(b), (2)(c), and (2)(f) are equivalent.
\end{proof}

An  obvious ``inverse'' equivalent of Theorem \ref{wsimG},  obtained by replacing $z$ with $1/z$ and $0$ with $\infty$,  shows that any asymptotic continued fraction expansion at $\infty$ with respect to a continued fraction that has nonzero terms $a_n(z)$ in $(1/z)\CC[1/z]$ and $b_n(z)$ in $\CC[1/z]\backslash (1/z)\CC[1/z]$ is equivalent to an  asymptotic expansion with respect to the asymptotic sequence $\left\{\frac{1}{z^n}\right\}$  at $\infty$.  Since Theorem \ref{wsimG} applies to any {\it associated continued fraction} \cite[p.\ 36]{cuyt},
 its ``inverse'' equivalent applies to any  Jacobi continued fraction.  Thus, Theorem \ref{wsimG} yields the equivalence of conditions (2)(a)--(d), (2)(f), and (2)(g) of Theorem \ref{wsimJ}.  Condition (2)(e), however, is unique to Jacobi continued fractions.

Theorem  \ref{wsimG} also applies to any {\it C-fraction},  {\it T-fraction}, or {\it M-fraction} \cite[pp.\ 35--38]{cuyt}.   Notably, a {\bf C-fraction}\index{C-fraction} is a continued fraction of the form
\begin{align}\label{CF}
 G(z) = a_0+\frac{a_1 z^{m_1}} {1 \, -} \ \frac{a_2z^{m_2}} {1 \, -}  \ \frac{a_3z^{m_3}} {1 \, -} \ \frac{a_4z^{m_4}} {1 \, -} \  \cdots,
\end{align}
where $a_n \in \CC$ and $m_n \in \ZZ_{> 0}$ for all  $n$.  Any C-fraction $G(z)$ converges formally to a  series $F(z) \in \CC[[z]]$.  Conversely, any  series $F(z) \in \CC[[z]]$ can be written as the limit of a unique  C-fraction $G(z)$, called the {\bf C-fraction expansion of $F(z)$}.\index{C-fraction expansion of a power series}    The C-fraction expansion $G(z)$ of $F(z)$ is terminating, i.e., $a_n = 0$ for some $n$, if and only if $F(z)$ is a rational function in $\CC(z)$.   For constructive proofs of these assertions, see \cite[Chapter V]{lor}.

A {\bf regular C-fraction}\index{regular C-fraction} is a C-fraction (\ref{CF})  with $m_n = 1$ for all $n$.   A {\bf modified C-fraction}\index{modified C-fraction} is a continued fraction of the form $G(1/z)$ for some C-fraction $G(z)$.    Theorem \ref{wsimG} and \cite[Chapter V Theorem 5]{lor} imply that any asymptotic expansion  of the form $f(z) \, \simeq \, \sum_{n = 0}^\infty c_nz^n \ (z \to 0)$, where $ \sum_{n = 0}^\infty c_n z^n \in \CC[[z]]\backslash \CC(z)$, has a unique equivalent formulation as an asymptotic C-fraction expansion  at $0$, and vice versa.  (Thus, any non-rational  real or complex function that is analytic at $0$ has a unique asymptotic C-fraction expansion at $0$.)  Equivalently, any asymptotic expansion of the form $f(z) \, \simeq \, \sum_{n = 0}^\infty \frac{c_n}{z^n} \ (z \to \infty)$ with $ \sum_{n = 0}^\infty c_n z^n  \in \CC[[z]]\backslash \CC(z)$ has a unique equivalent formulation as an asymptotic modified C-fraction expansion  at $\infty$, and vice versa.

\section{Prime counting counting functions}

Two asymptotic continued fraction expansions of $\PP(x)$ and $A(x) = \log x - \frac{1}{\PP(x)}$ follow  immediately from (\ref{asex2})  and the identity
\begin{align*}
\sum_{k = 0}^n k! z^{k+1} =  \frac{z}{1 \, -}  \ \frac{z} {1+z\, -} \ \frac{2z} {1+2z\, -} \   \frac{3z} {1+3z \, -} \ \frac{4z}{1+4z \, -} \  \cdots \ \frac{nz}{1+nz}, \quad \forall n \geq 0.
\end{align*}
The identity above can be verified easily by induction from the well-known recurrence relation  \cite[(1.4)]{wall}  for the numerator or denominator $C_n(z)$ of the  $n$th approximant  $u_n(z)$ of a continued fraction, which for the continued fraction
$$ \frac{z}{1 \, -}  \ \frac{z} {1+z\, -} \ \frac{2z} {1+2z\, -} \   \frac{3z} {1+3z \, -} \ \frac{4z}{1+4z \, -} \  \cdots \ = \  \cfrac{1}{\frac{1}{z} \, -}  \ \cfrac{\frac{1}{z}} {\frac{1}{z}+1\, -} \ \cfrac{\frac{2}{z}} {\frac{1}{z}+2\, -} \   \cfrac{\frac{3}{z}} {\frac{1}{z}+3 \, -} \ \cfrac{\frac{4}{z}}{\frac{1}{z}+4 \, -} \  \cdots$$
is the recurrence relation
$$C_{k+1}(z) = (1+kz)C_k(x) -kz C_{k-1}(z), \quad \forall k \geq 1.$$   It follows that  $u_n(z) = \sum_{k = 0}^n \frac{k!}{(\log x)^{k+1}}$  for all $n$, where $z = \frac{1}{\log x}$.  From the definiton of an asymptotic continued fraction expansion, then, we obtain the four restatements of (\ref{asex2}) provided in the following proposition.

\begin{proposition}\label{gencornew}
One has the following asymptotic continued fraction expansions.
\begin{enumerate}
\item $\PP(x)\,  \simeq \,  \cfrac{\frac{1}{\log x}}{1 \, -}  \ \cfrac{\frac{1}{\log x}} {1+\frac{1}{\log x}\, -} \ \cfrac{\frac{2}{\log x}} {1+\frac{2 }{\log x}\, -} \   \cfrac{\frac{3}{\log x}} {1+\frac{3}{\log x} \, -} \ \cfrac{\frac{4}{\log x}}{1+\frac{4}{\log x} \, -}  \  \cdots \ (x \to \infty)$.
\item $\PP(x)  \,  \simeq \, \displaystyle  \frac{1}{\log x \, -}  \ \frac{1\log x} {\log x+1 \, -} \ \frac{2\log x} {\log x+2 \, -} \   \frac{3\log x} {\log x+3\, -} \ \frac{4\log x}{\log x+4 \, -}  \  \cdots \  (x \to \infty).$ 
\item $A(x) \, \simeq \,  \cfrac{\frac{1}{\log x}} {1+\frac{1}{\log x}\, -} \ \cfrac{\frac{2}{\log x}} {1+\frac{2 }{\log x}\, -} \   \cfrac{\frac{3}{\log x}} {1+\frac{3}{\log x} \, -} \ \cfrac{\frac{4}{\log x}}{1+\frac{4}{\log x} \, -}  \  \cdots \ (x \to \infty)$.
\item $A(x) \, \simeq \, \displaystyle  \frac{1\log x} {\log x+1 \, -} \ \frac{2\log x} {\log x+2 \, -} \   \frac{3\log x} {\log x+3\, -} \ \frac{4\log x}{\log x+4 \, -}  \  \cdots \  (x \to \infty)$.
\end{enumerate}
\end{proposition}

It is known that the C-fraction expansion of the formal power series $\sum_{n = 0}^\infty n!z^{n+1}$ in $\CC[[z]]$ is given by
\begin{align}\label{stieltid}
\sum_{n = 0}^\infty n!z^{n+1}= \frac{z}  {1 \, -} \ \frac{z} {1 \, -}  \ \frac{z} {1 \, -} \ \frac{2z} {1 \, -}  \ \frac{2z} {1 \, -} \ \frac{3z} {1 \, -} \ \frac{3z} {1 \, -}  \  \cdots.
\end{align}
Indeed, this was first proved by Stieltjes in \cite[No.\ 57]{stie}, from Hankel matrix determinant expressions he obtained in \cite[No.\ 11]{stie} for the terms of a regular C-fraction in terms of its formal power series expansion in $\CC[[z]]$.    Together with Corollary \ref{wsim}, this yields the following theorem.

\begin{theorem}\label{maincontthm1}
One has the following asymptotic continued fraction expansions.
\begin{enumerate}
\item $\PP(x) \, \simeq \,  \cfrac{\frac{1}{\log x}}{1 \,-} \ \cfrac{\frac{1}{\log x}}{1 \,-}\  \cfrac{\frac{1}{\log x}}{1 \,-}\  \cfrac{\frac{2}{\log x}}{1 \,-}\  \cfrac{\frac{2}{\log x}}{1 \,-}\  \cfrac{\frac{3}{\log x}}{1 \,-} \  \cfrac{\frac{3}{\log x}}{1 \,-}\  \cfrac{\frac{4}{\log x}}{1 \,-}  \ \cfrac{\frac{4}{\log x}}{1 \,-}  \ \cdots \ (x \to \infty)$.
\item $\PP(x)  \,  \simeq \, \displaystyle \frac{1}{\log x -1\,-} \  \frac{1}{\log x  - 3 \,-}\  \frac{4}{\log x-5\,-}\  \frac{9}{\log x - 7 \,-} \ \frac{16}{\log x-9\,-}\   \cdots \  (x \to \infty)$.
\item $A(x) \,  \simeq \,   \cfrac{\frac{1}{\log x}}{1 \,-}\  \cfrac{\frac{1}{\log x}}{1 \,-}\  \cfrac{\frac{2}{\log x}}{1 \,-}\  \cfrac{\frac{2}{\log x}}{1 \,-}\  \cfrac{\frac{3}{\log x}}{1 \,-} \  \cfrac{\frac{3}{\log x}}{1 \,-}\  \cfrac{\frac{4}{\log x}}{1 \,-}  \ \cfrac{\frac{4}{\log x}}{1 \,-}  \ \cdots \ (x \to \infty)$.
\item $A(x)  \,  \simeq \, \displaystyle 1+  \frac{1}{\log x  - 3 \,-}\  \frac{4}{\log x-5\,-}\  \frac{9}{\log x - 7 \,-} \ \frac{16}{\log x-9\,-}\   \cdots \  (x \to \infty)$.
\end{enumerate}
\end{theorem}

\begin{proof}
Since $\PP(x)$  has the asymptotic expansion (\ref{asex2}), the theorem follows immediately from (\ref{stieltid})  and the equivalence of statements (2)(a)--(c) of Corollary \ref{wsim} (and Remark \ref{crem}(5)).
\end{proof}

As a corollary of Theorems \ref{maincontthm1} and \ref{wsimJ}, we obtain the following.

\begin{corollary}\label{maincor}
The best rational approximations of the function $\PP(e^x)$ are precisely the approximants $w_n(x)$ of the continued fraction $$\frac{1}{x-1 \,-} \  \frac{1}{x-3 \,-} \  \frac{4}{x-5 \,-}\  \frac{9}{x-7 \,-}\  \frac{16}{x-9 \,-} \  \cdots.$$  Moreover, for all $n \geq 0$, one has $\PP(e^x)- w_n(x) \sim \frac{(n!)^2}{x^{2n+1}}  \ (x \to \infty)$, and $w_n(x)$ is the unique rational function of rational degree at most $n$ such that $\PP(e^x) - w_n(x) = O\left( \frac{1}{x^{2n+1}}\right) \ (x \to \infty)$.  Furthermore, for all $n \geq 1$ and all $a \in \RR$, the function
\begin{align*}
\PP(e^x) -  \frac{1}{x-1 \,-} \  \frac{1}{x-3 \,-} \  \frac{4}{x-5 \,-} \ \cdots \   \frac{(n-2)^2}{x-(2n-3)\, -} \  \frac{(n-1)^2}{x -(2n-1)+a}
\end{align*}
is asymptotic to $\frac{((n-1)!)^2}{x^{2n}}a$ if $a \neq 0$, and to $\frac{(n!)^2}{x^{2n+1}}$ if $a = 0$.
\end{corollary}

Thus, for example,  one has
$$\PP(e^x)\sim \frac{1}{x} \ (x \to \infty),$$
$$\PP(e^x)- \frac{1}{x-1}\sim \frac{1}{x^{3}} \ (x \to \infty),$$
$$\PP(e^x)- \frac{x-3}{x^2-4x+2}\sim \frac{4}{x^{5}} \ (x \to \infty),$$
$$\PP(e^x)-\frac{x^2-8x+11}{x^3-9x^2+18x-6}\sim \frac{36}{x^{7}} \ (x \to \infty),$$
$$\PP(e^x)- \frac{x^3-15x^2+58x-50}{x^4-16x^3+72x^2-96x+24}\sim \frac{576}{x^{9}} \ (x \to \infty).$$
However, the convergence of the corresponding ratios is likely extremely slow.  For example,  the quantity
$$\left(\frac{\Ri(e^x)}{e^x}- \frac{x-3}{x^2-4x+2}\right)\frac{x^5}{4} \to 1$$ is within $\frac{1}{8}$ of its limit of $1$  for all $x \geq N$, but not until $N$ is greater, roughly, than $e^{113.33} \approx 1.605 \cdot 10^{49}$, which suggests that something similar may hold for the quantity
$$\left(\PP(e^x)- \frac{x-3}{x^2-4x+2}\right)\frac{x^5}{4} \to 1.$$   Based on such comparisons with $\frac{\Ri(e^x)}{e^x}$, it is likely preferable to express the error asymptotics precisely as they are expressed in the definition of asymptotic continued fraction expansions (using  \cite[(42.9)]{wall}),  as follows:
$$\PP(e^x)\sim \frac{1}{x-1} \ (x \to \infty),$$
$$\PP(e^x)- \frac{1}{x-1}\sim \frac{1}{(x-1)(x^2-4x+2)} \ (x \to \infty),$$
$$\PP(e^x)- \frac{x-3}{x^2-4x+2}\sim \frac{4}{(x^2-4x+2)(x^3-9x^2+18x-6)} \ (x \to \infty),$$
$$\PP(e^x)-\frac{x^2-8x+11}{x^3-9x^2+18x-6}\sim \frac{36}{(x^3-9x^2+18x-6)(x^4-16x^3+72x^2-96x+24)} \ (x \to \infty),$$
and so on.  The corresponding ratios appear to converge to $1$ much faster than those of their respective counterparts $\PP(e^x)- w_n(x) \sim \frac{(n!)^2}{x^{2n+1}}  \ (x \to \infty)$.  We leave it as an open problem to prove explicit inequalities based on any of the asymptotics above.

\begin{remark}[Best rational approximations of $\PP(e^x)$]\label{Laguerre} \
\begin{enumerate}
\item The numerator $P_n(x)$ and denominator  $Q_n(x)$ of the rational function $w_n(x)$ in Corollary \ref{maincor} are monic integer polynomials of degree $n-1$ and $n$, respectively, and  one has $w_{n+1}(x)-w_{n}(x) = \frac{(n!)^2}{Q_{n+1}(x)Q_{n}(x)}$  for all $n \geq 0$.     Laguerre showed in \cite{lagu} that the denominator $Q_n(x)$ of $w_n(x)$ is given by $Q_n(x) = \widehat{L}_n(x)$, where $$\widehat{L}_n(x) = (-1)^n n! L_n(x) = \sum_{k = 0}^n (-1)^k k!{n\choose k}^2 x^{n-k}$$
is the $n$th {\bf monic Laguerre polynomial}\index{monic Laguerre polynomials} and
$$L_n(x) = \sum_{k = 0}^n \frac{(-1)^k}{k!}{n\choose k} x^k = {}_{1}F_{1}(-n;1;x)$$
is the  {\bf $n$th Laguerre polynomial}. \index{Laguerre polynomials} The polynomial $L_n(x)$ is known to have $n$ distinct positive real roots, and its largest root lies between $3n-4$ and $4n+2$ \cite{sko}.  It follows that there is no fixed neighborhood of $\infty$ on which  all of the functions $w_n(x)-w_{n-1}(x)$ are defined.
\item From  \cite[{[1.14]}]{akh}, one can deduce that the numerator $P_n(x)$ of $w_n(x) = \frac{P_n(x)}{\widehat{L}_n(x)}$ is given by
$P_n(x) =  \sum_{k = 0}^{n-1} a_{n,k} x^{n-1-k},$
where $$a_{n,k} = \sum_{j = 0}^k  {(-1)^{j}}j!(k-j)! {n \choose j}^2=  (-1)^{k} k! {n \choose k}^2 {}_{3}F_{2}(1,1,-k; n-k+1, n-k+1; 1)$$
for all $n \geq 1$ and $1 \leq k \leq n$,
with explicit values $a_{n,0} = 1$,  $a_{n,1} = 1-n^2$, $a_{n,n-1} = (-1)^{n+1}n!H_n = s(n+1,2)$, and $a_{n,n-2} = (-1)^n n!((n+1)H_n-2n) =(-1)^n\langle \langle n,n-2 \rangle \rangle$, where $H_n = \sum_{k = 1}^n \frac{1}{k}$ denotes the $n$th harmonic number, $s(n,k)$ denotes the signed Stirling number of the first kind, and $\langle \langle n,k \rangle \rangle$ denotes the Eulerian number of the second kind.
\end{enumerate}
\end{remark}

 For $a, b \in \CC$, consider the {\it hypergeometric series} $${}_{2}F_{0}(a,b;;z) = \sum_{k = 0}^\infty \frac{(a)_k(b)_k}{n!}z^k$$
in $\CC[[z]]$, where $(x)_n = x(x+1)(x+2)\cdots(x+n-1)$ denotes the {\bf Pochhammer symbol}. \index{Pochhammer symbol} 
 In \cite[(89.5) and (92.2)]{wall} it is proved that ${}_{2}F_{0}(a,1;;z)$ for all $a \in \CC$ has the formal C-fraction expansion
$$\sum_{k = 0}^\infty (a)_k z^k = {}_{2}F_{0}(a,1;;z)  = \frac{1} {1 \, -} \ \frac{a z} {1 \, -}  \ \frac{1z} {1 \, -} \ \frac{(1+a)z} {1 \, -} \   \frac{2z} {1 \, -} \ \frac{(2+a)z} {1 \, -} \   \frac{3z} {1 \, -} \ \frac{(3+a)z} {1 \, -} \cdots.$$
We now use this, along with  Corollary \ref{wsim}, to prove the following.

\begin{theorem}\label{gentheorem}
For every nonnegative integer $n$ one has the asymptotic continued fraction expansions
$$p_n(x) \, \simeq \,  \cfrac{\frac{1}{\log x}}{1 \,-} \ \cfrac{\frac{1+n}{\log x}}{1 \,-}\  \cfrac{\frac{1}{\log x}}{1 \,-}\  \cfrac{\frac{2+n}{\log x}}{1 \,-}\  \cfrac{\frac{2}{\log x}}{1 \,-}\  \cfrac{\frac{3+n}{\log x}}{1 \,-} \  \cfrac{\frac{3}{\log x}}{1 \,-}\  \cfrac{\frac{4+n}{\log x}}{1 \,-}  \ \cfrac{\frac{4}{\log x}}{1 \,-}  \ \cdots \ (x \to \infty)$$
and
$$\PP_n(x)  \,  \simeq \, \frac{1}{\log x -1-n\,-} \  \frac{1(1+n)}{\log x - 3-n \,-}\  \frac{2(2+n)}{\log x-5-n \,-}\  \frac{3(3+n)}{\log x  - 7-n \,-}  \ \cdots \  (x \to \infty),$$
where  $\PP_n(x) = \frac{(\log x)^{n}}{n!} \left(\PP(x) - \sum_{k = 0}^{n-1} \frac{k!}{(\log x)^{k+1}}\right).$  Moreover, the same asymptotic continued fraction expansions hold for the function $l_n(x) = \frac{(\log x)^{n}}{n!} \left(\frac{\li(x)}{x} - \sum_{k = 0}^{n-1} \frac{k!}{(\log x)^{k+1}}\right).$
\end{theorem}

\begin{proof}
Formally in $\CC[[1/z]]$, one has
$$\sum_{k = 0}^\infty\frac{ (k+n)! }{n!}\frac{1}{z^{k+1}}  = \frac{1}{z}{}_{2}F_{0}(n+1,1;;1/z)  = \cfrac{\frac{1}{z}} {1 \, -} \ \cfrac{\frac{1+n}{z}} {1 \, -}  \ \cfrac{\frac{1}{z}} {1 \, -} \ \cfrac{\frac{2+n}{z}} {1 \, -} \   \cfrac{\frac{2}{z}} {1 \, -} \ \cfrac{\frac{3+n}{z}} {1 \, -} \   \cfrac{\frac{3}{z}} {1 \, -} \ \cfrac{\frac{4+n}{z}} {1 \, -} \cdots.$$
Also, from the asymptotic expansion (\ref{asex2})  of $\PP(x)$
we easily obtain the asymptotic expansion
\begin{align*}
\PP_n(e^x) \simeq \sum_{k = 0}^\infty \frac{(k+n)!}{n!} \frac{1}{x^{k+1}} \ (x \to \infty),
\end{align*}
and likewise from (\ref{asex}) we obtain the same asymptotic expansion for $l_n(e^x)$.
As a consequence, the theorem follows from  the equivalence of statements (2)(a)--(c) of Corollary \ref{wsim}.
\end{proof}

\begin{remark}[Best rational approximations of $\PP_n(e^x)$]
By \cite[pp.\ 95--96]{cuyt}, the monic denominator of the $k$th approximant of the Jacobi continued fraction in Theorem \ref{gentheorem} is equal to $\widehat{L}_k^{(n)}(\log x)$, where $\widehat{L}_k^{(n)}(z)$ denotes the {\it $k$th monic  generalized Laguerre  polynomial}.
\end{remark}

To provide further context for Theorem \ref{gentheorem}, we note that the fundamental asymptotic expansion  (\ref{asex2}) of $\PP(x)$ is by definition equivalent to  the asymptotic $\PP_n(x) \sim \frac{1}{\log x} \ (x \to \infty)$ holding for all $n \geq 0$, and the theorem for $n = 1$ yields the following.

\begin{corollary}\label{Plog}
One has the asymptotic continued fraction expansions
$$\PP(x)\log x  \,  \simeq \, 1+ \cfrac{\frac{1}{\log x}}{1 \,-} \ \cfrac{\frac{2}{\log x}}{1 \,-}\  \cfrac{\frac{1}{\log x}}{1 \,-}\  \cfrac{\frac{3}{\log x}}{1 \,-}\  \cfrac{\frac{2}{\log x}}{1 \,-}\  \cfrac{\frac{4}{\log x}}{1 \,-} \  \cfrac{\frac{3}{\log x}}{1 \,-}\  \cfrac{\frac{5}{\log x}}{1 \,-} \  \cfrac{\frac{4}{\log x}}{1 \,-}  \ \ \cdots \ (x \to \infty)$$
and
$$\PP(x)\log x  \,  \simeq \, 1+ \frac{1}{\log x -2\,-} \  \frac{1\cdot 2}{\log x - 4 \,-}\  \frac{2\cdot 3}{\log x-6 \,-}\  \frac{3 \cdot 4}{\log x - 8 \,-}  \  \frac{4 \cdot 5}{\log x - 10 \,-} \ \cdots \  (x \to \infty).$$
\end{corollary}

 By \cite[p.\ 243]{cuyt}, for all $a \in \CC$ the series $\frac{1}{z}{}_{2}F_{0}(a,1;;1/z) \in \CC[[1/z]]$ has  formal T-fraction expansion 
\begin{align*}
\frac{1}{z}{}_{2}F_{0}(a,1;;1/z) &  = \cfrac{\frac{1}{z}}{1+\frac{1-a}{z} \, -}  \ \cfrac{\frac{1}{z}} {1+\frac{2-a}{z}\, -} \ \cfrac{\frac{2}{z}} {1+\frac{3-a }{z}\, -} \   \cfrac{\frac{3}{z}} {1+\frac{4-a}{z} \, -} \ \cfrac{\frac{4}{z}}{1+\frac{5-a}{z} \, -}  \  \cdots  \\
&  = \frac{1}{z+1-a \, -}  \ \frac{1z} {z+2-a \, -} \ \frac{2z} {z+3-a \, -} \   \frac{3z} {z+4-a \, -} \ \frac{4z}{z+5-a \, -}  \  \cdots.
\end{align*}
Thus, from the ``inverse'' equivalent of Theorem  \ref{wsimG}, and as in the proof of Theorem \ref{gentheorem}, we obtain the following result,  which is a generalization of Proposition \ref{gencornew} and an analogue of Theorem \ref{gentheorem}.

\begin{theorem}\label{gentheoremnew}
For every nonnegative integer $n$ one has the asymptotic continued fraction expansion
$$\PP_n(x) \, \simeq \,  \cfrac{\frac{1}{\log x}}{1+\frac{-n}{\log x} \, -}  \ \cfrac{\frac{1}{\log x}} {1+\frac{1-n}{\log x}\, -} \ \cfrac{\frac{2}{\log x}} {1+\frac{2-n }{\log x}\, -} \   \cfrac{\frac{3}{\log x}} {1+\frac{3-n}{\log x} \, -} \ \cfrac{\frac{4}{\log x}}{1+\frac{4-n}{\log x} \, -}  \  \cdots \ (x \to \infty),$$
or, equivalently,
$$\PP_n(x)  \,  \simeq \,  \frac{1}{\log x-n \, -}  \ \frac{1\log x} {\log x+1 -n\, -} \ \frac{2\log x} {\log x+2-n \, -} \   \frac{3\log x} {\log x+3-n\, -}   \  \cdots \  (x \to \infty),$$
where  $\PP_n(x) = \frac{(\log x)^{n}}{n!} \left(\PP(x) - \sum_{k = 0}^{n-1} \frac{k!}{(\log x)^{k+1}}\right).$  Moreover, the same asymptotic continued fraction expansions hold for the function $l_n(x) = \frac{(\log x)^{n}}{n!} \left(\frac{\li(x)}{x} - \sum_{k = 0}^{n-1} \frac{k!}{(\log x)^{k+1}}\right).$
\end{theorem}

Note that the two  continued fractions in the theorem are equivalent, and, if $w_{n,k}(x)$ denotes the $k$th approximant, then by  Theorem  \ref{wsimG}   one has $\PP_n(x) - w_{n,k}(x) \sim \frac{k!}{(\log x)^k} \ (x \to \infty)$ for all $n$ and $k$.  For $n = 0$, the theorem yields Proposition \ref{gencornew}.  For $n = 1$, it yields the following.

\begin{corollary}\label{rkr}
One has the asymptotic continued fraction expansion
$$\PP(x)\log x \, \simeq \, 1+ \cfrac{\frac{1}{\log x}}{1-\frac{1}{\log x} \, -}  \ \cfrac{\frac{1}{\log x}} {1\, -} \ \cfrac{\frac{2}{\log x}} {1+\frac{1 }{\log x}\, -} \   \cfrac{\frac{3}{\log x}} {1+\frac{2}{\log x} \, -} \ \cfrac{\frac{4}{\log x}}{1+\frac{3}{\log x} \, -}  \  \cdots \ (x \to \infty),$$
or, equivalently,
$$\PP(x)\log x \, \simeq \, 1+ \frac{1}{\log x-1 \, -}  \ \frac{1\log x} {\log x \, -} \ \frac{2\log x} {\log x+1 \, -} \   \frac{3\log x} {\log x+2\, -} \ \frac{4\log x}{\log x+3 \, -}  \  \cdots \  (x \to \infty).$$
\end{corollary}

 Corollary \ref{rkr} also follows from (\ref{asex2}) and the fact, proved by induction, that the $n$th approximant of both continued fractions in the corollary is equal to $\frac{n \cdot n!}{(\log x - n)(\log x)^n} + \sum_{k = 0}^n \frac{k!}{(\log x)^k}$.

As a corollary of Theorem \ref{maincontthm1}, Corollary \ref{Plog}, and Lemmas \ref{asympprop} and \ref{asympprop2a}, we obtain the following.

\begin{proposition}\label{LP}
Let $L(x)$ be any function such that $L(x) = \log x + o ((\log x)^{-t}) \  (x \to \infty)$ for all $t > 0$, and let $P(x)$ be any function such that $P(x) = \frac{\li(x)}{x} + o ((\log x)^{-t}) \  (x \to \infty)$ for all $t > 0$.   For all $t \in \RR$, one has the following asymptotic continued fraction expansions at $\infty$.
\begin{enumerate}
\item  $\displaystyle P(e^t x) \,  \simeq \,  \cfrac{\frac{1}{L(x)+t}}{1 \,-} \ \cfrac{\frac{1}{L(x)+t}}{1 \,-}\  \cfrac{\frac{1}{L(x)+t}}{1 \,-}\  \cfrac{\frac{2}{L(x)+t}}{1 \,-}\  \cfrac{\frac{2}{L(x)+t}}{1 \,-}\  \cfrac{\frac{3}{L(x)+t}}{1 \,-} \  \cfrac{\frac{3}{L(x)+t}}{1 \,-}\  \cfrac{\frac{4}{L(x)+t}}{1 \,-}  \ \cfrac{\frac{4}{L(x)+t}}{1 \,-}  \ \cdots$.
\item  $\displaystyle P(e^t x)  \,  \simeq \, \frac{1}{L(x)+t -1\,-} \  \frac{1}{L(x)+t  - 3 \,-}\  \frac{4}{L(x)+t-5\,-}\  \frac{9}{L(x)+t - 7 \,-} \ \frac{16}{L(x)+t-9\,-}  \ \cdots$.
\item $\displaystyle {P(e^t x)(L(x)+t )}   \, \simeq \, 1+ \cfrac{\frac{1}{L(x)+t}}{1 \,-} \ \cfrac{\frac{2}{L(x)+t}}{1 \,-}\  \cfrac{\frac{1}{L(x)+t}}{1 \,-}\  \cfrac{\frac{3}{L(x)+t}}{1 \,-}\  \cfrac{\frac{2}{L(x)+t}}{1 \,-}\  \cfrac{\frac{4}{L(x)+t}}{1 \,-} \  \cfrac{\frac{3}{L(x)+t}}{1 \,-}\  \cfrac{\frac{5}{L(x)+t}}{1 \,-}  \ \cdots$.
\item $\displaystyle {P(e^t x)(L(x)+t) }   \, \simeq \,  1+\frac{1}{L(x)+t-2 \,-} \  \frac{1\cdot 2}{L(x)+t-4 \,-}\  \frac{2 \cdot 3}{L(x)+t-6 \,-}\  \frac{3 \cdot 4}{L(x)+t-8 \,-}  \ \cdots$.
\end{enumerate}
\end{proposition}

\begin{example}  The following provide examples of functions $L(x)$ satisfying the hypotheses of Proposition \ref{LP}.
\begin{enumerate}
\item   It is well known that $H_z -\gamma -\log z \sim \frac{1}{2z} \ (z \to \infty).$
Moreover, one has $H_z - \gamma = \Psi(z+1) = \Psi(z) + \frac{1}{z}$ for all $z \in \CC \backslash \{0,-1,-2,-3,\ldots \}$.  Thus,  one has
\begin{align}\label{harm}
H_x -\gamma = \log x+ o((\log x)^{-t}) \ (x \to \infty),
\end{align}
and likewise $$\Psi(x) = \log x + o((\log x)^{-t}) \ (x \to \infty),$$ for all $t > 0$.
\item  Mertens' third vc theorem can be improved to show that the function $L(x) = e^{-\gamma}\prod_{p \leq x}\left(1-\frac{1}{p}\right)^{-1}$  satisfies $L(x) \sim \log x + o ((\log x)^{-t}) \  (x \to \infty)$ for all $t > 0$: see Proposition \ref{mertensprop}.
\item By Propositions \ref{merttpropp} and \ref{mertprop1} in Section 10.1,  the functions $\sum_{p \leq x} \frac{\log p}{p}-B$ and $e^{-M} \prod_{p\leq x} e^{1/p}$ also satisfy the hypotheses of  Proposition \ref{LP} for $L(x)$.
\end{enumerate}
\end{example}

\begin{example}
Generalizing the obvious example $\frac{\pi(x)}{x}$,  Theorem \ref{arithsemigroup}  below provides examples $ \frac{\pi_G(x^{1/\delta})}{x}$ and $A\PP_G(x^{1/\delta})$ of  functions $P(x)$ that satisfiy the hypothesis $P(x) \sim \frac{\li(x)}{x} + o ((\log x)^{-t}) \  (x \to \infty)$, for all $t > 0$, of Proposition \ref{LP}.
\end{example}

From Proposition \ref{LP} and (\ref{harm}), we obtain the following.

\begin{corollary}\label{pas2}
One has the following asymptotic continued fraction expansions.
\begin{enumerate}
\item $\displaystyle \PP(x)\,  \simeq \,   \cfrac{\frac{1}{H_x-\gamma}}{1 \, -} \  \cfrac{\frac{1}{H_x-\gamma}}{1 \,-}\ \cfrac{\frac{1}{H_x-\gamma}}{1 \,-}\  \cfrac{\frac{2}{H_x-\gamma}}{1 \,-}\  \cfrac{\frac{2}{H_x-\gamma}}{1 \,-}  \ \cfrac{\frac{3}{H_x-\gamma}}{1 \,-}\  \cfrac{\frac{3}{H_x-\gamma}}{1 \,-} \  \cdots \ (x \to \infty).$
\item $\displaystyle  \PP(x) \, \simeq \,  \frac{1}{H_x-\gamma-1 \,-} \ \frac{1}{H_x-\gamma-3 \,-} \  \frac{4}{H_x-\gamma-5 \,-}\  \frac{9}{H_x-\gamma-7 \,-}\  \cdots \ (x \to \infty)$.
\item $\displaystyle \PP(e^\gamma x)\,  \simeq \,   \cfrac{\frac{1}{H_x}}{1 \, -} \  \cfrac{\frac{1}{H_x}}{1 \,-}\ \cfrac{\frac{1}{H_x}}{1 \,-}\  \cfrac{\frac{2}{H_x}}{1 \,-}\  \cfrac{\frac{2}{H_x}}{1 \,-}  \ \cfrac{\frac{3}{H_x}}{1 \,-}\  \cfrac{\frac{3}{H_x}}{1 \,-} \  \cdots \ (x \to \infty).$
\item $\displaystyle  \PP(e^\gamma x) \, \simeq \,  \frac{1}{H_x-1 \,-} \ \frac{1}{H_x-3 \,-} \  \frac{4}{H_x-5 \,-}\  \frac{9}{H_x-7 \,-}\  \frac{16}{H_x-9 \,-}\ \cdots \ (x \to \infty)$.
\end{enumerate}
\end{corollary}

All of the asymptotic continued fraction expansions of $\PP(x)$ we have proved in this section can be generalized to any  arithmetic semigroup $G$  satisfying the  equivalent conditions of  Theorem \ref{abstractPNT}.   (See Section 5.4 for the relevant definitions.)  Let $G$ be an arithmetic semigroup.  Recall that, for all $$\PP_G(x) = \frac{\pi_G(x)}{N_G(x)}$$ for all $x \geq 1$ denotes the probability that a randomly selected element of $G$ of norm less than or equal to $x$ is prime.   Using Theorem \ref{abstractPNT},  Lemma \ref{asympprop},  (\ref{asex}), and Corollary \ref{wsim},  we can generalize Theorem \ref{maincontthm1} and Corollary \ref{maincor} as follows.

\begin{theorem}\label{arithsemigroup}
Let $G$ be an arithmetic semigroup, and let $\delta > 0$.   Let $w_n(x)$ for any nonnegative integer $n$ denote the $n$th approximant of the continued fraction
$$\frac{1}{x-1 \,-} \  \frac{1}{x-3 \,-} \  \frac{4}{x-5 \,-}\  \frac{9}{x -7\,-} \ \frac{16}{x-9 \,-}\  \cdots.$$
\begin{enumerate}
\item The following conditions are equivalent. 
\begin{enumerate}
\item There is an $A > 0$ such that $N_G(x) = Ax^\delta+ O(x^\delta (\log x)^t) \ (x \to \infty)$ for all $t >0$.
\item $\pi_G(x) = \li(x^\delta) + O(x^\delta (\log x)^{-t}) \ (x \to \infty)$ for all  $t > 0$.
\item $\frac{\pi_G(x)}{x^\delta}$ has the asymptotic expansion
$\frac{\pi_G(x)}{x^\delta} \, \simeq \, \sum_{k = 0}^\infty \frac{k!}{(\delta \log x)^{k+1}} \ (x \to \infty).$
\item $\pi_G(x)$ has the asymptotic continued fraction expansion
 \begin{align*}
{\pi_G(x)} \, \simeq \, \cfrac{\frac{x^\delta}{\delta \log x}}{1 \,-} \  \cfrac{\frac{1}{\delta \log x}}{1 \,-} \  \cfrac{\frac{1}{\delta \log x}}{1 \,-}\  \cfrac{\frac{2}{\delta \log x}}{1 \,-}\  \cfrac{\frac{2}{\delta \log x}}{1 \,-} \  \cfrac{\frac{3}{\delta \log x}}{1 \,-}\  \cfrac{\frac{3}{\delta \log x}}{1 \,-} \ \cdots \ (x \to \infty).
\end{align*}
\item $\pi_G(x)$ has the asymptotic continued fraction expansion 
\begin{align*}
\pi_G(x) \, \simeq \, \frac{x^\delta}{\delta \log x-1 \,-} \  \frac{1}{\delta \log x-3 \,-} \  \frac{4}{\delta\log x-5 \,-}\  \frac{9}{\delta\log x-7 \,-} \ \frac{16}{\delta\log x-9 \,-}\  \cdots  \ (x \to \infty).
\end{align*}
\item $\frac{\pi_G(e^{x/\delta})}{e^x}-w_n(x) \sim \frac{(n!)^2}{x^{2n+1}} \ (x \to \infty)$ for all nonnegative integers $n$.
\item $\frac{\pi_G(e^{x/\delta})}{e^x}-w_n(x) =  O\left( \frac{1}{x^{2n+1}}\right) \ (x \to \infty)$ for all nonnegative integers $n$.
\item $w_n(x)$ for every nonnegative integer $n$ is the unique rational function $w(x) \in \CC(x)$ of rational degree at most $n$ such that $\frac{\pi_G(e^{x/\delta})}{e^x}-w(x) = O \left( \frac{1}{x^{2n+1}}\right) \ (x \to \infty)$. 
\end{enumerate}
\item Suppose  that the equivalent conditions (1)(a)--(g) hold.  Then $A \delta = {\operatorname{Res}_{s = \delta} \zeta_G(s)}$,  and one has the asymptotic continued fraction expansions
 \begin{align*}
\PP_G(x) \, \simeq \, \cfrac{\frac{1}{A\delta \log x}}{1 \,-} \  \cfrac{\frac{1}{\delta \log x}}{1 \,-} \  \cfrac{\frac{1}{\delta \log x}}{1 \,-}\  \cfrac{\frac{2}{\delta \log x}}{1 \,-}\  \cfrac{\frac{2}{\delta \log x}}{1 \,-} \  \cfrac{\frac{3}{\delta \log x}}{1 \,-}\  \cfrac{\frac{3}{\delta \log x}}{1 \,-} \ \cdots  \ (x \to \infty),
\end{align*}
and
\begin{align*}
\PP_G(x)\,  \simeq \, \frac{1/A}{\delta \log x-1 \,-} \  \frac{1}{\delta \log x-3 \,-} \  \frac{4}{\delta\log x-5 \,-}\  \frac{9}{\delta\log x-7 \,-} \ \frac{16}{\delta\log x-9 \,-}\  \cdots  \ (x \to \infty),
\end{align*}
Moreover, the $w_n(x)$ for all nonnegative integers $n$ are precisely the best rational approximations of the functions $\frac{\pi_G(e^{x/\delta})}{e^x}$ and $A\PP_G(e^{x/\delta})$.
\end{enumerate}
\end{theorem}

\begin{proof}
We may assume without loss of generality, by replacing $G$ with $G_\delta$, that $\delta = 1$.  By Theorem \ref{abstractPNT}, conditions (1)(a) and (1)(b) are equivalent, and, by Lemma \ref{asympprop} and  (\ref{asex}), conditions (1)(b) and (1)(c) are equivalent.   Moreover, conditions (1)(c)--(h) are equivalent by Corollary \ref{wsim} and the proof of  Theorem \ref{maincontthm1}.  This proves statement (1).  Assuming the hypothesis of statement (2) (with $\delta = 1$), condition (1)(b) then follows from    \cite[p.\ 152]{weg}, one has $A  = {\operatorname{Res}_{s = 1} \zeta_G(s)}$ by \cite[p.\ 304]{ny},  and one has $A\PP_G(x)=  \frac{\pi_G(x)}{N_G(x)/A} =   \frac{\li(x)}{x} + O((\log x)^{-t}) \ (x \to \infty)$ for all  $t > 0$.  Therefore, as with statement (1), statement (2) follows from  (\ref{asex}), Lemma \ref{asympprop}, and Corollary \ref{wsim}.
\end{proof}

Recall that an arithmetic semigroup $G$ satisfies {\bf Axiom A}  if there exist positive constants $A$ and $\delta$ and a nonnegative constant $\eta < \delta$ such that  $N_G(x) = Ax^\delta+O(x^\eta) \ (x \to \infty)$  \cite[Chapter 4]{kno}.  Clearly, if $G$ is an arithmetic semigroup satisfying Axiom A, then the hypothesis (\ref{Nhyp}) of Theorem \ref{abstractPNT} holds.  Thus, the various conditions in the theorem hold for any arithmetic semigroup satisfying Axiom A.

\section{Weighted prime counting functions}

With respect to the asymptotic sequence $\{\frac{1}{(\log x)^n}\}$, the functions $\pi(x)$, $\Pi(x)$, $\li(x)$, and $\Ri(x)$ all have the same asymptotic continued fraction expansions, as described by the previous section.    More generally, one has the following.

\begin{proposition}
Let $n$ be a positive integer, and let $f(x)$ be any of the following functions.
\begin{enumerate}
\item $\displaystyle \Pi(x)- \sum_{k = 1}^{n-1} \frac{1}{k}\pi(\sqrt[k]{x}) = \sum_{k = n}^{\infty} \frac{1}{k}\pi(\sqrt[k]{x})$.
\item $\displaystyle \li(x)- \sum_{k = 1}^{n-1} \frac{1}{k}\Ri(\sqrt[k]{x})$.
\item $\displaystyle \mu(n)\left(\pi(x)-\sum_{k = 0}^{n-1} \frac{\mu(k)}{k}\Pi(\sqrt[k]{x})\right) = \mu(n)\left(\sum_{k = n}^{\infty} \frac{\mu(k)}{k}\Pi(\sqrt[k]{x})\right)$.
\item $\displaystyle \mu(n)\left( \Ri(x)-\sum_{k = 1}^{n-1} \frac{\mu(k)}{k}\li(\sqrt[k]{x})\right) =\mu(n)\left( \sum_{k = n}^{\infty} \frac{\mu(k)}{k}\li(\sqrt[k]{x})\right)$.
\item $\displaystyle\frac{1}{n}\pi(\sqrt[n]{x})$.
\item $\displaystyle\frac{1}{n}\Pi(\sqrt[n]{x})$.
\item $\displaystyle\frac{1}{n}\Ri(\sqrt[n]{x})$.
\item $\displaystyle \frac{1}{n}\li(\sqrt[n]{x})$.
\item $\displaystyle \frac{1}{n}\pi^*_n(x)$, where $\displaystyle \pi^*_n(x) =  \sum_{k = n}^{\infty} \pi(\sqrt[k]{x}) = \sum_{k = n}^\infty \sum_{p^k \leq x} 1$.
\item $\displaystyle \sum_{k = n}^{\infty} \frac{1}{k}\pi(\sqrt[k]{x}) = \sum_{k = n}^\infty \sum_{p^k \leq x} \frac{1}{k}$.
\end{enumerate}
One has the asymptotic continued fraction expansions
$$f(x) \simeq \frac{\sqrt[n]{x}}{\log x \, -}  \ \frac{n \log x } {\log x+n\, -} \ \frac{2n \log x} {\log x+2n\, -} \   \frac{3n\log x} {\log x+3n \, -} \ \frac{4n \log x}{\log x+4n \, -}  \  \cdots \ (x \to \infty),$$
$$\displaystyle f(x) \simeq \cfrac{\frac{ \sqrt[n]{x}}{\log x}}{1 \,-} \  \cfrac{\frac{n}{\log x}}{1 \,-} \  \cfrac{\frac{n}{\log x}}{1 \,-}\  \cfrac{\frac{2n}{\log x}}{1 \,-}\  \cfrac{\frac{2n}{\log x}}{1 \,-} \  \cfrac{\frac{3n}{\log x}}{1 \,-}\  \cfrac{\frac{3n}{\log x}}{1 \,-} \ \cdots \ (x \to \infty),$$
$$\displaystyle f(x)  \simeq \frac{\sqrt[n]{x}}{\log x-n \,-} \  \frac{n^2}{\log x-3n \,-} \  \frac{(2n)^2}{\log x-5n \,-}\  \frac{(3n)^2}{\log x-7n \,-} \ \frac{(4n)^2}{\log x-9n \,-}\  \cdots \ (x \to \infty).$$
The best rational approximations of the function $e^{-x/n}f(e^x)$ are precisely the approximants $w_k(x)$ of the Jacobi continued fraction
$$\frac{1}{x-n \,-} \  \frac{n^2}{x-3n \,-} \  \frac{(2n)^2}{x-5n \,-}\  \frac{(3n)^2}{x-7n \,-} \  \frac{(4n)^2}{x-4n \,-} \ \cdots.$$
Moreover, for all $n, k \geq 0$, one has
$$e^{-x/n}f(e^x) - w_k(x) \sim \frac{(n^k k!)^2}{x^{2k+1}},$$ and $w_k(x)$ is the unique rational function of degree at most $k$ such that $e^{-x/n}f(e^x) - w_k(x) = O(\left( \frac{(n^k k!)^2}{x^{2k+1}}\right) \ (x \to \infty)$.
\end{proposition}

\begin{proof}
For the function $\frac{1}{n}\pi(\sqrt[n]{x})$, the conclusions of the  theorem follow from Theorem \ref{maincontthm1} and Corollary \ref{maincor} by substituting $\sqrt[n]{x}$ for $x$.   Moreover, each of the ten functions differs from  $\frac{1}{n}\pi(\sqrt[n]{x})$ by a function that is  $O(\pi(\sqrt[n+1]{x})) = O\left(\frac{\sqrt[n+1]{x}}{\log x}\right) = o(\sqrt[n]{x}(\log x)^{-t})$ for all $t > 0$.  Therefore, each of the ten functions has the same asymptotic continued fraction expansions as $\frac{1}{n}\pi(\sqrt[n]{x})$.  The theorem follows.
\end{proof}

Let $s \in \CC$.  Consider the function
$$\pi_s(x) = \sum_{p \leq x} p^s, \quad \forall x > 0\index[symbols]{.st P1@$\pi_s(x)$}$$
(where, of course, $\pi_0(x) = \pi(x)$).    
 By (\ref{abelrem}) (which follows from Abel's summation formula),  for any  positive integer  $N$ one has
$$\pi_s(N) = \sum_{p \leq N} p^s   =\pi(N)N^{s}-s\int_1^N \pi(x) x^{s-1}\, dx,$$
Moreover,  if $s \neq -1$, then
$$s\int \li(x) x^{s-1}\, dx =\li(x)x^{s}-\Ei((s+1)\log x)+C$$
 for some constant $C$ of integration (depending on $s$).   The prime number theorem  with error bound, then,  yields the following.

\begin{proposition}\label{pis}
For all $s \in \CC$ with $\operatorname{Re} s > -1$ and all $t > 0$, one has
$$\pi_s(x) = \Ei((s+1)\log x) + O \left(x^{\operatorname{Re} s+1}(\log x)^{-t} \right) \ (x \to \infty).$$
\end{proposition}

Propositions \ref{pis}, \ref{EIlem}, and \ref{asympprop},  the equivalence of statements (2)(a)--(c) of Corollary \ref{wsim},  and  (\ref{stieltid}),  then,  together yield the following.

\begin{proposition}\label{sprops}
Let $s \in \CC$ with $\operatorname{Re} s > -1$.  One has the following asymptotic expansions at $\infty$.
\begin{enumerate}
\item  $\displaystyle \pi_s(x) \simeq \sum_{k = 0}^\infty \frac{k!x^{s+1}}{((s+1)\log x)^{k+1}}.$
\item   $\displaystyle \pi_s(x) \simeq  \frac{\frac{1}{s+1}x^{s+1}}{\log x \, -}  \ \frac{\frac{1}{s+1}\log x} {\log x+\frac{1}{s+1} \, -} \ \frac{\frac{2}{s+1}\log x} {\log x+\frac{2}{s+1} \, -} \   \frac{\frac{3}{s+1}\log x} {\log x+\frac{3}{s+1}\, -} \ \frac{\frac{4}{s+1}\log x}{\log x+\frac{4}{s+1} \, -}  \ \cdots$.
\item $\displaystyle  \pi_s(x) \simeq \cfrac{\frac{x^{s+1}}{(s+1)\log x}}{1 \,-} \  \cfrac{\frac{1}{(s+1)\log x}}{1 \,-} \  \cfrac{\frac{1}{(s+1)\log x}}{1 \,-}\  \cfrac{\frac{2}{(s+1)\log x}}{1 \,-}\  \cfrac{\frac{2}{(s+1)\log x}}{1 \,-} \  \cfrac{\frac{3}{(s+1)\log x}}{1 \,-}\  \cfrac{\frac{3}{(s+1)\log x}}{1 \,-} \ \cdots$.
\item $ \displaystyle \pi_s(x)   \simeq \cfrac{x^{s+1}}{(s+1)\log x-1 \,-} \  \cfrac{1}{(s+1)\log x-3 \,-} \  \cfrac{4}{{(s+1)}\log x-5 \,-}\  \cfrac{9}{{(s+1)}\log x-7 \,-} \  \cdots$.
\end{enumerate}
\end{proposition}

Since $$\sum_{n \leq x} n^{s} = \frac{x^{s+1}}{s+1} + O(x^s) \ (x \to \infty)$$ for all $s \in \CC$ with $\operatorname{Re} s > -1$,  Proposition \ref{sprops} has the following corollary.

\begin{corollary} 
Let $s \in \CC$ with $\operatorname{Re} s > -1$.   Let $F(x)$ be either of the functions $\frac{\sum_{p \leq x} p^s}{\sum_{n \leq x} n^s}$ or  $\frac{s+1}{x^{s+1}}\sum_{p \leq x} p^s$.   One has the following asymptotic expansions at $\infty$.
\begin{enumerate}
\item  $\displaystyle F(x) \simeq \sum_{k = 0}^\infty \frac{(s+1)k!}{((s+1)\log x)^{k+1}}.$
\item   $\displaystyle F(x) \simeq  \frac{1}{\log x \, -}  \ \frac{\frac{1}{s+1}\log x} {\log x+\frac{1}{s+1} \, -} \ \frac{\frac{2}{s+1}\log x} {\log x+\frac{2}{s+1} \, -} \   \frac{\frac{3}{s+1}\log x} {\log x+\frac{3}{s+1}\, -} \ \frac{\frac{4}{s+1}\log x}{\log x+\frac{4}{s+1} \, -}  \ \cdots$.
\item $\displaystyle  F(x) \simeq \cfrac{\frac{1}{\log x}}{1 \,-} \  \cfrac{\frac{1}{(s+1)\log x}}{1 \,-} \  \cfrac{\frac{1}{(s+1)\log x}}{1 \,-}\  \cfrac{\frac{2}{(s+1)\log x}}{1 \,-}\  \cfrac{\frac{2}{(s+1)\log x}}{1 \,-} \  \cfrac{\frac{3}{(s+1)\log x}}{1 \,-}\  \cfrac{\frac{3}{(s+1)\log x}}{1 \,-} \ \cdots.$
\item $ \displaystyle F(x)   \simeq \cfrac{1}{\log x-\frac{1}{s+1} \,-} \  \cfrac{\left(\frac{1}{s+1}\right)^2}{\log x-\frac{3}{s+1} \,-} \  \cfrac{\left(\frac{2}{s+1}\right)^2}{\log x-\frac{5}{s+1} \,-}\  \cfrac{\left(\frac{3}{s+1}\right)^2}{\log x-\frac{7}{s+1} \,-} \  \cdots$.
\end{enumerate}
Moreover, the best rational approximations of  $F(e^x)$ are precisely the approximants $v_n(x) = (s+1)w_n((s+1)x)$ of the Jacobi continued fraction
$$\cfrac{1}{x-\frac{1}{s+1} \,-} \  \cfrac{\left(\frac{1}{s+1}\right)^2}{x-\frac{3}{s+1} \,-} \  \cfrac{\left(\frac{2}{s+1}\right)^2}{x-\frac{5}{s+1} \,-}\  \cfrac{\left(\frac{3}{s+1}\right)^2}{x-\frac{7}{s+1} \,-} \  \cdots,$$
where the $w_n(x)$ are the best rational approximations of $\PP(e^x)$, equal to the approximants of the Jacobi continued fraction
$$\frac{1}{x-1 \,-} \  \frac{1}{x-3 \,-} \  \frac{4}{x-5 \,-}\  \frac{9}{x-7 \,-}\  \frac{16}{x-9 \,-} \  \cdots.$$
Furthermore, for all nonnegative integers $k$ one has
$$F(e^x) - v_k(x) \sim \frac{(k!)^2}{(s+1)^{2k}x^{2k+1} }\ (x \to \infty).$$
\end{corollary}

Now, the identity
\begin{align*}
\sum_{k = 0}^n  \frac{(-1)^{k}}{(k+1)z^{k+1}}  =   \frac{1}{z \, + }  \ \frac{z}{2z-1\, +} \  \ \frac{4z}{3z-2\, +}  \ \frac{9z}{4z-3 \, +} \ \frac{16z}{5z-4 \, +} \  \cdots \ \frac{n^2 z}{(n+1)z-n}, \quad \forall n \geq 0.
\end{align*}
can be verified easily by induction from the well-known recurrence relation  \cite[(1.4)]{wall}  for the numerator or denominator $C_n(z)$ of the  $n$th approximant  $w_n(z)$ of a continued fraction, which for the continued fraction
$$    \frac{1}{z \, + }  \ \frac{z}{2z-1\, +} \  \ \frac{4z}{3z-2\, +}  \ \frac{9z}{4z-3 \, +} \ \frac{16z}{5z-4 \, +} \  \cdots$$
is the recurrence relation
$$C_{k+1}(z) = ((k+1)z-k)C_k(x) -k^2 z C_{k-1}(z), \quad \forall k \geq 1.$$  
It follows that
\begin{align*}
\sum_{k = 0}^\infty \frac{(-1)^{k}}{(k+1)z^{k+1}}  =   \frac{1}{z \, + }  \ \frac{z}{2z-1\, +} \  \ \frac{4z}{3z-2\, +}  \ \frac{9z}{4z-3 \, +} \ \frac{16z}{5z-4 \, +} \  \cdots
\end{align*}
formally in $\CC[[1/z]]$, while also 
\begin{align*}
\log \left(1+\frac{1}{z} \right) = \sum_{k = 0}^\infty \frac{(-1)^{k}}{(k+1)z^{k+1}}  =   \frac{1}{z \, + }  \ \frac{z}{2z-1\, +} \  \ \frac{4z}{3z-2\, +}  \ \frac{9z}{4z-3 \, +} \ \frac{16z}{5z-4 \, +} \  \cdots
\end{align*}
for all $z \in \CC$ with $|z| >1$.
Since also $\log \left(1+\frac{1}{z} \right)$ has the asymptotic expansion
\begin{align*}
\log \left(1+\frac{1}{z} \right) \simeq \sum_{k = 0}^\infty \frac{(-1)^{k}}{(k+1)z^{k+1}}  \ (z \to \infty),
\end{align*}
it also has the asymptotic continued fraction expansion
\begin{align}\label{logas0}
\log \left(1+\frac{1}{z} \right) \simeq   \frac{1}{z \, + }  \ \frac{z}{2z-1\, +} \  \ \frac{4z}{3z-2\, +}  \ \frac{9z}{4z-3 \, +} \ \frac{16z}{5z-4 \, +} \  \cdots \ (z \to \infty).
\end{align}

By \cite[(90.1) and (90.4)]{wall}, the formal power series $F(z) =  \sum_{k = 0}^\infty \frac{(-1)^{k}}{(k+1)z^{k+1}} $ also has formal continued fraction expansions
$$F(z) =  \frac{\frac{1}{z}}{1 \,+}\  \frac{\frac{1}{z}}{2 \,+}\  \frac{\frac{1}{z}}{3 \,+}\   \frac{\frac{4}{z}}{4 \,+}\  \frac{\frac{4}{z}}{5 \,+} \  \frac{\frac{9}{z}}{6 \,+} \  \frac{\frac{9}{z}}{7 \,+} \  \frac{\frac{16}{z}}{8 \,+} \  \frac{\frac{16}{z}}{9 \,+} \  \cdots $$
and
$$F(z) =  \frac{2}{2z+1 \,-}\  \frac{1}{6z+3 \,-}\  \frac{4}{10z +5 \,-}\   \frac{9}{14z +7\,-}\  \frac{16}{18z+9 \,-} \ \cdots$$
in $\CC[[1/z]]$, converging in the $(1/z)$-adic topology. Moreover, the continued fractions above converge to $\log\left( 1+ \frac{1}{z} \right)$ for all $z \in \CC\backslash[-1,0]$.  Thus, one has the asymptotic continued fraction expansions
\begin{align}\label{logas}
\log\left(1+\frac{1}{z}\right) \simeq \cfrac{\frac{1}{z}}{1 \,+}\  \cfrac{\frac{1}{z}}{2 \,+}\  \cfrac{\frac{1}{z}}{3 \,+}\   \cfrac{\frac{4}{z}}{4 \,+}\  \cfrac{\frac{4}{z}}{5 \,+} \  \cfrac{\frac{9}{z}}{6 \,+} \  \cfrac{\frac{9}{z}}{7 \,+} \  \cfrac{\frac{16}{z}}{8 \,+} \ \cdots \  (z \to \infty)
\end{align}
and
\begin{align}\label{logas2}
\log\left(1+\frac{1}{z}\right) \simeq \frac{2}{2z+1 \,-}\  \frac{1}{6z+3 \,-}\  \frac{4}{10z +5 \,-}\   \frac{9}{14z +7\,-}\  \frac{16}{18z+9 \,-} \ \cdots \  (z \to \infty).
\end{align}

Let $s \in \RR$ be nonzero and not equal to a prime.   Let
$$G(s) = -\lim_{x \to \infty} \left(\frac{1}{s}\sum_{p \leq x}\log \left(1-\frac{s}{p}\right) +\log \log x\right)\index[symbols]{.st P2@$G(s)$}$$
and 
$$H(s)   = -\frac{1}{s}\sum_p \left(\frac{1}{p} +\frac{1}{s}\log\left(1- \frac{s}{p} \right) \right)\index[symbols]{.st P3@$H(s)$}$$ 
for all $s$ such that the give limits exist.    The following estimates are well known for $s = 1$.

\begin{proposition}\label{mertensprop}
Let $s \in \RR$ be nonzero and not equal to a prime.  Then the limits $G(s)$ and $H(s)$ exist.   Moreover,  
for all $t > 0$, one has the following.
\begin{enumerate}
\item $\displaystyle -\frac{1}{s}\sum_{p \leq x}\log \left(1-\frac{s}{p}\right) = G(s) + \log \log x + o ((\log x)^{-t}) \  (x \to \infty)$.  
\item $\displaystyle - \sum_{p \leq x} \left(\frac{1}{p} +\frac{1}{s}\log \left(1-\frac{s}{p}\right) \right) = sH(s) + o ((\log x)^{-t}) \  (x \to \infty)$.
\item $\displaystyle  \prod_{p \leq x}\left(1-\frac{s}{p}\right)^{-1/s}  =  e^{G(s)}\log x +  o ((\log x)^{-t}) \  (x \to \infty)$.
\item $\displaystyle  \prod_{p \leq x}\left(1-\frac{s}{p}\right)^{-1}  =  e^{sG(s)}(\log x)^s +  o ((\log x)^{-t}) \  (x \to \infty)$.
\end{enumerate}
\end{proposition}

\begin{proof}
We prove (1), from which the other statements readily follow.  Since $s$ is neither zero or a prime, the sum $\sum_{p \leq x}\log \left(1-\frac{s}{p}\right)$ is finite for all $x > 0$.  Let $N = \max(2,\lfloor |s|\rfloor + 1)$.   By the series expansion
\begin{align}\label{logest} \log \left(1-\frac{s}{t}\right) = -\sum_{k = 1}^\infty \frac{s^k}{k t^k}, \quad \forall s,t \in \RR: |s|<|t|, \end{align}
one has
\begin{align*} \log \left(1-\frac{s}{t}\right) = -\frac{s}{t} + O \left ( \frac{1}{t^2} \right) \ (t \to \infty). \end{align*}
It follows that the function $F(u) = \log \left(1-\frac{s}{u}\right)$ satisfies the three necessary hypotheses of Landau's theorem \cite[p.\ 201--203]{land}, and therefore one has
$$\sum_{p \leq x}\log \left(1-\frac{s}{p}\right) = A(s)+ \int_{N}^x \frac{\log \left(1-\frac{s}{t}\right)}{\log t} dt + O \left((\log x)^{-u} \right)\ (x \to \infty)$$
for all $u > 0$, for some constant $A(s)$ depending on $s$.   Now, since $|t| > |s|$ for all $t \geq N$, from (\ref{logest}) it follows that
\begin{align*}
\int_N^x \frac{\log \left(1-\frac{s}{t}\right)}{\log t} dt   = B(s)-s\log \log x - \sum_{k = 1}^\infty \frac{s^{k+1}}{k+1}\li (x^{-k}) 
\end{align*}
for some constant $B(s)$ depending on $s$.  But also
$$0< - \li (1/x)  < \frac{1}{x \log x}, $$
for all $x> 1$
and therefore
$$\left|\sum_{k = 1}^\infty \frac{s^{k+1}}{k+1}{\li (x^{-k})}\right| \leq \sum_{k = 1}^\infty \frac{|s|^{k+1}}{k+1}{\li (x^{-k})} <  \sum_{k = 1}^\infty \frac{|s|^{k+1}}{k(k+1)x^k \log x} = O \left(\frac{1}{x \log x} \right) \ (x \to \infty)$$
for all $x > 1$.
Thus we have
\begin{align*}  \int_N^x \frac{\log \left(1-\frac{s}{t}\right)}{\log t} dt = B(s)-s\log \log x   +  O \left(\frac{1}{x \log x}\right)  \ (x \to \infty) \end{align*}
and therefore
\begin{align*}
\sum_{p \leq x}\log \left(1-\frac{s}{p}\right)  & = A(s)+ B(s) -s\log \log x   +    O \left( \frac{1}{x \log x}\right) +  O \left((\log x)^{-u} \right) \ (x \to \infty) \\
 &= - sG(s) -s\log \log x  + O \left((\log x)^{-u}\right)  \ (x \to \infty),
\end{align*}
where $G(s) = -\frac{1}{s}(A(s)+B(s))$.
\end{proof}

\begin{figure}[ht!]
\includegraphics[width=80mm]{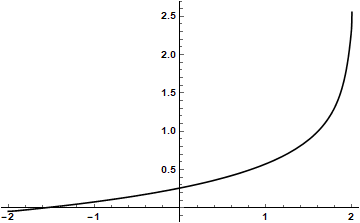}
\caption{\centering Approximation of $G(s)$ on $[-2,2)$ by the first 2000 terms of its Taylor series expansion}
\label{graphG}
\end{figure}

 One has
$$G(s) = M+ sH(s), \quad \forall s \in (-\infty,2),$$
and both $G(s)$ and $sH(s)$ are increasing on $(-\infty,2)$.
Moreover, one has
\begin{align*}
H(s)
 & = \sum_p \left(\frac{1}{2p^2} + \frac{s}{3p^3} + \frac{s^2}{4p^4} + \cdots \right) \\
 & =  \sum_{n = 0}^\infty \frac{P(n+2)}{n+2}s^n,
\end{align*}
provided that the given series converges absolutely.  In fact, for any $r \geq 0$ the sequence $ \frac{P(n+1)}{n+1}r^n$ converges monotonically to $0$ if and only if $r \leq 2$, so the radius of convergence of the series $ \sum_{n = 0}^\infty \frac{P(n+2)}{n+2}s^n$ for $s \in \CC$ is $2$, and the series converges on the entire disk $|s| \leq 2$ except at $s = 2$.  Thus,  the complex functions
$$G(s) =  M+\sum_{n = 1}^\infty \frac{P(n+1)}{n+1}s^n$$
and 
$$H(s) =  \sum_{n = 0}^\infty \frac{P(n+2)}{n+2}s^n$$
are defined on $\{s \in \CC: |s|\leq 2\}\backslash\{2\}$, and they are analytic on the open disk $\{s \in \CC: |s| < 2\}$, where convergence is absolute.  One has the following.
\begin{enumerate}
\item $G(0) =  M$.
\item $G(1) - G(0) = H = H(1)$.
\item $G(1) = \gamma = M+H$.
\item  $G(-1) = 0.07951536\ldots$. 
\item $G(-2) = 0.04723617\ldots$. 
\item $G(-1.59055202\ldots) = 0$.  
\item $G(2^-) = H(2^-) = \infty$.
\item  $G'(0)  = H(0)  =  \frac{1}{2} P(2)= 0.22612371 \ldots$. 
\end{enumerate}
In particular, the function $G(s)$ smoothly deforms the constant $M$ to the constant $\gamma = M+H$ over the interval $[0,1]$.   This  is yet another example, similar to Example \ref{stiec}(2), illustrating the fact that many well-known constants in number theory are special values of certain real or complex functions.  An approximation of the graph of $G(s)$ on $[-2,2)$ by the first 2000 terms of its Taylor series expansion is provided in Figure \ref{graphG}.

Next, we prove the following.

\begin{proposition}\label{recipthm}
Let $a > 1$.  One has the following asymptotic expansions at $ \infty$.
\begin{enumerate}
\item $\displaystyle  \sum_{x < p \leq ax} \frac{1}{p}\simeq  \sum_{n = 1}^\infty \frac{(-1)^{n-1}}{n(\log_a x)^{n}}$.
\item $\displaystyle  \sum_{x< p \leq ax} \frac{1}{p}\simeq   \frac{1}{\log_a x \, + }  \ \frac{\log_a x}{2\log_a x-1\, +} \  \ \frac{4\log_a x}{3\log_a x-2\, +}  \ \frac{9\log_a x}{4\log_a x-3 \, +} \ \frac{16\log_a x}{5\log_a x-4 \, +} \  \cdots$.
\item $\displaystyle \sum_{x< p \leq ax} \frac{1}{p}\simeq  \cfrac{\frac{1}{\log_a x}}{1 \,+}\  \cfrac{\frac{1}{\log_a x}}{2 \,+}\  \cfrac{\frac{1}{\log_a x}}{3 \,+}\   \cfrac{\frac{4}{\log_a x}}{4 \,+}\  \cfrac{\frac{4}{\log_a x}}{5 \,+} \  \cfrac{\frac{9}{\log_a x}}{6 \,+}\  \cfrac{\frac{9}{\log_a x}}{7 \,+}  \ \cfrac{\frac{16}{\log_a x}}{8 \,+}  \ \cfrac{\frac{16}{\log_a x}}{9 \,+}  \ \cdots$.
\item $\displaystyle \sum_{x< p \leq ax} \frac{1}{p}\simeq  \cfrac{2}{2\log_a x+1 \,-}\  \cfrac{1}{6\log_a x+3 \,-}\  \cfrac{4}{10\log_a x +5 \,-}\   \cfrac{9}{14\log_a x +7\,-}\ \cdots$.
\item $\displaystyle  \sum_{a^x < p \leq a^{x+1}} \frac{1}{p}\simeq  \sum_{n = 1}^\infty \frac{(-1)^{n-1}}{nx^{n}}$.
\item $\displaystyle  \sum_{a^x < p \leq a^{x+1}} \frac{1}{p}\simeq   \frac{1}{x \, + }  \ \frac{x}{2x-1\, +} \  \ \frac{4x}{3x-2\, +}  \ \frac{9x}{4x-3 \, +} \ \frac{16x}{5x-4 \, +} \  \cdots$.
\item $\displaystyle  \sum_{a^x < p \leq a^{x+1}} \frac{1}{p}\simeq  \cfrac{\frac{1}{x}}{1 \,+}\  \cfrac{\frac{1}{x}}{2 \,+}\  \cfrac{\frac{1}{x}}{3 \,+}\   \cfrac{\frac{4}{x}}{4 \,+}\  \cfrac{\frac{4}{x}}{5 \,+} \  \cfrac{\frac{9}{x}}{6 \,+}\  \cfrac{\frac{9}{x}}{7 \,+}  \ \cfrac{\frac{16}{x}}{8 \,+}  \ \cfrac{\frac{16}{x}}{9 \,+}  \ \cdots$.
\item $\displaystyle  \sum_{a^x < p \leq a^{x+1}} \frac{1}{p}\simeq  \cfrac{2}{2x+1 \,-}\  \cfrac{1}{6x+3 \,-}\  \cfrac{4}{10x +5 \,-}\   \cfrac{9}{14x +7\,-}\  \cfrac{16}{18x +9 \,-} \ \cdots$.
\end{enumerate}
\end{proposition}

\begin{proof}
Let $t = \log a > 0$.  One has the asymptotic expansion
\begin{align*}
\log (z+t)-\log z = \log \left(1+\frac{t}{z}\right) \, \simeq \, \sum_{k = 0}^\infty \frac{(-1)^{k}}{(k+1)(z/t)^{k+1}}  \ (z \to \infty)
\end{align*}
and therefore, letting $z =  \log x$, one has the asymptotic expansion
\begin{align*}
\log \log(ax)-\log \log x \,   \simeq \, \sum_{k = 0}^\infty \frac{(-1)^{k}}{(k+1)(\log_a x)^{k+1}} \ (x \to \infty)
\end{align*}
Since
\begin{align*}
\sum_{p \leq x}\frac{1}{p} = M+ \log \log x + O ((\log x)^{-t}) \ (x \to \infty),
\end{align*}
one has
$$\sum_{x < p \leq ax}\frac{1}{p}   = \sum_{p \leq ax}\frac{1}{p} - \sum_{p \leq x}\frac{1}{p} = \log \log(ax) -\log \log x + o ((\log x)^{-t}) \  (x \to \infty),$$ 
for all $t > 0$.  Therefore, the function $\sum_{x < p \leq ax}\frac{1}{p}$ has the same asymptotic expansion as $\log \log(ax)-\log \log x$.  This proves statement (1).  Statements (2)--(4) likewise follow from the asymptotic expansions (\ref{logas0}), (\ref{logas}),  and (\ref{logas2}).  Finally, statements (5)--(8) follow from statements (1)--(4) by substituting $a^x$ for $x$.
\end{proof}

It is noteworthy that the asymptotic expansions in the theorem of the function $\sum_{a^x < p \leq a^{x+1}} \frac{1}{p}$  do not depend on $a$.

By statements (1) and (2) of Proposition \ref{mertensprop}, for all $a > 1$ and all  nonzero $s <2$, the function  $$\log \prod_{x < p \leq ax} \left(1 - \frac{s}{p} \right)^{-1/s}  =- \frac{1}{s}\sum_{x < p \leq ax} \log \left(1-\frac{s}{p}\right)$$ has the same asymptotic expansions as the functions $\log \log (ax)-\log \log x$ and  $\sum_{x< p \leq ax} \frac{1}{p}$, as in Proposition \ref{recipthm}.  Thus, we have the following.

\begin{proposition}\label{lasf}
Let $a > 1$, and let $s \in \RR$ be nonzero and not equal to a prime.  The asymptotic expansions (1)--(4) of Proposition \ref{recipthm} also hold for the functions $\log(1 +\frac{1}{\log_a x}) = \log \log(ax)-
\log \log x$ and $$\log \prod_{x < p \leq ax} \left( 1-\frac{s}{p}\right)^{-1/s} = -\frac{1}{s}\sum_{x < p \leq ax} \log \left( 1-\frac{s}{p}\right),$$ 
while the asymptotic expansions (5)--(8)  of Proposition \ref{recipthm} also hold for the functions  
$\log(1 +\frac{1}{x})$ and $$\log \prod_{a^x < p \leq a^{x+1}} \left( 1-\frac{s}{p}\right)^{-1/s} = -\frac{1}{s}\sum_{a^x < p \leq a^{x+1}} \log \left( 1-\frac{s}{p}\right).$$
\end{proposition}

By \cite[(90.5)]{wall}, one has the asymptotic Jacobi continued fraction expansion
\begin{align}\label{uni2}
\log \left(\frac{z+1}{z-1}\right) =  \cfrac{2}{z \,-}\  \cfrac{\sfrac{1}{3}}{z \,-}\  \cfrac{\sfrac{4}{15}}{z \,-}\   \cfrac{\sfrac{9}{35}}{z \,-}\  \cfrac{\sfrac{16}{63}}{z  \,-} \  \cfrac{\sfrac{25}{99}}{z  \,-}\  \cfrac{\sfrac{36}{143}}{z  \,-} \cdots, \quad \forall z \in \CC\backslash [-1,1].
\end{align}
Substituting $2z+1$ for $z$,  we obtain the asymptotic Jacobi continued fraction expansion
$$\log \left( 1+ \frac{1}{z}\right) \, \simeq \,  \cfrac{1}{z+\sfrac{1}{2} \,-}\  \cfrac{\sfrac{1}{4\cdot 3}}{z+\sfrac{1}{2} \,-}\  \cfrac{\sfrac{4}{4\cdot 15}}{z+\sfrac{1}{2} \,-}\   \cfrac{\sfrac{9}{4\cdot 35}}{z+\sfrac{1}{2}\,-}\  \cfrac{\sfrac{16}{4\cdot 63}}{z+\sfrac{1}{2} \,-} \ \cdots \ (z \to \infty),$$
which is a restatement of (\ref{logas2}).  It is known that the denominator in the $n$th approximant of the Jacobi continued fraction in (\ref{uni2}) is the {\it $n$th Legendre polynomial} $P_n(z)$ \cite[p.\ 344]{wall}.  It follows that the denominator in the $n$th approximant of the continued fraction in the expansion of $\log(1+1/z)$ above is the polynomial $\widehat{P}_n(z) = P_n(2z+1)$.   The polynomial $\widehat{P}_n(z)$ is a degree $n$ polynomial with positive integer coefficients, given explicitly by
$$\widehat{P}_n(z)  = \sum_{k=0}^{n}{\binom {n}{k}}{\binom {n+k}{k}}z^{k}$$
for all $n$.  Applying Theorem \ref{wsimJ} and  Propositions  \ref{recipthm}(8) and \ref{lasf}, we obtain the following.

\begin{corollary}\label{ccc}
Let $a > 1$,  and let $s \in \RR$ be nonzero and not equal to a prime.  Let $F(x)$ denote any of the three functions $\log(1+1/x)$, $\sum_{a^x < p \leq a^{x+1}} \frac{1}{p}$, and $\log \prod_{a^x < p \leq a^{x+1}} (1 -s/p)^{-1/s}$.  Then $F(x)$ has the asymptotic Jacobi continued fraction expansion
\begin{align*}
 F(x) \, \simeq  \, \cfrac{2}{2x+1 \,-}\  \cfrac{1}{6x+3 \,-}\  \cfrac{4}{10x +5 \,-}\   \cfrac{9}{14x +7\,-}\  \cfrac{16}{18x +9 \,-} \ \cdots \ (x \to \infty),
\end{align*}
The best rational appoximations of the function $F(x)$ are precisely the approximants $w_n(x)$ of the given continued fraction for $n \geq 0$, which converge to $\log (1+1/x)$ for all $x\in \CC\backslash [-1,0]$ as $n \to \infty$.   Moreover, one has
 $$F(x)- w_n(x) \sim \frac{c_n}{x^{2n+1}}  \ (x \to \infty)$$
for all $n \geq 0$,
where $c_0 = 1$ and
$$c_n = \frac{1}{2^{2n}} \prod_{k = 1}^n \frac{k^2}{4k^2-1} = \frac{1}{(2n+2){2n+1 \choose n}{2n-1 \choose n}}$$ for all $n \geq 1$, and $w_n(x)$ is the unique rational function of degree at most $n$ such that $F(x)- w_n(x)=  O \left(\frac{1}{x^{2n+1}}\right)  \ (x \to \infty)$.
Furthermore, one has
$$w_n(x) = \sum_{k = 1}^n \frac{c_{k-1}}{P_k(2x+1)P_{k-1}(2x+1)}$$
for all $n \geq 0$, where $P_n(x)$ denotes the $n$th Legendre polynomial and
$$P_n(2x+1) = \sum_{k=0}^{n}{\binom {n}{k}}{\binom {n+k}{k}}x^{k}$$
for all $n \geq 0$.
\end{corollary}

\begin{remark}[Numerators in Corollary \ref{ccc}]
From \cite[{[1.14]}]{akh}, one can show that the numerator $R_n(x)$ of $w_n(x) = \frac{ R_n(x)}{\widehat{P}_n(x)}$ in Corollary \ref{ccc} is
$$R_n(x) =  \sum_{k = 1}^{n} a_{n,k} x^{k-1},$$
where $$a_{n,k} =  \sum_{j = k}^n \frac{(-1)^{j-k}}{j-k+1} {n \choose j}{n+j \choose j}   = {n+k \choose 2k}{2k \choose k} {}_{4}F_{3}(1,1,k-n, n+k+1; 2, k+1, k+1; 1)$$
for all $n \geq 1$ and $1 \leq k \leq n$.  
\end{remark}

\section{Measure-theoretic interpretation}

In this section we consider the meausure-theoretic aspects of Stieltes' theory of continued fractions  \cite{stie} and their consequences for asymptotic Stieltjes and Jacobi continued fraction expansions.  

Let $\mu$ be a  (positive) measure on $\RR$.  For all integers $k$, the {\bf $k$th moment of $\mu$}\index{moments $m_k(\mu)$ of a measure $\mu$} is the integral
$$m_k(\mu) = \int_{-\infty}^\infty t^k d \mu(t).$$  In its modern formulation, the {\bf Stieltjes moment problem},\index{Stieltjes moment problem} posed and motivated in \cite[No.\ 24]{stie} by Stieltjes in connection with his extensive theory of continued fractions \cite{stie}, is the problem of determining for which sequences $\{\mu_k\}_{k = 0}^\infty$ of real numbers there exists a Borel measure $\mu$ on $[0, \infty)$ such that $\mu_k = m_k(\mu)$ for all nonnegative integers $k$.  To solve this problem, Stieltjes introduced what we now call the {\bf Stieltjes transform of $\mu$},\index{Stieltjes transform $\SS_\mu$ of a measure $\mu$} which is the complex function
$$\SS_\mu(z) = \int_{-\infty}^\infty \frac{d \mu(t)}{z-t}.$$   
If the measure $\mu$ is finite, then $\SS_\mu(z)$ is analytic on $\CC\backslash \operatorname{supp} \mu$ with derivative  $\frac{d}{dz}\SS_\mu(z)  = -\int_{-\infty}^\infty \frac{d \mu(t)}{(z-t)^2}$.  Moreover, Stieltjes established in \cite{stie} the following remarkable result.   

\begin{theorem}[\cite{stie} {\cite[Chapter VII Theorems 3 and 4]{lor}} {\cite[Theorems 5.1.1 and 5.2.1]{cuyt}}]\label{mutheorem1}
Let $\{\mu_k\}_{k = 0}^\infty$ be a sequence of real numbers.   There exists a Borel measure $\mu$ on $[0,\infty)$ with infinite support (i.e., that is not a finite sum of point masses) such that $\mu_k = m_k(\mu)$ for all nonnegative integers $k$ if and only if there exists a Stieltjes continued fraction
\begin{align*}
G(z) = \cfrac{\frac{a_1}{z} }{1 \, -} \ \cfrac{\frac{a_2}{z}} {1 \, -}  \ \cfrac{\frac{a_3}{z}} {1 \, -} \ \cfrac{\frac{a_4}{z}} {1 \, -} \  \cdots \ = \ \cfrac{a_1} {z \, -} \ \cfrac{a_2} {1 \, -}  \ \cfrac{a_3} {z \, -} \ \cfrac{a_4} {1 \, -} \  \cdots,
\end{align*}
where $a_n \in \RR_{> 0}$ for all $n$, such that $G(z)$ converges $(1/z)$-adically  in $\CC[[1/z]]$ to the series  $\sum_{k = 0}^\infty\frac{ \mu_k}{ z^{k+1}}$.   The function
$\SS_\mu(z)$ is analytic on $\CC\backslash [0,\infty)$ and for all $\varepsilon > 0$ has the asymptotic expansion
$\SS_\mu(z) \simeq \sum_{k = 0}^\infty \frac{\mu_k}{z^{k+1}} \ (z \to \infty)$
over  ${\CC_\varepsilon} = \{z \in \CC: |\operatorname{Arg}(z)| \geq \varepsilon\}$.  Moreover, the Borel measure $\mu$ is unique if and only if the continued fraction  $G(z)$ converges for some $z \in \CC$, in which case it converges for all $z \in \CC \backslash [0,\infty)$ and  $\SS_\mu(z) = G(z)$ for all $z \in \CC\backslash [0,\infty)$.
\end{theorem}

Stieltjes' theorem above and Corollary \ref{wsim} have the following consequence for asymptotic Stieltjes continued fraction expansions.

\begin{theorem}\label{mutheorem}
Let $\mu$ be a Borel measure on $[0, \infty)$ with infinite support and finite moments, let $f$ be a complex function defined on an unbounded subset of $\CC$.    Then $f$ has the asymptotic expansion
$$f(z) \simeq \sum_{k = 0}^\infty \frac{m_k(\mu)}{z^{k+1}} \ (z \to \infty)$$
if and only if $f$ has the  asymptotic continued fraction expansion
$$f(z) \, \simeq \,  \cfrac{\frac{a_1}{z}} {1 \, -} \ \cfrac{\frac{a_2}{z}} {1 \, -}  \ \cfrac{\frac{a_3}{z}} {1 \, -} \ \cfrac{\frac{a_4}{z}} {1 \, -}  \  \cdots \ (z \to \infty),$$
where the $a_n \in \RR_{> 0}$ are as in Theorem \ref{mutheorem1}. 
\end{theorem}

\begin{proof}
The theorem  follows immediately Theorem \ref{mutheorem1} and the equivalence of statements (2)(a) and (2)(b) of Corollary  \ref{wsim}. 
\end{proof}

One also has the following analogues  of Theorems \ref{mutheorem1} and \ref{mutheorem}  for Borel  measures on $\RR$ and Jacobi continued fractions.  Theorem  \ref{mutheorem1bb}  is due to Hamburger \cite{ham}

\begin{theorem}[\cite{ham} {\cite[Theorems 5.1.3 and 5.2.3]{cuyt}}]\label{mutheorem1bb} 
Let $\{\mu_k\}_{k = 0}^\infty$ be a sequence of real numbers. 
There exists a Borel measure $\mu$ on $\RR$  with infinite support such that $\mu_k = m_k(\mu)$ for all nonnegative integers $k$ if and only if there exists a Jacobi continued fraction
\begin{align*}
G(z) = \frac{a_1}{z+b_1 \,-} \  \frac{a_2}{z+b_2 \,-} \  \frac{a_3}{z+b_3 \,-}  \ \cdots,
\end{align*}
where $a_n \in \RR_{> 0}$ and $b_n  \in \RR$ for all $n$, that converges $(1/z)$-adically  in $\CC[[1/z]]$ to the series $\sum_{k = 0}^\infty \frac{\mu_k }{z^{k+1}}$.  If these conditions hold, then
$\SS_\mu(z)$ is analytic on $\CC\backslash \RR$ and for all $\delta, \varepsilon > 0$ has the asymptotic expansion
$\SS_\mu(z) \simeq \sum_{k = 0}^\infty \frac{\mu_k}{z^{k+1}} \ (z \to \infty)$ over  $\CC_{\delta, \varepsilon} = \{z \in \CC : \delta \leq |\operatorname{Arg}(z)| \leq \pi- \varepsilon\}$, and one has $\SS_\mu(z) =  G(z)$ for all $z \in \CC\backslash \RR$.
\end{theorem}

\begin{theorem}\label{mutheorem1cc}
Let $\mu$ be a  Borel measure on $\RR$ with infinite support and finite moments, let $f$ be a complex function defined on an unbounded subset of $\CC$.   Then $f$ has the asymptotic expansion
$$f(z) \simeq \sum_{k = 0}^\infty \frac{m_k(\mu)}{z^{k+1}} \ (z \to \infty)$$
if and only if $f$ has the  asymptotic continued fraction expansion
$$f(z) \, \simeq \,\frac{a_1}{z+b_1 \,-} \  \frac{a_2}{z+b_2 \,-} \  \frac{a_3}{z+b_3 \,-} \  \cdots  \ (z \to \infty),$$
where the $a_n \in \RR_{> 0}$ and $b_n \in \RR$ are as in Theorem \ref{mutheorem1bb}.  
\end{theorem}

\begin{proof}
The theorem  follows immediately Theorem \ref{mutheorem1bb}  and the equivalence of statements (2)(a) and (2)(b) of Theorem \ref{wsimJ}. 
\end{proof}

We may use the results above  to provide an alternative measure-theoretic proof of Theorems \ref{maincontthm1} and \ref{gentheorem}.   Let $\gamma_0$ denote the probability measure on $[0,\infty)$ with density function $e^{-t}$, which is known as the {\bf exponential distribution with rate parameter $1$}. \index{exponential distribution $\gamma_0$ with rate parameter $1$}   The measure $\gamma_0$ is the unique Borel measure $\mu$ on $[0,\infty)$ that has $k$th moment $m_k(\mu) = \int_0^\infty t^k d\mu(t)$ equal to $k!$ for all $k \geq 0$.   Moreover, its density function $e^{-t}$ is the unique piecewise continuous function $\rho(t)$ on $[0,\infty)$ whose Mellin transform $\int_0^\infty t^{s-1}\rho(t)\, dt$ is equal to the gamma function $\Gamma(s)$. Furthermore, $\gamma_0$ is the unique Borel  measure $\mu$ on $[0,\infty)$ whose Stieltjes transform on $\CC\backslash [0,\infty)$ is equal to $-e^{-z}E_1(-z)$, where $E_1(z)$ is the exponential integral function.  In  \cite[No.\ 57]{stie}, Stieltjes proved that the Stieltjes transform $\SS_{\gamma_0}(z) = -e^{-z}E_1(-z)$ of $\gamma_0$ has the Stieltjes and Jacobi continued fraction expansions
\begin{align}
 -e^{-z}E_1(-z) & = \frac{1}{z \,-} \  \frac{1}{1 \,-} \  \frac{1}{z \,-}\  \frac{2}{1 \,-}\  \frac{2}{z \,-} \  \frac{3}{1 \,-}\  \frac{3}{z \,-} \ \cdots  \label{E11} \\
&  = \frac{1}{z-1 \,-} \  \frac{1}{z-3 \,-} \  \frac{4}{z-5 \,-}\  \frac{9}{z-7 \,-}\  \frac{16}{z-9 \,-} \  \cdots   \label{E12}
\end{align}
on $\CC \backslash [0,\infty)$.  

\begin{proof}[Proof of Theorem \ref{maincontthm1}]
Stieltjes proved in \cite[No.\ 57]{stie} that  $\gamma_0$ is the unique Borel measure associated as in Theorem \ref{mutheorem1} to the Stieltjes continued fraction (\ref{E11}).     The moments of  $\gamma_0$ are given by $m_k(\gamma_0) = k!$ for all $k$.  Therefore, by (\ref{asex2}), one has the asymptotic expansion $\PP(e^x) \simeq \sum_{k = 0}^{\infty} \frac{m_k(\gamma_0)}{x^{k+1}} \ (x \to \infty).$
The theorem therefore follows from Theorems \ref{mutheorem} and \ref{mutheorem1cc} and Remark \ref{crem}(5).
\end{proof}

We extend our analysis of the probability measure $\gamma_0$ to the probability measure $\gamma_n$ on $[0,\infty)$ with density function $\frac{t^{n}}{n!}e^{-t}$, as follows.  The measure $\gamma_n$ is the gamma distribution with  shape parameter  $n+1$  and rate parameter $1$ and is the $(n+1)$-fold convolution of the measure $\gamma_0$ with itself, where for all $a,b > 0$ the {\bf gamma distribution with shape parameter $a$ and and rate parameter $b$}\index{gamma distribution $\gamma_{a,b}$} is the  probability measure $\gamma_{a,b}$ on $[0,\infty)$ with density function $\frac{b^a}{\Gamma(a)}t^{a-1}e^{-bt}$, which has moments given by $m_k(\gamma_{a,b}) =  \frac{(a)_k}{b^{k}}$ for all $k$, where $(a)_k$ denotes the Pochhammer symbol.  By \cite[pp.\ 239--240]{cuyt}, the Stieltjes transform $\SS_{\gamma_{a,b}}(z)$ of $\gamma_{a,b}$ is equal to $-be^{-bz}E_{a}(-bz)$, where the {\bf exponential integral function}\index{exponential integral function $E_a(z)$} $E_a(z)$, for any $a \in \CC$, is the function
\begin{align*}
E_a(z) = z^{a-1}\int_{z}^\infty \frac{e^{-t}}{t^a} dt =  z^{a-1}\Gamma(1-a,z), \quad \forall z \in \CC \backslash (-\infty,0],
\end{align*}
where $$\Gamma(s,z) = \int_{z}^\infty t^{s-1}e^{-t} dt, \quad \forall s \in \CC, \quad \forall z \in \CC \backslash (-\infty,0],$$ denotes the {\bf upper incomplete gamma function},\index{upper incomplete gamma function $\Gamma(s,z)$} where the integrals are along any path of integration not crossing $(-\infty, 0]$ \cite[pp.\ 238 and 275]{cuyt}. 
Stieltjes proved in  \cite[No.\ 62]{stie} that the Stieltjes transform $\SS_{\gamma_{a,1}}(z) = -e^{-z}E_a(-z)$ of $\gamma_{a,1}$ has the continued fraction expansions
\begin{align}
-e^{-z}E_{a}(-z) & =  \cfrac{\frac{1}{z}}{1 \,-} \  \cfrac{\frac{a}{z}}{1 \,-}\  \cfrac{\frac{1}{z}}{1 \,-}\  \cfrac{\frac{1+a}{z}}{1 \,-}\  \cfrac{\frac{2}{z}}{1 \,-}\  \cfrac{\frac{2+a}{z}}{1 \,-} \ \cfrac{\frac{3}{z}}{1 \,-}\  \cfrac{\frac{3+a}{z}}{1 \,-} \ \cdots  \label{Eae}  \\
  &=  \frac{1}{z -a\,-} \  \frac{a}{z - 2-a \,-}\  \frac{2(1+a)}{z-4-a \,-}\  \frac{3(2+a)}{z -6-a \,-}  \ \frac{4(3+a)}{z -8-a \,-} \  \cdots  
\end{align}
on $\CC\backslash [0,\infty)$ (and  both continued fractions converge  formally in $\CC[[1/z]]$ to $(1/z){}_{2}F_{0}(a,1;;1/z)$).     It follows that  the Stieltjes transform $\SS_{\gamma_n}(z) = -e^{-z}E_{n+1}(-z)$ of the measure $\gamma_n = \gamma_{n+1,1}$ has the continued fraction expansions 
\begin{align}
-e^{-z}E_{n+1}(-z) & =  \cfrac{\frac{1}{z}}{1 \,-} \  \cfrac{\frac{1+n}{z}}{1 \,-}\  \cfrac{\frac{1}{z}}{1 \,-}\  \cfrac{\frac{2+n}{z}}{1 \,-}\  \cfrac{\frac{2}{z}}{1 \,-}\  \cfrac{\frac{3+n}{z}}{1 \,-} \ \cfrac{\frac{3}{z}}{1 \,-}\  \cfrac{\frac{4+n}{z}}{1 \,-} \ \cfrac{\frac{4}{z}}{1 \,-} \ \cdots \label{Enexpansion} \\
  &=  \frac{1}{z -1-n\,-} \  \frac{1(1+n)}{z -3-n \,-}\  \frac{2(2+n)}{z-5-n \,-}\  \frac{3(3+n)}{z -  7-n \,-}  \ \cdots \label{Enexpansion2}
\end{align}
on $\CC\backslash [0,\infty)$.

\begin{proof}[Proof of Theorem \ref{gentheorem}]
By  \cite[No.\ 62]{stie}, the measure $\gamma_{a,1}$ for any $a > 0$ is the unique Borel measure associated as in Theorem \ref{mutheorem1} to the  continued fraction (\ref{Eae}), and thus $\gamma_n$ to the  continued fraction (\ref{Enexpansion}).   The moments  of $\gamma_n$ are given by $m_k(\gamma_n) =  \frac{(k+n)!}{n!}$ for all $k$.  But then, from the asymptotic expansion (\ref{asex2}) of $\PP(x)$,
we easily obtain the asymptotic expansion $\PP_n(e^x) \simeq \sum_{k = 0}^\infty \frac{ m_k(\gamma_n) }{x^{k+1}} \ (x \to \infty)$, and likewise for $l_n(e^x)$.  The theorem therefore follows from Theorems \ref{mutheorem} and \ref{mutheorem1cc}.
\end{proof}

\begin{remark}[Relation to Hilbert transform of a measure]
We  provide further context for the functions $l_n(x)$ in  Theorem \ref{gentheorem}, as follows.  Let $\mu$ be a finite measure on $\RR$.  Let $\SS_\mu(x+0^+i) = \lim_{\varepsilon \to 0^+}  \SS_\mu(x+\varepsilon i)$  for all $x \in \RR$ such that the limit exists.  In fact the limit exists for all $x \in \RR$ outside a set of Lebesgue measure zero, and  its real part by definition is equal to $\pi$ times the {\bf Hilbert transform}\index{Hilbert transform ${\mathcal H}_\mu$ of a measure $\mu$}  ${\mathcal H}_\mu(x) = \frac{1}{\pi}  \operatorname{Re} \SS_\mu(x+0^+i)$ of the measure $\mu$ \cite{PSZ}. 
For all $a \in \CC$ and all $x <0$ one sets $E_a(x) := \lim_{\varepsilon \to 0^+} E_a(x+\varepsilon i)$, where also $\overline{E_{\overline{a}}(x)}  =  \lim_{\varepsilon \to 0^-} E_a(x+\varepsilon i)$.    By \cite[(14.1.9)]{cuyt}, one has
$-e^{-z}E_{n+1}(-z)  = \frac{z^n}{n!} \left(-e^{-z}E_{1}(-z)- \sum_{k = 0}^{n-1}\frac{k!}{z^{k+1}}\right)$
for all $z \in \CC \backslash \{0\}$,  which implies that
$$\SS_{\gamma_n}(x+0^+i)  = \lim_{\varepsilon \to 0^+}\left( -e^{-(x+\varepsilon i)}E_{n+1}(-(x+\varepsilon i)) \right) = l_n(e^x)- \pi i \frac{x^{n}e^{-x}}{n! },$$
and therefore $$l_n(e^x) =  \operatorname{Re}(-e^{-x}E_{n+1}(-x)) =  \operatorname{Re}\SS_{\gamma_n}(x+0^+i) = \pi {\mathcal H}_{\gamma_n}(x),$$
for all $x > 0$. It is clear, more generally, that
$\SS_{\gamma_{a,1}}(x+0^+i)  =  -e^{-x}\overline{E_a(-x)}$ and therefore
$\pi {\mathcal H}_{\gamma_{a,1}}(x) = \operatorname{Re}\SS_{\gamma_{a,1}}(x+0^+i) = \operatorname{Re}(-e^{-x}E_a(-x))$  for all $a > 0$ and all $x> 0$.
\end{remark}

We now apply the results above   to the function $\pi(ax)-\pi(bx)$ for $a > b > 0$.

\begin{example}
Let $s<t$ be real numbers.  Consider the measure $\mu$ on $[s,t]$ of density $e^{-u} \, du$.  The $n$th moment of $\mu$ is
$$m_n(\mu) = \int_{s}^t u^n e^{-u} \, du = \int_{s}^\infty u^n e^{-u} \, du -\int_{t}^\infty u^n e^{-u} \, du = e^{-s} r_n(s)-e^{-t} r_n(t),$$
where $r_n(X) = \sum_{k = 1}^n \frac{n!}{k!} X^k \in \ZZ[X]$.  Moreover, one has the asymptotic expansion
$$\frac{\li(e^{x-s})-\li(e^{x-t})}{e^x} \simeq \sum_{n = 0}^\infty \frac{m_n(\mu)}{x^{n+1}} \ (x \to \infty),$$
and the same expansion holds for the function $\frac{\pi(e^{x-s})-\pi(e^{x-t})}{e^x}$.
The Stieltjes transform of $\mu$ is
$${\mathcal S}_\mu(z) = -e^{-z}\left(E_1(-z+s)-E_1(-z+t)\right), \quad \forall z \in \CC\backslash[s,t],$$
and one has
$${\mathcal S}_\mu(x) = -e^{-x}\left(E_1(-x+s)-E_1(-x+t)\right)  = \frac{\li(e^{x-s})-\li(e^{x-t})}{e^x}, \quad \forall x \in \RR\backslash[s,t].$$
It follows that the Jacobi continued fraction expansion of ${\mathcal S}_\mu(x)$ provides  asymptotic continued fraction expansions of both $\frac{\li(e^{x-s})-\li(e^{x-t})}{e^x}$ and $\frac{\pi(e^{x-s})-\pi(e^{x-t})}{e^x}$ as $x \to \infty$.  These take the form
$${\mathcal S}_\mu(z)  =  \cfrac{b_1(s,t)}{z + a_1(s,t)\ -} \ \cfrac{b_2(s,t)}{z+a_2(s,t)\ -} \ \cfrac{b_3(s,t)}{z+a_3(s,t)\ -} \  \cdots, \quad \forall z \in \CC\backslash[s,t].$$
Under the obvious transformation, for $a > b > 0$ the expression above yields an asymptotic exansion of the form
$$\frac{\pi(ax)-\pi(bx)}{x} \simeq  \cfrac{d_1(a,b)}{\log x +c_1(a,b) \ -}\ \cfrac{ d_2(a,b)}{\log x +c_2(a,b) \ -} \ \cfrac{ d_3(a,b)}{\log x +c_3(a,b) \ -} \  \cdots  \ (x \to \infty).$$
Using Theorem \ref{wsimJ} and Mathematica, we compute the $d_1$, $c_1$, $d_2$, $c_2$, and $d_3$:
$$d_1(a,b) = a-b,$$
$$c_1(a,b) = -1+ \frac{a\log a-b\log b}{a-b},$$
$$d_2(a,b) = 1-ab\left(\frac{\log a-\log b}{a-b}\right)^2,$$
$$c_2(a,b) = -3+\frac{(a-b)^2(a \log a-b \log b)-ab(\log a-\log b)^2(2(a-b)+a \log b-b \log a)}{(a-b)((a-b)^2-ab(\log a -\log b)^2)},$$
$$d_3(a,b)  =4-ab\left(\frac{(\log a -\log b)(2(a-b)-(a+b)(\log a-\log b))}{(a-b)^2-ab(\log a -\log b)^2}\right)^2.$$
In particular, the first two best rational approximations of the function $\frac{\pi(ae^x)-\pi(be^x)}{e^x}$ are
$$w_1(x) = \cfrac{a-b}{x -1 + \frac{a\log a-b\log b}{a-b} } $$
and
$$w_2(x) =  \cfrac{a-b}{x -1 + \frac{a\log a-b\log b}{a-b} \ -}\ \cfrac{ 1-ab\left(\frac{\log a-\log b}{a-b}\right)^2}{x -3+\frac{(a-b)^2(a \log a-b \log b)-ab(\log a-\log b)^2(2(a-b)+a \log b-b \log a)}{(a-b)((a-b)^2-ab(\log a -\log b)^2)}},$$
which satisfy
$$\frac{\pi(ax)-\pi(bx)}{x} -  w_1(\log x) \sim \frac{(a-b)\left(1-ab\left(\frac{\log a-\log b}{a-b}\right)^2\right)}{(\log x)^3} \ (x \to \infty)$$
and
\begin{align*}
 \frac{\pi(ax)-\pi(bx)}{x} -w_2(\log x)  \sim \quad  \quad\quad\quad\quad\quad\quad\quad\quad\quad\quad\quad \quad\quad\quad\quad\quad\quad\quad\quad\quad\quad\quad \quad\quad\quad\\
 \quad \frac{(a-b)\left(1-ab\left(\frac{\log a-\log b}{a-b}\right)^2\right)\left( 4-ab\left(\frac{(\log a -\log b)(2(a-b)-(a+b)(\log a-\log b))}{(a-b)^2-ab(\log a -\log b)^2}\right)^2\right)}{(\log x)^5}\ (x \to \infty).
\end{align*}
We do not know a general expression for the $d_n(a, b)$ and $c_n(a,b)$.
For $a = 2$ and $b = 1$ we compute several more terms, where here $l = \log 2$:
$$d_1(2,1) = 1,$$
$$c_1(2,1) = -1+2 l = 0.3862943611\ldots,$$
$$d_2(2,1) = 1-2l^2 =  0.0390939721\ldots$$
$$c_2(2,1) = -3+ \frac{2l(-1+l)^2}{1-2l^2} = 0.3389190762\ldots$$
$$d_3(2,1) = 4-\frac{2l^2(-2+3l)^2}{(1-2l^2)^2} = 0.0321372040\ldots,$$
\begin{align*}
c_3(2,1)&  = -5+\frac{l(-2+4l-3l^2+2l^3)^2}{(1 - 2l^2)
   (2 - 12l^2 + 12l^3 -l^4)}   = 0.3453975675\ldots,
\end{align*}
\begin{align*}
d_4(2,1) & = 9-\frac{2l^2(6-18l+13l^2+l^4)^2}{(2-12l^2+12l^3-l^4)^2} = 0.0308982720\ldots,
\end{align*}
\begin{align*}
c_4(2,1) & = -7+\frac{l(12-36l+54l^2-74l^3+50l^4+l^6)^2}{
   (2 - 12l^2 + 12l^3 - l^4)(36 - 432l^2 + 
     864l^3 - 564l^4 + 72l^5 + 22l^6 + l^8)}\\
& = 0.3461893674\ldots,
\end{align*}
\begin{align*}
d_5(2,1) & = 16-\frac{2l^2(144-648l+936l^2-450l^3+120l^4-144l^5+8l^6-3l^7)^2}{(36-432l^2+864l^3-564l^4+72l^5+22l^6+l^8)^2} \\
  & = 0.0305081428\ldots.
\end{align*}
The computations above suggest that the sequences $c_n(2,1)$ and  $d_n(2,1)$ might be convergent.
\end{example}

\part{Applications of algebraic asymptotic analysis to number theory}

\chapter{The prime counting function $\pi(x)$ and related functions}

In this chapter, we apply the degree and logexponential degree formalisms to the study of the prime counting function $\pi(x)$ and various functions closely related to $\pi(x)$, including  the first and second Chebyshev functions $\vartheta(x)$ and $\psi(x)$,  and Riemann's prime counting function $\Pi(x)$.  
\section{The function $\li(x)-\pi(x)$}

Recall that $\Theta$ denotes the  {\bf Riemann constant}
$$\Theta = \sup\{\operatorname{Re} s : s \in \CC, \, \zeta(s) = 0\}.$$ 
Riemann proved in \cite{rie} that $\frac{1}{2} \leq \Theta \leq 1$ and  thus that the Riemann hypothesis holds if and only if $\Theta = \frac{1}{2}$.

By Corollary \ref{maindegree}, one has
$$\Theta = \deg(\li-\pi) = \deg (\id-\psi) = \deg(\id-\vartheta).$$
In particular,  one has $\frac{1}{2} \leq \deg(\li-\pi) \leq 1$, and the Riemann hypothesis holds if and only if $\deg(\li-\pi) = \frac{1}{2}$.
Thus, the problem of computing $\Theta$ and thereby settling the Riemann hypothesis  generalizes to the following even more difficult problem.

\begin{outstandingproblem}\label{mainproblem}
Compute $\dege(\li - \pi)$.
\end{outstandingproblem}

Motivated by the problem above,  let us define $$\Theta_k = \dege_k(\li-\pi) \in \overline{\RR}$$ for all $k$, or, equivalently, 
$$\dege(\li-\pi) = (\Theta_0, \Theta_1, \Theta_2, \ldots).$$
We call $\Theta_k$ the {\bf $k$th Riemann constant}.\index{Riemann constant $\Theta_k$}\index[symbols]{.f th@$\Theta_k$} Note, of course, that $\Theta = \Theta_0$ is the (zeroth) Riemann constant.   The Riemann constants provide fine-tuned information about the function $\li-\pi$.

Note that Corollaries \ref{montlem} and \ref{maindegree} imply that
$$\Theta_1 \leq 1.$$
 Littlewood's theorem (Theorem \ref{litt}) also has some bearing on Problem \ref{mainproblem}.   Indeed, from Littlewood's theorem it follows that $$\li(x)-\pi(x)  \neq o \left(\frac{\sqrt{x}\, \log \log \log x}{\log x} \right) \ (x \to \infty),$$
 and therefore, by Proposition \ref{oexpprop} (or by Corollary \ref{deglogprop}),  that
$$\dege(\li-\pi) \geq (\tfrac{1}{2},-1,0,1,0,0,0,\ldots).$$
It follows, then, that the Riemann hypothesis implies that $ \Theta_1 \geq -1$.    Therefore, given the Riemann hypothesis, $\Theta_1 = -1$ represents the best case scenario, given what is currently known, for what  $\Theta_1$ could be. 
 It is unknown  whether or not  \begin{align}\label{dem}
\li(x)-\pi(x) = O\left(\frac{ \sqrt{x}}{\log x}{\log \log \log x} \right) \ (x \to \infty),
\end{align}
but if this $O$ bound were to hold then one would have $ \dege(\li-\pi) =  (\frac{1}{2}, -1, 0, 1, 0, 0, 0, \ldots)$, which is therefore the current best case scenario for what $\dege(\li-\pi)$ could be.    The conjecture  (\ref{dem})  has  had at least two adherents: Stoll and Demichel  conjectured an explicit inequality \cite[(11)]{stoll} (see (\ref{dem1}))  that ``appears to hold over $10^{13}$ orders of magnitude'' and that implies (\ref{dem}).
However, in their efforts to find a reasonable fit of their data over an extensive range, they underestimate how  slowly the function $\log \log \log x$ grows.    

In 1979,  Montgomery provided motivation for a bold conjecture that implies the following.

\begin{conjecture}[{Montgomery \cite[Conjecture, p.\ 16]{mont1}}]\label{montconjec}
One has
\begin{align*}
0< \limsup_{x \to \infty} \frac{|\li(x)-\pi(x)|}{\sqrt{x}\, (\log \log \log x)^2/\log x} < \infty
\end{align*}
and therefore
$$\dege(\li-\pi)  = (\tfrac{1}{2},-1,0,2,0,0,0,\ldots).$$
\end{conjecture}

Thus, if either  Montgomery's conjecture or Stoll's and Demichel's conjecture  is true, then Littlewood's theorem is nearly the best of its kind.  See Section 14.2 for a  detailed discussion of these various conjectures regarding the function $\li(x)-\pi(x)$.

 The following proposition follows from Littlewood's theorem and a 1965 result of Grosswald \cite{gros}.

\begin{proposition}[{\cite{gros}}]\label{grossprop}
One has the following.
\begin{enumerate}
\item The  Riemann hypothesis is false if and only if $$\li(x)-\pi(x) = \displaystyle O\left(\frac{x^\Theta}{\log x}\right) \ (x \to \infty),$$ if and only if  $$x-\psi(x) = O(x^\Theta) \ (x \to \infty).$$
\item If there exists a zero of $\zeta(s)$ with real part equal to $\Theta$, then one has $$\li(x)-\pi(x) =\displaystyle \Omega_\pm\left(\frac{x^\Theta}{\log x}\right)\ (x \to \infty)$$ and $$x-\psi(x) = \Omega_\pm(x^\Theta) \ (x \to \infty).$$
\end{enumerate}
\end{proposition}

\begin{corollary}
If the Riemann hypothesis is false and yet there exists a zero of $\zeta(s)$ with real part equal to $\Theta$, then one has $\dege(\li-\pi) = (\Theta,-1,0,0,0,\ldots)$ and $\dege(\id -\psi) = (\Theta,0,0,0,\ldots)$ and $\Theta< 1$.
\end{corollary}

By Proposition \ref{grossprop} and Littlewood's theorem, we have the following.

\begin{corollary}
If $\Theta_1 > -1$, then the Riemann hypothesis is true, while if $\Theta_1 < -1$, then the Riemann hypothesis is false. 
\end{corollary}

Thus, any result that were to imply $\Theta_1 \neq -1$ would have to settle the Riemann hypothesis. 

The following theorem, which follows readily from results cited or proved earlier, summarizes much of what is known about $\dege (\li-\pi)$.

\begin{theorem}\label{riemannequivs}  One has the following.
\begin{enumerate}
\item $\Theta = \deg(\li - \pi)$.
\item $\Theta$ is the smallest real number $t$ such that $\pi(x) = \li(x) + O(x^t \log x) \ (x \to \infty)$.
\item One has $$\tfrac{1}{2} \leq \Theta \leq 1,$$ $$\Theta_1 \leq 1,$$
$$\dege(\li-\pi) \leq (\Theta,1,0,0,0,\ldots),$$
and
$$(\tfrac{1}{2},-1,0,1,0,0,0,\ldots) \leq  \dege(\li-\pi) \leq (1,-\infty,-\tfrac{3}{5},\tfrac{1}{5},0,0,0,\ldots).$$
\item The Riemann hypothesis is equivalent to $\deg (\li - \pi) = \frac{1}{2}$, and, more precisely, to
$$\dege(\li-\pi) \leq  (\tfrac{1}{2},1,0,0,0,\ldots)$$
and to
$$\dege(\li-\pi) >  (\Theta,-1,0,0,0,\ldots).$$
\item The Riemann hypothesis is false if and only if 
$$\dege(\li-\pi) \leq (\Theta,-1,0,0,0,\ldots).$$
\item If the Riemann hypothesis is true, then $-1 \leq \Theta_1 \leq 1.$
\item If the Riemann hypothesis is false, then $\Theta_1 \leq -1$.
\item The anti-Riemann hypothesis is equivalent to
$$\dege(\li-\pi) \geq (1,-\infty,-1,0,0,0,\ldots).$$
\item If there exists a zero of $\zeta(s)$ with real part $\Theta$, then the following conditions are equivalent.
\begin{enumerate}
\item The Riemann hypothesis is false.
\item $\frac{1}{2} < \Theta < 1$.
\item $\dege(\li-\pi) = (\Theta,-1,0,0,0,\ldots)$.
\item $\dege(\id -\psi) = (\Theta,0,0,0,\ldots)$.
\end{enumerate}
\end{enumerate}
\end{theorem}

 It is conceivable that improvements on the unconditional inequalities $\frac{1}{2} \leq \Theta \leq 1$ and $\Theta_1 \leq 1$ could be proved absent a proof or disproof of the Riemann hypothesis.    Thus, Problem \ref{mainproblem} specializes as follows.

\begin{outstandingproblem}\label{mainprob2}
Find unconditional bounds on $ \Theta$ and 
$$\Theta_1 = \dege_1 (\li - \pi) = \deg((\li(e^x)-\pi(e^x))e^{-\Theta x}) = \limsup_{x \to \infty} \frac{\log |\li(e^x)-\pi(e^x)|-\Theta x}{\log x}$$
stronger than $\frac{1}{2} \leq \Theta \leq 1$ and $\Theta_1 \leq 1$, respectively.
\end{outstandingproblem}

The following proposition, whose verification is trivial, shows how Conjectures \ref{eurekaconjecture} and \ref{eurekaconjecture2} relate to Problem  \ref{mainprob2}.

\begin{proposition}  One has the following.
\begin{enumerate}
\item Conjecture \ref{eurekaconjecture} is equivalent to  
$$\Theta = \tfrac{1}{2} \text{ and } \Theta_1 = -1,$$
and it implies that $\Theta_2 \geq 0$.
\item Conjecture \ref{eurekaconjecture2} is equivalent to  
$$\Theta= \tfrac{1}{2}, \ \Theta_1 = -1, \text{ and }\Theta_2 = 0,$$
and it implies that $\Theta_3 \geq 1$.
\end{enumerate}
\end{proposition}

Thus, Conjecture \ref{eurekaconjecture} says that the Riemann hypothesis holds and $\Theta_1$ is as small as it could be given the Riemann hypothesis and Littlewood's theorem.    Likewise, Conjecture \ref{eurekaconjecture2} says that $\Theta_2$ is as small as it could be given the Riemann hypothesis, Littlewood's theorem, and Conjecture \ref{eurekaconjecture}.    Each of these three conjectures, starting with the Riemann hypothesis, adds further confidence that Littlewood's theorem is  ``the best possible of its kind.'' To that very question, Hardy and Littlewood wrote  in 1916, ``We are naturally not prepared to express any very definite opinion on that point''  \cite[p.\ 127]{har1}.  Note, however, that Montgomery's Conjecture  \ref{montconjec} is far stronger than Conjectures \ref{eurekaconjecture} and \ref{eurekaconjecture2} and thus provides some independent theoretical support for the latter conjectures.

\begin{problem}
If possible, recast Conjectures \ref{eurekaconjecture} and \ref{eurekaconjecture2} as statements about the Riemann zeta function.
\end{problem}

\begin{remark}[Weaker versions of Littlewood's theorem]
As remarked by Littlewood  \cite{litt},  his Theorem \ref{litt} was a generalization of a prior result by M.\ E.\ Schmidt, namely, 
$$\pi(x)-\li(x)+\frac{1}{2}\li(\sqrt{x}) = \Omega_{\pm} \left( \frac{\sqrt{x}}{\log x} \right) \ (x \to \infty),$$
which by Corollary \ref{RAprop} is equivalent to 
$$\pi(x)-\Ri(x) = \Omega_{\pm} \left( \frac{\sqrt{x}}{\log x} \right) \ (x \to \infty)$$
and to
$$\pi(x)-\Ri(x) = \Omega_{\pm} \left( \li(x)-\Ri(x) \right) \ (x \to \infty).$$
A slightly different result, namely, 
$$\pi(x)-\li(x) = \Omega_{\pm} \left( \frac{\sqrt{x}}{\log x} \right) \ (x \to \infty),$$
is proved in \cite[Theorem 11.15]{bate}, and also follows from Littlewood's theorem.
\end{remark}

Now, we consider the anti-Riemann hypothesis.

\begin{proposition}\label{ARprop}
Write $$|\li(x)-\pi(x)| = xe^{-F(\log x)},$$
so that
$$F(x) = x-\log|\li(e^x)-\pi(e^x)|,$$
for all $x \geq 0$.  One has the folllowing.
\begin{enumerate}
\item $\displaystyle \lim_{x \to \infty} F(x)  = \infty$.
\item If the anti-Riemann hypothesis is false, then $\dege(1/F) = (-1,0,0,0,\ldots)$.
\item  If the anti-Riemann hypothesis is true, then $\dege(1/F) \geq (-1,0,0,0,\ldots)$, and one has
$$\dege(\li-\pi) = (1,-\infty,\deg(1/F), \dege_1(1/F), \dege_2(1/F), \ldots)$$
and $\Theta_2 = \deg (1/F) \in [-1,-\tfrac{3}{5}]$.
\end{enumerate}
\end{proposition}

\begin{proof}
By the definition of $F(x)$, one has
$$  \frac{-\log|\log |\li(e^x)-\pi(e^x)|-x|}{\log x} = \frac{-\log |F(x)|}{\log x}.$$
Since $\li(x)-\pi(x) = o(x) \ (x \to \infty)$, statement (1) holds, and one has
 $$\limsup_{x \to \infty} \frac{-\log|\log |\li(e^x)-\pi(e^x)|- x|}{\log x} = \deg (1/F).$$

Suppose that the anti-Riemann hypothesis is false, i.e.,  $\Theta < 1$.  Then, for some $t> 0$, one has $e^{-F(\log x)} \leq x^{-t}$, so that $0 \leq \frac{1}{F(x)} \leq \frac{1}{t x}$, for all $x \gg 0$, whence $\dege \frac{1}{F} \leq (-1,0,0,0,\ldots)$.   Suppose that $\dege \frac{1}{F}< (-1,0,0,0,\ldots)$, so that  $\lim_{x \to \infty} \frac{x}{F(x)} = 0$.  Then
$-F(\log x)\leq -\log x$, and therefore $$|\li(x)-\pi(x)| = xe^{-F(\log x)} \leq xe^{-\log x} = 1,$$
for all $x \gg 0$, which is a contradiction.  This proves (2).

If, on the other hand, the anti-Riemann hypothesis is true, then  $1/F = (\li-\pi)_{(2)}$, where  $f_{(k)}$ is defined as in the definition of $\dege f$.  Statement (3) follows.
\end{proof}

We call the constant $$\Theta_{-1} = \underline{\deg}(x-\log|\li(e^x)-\pi(e^x)|)= -\deg \frac{1}{x-\log|\li(e^x)-\pi(e^x)|}$$ the {\bf anti-Riemann constant}.\index{anti-Riemann constant $\Theta_{-1}$}\index[symbols]{.f te@$\Theta_{-1}$}    

\begin{outstandingproblem}
Compute the anti-Riemann constant $\Theta_{-1}$.  More generally, compute $ \underline{\dege}(x-\log|\li(e^x)-\pi(e^x)|) = -\dege\frac{1}{x-\log|\li(e^x)-\pi(e^x)|}$.
\end{outstandingproblem}

The following is a consequence of Littlewood's theorem.

\begin{proposition}
One has $\li(x) - \pi(x) = 0$ for an unbounded set of positive real numbers $x$.  Consequently,  one has $\underline{\dege}(\li-\pi) = (-\infty,-\infty,-\infty,\ldots)$ and  $\dege(x-\log|\li(e^x)-\pi(e^x)|) = (\infty,\infty, \infty,\ldots)$.  
\end{proposition}

Proposition \ref{ARprop} has the following corollary.

 \begin{corollary} One has the following.
 \begin{enumerate}
\item   $\Theta_{-1} \in [\tfrac{3}{5},1]$.  
\item  If the anti-Riemann hypothesis is false, then $\Theta_{-1} = 1$ and $$\underline{\dege}(x-\log|\li(e^x)-\pi(e^x)|) = (1,0,0,0,\ldots).$$
\item  If the anti-Riemann hypothesis is true, then $\Theta_{-1} =- \Theta_2$ and $$(1,-\infty,-1,0,0,0,\ldots) \leq \dege(\li-\pi) = (1,-\infty,\Theta_2,\Theta_3,\Theta_4, \ldots).$$
\end{enumerate}
\end{corollary}

 Thus, the anti-Riemann constant holds more interest if the anti-Riemann hypothesis is true.

\begin{proposition}\label{antiR}
The anti-Riemann constant $\Theta_{-1}$ is given by
\begin{align*}
\Theta_{-1} &  = -\inf\left\{t \in \RR: \li(x)-\pi(x) = O\left(xe^{-(\log x)^{-t}} \right) \ (x \to \infty)\right\}\\
 & = \sup\left\{t \in \RR: \li(x)-\pi(x) = O\left(xe^{-(\log x)^{t}} \right) \ (x \to \infty)\right\} \\
& =  \sup\left\{t \in \RR: x-\psi(x) = O\left(xe^{-(\log x)^{t}} \right) \ (x \to \infty)\right\} \\
& = \sup\left\{t \in \RR: x-\vartheta(x) = O\left(xe^{-(\log x)^{t}} \right) \ (x \to \infty)\right\}.
\end{align*}
Consequently, one has
$$\Theta_{-1}  = \underline{\deg} (x-\log|e^x-\psi(e^x)|) = \underline{\deg}(x-\log|e^x-\vartheta(e^x)|).$$
Moreover, $\Theta_{-1}$ is equal the supremum of all $\alpha \in (0,1]$ such that there exists some $t_0 > 0$ such that $\zeta(s) \neq 0$ for all $s = \sigma + it$ with $\sigma \geq 1-(\log t)^{1-1/\alpha}$  and $t \geq t_0$.
\end{proposition}

\begin{proof}
Let $t \in \RR$.  Then $$0\leq \frac{1}{x-\log|\li(e^x)-\pi(e^x)|} \leq x^t,$$ for all $x \gg 0$, is equivalent to  
$$|\li(x)-\pi(x)| \leq xe^{-(\log x)^{-t}}.$$
 The first two equations of the proposition follow.  Finally, the rest of the proposition follows from Theorem \ref{antiR2}. 
\end{proof}

We make the following conjecture.

\begin{conjecture}\label{arconst}
The anti-Riemann constant $\Theta_{-1}$ is equal to $1$.   Equivalently, one has $$\li(x)-\pi(x) = O\left(xe^{-(\log x)^t} \right) \ (x \to \infty)$$
for all $t <1$.  Equivalently still, for every $\beta > 0$ there exists some $t_0 > 0$ such that $\zeta(s) \neq 0$ for all $s = \sigma + it$ with $\sigma \geq 1-\frac{1}{(\log t)^\beta}$  and $t \geq t_0$.
\end{conjecture}

 Obviously, Conjecture \ref{arconst} is implied by the Riemann hypothesis.  More generally, the conjecture is true provided that the anti-Riemann hypothesis is false, and our bet is that the conjecture is true even if the anti-Riemann hypothesis is true.   It is thus conceivable that the conjecture could be proved true absent a proof or disproof of the anti-Riemann hypothesis.  However, a disproof of the conjecture would require a proof of the anti-Riemann hypothesis, and it is in that scenario where the constant $\Theta_{-1}$ holds the greatest interest.

The following are natural problems that arise when considering the various approximations of $\pi(x)$.

\begin{problem} \
\begin{enumerate}
\item By (\ref{asex2}) and Proposition \ref{degeprop}(2), one has
\begin{align*}
\inf\{\dege (\pi-f) : f \in \mathbb{L}\} & \leq  \inf\{(1,-n,0,0,0,\ldots): n \in \ZZ_{>0}\} \\
& = (1,-\infty,0,-1,0,0,0,\ldots).
\end{align*}
Does equality hold?
\item Let $\mathbb{H}$ be the universal Hardy adjoinable Hardy field, i.e., the Hardy field of all universally Hardian (germs of) functions,  as defined in Section 7.3. Since $\li \in  \mathbb{H}$,  one has
$$\inf\{\dege (\pi-f) : f \in \mathbb{H}\} \leq \dege(\li-\pi).$$
Does equality hold?
\item Is $\Ri$ Hardian?   Does $\Ri$ lie in $\mathbb{H}$?
\end{enumerate}
\end{problem}

Because of its importance in analytic number theory, we take the function $\li-\pi$ as a ``logexponential primitive,'' that is, we seek (in Sections 9.2, 9.3, and 10.1) to express $\dege f$ in terms of $\dege(\li-\pi)$ for as many functions $f$ of number-theoretic interest as possible.   In all such instances we encounter,   $\dege(\li-\pi)$ can also be expressed in terms of $\dege f$.

 Here is an example involving the harmonic numbers.    Let $t,a\in \RR$ with $a >0$.    It is straightforward to show that
$$\dege(\li-\pi) = \dege(\li(ax)-\pi(ax)) = \dege(\li(an)-\pi(an)).$$
Moreover, the author proved in 2021  \cite{ell2} that
 \begin{align}\label{betat}
 \frac{\li (e^t n)}{e^t} = \int_{\mu e^{-t}}^x \frac{du}{t+\log u} = \sum_{\mu e^{-t}  \leq k < n}   \frac{1}{H_k-\gamma+t} + \beta(t) - \frac{1+o(1)}{12n(\log n)^2} \ (n \to \infty)
 \end{align}
for some explicit constant $\beta(t) \in (0, \frac{1}{\log \mu})$ depending on $t$,  where $\mu$ is the Ramanujan--Soldner constant.   (The given sum above is a ``doubly discretized'' version of the logarithmic integral.)
From these results, we deduce the following.

\begin{proposition}\label{harmonicprime2}
For all $t \in \RR$, one has
$$\dege\left(\pi(e^{t} n) - e^{t}  \sum_{\mu e^{-t}  \leq k < n} \frac{1}{H_k -\gamma+t} \right) = \dege(\pi(e^{t} n)-\li(e^{t} n)) = \dege(\li-\pi).$$
\end{proposition}

\begin{remark}[Equivalents of the Riemann hypothesis using the harmonic numbers]\label{harmonicprime}
 The harmonic numbers appear in several reformulations of the Riemann hypothesis.  Most notably, in 2002, J.\ Lagarias proved \cite{lag} that the Riemann hypothesis is equivalent to
$$\sigma(n) < H_n + \exp(H_n) \log (H_n), \quad \forall  n > 1$$
and to
$$\sigma(n) < \exp(H_n) \log (H_n), \quad \forall  n > 60.$$
Let $t \in \RR$.  In  \cite{ell2},  the author used (\ref{betat}) to show that the Riemann constant $\Theta$ is the smallest real number $\alpha$ such that $$\pi(e^{t} n) = e^{t} \sum_{\mu e^{-t}  \leq k < n} \frac{1}{H_k -\gamma+t}  +  O\left(       n^\alpha H_n\right) \ (n \to \infty).$$
It follows that the Riemann hypothesis is equivalent to
$$\pi(e^{t} n) = e^{t}  \sum_{\mu e^{-t}  \leq k < n} \frac{1}{H_k -\gamma+t}  +  O\left(\sqrt{n} \, H_n\right) \ (n \to \infty).$$
Moreover,  the $O$ constant can be made explicit, in terms of $t$.  For example,  by \cite[Corollary 9]{ell2}, the Riemann hypothesis holds if and only if
$$\PP(e^\gamma n) =  \frac{1}{n} \sum_{k = 1}^{n} \frac{1}{H_k} + O \left(\frac{H_n}{ \sqrt{n}}\right) \ (n \to \infty),$$ 
if and only if
 $$\left|\PP(e^\gamma n) -\frac{1}{n} \sum_{k = 2}^{n-1}   \frac{1}{H_k}   \right| <  \frac{1}{33} \frac{H_n}{\sqrt{n}}+\frac{3}{2n}, \quad \forall n \geq 1.$$   
\end{remark}

\section{Riemann's functions $\Pi(x)$ and $\Ri(x)$}

An important application of the logexponential degree formalism is in relating $\dege(\li-\pi)$ to the logexponential degrees of other important number-theoretic functions.  In this section, we provide a simple application of the logexponential degree formalism to  Riemann's prime counting function $$\Pi(x) = \sum_{n = 1}^\infty \sum_{p^n \leq x} \frac{1}{n} = \sum_{n = 1}^\infty \frac{1}{n}\pi(x^{1/n})$$ and his remarkable approximation  $$\Ri(x)=\sum_{n=1}^\infty \frac{ \mu(n)}{n} \li(x^{1/n})$$ to both $\pi(x)$ and $\Pi(x)$. 

 One has the following.

\begin{proposition}\label{etaconj4}
One has  
$$\dege(\li-\pi)  = \dege(\Ri-\pi) = \dege(\li-\Pi) = \dege(\Ri-\Pi),$$
and
$$\dege(\li-\Ri) = \dege(\Pi-\pi)=  (\tfrac{1}{2},-1,0,0,0,\ldots) <\dege(\li-\pi).$$
\end{proposition}

\begin{proof}
Corollary \ref{elll} and Theorem \ref{riemannequivs}(3) imply that
$$\dege(\li-\Ri) = (\tfrac{1}{2},-1,0,0,0,\ldots) < (\tfrac{1}{2},-1,0,1,0,\ldots) \leq \dege(\li-\pi).$$ Therefore, by Theorem \ref{diffpropexp}(1), one has  $$\dege(\Ri-\pi) = \dege((\li-\pi)-(\li-\Ri)) = \dege(\li-\pi).$$  The other equalities are proved similarly.
\end{proof}

Since $\dege(\Ri-\pi) = \dege(\li-\pi)$, Riemann's approximation $\Ri(x)$, with respect to the logexponential degree formalism, is no better an approximation to the prime counting function $\pi(x)$ than is the  approximation $\li(x)$.

Recall that $\mathcal{E}(f,g)$ for any $f, g\in \RR^{\RR_\infty}$ denotes the smallest nonnegative integer $k$ such that $\dege_k f \neq \dege_k g$, where $\mathcal{E}(f,g) = \infty$ if no such integer $k$ exists, i.e., if $\dege f = \dege g$.  Likewise, $\mathcal{L}(f,g)$ denotes the smallest nonnegative integer $k$ such that $\degl_k f \neq \degl_k g$, where  $\mathcal{L}(f,g) = \infty$ if no such integer $k$ exists, i.e., if $\degl f = \degl g$. The following result is clear.

\begin{proposition}\label{etaconjcor}  One has the following.
\begin{enumerate}
\item $\mathcal{E}(\li-\pi,\li-\Ri) =  \mathcal{L}(\li-\pi,\li-\Ri) \leq 3$.
\item $\mathcal{E}(\li-\pi,\li-\Ri)$ is the smallest nonnegative integer $k$  such that $\dege_k(\li-\pi) > \dege_k(\li-\Ri)$ (or, alternatively, such that $\degl_k(\li-\pi) > \degl_k(\li-\Ri)$).
\item  $\mathcal{E}(\li-\pi,\li-\Ri)$ is the smallest nonnegative integer $k$ such that
$$\log \frac{|\li(e^x)-\pi(e^x)|}{ \li(e^x)-\Ri(e^x)} = \Omega_+(\log^{\circ k} x) \ (x \to \infty),$$ that is, such that
$$\limsup_{x \to \infty} \frac{\log |\li(e^x)-\pi(e^x)|-\log (\li(e^x)-\Ri(e^x)) }{\log^{\circ k} x} > 0.$$
\item The Riemann hypothesis is equivalent  to $\mathcal{E}(\li-\pi,\li-\Ri) = 1, 2,\text{or } 3$ 
and to $$\displaystyle \limsup_{x \to \infty} \frac{\log \frac{|\li(e^x)-\pi(e^x)|}{ \li(e^x)-\Ri(e^x)}}{x} = 0.$$
\item Conjecture  \ref{eurekaconjecture} is equivalent to $\mathcal{E}(\li-\pi,\li-\Ri) = 2 \text{ or } 3$ and to $$\displaystyle \limsup_{x \to \infty} \frac{\log \frac{|\li(e^x)-\pi(e^x)|}{ \li(e^x)-\Ri(e^x)}}{\log x} = 0.$$
\item Conjecture \ref{eurekaconjecture2} is equivalent to $\mathcal{E}(\li-\pi,\li-\Ri) = 3$
and to
$$\displaystyle \limsup_{x \to \infty} \frac{\log \frac{|\li(e^x)-\pi(e^x)|}{ \li(e^x)-\Ri(e^x)}}{\log \log x} = 0.$$
\item Montgomery's  Conjecture \ref{montconjec} implies $\mathcal{E}(\li-\pi,\li-\Ri) = 3$
 and 
$$\displaystyle \limsup_{x \to \infty} \frac{\log \frac{|\li(e^x)-\pi(e^x)|}{ \li(e^x)-\Ri(e^x)}}{\log \log \log x} = 2.$$
\end{enumerate}
\end{proposition}

Conjectures  \ref{eurekaconjecture} and \ref{eurekaconjecture2}, then,  are motivated by the following problem.

\begin{outstandingproblem}\label{lipip}
Is $\mathcal{E}(\li-\pi,\li-\Ri)$ equal to $0$, $1$, $2$, or $3$?
\end{outstandingproblem}

We remark that, in Proposition \ref{etaconjcor} and Problem \ref{lipip} above,  the function $\li-\Ri$ can  be replaced with $\frac{\sqrt{x}}{\log x}$.

\section{The first and second Chebyshev functions $\vartheta(x)$ and $\psi(x)$}

Recall  from Section 3.6 the first Chebyshev function  $$\vartheta(x) = \sum_{p \leq x} \log p$$  and the second Chebyshev function 
$$\psi(x) = \sum_{n = 1}^\infty \sum_{p^n \leq x} \log p = \sum_{n = 1}^\infty \vartheta(x^{1/n}),$$ 
which satisfy $$\vartheta(x) \sim \psi(x) \sim \pi(x) \log x \ (x \to \infty).$$
By Theorem \ref{psiestprop}, one has
$$x-\psi(x) = O(x^\Theta (\log x)^2) \ (x \to \infty),$$
and, by a 1914 result of  Littlewood  \cite{litt}, one has
$$x-\psi(x) = \Omega_{\pm} (\sqrt{x} \, \log \log \log x)  \ (x \to \infty).$$
Consequently, one has the following.

\begin{proposition}\label{etaconj3}
One has the following.
\begin{enumerate}
\item One has $$(\tfrac{1}{2},0,0,1,0,0,0,\ldots) \leq \dege (\id-\psi) \leq (\Theta,2,0,0,0,\ldots).$$
\item One has $$ \psi(x)-\vartheta(x) \sim \sqrt{x} \ (x \to \infty)$$
and therefore
$$\dege(\psi-\vartheta) =  (\tfrac{1}{2},0,0,0,\ldots).$$ 
\item One has 
$$\dege(\id-\psi) = \dege(\id -\vartheta).$$
\end{enumerate}
\end{proposition}

\begin{proof}
Statement (1) is follows from  the $O$ and $\Omega_{\pm}$ bounds on $x-\psi(x)$ noted earlier.  Statement (2)  follows from Propostion \ref{RAprop2} and the prime number theorem, and statement (3) follows from statements (1) and (2) and Theorem \ref{diffpropexp}(1).
\end{proof}

From Propositions \ref{pitheta} and  \ref{kregvar}, Theorems  \ref{diffpropexp} and \ref{infpropexp},  and Karamata's integral theorem (Theorem \ref{karam}), we deduce the following.

\begin{theorem}\label{lithetapsi}
Let $r$ be continuous on $[N,\infty)$ for some $N > 1$ and regularly varying of index $d$ (e.g., $r \in \mathbb{L}$ with $\deg r = d$) for some $d \in \RR$.  One has the following.
\begin{enumerate}
\item If $$\li(x)-\pi(x) = O(r(x)) \ (x \to \infty),$$
then
$$ \li(x)-\pi(x) -\frac{x-\vartheta(x)}{\log x} =  O\left( \frac{r(x) }{\log x}\right)\ (x \to \infty).$$ 
\item If $$x-\vartheta(x) = O(r(x)) \ (x \to \infty),$$
then
$$ \li(x)-\pi(x) -\frac{x-\vartheta(x)}{\log x} =  O\left( \frac{ r(x) }{(\log x)^2} \right)\ (x \to \infty).$$ 
\item One has $$\li(x)-\pi(x) = O(r(x)) \ (x \to \infty)$$
if and only if 
$$ x-\vartheta(x) = O(r(x) \log x) \ (x \to \infty).$$ 
\item One has $$\dege(\id-\vartheta) = \dege(\li-\pi)+ (0,1,0,0,0,\ldots)$$ and
$$\dege \left(\li-\pi - \frac{\id-\vartheta}{\log } \right)  \leq \dege(\li-\pi)+ (0,-1,0,0,0,\ldots).$$
\end{enumerate}
\end{theorem}

\begin{proof} 
Suppose that  $$\li(x)-\pi(x) = O( r(x)) \ (x \to \infty).$$
Then $d := \deg r  \geq \Theta \geq \frac{1}{2}$.  By Proposition \ref{pitheta} and Karamata's integral theorem, for all $x \geq N$ one has
\begin{align*}
\li(x)-\pi(x) -  \frac{x-\vartheta(x)}{\log x}  & = \frac{1}{\log x}\int_0^x \frac{\li(t)-\pi(t)}{t} \, dt  \\
 & = \frac{1}{\log x}\int_0^N \frac{\li(t)-\pi(t)}{t} \, dt + \frac{1}{\log x}\int_N^x \frac{\li(t)-\pi(t)}{t} \, dt  \\
& = O\left(\frac{1}{\log x}\right)  + O\left(\frac{1}{\log x} \int_N^x \frac{|r(t)|}{t} \, dt\right) \ (x \to \infty)  \\
& =  O\left( \frac{r(x) }{\log x} \right) \ (x \to \infty).
\end{align*}
This proves statement (1).  Statement (3) follows from statement (1) and Theorem \ref{infpropexp}.  It follows that $\deg(\id -\vartheta) = \deg(\li-\pi)$.   To prove statement (2),  suppose that  $$x-\vartheta(x)= O(r(x)) \ (x \to \infty).$$  Then $d := \deg r  \geq \Theta \geq \frac{1}{2}$ and, again by Proposition \ref{pitheta} and Karamata's integral theorem, for all $x \geq N$ one has
\begin{align*}
\li(x)-\pi(x) -\frac{x-\vartheta(x)}{\log x} & =  \li(N)-\pi(N) -\frac{N-\vartheta(N)}{\log N}+   \int_N^x \frac{t-\vartheta(t)}{t \log^2 t} \, dt \\
& = O(1) + O\left( \int_N^x \frac{|r(t)|}{t(\log t)^2} \, dt \right) \ (x \to \infty) \\
& =  O\left(\frac{r(x)}{(\log x)^2}\right)\ (x \to \infty).
\end{align*}
This proves statement (2).   Finally,   statement (3) follows from statements (1) and (2), which then together imply statement (4), by Theorem \ref{infpropexp} and Proposition \ref{aspropexoexp}(1).
\end{proof}

We  note the following.

\begin{proposition}[{cf.\  \cite[Theorem 13.2]{mont}}]\label{lipitheta}
One has
$$\li(x)-\pi(x) = \frac{x-\vartheta(x)}{\log x}+O\left(\frac{x^\Theta}{(\log x)^2} \right) \ (x \to \infty)$$
and therefore
\begin{align*}
\dege\left(\li-\pi - \frac{\id-\vartheta}{\log} \right) \leq (\Theta,-2,0,0,0,\ldots).
\end{align*}
\end{proposition}

\begin{proof}
The proof is a straightforward generalization of the proof of \cite[Theorem 13.2]{mont}, which assumes the Riemann hypothesis $\Theta = \frac{1}{2}$.
\end{proof}

The following proposition expresses $\dege(\psi(x)-\vartheta(x)-\sqrt{x})$ in terms of $\dege(\id-\psi) = \dege(\id-\vartheta)$.  Curiously, the expression depends on whether $\Theta <\tfrac{2}{3}$ or $\Theta \geq \tfrac{2}{3}$.

\begin{proposition}\label{lipithetap}  One has the following.
\begin{enumerate}
\item One has
$$\psi(x)-\vartheta(x)-\psi(\sqrt{x})    \sim x^{1/3}  \ (x \to \infty)$$
and therefore
$$\dege(\psi(x)-\vartheta(x)-\psi(\sqrt{x})) = (\tfrac{1}{3},0,0,0,\ldots).$$
\item One has
$$\dege(  \psi(x)-\vartheta(x)-\sqrt{x})   = \begin{cases}    (\tfrac{1}{3},0,0,0,\ldots) &  \text{if } \Theta < \tfrac{2}{3} \\
  \dege(\sqrt{x}-\psi(\sqrt{x})) \geq  (\tfrac{1}{3},0,0,1,0,0,0,\ldots) &  \text{if } \Theta \geq \tfrac{2}{3},
\end{cases}$$
where $$ \dege(\sqrt{x}-\psi(\sqrt{x}))   = \dege(\id -\psi)+(-\tfrac{\Theta}{2},0,0,0,\ldots).$$
\end{enumerate}
\end{proposition}

\begin{proof}
To prove statement (1), simply note that
\begin{align*}
\psi(x)-\vartheta(x)-\psi(\sqrt{x}) & =  \sum_{k = 2}^\infty \vartheta(x^{1/k})- \sum_{k = 1}^\infty \vartheta(x^{1/(2k)})  \\
& = \vartheta(x^{1/3}) + \vartheta(x^{1/5})+\vartheta(x^{1/7})+\vartheta(x^{1/9})+\cdots \\
& \sim  \vartheta(x^{1/3})  \ (x \to \infty).
\end{align*}
To prove statement (2), note first that
$$\psi(x)-\vartheta(x)-\sqrt{x} =(\psi(x)-\vartheta(x)-\psi(\sqrt{x}))- (\sqrt{x}-\psi(\sqrt{x})),$$
where, by Proposition \ref{etaconj3} and Theorem \ref{fgie00}, one has
\begin{align*}
(\tfrac{\Theta}{2},0,0,1,0,0,0,\ldots)  & \leq  \dege(\sqrt{x}-\psi(\sqrt{x}))  \\ &  = \dege(\id -\psi)+(-\tfrac{\Theta}{2},0,0,0,\ldots) \\ & \leq (\tfrac{\Theta}{2},2,0,0,0\ldots).
\end{align*}
It follows that
$$\dege(\psi(x)-\vartheta(x)-\sqrt{x}) = \dege(\sqrt{x}-\psi(\sqrt{x})) $$
provided that 
$$(\tfrac{\Theta}{2},0,0,1,0,0,0,\ldots) > (\tfrac{1}{3},0,0,0,\ldots),$$
that is, provided that $\Theta \geq \tfrac{2}{3}$, while, on the other hand, one has
$$\dege(\psi(x)-\vartheta(x)-\sqrt{x}) = \dege(\psi(x)-\vartheta(x)-\psi(\sqrt{x}))$$
provided that 
$$(\tfrac{\Theta}{2},2,0,0,0\ldots) < (\tfrac{1}{3},0,0,0,\ldots),$$
that is, provided that $\Theta < \tfrac{2}{3}$.  The proposition follows.
\end{proof}

By Propositions \ref{lipitheta} and \ref{lipithetap},  one has the following.

\begin{corollary}\label{lipithetapcor}
One has
$$\li(x)-\pi(x) = \frac{x-\psi(x)}{\log x} + O\left(\frac{\sqrt{x}}{\log x}+\frac{x^\Theta}{(\log x)^2} \right) \ (x \to \infty).$$
Equivalently, one has
$$\li(x)-\pi(x) = \frac{x-\psi(x)}{\log x} + O\left(\frac{x^\Theta}{\log x} \right) \ (x \to \infty),$$
while also
$$\li(x)-\pi(x) = \frac{x-\psi(x)}{\log x} + O\left(\frac{x^\Theta}{(\log x)^2} \right) \ (x \to \infty)$$
if the Riemann hypothesis is false.
\end{corollary}

Next, recall from Section 5.1  that the  Riemann--von Mangoldt explicit formula for $\psi_0(x)$   is
$$\psi_0(x) = x - \sum_\rho \frac{x^\rho}{\rho} - \log 2\pi- \frac{1}{2}\log(1-x^{-2}), \quad \forall x>1,$$
where the sum is over all nontrivial zeros $\rho$ of the Riemann zeta function $\zeta(s)$, taken in taken in order of increasing absolute value of the imaginary part and repeated to multiplicity.  It follows that 
\begin{align}\label{xminuspsi}
x - \psi(x) = \sum_\rho \frac{x^\rho}{\rho} + O(1) \ (x \to \infty),
\end{align} which, along with Theorem \ref{lithetapsi} and (\ref{psieq}), implies the following.

\begin{proposition}\label{Rprop}
One has $$\dege(\id - \vartheta) = \dege(\id-\psi) = \dege\left( \sum_\rho \frac{x^\rho}{\rho} \right),$$ where the sum is over all nontrivial zeros $\rho$ of the Riemann zeta function.   Let
$$G(x) =e^{-\Theta x} \sum_\rho \frac{e^{\rho x}}{\rho}, \quad \forall x >0.$$ Then one has
$$G(x) = 2\sum_{\rho: \, \operatorname{Im}  \rho> 0} \frac{e^{-(\Theta-\operatorname{Re} \rho)x}}{|\rho|} \cos (x \operatorname{Im} \rho - \arctan(\operatorname{Im} \rho /\operatorname{Re} \rho ))$$
for all $x > 0$, and
$$\dege(\id - \psi) =   \dege\left( \sum_\rho \frac{x^\rho}{\rho} \right) = (\Theta, \deg G, \dege_1 G, \dege_2 G, \dege_3 G, \ldots).$$
Thus, one has
$$\dege_1(\li-\pi) = \dege_1 (\id -\psi) -1= \deg G-1$$
and
$$\dege_k(\li-\pi) = \dege_k (\id -\psi) = \dege_{k-1} G$$
for all integers $k > 1$.
\end{proposition}

\begin{corollary}
 If the Riemann hypothesis holds, then one has
$$\dege(\li-\pi) = (\tfrac{1}{2},  \deg G-1, \dege_1 G, \dege_2 G, \dege_3 G, \ldots).$$ Consequently, Conjecture \ref{eurekaconjecture} is equivalent to the conjunction of the Riemann hypothesis and the statement $\deg G = 0$, while Conjecture \ref{eurekaconjecture2} is equivalent to the conjunction of the Riemann hypothesis and the statements $\deg G = 0$ and $\dege_1 G = 0$.
\end{corollary}

In Montgomery and Vaughan's book \cite{mont}, it is asserted that the conjecture \cite[(15.24)]{mont} of Monach and Montgomery  implies that 
\begin{align}\label{MMCa}
\limsup_{x \to \infty} \frac{x-\psi(x)}{\sqrt{x}\, (\log \log \log x)^2} \geq \frac{1}{2\pi}
\end{align}
and 
\begin{align}\label{MMCb}
\liminf_{x \to \infty} \frac{x-\psi(x)}{\sqrt{x}\, (\log \log \log x)^2} \leq -\frac{1}{2\pi},
\end{align}
and, ``in view of (13.48), it is plausible that equality holds'' \cite[pp.\ 483--484]{mont}.  Indeed, in \cite[Conjecture, p.\ 16]{mont1}, some 28 years earlier, Montgomery had conjectured that equalities hold, that is, that
\begin{align}\label{MMC2}
\limsup_{x \to \infty} \frac{x-\psi(x)}{\sqrt{x}\, (\log \log \log x)^2} = \frac{1}{2\pi} \ \mbox{ and } \  \liminf_{x \to \infty} \frac{x-\psi(x)}{\sqrt{x}\, (\log \log \log x)^2} = -\frac{1}{2\pi}.
\end{align}
The following result concerning Montgomery's conjecture follows from Proposition \ref{oexpprop} (or Corollary \ref{deglogprop}) and Corollary \ref{lipithetapcor}.

\begin{proposition}   Montgomery's conjecture (\ref{MMC2}) \cite[Conjecture, p.\ 16]{mont1} implies
that
\begin{align}\label{MMC3}
0< \limsup_{x \to \infty} \frac{|x-\psi(x)|}{\sqrt{x}\, (\log \log \log x)^2} < \infty.
\end{align}
Moreover,  the conjecture (\ref{MMC3}) is equivalent to each of the  following statements.
\begin{enumerate}
\item One has
$$x-\psi(x) = O(\sqrt{x}\, (\log \log \log x)^2) \ (x \to \infty)$$ and  $$x-\psi(x) \neq o(\sqrt{x}\, (\log \log \log x)^2) \ (x \to \infty).$$
\item One has
$$\li(x)-\pi(x) = O\left(\frac{ \sqrt{x}\, (\log \log \log x)^2}{\log x} \right) \ (x \to \infty)$$ and  $$\li(x)-\pi(x) \neq o\left(\frac{ \sqrt{x}\, (\log \log \log x)^2}{\log x} \right) \ (x \to \infty).$$
\item Conjecture \ref{montconjec} holds, that is, one has
\begin{align*}
0< \limsup_{x \to \infty} \frac{|\li(x)-\pi(x)|}{\sqrt{x}\, (\log \log \log x)^2/\log x} < \infty.
\end{align*}
\end{enumerate}
Moreover, each of the equivalent statements above implies
that
$$\dege(\id-\psi) =  (\tfrac{1}{2},0,0,2,0,0,0,\ldots),$$
which is equivalent to
$$\dege(\li-\pi)  = (\tfrac{1}{2},-1,0,2,0,0,0,\ldots).$$
\end{proposition}

\begin{remark}[The second Chebyshev function and the Farey sequences]
It is straightforward to check that
$$e^{\psi(x)} = \operatorname{lcm}(1,2,3,\ldots,\lfloor x \rfloor)$$ 
for all $x \geq 0$.  It is known \cite{lus}  that
$$\operatorname{lcm}(1,2,3,\ldots,n) =  \prod_{r \in F_{n}, \, 0< r < 1} 2 \sin \pi r  = 2\prod_{r \in F_n, \, 0< r < \frac{1}{2}} 4 \sin^2 \pi r ,$$
where $F_n$ for any positive integer $n$ denotes the {\it $n$th  Farey sequence}.  It follows that
$$\psi(n) = \log  \operatorname{lcm}(1,2,3,\ldots,n) = \sum_{r \in F_n, \, 0< r< 1} \log | 2\sin \pi r |  = \log 2+ \sum_{r \in F_n, \, 0< r< \frac{1}{2}} \log ( 4\sin^2 \pi r ).$$
for all positive integers $n$. Therefore, one has
$$\dege(\id - \psi) = \dege (n- \log  \operatorname{lcm}(1,2,3,\ldots,n)) =  \dege \left( n-\sum_{r \in F_n, \, 0< r< 1} \log | 2\sin \pi r |\right),$$
and, in particular, one has
$$\Theta =  \deg (\id- \log  \operatorname{lcm}(1,2,3,\ldots,n)) = \deg \left( n-\sum_{r \in F_n, \, 0< r< 1} \log | 2\sin \pi r |\right).$$
\end{remark}

\section{The function $\li(x)-\pi(x)$ on average}

In this section, which is largely based on D.\ R.\ Johnston's work \cite{johnston}, we study the average value $$A_c(x) = \frac{1}{x-c}\int_c^x (\li(t)-\pi(t))\, dt$$ of $\li(x)-\pi(x)$ on the interval $[c,x]$ for various $c \geq 0$, the most natural choices being $c = 1$, $c = 2$, and $c = \mu$, where $\mu$ is the Ramanujan--Soldner constant.    More generally, we are concerned with the function $\int_c^x \frac{\li(t)-\pi(t)}{t^s}\, dt$ for any fixed $s,c \in \RR$ with $c \geq 0$.

By \cite[Theorems 1.1 and 1.2]{johnston} and its proof (adapted appropriately), one has the following.
  
  \begin{theorem}[{cf.\  \cite[Theorems 1.1 and 1.2]{johnston}}]\label{johnstthm}
 Let $c \geq 1$.  Each of the following statements is equivalent to the Riemann hypothesis.
  \begin{enumerate}
  \item $\int_2^x (\li(t)-\pi(t)) \, dt > 0$ for all $x > 2$.
  \item  $\int_\mu^x (\li(t)-\pi(t))\, dt > 0$ for all $x > \mu$.
  \item  $\int_1^x (\li(t)-\pi(t))\, dt > 0$ for all $x > 2.498183\ldots$.
    \item  $\int_1^x (\li(t)-\pi(t))\, dt > 0$ for all $x \gg 0$.
  \item  $\int_2^x (t -\vartheta(t)) \, dt > 0$ for all $x > 2$.
    \item  $\int_1^x (t -\vartheta(t)) \, dt > 0$ for all $x > 1$.
        \item  $\int_1^x (t -\vartheta(t)) \, dt > 0$ for all $x \gg 0$.
                                  \item $\int_c^x (\li(t)-\pi(t)) \, dt  \neq \Omega_{-}(x^\kappa) \ (x \to \infty)$ for all $\kappa > \frac{3}{2}$.
      \item $\int_c^x (t -\vartheta(t)) \, dt  \neq \Omega_{-}(x^\kappa) \ (x \to \infty)$ for all $\kappa > \frac{3}{2}$.
          \item  For all $\varepsilon > 0$, one has $\int_c^x (\li(t)-\pi(t)) \, dt > (0.58 -\varepsilon)\frac{x^{3/2}}{\log x} $ for all $x \gg 0$.
           \item  There is an $M > 0$ such that $\int_c^x (\li(t)-\pi(t))\, dt > M\frac{x^{3/2}}{\log x}$ for all $x \gg 0$.
                    \item  For all $\varepsilon > 0$, one has $ \int_c^x (t -\vartheta(t)) \, dt > (0.57 -\varepsilon)x^{3/2}$ for all $x \gg 0$.
                                       
                                        \item  There is an $M > 0$ such that $ \int_c^x (t -\vartheta(t)) \, dt > Mx^{3/2}$ for all $x \gg 0$.
  \end{enumerate}  
  \end{theorem}
  
Note that statement (3) is equivalent to $A_1(x) > 0$  for all $x > 2.498183\ldots$, and statement (11) is equivalent to the existence of an $M > 0$ such that $A_c(x) > M\frac{\sqrt{x}}{\log x}$  for all $x \gg 0$ (or, equivalently,  such that $A_c(x) > M(\li(x)-\Ri(x))$ for all $x \gg 0$).  Thus, if the Riemann hypothesis is true, then the function $\li(x)-\pi(x)$ has an overwhelming positive bias.  

Building largely upon Johnston's work,  in 2023, C.\ Axler proved the following.

\begin{proposition}[{\cite[Corollary 1.3]{axl}}]
The Riemann hypothesis holds if and only if 
$$(\tfrac{2}{3}-\lambda_0)\frac{x^{3/2}}{\log x}<\int_2^x (\li(t)-\pi(t))\, dt < \frac{x^{3/2}}{\log x},$$
where $\lambda_0 = 0.04611764442151$, and
where the inequality on the right holds for all $x \geq  1200963$ and the  inequality on the left holds for all $x \geq 10$.
\end{proposition}

Thus,  assuming the Riemann hypothesis, various subtleties in the asymptotic behavior the function $\li(x)-\pi(x)$,  e.g., those described by Littlewood's theorem,  do not arise in the asymptotic behavior of the running average of  $\li(x)-\pi(x)$.

Since
$$\int_2^x(\li(t)-\Ri(t))\, dt \sim \frac{2}{3} \frac{x^{3/2}}{\log x} \ (x \to \infty)$$ 
 by Karamata's integral theorem, and in fact
 $$\int_2^x(\li(t)-\Ri(t))\, dt > \frac{2}{3} \frac{x^{3/2}}{\log x}$$ 
 for all $x \geq 13$, one has the following.

\begin{corollary}
The Riemann hypothesis holds if and only if, for every $\varepsilon > 0$, one has
$$-(\varepsilon+\lambda_0)\frac{x^{3/2}}{\log x}<\int_2^x (\Ri(t)-\pi(t))\, dt < \frac{1}{3}\frac{x^{3/2}}{\log x}$$
for all $x \gg 0$, where the inequality on the right holds for all $x \geq  1200963$.
\end{corollary}

\begin{theorem}\label{strongdeglemma2}  One has the following.
\begin{enumerate}
\item Let $F(x) = \int_\mu^x (\li(t)-\pi(t)) \, dt$.   For all $s \in \RR$ and all $x \geq  \mu$, one has
$$\int_\mu^x \frac{\li(t)-\pi(t)}{t^s} \, dt = \frac{F(x)}{x^s} +  s \int_\mu^x \frac{F(t)}{t^{s+1}} \, dt.$$
\item One has
\begin{align*}
\int_\mu^\infty \frac{\li(t)-\pi(t)}{t^{s+1}} \, dt & = \frac{1}{s}\left(\Log \frac{1}{s-1}- P(s)\right) -\int_1^\mu \frac{\li(t)}{t^{s+1}} \, dt \\ 
& = -\frac{1}{s}\left(\li(\mu^{1-s})+ P(s)\right),
\end{align*}
where the first equality holds for all $s > \Theta$, and the second equality holds for all $s > 1$.
\item If the Riemann hypothesis is true, then for all $s \in \RR$ one has
$$\int_\mu^x \frac{\li(t)-\pi(t)}{t^s} \, dt \geq \frac{1}{x^{s}}\int_\mu^x (\li(t)-\pi(t))\, dt>0$$
for all  $x> \mu$  and 
$$\int_\mu^x \frac{\li(t)-\pi(t)}{t^s} \, dt \geq \frac{1}{x^{s}}\int_\mu^x (\li(t)-\pi(t))\, dt > (0.58 -\varepsilon)\frac{x^{3/2-s}}{\log x}$$
for all $x \gg 0$.  
\item If the Riemann hypothesis is true, then  one has
$$\int_\mu^x \frac{\li(t)-\pi(t)}{t^s} \, dt  \asymp \frac{x^{3/2-s}}{\log x} \ (x \to \infty)$$
for all $s < \frac{3}{2}$ and
$$\int_x^\infty \frac{\li(t)-\pi(t)}{t^s} \, dt  \ll \frac{x^{3/2-s}}{\log x} \ (x \to \infty)$$
for all $s > \frac{3}{2}$.
\item If the Riemann hypothesis is false, then, for all $s,u \in \RR$ with $u< s < \Theta+1$,  one has
$$\int_\mu^x \frac{\li(t)-\pi(t)}{t^s} \, dt = \Omega_{-} (x^{\Theta+1-u}) \ (x \to \infty).$$  
\item Unconditionally, one has
$$\int_\mu^x \frac{\li(t)-\pi(t)}{t^s} \, dt = O(x^{\Theta+1-s}\log x) \ (x \to \infty)$$
for all $s< \Theta+1$ and
$$\int_x^\infty \frac{\li(t)-\pi(t)}{t^s} \, dt = O(x^{\Theta+1-s}\log x) \ (x \to \infty)$$
for all $s > \Theta+1$,
while
$$\int_\mu^x \frac{\li(t)-\pi(t)}{t^{\Theta+1}} \, dt = O((\log x)^2) \ (x \to \infty).$$
\end{enumerate}
\end{theorem}

\begin{proof}
Statement (1) follows from Riemann--Stieltjes integration by parts \cite[Section 1.1.3]{borg}.
By (\ref{zetaP}), one has
\begin{align*}
\Log \frac{1}{s-1}- P(s)=  s \int_1^\infty  \frac{\li(x)-\pi(x)}{x^{s+1}} \, dt,
\end{align*}
on $\{s \in \CC: \operatorname{Re} s > \Theta\}$,  where the left-hand side is analytically continued to that domain.
Moreover, one has
$$\int_1^\mu \frac{\li(t)}{t^{s+1}} \, dt = \frac{1}{s}\left(\li(\mu^{1-s})-\log(s-1)\right)$$
for all $s >1$.  Statement (2) follows. 
Statement (3)  follows from  statement (1) and Theorem \ref{johnstthm}, and statement (4) follows readily from statements (1)--(3).   Finally,  statement (5) follows from \cite[Lemma 2.8]{johnston} and the proof of \cite[Theorems 1.1 and 1.4]{johnston}, and statement (6) follows  from the fact that $\li(x)-\pi(x) = O(x^\Theta\log x) \ (x \to \infty)$.
\end{proof}

By Proposition \ref{pitheta},  one has
\begin{align}\label{pithetaeq}
\li(x)-\pi(x) -  \frac{x-\vartheta(x)}{\log x}   = \frac{1}{\log x}\int_0^x \frac{\li(t)-\pi(t)}{t} \, dt, \quad \forall x >0.
\end{align}
Consequently, one has the following.

\begin{corollary}
The following statements are equivalent to the Riemann hypothesis.
\begin{enumerate}
\item $\int_\mu^x \frac{\li(t)-\pi(t)}{t} \, dt \asymp  \frac{x^{1/2}}{\log x} \ (x \to \infty)$.
\item $\dege \int_\mu^x \frac{\li(t)-\pi(t)}{t} \, dt = (\frac{1}{2},-1,0,0,0,\ldots)$.
\item $\deg \int_\mu^x \frac{\li(t)-\pi(t)}{t} \, dt = \frac{1}{2}$.
\item $\li(x)-\pi(x) -  \frac{x-\vartheta(x)}{\log x}  \asymp  \frac{x^{1/2}}{(\log x)^2} \ (x \to \infty)$.
\item $\dege\left(\li-\pi -  \frac{\id-\vartheta}{\log} \right) = (\frac{1}{2},-2,0,0,0,\ldots)$.
\item $\deg\left(\li-\pi -  \frac{\id-\vartheta}{\log} \right) = \frac{1}{2}$.
\end{enumerate}
\end{corollary}

\begin{corollary}\label{strongdeglemma}
One has
 $$\deg \int_\mu^x (\li(t)-\pi(t))  \, dt = \Theta+1.$$
\end{corollary}

\begin{proof}
If the Riemann hypothesis is true, then by Theorem \ref{johnstthm} one has
$$\deg \int_\mu^x (\li(t)-\pi(t))  \, dt \geq \tfrac{3}{2} = \Theta+1.$$
On the other hand, if the Riemann hypothesis is false, then by  Theorem \ref{strongdeglemma2}(5) one has
$$\deg \int_\mu^x (\li(t)-\pi(t))  \, dt \geq \Theta+1.$$
Moreover,  by Theorem \ref{strongdeglemma2}(6), the reverse inequalities are unconditional.  
\end{proof}

\begin{theorem}\label{mainjohnston}
 For all $s \in \RR$, one has
  $$\deg  \int_\mu^x \frac{\li(t)-\pi(t)}{t^s} \, dt  = \max(\Theta+1-s,0).$$
Moreover,  for all $s < \Theta+1$, one has 
  $$\dege  \int_\mu^x \frac{\li(t)-\pi(t)}{t^s} \, dt  = \dege \int_\mu^x (\li(t)-\pi(t))  \, dt+(-s,0,0,0,\ldots),$$
  and, for all $s > \Theta+1$, one has
  $$\dege  \int_x^\infty \frac{\li(t)-\pi(t)}{t^s} \, dt  =  \dege \int_\mu^x (\li(t)-\pi(t))  \, dt+(-s,0,0,0,\ldots).$$
Furthermore, the function $$ \int_\mu^\infty \frac{\li(t)-\pi(t)}{t^s} \, dt \sim \frac{1}{s^2\log \mu} \mu^{1-s} \ (s \to \infty)$$  of $s$  is analytic, positive, and decreasing to $0$ on $(\Theta+1,\infty)$.
\end{theorem}

\begin{proof}
From the fact that $\li(\mu^{1-s}) \sim  -\frac{\mu}{(s-1)\log \mu} \mu^{-s}$ and $P(s) \sim 2^{-s}$,  one can verify, with assistance from Mathematica,  that the function $-\frac{1}{s}\left(\li(\mu^{1-s})+ P(s)\right)$ is positive and decreasing on $(1,\infty)$.  A separate computation shows that $\frac{1}{s}\left(\Log \frac{1}{s-1}- P(s)\right) -\int_1^\mu \frac{\li(t)}{t^{s+1}} \, dt$ is (defined and) positive and decreasing on $(\frac{1}{2},1]$.
In particular,  by Theorem \ref{strongdeglemma2}(2), one has $$\lim_{x \to \infty}\int_\mu^x \frac{\li(t)-\pi(t)}{t^s} \, dt >0,$$ whence $\deg \int_\mu^x \frac{\li(t)-\pi(t)}{t^s} \, dt = 0$, for all $s >\Theta+1$.   The theorem follows, then,  from Theorem \ref{strongdeglemma2},  Corollary \ref{strongdeglemma},  and Proposition \ref{strongdegg}.
\end{proof}

\begin{corollary}\label{strongdegree}
The Riemann hypothesis is equivalent to each of the following statements.
\begin{enumerate}
  \item $\deg \int_\mu^x \frac{\li(t)-\pi(t)}{t^s} \, dt = \max\left( \tfrac{3}{2}-s, 0\right)$ for all $ s \in \RR$.
  \item $\deg \int_\mu^x \frac{\li(t)-\pi(t)}{t^s} \, dt =  \tfrac{3}{2}-s$ for some $s \leq \frac{3}{2}$.
  \item $\deg \int_x^\infty \frac{\li(t)-\pi(t)}{t^s} \, dt =  \tfrac{3}{2}-s$ for all $s > \frac{3}{2}$.
    \item $\deg \int_x^\infty \frac{\li(t)-\pi(t)}{t^s} \, dt =  \tfrac{3}{2}-s$ for some $s >2$ (or for some $s > \Theta+1$).
            \item $\dege \int_\mu^x \frac{\li(t)-\pi(t)}{t^s} \, dt = (\frac{3}{2}-s,-1,0,0,0,\ldots)$  for all $s < \frac{3}{2}$.
    \item $\dege \int_\mu^x \frac{\li(t)-\pi(t)}{t^s} \, dt = (\frac{3}{2}-s,-1,0,0,0,\ldots)$  for some $s < \frac{3}{2}$.
            \item $\dege \int_x^\infty \frac{\li(t)-\pi(t)}{t^s} \, dt = (\frac{3}{2}-s,-1,0,0,0,\ldots)$  for all $s > \frac{3}{2}$.
        \item $\dege \int_x^\infty \frac{\li(t)-\pi(t)}{t^s} \, dt = (\frac{3}{2}-s,-1,0,0,0,\ldots)$ for some $s >2$ (or for some $s > \Theta+1$).
  \end{enumerate}
\end{corollary}

\begin{corollary}
One has
\begin{align*}
\dege\left(\li-\pi - \frac{\id-\vartheta}{\log} \right) &  = \dege  \int_1^x \frac{\li(t)-\pi(t)}{t} \, dt +(0,-1,0,0,0,\ldots) \\
& = \dege  \int_1^x (\li(t)-\pi(t)) \, dt + (-1,-1,0,0,\ldots),
\end{align*}
and therefore
$$\deg\left(\li-\pi - \frac{\id-\vartheta}{\log} \right) = \Theta.$$
\end{corollary}

\begin{remark}[The function $\int_c^x \frac{\li(t)-\pi(t)}{t^s} \, dt$]
 For all $s \in \RR \backslash\{\alpha\}$, one has
  $$\deg  \int_1^x \frac{\li(t)-\pi(t)}{t^s} \, dt  = \max(\Theta+1-s,0),$$
where $\alpha =  2.2983410255\ldots$ is the unique $s  \in \RR$ for which $s > \Theta+1$ and
$\int_1^\infty \frac{\li(t)-\pi(t)}{t^s} \, dt = 0$, or, equivalently, for which $s > 1$ and $P(s) = \log\frac{1}{s-1}$.
  Consequently, one has
    $$\deg  \int_1^x \frac{\li(t)-\pi(t)}{t^{\alpha}} \, dt = \deg  \int_x^\infty \frac{\li(t)-\pi(x)}{t^{\alpha}} \, dt  = \Theta+1-\alpha < 0.$$  Note also that the function $$\int_1^\infty \frac{\li(t)-\pi(t)}{t^s} \, dt = \Log \frac{1}{s-2}- P(s-1) \sim -\log s \ (s \to \infty)$$ is positive for $s \in (\Theta+1,\alpha)$, but negative for $s \in (\alpha,\infty)$, yet is still decreasing on $(\Theta+1,\infty)$ (and it assumes the limiting value $H$ at $2 \in (\Theta+1,\alpha)$).    The function
    $$\int_2^\infty \frac{\li(t)-\pi(t)}{t^{s}} \, dt \sim  \frac{\li(2)-1}{s}2^{1-s} \ (s \to \infty),$$
on the other hand,    is positive and decreasing on $(\Theta+1,\infty)$.  
\end{remark}

Given the results above, the following is clear.

\begin{proposition}\label{integralpiri}
One has the following.
\begin{enumerate}
\item Unconditionally, one has
$$(\Theta+1,-\infty,-1,0,0,0,\ldots) \leq \dege  \int_1^x (\li(t)-\pi(t)) \, dt \leq (\Theta+1,1,0,0,0,\ldots).$$
\item The Riemann hypothesis is equivalent to
$$\int_1^x (\li(t)-\pi(t)) \, dt \asymp \frac{x^{3/2}}{\log x} \ (x \to \infty)$$
and to
$$\dege  \int_1^x (\li(t)-\pi(t)) \, dt = (\tfrac{3}{2},-1,0,0,0,\ldots).$$
\item One has
$$\int_1^x(\li(t)-\Ri(t))\, dt \sim \frac{2}{3} \frac{x^{3/2}}{\log x} \ (x \to \infty),$$  and therefore
$$\dege  \int_1^x (\li(t)-\Ri(t)) \, dt  = (\tfrac{3}{2},-1,0,0,0,\ldots).$$
\item One has
 $$\dege  \int_1^x (\Ri(t)-\pi(t)) \, dt \leq \dege  \int_1^x (\li(t)-\pi(t)) \, dt.$$
\item The Riemann hypothesis is true if and only if
$$\int_1^x (\Ri(t)-\pi(t)) \, dt =  O\left( \frac{x^{3/2}}{\log x} \right) \ (x \to \infty),$$
if and only if
$$\dege \int_1^x (\Ri(t)-\pi(t)) \, dt   \leq (\tfrac{3}{2},-1,0,0,0,\ldots).$$
\item The Riemann hypothesis is false if and only if
$$\dege  \int_1^x (\li(t)-\pi(t)) \, dt >(\tfrac{3}{2},-1,0,0,0,\ldots),$$ 
if and only if
$$\dege  \int_1^x (\Ri(t)-\pi(t)) \, dt >(\tfrac{3}{2},-1,0,0,0,\ldots).$$
\item If the Riemann hypothesis is false, then one has
$$\dege  \int_1^x (\li(t)-\pi(t)) \, dt =  \dege  \int_1^x (\Ri(t)-\pi(t)) \, dt.$$
\end{enumerate}
\end{proposition}

\begin{outstandingproblem}\label{integralripi}
Is the inequality 
 $$\dege  \int_1^x (\Ri(t)-\pi(t)) \, dt \leq \dege  \int_1^x (\li(t)-\pi(t)) \, dt$$
strict?
\end{outstandingproblem}

Note that the answer is in the negative if the Riemann hypothesis is false.  Equivalently, if
 $$\dege  \int_1^x (\Ri(t)-\pi(t)) \, dt < \dege  \int_1^x (\li(t)-\pi(t)) \, dt,$$
 then the Riemann hypothesis is true.
 

\chapter{Summatory functions}

In this chapter, we study various arithmetic functions and their  corresponding summatory functions, including the Mertens function, the summatory Liouville function, the summatory divisor function, and the error terms in Mertens' theorems.

\section{Error terms in Mertens' theorems}

   Recall that  Mertens' second theorem states that the limit
$$M = \lim_{x \to \infty} \left( \sum_{p \leq x} \frac{1}{p} - \log \log x \right) = 0.261497212847\ldots,$$  
known as the Meissel--Mertens constant, exists.  Let
$$\Mert(x) =\sum_{p \leq x} \frac{1}{p}  -  \log \log x-M, \quad \forall x > 1.$$
The bound 
$$\int_x^\infty \frac{ \li(t)-\pi(t)}{t^2}\, dt = o\left((\log  x) e^{- A(\log x)^{3/5}(\log \log x)^{-1/5}}\right) \ (x \to \infty),$$
where $A  = 0.2098$, follows from Theorem \ref{bestPNT} and  Corollary \ref{karamcor}(2).   From this bound, along with Proposition \ref{psthm}, we deduce the following.

\begin{proposition}\label{psthmcor}
One has
$$M_2(x) =   O\left(e^{- C(\log x)^{3/5}(\log \log x)^{-1/5}}\right) \ (x \to \infty)$$ for all positive $C< 0.2098$.
\end{proposition}

Regarding lower bounds for $M_2(x)$, one has the following.

\begin{proposition}[{\cite{lay}}]\label{psthmcor22}
One has $$\Mert(x)  = \Omega_{\pm} \left(\frac{1}{\sqrt{x}\log x} \right)\ (x \to \infty),$$
and also, assuming the Riemann hypothesis,
$$\Mert(x) = \Omega_{\pm} \left(\frac{\log \log \log x}{\sqrt{x}\log x} \right)\ (x \to \infty).$$
\end{proposition}

 From the results above, we deduce the following.

\begin{proposition}\label{dege1p} One has the following.
\begin{enumerate}
\item $\displaystyle \dege \Mert \leq (0,-\infty,-\tfrac{3}{5},\tfrac{1}{5},0,0,0,\ldots).$
\item $\displaystyle \dege \Mert \geq (-\tfrac{1}{2},-1,0,0,0,\ldots).$
\item $\displaystyle\dege \Mert \geq (-\tfrac{1}{2},-1,0,1,0,0,0,\ldots)$, provided that the Riemann hypothesis holds.
\end{enumerate}
\end{proposition}

Recall  that Proposition \ref{psthm} relates the functions $\Mert$ and $\li-\pi$.  Using  that result, we prove the following.

\begin{theorem}\label{mertenssecond}
Let $r$ be continuous on $[N,\infty)$ for some $N > 1$ and regularly varying of index $d$ (e.g., $r \in \mathbb{L}$ with $\deg r = d$) for some $d \in \RR$.
\begin{enumerate}
\item If $d<1$, then one has $$\li(x)-\pi(x) = O(r(x)) \ (x \to \infty)$$
if and only if
$$ \Mert(x)=  O(x^{-1}  r(x)) \ (x \to \infty).$$ 
\item  If $r$ is Hardian (or $r \in \mathbb{L}$), $d = 1$,  and $\dege_1 r = -\infty$,   then each of the following conditions implies the next.
\begin{enumerate}
\item $\li(x)-\pi(x) = O(r(x)) \ (x \to \infty)$.
\item $ \Mert(x)=  o(x^{-1} r(x) \log x) \ (x \to \infty)$.
\item $ \Mert(x)=  O(x^{-1} r(x) \log x) \ (x \to \infty)$.
\item $\li(x)-\pi(x) =  o(r(x) (\log x)^2) \ (x \to \infty)$.
\end{enumerate}
\item One has $$\dege \Mert= \dege(\li-\pi) +(-1,0,0,0,\ldots).$$
\end{enumerate}
\end{theorem}

\begin{proof}
 Suppose first that $d <1$.   Suppose that  $$\li(x)-\pi(x) = O(r(x)) \ (x \to \infty),$$
 so that $d-2 \in [-\frac{3}{2},-1)$.   By Proposition \ref{psthm} and Karamata's integral theorem, for all $x \geq N$ one has
\begin{align*}
\Mert(x)& =  - \frac{\li(x)-\pi(x)}{x} + \int_x^\infty \frac{ \li(t)-\pi(t)}{t^2}\, dt  \\
& = O(x^{-1} r(x))  + O\left(\int_x^\infty  \frac{|r(t)|}{t^2} \, dt\right) \ (x \to \infty)  \\
& =  O(x^{-1} r(x)) \ (x \to \infty).
\end{align*}
A similar proof establishes the converse.  This proves statement (1).   The proof of (2) is similar but uses Corollary \ref{karamcor}(2) in place of  Karamata's integral theorem. Statement (3) then follows from statements (1) and (2) and Theorem \ref{infpropexp}.  Indeed,   if $\Theta < 1$, then statement (3) follows immediately from statement (1) and Theorem \ref{infpropexp}.   If,  on the other hand,  one has
$d = \Theta = 1$, then we may  restrict our infimum as in Theorem \ref{infpropexp} to functions $r \in \mathbb{L}$ such that $r(x) = O(xe^{-\sqrt{\log x}}) \ (x \to \infty)$,  so that $\dege_1 r  = -\infty$ and therefore $ \dege r(x) = \dege ( r(x) \log x) = \dege ( r(x)( \log x)^2)$ and $\dege (x^{-1} r(x)) = \dege ( x^{-1}r(x) \log x)$.  
\end{proof}

By the theorem above and the results of Section 9.1, we have the following.

\begin{corollary}
One has the following.
\begin{enumerate}
\item The Riemann hypothesis is equivalent to $$\deg \Mert  = -\tfrac{1}{2},$$
to $$ \Mert(x) = O(x^{-1/2}\log x) \ (x \to \infty),$$ and to
$$(-\tfrac{1}{2},-1,0,1,0,\ldots) \leq \dege \Mert \leq (-\tfrac{1}{2},1,0,0,0,\ldots).$$
\item The anti-Riemann hypothesis $\Theta = 1$ is equivalent to $\deg \Mert = 0$ and to
$$(0,-\infty,-1,0,0,0,\ldots) \leq \dege \Mert \leq (0,-\infty,-\tfrac{3}{5},\tfrac{1}{5},0,0,0,\ldots),$$
and it implies that  $\dege_1 \Mert = - \infty$ and also that
$$-\dege_2 \Mert = -\Theta_2 \in [\tfrac{3}{5},1]$$
is equal to the anti-Riemann constant.
\item The anti-Riemann constant is equal to $\underline{\deg}\,  \log|\Mert(e^x)|= -\deg \frac{1}{\log|\Mert(e^x)|}$.
\end{enumerate}
\end{corollary}

Now, let
$$N(x) = \sum_{p \leq x} \log\left(1- \frac{1}{p}\right) +\log\log x +\gamma, \quad \forall x > 0,$$
and
$$\Mertt(x) = e^{-\gamma} \prod_{p \leq x}\left(1-\frac{1}{p}\right)^{-1} - \log x, \quad \forall x > 0.$$
By \cite[Lemma 5]{lay} (or by \cite[p.\ 81]{prachar}),  one has
\begin{align*} \Mert(x) = -N(x)+ O\left(\frac{1}{x} \right) \ (x \to \infty).
\end{align*}
(In fact,  in Proposition \ref{merttpropp} below  we show that $\Mert(x) = -N(x)+ O\left(\frac{1}{x \log x} \right) \ (x \to \infty)$.)
It follows  that $\lim_{x \to \infty} N(x) = 0$ and $$\dege N = \dege M_2.$$
Consequently, one has
$$ (\log x)e^\gamma\prod_{p \leq x} \left(1- \frac{1}{p}\right)- 1= e^{N(x)} -1 \sim N(x) \ (x \to \infty),$$
and therefore
$$-\frac{M_3(x)}{(\log x)^2}  \sim e^\gamma \prod_{p \leq x} \left(1- \frac{1}{p}\right) -\frac{1}{\log x} \sim \frac{N(x)}{\log x}  \ (x \to \infty)$$
and
$$\Mertt(x) \sim -(\log x )N(x) \ (x \to \infty).$$
By the remarks above, one has the following.

\begin{corollary}
One has the following.
\begin{enumerate}
\item $\dege N = \dege \Mert$.
\item $\Mertt(x) \sim - (\log x)N(x) \ (x \to \infty)$.
\item $\dege \Mertt =  \dege M_2 + (0,1,0,0,0,\ldots) = \dege(\li-\pi) +(-1,1,0,0,\ldots)$.
\item $(-\tfrac{1}{2},0,0,0,\ldots) \leq \dege \Mertt \leq (0,-\infty,-\tfrac{3}{5},\tfrac{1}{5},0,0,0,\ldots).$
\item The Riemann hypothesis is equivalent to $\deg \Mertt =  -\frac{1}{2}$, to
 $$ \Mertt(x) = O(x^{-1/2}(\log x)^2) \ (x \to \infty),$$ and to
$$(-\tfrac{1}{2},0,0,1,0,0,0\ldots) \leq \dege \Mertt \leq (-\tfrac{1}{2},2,0,0,0,\ldots).$$
\item The anti-Riemann hypothesis $\Theta = 1$ is equivalent to $\deg \Mertt = 0$ and to
$$(0,-\infty,0,0,0,0,\ldots) \leq \dege \Mertt \leq (0,-\infty,-\tfrac{3}{5},\tfrac{1}{5},0,0,0,\ldots),$$
and it implies that  $\dege_1 \Mertt = - \infty$ and also that
$$-\dege_2 \Mertt =- \Theta_2 \in [\tfrac{3}{5},1]$$
is equal to the anti-Riemann constant.
\item The anti-Riemann constant is equal to $\underline{\deg}\, \log|\Mertt(e^x)|= -\deg \frac{1}{\log|\Mertt(e^x)|}$.
\end{enumerate}
\end{corollary}

Thus, we have expressed all of $\dege \Mert$, $\dege \Mertt$,  $\dege N$, and $\dege(\li-\pi)$ in terms of any one of them.    Note also  that, since $H_n-\gamma - \log n \sim \frac{1}{2n} \ (n \to \infty),$
where $H_n$ denotes the $n$th harmonic number, one also has
$$\dege \Mertt = \dege \left( e^{-\gamma} \prod_{p \leq n}\left(1-\frac{1}{p}\right)^{-1} - H_n +\gamma\right).$$

Another consequence of the results above is the following.

\begin{corollary}
Let $\alpha \in (0,1]$ and $\beta \in \RR$.    The equivalent conditions of Theorem \ref{antiR2} are  each equivalent to the existence of a constant $C > 0$ satisfying any of the following equivalent conditions.
\begin{enumerate}
\item $N(x)= O \left(e^{-C(\log x)^\alpha (\log \log x)^\beta} \right) \ (x \to \infty)$.
\item $\Mert(x)= O \left(e^{-C(\log x)^\alpha (\log \log x)^\beta} \right) \ (x \to \infty)$.
\item $\Mertt(x)= O \left((\log x )e^{-C(\log x)^\alpha (\log \log x)^\beta} \right) \ (x \to \infty)$.
\end{enumerate}
\end{corollary}

If Lamzouri's conjecture \cite[Conjecture 1.5]{lam} holds, then  one has
\begin{align}\label{manzc}
0< \limsup_{x \to \infty} \frac{|\Mertt(x)|}{(\log \log \log x)^2/\sqrt{x}} < \infty.
\end{align}
Moreover,  by Theorem \ref{mertenssecond}, the conjectural bound above is equivalent to Montgomery's conjecture (\ref{MMC3}).  Thus, we have the following.

\begin{corollary}  One has the following.
\begin{enumerate}
\item Conjecture \ref{eurekaconjecture} is equivalent to  
$$\deg \Mert = -\tfrac{1}{2} \text{ and }\dege_1\Mert = -1$$
and to
$$\deg \Mertt = -\tfrac{1}{2} \text{ and }\dege_1\Mertt = 0.$$
\item Conjecture \ref{eurekaconjecture2} is equivalent to  
$$\deg \Mert = -\tfrac{1}{2}, \ \dege_1 \Mert = -1, \text{ and } \dege_2 \Mert = 0$$
and to
$$\deg \Mertt = -\tfrac{1}{2}, \ \dege_1 \Mertt = 0, \text{ and } \dege_2 \Mertt = 0.$$
\item The conjectural statements
$$\dege (\li-\pi) = (\tfrac{1}{2}, -1, 0, 2, 0, 0, \ldots),$$
$$\dege \Mert = (-\tfrac{1}{2}, -1, 0, 2, 0, 0, \ldots),$$
and
$$\dege \Mertt = (-\tfrac{1}{2}, 0, 0, 2, 0, 0, \ldots),$$  
are all equivalent and are implied by Montgomery's conjecture (\ref{MMC3}) and by Lamzouri's conjecture (\ref{manzc}), both of which are equivalent to
\begin{align*}
0< \limsup_{x \to \infty} \frac{|\li(x)-\pi(x)|}{\sqrt{x}\, (\log \log \log x)^2/\log x} < \infty
\end{align*}
and to 
\begin{align*}
0< \limsup_{x \to \infty} \frac{|\Mert(x)|}{(\log \log \log x)^2/(\sqrt{x}\log x)} < \infty.
\end{align*}
\end{enumerate}
\end{corollary}

Let $s \in \RR$ be nonzero and not equal to a prime.  Recall the definitions of $G(s)$ and $H(s)$ from Section 8.5.  Let
$$N(x,s) = \frac{1}{s} \sum_{p \leq x} \log\left(1- \frac{s}{p}\right) +\log\log x +G(s), \quad \forall x > 0,$$
and
$$\Mertt(x,s)= -(\log x) \left( 1-e^{-N(x,s)}\right) = e^{-G(s)} \prod_{p \leq x}\left(1-\frac{s}{p}\right)^{-1/s} - \log x, \quad \forall x > 0.$$
Note that 
$$N(x) = N(x,1)$$
and $$\Mertt(x) = \Mertt(x,1)$$
for all $x > 0$.

\begin{proposition}\label{merttpropp}
Let $s \in [-2,2)\backslash \{0\}$.  
\begin{enumerate}
\item One has
$$\displaystyle M_2(x) =-\lim_{s \to 0} N(x,s) = -N(x,s) + O\left(\frac{1}{x \log x}\right)\ (x \to \infty)$$
and therefore 
$$\dege N(-,s) = \dege M_2.$$
\item One has
\begin{align*}
M_2(x)+N(x,s) & = sH(s) + \sum_{p \leq x} \left(\frac{1}{p} +\frac{1}{s}\log \left(1-\frac{s}{p}\right) \right) \\ &
 = \sum_{p>x} \sum_{k = 2}^\infty   \frac{s^{k-1}}{k p^k} \\ &  \sim \frac{s}{2x\log x} \  (x \to \infty)
\end{align*}
and therefore
\begin{align*}
\dege (M_2(x)+N(x,s) )  = (-1,-1,0,0,0,\ldots).
\end{align*}
\item One has
$$\Mertt(x,s)  \sim -(\log x )N(x,s) \ (x \to \infty)$$
and therefore
$$\dege \Mertt(-,s)   = \dege M_2 + (0,1,0,0,0,\ldots).$$
\item One has
\begin{align*}
 \lim_{s \to 0} M_3(x,s) & = e^{-M}\prod_{p\leq x} e^{1/p}-\log x \\ & =(\log x)\left( e^{-M_2(x)}-1 \right) \\ & \sim -(\log x)  M_2(x) \ (x \to \infty)
 \end{align*}
for all $x > 0$, and therefore
$$\dege\left(e^{-M}\prod_{p\leq x} e^{1/p}-\log x \right)   = \dege M_2 + (0,1,0,0,0,\ldots).$$
\end{enumerate}
\end{proposition}

\begin{proof}
One has
\begin{align*}
-\frac{1}{s}\sum_{p \leq x}\log \left(1-\frac{s}{p}\right) & = \frac{1}{s}\sum_{p \leq x} \sum_{k = 1}^\infty \frac{s^k}{kp^k} \\
&  = \sum_{p \leq x} \frac{1}{p} + \sum_p \sum_{k = 2}^\infty   \frac{s^{k-1}}{k p^k}  -\sum_{p>x} \sum_{k = 2}^\infty   \frac{s^{k-1}}{k p^k} \\
& = M_2(x) +\log \log x  + M+sH(s) - \sum_{p>x} \sum_{k = 2}^\infty   \frac{s^{k-1}}{k p^k} \\
& = M_2(x) +\log \log x  + G(s) -\sum_{p>x} \sum_{k = 2}^\infty   \frac{s^{k-1}}{k p^k},
\end{align*}
where 
$$\sum_{p>x} \sum_{k = 2}^\infty   \frac{s^{k-1}}{k p^k} \sim \frac{s}{2} \sum_{p>x}   \frac{1}{ p^2} \sim \frac{s}{2x \log x} \  (x \to \infty).$$
 Statements (1) and (2) follow.  
By (1),
one has
$$ (\log x)e^{G(s)} \prod_{p \leq x} \left(1- \frac{s}{p}\right)^{1/s}- 1= e^{N(x,s)} -1 \sim N(x,s) \ (x \to \infty),$$
and therefore
$$e^{G(s)} \prod_{p \leq x} \left(1- \frac{s}{p}\right)^{1/s} -\frac{1}{\log x} \sim \frac{N(x,s)}{\log x}  \ (x \to \infty)$$
and
$$\Mertt(x,s) = e^{-G(s)} \prod_{p \leq x} \left(1- \frac{s}{p}\right)^{-1/s}- \log x \sim -(\log x )N(x,s) \ (x \to \infty).$$
Statements (3) and (4) follow.
\end{proof}

Recall that Mertens' first theorem, in its stronger limit form, states that the limit
$$B = \lim_{x \to \infty} \left( \log x - \sum_{p \leq x} \frac{\log p}{p} \right) = 1.332582275733\ldots$$
exists.    All of the results above have analogues for the function
$$M_1(x) = \log x - \sum_{p \leq x} \frac{\log p}{p}-B.$$
The proofs are similar.  (See \cite[Lemmas 1 and 4]{lay}, for example.)  In particular, the following proposition is the analogue of Theorem \ref{mertenssecond}(3) for $M_1$.

\begin{proposition}\label{mertprop1}
 One has
$$\dege M_1 = \dege(\li-\pi)+(-1,1,0,0,0,\ldots).$$
\end{proposition}

Since one has
\begin{align*}
\dege \left(e^{-\gamma} \prod_{p \leq x}\left(1-\frac{1}{p}\right)^{-1}-\log x \right) & =
\dege \left(e^{-M} \prod_{p\leq x} e^{1/p}-\log x \right) \\ & = 
\dege  \left(\sum_{p \leq x} \frac{\log p}{p}+B-\log x \right) \\ & = 
\dege (\li-\pi) + (-1,1,0,0,0,\ldots),
\end{align*}
one might inquire as to the logexponential degrees
$$\dege \left(e^{-\gamma} \prod_{p \leq x}\left(1-\frac{1}{p}\right)^{-1}- e^{-M} \prod_{p\leq x} e^{1/p}\right) =  \dege\left(\prod_{p \leq x}\left(1-\frac{1}{p}\right)^{-1}-e^{H}\prod_{p\leq x} e^{1/p} \right),$$
$$\dege \left(e^{-\gamma} \prod_{p \leq x}\left(1-\frac{1}{p}\right)^{-1}- \sum_{p \leq x} \frac{\log p}{p}-B\right),$$
and
$$\dege \left(e^{-M} \prod_{p\leq x} e^{1/p}- \sum_{p \leq x} \frac{\log p}{p}-B\right),$$
all three of which are less than or equal to $\dege (\li-\pi) + (-1,1,0,0,0,\ldots)$,  and at least two of which are equal.  Note, for instance,  the following.

\begin{proposition}\label{mertprop2}
One has
\begin{align*}
e^{-\gamma}\prod_{p \leq x}\left(1-\frac{1}{p}\right)^{-1}-e^{-M}\prod_{p\leq x} e^{1/p} \sim  \frac{1}{2x} \ (x \to \infty)
\end{align*}
and therefore
$$\dege \left(e^{-\gamma} \prod_{p \leq x}\left(1-\frac{1}{p}\right)^{-1}- e^{-M} \prod_{p\leq x} e^{1/p}\right) =  (-1,0,0,0,\ldots).$$
\end{proposition}

\begin{proof}
One has
\begin{align*}
e^{-\gamma}\prod_{p \leq x}\left(1-\frac{1}{p}\right)^{-1}-e^{-M}\prod_{p\leq x} e^{1/p} & =e^{-\gamma} \prod_{p\leq x} e^{1/p} \left(\prod_{p \leq x}\left(1-\frac{1}{p}\right)^{-1}e^{-1/p}-e^H \right)  \\ & =
e^{-M}\prod_{p\leq x} e^{1/p} \left(\prod_{p > x}\left(1-\frac{1}{p}\right)^{-1}e^{-1/p}-1 \right)  \\
 & \sim - \log x \left(\sum_{p > x} \left(\frac{1}{p} + \log \left(1-\frac{1}{p}\right) \right) \right) \\
 & =  \log x \left(H + \sum_{p \leq x} \left(\frac{1}{p} + \log \left(1-\frac{1}{p}\right) \right)\right) \\
 & \sim  \frac{\log x}{2x\log x} \\
 & =  \frac{1}{2x}
\end{align*}
as $x \to \infty$.
\end{proof}

\begin{problem}
Compute $$\dege \left(e^{-\gamma} \prod_{p \leq x}\left(1-\frac{1}{p}\right)^{-1}- \sum_{p \leq x} \frac{\log p}{p}-B\right)$$
and
$$\dege \left(e^{-M} \prod_{p\leq x} e^{1/p}- \sum_{p \leq x} \frac{\log p}{p}-B\right).$$
\end{problem}

Finally, we note that the method we used to prove Theorems \ref{lithetapsi} and \ref{mertenssecond} can be used to generalize Proposition \ref{pis}, and its proof, as follows.

\begin{theorem}\label{pis2}
Let $s \in \CC$ with $\operatorname{Re} s > -\Theta$.   One has
$$\dege\left(\sum_{p\leq x} p^s - \Ei((s+1)\log x)\right) = \dege(\li-\pi)+(\operatorname{Re} s , 0,0,0,\ldots).$$
In particular, the Riemann hypothesis is equivalent to
$$\deg\left(\sum_{p\leq x} p^s - \Ei((s+1)\log x)\right) = \frac{1}{2}+\operatorname{Re} s.$$
\end{theorem}

\section{The Mertens function $M(x)$ and summatory Liouville function $L(x)$}

Recall that the Mertens function is the summatory function
$$M(x) = S_{\mu}(x) =   \sum_{n \leq x} \mu(n)$$
of the M\"obius function $\mu(n)$.   Closely related to the Mertens function is the
 {\bf summatory Liouville function}
$$L(x) = S_{\lambda}(x)=  \sum_{n \leq x} \lambda(n).$$
By Example \ref{direxample}(3), one has
\begin{align}\label{Dmu}
D_\mu(s)  = \frac{1}{\zeta(s)} =  s \int_1^\infty  \frac{M(x)}{x^{s+1}} \, dx
\end{align}
and
\begin{align}\label{Dlamb}D_\lambda(s) = \frac{\zeta(2s)}{\zeta(s)} =  s \int_1^\infty  \frac{L(x)}{x^{s+1}} \, dx
\end{align}
for all $s \in \CC$ with $\operatorname{Re} s > 1$.

\begin{theorem}\label{MLProp}
One has $$\deg M = \deg L = \Theta.$$ Equivalently, both (\ref{Dmu}) and (\ref{Dlamb}) hold for all $s \in \CC$ with $\operatorname{Re} s > \Theta$.
\end{theorem}

\begin{proof}
Let $\Theta' = \deg M$, which, by Theorem \ref{dirichlet}, is the abscissa  of convergence of the Dirichlet series $D_\mu(s) = \frac{1}{\zeta(s)}$.  If $D_\mu(s)$ were to converge at some $s_0 \in \CC$ with $\sigma_0 = \operatorname{Re} s_0 < \Theta$, then $D_\mu(s)$ would be analytic on $\{s \in \CC: \operatorname{Re} s > \sigma_0\}$, and thus so would $\frac{1}{\zeta(s)}$.  However,  the meromorphic function  $\frac{1}{\zeta(s)}$ has poles in $\{s \in \CC: \operatorname{Re} s > \sigma_0\}$, since $\zeta(s)$, by the definition of $\Theta$, has zeros in that domain.  Thus, the domain of convergence of $D_\mu(s)$ must be contained in $\{s \in \CC: \operatorname{Re} s \geq \Theta\}$, whence $\Theta' \geq \Theta$.  
The reverse inequality $\Theta' \leq \Theta$ is equivalent to the fact that, for all $t \geq \frac{1}{2}$, if $\zeta(s)$ has no zeros $s$ with $\operatorname{Re} s> t$, then $\deg M \leq t$.   The proof of the latter fact is a straightforward generalization of the proof of \cite[Theorem 4.16]{broughan} or  \cite[Section 12.1]{edw} (which assume $t = \frac{1}{2}$, i.e., the Riemann hypothesis).   The proof that $\deg L = \Theta$ is similar.  Alternatively,  one has $\deg L = \deg M$ by Proposition \ref{LMprop},  which is proved later in this section.
\end{proof}

The following corollary of the theorem above was first proved by Littlewood in 1912 \cite{litt0}.

\begin{corollary}[{\cite{litt0} \cite[Theorems 4.16 and 4.18]{broughan} \cite[Section 12.1]{edw}}]
Each of the following conditions is equivalent to the Riemann hypothesis.
\begin{enumerate}
\item  $D_\mu(s)$ has abscissa of convergence $\frac{1}{2}$.
\item $M(x)$ has degree $\frac{1}{2}$.
\item $M(x) = o(x^{\frac{1}{2}+\varepsilon})$ for all $\varepsilon > 0$.  
\item $D_\lambda(s)$ has abscissa of convergence $\frac{1}{2}$.
\item $L(x)$ has degree $\frac{1}{2}$.
\item $L(x) = o(x^{\frac{1}{2}+\varepsilon})$ for all $\varepsilon > 0$.  
\end{enumerate}
\end{corollary}

The following problem is as at least as difficult as Problem \ref{mainproblem}.

\begin{outstandingproblem}
Compute $\dege M$ and $\dege L$.
\end{outstandingproblem}  

It is known that $$M(x) = \Omega_{\pm}(x^{\Theta-\varepsilon}) \ (x \to \infty)$$ for all $\varepsilon > 0$, where $\Theta = \deg(\li-\pi)$ is the Riemann constant \cite[p.\ 90]{ing2}.  It is also known that $$M(x) = \Omega_{\pm}(\sqrt{x}) \ (x \to \infty)$$ 
and, more generally \cite{ng}, if $\zeta(s)$ has a  zero of order $m$, then  $$M(x) = \Omega_{\pm}(\sqrt{x}\, (\log x)^{m-1}) \ (x \to \infty).$$
It is widely suspected that $M(x) \neq O(\sqrt{x}) \ (x \to \infty)$ (which is known to follow from the conjecture that the positive imaginary parts of the zeros of $\zeta(s)$ are linearly independent over $\QQ$ \cite[Theorem A]{ing}), despite earlier conjectures of Stieltjes, Mertens, et.\ al., to the contrary.  To date, the best known explicit bounds are
$$\limsup_{x \to \infty} \frac{M(x)}{\sqrt{x}} > 1.826054$$ and $$\liminf_{x \to \infty} \frac{M(x)}{\sqrt{x}} < -1.837625,$$
and there exists a real number $x \leq e^{1.004\cdot 10^{33}} $ such that $|M(x)| > \sqrt{x}$, but there are no such $x$ less than or equal to $10^{16}$ \cite{hur} \cite{sao}.

Various bounds on $M(x)$ are directly related, via (\ref{Dmu}), to bounds concerning $\zeta(s)$ and its zeros.   The following result, for example,  is due to Walfisz and Allison.

\begin{theorem}[{\cite[Theorem]{allison}  \cite[Hilfsatz V.5.5]{wal}}]
Let $\alpha \in (0,1]$ and $\beta \in \RR$.    The equivalent conditions of Theorem \ref{antiR2} hold if and only if there exists a constant $C > 0$ such that
$$M(x)= O \left(x e^{-C(\log x)^\alpha (\log \log x)^\beta} \right) \ (x \to \infty).$$
\end{theorem}

\begin{proof}
By \cite[Hilfsatz V.5.5]{wal}, condition (1) of Theorem \ref{antiR2} implies the given $O$ bound for $M(x)$, and the reverse implication follows from \cite[Theorem]{allison}.
\end{proof}

From the prime number theorem with error bound and Proposition \ref{antiR},  we immediately obtain the following corollaries.

\begin{corollary}[{\cite[Satz V.5.3]{wal}}]\label{walf}
There exists a constant $C > 0$
\begin{align*}
M(x) = O\left(xe^{ - C(\log x)^{3/5}(\log \log x)^{-1/5}}\right) \ (x \to \infty).
\end{align*}
\end{corollary}

\begin{corollary}
The anti-Riemann constant $\Theta_{-1} = \underline{\deg}(x-\log|\li(e^x)-\pi(e^x)|)$ is given by
\begin{align*}
\Theta_{-1} &   = -\inf\left\{t \in \RR: M(x) = O\left(xe^{-(\log x)^{-t}} \right) \ (x \to \infty)\right\}\\
 & =  \sup\left\{t \in \RR: M(x) = O\left(xe^{-(\log x)^{t}} \right) \ (x \to \infty)\right\}  \\ & = \underline{\deg} ( x -\log |M(e^x)|).
\end{align*}
\end{corollary}

A 2008 result due to Soundararajan, Balazard, and de Roton implies that the Riemann hypothesis is equivalent to
\begin{align}\label{sbd}
M(x) = O\left (\sqrt{x} e^{(\log x)^{1/2} (\log \log x)^{5/2+\varepsilon}}\right) \ (x \to \infty)
\end{align}
for all $\varepsilon > 0$ \cite{bal}.  

By \cite{hum},  all of the $O$ and $\Omega_\pm$ bounds noted above for $M(x)$ also hold for the function $L(x)$.  Thus, we obtain the following.

\begin{theorem}\label{multzero}
 One has the following.
\begin{enumerate}
\item $\deg M = \deg L = \Theta.$
\item One has $$(\tfrac{1}{2},0,0,0,\ldots) \leq \dege M \leq (1,-\infty,-\tfrac{3}{5},\tfrac{1}{5},0,0,0,\ldots)$$
and
$$(\tfrac{1}{2},0,0,0,\ldots) \leq \dege L \leq (1,-\infty,-\tfrac{3}{5},\tfrac{1}{5},0,0,0,\ldots).$$
\item The Riemann hypothesis is equivalent to $\deg M = \frac{1}{2}$, to $\deg L = \frac{1}{2}$, to
$$\dege M \leq (\tfrac{1}{2},\infty,\tfrac{1}{2},\tfrac{5}{2},0,0,0,\ldots),$$
and to
$$\dege L \leq (\tfrac{1}{2},\infty,\tfrac{1}{2},\tfrac{5}{2},0,0,0,\ldots).$$
\item The anti-Riemann hypothesis  is equivalent to $\deg M= 1$, to $\deg L = 1$,  to 
$$\dege M \geq (1,-\infty,-1,0,0,0,\ldots),$$
and to
$$\dege L \geq (1,-\infty,-1,0,0,0,\ldots).$$
\item If $\zeta(s)$ has a  zero of order $m$, then $$\dege M \geq (\tfrac{1}{2},m-1,0,0,0,\ldots)$$ and
$$\dege L \geq (\tfrac{1}{2},m-1,0,0,0,\ldots).$$
\end{enumerate}
\end{theorem}

Statement (2) of the theorem represents the best case and worst case scenarios, given what is known, for what $\dege M$ and $\dege L$ could be.  Probably the truth is somewhere in between.  Relevant conjectures are discussed in Section 14.3.

The functions $M$ and $L$ are related  \cite[p.\ 552]{hum} via the identity
$$L(x) = \sum_{d^2 \leq x} M(x/d^2)$$
and therefore, by M\"obius inversion, 
$$M(x) = \sum_{d^2 \leq x} \mu(d) L(x/d^2).$$
The first of these follows easily from the more basic identity 
$$\lambda(n) = \sum_{d^2 | n} \mu(n/d^2).$$
Using these relations and the fact that $M(x) \neq o(\sqrt{x}) \ (x \to \infty)$, we deduce the following.

\begin{proposition}\label{LMprop} Let $r$ be an eventually positive and nondecreasing real function defined on a neighborhood of $\infty$.
\begin{enumerate}
\item  If  $M(x) = O(x^{1/2}r(x)) \ (x \to \infty)$, then $L(x) = O(x^{1/2}r(x)\log x) \ (x \to \infty)$.
\item  If  $L(x) = O(x^{1/2}r(x)) \ (x \to \infty)$, then $M(x) = O(x^{1/2}r(x)\log x) \ (x \to \infty)$.
\item  $\dege M+(0,-1,0,0,0,\ldots) \leq  \dege L  \leq \dege M +(0,1,0,0,0,\ldots)$.
\item  If the anti-Riemann hypothesis is true, or, more generally, if $\dege_1 M = -\infty$, then $\dege L = \dege M$.
\end{enumerate}
\end{proposition}

\begin{proof}
Suppose that $r$ is positive and nondecreasing on $[N,\infty)$, where $N \geq 1$, and that $|M(x)| \leq C x^{1/2} r(x)$ for all $x \geq N$, where $C > 0$.   Let $D = \sum_{k = 1}^{N-1} |M(k)|$. Then one has
\begin{align*}
|L(x)| & =\left| \sum_{d^2 \leq x} M(x/d^2) \right|\\
& = \left|\sum_{x/N < d^2 \leq x} M(x/d^2)+  \sum_{d^2 \leq x/N} M(x/d^2)  \right|
\\ & \leq  \sum_{1 \leq  x/d^2 < N} |M(x/d^2)|+  \sum_{d \leq \sqrt{x/N}}| M(x/d^2) |
\\ & \leq D x^{1/2}+  C \sum_{d \leq \sqrt{x/N}} \frac{x^{1/2}}{d} r(x/d^2)  
 \\ & \leq  D x^{1/2}+Cx^{1/2} r(x) \sum_{d \leq \sqrt{x/N}} \frac{1}{d}
\\& = O(x^{1/2}r(x)\log x) \ (x \to \infty).
\end{align*}
This proves statement (1), and a similar argument yields statement (2).

Now, let $q(x) = M(x)/x^{1/2}$.
Let $r \in  \mathbb{L}_{>0}$ with $|q(x)| \leq r(x)$ for all $x  \gg 0$.  Since  $q(x) \neq o(1) \ (x \to \infty)$, one has $r(x) \neq o(1) \ (x \to \infty)$,  and therefore $r'(x)$ is eventually nonnegative, whence $r$ is eventually nondecreasing.  From statement (1), then, it follows that $$\frac{L(x)}{x^{1/2} \log x} = O(r(x)) \ (x \to \infty).$$  Thus, by Theorem \ref{infpropexp}, one has
$$\dege \frac{L(x)}{x^{1/2} \log x} \leq \dege q = \dege \frac{M(x)}{x^{1/2}}$$
and therefore
$$\dege L + (0,-1,0,0,0,\ldots) = \dege \frac{L}{ \log } \leq \dege M.$$
A similar argument, interchanging the roles of $M$ and $L$, yields
$$\dege M + (0,-1,0,0,0,\ldots) =  \dege \frac{M}{ \log }  \leq  \dege L $$
and
$$\dege M \leq \dege L+ (0,1,0,0,0,\ldots).$$
Statements (3) and (4) follow. 
\end{proof}

Note that Proposition \ref{LMprop} implies that $\deg M = \deg L$.

The estimates in Proposition \ref{LMprop} are somewhat coarse.  It seems  likely that more is true, namely, that $\dege L  = \dege M$.  Thus we pose the following problem.

\begin{problem}
Express $\dege M$ in terms of $\dege L$ and vice versa.  In particular, are they equal?
\end{problem}  

The following problem is much more difficult and may require computing both $\dege M$ and $\dege (\li-\pi)$.

\begin{outstandingproblem}\label{Mprobpi}
Express $\dege M$ in terms of $\dege (\li-\pi)$ and vice versa.  
\end{outstandingproblem}

Various conjectures regarding the Mertens function $M$ are discussed in Section 14.3.

\section{The Dirichlet divisor problem and its various analogues}

Another famous problem that can be expressed using the degree formalism regards the {\bf summatory divisor function}\index{summatory divisor function}
$$S_d(x) = \sum_{n  \leq x} d(n),$$
where $d(n) = \sum_{d|n} 1$ is the divisor function.  Recall from Section 3.7 that 
\begin{align*}
S_d(x) = x \log x + (2 \gamma-1)x + O(x^{t}) \ (x \to \infty)
\end{align*}
for $t = \frac{1}{2}$.    The {\bf Dirichlet divisor problem}\index{Dirichlet divisor problem} is the problem of determining the infimum of all such $t$ for which the $O$ bound above holds, which is equal to the degree of the function $$D(x) = S_d(x) - x \log x - (2 \gamma-1)x.$$  Hardy proved in 1914 that $D(x) = \Omega_{\pm}(x^{1/4}) \ (x \to \infty)$ and therefore $\deg D \in [\frac{1}{4},\frac{1}{2}]$.   It is widely conjectured that $\deg D = \frac{1}{4}$.    The smallest known $O$ bound for $D$ at the writing of this text was  proved by Huxley in 2003 \cite{hux2}, namely, 
 $$D(x) = O(x^{131/416}(\log x)^{26947/8320}) \ (x \to \infty)$$ and therefore  $\deg D \in [\frac{1}{4},\frac{131}{416}]$ and
$$(\tfrac{1}{4},0,0,0,\ldots) \leq  \dege D \leq (\tfrac{131}{416},\tfrac{26947}{8320},0,0,0,\ldots).$$  
 The Dirichlet divisor problem thus generalizes as follows.

\begin{outstandingproblem}
Compute $\dege D$.
\end{outstandingproblem}

\begin{remark}[An alternative expression for $\dege D$]
Note that
$$S_d(x) = \sum_{k \leq x}  \sum_{r|k} d(r) = \sum_{n \leq x}\left\lfloor\frac{x}{n} \right\rfloor,$$
and therefore
$$S_d(x) +\sum_{n \leq x}\left\{\frac{x}{n} \right\}= xH_{\lfloor x\rfloor}= x \left(\log x+\gamma +O\left(\frac{1}{x}\right)\right) = x \log x + \gamma x + O(1) \ (x \to \infty).$$
It follows that
$$D(x)  =- \sum_{n \leq x}\left\{\frac{x}{n} \right\} + (1-\gamma)x + O(1) \ (x \to \infty),$$
and therefore
$$\dege D = \dege \left( \sum_{n \leq x}\left\{\frac{x}{n} \right\} - (1-\gamma)x \right).$$
An interesting consequence of this is that
$$\gamma = 1-\lim_{x \to \infty} \frac{1}{x} \sum_{n \leq x}\left\{\frac{x}{n} \right\}.$$
\end{remark}

Another famous problem, analogous to the Dirichlet divisor problem, regards the number $$S_{|\mu|}(x) = \sum_{n \leq x} |\mu(n)|$$ of squarefree positive integers less than or equal to $x$.   Recall from Example \ref{sumex0}(2) that
\begin{align*}
Q(x) = O(x^{1/2}) \ (x \to \infty),
\end{align*}
 where $$Q(x) = S_{|\mu|}(x) - \frac{6}{\pi^2} x.$$ 
By \cite[Theorem 1.3]{eve} one has
$Q(x) \neq o(x^{1/4}) \ (x \to \infty).$
Therefore $\frac{1}{4} \leq \deg Q \leq \frac{1}{2}$ and in fact
$$(\tfrac{1}{4} , 0,0,0,\ldots) \leq \dege Q \leq (\tfrac{1}{2} , 0,0,0,\ldots).$$
  Assuming the Riemann hypothesis, it is known that $\deg Q \leq \frac{11}{35}$  \cite{liu},
and this is the best upper bound conditional on the Riemann hypothesis to date.  As with $D$, it is plausible that $\deg Q = \frac{1}{4}$.

\begin{outstandingproblem}
Compute  $\dege Q$.
\end{outstandingproblem}

Let $$R(x) = S_\phi(x) - \frac{3}{\pi^2}x^2,$$
where $\phi$ denotes Euler's totient.
By an 1874 result of Mertens \cite{mertphi}, one has
$$R(x) =  O(x \log x) \ (x \to \infty).$$
The best known $O$ bound for $R(x)$ is due to H.-Q. Liu \cite{liu2}:
$$R(x) = O(x (\log x)^{2/3}(\log \log x)^{1/3}) \ (x \to \infty),$$
The best known lower bound, namely,
$$R(x) = \Omega_{\pm} (x\sqrt{\log \log x}) \ (x \to \infty),$$
is due to Montgomery \cite{mont0}.
Thus, one has
$$(1,0,\tfrac{1}{2},0,0,0,\ldots) \leq \dege R \leq (1,\tfrac{2}{3}, \tfrac{1}{3}, 0,0,0,\ldots).$$
Montgomery conjectured in \cite{mont0} that
\begin{align}\label{montR1}
R(x) = O(x \log \log x) \ (x \to \infty)
\end{align}
and 
\begin{align}\label{montR2}
R(x) = \Omega_{\pm} (x\log \log x) \ (x \to \infty),
\end{align}
hence that
$$\dege R = (1,0,1,0,0,0,\ldots).$$

\begin{outstandingproblem}
Compute  $\dege \left(S_\phi(x)-\frac{3}{\pi^2} x^2 \right)$.
\end{outstandingproblem}

By \cite[Theorem 1]{mont0},  the multiplicative arithmetic function
$$\rho(n) = \frac{\phi(n)}{n} = \prod_{p |n} \left(1-\frac{1}{p} \right),$$
where $\rho(n)$ for any positive integer $n$ is equal to the probability that a randomly chosen integer from $1$ to $n$ is relatively prime to $n$, 
satisfies
 $$S_\rho(x)-\frac{6}{\pi^2} x =  \frac{R(x)}{x} + e^{-c\sqrt{x}} \ (x \to \infty)$$
for some $c > 0$, where $R(x) = S_\phi(x) - \frac{3}{\pi^2}x^2$.  It follows that
$$\dege \left(S_\rho(x)-\frac{6}{\pi^2} x \right) = \dege R + (-1,0,0,0,\ldots).$$

Finally, recall from Example \ref{sumex}(2) that 
$$S_{2^\omega} (x) = \frac{1}{\zeta(2)}x \log x + \left(\frac{2\gamma-1}{\zeta(2)}- \frac{2\zeta'(2)}{\zeta^2(2)}\right) x+O(x^{1/2}\log x) \ (x \to \infty).$$
It is known, under the Riemann hypothesis, that the error term is $O(x^t)$ for all $t > \frac{1}{4}$.  It follows that
$$\dege \left(S_{2^\omega} (x) - \frac{1}{\zeta(2)}x \log x- \left(\frac{2\gamma-1}{\zeta(2)}- \frac{2\zeta'(2)}{\zeta^2(2)}\right)x \right) \leq (\tfrac{1}{2}, 1,0,0,0,\ldots),$$
and the given function has degree at most $\frac{1}{4}$ under the  Riemann hypothesis.

\begin{outstandingproblem}
Compute $\dege \left(S_{2^\omega} (x) - \frac{1}{\zeta(2)}x \log x- \left(\frac{2\gamma-1}{\zeta(2)}- \frac{2\zeta'(2)}{\zeta^2(2)}\right)x \right)$.
\end{outstandingproblem}

\section{The degree and logexponential degree of a summatory function}

This section is not used in later chapters.

In this section, we seek to relate the degree (resp., logexponential degree) of  an arithmetic function $f: \ZZ_{> 0} \longrightarrow \CC$ to that of its  summatory function $S_f$, or,  more generally, to that of the function
$$S_{f,0}(x) = \begin{cases} D_f(0)-S_f(x)  = \sum_{n > x} f(n) & \quad \text{if } D_f(0) = \sum_{n = 1}^\infty f(n) \text{ exists} \\
    S_f(x)  & \quad \text{otherwise}.
\end{cases}$$   A variety of examples is provided in Table  \ref{table1}, which extends Table  \ref{table1a}.  Note that, of the functions in the table, only the logexponential degrees of the functions $M$ and $L$ are unknown.  

For sufficiently nice arithemetic functions $f$,  one expects that  the degree of the running average $\frac{1}{x}S_f(x)$ of $f$ is equal to the degree of $f$,  i.e., that $\deg S_f = 1+\deg f$.  This expectation is corroborated by the  following elementary result.

\begin{table}[!htbp]
\scriptsize
  \caption{\centering Dirichlet series and summatory function of arithmetic functions $f(n)$}
\begin{tabular}{|l||l|l|l|l|l|} \hline
$f(n)$ &  $D_f(s)$ &  $S_{f,0}(x)$ & $S_{f,0}(x) \sim$ &  $\dege S_{f,0}$  & $\dege f$ \\  \hline\hline
$1$ & $\zeta(s)$  &  $\lfloor x \rfloor$ & $x$ &  $(1,0,0,0,\ldots)$  & $(0,0,0,\ldots)$  \\ \hline
$(-1)^{n-1}$ & $\left(1-2^{1-s}\right) \zeta(s)$  &  $\sum_{n\leq x} (-1)^{n-1}$  & $\sum_{n\leq x} (-1)^{n-1}$ &  $(0,0,0,0,\ldots)$  & $(0,0,0,\ldots)$  \\ \hline
$\chi_{\pp}(n)$ & $P(s)$  &  $\pi(x)$ &  $\li x$ &  $(1,-1,0,0,\ldots)$  &  $(0,0,0,\ldots)$  \\ \hline
$\chi_{\pp}(n)\log n$ & $-P'(s)$  &  $\vartheta(x)$ &  $x$ &  $(1,0,0,0,\ldots)$  & $(0,1,0,0,\ldots)$  \\ \hline
$\chi_{\pp}(n)\frac{\log n}{n}$ & $-P'(s+1)$  &  $ \sum_{p \leq x} \frac{\log p}{p}$ &  $\log x $ &  $(0,1,0,0,0,\ldots)$  & $(-1,1,0,0,\ldots)$  \\ \hline
$\Lambda(n)$ & $\frac{\zeta'(s)}{\zeta(s)}$ &  $\psi(x)$ & $x$  & $(1,0,0,0,\ldots)$  & $(0,1,0,0,\ldots)$  \\ \hline
$\frac{\Lambda(n)}{\log n}$ & $\log \zeta(s)$ &  $\Pi(x)$ & $\li(x)$  & $(1,-1,0,0,\ldots)$  &  $(0,0,0,\ldots)$ \\ \hline
$n^a$, $a > -1$ &  $\zeta(s-a)$ & $\sum_{n \leq x}n^a$ & $\frac{x^{a+1}}{a+1}$ & $(a+1,0,0,0,\ldots)$  &  $(a,0,0,0,\ldots)$ \\ \hline
$\frac{1}{n}$ &  $\zeta(s+1)$ & $H_{\lfloor x \rfloor}$ & $\log x$ & $(0,1,0,0,0,\ldots)$  & $(-1,0,0,0,\ldots)$  \\ \hline
$n^a$, $a < -1$ &  $\zeta(s-a)$ & $\sum_{n > x}n^a$ &  $-\frac{x^{a+1}}{a+1} $ &  $(a+1,0,0,0,\ldots)$  & $(a,0,0,0,\ldots)$  \\ \hline
$\chi_{\pp}(n)n^a$, & $P(s-a)$  &  $\sum_{p \leq x}p^a$ & $\li(x^{a+1})$ &   $(a+1,-1,0,0,\ldots)$  & $(a,0,0,0,\ldots)$  \\ 
   $ a > -1$&   &   &   &  &   \\  \hline
$\chi_{\pp}(n)\frac{1}{n}$ & $P(s+1)$  &  $\sum_{p \leq x}\frac{1}{p}$ &  $\log \log x$ &  $(0,0,1,0,0,0,\ldots)$  &   $(-1,0,0,0,\ldots)$ \\ \hline
$\chi_{\pp}(n)n^a$, & $P(s-a)$  &  $\sum_{p > x}p^a$ & $-\li(x^{a+1})$ & $(a+1,-1,0,0,\ldots)$  &  $(a,0,0,0,\ldots)$ \\ 
   $ a < -1$&   &   &   &  &   \\  \hline
$(\log n)^k$, $k \in \ZZ_{\geq 0}$ & $(-1)^k\zeta^{(k)}(s)$  &  $\sum_{n \leq x} (\log n)^k$ & $x (\log x)^k$ & $(1,k,0,0,\ldots)$  & $(0,k,0,0,0,\ldots)$  \\ \hline
$n^a(\log n)^k$, & $(-1)^k\zeta^{(k)}(s-a)$  &  $\sum_{n \leq x} n^a(\log n)^k$ & $\frac{ x^{a+1}(\log x)^k}{a+1}$ & $(a+1,k,0,0,\ldots)$  & $(a,k,0,0,0,\ldots)$  \\ 
   $ a > -1,  \ k \in \ZZ_{\geq 0}$ &   &   &   &  &   \\  \hline
$\frac{1}{n}(\log n)^k$, $k \in \ZZ_{\geq 0}$ & $(-1)^k\zeta^{(k)}(s+1)$  &  $\sum_{n \leq x} \frac{1}{n}(\log n)^k$ & $\frac{ (\log x)^{k+1}}{k+1}$ & $(0,k+1,0,0,\ldots)$  &  $(-1,k,0,0,0,\ldots)$  \\ \hline
$n^a(\log n)^k$, & $(-1)^k\zeta^{(k)}(s-a)$  &  $\sum_{n >x} n^a(\log n)^k $ & $-\frac{ x^{a+1}(\log x)^k}{a+1}$ & $(a+1,k,0,0,\ldots)$  & $(a,k,0,0,0,\ldots)$  \\
   $ a < -1,  \ k \in \ZZ_{\geq 0}$ &   &  &   &  &   \\  \hline
$\sigma_a(n)$, $a > 1$ & $\zeta(s)\zeta(s-a)$ &  $\sum_{n \leq x} \sigma_a(n)$ & $\frac{\zeta(a+1)x^{a+1}}{a+1}$ &  $(a+1,0,0,0,\ldots)$  & $(a,0,0,0,\ldots)$  \\ \hline
$\sigma(n) = \sigma_1(n)$ & $\zeta(s)\zeta(s-1)$ &  $\sum_{n \leq x} \sigma(n)$ & $\frac{\zeta(2)x^{2}}{2}$ &  $(2,0,0,0,\ldots)$  &  $(1,0,1,0,0,0,\ldots)$ \\ \hline
$\sigma_a(n)$, $0<a<1$ & $\zeta(s)\zeta(s-a)$ &  $\sum_{n \leq x} \sigma_a(n)$ & $\frac{\zeta(a+1)x^{a+1}}{a+1}$ &  $(a+1,0,0,0,\ldots)$  & $(a,\infty,1-a,-1,0,0\ldots)$  \\ \hline
$d(n) = \sigma_0(n)$ & $\zeta(s)^2$ &  $\sum_{n \leq x} d(n)$ & $x \log x$ &  $(1,1,0,0,\ldots)$  & $(0,\infty,1,-1,0,0,0,\ldots)$   \\ \hline
$\sigma_{a}(n)$, $a< 0$  & $\zeta(s)\zeta(s-a)$ &  $\sum_{n \leq x} \sigma_{a}(n)$ & $\zeta(1-a)x$ &  $(1,0,0,0,\ldots)$  & $\dege \sigma_{-a} + (a,0,0,\ldots)$  \\ \hline
$\omega(n)$ & $\zeta(s)P(s)$ &  $\sum_{n \leq x} \omega(n)$ & $x \log \log x$ &  $(1,0,1,0,0,0,\ldots)$  & $(0,1,-1,0,0,0,\ldots)$  \\ \hline
$\Omega(n)$ & $\zeta(s)\sum_{n = 1}^\infty P(ns)$ &  $\sum_{n \leq x} \Omega(n)$ & $x \log \log x $ & $(1,0,1,0,0,0,\ldots)$     & $(0,1,0,0,0,\ldots)$   \\  \hline
$\mu(n)$ & $\frac{1}{\zeta(s)}$ &  $M(x)$ & $M(x)$ &  $(\Theta,m_1,m_2,\ldots)$  &   $(0,0,0,\ldots)$ \\ \hline
$\lambda(n)$ & $\frac{\zeta(2s)}{\zeta(s)}$ &  $L(x)$ &  $L(x)$ & $(\Theta,l_1,l_2,\ldots)$  &  $(0,0,0,\ldots)$ \\ \hline
$\mu(n)^2$ & $\frac{\zeta(s)}{\zeta(2s)}$ &  $Q(x)$ & $\frac{6}{\pi^2} x$ &  $(1,0,0,0,\ldots)$  &  $(0,0,0,\ldots)$  \\ \hline
$\phi(n)$ & $\frac{\zeta(s-1)}{\zeta(s)}$ &  $\sum_{n \leq x} \phi(n)$ & $\frac{3}{\pi^2} x^2$ &  $(2,0,0,0,\ldots)$  &  $(1,0,0,0,\ldots)$ \\ \hline
$\frac{\phi(n)}{n}$ & $\frac{\zeta(s)}{\zeta(s+1)}$ &  $\sum_{n \leq x} \frac{ \phi(n)}{n}$ & $\frac{6}{\pi^2} x$ &  $(1,0,0,0,\ldots)$  &  $(0,0,0,\ldots)$ \\ \hline
$2^{\omega(n)}$ & $\frac{\zeta(s)^2}{\zeta(2s)}$ &  $\sum_{n \leq x} 2^{\omega(n)}$ &  $ \frac{6}{\pi^2}x \log x$ &  $(1,1,0,0,0,\ldots)$  &  $(0,\infty,1,-1,0,0,0,\ldots)$  \\  \hline
\end{tabular}\label{table1}
\end{table}

\begin{proposition}
Let $f$ be an arithmetic function, and let $N$ be a nonnegative integer.
\begin{enumerate}
\item If $\deg f \geq -1$, then $\deg S_f \leq 1+\deg f$.
\item If $\deg f < -1$, then $D_f(0) = \sum_{n = 1}^\infty f(n)$ converges absolutely and $\deg (D_f(0) -S_f) \leq 1+\deg f$.
\item If $\infty >\degl_k f \geq -1$ for all $k \leq N$, then $\degl_k S_f < 1+\degl_k f$ for the smallest $k \leq N$, if any, such that $\degl_k S_f \neq 1+\degl_k f$.
\item If $f$ is nonnegative, then $\dege S_f  \geq \dege f$.
\item If $f$ is nonnegative and nonincreasing, then $\dege S_f \geq \dege f \oplus(1,0,0,0,\ldots)$. 
\item If $f$ is nonnegative and nonincreasing and $\deg f \geq -1$, then $\deg S_f = 1+ \deg f$.
\end{enumerate}
\end{proposition}

\begin{proof}
Let $t > \deg f$, so that $f(n)  =  O(n^t) \ (x \to \infty)$, whence there exists a $C > 0$ such that $|f(n)| \leq Cn^t$ for all $n$.  Suppose that $\deg f  \geq -1$, so that $t > -1$.  Then, for all $x \geq 0$, one has
\begin{align*}
|S_f(x)|   = \left|\sum_{n \leq x} f(n)\right|  \leq  C \sum_{n \leq x} n^t  = O (x^{t+1}) \ (x \to \infty),
\end{align*}
so that $\deg S_f \leq t+1$ and thus $\deg S_f-1 \leq t$.  Taking the infimum over all $t > \deg f$, we see that $\deg S_f -1 \leq \deg f$.   
This proves statement (1).

Suppose now that $\deg f < t < -1$.  Absolute convergence of the sum $D_f(0)$ follows from convergence of the sum $\sum_{n =1}^\infty n^t = \zeta(-t)$.   One then has
\begin{align*}
|D_f(0) -S_f(x)|   = \left|\sum_{n > x} f(n)\right|  \leq  C \sum_{n > x} n^t  = O (x^{t+1}) \ (x \to \infty),
\end{align*}
so that $\deg (D_f(0) -S_f) \leq t+1$ and thus $\deg (D_f(0) -S_f) -1 \leq t$.  Taking the infimum over all $t$ as chosen, we see that $\deg (D_f(0) -S_f) -1 \leq \deg f$.   
This proves statement (2).

Statement (1) implies statement (3) for $N = 0$.  To prove (3) for $N = 1$,  we may suppose that $\deg S_f = 1+\deg f$, and $-1 \leq d:= \deg f  < \infty$, and also $\degl_1 f \geq -1$.  Let $f(x) = 0$ if $x$ is not a positive integer, and let $t > \degl_1 f$, so that $t > -1$ and there exists a $C > 0$ such that
$$|f(e^x)e^{-dx} | \leq C x^t$$ for all $x$. Then one has
\begin{align*}
\left|S_f(e^x)e^{-(d+1)x} \right| & \leq \sum_{n \leq e^x }\left|f(n)e^{-(d+1)x}\right|  \\ 
& \leq  \sum_{n \leq e^x} \left|f(n)n^{-(d+1)}\right| \\
& = \sum_{y = \log n \leq x} \left|f(e^y)e^{-(d+1)y}\right| \\
& \leq \sum_{y = \log n \leq x} \frac{Cy^t}{n}   \\
& = C\sum_{n \leq e^x} \frac{(\log n)^t}{n}   \\
&  = O (x^{t+1}) \ (x \to \infty),
\end{align*}
where the $O$ bound follows by comparison with the integral $\int_1^{e^x} \frac{(\log u)^t}{u}\, du = \frac{1}{t+1}x^{t+1}$.  Therefore $\degl_1 S_f \leq  1+ \degl_1 f$.  This proves statement (3) for $N = 1$.  The  statement for $N = 2$ then follows from the statement for $N = 1$, using the integral $\int_{e}^{e^{e^x}} \frac{(\log \log u)^t}{u \log u }\, du = \frac{1}{t+1}x^{t+1}$.  It is clear, then, how to proceed by induction.

To prove statement (4),  note that $S_f(n) \geq f(n) \geq 0$ for all $n$, so that $\dege S_f(x) = \dege S_f(n) \geq \dege f(n)$.  Similarly, to prove (5),  note that $S_f(n) \geq nf(n) \geq 0$ for all $n$, so that $\dege S_f(x) = \dege S_f(n) \geq \dege nf(n) =  \degl f \oplus (1,0,0,0,\ldots)$.  Finally, statement (6) follows from statements (1) and (5).
\end{proof}

\begin{remark}[Subcases of $\deg f \leq -1$]
Let $f$ be an arithmetic function.  The case $\deg f \leq -1$ splits into many subcases.  For example, one has the following.
\begin{enumerate}
\item  If $f(n) = O\left(\frac{1}{n} \right) \ (n \to \infty)$, then $S_f(x) = O(\log x) \ (x \to \infty)$ and thus $\dege S_f \leq (0,1,0,0,0,\ldots)$.
\item  If $f(n) = O\left(\frac{1}{n \log n} \right) \ (n \to \infty)$, then $S_f(x) = O(\log \log x) \ (x \to \infty)$ and thus $\dege S_f \leq (0,0,1,0,0,0,\ldots)$.
\item  If $f(n) = O\left(\frac{1}{n \log n \log \log n} \right) \ (n \to \infty)$, then $S_f(x) = O(\log \log \log x) \ (x \to \infty)$ and thus $\dege S_f \leq (0,0,0,1,0,0,0,\ldots)$.
\end{enumerate}
And so on.
\end{remark}

By Theorem \ref{hardintth} and Proposition \ref{iwanncor}, we have the following.

\begin{proposition}\label{iwannapp}
Let $f: [1,\infty) \longrightarrow \RR$ be any Hardian function of finite degree that is continuous and monotonic on $[M,\infty)$ for some $M\geq 1$,  with $\dege f \neq (-1,-1,-1,\ldots)$, and let $N$ be the smallest nonnegative integer such that $\dege_N f \neq -1$.    Then  $\int_M^x f(t) \, dt$ converges if and only if $\sum_{n = 1}^\infty f(n)$ converges, if and only if $\dege f < (-1,-1,-1,\ldots)$, if and only if  $\dege_N f < -1$.   If both diverge, then let
 $F_1(x) = \int_M^x f(t) \, dt$  and $S_1(x)  = S_f(x)$.    If both converge, then let $F_2(x) = \int_x^\infty f(t) \, dt$  and $S_2(x)  = \sum_{n > x} f(n)$.   In either case, one has
 $$S_i(x) = F_i(x) +C_i+ O(f(x)) = F_i(x) +o(F_i(x))  \sim F_i(x) \ (x \to \infty)$$
 for some constant $C_i$, where $C_i = 0$ if $i = 2$.   Consequently, one has
 $$\dege S_i  = \dege F_i = \dege f \oplus (1,1,1,\ldots,1,0,0,0,\ldots) \neq (0,0,0,\ldots),$$
where $(1,1,1,\ldots,1,0,0,0,\ldots)$ denotes the vector with $N+1$ $1$s followed by a  tail of $0$s.  Moreover, since $F_i$ is Hardian,   both $F_i$ and $S_i$ have exact logexponential degree.
\end{proposition}

 From the proposition above, we deduce the following.

\begin{proposition}
Let $f$ be an arithmetic function of finite degree with $\dege f \neq (-1,-1,-1,\ldots)$,  let $N$ be the smallest nonnegative integer such that $\dege_N f \neq -1$,  and let $(1,1,1,\ldots,1,0,0,0,\ldots)$ denote the vector with $N+1$ $1$s followed by a  tail of $0$s.   One has the following.
\begin{enumerate}
\item If $\dege f < (-1,-1,-1,\ldots)$, then $\sum_{n = 1}^\infty f(n)$ converges absolutely,  and one has
$$\dege \sum_{n >x} f(n) \leq \dege f \oplus (1,1,1,\ldots,1,0,0,0,\ldots) < (0,0,0,\ldots).$$
\item  If $\dege f > (-1,-1,-1,\ldots)$, then $$\dege S_f  \leq \dege f \oplus (1,1,1,\ldots,1,0,0,0,\ldots).$$
\item  If   $\underline{\dege}\,  f > (-1,-1,-1,\ldots)$, then $\sum_{n = 1}^\infty f(n)$ does not converge absolutely, and one has
$$\underline{\dege}\, S_{|f|}  \geq \underline{\dege}\, |f| \oplus (1,1,1,\ldots,1,0,0,0,\ldots)> (0,0,0,\ldots).$$
\item If $\dege f > (-1,-1,-1,\ldots)$ and $f$ is eventually positive and has exact logexponential degree,  then $S_f$ has exact logexponential degree
$$\dege S_f  = \dege f \oplus (1,1,1,\ldots,1,0,0,0,\ldots)> (0,0,0,\ldots).$$
\end{enumerate}
\end{proposition}

 \begin{proof}
 Suppose that $\dege f < (-1,-1,-1,\ldots)$.
 Let $r \in \mathbb{L}_{>0}$ with $\dege f < \dege r$.  In fact,  since $\dege_N f < -1$, there is such an $r$ with $\dege_k r = -1$ for all $k < N$ and $\dege_N f < \dege_N r < -1$.  It follows that  $f(n) = o(r(n)) \ (n \to \infty)$ and $\dege r < (-1,-1,-1,\ldots)$.    We assume without loss of generality that $r(x) > 0$ for all $x \geq 1$.  By Proposition \ref{iwannapp},  one has
 $$\dege \sum_{n > x} r(n) = \dege r \oplus (1,1,1,\ldots,1,0,0,0,\ldots) < (0,0,0,\ldots)$$
 and also $\sum_{n = 1}^\infty r(n)$ converges.    
 Since $f(n) = o(r(n)) \ (n \to \infty)$,  the sum $\sum_{n =1}^\infty f(n)$ is absolutely convergent and  one has
 $$ \sum_{n > x} f(n)  = o\left( \sum_{n > x} r(n) \right) \ (x \to \infty),$$
and therefore
  $$\dege \sum_{n > x} f(n) \leq \dege \sum_{n > x} r(n) =  \dege r \oplus (1,1,1,\ldots,1,0,0,0,\ldots).$$
  Taking the infimum over all $r$ as chosen, we see that
    $$\dege \sum_{n > x} f(n) \leq  \dege f \oplus (1,1,1,\ldots,1,0,0,0,\ldots).$$
    This proves (1).

    To prove (2),  suppose that $\dege f > (-1,-1,-1,\ldots)$.  Suppose first that $$\dege f  < (-1,-1,-1,\ldots,-1,\infty,1,0,0,0,\ldots),$$
    where $\infty$ is preceded by $N+1$ $-1$s.
By Proposition \ref{degerprop}, there is an $r \in \mathbb{L}_{>0}$ with 
$$\dege f  < \dege r  < (-1,-1,-1,\ldots,-1,\infty,1,0,0,0,\ldots).$$
It follows that  $f(n) = o(r(n)) \ (n \to \infty)$.   Moreover, by Proposition \ref{iwannapp},  one has
 $$\dege S_r = \dege r \oplus (1,1,1,\ldots,1,0,0,0,\ldots).$$
 It follows that
 $$ S_f(x)  = O\left( S_r(x) \right) \ (x \to \infty)$$ and therefore 
  $$\dege S_f \leq \dege S_r = \dege r \oplus (1,1,1,\ldots,1,0,0,0,\ldots).$$ 
  Taking the infimum over all $r$ as chosen, we see that
    $$\dege S_f \leq  \dege f \oplus (1,1,1,\ldots,1,0,0,0,\ldots).$$
    
    Suppose, on the other hand, that  $$\dege f  = (-1,-1,-1,\ldots,-1,\infty,1,0,0,0,\ldots).$$  We may let $r \in \mathbb{L}_{>0}$ with $$\dege r  = (-1,-1,-1,\ldots,-1,t,0,0,0,\ldots)> \dege f,$$
where $\dege_{N-1} r = t \in (-1,0)$.   
 By Proposition \ref{iwannapp},  one has 
 $$\dege S_r = \dege r \oplus (1,1,1,\ldots,1,0,0,0,\ldots) = (0,0,0,\ldots,0,t+1,0,0,0,\ldots).$$
Again, one has $f(n) = o(r(n)) \ (n \to \infty)$ and therefore $$ S_f(x)  = O\left( S_r(x) \right) \ (x \to \infty),$$ whence
  $$\dege S_f \leq \dege S_r = \dege r =(0,0,0,\ldots,0,t+1,0,0,0,\ldots).$$
  Taking the infimum over all $r$ as chosen, we see that
        $$\dege S_f \leq  (0,0,0,\ldots,0,0,\infty, 1,0,0,0,\ldots) = \dege f \oplus (1,1,1,\ldots,1,0,0,0,\ldots).$$
    This proves (2).
  
    Finally, suppose that  $\underline{\dege} \, f > (-1,-1,-1,\ldots)$.   By Proposition \ref{oexppropstrong}, there exists an $r \in \mathbb{L}_{>0}$ such that $\underline{\dege}\,  f > \dege r > (-1,-1,-1,\ldots)$.  By Proposition \ref{iwannapp}, the sum $\sum_{n = 1}^\infty r(n)$  diverges, and, by Proposition \ref{oexpprop}, one has $f \gg r$.  It follows, then,  that the sum $\sum_{n = 1}^\infty |f(n)|$ also diverges.  The remainder of the proof of (3) is similar to the proof of (1), and (4) follows from (2) and (3).  This completes the proof.
 \end{proof}

\chapter{The Riemann zeta function $\zeta(s)$}

In this chapter, we apply the degree and  logexponential degree formalisms to the study of the Riemann zeta function $\zeta(s)$ and its zeros.

\section{The Lindel\"of hypothesis and density hypothesis}

The famous {\bf Lindel\"of hypothesis}\index{Lindel\"of hypothesis}, conjectured by  Lindel\"of  in 1908 \cite{lind}, is the conjecture that $$\zeta(\tfrac{1}{2}+it) = o(t^{a}) \ (t \to \infty), \quad \forall a > 0.$$  For reasons described below, the conjecture is implied by the Riemann hypothesis.   For all  $\sigma \in \RR$, let $$\mu(\sigma) =  \inf\{a \in \RR: \zeta(\sigma+it) = o(t^a) \ (t \to \infty) \},$$
so that $$\mu(\sigma) = \deg \zeta(\sigma + it).$$ 
The Lindel\"of hypothesis  is  equivalent to $\mu(\frac{1}{2}) = 0$.  To date, the best known upper bound for $\mu(\tfrac{1}{2})$ is  $\mu(\tfrac{1}{2})\leq  \tfrac{13}{84}$, proved by Bourgain in 2017 \cite{bourgain}.   It is known that $\mu(\sigma) = 0$ for all $\sigma \geq 1$ and that the functional equation of $\zeta(s)$ implies that $$\mu(\sigma) = \mu(1 -\sigma) - \sigma + \tfrac{1}{2}$$ for all $\sigma$.  It is also known  that the function $\mu(\sigma)$ is nonincreasing, convex, and continuous \cite[p.\ 25]{ivic} \cite{lind}.   Consequently, one has the following.

\begin{proposition} One has the following.
\begin{enumerate}
\item $\mu(\sigma) = \mu(1 -\sigma) - \sigma + \tfrac{1}{2}$ for all $\sigma \in \RR$.
\item $\mu(\sigma)$ is a nonincreasing, convex, and continuous function of $\sigma$.
\item $\mu(\sigma) = 0$ for all $\sigma \geq 1$.
\item $\mu(\sigma) = \tfrac{1}{2}-\sigma$ for all $\sigma \leq 0$.
\item $\mu(\sigma) \leq \tfrac{1}{2}(1-\sigma)$ for all $\sigma \in [0,1]$.
\item $\mu(\tfrac{1}{2}) \leq \tfrac{1}{4}$. 
\end{enumerate}
\end{proposition}

\begin{figure}[ht!]
\includegraphics[width=70mm]{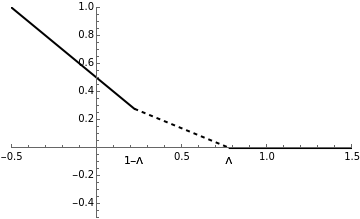}
\caption{\centering Graph of $\mu(\sigma)$ on $[-\tfrac{1}{2},\tfrac{3}{2}] \backslash (1-\Lambda,\Lambda)$, with upper bound $\mu(\sigma) \leq \tfrac{1}{2}(\Lambda-\sigma)$ on $[1-\Lambda,\Lambda]$}
   \label{lindelofg}
\end{figure}

Now, let$$\Lambda = \inf\{ \sigma \in \RR: \mu(\sigma) = 0\} = \inf\{ \sigma \in \RR:  \forall a > 0 \ (\zeta(\sigma+it) = o(t^{a}) \ (t \to \infty))\} \in [\tfrac{1}{2},1].$$
By the following corollary,  which is illustrated by Figure \ref{lindelofg}, the constant $\Lambda$ is an analogue of the Riemann constant $\Theta$ for the  Lindel\"of hypothesis.  As with the Riemann constant, the only known bounds on $\Lambda$ are $\tfrac{1}{2} \leq \Lambda  \leq 1$.    In particular,  no one has yet ruled out the {\bf anti-Lindel\"of hypothesis}\index{anti-Lindel\"of hypothesis} $\Lambda = 1$.

\begin{corollary} One has the folllowing.
\begin{enumerate}
\item  $\Lambda$ is the smallest real number such that $\mu(\sigma) = 0$ for all $\sigma \geq \Lambda$.
\item $1-\Lambda$ is the largest real number such that $\mu(\sigma) = \tfrac{1}{2}-\sigma$ for all $\sigma \leq \Lambda$.
\item One has $1-\Lambda = \sup\{\sigma \in \RR: \mu(\sigma) = \tfrac{1}{2}-\sigma\} \in [0,\tfrac{1}{2}]$.
\item The Lindel\"of hypothesis holds if and only if $\Lambda = \tfrac{1}{2}$, if and only if $1-\Lambda = \tfrac{1}{2}$.
\item $\mu(\sigma) \leq \tfrac{1}{2}(\Lambda-\sigma)$ if and only if $\sigma \in [1-\Lambda,\Lambda]$, for all $\sigma \in \RR$.
\item $\mu(\tfrac{1}{2})\leq  \tfrac{1}{2}(\Lambda-\tfrac{1}{2}) \leq \tfrac{1}{4}$.
\end{enumerate}
\end{corollary}

As a consequence, statements (1)--(3) of the following proposition are equivalent.

\begin{proposition}[{\cite{back} \cite[p.\ 45]{ivic}}]\label{lindelof}
The following conditions are equivalent.
\begin{enumerate}
\item The Lindel\"of hypothesis holds.
\item $\mu(\frac{1}{2}) = 0$, that is, $\deg \zeta(\frac{1}{2}+it) = 0$.
\item One has
\begin{align*}
\mu(\sigma)  = \begin{cases}  \frac{1}{2}-\sigma & \quad   \text{if } \sigma \leq \frac{1}{2}  \\
  0 & \quad  \text{if } \sigma \geq \frac{1}{2}.
\end{cases}
\end{align*}
\item $\deg \frac{1}{x} \int_0^x |\zeta(\frac{1}{2}+it)|^k \, dt = 0$ for all even positive integers $k$.
\item For every $\varepsilon > 0$, the number of zeros of $\zeta(s)$ with real part at least $\frac{1}{2}+\varepsilon$ and with imaginary part lying in $(T,T+1]$ is $o(\log T)$ as $T \to \infty$.
\end{enumerate}
\end{proposition}

The equivalence of the Lindel\"of hypothesis and statement (5) of the proposition was proved by R.\ Backlund in 1918--19 \cite{back}.  Moreover, the Riemann hypothesis states that there are no zeros in the region described in statement (5), for all $\varepsilon >0$ and all $T > 0$.  Thus, one has the following.

\begin{corollary}
The Riemann hypothesis implies the Lindel\"of hypothesis. 
\end{corollary}

 By the Riemann--von Mangoldt formula, the number  $N(T+1)-N(T)$ of zeros of $\zeta(s)$ with imaginary part lying in $(T,T+1]$ is asymptotic to $\log T$ as $T \to \infty$.  Thus, the Lindel\"of hypothesis is also equivalent to the statement that, for all $\varepsilon > 0$, the number of zeros of $\zeta(s)$ with real part lying in $[\frac{1}{2}-\varepsilon,\frac{1}{2}+\varepsilon]$ and with imaginary part lying in $(T,T+1]$ is asymptotic to $\log T$ as $T \to \infty$, and the Riemann hypothesis implies that this holds even for $\varepsilon = 0$.  
 
 
 The following proposition generalizes parts of Proposition \ref{lindelof}.  The proof follows that of  \cite{back}.   (See also https://mathoverflow.net/questions/486575/.)
 
\begin{proposition}[{cf.\ \cite{back} \cite[p.\ 45]{ivic}}]
The constant $\Lambda$ is the smallest constant $\kappa \in [\tfrac{1}{2},1]$ such the number of zeros of $\zeta(s)$ with real part at least $\kappa+\varepsilon$ and with imaginary part lying in $(T,T+1]$ is $o(\log T)$ as $T \to \infty$.  Consequently, one has $\Lambda \leq \Theta$.
\end{proposition}

The discussion above motivates the following problem.


\begin{outstandingproblem}\label{zetaorder}
Compute  $\dege \zeta(\sigma+it)$ for all $\sigma \in \RR$.
\end{outstandingproblem}

We first address the case where $\sigma = 1$.   In 1910, Bohr and Landau proved that
$$\zeta(1+it) \neq o(\log \log t) \ (t \to \infty),$$
while in 1912 Littlewood proved under the assumption of the Riemann hypothesis that
$$\zeta(1+it) = O(\log \log t \log \log \log t) \ (t \to \infty).$$
These results were later improved by Titchmarsh, who proved that 
$$e^\gamma \leq \limsup_{t \to \infty} \frac{|\zeta(1+it)|}{\log \log t} \leq 2 e^\gamma,$$
where the first inequality is unconditional and the second assumes the Riemann hypothesis \cite[Chapter VII]{tit}.  Since then one has $\zeta(1+it) = O(\log \log t) \ (t \to \infty)$
but $\zeta(1+it) \neq o(\log \log t) \ (t \to \infty)$, one has $$\dege \zeta(1+it) = (0,0,1,0,0,0,\ldots),$$ all assuming the Riemann hypothesis.  Thus a solution to Problem \ref{zetaorder} is known for $\sigma = 1$, on condition of the Riemann hypothesis.

The analogous problem for the values of $\zeta(s)$ on the critical line is more difficult.   In 1924, Littlewood proved \cite[Theorem 12]{litt24} that there exists a constant $A >0$ such that
$$\zeta(\tfrac{1}{2}+it) = O\left( \exp\left(A\frac{\log t}{\log \log t} \right) \right) \ (t \to \infty),$$
on condition of the Riemann hypothesis.
Moreover, in 2017, A.\ Bondarenko and K.\ Seip proved \cite{bonda} that, for all $\varepsilon > 0$, one has
$$\zeta(\tfrac{1}{2}+it) \geq \exp\left( \left(\tfrac{1}{2}-\varepsilon \right)\frac{\sqrt{\log t \log \log \log t}}{\sqrt{\log \log t}} \right)  \ (t \to \infty)$$
for an unbounded set of $t > 0$.  (Take $\beta = \frac{1}{2}$ and $c = \frac{1}{\sqrt{2}}-\varepsilon\sqrt{2}$ in \cite[Theorem 1]{bonda}.)   It follows that
$$(0,\infty, \tfrac{1}{2}, -\tfrac{1}{2},\tfrac{1}{2},0,0,0,\ldots) \leq \dege \zeta(\tfrac{1}{2}+it) \leq (0,\infty, 1, -1,0,0,0,\ldots),$$
where the upper bound is conditional on the Riemann hypothesis.
 In \cite{farm}, Farmer, Gonek, and Hughes  conjectured that 
\begin{align}\label{fghcb}\max_{t \in [0,T]} |\zeta(\tfrac{1}{2}+it)|  = \exp\left( (1+o(1)) \sqrt{\tfrac{1}{2} \log T \log \log T} \right) \ (T \to \infty).
\end{align}
If true, their conjecture would imply that
$$\dege \zeta(\tfrac{1}{2}+it) = \dege \max_{t \in [0,T]} |\zeta(\tfrac{1}{2}+it)|   = (0,\infty, \tfrac{1}{2}, \tfrac{1}{2},0,0,0,\ldots),$$
since the first equality holds by Proposition \ref{supprop}.

\begin{remark}[Power moments of $\zeta(\tfrac{1}{2}+it)$]
Regarding statement (4) of Proposition \ref{lindelof}, it is conjectured that, for every nonnegative integer $k$, there exists  a $c_k > 0$ such that
$$\frac{1}{x} \int_0^x |\zeta(\tfrac{1}{2}+it)|^{2k} \, dt \sim c_k (\log x)^{k^2} \ (x \to \infty).$$
This conjecture is known to hold for $k = 0,1,2$ with $c_0 = c_1 = 1$ and $c_2 = \frac{1}{2\pi^2}$.  It  is also known that
$$\frac{1}{x} \int_0^x |\zeta(\tfrac{1}{2}+it)|^{2k} \, dt  \gg (\log x)^{k^2} \ (x \to \infty)$$
on condition of the Riemann hypothesis.  For references, see \cite{conrey} \cite{heath0}.
\end{remark}

In  1937 \cite{ing0},  Ingham proved Theorem \ref{inghamtheorema} below.   The theorem establishes  a relationship between the zeros of the  Riemann zeta function $\zeta(s)$ and the order of growth of  $\zeta(s)$  on the critical line.  It is now customary, for all  $\sigma \in [0,1]$ and all $T > 0$,  to let $N(\sigma,T)$ denote the number of zeros $\rho$ of $\zeta(s)$ with $\operatorname{Re}\rho \geq \sigma$ and $0< \operatorname{Im}\rho \leq T$.  Note that $N(1,T) = 0$ for all $T > 0$ (i.e., there are no zeros of $\zeta(s)$ on the critical line), and the Riemann hypothesis is equivalent to $N(\sigma,T) = 0$ for all $\sigma > \frac{1}{2}$ and all $T > 0$.

\begin{theorem}[{\cite{ing0}}]\label{inghamtheorema}
 Suppose that $\zeta\left(\frac{1}{2}+it \right) = o(t^c)  \ (t \to \infty)$, where $c > 0$.   Then one has
$$N(\sigma,T) = O(T^{(2+4c)(1-\sigma)}(\log T)^5) \ (T \to \infty)$$
uniformly for $\frac{1}{2} \leq \sigma \leq 1$.
\end{theorem}

From Bourgain's bound $\mu(\tfrac{1}{2})\leq  \tfrac{13}{84}$, we deduce the following.

\begin{corollary}
For all $b> \tfrac{55}{21}$, one has
$$N(\sigma,T) = O(T^{b(1-\sigma)}) \ (T \to \infty)$$
uniformly for $\frac{1}{2} \leq \sigma \leq 1$ (where the $O$ constant depends only on $b$).  
\end{corollary}

In fact, one has the following result, proved  by Ingham in 1940 \cite{ing1}.

\begin{theorem}[{\cite{ing1} \cite[Theorem 19.9(B)]{tit}}]
One has
$$N(\sigma,T) = O(T^{3(1-\sigma)/(2-\sigma)}(\log T)^5) \ (T \to \infty)$$
uniformly for $\frac{1}{2} \leq \sigma \leq 1$.  Consequently, one has $\deg N(\sigma,T) < 1$ for $\sigma > \tfrac{1}{2}$.
\end{theorem}

Thus, since $\deg N(T) = 1$,  ``most'' of the zeros of $\zeta(s)$ are concentrated near the critical line.

 By Theorem \ref{inghamtheorema},  the  Lindel\"of hypothesis implies the following conjecture.

\begin{conjecture}\label{densityconjecture}
For all $b>2$, one has
$$N(\sigma,T) = O(T^{b(1-\sigma)}) \ (T \to \infty)$$
uniformly for $\frac{1}{2} \leq \sigma \leq 1$  (where the $O$ constant depends only on $b$).  
\end{conjecture}

Note that the $O$ bound  in the conjecture is automatic for $\sigma = \frac{1}{2}$ and $\sigma = 1$.

The following 1969 result of Montgomery implies that $\zeta(s)$ has relatively few zeros with real part in $[\tfrac{9}{10},1]$.

\begin{theorem}[{\cite[p.\ 348]{montg}}]
One has 
$$N(\sigma, T) = O(T^{2(1-\sigma)} (\log T)^{11})  \ (T \to \infty)$$
uniformly for $\frac{9}{10} \leq \sigma \leq 1$.
\end{theorem}

Consequently, Conjecture \ref{densityconjecture} is equivalent to the well-known {\bf density hypothesis}, \index{density hypothesis} stated below.

\begin{conjecture}[{Density hypothesis, e.g.,  \cite[(1.137)]{ivic}}]\label{densityhypothesis}
For every $\varepsilon > 0$, one has
$$N(\sigma,T) = O(T^{2(1-\sigma)+\varepsilon}) \ (T \to \infty)$$
uniformly for $\frac{1}{2} \leq \sigma \leq 1$  (where the $O$ constant depends only on $\varepsilon$). Consequently, one has $\deg  N(\sigma,T)  \leq 2(1-\sigma)$ for all $\sigma \in [\tfrac{1}{2},1]$.
\end{conjecture}

Thus, we have the following.

\begin{corollary}
The Lindel\"of hypothesis implies Conjectures \ref{densityhypothesis} (i.e., the density hypothesis) and \ref{densityconjecture}, which are equivalent.
\end{corollary} 

Note that some  related and well-known ``density conjectures,'' e.g., \cite[Density Conjecture, p.\ 249]{iwan}, are stronger than the density hypothesis.

\begin{outstandingproblem}
Determine $\deg N(\sigma,T)$, and, more generally,  $\dege N(\sigma,T)$, as a function of $\sigma$.
\end{outstandingproblem}

Let $N_0(T)$ denote the number of zeros $\rho$ of $\zeta(s)$ on the critical line with $0< \operatorname{Im}\rho \leq T$,  and let $N_0^*(T)$ denote the number of such zeros that are simple.
We also note the following.

\begin{theorem}[{\cite[Theorem 1]{conrey}}]
One has
$$\liminf_{T \to \infty} \frac{N_0(T)}{N(T)} \geq 0.4077$$
and
$$\liminf_{T \to \infty} \frac{N_0^*(T)}{N(T)} \geq 0.401.$$
In particular, more than $\frac{2}{5}$ of the nontrivial zeros of $\zeta(s)$ are simple and lie on the critical line.
\end{theorem}

\section{The nontrivial zeros of $\zeta(s)$}

In this section we assume the notation of Section 5.3.  We write $\rho_n = \sigma_n + i \gamma_n$ for the  $n$th nontrival zero of $\zeta(s)$ with positive imaginary part, where $0 < \gamma_1 \leq \gamma_2 \leq \gamma_3 \leq \cdots$ are the consecutive imaginary parts of the nontrivial zeros of $\zeta(s)$ with positive imaginary part, repeated to multiplicity, and where $\sigma_n \leq \sigma_{n+1}$ if $\gamma_n = \gamma_{n+1}$.    We also write $$\tau_n = \frac{\gamma_n}{2\pi}$$ and $$\widehat{\gamma}_n = \tau_n \log \frac{\tau_n}{e} + \frac{11}{8}$$
for all $n$.

We consider the following problem.

\begin{outstandingproblem}
Compute $\dege S$.
\end{outstandingproblem}

The $O$ and $\Omega$ results for the function $S(T)$ discussed in Section 5.3 yield the following.

\begin{proposition}\label{degeSbounds}
One has
\begin{align*}
(0,\tfrac{1}{3},-\tfrac{1}{3},0,0,0,\ldots) \leq \dege S \leq (0,1,0,0,0,\ldots) 
\end{align*}
and, assuming the Riemann hypothesis, also
\begin{align*}
(0,\tfrac{1}{2},-\tfrac{1}{2}, \tfrac{1}{2},0,0,0,\ldots) \leq \dege S \leq (0,1,-1,0,0,0,\ldots) 
\end{align*}
\end{proposition}

Recall that one has 
$$S(T) = O(\log T) \ (T \to \infty)$$
and
$$\int_0^T S(t) \, dt = O(\log T) \ (T \to \infty).$$
In 1924, Littlewood proved that
$$S(T) = o(\log T) \ (T \to \infty)$$
and
$$\int_0^T S(t) \, dt = o(\log T) \ (T \to \infty)$$
on condition of the  Lindel\"of hypothesis  \cite[Theorems 8 and 9]{litt24}.

\begin{proposition}\label{degeStau}  One has the following.
\begin{enumerate}
\item Let $r$ be Hardian and eventually positive (e.g.,  $r \in \mathbb{L}_{>0}$) with $r(t) = O(\log t) \ (t \to \infty)$.  One has
$$S(t) = O(r(t)) \ (t \to \infty)$$  
if and only if
$$S(\gamma_n) = O(r(\gamma_n)) \ (n \to \infty).$$  Likewise, one has
$$S(t) = o(r(t)) \ (t \to \infty)$$  if and only if
$$S(\gamma_n) = o(r(\gamma_n)) \ (n \to \infty).$$ 
Moreover, if any of the conditions above holds, or if $r$ is eventually increasing, then one has $\dege r(\gamma_n) = \dege r$.
\item $\dege S(\gamma_n) = \dege S$.
\end{enumerate}
\end{proposition}

\begin{proof}
Let $r \in \mathbb{L}_{>0}$ with $r(t) = O(\log t) \ (t \to \infty)$ and
$$S(\gamma_n)  = O(r(\gamma_n)) \ (n\to \infty).$$
 Given $t > 0$, let $n$ be such that
$$\gamma_n < t \leq \gamma_{n+1}.$$
It follows that
$$S(\gamma_{n+1})-\frac{1}{2} = S(\gamma_{n+1}^-) < S(t) < S(\gamma_{n}^+) = S(\gamma_n)+\frac{1}{2}$$
and therefore
$$ |S(t)| < \max(|S(\gamma_{n}^-)|,|S(\gamma_{n+1}^+) |) \leq \max(|S(\gamma_{n})|,|S(\gamma_{n+1}) |)+\frac{1}{2}.$$
Since then $\limsup_{n \to \infty} |S(\gamma_n)| = \infty$, it follows that $r(t)$ is eventually increasing without bound.   Since $ \frac{r(t)}{\log t}$ is Hardian,  one has
$$\lim_{t \to \infty} \frac{r(t)}{\log t} = \limsup_{t \to \infty} \frac{r(t)}{\log t} < \infty,$$ so that, by  L'H\^opital's rule and the fact that $r'$ is Hardian, one also has
$$\lim_{t \to \infty} \frac{r'(t)}{1/t} = \limsup_{t \to \infty} \frac{r'(t)}{1/t} < \infty.$$  Consequently, one has
$$r'(t) = O\left( \frac{1}{t}\right) \ (t \to \infty).$$  If follows from the mean value theorem, then, that
$$\lim_{n \to \infty} \frac{r(\gamma_{n+1})-r(\gamma_n)}{\gamma_{n+1}-\gamma_n} = 0,$$
so that, by (\ref{tau2}), one has
$$r(\gamma_{n+1})-r(\gamma_n) = o(\gamma_{n+1}-\gamma_n) = o\left( \frac{1}{\log\log\log n} \right) \ (n \to \infty)$$
and therefore
$$\lim_{n \to \infty} \frac{r(\gamma_n)}{r(\gamma_{n+1})} = 1.$$
It follows that there exists a $C> 1$  such that
\begin{align*}
\sup_{t \in [\gamma_n, \gamma_{n+1}]} \frac{|S(t)|}{r(t)}  & \leq  \max\left(\frac{\left|S(\gamma_{n})\right|}{r(\gamma_n)},\frac{\left|S(\gamma_{n+1}) \right|}{r(\gamma_n)}\right)+\frac{1}{2 r(\gamma_n)}  \\  & \leq C \max\left(\frac{\left|S(\gamma_{n})\right|}{r(\gamma_n)},\frac{\left|S(\gamma_{n+1}) \right|}{r(\gamma_{n+1})}\right) + o(1) \\
& = O(1) \ (n \to \infty)
\end{align*}
for all $n  \gg 0$, and therefore
$$S(t) = O(r(t)) \ (t \to \infty).$$  
The corresponding equivalence for $o$ follows from the inequalities above.  

Now, since $\gamma_n < n$ for all $n \gg 0$, one has $\dege r(\gamma_n) \leq \dege r$.  Since 
$$\gamma_n \sim  \frac{2\pi n}{\log n} \ (n \to \infty),$$ for any $C < 2 \pi$ one has
$$\gamma_n >  \frac{Cn}{\log n}$$
for all $n \gg 0$.  It follows that
$$\dege r(\gamma_n) \geq \dege r\left( \frac{Cn}{\log n}\right).$$
Moreover, by Proposition \ref{arithbridge2} and Theorem \ref{fgie}, one has
$$\dege r\left( \frac{Cn}{\log n}\right) = \dege r\left( \frac{Cx}{\log x}\right) = \dege r.$$
This proves statement (1).  Statement (2) then follows from statement (1) and  Theorem \ref{infpropexp}.
\end{proof}

The following result, which follows immediately from  Propositions \ref{lnnew}, \ref{sntn}, and \ref{degeStau}, provides an alternative expression for $\dege S$. 

\begin{proposition}\label{alss}
One has
$$\dege \left(n-\widehat{\gamma}_n \right) = \dege S.$$  Moreover, one has
$$n - \widehat{\gamma}_n \sim (r_n-\tau_n)\log n \ (n \to \infty),$$
and therefore
$$\dege (r_n-\tau_n) = \dege S + (0,-1,0,0,0,\ldots),$$
where 
$$r_n = \frac{n-11/8}{W((n-11/8)/e)}.$$
\end{proposition}

We also have the following immediate consequence of Theorem \ref{tauprop} and Propositions \ref{degeSbounds} and \ref{degeStau}.

\begin{theorem}\label{rndege}
Let $a \in \RR$, and let 
$$r_n = r_n(a) =  \frac{n-a}{W((n-a)/e)}.$$
Then 
\begin{align*}
\dege(r_n-\tau_n) = \dege\left( \frac{S(\gamma_n)+\frac{11}{8}-a}{\log r_n}  \right)   = \dege S + (0,-1,0,0,0,\ldots)
\end{align*}
and
\begin{align*}
\dege\left(r_n-\tau_n - \frac{S(\gamma_n)+\frac{11}{8}-a}{\log r_n} \right) \leq (-1,0,0,0,\ldots).
\end{align*} 
Moreover, one has
\begin{align*}
(0,-\tfrac{2}{3},-\tfrac{1}{3},0,0,0,\ldots) \leq \dege(r_n-\tau_n) \leq (0,0,0,0,0,\ldots) 
\end{align*}
and, assuming the Riemann hypothesis, also
\begin{align*}
(0,-\tfrac{1}{2},-\tfrac{1}{2}, \tfrac{1}{2},0,0,0,\ldots) \leq \dege(r_n-\tau_n) \leq (0,0,-1,0,0,0,\ldots).
\end{align*}
\end{theorem}

Proposition \ref{degeStau} (applied to $r(t) = \log t$)  and  the proof of Theorem \ref{tauprop} yield the following.

\begin{proposition}
Let $a \in \RR$, and let 
$$r_n = r_n(a) =  \frac{n-a}{W((n-a)/e)}.$$
The conjecture 
$$S(T) = o(\log T) \ (T \to \infty),$$
which is implied by the Lindel\"of hypothesis, is equivalent to 
$$S(\gamma_n) = o(\log \gamma_n) \ (n \to \infty)$$
and to 
$$\lim_{n \to \infty} \left(r_n -\tau_n \right) = 0$$
for some (or all) $a \in \RR$.  Moreover, the conjecture implies that
\begin{align*}
r_n -\tau_n = \frac{S(\gamma_n)+\frac{11}{8}-a }{\log r_n} + \frac{1+o(1)}{96\pi^2 n }  \ (n \to \infty)
\end{align*}
and therefore that
\begin{align*}
\dege \left(r_n -\tau_n - \frac{S(\gamma_n)+\frac{11}{8}-a }{\log r_n} \right) = (-1,0,0,0,\ldots).
\end{align*}
\end{proposition}

\begin{proof}
The first claim of the proposition follows from Proposition \ref{degeStau}, applied to $r(t) = \log t$, and Theorem \ref{tauprop}.  The second claim follows by observing that, if the conjecture holds, then, in the proof of Theorem \ref{tauprop}, the second $O$ in (\ref{lindelo}) can be replaced with $o$:
\begin{align*}
\tau_n-r_n = \left(\frac{a-\frac{11}{8} -S(\gamma_n) -\frac{1}{96 \pi^2 r_n}+ O \left( \frac{1}{r_n^3}\right)}{\log r_n}\right)\left(1+ o \left( \frac{1}{r_n \log r_n}\right) \right) \ (n \to \infty).
\end{align*}
This completes the proof.
\end{proof}

The {\bf Gram spacing}\index{Gram spacing}  $\frac{1}{\log \tau_n}$ is an approximation for the exact average spacing
$$\frac{1}{n}\sum_{k  = 1}^n (\tau_{k+1}-\tau_k) =  \frac{\tau_{n+2}-\tau_1}{n} \approx \frac{1}{\log \tau_n}.$$

\begin{proposition}\label{zzgap}
Let $a \in \RR$.  One has
\begin{align*}
\dege\left(\log \tau_n-1+W((n-a)/e)\right) & = \dege S + (-1,0,0,0,\ldots) \\ & \leq (-1,1,0,0,0,\ldots)
\end{align*}
and
\begin{align*}
\dege\left( \frac{1}{1+W((n-a)/e)} -\frac{1}{\log \tau_n} \right) & = \dege\left( \frac{1}{W((n-a)/e)} -\frac{1}{\log \tau_n-1} \right) \\
& = \dege S + (-1,-2,0,0,0,\ldots) \\ & \leq (-1,-1,0,0,0,\ldots).
\end{align*}
\end{proposition}

\begin{proof}
Since $\frac{x}{W(x)} = e^{W(x)}$, one has
$$r_n :=  \frac{n-a}{W((n-a)/e)} = e^{1+W((n-a)/e)}$$
and therefore
$$\log r_n = 1+W((n-a)/e).$$
Since $ \log x \sim x-1\ (x \to 1)$,  one has
$$\log \tau_n -\log r_n \sim \frac{\tau_n}{r_n}-1  \ (n \to \infty)$$
and therefore
\begin{align*}
 \frac{1}{1+W((n-a)/e)} -\frac{1}{\log \tau_n} &  = \frac{1}{\log r_n} -\frac{1}{\log \tau_n} \\
 &  \sim \frac{1}{(\log r_n)(\log \tau_n)}\left( \frac{\tau_n}{r_n}-1\right) \\
 &  \sim \frac{1}{n \log n}( \tau_n-r_n) \ (n \to \infty),
\end{align*}
whence by Theorem \ref{rndege} one has
\begin{align*}
\dege\left( \frac{1}{1+W((n-a)/e)} -\frac{1}{\log \tau_n} \right) & = \dege (\tau_n-r_n)+(-1,-1,0,0,0,\ldots) \\
& = \dege S+ (-1,-2,0,0,0,\ldots).
\end{align*}
The rest of proposition follows in similar fashion.
\end{proof}

Now we address the Riemann zeta zero gaps, and, specifically, the following problem.

\begin{outstandingproblem}
Compute $ \dege (\gamma_{n+1}-\gamma_n) = \dege (\tau_{n+1}-\tau_n)$.
\end{outstandingproblem}

By (\ref{tau1}), (\ref{tau2}), and (\ref{tau3}), one has the following.

\begin{proposition}\label{zzgap2}
One has
$$ (0,-1,0,0,0,\ldots) \leq \dege (\tau_{n+1}-\tau_n) \leq (0,0,0,-1,0,0,0,\ldots),$$
while also
$$\dege (\tau_{n+1}-\tau_n) \leq (0,0,-1,0,0,0,\ldots)$$
on condition of the Riemann hypothesis.  
\end{proposition}

Regarding the gaps $\widehat{\gamma}_{n+1}-\widehat{\gamma}_n$ between the normalized zeros $$\widehat{\gamma}_n = \tau_n \log \frac{\tau_n}{e}+\frac{11}{8},$$ in addition to the normalized spacings $$\delta_n = (\tau_{n+1}-\tau_n)\log \tau_n,$$ we note the following result, which follows immediately from Proposition \ref{deltal}.

\begin{proposition}\label{tnln}
One has
$$\widehat{\gamma}_{n+1}-\widehat{\gamma}_n \sim \delta_n \sim (\tau_{n+1}-\tau_n) \log n \ (n \to \infty),$$ and therefore
$$\dege(\widehat{\gamma}_{n+1}-\widehat{\gamma}_n) = \dege \delta_n = \dege (\tau_{n+1}-\tau_n) + (0,1,0,0,0,\ldots).$$
\end{proposition}

The following proposition concerns the Gram spacing $\frac{1}{\log \tau_n}$ and the exact average spacing 
$$\frac{1}{n}\sum_{k  = 1}^n (\tau_{k+1}-\tau_k) =  \frac{\tau_{n+2}-\tau_1}{n}.$$

\begin{proposition}
Let $a \in \RR$. 
\begin{enumerate}
\item One has 
\begin{align*}
 \dege\left (\frac{\tau_{n}}{n}- \frac{1}{\log \tau_n-1} \right) & =  \dege\left (\frac{\tau_{n+2}}{n}- \frac{1}{\log \tau_n-1} \right)  \\
  &  = \dege S + (-1,-1,0,0,0,\ldots) \\
   & \leq (-1,0,0,0,\ldots).
\end{align*}
and
$$\dege\left(\frac{\tau_{n+2}-\tau_1}{n}- \frac{1}{\log \tau_n-1}\right) \leq (-1,0,0,0,\ldots).$$
\item One has
$$\frac{\tau_{n+2}-\tau_1}{n}- \frac{1}{\log \tau_n} \sim  \frac{1}{(\log \tau_n)^2} \ (n \to \infty)$$
and therefore
$$\dege\left(\frac{\tau_{n+2}-\tau_1}{n}- \frac{1}{\log \tau_n}\right) = (0,-2,0,0,0,\ldots).$$
\item If the  Lindel\"of hypothesis holds, or, more generally, if
$$S(T) = o(\log T) \ (T \to \infty),$$
then one has
$$\frac{\tau_{n+2}-\tau_1}{n}- \frac{1}{\log \tau_n-1} \sim -\frac{\tau_1}{n} \ (n \to \infty)$$
and therefore
$$\dege\left(\frac{\tau_{n+2}-\tau_1}{n}- \frac{1}{\log \tau_n-1}\right) = (-1,0,0,0,\ldots).$$
\item If the Riemann hypothesis holds, or, more generally, if 
$$\dege S < (0,1,0,0,0,\ldots),$$
then one has
\begin{align*}
 \dege\left (\frac{\tau_{n}}{n}- \frac{1}{\log \tau_n-1} \right) & =  \dege\left (\frac{\tau_{n+2}}{n}- \frac{1}{\log \tau_n-1} \right) < (-1,0,0,0,\ldots).
\end{align*}
\end{enumerate}
\end{proposition}

\begin{proof}
Let $$r_n =  \frac{n}{W(n/e)}.$$
By the proof of Proposition \ref{zzgap}, one has
\begin{align*}
 \frac{1}{W(n/e)} -\frac{1}{\log \tau_n-1} \sim \frac{1}{n \log n}( \tau_n-r_n) \ (n \to \infty)
\end{align*}
and therefore
\begin{align*}
 r_n-\frac{n}{\log \tau_n-1} \sim \frac{1}{\log n}( \tau_n-r_n) \ (n \to \infty),
\end{align*}
whence also
\begin{align*}
 \tau_n-\frac{n}{\log \tau_n-1} & = (1+o(1) )\frac{1}{\log n}( \tau_n-r_n) +(\tau_n-r_n)   \ (n \to \infty) \\
 & \sim \tau_n-r_n \ (n \to \infty).
\end{align*}
Moreover, by Theorem \ref{tauprop}, one has
\begin{align*}
r_n-\tau_n =\frac{ S(\gamma_n)+\frac{11}{8}}{\log  r_n} + O \left( \frac{1}{n}\right)  \ (n \to \infty).
\end{align*}
Thus, one has
\begin{align*}
 \tau_n-\frac{n}{\log \tau_n-1} = -(1+o(1))\frac{ S(\gamma_n)+\frac{11}{8}}{\log  r_n} + O \left( \frac{1}{n}\right)  \ (n \to \infty),
\end{align*}
and therefore
\begin{align*}
 \tau_{n+2}-\frac{n+2}{\log \tau_{n+2}-1} & = -(1+o(1))\frac{ S(\gamma_{n+2})+\frac{11}{8}}{\log  r_{n+2}} + O \left( \frac{1}{n}\right)   \ (n \to \infty),
\end{align*}
whence
\begin{align}\label{ttthelp}
 \tau_{n+2}-\frac{n}{\log \tau_{n}-1} & = -(1+o(1))\frac{ S(\gamma_{n+2})+\frac{11}{8}}{\log  r_{n+2}} + O \left( \frac{1}{\log n}\right)   \ (n \to \infty).
\end{align}
It follows that
\begin{align*}
 \dege\left (\frac{\tau_{n}}{n}- \frac{1}{\log \tau_n-1} \right)  =   \dege \frac{ S(\gamma_{n})+\frac{11}{8}}{n\log r_n}   =   \dege S + (-1,-1,0,0,0,\ldots)
\end{align*}
and
\begin{align*}
 \dege\left (\frac{\tau_{n+2}}{n}- \frac{1}{\log \tau_n-1} \right)  =   \dege \frac{ S(\gamma_{n+2})+\frac{11}{8}}{n\log  r_{n+2}}   =   \dege S + (-1,-1,0,0,0,\ldots).
\end{align*}
This proves statement (1), and then statements (2)--(4) readily follow from statement (1) and (\ref{ttthelp}).
\end{proof}

A consequence of the proposition is that  $\frac{1}{\log \tau_n-1}$ is a more accurate estimate (in terms of logexponential degree)  for
 the exact average spacing $\frac{\tau_{n+2}-\tau_1}{n}$
than the Gram spacing $\frac{1}{\log \tau_n}$.  This is analogous to how $\frac{x}{\log x -1}$ is a more accurate estimate for $\pi(x)$ than the estimate $\frac{x}{\log x}$.  Similarly, $\frac{1}{\log \tau_n-1} - \frac{\tau_1}{n}$  is a more accurate estimate  for the exact average spacing $\frac{\tau_{n+2}-\tau_1}{n}$ than the estimate $\frac{1}{\log \tau_n-1}$, provided that 
$\dege S < (0,1,0,0,0,\ldots)$.

Assuming Montgomery's conjecture \cite[(17)]{mont3} (see also Conjecture \ref{normsp}) that
$$\limsup_{n \to \infty}\, (\tau_{n+1}-\tau_n) \log \tau_n = \infty,$$
that is, that
$$\tau_{n+1}-\tau_n \neq O\left(\frac{1}{\log n}\right) \ (n \to \infty),$$
one has
$$\tau_{n+1}-\tau_n -\frac{1}{\log \tau_n} \neq O\left(\frac{1}{\log n}\right) \  (n \to \infty),$$
and therefore
$$\dege\left(\tau_{n+1}-\tau_n -\frac{1}{\log \tau_n} \right) \geq (0,-1,0,0,0,\ldots),$$
or, equivalently, 
$$\dege\left(\delta_n- 1\right) \geq (0,0,0,0,0,\ldots),$$
where $\delta_n = (\tau_{n+1}-\tau_n)\log \tau_n$.  By Proposition \ref{zzgap}, then, Montgomery's conjecture above implies that
$$\dege\left(\tau_{n+1}-\tau_n -\frac{1}{\log \tau_n} \right)= \dege\left( \tau_{n+1}-\tau_n- \frac{1}{1+W((n-a)/e)} \right)$$
for all $a \in \RR$.  Assuming the stronger conjecture that
$$\dege (\tau_{n+1}-\tau_n) > (0,-1,0,0,0,\ldots),$$ 
one has an affirmative answer to the question posed in the following problem.

\begin{problem}\label{mainx}
Are all of
 $$\dege(\tau_{n+1}-\tau_n ),$$
$$\dege\left(\tau_{n+1}-\tau_n -\frac{1}{\log \tau_n} \right),$$ 
 $$\dege\left(\tau_{n+1}-\tau_n -\frac{1}{\log \tau_n-1} \right),$$
and 
 $$\dege\left(\tau_{n+1}-\tau_n -\frac{\tau_{n+2}-\tau_1}{n} \right)$$  equal to one another? 
\end{problem}

Various conjectures regarding the $\gamma_n$, $\tau_n$, and $\widehat{\gamma}_n$ are discussed in Section 14.4.

\chapter{Primes in intervals, the $n$th prime, and  the $n$th prime gap}

In this chapter, we use the degree formalism to study further functions related to the prime numbers, including the function $\pi(f(x))-\pi(g(x))$ for functions $f(x)$ and $g(x)$ satisfying various conditions,  the prime listing function $p_n$, the prime gap function $g_n = p_{n+1}-p_n$,  and the maximal prime gap function $G(x) = \max_{p_k \leq x} g_k$.

\section{Asymptotics for prime counts in intervals}

In this section we seek asymptotics for prime counts in intervals.

{\bf Legendre's conjecture}\index{Legendre's conjecture} states that there is a prime number between $n^2$ and $(n + 1)^2$, or, equivalently, that $\pi((n+1)^2)-\pi(n^2) \geq 1$, for every positive integer $n$.  
{\bf Oppermann's conjecture}\index{Opperman's conjecture} states that  there is at least one prime number between $n(n-1)$ and $n^2$ and at least another prime between $n^2$ and $n(n+1)$, or, equivalently,
that $\pi(n^2)-\pi(n(n-1)) \geq 1$ and $\pi(n(n+1))-\pi(n^2) \geq 1$, for every integer $n > 1$.   Opperman's conjecture implies that there are at least two primes between $n^2$ and $(n+1)^2$ for all positive integers $n$, namely, one between $n^2$ and $n(n+1)$ and another between $n(n+1)$ and $(n + 1)^2$.  Thus, Opperman's conjecture is stronger than Legendre's conjecture.   It is known (see Corollary \ref{lindel}) that the Lindel\"of hypothesis ``almost'' yields an ``eventual'' version of Opperman's conjecture in that it implies that
\begin{align}\label{tpi}
\pi(x +ax^{t}) - \pi(x) \sim \frac{ax^{t}}{\log x} \ (x \to \infty),
\end{align}
or, equivalently, that
\begin{align*}
\pi(x^2+a x^{2t}) - \pi(x^2) \sim \frac{ax^{2t}}{2\log x} \ (x \to \infty),
\end{align*}
for all $a \neq 0$ and all $t \in (\frac{1}{2},1)$.   If  (\ref{tpi}) were to hold for $t = \frac{1}{2}$ and $a = \pm 1$ (in fact, $a = 1$ suffices), then one could conclude that
$$\lim_{x \to \infty}|\pi(x^2 \pm x) - \pi(x^2)|  = \infty,$$
i.e., that for all $M > 0$ one has $|\pi(x^2 \pm x) - \pi(x^2)| \geq M$ for all sufficiently large $x$.   
 It is known  that (\ref{tpi}) holds unconditionally for all $a \neq 0$ and all $t \in [ \frac{7}{12},1)$: see \cite{hux} and \cite[Theorem]{heathbrown}.   From this it follows by a simple substitution that
$$\pi(x^d +a x) - \pi(x^d) \sim \frac{ax}{d\log x} \ (x \to \infty)$$
for all $d \in(1, 1+\frac{5}{7}]$, and likewise that
$$\pi((x+1)^d)-\pi(x^d)) \geq \pi(x^d + dx^{d-1}) - \pi(x^d) \sim  \frac{x^{d-1}}{\log x} \ (x \to \infty)$$
for all $d \geq 2+\frac{2}{5}$.  The latter implies that, for all $d \geq 2+\frac{2}{5}$ and for all sufficiently large $x$, there exists a prime between $x^d$ and $(x+1)^d$, and in fact
$$\lim_{x \to \infty} (\pi((x+1)^d)-\pi(x^d))=  \lim_{x \to \infty}  (\pi(x^d + x^{d-1}) - \pi(x^d)) = \infty.$$

In this section, we explore the following problem, which is motivated by the discussion above.

\begin{outstandingproblem}[{\cite{maier}}]\label{maierprob}
From Maier \cite{maier}: ``The following question is suggested by the prime number theorem: for which functions $\Phi$ is it true that 
$$\pi(x+\Phi(x))-\pi(x) \sim \frac{\Phi(x)}{\log x} \ (x \to \infty)?"$$
\end{outstandingproblem}

Further related problems include the following.

\begin{outstandingproblem}\label{Gmaxprob}
For which functions $h$ defined in a neighborhood of $\infty$ is there a prime greater than $x$ and less than or equal to $x+h(x)$ for all $x \gg 0$?  Equivalently,  for which functions $h$ defined in a neighborhood of $\infty$ does one have $\pi(x+h(x)) -\pi(x) \geq 1$ for all $x  \gg 0$?
\end{outstandingproblem}

\begin{outstandingproblem} Let $h$ be a particular function defined on $(N,\infty)$ for some $N> 0$.
\begin{enumerate}
\item Is there a prime between $x$ and $x+h(x)$ for all $x > N$?
\item Is there a prime between $x$ and $x+h(x)$ for all $x \gg 0$?  More generally, what are $\liminf_{x \to \infty} |\pi(x+h(x))-\pi(x)|$ and $\limsup_{x \to \infty} |\pi(x+h(x))-\pi(x)|$?
\item Is  it true that $\pi(x+h(x))-\pi(x) \sim \frac{h(x)}{\log x} \ (x \to \infty)$?
\item What is $\dege (\pi(x+h(x))-\pi(x))$? 
\end{enumerate}
\end{outstandingproblem}

The following is a useful result of Montgomery and Vaughan from 1973.

\begin{theorem}[{\cite[(1.12)]{mv}}]\label{mv}
For all $x \geq 0$ and  all $y > 1$, one has
$$0 \leq \pi(x+y) -\pi(x) < \frac{2 y}{\log y }.$$
\end{theorem}

\begin{corollary}
For all $ x> 1$ and all $t >0$, one has
$$\pi(x+x^t)-\pi(x) <  \frac{2x^t}{t\log x}.$$
\end{corollary}

Let us extend the domain of $\pi(x)$ to all of $\RR$ by setting $\pi(x) = 0$ for all $x < 0$.  Then we also have the following.

\begin{corollary}\label{mvcor}
For all $x,y\in \RR$, one has
$$| \pi(x+y) -\pi(x)|  \leq 2\li(|y|+\mu)+1,$$
where $\mu$ is the Ramanujan--Soldner constant.
\end{corollary}

\begin{lemma}\label{hloglemma}
Let $h$ be a real function defined on an unbounded subset of $\RR_{>0}$ with $\log (x+h(x)) \sim \log x \ (x \to \infty)$ (which holds if $h(x) \ll x \log x \ (x \to \infty)$ and $h(x)$ is eventually positive) and $h(x) \gg \frac{x}{(\log x)^a}\ (x \to \infty)$  for some $a > 0$.  Then one has
$$\pi(x+h(x))-\pi(x) \sim \li(x+h(x))-\li(x) \sim  \frac{h(x)}{\log x} \ (x \to \infty).$$
\end{lemma}

 \begin{proof}
By the prime number theorem with error bound, one has
\begin{align*}
\pi(x+h(x)) -\pi(x) -( \li(x+h(x))-\li(x) ) & \leq  | \pi(x+h(x)) - \li(x+h(x))|+ |\pi(x) -\li(x) | \\
& \ll \left|\frac{x+h(x)}{(\log(x+h(x)))^{n+a}} \right| + \left|\frac{x}{(\log x)^{n+a}}\right| \\
& \ll \frac{h(x)}{(\log x)^n}\ (x \to \infty)
\end{align*}
 for all $n > 0$.  Moreover, since $\li'(x) = \frac{1}{\log x}$ is decreasing  and $\log (x+h(x)) \sim \log x \ (x \to \infty)$,  by the mean value theorem one has 
 $$\li(x+h(x))-\li(x) \sim \frac{h(x)}{\log x} \ (x \to \infty).$$
 The lemma follows.
 \end{proof}

Our first main result concerning Problem \ref{maierprob} is as follows.

 \begin{theorem}\label{IVL}
 Let $t  \in (0,1)$, and suppose that 
 $$\pi(x+x^t) - \pi(x) \sim \frac{x^t}{\log x}  \ (x \to \infty).$$
 Let $h$ be an eventually positive real function whose domain is an unbounded subset of $\RR_{> 0}$,  and suppose that $\log (x+h(x)) \sim \log x \ (x \to \infty)$ and 
 $x^t|_{\dom h} = o(h(x)) \ (x \to \infty)$.   Then (on $\dom h$) one has
 $$\pi(x+h(x)) - \pi(x) \sim \frac{h(x)}{\log x}  \ (x \to \infty).$$
 \end{theorem}

 \begin{proof} 
  By  Lemma \ref{hloglemma}, we may suppose without loss of generality that $h(x) = o(x) \ (x \to \infty)$.
Let  $$T(x) = x+x^t$$
for all $x \geq 0$.  Also,  for all $x \in \dom h$,  let
 $$K(x) = \left \lfloor \frac{h(x)}{x^t} \right \rfloor.$$
 For all  $x \in \dom h$, one has
 $$0 \leq h(x)-K(x)x^t < x^t = o(h(x)) \ (x \to \infty),$$
 and therefore
 $$K(x)x^t \sim h(x) \ (x \to \infty).$$
 Note that, for all $x \geq 0$, one has
 $$0 \leq x^t = T(x)-x \leq T(T(x))-T(x) \leq T(T(T(x)))-T(T(x)) \leq \cdots,$$
so that
  $$0 \leq x+ kx^t \leq T^{\circ k}(x)$$
  for all positive integers $k$,
  whence
  $$0 \leq x+K(x)x^t  \leq T^{\circ K(x)}(x),$$
  for all $x \in \dom h$.   For all $x \in \dom h$, we may write
\begin{align}\label{pp1}
\pi(x+h(x))-\pi(x) = \sum_{1 \leq k \leq K(x)} (\pi(T^{\circ k}(x))-\pi(T^{\circ (k-1)}(x)) + \pi(x+h(x))-\pi(T^{\circ K(x)}(x)).
\end{align}
By our hypothesis on $\pi(x)$, for all fixed positive integers $k$ one has
 $$\pi(T^{\circ k}(x))-\pi(T^{\circ (k-1)}(x)) \sim \frac{T^{\circ (k-1)}(x)^t}{\log T^{\circ (k-1)}(x)} \sim \frac{x^t}{\log x} \ (x \to \infty),$$
which suggests that
\begin{align}\label{pp2}
 \sum_{1 \leq k \leq K(x)} (\pi(T^{\circ k}(x))-\pi(T^{\circ (k-1)}(x))  \sim \frac{K(x)x^t}{\log x} \sim \frac{h(x)}{\log x} \ (x \to \infty).
 \end{align}
 To verify (\ref{pp2}),  it suffices to show that
 $$T^{\circ (K(x)-1)}(x) \sim x \ (x \to \infty),$$
 or, equivalently, that 
  $$T^{\circ K(x)}(x) \sim x \ (x \to \infty).$$
To this end,  note that
\begin{align*}
1  &\leq \frac{T^{\circ K(x)}(x) }{x} & \\ 
  & = \frac{T(x)}{x} \cdot  \frac{T(T(x))}{T(x)}  \cdots \frac{T^{\circ K(x)}(x) }{T^{\circ (K(x)-1)}(x) }    \\  
  & =  \left( 1+ \frac{1}{x^{1-t}}  \right)  \left( 1+ \frac{1}{T(x)^{1-t}}  \right) \cdots \left( 1 + \frac{1}{T^{\circ (K(x)-1)}(x)^{1-t}}  \right)  \\
  & \leq \left( 1+ \frac{1}{x^{1-t}}  \right)^{K(x)} \\
  & \leq \left( 1+ \frac{1}{x^{1-t}}  \right)^{h(x)/x^t} \\
  & =  \left(\left( 1+ \frac{1}{x^{1-t}}  \right)^{x^{1-t}}\right)^{h(x)/x} \\
  & \leq e^{h(x)/x}  \\
  &\sim 1 \ (x \to \infty)
  \end{align*}
  for all $x \gg 0$ in $\dom h = \dom K$.
This proves (\ref{pp2}).  By (\ref{pp1}) and (\ref{pp2}), then,  it remains only to show that
 \begin{align}\label{NHT1}
  \pi(x+h(x))-\pi(T^{\circ K(x)}(x)) = o\left( \frac{h(x)}{\log x} \right) \ (x \to \infty).
  \end{align}

To prove (\ref{NHT1}), we show that the difference $x+h(x)-T^{\circ K(x)}(x)$ is bounded by $o(h(x))$ and then apply Corollary \ref{mvcor}. Note  first that
$$x+h(x)-T^{\circ K(x)}(x) \leq h(x)-K(x)x^t < x^t = o(h(x))$$
on $\dom h$.  We wish to bound the difference by $o(h(x))$ also from below.
By the binomial series expansion of $(x+x^t)^t$, one has
\begin{align}\label{NHT0}
T(T(x))-T(x)-(T(x)-x) = (x+x^t)^t-x^t \sim tx^{2t-1}  \ (x \to \infty)
\end{align}
unless $t = \frac{1}{2}$, for which one has
$$T(T(x))-T(x)-(T(x)-x) = (x+x^{1/2})^{1/2}-x^{1/2} \sim \tfrac{1}{2}  \ (x \to \infty).$$
Let us first assume $t < \frac{1}{2}$.   From  (\ref{NHT0}), there exists an $N> 0$ such that
\begin{align*}
T(T(x))-T(x)-(T(x)-x)  \leq  2tx^{2t-1}
\end{align*}
for all $x \geq N$.
Therefore, since $x \leq T(x) \leq T(T(x)) \leq \cdots$ for all $x$,  for all $k\geq 2$ we have
\begin{align*}
T^{\circ k}(x)-T^{\circ (k-1)}(x)-(T^{\circ (k-1)}(x)-T^{\circ (k-2)}(x)) \leq  2t(T^{\circ (k-2)}(x))^{2t-1} \leq  2tx^{2t-1}
\end{align*}
for all $k \geq 2$ and all $x$ such that $T^{\circ (k-2)}(x) \geq N$,  hence for all $x \geq N$.   Summing over $k$,  for all $x \geq N$ we then have
$$T^{\circ k}(x)-T^{\circ (k-1)}(x) \leq x^t+ 2(k-1)tx^{2t-1}$$
for all $k \geq 1$, whence, summing over $k$ from $1$ to $K(x)$, we find that
 $$T^{\circ K(x)}(x) \leq  x + K(x) x^t+ t K(x)^2x^{2t-1}$$
 for all $x \geq N$.
  It follows that
\begin{align*}
x+h(x)-T^{\circ K(x)}(x) & \geq  h(x)-K(x) x^t- t K(x)^2x^{2t-1} \\ 
& \geq - t K(x)^2x^{2t-1} \\
& = - t (K(x)x^t)^2x^{-1} \\
& \sim  - t\frac{h(x)}{x}   h(x) \\
& = o(h(x)) 
\end{align*}
as $x \to \infty$.
Therefore, one has
\begin{align}\label{NHT2}
x+h(x)-T^{\circ K(x)}(x)= o(h(x))
\end{align}
if $t <\frac{1}{2}$.    Suppose, on the other hand, that $t = \frac{1}{2}$.  Then there exists an $N$ such that
\begin{align*}
T(T(x))-T(x)-(T(x)-x)  < 1
\end{align*}
for all $x \geq N$.  It follows, for all $x \geq N$, that
$$T^{\circ k}(x)-T^{\circ (k-1)}(x) \leq x^t+ (k-1)$$
for all $k$, whence
 $$T^{\circ K(x)}(x) \leq x + K(x) x^t+ \tfrac{1}{2}  K(x)^2,$$
and therefore
\begin{align*}
x+h(x)-T^{\circ K(x)}(x)  & \geq h(x)-K(x) x^t- \tfrac{1}{2}  K(x)^2 \\ 
& = - \tfrac{1}{2} (K(x)x^t)^2x ^{-1}  \\
& = o(h(x)) 
\end{align*}
as $x \to \infty$.
This proves that (\ref{NHT2}) holds for all $t \leq \frac{1}{2}$.

To handle the case $t > \frac{1}{2}$,  the arguments above for $t \leq \frac{1}{2}$ yield
an $N > 0$ such that, for all $x \geq N$,  one has
\begin{align*}
T(T(x))-T(x)-(T(x)-x) = (x+x^t)^t-x^t \leq 2 tx^{2t-1},
\end{align*}
hence
\begin{align*}
T^{\circ k}(x)-T^{\circ (k-1)}(x)-(T^{\circ (k-1)}(x)-T^{\circ (k-2)}(x)) \leq  2t(T^{\circ (k-2)}(x))^{2t-1} \leq  2t(T^{\circ k}(x))^{2t-1}
\end{align*}
for all $k \geq 2$
hence
$$T^{\circ k}(x)-T^{\circ (k-1)}(x) \leq x^t+ 2t(k-1)(T^{\circ k}(x))^{2t-1},$$
and therefore
 $$T^{\circ K(x)}(x) \leq  x + K(x) x^t+ t K(x)^2(T^{\circ K(x)}(x))^{2t-1}.$$
  It follows that
\begin{align*}
x+h(x)-T^{\circ K(x)}(x) & \geq  h(x)-K(x) x^t- t K(x)^2(T^{\circ K(x)}(x))^{2t-1} \\
& \geq - t K(x)^2(T^{\circ K(x)}(x))^{2t-1}  \\
& = - t (K(x)x^t)^2x^{-1} \left( \frac{T^{\circ K(x)}(x)}{x} \right)^{2t-1}\\
& \sim  - t\frac{h(x)}{x}   h(x)1^{2t-1} \\
& = o(h(x)) 
\end{align*}
as $x \to \infty$.  This establishes (\ref{NHT2}) for all $t \in (0,1)$.

Finally,  we apply Corollary  \ref{mvcor} to obtain
 \begin{align*}
 |\pi(x+h(x))-\pi(T^{\circ K(x)}(x))|  \leq 2 \li ( |x+h(x)-T^{\circ K(x)}(x)| + \mu)+1
 \end{align*}
 for all $x \gg 0$ in  $\dom h$.   By (\ref{NHT2}), given any $\varepsilon > 0$,  for all $x \gg 0$ in $\dom h$ one has
 $$|x+h(x)-T^{\circ K(x)}(x)|  < \varepsilon h(x)$$
 and  $x^t < \varepsilon h(x)$, whence
  \begin{align*}
 |\pi(x+h(x))-\pi(T^{\circ K(x)}(x))| & < 2\li ( \varepsilon h(x) + \mu) +1 \\
& = (1+o(1)) \frac{2(\varepsilon h(x) + \mu)}{\log( \varepsilon h(x) + \mu)} \\
&  < \frac{(2+\varepsilon)\varepsilon h(x)}{t \log x}.
 \end{align*}
  This establishes    (\ref{NHT1}) and therefore completes the proof.
 \end{proof}

 \begin{corollary}
The set of all $t \in \RR$ such that
  $$\pi(x+x^t) - \pi(x) \sim \frac{x^t}{\log x}  \ (x \to \infty)$$
is a subinterval of $(0,1]$ containing $1$.
 \end{corollary}

Theorem \ref{IVL} and its corollary above motivate the following problem.

\begin{outstandingproblem}\label{Tproblem}
Compute the infimum of all $t \in \RR$ such that
 \begin{align}\label{tsim}
 \pi(x+x^t) - \pi(x) \sim \frac{x^t}{\log x}  \ (x \to \infty),
 \end{align}
  which,  by Theorem \ref{IVL}, is equal to the infimum of all $t\in \RR$ such that
$$\pi(x+h(x))-\pi(x) \sim \frac{h(x)}{\log x} \ (x \to \infty)$$
for any eventually positive real function $h$ whose domain is an unbounded subset of $\RR_{>0}$ and  for which $h(x) = O(x) \ (x \to \infty)$ and $x^t|_{ \dom h} = o(h(x)) \ (x \to \infty)$ (or, instead, for which $h(x) = O(x) \ (x \to \infty)$ and  $\underline{\deg}\, h >  t$).
\end{outstandingproblem}

Hoheisel,  in 1930, was the first to prove that (\ref{tsim}) holds for some $t < 1$, namely, for $t   = 1- \frac{1}{33000}$ \cite{hoh}.   Then, in 1936 \cite{chud} and 1937 \cite{ing0}, respectively, Chudakov and Ingham proved independent results that,  when combined  with  Hoheisel's results, established Theorem \ref{inghamtheorem} below.   Recall from Section 11.1 that, for all  $\sigma \in [0,1]$ and all $T > 0$,  we let $N(\sigma,T)$ denote the number of zeros $\rho$ of $\zeta(s)$ with $\operatorname{Re}\rho \geq \sigma$ and $0< \operatorname{Im}\rho \leq T$.  

\begin{theorem}[{\cite{hoh} \cite{chud} \cite{ing0}}]\label{inghamtheorem}
Suppose that $b > 0$ and $B \geq 0$ are constants such that
$$N(\sigma,T) = O(T^{b(1-\sigma)}(\log T)^B) \ (T \to \infty)$$
uniformly for $\frac{1}{2} \leq \sigma \leq 1$.
Let $\kappa = 1-  \frac{1}{b}$.  Then $b \geq 2$ and $\frac{1}{2} \leq \kappa < 1$.
\begin{enumerate}
\item For any $\theta \in (\kappa, 1]$, one has
\begin{align*}
\pi(x+x^t)-\pi(x) \sim \frac{x^t}{\log x} \sim \frac{ \psi(x+x^t)-\psi (x) }{\log x} \ (x \to \infty).
\end{align*}
\item  Suppose that $\zeta\left(\frac{1}{2}+it \right) = o(t^c)  \ (t \to \infty)$, where $c > 0$.   Then the hypotheses above are satisfied for $b = 2+4c$ and $B = 5$, and therefore for  $\kappa =  \frac{1+4c}{2+4c}$.
\end{enumerate}
\end{theorem}

By Theorems \ref{inghamtheorem} and  \ref{IVL}, we have the following corollary.

\begin{corollary}
Assume the hypothesis of Theorem \ref{inghamtheorem}(1) on $\kappa$.    Let $h$ be a real function defined on an unbounded subset of $\RR_{> 0}$, and suppose that $h$ is eventually positive,  $\log (x+h(x)) \sim \log x \ (x \to \infty)$, and $\underline{\deg}\, h > \kappa$.  Then one has
\begin{align*}
\pi(x+h(x))-\pi(x) \sim \frac{h(x)}{\log x} \ (x \to \infty).
\end{align*}
\end{corollary}

In 1972, Huxley showed \cite{hux} that the hypotheses of Theorem  \ref{inghamtheorem} are satisfied for $b = \frac{12}{5}$ and thus for $\kappa = \frac{7}{12}$.  (Currently this is the smallest known admissible value for $\kappa$.)  Consequently, one has the following.

\begin{corollary}[Huxley, et.\ al.]\label{huxley}
Let $h$ be a real function defined on an unbounded subset of $\RR_{> 0}$, and suppose that $h$ is eventually positive,  $\log (x+h(x)) \sim \log x \ (x \to \infty)$, and $\underline{\deg}\, h > \frac{7}{12}$. Then one has
$$\pi(x+h(x))-\pi(x) \sim \frac{h(x)}{\log x}  \ (x \to \infty).$$
\end{corollary}

In 1988, Heath-Brown proved a generalization  \cite[Theorem]{heathbrown} of Corollary \ref{huxley}  that implies that the conclusion of the corollary holds even if $\underline{\deg}\, h = \frac{7}{12}$,  i.e., it implies the following.

\begin{theorem}[{\cite{heathbrown}}]\label{hbthm}
Let $h$ be a real function defined on an unbounded subset of $\RR_{> 0}$, and suppose that $h$ is eventually positive,  $\log (x+h(x)) \sim \log x \ (x \to \infty)$, and $\underline{\deg}\, h \geq \frac{7}{12}$. Then one has
$$\pi(x+h(x))-\pi(x) \sim \frac{h(x)}{\log x} \ (x \to \infty).$$
\end{theorem}

It seems that \cite[Theorem]{heathbrown} is the strongest unconditional result regarding Problem \ref{maierprob} to date.

 Theorem \ref{inghamtheorem} also yields the following.

\begin{corollary}\label{lindel}
Suppose that the density hypothesis holds.   Let $h$ be any real function, with domain an unbounded subset of $\RR_{>0}$,  that is eventually positive and satisfies  $\log (x+h(x)) \sim \log x \ (x \to \infty)$ and $\underline{\deg}\, h > \frac{1}{2}$ (e.g., $h(x) = a x^d$, where $a > 0$ and $d \in (\frac{1}{2},1]$).  Then one has
$$\pi(x+h(x))-\pi(x) \sim \frac{h(x)}{\log x} \ (x \to \infty).$$
\end{corollary}

The Riemann hypothesis yields the conclusion of  Corollary \ref{lindel} for a slightly larger class of functions than does the Lindel\"of or density hypothesis, e.g., for $h(x) = \sqrt{x}\, (\log x)^d$ for any $d > 2$.  To verify this, and, more generally, to further address Problem \ref{maierprob}, it makes sense to compare $\pi(f(x))-\pi(g(x))$ with $\li(f(x))-\li(g(x))$ for various functions $f$ and $g$ whose limits at $\infty$ are $\infty$.  The following result provides some useful asymptotics for functions of the form  $\li(f(x))-\li(g(x))$.

\begin{proposition}\label{primesinlarge}
Let $f$ and $g$ be real functions defined on a neighborhood of $\infty$.  Suppose that $f(x) > 1$ and $g(x) > 1$ for all $x \gg 0$  and that  $\log f(x) \sim \log g(x) \ (x \to \infty)$.
\begin{enumerate}
\item One has
$$\li(f(x))-\li(g(x)) \sim \frac{f(x)-g(x)}{\log f(x)}\sim  \frac{f(x)-g(x)}{\log g(x)} \ (x \to \infty).$$
\item Suppose that $\lim_{x \to \infty} f(x) = \lim_{x \to \infty} g(x) = \infty$ and  $f(x) \sim g(x) \ (x \to \infty)$.  Then one has
\begin{align*}
\li(f(x))-\li(g(x)) - \frac{f(x)-g(x)}{ \log f(x)} & \sim - \frac{1}{2} \frac{f(x)-g(x)}{\log f(x)} \left( 1-\frac{\log f(x)}{\log g(x)} \right) \ (x \to \infty)
\end{align*}
and
\begin{align*}
\li(f(x))-\li(g(x)) - \frac{f(x)-g(x)}{ \log g(x)} & \sim- \frac{1}{2} \frac{f(x)-g(x)}{\log g(x)} \left( 1-\frac{\log g(x)}{\log f(x)} \right) \ (x \to \infty).
\end{align*}
\end{enumerate}
\end{proposition}

\begin{proof}
Since $\li'(x) = \frac{1}{\log x}$ is decreasing, one has
$$ 0 < \frac{y-x}{\log y} < \li(y)-\li(x) < \frac{y-x}{\log x}$$
for all $y > x > 1$.  Let $\varepsilon > 0$, and choose $N$ so that $\left|1-\frac{\log g(x)}{\log f(x)}\right|< \varepsilon$ and $\left|1-\frac{\log f(x)}{\log g(x)}\right| < \varepsilon$ for all $x \geq N$.   Let $x \geq N$.  If $f(x) \geq g(x)$, then
$$ 0 \leq \frac{f(x)-g(x)}{\log f(x)} \leq \li(f(x))-\li(g(x)) \leq \frac{f(x)-g(x)}{\log g(x)},$$
so that
$$ 0 \leq \frac{f(x)-g(x)}{\log g(x)}-\li(f(x))-\li(g(x)) \leq \frac{f(x)-g(x)}{\log g(x)}\left(1-\frac{\log g(x)}{\log f(x)}\right)\leq  \varepsilon \frac{f(x)-g(x)}{\log g(x)}.$$  
On the other hand, if $f(x) \leq g(x)$, then
$$ 0 \geq \frac{f(x)-g(x)}{\log f(x)}-\li(f(x))-\li(g(x)) \geq \frac{f(x)-g(x)}{\log f(x)}-\frac{f(x)-g(x)}{\log g(x)}$$
and also
$$ \frac{f(x)-g(x)}{\log g(x)}-  \frac{f(x)-g(x)}{\log f(x)} =   \frac{f(x)-g(x)}{\log g(x)} \left(1 -\frac{\log g(x)} {\log f(x)} \right),$$
so that, adding the latter equation to the former inequalities, we see that
$$-\varepsilon\frac{f(x)-g(x)}{\log g(x)}  \geq \frac{f(x)-g(x)}{\log g(x)} \left(1 -\frac{\log f(x)} {\log g(x)} \right)  \geq \frac{f(x)-g(x)}{\log g(x)}-\li(f(x))-\li(g(x)) \geq 0.$$
Thus, in either case, one has
$$\left| \frac{f(x)-g(x)}{\log g(x)}-\li(f(x))-\li(g(x)) \right|  \leq \varepsilon \left|\frac{f(x)-g(x)}{\log g(x)} \right|.$$
Therefore
$$\frac{f(x)-g(x)}{\log g(x)}-\li(f(x))-\li(g(x)) = o \left( \frac{f(x)-g(x)}{\log g(x)} \right) \ (x \to \infty).$$
Statement (1) follows.

Suppose now the hypotheses in statement (2).   Observe that either of the first two asymptotics in statement (2) implies the other.   We may assume without loss of generality that one does not have $f(x) = g(x)$ for all sufficiently large $x$.  Let $y > x > 1$, and let $$H(x,y) =  \frac{y-x}{2} \left(\frac{1}{\log x}-\frac{1}{\log y}\right)$$  and $$E(x,y) = \frac{y-x}{\log x} - (\li(y)-\li(x)) -H(x,y).$$
Since $\frac{1}{\log x}$ is decreasing and concave up and  $\frac{d^2}{dx^2} \frac{1}{\log x} = \frac{2+ \log x }{x^2 (\log x)^3}$ is decreasing on $(1,\infty)$, by the trapezoid rule and corresponding error bound one has
 $$0 < E(x,y) \leq \frac{(y-x)^3}{12}\frac{2+ \log x }{x^2 (\log x)^3}  = (y-x)\frac{(2+ \log x) (\frac{y}{x}-1)^2 }{12(\log x)^3}.$$
Note then that
$$0 < \frac{E(x,y)}{H(x,y)} \leq \frac{\frac{(2+ \log x) (\frac{y}{x}-1 )^2}{12(\log x)^3}}{\frac{1}{2} \left(\frac{1}{\log x}-\frac{1}{\log y}\right)} = \frac{1}{6} \frac{\log y}{\log x}\left(1+\frac{2}{\log x}\right)\frac{(\frac{y}{x}-1)^2}{\log(\frac{y}{x})}.$$
Consequently, if $f(x)> g(x)$, then one has
$$0 < \frac{E(g(x),f(x))}{H(g(x),f(x))} \leq\frac{1}{6} \frac{\log f(x)}{\log g(x)}\left(1+\frac{2}{\log g(x)}\right)\frac{\frac{f(x)}{g(x)}-1}{\log\left(\frac{f(x)}{g(x)}\right)}\left(\frac{f(x)}{g(x)}-1\right).$$
  Since also  $\lim_{x \to \infty} \frac{f(x)}{g(x)} = \lim_{x \to \infty} \frac{\log f(x)}{\log g(x)} = 1$ and $\lim_{x \to \infty} g(x) = \infty$ and   $\lim_{z \to 1} \frac{z-1}{\log z} = 1$, the upper bound above tends to  $\frac{1}{6}(1)(1)(1)(0) = 0$ as $x \to \infty$.
One has a similar bound if $f(x) < g(x)$, and therefore
$$\lim_{x \to \infty} \frac{E(g(x),f(x))}{H(g(x),f(x))} = 0,$$
where the limit is taken over the set of all $x$ such that $f(x) \neq g(x)$.  Also, if  $f(x) = g(x)$, then one has $E(g(x),f(x)) = H(g(x),f(x)) = 0$.  It follows, then, that
$$E(g(x),f(x)) = o(H(g(x),f(x)) \ (x \to \infty),$$
and therefore
$$\frac{f(x)-g(x)}{\log g(x)} - (\li(f(x))-\li(g(x)))- H(g(x),g(y)) =  E(g(x),f(x)) = o (H(g(x),f(x)) \ (x \to \infty).$$
Statement (2) follows.  
\end{proof}

The proposition below applies Proposition \ref{primesinlarge} above to yield, in more restricted circumstances, the conclusion $\pi(f(x))-\pi(g(x)) \sim \li(f(x))-\li(g(x)) \ (x \to \infty)$.

\begin{proposition}\label{big}
Let $r$, $f$, and $g$ be real functions defined on a neighborhood of $\infty$ and satisfying the following conditions.
\begin{enumerate}
\item $\displaystyle \pi(x)-\li(x) = o\left(\frac{r(x)}{\log x} \right) \ (x \to \infty).$
\item  $\displaystyle \lim_{x \to \infty} f(x) = \lim_{x \to \infty} g(x) = \infty$.
\item  $\displaystyle \log f(x) \sim \log g(x) \ (x \to \infty)$.
\item $\displaystyle r(f(x)) = O (f(x)-g(x)) \ (x \to \infty)$.
\item $\displaystyle r(g(x)) = O(f(x)-g(x)) \ (x \to \infty).$
\end{enumerate}
Then one has
$$\pi(f(x))-\pi(g(x)) \sim \li(f(x))-\li(g(x)) \sim \frac{f(x)-g(x)}{\log f(x)}\sim  \frac{f(x)-g(x)}{\log g(x)} \ (x \to \infty).$$
\end{proposition}

\begin{proof}
By Proposition \ref{primesinlarge}(1), one has
$$\li(f(x))-\li(g(x)) \sim \frac{f(x)-g(x)}{\log f(x)}\sim  \frac{f(x)-g(x)}{\log g(x)} \ (x \to \infty).$$
By hypotheses (1), (2), and (4), one has
$$\pi(f(x))-\li(f(x)) = o \left(\frac{r(f(x))}{\log f(x)} \right) = o \left(\frac{f(x)-g(x)}{\log f(x)}\right) = o(\li(f(x))-\li(g(x)))  \ (x \to \infty),$$ 
and likewise, by (1), (2), and (5), one has
$$\pi(g(x))-\li(g(x)) =  o \left(\frac{r(g(x))}{\log g(x)} \right)  = o \left(\frac{f(x)-g(x)}{\log g(x) }\right) = o(\li(f(x))-\li(g(x)))   \ (x \to \infty).$$ 
Therefore one has
\begin{align*}
\pi(f(x))-\pi(g(x))  & = \li(f(x))-\li(g(x)) + ( \pi(f(x))-\li(f(x)))-( \pi(g(x))-\li(g(x)) )   \\
  & = \li(f(x))-\li(g(x)) + o \left(\li(f(x))-\li(g(x))\right)  \ (x \to \infty).
\end{align*}
The proposition follows.
\end{proof}

\begin{corollary}\label{bigcor}
Let $r$ and $h$ be real functions defined on a neighborhood of $\infty$ and satisfying the following conditions.
\begin{enumerate}
\item $\displaystyle \pi(x)-\li(x) = o\left(\frac{r(x)}{\log x} \right) \ (x \to \infty).$
\item $\displaystyle \lim_{x \to \infty} (x+h(x)) = \infty$.
\item  $\displaystyle \log (x+h(x)) \sim \log x \ (x \to \infty)$.
\item $\displaystyle r(x) = O \left(h(x) \right) \ (x \to \infty).$
\item $\displaystyle r(x+h(x))= O(h(x)) \ (x \to \infty)$.
\end{enumerate}
Then one has
$$\pi(x+h(x))-\pi(x) \sim \li(x+h(x))-\li(x) \sim \frac{h(x)}{\log x} \ (x \to \infty).$$
Moreover, conditions (4) and (5) together are implied by the following conditions.
\begin{enumerate}
\item[(6)] $r(x) = O(h(x)) \ (x \to \infty)$.
\item[(7)] $r(x+h(x))  = O(r(x)) \ (x \to \infty)$.
\end{enumerate}
\end{corollary}

Applying Proposition \ref{regvarfun} and Lemma \ref{hloglemma},  we obtain the following.

\begin{corollary}\label{bigcorr}
Let $r$ and $h$ be real functions defined on a neighborhood of $\infty$ and satisfying the following conditions.
\begin{enumerate}
\item $r$ is eventually positive, monotonic, and regularly varying of finite index.
\item $\displaystyle \pi(x)-\li(x) = o\left(\frac{r(x)}{\log x} \right) \ (x \to \infty)$.
\item $\displaystyle r(x) = O \left(h(x) \right) \ (x \to \infty)$.
\item $\log (x+h(x)) \sim \log x \ (x \to \infty)$.
\end{enumerate}
Then one has
$$\pi(x+h(x))-\pi(x) \sim \li(x+h(x))-\li(x) \sim \frac{h(x)}{\log x} \ (x \to \infty).$$
\end{corollary}

\begin{proof}
By Lemma \ref{hloglemma}, one can assume without loss of generality that $h(x) = o(x) \ (x \to \infty)$,  i.e., that 
$x+h(x) \sim x \ (x \to \infty)$, and then the corollary follows from Corollary \ref{bigcor} and Proposition \ref{regvarfun}.
\end{proof}

The corollary above relates to $\dege(\li-\pi)$ as follows.

\begin{theorem}\label{1overh}
Let $h$ be a real function defined on a neighborhood of $\infty$ with $\log (x+h(x)) \sim \log x \ (x \to \infty)$ and $$\underline{\dege}\, h >  \dege (\li-\pi)+(0,1,0,0,0,\ldots).$$  Then one has
$$\pi(x+h(x))-\pi(x) \sim \li(x+h(x))-\li(x) \sim \frac{h(x)}{\log x} \ (x \to \infty).$$
\end{theorem}

\begin{proof}
By Proposition \ref{oexppropstrong},  since $\underline{\dege}\, h > \dege ((\li-\pi)\log)$,  there exists an $r \in \mathbb{L}_{>0}$ such that 
$$\li(x)-\pi(x) = o\left ( \frac{r(x)}{\log x}\right) \ (x \to \infty)$$
and
$$r(x)  = o(h(x)) \ (x \to \infty).$$
 Moreover,  $r$ is eventually positive and increasing, and, since $d:= \deg r \in [\frac{1}{2},1]$,    it is also regularly varying of index $d$, by Proposition \ref{kregvar}.   The theorem follows, then, from Corollary \ref{bigcorr}.
\end{proof}

Note that, by Littlewood's theorem, $\underline{\dege}\, h >  \dege (\li-\pi)+(0,1,0,0,0,\ldots)$ implies
$$\underline{\dege}\, h  > (\tfrac{1}{2},0,0,1,0,0,0,\ldots) = \dege(\sqrt{x}\, (\log \log \log x))$$
and therefore
$$\sqrt{x} \, (\log \log \log x) = o(h(x)) \ (x \to \infty).$$
This sets a lower limit on the applicability of Theorem \ref{1overh}.    In particular, it cannot yield an eventual version of Legendre's conjecture.  Nevertheless, the theorem has the following consequences.

\begin{corollary}\label{longgcor}
Let $h$ be a real function defined on a neighborhood of $\infty$ with $\log (x+h(x)) \sim \log x \ (x \to \infty)$, and let $a \in \RR^*$.  Each of the following conjectures, each equipped with conditions on $h$, implies that
$$\pi(x+h(x))-\pi(x) \sim \li(x+h(x))-\li(x) \sim \frac{h(x)}{\log x} \ (x \to \infty).$$
\begin{enumerate}
\item  $\Theta < 1$; and $\underline{\deg}\, h >  \Theta$ (e.g., $h(x) = ax^d$ with $d \in (\Theta, 1]$).  
\item  $\Theta < 1$;  $\underline{\deg}\, h =  \Theta$;  and $\underline{\deg}\, h >  \Theta_1+1$ (or $\underline{\deg}_1 h > 2$)  (e.g., $h(x) = ax^\Theta(\log x)^d$ with $d >  \Theta_1+1$).  
\item  The Riemann hypothesis holds; and $\underline{\deg}\, h > \frac{1}{2}$ (e.g., $h(x) = ax^d$ with $d \in (\frac{1}{2}, 1]$).  
\item  The Riemann hypothesis holds; $\underline{\deg}\,  h = \frac{1}{2}$; and $\underline{\dege}_1 h >  \Theta_1+1$ (or $\underline{\dege}_1 h > 2$) (e.g., $h(x) = \sqrt{x} \,  (\log x)^{d}$ with $d > \Theta_1+1$).   
\item Conjecture \ref{eurekaconjecture}  holds;   $\underline{\deg}\,  h = \frac{1}{2}$; and $\underline{\dege}_1 h  > 0$ (e.g., $h(x) = a\sqrt{x} \,  (\log x)^{\varepsilon}$ with $\varepsilon > 0$). 
\item Conjecture \ref{eurekaconjecture2} holds; $\underline{\deg}\,  h = \frac{1}{2}$; $\underline{\dege}_1 h  = 0$; and $\underline{\dege}_2 h  > 0$ (e.g., $h(x) = a\sqrt{x} \,  (\log  \log x)^{\varepsilon}$ with $\varepsilon > 0$).  
\item Conjecture \ref{eurekaconjecture2}  holds and $\Theta_3  < \infty$; $\underline{\deg}\,  h = \frac{1}{2}$; $\underline{\dege}_1 h   = \underline{\dege}_2 h  = 0$; and $\underline{\dege}_3 h > \Theta_3$ (e.g., $h(x) = a\sqrt{x} \,  (\log \log \log x)^{d}$ with $d > \Theta_3$).  
\end{enumerate}
\end{corollary}

\begin{corollary}
Suppose that Montgomery's conjecture (\ref{MMC3}) holds, and let $r \in \mathbb{L}$ with 
$\lim_{x \to \infty} r(x) = \infty$ and $r(x) \ll \sqrt{x}$.  Then, for all $a \neq 0$, one has
\begin{align*}
\pi(x + a \sqrt{x}\,  (\log \log \log x)^{2}r(x))-\pi(x)  \sim  \frac{a\sqrt{x}\, (\log \log \log x)^{2}r(x)  }{\log x}\ (x \to \infty).
\end{align*}
\end{corollary}

Thus, the negation of the anti-Riemann hypothesis, the Riemann hypothesis, Conjecture \ref{eurekaconjecture},  Conjecture \ref{eurekaconjecture2}, and Montgomery's conjecture (\ref{MMC3}) each get progressively closer to implying 
\begin{align*}
\pi(x + \sqrt{x})-\pi(x)  \sim \frac{\sqrt{x} }{\log x}  \ (x \to \infty),
\end{align*}
which would yield the eventual version of Opperman's conjecture,  i.e.,
$$\pi(n^2 + n) -\pi(n^2) >0 \text{ and } \pi(n^2)-\pi(n^2-n) > 0, \quad \forall n \gg 0.$$

Theorem \ref{1overh} motivates the following problems.

\begin{outstandingproblem}\label{dinf}
Compute the infimum, or smallest, ${\mathbf h}$ of all $\dd \in \prod_{n = 1}^{\infty*}\overline{\RR}$ such that
$$\pi(x+h(x))-\pi(x) \sim \frac{h(x)}{\log x} \ (x \to \infty)$$
for all eventually positive real functions $h$ defined on a neighborhood of $\infty$ with $h(x) = O(x) \ (x \to \infty)$ and $\underline{\dege}\, h >  \dd$.  
\end{outstandingproblem}

\begin{problem}\label{dpid}
If 
$$\pi(x+r(x))-\pi(x) \sim \frac{r(x)}{\log x} \ (x \to \infty),$$
where $r \in \mathbb{L}_{> 0}$ and $r(x)= o(x) \ (x \to \infty)$,  then does it follow that
$$\pi(x+h(x))-\pi(x) \sim \frac{h(x)}{\log x} \ (x \to \infty),$$
for all eventually positive real functions $h$ defined on a neighborhood of $\infty$ with $h(x) = O(x) \ (x \to \infty)$ and $\underline{\dege}\, h >  \dege r$? 
\end{problem}

The following theorem partially addesses Problem \ref{dpid}.  It also extends Theorem  \ref{IVL} to functions more general than $r(x) = x^t$ for $t  \in (0, \frac{1}{2})$ and, optimistically, might even apply to some functions $r$ of degree $0$.  The proof follows  that of  Theorem  \ref{IVL}, with some straightforward modifications, and is provided at the end of this section.

 \begin{theorem}\label{IVL2}
 Let $r$ and $h$ be real functions, where $r$ is defined on $\RR_{>0}$ and is eventually positive and increasing,  and where $h$ is defined and eventually positive on an unbounded subset of $\RR_{> 0}$.  Suppose that the following conditions hold.
\begin{enumerate}
\item $\pi(x+r(x)) - \pi(x) \sim \frac{r(x)}{\log x}  \ (x \to \infty)$.
\item  $1 \ll r(x)= o(x) \ (x \to \infty)$.
\item $r(x+r(x)) \sim r(x) \ (x \to \infty)$. 
\item The function $s(x) = r(x+r(x))-r(x)$ is eventually (positive and) nonincreasing.
\item $(r(x)\log x)|_{\dom h} = O(h(x)) \ (x \to \infty)$.
\item  $h(x) = o\left(\frac{r(x)^2}{s(x) \log x } \right) \ (x \to \infty)$.
\end{enumerate}
Then  one has
 $$\pi(x+h(x)) - \pi(x) \sim \frac{h(x)}{\log x}  \ (x \to \infty).$$
 \end{theorem}

 \begin{corollary}\label{IVL2cor}
 Let $r$ and $h$ be real functions, where $r$ is defined on $\RR_{>0}$ and is eventually positive, increasing, and twice differentiable with eventually nonincreasing derivative, and where $h$ is defined and eventually positive on an unbounded subset of $\RR_{> 0}$.  Suppose that the following conditions hold.
 \begin{enumerate}
 \item $\pi(x+r(x)) - \pi(x) \sim \frac{r(x)}{\log x}  \ (x \to \infty)$.
\item  $r(x) \gg 1$.
\item  $r''(x) = o(x^{-3/2}) \ (x \to \infty)$.
\item $(r(x)\log x)|_{\dom h} = O(h(x)) \ (x \to \infty)$.
\item  $ h(x) = o\left(\frac{r(x)}{r'(x) \log x } \right) \ (x \to \infty)$, which holds if $h(x) = o \left( \frac{x}{\log x} \right) \ (x \to \infty)$ and $\lim_{x \to \infty} \frac{xr'(x)}{r(x)}$ exists.
\end{enumerate} 
Then one has
 $$\pi(x+h(x)) - \pi(x) \sim \frac{h(x)}{\log x}  \ (x \to \infty).$$
 \end{corollary}
 
\begin{proof}
One has $\frac{r'(x)}{\frac{1}{\sqrt{x}}} = o(1)$ and $\frac{r(x)}{\sqrt{x}} = o(1)$,  and therefore
$r(x)r'(x) = o(1)$, as $x \to \infty$.    By the mean value theorem, the hypotheses  on $r$  imply that $$r(x+r(x))-r(x) \sim r(x)r'(x) = o(r(x)) \ (x \to \infty).$$
Therefore,   the function $s(x) = r(x+r(x))-r(x)$ is  $o(1)$ and is therefore eventually decreasing to $0$.  It follows that all of the conditions of  Theorem \ref{IVL2} hold.
\end{proof}

 \begin{corollary}\label{IVL2cor2}
 Let $r$ and $h$ be real functions, where $r$ is Hardian and eventually positive, and where $h$ is defined and eventually positive on an unbounded subset of $\RR_{> 0}$.  Suppose that the following conditions hold.
 \begin{enumerate}
 \item $\pi(x+r(x)) - \pi(x) \sim \frac{r(x)}{\log x}  \ (x \to \infty)$.
\item $1 \ll r(x) = o(\sqrt{x}) \ (x \to \infty)$.
\item $(r(x)\log x)|_{\dom h} = O(h(x)) \ (x \to \infty)$.
\item $\log (x+h(x)) \sim \log x \ (x \to \infty)$.
\end{enumerate} 
Then (on $\dom h$) one has
 $$\pi(x+h(x)) - \pi(x) \sim \frac{h(x)}{\log x}  \ (x \to \infty).$$
 \end{corollary}

\begin{proof}
Without loss of generality, we can assume that $r$ is defined on $\RR_{>0}$.   Since $r$ is Hardian and $r(x) = o(\sqrt{x})$, one has $r'(x) = o(x^{-1/2})$ and $r''(x) = o(x^{-3/2})$.   Since $r(x) \gg 1$ and $r'(x) = o(1)$, the function $r(x)$ is eventually increasing, and $r'(x)$ is eventually decreasing.   Moreover, the limit $\lim_{x \to \infty} \frac{xr'(x)}{r(x)}$ exists.   Thus,  by Corollary \ref{IVL2cor}, the conclusion of the corollary holds if  $h(x) = o \left( \frac{x}{\log x} \right) \ (x \to \infty)$, and thus, by Lemma \ref{hloglemma}, it holds more generally if $\log (x+h(x)) \sim \log x \ (x \to \infty)$.
\end{proof}

 Although there are no known functions $r$ for which conditions (1) and (2) of  Corollary \ref{IVL2cor2} both hold, the most interesting applications of the corollary,  reaching beyond those of Theorem \ref{IVL}, are in various hypothetical situations where $\deg r = 0$.  See Section 14.5 for a discussion of the corollary in relation to the functions $r(x) = e^{(\log x)^t}$ for $t \in (0,1)$.  The following negative result,  proved by Maier in 1985,  shows that the functions $r(x) = (\log x)^a$ for  all $a \in \RR$ provide a lower limit for applications of the corollary.

\begin{theorem}[{\cite{maier}}]\label{maierthm}
Let $h(x) = (\log x)^\lambda$, where $\lambda > 1$.  Then one has
$$\limsup_{x \to \infty} \frac{\pi(x+h(x))-\pi(x)}{h(x)/\log x} > 1$$
and
$$\liminf_{x \to \infty} \frac{\pi(x+h(x))-\pi(x)}{h(x)/\log x} < 1.$$
\end{theorem}

\begin{corollary}\label{maierthmcor}
 For any $a \in \RR$, the asymptotic 
 $$\pi(x+h(x)) - \pi(x) \sim \frac{rhx)}{\log x}  \ (x \to \infty)$$
fails for $h(x) = (\log x)^a$.
\end{corollary}

By Theorems \ref{hbthm}, \ref{1overh}, and \ref{maierthm},   the infimum ${\mathbf h}$ defined in Problem  \ref{dinf} satisifies
$$(0,\infty,0,1,0,0,0,\ldots) \leq {\mathbf h} \leq \min(\dege(\li-\pi) + (0,1,0,0,0,\ldots),  (\tfrac{7}{12},-\infty,0,-1,0,0,0,\ldots)).$$  Moreover, by Theorem \ref{IVL}, if
$$\pi(x+x^t)-\pi(x) \sim \frac{x^t}{\log x} \ (x \to \infty),$$ 
where $t \in (0,1]$, then ${\mathbf h} \leq (t,0,0,0,\ldots)$,
and  ${\mathbf h}_0$ is the infimum of all $t > 0$ for which the asymptotic above holds.    Since ${\mathbf h}_0 \in (0,\frac{7}{12}]$ is finite, one has $${\mathbf h} \leq ({\mathbf h}_0,\infty,1,0,0,0,\ldots).$$
 We have also seen that the density hypothesis (or the Lindel\"of hypothesis) implies that ${\mathbf h}_0 \leq \frac{1}{2}$, hence also  that ${\mathbf h} \leq (\frac{1}{2},\infty,1,0,0,0,\ldots)$.  Note also that 
 $${\mathbf h}_0 = 0 \  \Leftrightarrow  \ (0,\infty,0,1,0,0,0,\ldots) \leq {\mathbf h} \leq (0,\infty,1,0,0,0,\ldots).$$
 Thus,  by Corollaries \ref{IVL2cor2} and \ref{maierthmcor},  we have the following.
 
 \begin{corollary}\label{maierthmcor2}
For any Hardian function $r$ (or  any $r \in \mathbb{L}$) with
  $$\pi(x+r(x)) - \pi(x) \sim \frac{r(x)}{\log x}  \ (x \to \infty),$$ 
  one has $$\dege r \geq (0,\infty, 0,1,0,0,0,\ldots).$$
  Moreover, if there exists such an $r$ with $\deg r = 0$, then $\dege_1 r = \infty$ and $\dege_2 r \in [0,1]$, and one has  $$\pi(x+h(x)) - \pi(x) \sim \frac{h(x)}{\log x}  \ (x \to \infty)$$ 
  for any function $h$ defined on an unbounded subset of $\RR_{> 0}$ with $\log (x+h(x)) \sim \log x \ (x \to \infty)$
  and  $(r(x)\log x)|_{\dom h} = O(h(x)) \ (x \to \infty)$,  which holds, for example, if $\log (x+h(x)) \sim \log x \ (x \to \infty)$ and $\underline{\dege} \, h > \dege r$.    Thus,  if such an $r$ exists, then Problem \ref{dpid} has a positive answer, and the infimum ${\mathbf h}$ defined in Problem  \ref{dinf} satisfies ${\mathbf h} \leq \dege r$.
  \end{corollary}

  A related problem is to determine the infimum $\alpha \ (\leq {\mathbf h}_0)$ of all $t > 0$ such that
$$\pi(x+x^t)-\pi(x) \gg \frac{x^t}{\log x} \ (x \to \infty),$$
which, by Theorem \ref{mv}, is equivalent to
$$\pi(x+x^t)-\pi(x) \asymp \frac{x^t}{\log x} \ (x \to \infty).$$   The best unconditional result in this direction is that 
\begin{align}\label{bakerres}
\pi(x+x^{21/40})-\pi(x) \geq \frac{9}{100} \frac{ x^{21/40}}{\log x}
\end{align}
for all $x \gg 0$ \cite{baker}, and thus $\alpha \leq \frac{21}{40}$.  One is also interested in the infimum $\alpha' \ (\leq \alpha)$ of all $t > 0$ such that
$$\pi(x+x^t)-\pi(x)  > 0, \quad \forall x \gg 0.$$  
 The best unconditional result in this direction is 
also $\alpha' \leq  \frac{21}{40}$.

 \begin{proposition}
For any eventually positive real function $h$ defined on an unbounded subset of $\RR_{>0}$ with $h(x) \ll x$ and 
  $$\pi(x+h(x)) - \pi(x) \sim \frac{h(x)}{\log x}  \ (x \to \infty),$$ 
  one has $$\dege h \geq (0,\infty, 0,1,0,0,0,\ldots).$$
  Moreover, if there exists such a function, say $r$, that is Hardian with $\deg r = 0$, then $\dege_1 r = \infty$ and $\dege_2 r \in [0,1]$, and one has  $$\pi(x+h(x)) - \pi(x) \sim \frac{h(x)}{\log x}  \ (x \to \infty)$$ 
  for any function $h$ defined on an unbounded subset of $\RR_{> 0}$ with $\log (x+h(x)) \sim \log x \ (x \to \infty)$
  and  $(r(x)\log x)|_{\dom h} = O(h(x)) \ (x \to \infty)$,  which holds, for example, if $\log (x+h(x)) \sim \log x \ (x \to \infty)$ and $\underline{\dege} \, h > \dege r$.    Thus,  if such an $r$ exists, then Problem \ref{dpid} has a positive answer, and the infimum ${\mathbf h}$ defined in Problem  \ref{dinf} satisfies ${\mathbf h} \leq \dege r$.
  \end{proposition}
  
  \begin{proof}
Suppose that $h \gg 1$ and $\dege h < (0,\infty, 0,1,0,0,0,\ldots)$,  i.e., that $1 \ll h(x) \leq (\log x)^a$ for some $a > 0$.   The infimum of all such $a$ is $d:= \dege_1 h \in [0,\infty)$.   One has $r = h|_X$ for some unbounded subset $X$ of $\dom h$ with $r(e^x)$ of exact degree  $d$, so that $r(x) = (\log x)^{d+o(1)}$.
The remainder of the proposition then follows from Corollaries \ref{IVL2cor2} and \ref{maierthmcor}.
  \end{proof}

Next, we relate our results on $\pi(x)$  in this section  to the first and second Chebyshev functions.

\begin{proposition}\label{pithetapsi}
Let $h$ be any real function satisfying  $0< h(x) \leq Mx$ for all $x\geq N$, where $M,N > 0$.
\begin{enumerate}
\item  One has
$$\frac{\vartheta(x+h(x))-\vartheta(x)}{h(x)} =\frac {\pi(x+h(x))-\pi(x)}{h(x)/\log x} + R(x),$$
where
$$0 \leq R(x) \leq \frac{2+2M+o(1)}{\log x}$$
for all $x \geq N$.
\item One has  $$\pi(x+h(x))-\pi(x)  \sim \frac{h(x)}{\log x} \ (x \to \infty)$$
if and only if 
$$\vartheta(x+h(x))-\vartheta(x) \sim h(x) \ (x \to \infty).$$  More generally, the lim sup (resp., lim inf)
of $\frac {\pi(x+h(x))-\pi(x)}{h(x)/\log x}$ and $\frac{\vartheta(x+h(x))-\vartheta(x)}{h(x)}$ as $x \to \infty$ are equal.
 \item One has
$$\frac{\psi(x+h(x))-\psi(x)}{h(x)} = \frac{\vartheta(x+h(x))-\vartheta(x)}{h(x)} + S(x),$$
where
$$0 \leq S(x) \leq \frac{1}{2\sqrt{x}} + o\left(\frac{\sqrt{x}}{h(x)} \right)$$
for all $x \geq N$. 
 \item Suppose that $\sqrt{x} = O(h(x)) \ (x \to \infty)$.  Then one has
  $$\vartheta(x+h(x))- \vartheta(x) \sim h(x) \ (x \to \infty)$$
 if and only if
  $$\psi(x+h(x))- \psi(x) \sim h(x) \ (x \to \infty).$$  More generally, the lim sup (resp., lim inf)
of  $\frac{\vartheta(x+h(x))-\vartheta(x)}{h(x)}$ and $\frac{\psi(x+h(x))-\psi(x)}{h(x)}$  as $x \to \infty$ are equal.
\end{enumerate}
\end{proposition}

\begin{proof}
Let $x \geq  N$.
 By Proposition \ref{pitheta} (which follows from Abel's summation formula), one has
$$\vartheta(x) = \pi(x) \log x -\int_0^x \frac{\pi(t)}{t} \, dt$$
and therefore
\begin{align*}
R(x) h(x)=  \pi(x+h(x))\log\left(1+\frac{h(x)}{x}\right) +  \int_x^{x+h(x)} \frac{\pi(t)}{t} \, dt \geq 0.
\end{align*}
Note that
\begin{align*}
 0 \leq \int_x^{x+h(x)} \frac{\pi(t)}{t} \, dt &   \leq \pi(x+h(x))\int_x^{x+h(x)} \frac{1}{t} \, dt  \\
& = \pi(x+h(x))\log\left(1+\frac{h(x)}{x}\right)
\end{align*} 
and
\begin{align*}
0 \leq \pi(x+h(x))\log\left(1+\frac{h(x)}{x}\right) 
& \leq \pi(x+h(x))\frac{h(x)}{x}  \\
& \sim \frac{x+h(x)}{\log( x+h(x))}\frac{h(x)}{x} \\
& \sim \left( 1+ \frac{h(x)}{x} \right)\frac{h(x)}{\log x} \\
& \leq  (1+M) \frac{h(x) }{\log x}.
\end{align*} 
Therefore, one has
$$0 \leq R(x) h(x) \leq (2+2M+o(1))\frac{h(x) }{\log x}.$$
This proves (1), from which (2) immediately follows.   The proof of statements (3) and (4) is similar and is based on Proposition \ref{etaconj3}.
\end{proof}

Thus, once again,  $\sqrt{x}$ shows up as a natural barrier to understanding asymptotics for $\pi(x+h(x))-\pi(x)$, this time, in relation to  asymptotics for $\psi(x+h(x))- \psi(x)$.

The proposition, along with Theorem \ref{1overh},  yields the following.

\begin{corollary}\label{1overhhcor}
Let $h$ be an eventually positive real function defined on a neighborhood of $\infty$ with $h(x) = O(x) \ (x \to \infty)$ and $$\underline{\dege}\, h >  \dege (\li-\pi)+(0,1,0,0,0,\ldots).$$  Then one has
$$\pi(x+h(x))-\pi(x) \sim \frac{h(x)}{\log x} \ (x \to \infty),$$
  $$\vartheta(x+h(x))- \vartheta(x) \sim h(x) \ (x \to \infty),$$
and
  $$\psi(x+h(x))- \psi(x) \sim h(x) \ (x \to \infty).$$
\end{corollary}

Note also the following trivial result.

\begin{proposition}
Let $h$ be any eventually positive real function defined on $[N,\infty)$, where $N > 0$.    The following conditions are equivalent.
\begin{enumerate}
\item For all $x \geq N$, there is a prime $p$ with $x < p \leq x+h(x)$.
\item $\pi(x+h(x))-\pi(x) > 0$ for all $x \geq N$.
\item $\vartheta(x+h(x))-\vartheta(x)  > 0$ for all $x \geq N$.
\end{enumerate}
Moreover, the equivalent conditions above imply the following equivalent conditions.
\begin{enumerate}
\item[(4)]   For all $x \geq N$, there is a prime power $q>1$ with $x < q \leq x+h(x)$.
\item[(5)] $\psi(x+h(x))-\psi(x) > 0$ for all $x \geq N$.
\end{enumerate}
\end{proposition}

Next, we address asymptotics for prime counts in intervals of constant length.  Following \cite{hensley}, for any positive integer $n$, we let
$$\rho(n) = \limsup_{m \to \infty} \, (\pi(m+n)-\pi(m))$$
and
$$\rho_1(n) = \max_{m \geq n} \, (\pi(m+n)-\pi(m)),$$
and we let $\rho^*(n)$ denote the maximum length $k$ of any admissible $k$-tuple $(m_1,m_2,\ldots,m_k)$ of distinct integers $m_i$ with $0\leq m_i < n$ for all $i$, where a $k$-tuple $(m_1,m_2,\ldots,m_k)$ is {\bf admissible}\index{admissible $k$-tuple} if there are no primes $p$ for which the tuple includes every possible residue modulo $p$.  The {\bf prime $k$-tuples  conjecture}\index{prime $k$-tuples conjecture} states that,  for any admissible $k$-tuple $(m_1,m_2,\ldots,m_k)$ (e.g., $(0,2)$, $(0,2,6)$, $(0,4,6)$, and $(0,2,6,8)$), there are infinitely many positive integers $n$ such that the $k$-tuple $(n+m_1, n+m_2,\ldots,n+m_k)$ is a $k$-tuple of primes.  It is not difficult to show that
$$\rho(n) \leq \rho_1(n) \leq \rho^*(n)$$
for all $n$, and, furthermore, that the prime $k$-tuples conjecture implies that equalities hold for all $n$ \cite{hensley}.
Clearly one has $\rho^*(2) = 1$, and the {\bf twin prime conjecture},\index{twin prime conjecture} still unresolved to this day, states that there are infinitely many primes that differ by $2$, or, equivalently, that $\rho(2) = 1$.

Consider now the constant
$$L =  \limsup_{n \to \infty}\, \limsup_{m \to \infty} \frac{\pi(m+n)-\pi(m)}{\pi(n)} =  \limsup_{n \to \infty}  \frac{\rho(n)}{\pi(n)}.$$ 
By the  Theorem \ref{mv}, one has $L \leq 2$.   The best known unconditional lower bound for $L$ is the trivial lower bound $0$.  Assuming the prime $k$-tuples conjecture, however, one has $L \geq 1$.  To prove this, note first that $L \leq L'$, where
$$L' =  \limsup_{n \to \infty}  \frac{\rho^*(n)}{\pi(n)}.$$  
Assuming the prime $k$-tuples conjecture, one has $L = L'$.
Moreover, since, for any $k$, the $k$-tuple consisting of the first $k$ primes larger than $k$ is admissible, one has 
$$\rho^*(p_{\pi(k)+k+1}-p_{\pi(k)+1}+1) \geq k,$$
while also
$$p_{\pi(k)+k+1}-p_{\pi(k)+1} +1 \sim p_{\pi(k)+k+1} \sim p_k \ (k \to \infty)$$
and therefore
$$\pi(p_{\pi(k)+k+1}-p_{\pi(k)+1}+1) \sim \pi(p_k) = k \ (k \to \infty).$$ 
Thus, it follows that
$$L' \geq 1.$$
Assuming the prime $k$-tuples conjecture,  then, one has $$1 \leq L \leq 2.$$

It is known, unconditionally, that, for every $\varepsilon > 0$, one has
\begin{align}\label{hen}
\rho^*(n) \geq \pi(n) + (\log 2-\varepsilon)\lceil n/(\log n)^2 \rceil
\end{align}
for all $n \gg 0$ \cite{hensley}.  It follows that
$$  \liminf_{n \to \infty}  \frac{\rho^*(n)}{\pi(n)} \geq 1.$$  
Thus, assuming the prime $k$-tuples conjecture, one has $$1 \leq \liminf_{n \to \infty}  \frac{\rho(n)}{\pi(n)} \leq \limsup_{n \to \infty}  \frac{\rho(n)}{\pi(n)} = L \leq 2,$$
and one has $L = 1$ if and only if
$$\lim_{n \to \infty}\, \limsup_{m \to \infty} \frac{\pi(m+n)-\pi(m)}{\pi(n)} = 1,$$
if and only if
$$\rho(n) \sim \pi(n) \ (n \to \infty).$$  The author is unaware of any conjectures that determine the exact value of $L$.

\begin{remark}[The second Hardy--Littlewood conjecture]
The {\bf second Hardy--Littlewood conjecture}\index{second Hardy--Littlewood conjecture} states that
$$\pi(m+n) \leq \pi(m)+\pi(n)$$
for all integers $m,n>1$, or, equivalently, that the number of primes from $m + 1$ to $m + n$ is always less than or equal to the number of primes from $1$ to $n$.  By definition of $\rho_1(n)$, the conjecture is equivalent to $$\rho_1(n) \leq \pi(n)$$ for  all $n >1$.  Although no counterexamples to the conjecture are known, it is widely believed to be false, as it is known to be incompatible with the prime $k$-tuples conjecture  \cite{hensley}.  Indeed, by (\ref{hen}), if the prime $k$-tuples conjecture holds, then one has $\rho_1(n) = \rho(n) = \rho^*(n) > \pi(n)$ for all sufficiently large $n$.  To make matters worse, by definition of $\rho(n)$, if $\rho(n_0) > \pi(n_0)$ for some particular $n_0$, then $$\pi(m+n_0)-\pi(m) > \pi(n_0),$$
and thus
$$\pi(m+n_0) > \pi(m)+\pi(n_0),$$ for infinitely many positive integers $m$.  
\end{remark}

Finally, we provide a proof of Theorem \ref{IVL2}.
 
 \begin{proof}[Proof of Theorem \ref{IVL2}]
Let  $$T(x) = x+r(x)$$
for all $x \geq 0$.  Also,  for all $x \in \dom h$,  let
 $$K(x) = \left \lfloor \frac{h(x)}{r(x)} \right \rfloor.$$
 For all  $x \in \dom h$, one has
 $$0 \leq h(x)-K(x)r(x) < r(x) = o(h(x)) \ (x \to \infty),$$
 and therefore
 $$K(x)r(x) \sim h(x) \ (x \to \infty).$$
 Note that, for all $x > 0$, one has
 $$0 \leq r(x) = T(x)-x \leq T(T(x))-T(x) \leq T(T(T(x)))-T(T(x)) \leq \cdots,$$
so that
  $$0 \leq x+ kr(x) \leq T^{\circ k}(x)$$
  for all positive integers $k$,
  whence
  $$0 \leq x+K(x)r(x)  \leq T^{\circ K(x)}(x),$$
  for all $x \in \dom h$.   For all $x \in \dom h$, we may write
\begin{align}\label{pp1b}
\pi(x+h(x))-\pi(x) = \sum_{1 \leq k \leq K(x)} (\pi(T^{\circ k}(x))-\pi(T^{\circ (k-1)}(x)) + \pi(x+h(x))-\pi(T^{\circ K(x)}(x)).
\end{align}
By our hypotheses, for all fixed positive integers $k$ one has
 $$\pi(T^{\circ k}(x))-\pi(T^{\circ (k-1)}(x)) \sim \frac{r(T^{\circ (k-1)}(x))}{\log T^{\circ (k-1)}(x)} \sim \frac{r(x)}{\log x} \ (x \to \infty)$$
and
\begin{align}\label{pp2b}
 \sum_{1 \leq k \leq K(x)} (\pi(T^{\circ k}(x))-\pi(T^{\circ (k-1)}(x))  \sim \frac{K(x)r(x)}{\log x} \sim \frac{h(x)}{\log x} \ (x \to \infty).
 \end{align}
Indeed, to verify these two claims,  it suffices to show that
 $$r(T^{\circ (K(x)-1)}(x)) \sim r(x) \ (x \to \infty),$$
 or, equivalently, that 
  $$r(T^{\circ K(x)}(x)) \sim r(x) \ (x \to \infty).$$
To this end,  note that
\begin{align*}
1  &\leq \frac{r(T^{\circ K(x)}(x)) }{r(x)} & \\ 
  & = \frac{r(T(x))}{r(x)} \cdot  \frac{r(T(T(x)))}{r(T(x))}  \cdots \frac{r(T^{\circ K(x)}(x) )}{r(T^{\circ (K(x)-1)}(x) )}    \\  
  & =  \left( 1+ \frac{s(x)}{r(x)}  \right)  \left( 1+  \frac{s(T(x))}{r(T(x))}   \right) \cdots \left( 1 + \frac{s(T^{\circ (K(x)-1)}(x))}{r(T^{\circ (K(x)-1)}(x))}  \right)  \\
  & \leq \left( 1+ \frac{s(x)}{r(x)}  \right)^{K(x)} \\
  & \leq \left( 1+ \frac{s(x)}{r(x)}  \right)^{h(x)/r(x)} \\
  & =  \left(\left( 1+ \frac{s(x)}{r(x)}  \right)^{r(x)/s(x)}\right)^{h(x)s(x)/r(x)^2} \\
  & \leq e^{o(1/\log x)}  \\
  &\sim 1 \ (x \to \infty)
  \end{align*}
  for all $x \gg 0$ in $\dom h = \dom K$.
This proves (\ref{pp2b}).  By (\ref{pp1b}) and (\ref{pp2b}), then,  it remains only to show that
 \begin{align}\label{NHT1b}
  \pi(x+h(x))-\pi(T^{\circ K(x)}(x)) = o\left( \frac{h(x)}{\log x} \right) \ (x \to \infty).
  \end{align}

To prove (\ref{NHT1b}), we show that the difference $x+h(x)-T^{\circ K(x)}(x)$ is  $o\left( \frac{h(x)}{\log x} \right)$ and then apply Corollary \ref{mvcor}. Note  first that
$$x+h(x)-T^{\circ K(x)}(x) \leq h(x)-K(x)r(x) < r(x) = o\left( \frac{h(x)}{\log x} \right)$$
on $\dom h$.  We wish to bound the difference by $o\left( \frac{h(x)}{\log x} \right)$ also from below.
One has
\begin{align}\label{NHT0b}
T(T(x))-T(x)-(T(x)-x) = r(x+r(x))-r(x) = s(x)  \ (x \to \infty).
\end{align}
Let assume that $s(x)$  is nonincreasing for $x \geq N$.
Since $x \leq T(x) \leq T(T(x)) \leq \cdots$ for all $x$,  for all $k\geq 2$ we have
\begin{align*}
T^{\circ k}(x)-T^{\circ (k-1)}(x)-(T^{\circ (k-1)}(x)-T^{\circ (k-2)}(x)) =  s((T^{\circ (k-2)}(x))) \leq  s(x)
\end{align*}
for all $k \geq 2$ and all $x \geq N$.   Summing over $k$,  for all $x \geq N$ we then have
$$T^{\circ k}(x)-T^{\circ (k-1)}(x) \leq r(x)+ (k-1)s(x)$$
for all $k \geq 1$, whence, summing over $k$ from $1$ to $K(x)$, we find that
 $$T^{\circ K(x)}(x) \leq  x + K(x) r(x)+ K(x)^2s(x)$$
 for all $x \geq N$.
  It follows that
\begin{align*}
x+h(x)-T^{\circ K(x)}(x) & \geq  h(x)-K(x) r(x)- K(x)^2s(x) \\ 
& \geq -  K(x)^2 s(x) \\
& \sim  - \frac{h(x)^2}{r(x)^2}   s(x) \\
& = o\left( \frac{h(x)}{\log x} \right)
\end{align*}
as $x \to \infty$.
Therefore, one has
\begin{align}\label{NHT2b}
x+h(x)-T^{\circ K(x)}(x)= o\left( \frac{h(x)}{\log x} \right).
\end{align}
Next,  we apply Corollary  \ref{mvcor} to obtain
 \begin{align*}
 |\pi(x+h(x))-\pi(T^{\circ K(x)}(x))|  \leq 2 \li ( |x+h(x)-T^{\circ K(x)}(x)| + \mu)+1
 \end{align*}
 for all $x \gg 0$ in  $\dom h$.   Finally, by (\ref{NHT2b}),  for all $x \gg 0$ in $\dom h$, one has
 $$|x+h(x)-T^{\circ K(x)}(x)|  < \frac{ h(x)}{\log x},$$
whence
  \begin{align*}
 |\pi(x+h(x))-\pi(T^{\circ K(x)}(x))| & < 2\li (  h(x)/\log x + \mu) +1 \\
& = (2+o(1)) \frac{h(x)/\log x + \mu}{\log(  h(x)/\log x + \mu)} \\
&  \ll \frac{ h(x)/\log x}{\log ( r(x))} \\
&  = o\left(\frac{ h(x)}{\log x}\right).
 \end{align*}
  This establishes   (\ref{NHT1b}) and therefore completes the proof.
 \end{proof}

\section{The $n$th prime  $p_n$}

Let $p_n$ denote the $n$th prime for every positive integer $n$.  We call the function $p_n$ the {\bf prime listing function}.\index{prime listing function $p_n$} It is well known that the prime number theorem is equivalent to
$$p_n \sim n \log n  \ (n \to \infty).$$ 
The prime listing  function $p_n: \ZZ_{\geq 0} \longrightarrow \RR_{\geq 0}$ is a right inverse of the prime counting function $\pi(x): \RR_{\geq 0} \longrightarrow \ZZ_{\geq 0}$, where we define $p_0 = 0$:  one has $\pi(p_n) = n$ for all $n$.  Moreover, the reverse composition $p_{\pi(x)}: \RR_{\geq 0} \longrightarrow  \RR_{\geq 0}$ for any $x$   equals the largest prime less than or equal to $x$, and it is thus a  nondecreasing step function satisfying $p_{\pi(x)} \leq x$ for all $x$, where equality holds if and only if $x$ is a prime or $0$.  Likewise, $p_{\pi(x)+1}$ is the smallest prime greater than $x$, and so $p_{\pi(x)} \leq  x < p_{\pi(x)+1}$ for all $x$.   Note also that $\pi(x) =  n$  if and only if $p_n \leq x < p_{n+1}$.  Thus, one has $p_n \leq x$ if and only if  $n \leq \pi(x)$, for all  $x \in \RR_{\geq 0}$ and all $n \in \ZZ_{\geq 0}$.  In fancy terminology, we say that the pair  $p_n: \ZZ_{\geq 0} \longrightarrow \RR_{\geq 0}$ and $\pi(x): \RR_{\geq 0} \longrightarrow \ZZ_{\geq 0}$ of functions is a {\it monotone Galois injection} \cite{ore} of the poset  $(\ZZ_{\geq 0}, \leq)$ into  the poset  $(\RR_{\geq 0}, \leq)$, and, in particular, $p_n$ is a left adjoint of $\pi(x)$.  (Analogously, the same situation holds between  the identity function $\id: \ZZ \longrightarrow \RR$ and the floor function $\lfloor x \rfloor: \RR \longrightarrow \ZZ$.)    One most naturally extends  $p_n$ to a function $p_x : \RR_{\geq 0} \longrightarrow \RR_{\geq 0}$ by setting $$p_x = p_{\lceil x \rceil}\index[symbols]{.rt Pa@$p_x$}$$ for all $x \geq 0$.  Indeed, one then has
 $ p_x \leq y$ if and only if $x \leq \pi(y)$, for all $x, y \geq 0$. Thus,   the pair  $p_x$  and $\pi(x)$ is a {\it monotone Galois connection} \cite{ore},  now  from $(\RR_{\geq 0}, \leq)$ to $(\RR_{\geq 0}, \leq)$, where $\pi(p_x) = \lceil x \rceil$ and where $p_{\pi(x)}$ is as before.  (Analogously, the same situation holds between the ceiling function $\lceil x \rceil: \RR \longrightarrow \RR$ and the floor function $\lfloor x \rfloor: \RR \longrightarrow \RR$.)   

Since $\li(x)$ is  increasing on $(1,\infty)$ with range $\RR$, it has an increasing (two-sided) inverse $\li^{-1}: \RR \longrightarrow (1,\infty)$.  The derivative of $\li^{-1} x$ is the function $$\frac{d}{dx} \li^{-1}x = \log \li^{-1} x  = \Ei^{-1} x > 0,$$
which has inverse $\Ei x = \li(e^x): (0,\infty) \to \RR$.   The second derivative of $\li^{-1} x$  is also positive, satisfying
 $$0< \frac{d^2}{dx^2} \li^{-1}x = \frac{\log \li^{-1}x}{\li^{-1}x}  \leq \frac{1}{e}$$
for all $x$,
where the maximum value of $\frac{1}{e}$ is assumed at $x = \li(e) = 1.895117816355\ldots$, and thus the derivative of $\li^{-1}x$ has an inflection point at the point $(\li(e),1)$.
One also has
$$\frac{d}{dx} \li^{-1}x = \log \li^{-1} x \sim \log x  \ (x \to \infty)$$
and
$$ \frac{d^2}{dx^2} \li^{-1}x =  \frac{\log \li^{-1}x}{\li^{-1}x} \sim \frac{1}{x} \ (x \to \infty).$$
Note that, because $\li$ and $\li^{-1}$ are increasing, one has  $ \li^{-1}x \leq y$ if and only if  $x \leq \li y$, for all $x \in \RR$ and all $x \in (1,\infty)$, and the pair forms a monotone Galois isomorphism from  $(\RR, \leq)$ to $((1,\infty), \leq)$.

Because of these adjoint relationships, one might expect that, just as $\pi(x)$ and $\li(x)$ are closely related, so are their respective right adjoints $p_x$ and $\li^{-1} x$.
Indeed, it is known \cite{mas} that
\begin{align*}
p_n - \li^{-1}n = O \left( n e^{-c \sqrt{\log n}} \right)  \ (n \to \infty)
\end{align*}
for some $c > 0$, so that
$$ \dege(p_n-\li^{-1}n) \leq (1,-\infty,-\tfrac{1}{2},0,0,0,\ldots).$$
 It is also known \cite{rey} that $\deg(p_n-\li^{-1}n) \in [\frac{1}{2},1]$, and the Riemann hypothesis is equivalent to $\deg(p_n-\li^{-1}n) = \frac{1}{2}$ and to
$$|p_n - \li^{-1}n| < \frac{1}{\pi} \sqrt{n} \,  (\log n)^{5/2}, \quad \forall n \geq 11,$$
and therefore also to
$$\dege(p_n-\li^{-1}n) \leq (\tfrac{1}{2}, \tfrac{5}{2},0,0,0,\ldots).$$  

It is perhaps intuitively plausible that one can express $\dege(p_x-\li^{-1}x)$ in terms of $\dege(\li-\pi)$ and vice versa, as follows.

\begin{theorem}\label{pnthm}
Let $\Theta$ denote the Riemann constant $\deg(\li-\pi)$. One has
$$\dege(p_x - \li^{-1}x)  =  \dege(\li-\pi) + (0,\Theta+1,0,0,0,\ldots)$$ and
$$ \dege(\li(p_x)-\pi(p_x))  = \dege (\li-\pi) + (0, \Theta, 0 , 0 , 0 ,\ldots).$$
\end{theorem}

To prove the theorem, we first note the following lemma, which follows immediately from Proposition \ref{primesinlarge}.

\begin{lemma}\label{ali}
One has
$$\li(p_x)-x \sim \frac{p_x - \li^{-1}x}{\log \li^{-1}x}  \sim \frac{p_x - \li^{-1}x}{\log x}  \ (x \to \infty).$$
More precisely, one has
\begin{align*}
\li(p_x)-x - \frac{p_x - \li^{-1}x}{\log \li^{-1}x}  & \sim - \frac{1}{2} \frac{p_x-\li^{-1}x}{(\log x)^2} \log \left(\frac{p_x}{\li^{-1}x} \right) \ (x \to \infty) \\
 & \sim - \frac{1}{2} \frac{(p_x-\li^{-1}x)^2}{x(\log x)^3}  \ (x \to \infty).
\end{align*}
\end{lemma}

The following proposition implies that the two equations in Theorem \ref{pnthm} are equivalent.

\begin{proposition}\label{lipiprop}
One has
\begin{align*}
\dege(p_n-\li^{-1}n) = \dege(p_x-\li^{-1}x)  =  \dege(\li(p_x)-x)+(0,1,0,0,0,\ldots)
\end{align*}
and
\begin{align*}
\dege(\li(p_x)-x)  =  \dege(\li(p_x)-\pi(p_x)) =  \dege(\li(p_n)-\pi(p_n)).
\end{align*}
\end{proposition}

\begin{proof}
From Lemma \ref{ali}, it follows that
$$\dege(p_x-\li^{-1}x)  =  \dege(\li(p_x)-x)+(0,1,0,0,0,\ldots).$$
Since $$\li^{-1}(x+1)-\li^{-1}(x) \sim \log x  \ (x \to \infty),$$
the function $\li^{-1}\lceil x \rceil-\li^{-1}x$ has logexponential degree $(0,1,0,0,0,\ldots) <  \dege(p_n-\li^{-1}n)$, and therefore
$$\dege(p_x-\li^{-1}x) = \dege(p_x-\li^{-1}\lceil x \rceil) = \dege(p_n-\li^{-1}n)$$
(where the equality $\dege(p_x-\li^{-1}\lceil x \rceil) = \dege(p_n-\li^{-1}n)$ follows from Proposition \ref{arithb}).
Similarly, since  $\pi(p_x)-x = \lceil x \rceil-x$ has logexponential degree $(0,0,0,\ldots) < \dege(\li(p_x)-x)$, one has
$$ \dege(\li(p_x)-x) =  \dege(\li(p_x)-\pi(p_x)) = \dege(\li(p_n)-\pi(p_n)).$$
This completes the proof.
\end{proof}

Now, we prove Theorem \ref{pnthm}.

\begin{proof}[Proof of Theorem \ref{pnthm}]
By Proposition \ref{lipiprop}, it suffices to prove that
$$ \dege(\li(p_x)-\pi(p_x))  = \dege (\li-\pi) + (0, \Theta, 0 , 0 , 0 ,\ldots).$$

Let $r \in \mathbb{L}_{>0}$  with  $$\li(x)-\pi(x) = O(r(x)) \ (x \to \infty)$$ 
and $\deg r  = \Theta$. 
 Since $p_x \sim x \log x \ (x \to \infty)$ and $r$ has finite degree and exact logexponential degree, by Theorem \ref{fgie} one has
$$\dege r(p_x) = \dege r + (0,\Theta,0,0,0,\ldots).$$
Moreover, since $$\li(p_x)-\pi(p_x) = O(r(p_x)) \ (x \to \infty),$$ one has
$$\dege(\li(p_x)-\pi(p_x)) \leq \dege r(p_x) = \dege r + (0,\Theta,0,0,0,\ldots),$$
whence
$$\dege(\li(p_x)-\pi(p_x))- (0,\Theta,0,0,0,\ldots)\leq \dege r.$$
Taking the infimum over all $r$ as chosen, it follows from Theorem \ref{infpropexp}  that
$$\dege(\li(p_x)-\pi(p_x))- (0,\Theta,0,0,0,\ldots) \leq \dege(\li-\pi).$$

To prove the reverse inequality, let $r \in \mathbb{L}_{>0}$ with
$$|\li(p_x)-\pi(p_x)| \leq r(x), \quad \forall x \geq N,$$
where $N > 0$.  The function $r$ is eventually increasing, so we may suppose $r$ is increasing for $x \geq N$.  Let $x \geq 2p_N$.  Let $p_n = p_{\pi(x)}$ be the largest prime less than or equal to $x$, so that
$n \log n < p_n \leq x < p_{n+1} < 2p_n$ and therefore $x \geq p_n \geq p_N$ and thus $n \geq N$.  The inverse of the function $x \log x$  on $[1/e,\infty)$ is $L(x) = \frac{x}{W(x)} = e^{W(x)}$, where $W(x): [-1/e,\infty) \longrightarrow [-1,\infty)$ is the Lambert $W$ function, and the function $L(x)$ is increasing on $ [-1/e,\infty)$.   Thus one has $n < L(x)$.   Now, $\pi(x) = n$, and $\li(x)$ is increasing,  so that
$$|\li(x)-\pi(x)| \leq \max(|\li(p_n)-n|,|\li(p_{n+1})-n|).$$
If $|\li(x)-\pi(x)| \leq |\li(p_n)-n|$, then
$$|\li(x)-\pi(x)| \leq r(n).$$
On the other hand, if $|\li(x)-\pi(x)| \leq |\li(p_{n+1})-n|$, then
\begin{align*}
|\li(x)-\pi(x)| & \leq |\li(p_{n+1})-n|   \\
  & \leq  |\li(p_{n+1})-(n+1)|+1  \\
  & \leq r(n+1) + 1. 
\end{align*}
Thus, since $r$ is increasing for $x \geq N$, one has
$$|\li(x)-\pi(x)| \leq r(n+1) + 1 \leq r(L(x)+1)+1$$
for all $x \geq N$.   Since $L(x)+1 \sim x(\log x)^{-1} \ (x \to \infty)$,  it follows that
$$\Theta = \deg(\li-\pi) \leq \deg (r(L(x)+1)+1) = \deg r.$$
Taking the infimum over all $r$ as  chosen,  we see that
$$\Theta \leq \deg(\li(p_x)-\pi(p_x)) \leq \deg(\li-\pi) = \Theta,$$
whence equalities hold.   Thus,  for any $r$ as chosen earlier, we may assume additionally that $\deg r = \Theta$.  Then, by Theorem \ref{fgie},  one has
$$ \dege( r( L(x)+1) +1) = \dege r( L(x)+1)  = \dege r -(0,\Theta, 0,0,0,\ldots),$$
and therefore
$$\dege(\li-\pi)+(0,\Theta,0,0,0,\ldots)\leq  \dege( r( L(x)+1) +1)+(0,\Theta, 0,0,0,\ldots)  = \dege r.$$
Taking the infimum over all $r$ as  chosen, it follows from Theorem \ref{infpropexp} that
$$\dege(\li-\pi)+ (0, \Theta, 0 , 0 , 0 ,\ldots) \leq \dege(\li(p_x)-\pi(p_x)).$$
The theorem follows.
\end{proof}

\begin{corollary} One has the following.
\begin{enumerate}
\item
 One has $$\dege(p_x- \li^{-1}x ) \geq (\tfrac{1}{2},\Theta,0,1,0,0,\ldots) \geq  (\tfrac{1}{2},\tfrac{1}{2},0,1,0,0,\ldots)$$
and therefore
$$p_x- \li^{-1}x \neq O\left( \sqrt{x}\, (\log x)^\Theta (\log \log \log x)^t\right) \ (x \to \infty) \text{ if } t < 1$$
and
$$p_x- \li^{-1}x \neq O\left( \sqrt{x \log x} \, (\log \log \log x)^t \ (x \to \infty) \right) \text{ if }   t < 1.$$
\item The Riemann hypothesis is equivalent to 
$$p_x- \li^{-1}x = O( x^t) \ (x \to \infty) \text{ if (and only if) } t > \tfrac{1}{2}$$
and to
$$p_x- \li^{-1}x = O\left ( \sqrt{x}\,  (\log x)^{5/2}\right)\ (x \to \infty).$$
\item Conjecture \ref{eurekaconjecture}  is equivalent to
$$p_x- \li^{-1}x = O\left ( \sqrt{x}\,  (\log x)^{t}\right) \ (x \to \infty)  \text{ if (and only if) } t > \tfrac{1}{2}.$$
\item  Conjecture \ref{eurekaconjecture2}  is equivalent to
$$p_x- \li^{-1}x = O\left ( \sqrt{x \log x}\, (\log \log x)^{t}\right) \ (x \to \infty) \text{ if (and only if) } t > 0,$$
and it implies that $\degl_3(\li-\pi) \geq 1$ and that
$$p_x- \li^{-1}x = O\left ( \sqrt{x \log x}\, (\log \log  \log x)^{t} \right) \ (x \to \infty)$$ 
for all $t > \degl_3(\li-\pi)$ but for no $t< \dege=l_3(\li-\pi)$.
\item Montgomery's conjecture (\ref{MMC3}) implies that 
$$\dege(p_x- \li^{-1}x ) = (\tfrac{1}{2},\tfrac{1}{2}, 0,2,0,0,0,\ldots),$$
or, equivalently, that
$$p_x- \li^{-1}x = o\left( \sqrt{x \log x} \, (\log \log  \log x)^{2}r(x)\right) \ (x \to \infty)$$
for all $r \in \mathfrak{L}$ with $\lim_{x \to \infty} r(x) = \infty$ but for no $r \in \mathfrak{L}$ with $\lim_{x \to \infty} r(x) = 0$.
\item The anti-Riemann hypothesis is equivalent to $\deg(p_x- \li^{-1}x ) = 1$, 
to
$$p_x -\li^{-1} x = o(x^t) \ (x \to \infty) \text{ (if and) only if } t \geq 1,$$
and to
$$(1,-\infty,-1,0,0,0,\ldots) \leq \dege(p_x- \li^{-1}x ) \leq (1,-\infty,-\tfrac{3}{5},\tfrac{1}{5},0,0,0,\ldots).$$
Moreover, it implies that
$$\dege (p_x- \li^{-1}x )  = \dege(\li(p_x)-\pi(p_x))  =  \dege(\li-\pi).$$
\end{enumerate}
\end{corollary}

Note that
$$\Ri^{-1}x-\li^{-1}x \sim \sqrt{x \log x} \ (x \to \infty)$$
and thus
$$\dege (\Ri^{-1}x-\li^{-1}x) = (\tfrac{1}{2},\tfrac{1}{2},0,0,0,\ldots) < \dege(p_x- \li^{-1}x ),$$ 
whence
 $$\dege(p_x- \Ri^{-1}x ) = \dege(p_x- \li^{-1}x ).$$
 Thus, one can replace $\li^{-1}x$ with $\Ri^{-1}x$ in the results above.  Moreover,  the function
 $$\frac{p_x- \li^{-1}x }{\Ri^{-1}x-\li^{-1}x} =1+ \frac{p_x- \Ri^{-1}x }{\Ri^{-1}x-\li^{-1}x}$$
 is a natural normalization of the function $p_x- \li^{-1}x$.  Note,  then, that
 $$\dege\frac{p_x- \li^{-1}x }{\Ri^{-1}x-\li^{-1}x}  =  \dege \frac{\li(x)-\pi(x)}{\li(x)-\Ri(x)} + (0,\Theta-\tfrac{1}{2},0,0,0,\ldots),$$
 and thus the Riemann hypothesis or $\Theta_1 = -\infty$ holds if and only if
 $$\dege\frac{p_x- \li^{-1}x }{\Ri^{-1}x-\li^{-1}x}  =  \dege \frac{\li(x)-\pi(x)}{\li(x)-\Ri(x)},$$
 and the Riemann hypothesis holds if and only if either function has degree $0$.
 
The following result is in some ways more fine-tuned than  Theorem \ref{pnthm}, as it allows one to pass from known $O$ bounds of $\li(x)-\pi(x)$ to $O$ bounds of $p_n-\li^{-1}n$.

\begin{proposition}\label{rsprop}
Let $N \geq 1$, and let $r(n)$ and $s(n)$ be functions that are defined and positive for all integers $n > N$ such  that  $$-s(n)\leq \li(n)-\pi(n) \leq r(n)$$ for all $n > N$.
\begin{enumerate}
\item One has
\begin{align*}
-s(p_n) \log\li^{-1}n  <  p_n-\li^{-1}n < r(p_n) \log \li^{-1}(n+r(p_n))
\end{align*}
for all $n >  \pi(N)$.
\item  One has 
$$\inf_{n  >  \pi(N)} \frac{p_n - \li^{-1}n}{s(p_n) \log n} \geq -1.$$
Moreover, if $r(p_n) = O(n) \ (n \to \infty)$, or more generally if
 $\log  \left(1+\frac{r(p_n)}{n}\right) = o (\log n) \ (n \to \infty),$   then
$$\limsup_{n \to \infty} \frac{p_n - \li^{-1}n}{r(p_n)\log n} \leq 1,$$
so that, if also $r(n) \geq s(n)$ for all sufficiently large $n$, then
$$\limsup_{n \to \infty} \frac{|p_n - \li^{-1}n|}{r(p_n)\log n} \leq 1$$
and therefore
$$ p_n-\li^{-1}n = O(r(p_n) \log n) \ (n \to \infty).$$
\end{enumerate}
\end{proposition}

\begin{proof}
Note that $p_n > N$ holds if and only if $n > \pi(N)$.  Let $n > \pi(N)$.  Since then
$$\li(p_n)-n = \li(p_n)-\pi(p_n) \leq r(p_n),$$ one has
$$p_n \leq \li^{-1}(n+ r(p_n))$$
and therefore
$$p_n-\li^{-1}n \leq \li^{-1}(n+r(p_n)) -\li^{-1}n,$$
 Likewise, one has
$$p_n-\li^{-1}n \geq \li^{-1}(n-s(p_n)) -\li^{-1}n.$$
 Thus, we have
$$ -(\li^{-1}n -\li^{-1}(n-s(p_n))) \leq p_n-\li^{-1}n \leq \li^{-1}(n+r(p_n))-\li^{-1}n.$$
Now, since $\frac{d}{dn} \li^{-1}n =\log \li^{-1}n $ is positive and increasing, one has
\begin{align*}
p_n-\li^{-1}n & \leq \li^{-1}(n+r(p_n)) -\li^{-1}n \\ & =  \int_n^{n+r(p_n)} \log \li^{-1}(t) \, dt \\ & < r(p_n) \log \li^{-1}(n+r(p_n)).
\end{align*}
Likewise, one has
\begin{align*}
p_n-\li^{-1}n & \geq -(\li^{-1}n -\li^{-1}(n-s(p_n))) \\ & = - \int_{n-s(p_n)}^n \log \li^{-1}(t) \, dt \\ & > -s(p_n) \log \li^{-1}n.
\end{align*}
This proves statement (1).  To prove (2), note that
$$\frac{\log  (n+r(p_n))}{\log n}-1 =\frac{ \log  \left(1+\frac{r(p_n)}{n}\right)}{\log n} \to 0$$
as $n \to \infty$, so that
$$\log \li^{-1}(n+r(p_n)) \sim \log (n+r(p_n)) \sim \log n \ (n \to \infty).$$  Therefore the upper bound $ r(p_n) \log \li^{-1}(n+r(p_n))$ of $p_n-\li^{-1}n$ is asymptotic to $r(p_n) \log n$.   Likewise, the lower bound $ -s(p_n) \log \li^{-1}n$ of $p_n-\li^{-1}n$ is asymptotic to $-s(p_n) \log n$.  Statement (2) follows.
\end{proof}

\begin{corollary}
One has the following.
\begin{enumerate}
\item $p_x - \li^{-1}x =  O \left( x^\Theta( \log x)^{\Theta+1+\Theta_1+\varepsilon} \right) \ (x \to \infty)$ for all $\varepsilon > 0$.
\item  $p_x - \li^{-1}x =  O \left( x^\Theta( \log x)^{\Theta+2} \right) \ (x \to \infty)$.
\item  $p_x - \li^{-1}x = O\left(xe^{ - C(\log x)^{3/5}(\log \log x)^{-1/5}}\right)  \ (x \to \infty)$ for all positive $C  <0.2098$.   
\end{enumerate}
\end{corollary}

Table \ref{duality}  below describes an informal ``duality''  between the prime numbers and the zeros of the  Riemann zeta function that was touched upon in Section 5.3.   This duality is useful for translating problems about the primes to analogues for the zeros of $\zeta(s)$ and vice versa.  For example, the problem of estimating $\pi(x)$ (resp., the prime gaps $p_{n+1}-p_n$) translates to the problem of estimating $N(T)$ (resp., the gaps $\tau_{n+1}-\tau_n$), and vice versa.

\begin{table}[!htbp]
\footnotesize \caption{\centering Duality between the prime numbers and the zeros of the  Riemann zeta function}
\begin{tabular}{|l|l||l|} \hline
 $p_n$  & $\tau_n$ &  $n$   \\  \hline \hline
 $\pi_0(x)$ &  $N(2\pi x)$ & $\lfloor x \rfloor_0$   \\  \hline
   $\log x$ & $\frac{1}{\log x}$ & $1$  \\ \hline
  $\li(x) = \int_0^{n} \frac{dt}{\log t}$ & $x \log\frac{x}{e}+ \frac{7}{8} = \int_0^{x} \log t \, dt + \frac{7}{8} $  & $x = \int_0^x 1 \, dt$  \\ \hline
 $\Ri(x)$ & $1+\frac{1}{\pi}\theta(2 \pi x)$ & $x-\tfrac{1}{2}$ \\ \hline
 $\pi_0(x)-\Ri(x)$ & $S(2\pi x)$ &  $-\{x\}_0 +\tfrac{1}{2}$   \\  \hline
 $\pi_0(x) = \Ri(x) - \sum_{n=1}^\infty \sum_\rho  \frac{\mu(n)}{n}\Ei\left(\frac{\rho \log x}{n}\right)$ & $N(2\pi x) = 1+\frac{1}{\pi}\theta(2 \pi x)+S(2 \pi x)$ & $\lfloor x \rfloor_0 = x  -\frac{1}{2} +\sum_{n = 1}^\infty \frac{ \sin 2 \pi n x}{\pi n}$  \\ \hline
 $\li(p_n)$ & $\widehat{\gamma}_n =  \tau_n \log\frac{\tau_n}{e}+ \frac{11}{8}$ &  $n$   \\ \hline
$\Ri(p_n)$ & $1+\frac{1}{\pi}\theta(2\pi \tau_n)+\frac{1}{2}$ &  $n$  \\ \hline
$\li^{-1}x$ &   $r_x = \frac{x-7/8}{W((x-7/8)e)}$ & $x$  \\ \hline
$\pi_0(\li^{-1}x)$, $\vartheta_0(x)$, or $\psi_0(x)$ &   $N\left(2\pi\frac{x-7/8}{W((x-7/8)e)}\right)$ & $\lfloor x \rfloor_0$ \\ \hline
 $p_{x}$ &  $\tau_{\lceil x \rceil}$  & $\lceil x \rceil$  \\ \hline
 $p_{x}-\li^{-1}x$ &  $\tau_{\lceil x \rceil}-r_x$  & $\lceil x \rceil-x$  \\ \hline
 $p_{n+1}-p_n$ &  $\tau_{n+1}-\tau_n$ & $1$  \\ \hline
 $\li(p_{n+1})-\li(p_n)$ &  $\widehat{\gamma}_{n+1}-\widehat{\gamma}_n$ & $1$    \\ \hline
$\Ri(p_{n+1})-\Ri(p_n)$ &  $\frac{1}{\pi}\theta(2\pi \tau_{n+1})-\frac{1}{\pi}\theta(2\pi \tau_n)$  & $1$  \\ \hline
\end{tabular}\label{duality}
\end{table}

Just as the $\widehat{\gamma}_n \approx n$ and $1+\frac{1}{\pi}\theta(2\pi \tau_n) +\tfrac{1}{2} \approx n$ are natural normalizations of the $\gamma_n$, the $\li(p_n) \approx n$ and $\Ri(p_n) \approx n$ are natural normalizations of the $p_n$.  We note the following.

\begin{proposition}
One has
$$\li(p_n)-n \sim \frac{p_n-\li^{-1}n}{\log n} \ (n \to \infty)$$
and
$$\Ri(p_n)-n  \sim \frac{p_n-\Ri^{-1}n}{\log n} \ (n \to \infty).$$
Moreover, one has
\begin{align*}
\dege(\Ri(p_n)-n)  & = \dege (\li(p_n)-n ) \\
  & = \dege(p_n-\li^{-1}n)+ (0,-1,0,0,0,\ldots) \\
  & = \dege(\li-\pi)+(0,\Theta,0,0,0,\ldots)
\end{align*}
and
$$\dege (p_n-\Ri^{-1}n) = \dege(p_n-\li^{-1}n).$$ 
\end{proposition}

\begin{figure}[h!]
\includegraphics[width=70mm]{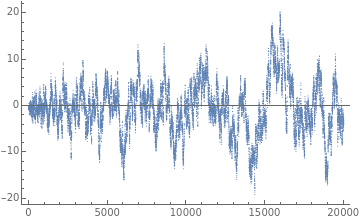}
    \caption{\centering Plot of $\Ri(p_n)-n$ for $n = 1,2,3,\ldots,20000$}
\label{ripnna}
\end{figure}

\begin{figure}[h!]
\includegraphics[width=70mm]{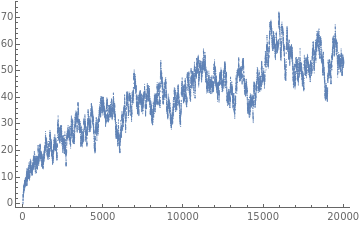}
    \caption{\centering Plot of $\li(p_n)-n$ for $n = 1,2,3,\ldots,20000$}
\label{lipnna}
\end{figure}

Figures \ref{ripnna} and \ref{lipnna} provide plots of $\Ri(p_n)-n = \Ri(p_n)-\pi(p_n)$ and $\li(p_n)-n = \li(p_n)-\pi(p_n)$, respectively, for $n = 1,2,3,\ldots,20000$.  Although $\Ri(p_n)$ appears to be much better approximated by $n$ than is $\li(p_n)$, both, with respect to the degree formalism, are equally well approximated by $n$ in the long run.

Recall that  $ p_x \leq y$ if and only if $x \leq \pi(y)$, for all $x, y \geq 0$.     Thus $p_{\li(x)} \leq y$ if and only if $\li(x) \leq \pi(y)$.  It is thus natural to consider the quantity $p_{\li(x)}\pi(x) - x \li(x)$.   Note that
$$p_{\li(x)}\pi(x) - x \li(x)  = E(x) + F(x),$$
where $$E(x) = p_{\li(x)}(\pi(x)-\li(x))$$ and
$$F(x) =  \li(x)(p_{\li (x)}-x),$$
and where $E(x) \geq 0$ if and only if $F(x) \leq 0$.  Moreover, one has the following.

\begin{proposition}\label{pppl}
One has 
$$\dege (p_{\li(x)}-x)= \dege(\li-\pi) + (0,1,0,0,0,\ldots).$$
It follows that
$$\dege (p_{\li(x)}(\pi(x)-\li(x))) = \dege(\li-\pi)+(1,0,0,0,\ldots) = \dege(\li(x)(p_{\li (x)}-x)),$$
and therefore
$$\dege(p_{\li(x)}\pi(x) - x \li(x)) \leq \dege(\li-\pi)+(1,0,0,0,\ldots).$$
Similar statements hold when one replaces $\li(x)$ by $\Ri(x)$.  
\end{proposition}

Thus, one has $\dege E = \dege F$.   It is conceivable, then, that $\dege(p_{\li(x)}\pi(x) - x \li(x))$ is strictly smaller than $\dege(\li-\pi)+(1,0,0,0,\ldots)$.   If true, that would show that the product $p_{\li(x)}\pi(x)$  is better approximated by $x \li(x)$ than its respective factors are by $x$ and $\li(x)$, separately.

Now,  for all $x> 0$, let
$$\bm{p}(x)  = \frac{p_x}{x}.$$
Note that $$\PP(p_x) = \frac{\pi(p_x)}{p_x} = \frac{\lceil x \rceil}{p_x} = \frac{\lceil x \rceil}{x}\frac{1}{\bm{p}(x)}$$
and therefore
$$\bm{p}(x) = \frac{\lceil x \rceil}{x}\frac{1}{\PP(p_x)}.$$
Note also that
$$\PP(x)-\frac{1}{\bm{p}(\li(x))} = \frac{1}{xp_{\li(x)}}(p_{\li(x)}\pi(x) - x \li(x))  = \frac{1}{x}(\pi(x)-\li(x))+ \frac{\li(x)}{xp_{\li (x)}}(p_{\li (x)}-x).$$
By Theorems \ref{pnthm} and \ref{fgie00}, one has the following.

\begin{corollary}
One has 
$$\dege  \left(\frac{1}{x}(\pi(x)-\li(x))\right) = \dege\left( \frac{\li(x)}{xp_{\li (x)}}(p_{\li (x)}-x)\right) = \dege(\li-\pi) + (-1,0,0,0,0,\ldots),$$
$$\dege \left(\PP(x)-\frac{1}{\bm{p}(\li(x))}  \right) \leq \dege\left(\PP(x)-\frac{\li(x)}{x}\right)  = \dege(\li-\pi) + (-1,0,0,0,0,\ldots),$$
and
$$\dege \left(\frac{1}{\PP(x)}-\bm{p}(\li(x)) \right) = \dege \left(\PP(x)-\frac{1}{\bm{p}(\li(x))}  \right).$$
Similar statements hold when one replaces $\li(x)$ by $\Ri(x)$.  
\end{corollary}

\begin{corollary}
One has
$$\dege\left(\pi(x)- \frac{x}{p_{\li(x)}}\li(x)\right) \leq \dege(\li-\pi)$$
and
$$\dege\left(\pi(x)- \frac{x}{p_{\Ri(x)}}\Ri(x)\right) \leq \dege(\li-\pi).$$
\end{corollary}

 Note, then, that $$\frac{x}{\bm{p}(\li(x))}= \frac{x}{p_{\li(x)}}\li(x)$$ is an excellent approximation for $\pi(x)$ (at least as good as $\li(x)$),  as is $$\frac{x}{\bm{p}(\Ri(x))} =  \frac{x}{p_{\Ri(x)}} \Ri(x)$$ (at least as good as $\Ri(x)$).  We do not know whether or not the inequalities in Proposition \ref{pppl} and its two corollaries  above are equalities.
 
 \begin{problem}
Express $\dege\left(\pi(x)- \frac{x}{p_{\li(x)}}\li(x)\right)$ and $\dege\left(\pi(x)- \frac{x}{p_{\Ri(x)}}\Ri(x)\right)$ in terms of $\dege(\li-\pi)$.  In particular, are they equal to each other, and is either equal to $\dege(\li-\pi)$?
 \end{problem}

\begin{figure}[ht!]
\includegraphics[width=80mm]{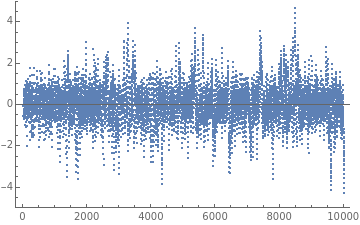}
\caption{Plot of $\pi(n)-\frac{n}{p_{\Ri(n)}} \Ri(n)$ for $n  = 1,2,3,\ldots, 10000$}
  \label{pima}
\end{figure}

\begin{figure}[ht!]
\includegraphics[width=80mm]{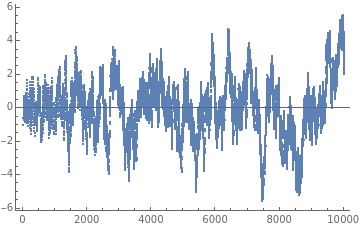}
\caption{Plot of $\pi(n)- \Ri(n)$ for $n  = 1,2,3,\ldots, 10000$}
  \label{pimb}
\end{figure}

\begin{figure}[ht!]
\includegraphics[width=80mm]{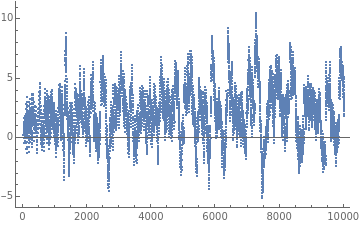}
\caption{Plot of $\pi(n)-\frac{n}{p_{\li(n)}} \li(n)$ for $n  = 2,3,4,\ldots, 10000$}
  \label{pimc}
\end{figure}

\begin{figure}[ht!]
\includegraphics[width=80mm]{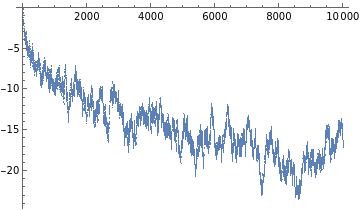}
\caption{Plot of $\pi(n)- \li(n)$ for $n  = 2,3,4,\ldots, 10000$}
  \label{pimd}
\end{figure}

See Figures \ref{pima}, \ref{pimb}, \ref{pimc}, and \ref{pimd} for plots of $\pi(n)-\frac{n}{p_{\Ri(n)}} \Ri(n)$,  $\pi(n)- \Ri(n)$, $\pi(n)-\frac{n}{p_{\li(n)}} \li(n)$, and $\pi(n)-\li(n)$, respectively.

\begin{remark}[Asymptotic continued fraction expansions of $p_x$]
Since $\frac{\li(x)}{x}$ has the asymptotic expansion (\ref{asex}), and since $\log \li^{-1}x = \Ei^{-1}x$,  the function $\frac{x}{\li^{-1}x}$ has the asymptotic expansion
\begin{align*}
\frac{x}{\li^{-1}x} \simeq \sum_{n = 0}^\infty {\frac {n!}{(\Ei^{-1}x)^{n+1}}} \ (x \to \infty)
\end{align*}
If follows that the asymptotic continued fraction expansions of $\PP(x)$ in Section 8.4
can be re-expressed as the following asymptotic continued fraction expansions of $\bm{p}(x) = \frac{p_x}{x}$.
\begin{enumerate}
\item $\bm{p}(x)  \,  \simeq \, \displaystyle  \Ei^{-1}x-  \ \frac{1\Ei^{-1}x} {\Ei^{-1}x+1 \, -} \ \frac{2\Ei^{-1}x} {\Ei^{-1}x+2 \, -} \   \frac{3\Ei^{-1}x} {\Ei^{-1}x+3\, -} \ \frac{4\Ei^{-1}x}{\Ei^{-1}x+4 \, -} \ \cdots \ (x \to \infty)$. 
\item $\bm{p}(x) \, \simeq \,  \cfrac{\Ei^{-1} x}{1 \,-} \ \cfrac{\frac{1}{\Ei^{-1} x}}{1 \,-}\  \cfrac{\frac{1}{\Ei^{-1} x}}{1 \,-}\  \cfrac{\frac{2}{\Ei^{-1} x}}{1 \,-}\  \cfrac{\frac{2}{\Ei^{-1} x}}{1 \,-}\  \cfrac{\frac{3}{\Ei^{-1} x}}{1 \,-} \  \cfrac{\frac{3}{\Ei^{-1} x}}{1 \,-}\  \cfrac{\frac{4}{\Ei^{-1} x}}{1 \,-}  \ \cfrac{\frac{4}{\Ei^{-1} x}}{1 \,-}  \ \cdots \ (x \to \infty)$.
\item $\bm{p}(x)  \,  \simeq \, \displaystyle \Ei^{-1}x -1 - \  \frac{1}{\Ei^{-1}x  - 3 \,-}\  \frac{4}{\Ei^{-1}x-5\,-}\  \frac{9}{\Ei^{-1}x - 7 \,-} \ \frac{16}{\Ei^{-1}x-9\,-} \ \cdots \ (x \to \infty)$.
\end{enumerate}
 Equivalently,  one has the following asymptotic continued fraction expansions of  $\bm{p}(\Ei x)$ at $\infty$.
\begin{enumerate}
\item $\bm{p}(\Ei x)  \,  \simeq \, \displaystyle  x-  \ \frac{1x} {x+1 \, -} \ \frac{2x} {x+2 \, -} \   \frac{3x} {x+3\, -} \ \frac{4x}{x+4 \, -}\ \cdots \ (x \to \infty)$. 
\item $\bm{p}(\Ei x) \, \simeq \,  \cfrac{x}{1 \,-} \ \cfrac{\frac{1}{x}}{1 \,-}\  \cfrac{\frac{1}{x}}{1 \,-}\  \cfrac{\frac{2}{x}}{1 \,-}\  \cfrac{\frac{2}{x}}{1 \,-}\  \cfrac{\frac{3}{x}}{1 \,-} \  \cfrac{\frac{3}{x}}{1 \,-}\  \cfrac{\frac{4}{x}}{1 \,-}  \ \cfrac{\frac{4}{x}}{1 \,-}  \ \cdots \ (x \to \infty)$.
\item $\bm{p}(\Ei x)  \,  \simeq \, \displaystyle x -1 - \  \frac{1}{x  - 3 \,-}\  \frac{4}{x-5\,-}\  \frac{9}{x - 7 \,-} \ \frac{16}{x-9\,-} \ \cdots \ (x \to \infty)$.
\end{enumerate}
Moreover,  the best rational approximations of the functions $\PP(e^x)$ and $\frac{1}{\bm{p}(\Ei x)}$ are one and the same.    The statements above hold also when $\Ei$ is replaced with $\ERi$.
\end{remark}

\begin{remark}[Largest order of an element in the symmetric group $S_n$]\label{landaufunction}
Note that the function $\li^{-1}x$ arises in connection with another equivalent of the Riemann hypothesis, concerning the largest order $g(n)$ of an element of the symmetric group $S_n$.  The function $g(n)$ is known as {\bf Landau's function} and is OEIS Sequence A000793.   Let $h(n) = \log g(n)-\sqrt{\li^{-1}n}$.  By \cite[Theorems 10.41 and 10.42]{broughan}, one has the following.
\begin{enumerate}
\item If $\Theta < 1$,  then $\dege h \leq (\frac{\Theta}{2},\frac{\Theta}{2}, 0,0,0,\ldots)$.
\item If $\Theta = 1$, then $\dege h \leq (\frac{1}{2},0,-\infty,-\frac{1}{2}, 0,0,0,\ldots)$.
\item If $\Theta > \frac{1}{2}$,  then $\dege h \geq (\frac{\Theta}{2},-\infty,-1,0,0,0,\ldots)$.
\item If $\zeta(s)$ has a zero with real part $\Theta$, then $\dege h = (\frac{\Theta}{2}, \frac{\Theta}{2},0,0,0,\ldots)$.
\item As a consequence of (1)--(4), one has $\deg h  = \frac{\Theta}{2}$.
\item The following conditions are equivalent to the Riemann hypothesis $\Theta = \frac{1}{2}$.
\begin{enumerate} 
\item $\deg h = \frac{1}{4}$.
\item $\dege h =  (\frac{1}{4},\frac{1}{4}, 0,0,0,\ldots)$.
\item $\log g(n)<\sqrt{\li^{-1}n}$ for all $n \gg 0$.
\end{enumerate}
\end{enumerate}
\end{remark}

\section{The $n$th prime gap $g_n = p_{n+1}-p_n$}

Far less is known about the {\bf prime gap function}\index{prime gap function $g_n$}\index[symbols]{.ru N@$g_n$}  $$g_n  = p_{n+1}-p_n$$ than is known about the prime listing function $p_n$.   It is easy to show that there are arbitrarily long sequences of consecutive composite numbers, or, equivalently, that
$$\limsup_{n \to \infty}  g_n = \infty.$$  The {\bf twin prime conjecture} states that there are infinitely many primes that differ by $2$, or, equivalently,
$$\liminf_{n \to \infty}  g_n = 2.$$ 
Until 2013, it was unknown whether or not $\liminf_{n \to \infty}  g_n$ is finite.   Until then, the best known bound on small prime gaps was one proved in 2007 by Goldson, Pintz, and Yildrim \cite{gold2}, namely, that
$$\liminf_{n \to \infty}  \frac{g_n}{\sqrt{\log n}\, (\log \log n)^2} < \infty$$
and therefore
$$\underline{\dege} \, g_n \leq (\tfrac{1}{2},2,0,0,0,\ldots).$$
 However, in 2013, Y.\ Zhang made a startling breakthrough towards settling the twin prime conjecture by proving that there are infinitely many primes that differ by at most 70000000 \cite{zhang}.  In the immediate aftermath of Zhang's proof, several mathematicians teamed up to improve Zhang's result, and, ultimately, by April 2014, proved that   there are infinitely many primes that differ by at most 246 \cite{poly}.  Equivalently, they proved that $$\liminf_{n \to \infty}  g_n  \leq 246.$$   As a consequence, one has
 $$\underline{\dege} \, g_n = (0,0,0,\ldots).$$  Here, we consider the following  problem.

\begin{outstandingproblem}
Compute $\dege g_n$.
\end{outstandingproblem}

The  {\bf Bertrand--Chebyshev theorem},\index{Bertrand--Chebyshev theorem} conjectured by Bertrand in 1845 and proved by Chebyshev in 1850 \cite{cheb1}, states that, for every integer $n >1$, there is at least one prime $p$ such that $n < p < 2n$.  Equivalently, one has
$$p_{n+1} < 2 p_n,$$
that is,
$$g_n < p_n,$$
for all positive integers $n$.   Many improvements on this inequality have since been proved, e.g.,  \cite[Proposition 5.4 and Corollary 5.5]{dus}.   One can do much better if one is concerned only with $O$ bounds rather than with explicit inequalities.  By (\ref{bakerres}), which is proved in \cite{baker},  and by Corollary \ref{simplepgcor} below, one has
$$g_n = O (p_n^{21/40}) \ (n \to \infty).$$
Since the prime number theorem is equivalent to 
$$p_n \sim n \log n  \ (n \to \infty),$$ the $O$ bound above is equivalent to
$$g_n = O ((n \log n)^{21/40}) \ (n \to \infty).$$
In \cite{cramer1}, Cram\'er proved on condition of the Riemann hypothesis that
$$g_n = O (p_n^{1/2}\log p_n) \ (n \to \infty),$$
or, equivalently,
$$g_n = O (n^{1/2}(\log n)^{3/2}) \ (n \to \infty).$$
In particular, one has
$$\dege g_n \leq (\tfrac{21}{40},\tfrac{21}{40},0,0,0,\ldots),$$
and one has
$$\dege g_n \leq (\tfrac{1}{2},\tfrac{3}{2},0,0,0,\ldots)$$
on condition of the Riemann hypothesis.  

Since $$p_{n+1} = 2+\sum_{k = 1}^n g_k, \quad \forall n \geq 1,$$
the prime number theorem implies that $$g_n \neq o(\log n) \ (n \to \infty).$$  Thus, an obvious lower bound for $\dege g_n$ is $$\dege g_n  \geq (0,1,0,0,0,\ldots).$$ The state-of-the-art lower bound for $g_n$ at the writing of this text improves the lower bound above by nearly a factor of $\log \log n$: in 2016, Ford, Green, Konyagin, Maynard, and Tao  proved \cite{fgktm} that
\begin{align}\label{fgkmt}
g_n \neq o\left(\frac{\log n \log \log n\log \log \log \log n}{\log \log \log n}\right) \ (n \to \infty).
\end{align}  This result, of course, implies that $$\dege g_n \geq (0,1,1,-1,1,0,0,0,\ldots).$$
Note that there is a long history that led up to the 2016 result (\ref{fgkmt}),  e.g.,  in 1931,  E.\ Westzynthius was the first to prove that
$g_n  \neq O(\log n) \ (n \to \infty)$, via a more precise result that implies $\dege g_n \geq (0,1,0,1,-1,0,0,0,\ldots)$ \cite{west}.

By the following simple proposition, which relates the prime gap function to the prime counting function, many of the results of Section 12.1 apply to the study of the function $g_n$.

\begin{proposition}\label{simplepg}
Let $h(x)$ be a nonnegative real function defined on $[p_N, \infty)$,  where $N$ is a positive integer.  If 
$$\pi(x+h(x)) -\pi(x) \geq 1, \quad \forall x \geq p_N,$$
then
$$h(p_n) \geq g_n, \quad \forall n \geq N.$$
Moreover, the converse holds if $h(x)$ is nondecreasing on $[p_N, \infty)$. 
\end{proposition}

\begin{proof}
Suppose that
$$\pi(x+h(x)) -\pi(x) \geq 1, \quad \forall x \geq p_N.$$
Then one has
$$\pi(p_n+h(p_n))-\pi(p_n) \geq 1,$$
and therefore
$$p_{n+1} \leq  p_n+h(p_n),$$
whence
$$g_n \leq h(p_n),$$
for all $n$ with $p_n \geq p_N$, that is, with $n \geq N$.  Conversely,  suppose that $h(x)$ is nondecreasing for all $x \geq p_N$ and 
$$g_n \leq h(p_n), \quad \forall n \geq N.$$
Then one has
$$\pi(p_n+h(p_n))-\pi(p_n) \geq 1$$
for all $n \geq N$.    Let $x \geq  p_N$.  Let $n = \pi(x)$, so that  $n \geq N$ and $x \geq p_n \geq p_N$, and therefore also
$$\pi(x+h(x)) -\pi(x) \geq \pi(p_n+h(p_n))-\pi(p_n) \geq 1.$$
The proposition follows.
\end{proof}

\begin{corollary}\label{simplepgcor}
One has
\begin{align*}
\deg g_n =  \inf\{t > 0: \forall x \gg 0 \ (\pi(x+x^t) -\pi(x) \geq 1 )\}.
\end{align*}
\end{corollary}

Corollaries \ref{simplepgcor} and \ref{lindel} yield the following.

\begin{corollary}\label{lindelb}
If the density hypothesis holds, then $\deg g_n \leq \tfrac{1}{2}$.
\end{corollary}

Thus, the density hypothesis yields the same upper bound of $\deg g_n$, namely, $\frac{1}{2}$, as does the Riemann hypothesis.  

The following conjecture was first made by Piltz in 1884 and is implied by many well-known conjectures regarding the prime gap function $g_n$.

\begin{conjecture}[{Piltz's conjecture \cite{piltz}}]\label{gapconj}
One has $g_n = O(n^t) \ (n \to \infty)$ for all $t > 0$.  Equivalently, one has $\deg g_n = 0$.   Equivalently still, for all $t > 0$ one has $\pi(x+x^t) -\pi(x) \geq 1$ for all $x \gg 0$.
\end{conjecture}

Now, let $G(x)$ denote the   {\bf maximal prime gap function}\index{maximal prime gap function  $G(x)$}\index[symbols]{.st P@$G(x)$}
$$G(x) = \max_{p_{k} \leq x}  g_k = \max_{k \leq \pi(x)}  g_k, \quad \forall x \geq 0,$$
\cite{gran} \cite{mayn} \cite{pintz}.     Note that, in the literature, one encounters the function $G(x)$  and slight variants thereof,  such as  the function
$$\widetilde{G}(x) = \max_{p_{k+1} \leq x}  g_k = \max_{k < \pi(x)}  g_k, \quad \forall x \geq 0,$$
 \cite{fgktm} and the function $\max_{p_{k+1} < x}  g_k$ \cite{wolf1},  but such variations do not affect the various  asymptotics proved and conjectured for these functions.  The function $G(x)$ is the most natural of these variants.  It satisfies $$G(p_{n}) = \max_{k \leq n} g_k$$ for all positive integers $n$, and it is a nondecreasing step function with jumps precisely at those primes $p_{n}$ such that ($g_n = G(p_{n})$ and) $g_k < g_n$ for all $k < n$.   By Proposition \ref{supprop},  one has
$$\dege g_n = \dege G(p_n),$$ and, dually, 
$$\dege G(x) = \dege g_{\pi(x)}.$$
Also, $G(x)$ is the smallest nondecreasing  function that is greater than or equal to $g_{\pi(x)}$.  Likewise, by Proposition \ref{simplepg}, it is the smallest nonnegative nondecreasing  function  $h(x)$ on $\RR_{\geq 0}$ such that $\pi(x+h(x)) - \pi(x) \geq 1$ for all $x \geq 2$.  Indeed,  one has the following.

\begin{corollary}\label{GmaxCor} If $h(x)$ is any function that is nonnegative and nondecreasing  on $\RR_{\geq 0}$,  then one has
$$\pi(x+h(x)) -\pi(x) \geq 1, \quad \forall x \geq 2,$$
if and only if $$h(x) \geq G(x), \quad \forall x \geq 0.$$  Likewise,  if $h(x)$ is any function that is  nonnegative and nondecreasing  on some neighborhood of $\infty$, then one has
$$\pi(x+h(x)) -\pi(x) \geq 1, \quad \forall x \gg 0,$$
if and only if $$h(x) \geq G(x), \quad \forall x \gg 0.$$  
\end{corollary}

Thus,  the function $G(x)$ acts as a very natual  bridge between the prime gap function $g_n$ and the prime counting function $\pi(x)$ and bears directly on Problem \ref{Gmaxprob}.
Note that the  2016 result (\ref{fgkmt}) in \cite{fgktm} is proved there in the following stronger form.

\begin{theorem}[{\cite{fgktm}}] One has
$$G(x) \gg \frac{\log x \log \log x\log \log \log \log x}{\log \log \log x} \ (x \to \infty).$$  
Equivalently,  if $h(x)$ is any function that is  nonnegative and nondecreasing  on some neighborhood of $\infty$ and satisfies
$$\pi(x+h(x)) -\pi(x) \geq 1, \quad \forall x \gg 0,$$
then one has
$$h(x) \gg  \frac{\log x \log \log x\log \log \log \log x}{\log \log \log x} \ (x \to \infty).$$  
\end{theorem}

Note that the lower bound 
$$G(x) > (1-\varepsilon) \log x, \quad \forall x \gg 0,$$
for any fixed $\varepsilon > 0$, follows from the prime number theorem.   This is because, if $G(x) \leq (1-\varepsilon ) \log x$ for some $x$, then
$$p_{\pi(x)+1} -2 = g_1 + g_2 + \cdots g_{\pi(x)} \leq (1-\varepsilon) \pi(x)\log x,$$
and, as $x \to \infty$,  one has $p_{\pi(x)+1} -2 \sim \pi(x)\log x$, so the two inequalities must fail for all $x \gg 0$.   In a sense, this is the sharpest lower bound that can be proved from the prime number theorem alone.   Note also that, again by the prime number theorem, the average prime gap satisfies
$$\frac{g_1 + g_2 + \cdots+ g_{n}}{n} = \frac{p_{n+1}-2}{n} \sim \log n \sim \log p_n \ (n \to \infty).$$

The increasing sequence $a(n)$ of all $N$ for which $g_k < g_N$ for all $k < N$ is OEIS Sequence  A005669 and 
is given by $$1, 2, 4, 9, 24, 30, 99, 154, 189, 217, 1183, 1831, \ldots.$$
The corresponding sequence 
$$2, 3, 7, 23, 89, 113, 523, 887, 1129, 1327, 9551, 15683, \ldots$$
of primes $p_{a(n)}$ is OEIS Sequence A002386,
and the corresponding sequence
$$1, 2, 4, 6, 8, 14, 18, 20, 22, 34, 36, 44, \ldots $$
of record prime gaps $g_{a(n)} = G(p_{a(n)})$ is OEIS Sequence A005250.   The primes $p_{a(n)}$ are precisely the points of discontinuity of the step function $G$,  and the respective values of $G$ at these points are the gaps $g_{a(n)}$.    In particular, the image of $G$ is equal to $ \{g_{a(n)}: n \geq 1\}\cup \{0\}$.
From $(p_{a(n)})^-$ to $p_{a(n)}$, the function $G$ jumps up by $g_{a(n)}-g_{a(n-1)}$, where the sequence $g_{a(n)}-g_{a(n-1)}$ for $n \geq 2$ is OEIS Sequence A053695 and is given by  $$1, 2, 2, 2, 6, 4, 2, 2, 12, 2, 8, 8,\ldots.$$

The following theorem summarizes what we can prove about $\dege g_n$ and $\dege G$ from the results noted above.

\begin{theorem}\label{GTM}   One has the following.
\begin{enumerate}
\item One has
\begin{align*}
\dege G & = \dege g_{\pi(x)} \\ 
& =  \dege g_n +(0,-\deg g_n, 0,0,0,\ldots) \\
 & = \inf\{\dege r: r \in \mathbb{L}_{>0}, \, \forall x \gg 0 \ (\pi(x+r(x)) -\pi(x) \geq 1 )\}.
\end{align*}
\item $(0,1,1,-1,1,0,0,0,\ldots) \leq \dege g_n \leq (\tfrac{21}{40},\tfrac{21}{40},0,0,0,\ldots).$
\item $(0,1,1,-1,1,0,0,0,\ldots) \leq \dege G \leq (\tfrac{21}{40},0,0,0,0,\ldots).$
\item If the Riemann hypothesis holds, then $\dege g_n \leq (\tfrac{1}{2},\tfrac{3}{2},0,0,0,\ldots)$.
\item If the Riemann hypothesis holds, then $\dege G \leq (\tfrac{1}{2},1,0,0,0,\ldots)$.
\item If the density hypothesis holds, then $\deg g_n = \deg G \leq \tfrac{1}{2}$.
\item Conjecture \ref{gapconj}, namely, that $\deg g_n = 0$, is equivalent to $\deg G = 0$, and it implies that $\dege g_n = \dege G$ and 
$$\dege_1 g_n = \dege_1 G = \inf\{t \in \RR: \forall x \gg 0\, ( \pi(x+(\log x)^t)-\pi(x) \geq 1 )\},$$
and therefore also that
$$\liminf_{x \to \infty}\,  (\pi(x+(\log x)^{t})-\pi(x)) = 
\left.
 \begin{cases}
   \geq 1 & \text{if } t> \dege_1 g_n \\
    0 & \text{if } t <\dege_1 g_n.
 \end{cases}
\right.$$ 
\end{enumerate}
\end{theorem}

To prove the theorem, we first establish a few lemmas.

\begin{lemma}
One has
$$\dege G = \dege g_{\pi(x)} = \dege g_n + (0,-\deg g_n, 0,0,0,\ldots).$$
\end{lemma}

\begin{proof}
By Proposition \ref{supprop}, one has
$$\dege  G =  \dege g_{\pi(x)}.$$
Let $r \in \mathbb{L}_{>0}$ with $g_n \leq r(n)$ for all $n \gg 0$.  Then one has $g_{\pi(x)} \leq r(\pi(x))$ for all $x \gg 0$,  so that
$$\dege g_{\pi(x)} \leq \dege r(\pi(x)) = \dege r + (0,-\deg r, 0,0,0,\ldots).$$
Taking the infimum over all $r$ as chosen, we find that
$$\dege g_{\pi(x)} \leq  \dege g_n + (0,-\deg g_n, 0,0,0,\ldots).$$
(To verify this,  there are two cases to consider: $\dege g_n < (\deg g_n, \infty, 1,0,0,0,\ldots)$ and $\dege g_n = (\deg g_n, \infty, 1,0,0,0,\ldots)$.)  Conversely, let $s \in \mathbb{L}_{>0}$ with $g_{\pi(x)} \leq s(x)$ for all $x \gg 0$.    Then one has $g_n = g_{\pi(p_n)} \leq s(p_n)$ for all $n \gg 0$,  so that
$$\dege g_{n} \leq \dege s(p_n) = \dege s + (0,\deg s, 0,0,0,\ldots).$$  
Taking the infimum over all $s$ as chosen, we find that
$$\dege g_{n} \leq  \dege g_{\pi(x)} + (0,\deg g_n, 0,0,0,\ldots).$$
This completes the proof.
\end{proof}

Replacing $\pi(x)$ with $\pi(x)-1$ in the proof above yields the following lemma.

\begin{lemma}\label{degeGcor0}
One has
$$\dege  \widetilde{G} = \dege g_{\pi(x)-1} = \dege g_n + (0,-\deg g_n, 0,0,0,\ldots) = \dege G.$$
\end{lemma}

\begin{lemma}\label{degeGcor}
One has
$$\dege g_{\pi(x)} = \inf\{\dege r: r \in \mathbb{L}_{>0}, \, \forall x \gg 0 \ (\pi(x+r(x)) -\pi(x) \geq 1 )\}.$$
\end{lemma}

\begin{proof}
Let $r \in \mathbb{L}_{>0}$.   By  Proposition \ref{simplepg},  one has
$$\pi(x+r(x)) -\pi(x) \geq 1, \quad \forall x \gg 0,$$
if and only if
$$g_n \leq r(p_n), \quad \forall n \gg 0,$$
if and only if
$$g_{\pi(x)} \leq r(p_{\pi(x)}), \quad \forall x \gg 0,$$
if and only if 
$$g_{\pi(x)} \leq r(x), \quad \forall x \gg 0,$$
since each of the four conditions implies that $r$ is eventually increasing without bound, and since the latter three conditions are trivially equivalent.   It follows that
\begin{align*}
 \dege g_{\pi(x)} & = \inf\{\dege r: r \in \mathbb{L}_{>0}, \, \forall x \gg 0 \ (g_{\pi(x)} \leq r(x)) \} \\
 & =  \inf\{\dege r: r \in \mathbb{L}_{>0}, \, \forall x \gg 0 \ (\pi(x+r(x)) -\pi(x) \geq 1 )\}.
\end{align*}
This completes the proof.
\end{proof}

Now, we prove Theorem \ref{GTM}.

\begin{proof}[Proof of Theorem \ref{GTM}]
Statement (1) follows from Lemmas \ref{degeGcor0} and \ref{degeGcor}.
We have already verified statements (2) and (4).   Finally, statements (3) and (5) follow immediately from statements (1), (2), and (4), and statements (6) and (7) follow from statement (1),  Lemma \ref{degeGcor}, and Corollary \ref{lindelb}.
\end{proof}

Now, let $$g(x) = x-p_{\pi(x)},  \quad \forall x \geq 0.$$  Note that $p_{\pi(x)}$ is the largest prime less than or equal to $x$ (or $0$ if $x < 2$),  and thus $g(x) \geq 0$ is the distance from $x$ to the prime closest to $x$ from below.    It follows that $g(x)$ increases linearly (with slope $1$) from $0$ to $g_n$  on $[p_n,p_{n+1})$,  that is,  $g(x) = x-p_n$ on $[p_n,p_{n+1})$.  Thus, one has
$$g(p_{n+1}^-) = g_n$$
for all positive integers $n$ and
$$g(x) \leq g_{\pi(x)}$$
for all $x \geq 2$.

\begin{proposition}
One has
$$\dege g(x) = \dege G.$$
\end{proposition}

\begin{proof}
Clearly one has
$$\dege g(x) \leq \dege g_{\pi(x)} = \dege G.$$
It remains only to show that $\dege g_{\pi(x)} \leq \dege g$.
Let $r \in \mathbb{L}_{>0}$ with $g(x) \leq r(x)$ for all $x \gg 0$.
Since
$$x-p_n = g(x) \leq r(x)$$
for all $x \in [p_n,p_{n+1})$ and $r$ is eventually continuous and increasing, one has
$$g_n = p_{n+1}-p_n \leq r(p_{n+1}),$$
for all $n \gg 0$.  Therefore, since $p_{\pi(x)} \leq x$ and thus $p_{\pi(x)+1}< 2x$ by the Bertrand--Chebyshev theorem, one has
$$g_{\pi(x)} \leq r(p_{\pi(x)+1}) < r(2x),$$
for all $x \gg 0$.   It follows that
$$\dege g_{\pi(x)} \leq \dege r(2x) = \dege r.$$
Taking the infimum over all $r$ as chosen, we find that $$\dege g_{\pi(x)} \leq \dege g(x),$$ as desired.
\end{proof}

Note that
$$\dege (x-p_{\li(x)})= \dege(\li-\pi) + (0,1,0,0,0,\ldots),$$
while
$$\dege (x-p_{\pi(x)})= \dege g = \dege G.$$
Note also that, since $\pi(p_x) = \lceil x \rceil$ for all $x \geq 0$, one has
$$\dege (x-p_{\pi(x)}) = \dege (\pi(p_x)-p_{\pi(x)}).$$

\begin{remark}[An alternatve maximal gap function]
Let 
$$\widetilde{g}_x =  \max_{k < x}  g_k, \quad \forall x \geq 0,$$
where $\widetilde{g}_x  = 0$ if $x \leq 1$.    Then one has
$$\widetilde{G}(x) = \widetilde{g}_{\pi(x)}, \quad \forall x \geq 0,$$
and, by Proposition \ref{supprop}, one has
$$\dege \widetilde{g}_x  = \dege g_n.$$
\end{remark}

\begin{remark}[Alternative expressions for $\dege g_n$]
Note that
$$\log \frac{p_{n+1}}{p_n} = \log \left(1+ \frac{g_n}{p_n} \right) \sim \frac{g_n}{p_n}  \ (n \to \infty),$$
so that
$$\dege \left(\log \frac{p_{n+1}}{p_n}\right) = \dege g_n + (-1,-1,0,0,0,\ldots).$$  Similarly, one has
$$\frac{\log p_{n+1}}{\log p_n}-1 = \frac{ \log \frac{p_{n+1}}{p_n}}{ \log p_n} \sim \frac{g_n}{p_n \log p_n}  \ (n \to \infty),$$
so that
$$\dege \left(\frac{\log p_{n+1}}{\log p_n}-1\right) = \dege g_n + (-1,-2,0,0,0,\ldots).$$ 
Thus, the results in this section can be recast in terms of the functions $\log \frac{p_{n+1}}{p_n}$ and $\frac{\log p_{n+1}}{\log p_n}-1$. 
\end{remark}

The results in this section can be recast in terms of the gaps $ \widehat{p}_{n+1}- \widehat{p}_n$ between the {\bf  normalized primes}\index{normalized prime $\widehat{p}_n$} $$\widehat{p}_n = \Ri(p_n) \approx n\index[symbols]{.rt Pb@$\widehat{p}_n$}$$
discussed  briefly in the previous section. 
Indeed, by the mean value theorem and the fact that
$$\Ri'(x) \sim \frac{1}{\log x} = \li'(x) \ (x \to \infty),$$ one has
$$\Ri(p_{n+1})-\Ri(p_n) \sim \li(p_{n+1})-\li(p_n) \sim \frac{g_n}{\log p_n} \ (n \to \infty),$$ so that
$$ \dege( \Ri(p_{n+1})-\Ri(p_n)) =\dege( \li(p_{n+1})-\li(p_n)) = \dege g_n+(0,-1,0,0,0,\ldots).$$  

Note that the $n$th average prime gap is given by
$$\frac{1}{n}\sum_{k = 1}^n g_k = \frac{p_{n+1}-2}{n} \sim \log n \sim \log p_n \ (n \to \infty).$$
 The quantity 
$$\widehat{g}_n = \frac{g_n}{\log p_n}\index[symbols]{.ru N1@$\widehat{g}_n$}$$
is called the {\bf $n$th normalized prime gap},\index{normalized prime gap $\widehat{g}_n$} or the {\bf merit}\index{merit of the prime gap $g_n$} of the prime gap $g_n$.   The following well-known conjecture is an analogue of the Montgomery pair correlation conjecture for the normalized prime gaps.

\begin{conjecture}[{\cite{cramer} \cite{sound}}]\label{poiss}
For all $0 \leq a \leq b$, one has
$$\frac{\# \left\{n \leq T: \widehat{g}_n \in (a,b] \right \} }{T} \sim  \int_a^b e^{-t} \,  dt  \  (T \to \infty).$$
\end{conjecture}

By \cite[Exercise 3]{sound}, Conjecture \ref{poiss} follows from the following conjecture.

\begin{conjecture}[{\cite{sound}}]\label{poissgen}
For all $\lambda > 0$ and all nonnegative integers $k$, one has
$$\frac{\# \{n \leq T: \pi(n+\lambda \log T)-\pi(n) = k \} }{T} \sim  \frac{\lambda^k}{k!} e^{-\lambda}  \  (T \to \infty).$$
\end{conjecture}

In 1976,   P.\ X.\ Gallagher proved that Conjecture \ref{poissgen} follows from a certain uniform version of the Hardy--Littlewood conjecture for prime $k$-tuples \cite{gall} \cite{sound}.

One can recast Conjecture \ref{poiss}  in terms of the gaps between the normalized primes $\widehat{p}_n = \Ri(p_n)$ as follows.

\begin{conjecture}\label{poiss2}
For all $0 \leq a \leq b$, one has 
$$\frac{\# \left\{n \leq T:   \widehat{p}_{n+1}- \widehat{p}_n \in (a,b] \right \} }{T} \sim  \int_a^b e^{-t} \,  dt  \  (T \to \infty).$$
Equivalently, one has
$$\frac{\# \left\{n \in \ZZ_{>0}:   \widehat{p}_n  \leq T, \, \widehat{p}_{n+1}- \widehat{p}_n \in (a,b] \right \} }{\# \left\{n \in \ZZ_{>0} : \widehat{p}_n \leq T \right\} } \sim  \int_a^b e^{-t} \,  dt  \  (T \to \infty),$$
where also
$$\# \left\{n \in \ZZ_{>0} : \widehat{p}_n \leq T \right\} = \pi(\Ri^{-1}(T)) \sim T \ (T \to \infty).$$
\end{conjecture}

Note that Conjecture \ref{poissgen} implies both Conjectures \ref{poiss} and \ref{poiss2}, which are likely equivalent.  Moreover, Conjecture \ref{poissgen} implies that, for all nonnegative integers $k$ and all $0 \leq a \leq b$, one has 
$$\frac{\# \left\{n \leq T:   \widehat{p}_{n+k}- \widehat{p}_n \in (a,b] \right \} }{T} \sim  \int_a^b \frac{t^k}{k!} e^{-t} \,  dt  \  (T \to \infty).$$  The distribution $ \frac{t^k}{k!} e^{-t}$ is the {\bf gamma distribution with shape parameter $k+1$ and rate parameter $1$}, which has mean $k+1$ and variance $k+1$ and converges to the normal distribution with mean $k+1$ and variance $k+1$ as $k \to \infty$.   Section 8.6 provides an application of the gamma distribution $ \frac{t^k}{k!} e^{-t}$ to the prime counting function.

Further conjectures regarding the prime gaps are discussed in Section 14.6.

\chapter{ Diophantine approximation and continued fractions}

{\it Diophantine approximation} is the study of the approximation of real numbers by rational numbers. In this chapter, we apply degree and logexponential degree to the theories of Diophantine approximation and continued fractions.    We assume some familiarity with the theory of regular continued fractions (e.g., as contained in \cite[Chapter 9]{leve}), but not with  Diophantine approximation more broadly.   

\section{The basics of Diophantine approximation}

In this chapter, we use the  notations $[a_0, a_1,a_2,\ldots,a_n]$  and $[a_0, a_1,a_2,\ldots]$ for the finite and infinite regular continued fractions
\begin{align*}
[a_0, a_1,a_2,\ldots,a_n] = a_0 +\cfrac{1}{a_1+\cfrac{1}{a_2+\cdots +\cfrac{1}{a_{n-1}+\cfrac{1}{a_{n}}}}} \in \QQ\index[symbols]{.t  Ha@$[a_0, a_1,a_2,\ldots,a_n]$}
\end{align*}
and
\begin{align*}
[a_0, a_1,a_2,\ldots] = a_0 +\cfrac{1}{a_1+\cfrac{1}{a_2+\cfrac{1}{a_3+\cdots}}} = \lim_{n \to \infty} [a_0, a_1,a_2,\ldots,a_n] \in \RR\backslash \QQ,\index[symbols]{.t  Hb@$[a_0, a_1,a_2,\ldots]$}
\end{align*}
respectively, where $a_0 \in \ZZ$ and $a_k \in \ZZ_{>0}$ for all $k> 0$. 
For each nonnegative integer $k$, the integer $a_k$ is called the {\bf $k$th term}, or {\bf  $k$th partial quotient}, of the finite or infinite regular continued fraction above, and its {\bf $k$th convergent} $[a_0, a_1,a_2,\ldots,a_k]$  is written $\frac{p_k}{q_k}$ for unique relatively prime integers $p_k$ and $q_k$ with $q_k > 0$, where $k \leq n$ in the finite case.    The integers $p_k$ and $q_k$ satisfy identical recurrence relations of order $2$, namely,
$$p_k = a_k p_{k-1}+ p_{k-2}$$
and
$$q_k = a_k q_{k-1}+ q_{k-2}$$
for all $k \geq 1$,  with initial conditions $p_{-1} = 1$, $p_0 = a_0$, and $q_{-1} = 0$, $q_0 = 1$, respectively.

  For any  $\alpha \in \RR\backslash \QQ$, there are unique integers $a_k = a_k(\alpha)$\index[symbols]{.u  a@$a_k(\alpha)$} with $a_k > 0$ for all $k > 0$ so that $\alpha = [a_0,a_1,a_2,\ldots]$.   The continued fraction $[a_0,a_1,a_2,\ldots]$ is called the {\bf regular continued fraction  (expansion)  of $\alpha$}.\index{regular continued fraction expansion} As with the terms $a_k$,  the $p_k= p_k(\alpha)$\index[symbols]{.u  b@$p_k(\alpha)$} and $q_k= q_k(\alpha)$\index[symbols]{.u  c@$q_k(\alpha)$} are uniquely determined by $\alpha$.    Any  $\alpha \in \QQ$, on the other hand, has exactly two finite regular continued fraction expansions,  where one has
$$\alpha = [a_0, a_1,a_2,\ldots,a_n] = [a_0, a_1,a_2,\ldots,a_{n-1},a_n-1,1]$$
 if $a_n > 1$ or $n = 0$.   Nevertheless,  the $p_k= p_k(\alpha)$\index[symbols]{.u  b@$p_k(\alpha)$} and $q_k= q_k(\alpha)$\index[symbols]{.u  c@$q_k(\alpha)$} are both uniquely determined by $\alpha \in \QQ$,  where,  for $\alpha = [a_0, a_1,a_2,\ldots,a_n]$, we let $p_k(\alpha) = p_n(\alpha)$ and $q_k(\alpha) = q_n(\alpha)$, so that $\alpha = \frac{p_k(\alpha)}{q_k(\alpha)}$, for all $k \geq n$.  

 Let  $$S: (-\infty,\infty]  \longrightarrow (1,\infty] $$ be given by
$$S(\alpha)=   \left.
 \begin{cases}
   \frac{1}{\alpha-\lfloor \alpha \rfloor} & \text{if } \alpha \notin \ZZ\cup \{ \infty\}   \\
   \infty  & \text{if } \alpha \in \ZZ\cup \{ \infty\}.\index[symbols]{.u  d@$S(\alpha)$}
 \end{cases}
\right.$$
We call the function $S$ the {\bf regular continued fraction operator}.\index{regular continued fraction operator $S$}
If we write
$$[a_0, a_1,a_2,\ldots,a_n] =[a_0, a_1,a_2,\ldots,a_n,\infty,\infty,\infty,\ldots],$$   
then the regular continued fraction expansion of any $\alpha \in \RR$ is given by 
\begin{align}\label{regcontfrac}
\alpha =  [\lfloor  \alpha \rfloor, \lfloor S(\alpha) \rfloor,    \lfloor S(S(\alpha)) \rfloor,    \lfloor S(S(S(\alpha)))\rfloor, \ldots],
\end{align}
where we set $\lfloor \infty \rfloor  = \infty$.    Thus, $\alpha$ is rational if and only if $S^{\circ n}(\alpha) = \infty$ for  some nonneagtive  integer $n$,  if and only if the regular continued fraction of $\alpha$  terminates in a tail of $\infty$s, and in fact (\ref{regcontfrac}) always yields the shorter of the two  regular continued fraction expansions of $\alpha$ in that case.   With these conventions, one has
$$a_n(\alpha) = \lfloor S^{\circ n}(\alpha) \rfloor$$
and 
$$ S^{\circ n}(\alpha) = [a_n(\alpha), a_{n+1}(\alpha), a_{n+2}(\alpha),\ldots]$$
for all $\alpha \in \RR$ and all nonnegative integers $n$.    Moreover, by \cite[Theorem 9.5]{leve}, one has
\begin{align}\label{cfineq00}
S^{\circ (n+1)}(\alpha) = -\frac{q_{n-1}(\alpha)\alpha-p_{n-1}(\alpha)}{q_{n}(\alpha)\alpha-p_{n}(\alpha)}
\end{align}
for all nonnegative integers $n$ for which $S^{\circ (n+1)}(\alpha) \neq \infty$, while also $S^{\circ (n+1)}(\alpha) = \infty$ if and only if $\alpha = \frac{p_n(\alpha)}{q_n(\alpha)}$.

Let $\alpha \in \RR$ and $n$ a positive integer with $S^{\circ (n+1)}(\alpha) \neq \infty$.  One has
$$a_n(\alpha) = \left\lfloor \frac{a_{n}(\alpha)q_{n-1}(\alpha)+q_{n-2}(\alpha)}{q_{n-1}(\alpha)}  \right\rfloor = \left\lfloor \frac{q_{n}(\alpha)}{q_{n-1}(\alpha)}  \right\rfloor,$$
and therefore
\begin{align}\label{cfineqa}
\frac{q_{n}(\alpha)}{q_{n-1}(\alpha)} = [a_n(\alpha),a_{n-1}(\alpha),\ldots,a_1(\alpha)].
\end{align}
By \cite[Theorem 9.9]{leve},  one has
\begin{align}\label{cfineq0}
 \alpha -\frac{p_n(\alpha)}{q_n(\alpha)}  = \frac{ (-1)^n}{q_n(\alpha) ( q_n(\alpha) S^{ \circ(n+1)} (\alpha)+q_{n-1}(\alpha))}
\end{align}
and therefore
\begin{align}\label{cfineq}
\frac{1}{2q_n(\alpha) q_{n+1}(\alpha) } < \frac{1}{q_n(\alpha) (q_{n+1}(\alpha)+q_{n}(\alpha) )} < \left| \alpha -\frac{p_n(\alpha)}{q_n(\alpha)} \right|< \frac{1}{q_n(\alpha)  q_{n+1}(\alpha)}.
\end{align}
The relations (\ref{cfineq00})--(\ref{cfineq}) above are fundamental to the theory of  continued fractions and are used on several occasions in this chapter.


Throughout this chapter, we let
$$\Phi = \frac{1+\sqrt{5}}{2} =  [1,1,1,\ldots]\index[symbols]{.u  k@$\Phi$}$$
denote the {\bf golden ratio}.  We also let $F_n$\index[symbols]{.u  l@$F_n$} denote the $n$th Fibonacci number,  where $F_0 = 0$ and $F_1 = 1$, so that $F_n \sim \frac{1}{\sqrt{5}} \Phi^n \ (n \to \infty)$ and $\log F_{n} \sim n \log \Phi \ (n \to \infty)$.   Note that, by the recurrence relation for the $q_n(\alpha)$ and the $F_n$, one has $q_n(\alpha) \geq F_{n+1} = q_n(\Phi)$ for any irrational $\alpha$ and for every nonnegative integer $n$.

Now, let $\alpha \in \RR$.   We let
$$\ord_1 \alpha = \# ((1,\alpha)\ZZ/\ZZ),\index[symbols]{.u  e@$\ord_1 \alpha$}$$
denote the {\bf order of $\alpha$ modulo $1$},\index{order modulo $1$}
where $(1,\alpha)\ZZ/\ZZ$ is the subgroup of the circle group $\RR/\ZZ$ generated by $\overline{\alpha} = \alpha+\ZZ$.  In other words,  $\ord_1  \alpha$ is the order of $\overline{\alpha}$ in the abelian group $\RR/\ZZ$.
This defines a map
$$\ord_1: \RR  \longrightarrow \ZZ_{>0} \cup \{\infty\}.$$
We also  let
$$\left|\alpha\right|_1 = \frac{1}{\ord_1 \alpha}\index[symbols]{.u  f@$\vert\alpha\vert_1$}$$
for all $\alpha \in \RR$, where $\frac{1}{\infty} = 0$.  
If $\frac{a}{b} \in \QQ$, where $a,b \in \ZZ$ are relatively prime and $b >0$, then one has
\begin{align*}
\ord_1 \frac{a}{b} = b
\end{align*}
and
\begin{align*}
\left|\frac{a}{b}\right|_1 = \frac{1}{b}.
\end{align*}
On the other hand,  one has $\ord_1 \alpha = \infty$ and $\left|\alpha\right|_1  = 0$
for all $\alpha \in \RR\backslash \QQ$.  The following proposition provides some basic properties of the maps $\ord_1$ and $|\!-\!|_1$.

\begin{proposition}
Let $r,s \in \QQ$ and $n,k \in \ZZ$ with $k > 0$. One has the following.
\begin{enumerate}
\item $\ord_1 r$ is the smallest $b \in \ZZ_{>0}$ such that $br \in \ZZ$.
\item $\ord_1 r = 1$ if and only if $|r|_1 = 1$, if and only if $r \in \ZZ$.
\item $|r|_1 = |r+n|_1 = |\{r\}|_1 = |-r|$.
\item If $r = \frac{a}{b}$, where $a, b \in \ZZ$ and $b \neq 0$, then $\ord_1 r =  \frac{|b|}{\gcd(a,b)}$ divides $b$ and $|r|_1 \geq \frac{1}{|b|}$.
\item $\ord_1(r+s) \mid \operatorname{lcm}( \ord_1 r, \ord_1 s)$.
\item $\ord_1(rs) \mid \ord_1 r \ord_1 s$.
\item $|r+s|_1 \geq |r|_1|s|_1$.
\item $|rs|_1 \geq |r|_1|s|_1$.
\item One has $\ord_1 r = k$ if and only if $\{r\} = \tfrac{a}{k}$ with $\gcd(a, k) = 1$ and $a \in [0,k-1]$.
\item There are exactly $\phi(k)$ elements $\overline{\alpha} = \alpha+\ZZ$ of $\QQ/\ZZ$  (or  of $\RR/\ZZ$)  with $\ord_1 \alpha  = k$.
\item For all $x \geq 1$, there are exactly $\sum_{k \leq x} \phi(k)$  elements $\overline{\alpha} = \alpha+\ZZ$ of $\RR/\ZZ$ (or of $\RR/\ZZ$) with $\ord_1 \alpha \leq x$.
\item For all $\alpha \in \RR$, one has
$$\left|\alpha\right|_1= \inf \, (1,\alpha)\ZZ \cap \RR_{>0} = \inf\left\{m+n \alpha: m,n \in \ZZ \text{ and } m+n \alpha  > 0\right \}.$$
\end{enumerate}
\end{proposition}

 For any $x \in \RR$, let $\operatorname{nint}(x)$\index[symbols]{.u  g@$\operatorname{nint}(x)$}  denote the nearest integer  to $x$, which is equal to the floor of $x$ or the ceiling of $x$, whichever is closer to $x$.  Of course, there is ambiguity here in the case where $x$ is half an odd integer: in that case, we assume that $\operatorname{nint}(x)$ is even.   Thus, for example, one has $\operatorname{nint}(1/2) = 0$ and $\operatorname{nint}(3/2) = 2$.  For all $x \in \RR$,  we let
$$\Vert x\Vert  = |x-\operatorname{nint}(x)| = \min\{|x-n|: n \in \ZZ\} = \min(\{x\},1-\{x\}) \in [0,\tfrac{1}{2}]\index[symbols]{.u  h@$\Vert\alpha\Vert$}$$ 
denote the distance from $x$ to the nearest integer.   The function $\Vert \!-\! \Vert: \RR \longrightarrow [0,\frac{1}{2}]$ is continuous, even, and periodic of period $1$ and equals $|x|$ on $[-\frac{1}{2},\frac{1}{2}]$ and thus has the Fourier series representation
 $$\Vert x\Vert  = \frac{1}{4}  - \frac{2}{\pi^2}\sum_{n = 0}^\infty \frac{\cos (( 2n+1)2\pi x)}{(2n+1)^2}, \quad \forall x \in \RR,$$  
 where $\frac{1}{4} $ is the average value of $\Vert x\Vert$ over any full period,  i.e.,  over any interval of length $1$ (or of any positive integer length).    Moreover,  any function $f:\RR \longrightarrow \overline{\RR}$ is even and periodic, with $1$ as a period, if and only if $f(x) = f(\Vert x\Vert)$ for all $x \in \RR$, if and only if $f$ factors through the map  $\Vert \!-\!\Vert$.   By contrast, a function $f:\RR \longrightarrow \overline{\RR}$ is periodic, with $1$ as a period, if and only if $f(x) = f(\{x\})$ for all $x \in \RR$, if and only if $f$ factors through the map  $\{ - \}$.   Note also that the map $d(x,y) = \Vert x-y\Vert $ is a metric on the circle group $\RR/\ZZ$ inducing its topology as a topological group (or as a one-dimensional compact  Lie group), and in fact the map  $\Vert \!-\! \Vert$ is a {\it norm} on the  group $\RR/\ZZ$, in the sense of \cite{fark}.    Moreover, one has $\Vert x \Vert    =  \Vert y \Vert$,  if and only if $\overline{x} = \pm\overline{y}$ in $\RR/\ZZ$, if and only if $x= \pm y + n$ for some $n \in \ZZ$.

Propositions \ref{vert1}--\ref{dirapproxthm} below are the first of several results in this chapter that reveal connections between the Diophantine approximation properties of a real number $\alpha$ and the properties of the arithmetic function $\Vert n\alpha \Vert$. 
The first two of these are easily verified.
 
 \begin{proposition}\label{vert1}
Let $\alpha \in \RR$, and let $n$ be  a positive integer.   One has
$$\frac{\Vert n \alpha \Vert }{n} = \left|  \alpha - \frac{\operatorname{nint}(n\alpha)}{n}\right| = \min \left\{| \alpha-r|: \text{$r \in\QQ$ and $\ord_1 r \mid n$}\right\},$$ where
  $$|n \alpha-\operatorname{nint}(n\alpha)|  = \Vert n\alpha \Vert  \leq \frac{1}{2},$$
or, equivalently,
 $$\left|\alpha-\frac{\operatorname{nint}(n\alpha)}{n} \right| = \frac{\Vert n \alpha \Vert }{n}  \leq \frac{1}{2n}.$$
If the inequalities above are strict, then $\frac{\operatorname{nint}(n\alpha)}{n}$ is the closest rational number $r$ to $\alpha$ with $\ord_1 r \mid n$.  Moreover, equalities hold if and only if $n\alpha \in 2\ZZ\pm \frac{1}{2}$, and, if equalities hold, then  $\alpha \in \QQ$, and both  $\frac{\operatorname{nint}(n\alpha)}{n}$ and $\frac{\operatorname{nint}(n\alpha) \mp 1}{n}$ are  the closest rational numbers $r$ to $\alpha$ with $\ord_1 r \mid n$, where the sign of $\mp 1$ is opposite that of $\pm \frac{1}{2}$ and is uniquely determined by $\alpha$.
\end{proposition}

\begin{proposition}
Let $\alpha \in \RR$.   Then $\alpha$ is irrational if and only if  the arithmetic function $\Vert n\alpha \Vert$ is injective, if and only if the image of the arithmetic function $\Vert n\alpha \Vert$ is dense in $[0,\frac{1}{2}]$.  Moreover,  the following conditions are equivalent.
\begin{enumerate}
\item  $\alpha$ is rational.
\item  The arithmetic function $\Vert n\alpha \Vert$ is periodic (of period $\ord_1 \alpha$).
\item The arithmetic function $\Vert n\alpha \Vert$ has zeros (precisely at the multiples of  $\ord_1 \alpha$).
\item  The image of the arithmetic function $\Vert n\alpha \Vert$ is finite (and equal to $|\alpha|_1 \ZZ \cap [0,\frac{1}{2}]$,  of cardinality $1+ \left \lfloor \frac{\ord_1 \alpha}{2}\right \rfloor$).
\end{enumerate}
\end{proposition}

Let $\alpha \in \RR$.  One says that a rational number $r$ is a {\bf best rational approximation of $\alpha$ (of the first kind)} if $|\alpha -r| < |\alpha - s|$  for all rational numbers $s \neq r$ with $\ord_1 s \leq \ord_1 r$ \cite[p.\ 21]{khin}.   By \cite[Theorems 15--17]{khin}, every convergent of the regular continued fraction expansion of $\alpha$ is a best rational approximation of $\alpha$, and, reversely, every best rational approximation of $\alpha$ is either a convergent or an {\it intermediate fraction} of the regular continued fraction expansion of $\alpha$.    In particular, the denominators of the best rational approximations of $\alpha$ form an increasing sequence that contains the sequence $q_n(\alpha)$ as a subsequence.    For $\alpha = \pi$, for example,  these sequences are OEIS Sequence A063673, given by $1, 4, 5, 6, 7, 57, 64, 71, 78, 85, 92, 99, 106, 113, 16604, 16717,\ldots$,  and OEIS Sequence A002486, given by  $1, 7, 106, 113, 33102, 33215, 66317,\ldots$, respectively.   

One says that a rational number $\frac{a}{b}$ with $a,b \in \ZZ$ and $b \neq 0$ is a {\bf best rational approximation of $\alpha$ of the second kind} if $|b\alpha -a| < |d\alpha - c|$  for all $c,d \in \ZZ$ with $0< |d| \leq |b|$  and $\frac{a}{b} \neq \frac{c}{d}$ (which implies that $\gcd(a,b) = 1$ and therefore $|b| = \ord_1 \frac{a}{b}$) \cite[p.\ 24]{khin}.      By Proposition \ref{denomschar} below, any best rational approximation of the second kind is a best rational approximation of the first kind.   By \cite[Theorems 16 and 17]{khin}, the convergents of the regular continued fraction expansion of an irrational number $\alpha$ are precisely the best rational approximations of $\alpha$ of the second kind.  

The following result follows from the definitions and remarks above.

\begin{proposition}\label{denomschar}
Let $\alpha \in \RR$ be irrational.  
\begin{enumerate}
\item The denominators $q_k(\alpha)$ of the convergents of the regular continued fraction expansion of $\alpha$ for $k \geq 0$ are precisely those positive integers $n$ such that  $\Vert n \alpha \Vert < \Vert m \alpha \Vert$ for all positive integers $m < n$.   Equivalently,  the points of discontinuity of the step function $\min\left \{\Vert n\alpha \Vert : n \leq x \right \}$ occur precisely at $x = q_k(\alpha)$ for $k \geq 0$.   Moreover, the $q_k(\alpha)$ are precisely the denominators of the best rational approximations of $\alpha$ of the second kind.
\item The denominators of the best rational approximations of $\alpha$ are precisely those positive integers $n$ such that  $\frac{\Vert n \alpha \Vert}{n} < \frac{ \Vert m \alpha \Vert}{m}$ for all positive integers $m < n$.   Equivalently,   the denominators of the best rational approximations of $\alpha$ are precisely the  points of discontinuity of the step function 
\begin{align*}
\mu_\alpha(x) & =  \min\left \{\frac{\Vert n\alpha \Vert }{n} : n \leq x \right \} \\
& =  \min\left\{ | \alpha-r|: \text{$r \in\QQ$ and $\ord_1 r \leq x$} \right\} \\
& =  \min\left\{ \left| \alpha-\frac{a}{b}\right|: \text{$a,b \in\ZZ$ and $1 \leq b \leq x$} \right\}
\end{align*}
on $[1,\infty)$.
\item The step functions in (1) and (2) are positive, nonincreasing with limit $0$, and continuous from the right.
\item For all positive integers $n$, one has
$$\operatorname{nint}(q_n(\alpha)\alpha)= p_n(\alpha)$$
and
$$\frac{\Vert q_n(\alpha)\alpha) \Vert}{q_n(\alpha)} = \left| \alpha- \frac{p_n(\alpha)}{q_n(\alpha)} \right|.$$
\end{enumerate}
\end{proposition}

 A general problem in Diophantine approximation is to find explicit bounds on the function $\mu_\alpha(x)$ for various irrationals or classes of irrationals $\alpha$ \cite{feld}.

The following theorem is known as {\bf Dirichlet's approximation theorem}\index{Dirichlet's approximation theorem} and  is a simple  consequence of the pigeonhole principle. 

\begin{proposition}[Dirichlet's approximation theorem]\label{dirapproxthm}
Let $\alpha \in \RR$.   For any real number $x \geq 1$, there exist integers $a$ and $b$ with $1\leq b \leq x$ such that
$$\left|b\alpha -a\right| < {\frac {1}{\lfloor x\rfloor}} \leq \frac{1}{b},$$
i.e., such that
$$\left|\alpha -\frac{a}{b}\right| < {\frac {1}{\lfloor x\rfloor b}} \leq \frac{1}{b^2}.$$ 
Equivalently, one has
$$\min\left \{\Vert n\alpha \Vert : n \leq x \right \}<   {\frac {1}{\lfloor x\rfloor}}$$
for all $x \geq 1$.  Consequently, one has $\Vert n \alpha \Vert < \frac{1}{n}$ for infinitely many positive integers $n$.
\end{proposition}

\begin{proof}
Consider the $\lfloor x \rfloor$ disjoint intervals $[k/\lfloor x \rfloor,(k+1)/\lfloor x \rfloor))$ of length $\frac{1}{\lfloor x \rfloor}$  for $k = [0,\lfloor x \rfloor-1]$, whose union is $[0,1)$, and the $\lfloor x \rfloor+1$ fractional parts $\{k\alpha\} \in [0,1)$ for $k \in [0,\lfloor x \rfloor]$.
By the pigeonhole principle, there exist integers $k>l$ in $[0,\lfloor x \rfloor]$ such that 
$$|b\alpha - a| =  |\{k\alpha \} -\{l\alpha\}| < \frac{1}{\lfloor x\rfloor},$$
where $b = k-l \in [1,\lfloor x \rfloor]$ and $a = \lfloor  k\alpha\rfloor  -\lfloor  l\alpha\rfloor$.
\end{proof}

\begin{corollary}\label{dircor}
Let $\alpha \in \RR$.   For any real number $x \geq 1$, there exists a rational number $r$ with $\ord_1 r \leq x$ and
$$\left|\alpha -r\right| < {\frac {1}{\lfloor x\rfloor }}|r|_1 \leq |r|_1^2.$$
\end{corollary}

The following well-known result illustrates why,  in studying the rational approxi­mations of a real number, one can often restrict one's attention to the convergents of its regular continued fraction expansion.   It is thus one of the main reasons that continued fractions are a central tool in the theory of Diophantine approximation.

\begin{theorem}[{\cite[Theorems 9.9 and 9.10]{leve} \cite[Theorem 18]{khin}}]\label{ccfa}
Let $\alpha \in \RR$.   For every convergent $\frac{p_n(\alpha)}{q_n(\alpha)}$  of the regular continued fraction expansion of $\alpha$,  one has
$$\left | \alpha -\frac{p_n(\alpha)}{q_n(\alpha)} \right| < \frac{1}{q_n(\alpha)^2},$$
and, if $n > 0$, then  at least one of $\frac{a}{b} = \frac{p_n(\alpha)}{q_n(\alpha)}$ and $\frac{a}{b} =\frac{p_{n-1}(\alpha)}{q_{n-1}(\alpha)}$ satisfies
$$\left | \alpha -\frac{a}{b} \right| < \frac{1}{2b^2}.$$
Conversely, if $a,b \in \ZZ$ with $b \neq 0$ and
$$\left | \alpha -\frac{a}{b} \right| < \frac{1}{2b^2},$$
then $\gcd(a,b) = 1$ and $\frac{a}{b}$ is a convergent of the regular continued fraction expansion of $\alpha$.
\end{theorem}

Let $f: \ZZ_{>0} \longrightarrow \RR_{\geq 0}$ be a nonnegative real-valued arithmetic function, and let $\alpha \in \RR$.     Let us write
$$\alpha \gg_1 f\index[symbols]{.u  i@$\gg_1$}$$
if there are only finitely many $r \in \QQ$ such that $$\left|\alpha -r \right| < f(\ord_1 r).$$
that is, if  one has
 $$\left|\alpha -r \right| \geq f(\ord_1 r)$$
for all but finitely many  $r \in \QQ$.
Let us also write
$$\alpha \ggg_1 f\index[symbols]{.u  j@$\ggg_1$}$$
if $$\left|\alpha -r \right| \geq f(\ord_1 r)$$
for all $r \in \QQ$ with $r \neq \alpha$, that is, if $\left|\alpha -r \right| < f(\ord_1 r)$ implies $r = \alpha$,  for all $r \in \QQ$.   Note that 
$\alpha \ggg_1 f$ implies $\alpha \gg_1 f$,  and both conditions are specific restrictions on how well $\alpha$ can be approximated by rational number.   More precisely,  the conditions $\alpha \gg_1 f$ and $\alpha \ggg_1 f$ place restrictions on how many rational numbers with denominator exactly $d$ are within $f(d)$ of $\alpha$,      for all $d \gg 0$,  or for all $d$ in total,  respectively.   Not only do these two notations allow us to shorten the statements of several definitions and results in this chapter,  but they also suggest various analogies between definitions and results in the theory of Diophantine approximation and those in algebraic asymptotic analysis.  They are meant, specifically, to be suggestive of the notation $f(x)\gg g(x) \ (x \to a)$.

Either by iterating Corollary \ref{dircor} or by applying Theorem \ref{ccfa}, one obtains the following.

\begin{corollary}\label{dirapproxcor}
Let $\alpha \in \RR$ be irrational.    Then $\alpha \not \gg_1  n^{-2}$, that is, there exist infinitely many rational numbers $r$ such that $|\alpha - r| <  |r|_1^2$.
\end{corollary}

Again, let $f: \ZZ_{>0} \longrightarrow \RR_{\geq 0}$ be a nonnegative real-valued arithmetic function, and let $\alpha \in \RR$.     We also write
$$\alpha \succ_1 f\index[symbols]{.u  ja@$\succ_1$}$$ if  there are only  finitely many pairs $(a,b)$ of  integers with $b > 0$ and
$$\left|\alpha-\frac{a}{b}\right|< f(b).$$   Note that $\alpha \succ_1 f$ implies $\alpha \gg_1 f$, but not conversely.    
  The following proposition provides some basic properties of the relations $ \gg_1$, $\ggg_1$, and $\succ_1$.

\begin{proposition}\label{basicrelations}
Let $\alpha \in \RR$, and let $f,g: \ZZ_{>0} \longrightarrow \RR_{\geq 0}$ be nonnegative real-valued  arithmetic functions.   One has the following.
\begin{enumerate}
\item For all $x \geq 1$, there are only finitely many pairs $(a,b)$ of integers $a$ and $b$ with 
$$\left|\alpha -\frac{a}{b} \right| < f(b)$$
and $1 \leq b \leq x$.   
\item For all $x \geq 1$, there are only finitely many $r \in \QQ$ such that
$$\left|\alpha -r \right| < f(\ord_1 r)$$
and $\ord_1 r \leq x$.   
\item  One has $\alpha \gg_1 f$ if and only if there are only finitely many pairs $(a,b)$ of relatively prime integers $a$ and $b$  with $b > 0$ and
$$\left|\alpha -\frac{a}{b} \right| < f(b),$$  if and only if
$$\left|\alpha -\frac{a}{b} \right| \geq f(b)$$
for all relatively prime integers $a$ and $b$  with $b > 0$ and $b \gg 0$, if and only if 
 $$\left|\alpha -r \right| \geq f(\ord_1 r)$$
 for all $r \in \QQ$ with $\ord_1 r \gg 0$.
 \item One has $\alpha \succ_1 f$ if and only if $$\left|\alpha-\frac{a}{b}\right| \geq f(b)$$ for all integers $a$ and $b$ with $b > 0$ and $b\gg 0$, if and only  if $\frac{\Vert n\alpha \Vert}{n} \geq f(n)$ for all   $n \gg 0$, if and only if there are only finitely many positive integers $n$ such that $\frac{\Vert n\alpha \Vert}{n} < f(n)$.
 \item For all $n \in \ZZ$, one has $\alpha \gg_1 f$ if and only if $\alpha +n \gg_1 f$, if and only if $-\alpha \gg_1 f$.
  \item For all $n \in \ZZ$, one has $\alpha \ggg_1 f$ if and only if $\alpha +n \ggg_1 f$, if and only if $-\alpha \ggg_1 f$.
   \item For all $n \in \ZZ$, one has $\alpha \succ_1 f$ if and only if $\alpha +n \succ_1 f$, if and only if $-\alpha \succ_1 f$.
 \item If $f(n) \geq g(n)$ for all $n \gg 0$, then $\alpha \gg_1 f$ implies $\alpha \gg_1 g$.
  \item If $f(n) \geq g(n)$ for all $n$, then $\alpha \ggg_1 f$ implies $\alpha \ggg_1 g$.
 \item If $f(n) \geq g(n)$ for all $n \gg 0$, then $\alpha \succ_1 f$ implies $\alpha \succ_1 g$.
    \item If $\alpha \ggg_1 f$, then  $\alpha \gg_1 f$.
     \item If $\alpha \succ_1 f$, then  $\alpha \gg_1 f$.
 \item If $\alpha \gg_1 f$, then there exists a $C > 0$ such that $\alpha \ggg_1 Cf$. 
 \item Suppose that $\alpha$ is irrational and that $\lim_{n \to \infty} f(n) = 0$.   Then  $\alpha \gg_1 f$  if and only if $\alpha \succ_1 f$.
\item Suppose that $\alpha$ is irrational, and let  $t  > 0$.  Then one has $\alpha \gg_1 n^{-t}$ if and only if $\alpha \succ_1 n^{-t}$, if and only if there are only finitely many pairs $(a,b)$ of integers $a$ and $b$  with $b> 0$ and
$$\left|\alpha -\frac{a}{b} \right| < \frac{1}{b^t},$$ if and only if $\frac{\Vert n\alpha \Vert}{n}  \geq \frac{1}{n^{t}}$ for all $n \gg 0$.
\item  Suppose that $f(n) > 0$ for all $n \gg 0$.  Then there are uncountably many $\beta \in \RR$ such that $\beta \not \gg_1 f$.
\end{enumerate}
\end{proposition}

\begin{proof}
For all positive integers $b$, there are only finitely many integers $a$ such that $|b\alpha-a|< bf(b)$, i.e., such that $|\alpha - a/b|<f(b)$.   Statements (1)--(4) follow, and statements (5)--(12) are trivial to check.  To prove (13), suppose that $\alpha \gg_1 f$, so that there exists an $N$  such that  $$\left|\alpha -r \right| \geq f(\ord_1 r)$$
 for all $r \in \QQ$ with $\ord_1 r > N$.  Let
 $$C = \min\{|\alpha-r|f(\ord_1 r)^{-1}: r \in \QQ\backslash\{\alpha\},\, \left|\alpha -r \right| < f(\ord_1 r), \text{ and } \ord_1 r \leq N\}.$$
 Then $C \in (0,1)$.    It follows that, for all $r \in \QQ$,  one has
 $$\left|\alpha -r \right| \geq f(\ord_1 r) \geq Cf(\ord_1 r)$$
if $\ord_1 r > N$,  while, if
 $\ord_1 r \leq N$ and   $$\left|\alpha -r \right| < Cf(\ord_1 r) < f(\ord_1 r),$$
 then
  $$\left|\alpha -r \right| f(\ord_1 r)^{-1} \geq C$$
  if $r \neq \alpha$, which is impossible.   Therefore, one has $\alpha \ggg_1 Cf$.   This proves (13).

    Next, we prove (14).  Suppose that $\alpha$ is irrational and $\alpha \gg_1 f$, and  let $(a,b)$ be any  of the finitely many pairs of relatively prime integers $a$ and $b$  with $b > 0$ and
$$\left|\alpha -\frac{a}{b} \right| < f(b).$$
Since $\left|\alpha -\frac{a}{b} \right|> 0$ and $\lim_{c \to \infty} f(cb) = 0$, one has
$$\left|\alpha -\frac{ca}{cb} \right|  = \left|\alpha -\frac{a}{b} \right| >  f(cb)$$
for all $c \gg 0$.
It follows that there are only finitely many pairs $(a',b') = (ca,cb)$ with $c > 0$ and
$$\left|\alpha -\frac{a'}{b'} \right| \leq  f(b')$$ corresponding to each pair $(a,b)$,
whence $\alpha \succ_1 f$.  Statements (14) and  (15) follow.  Finally, statement (16) follows from the proof of  \cite[Theorem 22]{khin} (which is a simple application of the theory of continued fractions).
\end{proof}

The following result,  known as {\bf Hurwitz's theorem},\index{Hurwitz's theorem}  was proved by Hurwitz in 1891 and is an optimized version of Corollary \ref{dirapproxcor}.

\begin{theorem}[{Hurwitz's theorem \cite{hurw}}]
Let $\alpha \in \RR$.  Then  $\alpha$ is irrational if and only if  there are infinitely many rational numbers $r$ such that
$$\left|\alpha -r \right| < \frac{1}{\sqrt{5}} \left|r\right|_1^2,$$ if and only if
$$\alpha \not \gg_1 \frac{1}{\sqrt{5}}n^{-2}.$$  
Moreover,  one has
$$\Phi \gg_1  t n^{-2}$$
for all $t < \frac{1}{\sqrt{5}}$, so the constant  $\frac{1}{\sqrt{5}}$ is optimal.
\end{theorem}

We also note the following result, which is complementary to Hurwitz's theorem.

\begin{proposition}\label{dioprop}
Let $\alpha \in \QQ$.   Then there are infinitely many $r \in \QQ$ such that
$$\alpha- r  = |\alpha|_1 |r|_1,$$
and therefore
$$\alpha  \not \gg_1  t n^{-1}$$
for all  $t > |\alpha|_1$.
Moreover, one has
$$\left|\alpha- r \right| \geq |\alpha|_1 |r|_1$$
for all $r \in \QQ\backslash \{\alpha\}$, that is, one has
$$\alpha \ggg_1 |\alpha|_1 n^{-1}.$$
\end{proposition}

\begin{proof}
Let $\alpha = \frac{a}{b}$, where $a$ and $b$ are relatively prime integers and $b > 0$.   For any integers $x$ and $y$ such that $ay-bx = 1$ and $y \neq 0$, one has
$$\alpha - r  = \frac{1}{by} = |\alpha|_1 |r|_1,$$
where $r = \frac{x}{y}$.  Since there are infinitely many (necessarily relatively prime) pairs $(x,y)$ of such integers,  the first statement follows.  To prove the second statement, let $r =  \frac{x}{y}$, where  $x$ and $y$ are relatively prime integers and $y > 0$.   Then one has
$$|\alpha - r| = \frac{|ay-bx|}{by} \geq \frac{1}{by} = |\alpha|_1 |r|_1$$
if $\alpha  \neq r$, that is, if $ay-bx \neq 0$.
\end{proof}

Here and henceforth,  {\it almost all real numbers} means all real numbers outside a set of Lebesgue measure zero.     The following remarkable result, proved by Koukoulopoulos  and Maynard  in 2019,   settled two important conjectures in the positive,  namely,  the {\bf Duffin--Schaeffer conjecture} and {\bf Catlin's conjecture}. \index{Duffin--Schaeffer conjecture}    As is usual with deep theorems, there is a long chain of partial results preceding it: for a survey, see \cite{koukou}.

\begin{theorem}[{\cite[Theorems 1 and 2]{kou2}}]\label{duffin}
Let  $f: \ZZ_{>0} \longrightarrow \RR_{\geq 0}$ be a nonnegative real-valued arithmetic function.   For all positive integers $n$, let
$$f^*(n) =  \sup\left\{ f(kn) : k \in \ZZ_{>0}\right\},$$ so that  $f^*$ is the smallest extended real-valued function greater than equal to $f$ such that $f^*(a) \geq f^*(b)$ if $a \mid b$.   One has the following.
\begin{enumerate} 
\item  If $\sum_{n = 1}^\infty \phi(n) f(n) < \infty$, then  $\alpha \gg_1 f$  for almost all $\alpha \in \RR$.  On the other hand, if 
$\sum_{n = 1}^\infty \phi(n) f(n) = \infty$, then $\alpha \not \gg_1 f$  for almost all $\alpha \in \RR$.
\item  If $\sum_{n = 1}^\infty \phi(n) f^*(n)  < \infty$, then $\alpha \succ_1 f$   for almost all $\alpha \in \RR$.  On the other hand, if $\sum_{n = 1}^\infty \phi(n) f^*(n)  = \infty$, then $\alpha \not \succ_1 f$  for almost all $\alpha \in \RR$.
\end{enumerate}
\end{theorem}

Note, of course, that if $f$ is nonincreasing, then $f^* = f$ and thus the convergence conditions in statements (1) and (2) of the theorem are equivalent.   Likewise, if $\lim_{n \to \infty} f(n) = 0$, then  $\alpha \gg_1 f$ is equivalent to $\alpha \succ_1 f$.

 \begin{remark}[Littlewood's conjecture]
{\bf Littlewood's conjecture},\index{Littlewood's conjecture} proposed by Littlewood around 1930, states that, for all $\alpha,\beta \in \RR$, one has
$$ \liminf _{n\to \infty }\ n\,\Vert n\alpha \Vert \,\Vert n\beta \Vert =0,$$
or, equivalently,  for every $\varepsilon > 0$, there exist infinitely many positive integers $n$ such that 
$$n\,\Vert n\alpha \Vert \,\Vert n\beta \Vert < \varepsilon.$$
 \end{remark}

 \begin{remark}[On the integral $\int_1^\infty \frac{\Vert x\Vert}{x^{s+1}}\, dx $]\label{vertremark}
For all $s \in  \CC$ with $\operatorname{Re} s > 2$, one has
 \begin{align*}
 \int_1^\infty \frac{\Vert x\Vert}{x^{s+1}}\, dx & = \frac{1}{s(s-1)} \sum_{n = 1}^\infty (n^{1-s} +(n+1)^{1-s} -2^{s}(2n+1)^{1-s}) \\
 & =  \frac{(4-2^{s})\zeta(s-1)+2^{s}-1}{s(s-1)}.
 \end{align*}
It follows that 
  \begin{align*}
 \int_1^\infty \frac{\Vert x\Vert}{x^{s+1}}\, dx  =  \frac{(4-2^{s})\zeta(s-1)+2^{s}-1}{s(s-1)}
 \end{align*}
 for all $s \in  \CC \backslash \{1,2\}$ with $\operatorname{Re} s > 0$, 
with limiting value
 $$\int_1^\infty \frac{\Vert x\Vert}{x^2}\, dx = \sum_{n = 1}^\infty \log \left( \frac{(2n+1)^2}{
(2n+1)^2-1}\right)  = \log \prod_{n = 1}^\infty \frac{(2n+1)^2}{
(2n+1)^2-1}  = \log \frac{4}{\pi}$$
 at $s = 1$, where $ \log \frac{4}{\pi}$ is the {\bf alternating Euler constant\index{alternating Euler constant}} \cite{sond}, and limiting value
  $$\int_1^\infty \frac{\Vert x\Vert}{x^3}\, dx = \frac{1}{2} \sum_{n = 1}^\infty \frac{1}{n(n+1)(2n+1)} =  \frac{3}{2}-2\log 2$$
at $s = 2$.    Reversely, one has
    \begin{align*}
(2-2^s)\zeta(s) =  \frac{s(s+1)}{2}\int_1^\infty \frac{\Vert x\Vert}{x^{s+2}}\, dx+2^{-1}-2^s
\end{align*}
for all $s \in \CC\backslash \{1\}$ with $\operatorname{Re} s > -1$, where the left hand side is entire.   It follows that 
$$\lim_{{s \to 0}\atop {\operatorname{Re} s > 0}}s \int_1^\infty \frac{\Vert x\Vert}{x^{s+1}}\, dx = -3\zeta(-1) = \frac{1}{4},$$
and therefore the function
\begin{align*}
s \int_1^\infty \frac{\Vert x\Vert}{x^{s+1}}\, dx =  \frac{(4-2^{s})\zeta(s-1)+2^{s}-1}{s-1}  
 \end{align*}
entends to an entire function with values $\frac{1}{4}$ at $s = 0,-1$  and $\log \frac{4}{\pi}$ at $s = 1$ (cf.\  (\ref{zetaexp})), and the function
\begin{align*}
 \int_1^\infty \frac{\Vert x\Vert-\frac{1}{4}}{x^{s+1}}\, dx  =  \int_1^\infty \frac{\Vert x\Vert}{x^{s+1}}\, dx  + \frac{3\zeta(-1)}{s}  
 \end{align*}
on $\operatorname{Re} s  > 0$ also extends to entire function with a zero at $s = -1$.  It also follows that the function $ \int_1^\infty \frac{\Vert x\Vert}{x^{s+1}}\, dx$ has meromorphic continuation to $\CC$ with a single (simple) pole at $s = 0$ with residue $-3\zeta(-1) = \frac{1}{4}$, which is the average value of $\Vert x\Vert$ over any full period.  This provides an interesting interpretation of the well-known identity $\zeta(-1) = -\frac{1}{12}$.  Moreover,  at $s = 0$  one has
$$ \int_1^\infty \frac{\Vert x\Vert-\frac{1}{4}}{x}\, dx  = -3(\zeta(-1)+\zeta'(-1)) -\frac{13}{12}\log 2 = 3\log A-\frac{13}{12}\log 2 = -0.00464601\ldots,$$  where $A$ is the {\it  Glaisher--Kinkelin constant}.
 \end{remark}

\section{Irrationality measure}

The {\bf irrationality measure}, {\bf approximation order}, or {\bf approximation exponent},\index{irrationality measure $\mu(\alpha)$}\index{approximation order $\mu(\alpha)$}\index{approximation exponent $\mu(\alpha)$} {\bf of} $\alpha \in \RR$, denoted $\mu(\alpha)$, is given by
\begin{align*}
\mu(\alpha)  = \inf\{t \in \RR: \alpha \gg_1 n^{-t}\}  = -\sup\{t \in \RR: \alpha \gg_1 n^{t}\} \in [0,\infty].\index[symbols]{.u  m@$\mu(\alpha)$} 
\end{align*}
Given the notation $\gg_1$ and its  resemblance (by design) to the asymptotic notation $\gg$, the definition above is reminiscent of the expression $$\underline{\deg} \, f  = \sup \{t \in \RR: f(x) \gg x^t  \ (x \to \infty)\}$$ in Proposition \ref{firstprop0lower}(3).    Corollary \ref{degmu} below shows that this resemblance is more than superficial.

By Hurwitz's theorem and  Proposition \ref{dioprop}, one has the following.

\begin{proposition}
Let $\alpha \in \RR$.   Then $\alpha$ is rational if and only if $\mu(\alpha) = 1$.  Moreover,  $\alpha$ is irrational if and only if $\mu(\alpha) \geq 2$.
\end{proposition}

The following proposition, whose proof is straightforward,  given Proposition \ref{basicrelations}, provides some basic properties of irrationality measure.

\begin{proposition}\label{muapprox}
Let $\alpha, t \in \RR$ with $t > 0$.  
\begin{enumerate}
\item  The map $\mu: \RR \longrightarrow [1,\infty]$ is even and periodic, with $1$ as a period.  Equivalently, one has  $\mu(\alpha+n) = \mu(\left\{\alpha\right\}) = \mu(\Vert \alpha\Vert ) = \mu(-\alpha) = \mu(\alpha)$ for all $n \in \ZZ$.
\item $\mu(\alpha)$ is the unique $T > 0$ such that $\alpha \gg_1 n^{-u}$ for all $u > T$ and $ \alpha \not \gg_1 n^{-u}$ for all $u < T$.
\item One has
\begin{align*} \mu(\alpha) = \sup\{t \in \RR: \alpha  \not \gg_1   n^{-t}\} = -\inf\{t \in \RR: \alpha  \not \gg_1   n^{t}\}.
\end{align*}
\item $t \geq \mu(\alpha)$ if and only if  $\alpha \gg_1  n^{-u}$ for all $u > t$.
\item $t< \mu(\alpha)$ if and only if  $\alpha  \not \gg_1   n^{-u}$ for some $u > t$.
\item If $t > \mu(\alpha)$, then $\alpha \gg_1 n^{-t}$,  and there exists a $C >0$ such that $\alpha \ggg_1 Cn^{-t}$.
\item If $t< \mu(\alpha)$, then $\alpha  \not \gg_1   n^{-t}$ and $\alpha \, {\not \ggg}_1 \, Cn^{-t}$ for all $C > 0$.
\end{enumerate}
Suppose that $\alpha$ is irrational.
\begin{enumerate}
\item[(8)] $t \geq \mu(\alpha)$ if and only if  $\alpha \succ_1  n^{-u}$ for all $u > t$.
\item[(9)] $t< \mu(\alpha)$ if and only if  $\alpha  \not \succ_1   n^{-u}$ for some $u > t$.
\item[(10)] $\mu(\alpha)$ is the unique $T > 0$ such that $\alpha \succ_1 n^{-u}$ for all $u > T$ and $ \alpha \not \succ_1 n^{-u}$ for all $u < T$.
\item[(11)] One has \begin{align*}
\mu(\alpha) & = \sup\{t \in \RR:  \alpha  \not \succ_1   n^{-t} \}   \\
& = \inf\{t \in \RR:  \alpha \succ_1 n^{-t}\} \\
 & = \inf\left\{t \in \RR: \forall n \gg 0\, \left(  \tfrac{n}{ \Vert n\alpha \Vert} \leq n^{t}\right) \right\}.
\end{align*}
\end{enumerate}
\end{proposition}

As an immediate consequence of statement (11) of the proposition, one has the following.

\begin{corollary}\label{degmu}
For all irrational numbers $\alpha$, one has
\begin{align*}
\mu(\alpha) = \deg \frac{n}{\Vert n \alpha \Vert} = 1+ \deg \frac{1}{\Vert n \alpha \Vert},
\end{align*}
where  $$\frac {\Vert n \alpha \Vert }{n}=  \min \left\{| \alpha-r|: \text{$r \in\QQ$ and $\ord_1 r \mid n$}\right\}$$
for all positive integers $n$.
\end{corollary}

Thus, $\mu(\alpha)$ for any irrational number $\alpha$ can be interpreted as the degree of the  arithmetic function $ \frac{n}{\Vert n \alpha \Vert}$.  In Section 13.4, we study, more generally, the logexponential degree $\dege \frac{n}{ \Vert n \alpha \Vert}$ of the arithmetic function $ \frac{n}{\Vert n \alpha \Vert}$.  

Applying Proposition \ref{supprop}, we obtain the following.

\begin{corollary}\label{degmu2}
For all irrational numbers $\alpha$, one has
\begin{align*}
\mu(\alpha) = \deg \frac{1}{\mu_\alpha(x)} = -\underline{\deg} \,  \mu_\alpha(x),
\end{align*}
where $\mu_\alpha(x)$ is the step function on $[1,\infty)$ defined in Proposition \ref{denomschar}
\end{corollary}

A real number $\alpha$ is said to be {\bf very well approximable}\index{very well approximable} if $\mu(\alpha) > 2$, or, equivalently, if $\alpha  \not \gg_1  n^{-t}$ for some $t > 2$.   By this measure, the larger $\mu(\alpha)$ is, the  ``more well approximable'' the number $\alpha$  is.    It is known (see Corollary \ref{2cor}) that the set of all very well approximable numbers has  (Lebesgue) measure $0$.  In other words, one has $\mu(\alpha) = 2$ for almost all real numbers $\alpha$, that is,  for all $\alpha \in \RR$ outside a set of measure $0$.

\begin{corollary}
Let $\alpha \in \RR$ be irrational.  The following conditions are equivalent.
\begin{enumerate}
\item $\alpha$ is not very well approximable, that  is, $\mu(\alpha) = 2$.
\item For every $t > 2$, one has $\alpha \gg_1 n^{-t}$.
\item For every  $t >2$, there exists a $C > 0$ such that $\alpha \ggg_1 Cn^{-t}$.
\item For every $t > 1$, one has $\Vert n\alpha \Vert < \tfrac{1}{n^{t}}$ for only finitely many positive integers $n$.
\item For every $t \in \RR$, one has $\Vert n\alpha \Vert < \tfrac{1}{n^{t}}$ for infinitely many positive integers $n$ if and only if $t \leq 1$.
\item $\deg \frac{1}{ \Vert n \alpha \Vert }= 1$.
\end{enumerate}
\end{corollary}

Note that, by Dirichlet's approximation theorem,  the ``if'' portion of statement (5) of the corollary is automatic  for any real number $\alpha$.  Thus,  the real numbers $\alpha$ with $\mu(\alpha) = 2$ are those real numbers for which the conclusion of Dirichlet's approximation theorem for $\alpha$ is more or less optimal---and this holds for almost all $\alpha \in \RR$.

A real number $\alpha$ is said to be {\bf Liouville}\index{Liouville number} if $\mu(\alpha) = \infty$.   Liouville numbers are  real numbers that are ``extremely well approximable.''   The following theorem, proved by Liouville in 1844,  implies that all Liouville numbers are transcendental.

\begin{theorem}[{Liouville \cite{liou} \cite{liou2}}]
If $\alpha$ is an algebraic number of degree at most $N$, then there exists a $C > 0$ such that 
$$\alpha \ggg_1 Cn^{-N},$$
and therefore $\mu(\alpha) \leq N$.
\end{theorem}

Continuing the work of Thue (1909), Siegel (1921),  and Dyson (1947),  in 1955,  Roth proved the following generalization of Liouville's theorem,  known as {\bf Roth's theorem}, or the {\bf Thue--Siegel--Roth theorem},\index{Roth's theorem}\index{Thue--Siegel--Roth theorem} for which he received the Fields Medal in 1958.

\begin{theorem}[{Roth's theorem \cite{roth}}]
One has $\mu(\alpha) = 2$ for every algebraic irrational number $\alpha$.    Equivalently,  every very well approximable real number is transcendental. 
\end{theorem}

Next,  let $f,  h \in \RR^{\RR_\infty}$, and suppose that $\dom h$ contains $[N,\infty) \cap \dom f$ for some $N >0$, and also    that $h(x) \neq 0$ and $h(x) \neq 1$  for all $x \gg 0$ in $\dom h$.   We define
$$\deg(f;h) = \limsup_{x \to \infty} \frac{\log |f(x)|}{\log |h(x)|},$$
which we call the {\bf degree of $f$ with respect to $h$}.\index{degree $\deg(f;h)$ of $f$ with respect to $h$}\index[symbols]{.g dz@$\deg(f;h)$} This generalizes the usual degree map $\deg(-) = \deg(-;\id)$.    We also define
$$\underline{\deg}(f;h) = \liminf_{x \to \infty} \frac{\log |f(x)|}{\log |h(x)|} = -\deg(f; 1/h),$$
which we call the {\bf lower degree of $f$ with respect to $h$}.\index{lower degree $\deg(f;h)$ of $f$ with respect to $h$}\index[symbols]{.g dz@$\underline{\deg}(f;h)$}
 Note that the function $\frac{\log |f(x)|}{\log |h(x)|}$
is the unique (germ of a) function   $L(x)$  such that $$|f(x)| = |h(x)|^{L(x)}$$ for all $x \gg 0$ in $\dom f$.  Also note that, if $h$ has finite and positive exact degree, then 
$$\deg(f;h) = \frac{\deg f}{\deg h},$$
while, if $h$ has finite and negative exact degree, then 
$$\deg(f;h) = \frac{\underline{\deg}\, f}{\deg h}.$$
In particular, one has 
$$\deg f =  \deg(f;\id)$$
and 
$$\underline{\deg}f =  \underline{\deg}(f;\id) = - \deg(f;1/\id).$$  In general, one has
$$\deg(f;h) = \deg |f(x)|^{\frac{\log x}{\log |h(x)|}}.$$
Moreover, if $|h(x)| > 1$ for all $x \gg 0$, then
$$\deg(f;h) = \inf\{t \in \RR: |f(x)| \ll |h(x)|^t \ (x \to \infty)\},$$
while,   if $0<|h(x)| < 1$ for all $x \gg 0$, then
$$\deg(f;h) = \inf\{t \in \RR: |f(x)| \gg |h(x)|^t \ (x \to \infty)\}.$$

The following result,  whose proof  in \cite{sondow} is a straightforward application of Theorem \ref{ccfa} and (\ref{cfineq}),    expresses $\mu(\alpha)$ in terms of the regular continued fraction expansion of $\alpha$ and provides motivation for introducing the degree notions defined above.

\begin{proposition}[{\cite[Theorem 1]{sondow}}]\label{sondprop}
For all irrational numbers $\alpha$, one has
\begin{align*}
\mu(\alpha) & = 1+  \limsup_{n \to \infty} \frac{\log q_{n+1}(\alpha)}{\log q_{n}(\alpha)} \\ & = 2+  \limsup_{n \to \infty} \frac{\log a_{n+1}(\alpha)}{\log q_{n}(\alpha)},
\end{align*}
or, equivalently, 
\begin{align*}
\mu(\alpha) & =\deg(q_n(\alpha) q_{n+1}(\alpha); q_n(\alpha)) \\ & =  1+\deg(q_{n+1}(\alpha); q_n(\alpha)) \\ & =  2+\deg(a_{n+1}(\alpha); q_n(\alpha)).
\end{align*}
\end{proposition}

\begin{corollary}\label{sondpropcor}
For all irrational numbers $\alpha$, one has
\begin{align*}
\mu(\alpha) = \deg\left(\left|\alpha-\frac{p_n(\alpha)}{q_n(\alpha)}\right|; \frac{1}{q_n(\alpha)}\right).
\end{align*}
\end{corollary}

\begin{proof}
This follows readily from the proposition and (\ref{cfineq}).
\end{proof}

\begin{corollary}
For all irrational numbers $\alpha$, one has $\mu(\alpha) = 2$ if and only if
$$\lim_{n \to \infty} \frac{\log q_{n+1}(\alpha) }{\log q_{n}(\alpha) } = 1,$$
if and only if
$$\log q_{n+1}(\alpha) \sim \log q_{n}(\alpha) \ (n \to \infty).$$
\end{corollary}

\begin{corollary}\label{mu0}
For all irrational numbers $\alpha$,  if $\deg a_n(\alpha) < \infty$,  or, more generally,  if $\deg(a_{n}(\alpha);F_n) = 0$ (or, equivalently,  if $\deg(a_{n}(\alpha);e^n) = 0$), then $\mu(\alpha) = 2$.
\end{corollary}

\begin{proof}
Since $F_n$ has exact degree $\infty$,  if  $\deg a_n(\alpha) < \infty$,  then
$$\deg(a_{n}(\alpha);F_n) =  \limsup_{n \to \infty}  \frac{\log a_{n}(\alpha)}{\log n}   \frac{\log n}{\log F_n} = 0.$$
Supposing that $\deg(a_{n}(\alpha);F_n) = 0$,  one has $\deg(a_{n+1}(\alpha);F_{n+1}) = 0$,  and therefore,  since also $q_n(\alpha) \geq F_{n+1}$ for all $n \geq 0$,  one has $\deg (a_{n+1}(\alpha),q_n(\alpha)) = 0$ and therefore
$$\mu(\alpha) =  2+\deg(a_{n+1}(\alpha); q_n(\alpha)) = 2.$$
This completes the proof.
\end{proof}

Recall that the degrees of many important number-theoretic functions are known only within some bound, e.g., one has $\Theta = \deg(\li-\pi) \in [\tfrac{1}{2},1]$.  By the same token, $\mu(\alpha) = \deg \frac{n}{\Vert n\alpha \Vert}$ is unknown for a great many important constants in  analysis.    For example,  the Euler--Mascheroni constant is not known to be irrational, i.e., one may have $\mu(\gamma) = 1$, for all we know.   For the constant $\log 2$, the best we can do thus far is $\mu(\log 2) \in [2,3.57455391]$ \cite{marc}, and,  likewise, $\mu(\pi) \in [2,7.103205334137\ldots]$ for the constant $\pi$ \cite{zei}.   Probably one has $\mu(\alpha) = 2$ for all three of these constants $\alpha$.    Since the Riemann hypothesis is equivalent to $\deg(\li-\pi) = \tfrac{1}{2}$, such conjectures can be regarded, loosely speaking,  as Diophantine approximation analogues of the Riemann hypothesis.

In 1937,  P.\  L\'evy proved the following.

\begin{theorem}[{L\'evy \cite{levy}}]\label{levythm}
 For almost all real numbers $\alpha$, that is, for all $\alpha \in \RR$ outside a set of (Lebesgue) measure $0$, one has
$$\lim_{n \to \infty} q_n(\alpha)^{1/n} = e^{\pi^2/(12 \log 2)},$$
or, equivalently, 
$$\log q_n(\alpha) \sim \frac{\pi^2}{12\log 2} n\ (n \to \infty),$$
where $ \frac{\pi^2}{12\log 2} = 1.186569110415\ldots$.
\end{theorem}

Consequently, one has the following.

\begin{corollary}\label{2cor}
One has
$$\log q_{n+1}(\alpha) \sim \log q_{n}(\alpha) \ (n \to \infty),$$
or, equivalently,  $\mu(\alpha) = 2$,
for almost all $\alpha \in \RR$.
\end{corollary}

By the following corollary, one cannot bound $\mu(\alpha+\beta)$ or $\mu(\alpha \beta)$ generally in terms of $\mu(\alpha)$ and $\mu(\beta)$.

\begin{corollary}
Let $\alpha_0 \in \RR$  be irrational.   Then, for almost all real numbers $\alpha$,  one has $\mu(\alpha) = \mu(\alpha_0-\alpha) = \mu(\alpha_0/\alpha) = 2$, and yet $\mu(\alpha+(\alpha_0-\alpha) )= \mu(\alpha (\alpha_0/\alpha)) = \mu(\alpha_0)$.
\end{corollary}

\begin{proof}
The sets $\{\alpha \in \RR: \mu(\alpha) \neq 2\}$,  $\{\alpha \in \RR: \mu(\alpha_0-\alpha) \neq 2\}$, and $\{\alpha \in \RR: \mu(\alpha_0/\alpha) \neq 2\}$ have measure $0$, and thus so does their union.  Therefore, almost all real numbers $\alpha$ do not lie in any of those three sets.
\end{proof}

Nevertheless, one has the following.

\begin{proposition}\label{muplusrat}
For all $\alpha \in \RR$ and all $r \in \QQ$, one has $\mu(\alpha+r) = \mu(\alpha)$.
\end{proposition}

\begin{proof}
Let $r = \frac{m}{n}$, where $m,n \in \ZZ$ are relatively prime and $n > 0$.  Suppose to obtain a contradiction that $\mu(\alpha+r) > \mu(\alpha)$.  Let $t, u \in \RR$ with $\mu(\alpha+r) > u > t > \mu(\alpha)$.    Since $t > \mu(\alpha)$,  one has $$\left| \alpha - \frac{a}{b} \right| \geq \frac{1}{b^t}$$ for all $a, b \in \ZZ$ with $b \gg 0$.   It follows that
$$\left| \alpha + r - \frac{a}{b} \right|  = \left| \alpha - \frac{na-mb}{nb}  \right|  \geq \frac{1}{(nb)^t} > \frac{1}{b^u}$$
for all $a, b \in \ZZ$ with $b \gg 0$.  But that contradicts $\mu(\alpha+r) > u$.
\end{proof}

For any $\alpha \in \RR$, we let
 $$\underline{\mu}(\alpha) = 1+\liminf_{n \to \infty} \frac{\log q_{n+1}(\alpha)}{\log q_{n}(\alpha)}\index[symbols]{.ua  n@$\underline{\mu}(\alpha)$}$$
if $\alpha$ is irrational,   and we let $\underline{\mu}(\alpha) = 1$  if $\alpha$ is rational.
We call $\underline{\mu}(\alpha)$ the {\bf lower irrationality measure}, {\bf lower approximation order}, or {\bf lower approximation exponent},\index{lower irrationality measure $\mu(\alpha)$}\index{lower approximation order $\mu(\alpha)$}\index{lower approximation exponent $\mu(\alpha)$} {\bf of} $\alpha$.  Note that the obvious analogue of Proposition \ref{sondprop} holds  for $\underline{\mu}(\alpha)$.   One also has
 $$\mu(\alpha) = \inf\left\{t \in \RR: \forall n \gg 0 \ \left| \alpha-\frac{p_n(\alpha)}{q_n(\alpha)} \right| \geq  \frac{1}{q_n(\alpha)^t} \right\},$$
 while
 $$\underline{\mu}(\alpha) = \sup\left\{t \in \RR: \forall n \gg 0 \ \left| \alpha-\frac{p_n(\alpha)}{q_n(\alpha)} \right| \leq  \frac{1}{q_n(\alpha)^t} \right\},$$
 for all irrationals $\alpha$.  Note that $\underline{\mu} (\RR\backslash \QQ) = [2,\infty)$.

Two real numbers are said to be {\bf continued-fraction equivalent}\index{continued-fraction equivalent} if they are both rational or if  they are both irrational and their regular continued fraction expansions have identical tails \cite[Section 9.6]{leve}.   In other words, $\alpha, \beta \in \RR$ are continued-fraction equivalent if and only if $S^{\circ n}(\alpha) = S^{\circ m}(\beta)$ for some nonnegative integers $n$ and $m$.   It is well known   that  $\alpha, \beta \in \RR$ are continued-fraction equivalent if and only if there exist integers $a,b,c,d$ such that $$\beta = \frac{a\alpha +b}{c\alpha + d}$$  and $ad-bc = \pm 1$ \cite[Theorem 9.14]{leve}.   In particular, for any $\alpha \in \RR$, the real numbers $\alpha+n$, and $\frac{1}{\alpha}$ if $\alpha \neq 0$, are continued-fraction equivalent to $\alpha$ for all $n \in \ZZ$.   Moreover, continued-fraction equivalence is the smallest equivalence relation on $\RR$  with that property, since then $\alpha$ is  equivalent  to $S(\alpha) =  \frac{1}{\alpha-\lfloor \alpha \rfloor}$ for all $\alpha \in \RR\backslash \ZZ$.  Thus, a function $f: \RR\backslash \QQ \longrightarrow \overline{\RR}$ satisfies $f(\alpha) = f(\beta)$ for all continued-fraction equivalent irrationals $\alpha, \beta$ if and only if $f(\alpha) = f(\alpha+1) = f(\frac{1}{\alpha})$ for all irrationals $\alpha$, if and only if $f(\alpha) = f(S(\alpha))$ for all irrationals $\alpha$, if and only if the function $f$ factors through the map $S: \RR\backslash \QQ \longrightarrow \RR \backslash \QQ$.  Moreover, since $-\alpha$ is continued-fraction equivalent to $\alpha$ for all $\alpha$, any  such function $f$ is even (and periodic, with $1$ as a period).

\begin{proposition}
For all $\alpha,\beta \in \RR$, if $\alpha$ and $\beta$  are continued-fraction equivalent, then $\mu(\alpha) = \mu(\beta)$ and $\underline{\mu}(\alpha) = \underline{\mu}(\beta)$.
\end{proposition}

\begin{proof}
We may assume without loss of generality that $\alpha$ and $\beta$ are irrational.   We prove the corollary for $\mu(\alpha)$; the proof for $\underline{\mu}(\alpha)$ is similar.  Since $\mu(\alpha) = \mu(\{\alpha\})$, it suffices to show that $\mu(\alpha) = \mu(\tfrac{1}{\alpha})$ if $\alpha \in (0,1)$ is irrational.   Note that, for any irrational $\alpha \in (0,1)$, one has
\begin{align}\label{qpq}
\frac{q_{n+1}(\alpha)}{q_n(\tfrac{1}{\alpha})}  = \frac{p_n(\tfrac{1}{\alpha})}{q_n(\tfrac{1}{\alpha})}  \to \frac{1}{\alpha} \quad \text{as } n \to \infty,
\end{align}
and therefore
$$\mu(\tfrac{1}{\alpha})= 1+\deg(q_{n+1}(\tfrac{1}{\alpha}); q_n(\tfrac{1}{\alpha})) =  1+\deg(q_{n+2}(\alpha); q_{n+1}(\alpha)) = \mu(\alpha).$$
This completes the proof.
\end{proof}

\begin{remark}[Irrationality base]
Let $\alpha\in \RR$ be irrational.   By \cite[Theorem 1]{sondow},   the  {\it irrationality base $\beta(\alpha) \in [1,\infty]$ of $\alpha$} \cite[Definition 3]{sondow} is equivalently given by
\begin{align*}
\beta(\alpha) & = \exp \deg (q_{n+1}(\alpha);  \exp q_{n}(\alpha))  \\
& = \exp \deg q_{n+1}(\alpha)^{\log n/q_{n}(\alpha)} \\ 
& = \exp \limsup_{n \to \infty} \frac{\log q_{n+1}(\alpha)}{q_{n}(\alpha)} \\ 
& = \limsup_{n \to \infty} q_{n+1}(\alpha)^{1/q_{n}(\alpha)}.
\end{align*}
Consequently, if $\mu(\alpha) < \infty$, then $\beta(\alpha)  = 1$ (so, if $\beta(\alpha) > 1$ then  $\mu(\alpha) = \infty$).  The larger the irrationality base of $\alpha$, the ``more extremely well approximable,''  or ``more Liouville,'' $\alpha$ is.
\end{remark}

\section{Markov and relativized Markov constants}

The {\bf Markov constant $M(\alpha)$}\index{Markov constant $M(\alpha)$} of $\alpha \in \RR$ is defined by
\begin{align*}
M(\alpha) & = \inf \left\{ t \in \RR_{>0}: \alpha  \gg_1 \frac{1}{t} n^{-2}\right\} \\ 
 & =  \sup \left\{ t \in \RR_{>0}: \alpha \,  {\not \gg}_1 \, \frac{1}{t} n^{-2}\right\}.\index[symbols]{.ub  o@$M(\alpha)$}
\end{align*}

\begin{lemma}
For all  real numbers $\alpha$,  the Markov constant $M(\alpha)$ is the unique $T \in [0,\infty]$ such that  $\alpha    \gg_1  \frac{1}{ t} n^{-2}$  for all $t > T$ and  $\alpha \,  {\not \gg}_1  \, \frac{1}{ t} n^{-2}$ for all positive $t< T$.
\end{lemma}

From Hurwitz's theorem and Proposition \ref{dioprop},  we obtain the following.

\begin{proposition}
For all $\alpha \in \RR$, one has $M(\alpha) = 0$ if and only if  $\alpha$ is rational,  and $M(\alpha) \geq M(\Phi) = \sqrt{5}$ if and only if $\alpha$ is irrational.
\end{proposition}

Thus,  for all irrational numbers $\alpha$, one has $M(\alpha) \geq M(\Phi) =\sqrt{5}$.  Moreover, equality holds if and only if $\alpha$ is continued-fraction equivalent to the golden ratio  $\Phi = [1,1,1,\ldots]$.  If equality fails,  that is,  if $\alpha$ is not continued-fraction equivalent to $\Phi$,   then one has $M(\alpha) \geq M(\sqrt{2}) = \sqrt{8}$.  Proceeding in this way, one can study the image of the function $M$ on $\RR\backslash \QQ$, which is known as the {\bf Lagrange spectrum},\index{Lagrange spectrum} and which has a remarkably intricate structure \cite{cus} \cite{matheus}.

A real number $\alpha$ is said to be  {\bf badly approximable}\index{badly approximable} if 
 $M(\alpha) < \infty$, and {\bf well approximable}\index{well approximable} if 
 $M(\alpha) = \infty$.  Under this measure,  the smaller $M(\alpha)$ is, the ``more badly approximable'' $\alpha$ is, while, the larger $M(\alpha)$ is, the ``more well approximable'' $\alpha$ is.    Thus, the rational numbers, somewhat ironically, are the ``most badly approximable'' numbers, followed by  the golden ratio $\Phi$ and all of its continued-fraction equivalents, followed by $\sqrt{2}$ and all of its equivalents.

Using (\ref{cfineqa})--(\ref{cfineq}) and Theorem \ref{ccfa},  one can express the Markov constant $M(\alpha)$ of an irrational number  $\alpha$ in terms of its regular continued fraction expansion,  and also in terms of the arithmetic function $\Vert n \alpha \Vert$, as follows.

\begin{theorem}[{\cite[Theorem 9.9 and the proof of Theorem 9.16]{leve}}]\label{markovconst}
For all irrational numbers $\alpha$, one has
\begin{align*}
M(\alpha )   = \limsup_{n \to \infty} \frac{1}{n \Vert n\alpha\Vert} =\limsup_{n\to \infty }  M_n(\alpha),
\end{align*}
where 
\begin{align*}
M_n(\alpha) & = \frac{1}{q_n(\alpha) \Vert q_n(\alpha)\alpha\Vert} \\
& = \frac{1}{ {q_n(\alpha)^2 \left | \alpha- \frac{p_n(\alpha)}{q_n(\alpha)} \right|}} \\
& = S^{\circ (n+1)}(\alpha)+ \frac{q_{n-1}(\alpha)}{q_n(\alpha)} \\
& = [a_{n+1},a_{n+2},a_{n+3},\ldots]+[0,a_{n},a_{n-1},\ldots, a_2, a_{1}]
\end{align*}
for all positive integers $n$,  and where $[a_0, a_1,a_2,\ldots]$ is the regular continued fraction expansion of $\alpha$.
Consequently, if $A(\alpha) = \limsup_{n \to \infty} a_n$, then  $M(\alpha) \in [A(\alpha),A(\alpha)+2]$,  and therefore  $A(\alpha)$ is finite if and only if  the $a_n$ are bounded,  if and only if   $M(\alpha) < \infty$.    Moreover, if $\alpha, \beta \in \RR$  are continued-fraction equivalent, then $M(\alpha) = M(\beta)$.
\end{theorem}

 \begin{corollary}\label{markcor}
 Let $\alpha \in \RR$ be irrational.  Then the following conditions are equivalent.
 \begin{enumerate}
 \item $\alpha$ is badly approximable.
  \item $\alpha \gg_1 Cn^{-2}$ for some $C > 0$.
 \item $\alpha \ggg_1 Cn^{-2}$ for some $C > 0$.
 \item The terms $a_n(\alpha)$ of the regular continued fraction expansion of $\alpha$ are bounded.
 \item  The sequence $\frac{q_{n+1}(\alpha)}{q_n(\alpha)}$ is bounded.
 \end{enumerate}
 Moreover,  if $B(\alpha) = \limsup_{n \to \infty}\frac{q_{n+1}(\alpha)}{q_n(\alpha)}$, then one has $M(\alpha) \in [B(\alpha),B(\alpha)+1]$.
 \end{corollary}

 \begin{corollary}
 If $\alpha \in \RR$ is badly approximable, then $\mu(\alpha) \leq 2$, i.e., $\alpha$ is not very well approximable.    Equivalently, if $\alpha$ is very well approximable, then $\alpha$ is well approximable.
\end{corollary}

 Although almost all real numbers are both well approximable and not very well approximable, all quadratic irrationals, that is, all real algebraic numbers of degree $2$, are badly approximable, as their regular continued fraction expansions are eventually periodic.   On the other hand, it has been conjectured that all algebraic irrationals of degree greater than $2$ are well approximable.   Unfortunately, however,  there is not a single algebraic irrational of degree greater than $2$ that is known to be well approximable,  or to be badly approximable!  Regarding the transcendental number $e$, since
$$e-1 = [1,1,2,1,1,4,1,1,6,1,1,8,\ldots],$$
one has $M(e) = \infty$ and
$\mu(e) = 2$,  that is, $e$ is well approximable but not  very well approximable, by Corollaries \ref{markcor} and \ref{mu0}.    It is conjectured that $\pi$ is also well approximable, since, for example,  all of the terms $a_n(\pi)$ have been computed for $n \leq 3 \cdot 10^{10}$, and the largest term in this range is $a_{15621034283}(\pi) = 52662113289$: see OEIS Sequences A033089 and A033090.  Similar computations suggest that the Euler--Mascheroni constant $\gamma$ is well approximable, and thus irrational: see OEIS Sequences A033091 and A033092.

  If $\alpha \in \RR$ is not a Liouville number,  that is, if $\mu(\alpha)$  is finite, 
then we let  
\begin{align*}
m(\alpha)  & = \inf \left\{ t \in \RR_{>0}: \alpha  \gg_1 \frac{1}{t} n^{-\mu(\alpha)}\right\} \\
& = \sup \left\{ t \in \RR_{>0}: \alpha \,  {\not \gg}_1  \, \frac{1}{t} n^{-\mu(\alpha)}\right\};\index[symbols]{.ub  p@$m(\alpha)$}
\end{align*}
otherwise,  we let $m(\alpha) = \infty$.    We call $m(\alpha)$ the {\bf relativized Markov constant of $\alpha$}.\index{relativized Markov constant $m(\alpha)$}

\begin{lemma}\label{marklem}
For all  real numbers $\alpha$ with $\mu(\alpha) < \infty$,  the relativized Markov constant $m(\alpha)$ is the unique $T \in [0,\infty]$ such that $\alpha    \gg_1  \frac{1}{ t} n^{-\mu(\alpha)}$  for all $t > T$ and
$\alpha \,  {\not \gg}_1  \, \frac{1}{ t} n^{-\mu(\alpha)}$ for all positive $t< T$.
\end{lemma}

For all irrational $\alpha$, one has
$$m(\alpha) \leq M(\alpha),$$ with equality if $\mu(\alpha) = 2$.   Moreover, one has $M(\alpha) = \infty$ if $\mu(\alpha)> 2$.    Thus, $m(\alpha)$  is a more refined invariant than $M(\alpha)$ is.   Note, however, that, unlike $M$,  the function $m$ is not invariant under continued-fraction equivalence, e.g., see Corollary \ref{minva} below.
Nevertheless,  both functions $M: \RR \longrightarrow [0,\infty]$ and $m: \RR \longrightarrow [0,\infty]$ are even and periodic of period $1$.   Moreover,  Proposition \ref{dioprop} implies the following.

\begin{proposition}\label{mmprop}
For all $\alpha \in \QQ$, one has
$$m(\alpha)  =\ord_1 \alpha = \frac{1}{|\alpha|_1}.$$
\end{proposition}

Thus,  $m(\alpha)$ for $\alpha \in \QQ$ provides an {\it analytic} expression for $\ord_1 \alpha$, which was defined algebraically.  Note also that, if $\alpha \in \QQ$ has (finite) regular continued fraction expansion $[a_0, a_1, a_2, \ldots, a_n]$, then $$m(\alpha) =\ord_1 \alpha = q_n(\alpha) = \max_k q_k(\alpha).$$

Theorem \ref{markovconst} generalizes as follows.

\begin{theorem}\label{markk}
For all irrational numbers $\alpha$ with $\mu(\alpha) < \infty$, one has
\begin{align*}
m(\alpha )  & = \limsup_{n \to \infty} \frac{1}{n^{\mu(\alpha)-1} \Vert n\alpha\Vert} =  \limsup_{n\to \infty } m_n(\alpha),
\end{align*}
where 
\begin{align*}
m_n(\alpha)  & = \frac{1}{q_n(\alpha)^{\mu(\alpha)-1} \Vert q_n(\alpha)\alpha\Vert}\\
& = \frac{1}{ {q_n(\alpha)^{\mu(\alpha)} \left | \alpha- \frac{p_n(\alpha)}{q_n(\alpha)} \right|}}  \\
& = q_n(\alpha)^{2-\mu(\alpha)}\left( S^{\circ (n+1)}(\alpha)+ \frac{q_{n-1}(\alpha)}{q_n(\alpha)}\right) \\
 & =  q_n(\alpha)^{2-\mu(\alpha)} ([a_{n+1},a_{n+2},a_{n+3},\ldots]+[0,a_{n},a_{n-1},\ldots, a_2, a_{1}])
\end{align*}
for all positive integers $n$, 
and where $[a_0, a_1,a_2,\ldots]$ is the regular continued fraction expansion of $\alpha$.    Moreover, if $\mu(\alpha) > 2$,  then
$$0< m_n(\alpha) - \frac{q_{n+1}(\alpha)}{q_n(\alpha)^{\mu(\alpha)-1}} < \frac{1}{q_n(\alpha)^{\mu(\alpha)-2}} = o(1)  \ (n \to \infty)$$
for all positive integers $n$ and therefore
\begin{align*}
m(\alpha )  = \limsup_{n \to \infty} \frac{q_{n+1}(\alpha)}{q_n(\alpha)^{\mu(\alpha)-1}} = \limsup_{n \to \infty} \frac{a_{n+1}(\alpha)}{q_n(\alpha)^{\mu(\alpha)-2}}.
\end{align*}
\end{theorem}

\begin{proof}
Let 
$$\kappa_n = \frac{1}{n^{\mu(\alpha)-1} \Vert n\alpha\Vert} $$
for all positive integers $n$.  Since
$$\left|\alpha-\frac{\operatorname{nint}(n\alpha)}{n} \right| = \frac{1}{\kappa_n} n^{-\mu(\alpha)}$$
for all $n$,  we conclude from Lemma \ref{marklem} that
$$m(\alpha) = \limsup_{n \to \infty} \kappa_n = \limsup_{n \to \infty}\frac{1}{n^{\mu(\alpha)-1} \Vert n\alpha\Vert}.$$
By Theorem \ref{markovconst},  to prove the rest of the  theorem, we may suppose without loss of generality that $\mu(\alpha)> 2$.  
By (\ref{cfineq0}), one has
\begin{align*}
 \left|\alpha -\frac{p_n(\alpha)}{q_n(\alpha)} \right| = \frac{ 1}{q_n(\alpha)^2 \left( S^{ \circ(n+1)} (\alpha)+\frac{q_{n-1}(\alpha)}{q_n(\alpha)} \right)} = \frac{1}{m_n(\alpha)}q_n(\alpha)^{-\mu(\alpha)}.
\end{align*}
From Theorem \ref{ccfa} and Lemma \ref{marklem}, it follows, readily, that $$m(\alpha)  = \limsup_{n \to \infty} m_n(\alpha).$$   Moreover, by (\ref{cfineq}), one has
\begin{align*}
q_n(\alpha)^{2-\mu(\alpha)}\frac{q_{n+1}(\alpha)}{q_n(\alpha)}< m_n(\alpha)  <q_n(\alpha)^{2-\mu(\alpha)}\left(1+\frac{q_{n+1}(\alpha)}{q_n(\alpha)}  \right)
\end{align*}
for all positive integers  $n$.   Since also
\begin{align*}
q_n(\alpha)^{2-\mu(\alpha)} \frac{q_{n+1}(\alpha) }{ q_n(\alpha)} &  =  q_n(\alpha)^{2-\mu(\alpha)}\left(a_{n+1}(\alpha) +  \frac{q_{n-1}(\alpha) }{ q_n(\alpha)}  \right)  \\
& =  q_n(\alpha)^{2-\mu(\alpha)}a_{n+1}(\alpha) + o(1) \ (n \to \infty),
\end{align*}
the theorem follows.
\end{proof}

Note that, by Corollary \ref{degmu}, the function $\frac{1}{n^{\mu(\alpha)-1} \Vert n\alpha\Vert}$ has degree $0$. 
Moreover,  the degree of both $m_n(\alpha)$ and $ \frac{q_{n+1}(\alpha)}{q_n(\alpha)^{\mu(\alpha)-1}}$ with respect to $q_n(\alpha)$ is also $0$.

\begin{remark}[Relativized Markov constant as an upper leading coefficient]
 It follows from Corollary \ref{degmu} and Theorem \ref{markk} that the relativized Markov constant $m(\alpha)$ of an irrational $\alpha$ with $\mu(\alpha) < \infty$ is precisely  the {\it upper leading coefficent} $$m(\alpha) = \operatorname{lc}^{+} \frac{n}{\Vert n \alpha \Vert} = \operatorname{lc}^{+} \frac{1}{\Vert n \alpha \Vert}$$ of the arithmetic functions $\frac{n}{\Vert n \alpha \Vert}$ and $\frac{1}{\Vert n \alpha \Vert}$, in the sense of Remark \ref{lc}. 
\end{remark}

\begin{corollary}\label{minva}
The function $m$ is even and periodic of period $1$. 
For any irrational $\alpha \in \RR$ with $\mu(\alpha) < \infty$, one has
$$m(\tfrac{1}{\alpha}) = |\alpha|^{2-\mu(\alpha)} m(\alpha)$$
and
$$m(S(\alpha)) = |S(\alpha)|^{\mu(\alpha)-2} m(\alpha),$$
and therefore $m(\tfrac{1}{\alpha})  = m(\alpha)$ if and only if  $m(S(\alpha)) = m(\alpha)$,  if and only if $\mu(\alpha) = 2$ or $m(\alpha) \in \{0,\infty\}$.
In particular, for any irrational $\alpha$ with $\mu(\alpha) < \infty$ and $0<m(\alpha) < \infty$, the irrationality measure $\mu(\alpha)$ of $\alpha$ can be recovered from $m(\alpha)$ and either of $m(\tfrac{1}{\alpha})$ or $m(S(\alpha))$ by $$\mu(\alpha) =  2+ \log_{|\alpha|}  \left(\frac{m(\alpha)}{m(\tfrac{1}{\alpha})} \right) = 2+ \log_{|S(\alpha)|}  \left(\frac{m(S(\alpha))}{m(\alpha)} \right).$$
\end{corollary}

\begin{corollary}
Let $\alpha$ be an irrational  number with $\mu(\alpha)<\infty$.   One has
$$ \frac{q_{n+1}(\alpha)}{q_{n}(\alpha)^{\mu(\alpha)-1}} \asymp 1 \ (n \to \infty)$$
if and only if
$$  \frac{\log q_{n+1}(\alpha)}{\log q_{n}(\alpha)} +1-\mu(\alpha) = O\left( \frac{1}{\log q_{n}(\alpha)} \right)  \ (n \to \infty).$$  Moreover, if the equivalent conditions above hold, then one has 
$$0<m(\alpha)< \infty,$$
$$  \frac{\log q_{n+1}(\alpha)}{\log q_{n}(\alpha)} +1-\mu(\alpha) = O\left( \frac{1}{n}\right) \ (n \to \infty),$$  
and
$$\mu(\alpha) = 1+   \lim_{n \to\infty} \frac{\log q_{n+1}(\alpha)}{\log q_{n}(\alpha)}.$$
\end{corollary}

\begin{proof}
The stated equivalence follows from the identity $$ \frac{\log q_{n+1}(\alpha)}{\log q_{n}(\alpha)} +1-\mu(\alpha)= \frac{\log \frac{ q_{n+1}(\alpha)}{q_{n}(\alpha)^{\mu(\alpha)-1}}}{\log q_{n}(\alpha)} $$ and the fact that, for any real function $f$, one has $f(x) \asymp 1 \ (x \to \infty)$ if and only if $\log |f(x)| = O(1) \ (x \to \infty)$.  The rest of the proposition then follows from Theorem \ref{markk} and the fact that $\log q_n(\alpha) \geq \log F_{n+1} \asymp n \ (n \to \infty)$.
\end{proof}

\begin{corollary}\label{badapp}
Let $\alpha$ be an irrational  number.  Then $\alpha$ is badly approximable if and only if 
$$ \frac{q_{n+1}(\alpha)}{q_{n}(\alpha)} \asymp 1 \ (n \to \infty),$$ if and only if
$$  \frac{\log q_{n+1}(\alpha)}{\log q_{n}(\alpha)} -1 \asymp \frac{1}{\log q_{n}(\alpha)} \ (n \to \infty).$$
Moreover,  if the equivalent conditions above hold,  then one has $\log q_n(\alpha) \asymp  n \ (n \to \infty)$ and therefore
$$ \frac{\log q_{n+1}(\alpha)}{\log q_{n}(\alpha)} -1  \asymp \frac{1}{n} \ (n \to \infty).$$
\end{corollary}

The following lemma is useful for constructing irrationals $\alpha$ with prescribed values of $\mu(\alpha) > 2$ and $m(\alpha)$.

\begin{lemma}\label{prescribeda}
Let $f: \ZZ_{>0} \longrightarrow \RR_{>0}$ be a positive real-valued arithmetic function.  Then there exists a unique (irrational) real  number $\alpha$ with $a_0(\alpha) = 0$ and $$a_{n+1}(\alpha) = \lceil f(q_n(\alpha)) \rceil$$
 for all nonnegative integers $n$.    Moreover, if  $\lim_{n \to \infty} f(n) = \infty$, then  one has
 $$\frac{q_{n+1}(\alpha) }{ q_n(\alpha)} \sim a_{n+1}(\alpha) \sim f(q_n(\alpha)) \ (n \to \infty).$$
\end{lemma}

\begin{proof}
The $p_n(\alpha)$ and $q_n(\alpha)$, and thus $\alpha = \lim_{n \to \infty} \frac{p_n(\alpha)}{q_n(\alpha)}$, are uniquely determined by the initial conditions $p_0(\alpha) = 0$, $q_0(\alpha) = 1$,  $p_1(\alpha) = 1$,   $q_1(\alpha) = \lceil f(1) \rceil$, and the recurrence relations
$$p_{n+1}(\alpha)  =\lceil f(q_n(\alpha)) \rceil p_{n}(\alpha) +  p_{n-1}(\alpha)$$
and
$$q_{n+1}(\alpha)  =\lceil f(q_n(\alpha)) \rceil q_{n}(\alpha) +  q_{n-1}(\alpha).$$
Suppose that $\lim_{n \to \infty} f(n) = \infty$.  Then one has
\begin{align*}
\frac{q_{n+1}(\alpha) }{ q_n(\alpha)} =  a_{n+1}(\alpha) + \frac{q_{n-1}(\alpha) }{ q_n(\alpha)}  \sim a_{n+1}(\alpha) \sim f(q_n(\alpha)) \ (x \to \infty).
\end{align*}
This completes the proof.
\end{proof}

By Theorem \ref{markk} and the lemma above, one has the following.

\begin{corollary}\label{mcor}
Let $\mu >2$ and $c > 0$.  
\begin{enumerate}
\item  Let $\alpha$ be the unique real number with  $a_0(\alpha) = 0$ and
 $$a_{n+1}(\alpha) = \lceil c q_n(\alpha)^{\mu-2} \rceil$$
 for all nonnegative integers $n$.  Then one has
 $$q_{n+1}(\alpha)  \sim c q_n(\alpha)^{\mu-1} \ (n \to \infty).$$
 Consequently, one has $\mu(\alpha) = \mu$ and $m(\alpha) = c$.
 \item Let $\alpha $ be the unique real number with  $a_0(\alpha) = 0$ and
 $$a_{n+1}(\alpha) = \lceil n q_n(\alpha)^{\mu-2} \rceil$$
 for all nonnegative integers $n$.  Then one has
 $$q_{n+1}(\alpha)  \sim n q_n(\alpha)^{\mu-1} \ (n \to \infty).$$
 Consequently, one has $\mu(\alpha) = \mu$ and $m(\alpha) = \infty$.
 \item Let $\alpha$ be the unique real number with  $a_0(\alpha) = 0$ and
 $$a_{n+1}(\alpha) = \lceil n^{-1} q_n(\alpha)^{\mu-2} \rceil$$
 for all nonnegative integers $n$.  Then one has
 $$q_{n+1}(\alpha)  \sim n^{-1} q_n(\alpha)^{\mu-1} \ (n \to \infty).$$
 Consequently, one has $\mu(\alpha) = \mu$ and $m(\alpha) = 0$.
 \end{enumerate}
\end{corollary}

\begin{corollary}\label{mmcor}
The image  $m(\mu^{-1}(1))$ of  the function $m$ on $\mu^{-1}(1) = \QQ$ is equal to $ \ZZ_{>0}$; the image $m(\mu^{-1}(2))$  of $m$ on $\mu^{-1}(2)$ is the Lagrange spectrum together with $\infty$; and, for all $t >2$, the image $m(\mu^{-1}(t))$  of $m$ on  $\mu^{-1}(t)$ is equal to $[0,\infty]$.
\end{corollary}

The marked contrast between the values of $m$ on $\mu^{-1}(t)$ for $ t > 2$ versus $t = 1,2$ is due to the fact that $t \in \RR$ is in the interior of $\im \mu$ if and only if $t > 2$.

\begin{problem}
 What is the image of the map
$\deg\left(\frac{\log q_{n+1}(-)}{\log q_{n}(-)} - 1\right)$ on $\mu^{-1}(2)$? 
\end{problem}

\section{Logexponential irrationality degree}

In this section, we introduce a  measure of the rational approximability of a real number,  which  we call {\it logexponential irrationality degree}, that is far more refined than irrationality measure, in the same way that logexponential degree is far more refined than degree.   Relevant to this generalization is the following conjecture, due to Lang,  which he proposed as a generalization of  Roth's theorem.

\begin{conjecture}[{Lang \cite[p.\ 98]{lang}}]\label{langconj}
Let $\alpha$ be an algebraic irrational number.  Then, for every $t > 1$, 
one has
$$\alpha \gg_1  n^{-2} (\log (n+1))^{-t}.$$
\end{conjecture}

We use $n+1$ above instead of $n$ to formulate Lang's conjecture so that the given function is defined at all positive integers.  For many purposes,  if $f$ is a function in $\mathbb{L}_{>0}$,  then one replace $f(x)$ by $f(x+c)$ for some sufficiently large $c >0$ to ensure that the resulting function is defined and positive on $[1,\infty)$.  

Lang's conjecture motivates the following definition:  for all $\alpha \in \RR$,  we let 
\begin{align*}
\pmb{\mu}( \alpha) & = \inf\{\dege f: f \in \mathbb{L}_{>0} \text{ and } \alpha \gg_1  1/f|_{\ZZ_{>0}}\} \\
& =  \inf\{\dege f: f \in \mathbb{L}_{>0} \text{ and } |\alpha -r| \geq 1/f(\ord_1 r) \text{ for all but finitely many } r \in \QQ\},\index[symbols]{.ua  o@$\pmb{\mu}(\alpha)$} 
\end{align*}
where the infima are computed in $\prod_{n = 0}^{ \infty*}\overline{\RR}$.  For all nonnegative integers $k$, we let $\pmb{\mu}_k (\alpha)$ denote the $k$th coordinate of $\pmb{\mu}( \alpha)$.  We call $\pmb{\mu}(\alpha)$ the {\bf logexponential irrationality degree of $\alpha$}.\index{logexponential irrationality degree $\pmb{\mu}(\alpha)$}   This notation and terminology are justified by the following proposition.

\begin{proposition}\label{refine}
For all $\alpha \in \RR$, one has
\begin{align*}
\mu( \alpha)= \pmb{\mu}_0 (\alpha) = \inf\{\deg f: f \in \mathbb{L}_{>0} \text{ and } \alpha \gg_1  1/f|_{\ZZ_{>0}}\}.
\end{align*}
\end{proposition} 

\begin{proof}
Since $\deg x^t = t$ and $x^t \in \mathbb{L}_{>0}$, it is  clear from the definitions of $\mu( \alpha)$ and $\pmb{\mu}( \alpha)$ that 
\begin{align*}
\mu( \alpha) \geq \inf\{\deg f: f \in \mathbb{L}_{>0} \text{ and } \alpha \gg_1  1/f|_{\ZZ_{>0}}\} = \pmb{\mu}_0 (\alpha).
\end{align*}
Suppose to obtain a contradiction that the inequality is strict.  Then there exists an $f \in  \mathbb{L}_{>0}$ such that  $\alpha \gg_1  1/f|_{\ZZ_{>0}}$ and $\mu( \alpha) > \deg f$.  We may let $t \in \RR$ with  $\mu( \alpha) > t > \deg f$.   Then,  by the definition of $\mu(\alpha)$,  one has $\alpha \not \gg_1  n^{-t}$.  But $f(x) = o(x^t) \ (x \to \infty)$, so $f(x) \leq x^t$, whence $1/f(n) \geq n^{-t}$,  for all $x, n \gg 0$, and therefore $\alpha \not \gg_1 1/f|_{\ZZ_{>0}}$, which is our desired contradiction.
\end{proof}

\begin{corollary}
Roth's theorem is equivalent to 
$$\pmb{\mu}( \alpha) \leq (2,\infty,1,0,0,0,\ldots)$$
for all real algebraic numbers $\alpha$.  Moreover,  Lang's conjecture, Conjecture \ref{langconj},  is equivalent to $\pmb{\mu}_1 (\alpha) \leq 1$, and thus also to
$$\pmb{\mu} (\alpha) \leq (2,1,\infty,1,0,0,0,\ldots),$$
for all real algebraic  numbers $\alpha$.  
\end{corollary}

 The following theorem provides several  more intuitive descriptions of logexponential irrationality degre, any of which we could have taken as the definition.

\begin{theorem}\label{mupropo}
Let $\alpha \in \RR$.  
\begin{enumerate}
\item  $\pmb{\mu} (\alpha)$  is the smallest $\dd \in \prod_{n = 0}^{ \infty*} \underline{\RR}$ such that, for any positive arithmetic function $f: \ZZ_{>0} \longrightarrow \RR_{>0}$ with $\dege f < -\dd$, one has $\alpha \gg_1 f$.   
\item  $\pmb{\mu} (\alpha)$  is the largest $\dd \in \prod_{n = 0}^{ \infty*} \underline{\RR}$ such that, for any positive arithmetic function $f: \ZZ_{>0} \longrightarrow \RR_{>0}$ with $ \underline{\dege}  \, f >  -\dd$, one has $\alpha \not \gg_1 f$. 
\item $\pmb{\mu} (\alpha)$  is the unique $\dd \in \prod_{n = 0}^{ \infty*} \underline{\RR}$ such that, for any positive arithmetic function $f: \ZZ_{>0} \longrightarrow \RR_{>0}$, one has $\alpha \gg_1 f$ if  $\dege f < -\dd$, and $\alpha  \not \gg_1 f$ if  $ \underline{\dege}\, f > -\dd$.  
\item $\pmb{\mu} (\alpha)$  is the unique $\dd \in \prod_{n = 0}^{ \infty*} \underline{\RR}$ such that, for any positive arithmetic function $f: \ZZ_{>0} \longrightarrow \RR_{>0}$, one has  $\dege f \geq -\dd$ if $\alpha \not \gg_1 f$, and  $ \underline{\dege} \, f \leq -\dd$ if $\alpha  \gg_1 f$.
\item One has
\begin{align*}
\pmb{\mu} (\alpha) & = -\inf\{ \dege f: f: \ZZ_{>0}  \longrightarrow \RR_{>0} \text{ and } \alpha \not \gg_1  f\} \\ & = -\sup\{\underline{\dege}\, f : f: \ZZ_{>0}  \longrightarrow \RR_{>0} \text{ and } \alpha \gg_1  f\}.
\end{align*}
\item  $\mu (\alpha)$  is the smallest $d \in   \underline{\RR}$ such that, for any positive arithmetic function $f: \ZZ_{>0} \longrightarrow \RR_{>0}$ with $\deg f < -d$, one has $\alpha \gg_1 f$.   
\item  $\mu (\alpha)$  is the largest $d \in   \underline{\RR}$ such that, for any positive arithmetic function $f: \ZZ_{>0} \longrightarrow \RR_{>0}$ with $ \underline{\deg} f >  -d$, one has $\alpha \not \gg_1 f$. 
\item $\mu (\alpha)$  is the unique $d \in \underline{\RR}$ such that, for any positive arithmetic function $f: \ZZ_{>0} \longrightarrow \RR_{>0}$, one has $\alpha \gg_1 f$ if  $\deg f < -d$, and $\alpha  \not \gg_1 f$ if  $ \underline{\deg}\, f > -d$.  
\item $\mu$  is the unique $d \in \underline{\RR}$ such that, for any positive arithmetic function $f: \ZZ_{>0} \longrightarrow \RR_{>0}$, one has  $\deg f \geq -d$ if $\alpha \not \gg_1 f$, and  $ \underline{\deg} \, f \leq -d$ if $\alpha  \gg_1 f$.
\item One has
\begin{align*}
\mu (\alpha) & = -\inf\{\deg f: f: \ZZ_{>0}  \longrightarrow \RR_{>0} \text{ and } \alpha \not \gg_1  f\} \\
& = -\sup\{ \underline{\deg}\, f: f: \ZZ_{>0}  \longrightarrow \RR_{>0} \text{ and } \alpha  \gg_1  f\}.
\end{align*}
\end{enumerate}  
\end{theorem}

\begin{proof}
 Suppose to obtain a contradiction that $\dege  f < - \pmb{\mu} (\alpha)$ and $\alpha \not \gg_1 f$ for some positive arithmetic function $f$.  Then, for any   $g \in \mathbb{L}_{>0}$ with $\dege g < -\dege f $ one has $ \dege  f < \dege(1/g)$ and therefore $f(x) = o(1/g(x)) \ (x \to \infty)$, whence  $\alpha \not \gg_1 1/g|_{\ZZ_{>0}}$.   By the definition of $\pmb{\mu} (\alpha)$, then, it follows that $\pmb{\mu} (\alpha) \geq -\dege f >  \pmb{\mu} (\alpha)$, which is a contradiction.  Therefore,   $\dege  f < - \pmb{\mu} (\alpha)$ implies $\alpha \gg_1 f$.    Suppose, on the other hand,  that $\dd < \pmb{\mu} (\alpha)$, where $\dd \in \prod_{n = 0}^{ \infty*}\overline{\RR}$.  Choose $r \in \mathbb{L}_{>0}$, defined and positive on $[1,\infty)$,  so that  $\dd< \dege r < \pmb{\mu} (\alpha)$.   Then $ \dege (1/r) < -\dd$  and also $\alpha \not \gg_1 1/r|_{\ZZ_{>0}}$ by definition of $\pmb{\mu} (\alpha)$, whence $\dege f < - \dd$ does not imply $\alpha \gg_1 f$.  This proves statement (1).

Now suppose to obtain a contradiction that $\underline{\dege} \, f > - \pmb{\mu} (\alpha)$ and $\alpha \gg_1 f$ for some positive arithmetic function $f$.  Then, for any   $g \in \mathbb{L}_{>0}$ with $\dege g > -\underline{\dege} \, f$ one has $ \dege(1/g) <  \underline{\dege} \, f$ and therefore $1/g(x) = o(f(x)) \ (x \to \infty)$, whence  $\alpha  \gg_1 1/g|_{\ZZ_{>0}}$.   By the definition of $\pmb{\mu} (\alpha)$, then, it follows that $\pmb{\mu} (\alpha) \leq -\underline{\dege} \, f <  \pmb{\mu} (\alpha)$, which is a contradiction.  Therefore,   $\underline{\dege} \, f > - \pmb{\mu} (\alpha)$ implies $\alpha \not \gg_1 f$.    Suppose, on the other hand,  that $\dd > \pmb{\mu} (\alpha)$, where $\dd \in \prod_{n = 0}^{ \infty*}\overline{\RR}$.  Choose $r \in \mathbb{L}_{>0}$, defined and positive on $[1,\infty)$,  so that  $\dd>\dege r > \pmb{\mu} (\alpha)$.   Then $ \dege (1/r) > -\dd$,  and also $\dege(1 /r) < -\pmb{\mu} (\alpha)$  and therefore $\alpha \gg_1 1/r|_{\ZZ_{>0}}$ by statement (1), whence $\underline{\dege} \, f > - \dd$ does not imply $\alpha \not \gg_1 f$.  This proves statement (2).
 
 Statements (3)--(5) follow immediately from statements (1) and (2).  Finally,  the proofs of statements (1)--(5) can be modified appropriately to yield statements (6)--(10).
\end{proof}

\begin{corollary}\label{pmbmucor}
Let $\alpha \in \RR$.   
\begin{enumerate}
\item  $\pmb{\mu}( \alpha)$ is the unique $\dd \in \prod_{n = 0}^{ \infty*}\overline{\RR}$ such that,   for all $f \in  \mathbb{L}_{>0}$ defined and positive on $[1,\infty)$, one has $ \alpha  \gg_1  1/f|_{\ZZ_{>0}}$ if $\dege f > \dd$, and $ \alpha \not \gg_1  1/f|_{\ZZ_{>0}}$  if $\dege f < \dd$.
Consequently, one has
\begin{align*}
\pmb{\mu}( \alpha) & = \inf\{\dege f: f \in \mathbb{L}_{>0} \text{ and } \alpha \gg_1  1/f|_{\ZZ_{>0}}\} \\
& = \sup\{\dege f: f \in \mathbb{L}_{>0} \text{ and } \alpha  \not \gg_1  1/f|_{\ZZ_{>0}}\}.
\end{align*}
\item Suppose that $\alpha$ is irrational.  Then $\pmb{\mu}( \alpha)$ is the unique $\dd \in \prod_{n = 0}^{ \infty*}\overline{\RR}$ such that, for all $f \in  \mathbb{L}_{>0}$  defined and positive on $[1,\infty)$, one has
$ \alpha  \succ_1  1/f|_{\ZZ_{>0}}$  if $\dege f > \dd$, and $ \alpha \not \succ_1 1/f|_{\ZZ_{>0}}$ if $\dege f < \dd$.     Consequently, one has
\begin{align*}
\pmb{\mu}( \alpha) & = \inf\{\dege f: f \in \mathbb{L}_{>0} \text{ and } \alpha \succ_1  1/f|_{\ZZ_{>0}}\} \\
& = \sup\{\dege f: f \in \mathbb{L}_{>0} \text{ and } \alpha  \not \succ_1  1/f|_{\ZZ_{>0}}\}.
\end{align*}
\end{enumerate}
\end{corollary}

By Theorem \ref{mupropo},  the relation $\pmb{\mu}( \alpha) < \pmb{\mu}(\beta)$ means in a precise sense that $\beta$ is ``better'' approximated by rational numbers than $\alpha$ is.  Indeed, we have the following corollary.

\begin{corollary}\label{orderedmu}
Let $\alpha, \beta \in \RR$.    One has $\pmb{\mu}(\alpha)< \pmb{\mu}( \beta)$ if and only there exist  $\dd < \ee \in  \prod_{n = 0}^{ \infty*}\overline{\RR}$ such that $\alpha \gg_1 f$ and $\beta \not \gg_1 f$  for all  positive arithmetic functions $f: \ZZ_{>0} \longrightarrow \RR_{>0}$ with $\dd \leq \underline{\dege} \,  f \leq \dege f  \leq  \ee$, if and only if there exist $r, s \in \mathbb{L}_{>0}$ such that $\dege r < \dege s$,  $\alpha \gg_1 s|_{\ZZ_{>0}}$, and $\beta \not \gg_1 r|_{\ZZ_{>0}}$.  Consequently, if $\pmb{\mu}( \alpha)< \pmb{\mu}(\beta)$,  then, for all  positive arithmetic functions $f: \ZZ_{>0} \longrightarrow \RR_{>0}$ of exact logexponential degree, if  $\beta \gg_1 f$, then $\alpha \gg_1 f$.
\end{corollary}

Note that the conditions equivalent to $\pmb{\mu}( \alpha) < \pmb{\mu}(\beta)$ expressed in the corollary above make no reference to $\pmb{\mu}$, much like the standard definition of $|A| \leq |B|$ for sets $A$ and $B$ makes no reference to the cardinality.  By  contrast,  it is not so easy to define the relations $\dege f < \dege g$ and $\dege f <  \underline{\dege} \, g$ without reference to $\dege$.

\begin{corollary}\label{simplemu}
Let $\alpha \in \RR$, and let $t = \mu(\alpha)$.   One has
$$\pmb{\mu}(\alpha)    \left.
 \begin{cases}
    \leq (t,0,0,0,\ldots) & \text{if $m(\alpha)  = 0$} \\
    =  (t,0,0,0,\ldots)  & \text{if $0  < m(\alpha) < \infty$} \\
      \geq   (t,0,0,0,\ldots)  & \text{if $ m(\alpha) = \infty$.}
 \end{cases}
\right.$$ 
Thus, if $\pmb{\mu}(\alpha) < (t,0,0,0,\ldots)$, then $m(\alpha)  = 0$, and if $\pmb{\mu}(\alpha) > (t,0,0,0,\ldots)$, then $m(\alpha)  = \infty$.
\end{corollary}

\begin{proof}
Suppose that $m(\alpha) <\infty$, so that
$$ \alpha  \gg_1 \frac{1}{\lambda} n^{-\mu(\alpha)}$$
for some  $\lambda >0$.   Then, by Corollary  \ref{pmbmucor}(1), one has $$\pmb{\mu} (\alpha) \leq \dege \lambda n^{\mu(\alpha)} = (\mu(\alpha),0,0,0,\ldots).$$
Likewise, if $m(\alpha) >0$, then
$$ \alpha \,  {\not \gg}_1  \, \frac{1}{\lambda} n^{-\mu(\alpha)}$$
 for some  $\lambda >0$, and therefore
 $$\pmb{\mu} (\alpha) \geq \dege \lambda n^{\mu(\alpha)} =(\mu(\alpha),0,0,0,\ldots).$$
The corollary follows.
\end{proof}

\begin{corollary}
Let $\alpha \in \RR$.   One has
$$\pmb{\mu}(\alpha) = (1,0,0,0,\ldots)$$
if and only if $\alpha \in \QQ$, and one has
$$\pmb{\mu}(\alpha) \geq (2,0,0,0,\ldots)$$
if and only if $\alpha \in \RR\backslash \QQ$, with equality if  $\alpha \in \RR\backslash \QQ$ is  badly approximable (e.g., if $\alpha$ is a quadratic irrational).
\end{corollary}

By Corollaries \ref{orderedmu}  and \ref{simplemu},  the real numbers  $\alpha$ for which $\pmb{\mu}(\alpha) = (2,0,0,0,\ldots)$ are in a precise  sense the  ``most irrational'' of all of the irrational numbers.  Of those, the badly approximable numbers are the worst, and they themselves are ordered in ``irrationality'' inversely by the Markov constant $M(\alpha)$.

This following problem currently is very far from being solved.  

\begin{outstandingproblem}
Compute  $\pmb{\mu}(\alpha)$ for any (or all) real algebraic numbers $\alpha$ of degree greater than $2$, e.g., for $\alpha = \sqrt[3]{2}$.
\end{outstandingproblem}

In fact, the current situation is rather bleak: for all we know, all algebraic numbers $\alpha$ could be badly approximable and therefore satisfy $\pmb{\mu}(\alpha) \leq  (2,0,0,0,\ldots)$.  

The following theorem expresses $\pmb{\mu}(\alpha)$ for any irrational number $\alpha$ in terms of the arithmetic functions $\Vert n \alpha \Vert$ and $$\frac{\Vert n \alpha \Vert }{n} = \min \left\{| \alpha-r|: \text{$r \in\QQ$ and $\ord_1 r \mid n$}\right\},$$

\begin{theorem}\label{pmbmu}
If $\alpha \in \RR$ is irrational, then one has 
\begin{align*}
\pmb{\mu}( \alpha) =  \dege \frac{n}{\Vert n \alpha \Vert} =
 \dege \frac{1}{\Vert n \alpha \Vert} + (1,0,0,0,\ldots).
 \end{align*}
\end{theorem}

\begin{proof}
By Proposition \ref{basicrelations}(4), the condition $\alpha \succ_1 1/ f|_{\ZZ_{>0}}$ is equivalent to 
$\frac{n}{\Vert n \alpha \Vert} \leq f(n) $ for all $n \gg 0$.   Therefore,  by Corollary \ref{pmbmucor}(2) and Theorem \ref{infpropexp},  one has
\begin{align*}
 \pmb{\mu}(\alpha) & = \inf\{\dege f: f \in \mathbb{L}_{>0} \text{ and } \alpha \succ_1 1/f|_{\ZZ_{>0}}\} \\
&  = \inf\left\{\dege f: f \in \mathbb{L}_{>0},  \, \forall n \gg 0 \, \left(  \frac{n}{\Vert n \alpha \Vert}  \leq f(n) \right) \right \} \\
& = \dege \frac{n}{\Vert n \alpha \Vert}.
\end{align*}
This completes the proof.
\end{proof}

Note that, since the image of $\Vert n \alpha \Vert$ is dense in $[0,\frac{1}{2}]$, one has $\dege \Vert n \alpha \Vert = (0,0,0,\ldots)$, whence $\underline{\dege}\, \frac{n}{ \Vert n \alpha \Vert} = (1,0,0,0,\ldots)$, for all irrationals $\alpha$.   

\begin{corollary}
Let $\alpha \in \RR$ be irrational.   One has
\begin{align*}
\pmb{\mu}(\alpha)  = \dege \frac{1}{\mu_\alpha(x)} = -\underline{\dege}\, \mu_\alpha(x),
\end{align*}
where 
\begin{align*}
\mu_\alpha(x) & =  \min\left \{\frac{\Vert n\alpha \Vert }{n} : n \leq x \right \} \\
& =  \min\left\{ | \alpha-r|: \text{$r \in\QQ$ and $\ord_1 r \leq x$} \right\} \\
& =  \min\left\{ \left| \alpha-\frac{a}{b}\right|: \text{$a,b \in\ZZ$ and $1 \leq b \leq x$} \right\} 
\end{align*}
and
\begin{align*}
\frac{1}{\mu_\alpha(x)}  & = \max\left \{\frac{n}{\Vert n\alpha \Vert } : n \leq x \right \} \\
& =  \max\left\{ \frac{1}{| \alpha-r|}: \text{$r \in\QQ$ and $\ord_1 r \leq x$} \right\} \\
& =  \max\left\{ \frac{1}{\left| \alpha-\frac{a}{b}\right|}: \text{$a,b \in\ZZ$ and $1 \leq b \leq x$} \right\}
\end{align*}
on $[1,\infty)$.
\end{corollary}

\begin{proof}
This follows from the theorem and Propositions \ref{denomschar}(2) and \ref{supprop}.
\end{proof}

The following theorem expresses the logexponential irrationality degree of any irrational number  $\alpha$ in terms of the regular continued fraction expansion of $\alpha$.

\begin{theorem}\label{pmbmutheorem}
Let $\alpha \in \RR$ be irrational,  let $D = \{q_n(\alpha): n \in \ZZ_{\geq 0}\}$, and let $F_1, F_2,G,H_1,H_2: D \longrightarrow [1,\infty)$ be defined, respectively, by
$$F_1(q_n(\alpha)) = \left| \alpha-\frac{p_n(\alpha)}{q_n(\alpha)}\right|^{-1},$$
$$F_2(q_n(\alpha)) = q_n(\alpha)q_{n+1}(\alpha),$$
$$G(q_n(\alpha)) = q_{n+1}(\alpha),$$
$$ H_1(q_n(\alpha)) = \frac{q_{n+1}(\alpha)}{q_n(\alpha)} = a_{n+1}(\alpha)+\frac{q_{n-1}(\alpha)}{q_n(\alpha)},$$
and
$$H_2(q_n(\alpha)) =  a_{n+1}(\alpha) =\left \lfloor \frac{q_{n+1}(\alpha)}{q_n(\alpha)} \right  \rfloor$$
for all nonnegative integers $n$.
Then one has
\begin{align*}
\pmb{\mu}(\alpha)  & =  \dege F_i  \\ & =  \dege  G+ (1,0,0,0,\ldots) \\ & 
= \dege H_i + (2,0,0,0,\ldots)
\end{align*}
for $i = 1,2$.  Moreover,  one has
$$M(\alpha) =  \limsup_{n \to \infty} \frac{F_1(q_n(\alpha))}{q_n(\alpha)^2} ,$$
and, if $\mu(\alpha) > 2$, then  one has 
\begin{align*}
m(\alpha) = \limsup_{n \to \infty} \frac{F_i(q_n(\alpha))}{q_n(\alpha)^{\mu(\alpha)}} = \limsup_{n \to \infty} \frac{G(q_n(\alpha))}{q_n(\alpha)^{\mu(\alpha)-1}}    =\limsup_{n \to \infty} \frac{H_i(q_n(\alpha))}{q_n(\alpha)^{\mu(\alpha)-2}}
\end{align*}
for $i = 1,2$.
\end{theorem}

\begin{proof}
Let $f  \in \mathbb{L}_{>0}$ be defined on $[1,\infty)$ with $f(n) \geq 2 n^2$ for all $n \gg 0$, and consider the condition that
$$\frac{n}{\Vert n \alpha \Vert} \leq f(n)$$ for all $n \gg 0$.
Note that,  for all $n \gg 0$, the inequality
$$\frac{n}{\Vert n \alpha \Vert} >  f(n)$$ implies
$$\frac{\Vert n \alpha \Vert}{n} < \frac{1}{2n^2}$$
and therefore,  by Theorem \ref{ccfa},  that $n$ is the denominator of a convergent of the regular continued fraction expansion of $\alpha$.
Thus,  by Proposition \ref{denomschar}(4), the condition under consideration is equivalent to
$$ F_1(q_k(\alpha)) \leq f(q_k(\alpha))$$ 
for all $k \gg 0$.  Also by Theorem \ref{ccfa}, one has $F_1(q_k(\alpha))> q_k(\alpha)^2$ for all $k$ and therefore
$$\dege F_1  \geq (2,0,0,0,\ldots).$$
By Theorem \ref{pmbmu}(3), Lemma \ref{difflemexp}, and Theorem \ref{infpropexp}, then, one has
\begin{align*}
\pmb{\mu}( \alpha) & =  \dege \frac{n}{\Vert n \alpha \Vert} \\
& =   \dege \max\left (\frac{n}{\Vert n \alpha \Vert},2n^2 \right) \\
&  = \inf\left\{\dege f: f \in \mathbb{L}_{>0},  \, \forall n \gg 0 \, \left( \max\left( \frac{n}{\Vert n \alpha \Vert},2n^2\right)  \leq f(n) \right) \right \} \\
& = \inf\{\dege f: f \in \mathbb{L}_{>0}, \, \forall k \gg 0 \, ( \max( F_1(q_k(\alpha)), 2q_k(\alpha)^2)) \leq f(q_k(\alpha)) \} \\
& = \dege \max (F_1(q_k(\alpha)), 2q_k(\alpha)^2) \\
&  = \dege F_1.
\end{align*}
Moreover, by (\ref{cfineq}), one has
$$F_2(q_k(\alpha)) < F_1(q_k(\alpha))< 2  F_2(q_k(\alpha))$$
for all $k$.    Therefore, one has
\begin{align*}
\dege F_1 & = \dege F_2 \\
& = \dege G + (1,0,0,0,\ldots) \\
& = \dege H_1 + (2,0,0,0,\ldots).
\end{align*}
Finally,  since $0 < H_1-H_2 < 1$, one has $$\dege H_1 = \dege H_2.$$
The rest of the theorem then follows from Theorem  \ref{markk}.
\end{proof}

One application of Theorem \ref{pmbmutheorem} is in determining $\pmb{\mu}(\alpha)$ for real constants $\alpha$, like $\alpha = e$,  whose partial quotients in its regular continued fraction expansion exhibit some simple patterns.

\begin{corollary} \label{ecorr}
One has  $\pmb{\mu}(e) = (2,1,-1,0,0,0, \ldots)$.
\end{corollary}

\begin{proof}
It is well known that $a_{n+1}(e) = 2n$ for all $n$ congruent to $1$ modulo $3$ and $a_{n+1}(e) = 1$ for all $n$ not congruent to $1$ modulo $3$.  Moreover,  according to OEIS Sequence A007677, one has 
$$q_{3k-2}(e) = {}_{1}F_1(-k;-2k;-1) \frac{(2k)! }{ k!} \sim  e^{-1/2}  \frac{(2k)! }{ k!}\ (k \to \infty)$$
 This completely determines the values of function
$$H_2(q_n(e)) =  a_{n+1}(e),$$
as defined in Theorem \ref{pmbmutheorem}, that are not equal to $1$.   Since
$$k =  \frac{k \log k}{\log k} \sim \frac{\log H(k)}{\log \log H(k)}  \ (k \to \infty),$$
where $H(k) =  \tfrac{(2k)! }{ k!}$,
one has
$$H_2(q_{3k-2}(e))  = a_{3k-1}(e) = 2(3k-2)  \sim \frac{6\log q_{3k-2}(e)}{\log \log q_{3k-2}(e)} \ (k \to \infty),$$
whence 
$$\dege H_2 = (0,1,-1,0,0,0,\ldots)$$
and therefore 
$$\pmb{\mu}(e) = \dege H_2 + (2,0,0,0,\ldots) =   (2,1,-1,0,0,0, \ldots).$$
This completes the proof.
\end{proof}

\begin{problem}
Compute $\pmb{\mu}(\alpha)$ for various constants $\alpha$, like $\alpha =\tan 1$,  for which an explicit formula for the partial quotients of $\alpha$ is known.
\end{problem}

Another application of Theorem \ref{pmbmutheorem} is in describing the image of function $\pmb{\mu}$: as one might expect, the only additional constraints on $\pmb{\mu}(\alpha)$ besides $\pmb{\mu}(\alpha) \in \prod_{n = 1}^{\infty *} \overline{\RR}$  are that $\pmb{\mu}(\alpha) = (1,0,0,0,\ldots)$ or $\pmb{\mu}(\alpha) \geq (2,0,0,0,\ldots)$.

\begin{corollary}\label{inversemu}
Let $f: \ZZ_{>0} \longrightarrow \RR_{>0}$ be a positive arithmetic function with  $\lim_{n \to \infty} f(n) = \infty$, and, per Lemma \ref{prescribeda}, let  $\alpha$ be the unique irrational number with $a_0(\alpha) = 0$ and $$a_{n+1}(\alpha) = \lceil f(q_n(\alpha)) \rceil$$
 for all nonnegative integers $n$, so that
  $$ a_{n+1}(\alpha) \sim f(q_n(\alpha)) \ (n \to \infty).$$ Then one has
 $$\pmb{\mu}( \alpha) =  \dege f|_D+ (2,0,0,0,\ldots),$$
 where $D = \{q_n(\alpha): n \in \ZZ_{\geq 0}\}$.  Consequently, for any $\dd \in \prod_{n = 1}^{\infty *} \overline{\RR}$ with $\dd \geq (2,0,0,0,\ldots)$, there exists an irrational number $\alpha$ with $\pmb{\mu}(\alpha) = \dd$, and therefore
  $$\pmb{\mu}(\QQ)  = \{(1,0,0,0,\ldots)\}$$
  and
 $$\pmb{\mu}(\RR \backslash \QQ) = \left\{\dd \in  \prod_{n = 1}^{\infty *} \overline{\RR}: \dd \geq (2,0,0,0,\ldots)\right\}.$$ 
 \end{corollary}

\begin{proof}
The last claim of the corollary follows from the first by using Theorem \ref{degeequiv} to choose, given any $\dd > (2,0,0,0,\ldots)$ in $\prod_{n = 1}^{\infty *} \overline{\RR}$,  some positive function $f$ on $\RR_{>0}$  of exact logexponential degree $\dd' = \dd - (2,0,0,0,\ldots) > (0,0,0,\ldots)$ with $\dege f = \dd'$,  so that  $\dege f|_D = \dege f = \dd'$ and $\lim_{n \to \infty} f(n) = \infty$, and therefore also $\pmb{\mu}(\alpha) = \dd$.
\end{proof}

In fact,  the only additional constraints on the sets $\pmb{\mu}(\mu^{-1}(t) \cap m^{-1}(c))$ for $t \in  \{1\} \cap [2,\infty)$ and $c \in [0,\infty]$ come from Corollaries \ref{mmcor} and \ref{simplemu}.

\begin{corollary}\label{muimage}
Let $t \in \{1\} \cap [2,\infty)$,  let $c \in [0,\infty]$,  and let
$$ \pmb{\operatorname{M}}_{t,c} = \pmb{\mu}(\mu^{-1}(t) \cap m^{-1}(c)).$$ Exactly one of the following conditions holds.
\begin{enumerate}
\item  $\pmb{\operatorname{M}}_{t,c}$  is empty.
\item $t  = 1$ and $c$ is a positive integer, $t = 2$ and $c$ lies in the Lagrange spectrum, or $t \in (2,\infty)$ and $c \in (0,\infty)$; in these cases, one has
$$\pmb{\operatorname{M}}_{t,c} = \{(t,0,0,0,\ldots)\}.$$
\item $t \in (2,\infty)$ and $c = 0$,  in which case
$$\pmb{\operatorname{M}}_{t,0} = \left\{\dd \in  \prod_{n = 1}^{\infty *} \overline{\RR}: (t,-\infty,-1,0,0,0\ldots)\leq  \dd \leq (t,0,0,0,\ldots)\right\}.$$
\item $t \in [2,\infty)$ and $c = \infty$, in which case
$$\pmb{\operatorname{M}}_{t,\infty} = \left\{\dd \in  \prod_{n = 1}^{\infty *} \overline{\RR}: (t,0,0,0\ldots)\leq  \dd \leq (t,\infty,1,0,0,0,\ldots)\right\}.$$ 
\end{enumerate}
\end{corollary}

\begin{proof}
This is a direct consequence of Corollaries \ref{mcor}, \ref{mmcor},  \ref{simplemu},  and \ref{inversemu}.
\end{proof}

Note that the ``rational approximability'' of a number $\alpha$  can be measured inversely by  the ordered pair $(\pmb{\mu}(\alpha) , m(\alpha))$ in the lexicographic ordering, where we set  $m(\alpha ) = \infty$ if  $\mu(\alpha) = \infty$.  We  can describe this ordering in further detail by first ordering the reals $\alpha$ by $\mu(\alpha)$, with the set $\mu^{-1}(2)$ then  further ordered by $\pmb{\mu}(\alpha)$ from  $(2,0,0,0,\ldots)$ to $(2,\infty,1,0,0,0,\ldots)$,  with each  set $\mu^{-1}(t)$ for $2< t < \infty$ further  ordered by $\pmb{\mu}(\alpha)$ from  $(t,-\infty,-1,0,0,0,\ldots)$ to $(t,\infty,1,0,0,0,\ldots)$,  with the set $\mu^{-1}(\infty)$ of all Liouville numbers further ordered by $\pmb{\mu}(\alpha)$ from  $(\infty, 0,1,0,0,0,\ldots)$ to $(\infty,\infty,\infty,\ldots)$,   with the set $\pmb{\mu}^{-1}((1,0,0,0,\ldots)) = \QQ$ further ordered by $m(\alpha) = \ord_1 \alpha \in \ZZ_{>0}$, with the set $\pmb{\mu}^{-1}((2,0,0,0,\ldots))$ further  ordered by $m(\alpha) = M(\alpha)$, and, finally, with each set $\pmb{\mu}^{-1}((t,0,0,0,\ldots))$ for $t \in  (2,\infty)$ further ordered by $m(\alpha) \in [0,\infty]$.  

The following proposition generalizes the corresponding results for the function $\mu$ proved in Section 13.2.

\begin{proposition}
If $\alpha, \beta \in \RR$ are continued-fraction equivalent, then $\pmb{\mu}(\alpha) = \pmb{\mu}(\beta)$.  Moreover, if $\alpha \in \RR$ and $r \in \RR$, then $\pmb{\mu}(\alpha+r) = \pmb{\mu}(\alpha)$.
\end{proposition}

\begin{proof}
We may suppose without loss of generality that $\alpha$ is irrational.  By (\ref{qpq}),  if $\alpha \in (0,1)$, then one has $$q_n(\tfrac{1}{\alpha}) \sim \alpha q_{n+1}(\alpha) \ (n \to \infty)$$ and therefore 
 $$\log^{\circ k} q_n(\tfrac{1}{\alpha}) \sim \log^{\circ k} q_{n+1}(\alpha) \ (n \to \infty)$$ for all positive integers $k$.    The first statement then follows readily from the identity $\pmb{\mu}(\alpha) = \dege G$ of Theorem \ref{pmbmutheorem}.   Let $r = \frac{m}{n}$, where $m,n \in \ZZ$ are relatively prime and $n > 0$.  Suppose to obtain a contradiction that $\pmb{\mu}(\alpha+r) > \pmb{\mu}(\alpha)$.  Then, by Proposition \ref{oexppropstrong}, there exist $f,g \in \mathbb{L}_{>0}$ defined on $[1,\infty)$ with $\pmb{\mu}(\alpha+r) > \dege f > \dege g  > \pmb{\mu}(\alpha)$.    Since $\dege g > \pmb{\mu}(\alpha)$,  one has $$\left| \alpha - \frac{a}{b} \right| \geq \frac{1}{g(b)}$$ for all $a, b \in \ZZ$ with $b \gg 0$.   Moreover,  by Corollary \ref{fgie2}, one has $\dege g(nx) = \dege g < \dege f$, and therefore $g(nx) < f(x)$ for all $x \gg 0$.  It follows that
$$\left| \alpha + r - \frac{a}{b} \right|  = \left| \alpha - \frac{na-mb}{nb}  \right|  \geq \frac{1}{g(nb)} > \frac{1}{f(b)}$$
for all $a, b \in \ZZ$ with $b \gg 0$, whence $\alpha \succ_1 1/f|_{\ZZ_{>0}}$.  But, by Theorem \ref{pmbmu}(2), that contradicts $\pmb{\mu}(\alpha+r) > \dege f$.
\end{proof}

Next,  using Proposition \ref{iwannapp} and the following early predecessor of Theorem \ref{duffin} proved by Khinchin in 1924, we show that $\pmb{\mu}(\alpha) = (2,1,1,1,\ldots)$ for almost all real numbers $\alpha$.

\begin{theorem}[{Khinchin \cite{khin2} \cite[Theorem 2.1]{koukou}}]\label{khinchin}
Let  $f: \ZZ_{>0} \longrightarrow \RR_{\geq 0}$ be a nonnegative real-valued arithmetic function.     If $\sum_{n = 1}^\infty n f(n) < \infty$, then  $\alpha \succ_1 f$  for almost all $\alpha \in \RR$.  On the other hand, if  $\sum_{n = 1}^\infty n f(n) = \infty$ and $n^2 f(n)$ is nonincreasing, then $\alpha \not \succ_1 f$  for almost all $\alpha \in \RR$.
\end{theorem}

\begin{theorem}\label{mualmostall}
Let  $F: [1,\infty) \longrightarrow \RR_{>0}$ be a positive continuous monotonic Hardian function of finite degree,  and let $f = F|_{\ZZ_{>0}}$.
\begin{enumerate}
\item If $\dege f > (2,1,1,1,\ldots)$, then one has $\alpha \succ_1 1/f$ for almost all $\alpha \in \RR$.
\item If $\dege f < (2,1,1,1,\ldots)$ and $\frac{n^2}{f(n)}$ is nonincreasing, then one has $\alpha  \not \succ_1 1/f$  for almost all $\alpha \in \RR$.
\item One has $\pmb{\mu}(\alpha) = (2,1,1,1,\ldots)$ for almost all $\alpha \in \RR$.
\end{enumerate}
\end{theorem}

\begin{proof}
Let $g = 1/f$.   If $\dege f > (2,1,1,1,\ldots)$, then  $\dege n g(n) < (-1,-1,-1,\ldots)$,  and therefore, by Proposition \ref{iwannapp}, the sum $\sum_{n = 1}^\infty n g(n)$ converges.  On the other hand,  if $\dege f < (2,1,1,1,\ldots)$, then $\dege n g(n) > (-1,-1,-1,\ldots)$,  and therefore, by Proposition \ref{iwannapp}, the sum $\sum_{n = 1}^\infty n g(n)$ diverges.   Thus,  statements (1) and (2) follow from Theorem \ref{khinchin}. Now, for any positive integer $n$, we may choose any positive continuous monotonic logarithmico-exponential functions $F_n(x)$ and $G_n(x)$ on $[1,\infty)$ of logexponential degree $(2,1,1,1,\ldots, 1,2,0,0,0,\ldots)$ and $(2,1,1,1,\ldots, 1,0,0,0,0,\ldots)$, respectively,  where the initial $2$ is followed by a sequence of $n$ $1$s, and where  $\frac{x^2}{G_n(x)}$ is  nonincreasing on $[1,\infty)$.   Then, for almost all real numbers $\alpha$, one has
$\alpha \succ_1 1/F_n|_{\ZZ_{>0}}$ and therefore $\pmb{\mu}(\alpha) \leq \dege F_n$.  Likewise, for almost all real numbers $\alpha$, one has $\alpha  \not \succ_1 1/G_n|_{\ZZ_{>0}}$ and therefore $\pmb{\mu}(\alpha) \geq \dege G_n$.    Since a countable union of measure $0$ sets is measure $0$,   it follows that, for all real numbers $\alpha$ outside a set of measure $0$,  one has
$$(2,1,1,1,\ldots, 1,0,0,0,0,\ldots) \leq \pmb{\mu}(\alpha) \leq (2,1,1,1,\ldots, 1,2,0,0,0,\ldots)$$
for any finite sequence of $1$s.   Statement (3) follows.
\end{proof}

As a frame of reference for Theorem \ref{mualmostall}, note the following.

\begin{proposition}\label{ordinaryreals}
Let $\alpha \in \RR$.
\begin{enumerate}
\item One has  $\pmb{\mu}(\alpha)\leq  (2,1,1,1,\ldots)$ if and only if for any positive arithmetic function $f$ with $\dege f  < (-2,-1,-1,-1,\ldots)$ one has $\alpha \gg_1  f$, if and only if for every nonnegative integer $k$  there exists a $c > 0$ such that 
\begin{align*}
\alpha  \gg_1 \ & (n+c)^{-2} (\log (n+c))^{-1} (\log^{\circ 2} (n+c))^{-1} \cdots (\log^{\circ k} (n+c))^{-1}\\ 
& \quad  \cdot (\log^{\circ( k+1)} (n+c))^{-2}.
\end{align*} 
\item One has  $\pmb{\mu}(\alpha) >  (2,1,1,1,\ldots)$ if and only if  there exists a positive arithmetic function $f$ with $\dege f  < (-2,-1,-1,-1,\ldots)$ and  $\alpha \not \gg_1  f$, if and only if there exists a nonnegative integer $k$ and a  $c > 0$ such that 
\begin{align*}
\alpha \not \gg_1  \ & (n+c)^{-2} (\log (n+c))^{-1} (\log^{\circ 2} (n+c))^{-1} \cdots (\log^{\circ k} (n+c))^{-1}\\ 
& \quad  \cdot (\log^{\circ (k+1)} (n+c))^{-2}.
\end{align*} 
\item One has $\pmb{\mu}(\alpha) \geq  (2,1,1,1,\ldots)$ if and only if for any positive arithmetic function $f$ with $\dege f  > (-2,-1,-1,-1,\ldots)$, one has $\alpha  \not \gg_1  f$, if and only if for every nonnegative integer $k$  there exists a $c > 0$ such that  $$\alpha  \not \gg_1  (n+c)^{-2} (\log (n+c))^{-1} (\log^{\circ 2} (n+c))^{-1} \cdots (\log^{\circ k} (n+c))^{-1}.$$
\item One has $\pmb{\mu}(\alpha) < (2,1,1,1,\ldots)$ if and only if there exists a positive arithmetic function $f$ with $\dege f > (-2,-1,-1,-1,\ldots)$ and  $\alpha \gg_1  f$, if and only if there exists a nonnegative integer $k$ and a $c > 0$ such that  $$\alpha  \gg_1  (n+c)^{-2} (\log (n+c))^{-1} (\log^{\circ 2} (n+c))^{-1} \cdots (\log^{\circ k} (n+c))^{-1}.$$
\end{enumerate}
\end{proposition}

We now use the results in this section to  motivate and formulate several conjectures concerning the rational approximation of  algebraic and transcendental numbers.      Let us say that a real number $\alpha$ is {\bf regularly approximable}  if $\pmb{\mu}(\alpha)= (2,1,1,1,\ldots)$, that is, if $\alpha$ satisfies conditions (1) and (3) of Proposition \ref{ordinaryreals}.   Loosely speaking, this means that Diophantine approximation is    unable to discern the distinction between $\alpha$ and a ``random'' real number, at least using logarthmico-exponential functions as benchmarks against the Diophantine approximation functions $\frac{\Vert n \alpha \Vert }{n}$, $\mu_\alpha(x)$, and any of the other functions shown in this section to determine $\pmb{\mu}(\alpha)$.   We  first conjecture that $\pi$, $\log 2$, and $\gamma$ are regularly approximable.

\begin{conjecture}\label{diophconj0}
One has $\pmb{\mu}(\alpha)= (2,1,1,1,\ldots)$ for the real constants $\alpha = \pi, \log 2, \gamma$.
\end{conjecture}

Morally, this conjecture should extend to many other constants ``arising in nature''  whose partial quotients exhibit no discernible pattern, like the constants $e^\gamma$, $M$, and $H$, but very much unlike the constants $e$ and $\tan 1 = [1,1,1,3,1,5,1,7,\ldots]$, for example.    Note that Khinchin proved, famously, that 
\begin{align}
\lim _{n\rightarrow \infty }\left(a_{1}(\alpha)a_{2}(\alpha) \cdots a_{n}(\alpha)\right)^{1/n}=K_{0}\label{khinlaw}
\end{align} for almost all real numbers $\alpha$, where $$K_{0}=\prod _{k=1}^{\infty }{\left(1+{1 \over k(k+2)}\right)}^{\log _{2}k} = 2.685452001065\ldots$$ is {\bf Khinchin's constant} \cite{khin}. \index{Khinchin's constant}   One might conjecture that Conjecture \ref{diophconj0} extends to all real constants ``arising in nature'' that obey  Khinchin's law (\ref{khinlaw}).  Of course,  any precise sufficient criteria for ``arising in nature'' would have to exclude constants designed precisely to satisfy Khinchin's law.

We also make the following conjecture.

\begin{conjecture}\label{diophconj4}
One has $\pmb{\mu}(\RR\backslash \overline{\QQ}) = \pmb{\mu}(\RR\backslash\QQ)$.  More generally,  for every  $\dd \in  \prod_{n = 1}^{\infty *} \overline{\RR}$ with $\dd \geq (2,0,0,0,\ldots)$, there exist uncountably many  transcendental real numbers $\alpha$ with $\pmb{\mu}(\alpha) = \dd$.  
\end{conjecture}

The {\bf Jarn\'ik--Besicovitch theorem}, \index{Jarn\'ik--Besicovitch theorem} proved independently by V.\ Jarn\'ik \cite{jarn} and A.\ S.\ Besicovitch \cite{bes} in  1929 and 1934, respectively, states that the Hausdorff dimension of the subset $\mu^{-1}(t)$ of $\RR$ is equal to $2/t$ for all $t \geq 2$.     This naturally leads to the following problem.  

\begin{problem} 
If possible, using \cite[Theorem 1]{beres} and/or the notion of {\it Hausdorff measure with respect to a dimension function (or gauge function)} \cite{beres} \cite{rog}, generalize the Jarn\'ik--Besicovitch theorem to the function $\pmb{\mu}$.  More broadly,  describe how the inverse image $\pmb{\mu}^{-1}(\dd)$ of the map  $\pmb{\mu}$ varies with $\dd \in  \prod_{n = 1}^{\infty *} \overline{\RR}$.  For instance, describe the subsets $\mu^{-1}(t )\cap \pmb{\mu}_1^{-1}(u)$ of $\RR$ for $t \in [2,\infty]$ and $u \in [-\infty,\infty]$ (which have Lebesgue measure zero if and only if $(t,u) \neq (2,1)$).
\end{problem}

Our main problem concerning real algebraic numbers is the following.

\begin{outstandingproblem}\label{algprob}
Compute $\pmb{\mu}(\RR \cap  \overline{\QQ})$.  Also,  describe each of  the $\pmb{\mu}$-equivalence classes $\pmb{\mu}^{-1}(\dd) \cap \overline{\QQ}$ of the real algebraic numbers for $\dd \in  \pmb{\mu}(\RR \cap  \overline{\QQ})$.
\end{outstandingproblem}

 Let us say that a real number $\alpha$ is {\bf strongly approximable} if $\pmb{\mu}(\alpha)> (2,1,1,1,\ldots)$, that is, if $\alpha$ satisfies condition (2) of Proposition \ref{ordinaryreals}.  Loosely speaking, this means that $\alpha$ can be rationally approximated better than any regularly approximable real number.    Note, for example, that all very well approximable numbers are strongly approximable, and all strongly approximable numbers are well approximable.    Given  Theorem \ref{mualmostall}, Roth's theorem that all very well approximable numbers are transcendental,  and Lang's conjecture,  it is natural to hypothesize the generalization of Lang's conjecture stating that all strongly approximable numbers are transcendental.

\begin{conjecture}\label{diophconja}
One has $\pmb{\mu}(\alpha) \leq (2,1,1,1,\ldots)$  for all real algebraic  numbers $\alpha$,  or, equivalently,  for any real number $\alpha$,  if $\pmb{\mu}(\alpha)> (2,1,1,1,\ldots)$, then $\alpha$ is transcendental.  
\end{conjecture}

Note that Conjecture \ref{diophconja} is equivalent to $$\sup \pmb{\mu}(\RR \cap \overline{\QQ}) \leq (2,1,1,1,\ldots).$$   For obvious reasons,  then, we are interested in  $\sup \pmb{\mu}(\RR \cap \overline{\QQ})  \in \prod_{n = 1}^{\infty*}\overline{\RR}$,  i.e.,  the supremum of $\pmb{\mu}(\alpha)$ over all real algebraic numbers $\alpha$.

\begin{outstandingproblem}\label{muprob}
Compute $\sup \pmb{\mu}(\RR \cap \overline{\QQ})$.
\end{outstandingproblem}

\begin{proposition}
One has $\sup \pmb{\mu}(\RR \cap \overline{\QQ})  \geq (2,0,0,0,\ldots)$.   Moreover,  Roth's theorem is equivalent to  $\sup \pmb{\mu}(\RR \cap \overline{\QQ})  \leq (2,\infty,1,0,0,0,\ldots)$, and Lang's conjecture is equivalent to $\sup \pmb{\mu}(\RR \cap \overline{\QQ})   \leq (2,1,\infty,1,0,0,0,\ldots)$.
\end{proposition}

Note also that,  if it were the case that 
$$\alpha  \not \gg_1  n^{-2} (\log (n+1))^{-1}$$
for some  algebraic numbers $\alpha$, that is, if the ``$t>1$'' in Lang's conjecture were optimal, then the lower bound of $(2,0,0,0,\ldots)$ would increase to $(2,1,0,0,\ldots)$.   On the other hand, if it were the case that $\sup \pmb{\mu}(\RR \cap \overline{\QQ}) < (2,1,1,1,\ldots)$,  then all algebraic  numbers would be substantially less approximable than the  regularly approximable real numbers.   All of this provides some rationale for the following generalization of Conjecture \ref{diophconja}.

\begin{conjecture}\label{diophconj2}
One has $\sup \pmb{\mu}(\RR \cap \overline{\QQ}) = (2,1,1,1,\ldots)$.   Equivalently,  Conjecture \ref{diophconja} holds and, for every nonnegative integer $k$, there exists an algebraic number $\alpha$ and a $c > 0$ such that $$\alpha  \not \gg_1  (n+c)^{-2} (\log (n+c))^{-1} (\log^{\circ 2} (n+c))^{-1} \cdots (\log^{\circ k} (n+c))^{-1}.$$
\end{conjecture}

Assuming Conjecture \ref{diophconja}, one might also suspect that either all real algebraic numbers of degree greater than $2$ are regularly approximable,  or none of them are.   Thus, Conjecture \ref{diophconja} motivates the following even stronger conjecture.

\begin{conjecture}\label{diophconj}
One of the following statements holds.
\begin{enumerate}
\item  $\pmb{\mu}(\alpha)= (2,1,1,1,\ldots)$  for all  real algebraic numbers $\alpha$ of degree greater than $2$. 
\item  $\pmb{\mu}(\alpha) < (2,1,1,1,\ldots)$  for all  real algebraic numbers $\alpha$. 
\end{enumerate}
\end{conjecture}

Scenario (1) of the conjecture seems to be the more likely of the two scenarios, if the regular continued fraction expansion of any real algebraic number of degree greater than $2$ is expected to be asymptotically indistinguishable from that of a ``random'' real number $\alpha$ and thus, for example,  to satisfy Khinchin's law.   In that case,  the set $\pmb{\mu}(\RR \cap \overline{\QQ})$ contains only three elements.   On the other hand, if both Conjecture \ref{diophconj2} and  scenario  (2) of Conjecture \ref{diophconj} hold, then $\pmb{\mu}(\RR \cap \overline{\QQ})$ is infinite with $(2,1,1,1,\ldots)$ as a limit point.   The  simplest  wager,  then, and perhaps the safest, is that scenario (1)  of Conjecture \ref{diophconj} holds.   Thus,   Conjectures  \ref{diophconj2} and \ref{diophconj} together motivate the ``optimized'' generalization of  Lang's conjecture below.

\begin{conjecture}\label{diophconj3}
One has  $\pmb{\mu}(\alpha)= (2,1,1,1,\ldots)$  for all  real algebraic numbers $\alpha$ of degree greater than $2$.   Equivalently,  an algebraic irrational is regularly approximable (if and) only if its partial quotients are not eventually periodic.
\end{conjecture}

This conjecture provides the simplest possible answer to Problem \ref{algprob}, second only to  the rather  unpalatable conjecture $\pmb{\mu}(\RR \cap  \overline{\QQ}) = \{(1,0,0,0,\ldots), (2,0,0,0,\ldots)\}$.    The conjecture states, equivalently, that all algebraic numbers of degree greater than $2$ exhibit the same logexponential Diophantine approximation profile as almost all real numbers.  The inequality $\leq$ implies a kind of rigidity: the approximation error cannot grow any faster.  Thus, equality represents Diophantine rigidity without structural uniqueness. A strict inequality would, by contrast, reveal a deeper trace of arithmetic origin within the logexponential approximation profile.   In any event,  positive or negative answers to our conjectures would reveal  more precisely how  logexponential irrationality degree  distinguishes real algebraic numbers of degree greater than $2$ from other irrational numbers of irrationality measure $2$.

\section{L\'evy and relative L\'evy constants}

For any irrational $\alpha$, the {\bf L\'evy constant $\lambda(\alpha)$ of $\alpha$}\index{Levy @L\'evy constant $\lambda(\alpha)$} is defined by
$$\lambda(\alpha) = \lim_{n \to \infty} \frac{1}{n} \log q_n(\alpha),\index[symbols]{.v  b@$\lambda(\alpha)$}$$ provided that the limit exists.    Note that, if  $\lambda(\alpha)$ exists, then $\mu(\alpha) =2$.   By Theorem \ref{levythm}, which was proved by  L\'evy in 1937 \cite{levy}, one has
$\lambda(\alpha) = \frac{\pi^2}{12 \log 2}$
 for almost all real numbers $\alpha$.    Naturally, we define
$$\overline{\lambda}(\alpha) = \limsup_{n \to \infty} \frac{1}{n} \log q_n(\alpha) = \deg(q_n(\alpha); e^n) = \log \limsup_{n \to \infty} q_n(\alpha)^{1/n}\index[symbols]{.v  c@$\overline{\lambda}(\alpha)$}\index{upper L\'evy constant $\overline{\lambda}(\alpha)$}$$ 
and
$$\underline{\lambda}(\alpha) = \liminf_{n \to \infty} \frac{1}{n} \log q_n(\alpha) = \underline{\deg}(q_n(\alpha); e^n)= \log \liminf_{n \to \infty} q_n(\alpha)^{1/n},\index[symbols]{.v  d@$\underline{\lambda}(\alpha)$}\index{lower L\'evy constant $\underline{\lambda}(\alpha)$}$$
which we call the {\bf upper and lower L\'evy constants of $\alpha$}, respectively,  and which, unlike $\lambda(\alpha)$,
both exist  (in $[\log \Phi,\infty]$) for any irrational $\alpha$.   Note the following.

\begin{proposition}
Let $\alpha \in \RR$ be irrational.
\begin{enumerate}
\item $\overline{\lambda}(\alpha) \geq \underline{\lambda}(\alpha) \geq \lambda(\Phi) = \log\Phi$.
\item $\lambda(\alpha)$ exists if and only if $\log q_n(\alpha) \sim cn \ (n \to \infty)$ for some $c > 0$.
\item $\overline{\lambda}(\alpha) < \infty$ if and only if $\log q_n(\alpha) \ll n \ (n \to \infty)$, if and only if $\log q_n(\alpha) \asymp n \ (n \to \infty)$.
\item  $\overline{\lambda}(\alpha) \leq \log B(\alpha)$, where  $B(\alpha) = \limsup_{n \to \infty} \frac{q_{n+1}(\alpha)}{q_n(\alpha)}$.
\item $\underline{\lambda}(\alpha) \geq \log \underline{B}(\alpha)$, where  $\underline{B}(\alpha) = \liminf_{n \to \infty} \frac{q_{n+1}(\alpha)}{q_n(\alpha)}$.
\item If $\alpha$ is badly approximable, then $\overline{\lambda}(\alpha) < \infty$. 
\end{enumerate}
\end{proposition}

\begin{proof}
Statement (1) follows from the fact that $q_n(\alpha) \geq q_n(\Phi) =  F_{n+1}$ for all $n$.  Statements (2) and (3) are clear.   If $\frac{q_{n+1}(\alpha)}{q_n(\alpha)}$ is bounded above,  say, by $M > 1$, for all $n \geq N$,  where $N$ is a positive integer, then $q_n(\alpha) \leq q_N(\alpha) M^{n-N}$, whence $\log q_n(\alpha) \leq (n-N) \log M+ \log q_N(\alpha)$,  for all $n \geq N$, so that $\overline{\lambda}(\alpha) \leq \log M$.  Statement (4)  follows, and statement (5) is proved similarly, with the inequalities reversed.
Finally, statement (6) follows from statements (3) and (4), or, alternatively, from Corollary  \ref{badapp}.
\end{proof}

It is known that the image of $\lambda$ on the set of all irrational numbers whose L\'evy constant exists is $[\log\Phi,\infty)$,  and in fact the image of $\lambda$ on the set of all transcendental real numbers whose L\'evy constant exists is also $[\log\Phi,\infty)$ \cite{baxa}.  Moreover,  $\lambda(\alpha)$ exists for every quadratic irrational $\alpha$, and the image of $\lambda$ on the set of all quadratic irrationals is dense in $[\log\Phi,\infty)$ \cite{wu2}. 

 The proof of \cite[Lemma 11]{baxa}, which is a simple application of (\ref{cfineq00}) and (\ref{cfineq0}), yields the following.

\begin{proposition}[{cf.\  \cite[Lemma 11]{baxa}}]
Let $\alpha,\beta \in \RR$ be irrational.   One has
$$\log q_{n}(\alpha) = \sum_{k = 1}^n \log S^{\circ k}(\alpha) + O(1) \ (n \to \infty).$$  Consequently, one has
$$\overline{\lambda}(\alpha) = \limsup_{n \to \infty}  \frac{1}{n} \sum_{k = 1}^n \log  S^{\circ k}(\alpha) $$
and
$$\underline{\lambda}(\alpha) = \liminf_{n \to \infty}  \frac{1}{n} \sum_{k = 1}^n \log S^{\circ k}(\alpha).$$
Moreover, if $\alpha$ and $\beta$ are continued-fraction equivalent, then $\overline{\lambda}(\alpha) = \overline{\lambda}(\beta)$ and $\underline{\lambda}(\alpha) = \underline{\lambda}(\beta)$.
\end{proposition}

For all irrationals $\alpha$ and $\beta$,  let us define 
$$\overline{\lambda}(\alpha; \beta) = \limsup_{n \to \infty} \frac{\log q_n(\alpha)}{\log q_n(\beta)} = \deg(q_n(\alpha);q_n(\beta)) \in [0,\infty]\index[symbols]{.v  e@$\overline{\lambda}(\alpha; \beta)$}\index{upper relative L\'evy constant $\overline{\lambda}(\alpha;\beta)$}$$ 
and 
$$\underline{\lambda}(\alpha; \beta)= \liminf_{n \to \infty} \frac{\log q_n(\alpha)}{\log q_n(\beta)} = \underline{\deg}(q_n(\alpha);q_n(\beta)) \in [0,\infty],\index[symbols]{.v  f@$\underline{\lambda}(\alpha; \beta)$}\index{lower relative L\'evy constant $\overline{\lambda}(\alpha;\beta)$}$$
which we call the {\bf upper and lower relative L\'evy constants of $\alpha$ with respect to $\beta$}, respectively.  
Propositions \ref{levyprop1} and \ref{levyprop2} below  express both  the upper and lower L\'evy constants and the irrationality measure in terms of the upper and lower relative L\'evy constants.

\begin{proposition}\label{levyprop1}  One has the following.
\begin{enumerate}
\item For all irrationals $\alpha$ and $\beta$, one has
 $$\overline{\lambda}(\alpha;\alpha) = \underline{\lambda}(\alpha;\alpha) = 1,$$
 $$\underline{\lambda}(\alpha; \beta) = \frac{1}{\overline{\lambda}(\beta; \alpha)},$$
 and also
  $$\underline{\lambda}(\alpha)  \leq \underline{\lambda}(\beta) \underline{\lambda}(\alpha; \beta) \leq  \overline{\lambda}(\beta) \overline{\lambda}(\alpha; \beta) \leq \overline{\lambda}(\alpha)$$
whenever the respective products are defined in  $[0,\infty]$.
\item For any irrational $\alpha$, one has
$$\overline{\lambda}(\alpha) =\lambda (\Phi) \overline{\lambda}(\alpha;\Phi)$$
and
$$\underline{\lambda}(\alpha) = \lambda (\Phi)\underline{\lambda}(\alpha;\Phi).$$
\item For any irrational $\beta$,  the L\'evy constant $\lambda(\beta)$ exists and equals $l > 0$ if and only if 
 $$\overline{\lambda}(\alpha) = l  \overline{\lambda}(\alpha; \beta) $$ 
for all irrationals $\alpha$, if and only if 
$$\underline{\lambda}(\alpha) = l \underline{\lambda}(\alpha; \beta)$$
for all irrationals $\alpha$.   
\item For all irrationals $\alpha$ and $\beta$, one has 
$$\overline{\lambda}(\alpha';\beta')  = \overline{\lambda}(\alpha;\beta) $$
and
$$\underline{\lambda}(\alpha';\beta') = \underline{\lambda}(\alpha;\beta)$$
if $\alpha'=  \alpha +m$ and $\beta' =  \beta +n$ for any $m,n \in \ZZ$.  
\item For  all irrationals $\alpha$ and $\beta$, one has
$$\overline{\lambda}(\alpha; \beta) = \limsup_{n \to \infty}  \frac{ \sum_{k = 1}^n \log  S^{\circ k}(\alpha) }{\sum_{k = 1}^n \log  S^{\circ k}(\beta) }$$
and
$$\underline{\lambda}(\alpha; \beta) = \liminf_{n \to \infty}  \frac{ \sum_{k = 1}^n \log  S^{\circ k}(\alpha) }{\sum_{k = 1}^n \log  S^{\circ k}(\beta)}.$$
\end{enumerate}
\end{proposition}

Finally, by (\ref{qpq}) and Proposition \ref{levyprop1}(4) above,  one has the following.
 
 \begin{proposition}\label{levyprop2}
 For any irrational $\alpha$, one has
 \begin{align*}
 \mu(\alpha)   = 1+  \overline{\lambda}\left(\tfrac{1}{\{\alpha\}}; \{\alpha\}\right) = 1+ \overline{\lambda}(S(\alpha); \alpha) 
 \end{align*}
 and
  \begin{align*}
 \underline{\mu}(\alpha)   = 1+ \underline{\lambda}\left(\tfrac{1}{\{\alpha\}} ; \{\alpha\}\right)   = 1+ \underline{\lambda}(S(\alpha); \alpha).
 \end{align*}
 \end{proposition}
 
Thus, the upper and lower  relative  L\'evy constants are a natural common generalization of both L\'evy constants and irrationality measure.

\section{Rates of convergence}

In this section, we discuss various ways to define  the ``rate of convergence'' of a  convergent sequence of real numbers.   Our primary concern is with the sequence of convergents of the regular continued fraction expansion of an irrational number.

Let $\{\alpha_n\}_n$ be a sequence of real numbers that converge to $\alpha$ but are eventually not equal to $\alpha$.   The {\bf $Q$-order of convergence $Q\{\alpha_n\}_n$ of $\{\alpha_n\}_n$}\index{$Q$-order of convergence $Q\{\alpha_n\}_n$} is defined by
\begin{align*}
 Q\{\alpha_n\}_n & =  \sup\{ c \in [1,\infty): |\alpha-\alpha_{n+1}| = O\left( | \alpha-\alpha_n|^c \right) \ (n \to \infty) \} \\ &= \liminf_{n \to \infty } \frac{\log|\alpha-\alpha_{n+1}|}{\log|\alpha-\alpha_{n}|} 
 \\ & = \underline{\deg} \left( |\alpha-\alpha_{n+1}|;|\alpha-\alpha_{n}|\right)\\ &  \in [1,\infty]\index[symbols]{.v  k@$Q\{\alpha_n\}_n$}
 \end{align*}
 \cite{potra}.  If $Q\{\alpha_n\}_n < \infty$, then we define the {\bf upper and lower reciprocal rates of convergence  $\overline{\rho}\{\alpha_n\}_n$ and  $\underline{\rho}\{\alpha_n\}_n$ of $\{\alpha_n\}_n$} to be
$$ \overline{\rho}\{\alpha_n\}_n  = \limsup_{n \to \infty } \frac{|\alpha-\alpha_{n+1}|}{|\alpha-\alpha_{n}|^{Q\{\alpha_n\}_n}} \in [0,\infty]\index[symbols]{.va  a@$\overline{\rho}\{\alpha_n\}_n$}\index{upper reciprocal rate of convergence  $\overline{\rho}\{\alpha_n\}_n$}$$
and
$$ \underline{\rho}\{\alpha_n\}_n = \liminf_{n \to \infty } \frac{|\alpha-\alpha_{n+1}|}{|\alpha-\alpha_{n}|^{Q\{\alpha_n\}_n}} \in [0,\infty],\index[symbols]{.va  b@$\underline{\rho}\{\alpha_n\}_n$}\index{lower reciprocal rate of convergence  $\underline{\rho}\{\alpha_n\}_n$}$$
respectively.  

\begin{example} 
If for some $c \in [1,\infty)$ one has
$$|\alpha-\alpha_{n+1}| \asymp |\alpha-\alpha_{n}|^c \ (n \to \infty),$$
then $$Q\{\alpha_n\}_n = \lim_{n \to \infty } \frac{\log|\alpha-\alpha_{n+1}|}{\log|\alpha-\alpha_{n}|}  =  c,$$ and if also
$$|\alpha-\alpha_{n+1}| \sim r|\alpha-\alpha_{n}|^c \ (n \to \infty),$$
for some $r \in (0,\infty)$,
then 
$$ \underline{\rho}\{\alpha_n\}_n =  \overline{\rho}\{\alpha_n\}_n = r.$$
In numerical analysis, if the $\sim$ condition above holds, then $\{\alpha_n\}_n$ is said to have  {\bf order of convergence $c$} and {\bf rate of convergence $r$}.
\end{example}

\begin{example}
Let $c,r,s \in \RR$ with $0< s < r \leq 1 \leq c$.   Let $a_0 = 1$, and let
$$a_{n+1}=   \left.
 \begin{cases}
    r a_n^c & \text{if $n$ is even} \\
    s a_n^c  & \text{if $n$ is odd}
 \end{cases}
\right.$$ 
for all nonnegative integers $n$.  Then one has $\lim_{n \to \infty} a_n = 0$,  $Q\{\alpha_n\}_n = c$, $ \overline{\rho}\{\alpha_n\}_n  = r$, and $ \underline{\rho}\{\alpha_n\}_n  = s$.
\end{example}

For any irrational $\alpha$,  we let $$ Q(\alpha) = Q\left\{\frac{p_n(\alpha)}{q_n(\alpha)}\right\}_n \in [1,\infty],\index[symbols]{.va  c@$Q(\alpha)$}$$
and, if $Q(\alpha) < \infty$, then we let
$$ \overline{\rho}(\alpha) =  \overline{\rho}\left\{\frac{p_n(\alpha)}{q_n(\alpha)}\right\}_n\in [0,\infty]\index[symbols]{.va  d@$\overline{\rho}(\alpha)$}$$
and
$$\underline{\rho}(\alpha) = \underline{\rho}\left\{\frac{p_n(\alpha)}{q_n(\alpha)}\right\}_n\in [0,\infty].\index[symbols]{.va  e@$\underline{\rho}(\alpha)$}$$
Since 
 $$\left |\alpha -\frac{p_n(\alpha)}{q_n(\alpha)} \right|  \asymp \frac{1}{q_n(\alpha)q_{n+1}(\alpha)} = \left |\frac{p_{n+1}(\alpha)}{q_{n+1}(\alpha)}-\frac{p_n(\alpha)}{q_n(\alpha)} \right| \ (n \to \infty),$$ one has the following.
 
 \begin{proposition}\label{Qmu}
 Let $\alpha \in \RR$ be irrational.   One has
 \begin{align*}
 Q(\alpha) & =\underline{ \deg} \left (  \left |\alpha -\frac{p_{n+1}(\alpha)}{q_{n+1}(\alpha)} \right| ;  \left |\alpha -\frac{p_{n}(\alpha)}{q_{n}(\alpha)} \right|  \right) \\ & =\underline{ \deg} \left (  \left |\frac{p_{n+1}(\alpha)}{q_{n+1}(\alpha)}-\frac{p_n(\alpha)}{q_n(\alpha)} \right| ;  \left |\frac{p_{n}(\alpha)}{q_{n}(\alpha)}-\frac{p_{n-1}(\alpha)}{q_{n-1}(\alpha)} \right|  \right)  \\ & = \underline{\deg}( q_n(\alpha)q_{n+1}(\alpha); q_{n-1}(\alpha)q_{n}(\alpha)) \\
 & =  \liminf_{n \to \infty} \frac{1+ \frac{\log q_{n+1}(\alpha)}{\log q_{n}(\alpha)}}{1+\frac{1}{\frac{\log q_{n}(\alpha)}{\log q_{n-1}(\alpha)}}}
 \end{align*}
 and therefore 
 $$\underline{\deg} (q_{n+1}(\alpha); q_n(\alpha)) =  \underline{\mu}(\alpha)-1  \leq Q(\alpha) \leq \mu(\alpha)-1 = \deg (q_{n+1}(\alpha); q_n(\alpha)).$$
 \end{proposition}
 
Let us define $\overline{Q}\{\alpha_n\}_n$ and $\overline{Q}(\alpha)$ as in the definitions of  $Q\{\alpha_n\}_n$ and $Q(\alpha)$, but using limits superior instead of limits inferior.   Then one has
\begin{align*}
 \overline{Q}\{\alpha_n\}_n & =  \inf\{ c \in [1,\infty): |\alpha-\alpha_{n+1}| \gg | \alpha-\alpha_n|^c \ (n \to \infty) \} \\ &= \limsup_{n \to \infty } \frac{\log|\alpha-\alpha_{n+1}|}{\log|\alpha-\alpha_{n}|} 
 \\ & = \deg \left( |\alpha-\alpha_{n+1}|;|\alpha-\alpha_{n}|\right)\\ &  \in [1,\infty],
 \end{align*}
 and the obvious analogue of Proposition \ref{Qmu} holds for $\overline{Q}(\alpha)$.  In particular, one has
\begin{align}\label{qmuineq}
 \underline{\mu}(\alpha)-1  \leq Q(\alpha) \leq \overline{Q}(\alpha) \leq \mu(\alpha)-1
 \end{align}
 for all irrationals $\alpha$.    Consequently,  one has the following.
 
 \begin{corollary}
 Let $\alpha \in \RR$, and suppose that $\mu(\alpha) = \underline{\mu}(\alpha)$, i.e.,  suppose that
 $$\mu(\alpha) =1+ \lim_{n \to \infty}  \frac{\log q_{n+1}(\alpha)}{\log q_{n}(\alpha)}.$$
 Then one has  
 \begin{align*}
 Q(\alpha) = \overline{Q}(\alpha)  = \mu(\alpha)-1.
 \end{align*}
 \end{corollary}

 \begin{corollary}
For all $\alpha \in \RR$ with $\mu(\alpha) = 2$, one  has $Q(\alpha) =  \overline{Q}(\alpha) = 1$, and therefore
\begin{align*}
 \overline{\rho}(\alpha)    =   \limsup_{n \to \infty } \frac{ \left|\alpha-\frac{p_{n+1}(\alpha)} {q_{n+1}(\alpha)}\right|}{ \left|\alpha-\frac{p_{n}(\alpha)} {q_{n}(\alpha)}\right|} \in [0,1]
 \end{align*}
  and
\begin{align*} \underline{\rho}(\alpha)  = \liminf_{n \to \infty } \frac{ \left|\alpha-\frac{p_{n+1}(\alpha)} {q_{n+1}(\alpha)}\right|}{ \left|\alpha-\frac{p_{n}(\alpha)} {q_{n}(\alpha)}\right|} \in [0,1].
\end{align*}
 \end{corollary}
 
 \begin{example}
One has $\mu(e) = 2$, $M(e) =  \overline{\lambda}(e) = \underline{\lambda}(e) = \infty$,  $Q(e)  = \overline{Q}(e) = 1$,  $ \overline{\rho}(e) = 1$, and $ \underline{\rho}(e) = 0$.
 \end{example}

The following example shows that all three inequalities in (\ref{qmuineq}) can be strict,  even given prescribed values $a$ of $\mu(\alpha)$ and $b$ of $\underline{\mu}(\alpha)$ with $2 < b < a < \infty$.
 
 \begin{example}
For any real numbers $a > b> 1$,  there is a unique real number $\alpha$ with $a_0(\alpha) = 0$ and 
$$a_{n+1}(\alpha) =   \left.
 \begin{cases}
    \lceil q_n(\alpha)^{a-1} \rceil & \text{if $n$ is even} \\
     \lceil q_n(\alpha)^{b-1} \rceil  & \text{if $n$ is odd}
 \end{cases}
\right.$$ 
for all nonnegative integers $n$,  and one has
$$ \underline{\mu}(\alpha)-1   = b < Q(\alpha)  = \frac{1+b}{1+\frac{1}{a}} <  \overline{Q}(\alpha) =  \frac{1+a}{1+\frac{1}{b}}  <  \mu(\alpha)-1 = a.$$
 \end{example}

 \begin{problem}
It follows from Proposition \ref{RProp} that $ \underline{\rho}(\alpha) = 0$ for almost all irrationals $\alpha$.  Does one also have $\overline{\rho}(\alpha) = 1$ for almost all irrationals $\alpha$?
 \end{problem}
 
For any irrational number $\alpha$,  the quantity
$$R(\alpha)    =   \limsup_{n \to \infty } \frac{ \left|\alpha-\frac{p_{n}(\alpha)} {q_{n}(\alpha)}\right|}{ \left|\alpha-\frac{p_{n+1}(\alpha)} {q_{n+1}(\alpha)}\right|} \in [1,\infty] \index[symbols]{.va  f@${R}(\alpha)$}$$
is yet another measure of the rational approximability of $\alpha$.   Equivalently,  $R(\alpha)$ is the infimum of all $r \geq 1$ such that
\begin{align*}
 \left|\alpha-\frac{p_{n}(\alpha)} {q_{n}(\alpha)}\right| \leq r \left|\alpha-\frac{p_{n+1}(\alpha)} {q_{n+1}(\alpha)}\right|, \quad \forall n \gg 0.
\end{align*}
If $Q(\alpha)$ is minimal, that is, if $Q(\alpha) = 1$ (or if $\mu(\alpha) = 2$), then one has
$$R(\alpha) = \frac{1}{ \underline{\rho}(\alpha)},$$
while, if $Q(\alpha) > 1$, then $R(\alpha ) = \infty$.  Thus,  the quantity $ \frac{1}{ \underline{\rho}(\alpha)}$ is a more refined invariant than $R(\alpha)$.  However,  unlike $Q(\alpha)$ and $R(\alpha)$,  the quantity $ \frac{1}{ \underline{\rho}(\alpha)}$  is not invariant under continued-fraction equivalence,  since, for example,  $ \underline{\rho}(\tfrac{1}{\alpha})$ need not equal $ \underline{\rho}(\alpha)$ if $Q(\alpha) > 1$.  Indeed,  for any irrational $\alpha$, the convergents of $\frac{1}{\alpha}$ are the reciprocals of the convergents of $\alpha$ with an index shift forward or back by $1$, and, moreover, if  $\alpha_n \neq 0$ for all $n$ and $\alpha = \lim_{n \to \infty} \alpha_n \neq 0$, then one has $$Q\left\{\frac{1}{\alpha_n}\right\}_n = Q\{\alpha_n\}_n,$$
while
$$\overline{\rho}\left\{\frac{1}{\alpha_n}\right\}_n = |\alpha|^{2Q\{\alpha_n\}_n-2}\overline{\rho}\{\alpha_n\}_n$$
and
$$\underline{\rho}\left\{\frac{1}{\alpha_n}\right\}_n =|\alpha|^{2Q\{\alpha_n\}_n-2} \underline{\rho}\{\alpha_n\}_n$$
if $ Q\{\alpha_n\}_n < \infty$.
Thus, one has the following.

 \begin{proposition}
The functions $\overline{\rho}$ and $\underline{\rho}$ are even and periodic of period $1$.
 Let $\alpha, \beta \in \RR$ be irrational.  If  $\alpha$ and $\beta$  are continued-fraction equivalent, then $Q(\alpha ) = Q(\beta)$ and $R(\alpha) = R(\beta)$.    Moreover,  if $Q(\alpha) < \infty$, then one has
 $$\overline{\rho}\left(\frac{1}{\alpha} \right) = |\alpha|^{2Q(\alpha)-2}\overline{\rho}(\alpha),$$
$$\underline{\rho}\left(\frac{1}{\alpha} \right)  =|\alpha|^{2Q(\alpha)-2} \underline{\rho}(\alpha).$$
 $$\overline{\rho}(S(\alpha)) = |S(\alpha)|^{2-2Q(\alpha)}\overline{\rho}(\alpha),$$
 and
 $$\underline{\rho}(S(\alpha))  =|S(\alpha)|^{2-2Q(\alpha)} \underline{\rho}(\alpha).$$
\end{proposition}

Note that the quantity
$$\underline{R}(\alpha) = \liminf_{n \to \infty } \frac{ \left|\alpha-\frac{p_{n}(\alpha)} {q_{n}(\alpha)}\right|}{ \left|\alpha-\frac{p_{n+1}(\alpha)} {q_{n+1}(\alpha)}\right|}  \in [1,\infty]\index[symbols]{.va  fa@$\underline{R}(\alpha)$}$$
which, equivalently,  is the supremum of all $r \geq 1$ such that
\begin{align*}
 \left|\alpha-\frac{p_{n}(\alpha)} {q_{n}(\alpha)}\right| \geq r \left|\alpha-\frac{p_{n+1}(\alpha)} {q_{n+1}(\alpha)}\right|, \quad \forall n \gg 0,
\end{align*}
 is also invariant under continued-fraction equivalence.  
 
 By (\ref{cfineq00})--(\ref{cfineq0}), one has the following.
 
\begin{proposition}\label{RProp}
 Let $\alpha \in \RR$ be  irrational with regular continued fraction expansion $[a_0,a_1,a_2,\ldots]$.   One has
 \begin{align*}
  \frac{ \left|\alpha-\frac{p_{n}(\alpha)} {q_{n}(\alpha)}\right|}{ \left|\alpha-\frac{p_{n+1}(\alpha)} {q_{n+1}(\alpha)}\right|}   & = \left(\frac{q_{n+1}(\alpha)}{q_n(\alpha)} \right)^2 \left(\frac{S^{\circ (n+2)}(\alpha)+ \frac{q_{n}(\alpha)}{q_{n+1}(\alpha)}}{S^{\circ (n+1)}(\alpha)+ \frac{q_{n-1}(\alpha)}{q_n(\alpha)}}\right)  \\ & =  \frac{q_{n+1}(\alpha)}{q_n(\alpha)} S^{\circ (n+2)}(\alpha) \\& = [a_{n+1},a_{n},\ldots,a_1][a_{n+2},a_{n+3},\ldots].
 \end{align*}
In particular, one has
$$R(\alpha) = \limsup_{n \to \infty} \, [a_n,a_{n-1},\ldots,a_1][a_{n+1},a_{n+2},\ldots]$$
and
$$\underline{R}(\alpha) = \liminf_{n \to \infty}  \, [a_n,a_{n-1},\ldots,a_1][a_{n+1},a_{n+2},\ldots].$$
\end{proposition}

Note also the following.

\begin{lemma}\label{intervalfrac}
Let $n,  m \in \ZZ_{\geq 0}$, let $a_0 \in \ZZ$,  and let $a,a_1,a_2,\ldots,a_n \in \ZZ_{>0}$. 
\begin{enumerate}
\item  The smallest rational number  $[b_0,b_1,\ldots,b_m]$ with $b_k \in [1,a]$  for all $k \leq 2n \leq m$ is $[1,a,1,a,\ldots,1,a,1]$, where the sequence $1,a$ occurs $n$ times.    
\item  The largest rational number $[b_0,b_1,\ldots,b_m]$ with
$b_k(\alpha) \in [1,a]$  for all $k \leq 2n+1 \leq m$ is $[a,1,a,1,\ldots,a,1]$,  where the sequence $a,1$ occurs $n+1$ times. 
\item The smallest irrational  number $\alpha$ with $a_k(\alpha) = a_k$  for all $k < n$ and $a_k(\alpha) \in [1,a]$  for all $k \geq n$ is 
$$\alpha = \left. \begin{cases} [a_0,a_1,\ldots,a_{n-1},\overline{1,a}] & \quad  \text{if }  n \text{ is even} \\
 [a_0,a_1,\ldots,a_{n-1},\overline{a,1}]  & \quad \text{if } n \text{ is odd}.
\end{cases}\right.$$
\item The largest irrational  number $\alpha$ with $a_k(\alpha) = a_k$  for all $k < n$ and $a_k(\alpha) \in [1,a]$  for all $k \geq n$ is 
$$\alpha = \left. \begin{cases} [a_0,a_1,\ldots,a_{n-1},\overline{a,1}] & \quad  \text{if }  n \text{ is even} \\
 [a_0,a_1,\ldots,a_{n-1},\overline{1,a}]  & \quad \text{if } n \text{ is odd}.
\end{cases}\right.$$
\end{enumerate}
 \end{lemma}

We use the proposition and lemma above to prove the following result.

 \begin{proposition}\label{imageofR}
 Let $\alpha \in \RR$ be irrational, and let $a$ and $b$ be positive integers.  One has the following.
 \begin{enumerate}
  \item $R([\overline{a,b}]) = [\overline{a,b}] [\overline{b,a}]$.
 \item $R(\Phi) = \Phi^2$.
 \item  $R([\overline{2}]) = R(\sqrt{2}) =  3+2\sqrt{2} = 5.828427\ldots$, where $[\overline{2}] = 1+\sqrt{2}$.
\item  $R( [\overline{2,1}]) = R( [\overline{1,2}] )  = [\overline{2,1}] [\overline{1,2}] =  2+\sqrt{3} = 3.732050\ldots$, where $[\overline{2,1}]= 1+\sqrt{3}$ and $[\overline{1,2}] = \frac{1+\sqrt{3}}{2}$.
 \item  If $A(\alpha) \geq a$, then $$R(\alpha) \geq [\overline{a,1}]+ \frac{1}{a} = \frac{a+ \sqrt{a(a+4)}}{2}+ \frac{1}{a} > a + \frac{2}{a}.$$
    \item The supremum of $R(\alpha)$ over all irrationals $\alpha$ with $A(\alpha) = a$ is equal to the limit $$\sup\{[\overline{a}]^2, [\overline{a,1,a}]^2, [\overline{a,1,a,1,a}]^2,\ldots\}  = [\overline{a,1}]^2$$
 from below.   In particular, $[\overline{2,1}]^2 =  2(2+\sqrt{3})$ is a limit point of $R(\RR\backslash\QQ)$ from below.
  \item If $A(\alpha) \geq 2$, then  $R(\alpha) \geq \frac{3}{2} +\sqrt{3} = 3.232050\ldots$.
 \item If $A(\alpha) \geq 3$, then  $R(\alpha) \geq \frac{11}{6} + \frac{\sqrt{21}}{2} = 4.214621\ldots$.
 \item $\min R(\RR\backslash \QQ) = R(\Phi) = \Phi^2$.
 \item $\min R((\RR\backslash \QQ)\backslash \{ R(\Phi)\} )= R( [\overline{2,1}]) = 2+\sqrt{3}$.
 \item  If $\alpha$ is not continued-fraction equivalent to $\Phi$ or $[\overline{2,1}]$, then  $R(\alpha) \geq \frac{51 + 19 \sqrt{3}}{22} = 3.814043\ldots > R([\overline{2,1}])$.
  \item $ R( [\overline{2,1}])$ is an isolated point of $R(\RR\backslash \QQ)$.
 \end{enumerate}
 \end{proposition}

\begin{proof}
Statements (1)--(4) are clear.  Let  $r_n = [1,a,1,a,\ldots,1,a,1]$, where the sequence $1,a$ occurs $n$ times. If $A(\alpha) = a$, then $a_n(\alpha) \in [1,a]$ for all $n \gg 0$, and therefore,  by Proposition \ref{RProp} and statements (1) and (3) of Lemma \ref{intervalfrac},  one has
$$R(\alpha) \geq r_n [a,a,1,a,1,a,\ldots]$$
for all $n$, and therefore
$$R(\alpha) \geq [1,a,1,a,\ldots][a,a,1,a,1,a,\ldots] =  [\overline{a,1}]+ \frac{1}{a}.$$
Statements (5) and (7)--(9) follow.   To prove (6),  let $\alpha_n = [\overline{a,1,a,1,\ldots,a}]$, where the sequence $a,1$ occurs $n$ times.
Then $R(\alpha_n) = [\overline{a,1,a,1,\ldots,a}]^2$, and $$\lim_{n \to \infty} R(\alpha_n) = [\overline{a,1}]^2.$$
Moreover, if $A(\alpha) = a$, then it is clear from  Lemma \ref{intervalfrac} that  $R(\alpha) \leq [\overline{a,1}]^2$.  Statement (6) follows.
To prove (10), suppose that $A(\alpha) =2$ and $\alpha$ is not continued-fraction equivalent to $[\overline{2,1}]$.  Then the regular continued fraction expansion of $\alpha$ contains infinitely many sequences of the form $2,2$ or infinitely many sequences of the form $2,1,1$.  In the former case, one has
$$R(\alpha) \geq [2,2,1,2,1,2,1,2,\ldots][2,2,1,2,1,2,1,2,\ldots]  = 3+ \frac{3\sqrt{3}}{2} = 5.598076\ldots > 2+\sqrt{3},$$
by Proposition \ref{RProp} and statements (2) and (3) of Lemma \ref{intervalfrac}.
Likewise, in the latter case, one has
 $$R(\alpha) \geq [2,2,1,2,1,2,1,2,\ldots][1,1,1,2,1,2,1,2,\ldots] = 2+\sqrt{3} = R([\overline{2,1}]).$$
 Statement (10) follows.   Moreover,  the inequality above is strict if the regular continued fraction of $\alpha$  contains infinitely many  sequences of the form $2,2$ or $2,1,1,2$ or $2,1,1,1,1$: in the latter two cases,  respectively, one has
  $$R(\alpha) \geq [2,2,1,2,1,2,1,2,\ldots][1,1,2,2,1,2,1,2,\ldots] = \frac{27 + 10 \sqrt{3}}{11} = 4.029137\ldots$$
and
    $$R(\alpha) \geq [2,2,1,2,1,2,1,2,\ldots][1,1,1,1,1,2,1,2,1,2\ldots] = \frac{51 + 19 \sqrt{3}}{22} = 3.814043\ldots.$$
Suppose, then, that the regular continued fraction of $\alpha$  contains only finitely many such sequences, and also that $\alpha$ is not continued-fraction equivalent to $[\overline{2,1}]$ or $[\overline{2,1,1,1}]$.   It then follows that the tail of  the regular continued fraction of $\alpha$  is the concatenation of blocks of the form $2,1$ and $2,1,1,1$, and with infinitely many of both, and therefore  with infinitely many blocks of the form $2,1,2,1,1,1$.  It follows that
 $$R(\alpha) \geq [2,1,2,2,1,2,1,2,1,2\ldots][1,1,1,2,1,2,1,2,\ldots]  =  \frac{87 + 31 \sqrt{3}}{33} = 4.263441\ldots.$$
Moreover,  one has
$$R([\overline{2,1,1,1}]) = [\overline{2,1,1,1}][\overline{1,1,1,2}] = \frac{11+4\sqrt{6}}{5} = 4.159591\ldots.$$
Statements (11) and (12) follow.    
\end{proof}

Next, for  any irrational $\alpha$,  we let
$$A(\alpha) = \limsup_{n \to \infty} a_n(\alpha),\index[symbols]{.va  g@$A(\alpha)$}$$
$$\underline{A}(\alpha) = \liminf_{n \to \infty} a_n(\alpha),\index[symbols]{.va  h@$\underline{A}(\alpha)$}$$
$$ B(\alpha)= \limsup_{n \to \infty} \frac{q_{n+1}(\alpha)}{q_{n}(\alpha)},\index[symbols]{.va  i@$B(\alpha)$}$$
$$ \underline{B}(\alpha)= \liminf_{n \to \infty} \frac{q_{n+1}(\alpha)}{q_{n}(\alpha)},\index[symbols]{.va  j@$\underline{B}(\alpha)$}$$
and $$C(\alpha) = \limsup_{n \to \infty} \frac{q_{n+2}(\alpha)}{q_{n}(\alpha)} = 1+\limsup_{n \to \infty} a_{n+2}(\alpha)\frac{q_{n+1}(\alpha)}{q_{n}(\alpha)}.\index[symbols]{.va  k@$C(\alpha)$}$$
Using (\ref{cfineq00})--(\ref{cfineq}) and the definitions above, one can readily verify the following.

\begin{proposition}\label{ABCR}
Let $\alpha \in \RR$ be irrational.  One has
$$1 \leq A(\alpha) \leq M(\alpha) \leq A(\alpha)+2,$$
$$\max\left(M(\alpha)-1,\Phi \right) \leq B(\alpha) \leq M(\alpha),$$
$$1\leq \underline{A}(\alpha)  \leq \underline{B}(\alpha) \leq B(\alpha)\leq A(\alpha)+1,$$
$$1 \leq \lceil B(\alpha)\rceil-1 \leq A(\alpha)\leq \lfloor B(\alpha)\rfloor,$$
 $$ \Phi^2 \leq \max\left(B(\alpha)\underline{B}(\alpha), A(\alpha)\underline{B}(\alpha)+1,\underline{A}(\alpha)B(\alpha)+1\right) \leq C(\alpha) \leq \min\left( B(\alpha)^2, A(\alpha)B(\alpha)+1 \right),$$ 
 and
$$\tfrac{1}{2}\Phi^2\leq \tfrac{1}{2}C(\alpha)\leq \frac{ C(\alpha)}{1+\frac{1}{\underline{B}(\alpha)}} \leq  R(\alpha) \leq  C(\alpha)+B(\alpha)\leq 2C(\alpha)-1.$$
Consequently, one has $M(\alpha)<\infty$ if and only if  $A(\alpha) < \infty$, if and only if $B(\alpha) < \infty$, if and only if $C(\alpha) < \infty$,  if and only if $R(\alpha)< \infty$, if and only if $\alpha$ is badly approximable, and any of these equivalent conditions implies that $\mu(\alpha) = 2$ and  $Q(\alpha) = 1$.  Moreover, if $\beta \in \RR$ is continued-fraction equivalent to $\alpha$, then $A(\alpha) = A(\beta)$, $\underline{A}(\alpha) = \underline{A}(\beta)$, $B(\alpha) = B(\beta)$, $\underline{B}(\alpha) = \underline{B}(\beta)$,  and $C(\alpha) = C(\beta)$.
 \end{proposition}

Let $\alpha \in \RR$ be irrational.  For each positive integer $k$, let
$$B_k(\alpha) = \limsup_{n \to \infty} \frac{q_{n+k}(\alpha)}{q_{n}(\alpha)} \in [1,\infty]\index[symbols]{.va  l@$B_k(\alpha)$}$$
and
$$\underline{B}_k(\alpha) = \liminf_{n \to \infty} \frac{q_{n+k}(\alpha)}{q_{n}(\alpha)} \in [1,\infty],\index[symbols]{.va  m@$\underline{B}_k(\alpha)$}$$
and let
$$B_0(\alpha) =  A(\alpha) \in  \ZZ_{>0}$$
and
$$\underline{B}_0(\alpha) =  \underline{A}(\alpha) \in \ZZ_{>0}.$$
Thus, one has
$$B(\alpha)  = B_1(\alpha),$$
$$\underline{B}(\alpha)  = \underline{B}_1(\alpha),$$
and
$$C(\alpha)  = B_2(\alpha).$$
Note that the functions $B_k$ and $\underline{B}_k$ are even and periodic of period $1$ for all nonnegative integers $k$.
The following proposition is clear.

\begin{proposition}
Let $\alpha \in \RR$ be irrational.  For all nonnegative integers $k \leq l$, one has
$$\underline{B}_k(\alpha)  \leq  B_k(\alpha)  \leq B_l(\alpha)$$
and
$$\underline{B}_k(\alpha)  \leq \underline{B}_l(\alpha).$$
Moreover, one has
$$B_k(\alpha) \leq B_1(\alpha)^k$$
for all nonnegative integers $k$, where equality holds if $\lim_{n \to \infty} \frac{q_{n+1}(\alpha)}{q_{n}(\alpha)}$ exists, i.e.,  if
$B_1(\alpha) = \underline{B}_1(\alpha)$.
More generally, one has
$$B_k(\alpha) \leq B_{k_1}(\alpha) B_{k_2}(\alpha)\cdots  B_{k_n}(\alpha)$$
for any partition $k = k_1 + k_2 +\ldots + k_n$  of $k$, with equality if $B_{k_i}(\alpha) = \underline{B}_{k_i}(\alpha)$ for all $i<n$.    Consequently,   $\alpha$ is badly approximable if and only if $B_k(\alpha) < \infty$ for some $k$, if and only if $B_k(\alpha) < \infty$ for all $k$.
\end{proposition}

\begin{example}
One has $B_k(e) = R(e)= \infty$ for all $k \geq 0$, while $\underline{B}_0(e) = \underline{B}_1(e)  = \underline{R}(e)= 1$,  $\underline{B}_2(e) = 2$, and  $\underline{B}_k(e) = \infty$ for all $k \geq 3$.
\end{example}

For the  remainder of this section, we attempt to address the following problem.

\begin{problem}\label{imageprob}
  Describe the image of  the functions $B$, $\underline{B}$, $C$, and, more generally, of $B_k$ and $\underline{B}_k$ for all nonnegative integers $k$.    Do the same for the functions $R$ and $\underline{R}$.   
\end{problem}

For all nonnegative integers $k$,  let $$\iota_k = \inf  B_k(\RR\backslash \QQ).$$ 
One has  $$\iota_0 = B_0(\Phi) = A(\Phi) =1,$$
$$\iota_1 = B_1(\Phi) = B(\Phi) = \Phi,$$
and
$$\iota_2 = B_2(\Phi) = C(\Phi) = \Phi^2.$$
In fact, one has $\iota_k = B_k(\Phi) =  \Phi^k$ for all $k$, which we prove below.   Our proof rests on the following result, which provides a recurrence relation for the quantities  $\frac{q_{n+k}(\alpha)}{q_{n}(\alpha)}$, for which (\ref{cfineqa}) is a base case.

\begin{proposition}\label{anrn}
Let $\alpha \in \RR$  be irrational with regular continued fraction expansion $\alpha = [a_0, a_1,a_2,\ldots]$.  For all nonegative integers $n$ and $k$,  let
$$r_{n,k} = \frac{q_{n+k}(\alpha)}{q_{n}(\alpha)}.$$
For all nonnegative integers $n$ and $k > 1$ and $l < k-1$,  one has
$$r_{n,k} = p_{l}( [a_{n+k}, a_{n+k-1},\ldots,a_1])r_{n,k-l-1} +  p_{l-1}( [a_{n+k}, a_{n+k-1},\ldots,a_1])r_{n,k-l-2},$$
where $$r_{n,0} = 1,$$ 
$$r_{n,1} = \frac{q_{n+1}(\alpha)}{q_{n}(\alpha)} = [a_{n+1}, a_{n},\ldots,a_1],$$
and
$$r_{n,2} = a_{n+2} r_{n,1}+1.$$
Consequently, if $\beta \in \RR$ is continued-fraction equivalent to $\alpha$, then $B_k(\alpha) = B_k(\beta)$ and $\underline{B}_k(\alpha) = \underline{ B}_k(\beta)$ for all nonnegative integers $k$.  
\end{proposition}

\begin{proof}
The proof, though tedious, is by induction.
\end{proof}

\begin{proposition}\label{ABk}
One has
$$\min B_k(\RR\backslash \QQ) = B_k(\Phi) = \Phi^k$$
for all nonnegative integers $k$.  Equivalently, one has $$B_k(\alpha) \geq B_k(\Phi) = \Phi^k$$ for all irrationals $\alpha$ and all nonnegative integers $k$.   Let $k$ be a positive integer.   One has
$$B_k(\alpha) \geq F_{k} A(\alpha) +  F_{k-1},$$ 
where $F_k$ denotes the $k$th Fibonacci number.  Consequently, if 
$\alpha$ is not continued-fraction equivalent to $\Phi$, then 
$$B_k(\alpha) \geq F_{k+2} = 2F_{k}  +  F_{k-1}  > \Phi F_{k}  +  F_{k-1} = \Phi^k.$$ 
In particular, $\Phi^k$ is an isolated point of $B_k(\RR \backslash \QQ)$.  
\end{proposition}

\begin{proof}
Let $\alpha = [a_0, a_1,a_2,\ldots]$ be irrational.
For all positive integers $n$, let $$\beta_n = \frac{q_n(\alpha)}{q_{n-1}(\alpha)} = [a_n, a_{n-1},\ldots,a_1].$$  
By Proposition \ref{anrn},  for all $n$ and all $k  > 1$, one has
\begin{align*}
\frac{q_{n+k}(\alpha)}{q_n(\alpha)} & = p_{k-2}(\beta_{n+k}) \beta_{n+1} +  p_{k-3}(\beta_{n+k}) \\
& \geq F_{k} \beta_{n+1} + F_{k-1} \\
& >F_{k} a_{n+1} + F_{k-1} 
\end{align*}
and therefore
$$B_k(\alpha) \geq  F_{k}A(\alpha)  + F_{k-1}.$$
Moreover, since $B(\alpha) \geq A(\alpha)$, the inequality above also holds for $k = 1$.
This completes the proof.
\end{proof}

By the results in this chapter, one has the following result, which provides several ``universal properties'' of the golden ratio $\Phi$.

\begin{theorem}
On $\RR\backslash \QQ$, the functions $f = \mu$, $\underline{\mu}$, $M$,  $\overline{\lambda}$,  $\underline{\lambda}$,  $Q$,  $R$, and $B_k$ for all $k$ are all minimized at the golden ratio $\Phi$.  Moreover, $f(\Phi)$ is an isolated point of $\im f$ (and $f$ is minimized precisely at the real numbers that are continued-fraction equivalent to $\Phi$) precisely for $f = M$,  $R$, $B_k$.
\end{theorem}

Note that $q_k$, for any fixed $k$, assumes its minimal value $F_{k+1}$ precisely at $\frac{F_{k+2}}{F_{k+1}} +n = [n+1,1,1,\ldots,1,2]$ for $n \in \ZZ$ and at those real numbers $\alpha$ with $a_j(\alpha) = 1$ for all positive integers $j \leq k$,  and thus $q_k$ is minimized for all $k$ simultaneously  precisely at the numbers $\Phi+n$ for $n \in \ZZ$.

Problem \ref{imageprob} specializes as follows.

\begin{problem}
Let $k$ be a positive integer.  Compute $I_{k,1} = \inf\left(B_k(\RR\backslash \QQ)\backslash\{\Phi^{k}\} \right)$.  Compute $I_{k,2} = \inf\left(B_k(\RR\backslash \QQ)\backslash\{\Phi^{k},I_{k,2}\} \right)$.  And so on.    Address the analogous problems for the functions $\underline{B}_k$, $R$, and $\underline{R}$.
\end{problem}

Proposition \ref{imageofR} is a partial solution to  the problem above for the function $R$.
The next two propositions are analogous results for the functions $B= B_1$ and $C = B_2$, respectively.

 \begin{proposition}
 Let $\alpha \in \RR$ be irrational, and let $a$ be a positive integer.  One has the following.
 \begin{enumerate}
 \item $\min B(\RR\backslash \QQ) = B(\Phi) = \Phi$.
 \item  $B([\overline{2}]) = B(\sqrt{2}) = [\overline{2}]$, where $[\overline{2}] = 1+\sqrt{2} = 2.414213\ldots$.
 \item  If $A(\alpha) \geq a$, then $$B(\alpha) \geq [a,\overline{a,1}] =a+ \sqrt{\frac{1}{4}+\frac{1}{a}}- \frac{1}{2}> a+ \frac{1}{a+1}.$$ 
  \item If $\alpha$ is not continued-fraction equivalent to $\Phi$ or $[\overline{2}]$, then $B(\alpha) \geq 2 +\sqrt{\frac{1}{3}} = 2.577350\ldots > B([\overline{2}])$ 
 \item $\min B((\RR\backslash \QQ)\backslash \{ \Phi\} )= B([\overline{2}])$,  and $B( [\overline{2}])$ is an isolated point of 
 $B(\RR\backslash \QQ)$.
 \end{enumerate}
 \end{proposition}

\begin{proof}
Statement (1) and (2) are clear.   Let $r_n = [a,a,1,a,1,\ldots,a,1]$, where the sequence $a,1$ occurs $n$ times.
If $A(\alpha) = a$, then $a_n(\alpha) \in [1,a]$ for all $n \gg 0$, and therefore,  by (\ref{cfineqa}) and Lemma \ref{intervalfrac}(2),  one has
$$B(\alpha) \geq \limsup_{n \to \infty} r_n =[a,a,1,a,1,a,1,\ldots] = a+ \sqrt{\frac{1}{4}+\frac{1}{a}}-\frac{1}{2}.$$
Statement (3) follows.  Finally, to prove (4) and (5),  we may suppose that $A(\alpha) =2$ and $\alpha$ is not continued-fraction equivalent to $[\overline{2}]$.  Then the regular continued fraction expansion of $\alpha$ contains infinitely many sequences of the form $1,2$, and therefore
$$B(\alpha) \geq [2,1,1,2,1,2,1,2,\ldots]  = 2 +\sqrt{\frac{1}{3}}  > 1+\sqrt{2}.$$
Statements (4) and (5) follow.
\end{proof}

 \begin{proposition}
 Let $\alpha \in \RR$ be irrational, and let $a$ be a positive integer.  One has the following.
 \begin{enumerate}
 \item $\min C(\RR\backslash \QQ) = C(\Phi) = \Phi^2$.
 \item  $C([\overline{2,1}]) = C(\sqrt{3}) = [\overline{2,1}][\overline{1,2}] = 2+\sqrt{3} = 3.732050\ldots$, where $[\overline{2,1}] = 1+\sqrt{3}$.
 \item  If $A(\alpha) \geq a$, then $$C(\alpha) \geq 1+a[\overline{1,a}] = 1+[\overline{a,1}] = 1+ \frac{a+\sqrt{a(a+4)}}{2}> a+1+ \frac{a}{a+1}.$$ 
  \item If $\alpha$ is not continued-fraction equivalent to $\Phi$ or to $[\overline{2,1}]$, then $C(\alpha) \geq 3 +\sqrt{\frac{4}{3}} = 4.154700\ldots> C([\overline{2,1}])$ 
 \item $\min C((\RR\backslash \QQ)\backslash \{ \Phi\} )= C([\overline{2,1}])$,  and $C( [\overline{2,1}])$ is an isolated point of 
 $C(\RR\backslash \QQ)$.
 \end{enumerate}
 \end{proposition}

\begin{proof}
Statements (1) and (2) are clear, and statement (3) follows from Lemma \ref{intervalfrac}(1) and the  fact that $$C(\alpha) = 1+\limsup_{n \to \infty} a_{n+1}(\alpha)[a_{n}(\alpha),a_{n-1}(\alpha),\ldots,a_1(\alpha)].$$   Thus,  to prove (4) and (5), we may suppose that  $A(\alpha) =2$ and $\alpha$ is not continued-fraction equivalent to $[\overline{2,1}]$.  Then the regular continued fraction expansion of $\alpha$ contains infinitely many sequences of the form $2,2$ or infinitely many sequences of the form $1,1,2$.  By Lemma \ref{intervalfrac}(2), in the former case, one has
$$C(\alpha) \geq 1+ 2[2,2,1,2,1,2\ldots]  = 4+\sqrt{3} = 4.732050\ldots.$$
In the latter case, one has
 $$C(\alpha) \geq 1+2[1,1,1,2,1,2,\ldots] = 3+\sqrt{\frac{4}{3}} = 4.154700\ldots.$$
Statements (4) and (5) follow.
\end{proof}

Various techniques for studying the Lagrange spectrum, as in \cite{cus},  might allow one to further the study of the images of the functions $R$ and $B_k$ for positive integers $k$.   

\section{Quadratic irrationals}

In this section, we compute the various invariants defined in this chapter for the quadratic irrationals.
As is well known,  the quadratic irrationals are precisely those irrational numbers whose regular continued fraction  expansion is eventually periodic.  For any quadratic irrational $\alpha$, we let $\alpha^*$\index[symbols]{.va  n@$\alpha^*$} denote the conjugate of $\alpha$, so that
 $N_{\QQ(\alpha)/\QQ}(\beta) = \beta\beta^* \in \QQ$ for all $\beta \in \QQ(\alpha)$.
 Note that,  by \cite[Chapter 3 Satz 9]{perron}, if $$\alpha =  [\overline{a_{1},a_{2},\ldots, a_n}]$$ is purely periodic, then
 $$-\frac{1}{\alpha^*} =  [\overline{a_{n},a_{n-1},\ldots, a_1}],$$
 and therefore
 $$N_{\QQ(\alpha)/\QQ}(\alpha) = \alpha\alpha^* = -\frac{ [\overline{a_{1},a_{2},\ldots, a_n}]}{[\overline{a_{n},a_{n-1},\ldots, a_1}]} \in \QQ.$$
 One also has
\begin{align*}
N_{\QQ(\alpha)/\QQ} (\alpha) =- \frac{[a_1, a_2,\ldots,a_{n-1}]}{[a_n, a_{n-1},\ldots,a_{2}]}.
\end{align*}
We have the following result, whose proof is straightforward.
 
 \begin{theorem}\label{MRR}
 Let $\alpha$ be a quadratic irrational,  and let  $\gamma = [\overline{a_1, a_2,\ldots,a_n}]$ be the tail of  its regular continued fraction expansion.
For all  positive integers $k \leq n$ and all nonnegative integers $r$, let
 $$\alpha_k = \alpha_{rn+k}= [\overline{a_{k}, a_{k+1},\ldots, a_n,a_1,\ldots,a_{k-1}}] = S^{\circ (rn+k-1)}(\gamma)$$
 and
 $$\alpha_{-k} = \alpha_{-rn-k} = [\overline{a_{k-1}, a_{k-2},\ldots, a_1,a_n,\ldots,a_{k}}].$$
 For all $k$, one has
 $$\alpha_{-k} = -\frac{1}{\alpha_k^*}$$
 and therefore
  $$\alpha_k + \frac{1}{\alpha_{-k}} = \alpha_k -\alpha_k^*,$$
  and one also has
  $$\beta_k : = (\alpha_{k+1}^*)^2  \left( \frac{\alpha_k- \alpha_{k}^* }{\alpha_{k+1}- \alpha_{k+1}^* } \right) = -\frac{ \alpha_{k+1}^*}{ \alpha_{k+1}} = \frac{1}{\alpha_{k+1}\alpha_{-k-1}}= -\frac{N_{\QQ(\alpha)/\QQ}(\alpha_{k+1})}{\alpha_{k+1}^2}.$$ 
  Moreover, one has
 $$M(\alpha) = \max_{1\leq k \leq n} \left(\alpha_k + \frac{1}{ \alpha_{-k}} \right) =  \max_{1\leq k \leq n} \left(\alpha_k- \alpha_{k}^* \right),$$
\begin{align*}
 \overline{\rho}(\alpha)   = \max_{1\leq k \leq n} \beta_k = \frac{1}{\min_{1\leq k \leq n}(\alpha_{k}\alpha_{-k})} \in (0,1),
 \end{align*}
 \begin{align*}
 \underline{\rho}(\alpha) = \min_{1\leq k \leq n} \beta_k= \frac{1}{\max_{1\leq k \leq n}(\alpha_{k}\alpha_{-k})} \in (0,1),
 \end{align*}
  \begin{align*}
 R(\alpha) = \frac{1}{\underline{\rho}(\alpha)} = \underline{\rho}(\alpha)^* = \max_{1\leq k \leq n}(\alpha_{k}\alpha_{-k}) \in (1,\infty),
 \end{align*}
 and
   \begin{align*}
 \underline{R}(\alpha) = \frac{1}{\overline{\rho}(\alpha)} = \overline{\rho}(\alpha)^* = \min_{1\leq k \leq n}(\alpha_{k}\alpha_{-k}) \in (1,\infty).
 \end{align*}
  For all positive integers $k$, the limit points of the set $\left\{\frac{q_{l+k}(\alpha)}{q_l(\alpha)}: l \in \ZZ_{>0}\right\}$ are
 $$\alpha_{-1}\alpha_{-2}\cdots \alpha_{-k},  \, \alpha_{-2}\alpha_{-3}\cdots \alpha_{-k-1},  \, \ldots,  \, \alpha_{-n}\alpha_{-1}\alpha_{-2}\cdots \alpha_{-k+1},$$ or, equivalently, 
  $$\frac{(-1)^k}{(\alpha_{1}\alpha_{2}\cdots \alpha_{k})^*},   \frac{(-1)^k}{(\alpha_{2}\alpha_{3}\cdots \alpha_{k+1})^*},   \ldots,  \frac{(-1)^k}{ (\alpha_{n}\alpha_{1}\alpha_{2}\cdots \alpha_{k-1})^*}.$$
 More precisely, one has
 $$\lim_{r \to \infty} \frac{q_{rn+l+k}(\alpha)}{q_{rn+l}(\alpha)} = \alpha_{-l-2}\alpha_{-l-3}\cdots \alpha_{-l-k-1} =  \frac{(-1)^k}{( \alpha_{l+2}\alpha_{l+3} \cdots \alpha_{l+k+1})^*}$$
 for all nonnegative integers $k,l$.
  Consequently, one has
  $$B_k(\alpha) = \max(\alpha_{-1}\alpha_{-2}\cdots \alpha_{-k},  \, \alpha_{-2}\alpha_{-3}\cdots \alpha_{-k-1},  \, \ldots,  \, \alpha_{-n}\alpha_{-1}\alpha_{-2}\cdots \alpha_{-k+1})$$
  and
    $$\underline{B}_k(\alpha) = \min(\alpha_{-1}\alpha_{-2}\cdots \alpha_{-k},  \, \alpha_{-2}\alpha_{-3}\cdots \alpha_{-k-1},  \, \ldots,  \, \alpha_{-n}\alpha_{-1}\alpha_{-2}\cdots \alpha_{-k+1})$$
    for all positive integers $k$.
 \end{theorem}

Note that each of the numbers  $\frac{\alpha_k- \alpha_{k}^* }{\alpha_{k+1}- \alpha_{k+1}^* } $ in the theorem is rational, each of the numbers $\beta_k$
generates $\QQ(\alpha)$, in that $\QQ(\beta_k) = \QQ(\alpha)$, and each $\beta_k$ has norm $$N_{\QQ(\alpha)/\QQ} (\beta_k) = 1,$$
and therefore  $N_{\QQ(\alpha)/\QQ} (\overline{\rho}(\alpha)) = N_{\QQ(\alpha)/\QQ} (\underline{\rho}(\alpha))  = 1$.  
Note  also that  the identity
$$\alpha_{k+1}\alpha_{k+1}^* \left( \frac{\alpha_k- \alpha_{k}^* }{\alpha_{k+1}- \alpha_{k+1}^* } \right)  = -1,$$
and therefore the identity
$$ \beta_k = -\frac{ \alpha_{k+1}^*}{ \alpha_{k+1}} = \frac{1}{\alpha_{k+1}\alpha_{-k-1}} = -\frac{N_{\QQ(\alpha)/\QQ}(\alpha_{k+1})}{\alpha_{k+1}^2}$$
of the theorem,
follows from \cite[Section 23 Equation (5)]{perron}.

\begin{corollary}
 Let $\alpha$ be a quadratic irrational whose regular continued fraction expansion has tail $[\overline{a,b}]$ of period length $1$ or $2$, where $a$ and $b$ are positive integers with $a \geq b$.   One has
 $$M(\alpha) = [\overline{a,b}]+\frac{1}{[\overline{b,a}]},$$
 $$ \overline{\rho}(\alpha)   =  \underline{\rho}(\alpha) =  \frac{1}{[\overline{a,b}][\overline{b,a}]} \in (0,1),$$
  $$R(\alpha) = \underline{R}(\alpha) = C(\alpha) =   [\overline{a,b}][\overline{b,a}] \in (1,\infty),$$
    $$B_k(\alpha) =  [\overline{a,b}]^{\lceil k/2 \rceil}[\overline{b,a}]^{\lfloor k/2 \rfloor} \in (1,\infty),$$
    and
    $$\underline{B}_k(\alpha) =  [\overline{a,b}]^{\lfloor k/2 \rfloor}[\overline{b,a}]^{\lceil k/2 \rceil} \in (1,\infty)$$
    for all positive integers $k$,
 where $$[\overline{a,b}] = \frac{ab+\sqrt{ab(ab+4)}}{2b},$$
 $$[\overline{b,a}] = \frac{ab+\sqrt{ab(ab+4)}}{2a},$$
 $$\QQ(\alpha) = \QQ(\sqrt{ab(ab+4)}),$$
and
$$N_{\QQ(\alpha)/\QQ} ([\overline{a,b}])  = \frac{1}{N_{\QQ(\alpha)/\QQ} ([\overline{b,a}])}= -\frac{a}{b}.$$
\end{corollary}

\begin{proposition}
 Let $d$ be a positive integer that is not a square.   The regular continued fraction of  $\sqrt{d}$ is of the form
$$\sqrt{d} = [a_0, \overline{a_1,a_2,\ldots, a_{n-1},2a_0}]$$
with $a_k \leq a_0 = \lfloor \sqrt{d} \rfloor$ for all  $k \leq n-1$.
Consequently, 
$$a_0+\sqrt{d} = [\overline{2a_0,a_{1},a_{2},\ldots, a_{n-1}}]$$
is purely periodic, and one has
$$M(\sqrt{d}) =  2 \sqrt{d}.$$
Moreover, for each $k$ one may write
$$\alpha_{k+1} = \frac{A_k+\sqrt{d}}{B_k},$$
where $A_k$ and $B_k$ are positive integers  with $A_k \leq a_0 < \sqrt{d}$, and therefore
$$\beta_k = -\frac{A_k-\sqrt{d}}{A_k+\sqrt{d}} = 1-\frac{2A_k}{A_k+\sqrt{d}},$$
where  $\alpha= \sqrt{d}$ and the $\alpha_k$ and $\beta_k$ are defined as in Theorem \ref{MRR}.  Consequently,  one has
$$\underline{\rho}(\sqrt{d})= \frac{\sqrt{d}-a_0 }{\sqrt{d}+a_0} \leq \frac{\sqrt{d}-A_j }{\sqrt{d}+A_j}= \overline{\rho}(\sqrt{d}) \leq \frac{\sqrt{d}-1}{\sqrt{d}+1},$$
where $j$ is chosen so that $A_j = \min_{1\leq k \leq n} A_k$, and where also
\begin{align*}
R(\sqrt{d})  = \frac{1}{\underline{\rho}(\sqrt{d})} & = [\overline{2a_0,a_{1},a_{2},\ldots, a_{n-1}}][\overline{a_1,a_2,\ldots, a_{n-1},2a_0}]  \\
 &  =  [\overline{2a_0,a_{n-1},a_{n-2},\ldots, a_{1}}][\overline{a_{n-1},a_{n-2},\ldots, a_{1},2a_0}] \\
 & = \lim_{r \to \infty}  \frac{q_{rn}(\sqrt{d})}{q_{rn-1}(\sqrt{d})} \frac{q_{rn-1}(\sqrt{d})}{q_{rn-2}(\sqrt{d})} \\
 & =  \lim_{r \to \infty} \frac{q_{rn}(\sqrt{d})}{q_{rn-2}(\sqrt{d})} \\
 & = C(\sqrt{d}).
 \end{align*}
\end{proposition}

\begin{proof}
The first statement is well known (e.g.,  \cite[{Section 25}]{perron}), and then the second statement follows from Theorem \ref{MRR}.  Moreover, the third statement follows from \cite[{Sections 23--25}]{perron}, and then the last statement follows from Theorem \ref{MRR}.  
\end{proof}

Unfortunately,  we do not know a formula for $A_j$ in terms of $d$.  

Regarding the L\'evy constant of a quadratic irrational,  one has the following.

 \begin{theorem}[{\cite[Proposition 2.1]{faivre} \cite[p.\ 197]{smith}}]
 Let $\alpha$ be a quadratic irrational,  and let  $[\overline{a_1, a_2,\ldots,a_n}]$ be the tail of  its regular continued fraction expansion.
For all  positive integers $k \leq n$, let
 $$\alpha_k = [\overline{a_{k}, a_{k+1},\ldots, a_n,a_1,\ldots,a_{k-1}}] = S^{\circ (k-1)}([\overline{a_1, a_2,\ldots,a_n}]).$$
 One has
 $$\lambda(\alpha) = \frac{1}{n} \sum_{k = 1}^n \log \alpha_k,$$
 or, equivalently, 
 $$\exp \lambda(\alpha) = \lim_{n \to \infty} q_n(\alpha)^{1/n}  =  \left(\prod_{k = 1}^n \alpha_k\right)^{1/n}.$$
 Let $d$ be a positive integer that is not a square, and let $\alpha = \sqrt{d}$.
  Suppose that $n = N(\sqrt{d})$ is minimal,  so that $N(\sqrt{d})$ is the length of the period of the regular continued fraction expansion of $\sqrt{d}$.  Then one has
$$ \eta_0(\sqrt{d}) = (\exp \lambda(\alpha))^{N(\sqrt{d})} = \prod_{k = 1}^{N(\sqrt{d})} \alpha_k,$$
or, equivalently, 
 $$\lambda(\sqrt{d}) = \frac{1}{N(\sqrt{d})}\log \eta_0(\sqrt{d}),$$
where $ \eta_0(\sqrt{d})$ is the fundamental unit of $\ZZ[\sqrt{d}]$, that is,  where  $ \eta_0(\sqrt{d}) $  is the solution $\varepsilon \in \ZZ[\sqrt{d}]$ to $N_{\QQ(\sqrt{d})/\QQ}(\varepsilon) = \pm 1$ (or $N_{\QQ(\sqrt{d})/\QQ}(\varepsilon) = (-1)^{N(\sqrt{d})}$) with both $\varepsilon+\varepsilon^* > 0$ and $\varepsilon-\varepsilon^* > 0$ minimal.
 \end{theorem} 
 
 The {\it regulator} $R_K$ of a real quadratic number field $K = \QQ(\sqrt{d})$,  where $d >1$ is an integer that is not a square,  is $\log \eta_K$, where $\eta_K$ is the {\it fundamental unit} of $K$,  or, equivalently, where  $\varepsilon = 2\eta_K$ is the solution   $\varepsilon \in \ZZ[\sqrt{\Delta}]$ of $N_{K/Q}(\varepsilon) = \pm 4$ (or of $N_{K/Q}(\varepsilon) =(-1)^{N(\sqrt{d})}4$) with both $\varepsilon+\varepsilon^* > 0$ and $\varepsilon-\varepsilon^* > 0$ minimal,  where $\Delta$ is the {\it discriminant} of $K$  \cite[p.\ 15 and pp.\ 42--43]{neukirch}.  Note that, if $d >1$ is squarefree and not congruent to $1$ modulo $4$, then the fundamental unit of $\QQ(\sqrt{d})$ is the same as the fundamental unit of $\ZZ[\sqrt{d}]$ as defined in the proposition.  Thus, we have the following corollary.
 
\begin{corollary}\label{reglam}
Let $K = \QQ(\sqrt{d})$ be a real quadratic number field,  where $d \not \equiv 1 \; (\bmod \; 4)$ is a squarefree positive integer.   Then
$\eta_K = \eta_0(\sqrt{d})$, and the regulator of $K$ is given by
$$R_K = \log \eta_K = N(\sqrt{d})\lambda(\sqrt{d}),$$ and thus
$$\lambda(\sqrt{d}) = \frac{1}{ N(\sqrt{d})}R_K,$$
 where $N(\sqrt{d})$ is the length of the period of the regular continued fraction of $\sqrt{d}$.
\end{corollary}

\begin{example} One has the following.
\begin{enumerate}
\item  $\sqrt{2} = [1,\overline{2}]$, and therefore  $$M(\sqrt{2}) = [\overline{2}]+\frac{1}{[\overline{2}]}= 2 \sqrt{2},$$
$$\overline{\rho}(\sqrt{2}) = \underline{\rho}(\sqrt{2}) = \frac{1}{[\overline{2}]^2}= \frac{ \sqrt{2}-1}{\sqrt{2}+1} = 3-2\sqrt{2},$$
$$R(\sqrt{2}) = \underline{R}(\sqrt{2})=  3+2\sqrt{2},$$
$$\eta_0(\sqrt{2}) = [\overline{2}] = 1+\sqrt{2},$$
and
$$\lambda(\sqrt{2}) = \log (1+\sqrt{2}) =  R_{\QQ(\sqrt{2})}.$$
\item  $\sqrt{3} = [1,\overline{1,2}]$, and therefore  $$M(\sqrt{3}) = [\overline{2,1}]+\frac{1}{[\overline{1,2}]}= 2 \sqrt{3},$$
$$\overline{\rho}(\sqrt{3}) = \underline{\rho}(\sqrt{3}) = \frac{1}{[\overline{1,2}][\overline{2,1}]} = \frac{ \sqrt{3}-1}{\sqrt{3}+1}= 2-\sqrt{3},$$
$$R(\sqrt{3}) = \underline{R}(\sqrt{3}) = 2+\sqrt{3},$$
$$\eta_0(\sqrt{3}) = [\overline{1,2}][\overline{2,1}]= 2+\sqrt{3},$$
and
$$\lambda(\sqrt{3}) = \frac{1}{2} \log (2+\sqrt{3}) =  \frac{1}{2} R_{\QQ(\sqrt{3})}.$$
Note that $R(\sqrt{2}) > R(\sqrt{3})$, despite the fact that $M(\sqrt{2} )< M(\sqrt{3})$.
\item  $\sqrt{41} = [6,\overline{2,2,12}]$, and therefore  $$M(\sqrt{41}) = [\overline{12,2,2}]+\frac{1}{[\overline{2,2,12}]}= 2 \sqrt{41},$$
$$\overline{\rho}(\sqrt{41}) = \beta_1 =  \frac{1}{[\overline{2,12,2}]^2}= \frac{ \sqrt{41}-4}{\sqrt{41}+4}  = 0.231000\ldots,$$
$$\underline{\rho}(\sqrt{41}) = \beta_2 = \beta_3=  \frac{1}{[\overline{12,2,2}][\overline{2,2,12}]} = \frac{ \sqrt{41}-6}{\sqrt{41}+6} = 0.032501\ldots,$$
$$\eta_0(\sqrt{41}) = [\overline{2,2,12}][\overline{2,12,2}][\overline{12,2,2}]= 32+5\sqrt{41} = \eta_{\QQ(\sqrt{41})},$$
and
$$\lambda(\sqrt{41}) = \frac{1}{3} \log (32+5\sqrt{41}) = \frac{1}{3} R_{\QQ(\sqrt{41})}.$$
Note that $41$ is the smallest positive integer $d$ such that the continued fraction of  $\sqrt{d}$ has period length $3$.  This particular example shows that the inequality $\overline{\rho}(\sqrt{d}) \leq \frac{\sqrt{d}-1}{\sqrt{d}+1}$ can be strict.  Moreover,   $\eta_0(\sqrt{41})$ is equal to the fundamental unit  $ \eta_{\QQ(\sqrt{41})}$ of $\QQ(\sqrt{41})$, despite the fact that $41 \equiv 1 \; (\bmod\; 4)$ and that the period of $$\frac{1+\sqrt{41}}{2} = [3,\overline{1,2,2,1,5}]$$ is $5$, not $3$,  and one has
$$\lambda\left(\frac{1+\sqrt{41}}{2}\right)=  \frac{1}{5}\log \eta_{\QQ(\sqrt{41})} =  \frac{1}{5} R_{\QQ(\sqrt{41})} = \frac{3}{5} \lambda(\sqrt{41}).$$  This does not hold, for example, for $\QQ(\sqrt{d} )$ for $d = 5$ or $d = 13$,  i.e., in those cases, the fundamental unit of $\QQ(\sqrt{d} )$  is not equal to $\eta_0(\sqrt{d})$.
\end{enumerate}
\end{example}

 \begin{conjecture}
   Let $d \equiv 1 \;(\bmod\; 4)$ be a positive integer that is not a square, and let $\alpha = \frac{1+\sqrt{d}}{2}$.
 Then one has
$$(\exp \lambda(\alpha))^{N(\alpha)}  =  (\eta_{\QQ(\alpha)})^{K(\alpha)},$$
or, equivalently, 
 $$N(\alpha) \lambda(\alpha) =K(\alpha)\log \eta_{\QQ(\alpha)} =  K(\alpha)R_{\QQ(\alpha)},$$
 for some positive integer $K(\alpha)$, 
where $ \eta_{\QQ(\alpha)} $ is the fundamental unit and $R_{\QQ(\alpha)}$ is the regulator of $\QQ(\alpha)$.   Moreover, if $d$ is squarefree, then $K(\alpha) = 1$.
 \end{conjecture} 
 
We are unable to make a conjecture on the value of $K(\alpha) = \frac{N(\alpha) \lambda(\alpha)}{R_{\QQ(\alpha)}}$ beyond the conjecture that it is a positive integer depending on $\alpha$.  Note that Corollary \ref{reglam}  states that  the conclusion of the conjecture holds for $\alpha = \sqrt{d}$  for $d \not \equiv 1 \;(\bmod\; 4)$ a positive integer that is not a square, where $K(\alpha)= 1$.

 \begin{example}
$d = 621 =3^3 \cdot 23$ is the smallest nonsquarefree nonsquare congruent to $1$ modulo $4$ such that $K(\alpha)$ (exists and) is equal to $1$: one has $\alpha:= \frac{1+\sqrt{621}}{2} = [12,\overline{1,23}]$,
 $$ (\exp \lambda(\alpha))^{2} =  [\overline{1,23}][\overline{23,1}] = \frac{25 +  \sqrt{621}}{2} = \eta_{\QQ(\alpha)},$$
 $$\sqrt{621} = [24, \overline{1, 11, 2, 11, 1, 48}],$$
 $$\eta_0(\sqrt{621}) =(\eta_{\QQ(\alpha)})^3,$$
 and
 $$\lambda(\sqrt{621})= \frac{1}{6} \log \eta_0(\sqrt{621}) = \frac{1}{2} \log \eta_{\QQ(\alpha)} = \lambda(\alpha).$$
 In particular, one has $\lambda(\sqrt{621}) =  \lambda(\alpha)$, even though $\eta_0(\sqrt{621}) \neq \eta_{\QQ(\alpha)}$.
 \end{example}

 Finally, we  note the following.
 
 \begin{proposition}
Let $Q$ denote the set of all quadratic irrationals.  Then $\Phi^3 = 4.236067\ldots$ is a limit point from above of the set of all limit points of the set of all limit points of $R(Q)$, and $\Phi^4  = 6.854101\ldots$ is also a limit point of $R(Q)$.
 \end{proposition}
 
 \begin{proof}
 For each nonnegative integer $n$, let
 $$\alpha_n = [\overline{2,1,1,1,\ldots,1}],$$
 where the $2$ is followed by a sequence of $n$ $1$s.
 One has
 $$\lim_{n \to \infty} \alpha_n =  \Phi^2$$
 and
 $$R(\alpha_n) = [\overline{2,1,1,1,\ldots,1}][\overline{1,1,1,\ldots,1,2}],$$
 and therefore
 $$\lim_{n \to \infty} R(\alpha_n) = [2,\overline{1}][\overline{1}] = (1+\Phi)\Phi = \Phi^3 = \Phi^2+\Phi= 2+\sqrt{5} = 4.236067\ldots,$$
 where
 $$R(\alpha_{n+1}) < R(\alpha_{n+3})< \Phi^3 <R(\alpha_{n+2}) < R(\alpha_{n})$$
 for all even integers $n$.   Thus,  the number $\Phi^3 = \Phi^2+\Phi$ is a limit point of $R(Q)$ both from above and below.
 Repeating the argument above but with
 $$\alpha_n = [\overline{2,1,1,1,\ldots,1,2}],$$ we see that
 $$[2,\overline{1}][2,\overline{1}] = \Phi^4$$
 is also a limit point of $R(Q)$ both from above and below.

Considering the sequence
$$\beta_n = [\overline{2,1,2,1,1,1,\ldots,1}]$$
of  $2,1,2$ followed by $n$ $1$s,  one sees that
$$R(\beta_n) = \beta_n [\overline{1,1,1,\ldots,1,2,1,2}],$$
which converge to
$$\lim_{n \to \infty} R(\beta_n) = [2,1,2,\overline{1}][\overline{1}] =  \frac{13\Phi+1}{5} = 4.406888\ldots$$
from above and below.
 Repeating this argument, we see that
 $$[2,1,2,1,2\ldots,1,2,\overline{1}][\overline{1}] $$
 is a limit point from above and below, for all finite sequences $2,1,2,1\ldots,2$.  These limit points themselves converge from below to the limit
 $$[\overline{2,1}][\overline{1}] = (1+\sqrt{3})\Phi = 4.420551\ldots.$$
 
Further limit points from above and below are
 $$[2,1,1,1,2,\overline{1}][\overline{1}]  = \frac{83+81\sqrt{5}}{62} = 4.260024\ldots.$$
  $$[2,1,1,1,2,1,1,1,2,\overline{1}][\overline{1}]  =\frac{3912 + 3911 \sqrt{5}}{2971} = 4.260269\ldots,$$
  and so on, which themselves converge from below to the limit
    $$[\overline{2,1,1,1}][\overline{1}]  =\frac{3 + 2\sqrt{6}}{3}\Phi = 4.260272\ldots.$$
   By the same token, each of the values
    $$[\overline{2,1,1,1,1,1}][\overline{1}], [\overline{2,1,1,1,1,1,1,1}][\overline{1}],[\overline{2,1,1,1,1,1,1,1,1,1}][\overline{1}],\ldots$$
is a limit from below of limit points from above and below,  and those values converge from above to
   $$[2,\overline{1}][\overline{1}] = (1+\Phi)\Phi = \Phi^3 = 4.236067\ldots.$$
   Thus,  $\Phi^3$  is a limit point from above of the set of limit points of the set of limit points of $R(Q)$.  
  \end{proof}

 \begin{problem}
What is the smallest limit point of  $R(\RR\backslash \QQ)$?  Is it equal to $\Phi^3$? 
 \end{problem}

\chapter{Conjectures}

Throughout this book, we have considered various conjectures  in analytic number theory---some well-known and some new---in relation to the degree formalisms.  See, for example,  Examples \ref{degexamples}(1)--(8),  which include the abc conjecture,  Vinogradov's conjecture,  and conjectures regarding the Dirichlet divisor problem, the Gauss circle problem, the least primitive root modulo a prime $p$, and the least prime in an arithmetic progression.  Further degree and logexponential degree conjectures are listed in Tables \ref{tab1}--\ref{tab1b}.   In this chapter, we discuss several of these and other conjectures,  accompanied by modest numerical and graphical support.  All numerical tables referred to in this chapter appear at the end of the chapter, in Section 14.7.

\section{Iterated logarithmic degree statistics}

Let $f \in \RR^{\RR_\infty}$.  Writing $$|f(x)| = x^{L_f(x)}$$ for a unique function $L_f$, one has  
$$L_f(x) = \frac{\log |f(x)|}{\log x}$$\index[symbols]{.g k@$L_f(x)$}
and  therefore
$$\deg f = \limsup_{x \to \infty} L_f(x).$$  
Computing (resp., graphing) the function $L_f(x)$ or the function $L_f(e^x) = \frac{\log |f(e^x)|}{x}$ on as large and as wide a domain as possible is useful for obtaining numerical (resp., graphical) support for a conjectural value of $\deg f$.  The main advantage of the statistic $L_f(x)$ over the statistic $\frac{|f(x)|}{x^t}$ for any $t$ is that it does not depend on a prior guess for $\deg f$.   However, if one knows already what $\deg f$ is, and if that degree is finite, then obviously the statistic $\frac{|f(x)|}{x^{\varepsilon+\deg f}}$ for any $\varepsilon \geq 0$ provides further relevant information.  On the other hand, for  functions $f$ where $\deg f$ is not known, such as $f(x) = \li(x)-\pi(x)$ and $f(x) = M(x)$, the statistic $L_f(x)$ assists one in making an educated guess for $\deg f$, after which one may then move on to computing $\frac{|f(x)|}{x^{\deg f}}$.

\begin{example}
Figure \ref{eureka100} provides graphs of the functions $$L_{\li-\pi}(e^x) = \frac{\log|\li(e^x)-\pi(e^x)|}{x}$$
and
$$L_{\li-\Ri}(e^x) = \frac{\log|\li(e^x)-\Ri(e^x)|}{x},$$
which are just the graphs of the functions $L_{\li-\pi}(x)$ and $L_{\li-\Ri}(x)$  but on a lin-log scale.  Since $\li(x)-\Ri(x) \sim \frac{\sqrt{x}}{\log x} \ (x \to \infty)$, one has
$$\lim_{x \to \infty} L_{\li-\Ri}(x) = \lim_{x \to \infty} L_{\li-\Ri}(e^x) = \tfrac{1}{2}.$$
The Riemann constant $\Theta = \deg(\li-\pi)$ is exactly the lim sup of the blue curve in Figure \ref{eureka100} as $x \to \infty$, that is, one has
\begin{align*}
\Theta & = \limsup_{x \to \infty}  L_{\li-\pi}(x) = \limsup_{x \to \infty} L_{\li-\pi}(e^x).
\end{align*}
\end{example}

\begin{figure}[h!]
\includegraphics[width=70mm]{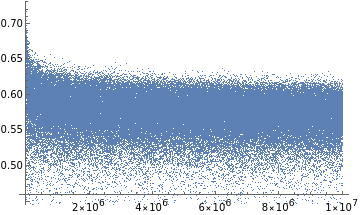}
    \caption{\centering Plot of $L_{H}(n)$ for the multiples $n$ of $100$ less than $10^7$}
\label{Ldegree}
\end{figure}

\begin{figure}[h!]
\includegraphics[width=70mm]{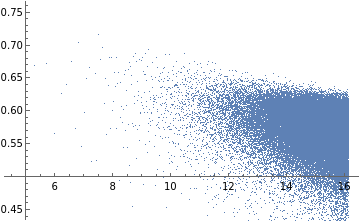}
    \caption{\centering Lin-log plot of $L_{H}(n)$ for multiples $n$ of $100$}
\label{LdegreeB}
\end{figure}

\begin{figure}[h!]
\includegraphics[width=70mm]{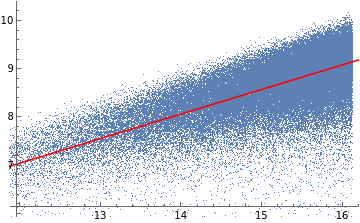}
    \caption{\centering Log-log plot of $|H(n)|$ for multiples $n = e^r$ of $100$, along with best-fit line $0.949479+0.508874 r$}
\label{Ldegree2}
\end{figure}

\begin{figure}[h!]
\includegraphics[width=70mm]{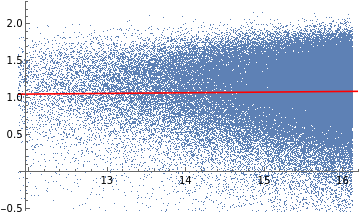}
    \caption{\centering Log-log plot of $|H(n)|$, with $0.5r$ subtracted, for multiples $n = e^r$ of $100$, along with best-fit line $0.949479+0.008874 r$}
\label{Ldegree3}
\end{figure}

\begin{example}
An example, based on the Gauss circle problem in Example \ref{degexamples}, is shown in   Figure \ref{Ldegree}, which provides a plot of $$L_H(n) = \frac{\log |H(n)|}{\log n}$$ for the multiples $n$ of $100$ less than $10^7$, where $H(r) = L(r)-\pi r^2$.  One might be tempted to conjecture on the basis of Figure \ref{Ldegree} that $\deg H > \frac{1}{2}$.  However, Figure \ref{LdegreeB}, which provides a lin-log plot of $L_H(n)$,  reveals more clearly a marked decrease in the local peaks.  Moreover, the conjecture  $\deg H =  \frac{1}{2}$ is supported by the log-log plot  of $|H(n)|$ shown in Figure \ref{Ldegree2}, which  ``linearizes'' the given data and has best fit line $0.949479+0.508874r$.  (This  suggestion and computation is due to Aditya Baireddy.)  Figure \ref{Ldegree2}, equivalently, is a (lin-lin) plot of the function $\log|H(e^r)|$, which is just the function $L_H(e^r)$ multiplied by $r$, for multiples $n = e^r$ of $100$.  When $0.5 r$ is subtracted from the data in  Figure \ref{Ldegree2}, we obtain  Figure \ref{Ldegree3}; that is, Figure \ref{Ldegree3} is a plot of the function $\log|H(e^r)|-0.5r = (L_H(e^r)-0.5)r$ for multiples $e^r$ of $100$.  Plotting that function divided by $r$ would yield (a portion of) the lin-log plot of $L_H(n)$ in Figure \ref{LdegreeB} but with $0.5$ subtracted.  
Moreover, one has $\deg H = \frac{1}{2}$ if and only if the function $\log|H(e^r)|-0.5r$ plotted in Figure \ref{Ldegree3}, divided by $r$, has a lim sup of $0$ as $r \to \infty$ (where $e^r$ ranges over all positive integers).  Both of these conditions hold if the function $\log|H(e^r)|-0.5r$ is bounded above on $\RR_{>0}$, or, equivalently, if $H(n) = O(\sqrt{n})  \ (n \to \infty)$.

The relationship between $L_H(n)$ and the log-log plot of $|H(n)|$ generalizes to any function $f \in \RR^{\RR_\infty}$.  Specifically, a plot of $L_f(e^x)$ can be obtained from a log-log plot of $|f(x)|$, i.e., from a plot of $\log |f(e^x)|$, simply  by dividing by $x$. 
\end{example}

We extend the definition of $L_f(x)$ to the iterated logarithmic degrees $\degl_N f$ of order $N > 0$  as follows.  Suppose that $\degl_i f \in \RR$ for $i < N$.  Given conjectural values $d_i \in \RR$ for $\degl_i f$ for $i < N$, we set
\begin{align*}
L_f[d_0,d_1,\ldots,d_{N-1}](x) =  \frac{\log |f(x)|-d_0\log x- d_1 \log^{\circ 2}x - \cdots - d_{N-1} \log^{\circ (N)}x}{\log^{\circ (N+1)} x},
\end{align*}\index[symbols]{.g m@$L_f[d_0,d_1,\ldots,d_{N-1}](x)$}
so that
$$|f(x)| = x^{d_0}(\log x)^{d_1}(\log  \log x)^{d_2}\ldots (\log^{\circ(N-1)}x)^{d_{N-1}}  (\log^{\circ N}x)^{L_f[d_0,d_1, \ldots, d_{N-1}](x)}$$
for all $x > \exp^{\circ N}(0)$ in $\dom f$.
If in fact $d_k = \degl_k f$ for $k < N$, then one has
$$\degl_N f = \limsup_{x \to \infty} L_f[d_0,d_1, \ldots, d_{N-1}](x),$$ while if $d_i \neq \degl_i f$ for some $i < N$, then $\limsup_{x \to \infty} L_f[d_0,d_1, \ldots, d_{N-1}](x) = \pm \infty$.  Thus, if the conjectural values $d_k \in \RR$ for $\degl_k f$ for $i < N$ are correct, then the statistic $L_f[d_0,d_1,\ldots,d_{N-1}](x)$ provides information regarding $\degl_N f$.  Note, however, that this statistic, unlike $L_f(x)$, is dependent on choice of bases for the various logs.  Our convention is to use base $e$ for all of them.

\begin{remark}[Renormalizations of the iterated logarithmic degree statistics]\label{Lnorm}
Suppose that in the definition of $L_f[d_0,d_1,\ldots,d_{N-1}](x)$ we replace the sequence of functions $\log^{\circ k}x$ with a sequence $l = \{l_k(x)\}$ of eventually positive functions with $l_k(x) \asymp \log^{\circ k}x \ (x \to \infty)$ for all $k$. In other words, for such a sequence $l$,  consider
\begin{align*}
L_f[l; d_0,d_1,\ldots,d_{N-1}](x) =  \frac{\log |f(x)|-d_0\log l_0(x)- d_1 \log l_1(x) - \cdots - d_{N-1} \log l_{N-1}(x)}{\log l_N(x)},
\end{align*}
so that
$$|f(x)| = l_0(x)^{d_0}l_1(x)^{d_1}l_2(x)^{d_2}\ldots l_{N-1}(x) ^{d_{N-1}} l_N(x)^{L_f[l;d_0,d_1, \ldots, d_{N-1}](x)}$$
for all sufficiently large $x$.  Then the relationship  $$\degl_N f = \limsup_{x \to \infty} L_f[l;d_0,d_1, \ldots, d_{N-1}](x)$$ still holds, again provided that  $d_k = \degl_k f$ for $i < N$.   The statistic $L_f[l; d_0,d_1,\ldots,d_{N-1}](x)$  is a renormalization of the statistic $L_f[d_0,d_1,\ldots,d_{N-1}](x)$ with respect to the sequence $l$.  In this text we employ the statistic $L_f[d_0,d_1,\ldots,d_{N-1}](x)$.  Although there are other natural choices for the sequence $l$ besides $\log^{\circ k} x$, it is important to avoid choosing a specific $l$ for the purpose of manipulating the data to fit a desired outcome. 
\end{remark}

It is not uncommon, given a relevant function $f(x)$, for someone to have computed upper bounds of $|f(x)|$ on various intervals, often of the form $[b^n, b^{n+1}]$ for some  fixed $b > 1$, in order to investigate numerically the asymptotic behavior of $f(x)$.  The following proposition provides a bridge between such bounds and the statistics  $L_f[d_0,d_1,\ldots,d_{N-1}](x)$.

\begin{proposition}\label{deltaprop} 
Let $f(x)$ be piecewise bounded and continuous function on $(0, \infty)$ such that $f(x^-) \leq f(x) \leq f(x^+)$ or $f(x^+) \leq f(x) \leq f(x^-)$ for all $x > 0$, and suppose that $\deg f < \infty$.  Let $b > 1$, and let $n$ be a nonnegative integer.  Let $M_n$ denote the supremum  of $|f(x)|$ on the interval $[b^n,b^{n+1}]$,  attained at or approached to the immediate left or right of some $a_n \in [b^n,b^{n+1}]$.  One has the following.
\begin{enumerate}
\item $\displaystyle \sup_{x \in [b^n,b^{n+1}]} \frac{\log |f(x)|}{\log x} \leq \frac{\log M_n}{\log a_n}+ \frac{|\log M_n|}{n(n+1)\log b}.$
\item $\displaystyle \deg f = \limsup_{n \to \infty} \frac{\log M_n}{\log a_n}$.
\item $\displaystyle \sup_{x \in [b^n,b^{n+1}]} \frac{\log |f(x)|}{\log \log x} \leq \frac{\log M_n}{\log \log a_n}+ \frac{|\log M_n|}{\log(n\log b )}-\frac{|\log M_n|}{\log ((n+1)\log b)}$ if $n  >  \frac{1}{\log b}$. 
\item  $\displaystyle \degl_1 f = \limsup_{n \to \infty} \frac{\log M_n}{\log \log a_n}$ if $\deg f = 0$.
\item $\displaystyle \sup_{x \in [b^n,b^{n+1}]} \frac{\log |f(x)|}{\log \log \log x} \leq \frac{\log M_n}{\log \log \log a_n}+ \frac{|\log M_n|}{\log \log(n\log b )}-\frac{|\log M_n|}{\log \log ((n+1)\log b)} $ if $n 
> \frac{e}{\log b}$. 
\item  $\displaystyle \degl_2 f = \limsup_{n \to \infty} \frac{\log M_n}{\log \log a_n}$ if $\deg f  = \degl_1 f  = 0$.
\end{enumerate}
\end{proposition}

\begin{proof}
Let $D(x) = \log |f(x)|,$ so that $L_n = \log M_n$ is the supremum of $D(x)$ on $[b^n,b^{n+1}]$, and where $L_n = D(a_n)$, $L_n  =  D(a_n^+)$, or $L_n  =D(a_n^-)$.     We assume that $L_n$ is positive; the proof when $L_n$ is negative is similar.  Let $x \in [b^n,b^{n+1}]$.  If $x \geq a_n$, then 
$$\frac{D(x)}{\log x} \leq \frac{L_n}{\log a_n}.$$
On the other hand, if $x < a_n$, then
$$ \frac{D(x)}{\log x} -\frac{L_n}{\log a_n} \leq  L_n\left(\frac{1}{\log x}-\frac{1}{\log a_n}\right)  \leq \frac{L_n}{n(n+1)\log b}.$$
It follows that
$$\frac{D(x)}{\log x} \leq \frac{L_n}{\log a_n}+ \frac{L_n}{n(n+1)\log b} $$
for all $x \in [b^n,b^{n+1}]$, 
which implies statement (2).    Moreover, since $f(x) = o(x^t) \ (x \to \infty)$ for some $t > 0$, one has $D(x) \leq t\log x$ for sufficiently large $x$, and therefore
$L_n \leq t\log(b^{n+1}) = t(n+1)\log b$ for sufficiently large $n$, so that
$\lim_{n \to \infty} \frac{L_n}{n(n+1)\log b}  =  0.$
It follows that
$$\deg f = \limsup_{x \to \infty} \frac{D(x)}{\log x} \leq \limsup_{n \to \infty} \left( \frac{L_n}{\log a_n}+ \frac{L_n}{n(n+1)\log b}\right) = \limsup_{n \to \infty}\frac{L_n}{\log a_n}.$$
Moreover,  since $\frac{L_n}{\log a_n} $ is equal to $D(a_n)$, $D(a_n^+)$, or $D(a_n^-)$, one  has $\sup_{x \in [b^n,b^{n+1}]} \frac{D(x)}{\log x} \geq \frac{L_n}{\log a_n}$ and therefore $$\deg f = \limsup_{x \to \infty} \frac{D(x)}{\log x} \geq  \limsup_{n \to \infty}\frac{L_n}{\log a_n}.$$  
This proves statement (2).

A similar argument yields
$$\frac{D(x)}{\log \log x} \leq \frac{L_n}{\log \log a_n}+ \frac{L_n\log(1+1/n)}{\log(n\log b)\log ((n+1)\log b)}\leq \frac{L_n}{\log \log a_n}+ \frac{L_n}{n\log(n\log b)\log ((n+1)\log b)}$$
for all $x \in [b^n,b^{n+1}]$, provided that $\log(n \log b )> 0$, where
$\lim_{n \to \infty} \frac{L_n}{n\log(n\log b)\log ((n+1)\log b)}= 0,$
and therefore
$$ \limsup_{x \to \infty} \frac{D(x)}{\log \log x}  = \limsup_{n \to \infty}\frac{L_n}{\log \log a_n}.$$
If $\deg f = 0$, then the lim sup on the left above is precisely $\degl_1 f$.  This yields statements (3) and (4).

Finally, one has 
\begin{align*}
\frac{D(x)}{\log \log \log x} & \leq \frac{L_n}{\log \log \log a_n}+ \frac{L_n\log\left(\frac{\log ((n+1)\log b)}{\log (n\log b)} \right)}{\log \log(n\log b)\log\log ((n+1)\log b)} \\ & \leq \frac{L_n}{\log \log \log a_n}+ \frac{L_n}{n\log n\log\log(n\log b)\log \log ((n+1)\log b)} 
\end{align*}
for all $x \in [b^n,b^{n+1}]$, provided that $\log \log (n \log b) > 0$, where 
$$\lim_{n \to \infty} \frac{L_n}{n\log n\log\log(n\log b)\log \log ((n+1)\log b)}  = 0,$$
and therefore
$$ \limsup_{x \to \infty} \frac{D(x)}{\log \log \log x}  = \limsup_{n \to \infty}\frac{L_n}{\log \log \log a_n}.$$
If $\deg f = 0$ and $\degl_1 f = 0$, then the lim sup on the left above is precisely $\degl_2 f$.    This yields statements (5) and (6).
\end{proof}

\section{Conjectures concerning $\li(x)-\pi(x)$}

 In this section, we use the iterated logarithmic degree statistics $L_f[d_0,d_1, d_2,\ldots,d_{N-1}](x)$ described in  the previous section to provide some modest evidence for Conjectures \ref{eurekaconjecture} and \ref{eurekaconjecture2} regarding the error $\li(x)-\pi(x)$ in the approximation $\pi(x) \approx \li(x)$.

\begin{figure}[h!]
\includegraphics[width=70mm]{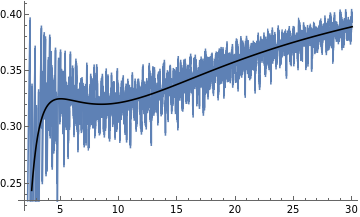}
    \caption{\centering Graphs of $L_{\li-\pi}(e^x)$ and $L_{\li-\Ri}(e^x)$}
\label{eureka1}
\end{figure}

 Figure \ref{eureka1}  (which is the same as Figure \ref{eureka100} of the Preface) provides a graph of $L_{\li-\pi}(e^x)$ and $L_{\li-\Ri}(e^x)$ on the interval $[2,30]$.    Notice that the smooth graph of $L_{\li-\Ri}(e^x)$ appears to track very closely the ``center'' of the fluctuating graph of $L_{\li-\pi}(e^x)$.  Since $\deg(\li-\Ri) = \lim_{x \to \infty} L_{\li-\Ri}(e^x) = \frac{1}{2}$, the Riemann hypothesis holds if and only if the fluctuations of the blue curve above the  black curve tend to $0$ as $x \to \infty$.  We venture that $L_{\li-\pi}(x)$ is one of the most natural statistics for measuring the conformity or nonconformity of $\li(x)-\pi(x)$ to the Riemann hypothesis $O$ bound (\ref{vk}).

\begin{figure}[h!]
\includegraphics[width=70mm]{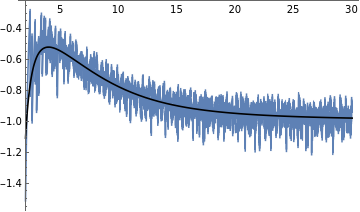}
  \caption{\centering Graphs of $L_{\li-\pi}[\frac{1}{2}](e^x)$ and $L_{\li-\Ri}[\frac{1}{2}](e^x)$}
 \label{eureka2}
\end{figure}

\begin{figure}[h!]
\includegraphics[width=70mm]{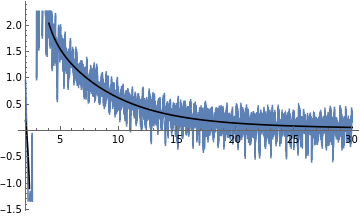}
    \caption{\centering Graphs of $L_{\li-\pi}[\frac{1}{2},-1](e^x)$ and $L_{\li-\Ri}[\frac{1}{2},-1](e^x)$}
   \label{eureka3}
\end{figure}

Figure \ref{eureka2} provides the graphs of $L_{\li-\pi}[\frac{1}{2}](e^x)$ and $L_{\li-\Ri}[\frac{1}{2}](e^x)$, while Figure \ref{eureka3} provides the graphs of $L_{\li-\pi}[\frac{1}{2},-1](e^x)$ and $L_{\li-\Ri}[\frac{1}{2},-1](e^x)$, all on the interval $[2,30]$.  Notice again that the graph of $L_{\li-\Ri}[\frac{1}{2}](e^x)$ appears to track the ``center'' of the fluctuating graph of $L_{\li-\pi}[\frac{1}{2}](e^x)$, and likewise for $L_{\li-\Ri}[\frac{1}{2},-1](e^x)$ and $L_{\li-\pi}[\frac{1}{2},-1](e^x)$.  Since $\degl_1(\li-\Ri) = \lim_{x \to \infty} L_{\li-\Ri}[\frac{1}{2}](e^x) = -1$ and $\degl_2(\li-\Ri) = \lim_{x \to \infty} L_{\li-\Ri}[\frac{1}{2},-1](e^x) = 0$, Conjectures \ref{eurekaconjecture} and \ref{eurekaconjecture2} hold, respectively, if the fluctuations of the blue curves above the respective black curves tend to $0$ as $x \to \infty$.

\begin{figure}[h!]
\includegraphics[width=70mm]{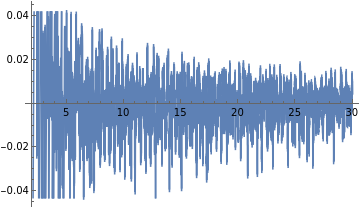}
    \caption{\centering Graph of $L_{f}(e^x)$, where $f = \frac{\li-\pi}{\li-\Ri}$}
\label{eureka4}
\end{figure}

\begin{figure}[h!]
\includegraphics[width=70mm]{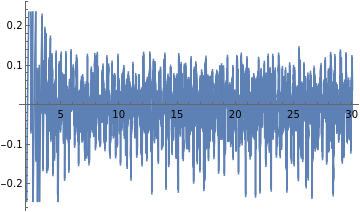}
  \caption{\centering Graph of $L_{f}[0](e^x)$, where $f = \frac{\li-\pi}{\li-\Ri}$}
 \label{eureka5}
\end{figure}

\begin{figure}[h!]
\includegraphics[width=70mm]{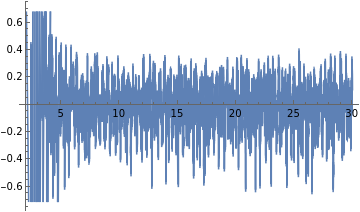}
    \caption{\centering Graph of $L_{f}[0,0](e^x)$, where $f = \frac{\li-\pi}{\li-\Ri}$}
   \label{eureka6}
\end{figure}

Figures \ref{eureka4}, \ref{eureka5}, and \ref{eureka6} plot the difference between the blue curve and the black curve in Figures \ref{eureka1}, \ref{eureka2}, and \ref{eureka3}, respectively; that is, they provide graphs of the functions 
$$L_{f}(e^x) = \frac{\log \frac{|\li(e^x)-\pi(e^x)|}{\li(e^x)-\Ri(e^x)}}{x},$$ $$L_{f}[0](e^x) = \frac{\log \frac{|\li(e^x)-\pi(e^x)|}{\li(e^x)-\Ri(e^x)}}{\log x},$$ and $$L_{f}[0,0](e^x) = \frac{\log \frac{|\li(e^x)-\pi(e^x)|}{\li(e^x)-\Ri(e^x)}}{\log \log x},$$ respectively, where $f = \frac{\li-\pi}{\li-\Ri}$ is the function $\li-\pi$ normalized by $\frac{1}{\li-\Ri}$.  By Proposition \ref{etaconjcor}, the Riemann hypothesis, Conjecture \ref{eurekaconjecture}, and Conjecture \ref{eurekaconjecture2} hold, respectively, if and only if the lim sup as $x \to \infty$ of these functions is $0$.

\begin{figure}[h!]
\includegraphics[width=70mm]{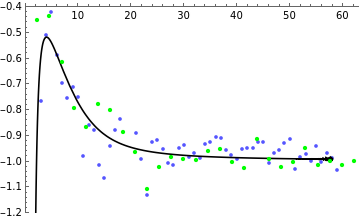}
\caption{\centering Plot of $L_{\li-\pi}[\frac{1}{2}](e^x)$ for $x = 3,4,5,\ldots, 59$ (in blue) and $x =1 \log 10,2\log 10,\ldots, 27 \log 10$ (in green), and graph of $L_{\li-\Ri}[\frac{1}{2}](e^x)$} 
   \label{dotfig1}
\end{figure}

\begin{figure}[h!]
\includegraphics[width=70mm]{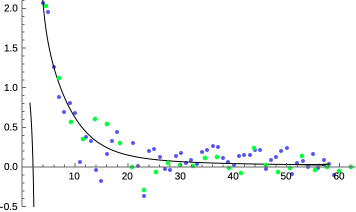}
\caption{\centering Plot of $L_{\li-\pi}[\frac{1}{2},\!-1](e^x)$ for $x = 3,4,5,\ldots, 59$ (in blue) and $x =2 \log 10,3\log 10,\ldots, 27 \log 10$ (in green), and graph of $L_{\li-\Ri}[\frac{1}{2},\!-1](e^x)$} 
   \label{dotfig2}
\end{figure}

One can provide further evidence for or against these conjectures  by estimating $\li(x)-\pi(x)$ for as many and as large values of $x$ as possible.  Such values that are already known are at our immediate disposal.   Table \ref{tab2} lists the values of $L_{\li-\pi}(x),$ $L_{\li-\pi}[\frac{1}{2}](x),$ and $L_{\li-\pi}[\frac{1}{2},-1](x)$ for $x = 10^n$, $n = 1,2,\ldots, 29$, based on known exact values of $\pi(10^n)$; Table \ref{tab3} lists their approximate values for $x = 10^n$, $n = 26,27,\ldots, 50$, based on the approximations of $\pi(10^n)$ provided in \cite{plan}; and Table \ref{tab4}  lists the their values for $x = e^n$, $n = 1,2,\ldots, 59$, based on known exact values of $\pi(e^n)$.   For comparison, we also provide the corresponding values for the function $ \li(x)-\Ri(x)$.   Figure  \ref{dotfig1} (resp., Figure \ref{dotfig2}) provides a plot of most of the values of $L_{\li-\pi}[\frac{1}{2}](x)$ (resp., $L_{\li-\pi}[\frac{1}{2},-1](x)$) given in Tables \ref{tab2} and \ref{tab4}, respectively, as a function of $e^x$, and compares them with the graph of $L_{\li-\Ri}[\frac{1}{2}](e^x)$ (resp., $L_{\li-\Ri}[\frac{1}{2},-1](e^x)$).

A more  informative approach is made possible by Proposition \ref{deltaprop} and the data in \cite{kul}.  Let
$$\Delta(x) = \left(\pi_0(x)- \Ri(x)  +\frac{1}{\log x} - \frac{1}{\pi} \arctan \frac{\pi}{\log x}\right) \frac{\log x}{\sqrt{x}}.$$
Since $\dege \Delta = \dege(\li-\pi)+(-\tfrac{1}{2},1,0,0,0\ldots)$,  Conjecture \ref{eurekaconjecture} is equivalent to $\dege \Delta  = (0,0,\Theta_2,\Theta_3,\ldots).$  Proposition \ref{deltaprop} implies the following.

\begin{corollary}\label{deltapropcor} 
Let $b > 1$, and let $n$ be a nonnegative integer.  Let $f(x)$ be any real function on $[1,\infty)$ such that $\degl f = \degl(\li-\pi)+(-\tfrac{1}{2},1,0,0,0\ldots)$ (e.g., $f(x) = \Delta(x)$).    Let $M_n$ denote the supremum  of $|f(x)|$ on the interval $[b^n,b^{n+1}]$, attained at or approached to the immediate left or right of some  $a_n \in [b^n,b^{n+1}]$.  One has the following.
\begin{enumerate}
\item $\displaystyle \sup_{x \in [b^n,b^{n+1}]} \frac{\log |f(x)|}{\log x} \leq \frac{\log M_n}{\log a_n}+ \frac{|\log M_n|}{n(n+1)\log b}.$
\item $\displaystyle \deg f = \limsup_{n \to \infty} \frac{\log M_n}{\log a_n}$.
\item $\displaystyle \sup_{x \in [b^n,b^{n+1}]} \frac{\log |f(x)|}{\log \log x} \leq \frac{\log M_n}{\log \log a_n}+ \frac{|\log M_n|}{\log(n\log b )}-\frac{|\log M_n|}{\log ((n+1)\log b)}.$ 
\item  $\displaystyle \degl_1 f = \limsup_{n \to \infty} \frac{\log M_n}{\log \log a_n}$ if the Riemann hypothesis holds.
\item $\displaystyle \sup_{x \in [b^n,b^{n+1}]} \frac{\log |f(x)|}{\log \log \log x} \leq \frac{\log M_n}{\log \log \log a_n}+ \frac{|\log M_n|}{\log \log(n\log b )}-\frac{|\log M_n|}{\log \log ((n+1)\log b)} $ if $n > 1$. 
\item  $\displaystyle \degl_1 f = \limsup_{n \to \infty} \frac{\log M_n}{\log \log a_n}$ if the Riemann hypothesis holds and $\degl_1(\li-\pi) = -1$.
\end{enumerate}
\end{corollary}

\begin{figure}[h!]
\includegraphics[width=70mm]{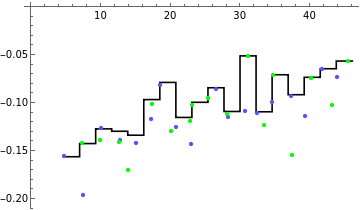}
\caption{\centering Graph of upper bound of $L_{\Delta}[0](e^x)$ on $[\log(10^{n}), \log (10^{n+1})]$ for $n = 2,3,\ldots,19$, and plot of $L_{\Delta}[0](e^x)$ for the suprema (in blue) and infima (in green) of $\Delta(x)$ on these intervals} 
\label{dotg1}
\end{figure}

\begin{figure}[h!]
\includegraphics[width=70mm]{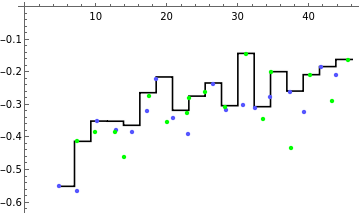}
\caption{\centering Graph of upper bound of $L_{\Delta}[0,0](e^x)$ on $[\log(10^{n}), \log (10^{n+1})]$ for $n = 2,3,\ldots,19$, and plot of $L_{\Delta}[0,0](e^x)$ for the suprema (in blue) and infima (in green) of $\Delta(x)$ on these intervals} 
\label{dotg2}
\end{figure}

The data in \cite{kul} provides, for $n = 0,1,2,\ldots,19$, the supremum $M_n^+$ and infimum $-M_n^-$ of $\Delta(k)$  on $[10^n,10^{n+1}]$, where the extrema are attained, say, at $b_n^+$ and $c_n^-$, so that $M_n = \max(M_n^+,M_n^-)$ is the supremum of $|\Delta(k)|$ on the given interval, attained at $a_n$, which is either $b_n^+$ or $c_n^-$.   This data is reproduced in Table \ref{tab5} (resp., Table \ref{tab6}), where we compute $L_\Delta(b_n^+) = \frac{\log M_n^+}{\log b_n}$,   $L_\Delta[0](b_n^+) = \frac{\log M_n^+}{\log  \log b_n}$,  $L_\Delta[0, 0](b_n^+) =  \frac{\log M_n^+}{\log  \log \log b_n}$ (resp., $L_\Delta(c_n^-) = \frac{\log M_n^-}{\log c_n}$,   $L_\Delta[0](c_n^-) = \frac{\log M_n^-}{\log  \log c_n}$,  $L_\Delta[0, 0](c_n^-) =  \frac{\log M_n^+}{\log  \log \log c_n}$).   Corollary \ref{deltapropcor}  implies that $\limsup_{n \to \infty} L_\Delta(a_n) = \deg \Delta$, while  $\limsup_{n \to \infty}L_\Delta[0](a_n)= \degl_1 \Delta$ if the Riemann hypothesis holds, and $\limsup_{n \to \infty} L_\Delta[0,0](a_n)= \degl_2 \Delta$ if the Riemann hypothesis holds and $\degl(\li -\pi)  = -1$.  Table \ref{tab7}  provides upper bounds of $L_\Delta(x) $, $L_\Delta[0](x)$, and  $L_\Delta[0, 0](x)$ on $[10^n,10^{n+1}]$.  For $n = 0,1,2,3,4$ these upper bounds were computed directly and are tight, while for $n = 5,6,7,\ldots, 19$ they were computed using Tables \ref{tab5} and \ref{tab6} and the bounds provided by Corollary \ref{deltapropcor}.  Figures \ref{dotg1} and \ref{dotg2} provide a graph of these upper bounds, along with a plot of  $L_\Delta[0](e^x)$ and $L_\Delta[0,0](e^x)$ for the supremum and infimum of $\Delta(x)$, on the intervals $[\log(10^n),\log(10^{n+1})]$ for $n = 2,3,4,\ldots,19$.  Since all of these upper bounds are negative, this might be construed as mild evidence for Conjectures \ref{eurekaconjecture} and \ref{eurekaconjecture2} that $\deg \Delta = \degl_1 \Delta = \degl_2 \Delta = 0$.

However, the data in Tables \ref{tab5}--\ref{tab7} is somewhat misleading: by Littlewood's theorem one has
$$\Delta(x)  =  \Omega_{\pm}(\log \log \log x) \ (x \to \infty),$$
and therefore
 $$\limsup_{x \to \infty} \Delta(x) = \limsup_{n \to\infty} M_n^+ = \infty$$ and  $$\liminf_{x \to \infty} \Delta(x) = -\limsup_{n \to\infty} M_n^- = -\infty,$$ so the suprema and infima observed in the three tables are much smaller than is expected for large $n$.  To compensate for this inadequacy of the data, one can settle for {\it approximations} of $\Delta(x)$ or of related functions of equal logarithmic degree.  In \cite{stoll}, Stoll and Demichel use a subtantial (but varying) number of nontrivial zeros of $\zeta(s)$ to approximate the function
$$V(x) = \frac{ \Ri(x)-\pi(x) }{\li(x)-\Ri(x)},$$
which satisfies $$V(x) =   \Omega_{\pm}(\log \log \log x) \ (x \to \infty)$$
(as does $\Delta(x)$) and
 $$\dege V =  \dege \Delta =  \dege(\li-\pi) + (-\tfrac{1}{2},1,0,0,0,\ldots).$$
Note that the red graph in Figure \ref{LiRi} is a graph of the function $V(e^x)$ on the interval $[2,30]$, and, as observed in Section 1.3, one has $V(x) < -1$ if and only if $\li(x)<\pi(x)$, for any $x \geq \mu$. For $x> 10^{20}$, Stoll and Demichel used the following sum of cosines  to approximate $V(x)$, using the first $N \geq 10^5$ zeros  $\rho_k = \frac{1}{2}+i\gamma_k$ of $\zeta(s)$ with positive imaginary part $\gamma_k$ and real part $\frac{1}{2}$:
$$V(x) \approx -\Delta(x) \approx  \frac{ \sum_\rho \Ei(\rho \log x)}{\sqrt{x}/\log x}\approx  \sum_{k = 1}^N 2\operatorname{Re}\left( \frac{1}{\rho} e^{i\gamma_k\log x }\right) =  \sum_{k = 1}^N \frac{2 \cos(\gamma_k \log x - \operatorname{Arg} \rho_k)}{|\rho_k|}.$$
To reduce otherwise unweildy computational time, they chose $N$ to vary depending on their degree of interest in determining a more accurate estimation of the location and value for various ``significant'' peak values of $V(x)$. 
In private communication with the author,  Stoll defined ``significant'' as ``at least exceeding 0.90 times the previous $+V(x)$ or $-V(x)$ maximum/minimum value.''

  One observation to be made is that the positive values of $V(x)$ in Table \ref{tab8} increase as $x$ increases, while the negative values of $V(x)$ increase in absolute value as $x$ increases with but three exceptions, at $x = 1.592776 \cdot 10^{1165}$ and in last two rows of the table.  However,  Stoll informed the author \cite{stolll} that the value $x =1.592776 \cdot 10^{1165}$ in the table was included for historical reasons, as it was first discovered by R.\ S.\ Lehman in 1966 in the form $e^{2682.9768}$ \cite{lehm} (where Lehman showed that there exist at least $10^{500}$ successive integers $x$ between $1.53 \cdot 10^{1165}$ and $1.65 \cdot 10^{1165}$ such that $\li(x)-\pi(x) < 0$), and the last two were included because their locations  were so much larger than prior locations and they demonstrated that the range investigated indeed climbed all the way to $10^{10^{13}}$.  With those three rows removed, the positive (resp. negative) values of $V(x)$ are record high values of $V(x)$ for the size of $x$, as far as is currently known.  Stoll has expressed in private communication the much stronger claim that, for each of these (approximate) values of $x$ there is no larger positive (resp., larger negative) value of $V(x)$ for any smaller $x$, but this claim is difficult to verify.  If true, it would imply, for example, that the true value of Skewes' number is close to $1.397162914 \cdot 10^{316}$.

\begin{remark}[The search for Skewes' number]
Skewes' number is equal to infimum of all $x \geq \mu$ such that $V(x) < -1$.  One of the purposes of Stoll's and Demichel's  $V(x)$ data was to try to find an upper bound for  Skewes' number that is smaller than the smallest verified upper bound, which, currently, is $e^{727.951336108} \leq  1.3971671494 \cdot 10^{316}$.  They carried out extensive searches at the two of the ``most promising'' regions indicated by Bays' and Hudson’s extensive sweep up to $10^{400}$ \cite{bays}.  The first was in the vicinity of  $x =1.258585341\cdot 10^{179}$, and the second in the vicinity of  $x =1.3379303302214\cdot 10^{190}$.  Even using the first $5 \cdot 10^9$ nontrivial zeros of $\zeta(s)$, they found no value of $V(x)$ less than than $-0.969$.  In both cases mentioned,  Demichel was able to show that including even more zeros of $\zeta(s)$ would not make the values approach $-1.0$.  It was their conclusion that Skewes' number was likely to be in the vicinity of  $1.397162914025 \cdot 10^{316} \approx e^{727.9513330766}$.  However, they did find a value of $V(x)$ almost equal to $+1.0$ in the vicinity of $x = 1.90987566088001\cdot 10^{215}$, with $V(x) \approx 0.999$ using $5\cdot 10^9$ zeros of $\zeta(s)$, which indicates a great shortage of primes in that region, not the great excess of primes needed for $\pi(x)$ to overtake $\li(x)$ as it does in the vicinity of Skewes’ number.
\end{remark}

Using the data in \cite[Table 3]{stoll} discussed above, in Table \ref{tab8} we provide the corresponding approximations of  $L_V(x)$, $L_V[0](x)$, $L_V[0,0](x)$, and $L_V[0,0,0](x)$ at these peak $V(x)$ points.  Once $x$ surpasses the proposed Skewes' number $1.397162914 \cdot 10^{316}$, these four values become positive at the given extreme points.  Also, the values $L_V(x)$  at these points increase monotonically to the maximum value $0.0000568$ at $6..5769904020884 \cdot 10^{370}$, and thereafter they decrease monotonically to $2.42 \cdot 10^{-14}$.  This provides some evidence for the Riemann hypothesis that $\deg V = 0$.  Moreover, the values of $L_V[0](x)$ at these points do not surpass $0.0167$, which might be construed as mild evidence for Conjecture \ref{eurekaconjecture} that $\degl_1 V = 1+ \degl_1(\li -\pi) = 0$.  However, the values of $L_V[0,0](x)$ at these points increase nearly monotonically from $-0.111$ to the seventh to last value of $0.121$, where the last seven values are $0.121$, $0.127$, $0.123$, $0.122$, $0.130$, $0.122$, and $0.121$. This suggests that the $V(x)$ peak data does not provide evidence that $\degl_2 V = \degl_2(\li -\pi) = 0$.

\begin{figure}[h!]
\includegraphics[width=70mm]{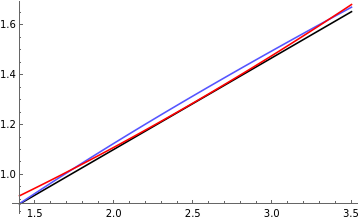}
    \caption{\centering Graph of $\frac{1}{e}(1+\log \log \log x)$ (in black) and upper bounds  $\frac{12}{19}(\frac{1}{6}+\log \log \log x)^{3/4}$ (in blue) and $\frac{1}{38}(\frac{9}{2}+\log \log \log x)^2$ (in red) on $[e^{e^{e^{1.4}}},e^{e^{e^{3.5}}}]$, as functions of $\log \log \log x$}
\label{stollconj}
\end{figure}

In  \cite[(11), p.\ 2384]{stoll}, Stoll and Demichel conjecture that
 \begin{align}\label{dem1}
|V(x)| < \frac{1}{e}\left(1+\log \log \log x\right)+ \frac{4}{3(\log x)^2}, \quad \forall x \geq 3,
\end{align}
and, since ``the third term is only required for $x \leq 10$ and is dropped for simplicity in later analysis,'' also that
 \begin{align*}
|V(x)| < \frac{1}{e}\left(1+\log \log \log x\right),  \quad \forall x >10.
\end{align*}
Either of these conjectures implies (\ref{dem}), which in turn implies the ``best case'' scenario $\degl (\li-\pi) = (\tfrac{1}{2},-1,0,1,0,0,0,\ldots)$.  However, the rather bold conjecture that $\degl (\li-\pi) = (\tfrac{1}{2},-1,0,1,0,0,0,\ldots)$ is incompatible with the conjectures (\ref{MMCa}) and (\ref{MMCb}) of Monach and Montgomery (and also with Montgomery's conjecture (\ref{MMC3})).  Stoll and Demichel  cite  \cite[Figure 1]{stoll} as evidence for their conjectures,  which is a plot of the $V(x)$ peak values in \cite[Table 3]{stoll} (which are reproduced in Table \ref{tab8})  against $\log \log \log x$ and shows an approximate linear relationship with $ \log \log \log x $ at those particular points, which satisfy $1.4535< \log \log \log x < 3.4168$.  They also note that ``the slope of  $1/e$ was derived empirically.''  However, several nonlinear functions of $\log \log \log x$ appear to be approximately linear  in the range  considered while also bounding $\frac{1}{e}(1+\log \log \log x)$ from above, e.g., the functions $\frac{12}{19}(\frac{1}{6}+\log \log \log x)^{3/4}$ and $\frac{1}{38}(\frac{9}{2}+\log \log \log x)^2$ pictured in Figure \ref{stollconj}; but replacing  $\frac{1}{e}(1+\log \log \log x)$ in their conjecture with the first function    would contradict Littlewood's theorem, while replacing it with the second function would make their conjecture compatible with the conjectures of Monach and Montgomery.  One can also add to these adjusted upper bounds  a  corrective term of the form $\frac{C}{(\log x)^d}$ so that the resulting sum would also bound the known peak $V(x)$ values for $x < 10^{20}$ as listed in \cite[Table 2]{stoll}.   However, none of this could be regarded as evidence that $|V(x)|$ is so bounded for all $x$.  Given the approximations  of $L_V[0,0](x)$ and $L_V[0,0,0](x)$ at the approximate peak values of $V(x)$ listed in Table \ref{tab8}, it seems far more likely that the true value of  $\degl_3 V = \degl_3 (\li-\pi)$ is not encroached upon until $\log \log \log x$ is substantially greater than $3.4168$.

 Generally, conjectures concerning $\li(x)-\pi(x)$ as it relates to  $\log \log \log x$ that are based on numerical and graphical data alone should be viewed with healthy skepticism.    As we have argued, the peak $V(x)$ data does not provide evidence for the conjecture that $\degl_2(\li-\pi) = 0$, and this makes it very difficult for us to accept the claim that the same data provides evidence for the conjecture that  $\degl (\li-\pi) = (\tfrac{1}{2},-1,0,1,0,0,0,\ldots)$.  At least for the forseeable future, competing conjectural values for $\degl_3 (\li-\pi)$ likely cannot be adjudicated by numerical computations alone.  In particular, Stoll's and Demichel's data cannot be regarded as evidence for Conjecture (\ref{dem1}), or even as evidence against the conjectures of Monach and Montgomery.  The statistic $L_f[d_0,d_1,\ldots,d_{N-1}]$  (specifically, $L_V[0,0]$ and $L_V[0,0,0]$) here serves a warning against the use of  curve-fitting to predict the long-term behavior of functions related to the prime numbers.

Because the data and graphs discussed in this section provide modest support for the conjecture that $\degl_1(\li-\pi) = -1$ but  less support for $\degl_2 (\li -\pi) = 0$, we assert Conjecture \ref{eurekaconjecture2} to a lesser degree of confidence than we do Conjecture \ref{eurekaconjecture}.    Nevertheless, both conjectures have theoretical appeal for the reasons described in Sections 9.1 and 9.2,  both are implied by the conjectures of Monach and Montgomery and by the conjectures of Stoll and Demichel, and no existing conjecture of which we are aware that is consistent with the Riemann hypothesis implies that either conjecture is false.

\section{Conjectures concerning the Mertens function}

Various well-known conjectures regarding the order of growth of the Mertens function $M(x) = \sum_{n \leq x} \mu(n)$ also have bearing on the Riemann hypothesis, since $\deg M = \deg(\li-\pi)$.   In \cite{ng} Ng provides  some theoretical evidence for an unpublished conjecture by Gonek that 
\begin{align}\label{gonek}
0 < \limsup_{x \to \infty} \frac{M(x)}{\sqrt{x}\, (\log \log \log x)^{5/4}} = \limsup_{x \to \infty}  \frac{-M(x)}{\sqrt{x} \, (\log \log \log x)^{5/4}} < \infty.
\end{align}
Previously, Good and Churchhouse \cite{good} and L\'evy \cite{lev} had made conjectures implying that
\begin{align}\label{good}
0 < \limsup_{x \to \infty} \frac{|M(x)|}{\sqrt{x \log \log x}} < \infty.
\end{align}

\begin{proposition}
Conjecture (\ref{gonek}) of Gonek and Ng implies that $$\degl M = (\tfrac{1}{2},0,0,\tfrac{5}{4},0,0,0,\ldots).$$ Conjecture  (\ref{good}) of Good and Churchhouse and L\'evy implies that $$\degl M = (\tfrac{1}{2}, 0, \tfrac{1}{2}, 0, 0, 0, \ldots).$$
\end{proposition}

Ng provides an explanation for the exponent $\frac{5}{4}$ on $\log \log \log x$ in his and Gonek's conjecture (\ref{gonek}) and the exponent $2$ on $\log \log \log x$ in Montgomery's conjecture (\ref{MMC2}): ``The difference between these cases is due directly to the different discrete moments of $\sum_{\gamma \leq T} \frac{1}{|\rho|} \asymp (\log T)^2$  and $\sum_{\gamma \leq T} \frac{1}{|\rho \zeta'(\rho)|} \asymp (\log T)^{\frac{5}{4}}$, where the second inequality is currently conjectural'' \cite{ng}.   Thus, Ng implies that his conjecture more generally is that $\degl M = (\tfrac{1}{2},0,0,d_3, 0,0,0,\ldots),$ where
$d_3 = \degl_1 \sum_{0<\operatorname{Im} \rho \leq T} \frac{1}{|\rho \zeta'(\rho)|}.$
Ng also proves  in \cite{ng} that, if the Riemann hypothesis holds, if all of the nontrivial zeros $\rho$ of $\zeta(s)$ are simple, and if 
\begin{align}\label{rm1}
\sum_{0<\operatorname{Im} \rho \leq T} \frac{1}{|\zeta'(\rho)|^2} = O(T) \ (T \to \infty),
\end{align}
then one has $M(x) = O\left(\sqrt{x}\, (\log x)^{3/2}\right) \ (x \to \infty)$
and therefore $\degl_1 M \leq \frac{3}{2}$; likewise, if also 
\begin{align}\label{rm2}
\sum_{0<\operatorname{Im} \rho \leq T} \frac{1}{|\zeta'(\rho)|} = O(T(\log x)^{1/4}) \ (T \to \infty),
\end{align}  then $M(x) = O\left(\sqrt{x}\, (\log x)^{5/4}\right) \ (x \to \infty)$
and therefore $\degl_1 M \leq \frac{5}{4}$.  Moreover, both (\ref{rm1}) and (\ref{rm2}) follow from \cite[Conjecture 2]{HKN}, which is motivated by a random matrix theory model of the asymptotics of the discrete moments of $\zeta'(s)$.

In this section we provide some evidence for Conjecture \ref{Mconjecture}.

By Theorem \ref{multzero}, we have the following.

\begin{proposition}
If Conjecture \ref{Mconjecture} holds, then the Riemann hypothesis holds, all of the zeros of the Riemann zeta function are simple, and $\degl_2 M \geq 0$.
\end{proposition}

\begin{figure}[ht!]
\includegraphics[width=70mm]{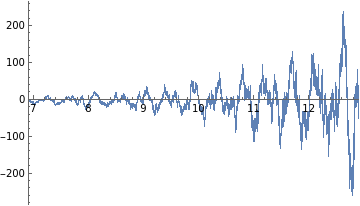}
\caption{\centering Plot of $M(e^x)$ for integer values of  $e^x \in [1000,400000]$} 
\label{Mfig1}
\end{figure}

\begin{figure}[ht!]
\includegraphics[width=70mm]{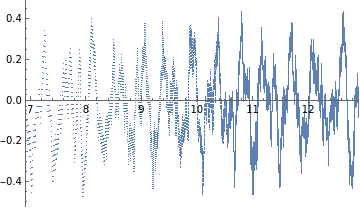}
\caption{\centering Plot of $M(e^x)e^{-x/2}$ for integer values of  $e^x \in [1000,400000]$} 
\label{Mfig2}
\end{figure}

\begin{figure}[ht!]
\includegraphics[width=70mm]{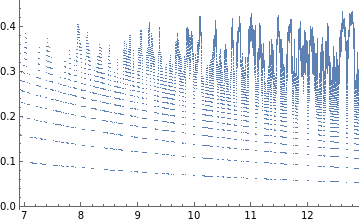}
\caption{\centering Plot of $L_M(e^x)$ for integer values of  $e^x \in [1000,400000]$} 
\label{Mfig3}
\end{figure}

\begin{figure}[ht!]
\includegraphics[width=70mm]{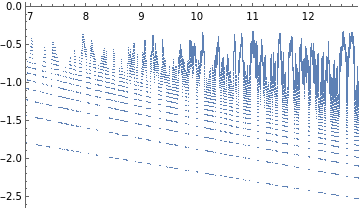}
\caption{\centering Plot of $L_M[\frac{1}{2}](e^x)$ for integer values of  $e^x \in [1000,400000]$} 
\label{Mfig4}
\end{figure}

\begin{figure}[ht!]
\includegraphics[width=70mm]{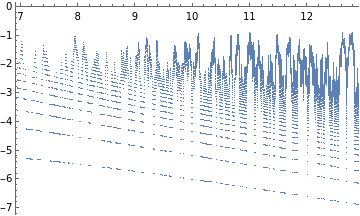}
\caption{\centering Plot of $L_M[\frac{1}{2},0](e^x)$ for integer values of  $e^x \in [1000,400000]$} 
\label{Mfig5}
\end{figure}

Conjecture \ref{Mconjecture}  is implied by the conjecture (\ref{gonek}) of Gonek and Ng and by the conjectures of Good, Churchhouse, and L\'evy and is also modestly supported, with some caveats, by numerical data.  Recall from Theorem \ref{multzero} that the Riemann hypothesis is equivalent to $\deg M = \frac{1}{2}$ and implies $\degl_1 M \geq 0$.  Table \ref{tab9}  in Section 14.7 computes  $L_{M}(x) = \frac{\log |M(x)|}{\log x}$, $L_{M}[\frac{1}{2}](x) = \frac{\log |M(x)|-\frac{1}{2}\log x}{\log \log x}$, and $L_{M}[\frac{1}{2},0](x) =  \frac{\log |M(x)|-\frac{1}{2}\log x}{\log \log \log x}$ for the peak values of $M(x)$  provided in \cite{hur} \cite{kuz}.  The six largest values for each statistic are listed in bold.   
Figures \ref{Mfig1}--\ref{Mfig5} provide, respectively, a plot of the functions $M(e^x)$, $M(e^x)e^{-x/2}$, $L_M(e^x)$, $L_M[\frac{1}{2}](e^x)$, and $L_M[\frac{1}{2},0](e^x)$ for integer values of $e^x$ in $[1000,400000]$.

Although Table \ref{tab9} provides modest numerical support for Conjecture \ref{Mconjecture}, all known data on $M(x)$ can be deceptive since the smallest counterexample to inequality  $|M(x)| < \sqrt{x}$ (known as the {\bf Mertens conjecture})\index{Mertens conjecture} is currently only known to be somewhere between $10^{16} \approx e^{36.84}$ and $e^{1.59 \cdot 10^{40}}$ \cite{hur} \cite{kot2}, yet the inequality is known to fail for arbitrarily large values of $x$ since $\limsup_{x \to \infty} \frac{M(x)}{\sqrt{x}} > 1.826054$ and $\liminf_{x \to \infty} \frac{M(x)}{\sqrt{x}} < -1.837625$ \cite{hur} \cite{odl}.      Since for any $x > e$ the  conditions $|M(x)| < \sqrt{x}$,  $L_{M}(x) < \frac{1}{2}$, and  $L_{M}[\frac{1}{2}](x) < 0$ are equivalent, all known data on $M(x)$  naively supports Conjecture \ref{Mconjecture}.  Data on $M(x)$ that would support the negation of Conjecture \ref{Mconjecture} is very difficult to obtain as it would require at the very least estimates of $\log \frac{|M(x)|}{\sqrt{x}}$ and $\log  \log x$ for a sampling of values of $x$ with $|M(x)| > \sqrt{x}$.    Currently, the largest known explicit value of $\frac{|M(x)|}{\sqrt{x}}$ is $0.58576768\ldots$, which occurs at $x = 11609864264058592345$ \cite{kuz}, and for this $x$ the values of $L_{M}(x)$, $L_{M}[\frac{1}{2}](x)$, and $L_{M}[\frac{1}{2},0](x)$ are the largest that we can extract from values of $M(x)$ in the literature.  No one, as far as we can surmise, has undertaken the project of bounding $\frac{|M(x)|}{\sqrt{x}}$ on consecutive intervals.    The  closest result to this effect we can find in the literature is that $$-0.5247975\ldots \leq \frac{M(x)}{\sqrt{x}} \leq 0.5705908\ldots, \quad \forall x\in [201, 10^{16}],$$
where the  positive and negative extreme values of $\frac{M(x)}{\sqrt{x}}$ are assumed at $x = 7766842813$ (where $M(x) = 50286$) and $x = 71578936427177$, respectively \cite[Section 6.1]{hur}.   From this result and Proposition  \ref{deltaprop} applied to $f(x) = \frac{M(x)}{\sqrt{x}}$, $b = 10^8$, and $n = 1$, it follows that
$$L_M(x) < 0.490592, \ \ L_M[\tfrac{1}{2},0](x) < -0.142501, \ \  L_M[\tfrac{1}{2},0](x) < -0.405043, \quad \forall x \in [10^{8},10^{16}],$$
where  also $L_M(351246529829131) = 0.480335\ldots$, $L_M[\tfrac{1}{2},0](7766842813) = -0.179513\ldots$, and $L_M[\tfrac{1}{2},0](7766842813) = -0.492341\ldots$ are the largest values  observed in Table  \ref{tab9} on $[10^8,10^{16}]$.  Table \ref{tab10} provides the supremum of $\frac{M(x)}{\sqrt{x}}$, $L_M(x)$, $L_M[\frac{1}{2}](x)$, and  $L_M[\frac{1}{2},0](x)$ on $[10^{n},10^{n+1}]$ for $n = 0,1,2,\ldots,7$, which we computed directly using Mathematica.  This data allows us to conclude that the bounds above for $x \in [10^{8},10^{16}]$ hold more generally for all $x \in (e^e,10^{16}]$.  Notably, the data also shows that the supremum of $\frac{M(x)}{\sqrt{x}}$  on $[10^n,10^{n+1}]$ decreases monotonically from $n = 0$ to $n = 6$, while the supremum of $L_M(x)$ on $[10^n,10^{n+1}]$ increases monotonically from 
$n = 2$ to $n = 6$.  This reinforces the fact that the statistics $\frac{M(x)}{\sqrt{x}}$ and $L_M(x)$  both reveal important information that the other  may not.

Next, we examine \cite[Tables 2 and 3]{kuz}, which approximate, for $x \in [10^{14}, e^{10^{15}}]$, the {\bf IL$q^+$ values}  (resp., {\bf IL$q^-$ values})\index{IL$q^\pm$ values}  of $q(x) = \frac{M(x)}{\sqrt{x}}$, which, roughly, are the values of $q(k)$ that are largest positive (resp., largest negative) among the positive integers $k \geq 10^4$ that have the same or fewer digits than $k$  \cite[Definition 3.1]{kot}.  (``IL'' stands for ``increasingly large.'')  There are exactly seven IL$q$ (i.e., IL$q^+$/IL$q^-$) values of $q(x)$ for $x \in [10^4,10^{14}]$, given precisely in \cite[p.\ 475]{kot}.   Assuming the Riemann hypothesis and that all zeros of $\zeta(s)$ are simple, one has
$${q(x)} = \sum_{k = 1}^\infty 2\operatorname{Re}\left(\frac{x^{\rho_k-1/2}}{\rho_k\zeta'(\rho_k)}\right) + O(x^{-1/2}) = \sum_{k = 1}^\infty 2|\beta_k| \cos (\gamma_k \log x + \arg \beta_k) + O(x^{-1/2}) \ (x \to \infty),$$
where $\gamma_k$ is the imaginary part of the $k$th nontrivial zero $\rho_k$ of $\zeta(s)$ with positive imaginary part,  where $\beta_k = \frac{1}{\rho_k \zeta'(\rho_k)}$, and where the terms are grouped appropriately \cite[(3)]{hur}.  Kotnik and van de Lune and Kuznetsov obtained the  approximations of the IL$q$ values of $q(x)$ listed in   \cite[Table 4]{kot} \cite[Tables 2 and 3]{kuz} by truncating the series  above after the first $10^6$ terms.  More precisely, they computed the IL$q_{10^6}$ values of $q_{10^6}(x)$, where 
$$q_K(x) =\sum_{k = 1}^K 2\operatorname{Re}\left(\frac{x^{\rho_k-1/2}}{\rho_k\zeta'(\rho_k)}\right) =  \sum_{k = 1}^K 2|\beta_k| \cos (\gamma_k \log x + \arg \beta_k).$$
Kotnik and van de Lune provide some evidence that the data in their table approximates the IL$q$ values for $x$ in the range $10^4 \leq x \leq 10^{10^{10}}$, thus yielding some values of $x$ where $q(x)$ is likely to be large positive or large negative for the size of $x$.  Table \ref{tab11} computes $L_f(x) \approx L_q(x),$ $L_f[0](x) \approx L_q[0](x),$ and $L_f[0,0](x) \approx L_q[0,0](x),$  where $f(x) = q_{10^6}(x) \approx q(x)$, for the IL$q_{10^6}$ values of $q_{10^6}(x)$ for $x$ in the range  $10^{14} \leq x \leq e^{10^{15}}$ in \cite[Tables 2 and 3]{kuz}, which is an extension of \cite[Table 4]{kot} to greater precision and larger domain. 
(Note that the explicit formula for $M(x)$ can be generalized to account for the possibility that the Riemann hypothesis is false or that not all the zeros of $\zeta(s)$ are simple \cite[Theorem]{bartz}, and in either of those instances the data in Table \ref{tab11} still has some significance.)  Since $|q_{10^6}(x)| < 1$ for all IL$q_{10^6}$ values in Table \ref{tab11}, one has $L_{q_{10^6}}(x) < 0$, $L_{q_{10^6}}[0](x) < 0$, $L_{q_{10^6}}[0,0](x) < 0$, and $L_{q_{10^6}}[0,0,0](x) < 0$ for all of these values.  Moreover, the values of $L_{q_{10^6}}[0](x)$ are  $L_{q_{10^6}}[0,0](x)$ in the two tables are quite small in absolute value, while $\log \log x$ reaches the value of $33.857$.  Thus the data provides evidence, albeit limited, for the conjecture $\degl_1 M = 0$ (and perhaps also  for the conjecture $\degl_2 M = 0$, although this is debated below).

Based on the IL$q_{10^6}$ values and the apparent linear relationship with $\sqrt{\log \log \log x}$ shown in \cite[Figure 4]{kot}, Kotnik and van de Lune claim that it is not clear how to reconcile the conjectures of Good, Churchhouse, and L\'evy with the data and suggest that the exponent $\frac{5}{4}$ of $\log \log \log x$ in the conjecture of Gonek and Ng may be an overestimate and appears to be closer or equal to $\frac{1}{2}$.  They make the more conservative conjecture \cite[Conjecture 5.1]{kot} that 
\begin{align}\label{kotc}
M(x) = \Omega_{\pm}(\sqrt{x \log \log \log x}) \ (x \to \infty),
\end{align}
which is given further graphical support by Kotnik and te Riele in \cite[Figure 3]{kot2}.   However, since  $\sqrt{\log \log \log x} < 2.36$  and $L_{q_{10^6}}[0,0,0](x) < 0$ for all values of $x$ considered in  \cite{kot} \cite{kot2}, it is not clear how their data  and graphs can support such claims involving $\log \log \log x$.  Indeed, their evidence at the very least suffers from the same problem that Stoll's and Demichel's evidence for the inequality (\ref{dem1}) has, namely, that the graphs of many other functions of a given power of $\log \log \log x$---even some involving $\log \log x$---appear very similar over the domain considered.  We can provide a more fine-tuned analysis by approximating the  IL$q_{10^6}$ values with a function more general than simply a constant times a power of $\log \log \log x$.

\begin{figure}[h!]
\includegraphics[width=70mm]{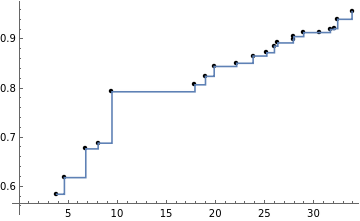}
    \caption{\centering Plot of $|q_{10^6}(e^{e^x})|$ for the IL$|q_{10^6}|$ values $e^{e^x}$, and graph of  implied approximation for $\sup_{t \in [10^2,x] }|q_{10^6}(e^{e^t})|$}
\label{Tq}
\end{figure}

\begin{figure}[h!]
\includegraphics[width=70mm]{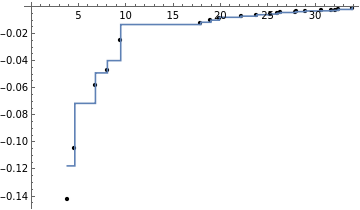}
    \caption{\centering Plot of $L_{q_{10^6}}[0](e^{e^x})$ for the IL$|q_{10^6}|$ values $e^{e^x}$, and graph of  implied upper bound for $\sup_{t \in [10^2,x] }L_{q_{10^6}}[0](e^{e^x})$}
\label{Tq0}
\end{figure}

\begin{figure}[h!]
\includegraphics[width=70mm]{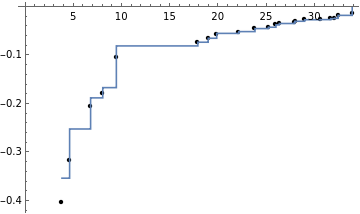}
    \caption{\centering Plot of $L_{q_{10^6}}[0,0](e^{e^x})$ for the IL$|q_{10^6}|$ values $e^{e^x}$, and graph of implied  upper bound for $\sup_{t \in [10^2,x] }L_{q_{10^6}}[0,0](e^{e^x})$}
\label{Tq00}
\end{figure}

First, we define the {\bf IL$|q|$ values}\index{IL$\vert q \vert$ values}  of $|q(x)|$ to be the values of $|q(x)|$ that are  the largest in absolute value among the positive integers $k \geq 10^4$ that have the same or fewer digits than $k$, and we define the {\bf IL$|q_{10^6}|$ values}  of $|q_{10^6}(x)|$ similarly.  One obtains a table of IL$|q_{10^6}|$ values by pruning Table \ref{tab11} of 
the IL$q_{10^6}$ values appropriately.  The result is Table \ref{tab12}, the values of which are plotted in Figure \ref{Tq}.   Assuming  that Table \ref{tab12} is not missing any IL$q_{10^6}$ values in the domain considered, Figure \ref{Tq} provides an  plot of $|q_{10^6}(e^{e^x})|$ for the IL$|q_{10^6}|$ values $e^{e^x}$, along with an approximate graph of the step function $\sup_{t \in [10^2,x]} |q_{10^6}(e^{e^x})|$.  (Any missing  IL$|q_{10^6}|$ values would create intermediate jumps in the graph.)  Presumably, then, the graph is a rough approximation of the graph of the step function $\sup_{t \in [10^2,x]} |q(e^{e^x})|$.   Figures \ref{Tq0} and \ref{Tq00} provide, repsectively, a  plot of $L_{q_{10^6}}[0](e^{e^x})$ and $L_{q_{10^6}}[0,0](e^{e^x})$ for the IL$|q_{10^6}|$ values $e^{e^x}$, along with implied approximate upper bounds for the step functions $\sup_{t \in [10^2,x]} L_{q_{10^6}}[0](e^{e^x})|$ and $\sup_{t \in [10^2,x]} L_{q_{10^6}}[0,0](e^{e^x})|$.

\begin{figure}[h!]
\includegraphics[width=70mm]{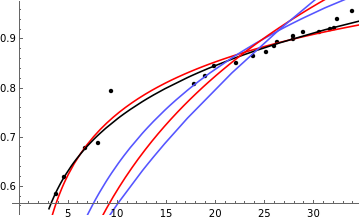}
    \caption{\centering Plot of $|q_{10^6}(e^{e^x})|$ for the IL$|q_{10^6}|$ values $e^{e^x}$, and graph of best fit curve $Ax^{a_2}(\log x)^{a_3}$, with no restrictions (in black);  with $a_2 = 0$, $a_3 = \frac{1}{2}$ (in red); with $a_2 = 0$, $a_3 = 1$ (in blue); with $a_2 = 0$, $a_3 = \frac{5}{4}$ (in red);  and with $a_2 = \frac{1}{2}$, $a_3 = 0$ (in blue)}
\label{BestfitM0}
\end{figure}

To analyze the IL$q_{10^6}$ data more precisely, we use  NonlinearModelFit in Mathematica to find some best fit curves to the absolute IL$|q_{10^6}|$ values $|q_{10^{6}}(e^{e^x})|$, expressed as a function of $e^{e^x}$.   We use the  IL$|q_{10^6}|$ values rather than the  IL$q_{10^6}$ values because the latter require two curves, not one.   As a disclaimer, we note that all parameters for the best fit curves discussed below were computed with NonlinearModelFit and may not be actual globally optimizing parameters.  
Let $$f(x) =Ae^{a_1 x}x^{a_2}( \log x)^{a_3}.$$
 Using NonlinearModelFit, we find that the best fit curve of the form $f(x)$ above to  $|q_{10^{6}}(e^{e^x})|$ at the IL$|q_{10^6}|$ values $e^{e^x}$ has
$$a_1 = 0.0116957, \quad a_2 = -0.398891,\quad a_3 = 1.04372, \quad A = 0.69784,$$
with an adjusted $R^2$ of 0.999601.   Since $a_1$ is small, $a_1 = 0$ is close to optimal, so we use  NonlinearModelFit again to find that the best fit curve of the form above with  $a_1 = 0$ to the given data has 
$$a_1 = 0, \quad a_2 = 0.132648, \quad a_3 =0.170467,  \quad A= 0.472356$$
with an adjusted $R^2$ of 0.999516.  Thus $a_1 = 0$ is close to optimal.   Assuming the Riemann hypothesis $\deg q = 0$, this provides some evidence for Conjecture \ref{Mconjecture} that $\degl_1 q = \degl_1 M= 0$. 
However, the data does not really suggest a value for $\degl_2 q$.       Figure \ref{BestfitM0} shows a plot of $|q_{10^6}(e^{e^x})|$ for the IL$|q_{10^6}|$ values $e^{e^x}$, along with graphs of the best fit curve $Ax^{a_2}(\log x)^{a_3}$, with no restrictions (in black);  with $a_2 = 0$, $a_3 = \frac{1}{2}$ (in red); with  $a_2 = 0$, $a_3 = 1$ (in blue);  with $a_2 = 0$, $a_3 = \frac{5}{4}$ (in red);  and with $a_2 = \frac{1}{2}$, $a_3 = 0$ (in blue).  Since among these the red curve $\sqrt{\log x}$ with $a_2 = 0$, $a_3 = \frac{1}{2}$ is closest to the best fitting curve in black, this appears to give some further rationale---but not evidence, in our view---for the Kotnik and van de Lune conjecture (\ref{kotc}) \cite[Conjecture 5.1]{kot} and their claim that $\degl_3 M$ appears to be closer to $\frac{1}{2}$ than $\frac{5}{4}$.

\begin{figure}[h!]
\includegraphics[width=70mm]{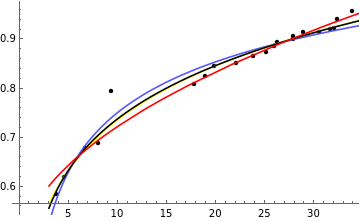}
    \caption{\centering Plot of $|q_{10^6}(e^{e^x})|$ for the IL$|q_{10^6}|$ values  $e^{e^x}$,  and graph of best fit curve $B+Ax^{a_2}(\log x)^{a_3}$, with no constraints (in black); with $a_2 = 0$, $a_3 = \frac{5}{4}$ (in yellow); with $a_2 = 0$, $a_3 = \frac{1}{2}$ (in blue); and with $a_2 = \frac{1}{2}$, $a_3 = 0$  (in red)}
\label{BestfitM}
\end{figure}

 In  order to probe $\degl_2 q$ and $\degl_3 q$ more closely, henceforth we  assume $\deg q = \degl_1 q = 0$ and set $a_1 = 0$, and we introduce an additive constant $B$ so that $$f(x) =B+Ax^{a_2}(\log x)^{a_3}.$$  The  extra parameter $B$ serves to enhance our curve-fitting capability, but its introduction, admittedly, is somewhat {\it ad hoc}. The parameters for the best fit curves $f(x)$ to $|q_{10^{6}}(e^{e^x})|$ at the IL$|q_{10^6}|$ values  $e^{e^x}$ under various constraints are provided in  Table \ref{tablebestfit} below.    Subject to no constraints, the best fit curve has
$$a_2 = 0.140224, \quad a_3 = 0.28531, \quad	A  = 0.329144, \quad B = 0.161942,$$
with an adjusted $R^2$ of 0.999488.   Subject to the constraint $a_2 = 0$, the best fit curve has
$$a_2 = 0, \quad  a_3  = 1.11312,\quad A = 0.126778, \quad B = 0.418209,$$
with an adjusted $R^2$ of 0.999502.  In fact, the two curves above are remarkably close to one another, and thus the latter curve is close to optimal.  In the latter curve, $a_3$ is somewhat close to the value $\frac{5}{4}$ for $\degl_3 q$ conjectured by Gonek and Ng, but  perhaps the large size of $B$ relative to $f(x)$ makes this less substantial evidence for their conjecture than  might appear at first sight.  Subject to the constraints $a_3 = 0$ and $B \geq 0$, the best fit curve has
$$a_2 = 0.202478, \quad  a_3  = 0, \quad A = 0.459861, \quad B = 0,$$
with  an adjusted $R^2$ of 0.999491, where now $a_2$ is  close to $\frac{1}{5}$, not $\frac{1}{2}$ as in the Good, Churchhouse and L\'evy conjectures. As shown in Figure \ref{BestfitM}, the best fit curves with $a_2 = 0$, $a_3 = \frac{5}{4}$ (in yellow), $a_2 = 0$, $a_3 = \frac{1}{2}$ (in blue), and $a_2 = \frac{1}{2}$, $a_3 = 0$ (in red) all fit the data well, although the  Gonek and Ng curve (in yellow) comes closest to the best fit curve of the form $f(x)$ with no constraints (in black), and in fact the black and yellow curves appear to be almost indistinguishable.  Thus, with the extra additive parameter $B$ added, the conjecture of Gonek and Ng  is favored over the other conjectures discussed, whereas, as we have seen, without the parameter $B$ (i.e., with $B = 0$), the claim of Kotnik and van de Lune that  $\degl_2 q = 0$ and $\degl_3 q$ appears close or equal to $\frac{1}{2}$ is favored over the other conjectures.

\begin{table}[!htbp]
\caption{\centering  Best fit curves $B+Ax^{a_2}(\log x)^{a_3}$ to $|q_{10^{6}}(e^{e^x})|$ at the IL$|q_{10^6}|$ values $e^{e^x}$,  under various constraints }
\footnotesize
\begin{tabular}{|l||r|r|r|r||r|} \hline
Constraints &   $a_2$ & $a_3$  & $A$ & $B$ & Adjusted $R^2$ \\ \hline\hline
$B = 0$  &   0.132648 & 0.170467 &  0.472356 & 0 & 0.999516 \\ \hline
$a_2 = 0$, $B = 0$  &   0 & 0.490714 &  0.498901 & 0 & 0.999458 \\ \hline
$a_2 = 0$, $a_3 =  {\tfrac{1}{2}}$, $B = 0$ & $0$ & $\tfrac{1}{2}$ &  $0.493834$ & 0 & 0.999481  \\ \hline
$a_2 = 0$, $a_3 = 1$, $B = 0$	&		0	&	1 & 0.280583	&   0 &  0.989116	\\ \hline 
$a_2 = 0$, $a_3 = {\tfrac{5}{4}}$, $B = 0$		& 	0	&	$\tfrac{5}{4}$ &	0.209703  &  0 &  0.979111	\\ \hline
$a_3 = 0$, $B = 0$  &   0.202478 & 0 &  0.459861 & 0 & 0.999517 \\ \hline
$a_2 = {\tfrac{1}{2}}$, $a_3 = 0$, $B = 0$		& $\tfrac{1}{2}$	&	$0$  &	0.178904	&  0 & 0.976633  \\ \hline
${a_2 = 0}$	&		0	&	1.11312 & 0.126778	&   0.418209 &  0.999502	\\ \hline 
$a_3 = \tfrac{1}{2}$ & $0.123669$ & $\tfrac{1}{2}$ &  $0.221761$ & 0.291093 & 0.999515  \\ \hline
${a_2 = 0}$, ${a_3 = {\tfrac{1}{2}}}$ & $0$ & $\tfrac{1}{2}$ &  $0.483888$ & 0.0171221 & 0.999460  \\ \hline

$a_3 = 1$		& 	$0.039753$	&	$1$ &	0.128462  & 0.41373 & 0.999506 \\ \hline
$a_2 = 0$, $a_3 = 1$	&		0	&	1 & 0.156217	&   0.381443 &  0.999526	\\ \hline 

$a_3 = {\tfrac{5}{4}}$		& 	$-0.013617$	&	$\tfrac{5}{4}$ &	0.10617  & 0.445702 & 0.999499 \\ \hline
${a_2 = 0}$, ${a_3 = {\tfrac{5}{4}}}$		& 	0	&	$\tfrac{5}{4}$ &	0.0997037  & 0.453741&  0.999525 	\\ \hline

none		& 	0.140224	&	$0.28531$	&	0.329144 & 0.161942 &  0.999488 \\ \hline
$a_3 = 0$  & 0.0784993 & 0 & 1.63255 & $-1.21821$ & 0.999506 \\ \hline
${a_3 = 0}$, ${B \geq 0}$  &   0.202478 & 0 &  0.459861 & 0 & 0.999491 \\ \hline

$a_2 = \tfrac{1}{2}$ 		& 	$\tfrac{1}{2}$	& $-0.0588966$	& 0.0941453 & 0.439545 &  0.999338 \\ \hline
${a_2 = {\tfrac{1}{2}}}$, ${a_3 = 0}$		& $\tfrac{1}{2}$	&	$0$  &	0.0848502	& 0.454631&  0.999372  \\ \hline
\end{tabular}
\end{table}\label{tablebestfit}

More variations on this theme can be found by instead finding best-fit curves of the form
$$f(x) = A (x^{a_2}+C_2) (( \log x)^{a_3}+C_3)$$
or of the form
$$f(x) = A (x+c_2)^{a_2} (\log (x+c_3))^{a_3}$$
under various constraints.  More generally, as indicated in Remark \ref{Lnorm}, one can use different normalizations of the sequence of iterated logarithms, the danger here being that choices can be made to favor one conjectural value of $\degl_2 M$ or $\degl_3 M$ over another.  All of this suggests that, with the data currently at hand, there is no way to adjudicate between various conjectural values of $\degl_2 M$ and $\degl_3 M$, including but not limited to those implied by the conjectures of Gonek and Ng, of Good, Churchhouse and L\'evy, and of Kotik and van de Lune.

For the reasons described above, we limited Conjecture \ref{Mconjecture} to the statement that $\deg M = \frac{1}{2}$ and $\degl_1 M  = 0$ and have refrained from making a conjecture regarding $\degl_2 M$, although we are somewhat sympathetic to the conjecture that $\degl_2 M = 0$.  As with $\li(x)-\pi(x)$, we suspect that statements involving $M(x)$ and $\log \log \log x$, including any purported bounds on $\degl_3 M$ and possibly even $\degl_2 M$, besides those proven rigorously,  will not have sufficient numerical support for the forseeable future.  Granted, say, any analogue $M(x) = \Omega_{\pm}(\sqrt{x} \, (\log \log \log x)^d) \ (x \to \infty)$ for some $d > 0$ of Littlewood's result $x-\psi(x) = \Omega_{\pm}(\sqrt{x} \, \log \log \log x) \ (x \to \infty)$   would be a substantial improvement upon $M(x) = \Omega_{\pm}(\sqrt{x}) \ (x \to \infty)$, as it would imply that $\degl M \neq (\frac{1}{2},0,0,0,\ldots)$ and $M(x) \neq O( \sqrt{x}) \ (x \to \infty)$ (and also that $\degl_3 M \geq d$ if $\degl_1 M = \degl_2 M = 0$).    However, currently the ``best case'' scenario that is consistent with known results is that $M(x) =O(\sqrt{x}) \ (x \to \infty)$, or more generally that $\degl M = (\frac{1}{2},0,0,0,\ldots)$, and the discussion above suggests that  known numerical data provides little to no evidence that these are false.  Notwithstanding the fact that  $M(x) \neq O(\sqrt{x}) \ (x \to \infty)$ is implied by the conjecture that the imaginary parts  of the zeros of  $\zeta(s)$ are linearly independent \cite[Theorem A]{ing}, and notwithstanding the fact that most experts seem to find $M(x) \neq O( \sqrt{x}) \ (x \to \infty)$ and either $\degl_2 M> 0$ or $\degl_3 M > 0$ to be more plausible than the alternatives, there is little numerical evidence to support any of  these conjectures.

On heuristic grounds, on the other hand, the similarities between the functions $M(x)$ and $x-\psi(x)$ are strong enough that one might expect $\degl M$ and $\degl(\id-\psi)$ to take on similar forms.  Just as the explicit formula for $\psi(x)$ yields $x - \psi(x) = \sum_{\rho} \frac{x^{\rho}}{\rho} + O(1) \ (x \to \infty)$ and $\dege (\id -\psi) = \dege \left(\sum_{\rho} \frac{x^{\rho}}{\rho}\right)$, the analogous explicit formula for $M(x)$ \cite[(3)]{hur} \cite{tit},  assuming the Riemann hypothesis and that all zeros of $\zeta(s)$ are simple, yields 
${M(x)} = \sum_{\rho} \frac{x^{\rho}}{\rho \zeta'(\rho)} + O(1) \ (x \to \infty)$
and therefore 
$\dege M = \dege \left(\sum_{\rho} \frac{x^{\rho}}{\rho \zeta'(\rho)}\right)$.
Moreover, under the Riemann hypothesis, the functions $(e^x-\psi(e^x))e^{-x/2}$ and $(\li(e^x)-\pi(e^x))xe^{-x/2}$ both possess a limiting distribution, and under the additional hypothesis \cite[(1.17)]{ANS} so does the function $M(e^x)e^{-x/2}$ \cite{ANS}.  Under these hypotheses, it is not  unreasonable to contemplate the possibility that their respective logexponential degrees all take on the same form, which  is $(0,0,d,0,0,0,\ldots)$ if we are also to believe something like the conjecture (\ref{MMC3}) of Montgomery or the conjecture (\ref{gonek}) of Gonek and Ng is true.  Their work thus provides motivation for the claims that
$$\dege(\id-\psi) = (\tfrac{1}{2},0,0,\Theta_3, 0, 0, 0, \ldots)$$
and
$$\dege M = (\tfrac{1}{2},0,0,d_3, 0, 0, 0, \ldots),$$
where $\Theta_3 \geq 1$ and $d_3 \geq 0$.   Since we have conjectured that $\Theta_2 = 0$, albeit with much less confidence than we have conjectured that $\Theta_1 = -1$, the comparisons between $M$ and $\id -\psi$ noted above lead us also  to be somewhat sympathetic to the conjecture that $\dege_2 M = 0$.  However, we do not formally assert $\dege_2 M = 0$ as a conjecture as it lacks sufficient numerical evidence and  contradicts the conjectures of Good, Churchhouse, and L\'evy.  By contrast, the conjecture $\dege_2 (\li-\pi)= 0$, though also lacking substantial numerical support, does not contradict any conjectures of which we are aware and is implied by the conjectures of Montgomery and of Stoll and Demichel.

Further comparisons between  $M$ and $\id -\psi$ led Kotnik and van de Lune to speculate that, ``if all this is taken to account, it would prehaps not be too surprising if $q(x) = o((x-\psi(x))/\sqrt{x})$'' \cite[p.\ 480]{kot}.  This speculation seems to be overly optimistic, however, as it would require for all $N \geq 0$   that the inequality  $|x-\psi(x) | \leq N$ imply $M(\lfloor x \rfloor) = 0$ for all sufficiently large $x$.  Nevertheless,  it might at least suggest the more conservative claim that $\dege M \leq \dege(\id-\psi)$, or perhaps even that $\dege M < \dege(\id-\psi)$.

In a different vein, note that
$$\psi(x) = \sum_{n \leq x} \Lambda(n) =  \sum_{n \leq x} \sum_{ab = n} \mu(a) \log b  = \sum_{ab \leq x} \mu(a) \log b =  \sum_{n \leq x} \log (n) M\left(\frac{x}{n}\right),$$
where $\Lambda(n) = \sum_{ab = n} \mu(a) \log b$ is the von Mangoldt function (which assumes the value $\log p$ if $n > 1$ is a power of a prime $p$, and $0$ otherwise).   By M\"obius inversion, one has
$$M(x) = \sum_{n \leq x} \mu(n) \log(n) \psi\left(\frac{x}{n}\right).$$
This motivates the following.

\begin{outstandingproblem}\label{MLprob2}
Assuming  hypotheses regarding the zeros of $\zeta(s)$, can one relate $\dege M$ and $\dege(\id-\psi)$ to each another?  If so, then how?  In particular, does either Conjecture \ref{eurekaconjecture} or Conjecture \ref{Mconjecture} imply, or contradict, the other?  Is it true that  $\dege M \leq \dege(\id-\psi)$?  If so, is it also true that $\dege M < \dege(\id-\psi)$?

\end{outstandingproblem}

In \cite[p.\ 578]{hum}, Humphries  suggests that the same reasoning that leads to the conjecture  (\ref{gonek}) of Gonek and Ng leads to the conjecture that
$$0< \limsup_{x \to \infty} \frac{|L(x)|}{\sqrt{x}\, (\log \log \log x)^{5/4}}   < \infty,$$
which if true would yield $\dege L = (\frac{1}{2}, 0, 0, \frac{5}{4}, 0, 0, \ldots)$.  Since also $L(x) = \sum_{b^2 \leq x} M(x/b^2)$ for all $x > 0$, and since $L(x)$, assuming the Riemann hypothesis and that the zeros of $\zeta(s)$ are simple, has an explicit formula similar to that of $M(x)$ but with central term $\frac{\sqrt{x}}{\zeta(1/2)}$ instead of $-2$ \cite{fawaz},  it seems reasonable to  conjecture the following.

\begin{conjecture}
One has $\dege L  = \dege \left(L(x)-\frac{\sqrt{x}}{\zeta(1/2)}\right) = \dege M$.
\end{conjecture}

\section{Conjectures concerning the Riemann zeta function}

We assume the notation of Sections 5.3 and 11.2. 

The constant $$\dege_1 S \in [\tfrac{1}{3},1]$$
is a sort of dual to the Riemann constant $\Theta = \deg(\li-\pi)  \in [\tfrac{1}{2},1]$.  Assuming the Riemann hypothesis, one has $\dege_1 S \in  [\tfrac{1}{2},1]$.

In \cite{farm}, Farmer, Gonek, and Hughes  conjectured that 
\begin{align}\label{fghc}
\limsup_{T \to \infty} \frac{|S(T)|}{\sqrt{\log T \log \log T}} = \frac{1}{\pi \sqrt{2}}.
\end{align}
If the conjecture above is true, or, more generally, if the given lim sup is finite and positive, then one has
$$\dege S = (0,\tfrac{1}{2},\tfrac{1}{2}, 0, 0, 0, \ldots),$$
or, equivalently, 
$$\dege(2\pi r_n(a)-\gamma_n) = (0,-\tfrac{1}{2},\tfrac{1}{2}, 0, 0, 0, \ldots)$$
for all (or some) $a \in \RR$.  In particular, their conjecture implies the following much weaker conjecture regarding the constant $\dege_1 S$.

\begin{conjecture}\label{tconj2} Let $a \in \RR$,  and let 
$$r_n(a) =  \frac{n-a}{W((n-a)/e)}.$$
One has
$$\dege_1 S = \tfrac{1}{2},$$
or, equivalently, 
$$\dege_1 \left(n-\widehat{\gamma}_n \right) =\tfrac{1}{2},$$
or, equivalently still,
$$\dege_1 (2\pi r_n(a)-\gamma_n) = -\tfrac{1}{2}.$$
\end{conjecture}

\begin{figure}[ht!]
\includegraphics[width=80mm]{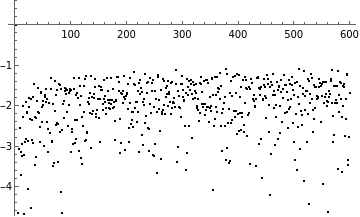}
\caption{\centering Plot of $\frac{\log\left|\frac{n-11/8}{W((n-11/8)/e)}-\frac{\gamma_n}{2\pi}\right|}{\log \log n}$ for $n = 1,2,3,\ldots,600$}
  \label{ZetaW}
\end{figure}

\begin{figure}[ht!]
\includegraphics[width=80mm]{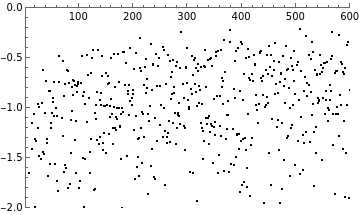}
\caption{\centering Plot of $\frac{\log\left|n-\widehat{\gamma}_n\right|}{\log \log n}$ for $n = 1,2,3,\ldots,600$}
\label{gammahat}
\end{figure}

Conjecture \ref{tconj2} states, loosely speaking, that the logexponential degree of order $1$ of each of the functions $S(T)$, $n-\widehat{\gamma}_n$, and $2\pi r_n(a)-\gamma_n$ is as small as the Riemann hypothesis will allow, as far as we currently know.  Figures \ref{ZetaW} and  \ref{gammahat}  provide plots of $\frac{\log\left|r_n\left(\frac{11}{8}\right)-\frac{\gamma_n}{2\pi}\right|}{\log \log n}$ and $\frac{\log |n-\widehat{\gamma}_n|}{\log \log n}$, respectively, for $n = 1,2,3,\ldots, 600$.  Using Mathematica, we computed that the maximum value of $\frac{\log |n-\widehat{\gamma}_n|}{\log \log n}$ on $[1,10000]$ is $-0.024009\ldots$ and occurs at $n = 8571$.  These computations provide very modest support for Conjecture \ref{tconj2}.

We now consider the spacings $\tau_{n+1}-\tau_n$, where $\tau_n = \frac{ \gamma_n}{2\pi}$.
Recall from Proposition \ref{zzgap2} that
 $$ (0,-1,0,0,0,\ldots) \leq \dege (\tau_{n+1}-\tau_n) \leq (0,0,0,-1,0,0,0,\ldots).$$
 In particular, one has $\deg (\tau_{n+1}-\tau_n)  = 0$ and $\dege_1 (\tau_{n+1}-\tau_n)   \in [-1,0]$.
We make the following conjecture.

\begin{conjecture}\label{gapc}
One has
$$\dege_1(\tau_{n+1}-\tau_n) >- 1.$$
Equivalently, there exists a $t \in (-1,0)$ such that
$$\tau_{n+1}-\tau_n \neq O((\log x)^t) \ (x \to \infty).$$
\end{conjecture}

\begin{figure}[ht!]
\includegraphics[width=80mm]{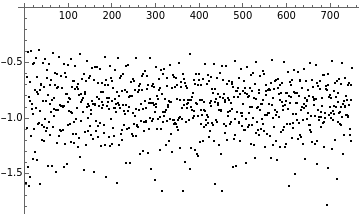}
\caption{\centering Plot of $\frac{\log(\tau_{n+1}-\tau_n)}{\log \log n}$ for $n = 3,4,5,\ldots,750$}
   \label{taugaps}
\end{figure}

 Figure \ref{taugaps} provides a plot of $\frac{\log(\tau_{n+1}-\tau_n)}{\log \log n}$ for $n = 3,4,5,\ldots,750$.

By Proposition \ref{tnln}, one has
$$\dege(\widehat{\gamma}_{n+1}-\widehat{\gamma}_n)  = \dege \delta_n =  \dege( \tau_{n+1}-\tau_n) +(0,1,0,0,0,\ldots).$$
Thus, Conjecture \ref{gapc} is equivalent to 
$$\dege_1(\widehat{\gamma}_{n+1}-\widehat{\gamma}_n) > 0$$
and to
$$\dege_1 \delta_n > 0,$$
where $\delta_n =   (\tau_{n+1}-\tau_n) \log \tau_n$ is the $n$th normalized spacing.   Moreover, Conjecture \ref{gapc} immediately implies Montgomery's conjecture 
\begin{align}\label{mcon}
\limsup_{n \to \infty}\, (\tau_{n+1}-\tau_n) \log \tau_n = \infty
\end{align}
\cite[p.\ 185]{mont3}, which is equivalent to
$$\tau_{n+1}-\tau_n \neq O\left(\frac{1}{\log n}\right) \ (n \to \infty).$$
It also implies an affirmative answer to the question posed in Problem \ref{mainx}.  Note,  however, that a more modest conjecture than Conjecture \ref{gapc}, namely,
$$\dege(\tau_{n+1}-\tau_n) > (0,-1,0,0,0,\ldots),$$
is  still strong enough  both to imply Montgomery's conjecture (\ref{mcon}) and  to answer the question posed in Problem \ref{mainx} affirmatively.

Much stronger conjectures have been made regarding the spacings $\tau_{n+1}-\tau_n$.
According to \cite{bui}, the {\it Gaussian Unitary Ensemble (GUE) hypothesis} suggests that the spacings should get as large as a constant times $\frac{1}{\sqrt{\log \tau_n}}$ for infinitely many $n$, that is, that
$$\tau_{n+1}-\tau_n = \Omega_+ \left(\frac{1}{\sqrt{\log \tau_n}} \right) \ (n \to \infty).$$
This conjecture implies that
$$\dege  (\tau_{n+1}-\tau_n) \geq (0,-\tfrac{1}{2},0,0,0,\ldots).$$
Moreover, a 2013 suggestion (if not conjecture) of  G.\ B.\ Arous and P.\ Bourgade \cite{arous}, namely, that
$$\limsup_{n \to \infty} \, (\gamma_{n+1}-\gamma_n) \sqrt{\frac{\log \gamma_n}{32}} = 1,$$
or even the weaker claim that
\begin{align}\label{arousbour}
0< \limsup_{n \to \infty} \, \frac{\gamma_{n+1}-\gamma_n} {(\log \gamma_n)^{-1/2}} < \infty,
\end{align}
implies that
$$\dege  (\tau_{n+1}-\tau_n) = (0,-\tfrac{1}{2},0,0,0,\ldots).$$

\begin{figure}[ht!]
\includegraphics[width=80mm]{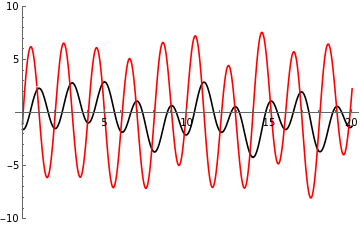}
\caption{\centering Graph of $\widehat{Z}(x) = Z\left(2\pi  \frac{x-7/8}{W((x-7/8)/e)}\right)$ (in black) and its derivative (in red) on $[0,20]$}
   \label{Zhata}
\end{figure}

\begin{figure}[ht!]
\includegraphics[width=80mm]{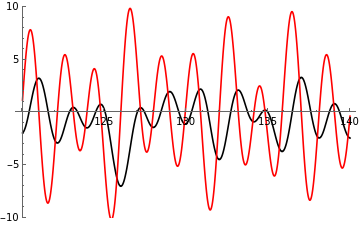}
\caption{\centering Graph of  $\widehat{Z}(x) = Z\left(2\pi \frac{x-7/8}{W((x-7/8)/e)}\right)$  (in black) and its derivative (in red)  on $[120,140]$}
   \label{Zhatb}
\end{figure}

\begin{figure}[ht!]
\includegraphics[width=80mm]{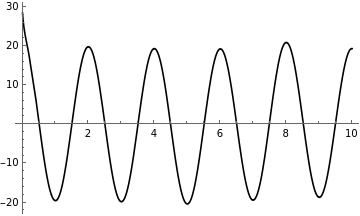}
\caption{\centering Graph of the second derivative of $\widehat{Z}(x)$ on $[0,10]$}
   \label{ZDP}
\end{figure}

\begin{figure}[ht!]
\includegraphics[width=80mm]{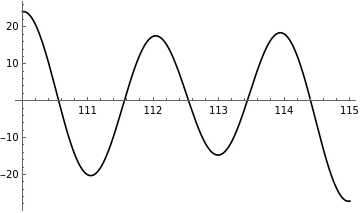}
\caption{\centering Graph of the second derivative of $\widehat{Z}(x)$ on  $[110,115]$}
   \label{ZDP2}
\end{figure}

Recall from Section 5.3 that we defined the normalization
$$\widehat{Z}(x) = Z\left(2\pi  \frac{x-7/8}{W((x-7/8)/e)}\right)$$
of the Riemann--Siegel $Z$ function $Z(T)$,
where $W$ is the  Lambert $W$ function.  
See Figures \ref{Zhata} and \ref{Zhatb} for a graph of the function $\widehat{Z}(x)$ and its derivative on $[0,20]$ and $[120,140]$, respectively.  Notice that the zeros of the derivative of $\widehat{Z}(x)$ in $[0,\infty)$ are very close to the nonnegative integers.   For example, the first zero is at $x =0.039096\ldots$ and the second at $x = 1.005778\ldots$, while the third is at $x = 1.986098\ldots$. See Figures  \ref{ZDP} and  \ref{ZDP2} for  a graph of the second derivative of the function $\widehat{Z}(x)$  on $[0,10]$ and $[110,115]$, respectively.  Examining up to the eighth derivative suggests a pattern: the $(n+N_k)$th largest zero of  the $k$th derivative of $\widehat{Z}(x)$ appears to be approximately $n-\frac{1}{2}$ if $k$ is even and $n$ if $k$ is odd, where $N_k$ is a small integer, $n$ is a positive integer, and $n+N_k \geq 1$ (where $N_k = 0,1,0,0,1,1,0,0,-1$ for $k = 0,1,2,\ldots,8$).  Regarding the case $k = 0$, we have seen that the $n$th largest zero of $\widehat{Z}(x)$, assuming the Riemann hypothesis, is $\widehat{\gamma}_n-\frac{1}{2}$, which satisfies
$$n-\tfrac{1}{2}-(\widehat{\gamma}_n-\tfrac{1}{2})  = S(\gamma_n) + (1+o(1)) \frac{\log n}{96 \pi^2 n}  \ (n \to \infty).$$

\begin{conjecture}\label{BIG}  For all positive integers $n$ and all nonnegative integers $k$, let $\widehat{\gamma}_{n,k}$ denote the $n$th largest zero of the $k$th derivative of $\widehat{Z}(x)$ on $(0,\infty)$.   For all nonnegative integers $k$, one has
$$\deg (n-  \widehat{\gamma}_{n,k}) = 0$$ and
$$\# \left\{\widehat{\gamma}_{n,k}: 0 < \widehat{\gamma}_{n,k} \leq T \right\} \sim T \ (T \to \infty).$$
\end{conjecture}

Regarding the function
$$\widehat{A}(x) = \frac{1}{\pi}\operatorname{Arg}\zeta\left(\frac{1}{2}+2\pi i \frac{x-7/8}{W((x-7/8)/e)}\right)$$
discussed in Section 5.3, we make the following conjecture.

\begin{conjecture}\label{BIG2}
 Let $-1\leq a \leq b \leq 1$.  For some $\beta \geq 1$, the function  $\widehat{A}(x)$ has values in $[a,b]\subseteq [-1,1]$  asymptotically as
$$\frac{\Gamma(\beta+\tfrac{1}{2})}{\Gamma(\beta)\Gamma(\tfrac{1}{2})}\int_a^b (1-t^2)^{\beta-1} \, dt,$$
which is known as the {\bf beta distribution on $[-1,1]$ with shape parameters $\beta$ and $\beta$}.\index{beta distribution}
In other words, one has
$$\frac{1}{T}\operatorname{meas}\{x \in [0,T]: \widehat{A}(x) \in [a,b] \} \sim \frac{\Gamma(\beta+\tfrac{1}{2})}{\Gamma(\beta)\Gamma(\tfrac{1}{2})}\int_a^b (1-t^2)^{\beta-1} \, dt \ (T \to \infty),$$
where $\operatorname{meas}$ denotes Lebesgue measure.
\end{conjecture}

\begin{remark}[The beta distribution]
The beta distribution on $[-1,1]$ with shape parameters $1$ and $1$ is just the uniform distribution.  With shape parameters $\frac{3}{2}$ and $\frac{3}{2}$, the beta distribution on $[-1,1]$ is given by
$$\frac{2}{\pi}\int_a^b \sqrt{1-t^2} \, dt =\left. \frac{1}{\pi} \left(  \arcsin t+ t\sqrt{1-t^2} \right) \right|_a^b,$$ known as the {\bf Wigner semicircle distribution on $[-1,1]$}.\index{Wigner semicircle distribution}  
The functions $\sin x$ and $\cos x$ have values in $[a,b] \subseteq [-1,1]$ asymptotically as $$\frac{1}{\pi}\int_a^b \frac{1}{\sqrt{1-t^2}} \, dt =\left. \frac{1}{\pi}  \arcsin t \right|_a^b,$$
which is the beta distribution on $[-1,1]$ with shape parameters $\frac{1}{2}$ and $\frac{1}{2}$. Because the functions $\sin x$ and $\cos x$ have derivative $0$ at their extreme values $1$ and $-1$, the density function $\frac{1}{\pi \sqrt{1-t^2}}$ is infinite there.  By contrast, the density function $\frac{2}{\pi} \sqrt{1-t^2}$ is zero at $1$ and $-1$, as is the density function of the beta distribution on $[-1,1]$ with shape parameters $\beta$ and $\beta$, for all $\beta > 1$ (which has infinite derivative at $1$ and $-1$ for $\beta \in (1,2)$).  

Remarkably, the beta distribution on $[0,1]$ with shape parameters   $\frac{1}{2}$ and $\frac{1}{2}$ arises in the study of the divisor function $d(n)$: by \cite[Th\'eor\`eme DDT]{ddt}, for all $x \in [0,1]$, one has
$$\lim_{T \to \infty} \frac{1}{T} \sum_{n \leq T} \frac{\#\{d|n: d \leq n^x\}}{d(n)} =\frac{1}{\pi} \int_0^x  \frac{1}{\sqrt{t(1-t)}} \, dt = \frac{2}{\pi} \arcsin \sqrt{x}.$$
In other words, for all $x \in [0,1]$, the mean value of the probability that a randomly chosen divisor of a positive integer $n$ is less than or equal to $n^x$ is $\frac{2}{\pi} \arcsin \sqrt{x}$.  (For $x = 0,\frac{1}{4},\frac{1}{2},\frac{3}{4},1$, the mean values are $0,\frac{1}{3},\frac{1}{2},\frac{2}{3},1$, respectively.)
\end{remark}

 An obvious question, if the values  of $\widehat{A}(x)$ are indeed beta distributed, is what the exact value of the shape parameter $\beta$ is.  (Note that symmetry suggests that the two shape parameters are equal.)    In the most extreme scenario, the functions $S(T)$ and $R(T)$ randomize, or ``even out,'' the values  of $\widehat{A}(x)$, so that they are no longer centered around $0$ but rather  are uniformly distributed in $[-1,1]$,  in the long run.  If this holds, which is a distinct possibility, then one has $\beta = 1$.   Relevent here is the open problem  \cite[p.\ 40]{reyna} as to whether or not the sequence $\widehat{\gamma}_n$ of normalized zeros  is uniformly distributed modulo $1$.  Note that the three sequences $\widehat{\gamma}_n$, $n-\widehat{\gamma}_n$, and $S(\gamma_n)$ are all uniformly distributed modulo $1$, or none of them is.  Should either (hence both) of the sequences $n-\widehat{\gamma}_n$ and $S(\gamma_n)$ be uniformly distributed modulo $2$ (hence also modulo $1$), then, given that the sequences $n-\widehat{\gamma}_n$ and $S(\gamma_n) \text{ mod }  2$ determine the locations of the  approximately linear portions of the graph $\widehat{A}(x)$ as in Proposition \ref{ahatprop}, that is likely enough to imply that the values of $\widehat{A}(x)$ are uniformly distributed in $[-1,1]$.  However, the claim that the given sequences are uniformly distrubuted modulo $1$ has yet to receive much support \cite{reyna}.

\begin{figure}[h!]
\includegraphics[width=60mm]{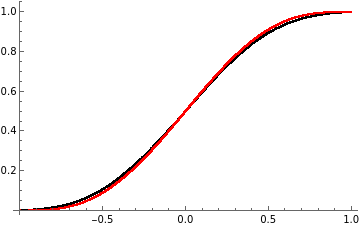}
    \caption{\centering Cumulative beta distrubution on $[-1,1]$ with shape parameter $B = 3.381750264764\ldots$ (in red), and cumulative distribution of values of $\widehat{A}(x)$ (in black) for $x \in [0,1000]$, with $10^5$ values sampled randomly (uniformly)}
\label{Ahat354}
\end{figure}

\begin{figure}[h!]
\includegraphics[width=60mm]{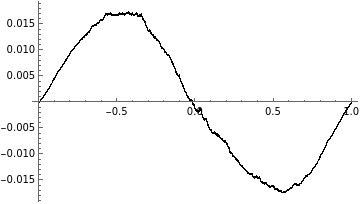}
    \caption{\centering Difference between two plots in Figure \ref{Ahat354} }
\label{Ahat354b}
\end{figure}

\begin{figure}[h!]
\includegraphics[width=60mm]{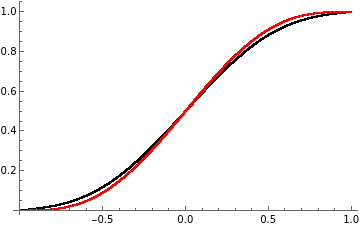}
    \caption{\centering Cumulative beta distrubution on $[-1,1]$ with shape parameter $B = 3.381750264764\ldots$ (in red), and cumulative distribution of values of $\widehat{A}(x)$ (in black) for $x \in [0,10000]$, with $10^6$ values sampled randomly (uniformly)}
\label{Ahat464}
\end{figure}

\begin{figure}[h!]
\includegraphics[width=60mm]{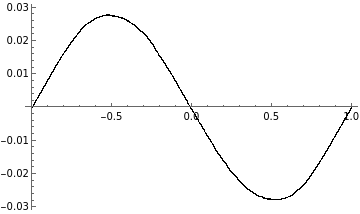}
    \caption{\centering Difference between two plots in Figure \ref{Ahat464} }
\label{Ahat464b}
\end{figure}

 In the opposite direction,  the heuristic evidence described below, along with data collected on the intervals $[0,10^2]$, $[0,10^3]$, and $[0,10^4]$, suggest that $\beta$ is less than the unique positive real number $B$ satisfying $$\frac{\Gamma(B+\tfrac{1}{2})}{\Gamma(B)\Gamma(\tfrac{1}{2})} = 1,$$ which is $$B = 3.381750264764\ldots.$$   Here, $\frac{\Gamma(\beta+\frac{1}{2})}{\Gamma(\beta)\Gamma(\frac{1}{2})}$ is the slope of the cumulative beta distribution, with shape parameters $\beta$ and $\beta$, at $0$.  Figure \ref{Ahat354} compares the beta distribution with shape parameters $B$ and $B$ to the observed distribution of values of $\widehat{A}(x)$ on $[0,1000]$, the latter of which was approximated by selecting $10^5$ values in $[0,10^3]$ randomly (uniformly) and plotting the proportion of these values that lie in the interval $I_k = [-1,-1+2 \cdot 10^{-4}k]$ for $k = 1,2,3,\ldots, 10^4$.  The slope of the observed distribution at $0$ is about $1.05$.  Figure \ref{Ahat354b} shows the difference between the two plots in  Figure \ref{Ahat354}.  Figure \ref{Ahat354} compares the same beta distribution  to the observed distribution of values of $\widehat{A}(x)$ on $[0,10^4]$, which was approximated by selecting $10^6$ values in $[0,10^4]$ randomly (uniformly) and plotting the proportion of these values that lie in the  intervals $I_k$ as above.  The slope of the observed distribution at $0$ is about $0.985$.   Figure \ref{Ahat464b} shows the difference between the two plots in  Figure \ref{Ahat464}.   These observations suggest that the values of $\widehat{A}(x)$ on $[0,T]$ are approximately beta distributed with shape parameter $\beta_T$ for some optimal $\beta_T$, which can be closely approximated via the slope method noted above.  The $\beta_T$ tend to decrease as $T \to \infty$, since the values of $\widehat{A}(x)$ apparently become less concentrated at $0$, and our conjecture is that the distribution of values converges to the beta distribution with shape parameter $\beta = \lim_{T \to \infty} \beta_T$.    The following analysis (which the author owes in part to Aditya Baireddy) suggests that $\beta \in [1,B)$.    Since $\widehat{A}(x)$ is piecewise approximately linear of slope $-1$, the inverse image of $[0,\varepsilon]$, intersected with $[0,T]$, is a union of finitely many disjoint intervals, all of which have length asymptotic to $\varepsilon$ as $\varepsilon \to 0^+$.   Moreover, one has $\widehat{A}(x) = 0$ only at the $\widehat{\gamma}_n$ and at the numbers $$\widehat{G}_n = \tfrac{G_n}{2\pi} \log \tfrac{G_n}{2\pi e}+\tfrac{7}{8} \approx n$$ for those Gram points $$G_n = \theta^{-1}((n-1)\pi)$$ that are {\it good} Gram points (i.e., for those $G_n$ at which $\zeta\left(\frac{1}{2}+it\right)$ is real and positive and thus has argument $0$).  Say, the number of such normalized good Gram points $\widehat{G}_n$ in $[0,T]$ is $G(T)$.  Then the derivative of the distribution of values of $\widehat{A}(x)$ on $[0,T]$ at $x = 0^+$ is approximately $\frac{G(T)}{T}$.  In the limit as $T \to \infty$, we  expect that the derivative of the distribution of values of $\widehat{A}(x)$ at $x = 0^+$ is at most the lower density $\liminf_{T\to\infty} \frac{G(T)}{T}< 1$ of all Gram points that are good.    Thus we expect that $\beta$ satisfies $$\frac{\Gamma(\beta+\frac{1}{2})}{\Gamma(\beta)\Gamma(\frac{1}{2})} \leq \liminf_{T\to\infty} \frac{G(T)}{T} < 1.$$   Note that $\beta > 1$ would require that the slope $\frac{\Gamma(\beta+\frac{1}{2})}{\Gamma(\beta)\Gamma(\frac{1}{2})}(1-t^2)^{\beta-1}$ of the cumulative distribution at $t = 1^-$ and $t = -1^+$ be $0$, signifying in that case   that the Gram points that are bad (i.e., not good) have a negligible effect on the distribution of values of $\widehat{A}(x)$, and thus implying a strong asymmetry between  the good Gram points and the bad Gram points.   It could be, then, the most compelling guess for $\beta$ is $\beta = 1$, after all.  Alternatively, Conjecture \ref{BIG2} is rendered false altogether if the values of $\widehat{A}(x)$ are nonuniformly distributed in $[-1,1]$ according to some density function that is nonzero in a neighborhood of $1$ and $-1$.

\begin{remark}[Alternative normalizations]
One can formulate analogues of the conjectures in this section for the alternative normalizations $\widecheck{\gamma}_n$, $\widecheck{Z}(T)$, and $\widecheck{A}(T)$ discussed in Remark \ref{altnorm}.   The conjectures in this section should be equivalent to their renormalized counterparts, with the possible exception of Conjecture \ref{BIG}, since it is not clear how to compare the zeros of the $n$th derivative of $\widehat{Z}(x)$ with the zeros of the $n$th derivative of $\widecheck{Z}(x)$.  Thus,  Conjecture \ref{BIG} ought to be replaced with its analogue for  $\widecheck{Z}(x)$, if the two conjectures are not equivalent.
\end{remark}

\section{Conjectures on asymptotics for prime counts in intervals}

Expanding on Section 12.1,  this section discusses various conjectures on asymptotics for  prime counts in intervals.   

The following conjecture is widely supported, for example, by Richards in \cite{rich}, by Heath-Brown in \cite{heathbrown}, and by Granville in \cite[p.\ 393]{gran0}.

\begin{conjecture}\label{piliconj}
For all $t  \in (0,1)$, one has
$$\pi(x+x^t)-\pi(x) \sim \frac{x^t}{\log x} \ (x \to \infty).$$ 
\end{conjecture}

\begin{figure}[h!]
\includegraphics[width=70mm]{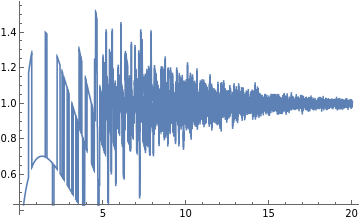}
    \caption{\centering Graph of $ \frac{\pi(e^x+e^{0.6x})-\pi(e^x)}{\li(e^x+e^{0.6x})-\li(e^x)}$ on $[0,20]$}
\label{pili1}
\end{figure}

\begin{figure}[h!]
\includegraphics[width=70mm]{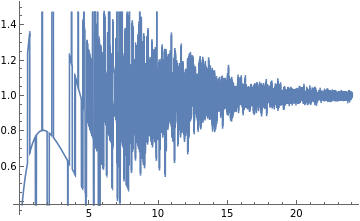}
    \caption{\centering Graph of $ \frac{\pi(e^x+e^{0.5x})-\pi(e^x)}{\li(e^x+e^{0.5x})-\li(e^x)}$ on $[0,24]$}
\label{pili2}
\end{figure}

\begin{figure}[h!]
\includegraphics[width=70mm]{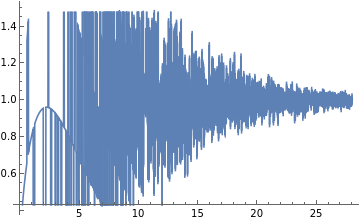}
    \caption{\centering Graph of $ \frac{\pi(e^x+e^{0.4x})-\pi(e^x)}{\li(e^x+e^{0.4x})-\li(e^x)}$ on $[0,28]$}
\label{pili3}
\end{figure}

\begin{figure}[h!]
\includegraphics[width=70mm]{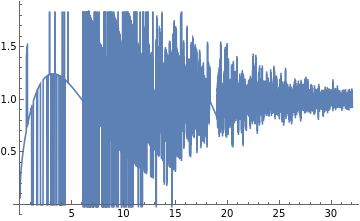}
    \caption{\centering Graph of $ \frac{\pi(e^x+e^{0.3x})-\pi(e^x)}{\li(e^x+e^{0.3x})-\li(e^x)}$ on $[0,32]$}
\label{pili4}
\end{figure}

If Conjecture \ref{piliconj} holds, then, for all $t  >0$, one has $\pi(p_n+p_n^t) > \pi(p_n)$, whence $g_n < p_n^t$, for all $n \gg 0$, and therefore $\deg g_n = 0$.   For all $t >0$, let 
$$F_t(x) = \frac{\pi(x+x^t)-\pi(x)}{\li(x+x^t)-\li(x)}.$$  Conjecture \ref{piliconj}  is equivalent to
\begin{align}\label{Ftconj}
\lim_{x \to \infty} F_t(x)   = 1
\end{align}
for all $t > 0$.  Note that, by the prime number theorem, (\ref{Ftconj}) holds for all $t \geq 1$, and thus, 
by Corollary \ref{lindel} and \cite[Theorem]{heathbrown},  it holds for all $t  \geq \frac{7}{12}$.  Mild evidence for the conjecture can be obtained by graphing the function $F_t(e^x)$ for successively smaller and smaller values of $t >0$ and observing that the graph of $F_t(e^x)$ for smaller values of $t$ seems to be  a ``stretched out'' version of the graph of $F_t(e^x)$ for larger values of $t$, all apparently having the same qualitative behavior as $x$ increases.  Since $F_t(e^x) \to 1$ for all $t \geq \frac{7}{12} = 0.58333\ldots$, this observation suggests that $F_t(e^x) \to 1$ for all $t > 0$.  See Figures \ref{pili1}, \ref{pili2}, \ref{pili3} and \ref{pili4} for graphs of the function $F_t(e^x)$ for $t = 0.6, 0.5, 0.4, 0.3$, respectively.

Note that Theorem  \ref{maierthm} and Conjecture \ref{piliconj} leave open the possibility that there might exist an unbounded function $h(x)> 0$ of degree $0$ (e.g., $h(x) = e^{\sqrt{\log x}}$) such that 
$$\pi(x+h(x))-\pi(x) \sim \frac{h(x)}{\log x} \ (x \to \infty).$$ 
An analysis similar to that above suggests that the asymptotic above holds for  $h(x) = e^{(\log x)^t}$ for any $t \in (0,1]$.

Recall that, in Problem \ref{dinf}, we let ${\mathbf h}$ denote the smallest $\dd \in \prod_{n = 1}^{\infty*}\overline{\RR}$ such that
$$\pi(x+h(x))-\pi(x) \sim \frac{h(x)}{\log x} \ (x \to \infty)$$
for all eventually positive real functions $h$ defined on a neighborhood of $\infty$ with $h(x) = O(x) \ (x \to \infty)$ and $\underline{\dege}\, h >  \dd$.  It follows from results in Section 12.1 that Conjecture \ref{piliconj} is equivalent to ${\mathbf h}_0 = 0$ and to
$$(0,\infty,0,1,0,0,0,\ldots) \leq {\mathbf h} \leq (0,\infty,1,0,0,0,\ldots),$$
where the lower bound is unconditional.  In particular,  Conjecture \ref{piliconj} is equivalent to the conjunction of ${\mathbf h}_0 = 0$, ${\mathbf h}_1 = \infty$, and ${\mathbf h}_2 \in [0,1]$.  It would be highly unusual under these circumstances if ${\mathbf h}_2$ turned out to be some constant in $(0,1)$.  Thus, it is natural to conjecture that ${\mathbf h}_2 = 0$ or ${\mathbf h}_2 = 1$.  In fact, if our suggestion above about the functions  $h(x) = e^{(\log x)^t}$ holds true, then by Corollary \ref{IVL2cor2} it also holds  true for any function $h$ with  $e^{(\log x)^t} \log x= O(h(x)) \ (x \to \infty)$ for some $t \in (0,1)$, which would then imply
$${\mathbf h} \leq (0,\infty,t,0,0,0,\ldots),$$
for all $t \in (0,1)$,  whence ${\mathbf h}_2 = 0$.  Thus, we have the following.

\begin{proposition}\label{e2is0}
Let ${\mathbf h}$ be defined as in Problem \ref{dinf}.   If ${\mathbf h}_0 = 0$, then ${\mathbf h}_1 = \infty$, and ${\mathbf h}_2$ is the infimum of all $t \in (0,1]$  such that 
$$\pi\left(x+e^{(\log x)^t }\right)-\pi(x) \sim \frac{e^{(\log x)^t }}{\log x} \ (x \to \infty).$$
Moreover,  the following conditions are equivalent.
\begin{enumerate}
\item ${\mathbf h}_0 = 0$ and ${\mathbf h}_2 = 0$.
\item $(0,\infty,0,1,0,0,0,\ldots) \leq {\mathbf h} \leq (0,\infty,0,\infty,1,0,0,0,\ldots)$.
\item For all $t \in (0,1)$ and any eventually positive real function $h$ defined in a neighborhood of $\infty$ with $h(x) = O(x) \ (x \to \infty)$ and $e^{(\log x)^t }= O(h(x)) \ (x \to \infty)$, one has
$$\pi(x+h(x))-\pi(x) \sim \frac{h(x)}{\log x} \ (x \to \infty).$$
\item For all $t \in (0,1]$, 
one has
$$\pi\left(x+e^{(\log x)^t }\right)-\pi(x) \sim \frac{e^{(\log x)^t }}{\log x} \ (x \to \infty).$$
\end{enumerate}
\end{proposition}

By Corollary \ref{maierthmcor2}, since $$\dege\, (\log x)^{\log^{\circ n}x} = (0,\infty,0,1,0,0,0,\ldots,0,1,0,0,0,\ldots)$$
for any $n \geq 3$, where the last $1$ is preceded by $n-3$ $0$s, one has the following.

\begin{proposition}
One has  ${\mathbf h} = (0,\infty,0,1,0,0,0,\ldots)$ if and only if
$$\pi(x+(\log x)^{\log^{\circ n}x}))-\pi(x) \sim \frac{(\log x)^{\log^{\circ n}x}}{\log x} \ (x \to \infty)$$
for all integers $n \geq 2$.
\end{proposition}

We make the following conjecture.

\begin{conjecture}
The equivalent conditions of Proposition \ref{e2is0} hold.
\end{conjecture}

If one suspects that Maier's result,  Theorem \ref{maierthm} \cite{maier}, is close to optimal, then one  might conjecture even more boldly that the unconditional lower bound $(0,\infty,0,1,0,0,0,\ldots)$ for ${\mathbf h}$ is an equality.     Such a conjecture can be motivated as follows.   Let  $h$ be any real function with $1 \ll h(x) = o(x) \ (x \to \infty)$.  Let
$$U(h) = \limsup_{x \to \infty} \frac{\pi(x+h(x))-\pi(x)}{h(x)/\log x}.$$
and
$$L(h) = \liminf_{x \to \infty} \frac{\pi(x+h(x))-\pi(x)}{h(x)/\log x}.$$
Of course,  one has $L(h) \leq U(h)$, and equality holds if and only if the asymptotic
$$\pi(x+h(x))-\pi(x) \sim \frac{h(x)}{\log x} \ (x \to \infty)$$
holds.
By Maier's theorem, one has
$$  L((\log x)^a) < 1 < U((\log x)^a)$$
for all $a > 1$.  It has been conjectured  (for example,  in \cite[(5)]{granl}) that 
$$\lim_{a \to \infty}  L((\log x)^a) = 1 = \lim_{a \to \infty}  U((\log x)^a),$$
or, equivalently,  that, for every $\varepsilon > 0$, there exists an $N>1$, and for each $a \geq N$ there exists an $N(a) > 1$,  such that 
$$1-\varepsilon \leq \frac{\pi(x+(\log x)^a)-\pi(x)}{(\log x)^{a}/\log x} \leq 1+\varepsilon$$
for all $x \geq N(a)$.  Consider the ``uniform'' generalization of the conjecture above stated in statement (4) of the following proposition.

\begin{proposition}\label{granconj0}
The following statements are equivalent.
\begin{enumerate}
\item For any real function $h$ defined and positive on some unbounded subset of $\RR_{>0}$ with $h(x) = O(x)$, if
$$\underline{\dege} \, h \geq (0,\infty,0,1,0,0,0,\ldots),$$
i.e., if $h(x) \gg (\log x)^a$ for all $a >0$, 
then
$$\pi(x+h(x))-\pi(x) \sim \frac{h(x)}{\log x}  \ (x \to \infty).$$
\item  For any real function $a$ defined and positive on some unbounded subset of $\RR_{>0}$,  if $1 \ll a(x) \leq  \frac{\log x}{\log \log x}$ for all $x \gg 0$,
then 
$$\pi(x+(\log x)^{a(x)})-\pi(x) \sim \frac{(\log x)^{a(x)}}{\log x}  \ (x \to \infty).$$
\item For every $\varepsilon > 0$, there exists an $N>1$ such that
$$1-\varepsilon \leq \frac{\pi(x+y)-\pi(x)}{y/\log x} \leq 1+\varepsilon$$
for all $x \geq N$ and all $y > 0$ with $(\log x)^N \leq y \leq x$.
\item For every $\varepsilon > 0$, there exists an $N>1$ such that
$$1-\varepsilon \leq \frac{\pi(x+(\log x)^a)-\pi(x)}{(\log x)^{a}/\log x} \leq 1+\varepsilon$$
for all $x \geq N$ and all $a \geq N$ with $a \leq  \frac{\log x}{\log \log x}$.
\item One has
$$\sup_{(\log x)^a \leq y \leq x} \left| \frac{\pi(x+y)-\pi(x)}{y/\log x}-1 \right| \to 0$$
as $\min(a,x) \to \infty$.
\end{enumerate}
\end{proposition}

\begin{outstandingproblem}\label{AConjec} 
Are the equivalent statements of Proposition \ref{granconj0} true?
\end{outstandingproblem}

\begin{proposition}\label{granconj}
Suppose that  Problem \ref{AConjec} has an affirmative answer, i.e.,  suppose that the equivalent statements of Proposition \ref{granconj0} are true.    One then has the following.
\begin{enumerate}
\item $\lim_{A \to \infty}  L((\log x)^A) = 1 = \lim_{A \to \infty}  U((\log x)^A)$.
\item There exists an eventually positive infinitely differentiable function $h$ with exact logexponential degree $(0,\infty,0,1,0,0,0,\ldots)$,  and,  for any such function $h$, one has
$$\pi(x+h(x))-\pi(x) \sim \frac{h(x)}{\log x}  \ (x \to \infty).$$
\item Let ${\mathbf h}$ be defined as in Problem \ref{dinf}.    One has $\mathbf{h} = (0,\infty,0,1,0,0,0,\ldots)$.
\end{enumerate}
 \end{proposition}

Note that an affirmative answer to Problem \ref{AConjec}, along with Corollary \ref{maierthmcor2}, which asserts the converse when $h$ is Hardian,  would yield a near complete answer to Problem \ref{maierprob}.   In 2019, Granville wrote \cite[p.\ 393]{gran0} that affirmative answer to Problem \ref{AConjec} is ``plausible.''  Specifically, he wrote that,  for $0 < y \leq x$,  ``it is plausible that (5) [viz., $\pi(x+y)-\pi(x) \sim \frac{y}{x}$] holds uniformly if  $\log y/\log \log x \to \infty$ as $x \to \infty$,'' and,  ``we conjecture,  presumably safely,  that (4) and (5) hold uniformly when $y > x^\varepsilon$.''    Thus,  it is tempting, though likely premature, to conjecture an affirmative answer to Problem \ref{AConjec}.

\section{Conjectures on prime gaps}

We assume the notation of Section 12.3.

 {\bf Cram\'er's conjecture}\index{Cram\'er's conjecture} states that
\begin{align}\label{cramer0}
g_n=O((\log n)^2) \ (n \to \infty).
\end{align}
The conjecture thus implies that $\deg g_n = 0$ and
$$1 \leq \dege_1 g_n \leq 2$$
and, more precisely,
$$(0,1,1,-1,1,0,0,0,\ldots) \leq \dege g_n \leq (0,2,0,0,0,\ldots).$$
Cram\'er offered the conjecture as being ``suggested'' by a ``heuristic method founded on probability arguments'' \cite{cramer}.
He  noted that his  heuristic model also suggests that  ``some  similar relation'' to the relation
\begin{align}\label{cramer}
 \limsup_{n\rightarrow\infty} \frac{g_n}{(\log p_n)^2}  = 1
\end{align} 
``may hold'' ``with a probability = $1$.'' 
However,  it quite possibly was not Cram\'er's intent to conjecture the full strength of (\ref{cramer}).  In \cite{gran}, Granville interprets Cram\'er as suggesting the stronger conjecture
\begin{align}\label{cramersim}
G(x) \sim (\log x)^2 \ (x \to \infty).
\end{align}
Despite his more stringent interpretation of Cram\'er's suggestions, Granville proposes a worthy revision to  Cram\'er's model that  suggests  that $G(x)$ is asymptotically greater than $2e^{-\gamma}(\log x)^2$, which if true would  invalidate (\ref{cramersim}).   Wolf, on the other hand, has provided further heuristic evidence for the conjecture (\ref{cramersim}) \cite{wolf1}.

At least the statement $G(x) \asymp (\log x)^2 \ (x \to \infty)$
is consistent with Granville's claims and all of the conjectures above and is enough to yield  both  (\ref{cramer}) 
and
$$\dege G  = \dege  g_n = (0,2,0,0,0,\ldots).$$  
However, Pintz suggests in \cite[p.\ 286]{pintz} that 
\begin{align}\label{pintzcon}
G(x) = O ((\log x)^t) \ (x \to \infty), \quad \forall t > 2,
\end{align}
and  speculates that one might have $G(x) \neq O ((\log x)^2) \ (x \to \infty)$.
At the very least, $\deg g_n = \deg G = 0$ and $\dege_1 g_n  =  \dege_1 G  \leq 2$ are implied by all related conjectures we can find in the literature.  
Thus, the following conjecture is an extension of (Piltz's) Conjecture \ref{gapconj} that is implied by each of the conjectures discussed above.

\begin{conjecture}\label{gapconj2}
Pintz's suggestion (\ref{pintzcon}) holds.
\end{conjecture}

The following result follows from Theorem \ref{GTM}.

\begin{proposition} 
Each of the following statements is equivalent to (\ref{pintzcon}), i.e., to Conjecture \ref{gapconj2}.
\begin{enumerate}
\item $\deg g_n = 0$ and $\dege_1 g_n \leq 2$.  
\item $\deg G = 0$ and $\dege_1 G \leq 2$.  
\item $g_n = O((\log n)^t) \ (n \to \infty)$ for all $t > 2$.
\item For all $t > 2$, one has $\pi(x+(\log x)^t))-\pi(x) \geq 1$ for all $x \gg 0$.
\item $(0,1,1,-1,1,0,0,0,\ldots) \leq \dege g_n \leq (0,2,\infty,1,0,0,0,\ldots).$
\item $(0,1,1,-1,1,0,0,0,\ldots) \leq \dege G \leq (0,2,\infty,1,0,0,0,\ldots).$
\end{enumerate}
\end{proposition}

One might also suspect that $\dege_1 g_n = 2$, which is an analogue of the conjecture $\dege_1(\gamma_{n+1}-\gamma_n) = -\frac{1}{2}$ discussed briefly in the previous section.   Indeed,  in 2018, Maynard conjectured the following.

\begin{conjecture}[Maynard {\cite[Conjecture 2]{mayn}}]\label{mayncram}
One has
$$G(x) = (\log x)^{2+o(1)} \ (x \to \infty).$$
\end{conjecture}

  Note the following.

\begin{proposition}
Each of the following statements is equivalent to Maynard's conjecture  $G(x) = (\log x)^{2+o(1)} \ (x \to \infty)$.
\begin{enumerate}
\item $G(e^x)$ has exact degree $2$.
\item $\lim_{x \to \infty} \frac{ \log G(x)}{\log \log x} = 2$.
\item  $G(x)$ has (exact) degree $0$ and exact logexponential degree $2$ of order $1$.  
\item $(0,2,-\infty,-1,0,0,0,\ldots) \leq \underline{\dege} \, G \leq\dege G \leq (0,2,\infty,1,0,0,0,\ldots)$.
\item  $G(x) \ll (\log x)^t$ for all $t > 2$ and $G(x) \gg (\log x)^t$ for all $t < 2$.
\item $(\log x)^{2-\varepsilon} \ll G(x) \ll (\log x)^{2+\varepsilon}$ for all $\varepsilon > 0$.
\item For any  eventually nonnegative nondecreasing function  $h(x)$ defined on a neighborhood of $\infty$,  if $h(x) \gg (\log x)^t$ for some $t > 2$, then
 $$\pi(x+ h(x))-\pi(x)) \geq 1, \quad \forall x \gg 0,$$
and, conversely, if
 $$\pi(x+ h(x))-\pi(x)) \geq 1, \quad \forall x \gg 0,$$
 then $h(x) \gg (\log x)^t$ for some $t < 2$.
\end{enumerate}
\end{proposition}

\begin{figure}[h!]
\includegraphics[width=70mm]{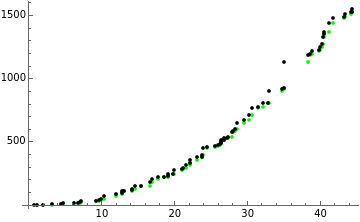}
    \caption{\centering Plot of $G(e^x)$ (in black) and $G((e^x)^-)$ (in green) at each (prime integer) value $n$ of $e^x$ where $G(e^x)$ has a discontinuity,  i.e., where $g_{\pi(n)}$ is a maximal prime gap, for the first $80$ maximal prime gaps}
\label{Gx}
\end{figure}

\begin{figure}[h!]
\includegraphics[width=70mm]{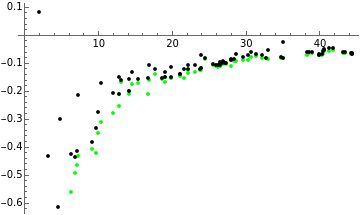}
    \caption{\centering Plot of $\frac{\log G(e^x)}{\log x}-2$ (in black) and $\frac{\log G((e^x)^-)}{\log x}-2$ (in green) at each  value $n$ of $e^x$ where $G(e^x)$ has a discontinuity,  for the first $80$ maximal prime gaps $g_{\pi(n)}$, excluding $n = 2,3$}
\label{LGx}
\end{figure}

\begin{figure}[h!]
\includegraphics[width=70mm]{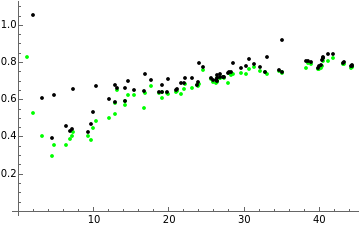}
    \caption{\centering Plot of $\frac{G(e^x)}{x^2}$ (in black) and $\frac{G((e^x)^-)}{x^2}$  (in green) at each value $n$ of $e^x$ where $G(e^x)$ has a discontinuity,  for the first $80$ maximal prime gaps $g_{\pi(n)}$, excluding $n = 2,3$}
\label{SGx}
\end{figure}

\begin{figure}[h!]
\includegraphics[width=70mm]{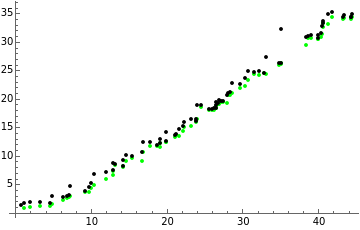}
    \caption{\centering Plot of $ \frac{G(e^x)}{x} = \widehat{g}_{\pi(n)}$ (in black) and $ \frac{G((e^x)^-)}{x}$ (in green) at each value $n$ of $e^x$ where $G(e^x)$ has a discontinuity,  for the first $80$ maximal prime gaps $g_{\pi(n)}$}
\label{MGx}
\end{figure}

To test Maynard's conjecture,  one most naturally plots and samples the statistic $\frac{ \log G(e^x)}{\log x}$ and checks to see if it appears to have a limit of $2$ as $x \to \infty$.  Such an analysis is a bit more straightforward than that of the function $\frac{G(x)}{(\log x)^2}$.    In particular, Maynard's conjecture  sidesteps  the various subtleties involved in conjecturing appropriate bounds on $\frac{G(x)}{(\log x)^2}$.    See Figures \ref{Gx},  \ref{LGx}, and \ref{SGx} for plots of $G(e^x)$, $\frac{ \log G(e^x)}{\log x}-2$, and $\frac{G(e^x)}{x^2}$,  respectively,  in black, and their respective limits from the left in green, at each (prime integer) value $n$ of $e^x$ where $G(e^x)$ has a discontinuity,  i.e.,  where $g_{\pi(n)}$  is a  (record) maximal  prime gap, for the first $80$ maximal prime gaps.  Note that the step function $G(e^x)$ is constant in between each such point.  The 80 $x$-coordinates in these plots are just $\log a(n)$, where $a(n)$ is OEIS Sequence A002386, while the $y$-coordinates for the points in black, respectively, are $b(n)$,  $\frac{\log b(n)}{\log \log a(n)}-2$,  and  $\frac{b(n)}{(\log a(n))^2}$, where $b(n)$ is OEIS Sequence A005250, for $n = 1,2,3,\ldots,80$.  The largest (i.e., the $80$th) of these maximal prime gaps occurs at the $20$-digit prime $$p_{423731791997205041} =  18361375334787046697,$$
for which the gap between the next prime is $g_{423731791997205041} = 1550$.  To date this is the largest known maximal prime gap, which was found by B.\ Nyman in 2018 and verified  by analyzing the prime gaps for all primes less than $2^{64} = 18446744073709551616$.  
The clear outlier in these three plots occurs at $\log p = 35.065386\ldots$, where $p$ is the prime number $$p_{49749629143526} =  1693182318746371,$$ for which the value of $$\frac{g_{\pi(p)}}{(\log p)^2}  = \frac{1132}{(\log( 1693182318746371))^2} = 0.920638\ldots$$
is the largest known, for $p > 7$.  The merit $\widehat{g}_{\pi(p)} = \frac{g_{\pi(p)}}{\log p}$, or normalized prime gap, of this prime gap of $1132$ is $32.282547\ldots$.  As of September 2022, the largest known merit of a prime gap, and the first found with merit over 40 (discovered by the Gapcoin network), is $41.938783\ldots$,  for a particular $87$-digit prime $p$ with $\log p = 199. 099717\ldots$ and with a prime gap of $8350$.  See Figure \ref{MGx} for a corresponding  plot of $\frac{G(e^x)}{x}$, i.e., a lin-log plot of the merits of the first $80$ maximal prime gaps.

Assuming that $\deg G = 0$,  the quantity $\dege_1 G$ is the lim sup of the points in black in Figure \ref{LGx}, while  $\underline{\dege}_1 G$ is the lim inf of the points in green in Figure \ref{LGx}, as $x \to \infty$.  Thus,  Maynard's conjecture holds if and only if all points, both black and green, tend to $0$ as $x \to \infty$.  In fact, by Proposition \ref{stepdeg}, one has $\dege G = \dege G|_X$ and $\underline{\dege} \, G = \underline{\dege} \, G(x^-)|_X$, where $X$ is the set of all points of discontinuity of $G$, i.e., where $X$ is the set of all primes $p_n$ for which $g_{n}$ is a maximal prime gap, so that $G(p_n) = g_{n}$ and $G(p_n^-) = \max_{k <n} g_k$ is the previous maximal prime gap, for each prime $p_n$ in $X$.    
 Thus,  one can restrict one's study of $G(x)$ to its points of discontinuity,  as we did in Figures \ref{Gx}--\ref{MGx}, at least for  the purpose of investigating  $\dege G$ and $\underline{\dege}\,  G$.   Recall that
$$\widetilde{G}(x) = \max_{p_{k+1} \leq x}  g_k = \max_{k < \pi(x)}  g_k, \quad \forall x \geq 0,$$
is a slight variant of $G(x)$ employed in place of $G(x)$ in some works, e.g.,  \cite{fgktm}.   
Note that
$$G(x^-) = \widetilde{G}(x), \quad \forall x \in X.$$   Thus, one has  $$\dege G = \dege \widetilde{G} \geq \dege \widetilde{G}|_X \geq \underline{\dege} \, \widetilde{G}|_X = \underline{\dege} \, G.$$
 Since the function $G(x)$ is relatively ``stable,''  e.g., it is a nonincreasing step function,  in light of Maynard's conjecture, it is plausible that $G(x)$ has exact logexponential degree to all orders, which holds if and only if $\dege G = \underline{\dege}\, \widetilde{G}|_X$.   This holds, for example, if $G(x) \asymp r(x)\ (x \to \infty)$ for some Hardian function $r(x)$.

 \section{Gemeralizations of the abc conjecture}

The {\bf radical} of a positive integer $n$, denoted $\operatorname{rad}(n)$,  is the product $\prod_{p |n} p $ of the distinct prime factors of $n$.    It is so called because the radical of the ideal $n\ZZ$ in the ring $\ZZ$ is equal to $\operatorname{rad}(n) \ZZ$.   It is clear that $\operatorname{rad}(n)$ is the largest squarefree divisor of $n$, and thus $\operatorname{rad}(n) \leq n$, with equality if and only if $n$ is squarefree.    It follows that $\dege  \operatorname{rad}(n) = (1,0,0,0,\ldots)$.    Since $\operatorname{rad}(n) \geq 1$ for all $n$ and $\operatorname{rad}(p^k) = p$ is bounded for any fixed prime $p$, one has $\underline{\dege}  \operatorname{rad}(n)  = (0,0,0,\ldots)$.    Note that $\operatorname{rad}(n)$  is a  multiplicative arithmetic function with  Dirichlet series
$$\sum _{n=1}^{\infty }{\frac {\operatorname {rad} (n)}{n^{s}}} =  \prod _{p}\left(1+{\frac {p^{1-s}}{1-p^{-s}}}\right)$$ having abscissa of convergence and absolute convergence equal to $2$.  Note also that
$$\frac{1}{x} \sum_{n \leq x} \operatorname{rad}(n) =  Cx + O(x^{1/2}\log x) \ (x \to \infty),$$
where 
$$C = \frac{1}{2} \prod_p \left(1-\frac{1}{p(p+1)} \right) = 0.352221\ldots,$$
and therefore $\operatorname{rad}(n)$ has average order $2Cn$.

The {\bf abc conjecture}\index{abc conjecture} states that, for every $t > 1$, there exist only finitely many triples $(a, b, c)$ of mutually prime positive integers such that  $a + b =  c>\operatorname {rad} (abc)^t$.   Let $\operatorname{ABC}(n)$ denote the maximal value of $c = a + b $ over all relatively prime positive integers $a, b$ such that $\operatorname{rad}(abc) = n$; that is, let
\[
\operatorname{ABC}(n) = \max\{a + b : a, b \in \mathbb{Z}_{>0},\ \gcd(a, b) = 1,\ \operatorname{rad}(ab(a+b)) = n\}.
\]
It is known that $\operatorname{ABC}(n)  < \infty$ for all $n$.  
Moreover, one has
\[ \dege \operatorname{ABC}  = \inf \left\{\dege r : \begin{array}{l}
r \in \mathbb{L} \text{ and only finitely many } (a, b, c) \in 
 \mathbb{Z}_{>0}^3 \\
 \text{satisfy } a + b = c > r(\operatorname{rad}(abc)) \text{ and} \gcd(a,b) = 1.
\end{array}\right\},
\]

\begin{proposition}
The abc conjecture is equivalent to $\deg \operatorname{ABC} = 1$ and thus also to 
\[
 \dege \operatorname{ABC}  \leq (1, \infty, 1, 0, 0, 0, \ldots).
\]
\end{proposition}

To date, the tightest known bounds on $\dege  \operatorname{ABC}$ are
\[
(1, 0, \infty, \tfrac{1}{2}, -1, 0, 0, 0, \ldots) \leq \dege \operatorname{ABC}   \leq (\infty, \tfrac{1}{3}, 3, 0, 0, 0, \ldots).
\]
Here, the lower bound follows from the fact that,  for all $k < 4$, there exist infinitely many  mutually prime triples $(a,b,c)$ such that $a+b = c$ and $$c>\operatorname {rad} (abc)\exp {\left(k{\sqrt {\log c}}/\log \log c\right)}$$
(Stewart and Tijdeman (1986)), and the upper bound follows from the existence of a $K>0$ such that
$$c<\exp {\left(K\operatorname {rad} (abc)^{\frac {1}{3}}\left(\log(\operatorname {rad} (abc)\right)^{3}\right)}$$
for all mutually prime triples  $(a,b,c)$ such that $a+b = c$ (Stewart and Yu (2001)).
A more precise conjecture of van Frankenhuysen (1995) for upper and lower bounds,  later refined by Robert,  Stewart and Tenenbaum (2014),  and supported by heuristics and numerical evidence,  implies that
\begin{align}\label{vf}
 \dege \operatorname{ABC} = (1, 0, \infty, \tfrac{1}{2}, -\tfrac{1}{2}, 0, 0, 0, \ldots).
\end{align}
This leads us to conjecture,  more modestly,  and in analogy with $\li(x)-\pi(x)$ and $M(x)$, the following.

\begin{conjecture}\label{abcextension}
One has $\deg \operatorname{ABC}(n) = 1$ and $\dege_1  \operatorname{ABC}(n) = 0$,  that is, 
$$\operatorname{ABC}(n) = O(n (\log n)^t) \ (n \to \infty)$$
if (and only if) $t > 0$.  Equivalently,   one has
$$\dege  \operatorname{ABC}(n) \leq (1, 0, \infty, 1, 0, 0, 0, \ldots).$$
or, equivalently still, for any $\varepsilon > 0$,  there exist only finitely many triples $(a,b,c)$ of mutually (or pairwise) relatively prime positive integers with $a+b = c > \operatorname{rad}(abc) \cdot( \log \operatorname{rad}(abc))^\varepsilon$.  
\end{conjecture}

Now, let
\begin{align*}
\underline{\operatorname{ABC}}(n) & =  \min\{ \operatorname{rad}(abc): a, b \in \mathbb{Z}_{>0},\ \gcd(a, b) = 1,\  c = a+b = n\}  \\ 
& \leq  \operatorname{rad}(n) \operatorname{rad}(n-1).
\end{align*}
This function is much more easily computed than $\operatorname{ABC}(n)$, since for any $n$ it is easy to enumerate all relatively prime pairs $(a,b)$ with $a+b = n$.   Two equivalent generalizations of the abc conjecture are to determine $\dd = \dege \operatorname{ABC}$ and
 $\dd' = \underline{\dege} \, \underline{\operatorname{ABC}}$, which, by Proposition \ref{abcp}(1) below, are related to one another as in Proposition \ref{thm:ledege_inverse}.   
 
\begin{proposition}\label{abcp}
One has the following.
\begin{enumerate}
\item One has
\[ \underline{\dege} \,  \underline{\operatorname{ABC}}  = \sup \left\{\dege r : \begin{array}{l}
r \in \mathbb{L} \text{ and only finitely many } (a, b) \in 
 \mathbb{Z}_{>0}^2 \\
 \text{satisfy } r(a+b)  > \operatorname{rad}(ab(a+b)) \text{ and} \gcd(a,b) = 1.
\end{array}\right\}
\]
\item  The abc conjecture is  equivalent to $\underline{\deg} \, \underline{\operatorname{ABC}}(n) = 1$.
\item  One has
\[
(-\infty, -\tfrac{1}{3}, -3, 0, 0, 0, \ldots)  \leq \underline{\dege} \, \underline{\operatorname{ABC}}   \leq (1, 0, -\infty, -\tfrac{1}{2}, 1, 0, 0, 0, \ldots),
\]
\item The conjecture (\ref{vf}), which is a consequence of a 1995 conjecture of van Frankenhuysen,  is equivalent to
$$\underline{\dege} \,  \underline{\operatorname{ABC}}(n)  = (1, 0, -\infty, -\tfrac{1}{2}, \tfrac{1}{2}, 0, 0, 0, \ldots).$$
\item   Conjecture \ref{abcextension} is equivalent to each of the following statements.
\begin{enumerate} 
\item $\underline{\deg}  \, \underline{\operatorname{{ABC}}}(n) = 1$ and $\underline{\dege}_1   \operatorname{\underline{{ABC}}}(n) = 0$.
\item  One has $$\underline{\operatorname{ABC}}(n) \gg n (\log n)^t \ (n \to \infty)$$
if (and only if) $t > 0$.
\item One has
$$\underline{\dege} \,  \underline{\operatorname{{ABC}}}(n) \geq (1, 0, -\infty, -1, 0, 0, 0, \ldots).$$
\item For all $\varepsilon > 0$,  there exist only finitely many pairs $(a,b)$ of relatively prime positive integers with $\frac{a+b}{(\log (a+b))^\varepsilon} > \operatorname{rad}(ab(a+b))$.   
\end{enumerate}
\end{enumerate}
\end{proposition}

Regarding $\dege \underline{\operatorname{ABC}}(n)$ and  $\underline{\dege }  \operatorname{ABC}(n)$, we make the following conjecture.

\begin{conjecture}\label{lowerABC}
There exist infinitely many positive integers $n$ such that $\underline{\operatorname{ABC}}(n) = n(n-1)$. 
\end{conjecture}

\begin{proposition}
Conjecture \ref{lowerABC} implies that $$\dege  \underline{\operatorname{ABC}}(n) = (2,0,0,0,\ldots)$$ and $$\underline{\dege}   \operatorname{ABC} = (\tfrac{1}{2},0,0,0,\ldots).$$
\end{proposition}

As evidence for Conjecture \ref{lowerABC}, note that it is a well-known classical result that the natural density of all integers $n$ such that both $n$ and $n-1$ are squarefree exists and is equal to
$$ \prod_p \left(1 - \frac{2}{p^2} \right) = 0.322634\ldots.$$
To guarantee  $\underline{\operatorname{ABC}}(n) = n(n-1)$,  one can impose a sole additional condition on $n$, namely, that $\operatorname{rad}(ab)$ be minimized over all relatively prime pairs $(a,b)$ such that $a+b = n$ by the pair $(1,n-1)$.   A simple search reveals that there are many such squarefree integers $n$, e.g., $n = 10542$ and $n = 20010$, leading us to conjecture that there are infinitely many such $n$, and thus that Conjecture \ref{lowerABC} holds.  However, without further investigation, it is not clear whether the density of such integers $n$ ought to (exist and) be positive or zero.

\section{Tables}

This section contains the tables discussed earlier in this chapter.

\begin{table}[h!]
\caption{\centering $f(x) = \li(x)-\pi(x)$ and $g(x) = \li(x)-\Ri(x)$ for $x = 10^n$ (https://oeis.org/A006880)}
\scriptsize
\begin{tabular}{|r|r|r|r|r|r|r|r|r|} \hline
$n$  & $\li(x)-\pi(x)$ & $\li(x)-\Ri(x)$ & $L_f(x)$ & $L_g(x)$ & $L_f[\frac{1}{2}](x)$ & $L_g[\frac{1}{2}](x)$ & $L_f[\frac{1}{2},\! \!-1](x)$ & $L_g[\frac{1}{2},\! \!-1](x)$ \\\hline \hline
1  &  2.165600 		&	1.601016    & 0.3356	&  0.2044	&  $-$0.4539	& $-$0.8161	&  $-$2.5095  & $-$0.8451  \\\hline
2  &   5.126142		&	4.464508  & 0.3549	&  0.3249	&  $-$0.4376 &  $-$0.5280	&  2.0286  & 1.7022 \\\hline
3  &  9.609658		&	9.250212   & 0.3276  &  0.3221	&  $-$0.6163	&  $-$0.6360 	&  1.1254   &1.0676  \\\hline
4  &  17.13722		&	19.20600  &	 0.3085 &  0.3209	&  $-$0.7944 	&  $-$0.7431	&  0.5722  & 0.7151 \\\hline
5  &  37.80900		&	42.37726  &	 0.3155 &  0.3254	&  $-$0.8692	&  $-$0.8225	&  0.3577  & 0.4853  \\\hline
6  &  129.5492		&	100.1497   &	 0.3521 &  0.3334	&  $-$0.7783	&  $-$0.8763	&  0.6030 & 0.3363  \\\hline
7  &  339.4050		&	250.9575  &	 0.3615 &  0.3428	&  $-$0.8028	&  $-$0.9114	& 0.5361   &  0.2408  \\\hline
8  &  754.3754		&	657.5081   &	 0.3597 &  0.3522	&  $-$0.8871	& $-$0.9342	& 0.3077   & 0.1792   \\\hline
9	 &	1700.957 &	1779.529 &	0.3590 &	0.3611 &	$-$0.9642 	& $-$0.9493 &	0.0051 &	0.0503 \\\hline
10	 &	3103.587 &	4931.280 &	0.3492 &	0.3693 &	$-$1.1071	& $-$0.9595 &	$-$0.2939 &	0.1111 \\\hline
11	&	11587.62 &	13905.99 &	0.3695 &	0.3767 &	$-$1.0231 &	$-$0.9667 &	$-$0.0636	& 0.0919 \\\hline
12	 &	38262.80 &	39738.58 &	0.3819 &	0.3833 &	$-$0.9832 &	$-$0.9718 &	0.0464 &	0.0779 \\\hline
13	 &	108971.05 &	114744.22 &	0.3875 &	0.3892 &	$-$0.9909 &	$-$0.9757 &	0.0254 &	0.0675 \\\hline
14	&	314889.95 &	334090.26 &	0.3927 &	0.3946 &	$-$0.9957 &	$-$0.9786 &	0.0120 &	0.0596 \\\hline
15	 &	1052618.58 &	979400.65 &	0.4015 &	0.3994 &	$-$0.9606 &	$-$0.9810 &	0.1103	& 0.0533 \\\hline
16	 &	3214631.79 &	2887579.61 &	0.4067 &	0.4038  &	$-$0.9531 &	$-$0.9828 &	0.1319 &	0.0482 \\\hline
17	 &	7956588.78 &	8554843.77 &	0.4059 &	0.4078	& $-$1.0041	& $-$0.9844 &	$-$0.0117 &	0.0441 \\\hline
18	 &	21949555.02 &	25450920.62 &	0.4079 &	0.4114	& $-$1.0254 &	$-$0.9857 &	$-$0.0719 &	0.0406  \\\hline
19		& 99877775.22 &	75993441.99 &	0.4210 &	0.4148 &	$-$0.9144 &	$-$0.9867 &	0.2433 &	0.0377 \\\hline
20	 &	222744643.55	 & 227636468.15 &	0.4174	& 0.4179 &	$-$0.9934 &	$-$0.9877 &	0.0190 &	0.0351  \\\hline
21	 &	597394254.33	& 683826458.60 &	0.4179 	& 0.4207 &	$-$1.0233 &	$-$0.9885 &	$-$0.0668 &	0.0329 \\\hline
22	 &	1932355208.15 &	2059487872.99 &	0.4221 &	0.4234 &	$-$1.0054 &	$-$0.9892 &	$-$0.0156 &	0.0310 \\\hline
23	 &	7250186215.78 &	6216886362.70 &	0.4287 &	0.4258  &	$-$0.9511  &	$-$0.9898 &	0.1408 &	0.0293  \\\hline
24		& 17146907278.15 &	18805896997.35	& 0.4264  &	0.4281 &	$-$1.0134  &	$-$0.9904  &	$-$0.0388 &	0.0277  \\\hline
25 &	55160980939.38	& 56995765653.55 &	0.4297 &	0.4302 	& $-$0.9990  &	$-$0.9909  &	0.0029  &	0.0263  \\\hline
26	 &  155891678120.79 &	173041013576.57 &	0.4305  &	0.4322 &	$-$1.0169 &	$-$0.9914 &	$-$0.0490 &	0.0251 \\\hline
27	&	508666658006.34 &	526202145940.76 &	0.4336 &	0.4341	& $-$1.0000 &	$-$0.9918 &	0.0000 &	0.0239 \\\hline
28  &  1427745660373.74  & 1602506180201.17   &  0.4341	& 	0.4359 & $-1.0199$	& $-0.9922$  &   $-0.0580$  & $0.0229$ \\\hline
29  &  4551193622464.34  & 4887005927157.18   &  0.4365	& 	0.4376 & $-1.0095$	& $-0.9925$  &   $-0.0277$  & $0.0219$ \\\hline
\end{tabular}\label{tab2}
\end{table}

\begin{table}[h!]
\caption{\centering Approximations of $f(x) = \li(x)-\pi(x)$ and $g(x) = \li(x)-\Ri(x)$ for $x = 10^n$ \cite{plan}}
\scriptsize 
\begin{tabular}{|r|r|r|r|r|r|r|r|r|} \hline 
$n$  & $\li(x)-\pi(x) \approx$   & $\li(x)-\Ri(x) \approx$ &  \ \ \ $L_f(x)$  & $L_g(x)$  & $L_{f}[\frac{1}{2}](x)$ &  $L_{g}[\frac{1}{2}](x)$  & $L_{f}[\frac{1}{2},\! \! -1]( x)$  & $L_{f}[\frac{1}{2},\! \! -1]( x)$ \\  \hline \hline

${26}$ & $1.5603 \cdot 10^{11}$  &  $1.7304\cdot 10^{11}$  & 0.4305 & 0.4322 & $-1.0167$ & $-$0.9914 &  $-0.0484$ &   0.0251 \\ \hline
${27}$ & $5.0811 \cdot 10^{11}$  &  $ 5.2620 \cdot 10^{11}$  & 0.4336 & 0.4341 & $-1.0003$ &$-$0.9918 & $-0.0008$ &  0.0239 \\ \hline
${28}$ &  $1.4271 \cdot 10^{12}$ &  $1.6025 \cdot 10^{12}$   & 0.4341  &  0.4359  &  $-$1.0200  &   $-$0.9922   &  $-$0.0584  & 0.0229  \\ \hline 
${29}$ &  $4.5533  \cdot 10^{12}$ &  $4.8870 \cdot 10^{12}$  &  0.4365   &  0.4376  &  $-$1.0093  &   $-$0.9925   &  $-$0.0274   &   0.0219  \\\hline 
${30}$ & $1.8779 \cdot 10^{13}$ &  $1.4922 \cdot 10^{13}$  &  0.4425  &  0.4391 &  $-$0.9386  &   $-$0.9928   &     0.1803   &   0.0210  \\\hline 
${31}$ &  $3.9314 \cdot 10^{13}$ &  $4.5620 \cdot 10^{13}$  &    0.4385  & 0.4406   &  $-$1.0280  &   $-$0.9931   &   $-$0.0823  &  0.0202  \\ \hline 
${32}$ & $9.0918 \cdot 10^{13}$  &   $1.3962 \cdot 10^{14}$  &  0.4362  & 0.4420   &  $-$1.0932   &   $-$0.9934   &  $-$0.2747  &   0.0194  \\\hline 
${33}$ &  $3.9709 \cdot 10^{14}$ &   $4.2775 \cdot 10^{14}$  &  0.4424  &  0.4434   &  $-$1.0108   &   $-$0.9937   &  $-$0.0320   &   0.0187  \\\hline 
${34}$ & $1.9605 \cdot 10^{15}$  &   $1.3118 \cdot 10^{15}$  &  0.4498   &  0.4446   &   $-$0.9017   &   $-$0.9939   &  0.2909   &   0.0181  \\\hline 
${35}$ &  $4.8963 \cdot 10^{15}$ &    $4.0265 \cdot 10^{15}$ &  0.4483   &  0.4459   &  $-$0.9496   &   $-$0.9941   &     0.1497  &   0.0175  \\\hline 
${36}$ & $7.4445 \cdot 10^{15}$  &   $ 1.2370 \cdot 10^{16}$ &  0.4409   &  0.4470  &   $-$1.1093  &   $-$0.9943   &   $-$0.3249    &   0.0169  \\\hline 
${37}$ &  $2.2485 \cdot 10^{16}$ &   $3.8033 \cdot 10^{16}$  & 0.4419   &  0.4481   &  $-$1.1128   &   $-$0.9945   &   $-$0.3360  &   0.0163  \\\hline 
${38}$ &  $1.1575 \cdot 10^{17}$ &    $1.1703 \cdot 10^{17}$  &  0.4490    &  0.4492   &  $-$0.9972   &   $-$0.9947  &   0.0085  &   0.0158  \\\hline 
${39}$ &  $3.6121 \cdot 10^{17}$ &   $3.6036 \cdot 10^{17}$  &  0.4502   &  0.4502   &  $-$0.9943    &   $-$0.9949   &  0.0169  &   0.0153  \\\hline 
${40}$ & $9.5267 \cdot 10^{17}$ &   $1.1104 \cdot 10^{18}$ &   0.4494   &  0.4511   &  $-$1.0291   &   $-$0.9950   &  $-$0.0866     &   0.0149  \\\hline 
${41}$ & $ 2.8478 \cdot 10^{18}$ &   $3.4238 \cdot 10^{18}$   &  0.4501    &  0.4521   &  $-$1.0359   &   $-$0.9952   &  $-$0.1072  &   0.0145  \\\hline 
${42}$ & $ 9.0106 \cdot 10^{18}$&    $ 1.0564 \cdot 10^{19}$  &  0.4513  &  0.4529   &   $-$1.0301   &   $-$0.9953   &   $-$0.0906   &   0.0141  \\\hline 
${43}$ & $3.3434 \cdot 10^{19}$  &    $3.2612 \cdot 10^{19}$ &   0.4541  &  0.4538   &  $-$0.9900    &   $-$0.9955   &   0.0300  &   0.0137  \\\hline 
${44}$ & $9.1921 \cdot 10^{19}$ &   $1.0073 \cdot 10^{20}$ &   0.4537   &  0.4546   &  $-$1.0154   &   $-$0.9956   &   $-$0.0465   &   0.0133  \\\hline 
${45}$ & $2.2524 \cdot 10^{20}$ &    $3.1132 \cdot 10^{20}$  &  0.4523  &  0.4554   &  $-$1.0655   &   $-$0.9957   &  $-$0.1979  &   0.0130  \\\hline 
${46}$ & $1.0325 \cdot 10^{21}$ &   $9.6266 \cdot 10^{20}$  &  0.4568   &  0.4562   &  $-$0.9808   &   $-$0.9958   &    0.0581    &   0.0126  \\\hline 
${47}$ & $3.6107 \cdot 10^{21}$ &   $2.9782 \cdot 10^{21}$ &   0.4587  &  0.4569   &  $-$0.9548   &   $-$0.9959   &   0.1370  &  0.0123  \\\hline 
${48}$ & $6.5485 \cdot 10^{21}$ &   $9.2178 \cdot 10^{21}$  &  $0.4545$   &  0.4576  &  $-$1.0687    &   $-$0.9960   &  $-$0.2088   &  0.0120  \\ \hline 
${49}$ & $2.2397 \cdot 10^{22}$ &   $2.8543 \cdot 10^{22}$ &  $0.4561$    &  0.4583   &   $-$1.0475   &   $-$0.9961 &  $-$0.1444   &   0.0117  \\\hline 
${50}$ & $8.2123 \cdot 10^{22}$ &    $8.8423 \cdot 10^{22}$  &  $0.4583$  &  0.4589   &  $-$1.0118   &   $-$0.9962 &   $-$0.0360   &  0.0115 \\ \hline
\end{tabular}\label{tab3}
\end{table}

\begin{table}[ht!]
\caption{\centering $f(x) = \li(x)-\pi(x)$ and $g(x) = \li(x)-\Ri(x)$ for $x = e^n$ (https://oeis.org/A040014)}
\tiny
\begin{tabular}{|r|r|r|r|r|r|r|r|r|}\hline
$n$  & $\li(x)-\pi(x)$ & $\li(x)-\Ri(x)$ & $L_f(x)$ & $L_g(x)$ & $L_f[\frac{1}{2}](x)$ & $L_g[\frac{1}{2}](x)$ & $L_f[\frac{1}{2},\! \!-1](x)$ & $L_g[\frac{1}{2},\! \!-1](x)$ \\\hline \hline
1    & 0.895118      & 0.015939     & $-$0.1108 & $-$4.1390 & ---   & ---   & ---      & ---     \\ \hline
2   & 0.954234     & 1.263390      & -0.0234 & 0.1168  & $-$1.5103  & $-$1.1054 & 0.9650   & 0.1993 \\ \hline
3 & 1.933833     & 2.382997      & 0.2198 & 0.2895  & $-$0.7650 & $-$0.5749  & 2.7445   & 4.9653 \\ \hline
4   & 3.630874      & 3.608914      & 0.3224  & 0.3209  & $-$0.5125 & $-$0.5169 & 2.0689  & 2.0503 \\ \hline
5   & 6.185275       & 5.088573      &  0.3644  & 0.3254  & $-$0.4212 & $-$0.5424 & 1.9576  & 1.5475 \\ \hline
6  & 6.989762      & 6.984925      & 0.3241  & 0.3240  & $-$0.5891 & $-$0.5895 & 1.2624  & 1.2612 \\ \hline
7    & 8.504743     & 9.518112      & 0.3058  & 0.3219  & $-$0.6986 & $-$0.6407 & 0.8810 & 1.0501 \\ \hline
8   & 11.37990     & 13.00712     & 0.3040  & 0.3207   & $-$0.7541 & $-$0.6899 & 0.6984  & 0.8809 \\ \hline
9   & 18.87829    & 17.92737    & 0.3264  & 0.3207  & $-$0.7109 & $-$0.7344 & 0.8070  & 0.7413 \\ \hline
10 & 26.22898      & 24.99738     & 0.3267  & 0.3219  & $-$0.7527 & $-$0.7736  & 0.6828  & 0.6251 \\ \hline
11  & 23.40637     & 35.31241     & 0.2866  & 0.3240  & $-$0.9788  & $-$0.8073  & 0.0582  & 0.5284 \\ \hline
12  & 47.53267     & 50.55248    & 0.3218  & 0.3269 & $-$0.8606  & $-$0.8358 & 0.3805  & 0.4481 \\ \hline
13  & 69.68849     & 73.30739     & 0.3265  & 0.3304  & $-$0.8795 & $-$0.8598 & 0.3280  & 0.3818 \\ \hline
14  & 75.51363    & 107.58650   & 0.3089  & 0.3342  & $-$1.0139 & $-$0.8797 & $-$0.0377 & 0.3270 \\ \hline
15  & 100.8525   & 159.6199   & 0.3076  & 0.3382  & $-$1.0658  & $-$0.8963 & $-$0.1790 & 0.2819 \\ \hline
16  & 219.9987    & 239.1208    & 0.3371  & 0.3423 & $-$0.9401 & $-$0.9100 & 0.1630  & 0.2447 \\ \hline
17 & 404.8940     & 361.2780     & 0.3532   & 0.3464  & $-$0.8811  & $-$0.9213 & 0.3234  & 0.2140 \\ \hline
18  & 718.3306    & 549.9068    & 0.3654  & 0.3505  & $-$0.8383 & $-$0.9308 & 0.4403   & 0.1889 \\ \hline
19  & 561.2510   & 842.4401   & 0.3332  & 0.3545  & $-$1.0766  & $-$0.9386 & $-$0.2087  & 0.1674 \\ \hline
20  & 1090.664   & 1297.847   & 0.3497  & 0.3584  & $-$1.0032  & $-$0.9452 & $-$0.0089 & 0.1496 \\ \hline
21  & 2409.356  & 2009.211  & 0.3708  & 0.3622  & $-$0.8911 & $-$0.9507 & 0.2979   & 0.1347  \\ \hline
22 & 2774.300   & 3123.742   & 0.3604  & 0.3658  & $-$0.9938 & $-$0.9554 & 0.0170  & 0.1221 \\ \hline
23  & 2826.830   & 4874.654   & 0.3455  & 0.3692  & $-$1.1332 & $-$0.9594 & $-$0.3654  & 0.1114 \\ \hline
24 & 8547.849  & 7631.985  & 0.3772  & 0.3725  & $-$0.9272 & $-$0.9628 & 0.2002  & 0.1022 \\ \hline
25  & 14041.53       &  11983.74  & 0.3820  & 0.3757  & $-$0.9165 & $-$0.9658 & 0.2298  & 0.0942   \\ \hline
26  & 19730.90      &   18865.54 & 0.3804  & 0.3787  & $-$0.9546 & $-$0.9683   & 0.1253  & 0.0874 \\ \hline
27  & 26318.88       & 29767.99  & 0.3770  & 0.3815  & $-$1.0079 & $-$0.9706 & $-$0.0219 & 0.0814 \\ \hline
28  & 40913.31       & 47068.67 & 0.3793  & 0.3843  & $-$1.0146 & $-$0.9725 & $-$0.0404 & 0.0761 \\ \hline
29  & 81089.41       & 74564.04 & 0.3898  & 0.3869  & $-$0.9493 & $-$0.9743 & 0.1405 & 0.0714 \\ \hline
30   & 134842.27      & 118322.37 & 0.3937 & 0.3894  & $-$0.9374 & $-$0.9758 & 0.1740  & 0.0673 \\ \hline
31  & 184278.39      & 188052.43 & 0.3911 & 0.3918  & $-$0.9831 & $-$0.9772 & 0.0472  & 0.0636 \\ \hline
32  & 310474.79      & 299301.41 & 0.3952  & 0.3940  & $-$0.9678 & $-$0.9784 & 0.0898  & 0.0603  \\ \hline
33   & 462754.50      & 476986.81 & 0.3953 & 0.3962  & $-$0.9881 & $-$0.9795 & 0.0331  & 0.0573  \\ \hline
34   & 902506.27     & 761075.61 & 0.4033 & 0.3983   & $-$0.9321 & $-$0.9805 & 0.1899  & 0.0546 \\ \hline
35  & 1483541.60     & 1215721.63 & 0.4060  & 0.4003  & $-$0.9253 & $-$0.9814  & 0.2091 & 0.0522 \\ \hline
36  & 2551387.77     & 1943975.49 & 0.4098  & 0.4022  & $-$0.9063 & $-$0.9822 & 0.2630  & 0.0500 \\ \hline
37  & 4051058.40     & 3111478.72 & 0.4112  & 0.4041    & $-$0.9099 & $-$0.9830 & 0.2534  & 0.0479 \\ \hline
38   & 5477912.65     & 4984649.34 & 0.4083  & 0.4058  & $-$0.9577 & $-$0.9837 & 0.1191 & 0.0460\\ \hline
39  & 8227929.60     & 7992252.66    & 0.4083  & 0.4075  & $-$0.9764 & $-$0.9843 & 0.0667 & 0.0443 \\ \hline
40   & 12547993.58    & 12824726.72    & 0.4086 & 0.4092  & $-$0.9908 & $-$0.9849 & 0.0260  & 0.0427  \\ \hline
41  & 23321012.11   & 20594486.91    & 0.4138  & 0.4107  & $-$0.9519 & $-$0.9854 & 0.1360 & 0.0412 \\ \hline
42   & 38092380.59    & 33094760.86    & 0.4156  & 0.4123  & $-$0.9483 & $-$0.9859 & 0.1465  & 0.0399 \\ \hline
43  & 61485600.03    & 53217692.13    & 0.4171   & 0.4137  & $-$0.9480 & $-$0.9864  & 0.1476 & 0.0386 \\ \hline
44  & 108071721.51  & 85629989.66    & 0.4204  & 0.4151  & $-$0.9253 & $-$0.9869  & 0.2123 & 0.0374 \\ \hline
45   & 174585516.15   & 137865194.52   & 0.4217  & 0.4165  & $-$0.9252 & $-$0.9873 & 0.2129   & 0.0362 \\ \hline
46    & 204498265.00   & 222090120.72   & 0.4160  & 0.4178  & $-$1.0092  & $-$0.9877  & $-$0.0263 & 0.0352  \\ \hline
47   & 385623847.54  & 357962620.61 & 0.4206  & 0.4191   & $-$0.9687 & $-$0.9880 & 0.0894  & 0.0342 \\ \hline
48  & 651043126.55  & 577256516.55   & 0.4228  & 0.4203  & $-$0.9573 & $-$0.9884 & 0.1221 & 0.0332  \\ \hline
49   & 1164251352.07 & 931348981.39   & 0.4260   & 0.4215  & $-$0.9314 & $-$0.9887 & 0.1966 & 0.0324 \\ \hline
50    & 2005203489.62 & 1503345912.12  & 0.4284 & 0.4226    & $-$0.9154 & $-$0.9890 & 0.2427  & 0.0315 \\ \hline
51     & 2068331099.91  & 2427724401.02  & 0.4206  & 0.4237   & $-$1.0301  & $-$0.9893 & $-$0.0863  & 0.0307 \\ \hline
52   & 4005009209.52  & 3922162119.58  & 0.4252 & 0.4248  & $-$0.9843 & $-$0.9896 & 0.0452  & 0.0300 \\ \hline
53   & 6781992345.02 & 6339129875.41  & 0.4271  & 0.4258   & $-$0.9728 & $-$0.9898 & 0.0782  & 0.0292 \\ \hline
54    & 9918796068.45  & 10249542084.52 & 0.4263  & 0.4269  & $-$0.9983 & $-$0.9901 & 0.0048  & 0.0285 \\ \hline
55   & 20087997222.46 & 16578421917.45   & 0.4313  & 0.4278  & $-$0.9424 & $-$0.9903 & 0.1662  & 0.0279 \\ \hline
56 & 25498675636.03 & 26824982492.19 & 0.4279   & 0.4288  & $-$1.0032 & $-$0.9906 & $-$0.0092 & 0.0273 \\ \hline
57  & 47560462482.65 & 43419742391.88  & 0.4313  & 0.4297  & $-$0.9683 & $-$0.9908 & 0.0919  & 0.0267  \\ \hline
58   & 71057460397.51 & 70304144135.27 & 0.4308  & 0.4306  & $-$0.9884 & $-$0.9910 & 0.0337  & 0.0261 \\ \hline
59   & 94758718695.13 & 113871533611.46 & 0.4284  & 0.4315  & $-$1.0362 & $-$0.9912 & $-$0.1052 & 0.0255 \\ \hline
\end{tabular}\label{tab4}
\end{table}

\begin{table}[h!]
\caption{\centering Sup of $\Delta(x) = \left(\pi_0(x)- \Ri(x)  +\frac{1}{\log x} - \frac{1}{\pi} \arctan \frac{\pi}{\log x}\right) \frac{\log x}{\sqrt{x}}$ on $I_n = [10^n,10^{n+1}]$ \cite{kul}}
\scriptsize
\begin{tabular}{|r|r|r|r|r|r|r|} \hline
$n$  & $b_n$ & $\sup_{I_n} \Delta = \Delta(b_n^+)$ & $L_\Delta(b_n^+)$   & $L_\Delta[0](b_n^+)$ &  $L_\Delta[0, \! 0](b_n^+)$  \\\hline \hline
$0$ 		&	1			&1.0000000000		&	$-2.0000$	&	$0.0000$	& ---	\\ \hline
$1$ 		&	19			&0.5607597113		&	$-$0.1965	&	 $-$0.5357	&	$-$7.5237	\\ \hline
$2$ 		&	113			&0.7848341482		&	$-$0.05125	&	$-$0.1560	&	$-$0.5501	\\ \hline
$3$ 		&	1627			&0.6754517455		&	$-$0.05306	&	$-$0.1961	&	$-$0.5658	\\ \hline
$4$ 		&	24137		&0.7457431860		&	$-$0.02907	&	$-$0.1269	&	$-$0.3501	\\ \hline
$5$ 		&	355111		&0.7008073861		&	$-$0.02782	&	$-$0.1395	&	$-$0.3801	\\ \hline
$6$ 		&	3445943		&0.6809987397		&	$-$0.02552	&	$-$0.1417	&	$-$0.3851	\\ \hline
$7$ 		&	30909673	&0.7157292126		&	$-$0.01939	&	$-$0.1174	&	$-$0.3196	\\ \hline
$8$ 		&	110102617	&0.7878100197		&	$-$0.01288	&	$-$0.08171	&	$-$0.2227	\\ \hline
$9$ 		&	1110072773	&0.6833192028		&	$-$0.01828	&	$-$0.1254	&	$-$0.3429	\\ \hline
$10$ 	&	10016844407	&0.6386706267		&	$-$0.01947	&	$-$0.1429	&	$-$0.3922	\\ \hline
$11$ 	&	330957852107		&0.7533813432	&	$-$0.01068	&	$-$0.08639 	&	$-$0.2385		\\ \hline
$12$ 	&	2047388353069		&0.6808028098	&	$-$0.01356	&	$-$0.1150	&	$-$0.3185	\\ \hline
$13$ 	&	21105695997889		&0.6896466780	&	$-$0.01211	&	$-$0.1085	&	$-$0.3019	\\ \hline
$14$ 	&	117396942462053		&0.6789107425	&	$-$0.01195	&	$-$0.1113	&	$-$0.3107	\\ \hline
$15$ 	&	1047930291039067	&0.7042622330	&	$-$0.01014	&	$-$0.09894	&	$-$0.2771	\\ \hline
$16$ 	&	16452596773450399	&0.7144542025	&	$-$0.009005	&	$-$0.09288	&	$-$0.2614	\\ \hline
$17$ 	&	125546149553907317	&0.6572554320	&	$-$0.01066	&	$-$0.1143	&	$-0.3226$	\\ \hline
$18$ 	&	1325005986250807813		&0.7839983342	&	$-$0.005832	&	$-$0.06522	&	$-$0.1848	\\ \hline
$19$ 	&	11538454954199984761	&0.7574646817	&	$-$0.006329	&	$-$0.07346	&	$-$0.2088	\\ \hline
\end{tabular}\label{tab5}
\end{table}

\begin{table}[h!]
\caption{\centering Inf of $\Delta(x) = \left(\pi_0(x)- \Ri(x)  +\frac{1}{\log x} - \frac{1}{\pi} \arctan \frac{\pi}{\log x}\right) \frac{\log x}{\sqrt{x}}$ on $I_n = [10^n,10^{n+1}]$ \cite{kul}}
\scriptsize
\begin{tabular}{|r|r|r|r|r|r|r|} \hline
$n$  & $c_n$ & $\inf_{I_n} \Delta = \Delta(c_n^-)$ & $L_\Delta(c_n^-)$   & $L_\Delta[0](c_n^-)$ &  $L_\Delta[0, \! 0](c_n^-)$  \\\hline \hline
$0$ 		&	5			& $-$0.3952461978		&	$-$0.5768	&	$-$1.9506	&	1.2500	\\ \hline
$1$ 		&	11			& $-$0.5492343329		&	$-$0.2499	&	$-$0.6852	&	4.4719	\\ \hline
$2$ 		&	223			& $-$0.6051733874		&	$-$0.09289	&	$-$0.2976	&	$-$0.9596	\\ \hline
$3$ 		&	1423			& $-$0.7542604400		&	$-$0.03884	&	$-$0.1423	&	$-$0.4121	\\ \hline
$4$ 		&	19373		& $-$0.7278356754		&	$-$0.03218	&	$-$0.1387	&	$-$0.3835	\\ \hline
$5$ 		&	302831		& $-$0.6995719492		&	$-$0.02831	&	$-$0.1409	&	$-$0.3840	\\ \hline
$6$ 		&	1090697		& $-$0.6389660809		&	$-$0.03222	&	$-$0.1702	&	$-$0.4628	\\ \hline
$7$ 		&	36917099	& $-$0.7489165055		&	$-$0.01659	&	$-$0.1012	&	$-$0.2753	\\ \hline
$8$ 		&	516128797	& $-$0.6775687236		&	$-$0.01940	&	$-$0.1298	&	$-$0.3544	\\ \hline
$9$ 		&	7712599823	& $-$0.6889577485		&	$-$0.01637	&	$-$0.1192	&	$-$0.3270	\\ \hline
$10$ 	&	11467849447	& $-$0.7251609705		&	$-$0.01387	&	$-$0.1023	&	$-$0.2807	\\ \hline
$11$ 	&	110486344211	& $-$0.7355462679		&	$-$0.01208	&	$-$0.09492	&	$-$0.2616	\\ \hline
$12$ 	&	1635820377397		& $-$0.6892596608		&	$-$0.01323	&	$-$0.1115	&	$-$0.3088		\\ \hline
$13$ 	&	36219717668609		& $-$0.8360329846	&	$-$0.005736	&	$-$0.05204	&$-$0.1449		\\ \hline
$14$ 	&	348323506633621		& $-$0.6494959371	&	$-$0.01289	&	$-$0.1229	&$-$0.3436		\\ \hline
$15$ 	&	1212562524413153		& $-$0.7750460589	&	$-$0.007337	&	$-$0.07183	&	$-$0.2012	\\ \hline
$16$ 	&	18019655286689201	& $-$0.5710665212	&	$-$0.01497	&	$-$0.1547	&	$-$0.4353 \\ \hline
$17$ 	&	266175790131587543	& $-$0.7599282036	&	$-$0.006842	& $-$0.07436	&$-$0.2102		\\ \hline
$18$ 	&	5805523423155128399	& $-$0.6804259482	&	$-$0.008912	&	$-$0.1022 &	$-$0.2904	\\ \hline
$19$ 	&	55496658217283199013	& $-$0.8042730098		&	$-$0.004791	&	$-$0.05707	&	$-$0.1626	\\ \hline
\end{tabular}\label{tab6}
\end{table}

\begin{table}[h!]
\centering
\caption{\centering Sup of $|\Delta(x)|$  and upper bounds of $L_\Delta(x) $, $L_\Delta[0](x)$, and  $L_\Delta[0, \!0](x)$ on $I_n = [10^n,10^{n+1}]$ \cite{kul}}
\scriptsize
\bigskip
\begin{tabular}{|r|r|r|r|r|r|r|} \hline
$n$  & $a_n$ & $\sup_{I_n} |\Delta| = |\Delta(a_n)|$ & $L_\Delta(x) \leq$   & $L_\Delta[0](x) \leq$ &  $L_\Delta[0, \! 0](x) \leq$  \\\hline \hline
$0$ 		&	$1^+$			&1.0000000000		&	$-0.3906$	&	$\infty$	&	$\infty$	\\ \hline
$1$ 		&	$19^+$			&0.5607597113		&	$-$0.1618	 
&	 $-$0.4710 
	&	$\infty$\\ \hline
$2$ 		&	$113^+$			&0.7848341482		&	$-$0.05125 
&	-$0.1559$
 &	$-$0.5501  
	\\ \hline
$3$ 		&	$1423^-$		& 0.7542604400 	&	$-$0.03884	
&	$-$0.1422		
&	$-$0.4121	
	\\ \hline
$4$ 		&	$24137^+	$	&0.7457431860		
&	$-$0.02907	&	$-$0.1269	&	$-$0.3500 \\ \hline
$5$ 		&	$355111^+$	&0.7008073861		&	$-$0.02267	&	$-$0.1294	&	$-$0.3504	\\ \hline
$6$ 		&	$3445943^+$		&0.6809987397		&	$-$0.02155	&	$-$0.1335		&	$-$0.3629	\\ \hline
$7$ 		&	$36917099^-$	& 0.7489165055		&	$-$0.01435	&	$-$0.09640	&	$-$0.2629	\\ \hline
$8$ 		&	$110102617^+$	&0.7878100197		&	$-$0.01144	&	$-$0.07853	&	$-$0.2146	\\ \hline
$9$ 		&	$7712599823^-$	& 0.6889577485		&	$-$0.01456	&	$-$0.1150	&	$-$0.3169	\\ \hline
 $10$ 	&	$11467849447^-$	& 0.7251609705		&	$-$0.01260	&	$-$0.09924	&	$-$0.2734	\\ \hline
$11$ 	&	$330957852107^+$		&0.7533813432	&	$-$0.009744	&	$-$0.08408 	&	$-$0.2331		\\ \hline
$12$ 	&	$1635820377397^-$		& 0.6892596608		&	$-$0.01219	&	$-$0.1088	&	$-$0.3027		\\ \hline
$13$ 	&	$36219717668609^-$		& 0.8360329846	&	$-$0.005308	&	$-$0.05091	&$-$0.1423		\\ \hline
$14$ 	&	$117396942462053^+$		&0.6789107425	&	$-$0.01115	&	$-$0.1091	&	$-$0.3058	\\ \hline
$15$ 	&	$1212562524413153	^-$	& 0.7750460589	&	$-$0.006876	&	$-$0.07054	&	$-$0.1984	\\ \hline
$16$ 	&	$16452596773450399^+$	&0.7144542025	&	$-$0.008468	&	$-$0.09134	&	$-$0.2579	\\ \hline
$17$ 	&	$266175790131587543^-$	& 0.7599282036	&	$-$0.006452	& $-$0.07321	&$-$0.2076		\\ \hline
$18$ 	&	$1325005986250807813^+$		&0.7839983342	&	$-$0.005522	&	$-$0.06428	&	$-$0.1828	\\ \hline
$19$ 	&	$55496658217283199013^-$	&  0.8042730098		&	$-$0.004542	&	$-$0.05629	&	$-$0.1609	\\ \hline
\end{tabular}\label{tab7}
\end{table}

\begin{table}[h!]
\caption{\centering Prominent peak values of $\sum_{k = 1}^N \frac{2 \cos(\gamma_k \log x - \operatorname{Arg} \rho_k)}{|\rho_k|} \approx V(x) = \frac{ \Ri(x)-\pi(x) }{\li(x)-\Ri(x)}$  for $10^{20} <  x < 10^{10^{13}}$  \cite[Table 3]{stoll}}
\scriptsize
\begin{tabular}{|r|r|r|r|r|r|r|} \hline
$x$ & $V(x) \approx $ & $L_V(x)\approx$   & $L_V[0](x) \approx$ &  $L_V[0,\! 0](x)\approx$ &  $L_V[0,\! 0,\! 0](x)\approx $  \\\hline \hline
$2.08963594805312 \cdot 10^{31}$	&	0.85111265	&	$-$0.00224		&	 $-$0.0377	& 	$-$0.111	&  $-$0.431   \\ \hline
$4.489584432 \cdot 10^{41}$			&	$-$0.92040123	&	$-$0.000865	&	$-$0.0182	& 	$-$0.0546	&    $-$0.199  \\ \hline
$5.38952511422773 \cdot 10^{84}$		&	$-$0.93251498	&	$-$0.000358	&	$-$0.0132	& 	$-$0.0420	&   $-$0.137  \\ \hline
$1.25855876915892 \cdot 10^{179}$	&	$-$0.95469600	&	$-$0.000112		&	 $-$0.00770	& 	$-$0.0258	&  $-$0.0792   \\ \hline
$1.33793033022140 \cdot 10^{190}$	&	$-$0.96875484	&	$-$0.0000725	&	$-$0.00522	& 	$-$0.0176	&    $-$0.0537  \\ \hline
$1.9098756608800 \cdot 10^{215}$	&	0.99896402	&	$-2.091\cdot10^{-6}$	& $-$0.000167		& 	 $-$0.000568	&   $-$0.00172   \\ \hline
$1.397162914 \cdot 10^{316}$			&	$-$1.000	&	0	&	0	& 	0	&    0 	 \\ \hline
$1.397166708 \cdot 10^{316}$			&	$-$1.002	&	$2.74\cdot 10^{-6}$		&	 0.000303	& 	0.00106	&    0.00315 \\ \hline
$1.398215 \cdot 10^{316}$			&	$-$1.005	&	$6.85\cdot 10^{-6}$	&	 0.000757	& 		0.00265 &   0.00786  \\ \hline
$1.39822803460925 \cdot 10^{316}$	&	$-$1.01770266	&	0.0000241		&	0.0026627	& 	0.00931	&   0.0277  \\ \hline
$6.5769904020884 \cdot 10^{370}$	&	$-$1.04965574	&	0.0000568	&	0.00718	& 	0.0254	&    0.0746  \\ \hline
$2.56771067186843 \cdot 10^{807}$	&	1.04307053	&	0.0000227	&	0.00560	& 	0.0209	&    0.0600 \\ \hline
$1.592776 \cdot 10^{1165}$			&	$-$1.042	&	0.0000153	&	0.00521	& 	0.0199	&	0.0567     \\ \hline
$2.76934778827004 \cdot 10^{4943}$		&	$-$1.11199752	&	$9.33 \cdot 10^{-6}$	&	0.0114	& 	0.0475	&  0.132   \\ \hline
$1.303324322225356 \cdot 10^{651157}$	&	$-$1.26421901	&	$1.56 \cdot 10^{-7}$	&	0.0165	& 	0.0883		&  0.240   \\ \hline
$1.1952052716413331 \cdot 10^{1748085}$	&	$-$1.28300855	&	$6.19 \cdot 10^{-8}$	&	0.0164	& 		0.0916	&   0.249  \\ \hline
$1.382217859037819 \cdot 10^{6359808}$	&	$-$1.2851	&	$1.71 \cdot 10^{-8}$		&	0.0152	& 	0.0895	&   0.243  \\ \hline
$4.63188213923253 \cdot 10^{25462014}$	&	1.30691350	&	$4.57 \cdot 10^{-9}$	&	0.0150	& 	0.0928	&   0.253  \\ \hline
$6.5722 \cdot 10^{30802655}$			&	$-$1.332	&	$4.04 \cdot 10^{-9}$	&	0.0159	& 	0.0990	&  0.270   \\ \hline
$1.49741610590218 \cdot 10^{140451912}$	&	1.34051829	&	$9.06 \cdot 10^{-10}$	&	0.0150	& 		0.0985  &  0.269   \\ \hline
$1.14937434741219 \cdot 10^{263854435}$	&	1.36190288	&	$5.08 \cdot 10^{-10}$	&	0.0153	& 	0.103	&   0.281  \\ \hline
$3.3380935222734 \cdot 10^{1048348162}$	&	$-$1.36843786	&	$1.30 \cdot 10^{-10}$	&	0.0145	& 		0.102	&  0.279   \\ \hline
$5.60458099779389 \cdot 10^{1715984906}$		&	1.41803522	&	$8.84 \cdot 10^{-11}$	&	0.0158	& 	0.113	&    0.309 \\ \hline
$4.66560656850816 \cdot 10^{3926214944}$	&	1.45919474	&	$4.18 \cdot 10^{-11}$		&	0.0165	& 	0.121	&   0.331  \\ \hline
$5.51920707045967 \cdot 10^{16251858572}$	&	$-$1.50218678	&	$1.09 \cdot 10^{-11}$	&	0.0167	& 	0.127	&  0.351   \\ \hline
$7.48642968663951 \cdot 10^{1741296446971}$	&	1.51320794	&	$1.03 \cdot 10^{-13}$	&	0.0142	& 	0.123	&   0.341  \\ \hline
$1.26735548862278 \cdot 10^{2176190682786}$	&	$-$1.51154647	&	$8.24 \cdot 10^{-14}$	&	0.0141	& 	0.122	&    0.340 \\ \hline
$4.9055169754782 \cdot 10^{3648403832739}$	&	$-$1.55246992	&	$5.24 \cdot 10^{-14}$	&	0.0148	& 	 0.130	&   0.360  \\ \hline
$2.9778893050643 \cdot 10^{4026762218841}$	&	$-$1.51359970	&	$4.47 \cdot 10^{-14}$	&	0.0139 & 	0.122	&    0.339 \\ \hline
$2.30170547097451 \cdot 10^{7424180200481}$	&	$-$1.51177469	&	$2.42 \cdot 10^{-14}$	&	0.0136& 	0.121	&   0.336  \\ \hline
\end{tabular}\label{tab8}
\end{table}

\begin{table}[h!]
\caption{\centering Some prominent peak values of $M(n)$ \cite{hur} \cite{kuz}} 
\scriptsize 
\begin{tabular}{|r|r|r|r|r|r|} \hline 
$n$  & $M(n)$  & $\frac{M(n)}{\sqrt{n}}$  & \ \  $L_M(n)$  &  $L_{M}[\frac{1}{2}](n)$ & $L_{M}[\frac{1}{2},\!  0](n)$  \\  \hline \hline
179919749		&	$-$6226			&	$-$0.4642			&	0.4596		&	$-$0.2606		&  $-$0.7106	\\ \hline
6631245058		&	$-$31206	&	$-$0.3832	&	0.4576	&	$-$0.3076	     &  $-$0.8433   \\ \hline
7766842813		&	50286		&	{\bf 0.5706}	&	0.4754	&	{\bf $-$0.1795}	     &  {\bf $-$0.4923}   \\ \hline
15578669387		&	$-$51116		&	$-$0.4029	&	0.4620 	&	$-$0.2829	     &    $-$0.7768 \\ \hline
19890188718		&	60442		&	0.4286	&	0.4643	&	$-$0.2676	     &   $-$0.7352  \\ \hline
22867694771		&	$-$62880	&	$-$0.4158	&	0.4632	&	$-$0.2767 	     &    $-$0.7602 \\ \hline
38066335279		&	$-$81220	&	$-$0.4163	&	0.4640	&	$-$0.2745 	     &    $-$0.7549 \\ \hline
48638777062		&	76946		&	0.3489	&	0.4572	&	$-$0.3287	     &   $-$0.9045  \\ \hline
56808201767		&	$-$87995	&	$-$0.3692	&	0.4597	&	$-$0.3105	     &    $-$0.8545 \\ \hline
101246135617		&	$-$129332	&	$-$0.4065	&	0.4645	&	$-$0.2785	     &  $-$0.7673    \\ \hline
108924543546	&	170358		&	{\bf 0.5162}	&	0.4740	&	{\bf $-$0.2044}	     &  {\bf  $-$0.5632}  \\ \hline
148491117087  	&	$-$131461	&	$-$0.3412	&	0.4582 	&	$-$0.3312	     &   $-$0.9130  \\ \hline
217309283735	&	$-$190936	&	$-$0.4096 	&	0.4658	&	$-$0.2736	     &    $-$0.7549 \\ \hline
297193839495	&	207478		&	0.3806 	&	0.4634	&	$-$0.2951	     &     $-$0.8145	\\ \hline
330508686218	&	$-$294816	&	{\bf $-$0.5128}	&	 0.4748	&	{\bf $-$0.2037}	     &    {\bf $-$0.5625} \\ \hline
402027514338	&	271498		&	0.4282	&	0.4683	&	$-$0.2582	     &    $-$0.7131 \\ \hline
661066575037	&	331302		&	0.4075	&	0.4670	&	$-$0.2717	     &    $-$0.7512 \\ \hline
1440355022306	&	$-$368527	&	$-$0.3071	&	0.4578	&	$-$0.3543	     &   $-$0.9810  \\ \hline
1653435193541	&	546666		&	0.4251	&	0.4696	&	$-$0.2563	     &   $-$0.7098  \\ \hline
2087416003490	&	$-$625681	&	$-$0.4331	&	0.4705	&	$-$0.2502	     &    $-$0.6930 \\ \hline
2343412610499	&	594442		&	0.3883	&	0.4668	&	$-$0.2824	     &  $-$0.7826   \\ \hline
3270926424607	&	$-$635558	&	$-$0.3514	&	0.4637	&	$-$0.3112	     &    $-$0.8627 \\ \hline
4098484181477	&	780932		&	0.3857	&	0.4672	&	$-$0.2828	     &    $-$0.7843 \\ \hline
5191164528277	&	$-$668864	&	$-$0.2936	&	0.4581	&	$-$0.3630	     &  $-$1.00716   \\ \hline
5197159385733	&	$-$689688	&	$-$0.3025	&	0.4592	&	$-$0.3540	     &   $-$0.9824  \\ \hline
10236053505745	&	1451233		&	0.4536	&	0.4736	&	$-$0.2325	     &   $-$0.6460  \\ \hline
21035055623987	&	$-$1740201	&	$-$0.3794	&	0.4684	&	$-$0.2831	     &   $-$0.7875  \\ \hline
21036453134939	&	$-$1745524	&	0.3806	&	0.4685	&	$-$0.2822	     &     $-$0.7850  \\ \hline
23431878209318	&	1903157		&	0.3932	&	0.4697	&	$-$0.2724	     &    $-$0.7579  \\ \hline
30501639884098	&	$-$1930205	&	$-$0.3495	&	0.4661	&	$-$0.3060	     &    $-$0.8518  \\ \hline
36161703948239	&	2727852		&	0.4536	&	0.4747	&	$-$0.2297	     &    $-$0.6397  \\ \hline
36213976311781	&	2783777		&	0.4626 &		 0.4753	&	$-$0.2240	     &   $-$0.6238  \\ \hline 
71578936427177	&	$-$4440015	&	{\bf $-$0.5248}	&	{\bf 0.4798}	& {\bf $-$0.1862}	     &   {\bf $-$0.5191}  \\ \hline
146734769129449	&	3733097		&	0.3082	&	0.4639	&	$-$0.3378	     &    $-$0.9428 \\ \hline
175688234263439	&	$-$5684793	&	$-$0.4289	&	0.4742	&	$-$0.2425	     &    $-$0.6772  \\ \hline
212132789199869	&	5491769		&	0.3771	&	0.4639	&	$-$0.3378	     &     $-$0.7792 \\ \hline 
212137538048059	&	5505045		&	0.3780	&	0.4705	&	$-$0.2783	     &   $-$0.7773  \\ \hline
304648719069787	&	$-$5757490	&	$-$0.3299	&	0.4667	&    $-$0.3162	     &   $-$0.8839  \\ \hline
351246529829131	&	9699950		&	{\bf 0.5176}	&	{\bf 0.4803}	&	{\bf $-$0.1876}	     &   {\bf $-$0.5244}  \\ \hline
1050365365851491&	$-$13728339	&	$-$0.4236	&	 0.4752	&	$-$0.2424	     &   $-$0.6790 \\ \hline
1211876202620741	&	16390637	&	0.4708	&	{\bf 0.4783}	&	$-$0.2123	     &    $-$0.5949 \\ \hline
2458719908828794	&	$-$20362905	&	$-$0.4107	&	0.4749	&	$-$0.2494	     &    $-$0.6997 \\ \hline
3295555617962269	&	18781262		&	0.3272	&	0.4687	&	$-$0.3124	     &    $-$0.8768  \\ \hline
3664310064219561	&	$-$23089949	&	$-$0.3814	&	0.4731	&	$-$0.2693	     &    $-$0.7559 \\ \hline
4892214197703689	&	24133331		&	0.3450	&	0.4705	&	$-$0.2967	     &   $-$0.8331 \\ \hline
6287915599821430	&	$-$35629003	&	$-$0.4493	&	{\bf 0.4780}	&	$-$0.2226	     &  $-$0.6254   \\ \hline
7332940231978758	&	40371499	&	0.4715	&	{\bf 0.4794}	&	$-$0.2090	     &    $-$0.5873 \\ \hline
11609864264058592345 	&	$-$1995900927  & {\bf $-$0.5857}  &	{\bf 0.4878}	&	{\bf$- $0.1414}	& {\bf $-$0.4021}	 \\ \hline
\end{tabular}\label{tab9}
\end{table}

\begin{table}[h!]
\caption{\centering Supremum of $\frac{|M(x)|}{\sqrt{x}}$, $L_M(x)$,  $L_{M}[\frac{1}{2}](x)$, and $L_{M}[\frac{1}{2},\!  0](x)$ on intervals  $I$}
\scriptsize 
\begin{tabular}{|r|r|r|r|r|r|} \hline 
$I$    &   $\sup_{I} \!\! \frac{|M(x)|}{\sqrt{x} }= \frac{|M(a_I)|}{\sqrt{a_I}} $& \ \   $\sup_I \!L_M(x)$  &  $\sup_I \!L_{M}[\frac{1}{2}](x)$ & $\sup_I \!L_{M}[\frac{1}{2},\!  0](x)  $  \\  \hline \hline
$[10^0,10^1]$ 	&	$1.00000000\ldots$ \hfill at $1$	& $0.430676$ \hfill at $5$	& 	$  
\infty$ at $e^-$ & ---	\\ \hline
$[10^1,10^2]$ 			&	$0.83205029\ldots$ \hfill at 13	&	$0.428317\ldots$ \hfill at 13	&	$-0.195195\ldots$	\hfill at 13 & 	$\infty$ \hfill at $(e^e)^-$ \\ \hline
$[10^2,10^3]$ 		&	$0.56710496\ldots$	\hfill at 199 &	$0.392843\ldots$ \hfill at 199 &	$-0.340372\ldots$	\hfill at 199 & 	$-1.110672\ldots$ \hfill at 199 \\ \hline
$[10^3,10^4]$ &	$0.47220269\ldots$	 \hfill at 2803 & $0.408988\ldots$ \hfill at 9861	&	$-0.362185\ldots$ \hfill at 2803 & $-1.030161\ldots$\hfill at 2803	\\ \hline
$[10^4,10^5]$ 		&	$0.46297703\ldots$ \hfill at 24185	&	 $0.429956\ldots$ \hfill at 59577	&		$-0.321227\ldots$ \hfill at 59577&  	$-0.880743\ldots$ \hfill at 59577\\ \hline
$[10^5,10^6]$ 		&	$0.43777620\ldots$ \hfill at 300551	&	 $0.434510\ldots$ \hfill at 300551	&		$-0.325888\ldots$ \hfill at 300551 &  	$-0.888129\ldots$ \hfill at 300551\\ \hline
$[10^6,10^7]$ 		&	$0.41825475\ldots$ \hfill at 1066854	&	 $0.444134\ldots$ \hfill at 6481601	&		$-0.318314\ldots$ \hfill at 6481601 &  	$-0.865335\ldots$ \hfill at 6481601\\ \hline
$[10^7,10^8]$ 		&	$0.46272869\ldots$ \hfill at 30919091	&	 $0.455318\ldots$ \hfill at 30919091	&		$-0.270616\ldots$ \hfill at 30919091 &  	$-0.736381\ldots$ \hfill at 30919091\\ \hline
$[33,10^{16}]$		&	$0.57059088\ldots$ \hfill at 7766842813	&	$\leq 0.490592$ \hfill \	&	$\leq-0.142501$ \hfill  \ & $\leq-0.405043$\hfill \	\\ \hline
\end{tabular}\label{tab10}
\end{table}

\begin{table}[h!]
\caption{\centering Increasingly large values of $q_{10^6}(x) \approx M(x)/\sqrt{x}$ (approximate IL$q$ values) \cite[Tables 2 and 3]{kuz}} 
\scriptsize
\begin{tabular}{|r|r|r|r|r|r|} \hline 
$\log x$  & $q_{10^6}(x)$  & $L_{q_{10^6}}(x)$  &  $L_{q_{10^6}}[0](x)$ & $L_{q_{10^6}}[0,  0](x)$  & $L_{q_{10^6}}[0, \! 0,\!0](x)$ \\  \hline \hline
43.898387   & $-0.58478$  & $-0.0122$  & $-0.142$   &   $-0.403$ &  $-$1.880  \\  \hline 
97.526523 & 0.61863  &  $-0.00492$ & $-0.105$   & $-0.316$  &  $-$1.144  \\  \hline 
140.373835  & $-0.58940$  &  $-0.00377$   &  $-0.107$  &  $-0.331$  &  $-$1.127  \\  \hline 
853.851589  & $-0.67715$  & $-0.000455$   & $-0.0578$   &  $-0.204$  &    $-$0.603 \\  \hline 
984.282019 & 0.62512 & $-0.000477$     &  $-0.0682$  &  $-0.243$ &   $-$0.714 \\  \hline 
1625.698493 & 0.62687 &  $-0.000287$  &   $-0.0632$  & $-0.233$  &  $-$0.673  \\  \hline 
1957.803133 & 0.62849 & $-0.000239$ &  $-0.0616$  & $-0.230$  &  $-$0.658  \\  \hline 
2709.485814  & 0.65467 & $-0.000156$ &  $-0.0536$  & $-0.205$  &  $-$0.583   \\  \hline 
2794.384965 & 0.65955 & $-0.000149$  &  $-0.0524$  &  $-0.201$ &  $-$0.572  \\  \hline 
3005.762748  & $-0.68878$  & $-0.000124$   &   $-0.0466$ & $-0.179$   &  $-$0.509  \\  \hline 
12277.362671 & 0.79344 &  $-0.0000188$ &  $-0.0246$  &  $-0.103$ &   $-$0.287  \\  \hline 
102494.024866 &  $-0.69580$   &  $-3.54 \cdot 10^{-6}$  &  $-0.0314$  &  $-0.148$  &  $-$0.406  \\  \hline 
150020.464414  &  $-0.70773$  &  $-2.30 \cdot 10^{-6}$  & $-0.0290$   &  $-0.139$  &  $-$0.381  \\  \hline 
178259.151801   &  $-0.71541$  &  $-1.88 \cdot 10^{-6}$  &  $-0.0277$  &  $-0.134$  &   $-$0.367 \\  \hline 
203205.659988  &  $-0.74083$  &  $-1.48 \cdot 10^{-6}$  &  $-0.0245$  &  $-0.120$  &   $-$0.327 \\  \hline 
860440.495719 &  $-0.75254$  &    $-3.30 \cdot 10^{-7}$  &  $-0.0208$  & $-0.109$   &  $-$0.296  \\  \hline 
1365643.292004  &  $-0.75927$  &  $-2.02 \cdot 10^{-7}$  & $-0.0195$   &   $-0.104$ &  $-$0.283  \\  \hline 
2765781.628095  &  $-0.76041$  &  $-9.90 \cdot 10^{-8}$  &  $-0.0185$ &  $-0.102$  &   $-$0.276 \\  \hline 
7078384.260482  &  $-0.76879$  &  $-3.71 \cdot 10^{-8}$  &  $-0.0167$  &   $-0.0953$ &  $-$0.259  \\  \hline 
13670267.747472 &   $-0.78505$   & $-1.77 \cdot 10^{-8}$   &  $-0.0147$  & $-0.0865$   &  $-$0.235  \\  \hline 
19371574.223934  &  $-0.78747$  & $-1.23 \cdot 10^{-8}$   &   $-0.0142$  & $-0.0847$   & $-$0.230   \\  \hline 
57334128.09084  &   $-0.80765$  & $-3.73 \cdot 10^{-9}$   & $-0.0120$   &  $-0.0741$   &  $-$0.202  \\  \hline 
86458087.131684 & 0.80443 & $-2.52 \cdot 10^{-9}$   &   $-0.0119$ &   $-0.0749$ & $-$0.204  \\  \hline
167211796.14902   &   $-0.82488$  &  $-1.15 \cdot 10^{-9}$  &  $-0.0102$  & $-0.0655$    &    $-$0.178	\\  \hline 
249548703.533702 & 0.81472  & $-8.21 \cdot 10^{-10}$ & $-0.0106$  & $-0.0692$  &    $-$0.189 \\  \hline 

405441986.398094  &   $-0.84497$  & $-4.15 \cdot 10^{-10}$    &  $-0.00850$  &   $-0.0564$ &  $-$0.155  \\  \hline 
1467573228.501077  & 0.81702& $-1.38 \cdot 10^{-10}$ &  $-0.00957$  & $-0.0663$  &   $-$0.181 \\  \hline  
2500922487.505913  & 0.82884   &   $-7.51 \cdot 10^{-11}$ &  $-0.00868$  &   $-0.0611$ &   $-$0.167 \\  \hline 
3847517705.646364 & 0.83682 & $-4.63 \cdot 10^{-11}$ &  $-0.00807$  &  $-0.0576$ &   $-$0.158 \\  \hline 
4016980126.87193  &  $-0.85146$  &  $-4.00 \cdot 10^{-11}$  &  $-0.00727$  &   $-0.0519$ &    $-$0.142 \\  \hline

10407545552.85608 & 0.842485 & $-1.65 \cdot 10^{-11}$ & $-0.00743$  &  $-0.0546$ &   $-$0.150 \\  \hline 
21334043144.02927	& 0.86622	& $-6.73 \cdot 10^{-12}$	&$-$0.00604	& $-$0.0453	&   $-$0.125\\  \hline 
86339883457.03526	&	$-$0.87296  &	$-1.57 \cdot 10^{-12}$	& $-$0.00540	&	$-$0.0421  &   $-$0.116 \\  \hline 
187114096628.77484	& 0.88636	&$-6.45 \cdot 10^{-13}$	&$-$0.00465	&  $-$0.0370	& $-$0.102   \\  \hline 
264421251554.46918	&	$-$0.89290  &	$-4.28 \cdot 10^{-13}$	&	$-$0.00431 &	$-$0.0346  &  $-$0.0956  \\  \hline 
1278282683343.76520	&	$-$0.89907  &	$-8.32 \cdot 10^{-14}$	& $-$0.00382	& $-$0.0320	&  $-$0.0885  \\  \hline 
1354137181464.62097 	& 0.90578	&	$-7.31 \cdot 10^{-14}$ & $-$0.00354 &	$-$0.0297 &  $-$0.0858  \\  \hline 
3680547202477.03623	&	$-$0.91392  &	$-2.45 \cdot 10^{-14}$	&	$-$0.00311 & $-$0.0267	&  $-$0.0742  \\  \hline 
6984497047106.74600	& 0.90744 &	 $-1.39 \cdot 10^{-14}$ &	$-$0.00328  &	$-$0.0287 &   $-$0.0796 \\  \hline 
18747209824980.73961	&	$-$0.91405 &	$-4.79 \cdot 10^{-15}$	& $-$0.00294	& $-$0.0263	& $-$0.0731   \\  \hline 
55714219637174.49540  &	$-$0.92025  &	$-1.49 \cdot 10^{-15}$	&	$-$0.00263 & $-$0.0241	&    $-$0.0670 \\  \hline 

84594507546024.46719	& 0.92208	&$-9.59 \cdot 10^{-16}$ & $-$0.00253	& $-$0.0234 &  $-$0.0652  \\  \hline 
117239588213313.90075	& 0.94102	& $-5.19 \cdot 10^{-16}$	&	$-$0.00188 &	$-$0.0175 &    $-$0.0488 \\  \hline 
117892199597999.02070	&   $-$0.92078 	&	$-7.00 \cdot 10^{-16}$   &  $-$0.00255	&  $-$0.0237	&   $-$0.0662 \\  \hline 
143697547951999.01914	&	$-$0.93038 &		$-5.02 \cdot 10^{-16}$	&	$-$0.00221  &	$-$0.0207  &   $-$0.0578 \\  \hline 
258592103887306.71643	&	$-$0.94097 &		$-2.35 \cdot 10^{-16}$	& 	$-$0.00183  &	$-$0.0174  &   $-$0.0485 \\  \hline 
505863698785929.24318 	&	$-$0.95652  &		$-8.79 \cdot 10^{-17}$	&	$-$0.00131  & $-$0.0126	 &   $-$0.0353  \\  \hline 
\end{tabular}\label{tab11}
\end{table}

\begin{table}[h!]
\caption{\centering Increasingly large values of $|q_{10^6}(x)| \approx |M(x)|/\sqrt{x}$ (approximate IL$|q|$ values) \cite[Tables 2 and 3]{kuz}} 
\scriptsize 
\begin{tabular}{|r|r|r|r|r|r|r|} \hline 
$\log x$  & $r_k = \lfloor\log_{10}x \rfloor$ (so $|q_{10^6}(x)| $  & $q_{10^6}(x)$  & $L_{q_{10^6}}(x)$  &  $L_{q_{10^6}}[0](x)$ & $L_{q_{10^6}}[0,  0](x)$  & $L_{q_{10^6}}[0, \! 0,\!0](x)$ \\  
 & $\approx \sup_{ [10^2,10^{r_{k+1}}]}|q_{10^6}(t)|)$  &   &  &   &   &  \\  \hline \hline
43.898387       &  19   &    $-0.58478$  & $-0.0122$  & $-0.142$   &   $-0.403$ &  $-$1.880  \\  \hline 
97.526523     &   42  &    0.61863  &  $-0.00492$ & $-0.105$   & $-0.316$  &  $-$1.144  \\  \hline 
853.851589      &  370   &    $-0.67715$  & $-0.000455$   & $-0.0578$   &  $-0.204$  &    $-$0.603 \\  \hline 
3005.762748      &  1305   &    $-0.68878$  & $-0.000124$   &   $-0.0466$ & $-0.179$   &  $-$0.509  \\  \hline 
12277.362671     &   5331  &    0.79344 &  $-0.0000188$ &  $-0.0246$  &  $-0.103$ &   $-$0.287  \\  \hline 
57334128.09084      &  24899895   &      $-0.80765$  & $-3.73 \cdot 10^{-9}$   & $-0.0120$   &  $-0.0741$   &  $-$0.202  \\  \hline 
167211796.14902       &   72619160  &      $-0.82488$  &  $-1.15 \cdot 10^{-9}$  &  $-0.0102$  & $-0.0655$    &    $-$0.178	\\  \hline 
405441986.398094      &  
176081217   &      $-0.84497$  & $-4.15 \cdot 10^{-10}$    &  $-0.00850$  &   $-0.0564$ &  $-$0.155  \\  \hline 
4016980126.87193      &   1744552303  &     $-0.85146$  &  $-4.00 \cdot 10^{-11}$  &  $-0.00727$  &   $-0.0519$ &    $-$0.142 \\  \hline 
21334043144.02927	     &  9265257214   &     0.86622	& $-6.73 \cdot 10^{-12}$	&$-$0.00604	& $-$0.0453	&   $-$0.125\\  \hline 
86339883457.03526	       &  37496934953   &    	$-$0.87296  &	$-1.57 \cdot 10^{-12}$	& $-$0.00540	&	$-$0.0421  &   $-$0.116 \\  \hline 
187114096628.77484	     & 81262619652    &     0.88636	&$-6.45 \cdot 10^{-13}$	&$-$0.00465	&  $-$0.0370	& $-$0.102   \\  \hline 
264421251554.46918	     &  114836690448   &    	$-$0.89290  &	$-4.28 \cdot 10^{-13}$	&	$-$0.00431 &	$-$0.0346  &  $-$0.0956  \\  \hline 
1278282683343.76520	     & 555151115688    &    	$-$0.89907  &	$-8.32 \cdot 10^{-14}$	& $-$0.00382	& $-$0.0320	&  $-$0.0885  \\  \hline 
1354137181464.62097 	     &   588094305650  &     0.90578	&	$-7.31 \cdot 10^{-14}$ & $-$0.00354 &	$-$0.0297 &  $-$0.0858  \\  \hline 
3680547202477.03623	     &   1598441340420  &    	$-$0.91392  &	$-2.45 \cdot 10^{-14}$	&	$-$0.00311 & $-$0.0267	&  $-$0.0742  \\  \hline 
18747209824980.73961	     &   8141809778071  &    	$-$0.91405 &	$-4.79 \cdot 10^{-15}$	& $-$0.00294	& $-$0.0263	& $-$0.0731   \\  \hline 
55714219637174.49540       & 24196378151970    &    	$-$0.92025  &	$-1.49 \cdot 10^{-15}$	&	$-$0.00263 & $-$0.0241	&    $-$0.0670 \\  \hline 
84594507546024.46719	     &  36738927826561   &     0.92208	&$-9.59 \cdot 10^{-16}$ & $-$0.00253	& $-$0.0234 &  $-$0.0652  \\  \hline 
117239588213313.90075	     &   50916506221651  &     0.94102	& $-5.19 \cdot 10^{-16}$	&	$-$0.00188 &	$-$0.0175 &    $-$0.0488 \\  \hline 
505863698785929.24318 	     &  219693812977897   &    	$-$0.95652  &		$-8.79 \cdot 10^{-17}$	&	$-$0.00131  & $-$0.0126	 &   $-$0.0353  \\  \hline 
\end{tabular}\label{tab12}
\end{table}

\printindex[symbols]
  \printindex

\end{document}